\documentclass[a4paper, 10pt]{report}

\usepackage{morewrites}

\usepackage[top=2.7cm, bottom=2.5cm, left=2.5cm, right=2.5cm]{geometry} 
\usepackage[main=english,french]{babel} 
\usepackage[utf8]{inputenc} 
\usepackage[T1]{fontenc} 
\usepackage{lmodern} 
\usepackage{graphicx} 
\usepackage{amsmath, amssymb, mathrsfs, mathtools, cancel} 
\usepackage{xcolor} 

\usepackage{environ} 
\usepackage{fancybox} 
\usepackage[thmmarks, amsmath]{ntheorem} 
\usepackage[framemethod=TikZ]{mdframed} 
\usepackage{tablefootnote} 
\usepackage{aliascnt} 

\usepackage{tikz}
\usepackage{pgfplots, pgfplotstable}
\usepackage[justification=centering]{caption} 
\usepackage{enumitem} 
\usepackage{longtable}
\usepackage{multirow} 
\usepackage[all]{xy} 

\usepackage{multicol} 
\usepackage{slashbox} 
\usepackage{contour} 
\usepackage[ntheorem]{empheq} 
\usepackage[section]{placeins} 
\usepackage{subcaption} 
\usepackage{appendix} 
\usepackage[nottoc]{tocbibind} 
\usepackage{minitoc} 
\usepackage{xparse} 
\usepackage[intoc,refpage]{nomencl} 
\usepackage{fancyvrb} 
\usepackage{fancyhdr} 
\usepackage{xpatch} 
\usepackage[ruled,linesnumbered,algochapter]{algorithm2e} 
\usepackage{listings} 
\usepackage{everypage} 
\usepackage{tabularray} 
\usepackage{pdfpages} 
\usepackage{CJKutf8} 

\usetikzlibrary{patterns.meta} 
\usetikzlibrary{arrows.meta} 
\usetikzlibrary{decorations.pathreplacing} 
\usetikzlibrary{fadings} 
\usetikzlibrary{babel} 
\pgfplotsset{compat=1.9}
\usepgfplotslibrary{colorbrewer}

\usetikzlibrary{external}
\tikzset{external/mode=graphics if exists}
\tikzset{external/export=false}

\usepackage[style=alphabetic,maxbibnames=99,datamodel=software]{biblatex} 
\usepackage{software-biblatex}
\AtBeginBibliography{\small} 
\usepackage[hidelinks, colorlinks=true, linkcolor=blue, urlcolor=cyan, citecolor=green!70!black]{hyperref} 
\usepackage[capitalise,noabbrev]{cleveref}

\catcode`\:=13 

\makeatletter 

\toksdef\list@toka=0 \toksdef\list@tokb=2 
\newcount\list@count

\def\newlist#1{\gdef#1{}}

\long\def\appendlistitem#1\to#2{%
	\list@toka={\list@item{#1}}%
	\list@tokb=\expandafter{#2}%
	\xdef#2{\the\list@tokb \the\list@toka}}
\long\def\leftappendlistitem#1\to#2{%
	\list@toka=\expandafter{#2}%
	\list@tokb={\list@item{#1}}%
	\xdef#2{\the\list@tokb \the\list@toka}}
\long\def\eappendlistitem#1\to#2{%
	\list@tokb=\expandafter{#2}%
	\xdef#2{\the\list@tokb \noexpand\list@item{#1}}}

\def\listsize#1\to#2{%
	#2=0\relax
	\long\def\list@item##1{\advance#2 by 1}%
	#1}
\def\printlistsize#1{\listsize#1\to\list@count \relax \the\list@count}

\def\getlistitem#1\of#2{
	\list@count=#1\relax
	\long\gdef\list@item##1{%
		\advance\list@count by -1\relax
		\ifnum \list@count=0 ##1\fi}%
	#2}
\def\putlistitem#1\of#2\in#3{
	\list@count=#1\relax
	\long\gdef\list@item##1{%
		\advance\list@count by -1\relax
		\ifnum \list@count=0 \gdef#3{##1}\fi}%
	#2}

\def\list@getsep sep={}
\newcommand{\printlist}[2][sep=]{%
	\edef\list@sep{\list@getsep#1}%
	\list@count=0\relax
	\long\def\list@item##1{%
		\ifnum \list@count=0 \else \list@sep \fi
		\advance\list@count by 1\relax
		##1}%
	#2}

\long\def\@deflist#1=#2{%
	\begingroup
	\gdef#1{}%
	\expandafter\def\expandafter\@@deflist\expandafter##\expandafter1\list@sep##2\@@deflist{
		\appendlistitem{##1}\to#1%
		\if\relax\detokenize{##2}\relax
		\else \@@deflist##2\@@deflist \fi}
	\def\list@temp##1{\@@deflist#2##1\@@deflist}%
	\expandafter\list@temp\expandafter{\list@sep}%
	\endgroup}
\newcommand{\deflist}[1][sep=,]{%
	\edef\list@sep{\list@getsep#1}%
	\@deflist}

\def\applytolistitems#1\in#2{
	\let\list@copy=#2%
	\def#2{}
	\def\list@item##1{
		\list@toka=\expandafter{#1{##1}}%
		\expandafter\appendlistitem\expandafter{\the\list@toka}\to#2}%
	\list@copy
}

\setlist{itemsep=0pt, topsep=3pt} 
\numberwithin{equation}{section}

\newcommand{\longCaption}[2]{\caption[#1]{#1#2}} 

\newcommand{\red}[1]{\textcolor{red!93!black}{#1}}
\definecolor{couleurNePasRelire}{rgb}{0.15, 0.25, 0.7}

\let\old@underline\underline
\renewcommand{\underline}[1]{\old@underline{#1\vrule width 0pt depth 0.25ex}}

\crefname{equation}{}{} 
\Crefname{equation}{Equation}{Equations} 

\newcommand{\inlineEqNb}{\refstepcounter{equation}\textup{\tagform@{\theequation}}}

\DeclareCaptionLabelFormat{lc}{\textup{#1}~#2}
\tikzifexternalizing{}{%
\dominitoc 
\AtBeginDocument{
	\let\old@chapter\chapter
	\RenewDocumentCommand{\chapter}{ s O{#3} m}{
		\IfBooleanTF{#1}{\old@chapter*{#3}}{\old@chapter[#2]{#3}\minitoc}}
}}

\makenomenclature 
\renewcommand{\nompreamble}{The list below enumerates in alphabetical order the notations that are used in multiple sections. For more details, see the page of definition indicated on the right.} 
\nomenclature[23]{$\mathbb{N}, \mathbb{N}^*$}{Set of natural integers including, resp.\ excluding, 0}

\fancypagestyle{plain}{
    \fancyhf{}
    
    \fancyfoot[C]{--\ \thepage\ --}}

\newcommand{\spewfootnotes}{%
\tfn@tablefootnoteprintout%
\global\let\tfn@tablefootnoteprintout\relax%
\gdef\tfn@fnt{0}%
}

\g@addto@macro\@floatboxreset\centering 
\let\ftype@table\ftype@figure 

\def\marksizedefaut{2pt} 
\newcommand{\creercyclepgfplotsregression}[3][sep=,]{
	\long\def\list@item##1{\@temp##1\relax}
	\long\def\@temp##1 + ##2\relax{
		opacity=0\\
		color=##1,dotted,mark size=\marksizedefaut,mark=##2,every mark/.append style={solid}\\
		color=##1\\}
	\deflist[#1]\liste={#3}
	\edef\liste{\liste}
	\def\@temp##1{\pgfplotscreateplotcyclelist{#2}{##1}}
	\expandafter\@temp\expandafter{\liste}
	\def\marksizedefaut{2pt} 
}

\definecolor{Dark2-cyan}{RGB}{27,158,119}
\definecolor{Dark2-orange}{RGB}{217,95,2}
\definecolor{Dark2-purple}{RGB}{117,112,179}
\definecolor{Dark2-magenta}{RGB}{231,41,138}
\definecolor{Dark2-green}{RGB}{102,166,30}
\definecolor{Dark2-yellow}{RGB}{230,171,2}
\definecolor{Dark2-kaki}{RGB}{166,118,29}
\definecolor{Dark2-grey}{RGB}{102,102,102}
\definecolor{Dark2-blue}{RGB}{0,80,180}
\definecolor{Dark2-brown}{RGB}{160,90,20}

\def\marksizedefaut{2.5pt}
\creercyclepgfplotsregression{Dark2_avec_regressions}{%
	Dark2-cyan + *,%
	Dark2-orange + square*,%
	Dark2-purple + triangle*,%
	Dark2-magenta + halfsquare*,%
	Dark2-green + pentagon*,%
	Dark2-yellow + halfcircle*,%
	Dark2-kaki + halfdiamond*,%
	Dark2-grey + otimes,%
	Dark2-blue + *,%
	Dark2-brown + square*}

\def\marksizedefaut{2.5pt}
\creercyclepgfplotsregression[sep=;]{Dark2_CL_2.5}{%
	Dark2-cyan + * ;%
	Dark2-orange + triangle* ;%
	Dark2-purple + pentagon* ;%
	Dark2-magenta + star, mark size=3pt ;%
	Dark2-green + square* ;%
	Dark2-yellow + halfsquare* ;%
	Dark2-kaki + x, mark size=3pt ;%
	Dark2-grey + square* ;%
	Dark2-brown + halfsquare* ;%
	Dark2-blue + +, mark size=3pt}

\def\marksizedefaut{2.5pt}
\creercyclepgfplotsregression{Dark2_CL_4a}{%
	Dark2-cyan + *,%
	Dark2-orange + square*,%
	Dark2-purple + triangle*,%
	Dark2-yellow + pentagon*,%
	Dark2-green + halfcircle*,%
	Dark2-magenta + star,%
	Dark2-kaki + otimes*,%
	Dark2-kaki + otimes*,%
	Dark2-grey + x,%
	Dark2-grey + x,%
	Dark2-grey + x}
\def\marksizedefaut{2.5pt}
\creercyclepgfplotsregression{Dark2_CL_4b}{%
	Dark2-grey + x,%
	Dark2-blue + *,%
	Dark2-brown + triangle*,%
	Dark2-brown + triangle*,%
	Dark2-orange + halfdiamond*,%
	Dark2-orange + halfdiamond*,%
	Dark2-cyan + otimes}

\def\largeurgrapheCL{75mm}

\newenvironment{axisCL}[1]{
	\catcode`\:=13 
	\gdef\num@fichier{figures/erreurs_CL/#1}
	\gdef\num@nopremierCL{-}
	\newlist\num@listepentes
	\begin{axis}[
		height=\largeurgrapheCL,
		xlabel={$\varepsilon$},
		ylabel={$\displaystyle \bigg\|u_\varepsilon -\sum_{p'<p,\,\ell} \varepsilon^{p'} \ln^\ell\!\varepsilon \, u_{p,\ell} \bigg\|_{H^1} / \.\|u_{0,0}\|_{H^1}$},
		enlargelimits=0.045,
    	xtick={0.04,0.065,0.1,0.16,0.25},
		xmajorgrids=true,
		ymajorgrids=true,
		xmode=log,
		ymode=log,
		log ticks with fixed point, 
		legend style={
			at={(1.05,0.5)},
			anchor=west,
			legend columns=1,}]
}{
	\addplot[opacity=0]
		table[y=errH1_CL_\num@nopremierCL] {\num@fichier};
		\addlegendentry[text height=0ex]{}
	\end{axis}
}

\newenvironment{axisME}[3][relative $H^1$~norm on ]{
	\catcode`\:=13 
	\gdef\num@fichier{figures/erreurs_ME/#2}
	\newlist\num@listepentes
	\begin{axis}[
		height=60mm,
		xlabel={$\varepsilon$},
		ylabel={#1#3},
		enlargelimits=0.045,
    	xtick={0.04,0.065,0.1,0.16,0.25},
		xmajorgrids=true,
		ymajorgrids=true,
		xmode=log,
		ymode=log,
		log ticks with fixed point, 
		legend style={
			at={(1.05,0.5)},
			anchor=west,
			legend columns=1,}]
}{
	\addplot[opacity=0]
		table[y=errH1_ME_3] {\num@fichier};
		\addlegendentry[text height=0ex]{}
	\end{axis}
}

\def\nameslope{slope}
\def\num@boitepCL#1{\begingroup
    \setbox0=\hbox{$p=0+0\pst\approx0.00$}%
    \hbox to \wd0{#1\hss}%
    \endgroup}
\gdef\num@boitepenteCL#1{\hphantom{\quad\;\nameslope\ 0.0}}
\def\grapheCLpcourt{%
	\def\num@boitepCL##1{\begingroup
    \setbox0=\hbox{$p=0$}%
    \hbox to \wd0{##1\hss}%
    \endgroup}}
\def\num@boiteuME#1{\begingroup
    \setbox0=\hbox{$u_\peps - u^{(3)}$}%
    \hbox to \wd0{#1\hss}%
    \endgroup}
\def\num@boiteuMEintro#1{\begingroup
    \setbox0=\hbox{$m=3$}%
    \hbox to \wd0{#1\hss}%
    \endgroup}
\def\num@boitepenteME#1{(\nameslope\ \pgfmathprintnumber[precision=1]{#1})}
\def\@courbeCLME[y=#1, reg=#2:#3, nom=#4, comms=#5#6]{
    \pgfplotstablesave[
        skip rows between index={0}{#2},
        skip rows between index={#3}{10}
    ]{\num@fichier}{tempo.dat}

    \addplot[opacity=0]
        table[y=#1] {\num@fichier};
        \addlegendentry[text height=0ex]{}
    \addplot
        table[y=#1] {\num@fichier};
        \addlegendentry{}
    \addplot
        table[y={create col/linear regression={y=#1}}] {tempo.dat};
        \eappendlistitem{\pgfplotstableregressiona}\to{\num@listepentes}
    	\pgfmathparse{\pgfplotstableregressiona==0 ? 0 : 1}%
		\ifnum\pgfmathresult=1
        	\addlegendentryexpanded[text height=0ex]{#5{#4} #6{\pgfplotstableregressiona}}
		\else \addlegendentryexpanded[text height=0ex]{#5{#4}} \fi
}
\def\courbeCL[no=#1, p=#2, reg=#3:#4]{
	\if-\num@nopremierCL \gdef\num@nopremierCL{#1}\fi
	\@courbeCLME[y=errH1_CL_#1, reg=#3:#4, nom={$p=#2$}, comms=\num@boitepCL\num@boitepenteCL]}
\def\courbeME[y=#1, fct=$#2$, reg=#3:#4]{
	\@courbeCLME[y=#1, reg=#3:#4, nom={$u_\peps-\noexpand#2$}, comms=\num@boiteuME\num@boitepenteME]}
\def\courbeMEintro[y=#1, m=#2, reg=#3:#4]{
	\@courbeCLME[y=#1, reg=#3:#4, nom={$m=#2$}, comms=\num@boiteuMEintro\num@boitepenteME]}

\def\nodepenteCL[nodbt=#1, nb=#2]{
	\begingroup
	\pgfmathsetmacro{\dbt}{#1 - \num@nopremierCL + 1}
	\pgfmathsetmacro{\fin}{#1 - \num@nopremierCL + #2}
	\pgfmathsetmacro{\num@pentemax}{-1000}
	\pgfmathsetmacro{\num@pentemin}{1000}
	\foreach \i in {\dbt,...,\fin} {
		\putlistitem\i\of\num@listepentes\in\@tempo
		\pgfmathsetmacro{\num@pentemax}{max(\num@pentemax, \@tempo)}
		\pgfmathsetmacro{\num@pentemin}{min(\num@pentemin, \@tempo)}
		\global\let\num@pentemax=\num@pentemax
		\global\let\num@pentemin=\num@pentemin}
	\endgroup
	\pgfmathsetmacro{\num@pentemoy}{(\num@pentemax+\num@pentemin)/2}
	\pgfmathsetmacro{\num@pentemoy}{floor(\num@pentemoy*10+0.5)/10}
	\pgfmathsetmacro{\num@errpente}{(\num@pentemax-\num@pentemin)/2}
	\listsize\num@listepentes\to\list@count
    \draw ({\largeurgrapheCL*1.05 + 1.25em},{\largeurgrapheCL*0.475 + (\the\list@count/2 - (#1-\num@nopremierCL+1) - #2/2)*3.82ex}) node[right] {\num@boitepCL{}$
        \if1#2 \big\}
        \else\if2#2 \Bigg\}
        \else\if3#2 \mkern-2mu\left. \vrule height5.9ex width0pt \right\} 
        \else\if4#2 \mkern-2mu\left. \vrule height7.6ex width0pt \right\} 
        \fi\fi\fi\fi
    \mkern-1mu$ \nameslope\ \pgfmathprintnumber{\num@pentemoy}%
    \pgfmathparse{\num@errpente > 0.099 ? 1 : 0}%
	\ifnum\pgfmathresult=1 ${}\pm{}$\pgfmathprintnumber[fixed]{\num@errpente}\fi
    };
}

\def\clbr@ecrirenb#1#2{%
	\footnotesize\hskip-1pt\relax
    \pgfmathparse{abs(#1) < 1e-4 ? 1 : 0}%
    \ifnum\pgfmathresult=1 \pgfmathprintnumber{#1}%
    \else
        \pgfmathsetmacro{\puis}{floor(log10(abs(#1)+3*1e-4))}%
        \pgfmathsetmacro{\val}{#1 * 10^(-\puis)}%
        \pgfmathsetmacro{\puis}{\puis+(#2)}%
        \pgfmathparse{abs(\puis) < 1.1 ? 1 : 0}%
        \ifnum\pgfmathresult=1
        	\pgfmathsetmacro{\val}{\val*10^(\puis)}%
			\pgfmathprintnumber[fixed]{\val}%
        \else
			\pgfmathprintnumber{\val}$ {} \cdot 10^{\text{\pgfmathprintnumber{\puis}}}$%
    \fi\fi\strut
}
\def\colorbarheight{3cm}
\def\colorbarticks{3.5}
\def\colorbarticks{2.5}
\def\colorbar[min=#1, max=#2, 1e#3]{%
    \begingroup
    \pgfmathsetmacro{\pasVise}{(#2-#1)/(\colorbarticks-1.5)}%
    \pgfmathsetmacro{\powTen}{10^(floor(log10(\pasVise) + 1e-5))}%
    \pgfmathsetmacro{\mantissa}{\pasVise/\powTen}%
    \pgfmathtruncatemacro{\mantissa}{\mantissa < 1.5 ? 1 : (\mantissa < 3.5 ? 2 : (\mantissa < 7.5 ? 5 : 10))}%
    \pgfmathsetmacro{\pas}{\mantissa*\powTen}
    \pgfmathsetmacro{\marge}{2ex/\colorbarheight * (#2-#1)/\pas}%
    \pgfmathsetmacro{\intA}{ceil(#1/\pas + \marge)}%
    \pgfmathsetmacro{\intB}{floor(#2/\pas - \marge)}%
    \pgfmathsetmacro{\nbTicks}{\intB - \intA}%
    \pgfmathsetlengthmacro{\yUnit}{\colorbarheight/(#2-#1)}%
    \raisebox{-.5\height}{\begin{tikzpicture}[y=\yUnit]
		\draw (0,#1) node[above left,inner sep=0,outer sep=0]{\includegraphics[height=\colorbarheight]{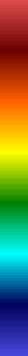}};
        \foreach \i in {0,...,\nbTicks} {
            \pgfmathsetmacro{\tickPos}{(\intA+\i)*\pas}%
            \draw (-0.12,\tickPos) -- (0.12,\tickPos)  node[right] {\clbr@ecrirenb{\tickPos}{#3}};
        }
        \draw (-0.12,#1) -- (0.12,#1) node[right] {\clbr@ecrirenb{#1}{#3}};
        \draw (-0.12,#2) -- (0.12,#2) node[right] {\clbr@ecrirenb{#2}{#3}};
    \end{tikzpicture}}%
    \endgroup
}

\newlength{\reducBoite}
\newcommand{\fbeq}[1][1]{\setlength{\reducBoite}{-2.3ex}\vspace{#1\reducBoite}} 
\newcommand{\fbi}[1][1]{\setlength{\reducBoite}{-1ex}\vspace{#1\reducBoite}} 

\NewDocumentCommand{\@clonerCompteur}{mmO{}}{
    \newaliascnt{#1}{paragraphe}
    \expandafter\renewcommand\csname the#1\endcsname{\thesection.\arabic{#1}}
	\crefname{#1}{#2}{\ifx\relax#3\relax #2s\else #3\fi}
	\Crefname{#1}{#2}{\ifx\relax#3\relax #2s\else #3\fi}}
\@clonerCompteur{definition}{Definition}
\@clonerCompteur{lemme}{Lemma}
\@clonerCompteur{proposition}{Proposition}
\@clonerCompteur{corollaire}{Corollary}[Corollaries]
\@clonerCompteur{theoreme}{Theorem}
\@clonerCompteur{remarque}{Remark}
\@clonerCompteur{exemple}{Example}
\@clonerCompteur{notation}{Notation}
\@clonerCompteur{hypothese}{Assumption}
\@clonerCompteur{ansatz}{Ansatz}[Ansatze]

\cornersize*{12pt}
\def\width@boite{\textwidth-12pt}

\newcommand{\avant@boite}{%
	\medskip\smallskip
	\begingroup
	\def\footnote{\tablefootnote}
	\setlength{\fboxrule}{0.5pt}%
	\setlength{\fboxsep}{5.5pt}}

\newcommand{\dbt@boite}[1]{
	\vspace{1pt}%
	\textbf{#1}%
	\dbt@@boite#1\dbt@@boite}
\def\deuxPointsEst#1{
	\def\dbt@@boite##1#1##2\dbt@@boite{
		\ifx\relax##2\relax \else
		{\def\soustitre{##2}\def\espace{ }%
		\if\soustitre\espace \else\vadjust{\kern1pt\relax}\par\fi}\fi}}
\deuxPointsEst{:}

\newcommand{\fin@boite}{\vskip1pt\relax}

\newcommand{\apres@boite}{%
	\spewfootnotes
	\endgroup
	\medbreak\smallskip}

\newcommand{\avant@par@sansskip}{%
	\medskip\smallskip
	\begingroup
	\def\footnote{\tablefootnote}}
\newcommand{\avant@par@avecskip}{%
	\avant@par@sansskip
	\def\\{\par\medskip}}
\newcommand{\avant@par}{\@ifstar{\avant@par@sansskip}{\avant@par@avecskip}}

\newcommand{\apres@par}{%
	\spewfootnotes
	\endgroup
	\medbreak\smallskip}

\newcommand{\boite@ronde}[2]{%
	\avant@boite
	\vskip-1pt\relax
	\ovalbox{\parbox{\width@boite}{\dbt@boite{#1}#2\fin@boite}}%
	\apres@boite}
\newcommand{\boite@carree}[2]{%
	\avant@boite
	\vskip-1pt\relax
	\fbox{\parbox{\width@boite}{\dbt@boite{#1}#2\fin@boite}}%
	\apres@boite}
\NewEnviron{boiteronde}[1]{\boite@ronde{#1}{\BODY}}
\NewEnviron{boitecarree}[1]{\boite@carree{#1}{\BODY}}

\mdfdefinestyle{styleboiteronde}{roundcorner=6pt, linewidth=0.5pt, innerleftmargin=5.5pt, innerrightmargin=5.5pt, innertopmargin=5.5pt, innerbottommargin=5.5pt, skipabove=0pt, skipbelow=0pt}
\newenvironment{boiterondesecable}[1]{%
	\avant@boite
	\begin{mdframed}[style=styleboiteronde]\dbt@boite{#1}\ignorespaces}
	{\fin@boite\end{mdframed}\apres@boite}

\mdfdefinestyle{styleboitecarre}{linewidth=0.5pt, innerleftmargin=5.5pt, innerrightmargin=5.5pt, innertopmargin=5.5pt, innerbottommargin=5.5pt, skipabove=0pt, skipbelow=0pt, needspace=4\baselineskip}
\newenvironment{boitecarreesecable}[1]{%
	\avant@boite
	\begin{mdframed}[style=styleboitecarre]\dbt@boite{#1}\ignorespaces}
	{\end{mdframed}\apres@boite}

\newcommand{\grosse@boite@ronde}[2]{%
	\avant@boite
	\setlength{\fboxrule}{1.2pt}\setlength{\fboxsep}{1.3pt}\cornersize*{17pt}%
	\noindent\Ovalbox{%
	\setlength{\fboxrule}{0.4pt}\setlength{\fboxsep}{5.5pt}\cornersize*{4mm}%
	\ovalbox{\parbox{0.965\textwidth}{\dbt@boite{#1}#2\fin@boite}}}%
	\apres@boite} 
\newcommand{\grosse@boite@carree}[2]{%
	\avant@boite
	\setlength{\fboxrule}{0.5pt}%
	\noindent\doublebox{\parbox{0.965\textwidth}{\dbt@boite{#1}#2\fin@boite}}%
	\apres@boite}
\NewEnviron{grosseboiteronde}[1]{\grosse@boite@ronde{#1}{\BODY}}
\NewEnviron{grosseboitecarree}[1]{\grosse@boite@carree{#1}{\BODY}}
	
\newcommand{\paragraphe@sans@boite}[2]{\avant@par\textbf{#1}#2\apres@par}
\newenvironment{paragraphesansboite}[1]{\avant@par\textbf{#1}}{\apres@par}
\newenvironment{textesansboite}{\avant@par}{\apres@par}

\NewEnviron{definition}[1]{%
	\refstepcounter{definition}%
	\boite@ronde{Definition \thedefinition: #1}{\BODY}}
\NewEnviron{proposition}[1]{%
	\refstepcounter{proposition}%
	\boite@carree{Proposition \theproposition: #1}{\BODY}}
\NewEnviron{lemme}[1]{%
	\refstepcounter{lemme}%
	\boite@carree{Lemma \thelemme: #1}{\BODY}}
\NewEnviron{corollaire}[1]{%
	\refstepcounter{corollaire}%
	\boite@carree{Corollary \thecorollaire: #1}{\BODY}}
\NewEnviron{theoreme}[1]{%
	\refstepcounter{theoreme}%
	\boite@carree{Theorem \thetheoreme: #1}{\BODY}}

\newenvironment{definitionsecable}[1]{%
	\refstepcounter{definition}%
	\begin{boiterondesecable}{Definition \thedefinition: #1}}{%
	\end{boiterondesecable}} 
\newenvironment{propositionsecable}[1]{%
	\refstepcounter{proposition}%
	\begin{boitecarreesecable}{Proposition \theproposition: #1}}{%
	\end{boitecarreesecable}} 
\newenvironment{lemmesecable}[1]{%
	\refstepcounter{lemme}%
	\begin{boitecarreesecable}{Lemma \thelemme: #1}}{%
	\end{boitecarreesecable}}

\newenvironment{theoremesecable}[1]{%
	\refstepcounter{theoreme}%
	\begin{boitecarreesecable}{Theorem \thetheoreme: #1}}{%
	\end{boitecarreesecable}}
		
\NewEnviron{grandtheoreme}[1]{
	\refstepcounter{theoreme}\groose@boite@carree{Theorem \thetheoreme: #1}{\BODY}}

\NewEnviron{definitionsansboite}[1]{%
	\refstepcounter{definition}%
	\paragraphe@sans@boite{Definition \thedefinition: #1}{\BODY}}
\NewEnviron{propositionsansboite}[1]{%
	\refstepcounter{proposition}%
	\paragraphe@sans@boite{Proposition \theproposition: #1}{\BODY}}
\NewEnviron{lemmesansboite}[1]{%
	\refstepcounter{lemme}%
	\paragraphe@sans@boite{Lemma \thelemme: #1}{\BODY}}
\NewEnviron{corollairesansboite}[1]{%
	\refstepcounter{corollaire}%
	\paragraphe@sans@boite{Corollary \thecorollaire: #1}{\BODY}}
\NewEnviron{theoremesansboite}[1]{%
	\refstepcounter{theoreme}%
	\paragraphe@sans@boite{Theorem \thetheoreme: #1}{\BODY}}
\NewEnviron{remarquenumerotee}{%
	\refstepcounter{remarque}%
	\paragraphe@sans@boite{Remark \theremarque:\ }{\BODY}}
\NewEnviron{exemplenumerote}{%
	\refstepcounter{exemple}%
	\paragraphe@sans@boite{Example \theexemple:\ }{\BODY}}

\NewEnviron{remarque}{\paragraphe@sans@boite{Remark:\ }{\BODY}}
\NewEnviron{exemple}{\paragraphe@sans@boite{Example:\ }{\BODY}}
\NewEnviron{remarques}{\paragraphe@sans@boite{Remarks:\ }{\BODY}}
\NewEnviron{exemples}{\paragraphe@sans@boite{Examples:\ }{\BODY}}

\theorempreskip{0pt}
\theorempostskip{0pt}
\theoremstyle{empty}
\theorembodyfont{\normalfont}
\theoremsymbol{\ensuremath{\scriptstyle\Box}}
\newtheorem{texteavecQED}{}

\definecolor{grisfonce}{gray}{0.3}
\newenvironment{demo}[1][Proof]{
	\avant@par*
	\begingroup\color{grisfonce}%
	\begin{texteavecQED}\underline{#1:}\ %
}{%
	\end{texteavecQED}%
	\endgroup
	\apres@par}

\newenvironment{demobis}[1][Proof]{
	\avant@par*
	\begingroup\color{grisfonce}%
	\underline{#1:}\ %
	}{\setbox0=\hbox to \hsize{\hfill\scriptsize$\Box$}
	\vadjust{{\color{grisfonce}\kern-\ht0\box0}}%
	\endgroup
	\apres@par}

\mdfdefinestyle{stylesousdemo}{linewidth=0.6pt, topline=false, rightline=false, bottomline=false, leftmargin=2em, innerleftmargin=5.5pt, innerrightmargin=0pt, innertopmargin=1pt, innerbottommargin=0pt, skipabove=2pt, skipbelow=0pt, fontcolor=grisfonce}
\def\sousdemo@crochet{\hbox{\kern -5.5pt \vrule width 1.5mm height 0.6pt}}
\newenvironment{sousdemo}{
	\begin{mdframed}[style=stylesousdemo]%
	}{\ifvmode \kern -1ex \sousdemo@crochet%
	\else \vadjust{\kern 1pt \sousdemo@crochet}\fi
	\end{mdframed}}

\let\Smallrightarrow\Rightarrow
\let\Smallleftarrow\Leftarrow

\renewcommand{\Rightarrow}{\mathrel{\mbox{$=\hspace{-6.5pt}\Smallrightarrow$}}}
\renewcommand{\Leftarrow}{\mathrel{\mbox{$\Smallleftarrow\hspace{-6.5pt}=$}}}
\renewcommand{\Leftrightarrow}{\mathrel{\mbox{$\Smallleftarrow\hspace{-6pt}\Smallrightarrow$}}}

\let\old@cr\relax
\def\def@cr{
	\ifx\old@cr\relax \let\old@cr\\ \fi
	\def\\{\vspace{\jot}\vspace{1pt}\old@cr}}
\newenvironment{alignjustif}[1][\qquad]{
	\[\hspace*{3\fill}
	\begin{array}{@{}>{\displaystyle}r@{}>{\displaystyle\def@cr{}}l@{#1}>{\def@cr}l}
	}{\end{array}
	\let\old@cr\relax
	\hspace*{2\fill}\]}

\def\*{\discretionary{\thinspace\the\textfont2\char1}{}{}} 

\let\ancien@point\.
\def\.{\ifmmode \mkern1.5mu \let\@next\relax \else \let\@next\ancien@point \fi \@next}

\begingroup
	\catcode`\'=\active
	\global\let\old@prime='
	\gdef'{\@ifnextchar,{\old@prime\prime@virg}{\@ifnextchar\!{\prime@excl}{\old@prime}}} 
\endgroup
\def\prime@excl#1\!{'}
\def\prime@virg#1,{\mathchoice{\mkern-0.5mu ,}{\mkern-0.5mu ,}{\mkern-2.7mu , \mkern1mu}{,}}

\def\protectsuperscript#1{\@ifnextchar^{(#1)}{#1}}
\def\open@subscr#1_#2_#3{#1_{#2#3}} 
\def\opensubscript#1{\@ifnextchar_{\open@subscr#1}{#1}} 

\newcommand{\setalignbox}[1]{%
	\def\alignbox##1{\begingroup
		\setbox0=\hbox{#1}%
		\ifmmode\hbox to \wd0{$##1$\hss}%
		\else\hbox to \wd0{##1\hss}\fi\endgroup}%
	\def\ralignbox##1{\begingroup
		\setbox0=\hbox{#1}%
		\ifmmode\hbox to \wd0{\hss$##1$}%
		\else\hbox to \wd0{\hss ##1}\fi\endgroup}}
\newcommand{\disablealignbox}{\def\alignbox##1{##1}}

\def\XXint#1#2#3{{\setbox0=\hbox{$#1{#2#3}{\int}$ }
\vcenter{\hbox{$#2#3$ }}\kern-.6\wd0}}

\newcommand{\couche}{\Lambda}
\let\layer\couche
\newcommand{\bordgauche}{\Sigma_\Omega}
\newcommand{\borddroit}{\Sigma_\Lambda}
\newcommand{\domsA}{\{\Omega,\Gamma,\couche,\Pi\}}
\newcommand{\Rcoin}{R_{\mathrm{c}}}
\newcommand{\fsource}{f_{\mathrm{s}}}
\newcommand{\Rsource}{R_{\mathrm{s}}}

\newcommand{\pst}{\frac\pi\Theta}
\newcommand{\txtpst}{\tfrac\pi\Theta}
\newcommand{\smallpst}{\pi/\Theta}
\newcommand{\psttab}{\tau}
\newcommand{\TZ}{\@ifstar{\pst\mathbb{Z}^*}{\pst\mathbb{Z}}}
\newcommand{\txtTZ}{\@ifstar{\txtpst\mathbb{Z}^*}{\txtpst\mathbb{Z}}}
\newcommand{\TN}{\pst\mathbb{N}^*}
\newcommand{\txtTN}{\txtpst\mathbb{N}^*}
\newcommand{\NN}{\mathbb{P}}

\let\peps\varepsilonLatex 
\let\feps\epsilon 
\newcommand{\powelne}{\varepsilon^p \ln^\ell \!\varepsilon}
\newcommand{\pgcd}[1]{\mathrm{gcd}(#1)}

\newcommand{\Rbasic}[2]{R_{#1}^{\raisebox{1pt}{$\scriptscriptstyle #2$}}}
\newcommand{\ROmegaPi}{R_{\Delta}}			
\newcommand{\RGammaPi}{R_{\mathrm{N}}}		
\newcommand{\RcouchePi}{R_{\partial_Y^2}}	
\newcommand{\ROmegaOmega}{\Rbasic{\Delta}{\Omega}}
\newcommand{\RGammaOmega}{\Rbasic{\mathrm{D}}{\Omega}}
\newcommand{\Aperp}{\mathcal{A}_{\perp}}

\newcommand{\phiOmega}{\phi^{\raisebox{1pt}{$\scriptscriptstyle\Omega$}}}
	
\newcommand{\bfu}{\mathbf{u}}
\newcommand{\bfv}{\mathbf{v}}
\newcommand{\bfr}{\mathbf{r}}
\newcommand{\bfx}{\mathbf{x}}
\newcommand{\bfue}{\mathbf{u}^\peps}
\newcommand{\Abarre}{\overline{\mathcal{A}}}
\newcommand{\Amacro}{\mathcal{A}^+\hspace{-0.15ex}}
\newcommand{\Amicro}{\mathcal{A}^-\hspace{-0.15ex}}
\newcommand{\Amaicro}{\mathcal{A}^\pm\hspace{-0.15ex}}
\newcommand{\Aeps}{\mathcal{A}_\feps}
\newcommand{\Aemacro}{\mathcal{A}^+_\feps}
\newcommand{\Aemicro}{\mathcal{A}^-_\feps}
\newcommand{\Aemaicro}{\mathcal{A}^\pm_\feps}
\newcommand{\Aefini}{\mathcal{A}_\feps^{\mathrm{fin}}}
\newcommand{\RIma}[1]{L_{#1}} 

\contourlength{0.05ex}
\DeclareMathOperator*{\fsum}{\mathchoice
	{\contour*{.}{$\displaystyle\sum$}}
	{\contour*{.}{$\textstyle\sum$}}%
	{\contour*{.}{$\scriptstyle\sum$}}%
	{\contour*{.}{$\scriptscriptstyle\sum$}}}
\newcommand{\fsumn}{\fsum_{n=0}^\infty}
\newcommand{\fsumd}{\fsum_{d\in\mathbb{R}}}
\newcommand{\fsump}{\fsum_{p\in\mathbb{R}}}
\newcommand{\fsumdp}{\fsum_{(d,p)\in\mathbb{R}^2}}
\newcommand{\FCX}{\mathbb{C}[[X]]}
\newcommand{\ord}{\operatorname{ord}}
\newcommand{\isfud}{\tilde{\iota}} 

\newcommand{\bfuDA}{\bfu^0}
\newcommand{\SDA}{S^\infty}
\newcommand{\uDA}{u^0}
\newcommand{\ccDA}{\mathcal{S}^\infty} 
\newcommand{\sumdmacro}{\!\!\!\fsum_{d\in\TZ*\cap(\NN-p)}\!\!\!\!}
\newcommand{\sumdmicro}{\!\!\!\fsum_{d\in\TZ*\cap(p-\NN)}\!\!\!\!}
\newcommand{\Rmacro}{\protectsuperscript{\mathcal{R}^+_\feps}}
\newcommand{\Rmicro}{\protectsuperscript{\mathcal{R}^-_\feps}}
\newcommand{\Rmaicro}{\mathcal{R}^\pm_\feps}

\def\@conemaicro#1#2{\mathchoice{K^{#1}}{K^{#1}}{K^{\mkern-#2mu #1}}{K^{#1}}}
\newcommand{\conemacro}{\@conemaicro{+}{1}}
\newcommand{\conemicro}{\@conemaicro{-}{1.5}}
\newcommand{\conemaicro}{\@conemaicro{\pm}{1}}
\newcommand{\Eps}{\mathcal{H}_\feps} 
\newcommand{\preEps}{\widetilde{\mathcal{H}}_\feps}
\newcommand{\pEps}{h} 
\newcommand{\evaleps}[1][\peps]{\mathrm{eval}_{#1}}
\newcommand{\projphi}{\pi_\sigma}
\newcommand{\cmacro}{c^{\bfu\leftarrow S}}
\newcommand{\cmicro}{c^{S\leftarrow\bfu}}
\newcommand{\cfinal}{c^{\bfu\leftarrow\bfu}}
\newcommand{\cfinalbis}{\tilde{c}^{\bfu\leftarrow\bfu}}
\newcommand{\intvl}[1]{[\![#1]\!]_{\NN\!,\mathbb{N}}}

\newcommand{\ptronc}{P}
\newcommand{\Vmacro}{\mathtt{V}^{\mathtt{u}}}
\newcommand{\Vmicro}{\mathtt{V}^{\mathtt{S}}}
\newcommand{\Pmacro}{\mathtt{P}}
\newcommand{\Pmicro}{\mathtt{Q}}
\newcommand{\Pmacroo}[2]{\Pmacro^{#1}_{#2}}
\newcommand{\Pmicroo}[2]{\Pmicro^{#1}_{#2}}
\newcommand{\smacro}{\mathtt{u}}
\newcommand{\smicro}{\mathtt{S}}
\newcommand{\sMacro}{\sigma(\bfu)}
\newcommand{\sMicro}{\sigma(S)}
\newcommand{\Vcoin}{V}
\newcommand{\crown}{C_\eta}

\newcommand{\mindeg}{\operatorname{deg_{min}}}
\newcommand{\maxdeg}{\operatorname{deg_{max}}}
\newcommand{\khiz}{\chi_0}
\newcommand{\khii}{\chi_\infty}
\newcommand{\khif}{\chi_f} 

\newcommand{\prodH}[1]{\mathcal{H}^{#1}_\times}
\newcommand{\Res}{\operatorname{Res}} 

\newcommand{\sbd@textstyle}[1]{
	\def\@sbdHmarge{3mu} \def\@sbdVmarge{0.2ex} \def\@sbdVthickness{-0.05ex}
	\mkern\@sbdHmarge{\overbracket[0.2pt][0pt]{
	\raisebox{0pt}[\height+\@sbdVthickness]{$\overbracket[0.2pt][0pt]{\mkern-\@sbdHmarge
	\raisebox{0pt}[\height+\@sbdVmarge]{$\textstyle #1$}
	\mkern-\@sbdHmarge}$}}}\mkern\@sbdHmarge}
\newcommand{\sbd@scriptstyle}[1]{
	\def\@sbdHmarge{2mu} \def\@sbdVmarge{-0.05ex} \def\@sbdVthickness{-0.05ex}
	\mkern\@sbdHmarge{\overbracket[0.2pt][0pt]{
	\raisebox{0pt}[\height+\@sbdVthickness]{$\overbracket[0.2pt][0pt]{\mkern-\@sbdHmarge
	\raisebox{0pt}[\height+\@sbdVmarge]{$\scriptstyle #1$}
	\mkern-\@sbdHmarge}$}}}\mkern\@sbdHmarge}
\newcommand{\sbd}[1]{\mathchoice{\sbd@textstyle{#1}}{\sbd@textstyle{#1}}{\sbd@scriptstyle{#1}}{\sbd@scriptstyle{#1}}}

\newcommand{\BesselGamma}{\mathcal{J}}
\newcommand{\BesselOmega}{\mathcal{J}^\Omega}
\newcommand{\BesselYOmega}{\mathcal{Y}^\Omega}%
\newcommand{\BesselYGamma}{\mathcal{Y}}%
\newcommand{\mcla}{\Phi^0} 
\newcommand{\mcll}{\Phi^{\rm lim}} 

\newcommand{\ZTZ}{\frac1a\mathbb{Z}}
\def\@AGprecis#1,#2\@AGprecis{\mathcal{A}_{#1}^{#2}(\Gamma)}
\def\@AOprecis#1,#2,#3\@AOprecis{\mathcal{A}_{#1,#3}^{#2}(\Omega)}
\def\@AMprecis#1,#2\@AMprecis{\mathcal{A}_{#2}^{+,#1}}
\newcommand{\AGprecis}[1]{\@AGprecis#1\@AGprecis}
\newcommand{\AOprecis}[1]{\@AOprecis#1\@AOprecis}
\newcommand{\AMprecis}[1]{\@AMprecis#1\@AMprecis} 
\newcommand{\normg}[1]{\|#1\|_\gamma}
\newcommand{\eAMP}{s} 
\newcommand{\Clap}[1]{C^\Delta_{#1}}
\newcommand{\eAMPmin}[1]{s_{#1}}

\newcommand{\nbcs}{N}
\newcommand{\fcs}[1][n]{\opensubscript{\varphi_{#1}}} 

\newcommand{\rmicro}{R^-}
\newcommand{\khiL}{\chi_L}
\newcommand{\evcs}{V_L}
\newcommand{\ounl}{\Omega_1^L}
\newcommand{\uL}{u_L}

\newcommand{\pef}{\mathfrak{h}}
\newcommand{\pcs}[1]{\mathfrak{s}_{#1}}
\newcommand{\SV}{S^V}
\newcommand{\SVL}{S^V_L}
\newcommand{\Srest}{S_{\mathrm{r}}}
\newcommand{\snL}[2]{s_{#1}^{#2}}
\newcommand{\psinL}[1]{\psi_n^{#1}}
\newcommand{\xcoin}{x_{\mathrm{c}}}
\newcommand{\phicoin}{\widetilde{\phi}}

\newcommand{\opF}{\mathbf{F}}
\newcommand{\uA}{\tilde{u}}
\newcommand{\uV}{u^V}
\newcommand{\Vh}{V^h}
\newcommand{\Vhg}{V^h_{\GME}}
\newcommand{\uh}{u^h_{p,\ell}}
\newcommand{\duh}[1]{u^h_{p,\ell,#1}}

\let\eff\circ
\def\Nvirg{N\mkern-2mu,\mkern1mu}
\newcommand{\RME}{\protectsuperscript{R^\eff}} 
\newcommand{\OME}{\Omega^\eff} 
\newcommand{\GME}{\Gamma^\eff} 
\newcommand{\bordgaucheME}{\bordgauche^\eff} 
\newcommand{\SME}{\Sigma^\eff} 
\newcommand{\PME}{P^\eff} 
\newcommand{\closedBall}{\overline{B}}

\def\@getcharcode#1 #2 #3\relax{`#3 }
\def\@nextcharcode{\expandafter\@getcharcode\meaning\@nextchar\relax}
\def\@eatnextchar#1{\@text}
\def\do@use@m{\let\@m=m \let\@donext=\@text
    \ifcat0\noexpand\@nextchar \ifnum\@nextcharcode>47 \ifnum\@nextcharcode<58 
	\let\@m=\@nextchar \let\@donext=\@eatnextchar \fi\fi\fi
    \@donext}
\def\use@m#1{\def\@text{#1} \futurelet\@nextchar\do@use@m}
\def\use@eps#1{
	\@ifnextchar0{\def\@text{#1} \let\@eps=0 \@eatnextchar}%
	{\let\@eps=\peps #1}}

\newcommand{\VME}{\use@m{H_{\@m}}}
\newcommand{\sme}{\use@m{\opensubscript{u^{(\@m)}_\peps}}}
\newcommand{\opS}{\use@m{\mathbf{S}^{(\@m)}_{\Nvirg\peps}}} 
\newcommand{\opSref}{\use@m{\mathbf{S}^{(\@m)}_{\Nvirg 0}}}
\newcommand{\opP}{\use@m{\mathbf{P}^{(\@m)}_{\!\Nvirg\peps}}}
\newcommand{\opD}{\use@m{\mathbf{D}^{(\@m)}_{\peps}}}
\newcommand{\scl}{\use@m{\opensubscript{u^{\rm f{}f}_{\@m,\peps}}}} 
\newcommand{\sclDA}{u^{\rm f{}f\mkern-0.5mu,\mkern0.5mu 0}_{m,\peps}}
\newcommand{\mcle}{\Phi^{N}}
\newcommand{\cme}{c^{-\leftarrow+}}
\newcommand{\cmebis}{\tilde{c}^{-\leftarrow+}}
\newcommand{\mse}{\use@eps{\Phi^{\Nvirg\@eps}}}

\newcommand{\mme}{\use@eps{\Psi^{\Nvirg\@eps}}}
\newcommand{\mmen}{\use@eps{\widetilde{\Psi}^{\Nvirg\@eps}}}
\newcommand{\mmeo}[1]{\opensubscript{\Psi^{\Nvirg\peps}_{#1,\perp}}}

\newcommand{\gram}{\use@eps{G_{\Nvirg\@eps}}}
\newcommand{\gramn}{\use@eps{\widetilde{G}_{\Nvirg\@eps}}}

\newcommand{\CdtnME}[2][\Rsource]{C^{\rm DtN\mkern-0.5mu,\mkern0.5mu #2}_{\mkern-1mu #1,N}}
\newcommand{\CdtnMEprecis}[1]{\widetilde{C}^{\rm DtN\mkern-0.5mu,\mkern0.5mu #1}_{\mkern-1mu R}}
\newcommand{\CstabME}{C^{\rm stab}_{\mkern-1mu\Rsource,N}}
\newcommand{\CstabMEprecis}{\widetilde{C}^{\rm stab}_{\mkern-1mu\Rsource}}
\newcommand{\CcontME}{C^{\rm stab}_{\rm s}}
\newcommand{\domi}[1]{\langle #1\rangle}

\newcommand{\eme}{\use@m{\opensubscript{w_\peps}}} 
\newcommand{\emeg}{\use@m{\opensubscript{w_{\peps,\GME\mkern-1.5mu}}}}
\newcommand{\emes}{\use@m{\opensubscript{w_{\peps,\SME\mkern-1.5mu}}}}
\newcommand{\snhg}[1]{\lvert#1\rvert_{\dot{H}^1(\GME)}}
\def\snh#1(#2){\lvert#1\rvert_{\dot{H}^1(#2)}}
\newcommand{\HME}[1]{H_{\rm sc}^{#1}}
\NewDocumentCommand{\errmdl}{O{0} O{m}}{e_{\rm model}^{#2,#1}}
\newcommand{\errdtn}[1][0]{e_{\rm dtn}^{R,#1}}
\def\diff{D}
\def\varlibre{\.\vcenter{\hbox{\scalebox{0.4}{$\bullet$}}} \.}

\newcommand{\Cplusplus}{C\raisebox{0.18\height}{\scalebox{0.9}{++}}}

\makeatother 

\begin{document}
\tikzifexternalizing{}{\includepdf{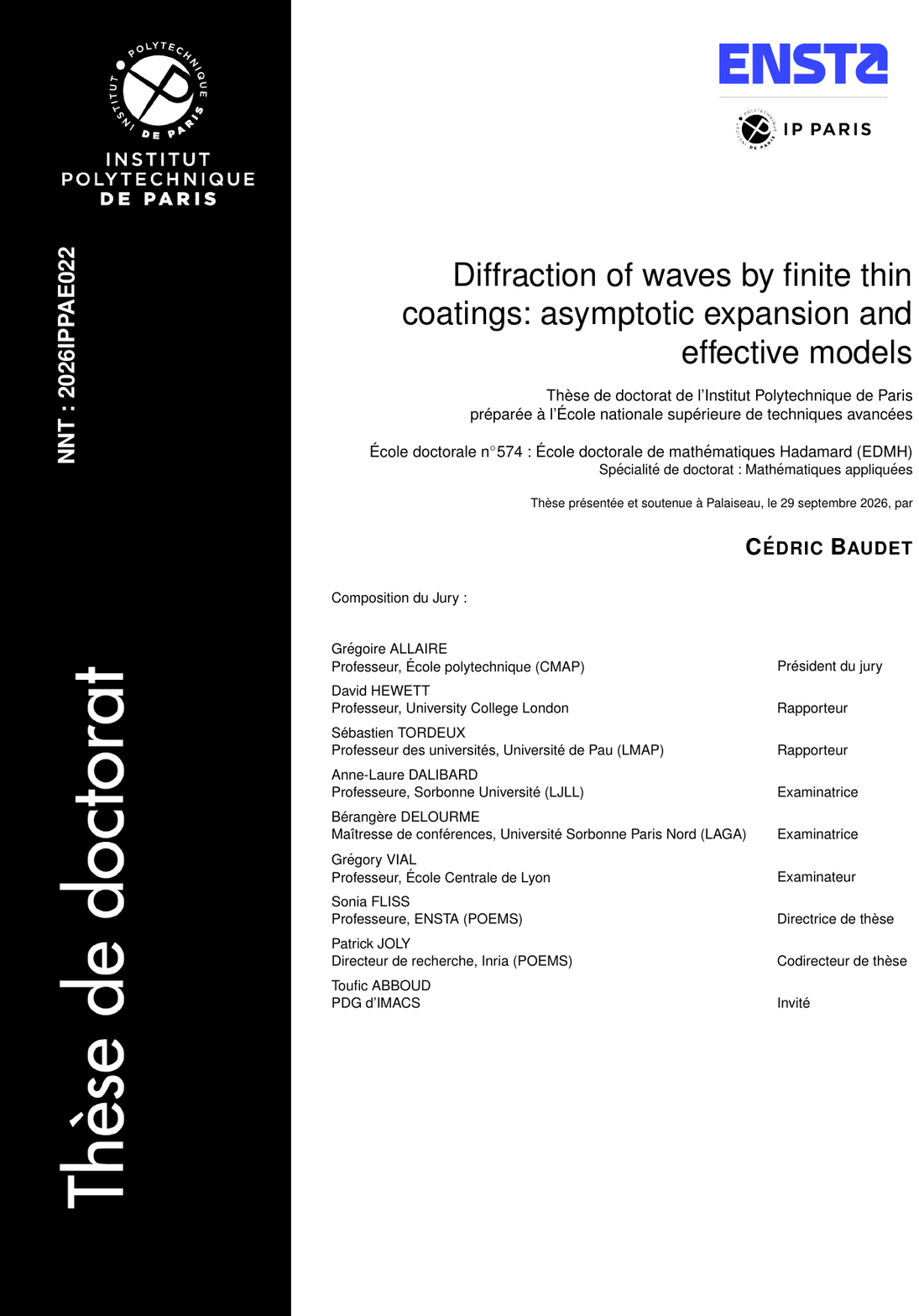}}
\setcounter{page}{1}


\chapter*{Remerciements}
\begin{otherlanguage}{french}
\setlength\parindent{2em}%

\hskip\parindent Je souhaite tout d'abord remercier mes encadrants, Sonia et Patrick. \c C'a été un plaisir de travailler avec vous durant ces quatre années et demie. Merci pour votre gentillesse, votre humour, votre patience et votre accompagnement. Vous avez à la fois toujours été disponibles quand j'avais besoin d'une réunion ou d'un conseil, et vous m'avez laissé beaucoup d'autonomie dans mon travail. Les deux pourraient s'opposer mais ça n'a pas été le cas. Quand j'ai proposé des idées qui vous étaient totalement étrangères en début de thèse, vous auriez pu les rejeter mais vous avez fait beaucoup d'efforts pour comprendre mon point de vue. Quand j'ai écrit mon premier article, nous aurions pu être les trois auteurs comme cela semble la norme dans la communauté, mais vous m'avez proposé de publier en mon nom seul ce que j'avais fait. Et malgré cela, votre investissement n'a pas diminué et il a considérablement amélioré ma manière de présenter mes travaux à cette communauté. Merci aussi Sonia pour ta patience quand je faisais de la résistance aux propositions de modifications. Et merci Patrick pour tes blagues toujours bienvenues (je retiens aussi le ``principe de conservation de l'emmerdement'' comme un théorème de la vie souvent vérifié !) et pour les invitations aux concerts que tu as organisés. En un mot, votre encadrement a été altruiste et j'ai eu beaucoup de chance d'en profiter.

Je remercie ensuite le jury d'avoir accepté d'évaluer mon travail. Merci du temps passé à étudier le manuscrit, malgré ses nombreux passages techniques, et d'avoir tous fait le déplacement, en particulier pour ceux --- Sébastien, David et Grégory --- qui venaient de loin. J'ai été très heureux de retrouver David et Bérangère que j'avais déjà eu le plaisir de rencontrer à plusieurs reprises. Thank you David for doing your best to offer me a postdoc position at UCL (which did not succeed for funding reasons). Merci Bérangère pour les discussions scientifiques sur le sujet que nous avons eues dès le début de thèse. Je crois que tu es la personne que j'ai le plus citée dans ma thèse, tu étais donc une incontournable du jury ! Merci Grégory Vial pour votre gentillesse et pour l'enthousiasme que vous avez exprimé pour mes travaux ; venant de l'auteur d'une thèse de référence sur ce sujet, cela me touche d'autant plus. Merci Toufic pour les échanges que nous avons eus depuis quatre ans dans le cadre du projet HyBox et pour le grand intérêt que vous avez montré pour ma thèse. \c C'a été un plaisir de collaborer avec vous, notamment car ça donnait du sens au travail que je faisais, et j'espère que nous aurons l'occasion de poursuivre ce projet. Je profite de ces remerciements pour réitérer mon admiration pour vos talents pianistiques que j'ai eu la chance d'entendre par deux fois ! Je remercie à nouveau Sébastien d'avoir eu le courage de faire le trajet depuis Pau et pour les discussions que nous avons eues à l'UMA et à Montréal. Enfin j'ai eu la chance de rencontrer et d'échanger à l'occasion de cette soutenance avec deux éminents chercheurs --- Grégoire Allaire et Anne-Laure Dalibard --- que je ne connaissais alors que de réputation. C'est un honneur de vous compter dans mon jury et ç'a été un plaisir de pouvoir échanger avec vous lors de la soutenance.

Je remercie les membres permanents du labo pour leur accueil chaleureux et pour leur accessibilité : quand je suis arrivé en stage de M2, mes professeurs sont devenus mes collègues, les liens hiérarchiques se sont relâchés et cela a ouvert la voie à de nombreux échanges, notamment pour en apprendre plus sur le métier de chercheur. Je remercie Nicolas et \'Eric pour leur aide avec XLiFE++ et leur réactivité quand j'avais besoin d'eux. \'Eric a notamment ajouté une fonctionnalité à XLiFE++ spécialement pour me permettre d'implémenter mon modèle effectif d'ordre 5, et Nicolas m'a aidé dans la mise en place de mon dispatcheur de tâches. Je remercie Laurent Bourgeois, Christophe Hazard et Laure Giovangigli de m'avoir fait confiance pour assurer des travaux dirigés dans leurs cours. \c C'a été une expérience très intéressante et j'ai beaucoup aimé travailler avec vous. Merci à Corinne, sur qui repose tout l'administratif du labo. Ton efficacité dans les missions et tous les tracas administratifs n'a d'égal que ta convivialité quand on vient te voir. Je remercie aussi \'Eliane pour son aide dans les démarches avec l'école doctorale ; c'est une chance d'avoir quelqu'un dans le labo à qui l'on peut demander directement des précisions quand on est perdu avec adum.fr ou la procédure de soutenance. Mention spéciale à Anne-Sophie qui m'a gentiment raccompagné en voiture chez moi un jour où j'avais fait une allergie et où j'étais trop malade pour rentrer par le RER. J'en profite pour remercier tous ceux qui m'ont porté assistance lorsque des allergies m'ont fait changer de forme et de couleur. Merci enfin à Jean-François pour sa jovialité, Axel et Patrick C. pour leur gentillesse, Frédéric pour son soutien, Zoïs, Luiz et Pierre pour leur bonne humeur à la salle de sport, Jérôme pour l'organisation des ``psaumes'', etc.

J'aimerais ajouter quelques mots pour mes anciens professeurs, même s'ils ne liront pas forcément ces remerciements, car ils ont participé à mon goût pour les mathématiques ou ont nourri ma curiosité intellectuelle. Je pense à mes professeurs de l'ENS, notamment Isabelle Gallagher pour son incroyable cours d'analyse fonctionnelle, Bertrand Maury pour son cours de modélisation riche et varié, et feu Nicolas Bergeron pour son enthousiasme communicatif en cours d'analyse complexe. Je pense ensuite à mes professeurs de prépa à Louis-le-Grand, en particulier Alain Troesch qui a construit les fondements de ma formation en mathématiques, ainsi que Jean-François Logeais dont j'admire le talent pour rendre la physique rigoureuse et élégante, à tel point que j'ai failli entrer au département de physique de l'ENS grâce à lui. Plus loin dans le temps, mais tout aussi important, viennent mes professeures de lycée Delphine Péron et Sophie Dumetz qui m'ont encouragé à viser les grandes prépas et m'ont donné une expérience initiatique de la recherche en me permettant de participer au TFJM$^2$.

Je termine cette thèse riche de nombreuses nouvelles amitiés parmi mes collègues doctorants. Je pense tout d'abord à mes voisins de bureau Sarah et Timothée. C'était une vraie joie de venir tous les matins pour vous retrouver ; le trajet depuis Paris en valait largement la peine. Merci pour votre bienveillance, votre bonne humeur et pour les innombrables fois où on a ri ensemble. Je garderai ces souvenirs précieusement. Merci aussi pour l'entraide et le soutien mutuel dans notre bureau, qui font toute la différence dans les moments difficiles. Merci Antonin pour l'organisation des nombreuses sorties karaoké et des cadeaux de soutenance. Merci aussi pour tous les bons moments que j'ai passés avec vous trois ainsi qu'avec Fabien à la salle de sport. J'avais envie de reprendre le sport au début de la thèse, mais je n'osais pas franchir les portes d'une salle de sport. \^Etre accompagné m'a beaucoup aidé à franchir le pas et à m'y investir avec constance, et ça m'a fait énormément de bien. J'en profite pour remercier notre coach Fabrice pour son accompagnement bienveillant. 
\begin{CJK*}{UTF8}{gbsn}
谢谢赖蒽，我准备去你美丽的国家的时候你帮了我很多。
\end{CJK*}
Merci Quentin pour ton accueil bienveillant qui a joué pour beaucoup dans mon intégration au labo quand je suis arrivé. Merci Louise pour ta bonne humeur quand nous partagions le même bureau. Je remercie aussi tous les autres doctorants que j'ai eu le plaisir de côtoyer au cours de mes années à l'UMA: Natalia, Alex, Morgane, Pierre A., Pierre B., Thibaut, Alice et Laura entre autres.

Je remercie également mes amis de la prépa et de l'ENS, notamment Philippe, Martin, Héloïse, Marc-Olivier et Charles. Merci à ceux qui sont venus le jour de la soutenance, et j'ai hâte de revoir bientôt ceux qui n'ont pas pu faire le déplacement. Je remercie Didier et \'Evelyne d'être venus depuis Tours pour assister à ma soutenance, ça m'a fait très plaisir. J'aimerais de plus remercier Claude Dumurgier, qui m'a aidé à deux reprises quand j'avais besoin d'un logement à Paris. Il était un propriétaire pas comme les autres, qui fonctionnait uniquement à la confiance, pas aux dossiers de location à rallonge. Merci également à Sath pour les moments passés ensemble et son soutien.

Enfin un immense merci à mes parents, sans qui tout cela n'aurait pas été possible. Vous avez mis toutes les chances de mon côté pour que je réussisse dans mes études, de la petite école jusqu'au doctorat. Vous m'avez inscrit dans un établissement privé pour m'offrir l'enseignement de la meilleure qualité possible malgré le coût que ça représentait pour une famille de la classe moyenne. Vous m'avez soutenu et aidé quand je suis parti à 17 ans pour aller étudier à Paris. Vous m'avez aussi fait voyager avec vous à chaque fois que vous le pouviez : aux \'Etats-Unis alors que j'étais encore tout petit, au Maroc, à Londres quand je n'avais pas pu y aller avec le collège pour quand même me faire découvrir l'Angleterre, en Italie, etc. Vous m'avez fait goûté à tant de délicieux plats faits maison : brioches, shortbreads, cannelés, croissants, blanquette, fondant au chocolat, langue de b\oe{}uf, et la liste pourrait encore être longue... Et vous m'avez globalement donné tellement d'amour depuis que je suis enfant que je n'ai qu'un mot : merci !
\end{otherlanguage}
\tableofcontents


\let\varepsilon=\peps
\chapter{Introduction in English}
\label{chap: intro}
\section{Motivations and setting of the problem}
\label{sec: motivations et pb etudie}

\begin{textesansboite}
Problems that involve thin layers appear in many areas, from composite materials engineering \cite{FakKha14} to biology \cite{CauHadLiNgu16}, including elasticity \cite{BelGeyKra16, FriJamMul02}, fluid mechanics \cite{RicRuyVil16, AchPirVal98, JagMik01, DalGer17} and electrochemistry \cite{VerMaiNov10}. Applications are especially numerous in electromagnetism: let us mention among others thin dielectric layers \cite{HadJia15, Mak08, SenVol95}, ferromagnetic films \cite{AmmHalHam00, HadJol01}, skin effects \cite{ZutKno05} and the Faraday cage effect \cite{DelHew20}. All these situations are numerically challenging because they require finely meshing the thin structure, which is very costly when its thickness is very small compared to the wavelength and the size of the objects.\\

To avoid this issue, one can build an approximation of the exact solution that is cheaper to compute, especially by not involving a fine mesh of the coating. There are two main ways to achieve this: either by building an asymptotic expansion of the exact solution when the thickness of the layer tends to zero, or by designing effective models that replace the layer with approximate effective boundary conditions. The latter is in general a better option because it directly provides an approximate solution, whereas the former requires building and summing several terms of the asymptotic expansion. However, asymptotic expansions are often a pathway towards effective models. In this work, we study these two approaches, one after the other.\\

As explained in \cref{sec: etat de l'art}, the case of objects entirely covered by a thin layer has already been extensively studied in the literature. On~the contrary, research on partially covering coatings is still in its early stages, both for asymptotic expansions and effective models. In~this work, we consider these latter situations. We consider a two-dimensional case and we~focus on the vicinity of the edge of the layer, which is modeled by an academic setting with a semi-infinite straight layer. In addition, we assume that the object may have a corner at that point, of arbitrary angle. We also consider a local perturbation of the geometry at the corner, at the scale of the layer's thickness.
\end{textesansboite}

\begin{textesansboite}
To describe the geometry of the problem, let us introduce the angle $\Theta  \in (0,2\pi )$ and the domains $\Omega  := \{(r \cos \theta ,r \sin \theta ) \mid r\in \mathbb{R}_+^*,\theta \in (0,\Theta )\}$, $\couche := \mathbb{R}_+^* \times (-1,0)$, $\Gamma  = \mathbb{R}_+^* \times \{0\}$, $\bordgauche=\{(r \cos \Theta ,r \sin \Theta ) \mid r\in \mathbb{R}_+^*\}$ and $\borddroit := \mathbb{R}_+^* \times \{-1\}$. All these sets are represented in \cref{fig:Omega et couche}. Then we consider $\Omega _{1} \subset \mathbb{R}^2 $ be an open set that coincides with $\Omega \cup \Gamma \cup \couche$ outside of the disc $B(0,\Rcoin)$ for some $\Rcoin \in \mathbb{R}_+^*$. In addition, let $\mu ,\rho  \in L^{\infty } (\Omega _{1} )$ be two functions greater than a positive constant (ellipticity assumption), and equal to the constants $\mu _{0}$ and $\rho _{0}$ in $\Omega  \setminus B(0,\Rcoin)$ and to the constants $\mu _{1}$ and $\rho _{1}$ in $\couche \setminus B(0,\Rcoin)$. See \cref{fig:Omega1,fig:autres exemples de Omega1} for different configurations.\\
\nomenclature[42]{$\Theta $}{Angle of the sector $\Omega $}
\nomenclature[59]{$\Omega $}{Sector of angle $\Theta $}
\nomenclature[46]{$\couche$}{Layer}
\nomenclature[53-46]{$\bordgauche$}{Outer edge of the sector $\Omega $}
\nomenclature[53-46]{$\borddroit$}{Outer edge of the layer}
\nomenclature[38a]{$\Gamma $}{Interface between the layer and the rest of the domain}
\nomenclature[59-01]{$\Omega _{1}$}{Unscaled version of the total domain (i.e.\ with a layer of thickness 1)}
\nomenclature[27-12]{$\Rcoin$}{Radius outside which $\Omega _{1}$ coincides with the union of a sector and a rectangular layer}
\nomenclature[47-52]{$\mu ,\rho $}{Functions that describe the properties of the medium}
\nomenclature[47-52-00]{$\mu _{0} ,\rho _{0} $}{Values of $\mu $ and $\rho $ outside the layer}
\nomenclature[47-52-01]{$\mu _{1} ,\rho _{1} $}{Values of $\mu $ and $\rho $ in the layer}

\begin{figure}[h]
\def\widthFig{2.3} 
\def\heightFig{2} 
\def\ThetaFig{205} 
\begin{tikzpicture}[scale=0.9]
	\coordinate (finLambda) at (-\widthFig, {-\widthFig*tan(\ThetaFig-180)});
	\fill [pattern={Lines[angle=45,line width=.2pt]}, pattern color=white!70!gray] (0,0) -- (\widthFig,0) -- (\widthFig,\heightFig) -- (-\widthFig,\heightFig) -- (finLambda) -- cycle;
	\draw [thick] (finLambda) -- (0,0) -- (\widthFig,0);
	\draw [>=stealth,->] (-\widthFig-0.4,0) -- (\widthFig+0.4,0) node [right] {$x$};
	\draw [>=stealth,->] (0,{-\widthFig*tan(\ThetaFig-180)}) -- (0,{\heightFig+0.3}) node [right] {$y$};
	\draw (-\widthFig*2/5,\heightFig/2) node {$\Omega$};
	\draw [->,>=latex] (0.35,0) arc (0:\ThetaFig-5:0.35) node[near start,above] {$\Theta$} -- (\ThetaFig:0.36);
	\draw [<-,>=latex] (-\widthFig/2, {-\widthFig/2*tan(\ThetaFig-180)}) to[bend right=20] ++(0.2,-0.45) node[right]{$\!\bordgauche$};
	\draw [<-,>=latex] (\widthFig/2,0) to[bend right=20] ++(0.15,-0.5) node[below]{$\ \ \Gamma$};
\end{tikzpicture}
\hspace{-8mm}
\def\epsFig{0.8} 
\def\widthFig{3.5}
\begin{tikzpicture}
	\fill [color=white!70!gray] (0,0) rectangle (\widthFig,-\epsFig);
	\draw [thick] (\widthFig,-\epsFig) -- (0,-\epsFig) -- (0,0) -- (\widthFig,0);
	\draw (-2.5,0) node{\ }; 
	\draw [>=stealth,->] (-0.4,0) -- (\widthFig+0.4,0) node [right] {$x$};
	\draw [>=stealth,->] (0,{-\epsFig-0.5}) -- (0,0.5) node [above] {$Y$};
	\draw (\widthFig/2,-\epsFig/2) node {layer $\couche$};
	\draw [<->,>=latex] (\widthFig*0.85,0) -- node[midway,right]{1} (\widthFig*0.85,-\epsFig);
	\draw [<-,>=latex] (\widthFig/2,0) to[bend left=20] ++(0.15,0.5) node[above]{$\ \ \Gamma$};
	\draw [<-,>=latex] (\widthFig/2,-\epsFig) to[bend right=20] ++(0.15,-0.5) node[below]{$\ \ \borddroit$};
\end{tikzpicture}\par
\caption[The domains $\Omega$ and $\couche$]{The domains $\Omega$ (on the left) and $\couche$ (on the right)}
\label{fig:Omega et couche}
\end{figure}

\begin{figure}[h]
\def\epsFig{0.8} 
\def\widthFig{3.5} 
\def\heightFig{2.3} 
\def\ThetaFig{210} 
\def\RCoinFig{(\epsFig/sin(\ThetaFig-180)+0.2)}
\begin{tikzpicture}[scale=0.9]
	\coordinate (finLambda) at (-\widthFig, {-\widthFig*tan(\ThetaFig-180)});
	\coordinate (sommetCouche) at ({-\epsFig*cot(\ThetaFig-180)}, -\epsFig);
	\fill [color=white!70!gray] (0,0) -- (\widthFig,0) -- (\widthFig,-\epsFig) -- (sommetCouche) -- cycle;
	\fill [pattern={Lines[angle=45,line width=.2pt]}, pattern color=white!70!gray] (0,0) -- (\widthFig,0) -- (\widthFig,\heightFig) -- (-\widthFig,\heightFig) -- (finLambda) -- cycle;
	\draw [thick] (finLambda) -- (0,0) -- (\widthFig,0);
	\draw [thick] (sommetCouche) -- (\widthFig,-\epsFig);
	\draw (-\widthFig-1.6,0) node {\ }; 
	\draw [>=stealth,->] (-\widthFig-0.4,0) -- (\widthFig+0.4,0) node [right] {$X$};
	\draw [>=stealth,->] (0,{-\RCoinFig-0.3}) -- (0,{\heightFig+0.3}) node [right] {$Y$};
	\draw (\widthFig*0.27,-\epsFig/2) node {layer};
	\draw (\widthFig*0.4,{-\RCoinFig-0.1}) node {obstacle};
	\draw (-\widthFig*0.7,\heightFig*0.6) node {$(\mu_0,\rho_0)$};
	\draw [->,>=latex] (0.35,0) arc (0:\ThetaFig-5:0.35) node[near start,above] {$\Theta$} -- (\ThetaFig:0.36);
	\draw [decorate, decoration={brace,amplitude=10pt}, xshift=4pt] (\widthFig+0.7,\heightFig) -- (\widthFig+0.7,-\epsFig) node [midway,right,xshift=10pt] {$\Omega_1$};
	\draw [<-,>=latex] (\widthFig*2/3,-\epsFig*0.6) to[bend right=20] ++(0.2,-0.6) node[below]{$\ \ (\mu_1,\rho_1)$};
	\draw [<->,>=latex] (\widthFig*0.85,0) -- node[midway,right]{1} (\widthFig*0.85,-\epsFig);
	\draw [dashed] (0,0) circle ({\RCoinFig});
	\draw [<->,>=latex] (0,0) -- node[near end, below] {$\Rcoin$} ({\RCoinFig*cos(150)}, {\RCoinFig*sin(150)});
\end{tikzpicture}\par
\caption{The domain $\Omega_1$ with a configuration example of $\Omega_1 \cap B(0,\Rcoin)$}
\label{fig:Omega1}
\end{figure}
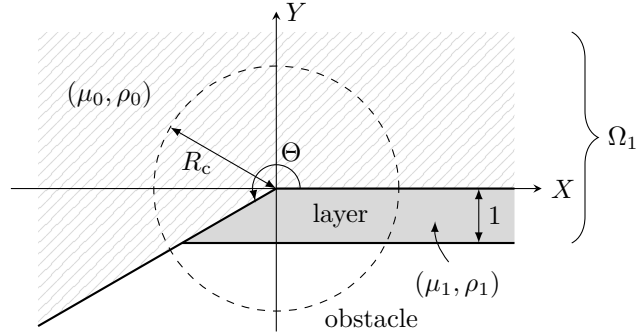

\begin{figure}[h]
\def\epsFig{0.8} 
\def\widthFig{2.4} 
\def\heightFig{2.1} 
\def\ThetaFig{210} 
\def\RCoinFig{(\epsFig/sin(\ThetaFig-180)+0.2)}
\begin{tikzpicture}[scale=0.9]
	\coordinate (finLambda) at (-\widthFig, {-\widthFig*tan(\ThetaFig-180)});
	\coordinate (sommetJoint) at ({-\epsFig*cot(\ThetaFig-180)}, -\epsFig);
	\coordinate (sommetCoucheHaut) at ({-\epsFig*cot(\ThetaFig-180)+\epsFig}, 0);
	\coordinate (sommetCoucheBas) at ({-\epsFig*cot(\ThetaFig-180)+\epsFig}, -\epsFig);
	\fill [pattern={Lines[angle=45,line width=.2pt]}, pattern color=white!70!gray] (0,0) -- (\widthFig,0) -- (\widthFig,\heightFig) -- (-\widthFig,\heightFig) -- (finLambda) -- cycle;
	\fill [color=white!70!gray] (sommetCoucheBas) -- (sommetCoucheHaut) -- (\widthFig,0) -- (\widthFig,-\epsFig) -- cycle; 
	\fill [color=white!50!gray] (sommetJoint) arc (180:90:\epsFig) -- (sommetCoucheBas) -- cycle; 
	\draw [thick] (finLambda) -- (sommetJoint) arc (180:90:\epsFig) -- (\widthFig,0);
	\draw [thin,dashed] (0,0) -- (sommetJoint);
	\draw [thick] (sommetJoint) -- (\widthFig,-\epsFig);
	\draw [thick] (sommetCoucheBas) -- (sommetCoucheHaut);
	\draw [>=stealth,->] (-\widthFig-0.3,0) -- (\widthFig+0.3,0);
	\draw [>=stealth,->] (0,{-\RCoinFig-0.2}) -- (0,{\heightFig+0.2});
	\draw [dashed] (0,0) circle ({\RCoinFig});
	\draw [<->,>=latex] (0,0) -- node[near start, right] {\;$\Rcoin$} ({\RCoinFig*cos(50)}, {\RCoinFig*sin(50)});
	\draw [<-,>=latex] (-\epsFig*1.05,-\epsFig*0.4) to[bend left=20] ++(-0.15,0.5) node[above]{a joint};
	\draw [<-,>=latex] (-\epsFig*0.95,-\epsFig*0.7) to[bend right=20] ++(0.03,-0.5) node[below]{\ \ $(\mu_2 ,\rho_2)$};
\end{tikzpicture}%
\nolinebreak\hspace{10mm}
\def\ThetaFig{150} 
\def\RCoinFig{(\epsFig/sin(180-\ThetaFig)+0.2)}
\begin{tikzpicture}[scale=0.9]
	\coordinate (finLambda) at (-\widthFig, {-\widthFig*tan(\ThetaFig-180)});
	\coordinate (sommetCoucheBas) at ({-\epsFig*cot(\ThetaFig-180)}, -\epsFig);
	\coordinate (sommetCoucheHaut) at (\ThetaFig-90:\epsFig);
	\coordinate (sommetCoucheMilieu) at ({\epsFig*tan(\ThetaFig-90)}, 0);
	\fill [pattern={Lines[angle=45,line width=.2pt]}, pattern color=white!70!gray] (0,0) -- (\widthFig,0) -- (\widthFig,\heightFig) -- (-\widthFig,\heightFig) -- (finLambda) -- cycle;
	\fill [color=white!70!gray] (0,0) -- (sommetCoucheHaut) -- (sommetCoucheMilieu) -- (\widthFig,0) -- (\widthFig,-\epsFig) -- (sommetCoucheBas) -- cycle;
	\draw [thick] (finLambda) -- (sommetCoucheBas) -- (\widthFig,-\epsFig);
	\draw [thick] (0,0) -- (sommetCoucheHaut) -- (sommetCoucheMilieu) -- (\widthFig,0);
	\draw [>=stealth,->] (-\widthFig-0.3,0) -- (\widthFig+0.3,0);
	\draw [>=stealth,->] (0,{-\RCoinFig-0.2}) -- (0,{\heightFig+0.2});
	\draw [->,>=latex] (0.35,0) arc (0:\ThetaFig-5:0.35) node[near end,above] {$\Theta$\;} -- (\ThetaFig:0.36);
	\draw [dashed] (0,0) circle ({\RCoinFig});
	\draw [<->,>=latex] (0,0) -- node[near end, right] {\ \;$\Rcoin$} ({\RCoinFig*cos(210)}, {\RCoinFig*sin(210)});
\end{tikzpicture}\par
\caption{Other configuration examples of $\Omega_1 \cap B(0,\Rcoin)$ for different values of $\Theta$}
\label{fig:autres exemples de Omega1}
\end{figure}
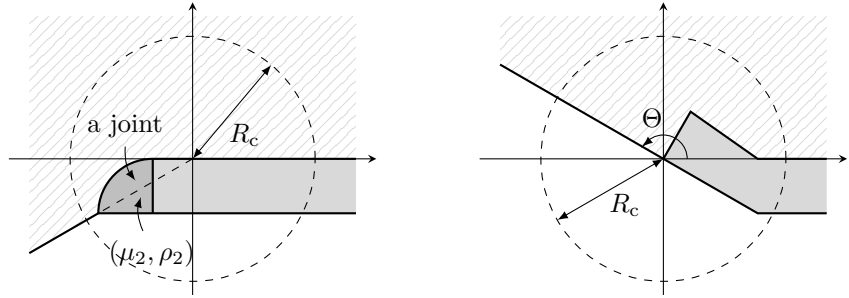

\nomenclature[40a]{$\peps$}{Thickness of the layer\nomnorefpage}
\nomenclature[33-34]{$X,Y$}{Scaled coordinates $( \frac{x}{\peps} , \frac{y}{\peps} )$}
\nomenclature[59-40]{$\Omega _{\peps}$}{Scaled version of the total domain (i.e.\ with a layer of thickness $\varepsilon $)}
\nomenclature[47-52-40]{$\mu _{\peps} ,\rho _{\peps}$}{Scaled version of the functions $\mu $ and $\rho $}
\nomenclature[15-28]{$\fsource$}{Source term of the problem satisfied by $u_{\peps}$}
\nomenclature[59]{$\omega $}{Pulsation ($2\pi $ times the frequency)}
\nomenclature[30-40]{$u_{\peps}$}{Exact solution}
Let $\varepsilon >0$ be the thickness of the layer. The physical domain is given by
\[\Omega _{\varepsilon } := \left\{(x,y)\in \mathbb{R}^2  \ \Big|\ \left( \frac{x}{\varepsilon } , \frac{y}{\varepsilon } \right) \in  \Omega _{1} \right\}.\]
Note that $\Omega _{\varepsilon }$ coincides with $\Omega _{1}$ when $\varepsilon =1$. We introduce the scaled coefficients $\mu _{\varepsilon } :(x,y)\in \Omega _{\varepsilon } \mapsto  \mu ( \frac{x}{\varepsilon } , \frac{y}{\varepsilon } )$ and $\rho _{\varepsilon } :(x,y)\in \Omega _{\varepsilon } \mapsto  \rho ( \frac{x}{\varepsilon } , \frac{y}{\varepsilon } )$, and the scaled variables $X := \frac{x}{\varepsilon }$ and $Y:= \frac{y}{\varepsilon }$. Let $\omega \in \mathbb{C} \setminus \mathbb{R}$ and $\fsource\in H^{-1} (\Omega )$ be a source term s.t.\ $\operatorname{dist}(\operatorname{supp}(\fsource),\Gamma )>0$. We denote $u_{\varepsilon }$ the unique solution in $H^{1} _{0} (\Omega _{\varepsilon } )$ of 
\begin{equation}{\label{eq: pb etudie}}
\left\{\begin{array}{r@{}ll}
    \operatorname{div}(\mu _{\varepsilon } \nabla  u_{\varepsilon } )+\omega ^2 \rho _{\varepsilon } u_{\varepsilon } &{}= \fsource &\text{ in }\Omega _{\varepsilon } \\[1pt]
    u_{\varepsilon } &{}= 0 &\text{ on }\partial \Omega _{\varepsilon }
\end{array}\right.
\end{equation}
Since $\Im(\omega ) \neq 0$, this problem is well-posed in $H^{1} _{0} (\Omega _{\varepsilon } )$ (it suffices to use the Lax-Milgram theorem) with a stability constant independent of $\varepsilon $ :
\begin{equation}{\label{eq: uniforme coercivite}}
\exists C>0,\;\forall \fsource \in H^{-1} (\Omega ), \;\forall \varepsilon >0, \quad  \|u_{\varepsilon } \|_{H^{1} (\Omega _{\varepsilon } )} \leqslant  C \.\|\fsource\|_{H^{-1} (\Omega )}
\end{equation}
\end{textesansboite}

\begin{remarque}
$\Im(\omega ) \neq 0$ is a technical assumption that is linked to the fact that we consider an academic setting with a semi-infinite layer. In this situation, $\Im(\omega )=0$ would require designing a specific radiation condition, which is a difficult question. For the more realistic case of a bounded obstacle partially covered by a thin layer, a classical Sommerfeld radiation condition can be used. Note that assuming $\Im(\omega ) \neq 0$ has little impact on the asymptotic analysis.
\end{remarque}

\begin{paragraphesansboite}{Notations:}
Throughout this document, we use the following notations:
\begin{itemize}
    \item $(x,y)$ denotes the cartesian coordinates, and $(r,\theta )$ denotes the polar coordinates with $\theta  \in [0,2\pi )$,
    \item for~any $r>0$, $B(0,r)$ is the disc of $\mathbb{R}^2 $ of radius $r$ centered at $(0,0)$,
    \item $k_{i} := \omega  \sqrt{\rho _{i} /\mu _{i} }$ with $i=0,1$ is the wave number in $\Omega $, resp.\ in the layer,
    \item $\NN := \mathbb{N}+ \frac{\pi }{\Theta } \mathbb{N}$,
	\item $\Rsource := \operatorname{dist}(0,\operatorname{supp}\fsource)$, which is positive by assumption on $\fsource$,
	\item $\bfr := r$ in $\Omega $ and $\bfr:=x$ in $\couche$,
    \item $\alpha  := e^{-\mathrm{i}\Theta }$.
\end{itemize}
In addition, we use the following conventions:
\begin{itemize}
    \item $\mathbb{N}:=\{0,1,2,3,\dots \}$,
    \item $\mathbb{N}^* :=\mathbb{N} \setminus \{0\}$, $\mathbb{R}^* :=\mathbb{R} \setminus \{0\}$ and likewise for $\mathbb{Z}^*$, $\mathbb{Q}^*$ and $\mathbb{C}^*$,
    \item $\mathbb{R}_+ :=[0,\infty )$, $\mathbb{R}_+^* :=(0,\infty )$ and likewise for $\mathbb{R}_-$ and $\mathbb{R}_-^*$.
	\item For any $n,m\in \mathbb{N}^*$ and any open set $O\subset \mathbb{R}^{n}$, $H^{m} _{\mathrm{loc}} (O)$ denotes the set of functions that are $H^{m}$ on every open subset of $O$ whose closure is a compact subset of $O$. Similarly for $L^2 _{\mathrm{loc}}$.
\end{itemize}
If needed, a more complete list of the notations can be found at the end of this document, see page~\pageref{nomenclature}.
\end{paragraphesansboite}
\nomenclature[25]{$\NN$}{$\mathbb{N}+ \frac{\pi }{\Theta } \mathbb{N}$\nomnorefpage}
\nomenclature[27-43]{$r,\theta $}{Polar coordinates with $\theta  \in [0,2\pi )$\nomnorefpage}
\nomenclature[36]{$\alpha $}{$e^{-\mathrm{i}\Theta }$\nomnorefpage}
\nomenclature[20-00-01]{$k_{0} ,k_{1}$}{Wave number outside, resp.\ inside, the layer: $\forall i\in \{0,1\},\ k_{i} := \omega  \sqrt{\rho _{i} /\mu _{i} }$\nomnorefpage}
\nomenclature[17-21-24]{$H^{k} _{\mathrm{loc}} (O)$}{The set of functions that are $H^{k}$ on every open set $O'$ s.t. $\overline{O'}$ is a compact subset of $O$\nomnorefpage}
\nomenclature[21-21-24]{$L^2 _{\mathrm{loc}} (O)$}{The set of functions that are $L^2 $ on every compact subset of $O$\nomnorefpage}

\section{State of the art}
\label{sec: etat de l'art}

	\tikzset{external/export=true}%
	\tikzsetfigurename{fig-couche_infinie-}%
	\begin{figure}[h]

\begin{subfigure}[c]{0.5\textwidth}
\centering
\def\epsFig{0.7} 
\def\widthFig{2.9} 
\def\heightFig{2.8} 
\begin{tikzpicture}[scale=0.9]
	\fill [color=white!70!gray] (-\widthFig,-\epsFig) rectangle (\widthFig,0);
	\fill [pattern={Lines[angle=45,line width=.2pt]}, pattern color=white!70!gray] (-\widthFig,0) rectangle (\widthFig,\heightFig);
	\draw [thick] (-\widthFig,0) -- (\widthFig,0);
	\draw [thick] (-\widthFig,-\epsFig) -- (\widthFig,-\epsFig);
	\draw (0,\heightFig/2) node {\parbox{25mm}{\centering exterior domain\\ of propagation}};
	\draw (0,-\epsFig/2) node {layer};
	\draw (0,-\epsFig-0.5) node {obstacle};
	\draw [<->,>=latex] (\widthFig*0.75,0) -- node[midway,right]{$\peps$} (\widthFig*0.75,-\epsFig);
\end{tikzpicture}
\end{subfigure}%
\begin{subfigure}[c]{0.5\textwidth}
\centering
\def\epsFig{0.8} 
\def\widthFig{3.5} 
\def\heightFig{2.55} 
\def\cycleExt{(-0.7,4.5) (-0.3,3) (1,2) (3,2) (5.2,1.7) (7,2.8) (7.8,5) (7,7) (5,7.5) (3,6.8) (1,6.9) (-0.3,6)}
\def\cycleInt{(-1.15,4.5) (-0.65,3) (0.7,2.15) (3,2.15) (5.2,1.85) (7.05,2.8) (8,5) (7.25,7.15) (5.1,7.6) (3,6.8) (0.6,6.95) (-0.8,6)}
\def\tension{0.7}
\begin{tikzpicture}[scale=0.9]
	\fill [pattern={Lines[angle=45,line width=.2pt]}, pattern color=white!70!gray] (-\widthFig,-\heightFig) rectangle (\widthFig,\heightFig);
	\draw (0,-\heightFig+0.4) node {exterior domain of propagation};
	\begin{scope}[shift={(-2.1,-2.45)}, scale=0.6]
		\fill [color=white!70!gray] plot [smooth cycle, tension=\tension] coordinates {\cycleExt} ;
		\draw [thick] plot [smooth cycle, tension=\tension] coordinates {\cycleExt} ;
		\begin{scope}[shift={(0.8,0.8)}, scale=0.8]
			\fill [color=white] plot [smooth cycle, tension=\tension] coordinates {\cycleInt} ;
			\draw [thick] plot [smooth cycle, tension=\tension] coordinates {\cycleInt} ;
		\end{scope}
		\draw (3.5,4.5) node {obstacle};
		\draw [<-,>=latex] (3.5,2.12) to[bend left=20] ++(0.5,0.8) node[right]{\!layer};
		\draw [<->,>=latex] (5.3,6.9) -- node[midway,right]{$\peps$} (5.3,7.5);
	\end{scope}
\end{tikzpicture}
\end{subfigure}

\caption{An infinite planar layer (on the left) and a smooth curved layer (on the right)}
\label{fig: couche infinie}
\end{figure}
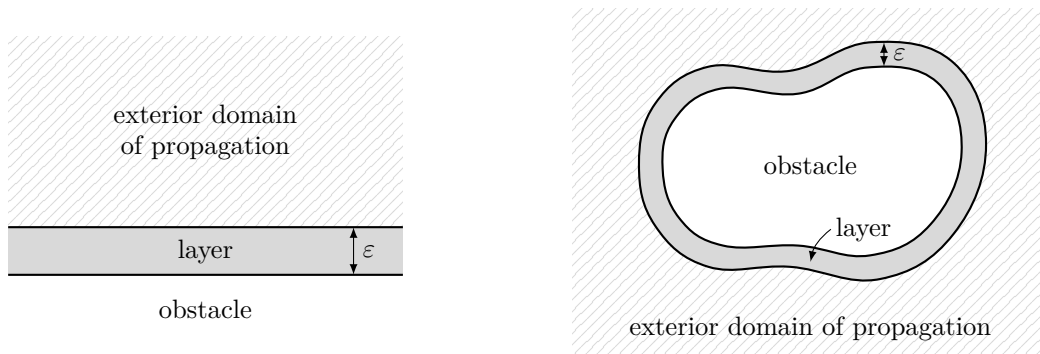%
	\tikzset{external/export=false}%
	\tikzsetfigurename{main-figure}

\begin{textesansboite}
Infinite planar layers and smooth curved layers (see \cref{fig: couche infinie}) were studied during the 90s in \cite{BenLem96, EngNed93}. Their method is to look for an asymptotic expansion of the exact solution that has the following form, denoting $s$ the curvilinear coordinate along the layer and $\nu $ the coordinate along its normal (e.g.\ $(s,\nu )=(x,y)$ for a straight layer),
\begin{equation}{\label{eq: ansatz couche infinie homogene}}
\left\{\begin{array}{ll}
    \displaystyle \sum_{n=0}^\infty  \varepsilon ^{n} \.u_{n} (s,\nu ) & \text{outside the layer }  (\nu >0) \vspace{1pt}\\
    \displaystyle \sum_{n=0}^\infty  \varepsilon ^{n} \.U_{n} \Big(s, \frac{\nu }{\varepsilon } \Big) & \text{inside the layer } (-\varepsilon <\nu <0).
\end{array}\right.
\end{equation}
The variable $\frac{\nu }{\varepsilon }$ in $U_{n}$ allows to describe the behavior of the solution at the microscopic scale in the layer. The terms of this expansion can be computed by induction w.r.t.\ $n$. Then \cite{BenLem96, EngNed93} derive effective models that replace the layer by an effective boundary condition. The simplest one is the limit problem, which imposes the same boundary condition as the original problem, simply without layer.
Its solution is $u_{0}$, the first term of the expansion. It converges to the exact solution outside the layer in $O(\varepsilon )$ for the $H^{1}$~norm.
This model is~not satisfactory for applications that are interested in studying the effects of the layer, as it completely ignores the thin structure.
More interestingly, \cite{BenLem96, EngNed93} propose higher order effective models. For instance, if one considers an infinite planar layer with a Dirichlet condition, these models involve a boundary condition of Robin, resp.\ Ventcel, type that can be written with our notations as
\begin{equation}{\label{eq: conditions effectives classiques intro}}
    \mu _{0} \partial _{n} v_{\varepsilon } = -\frac{\mu _{1} }{\varepsilon } v_{\varepsilon } , \qquad  \text{resp.}  \qquad  \mu _{0} \partial _{n} v_{\varepsilon } = -\frac{\mu _{1} }{\varepsilon } v_{\varepsilon } + \frac{\varepsilon }{3} (\mu _{1} \partial _{x} ^2 +\omega ^2 \rho _{1} )v_{\varepsilon } .
\end{equation}
This provides a model that converges at order $O(\varepsilon ^{3} )$, resp.\ $O(\varepsilon ^{5} )$, outside the layer. In the simple case of a planar infinite layer, these models can be directly obtained by expressing the Dirichlet-to-Neumann map associated to the layer using Fourier transformation, and doing a Taylor expansion w.r.t.\ $\varepsilon $ of the obtained expression.\\

\def\highlight#1{\emph{#1}}%
Those results were later extended to heterogeneous and periodic layers in \cite{AmmHe97, AmmHe98, AmmLat98} and more recently in \cite{MadVal07, DelHadJol12, DelCla13, Del14, DelHew20, Bou23}. In this case, the variations of the solution at the microscopic scale are not confined in the coating, like in the homogeneous setting, but leak in a thin boundary layer outside the physical coating. Thus, building an asymptotic expansion requires more involved techniques. These techniques are also useful to study layers in corner, so let us introduce them. They are two main types of methods: \highlight{multiscale asymptotic expansions} and \highlight{matched asymptotic expansions}. See \cite{VisLyu57, Naz81a, Naz81b} and \cite{Van64, Ili92, Eck94} respectively for a general presentation, and \cite{TorViaDau06} for a comparison of these methods in the context of a corner without layer. These methods both involve two types of fields: ``far fields'' that depend on the macroscopic scale described by $(x,y)$, and ``near fields'' that depend on the microscopic scale described by $( \frac{x}{\varepsilon } , \frac{y}{\varepsilon } )$. In multiscale expansions, far and near fields are defined and coexist in the whole domain (of a significant part of it) and the near fields tend to 0 towards infinity so that they describe the boundary layer effect in the vicinity of the layer, called the ``near zone''. In contrast, matched asymptotic expansions involve near fields only in the near zone and far fields only in the ``far zone'' (the rest of the domain), and the near and far fields have to coincide asymptotically in an intermediate zone called the ``matching zone''. This zone is located at distance $y\in (\varepsilon ^{a} ,2\mkern1mu \varepsilon ^{a} )$ of the layer with $a\in (0,1)$ fixed, so that, seen from the far fields, it tends to the interface with the layer, but seen from the near fields it tends to infinity, as $\varepsilon ^{a} /\varepsilon  \rightarrow \infty $ when $\varepsilon \rightarrow 0$. Hence matching asymptotic expansions involve a ``matching assumption'' that links the behaviors near the layer of the far fields and the behaviors at infinity of the near fields. In this case, near fields typically diverge at infinity, which is an important difference with multiscale expansions. Let us mention that these techniques have also been applied to related problems, including 3D settings \cite{Cham16}, random layers \cite{BouFliGio26} and boundary layers in homogenization \cite{AllAma99, Vin16, Ben21}.\\

Problems involving finite layers or layers in a corner present additional difficulties due to the presence of singular behaviors near the corner (or the edge of a flat finite layer). For instance, if the exact solution has a Dirichlet condition, then the function $u_{1}$ in \cref{eq: ansatz couche infinie homogene} satisfies the Dirichlet boundary condition $u_{1} = - \frac{\mu _{0} }{\mu _{1}} \partial _{n} u_{0}$. For an infinite layer, one can show that $u_{0}$ is smooth near the layer, so $\partial _{n} u_{0} \in H^{1/2}$. However, in a corner of angle $\Theta $, one has \[u_{0} (r,\theta )= c\,r^{\pst} \sin \!\Big(\pst \theta \Big) + u_{0} ^{\mathrm{reg}} \qquad  \text{with }  c\in \mathbb{C} \text{ and } u_{0} ^{\mathrm{reg}} \in H^2  .\]
The function $r^{\pst} \sin (\pst \theta )$ is called a ``corner singularity'' because it is not regular at the corner (e.g.\ not $H^2 $ when $\Theta >\pi $). Hence, one has $\partial _{n} u_{0} = -\pst c\.r^{\pst-1} + \partial _{n} u_{0} ^{\mathrm{reg}}$, which does~not belong to $H^{1/2}$ when $\pst<1$, i.e.\ $\Theta >\pi $ (and $c\neq 0$). In this case, one cannot look for $u_{1}$ in $H^{1}$ and has to consider a far field $u_{1}$ that diverges at the corner. For a smaller angle $\Theta $, $u_{1}$ is well-defined in $H^{1}$, but a similar problem arises for some far field $u_{n}$ with $n>1$.\\

These problems involving corners were first studied for Poisson equation with an arbitrary angle in \cite{CalCosDauVia06, Via03} (we cite together papers and associated PhD theses). This work provides an asymptotic expansion at any order using the method of multiscale asymptotic expansions. The obtained expansion has the form
\[\sum _{p\in \mathbb{N}+\pst \mathbb{N}} \.\sum _{\ell =0} ^{\infty } \.\varepsilon ^{p} \ln ^{\ell } \!\varepsilon  \cdot  u_{p,\ell } (x,y)\]
far from the corner. Note that it involves non-integer powers of $\varepsilon $ in $\NN:=\mathbb{N}+\pst \mathbb{N}$ and integer powers of $\ln \varepsilon $, which are linked to the corner singularities. This extends a previous work \cite{CosDau96} that studies a Poisson problem without layer, but with a mixed boundary condition of Neumann $\partial _{n} u_{\varepsilon } =0$ and Robin $u_{\varepsilon } +\varepsilon  \.\partial _{n} u_{\varepsilon } =0$. The expansion in this case has similarities (linked to the fact that this Robin condition is one of the effective conditions of \cref{eq: conditions effectives classiques intro}), but is less technical to obtain. These asymptotic expansions at any order were generalized to periodic layers in \cite{DelSchSem16}, still for Poisson's problems. This work uses at the same time a multiscale asymptotic expansion around the layer and a matched asymptotic expansion around the corner to handle the complex coupling of these two types of small structures.\\

In comparison, Helmholtz problems involving finite layers not only present the same difficulties, but also lead to more complex singularities. Indeed, the corner singularities that appear at any order cannot be organized into a family of inductive sequences, like in Poisson problems, because they are related to one another through a tree of dependencies (more details in \cref{rq: utilite des operateurs R}). Hence, building an asymptotic expansion at any order and for any angle is very challenging. This may explain a posteriori why previous asymptotic expansions for the Helmholtz equation were limited to the order 2, which avoids the issue by explicitly describing the first singularities. \cite{BenMakTor11, Mak08} proposes an expansion at order 2 for a homogeneous layer with a flat angle ($\Theta =\pi $), while \cite{SemDelSch18} generalizes this to a periodic layer with an arbitrary corner angle. Both use matched asymptotic expansions. These expansions of limited order are difficult to justify rigorously because one needs more terms of the expansion in order to prove an optimal error estimate at a given order. \cite{BenMakTor11, Mak08, SemDelSch18} either prove a sub-optimal error estimate or assume that the following terms can be build like the previous one (which seems hasty to us, given the rapid increase in complexity of the expansion w.r.t.\ the order).\\

These works on finite layers present an increasing level of complexity in calculations, due to the raising intricacy of the corner singularities. We propose to remedy this thanks to a system of algebraic calculus tools that replaces cumbersome calculations with abstract generalizable ones. This allows us to treat a more complex case -- Helmholtz's equation with a finite layer at any order for any angle -- while providing more explicit matching conditions. Besides, let us mention that \cite{JolTor06, JolTor08, Tor04} proved an asymptotic expansion  at any order for Helmholtz's equation in the related case of a half plane with a thin slot, using matched asymptotic expansion. Hence the complexity that we face here seems quite specific to the combination of the Helmholtz equation with a finite layer.\\

As explained above, building asymptotic expansions is much harder for finite layers than infinite ones, but this~is even more true for effective models. In fact, high order effective models for finite layers remained an open question, even in simpler cases like Poisson's equation or specific angles like $\Theta =\pi $, and despite a substantial literature on asymptotic expansions (see above).
A naive approach consists in imposing the effective boundary conditions \cref{eq: conditions effectives classiques intro} on $\Gamma $ entirely, i.e.\ up to the corner. For the Robin condition for instance, this writes as:
\[\left\{\begin{array}{rl}
	\mu _{0} \Delta  v_{\varepsilon } +\omega ^2 \rho _{0} v_{\varepsilon } = \fsource &\text{ in }\Omega \\[2pt]
	\mu _{0} \partial _{n} v_{\varepsilon } + \frac{\mu _{1} }{\varepsilon } v_{\varepsilon } =0&\text{ on }\Gamma .
\end{array}\right.\]
\cite{Via05,Via03} proved for this model an optimal rate of convergence of $O(\varepsilon ^{\min (3,\pst)} )$ in the total domain (and we think one has $O(\varepsilon ^{\min (3,2\pst)} )$ far from the corner). This deterioration is due to the presence of the corner singularities and is worse for large angles because the larger the angle, the less regular the singularities like $r^{\pst} \sin (\pst \theta )$. In addition, the Ventcel condition of \cref{eq: conditions effectives classiques intro} involves a second order tangential derivative, so it requires an additional boundary condition at the corner. The choice of the appropriate condition is not obvious at all in general.\\

\cite{AuvVia19,AuvVia18,Auv18} proposed two effective models for Poisson's equation with $\Theta =\pi $. These models have order 1.5 for the $L^2 $~norm and 0.5 for the $H^{1}$~norm. The first one considers on $\Gamma $ the Robin condition $\mu _{0} \partial _{n} u=-\frac{\mu _{1} }{\varepsilon } u$ of \cref{eq: conditions effectives classiques intro} (with our notations), but with a coefficient that is locally perturbed near the corner at scale $\varepsilon $. The perturbation is chosen by solving an optimisation problem that depends on the first near field. Unfortunately, computing the solution of this model requires finely meshing the domain near the corner. In their second effective model, the condition $\mu _{0} \partial _{n} u=-\frac{\mu _{1} }{\varepsilon } u +\varepsilon \beta  \partial _{x} ^2  u$ is imposed on $\Gamma $ with a coefficient $\beta $ that is again chosen through an optimisation problem. So this is a Ventcel condition, but different than \cref{eq: conditions effectives classiques intro}. Theoretical justifications were given to the first model but not the second one. These two models bring limited improvement for the $L^2 $~convergence compared to the naive models. Moreover their lack of physical meaning make them difficult to generalize or improve.
\end{textesansboite}

\section{Summary of the thesis}

\begin{textesansboite}
Let us present the main results and contributions of each chapter. The thesis can be divided into three main parts:
\begin{itemize}
	\item in \crefrange{chap: SF}{chap: construction du DA}, we build an asymptotic expansion at any order of $u_{\varepsilon }$, the solution of \cref{eq: pb etudie},
	\item in \crefrange{chap: complements sur le DA}{chap: validation numerique du DA}, we give additional properties of the expansion and we verify numerically the validity of the expansion,
	\item in \cref{chap: ME}, we deduce high order effective models, of order~3 and 5 respectively.
\end{itemize}
The results of \cref{chap: raccord algebrique,chap: construction du DA} are published in \cite{Bau26}. Several other articles are in preparation about \crefrange{chap: complements sur le DA}{chap: ME}. The results of this thesis were presented in the following conferences: Waves 2024 (Berlin), Aspect 2024 (Metz), ICOSAHOM 2025 (Montréal), CANUM 2026 (Bretagne), Waves 2026 (Montréal). In addition, part of the code written during this thesis has been published as a \Cplusplus{} library, see \cite{OPTEC}.\\

To begin with, the next section introduces the method used for the asymptotic expansion.
\end{textesansboite}

\subsection{Presentation of the matched asymptotic expansion}
\label{sec: ansatz du dvlpmt raccorde}

\begin{figure}[h]
\def\epsFig{0.6} 
\def\widthFig{3.6} 
\def\heightFig{2.7} 
\def\depthFig{1.8} 
\def\ThetaFig{225} 
\def\racInt{0.95} 
\def\racExt{1.7} 
\begin{tikzpicture}[scale=0.9]
	\coordinate (finLambda) at ({-\depthFig*cot(\ThetaFig)}, -\depthFig);
	\coordinate (sommetCouche) at ({-\epsFig*cot(\ThetaFig)}, -\epsFig);
	\draw (-\widthFig-1.8,0) node {\ }; 
	\fill [color=white!85!gray] (0,0) -- (\widthFig,0) -- (\widthFig,-\epsFig) -- (sommetCouche) -- cycle; 
	\fill [color=cyan, opacity=0.3] (\racExt,0) arc (0:\ThetaFig:\racExt) -- (finLambda) -- (-\widthFig,-\depthFig) -- (-\widthFig,\heightFig) -- (\widthFig,\heightFig) -- (\widthFig,0)  -- cycle; 
	\fill [color=green, opacity=0.3] (\racExt,0) rectangle (\widthFig,-\epsFig); 
	\fill [color=red, opacity=0.3] (\racInt,0) arc (0:\ThetaFig:\racInt) -- (sommetCouche) -- (\racInt,-\epsFig) -- cycle; 
	\foreach \i in {5,15,...,95}{
		\def\Rint{(\racInt*(100-\i+5)+\racExt*(\i-5))/100}
		\def\Rext{(\racInt*(100-\i-5)+\racExt*(\i+5))/100}
		\fill [color=cyan!\i !red, opacity=0.3] ({\Rint},0) -- ({\Rext},0) arc (0:\ThetaFig:{\Rext}) -- ({\Rint*cos(\ThetaFig)},{\Rint*sin(\ThetaFig)}) arc (\ThetaFig:0:{\Rint}) -- cycle; 
		\fill [color=green!\i !red, opacity=0.3] ({\Rint},0) rectangle ({\Rext},-\epsFig); 
	}
	\draw [thick] (finLambda) -- (sommetCouche) -- (\widthFig,-\epsFig);
	\draw (sommetCouche) -- (0,0) -- (\widthFig,0);
	\draw [thick,dashed] (\racInt,-\epsFig) -- (\racInt,0) arc (0:\ThetaFig:\racInt);
	\draw [thick,dashed] (\racExt,-\epsFig) -- (\racExt,0) arc (0:\ThetaFig:\racExt);
	\draw [thick,dashed] (\racExt,0) -- (\widthFig,0);
	\draw [>=stealth,->] (-\widthFig-0.4,0) -- (\widthFig+0.4,0) node [right] {$x$};
	\draw [>=stealth,->] (0,-\depthFig) -- (0,{\heightFig+0.3}) node [right] {$y$};
	\draw [decorate, decoration={brace,amplitude=10pt}, xshift=4pt] (\widthFig+0.9,\heightFig) -- (\widthFig+0.9,-\epsFig) node [midway,right,xshift=10pt] {$\Omega_\varepsilon$};
	\draw [color=black!50!blue] (-\widthFig*0.7,\heightFig*0.63) node {far zone};
	\draw [color=black!60!green,<-,>=latex] (\widthFig*0.77,-\epsFig*0.6) to[bend right=20] ++(-0.05,-0.6) node[below]{layer zone};
	\draw [color=black!50!purple] (0,\racInt*0.4) node {\parbox{20mm}{\centering corner\\ zone}};
	\draw [<->,>=latex] (\widthFig*0.9,0) -- node[midway,right]{$\varepsilon$} (\widthFig*0.9,-\epsFig);
	\draw [<->,>=latex] (0,{-\epsFig-0.1}) -- node[midway,below]{$\sqrt\varepsilon$} ({\racInt-0.1},{-\epsFig-0.1});
	\draw [<->,>=latex] (\racInt,{-\epsFig-0.1}) -- node[midway,below]{$\sqrt\varepsilon$} (\racExt,{-\epsFig-0.1});
	\draw [<-,>=latex] ({(\racInt+\racExt)/2*cos(70)},{(\racInt+\racExt)/2*sin(70)}) to[bend left=20] ++(0.3,0.6) node[right]{\!matching zone};
\end{tikzpicture}\par
\caption{Zones of the matched asymptotic expansion}
\label{fig: zones du dvlpmt raccorde}
\end{figure}
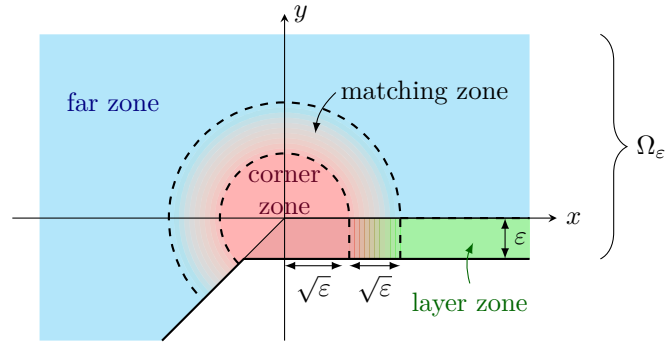

\begin{textesansboite}
In this work, we choose the method of \emph{matched asymptotic expansions} rather than multiscale expansions. Indeed, it provides a more intrinsic construction of the far fields (independent of a truncation function, see \cite{TorViaDau06}), which seemed to us more promising to derive effective models. To take into account the different behaviors of the solution in the layer, near the corner and far from the corner and the layer, we divide $\Omega _{\varepsilon }$ into three zones, illustrated in \cref{fig: zones du dvlpmt raccorde}. In each zone we postulate an asymptotic expansion in powers of $\varepsilon $ and $\ln \varepsilon $, called ``ansatz''.
\end{textesansboite}

\begin{textesansboite}
We recall that the geometry is defined in \cref{sec: motivations et pb etudie}. Let us assume that for~any $(p,\ell )\in \NN\times \mathbb{N}$ there exist three functions independent of $\varepsilon $ – namely $u_{p,\ell }$ defined on $\Omega $ called ``far field'', $U_{p,\ell }$ defined on $\couche$ called ``layer field'' and $S_{p,\ell }$ defined on the rescaled domain $\Omega _{1}$ called ``corner field'' – such that $u_{\varepsilon }$ is formally written as:
\begin{itemize}
	\item $u_{\varepsilon } (x,y) =\displaystyle \sum _{p\in \NN} \sum _{\ell \in \mathbb{N}} \varepsilon ^{p} \ln ^{\ell } \!\varepsilon \;u_{p,\ell } (x,y)$ when $r:=\|(x,y)\|\geqslant \sqrt{\varepsilon }$ and $(x,y)\not\in \mathbb{R}_+\times (-\varepsilon ,0)$ (far zone),
	\item $u_{\varepsilon } (x,y) =\displaystyle \sum _{p\in \NN} \sum _{\ell \in \mathbb{N}} \varepsilon ^{p} \ln ^{\ell } \!\varepsilon \;U_{p,\ell } \Big(x, \frac{y}{\varepsilon } \Big)$ when $x\geqslant \sqrt{\varepsilon }$ and $y\in (-\varepsilon ,0)$ (layer zone),
	\item $u_{\varepsilon } (x,y) =\displaystyle \sum _{p\in \NN} \sum _{\ell \in \mathbb{N}} \varepsilon ^{p} \ln ^{\ell } \!\varepsilon \;S_{p,\ell } \Big( \frac{x}{\varepsilon } , \frac{y}{\varepsilon } \Big)$ when $r\leqslant 2\sqrt{\varepsilon }$ (corner zone).
\end{itemize}\fbi
\end{textesansboite}
\nomenclature[30a]{$u_{p,\ell }$}{Far field}
\nomenclature[30b]{$U_{p,\ell }$}{Layer field}
\nomenclature[28]{$S_{p,\ell }$}{Corner field}

\begin{textesansboite}
The presence of integer powers of $\varepsilon $ is entirely classical in asymptotic analysis. Integer powers of $\varepsilon ^{\pi /\Theta }$ and $\ln \varepsilon $ can be found in other asymptotic expansions involving corners, see \cite{CalCosDauVia06, DelSchSem16, SemDelSch18} and \cref{sec: etat de l'art}. In addition, we will see in \cref{borne sur l} that we have:
\[\forall p\in \NN,\; \exists \ell _{p} \in \mathbb{N},\; \forall \ell >\ell _{p} ,\qquad  u_{p,\ell } = 0 \text{ and } U_{p,\ell } = 0 \text{ and } S_{p,\ell } = 0.\]\\

Note that we treat the layer in the same way as \cref{eq: ansatz couche infinie homogene} in \cite{BenLem96,EngNed93}, which is possible because the layer is homogeneous. It allows us to focus on the matching asymptotic expansion around the corner. The corner fields $S_{p,\ell }$ are, in our context, the ``near fields'' of the matched expansion.\\

\Cref{chap: raccord algebrique,chap: construction du DA} show how to build the fields $u_{p,\ell }$, $U_{p,\ell }$ and $S_{p,\ell }$. The construction is inductive w.r.t.\ $p$. For each $p$, all the fields $u_{p,\ell }$, $U_{p,\ell }$ and $S_{p,\ell }$ are built at once, but there~is only a finite number of these fields that do~not vanish.\\

The first and easiest step of the construction is to obtain the volume equations and boundary conditions for the various fields. Injecting the above ansatz in the Helmholtz equation \cref{eq: pb etudie}, using that $\partial _{x} ^2  \big[\varphi ( \frac{x}{\varepsilon } )\big]=\varepsilon ^{-2} [\partial _{X} ^2 \varphi ](X = \frac{x}{\varepsilon } )$ and $\partial _{y} ^2  \big[\varphi ( \frac{y}{\varepsilon } )\big]=\varepsilon ^{-2} [\partial _{Y} ^2 \varphi ](Y = \frac{y}{\varepsilon } )$ for~any function $\varphi $, and formally identifying the powers of $\varepsilon $ and $\ln \varepsilon $, one can easily derive the following equations.
\begin{equation}{\label{eq: eq vol champs lointains}}
	\left\{\begin{array}{r@{\;=\;}l@{\;}l}
	\mu _{0} \Delta  u_{p,\ell } +\omega ^2 \rho _{0} u_{p,\ell } & f\,\delta _{p,0} \,\delta _{\ell ,0} & \text{ in }  \Omega  \\[0.5mm]
	u_{p,\ell } & U_{p,\ell } & \text{ on }  \Gamma \\[0.5mm]
	u_{p,\ell } & 0 & \text{ on }  \bordgauche
	\end{array}\right. 
\end{equation}
\begin{equation}{\label{eq: eq vol champs de couche}}
	\left\{\begin{array}{r@{\;=\;}l@{\;}l}
	\mu _{1} \partial _{Y} ^2  U_{p,\ell } & -(\mu _{1} \partial _{x} ^2 +\omega ^2 \rho _{1} )U_{p-2,\ell } & \text{ in }  \couche \\[0.5mm]
	\mu _{1} \partial _{Y} U_{p,\ell } & \mu _{0} \,\partial _{y} u_{p-1,\ell } & \text{ on }  \Gamma  \\[0.5mm]
	U_{p,\ell } & 0 & \text{ on }  \borddroit
	\end{array}\right.
\end{equation}
\begin{equation}{\label{eq: eq vol champs de coin}}
	\left\{\begin{array}{r@{\;=\;}ll}
	\operatorname{div}(\mu  \nabla  S_{p,\ell } ) & -\omega ^2 \rho  \,S_{p-2,\ell } & \text{ in }  \Omega _{1} \\[0.5mm]
	S_{p,\ell } & 0 & \text{ on }  \partial  \Omega _{1}
	\end{array}\right.
\end{equation}
Here we denote by convention $u_{p,\ell } =0$, $U_{p,\ell } =0$ and $S_{p,\ell } =0$ for~any $(p,\ell )\in (\mathbb{R} \setminus \NN) \times \mathbb{N}$, and $\delta _{i,j} :=1$ if $i=j$ and 0 if not.
\end{textesansboite}

\begin{remarques}
\begin{itemize}
	\item The condition $\mu _{1} \partial _{Y} U_{p,\ell } = \mu _{0} \partial _{y} u_{p-1,\ell }$ is included in the problem \cref{eq: eq vol champs de couche} satisfied by $U_{p,\ell }$ whereas $u_{p,\ell } =U_{p,\ell }$ is included in the problem \cref{eq: eq vol champs lointains} satisfied by $u_{p,\ell }$ so~that the construction is inductive w.r.t.\ $p$: $u_{p-1,\ell }$ allows to build $U_{p,\ell }$, which allows to build $u_{p,\ell }$.
	\item The layer field $U_{p,\ell }$ satisfies a 1D problem w.r.t.\ $Y$, the variable $x$ playing the role of a parameter.\\
\end{itemize}\fbi
\end{remarques}

\begin{textesansboite}
\Crefrange{eq: eq vol champs lointains}{eq: eq vol champs de coin} would be sufficient to uniquely define the fields, if they were in their natural variational spaces (e.g.\ $H^{1} (\Omega )$ for $u_{p,\ell }$). But we need to take into account a \emph{matching assumption}: the far and corner fields must coincide asymptotically in the intersection of the far and corner zones, and similarly for the layer and corner fields. These intersections form the matching zone (see \cref{fig: zones du dvlpmt raccorde}). Given that $\varepsilon \rightarrow 0$ and $\sqrt{\varepsilon } /\varepsilon  \rightarrow \infty $, this zone tends to $(0,0)$ w.r.t.\ the far and layer fields, but it tends to infinity w.r.t.\ the corner fields. Thus, the matching assumption formally links the asymptotic behavior of the far and layer fields at the corner to the one of the corner fields at infinity:
\begin{equation}{\label{eq: hypothese de raccord}}
\left\{\begin{array}{r@{}ll}
	\displaystyle \sum _{p,\ell } \varepsilon ^{p} \ln ^{\ell } \!\varepsilon \;u_{p,\ell } (x,y)\. &{}\approx \.\displaystyle \sum _{p,\ell } \varepsilon ^{p} \ln ^{\ell } \!\varepsilon \;S_{p,\ell } \Big( \frac{x}{\varepsilon } , \frac{y}{\varepsilon } \Big) &\text{ in }\Omega  \text{ when }  r\rightarrow 0 \text{ and } \frac{r}{\varepsilon } \rightarrow \infty  \vrule width 0pt height 3.5ex \vspace{2pt}\\
	\displaystyle \sum _{p,\ell } \varepsilon ^{p} \ln ^{\ell } \!\varepsilon \;U_{p,\ell } \Big(x, \frac{y}{\varepsilon } \Big) &{}\approx \.\displaystyle \sum _{p,\ell } \varepsilon ^{p} \ln ^{\ell } \!\varepsilon \;S_{p,\ell } \Big( \frac{x}{\varepsilon } , \frac{y}{\varepsilon } \Big) &\text{ in }\couche \text{ when }  x\rightarrow 0 \text{ and } \frac{x}{\varepsilon } \rightarrow \infty .
\end{array}\right.
\end{equation}
The meaning of these ``$\approx $'' will be made more rigorous later (see \cref{def: raccord algebrique}).\\

We will see that each far field $u_{p,\ell }$ has an asymptotic expansion w.r.t.\ $r$ at the corner that is roughly speaking a sum of powers of $r$, some of which are positive (like in a Taylor expansion e.g.). The matching assumption implies that these positive powers of $r$ also appear in the asymptotic expansions of the corner fields at infinity, where they are divergent. We call them singularities for the corner fields. Conversely, the asymptotic expansions of the corner fields at infinity contain negative powers of $r$ corresponding to the decay of the variational part, and these powers must be found in the far fields, which corresponds to singularities at the corner. Thus the fields cannot be searched for in their natural variational spaces. To derive problems that define uniquely the fields, one must specify their singular part, i.e.\ their non-variational part.\\

The fields' singular parts are chosen so that the matching assumption \cref{eq: hypothese de raccord} is satisfied. The singular parts of any field $u_{p,\ell }$ or $U_{p,\ell }$ depends on the behaviors at infinity of the previous corner fields, i.e.\ the $S_{p',\ell '}$ with $p'<p$. Conversely, the singular part of any $S_{p,\ell }$ depends on the behaviors at the corner of the previous far and layer fields. Like in previous works \cite{BenMakTor11,DelSchSem16,SemDelSch18}, choosing these singular parts can be reduced to fixing properly some coefficients involved in the singular behaviors, according to some relations that we call \emph{matching conditions}. Deriving these conditions is quite straightforward for infinite periodic layers, but is much more difficult in corners, and particularly in our case due to the exceptional complexity of the corner singularities at any order. This led us to introduce new and more efficient tools to handle the matching process without being tied to the intricacy of the singularities. This work spans \cref{chap: SF,chap: raccord algebrique} and is presented in the next section.
\end{textesansboite}

\subsection{\texorpdfstring{\Cref{chap: SF,chap: raccord algebrique}}{Chapters 3 and 4}: the matching conditions}
\label{sec: resume raccord en intro}

\begin{textesansboite}
\phantomsection\label{page: definition de Pi}
To perform the matching of the corner fields with the far and the layer fields at the same time, we merge the latter two into a single field denoted $\bfu_{p,\ell }$ and called ``far-and-layer field''. It~is defined on a new domain $\Pi $, defined as follows.
\begin{itemize}
	\item If $\Theta \leqslant  \frac{3\pi }{2}$, then $\Pi  := \Omega \sqcup \Gamma \sqcup \couche$ (disjoint union) and it~is an open subset of $\mathbb{R}^2 $.
	\item If $\Theta > \frac{3\pi }{2}$, then $\Omega $ and $\couche$ intersect as subsets of $\mathbb{R}^2 $, so the previous definition is~not valid anymore (see \cref{fig:Pi}). Thus, we define $\Pi $ as the disjoint gluing of $\Omega $ and $\couche$ on $\Gamma $ (which is a flat Riemannian manifold).
\end{itemize}\fbi
\end{textesansboite}
\nomenclature[51]{$\Pi $}{$\Omega \sqcup \Gamma \sqcup \couche$ (with an exception if $\Theta > \frac{3\pi }{2}$)}

\begin{figure}[h]
\def\widthFig{2.2} 
\def\heightFig{1.8} 
\def\depthFig{-1.5} 
\def\epsFig{0.8} 
\def\ThetaFig{215} 
\begin{tikzpicture}
	\coordinate (finLambda) at (-\widthFig, {-\widthFig*tan(\ThetaFig-180)});
	\fill [pattern={Lines[angle=45,line width=.2pt]}, pattern color=white!70!gray] (0,0) -- (\widthFig,0) -- (\widthFig,\heightFig) -- (-\widthFig,\heightFig) -- (finLambda) -- cycle;
	\fill [color=white!70!gray] (0,0) rectangle (\widthFig,-\epsFig);
	\draw [thick] (finLambda) -- (0,0) -- (\widthFig,0);
	\draw [thick] (\widthFig,-\epsFig) -- (0,-\epsFig) -- (0,0);
	\draw [>=stealth,->] (-\widthFig-0.3,0) -- (\widthFig+0.3,0);
	\draw [>=stealth,->] (0,{\depthFig-0.2}) -- (0,{\heightFig+0.2});
	\draw (-\widthFig*2/5,\heightFig/2) node {$\Omega$};
	\draw [->,>=latex] (0.35,0) arc (0:\ThetaFig-5:0.35) node[near start,above] {$\Theta$} -- (\ThetaFig:0.36);
	\draw [<-,>=latex] (-\widthFig/2, {-\widthFig/2*tan(\ThetaFig-180)}) to[bend right=20] ++(0.3,-0.4) node[right]{$\bordgauche$};
	\draw (\widthFig/2,-\epsFig/2) node {layer $\couche$};
	\draw [<->,>=latex] (\widthFig*0.85,0) -- node[midway,right]{1} (\widthFig*0.85,-\epsFig);
	\draw [<-,>=latex] (\widthFig/2,0) to[bend left=20] ++(0.15,0.5) node[above]{$\ \ \Gamma$};
	\draw [<-,>=latex] (\widthFig/2,-\epsFig) to[bend right=20] ++(0.15,-0.5) node[below]{$\ \ \borddroit$};
\end{tikzpicture}%
\hspace{10mm}\nolinebreak%
\def\ThetaFig{310}%
\begin{tikzpicture}
	\coordinate (finLambda) at ({\depthFig*cot(\ThetaFig)}, \depthFig);
	\fill [color=white!70!gray] (0,0) rectangle (\widthFig,-\epsFig);
	\fill [pattern={Lines[angle=45,line width=.2pt]}, pattern color=white!70!gray] (0,0) -- (\widthFig,0) -- (\widthFig,\heightFig) -- (-\widthFig,\heightFig) -- (-\widthFig,\depthFig) -- (finLambda) -- cycle;
	\fill [pattern={Hatch[angle=45,line width=1pt]}, pattern color=red!70!white] ({-\epsFig*cot(\ThetaFig)}, -\epsFig)  -- (0,0) -- (0,-\epsFig) -- cycle;
	\draw [thick] (finLambda) -- (0,0) -- (\widthFig,0);
	\draw [thick] (\widthFig,-\epsFig) -- (0,-\epsFig) -- (0,0);
	\draw [>=stealth,->] (-\widthFig-0.3,0) -- (\widthFig+0.3,0);
	\draw [>=stealth,->] (0,{\depthFig-0.2}) -- (0,{\heightFig+0.2});
	\draw (-\widthFig*2/5,\heightFig/2) node {$\Omega$};
	\draw [->,>=latex] (0.35,0) arc (0:\ThetaFig-5:0.35) node[midway,left] {$\Theta$} -- (\ThetaFig:0.36);
	\draw ({\widthFig*0.6},-\epsFig/2) node {layer $\couche$};
	\draw [<-,>=latex] ({\epsFig/3},{-\epsFig*2/3}) to[bend left=20] ++(-0.5,-0.5) node[left]{$\Omega\cap\couche$};
\end{tikzpicture}\par
\longCaption{The domain $\Pi $}{. It~is equal to $\Omega \sqcup \Gamma \sqcup \couche$ when $\Theta \leqslant  \frac{3\pi }{2}$ (open subset of $\mathbb{R}^2 $) and it~is a flat Riemannian manifold when $\Theta > \frac{3\pi }{2}$ ($\Omega \cap \couche \neq  \varnothing $).}
\label{fig:Pi}
\end{figure}

\begin{textesansboite}
For~any $(p,\ell )\in \NN\times \mathbb{N}$, we define the \emph{far-and-layer field} $\bfu_{p,\ell }$ and the generalized radial variable $\bfr$ by
\begin{equation}{\label{eq: def bfu bfr}}
    \bfu_{p,\ell } := \left\{\begin{array}{@{\;}l@{\,}l}
        u_{p,\ell } &\text{ in }\Omega  \\[0.5mm]
        U_{p,\ell } &\text{ in }\couche
    \end{array}\right. \qquad  \text{ and } \qquad  \bfr := \left\{\begin{array}{@{\;}l@{\,}l}
        r &\text{ in }\Omega  \\[0.5mm]
        x &\text{ in }\couche
    \end{array}\right.
\end{equation}
\nomenclature[30b]{$\bfu_{p,\ell }$}{Far-and-layer field}
\nomenclature[27b]{$\bfr$}{Generalized radial variable}

A straightforward reformulation of \cref{eq: eq vol champs lointains} and \cref{eq: eq vol champs de couche} gives for~any $(p,\ell )\in \NN\times \mathbb{N}$: \begin{equation}{\label{eq: eq vol de bfu p l}}
\left\{\begin{array}{r@{\;=\;}ll}
	\mu _{0} \Delta  \bfu_{p,\ell } +\omega ^2 \rho _{0} \bfu_{p,\ell } & f\,\delta _{p,0} \,\delta _{\ell ,0} &\text{ in }\Omega  \\[0.5mm]
	\mu _{1} \partial _{Y} ^2  \bfu_{p,\ell } & -(\mu _{1} \partial _{x} ^2 +\omega ^2 \rho _{1} )\bfu_{p-2,\ell } &\text{ in }\couche \\[0.5mm]
	\mu _{1} \partial _{Y} \bfu_{p,\ell |Y=0^{-}} & \mu _{0} \,\partial _{y} \bfu_{p-1,\ell |y=0^{+}}  &\text{ on }\Gamma  \\[0.5mm]
	\bfu_{p,\ell |y=0^{+}} - \bfu_{p,\ell |Y=0^{-}} & 0 &\text{ on }\Gamma \\[0.5mm]
	\bfu_{p,\ell } & 0 &\text{ on }\bordgauche \cup \borddroit
\end{array}\right.\end{equation}
where by convention $\bfu_{p,\ell } =0$ when $p\in \mathbb{R} \setminus \NN$.\\

Let us give some starting point ideas to dive into this part. The matching assumption links the asymptotic behaviors of $\sum _{p,\ell } \varepsilon ^{p} \ln ^{\ell } \!\varepsilon \, \bfu_{p,\ell }$ when $\bfr\rightarrow 0$ and $\sum _{p,\ell } \varepsilon ^{p} \ln ^{\ell } \!\varepsilon \, S_{p,\ell }$ when $\bfr\rightarrow \infty $. So~we can begin with a look at the asymptotic of $\bfu_{p,\ell }$ when $\bfr\rightarrow 0$. Let focus on $\Omega $ because it~is the most interesting part. First, by \cref{eq: eq vol de bfu p l}, $u_{0,0}$ satisfies $\mu _{0} \Delta u_{0,0} +\omega ^2  \rho _{0} u_{0,0} =0$ in the vicinity of the corner in $\Omega $, with homogeneous Dirichlet condition on $\partial \Omega $. So using separation of variables, it~is easy to show that we have in a vicinity of the corner:
\begin{equation}{\label{eq: decomposition de u 00 en fcts de Bessel}}
	u_{0,0} (r,\theta ) = \sum _{d\in \TN} s_{d} (\bfu_{0,0} ) \,J_{d} (k_{0} r) \.\sin (d\theta ) = \sum _{d\in \TN} s_{d} (\bfu_{0,0} ) \,\Bigg(\sum_{n=0}^\infty  j_{d,n} \.r^{d+2n} \Bigg) \sin (d\theta )
\end{equation}
with $s_{d} (\bfu_{0,0} )$ and $j_{d,n}$ some coefficients in $\mathbb{C}$, $J_{d}$ the Bessel functions of the first kind and $k_{0} := \omega  \sqrt{\rho _{0} /\mu _{0} }$. Then, we will see that \cref{eq: eq vol de bfu p l} implies that $u_{1,0}$ satisfies the Helmholtz equation with the condition $\mu _{1} u_{1,0} = \mu _{0} \.\partial _{y} u_{0,0|y=0^{+}}$ on $\Gamma $. Using \cref{eq: decomposition de u 00 en fcts de Bessel}, one can show that there exist some functions $f_{q,i}$ and coefficients $s_{d} (\bfu_{1,0} )$ s.t.\ we have close to the corner:
\[u_{1,0} = \sum _{d\in \NN} r^{d-1} (f_{d,0} (\theta )+\ln (r)\.f_{d,1} (\theta ))+ \sum _{d\in \TN} s_{d} (\bfu_{1,0} ) \,J_{d} (k_{0} r) \.\sin (d\theta ).\]
The first sum above is a particular solution of the Helmholtz equation that has trace $\frac{\mu _{0} }{\mu _{1}} \,\partial _{y} u_{0|y=0^{+}}$ on $\Gamma $, and the second one is a homogeneous solution.\\

More generally, for any $(p,\ell ) \in \NN\times \mathbb{N}$, the behavior of $\bfu_{p,\ell }$ when $\bfr\rightarrow 0$ (resp.\ $S_{p,\ell }$ when $\bfr\rightarrow \infty $) has the following form:
\[\bfu_{p,\ell } \underset{\bfr\rightarrow 0}\approx \left\{\begin{array}{@{\;}l@{\;}l@{}}
	\displaystyle \sum_{n=0}^\infty  r^{d_{n}} \ln ^{j_{n}} \!r \cdot  g^{0} _{n,1} (\theta ) &\text{ in }\Omega  \\[2pt]
	\displaystyle \sum_{n=0}^\infty  r^{d_{n}} \ln ^{j_{n}} \!r \cdot  g^{0} _{n,2} (Y) &\text{ in }\couche,
\end{array}\right.  \qquad  \text{resp.}  \qquad  
S_{p,\ell } \underset{\bfr\rightarrow \infty }\approx \left\{\begin{array}{@{\;}l@{\;}l@{}}
	\displaystyle \sum_{n=0}^\infty  r^{d'_{n}} \ln ^{j'_{n}} \!r \cdot  g^{\infty } _{n,1} (\theta ) &\text{ in }\Omega  \\[2pt]
	\displaystyle \sum_{n=0}^\infty  r^{d'_{n}} \ln ^{j'_{n}} \!r \cdot  g^{\infty } _{n,2} (Y) &\text{ in }\couche
\end{array}\right.\]
where $(d_{n} )$ is an increasing sequence in $\mathbb{Z}+\TZ$ (resp.\ $(d'_{n} )$ is a decreasing sequence in $\mathbb{Z}+\TZ$)  and $j_{n} ,j'_{n} \in \mathbb{N}$. For now it is not clear at all that the far-and-layer fields and the corner fields can match. This depends on the functions $g^{0} _{n,i}$ and $g^{\infty } _{n,i}$ as, if these functions were too different, the fields would~not be able to match. Therefore, we need to describe the $g^{0} _{n,i}$ and $g^{\infty } _{n,i}$. In the literature, they are usually built inductively, by~solving 1D problems in $\theta $ and $Y$ \cite{CalCosDauVia06, DelSchSem16, SemDelSch18}. However, due to the inhomogeneity of the Helmholtz operator $\Delta +k_{0} ^2 $, the $g_{n,i}$ in our case do~not depend on each other only through linear sequences, but rather through a tree of dependencies (see \cref{rq: utilite des operateurs R}).\\

To find a way to handle this complexity, we looked back at the simpler case of an infinite periodic layer \cite{AmmHe97, AmmHe98, AmmLat98}. In this situation, the near zone is the layer and its vicinity, and the behaviors of the far fields near the layer are simply Taylor expansions w.r.t.\ $y$, the distance to the layer. Similarly, the behaviors at infinity of the near fields are also polynomial w.r.t.\ $y$. Thus the behavior of any field can be decomposed on a very simple basis, given by the powers $y^{n}$ with $n\in \mathbb{N}$, and 
the matching assumption reduces to matching the components of the fields' behaviors on each $y^{n}$.\\

In order to try reproducing the simplicity of this case in our much more complex setting, the first original idea of this thesis is to decompose the fields' behaviors on a well-chosen \emph{basis of functions}. Thus, only the coefficients of the decompositions on the basis need to be matched. Moreover, we provide explicit expressions of those basis functions (\cref{def: espaces A}), of the form:
\begin{equation}{\label{eq: def A dans intro}}
	\left\{\begin{array}{ll}
	\Im[(\alpha z)^{q} \,\overline{\alpha z}^{k} \log (\alpha z)^{j} ] &\text{ in }\Omega ,\; \text{ with }  z:=x+\mathrm{i}y \text{ and } \alpha :=e^{-\mathrm{i}\Theta } \\[2pt]
	x^{d} \,Q(\ln x,Y) &\text{ in }\couche,\; \text{ with }  Q \in \mathbb{C}[T,Y].
	\end{array}\right.
\end{equation}
This expression is not obvious at all and will be justified in \cref{sec: singularites quasi homogenes}. We denote $\mathcal{A}(\Omega )$, resp.\ $\mathcal{A}(\couche)$, the vector space generated by the basis functions on $\Omega $, resp.\ $\couche$. What we call ``corner singularities'' (or simply ``singularities'') are the elements of these spaces. For example, $\mathcal{A}(\Omega )$ contains the corner singularities $r^{d+2n} \sin (d\theta )$ already mentioned in \cref{eq: decomposition de u 00 en fcts de Bessel} because, for any $(d,n)\in \TZ*\times \mathbb{N}$, we have $r^{d+2n} \sin (d\theta )=\pm \Im[r^{d+n} e^{\mathrm{i}(d+n)(\theta -\Theta )} \cdot  r^{n} e^{-\mathrm{i}n(\theta -\Theta )} ] = \pm \Im[(\alpha z)^{d+n} (\alpha z)^{n} ]$ with $z=re^{\mathrm{i}\theta }$. \Cref{eq: def A dans intro} extends a formula ``$\Re[(-z)^{q} \log (-z)^{j} ]$'' that is briefly mentioned in \cite{CosDau96} in a simpler case without layer, with $\Theta =\pi $ and for the Poisson equation.\\

\Cref{eq: def A dans intro} is an appropriate definition of the basis functions because it allows to define some key \emph{operators} on the spaces $\mathcal{A}$, that are denoted $\ROmegaPi$, $\RGammaPi$ and $\RcouchePi$ (\cref{operateurs R}). These operators are used to describe how singularities are linked to each other and constitute the second original idea of this part. They are defined from the solution of a Poisson-like canonical problem in $\Pi $, given by \cref{eq: pb resolu dans A}. They can be composed, added and factorized in all possible ways, which, in contrast to linear sequences of functions, offers sufficient flexibility to describe all singularities. An important result is that all singularities that appear in the fields' behaviors are generated, through these operators, by the fundamental singularities $\phi_{d} :=r^{d} \sin (d \theta )$, $d\in \TZ*$ (which are part of the basis). More precisely, \cref{champs totaux en fct de sigma} shows the following decompositions (written between quotes for now as we have~not yet defined some of the involved objects).

\begin{equation}{\label{eq: champs totaux en fct de sigma en intro}}
	\fboxsep=4pt
	\fbox{$ \begin{array}{@{\,}l@{\,}}
	\text{``\,$\displaystyle \sum _{p,\ell } \powelne \,\bfuDA_{p,\ell } = \sum_{n=0}^\infty  R_{\varepsilon } ^{n} \bigg(\sum _{p,\ell } \powelne\, \sum _{\mathclap{d\in \TZ*}} \sigma _{d} (\bfuDA_{p,\ell } )\,\phi_{d} \bigg) $\,''}  \vspace{\jot}\\
	\text{``\,$\displaystyle \sum _{p,\ell } \powelne \,\SDA_{p,\ell } = \sum_{n=0}^\infty  \widetilde{R}_{\varepsilon } ^{n} \bigg(\sum _{p,\ell } \powelne\, \sum _{\mathclap{d\in \TZ*}} \sigma _{d} (\SDA_{p,\ell } )\,\phi_{d} \bigg) $\,''}  
	\end{array} $}
\end{equation}
Here $\bfuDA_{p,\ell }$ (resp.\ $\SDA_{p,\ell }$) is an object that describes the behavior when $\bfr\rightarrow 0$ of $\bfu_{p,\ell }$ (resp.\ when $\bfr\rightarrow \infty $ of $S_{p,\ell }$) and is explained below. The numbers $\sigma _{d} (\bfuDA_{p,\ell } )$ and $\sigma _{d} (\SDA_{p,\ell } )$ denote the coefficients of the fields' behaviors on the basis functions $\phi_{d}$. In addition, $R_{\varepsilon }$ and $\widetilde{R} _{\varepsilon }$ are operators acting on corner singularities that are defined from the operators $\ROmegaPi$, $\RGammaPi$ and $\RcouchePi$, and powers of $\varepsilon $. It follows from \cref{eq: champs totaux en fct de sigma en intro} that, for the far-and-layers fields and corner fields to match, it~is sufficient to fix properly the coefficients $\sigma _{d} (\bfuDA_{p,\ell } )$ and $\sigma _{d} (\SDA_{p,\ell } )$. More precisely only the coefficients involved in the non-variational behaviors of the fields need to be imposed, i.e.\ the $\sigma _{d} (\bfuDA_{p,\ell } )$ with $d<0$ and the $\sigma _{d} (\SDA_{p,\ell } )$ with $d>0$.\\

The third main idea is to use \emph{formal series} to focus on the calculative nature of the matching, without having to deal with sum truncations and asymptotic remainders like $o(r^{d} )$. Formal series are usually used in asymptotic analysis in an intuitive way like we did in \cref{sec: ansatz du dvlpmt raccorde}, or simply as a notation (see \cite{Won01}). However, here we provide rigorous foundations to formal series in order to use them as a proof tool and to justify non trivial operations on them, such as applying operators, and even formal series of operators like in \cref{eq: champs totaux en fct de sigma en intro}. For instance, for~any $(p,\ell )$, the object $\bfuDA_{p,\ell }$ mentioned above is a formal series that describes the asymptotic expansion of $\bfu_{p,\ell }$ when $\bfr\rightarrow 0$. It is roughly speaking defined as the sequence $(v_{n} )_{n\in \mathbb{N}}$ of terms of this expansion. This means that, for any $n\in \mathbb{N}$, $v_{n}$ behaves like $\bfr^{d_{n}}$ for some increasing sequence $(d_{n} )$, and we have $\bfu_{p,\ell } =\sum_{k=0}^n v_{k} +o(\bfr^{d_{n}} )$. Similarly for $\SDA_{p,\ell }$ w.r.t.\ $S_{p,\ell }$ when $\bfr\rightarrow \infty $. The exact definition of formal series is a bit less straightforward, but more flexible so that it permits summing series and applying operators on them.\\

We introduce several spaces of formal series, in which we show the existence and uniqueness of the solution of some Poisson and Helmholtz problems, with an expression of the solution (\cref{pb de Helmholtz dans A Pi,pb de Poisson dans A Pi,pb de Helmholtz dans A Omega,pb de Helmholtz dans Aemacro,pb de Helmholtz dans Aemicro}). This allows us to express the fields' behaviors w.r.t.\ $\bfr$ and to prove \cref{eq: champs totaux en fct de sigma en intro}, as a consequence of the equations satisfied by ``$\sum _{p,\ell } \powelne \,\bfu_{p,\ell }$'' and ``$\sum _{p,\ell } \powelne \,S_{p,\ell }$''.\\

The introduction of \cref{sec: C[[X]]} draws a parallel between formal series and distributions. Both are generalizations of classical objects (converging series, resp.\ $L^{1} _{\mathrm{loc}}$ functions) that offer more flexibility on some operations (e.g.\ derivation, resp.\ infinite sums), allowing to solve problems that do~not have a well-defined solution in the classical framework. Let us mention that \cite{Fao02} also used an approach with formal series to study a problem with a thin layer, but with more elementary algebraic tools.\\

These three ideas -- explicit basis singularities, operators on singularities, and formal series -- form a system of powerful calculus tools that are suitable to face the sheer complexity of the matching. It allows us to prove the following theorem, that provides the matching conditions and is the main result of \cref{chap: raccord algebrique}. One of the originalities of our method is that these conditions are obtained at all orders at once, unlike the step-by-step process that is often used for corner problems \cite{CalCosDauVia06,BenMakTor11,SemDelSch18,CosDau96}.
\end{textesansboite}

\begin{boitecarree}{\cref{condition de raccord}:}
\ There exist two families of constants $(\cmacro_{d,d',p,\ell } )$ and $(\cmicro_{d,d',p,\ell } )$ s.t.\ the matching assumption of \cref{eq: hypothese de raccord} is equivalent to the following set of equations called the \emph{matching conditions}:\vspace{-1mm}
\[\!\left\{\begin{array}{ll}
	\! \forall (d,p,\ell )\in (-\TN)\times \NN\times \mathbb{N},
	&\displaystyle \sigma _{d} (\bfuDA_{p,\ell } ) = \! \sum _{p'\in \intvl{0,p+d}} \sum _{\substack{d'\in \TZ*\\ p'-d'\in \intvl{0,p}}} \! \sum _{\ell '=0} ^{\ell } \,\cmacro_{d,d',p-p',\ell -\ell '} \cdot  \sigma _{d'} (\SDA_{p',\ell '} ) \\[-1mm]
	\! \forall (d,p,\ell )\in \TN\times \NN\times \mathbb{N},
	&\displaystyle \sigma _{d} (\SDA_{p,\ell } ) = \! \sum _{p'\in \intvl{0,p-d}} \sum _{\substack{d'\in \TZ*\\ p'+d'\in \intvl{0,p}}} \! \sum _{\ell '=0} ^{\ell } \,\cmicro_{d,d',p-p',\ell -\ell '} \cdot  \sigma _{d'} (\bfuDA_{p',\ell '} )
\end{array}\right.\]
with $\intvl{0,a}:=\{b\in [0,a] \mid b\in \NN \text{ and } a-b\in \mathbb{N}\}$ for~any $a\in \mathbb{R}$.
\end{boitecarree}

\begin{textesansboite}
%
%
We call \emph{matching coefficients} the numbers $\cmacro_{d,d',p,\ell }$ and $\cmicro_{d,d',p,\ell }$. They depend only on $\Theta $ and $(\mu _{0} ,\mu _{1} ,\rho _{0} ,\rho _{1} ,\allowbreak \omega )$ ; in particular they do~not depend on the source term $\fsource$, nor on the local perturbation of the geometry at the corner. In addition, these coefficients can be obtained through computations on corner singularities involving the operators $\ROmegaPi$, $\RGammaPi$ and $\RcouchePi$. Thanks to the knowledge of a basis of the singularities and using that the operators $\ROmegaPi$, $\RGammaPi$ and $\RcouchePi$ have explicit expressions on this basis, we were able to develop a \emph{symbolic computation module} specifically designed to handle these computations (\cref{sec: implementation des espaces A}). This allows to compute the matching coefficients, as well as the corner singularities appearing in the fields' behaviors, both exactly and very quickly (typically in a fraction of a second). This last point is an original contribution compared to the existing literature.\\

Note that the matching conditions in \cref{condition de raccord} have a convolutive structure w.r.t.\ $p$ and $\ell $, which is important later for the effective models. Moreover, the complexity of the matching conditions is partly due to the fact that we consider Helmholtz's equation. For Poisson's equation, we would have double sums instead of triple sums (because the index $d'$ of the second sum can be expressed from the index $p'$).\\

\Cref{chap: raccord algebrique} is one of the most original parts of this thesis. It not only derives the matching conditions, but also sets up tools that are used throughout the following chapters. In addition, our method uses abstract non case-specific calculations that make it suitable for generalizations to other situations. An insight of these extensions is given in \cref{sec: resume du raccord} and future works will explore this topic in greater depth.\\

Our approach uses a type of formal series that we did~not find in the literature and is presented in \cref{chap: SF}. It involves a summation index that can be any real number, in order to treat the non-integer powers of $\varepsilon $ and $r$ that appear in the expansion. 
\Cref{chap: SF} mostly defines notations, such as the formal sum symbol $\fsum$ (a bold $\sum $), and provides calculus rules on which infinite formal sums are allowed, commutation properties, and theorems on operator inversion (used to solve PDEs among formal series). It is a gradual introduction meant to familiarize the reader with manipulations that are extensively used in the following chapters.
\end{textesansboite}

\subsection{\texorpdfstring{\Cref{chap: construction du DA}}{Chapter 5}: construction of the expansion}

\begin{textesansboite}
In \cref{chap: construction du DA}, we build the terms of the asymptotic expansion of $u_{\varepsilon }$, and justify the expansion by error estimates.\\

First, we express the layer fields $U_{p,\ell }$ in terms of the far fields (\cref{expression de Upl}). Then we introduce model problems to build each far field and corner field, and we prove that they possess a unique solution. For instance, for the far fields, this is given by the theorem below. These problems are non-variational as the fields have singular behaviors (at the corner for far fields and at infinity for corner fields).
\end{textesansboite}

\begin{boitecarree}{\cref{cadre fonctionnel des champs lointains}:}
\ Let us denote $\mathcal{A}(\Gamma ) := \{\varphi _{|\Gamma } \mid \varphi \in \mathcal{A}(\Omega )\}$, where we recall that $\mathcal{A}(\Omega )$ is the space of corner singularities on $\Omega $ introduced in \cref{chap: raccord algebrique}. Let $\khiz$ be a smooth cutoff function localized around the corner. Let $f\in (H^{1} (\Omega ))'$, $g\in H^{1/2} _{00} (\Gamma )+\khiz \mathcal{A}(\Gamma )$ and $(s_{d} )\in \mathbb{C}^{-\TN}$ with finite support. The following system has a unique solution in $H^{1} (\Omega )+\khiz \mathcal{A}(\Omega )$.
\begin{equation}{\label{eq: pb modele de CL dans intro}}
	\left\{\begin{array}{r@{\;=\;}ll}
	\mu _{0} \Delta  u +\omega ^2 \rho _{0} u & f &\text{ in }\Omega  \\[0.5mm]
	u & g &\text{ on }\Gamma  \\[0.5mm]
	u & 0 &\text{ on }\bordgauche\\[0.5mm]
	\sigma _{d} (u)& s_{d} & \ \forall d\in -\TN
	\end{array}\right.
\end{equation}\fbeq
\end{boitecarree}

\begin{textesansboite}
The component of $u$ in $\khiz \mathcal{A}(\Omega )$ contains its singular behavior. Besides, \cref{cadre fonctionnel pour les champs de coin} provides a similar result for the corner fields, replacing the Helmholtz equation by Poisson and $d\in -\TN$ by $d\in \TN$.\\

These theorems show that the coefficients $\sigma _{d} (\bfuDA_{p,\ell } )$ with $d\in -\TN$ and $\sigma _{d} (\SDA_{p,\ell } )$ with $d\in \TN$ can be imposed (e.g.\ by the last line of \cref{eq: pb modele de CL dans intro}). For the other values of $d$, these coefficients, that appear in the right-hand side of the matching conditions, belong to the variational part of the fields and cannot be fixed. See \cref{tab: role de sigma d pour CL et CP}. The variational behaviors are studied by proving asymptotic expansions w.r.t.\ $r$ at any order (\cref{DA des champs lointains concatenes,DA des champs de coin}). This uses the Kondrat'ev theory, involving weighted spaces and the Laplace transform w.r.t.\ $\ln r$ (i.e.\ the Mellin transform w.r.t.\ $r$). This is classical, except that we express the expansion in terms of the corner singularities of the spaces $\mathcal{A}$.
\end{textesansboite}

\begin{table}[h]
\renewcommand{\arraystretch}{1.4}
\definecolor{tabvaria}{rgb}{0.83, 1, 0.83} 
\definecolor{tabsing}{rgb}{1, 0.85, 0.85} 
\caption{Summary of the variational/non-variational status of $\phi_{d}$ and the role of the coefficients $\sigma _{d}$ for the various fields}
\label{tab: role de sigma d pour CL et CP}
\begin{tblr}{Q[c,m]Q[c,m]Q[c,m]}
\cline{2-3}
& Far fields ($\bfr\rightarrow 0$) & Corner fields ($\bfr\rightarrow \infty $) \\
\hline
$d\in \TN$ 
& \SetCell{tabvaria} {$\phi_{d}$ variational,\\ $\sigma _{d} (\bfuDA_{p,\ell } )$ computed from $u_{p,\ell }$ \footnotesize (chap.\,\ref{chap: validation numerique du DA})}
& \SetCell{tabsing} {$\phi_{d}$ non-variational,\\ $\sigma _{d} (\SDA_{p,\ell } )$ fixed by matching condition} \\
\hline
$d\in -\TN$ 
& \SetCell{tabsing} {$\phi_{d}$ non-variational,\\ $\sigma _{d} (\bfuDA_{p,\ell } )$ fixed by matching condition}
& \SetCell{tabvaria} {$\phi_{d}$ variational,\\ \vphantom{$\bfuDA$}$\sigma _{d} (\SDA_{p,\ell } )$ computed from $S_{p,\ell }$ \footnotesize (chap.\,\ref{chap: validation numerique du DA})} \\
\hline
\end{tblr}
\end{table}

\begin{textesansboite}
Gathering all the equations that the fields must satisfy (i.e.\ \crefrange{eq: eq vol champs lointains}{eq: eq vol champs de coin} and the matching conditions) and using the results of this chapter, we prove the existence and uniqueness of the fields $u_{p,\ell }$, $U_{p,\ell }$ and $S_{p,\ell }$ (\cref{validation construction et DA des champs}). Moreover we provide a way to build directly the far fields, by eliminating the layer and corner fields, given by the following theorem.
\end{textesansboite}

\begin{boitecarree}{\cref{construction directe des upl}:}
\ Let $(T_{n} )_{n\in \mathbb{N}}$ be the sequence of Taylor coefficients of the tangent. $(u_{p,\ell } )_{p,\ell }$ is the unique family of $H^{1} (\Omega )+\khiz \mathcal{A}(\Omega )$ s.t.\ for~any $(p,\ell )\in \NN\times \mathbb{N}$:
\begin{equation}{\label{eq: construction directe des upl en intro}}
	\!\!\left\{\begin{array}{r@{\;=\;}ll@{}}
	\mu _{0} \Delta  u_{p,\ell } +\omega ^2 \rho _{0} u_{p,\ell } & \fsource \,\delta _{p,0} \,\delta _{\ell ,0} &\text{ in }\Omega  \\[0.5mm]
	u_{p,\ell } & 0 &\text{ on }\bordgauche\\[0.5mm]
	u_{p,\ell } &\displaystyle \sum_{n=0}^\infty  \frac{\mu _{0} }{\mu _{1}} \,T_{n} \cdot  (\partial _{x} ^2 +k_{1} ^2 )^{n} \partial _{y} u_{p-1-2n,\ell } &\text{ on }\Gamma  \\[-1.5mm]
	\sigma _{d} (u_{p,\ell } )&\displaystyle \!\! \sum _{\substack{p'\in \NN\\ p-p'\in \NN+ \frac{2\pi }{\Theta } }} \; \sum _{\substack{d'\in \TZ*\\ p'+d'\in \intvl{0,p+d}}} \! \sum _{\ell '=0} ^{\ell } \cfinal_{d,d',p-p',\ell -\ell '} \cdot \sigma _{d'} (\bfuDA_{p',\ell '} ) & \ \forall d\in -\TN
	\end{array}\right.\!\!\!
\end{equation}\fbeq
\end{boitecarree}

\begin{textesansboite}
In \cref{eq: construction directe des upl en intro}, the coefficients $\cfinal$ are defined explicitly from the matching coefficients $\cmacro$ and $\cmicro$ appearing in \cref{condition de raccord}, and from new coefficients $\sigma _{q} (\ccDA_{d,n} )$ with $q\in -\TN$, see \cref{eq: def cfinal}. The functions $\mathcal{S}_{d,n}$ are a basis of the corner fields called \emph{corner profiles} (\cref{def: correcteurs de coin}) that do~not depend on the far fields, and the $\sigma _{q} (\ccDA_{d,n} )$ describe their variational behavior at infinity.\\

Finally, the asymptotic expansion is rigorously justified by the following sharp error estimate.
\end{textesansboite}

\begin{boitecarree}{\cref{DA de SE avec CL et delta fixe}:}
\ For~any $\ptronc \geqslant 0$ and $\delta >0$, we have:
\[\bigg\|u_{\varepsilon } -\sum _{p\in \NN\mkern-1mu,\,p\leqslant \ptronc} \sum _{\ell =0} ^{\ell _{p}} \varepsilon ^{p} \ln ^{\ell } \!\varepsilon \, u_{p,\ell } \bigg\|_{H^{1} (\Omega  \setminus \closedBall(0,\delta ))} = o(\varepsilon ^{\ptronc} ) \qquad  \text{when }  \varepsilon \rightarrow 0.\]\fbeq
\end{boitecarree}

\begin{textesansboite}
This concludes the construction of the expansion, that spans \crefrange{chap: SF}{chap: construction du DA}.
\end{textesansboite}

\subsection{\texorpdfstring{\Cref{chap: complements sur le DA,chap: decompo modale CL}}{Chapter 6 and 7}: additional properties of the expansion}

\begin{textesansboite}
\Cref{chap: complements sur le DA,chap: decompo modale CL} establish miscellaneous additional properties of the fields, especially on the far fields. The most remarkable one is a substantial simplification of the matching conditions that happens in ``most cases'' (\cref{simplifications cas irrationnel,simplifications cas rationnel}). We prove that, if $\frac{\Theta }{\pi }$ is irrational or if it is rational and $p\in \NN$ is smaller than the denominator of $\frac{\Theta }{\pi }$ in irreducible form, we have:
\begin{itemize}
	\item the expansion of $u_{\varepsilon }$ up to order $p$ contains no $\ln \varepsilon $,
	\item the matching conditions become one-to-one relations of the form $\sigma _{d} (\bfuDA_{p,0} )=\sigma _{d} (\SDA_{p+d,0} )$,
	\item the corner condition in \cref{construction directe des upl} becomes
	\begin{equation}{\label{eq: condition de coin CL simplifiee en intro}}
		\sigma _{d} (u_{p,0} )= \sum _{\substack{n\in \mathbb{N}\\ 2n\leqslant p+d}} \sum _{\substack{d'\in \TN\\ d'\leqslant p+d-2n}} \omega ^{2n} \sigma _{d} (\ccDA_{d',n} )\cdot \sigma _{d'} (\bfuDA_{p-d'+d-2n,0} ).
	\end{equation}\fbeq
\end{itemize}
This is extended to a few other relevant special cases in Appendix~\ref{sec: annexe des complements sur le DA}. Thus, at a given order, e.g.\ 3 or 5 like in the effective models, only a small number of angles fall outside these simplifications, namely $\frac{2\pi }{3}$, $\pi $, $\frac{4\pi }{3}$ and $\frac{5\pi }{3}$ at order~5. We emphasize that the approach of \cref{chap: raccord algebrique} does~not rely on these simplifications, which allows us to treat any angle at any order.\\

To prove those simplifications, we express some of the far fields' singularities with Bessel functions. First, we show the following \emph{formal} decomposition of the far fields' behavior at the corner (\cref{expression uDA avec mcl})
\begin{equation}{\label{eq: expression uDA avec mcl en intro}}
    \fbox{$ \,\text{``$\.\displaystyle u_{p,\ell } \, \underset{r\rightarrow 0}\approx \, \sum _{n=0} ^{\lfloor {p}\rfloor } \.\sum _{d\in \TZ*\!,\,d\geqslant n-p} \sigma _{d} (\bfuDA_{p-n,\ell } ) \,\mcla_{d,n} \.$''} \, $}
\end{equation}
where the $\mcla_{d,n}$'s are formal series of corner singularities that formally satisfy for~any $d\in \TZ*$
\begin{equation}{\label{eq: pb verifie par mcl qd n=0 en intro}}
    \left\{\begin{array}{r@{\;=\;}ll}
    (\Delta +k_{0} ^2 )\mcla_{d,0} & 0 &\text{ in }\Omega  \\[3pt]
    \mcla_{d,0} & 0 &\text{ on }\Gamma  \\[3pt]
    \sigma _{q} (\mcla_{d,0} ) & \delta _{d,q} &\ \forall q\in \TZ*
    \end{array}\right.
\end{equation}
and, for~any $n\in \mathbb{N}^*$,
\begin{equation}{\label{eq: pb verifie par mcl qd n>0 en intro}}
    \left\{\begin{array}{r@{\;=\;}ll}
    (\Delta +k_{0} ^2 )\mcla_{d,n} & 0 &\text{ in }\Omega  \\[2pt]
    \mcla_{d,n} & \sum\limits_{m=0}^\infty  \frac{\mu _{0} }{\mu _{1}} \,T_{m} \cdot  (\partial _{x} ^2 +k_{1} ^2 )^{m} \partial _{y} \mcla_{d,n-1-2m} &\text{ on }\Gamma  \vspace{2pt}\\
    \sigma _{q} (\mcla_{d,n} ) & 0 &\ \forall q\in \TZ*
    \end{array}\right.
\end{equation}
We prove that, in a lot of cases, the series $\mcla_{d,n}$ converges to a finite linear combinations of Bessel functions (\cref{expression sing heterog CL}). For instance, when $d\not\in -\mathbb{N}$, $\mcla_{d,0}$ converges to $J_{d} (k_{0} r) \sin (d\theta )$ times a constant, where $J_{d}$ is a Bessel function of the first kind. This is coherent with \cref{eq: decomposition de u 00 en fcts de Bessel}. Since Bessel functions of the first kind do~not involve $\ln r$, neither do the far fields in most cases, which implies that there~is no $\ln \varepsilon $ and the other simplifications ``$\sigma _{d} (\bfuDA_{p,0} )=\sigma _{d} (\SDA_{p+d,0} )$'' and \cref{eq: condition de coin CL simplifiee en intro}. Note that the simplifications when $\pst \not\in \mathbb{Q}$ were already known for Poisson's equation \cite{CalCosDauVia06}, but they are significantly less obvious to prove for Helmholtz.\\

In \cref{sec: exemples de structures de DA}, we explicitly write down the matching conditions for a few examples of angle, and we plot the behavior w.r.t.\ the expansion's order of a few relevant quantities like the number of far fields. These data were generated by our symbolic computation module mentioned earlier. They show the announced simplifications and the complexity of the other cases. We use the notion of \emph{genericity} to characterize fields or coefficients that do~not vanish except for exceptional values of $\omega $, $\mu _{0}$, $\mu _{1}$, etc. (\cref{def: genericite}).\\

Next, in \cref{chap: decompo modale CL}, we prove that the formal series $\mcla_{d,n}$ always converges, even if not towards Bessel functions (\cref{convergence des mcl}). We denote $\mcll_{d,n}$ the limit. We show that the modal decomposition of the far fields in \cref{eq: expression uDA avec mcl en intro} actually holds in the classical sense, i.e.\ converges.
\end{textesansboite}

\begin{boitecarree}{\cref{decompo modale des CL}:}
\ Let $\Rsource:=\operatorname{dist}(0,\operatorname{supp}\fsource)>0$. Let $p\in \NN$ and $\ell \in \mathbb{N}$. On $\Omega  \cap  B(0,\Rsource)$, the series below converge and the following equality holds:
\[u_{p,\ell } \,=\, \sum _{n=0} ^{\lfloor {p}\rfloor } \sum _{\substack{d\in \TZ*\\ d\geqslant -p+n}} \sigma _{d} (\bfuDA_{p-n,\ell } ) \,\mcll_{d,n} .\]
In addition, for~any $R\in (0,\Rsource)$, the convergence is uniform on $\Omega  \cap  B(0,R)$ and there exists a constant $C_{R} >0$ s.t., for~any $N\in \mathbb{N}$, $r\in (0,R)$ and $\theta \in [0,\Theta ]$, we have
\begin{equation}{\label{eq: vitesse de cvg decompo modale des CL en intro}}
	\Bigg| u_{p,\ell } - \sum _{n=0} ^{\lfloor {p}\rfloor } \sum _{\substack{d\in \TZ*\\ -p+n\leqslant d\leqslant N\pst}} \sigma _{d} (\uDA_{p-n,\ell } ) \,\mcll_{d,n} \Bigg| \leqslant  C_{R} \left( \frac{r}{R} \right)^{(N+1)\pst-\lfloor {p}\rfloor } (1+|{\ln r}|)^{2p+1} 
\end{equation}
and similarly for the gradient.
\end{boitecarree}

\begin{textesansboite}
In addition, we prove that the error estimate \cref{eq: vitesse de cvg decompo modale des CL en intro} is preserved when replacing $\mcll_{d,n}$ by a partial sum of the series $\mcla_{d,n}$ that contains enough terms (\cref{decompo en modes tronques des CL}). This modal decomposition helps to justify the numerical computation of the far fields in \cref{chap: validation numerique du DA}, and is essential for the effective models of \cref{chap: ME}.\\

To prove the results of \cref{chap: decompo modale CL}, in particular \cref{decompo modale des CL} above, we introduce a norm on corner singularities that is equal to the sum of the coefficients on the singularities' basis (\cref{def: norme sur A_d}). Then we define related weighted norms on formal series of singularities, where the weight depends on the homogeneity degree of the singularities (i.e.\ the power of $\bfr$), see \cref{def: norme sur Amacro}. They are such that a formal series has a finite norm if and only if it converges, and the weight tells us about the speed of convergence. We show that the operators on singularities introduced in \cref{chap: raccord algebrique} are continuous for these norms on well-chosen subspaces of singularities (\cref{norme des op R sur series formelles}), which allows us to bound the norms of $\mcla_{d,n}$ (\cref{norme SF de mcl cas rationnel}). This yields that the $\mcla_{d,n}$ converge, bounds their limit, and implies \cref{decompo modale des CL} by sum (in short, because the details are quite technical). To our knowledge, the method and tools used here do~not exist in the literature.
\end{textesansboite}

\subsection{\texorpdfstring{\Cref{chap: validation numerique du DA}}{Chapter 8}: numerical validation}

\begin{textesansboite}
In this chapter, we illustrate numerically the convergence of the asymptotic expansion at various orders and for several angles. This leads us to designing specific numerical methods to compute the fields, that remain relevant for the effective models of \cref{chap: ME}. The computation of the expansion has four parts in general.
\begin{enumerate}
	\item The computation of the corner profiles $\mathcal{S}_{d,k}$, that are defined in $\Omega _{1}$ and satisfy:
	\[\left\{\begin{array}{r@{\;=\;}l@{\;}l}
		\operatorname{div}(\mu  \nabla  \mathcal{S}_{d,0} ) & 0 &\text{ in }\Omega _{1} \\[0.5mm]
		\mathcal{S}_{d,0} & 0 &\text{ on }\partial \Omega _{1} \\[0.5mm]
		\sigma _{q} (\mathcal{S}_{d,0} )& \delta _{q,d} & \ \forall q\in \TN
	\end{array}\right. \quad  \text{ and } \quad  \ \forall k\in \mathbb{N}^*,\; \left\{\begin{array}{r@{\;=\;}l@{\;}l}
		\operatorname{div}(\mu  \nabla  \mathcal{S}_{d,k} ) & -\rho \,\mathcal{S}_{d,k-1} &\text{ in }\Omega _{1} \\[0.5mm]
		\mathcal{S}_{d,k} & 0 &\text{ on }\partial \Omega _{1} \\[0.5mm]
		\sigma _{q} (\mathcal{S}_{d,k} )& 0 & \ \forall q\in \TN
	\end{array}\right.\]
	From the corner profiles, we only use the coefficients $\sigma _{q} (\ccDA_{d,k} )$ with $q\in -\TN$, which describe the variational behavior at infinity of $\mathcal{S}_{d,k}$. These numbers contain information about the local perturbation of the geometry at the corner, whereas we recall that the matching coefficients do~not see it.
	\item The computation of the matching coefficients $\cmicro_{d,d',p,\ell }$ and $\cmacro_{d,d',p,\ell }$. It~is independent from part~1 and done by our symbolic computation module.
	\item The computation of the coefficients $\cfinal_{d,d',p,\ell }$ from the numbers $\sigma _{q} (\ccDA_{d,k} )$ and the matching coefficients $\cmacro$ and $\cmicro$.
	\item Finally, the computation of the far fields $u_{p,\ell }$ inductively w.r.t.\ $p$, according to the equations of \cref{construction directe des upl}. One needs also to compute some of the coefficients $\sigma _{d} (\uDA_{p,\ell } )$ with $d\in \TN$, as they are involved in the construction of other far fields (when $d\in -\TN$, $\sigma _{d} (\uDA_{p,\ell } )$ is already known, see \cref{tab: role de sigma d pour CL et CP}).
\end{enumerate}
Only steps 1 and 4 require heavy computations as they involve solving a PDE with no explicit solution. Moreover one needs to extract asymptotic coefficients of the PDE's solutions, so special numerical methods are needed.
\end{textesansboite}

\begin{textesansboite}
For the corner profiles, we introduce in \cref{sec: calcul numerique des CC} an enriched Galerkin method that involves so called \emph{complements at infinity}. First, we decompose the corner profile as a singular part, made of explicitly known corner singularities, plus an unknown variational part $\mathcal{S}^{V} _{d,k}$, that must be numerically computed. We know that, for a given $\nbcs \in \mathbb{N}$, we have
\[\mathcal{S}^{V} _{d,k} = \sum _{n=1} ^{\nbcs} \sigma _{-n\pst} (\ccDA_{d,k} ) \,\varphi ^{\nbcs} _{n} +O(r^{-(\nbcs+1)\pst} ) \quad \text{when } r\rightarrow \infty .\]
The $\varphi ^{\nbcs} _{n}$ are explicit corner singularities, while the $\sigma _{-n\pst} (\ccDA_{d,k} )$ are the coefficients we are looking for. Then, we define a bounded subdomain $\Omega _{L} :=\Omega \cap \{\bfr<L\}$ and a smooth cutoff function $\chi _{L}$ localized in $\{\bfr>L-1\}$. We look for an approximation of $\mathcal{S}^{V} _{d,k}$ in the Galerkin space
\[\big\{u:\Omega _{1} \rightarrow \mathbb{C} \mid u_{|\Omega _{L}} \in H^{1} _{0} (\Omega _{L} ) \text{ and } u_{|\Omega _{1} \setminus \ounl} =0 \big\} \oplus  \operatorname{Span}(\{ \khiL \varphi ^{\nbcs} _{n} \mid n\in [\![1,\nbcs]\!]\}).\]
The functions $\khiL \varphi ^{\nbcs} _{n}$ are called \emph{complements at infinity}. In practice, the $H^{1} _{0}$-part is treated with finite elements.\\

We show that the approximate solution converges to $\mathcal{S}^{V} _{d,k}$ in $O\big(L^{-(\nbcs+1) \pst + \frac{1}{2} }\big)$ when $L\rightarrow \infty $ (\cref{estimate of uL-SV}). 
Moreover, we prove that, for~any $n\in [\![1,N]\!]$, the coordinate of the approximate solution on $\varphi ^{\nbcs} _{n}$ converges to $\sigma _{-n\pst} (\ccDA_{d,k} )$ in $O\big(L^{(n-\nbcs-1)\pst + \frac{1}{2} }\big)$ (\cref{convergence des CCC}). Hence, this is a high order method, where the larger $N$, the better the convergence rate. In addition, $\sigma _{-n\pst} (\ccDA_{d,k} )$ is better approximated for small values of $n$ than for the large ones.\\

Numerical results are provided in \cref{sec: validation numerique complements singuliers}. They show that the error made on the coefficients can be very small (from roughly $10^{-5}$ to $10^{-10}$ in relative error) even for a domain $\Omega _{L}$ of reasonable size (e.g.\ $L=10$).\\

This method has similarities with the singular function method \cite{StrFix73, Loh98, LubNic94, CiaHe03}, that is used to address the lack of regularity of solutions of elliptic PDEs in corners. The originality here is to take several complements, which requires a different technique to prove the convergence.
\end{textesansboite}

\begin{figure}[h]
\def\widthFig{3.5} 
\def\heightFig{2.9} 
\def\ThetaFig{205} 
\def\RMEFig{1.3} 
\begin{tikzpicture}[scale=0.9]
	\coordinate (finLambda) at (-\widthFig, {-\widthFig*tan(\ThetaFig-180)});
	\coordinate (dbtLambda) at ({\RMEFig*cos(\ThetaFig)}, {\RMEFig*sin(\ThetaFig)});
	\fill [pattern={Lines[angle=45,line width=.2pt]}, pattern color=white!70!gray] (dbtLambda) arc (\ThetaFig:0:\RMEFig) -- (\widthFig,0) -- (\widthFig,\heightFig) -- (-\widthFig,\heightFig) -- (finLambda) -- cycle;
	\draw [thick] (finLambda) -- (dbtLambda) arc (\ThetaFig:0:\RMEFig) -- (\widthFig,0);
	\draw [dashed] (dbtLambda) -- (0,0);
	\draw [>=stealth,->] (-\widthFig-0.4,0) -- (\widthFig+0.4,0) node [right] {$x$};
	\draw [>=stealth,->] (0,{-\widthFig*tan(\ThetaFig-180)}) -- (0,{\heightFig+0.3}) node [right] {$y$};
	\draw (-\widthFig*0.5,\heightFig*0.6) node {$\OME$};
	\draw [->,>=latex] (0.35,0) arc (0:\ThetaFig-5:0.35) node[near start,above] {$\Theta$} -- (\ThetaFig:0.36);
	\draw [<-,>=latex] (-\widthFig*0.6, {-\widthFig*0.6*tan(\ThetaFig-180)}) to[bend right=20] ++(0.4,-0.4) node[right]{$\!\bordgaucheME$};
	\draw [<-,>=latex] (\widthFig*0.65,0) to[bend right=20] ++(0.15,-0.55) node[below]{$\ \ \GME$};
	\draw (\RMEFig,0) node{$\bullet$};
	\draw [<-,>=latex] (\RMEFig+0.05,0.05) to[bend left=10] ++(0.45,0.45) node[right]{\raisebox{1.5ex}{$\!\PME$}};
	\draw [<-,>=latex] ({\RMEFig*cos(55)}, {\RMEFig*sin(55)}) to[bend left=20] ++(0.4,0.45) node[right]{\raisebox{1ex}{$\SME$}};
	\draw [<->,>=latex] (0,-0.2) -- node[midway, below] {$\RME$} (\RMEFig, -0.2);
\end{tikzpicture}\par
\caption{The domains used to compute the far fields and for the effective models}
\label{fig: domaine ME en intro}
\end{figure}
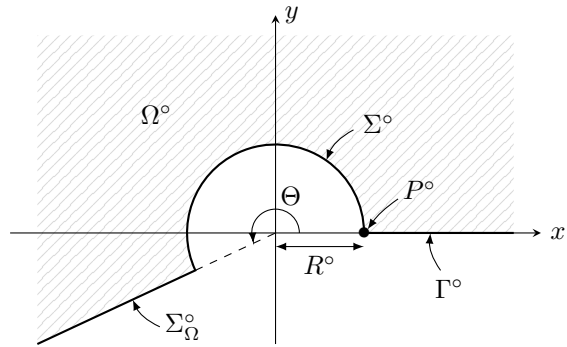

\begin{textesansboite}
Next, we compute the far fields according to the equations \cref{eq: construction directe des upl en intro} of \cref{construction directe des upl} (\cref{sec: numerique des CL}). The first difficulty is to handle the singular behavior of the far field at the corner. We address it by introducing an artificial boundary $\SME$ around the corner (see \cref{fig: domaine ME en intro}) and imposing a non-homogeneous Dirichlet-to-Neumann (DtN) condition of the form
\[\partial _{r} u_{p,\ell } -\opF(u_{p,\ell } )=g \quad \ \text{ on }\SME.\]
Here $\opF$ is the standard DtN operator for a Helmholtz problem in a Dirichlet corner. The right-hand side $g$ is~not standard and depends on the previous far fields and takes into account the condition on $\sigma _{d} (u_{p,\ell } )$ when $d\in -\TN$ of \cref{construction directe des upl}. This method is more convenient than the singular function method in this case, because it addresses three problems at once: the fact that $u_{p,\ell }$ is non variational, the lack of regularity at the corner of its variational part, and the need to compute $\sigma _{d} (\bfuDA_{p,\ell } )$ for $d\in \TN$.\\

An additional difficulty here is to treat the Dirichlet condition $u_{p,\ell } = \sum_{n=0}^\infty  \frac{\mu _{0} }{\mu _{1}} \,T_{n} \cdot  (\partial _{x} ^2 +k_{1} ^2 )^{n} \partial _{y} u_{p-1-2n,\ell }$ on $\Gamma $. It cannot be computed directly as derivatives of finite elements are poorly regular, so they cannot be a Dirichlet condition for a finite element unknown. We handle this by projecting the derivatives on the finite element space and regularizing the derivatives.\\

We provide numerical results that validate the asymptotic expansion up to order~5 (\cref{sec: resultats numeriques CL}). As~far as we know, this is the first time that an asymptotic expansion is computed numerically further than order 2 for a problem with a corner and a thin layer. This is possible thanks to the generality of our method for corner profiles and of the theoretical construction in \cref{chap: raccord algebrique,chap: construction du DA}.\\

The drawback of approaching the exact solution by a sum of far fields is that it may require computing a lot of far fields (e.g.\ 12 at order~5 for $\Theta = \frac{3\pi }{2}$) and the inductive construction of the far fields induces accumulated numerical errors. The effective models of the next chapter do~not present this issue.
\end{textesansboite}

\subsection{\texorpdfstring{\Cref{chap: ME}}{Chapter 9}: effective models}

\begin{textesansboite}
In \cref{chap: ME}, we propose two high order effective models, of respective order~3 and 5 w.r.t.\ $\varepsilon $. They permit computing an approximate solution without finely meshing the layer, nor the vicinity of the corner. By~taking $\omega =0$, our work also provides high order models for the Poisson equation (up to replacing $H^{1}$ by the appropriate variational space), which did~not exist in the literature either. The effective model of order $m\in \{3,5\}$ is designed so that its solution approaches $\sum _{p<m,\ell } \powelne \,u_{p,\ell }$ (which is an approximation of $u_{\varepsilon }$ at order $m$) and takes advantage of what we know on the far fields.\\

The model of order~3, resp.\ 5, is defined on the domain $\OME$ of \cref{fig: domaine ME en intro} according to the following set of equations resp.
\begin{equation}{\label{eq: ME3 en intro}}
	\text{find\ \,}  u\in H^{1} (\OME) \quad  \text{s.t.}  \quad \ \left\{ \begin{array}{r@{\,=\,}ll}
		\mu _{0} \Delta  u+\omega ^2 \rho _{0} u & \fsource &\text{ in }\OME\\[0.5mm]
		\mu _{0} \partial _{n} u + \frac{\mu _{1} }{\varepsilon } u & 0 &\text{ on }\GME\\[0.5mm]
		u & 0 &\text{ on }\bordgaucheME\\[0.5mm]
		\partial _{r} u - \opS3 (u_{|\SME} ) & 0 &\text{ on }\SME 
	\end{array}\right. \qquad 
\end{equation}
\begin{equation}{\label{eq: ME5 en intro}}
	\text{find\ \,}  u\in  H^{1} (\OME)\cap  H^{1} (\GME) \quad \! \text{s.t.}  \quad  \left\{ \begin{array}{r@{\,=\,}ll}
		\mu _{0} \Delta  u+\omega ^2 \rho _{0} u & \fsource &\text{ in }\OME\\[0.5mm]
		\mu _{0} \,\partial _{n} u + \frac{\mu _{1} }{\varepsilon } u - \frac{\varepsilon }{3} (\mu _{1} \partial _{x} ^2 +\omega ^2 \rho _{1} )u & 0 &\text{ on }\GME\\[0.5mm]
		u & 0 &\text{ on }\bordgaucheME\\[0.5mm]
		\partial _{r} u - \opS5 (u_{|\SME} ) & 0 &\text{ on }\SME \\[0.5mm]
		\partial _{x} u - \opP5 (u_{|\SME} ) & 0 & \text{ at }  \PME
	\end{array}\right. \qquad 
\end{equation}
We impose on $\GME$ the same effective boundary conditions as for the infinite layer, see \cref{eq: conditions effectives classiques intro}. Moreover, the novelty of our models is to replace the corner by a DtN condition on the artificial boundary $\SME$. In the model of order~5 in \cref{eq: ME5 en intro}, we also impose a non-classical DtN condition at the point $\PME$ (the tip of $\GME$, see \cref{fig: domaine ME en intro}), which is needed due to the presence of a $\partial _{x} ^2 $ in the boundary condition on $\GME$.\\

The DtN operators $\opS$ and $\opP$ take into account the singular behavior near the corner of the far fields as well as the corner profiles. They are based on an approximate modal decomposition of $\scl := \sum _{p<m,\ell } \powelne \,u_{p,\ell }$ at order $\varepsilon ^{m}$. This decomposition is truncated, keeping only the first $N$ modes, where $N$ is a parameter of the method. $\opS$ and $\opP$ are non-local operators and have a finite rank $N$.\\

A preliminary modal decomposition of $\scl$ is obtained by summing the modal decompositions of the far fields from \cref{chap: decompo modale CL} (\cref{decompo modale SE 1}). A simplified version is
\[\begin{array}{c}
	\displaystyle \scl = \sum _{d\in \TZ,\,d\geqslant -m} \sigma _{d} (\sclDA) \cdot \Phi ^{\varepsilon } _{d} \;+\, O(\varepsilon ^{m} \ln ^{\ell _{m}} \!\varepsilon ) \\
	\displaystyle \text{with } \; \ell _{m} \in \mathbb{N} ,\quad \ \Phi ^{\varepsilon } _{d} := \sum _{n=0} ^{m-1} \varepsilon ^{n} \.\mcll_{d,n} 
	 \quad \  \text{ and } \quad \ \forall d\in \txtTZ*,\ \sigma _{d} (\sclDA) := \! \sum _{p\in \NN \cap [0,m)} \sum _{\ell \in \mathbb{N}} \powelne \cdot  \sigma _{d} (\bfuDA_{p,\ell } ) \in \mathbb{C} .
\end{array}\]
The corner condition on the coefficients $\sigma _{d} (\bfuDA_{p,\ell } )$ in \cref{construction directe des upl} induces constraints on the numbers $\sigma _{d} (\sclDA)$ that must be taken into account. These constraints indirectly involve the coefficients at infinity $\sigma _{q} (\ccDA_{d',k} )$ of the corner profile, so they contain information on the local perturbation of the geometry at the corner. Hence we deduce a second decomposition (\cref{decompo modale SE 2})
\begin{equation}{\label{eq: decompo modale SE 2 en intro}}
	\scl = \sum_{n=0}^\infty  \sigma _{n\pst} (\sclDA) \cdot \Psi ^{\varepsilon } _{n} \;+\, O(\varepsilon ^{m} \ln ^{\ell _{m}} \!\varepsilon )
\end{equation}
where the modes $\Psi ^{\varepsilon } _{n}$ are linear combinations of the $\Phi ^{\varepsilon } _{d}$. The modes $(\Psi ^{\varepsilon } _{1} ,\dots ,\Psi ^{\varepsilon } _{N} )$ are the ones used to build the DtN operators $\opS$ and $\opP5$, or more precisely their Gram-Schmidt orthonormalization, as they are not orthogonal in general.
\\

We prove that, for $\varepsilon $ small enough, the effective models are well-posed and uniformly stable w.r.t.\ $\varepsilon $ (\cref{caractere bien pose des ME}). The strategy is to prove the coercivity of the problem by showing that the main part of the variational formulation (involving the Helmholtz equation and the condition on $\GME$) is coercive, that $\opSref$ is positive for $\varepsilon =0$, and that $\opS-\opSref$ and $\opP5$ have small norms when $\varepsilon $ is small enough. We use that, for $\varepsilon =0$, the modes $\Psi ^{0} _{n}$ are proportional to $J_{n\pst} (k_{0} r) \sin (n\pst \Theta )$, i.e.\ the modes of Helmholtz's equation in a Dirichlet corner without layer. In addition, $\|\opS-\opSref\|$ is controlled thanks to the continuity w.r.t.\ $\varepsilon $ of the $\Psi ^{\varepsilon } _{n}$. This is made more precise thanks to fine estimates of the functions $\mcll_{d,n}$ in \cref{chap: decompo modale CL}. We obtain the following stability condition, where $C$ is a constant independent of $(\varepsilon ,\RME,N)$:
\[\varepsilon  \leqslant  C \. \frac{\RME}{N^{3} \.(1+|{\ln \RME}|)^{3} } \ \ \text{if } m=3 \qquad  \text{ and } \qquad  \varepsilon  \leqslant  C \. \frac{\RME^{4/3} }{N^{3} \.(1+|{\ln \RME}|)^{2} } \ \ \text{if } m=5.\]\\

Moreover, denoting $\sme$ the solution of the effective model of order $m$, we prove the following error estimate (\cref{cvg des ME})
\begin{equation}{\label{eq: cvg des ME en intro}}
	\|\sme-u_{\varepsilon } \|_{H^{1} (\OME)} \leqslant  C \left[ \bigg( \frac{\varepsilon }{\RME} \bigg)^{m} (1+|{\ln \varepsilon }|)^{\ell _{m}} \.N^{3} + \left( \frac{\RME}{R} \right)^{N\pst} \right]
\end{equation}
where $C$ is independent of $(\varepsilon ,\RME,N)$ and $\ell _{m} \in [\![0,4]\!]$ depends on $\Theta $. The second term in \cref{eq: cvg des ME en intro} is due to the error made by the truncation of \cref{eq: decompo modale SE 2 en intro}, which we control thanks to \cref{eq: vitesse de cvg decompo modale des CL en intro}. In particular, taking e.g.\ $N:=\lceil {m \frac{\Theta }{\pi } \log _{2} \varepsilon }\rceil $ and $\RME<\Rsource/2$ fixed, we obtain 
\[\|\sme-u_{\varepsilon } \|_{H^{1} (\OME)} =O(\varepsilon ^{m} \ln ^{\ell _{m}} \!\varepsilon  ) \qquad  \text{when }  \varepsilon \rightarrow 0\]
with $\ell _{m} =0$ for most angles.\\

\Cref{eq: cvg des ME en intro} is proven in two steps. First, we obtain an error in $\varepsilon ^{m-1}$ due to the replacement of the Dirichlet condition by $\mu _{0} \partial _{n} \sme3 = \red{\varepsilon ^{-1} } \mu _{1} \sme3$\footnote{The red font is used in this document to draw attention to some important elements or details that can easily be missed.} (and similarly for $m=5$). Then, like for penalisation methods, the rate of convergence can be improved provided that $\partial _{n} \sme$ is regular enough. We use arguments of elliptic regularity, that are made more complicated due to the corner of the domain $\OME$ at $\PME$.\\

Numerical results are provided and show the efficiency of our models (\cref{sec: numerique ME}). See for example \cref{fig: ME angle pi SE intro,fig: ME angle pi err intro,fig: graphe err ME pi intro} for the case $\Theta =\pi$, an angle that is particularly relevant and presents one of the most complicated asymptotic expansions. \Cref{sec: numerique ME} also provides numerical results for other relevant angles.
\end{textesansboite}

	\tikzset{external/export=true}%
	\tikzsetfigurename{fig-ME_pi_intro-}%
	\begingroup
\def\largeurimages{53mm}
\def\colorbarheight{20mm}
\def\mycolorbar[#1]{\,\colorbar[#1]\hspace*{-7mm}}
\def\espaceapressubcaption{\hspace*{0mm}}
\def\boite#1{\raisebox{-.45\height}{\vbox{\vskip2mm\relax \hbox{#1}}}}%

\begin{figure}[h]
    \boite{%
    \dimen0=\largeurimages\relax
    \dimen0=1.375\dimen0\relax
    \edef\tempo{\the\dimen0}%
    \begin{tikzpicture}[x=\largeurimages, y=\largeurimages]
        \clip (-0.5,{-0.5*0.1/4}) -- (0.5,{-0.5*0.1/4}) -- (0.5,0) arc (0:180:0.5) -- cycle;
        \draw (0,0.335) node{\includegraphics[width=\tempo]{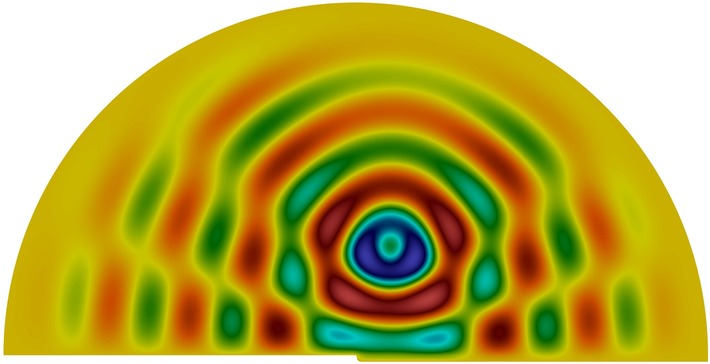}};
    \end{tikzpicture}}
    \mycolorbar[min=-2.1, max=1.4, 1e-2]

    \caption{Real part of the exact solution $u_\peps$ for $\Theta=\pi$ and $\peps=0.1$}
    \label{fig: ME angle pi SE intro}
\end{figure}

\begin{figure}[h]
    \begin{subfigure}[c]{0.5\textwidth}
        \centering
        \boite{\includegraphics[width=\largeurimages]{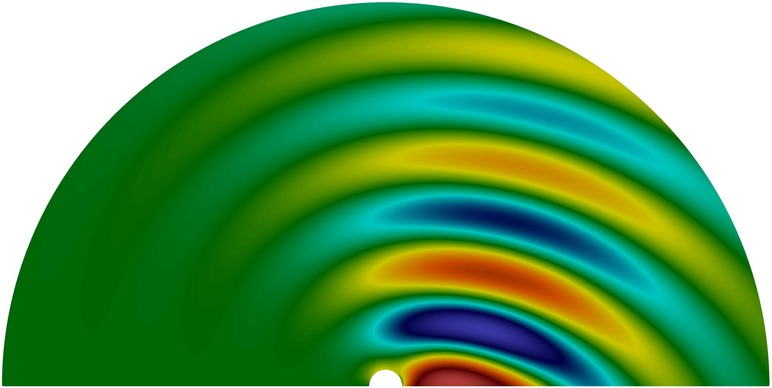}}
        \mycolorbar[min=-3.4, max=4.7, 1e-4]
        \subcaption{\normalsize $m=3$\espaceapressubcaption}
    \end{subfigure}%
    \begin{subfigure}[c]{0.5\textwidth}
        \centering
        \boite{\includegraphics[width=\largeurimages]{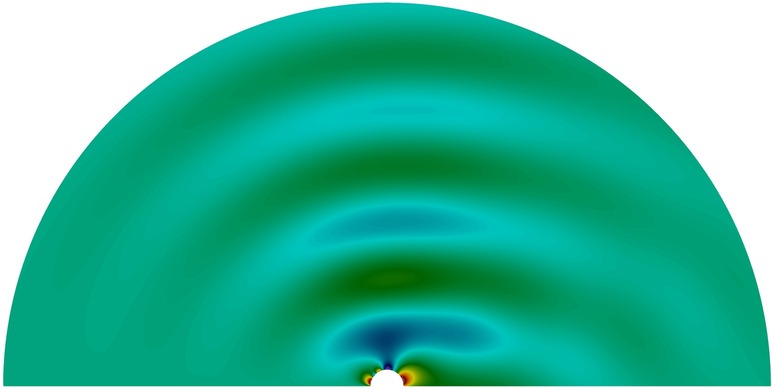}}
        \mycolorbar[min=-5, max=10, 1e-5]
        \subcaption{\normalsize $m=5$\espaceapressubcaption}
    \end{subfigure}

    \caption{Real part of the difference $\protect\sme -u_\peps$ between the solution of the effective models and the exact solution for $\Theta=\pi$ and $\peps=0.1$}
    \label{fig: ME angle pi err intro}
\end{figure}

\endgroup%
	\tikzset{external/export=false}%
	\tikzsetfigurename{main-figure}
	\tikzset{external/export=true}%
	\tikzsetfigurename{fig-erreurs_ME_pi_intro-}%
	\pgfplotsset{
    every axis/.append style={
        thick,
    },
    cycle list name=Dark2_avec_regressions,
    /pgfplots/table/x=epsilon,
    /pgfplots/table/y=nbCC,
    legend cell align=left,
}

\begin{figure}
\begin{tikzpicture}
\begin{axisME}[relative error $\dfrac{\|\sme-u_\peps\|_{H^1 (\OME)}}{\|u_{0,0}\|_{H^1 (\OME)}}$]{err_pi_tot.dat}{}
\courbeMEintro[y=errH1_ME_3, m=3, reg=0:5]
\courbeMEintro[y=errH1_ME_5, m=5, reg=1:5]
\end{axisME}
\end{tikzpicture}
\longCaption{Numerical errors of the effective models for $\Theta=\pi$}{ (the solid lines represent the linear regressions)}
\label{fig: graphe err ME pi intro}
\end{figure}%
	\tikzset{external/export=false}%
	\tikzsetfigurename{main-figure}

\chapter{Introduction en français}
\label{chap: intro fr}
\begin{otherlanguage}{french}
\deuxPointsEst{:}

\section{Motivations et description du problème}
\label{sec: motivations et pb etudie fr}

\begin{textesansboite}
Les problèmes comportant des couches minces apparaissent dans de nombreux domaines allant de l'ingénierie des matériaux composites \cite{FakKha14} à la biologie \cite{CauHadLiNgu16}, en passant par l'élasticité \cite{BelGeyKra16, FriJamMul02}, la mécanique des fluides \cite{RicRuyVil16, AchPirVal98, JagMik01, DalGer17} et l'électrochimie \cite{VerMaiNov10}. Les applications sont particulièrement nombreuses en électromagnétisme ; citons notamment les couches diélectriques minces \cite{HadJia15, Mak08, SenVol95}, les films ferromagnétiques \cite{AmmHalHam00, HadJol01}, les effets de peau \cite{ZutKno05} et l'effet de cage de Faraday \cite{DelHew20}. Toutes ces situations posent des défis numériques, car elles nécessitent un maillage fin de la structure mince, ce qui est très coûteux lorsque son épaisseur est très faible par rapport à la longueur d'onde et à la taille des objets.\\

Pour éviter ce problème, on peut construire une approximation de la solution exacte dont le calcul est moins coûteux, notamment en évitant d'utiliser un maillage fin de la couche. Il existe deux méthodes principales pour y parvenir : soit en construisant un développement asymptotique de la solution exacte lorsque l'épaisseur de la couche tend vers zéro, soit en concevant des modèles effectifs qui remplacent la couche par des conditions aux limites effectives approchées. Cette dernière option est généralement préférable car elle fournit directement une solution approchée, tandis que la première nécessite de construire et de sommer plusieurs termes du développement asymptotique. Cependant, les développements asymptotiques constituent souvent une première étape vers les modèles effectifs. Dans ce travail, nous étudions ces deux approches l'une après l'autre.\\

Comme expliqué dans la section~\ref{sec: etat de l'art fr}, le cas d'objets entièrement recouverts par une couche mince a déjà été étudié de manière approfondie dans la littérature. En revanche, la recherche en est encore à ses débuts en ce qui concerne les revêtements partiels, tant sur les développements asymptotiques que sur les modèles effectifs. Dans ce travail, nous nous intéressons à ces dernières situations. Nous considérons un cas bidimensionnel et nous nous concentrons sur le voisinage du bord de la couche, qui est modélisé par un problème académique comportant une couche rectiligne semi-infinie. De plus, nous supposons que l'objet peut présenter un coin à cet endroit, d'angle arbitraire. Nous considérons également une perturbation locale de la géométrie au niveau du coin, à l'échelle de l'épaisseur de la couche.
\end{textesansboite}

\begin{textesansboite}
Afin de décrire la géométrie, nous introduisons l'angle $\Theta  \in {]0,2\pi [}$ et les domaines $\Omega  := \{(r \cos \theta ,r \sin \theta ) \mid r\in \mathbb{R}_+^*,\theta \in {]0,\Theta [} \}$, $\couche := \mathbb{R}_+^* \times {]{-1},0[}$, $\Gamma  = \mathbb{R}_+^* \times \{0\}$, $\bordgauche=\{(r \cos \Theta ,r \sin \Theta ) \mid r\in \mathbb{R}_+^*\}$ et $\borddroit := \mathbb{R}_+^* \times \{-1\}$. Tous ces ensembles sont représentés dans la figure~\ref{fig:Omega et couche fr}. Puis nous considérons un ouvert $\Omega _{1} \subset \mathbb{R}^2 $ qui coïncide avec $\Omega \cup \Gamma \cup \couche$ en dehors du disque $B(0,\Rcoin)$ pour un certain $\Rcoin \in \mathbb{R}_+^*$. De plus, nous nous donnons deux fonctions $\mu ,\rho  \in L^{\infty } (\Omega _{1} )$ minorées par une constante positive (hypothèse d'ellipticité) et égales aux constantes $\mu _{0}$ et $\rho _{0}$ dans $\Omega  \setminus B(0,\Rcoin)$ et aux constantes $\mu _{1}$ et $\rho _{1}$ dans $\couche \setminus B(0,\Rcoin)$. Voir les figures \ref{fig:Omega1 fr} et~\ref{fig:autres exemples de Omega1 fr} pour différentes configurations.\\

\begin{figure}[h]
\def\widthFig{2.3} 
\def\heightFig{2} 
\def\ThetaFig{205} 
\begin{tikzpicture}[scale=0.9]
	\coordinate (finLambda) at (-\widthFig, {-\widthFig*tan(\ThetaFig-180)});
	\fill [pattern={Lines[angle=45,line width=.2pt]}, pattern color=white!70!gray] (0,0) -- (\widthFig,0) -- (\widthFig,\heightFig) -- (-\widthFig,\heightFig) -- (finLambda) -- cycle;
	\draw [thick] (finLambda) -- (0,0) -- (\widthFig,0);
	\draw [>=stealth,->] (-\widthFig-0.4,0) -- (\widthFig+0.4,0) node [right] {$x$};
	\draw [>=stealth,->] (0,{-\widthFig*tan(\ThetaFig-180)}) -- (0,{\heightFig+0.3}) node [right] {$y$};
	\draw (-\widthFig*2/5,\heightFig/2) node {$\Omega$};
	\draw [->,>=latex] (0.35,0) arc (0:\ThetaFig-5:0.35) node[near start,above] {$\Theta$} -- (\ThetaFig:0.36);
	\draw [<-,>=latex] (-\widthFig/2, {-\widthFig/2*tan(\ThetaFig-180)}) to[bend right=20] ++(0.2,-0.45) node[right]{$\!\bordgauche$};
	\draw [<-,>=latex] (\widthFig/2,0) to[bend right=20] ++(0.15,-0.5) node[below]{$\ \ \Gamma$};
\end{tikzpicture}
\hspace{-8mm}
\def\epsFig{0.8} 
\def\widthFig{3.5}
\begin{tikzpicture}
	\fill [color=white!70!gray] (0,0) rectangle (\widthFig,-\epsFig);
	\draw [thick] (\widthFig,-\epsFig) -- (0,-\epsFig) -- (0,0) -- (\widthFig,0);
	\draw (-2.5,0) node{\ }; 
	\draw [>=stealth,->] (-0.4,0) -- (\widthFig+0.4,0) node [right] {$x$};
	\draw [>=stealth,->] (0,{-\epsFig-0.5}) -- (0,0.5) node [above] {$Y$};
	\draw (\widthFig/2,-\epsFig/2) node {couche $\couche$};
	\draw [<->,>=latex] (\widthFig*0.85,0) -- node[midway,right]{1} (\widthFig*0.85,-\epsFig);
	\draw [<-,>=latex] (\widthFig/2,0) to[bend left=20] ++(0.15,0.5) node[above]{$\ \ \Gamma$};
	\draw [<-,>=latex] (\widthFig/2,-\epsFig) to[bend right=20] ++(0.15,-0.5) node[below]{$\ \ \borddroit$};
\end{tikzpicture}\par
\caption[Les domaines $\Omega$ et $\couche$]{Les domaines $\Omega$ (à gauche) et $\couche$ (à droite)}
\label{fig:Omega et couche fr}
\end{figure}

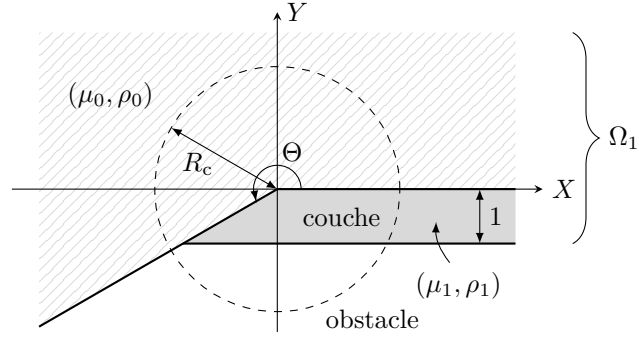
\begin{figure}[h]
\def\epsFig{0.8} 
\def\widthFig{3.5} 
\def\heightFig{2.3} 
\def\ThetaFig{210} 
\def\RCoinFig{(\epsFig/sin(\ThetaFig-180)+0.2)}
\begin{tikzpicture}[scale=0.9]
	\coordinate (finLambda) at (-\widthFig, {-\widthFig*tan(\ThetaFig-180)});
	\coordinate (sommetCouche) at ({-\epsFig*cot(\ThetaFig-180)}, -\epsFig);
	\fill [color=white!70!gray] (0,0) -- (\widthFig,0) -- (\widthFig,-\epsFig) -- (sommetCouche) -- cycle;
	\fill [pattern={Lines[angle=45,line width=.2pt]}, pattern color=white!70!gray] (0,0) -- (\widthFig,0) -- (\widthFig,\heightFig) -- (-\widthFig,\heightFig) -- (finLambda) -- cycle;
	\draw [thick] (finLambda) -- (0,0) -- (\widthFig,0);
	\draw [thick] (sommetCouche) -- (\widthFig,-\epsFig);
	\draw (-\widthFig-1.6,0) node {\ }; 
	\draw [>=stealth,->] (-\widthFig-0.4,0) -- (\widthFig+0.4,0) node [right] {$X$};
	\draw [>=stealth,->] (0,{-\RCoinFig-0.3}) -- (0,{\heightFig+0.3}) node [right] {$Y$};
	\draw (\widthFig*0.27,-\epsFig/2) node {couche};
	\draw (\widthFig*0.4,{-\RCoinFig-0.1}) node {obstacle};
	\draw (-\widthFig*0.7,\heightFig*0.6) node {$(\mu_0,\rho_0)$};
	\draw [->,>=latex] (0.35,0) arc (0:\ThetaFig-5:0.35) node[near start,above] {$\Theta$} -- (\ThetaFig:0.36);
	\draw [decorate, decoration={brace,amplitude=10pt}, xshift=4pt] (\widthFig+0.7,\heightFig) -- (\widthFig+0.7,-\epsFig) node [midway,right,xshift=10pt] {$\Omega_1$};
	\draw [<-,>=latex] (\widthFig*2/3,-\epsFig*0.6) to[bend right=20] ++(0.2,-0.6) node[below]{$\ \ (\mu_1,\rho_1)$};
	\draw [<->,>=latex] (\widthFig*0.85,0) -- node[midway,right]{1} (\widthFig*0.85,-\epsFig);
	\draw [dashed] (0,0) circle ({\RCoinFig});
	\draw [<->,>=latex] (0,0) -- node[near end, below] {$\Rcoin$} ({\RCoinFig*cos(150)}, {\RCoinFig*sin(150)});
\end{tikzpicture}\par
\caption{Le domaine $\Omega_1$ avec un exemple de configuration de $\Omega_1 \cap B(0,\Rcoin)$}
\label{fig:Omega1 fr}
\end{figure}

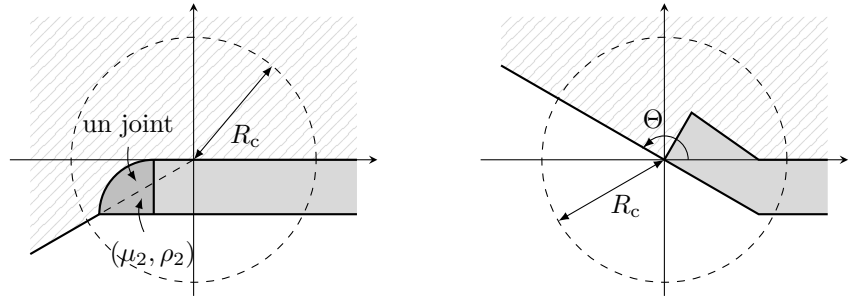
\begin{figure}[h]
\def\epsFig{0.8} 
\def\widthFig{2.4} 
\def\heightFig{2.1} 
\def\ThetaFig{210} 
\def\RCoinFig{(\epsFig/sin(\ThetaFig-180)+0.2)}
\begin{tikzpicture}[scale=0.9]
	\coordinate (finLambda) at (-\widthFig, {-\widthFig*tan(\ThetaFig-180)});
	\coordinate (sommetJoint) at ({-\epsFig*cot(\ThetaFig-180)}, -\epsFig);
	\coordinate (sommetCoucheHaut) at ({-\epsFig*cot(\ThetaFig-180)+\epsFig}, 0);
	\coordinate (sommetCoucheBas) at ({-\epsFig*cot(\ThetaFig-180)+\epsFig}, -\epsFig);
	\fill [pattern={Lines[angle=45,line width=.2pt]}, pattern color=white!70!gray] (0,0) -- (\widthFig,0) -- (\widthFig,\heightFig) -- (-\widthFig,\heightFig) -- (finLambda) -- cycle;
	\fill [color=white!70!gray] (sommetCoucheBas) -- (sommetCoucheHaut) -- (\widthFig,0) -- (\widthFig,-\epsFig) -- cycle; 
	\fill [color=white!50!gray] (sommetJoint) arc (180:90:\epsFig) -- (sommetCoucheBas) -- cycle; 
	\draw [thick] (finLambda) -- (sommetJoint) arc (180:90:\epsFig) -- (\widthFig,0);
	\draw [thin,dashed] (0,0) -- (sommetJoint);
	\draw [thick] (sommetJoint) -- (\widthFig,-\epsFig);
	\draw [thick] (sommetCoucheBas) -- (sommetCoucheHaut);
	\draw [>=stealth,->] (-\widthFig-0.3,0) -- (\widthFig+0.3,0);
	\draw [>=stealth,->] (0,{-\RCoinFig-0.2}) -- (0,{\heightFig+0.2});
	\draw [dashed] (0,0) circle ({\RCoinFig});
	\draw [<->,>=latex] (0,0) -- node[near start, right] {\;$\Rcoin$} ({\RCoinFig*cos(50)}, {\RCoinFig*sin(50)});
	\draw [<-,>=latex] (-\epsFig*1.05,-\epsFig*0.4) to[bend left=20] ++(-0.12,0.51) node[above]{un joint};
	\draw [<-,>=latex] (-\epsFig*0.95,-\epsFig*0.7) to[bend right=20] ++(0.03,-0.5) node[below]{\ \ $(\mu_2 ,\rho_2)$};
\end{tikzpicture}%
\nolinebreak\hspace{10mm}
\def\ThetaFig{150} 
\def\RCoinFig{(\epsFig/sin(180-\ThetaFig)+0.2)}
\begin{tikzpicture}[scale=0.9]
	\coordinate (finLambda) at (-\widthFig, {-\widthFig*tan(\ThetaFig-180)});
	\coordinate (sommetCoucheBas) at ({-\epsFig*cot(\ThetaFig-180)}, -\epsFig);
	\coordinate (sommetCoucheHaut) at (\ThetaFig-90:\epsFig);
	\coordinate (sommetCoucheMilieu) at ({\epsFig*tan(\ThetaFig-90)}, 0);
	\fill [pattern={Lines[angle=45,line width=.2pt]}, pattern color=white!70!gray] (0,0) -- (\widthFig,0) -- (\widthFig,\heightFig) -- (-\widthFig,\heightFig) -- (finLambda) -- cycle;
	\fill [color=white!70!gray] (0,0) -- (sommetCoucheHaut) -- (sommetCoucheMilieu) -- (\widthFig,0) -- (\widthFig,-\epsFig) -- (sommetCoucheBas) -- cycle;
	\draw [thick] (finLambda) -- (sommetCoucheBas) -- (\widthFig,-\epsFig);
	\draw [thick] (0,0) -- (sommetCoucheHaut) -- (sommetCoucheMilieu) -- (\widthFig,0);
	\draw [>=stealth,->] (-\widthFig-0.3,0) -- (\widthFig+0.3,0);
	\draw [>=stealth,->] (0,{-\RCoinFig-0.2}) -- (0,{\heightFig+0.2});
	\draw [->,>=latex] (0.35,0) arc (0:\ThetaFig-5:0.35) node[near end,above] {$\Theta$\;} -- (\ThetaFig:0.36);
	\draw [dashed] (0,0) circle ({\RCoinFig});
	\draw [<->,>=latex] (0,0) -- node[near end, right] {\ \;$\Rcoin$} ({\RCoinFig*cos(210)}, {\RCoinFig*sin(210)});
\end{tikzpicture}\par
\caption{Autres exemples de configuration de $\Omega_1 \cap B(0,\Rcoin)$ pour différentes valeurs de $\Theta$}
\label{fig:autres exemples de Omega1 fr}
\end{figure}

Soit $\varepsilon >0$ l'épaisseur de la couche. Le domaine physique est donné par
\[\Omega _{\varepsilon } := \left\{(x,y)\in \mathbb{R}^2  \ \Big|\ \left( \frac{x}{\varepsilon } , \frac{y}{\varepsilon } \right) \in  \Omega _{1} \right\}.\]
Notez que $\Omega _{\varepsilon }$ coïncide avec $\Omega _{1}$ quand $\varepsilon =1$. Nous introduisons les coefficients mis à l'échelle $\mu _{\varepsilon } :(x,y)\in \Omega _{\varepsilon } \mapsto  \mu ( \frac{x}{\varepsilon } , \frac{y}{\varepsilon } )$ et $\rho _{\varepsilon } :(x,y)\in \Omega _{\varepsilon } \mapsto  \rho ( \frac{x}{\varepsilon } , \frac{y}{\varepsilon } )$ ainsi que les variables mises à l'échelle $X := \frac{x}{\varepsilon }$ et $Y:= \frac{y}{\varepsilon }$. Soient $\omega \in \mathbb{C} \setminus \mathbb{R}$ et $\fsource\in H^{-1} (\Omega )$ un terme source tel que $\operatorname{dist}(\operatorname{supp}(\fsource),\Gamma )>0$. Nous notons $u_{\varepsilon }$ l'unique solution dans $H^{1} _{0} (\Omega _{\varepsilon } )$ de 
\begin{equation}{\label{eq: pb etudie fr}}
\left\{\begin{array}{r@{}ll}
    \operatorname{div}(\mu _{\varepsilon } \nabla  u_{\varepsilon } )+\omega ^2 \rho _{\varepsilon } u_{\varepsilon } &{}= \fsource & \text{ dans }  \Omega _{\varepsilon } \\[1pt]
    u_{\varepsilon } &{}= 0 & \text{ sur }  \partial \Omega _{\varepsilon }
\end{array}\right.
\end{equation}
Puisque $\Im(\omega ) \neq 0$, ce problème est bien posé dans $H^{1} _{0} (\Omega _{\varepsilon } )$ (il suffit d'utiliser le théorème de Lax-Milgram) avec une constante de stabilité indépendante de $\varepsilon $ :
\begin{equation}{\label{eq: uniforme coercivite fr}}
\exists C>0,\;\forall \fsource \in H^{-1} (\Omega ), \;\forall \varepsilon >0, \quad  \|u_{\varepsilon } \|_{H^{1} (\Omega _{\varepsilon } )} \leqslant  C \.\|\fsource\|_{H^{-1} (\Omega )}
\end{equation}
\end{textesansboite}

\begin{paragraphesansboite}{Remarque :}
$\Im(\omega ) \neq 0$ est une hypothèse technique liée au fait que nous considérons un cadre académique comportant une couche semi-infinie. Dans ce cas, $\Im(\omega ) = 0$ nécessiterait de définir une condition de rayonnement spécifique, ce qui est difficile. Dans le cas plus réaliste d'un obstacle borné partiellement recouvert par une couche mince, il est possible d'utiliser une condition de rayonnement classique de Sommerfeld. Notez que l'hypothèse $\Im(\omega ) \neq 0$ n'a que peu d'impact sur l'analyse asymptotique.
\end{paragraphesansboite}

\begin{paragraphesansboite}{Notations :}
Tout au long de ce document, nous utilisons les notations suivantes :
\begin{itemize}
    \item $(x,y)$ désigne les coordonnées cartésiennes et $(r,\theta )$ désigne les coordonnées polaires avec $\theta  \in [0,2\pi [$,
    \item pour tout $r>0$, $B(0,r)$ est le disque de $\mathbb{R}^2 $ de centre $(0,0)$ et de rayon $r$,
    \item $k_{i} := \omega  \sqrt{\rho _{i} /\mu _{i} }$ avec $i=0,1$ est le nombre d'onde dans $\Omega $, resp.\ dans la couche,
    \item $\NN := \mathbb{N}+ \frac{\pi }{\Theta } \mathbb{N}$,
	\item $\Rsource := \operatorname{dist}(0,\operatorname{supp}\fsource)$, qui est strictement positif par hypothèse sur $\fsource$,
	\item $\bfr := r$ dans $\Omega $ et $\bfr:=x$ dans $\couche$,
    \item $\alpha  := e^{-\mathrm{i}\Theta }$.
\end{itemize}
De plus, nous utilisons les conventions suivantes :
\begin{itemize}
	\item les intervalles ouverts sont notés $]a,b[$ dans la partie en français (i.e. ce chapitre) mais $(a,b)$ partout ailleurs, et de même pour les intervalles semi-ouverts,
    \item $\mathbb{N}:=\{0,1,2,3,\dots \}$,
    \item $\mathbb{N}^* :=\mathbb{N} \setminus \{0\}$, $\mathbb{R}^* :=\mathbb{R} \setminus \{0\}$ et de même pour $\mathbb{Z}^*$, $\mathbb{Q}^*$ et $\mathbb{C}^*$,
    \item $\mathbb{R}_+ :=[0,\infty [$, $\mathbb{R}_+^* :={]0,\infty [}$ et de même pour $\mathbb{R}_-$ et $\mathbb{R}_-^*$.
	\item Pour tous $n,m\in \mathbb{N}^*$ et tout ouvert $O\subset \mathbb{R}^{n}$, $H^{m} _{\mathrm{loc}} (O)$ désigne l'ensemble des fonctions qui sont $H^{m}$ sur tout ouvert $O'\subset O$ tel que $\overline{O'}$ est compact dans $O$. De même pour $L^2 _{\mathrm{loc}}$.
\end{itemize}
Si besoin, une liste plus complète des notations se trouve à la fin de ce document page~\pageref{nomenclature}.
\end{paragraphesansboite}

\section{\'Etat de l'art}
\label{sec: etat de l'art fr}

	\tikzset{external/export=true}%
	\tikzsetfigurename{fig-couche_infinie_fr-}%
	\begin{figure}[h]

\begin{subfigure}[c]{0.5\textwidth}
\centering
\def\epsFig{0.7} 
\def\widthFig{2.9} 
\def\heightFig{2.8} 
\begin{tikzpicture}[scale=0.9]
	\fill [color=white!70!gray] (-\widthFig,-\epsFig) rectangle (\widthFig,0);
	\fill [pattern={Lines[angle=45,line width=.2pt]}, pattern color=white!70!gray] (-\widthFig,0) rectangle (\widthFig,\heightFig);
	\draw [thick] (-\widthFig,0) -- (\widthFig,0);
	\draw [thick] (-\widthFig,-\epsFig) -- (\widthFig,-\epsFig);
	\draw (0,\heightFig/2) node {\parbox{30mm}{\centering domaine extérieur\\ de propagation}};
	\draw (0,-\epsFig/2) node {couche};
	\draw (0,-\epsFig-0.5) node {obstacle};
	\draw [<->,>=latex] (\widthFig*0.75,0) -- node[midway,right]{$\peps$} (\widthFig*0.75,-\epsFig);
\end{tikzpicture}
\end{subfigure}%
\begin{subfigure}[c]{0.5\textwidth}
\centering
\def\epsFig{0.8} 
\def\widthFig{3.5} 
\def\heightFig{2.55} 
\def\cycleExt{(-0.7,4.5) (-0.3,3) (1,2) (3,2) (5.2,1.7) (7,2.8) (7.8,5) (7,7) (5,7.5) (3,6.8) (1,6.9) (-0.3,6)}
\def\cycleInt{(-1.15,4.5) (-0.65,3) (0.7,2.15) (3,2.15) (5.2,1.85) (7.05,2.8) (8,5) (7.25,7.15) (5.1,7.6) (3,6.8) (0.6,6.95) (-0.8,6)}
\def\tension{0.7}
\begin{tikzpicture}[scale=0.9]
	\fill [pattern={Lines[angle=45,line width=.2pt]}, pattern color=white!70!gray] (-\widthFig,-\heightFig) rectangle (\widthFig,\heightFig);
	\draw (0,-\heightFig+0.4) node {domaine extérieur de propagation};
	\begin{scope}[shift={(-2.1,-2.45)}, scale=0.6]
		\fill [color=white!70!gray] plot [smooth cycle, tension=\tension] coordinates {\cycleExt} ;
		\draw [thick] plot [smooth cycle, tension=\tension] coordinates {\cycleExt} ;
		\begin{scope}[shift={(0.8,0.8)}, scale=0.8]
			\fill [color=white] plot [smooth cycle, tension=\tension] coordinates {\cycleInt} ;
			\draw [thick] plot [smooth cycle, tension=\tension] coordinates {\cycleInt} ;
		\end{scope}
		\draw (3.5,4.5) node {obstacle};
		\draw [<-,>=latex] (3.5,2.12) to[bend left=20] ++(0.4,0.85) node[right]{\raisebox{2mm}{\!couche}};
		\draw [<->,>=latex] (5.3,6.9) -- node[midway,right]{$\peps$} (5.3,7.5);
	\end{scope}
\end{tikzpicture}
\end{subfigure}

\caption{Une couche plane infinie (à gauche) et une couche courbe lisse (à droite)}
\label{fig: couche infinie fr}
\end{figure}
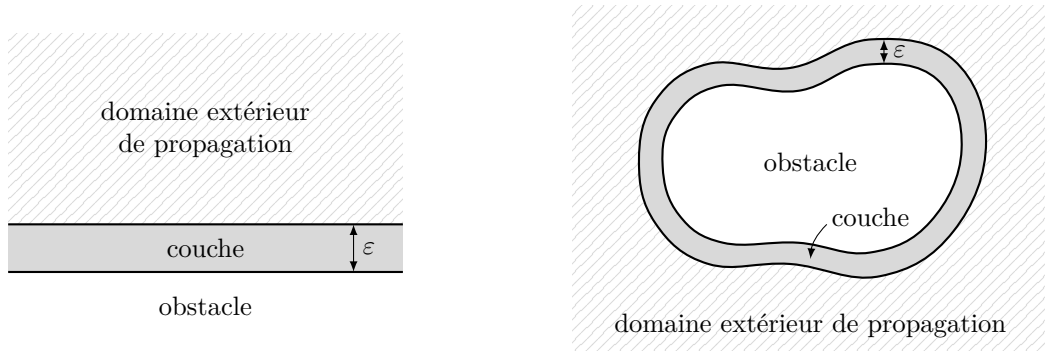%
	\tikzset{external/export=false}%
	\tikzsetfigurename{main-figure}

\begin{textesansboite}
Les couches planes infinies et les couches courbes lisses (voir la figure~\ref{fig: couche infinie fr}) ont été étudiées dans les années 90 par \cite{BenLem96, EngNed93}. Leur méthode consiste à rechercher un développement asymptotique de la solution exacte sous la forme suivante, où $s$ désigne l'abscisse curviligne le long de la couche et $\nu $ est la coordonnée normale associée (par exemple, $(s,\nu )=(x,y)$ pour une couche droite),
\begin{equation}{\label{eq: ansatz couche infinie homogene fr}}
	\left\{\begin{array}{ll}
		\displaystyle \sum_{n=0}^\infty  \varepsilon ^{n} \.u_{n} (s,\nu ) & \text{en dehors de la couche }  (\nu >0) \vspace{1pt}\\
		\displaystyle \sum_{n=0}^\infty  \varepsilon ^{n} \.U_{n} \Big(s, \frac{\nu }{\varepsilon } \Big) & \text{dans la couche } (-\varepsilon <\nu <0).
	\end{array}\right.
\end{equation}
La variable $\frac{\nu }{\varepsilon }$ dans $U_{n}$ permet de décrire le comportement de la solution à l'échelle microscopique dans la couche. Les termes de ce développement peuvent être calculés par récurrence par rapport à $n$. Ensuite, \cite{BenLem96, EngNed93} déduisent des modèles effectifs qui remplacent la couche par une condition au bord effective. Le plus simple est le problème limite, qui impose la même condition au bord que le problème initial, mais sans la couche. Sa solution est $u_{0}$, le premier terme du développement. Elle converge vers la solution exacte en dehors de la couche en $O(\varepsilon )$ pour la norme $H^{1}$. Ce modèle n'est pas satisfaisant pour les applications qui s'intéressent aux effets de la couche, car il ignore complètement cette structure mince.
De manière plus intéressante, \cite{BenLem96, EngNed93} proposent des modèles effectifs d'ordre supérieur. Par exemple, si l'on considère une couche plane infinie avec une condition de Dirichlet, ces modèles contiennent une condition au bord de type Robin, resp. Ventcel, qui peut s'écrire avec nos notations comme suit :
\begin{equation}{\label{eq: conditions effectives classiques intro fr}}
    \mu _{0} \partial _{n} v_{\varepsilon } = -\frac{\mu _{1} }{\varepsilon } v_{\varepsilon } , \qquad  \text{resp.}  \qquad  \mu _{0} \partial _{n} v_{\varepsilon } = -\frac{\mu _{1} }{\varepsilon } v_{\varepsilon } + \frac{\varepsilon }{3} (\mu _{1} \partial _{x} ^2 +\omega ^2 \rho _{1} )v_{\varepsilon } .
\end{equation}
On obtient ainsi un modèle qui converge d'ordre $O(\varepsilon ^{3} )$, resp. $O(\varepsilon ^{5} )$, en dehors de la couche. Dans le cas simple d'une couche infinie plane, ces modèles peuvent être obtenus directement en exprimant l'opérateur de Dirichlet-to-Neumann associé à la couche à l'aide de la transformée de Fourier, puis en effectuant un développement limité par rapport à $\varepsilon $ de l'expression obtenue.\\

\def\highlight#1{\emph{#1}}%
Ces résultats ont ensuite été étendus aux couches hétérogènes et périodiques dans \cite{AmmHe97, AmmHe98, AmmLat98} et, plus récemment, dans \cite{MadVal07, DelHadJol12, DelCla13, Del14, DelHew20, Bou23}. Dans ce cas, les variations de la solution à l'échelle microscopique ne sont pas confinées dans le revêtement, comme dans le cas homogène, mais débordent dans une fine couche limite à l'extérieur du revêtement physique. Ainsi, la construction d'un développement asymptotique requiert des techniques plus complexes. Ces techniques sont également utiles pour étudier les couches dans des coins ; nous allons donc les présenter. Il existe deux grands types de méthodes : les \highlight{développements asymptotiques multi-échelle} et les \highlight{développements asymptotiques raccordés}. Voir respectivement \cite{VisLyu57, Naz81a, Naz81b} et \cite{Van64, Ili92, Eck94} pour une présentation générale, et \cite{TorViaDau06} pour une comparaison de ces méthodes dans le contexte d'un coin sans couche. Ces deux méthodes impliquent deux types de champs : les \og champs lointains \fg{}, qui dépendent de l'échelle macroscopique décrite par $(x,y)$, et les \og champs proches \fg{}, qui dépendent de l'échelle microscopique décrite par $( \frac{x}{\varepsilon } , \frac{y}{\varepsilon } )$. Dans les développements multi-échelle, les champs lointains et proches sont définis et coexistent sur l'ensemble du domaine (ou sur une partie significative de celui-ci) et les champs proches tendent vers 0 à l'infini, de sorte qu'ils décrivent l'effet de couche limite dans le voisinage de la couche, appelé \og zone proche \fg{}. En revanche, les développements asymptotiques raccordés ne font intervenir les champs proches que dans la zone proche et les champs lointains que dans la \og zone lointaine \fg{} (le reste du domaine), et les champs proches et lointains doivent coïncider asymptotiquement dans une zone intermédiaire appelée \og zone de raccord \fg{}. Cette zone est située à une distance $y\in {]\varepsilon ^{a} ,2\mkern1mu \varepsilon ^{a} [}$ de la couche, avec $a\in {]0,1[}$ fixe, de sorte que, du point de vue des champs lointains, elle tend vers l'interface avec la couche, mais du point de vue des champs proches, elle tend vers l'infini, puisque $\varepsilon ^{a} /\varepsilon  \rightarrow \infty $ lorsque $\varepsilon \rightarrow 0$. Par conséquent, les développements asymptotiques raccordés contiennent une \og hypothèse de raccord \fg{} qui relie le comportement au voisinage de la couche des champs lointains et celui à l'infini des champs proches. Dans ce cas, les champs proches divergent généralement à l'infini, ce qui constitue une différence importante par rapport aux développements multi-échelle. Notons que ces techniques ont également été appliquées à des problèmes connexes, notamment les configurations tridimensionnelles \cite{Cham16}, les couches aléatoires \cite{BouFliGio26} et les couches limites en homogénéisation \cite{AllAma99, Vin16, Ben21}.\\

Les problèmes impliquant des couches finies ou des couches dans un coin posent des difficultés supplémentaires en raison de la présence de comportements singuliers à proximité du coin (ou de l'extrémité d'une couche finie plane). Par exemple, si la solution exacte possède une condition de Dirichlet, alors la fonction $u_{1}$ dans \eqref{eq: ansatz couche infinie homogene fr} satisfait la condition au bord de Dirichlet $u_{1} = - \frac{\mu _{0} }{\mu _{1}} \partial _{n} u_{0}$. Pour une couche infinie, on peut montrer que $u_{0}$ est régulière au voisinage de la couche, donc $\partial _{n} u_{0} \in  H^{1/2}$. Cependant, dans un coin d'angle $\Theta $, on a
\[u_{0} (r,\theta )= c\,r^{\pst} \sin \!\Big(\pst \theta \Big) + u_{0} ^{\mathrm{reg}} \qquad  \text{avec }  c\in \mathbb{C} \text{ et } u_{0} ^{\mathrm{reg}} \in H^2  .\]
La fonction $r^{\pst} \sin(\pst \theta )$ est appelée \og singularité de coin \fg{} car elle n'est pas régulière au niveau du coin (par exemple, elle n'est pas dans $H^2 $ lorsque $\Theta  > \pi$). On a donc $\partial _{n} u_{0} = -\pst c\.r^{\pst-1} + \partial _{n} u_{0} ^{\mathrm{reg}}$, qui n'appartient pas à $H^{1/2}$ lorsque $\pst<1$, c'est-à-dire lorsque $\Theta >\pi $ (et $c\neq 0$). Dans ce cas, on ne peut pas chercher $u_{1}$ dans $H^{1}$ et il faut considérer un champ lointain $u_{1}$ qui diverge au coin. Pour un angle $\Theta $ plus petit, $u_{1}$ est bien défini dans $H^{1}$, mais un problème similaire se pose pour certains champs lointains $u_{n}$ avec $n > 1$.\\

Ces problèmes de coins ont été étudiés pour la première fois dans le cadre de l'équation de Poisson avec un angle arbitraire dans \cite{CalCosDauVia06, Via03} (nous citons ensemble les articles et les thèses associées). Ce~travail construit un développement asymptotique à tout ordre à l'aide de la méthode des développements asymptotiques multi-échelle. Le développement obtenu est de la forme
\[\sum _{p\in \mathbb{N}+\pst \mathbb{N}} \.\sum _{\ell =0} ^{\infty } \.\varepsilon ^{p} \ln ^{\ell } \!\varepsilon  \cdot  u_{p,\ell } (x,y)\]
loin du coin. Notez qu'apparaissent des puissances non entières de $\varepsilon $ dans $\NN:=\mathbb{N}+\pst \mathbb{N}$ et des puissances entières de $\ln \varepsilon $, qui sont liées aux singularités de coin. Cela étend un travail antérieur \cite{CosDau96} qui étudie un problème de Poisson sans couche, mais avec une condition aux limites mixte de type Neumann $\partial _{n} u_{\varepsilon } =0$ et Robin $u_{\varepsilon } +\varepsilon  \.\partial _{n} u_{\varepsilon } =0$. Le développement dans ce cas présente des similitudes --- liées au fait que cette condition de Robin est l'une des conditions effectives de \eqref{eq: conditions effectives classiques intro fr} --- mais son obtention est moins technique. Ces développements asymptotiques à tout ordre ont été généralisés aux couches périodiques dans \cite{DelSchSem16}, toujours pour l'équation de Poisson. Ce travail utilise à la fois un développement asymptotique multi-échelle autour de la couche et un développement asymptotique raccordé autour de l'angle pour traiter le couplage complexe de ces deux types de petites structures.\\

En comparaison, les problèmes de Helmholtz impliquant des couches finies présentent non seulement les mêmes difficultés, mais conduisent également à des singularités plus complexes. En effet, les singularités de coin qui apparaissent à tous les ordres ne peuvent pas être organisées en une famille de suites récurrentes, comme dans les problèmes de Poisson, car elles sont liées les unes aux autres par un arbre de dépendances (pour plus de détails, voir la remarque~\ref{rq: utilite des operateurs R}). Par conséquent, construire un développement asymptotique à tout ordre et pour tout angle est très difficile. Ceci explique peut-être a posteriori pourquoi les développements asymptotiques antérieurs pour l'équation de Helmholtz se limitaient à l'ordre 2, ce qui permet d'éviter le problème en décrivant explicitement les premières singularités. \cite{BenMakTor11, Mak08} propose un développement d'ordre 2 pour une couche homogène et un angle plat ($\Theta =\pi $), tandis que \cite{SemDelSch18} généralise ce résultat à une couche périodique avec un angle arbitraire. Ces deux travaux utilisent des développements asymptotiques raccordés. Ces développements d'ordre limité sont difficiles à justifier rigoureusement, car prouver une estimation d'erreur optimale à un ordre donné nécessite davantage de termes du développement. \cite{BenMakTor11, Mak08, SemDelSch18} prouvent soit une estimation d'erreur sous-optimale, soit supposent que les termes suivants peuvent être construits comme les précédents (ce qui nous semble hâtif, compte tenu de la croissance rapide de la complexité du développement par rapport à l'ordre).\\

Ces travaux sur les couches finies présentent une complexité calculatoire croissante en raison de la complexité grandissante des singularités de coin. Nous proposons d'y remédier grâce à un ensemble d'outils calculatoires algébriques qui remplace des calculs fastidieux par des calculs abstraits et généralisables. Cela nous permet de traiter un cas plus complexe --- l'équation de Helmholtz avec une couche finie à tout ordre pour tout angle --- tout en obtenant des conditions de raccord plus explicites. Par ailleurs, notons que \cite{JolTor06, JolTor08, Tor04} ont démontré un développement asymptotique à tout ordre pour l'équation de Helmholtz dans le cas apparenté d'un demi-plan muni d'une fente mince, à l'aide d'un développement asymptotique raccordé. La complexité à laquelle nous sommes confrontés ici semble donc spécifique à l'association de l'équation de Helmholtz et d'une couche finie.\\

Comme expliqué ci-dessus, construire des développements asymptotiques est bien plus difficile pour les couches finies que pour les couches infinies, mais cela est encore plus vrai pour les modèles effectifs. En~effet, les modèles effectifs d'ordre élevé pour les couches finies restaient une question ouverte, même dans des cas plus simples comme l'équation de Poisson ou pour des angles spécifiques tels que $\Theta =\pi $, et ce malgré une littérature substantielle sur les développements asymptotiques (voir ci-dessus). Une approche naïve consiste à imposer les conditions aux limites effectives \eqref{eq: conditions effectives classiques intro fr} sur la totalité de $\Gamma $, c'est-à-dire jusqu'au coin. Pour la condition de Robin, par exemple, cela s'écrit ainsi :
\[\left\{\begin{array}{rl}
	\mu _{0} \Delta  v_{\varepsilon } +\omega ^2 \rho _{0} v_{\varepsilon } = \fsource & \text{ dans }  \Omega \\[2pt]
	\mu _{0} \partial _{n} v_{\varepsilon } + \frac{\mu _{1} }{\varepsilon } v_{\varepsilon } =0& \text{ sur }  \Gamma .
\end{array}\right.\]
\cite{Via05,Via03} a démontré pour ce modèle un taux de convergence optimal de $O(\varepsilon ^{\min (3,\pst)} )$ sur le domaine total (et nous pensons que l'on a $O(\varepsilon ^{\min (3,2\pst)} )$ loin du coin). Cette détérioration est due à la présence des singularités de coin et est plus importante pour les grands angles, car plus l'angle est grand, moins les singularités telles que $r^{\pst} \sin (\pst \theta )$ sont régulières. De plus, la condition de Ventcel de \eqref{eq: conditions effectives classiques intro fr} implique une dérivée tangentielle seconde, ce qui nécessite une condition aux limites supplémentaire au niveau du coin. Le choix de la condition appropriée n'est en général pas du tout évident.\\

\cite{AuvVia19,AuvVia18,Auv18} a proposé deux modèles effectifs pour l'équation de Poisson avec $\Theta =\pi $. Ces modèles sont d'ordre 1,5 pour la norme $L^2 $ et 0,5 pour la norme $H^{1}$. Le premier considère sur $\Gamma $ la condition de Robin $\mu _{0} \partial _{n} u = -\frac{\mu _{1} }{\varepsilon } u$ de \eqref{eq: conditions effectives classiques intro fr} (selon nos notations), mais avec un coefficient localement perturbé près du coin à l'échelle $\varepsilon $. La perturbation est choisie en résolvant un problème d'optimisation qui dépend du premier champ proche. Malheureusement, le calcul de la solution de ce modèle nécessite de mailler finement le domaine près du coin. Dans leur second modèle effectif, la condition $\mu _{0} \partial _{n} u=-\frac{\mu _{1} }{\varepsilon } u +\varepsilon \beta  \partial _{x} ^2  u$ est imposée sur $\Gamma $ avec un coefficient $\beta $ qui est à nouveau choisi via un problème d'optimisation. Il s'agit donc d'une condition de Ventcel, mais différente de \eqref{eq: conditions effectives classiques intro fr}. Des justifications théoriques ont été fournies pour le premier modèle, mais pas pour le second. Ces deux modèles n'apportent qu'une amélioration limitée de la convergence $L^2 $ par rapport aux modèles naïfs. De plus, leur absence de signification physique les rend difficiles à généraliser ou à améliorer.
\end{textesansboite}

\section{Résumé de la thèse}

\begin{textesansboite}
Dans cette section, nous présentons les principaux résultats et contributions de chaque chapitre. La thèse peut être divisée en trois parties principales :
\begin{itemize}
    \item dans les chapitres \ref{chap: SF} à~\ref{chap: construction du DA}, nous construisons un développement asymptotique à tout ordre de $u_{\varepsilon }$, la solution de \eqref{eq: pb etudie fr},
	\item dans les chapitres \ref{chap: complements sur le DA} à~\ref{chap: validation numerique du DA}, nous donnons des propriétés supplémentaires de ce développement et nous vérifions numériquement sa validité,
    \item dans le chapitre~\ref{chap: ME}, nous en déduisons des modèles effectifs d'ordre élevé, resp. d'ordre~3 et 5.
\end{itemize}
Les résultats des chapitres \ref{chap: raccord algebrique} et~\ref{chap: construction du DA} ont été publiés dans \cite{Bau26}. Plusieurs autres articles sont en préparation concernant les chapitres \ref{chap: complements sur le DA} à~\ref{chap: ME}. Les résultats de cette thèse ont été présentés lors des conférences suivantes : Waves 2024 (Berlin), Aspect 2024 (Metz), ICOSAHOM 2025 (Montréal), CANUM 2026 (Bretagne), Waves 2026 (Montréal). Par ailleurs, une partie du code écrit durant cette thèse a été publié sous forme de librairie \Cplusplus{} \cite{OPTEC}.\\

Tout d'abord, la section suivante présente la méthode utilisée pour le développement asymptotique.
\end{textesansboite}

\subsection{Présentation du développement asymptotique raccordé}
\label{sec: ansatz du dvlpmt raccorde fr}

\begin{textesansboite}
Dans ce travail, nous avons choisi la méthode des \emph{développements asymptotiques raccordés} plutôt que celle des développements multi-échelle. En effet, elle offre une construction plus intrinsèque des champs lointains (indépendante d'une fonction de troncature, voir \cite{TorViaDau06}), ce qui nous a semblé plus prometteur pour obtenir des modèles effectifs. Afin de prendre en compte les différents comportements de la solution dans la couche, à proximité du coin et loin du coin et de la couche, nous divisons $\Omega _{\varepsilon }$ en trois zones, illustrées dans la figure~\ref{fig: zones du dvlpmt raccorde fr}. Dans chaque zone, nous postulons un développement asymptotique en puissances de $\varepsilon $ et de $\ln \varepsilon $, appelé \og ansatz \fg{}.
\end{textesansboite}

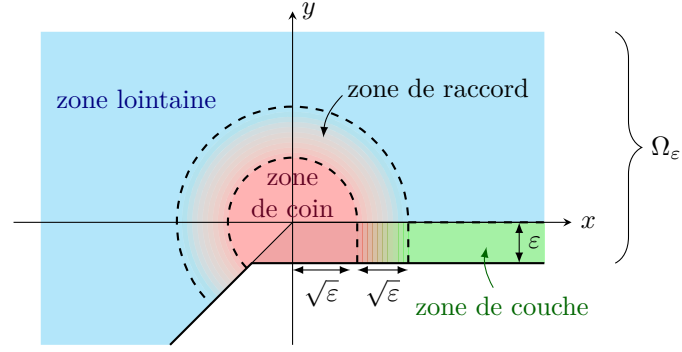
\begin{figure}[h]
\def\epsFig{0.6} 
\def\widthFig{3.7} 
\def\heightFig{2.8} 
\def\depthFig{1.8} 
\def\ThetaFig{225} 
\def\racInt{0.95} 
\def\racExt{1.7} 
\begin{tikzpicture}[scale=0.9]
	\coordinate (finLambda) at ({-\depthFig*cot(\ThetaFig)}, -\depthFig);
	\coordinate (sommetCouche) at ({-\epsFig*cot(\ThetaFig)}, -\epsFig);
	\draw (-\widthFig-1.8,0) node {\ }; 
	\fill [color=white!85!gray] (0,0) -- (\widthFig,0) -- (\widthFig,-\epsFig) -- (sommetCouche) -- cycle; 
	\fill [color=cyan, opacity=0.3] (\racExt,0) arc (0:\ThetaFig:\racExt) -- (finLambda) -- (-\widthFig,-\depthFig) -- (-\widthFig,\heightFig) -- (\widthFig,\heightFig) -- (\widthFig,0)  -- cycle; 
	\fill [color=green, opacity=0.3] (\racExt,0) rectangle (\widthFig,-\epsFig); 
	\fill [color=red, opacity=0.3] (\racInt,0) arc (0:\ThetaFig:\racInt) -- (sommetCouche) -- (\racInt,-\epsFig) -- cycle; 
	\foreach \i in {5,15,...,95}{
		\def\Rint{(\racInt*(100-\i+5)+\racExt*(\i-5))/100}
		\def\Rext{(\racInt*(100-\i-5)+\racExt*(\i+5))/100}
		\fill [color=cyan!\i !red, opacity=0.3] ({\Rint},0) -- ({\Rext},0) arc (0:\ThetaFig:{\Rext}) -- ({\Rint*cos(\ThetaFig)},{\Rint*sin(\ThetaFig)}) arc (\ThetaFig:0:{\Rint}) -- cycle; 
		\fill [color=green!\i !red, opacity=0.3] ({\Rint},0) rectangle ({\Rext},-\epsFig); 
	}
	\draw [thick] (finLambda) -- (sommetCouche) -- (\widthFig,-\epsFig);
	\draw (sommetCouche) -- (0,0) -- (\widthFig,0);
	\draw [thick,dashed] (\racInt,-\epsFig) -- (\racInt,0) arc (0:\ThetaFig:\racInt);
	\draw [thick,dashed] (\racExt,-\epsFig) -- (\racExt,0) arc (0:\ThetaFig:\racExt);
	\draw [thick,dashed] (\racExt,0) -- (\widthFig,0);
	\draw [>=stealth,->] (-\widthFig-0.4,0) -- (\widthFig+0.4,0) node [right] {$x$};
	\draw [>=stealth,->] (0,-\depthFig) -- (0,{\heightFig+0.3}) node [right] {$y$};
	\draw [decorate, decoration={brace,amplitude=10pt}, xshift=4pt] (\widthFig+0.9,\heightFig) -- (\widthFig+0.9,-\epsFig) node [midway,right,xshift=10pt] {$\Omega_\varepsilon$};
	\draw [color=black!50!blue] (-\widthFig*0.62,\heightFig*0.65) node {zone lointaine};
	\draw [color=black!60!green,<-,>=latex] (\widthFig*0.77,-\epsFig*0.6) to[bend right=20] ++(0.1,-0.6) node[below]{\ \ zone de couche};
	\draw [color=black!50!purple] (0,\racInt*0.42) node {\parbox{20mm}{\centering zone\\ de coin}};
	\draw [<->,>=latex] (\widthFig*0.9,0) -- node[midway,right]{$\varepsilon$} (\widthFig*0.9,-\epsFig);
	\draw [<->,>=latex] (0,{-\epsFig-0.1}) -- node[midway,below]{$\sqrt\varepsilon$} ({\racInt-0.1},{-\epsFig-0.1});
	\draw [<->,>=latex] (\racInt,{-\epsFig-0.1}) -- node[midway,below]{$\sqrt\varepsilon$} (\racExt,{-\epsFig-0.1});
	\draw [<-,>=latex] ({(\racInt+\racExt)/2*cos(70)},{(\racInt+\racExt)/2*sin(70)}) to[bend left=20] ++(0.28,0.62) node[right]{\raisebox{2mm}{\!zone de raccord}};
\end{tikzpicture}\par
\caption{Zones du développement asymptotique raccordé}
\label{fig: zones du dvlpmt raccorde fr}
\end{figure}

\begin{textesansboite}
Nous rappelons que la géométrie est définie en section~\ref{sec: motivations et pb etudie fr}. Supposons que, pour tout $(p, \ell ) \in \NN \times \mathbb{N}$, il existe trois fonctions indépendantes de $\varepsilon $ --- à savoir $u_{p,\ell }$ définie sur $\Omega $ et appelée \og champ lointain \fg{}, $U_{p,\ell }$ définie sur $\couche$ et appelée \og champ de couche \fg{}, et $S_{p,\ell }$ définie sur le domaine mis à l'échelle $\Omega _{1}$ et appelée \og champ de coin \fg{} --- telles que $u_{\varepsilon }$ s'écrit formellement ainsi :
\begin{itemize}
	\item $u_{\varepsilon } (x,y) =\displaystyle \sum _{p\in \NN} \sum _{\ell \in \mathbb{N}} \varepsilon ^{p} \ln ^{\ell } \!\varepsilon \;u_{p,\ell } (x,y)$ quand $r:=\|(x,y)\|\geqslant \sqrt{\varepsilon }$ et $(x,y)\not\in \mathbb{R}_+\times {]-\varepsilon ,0[}$ (zone lointaine),
	\item $u_{\varepsilon } (x,y) =\displaystyle \sum _{p\in \NN} \sum _{\ell \in \mathbb{N}} \varepsilon ^{p} \ln ^{\ell } \!\varepsilon \;U_{p,\ell } \Big(x, \frac{y}{\varepsilon } \Big)$ quand $x\geqslant \sqrt{\varepsilon }$ et $y\in {]-\varepsilon ,0[}$ (zone de couche),
	\item $u_{\varepsilon } (x,y) =\displaystyle \sum _{p\in \NN} \sum _{\ell \in \mathbb{N}} \varepsilon ^{p} \ln ^{\ell } \!\varepsilon \;S_{p,\ell } \Big( \frac{x}{\varepsilon } , \frac{y}{\varepsilon } \Big)$ quand $r\leqslant 2\sqrt{\varepsilon }$ (zone de coin).
\end{itemize}\fbi
\end{textesansboite}

\begin{textesansboite}
La présence de puissances entières de $\varepsilon $ est tout à fait classique en analyse asymptotique. Quant aux puissances entières de $\varepsilon ^{\pi /\Theta }$ et de $\ln \varepsilon $, elles sont présentes dans d'autres développements asymptotiques impliquant des coins, voir \cite{CalCosDauVia06, DelSchSem16, SemDelSch18} et la section~\ref{sec: etat de l'art fr}. De plus, nous verrons dans la proposition~\ref{borne sur l} que nous avons :
\[\forall p\in \NN,\; \exists \ell _{p} \in \mathbb{N},\; \forall \ell >\ell _{p} ,\qquad  u_{p,\ell } = 0 \text{ et } U_{p,\ell } = 0 \text{ et } S_{p,\ell } = 0.\]\\

Notons que nous traitons la couche de la même manière que \eqref{eq: ansatz couche infinie homogene fr} dans \cite{BenLem96,EngNed93}, ce qui est possible car la couche est homogène. Cela nous permet de nous concentrer sur le développement asymptotique raccordé au niveau du coin. Les champs de coin $S_{p,\ell }$ sont, dans notre contexte, les \og champs proches \fg{} du développement raccordé.\\

Les chapitres \ref{chap: raccord algebrique} et~\ref{chap: construction du DA} montrent comment construire les champs $u_{p,\ell }$, $U_{p,\ell }$ et $S_{p,\ell }$. La construction est récursive par rapport à $p$. Pour chaque $p$, tous les champs $u_{p,\ell }$, $U_{p,\ell }$ et $S_{p,\ell }$ sont construits en même temps, mais seul un nombre fini d'entre eux est non nul.\\

La première étape, qui est aussi la plus simple, consiste à déterminer les équations volumiques et les conditions aux limites des différents champs. En injectant l'ansatz ci-dessus dans l'équation de Helmholtz \eqref{eq: pb etudie fr}, en utilisant le fait que $\partial _{x} ^2  \big[\varphi ( \frac{x}{\varepsilon } )\big]=\varepsilon ^{-2} [\partial _{X} ^2 \varphi ](X = \frac{x}{\varepsilon } )$ et $\partial _{y} ^2  \big[\varphi ( \frac{y}{\varepsilon } )\big]=\varepsilon ^{-2} [\partial _{Y} ^2 \varphi ](Y = \frac{y}{\varepsilon } )$ pour toute fonction $\varphi $, et en identifiant formellement les puissances de $\varepsilon $ et de $\ln \varepsilon $, on peut facilement déduire les équations suivantes.
\begin{equation}{\label{eq: eq vol champs lointains fr}}
	\left\{\begin{array}{r@{\;=\;}l@{\;}l}
	\mu _{0} \Delta  u_{p,\ell } +\omega ^2 \rho _{0} u_{p,\ell } & f\,\delta _{p,0} \,\delta _{\ell ,0} & \text{ dans }  \Omega  \\[0.5mm]
	u_{p,\ell } & U_{p,\ell } & \text{ sur }  \Gamma \\[0.5mm]
	u_{p,\ell } & 0 & \text{ sur }  \bordgauche
	\end{array}\right. 
\end{equation}
\begin{equation}{\label{eq: eq vol champs de couche fr}}
	\left\{\begin{array}{r@{\;=\;}l@{\;}l}
	\mu _{1} \partial _{Y} ^2  U_{p,\ell } & -(\mu _{1} \partial _{x} ^2 +\omega ^2 \rho _{1} )U_{p-2,\ell } & \text{ dans }  \couche \\[0.5mm]
	\mu _{1} \partial _{Y} U_{p,\ell } & \mu _{0} \,\partial _{y} u_{p-1,\ell } & \text{ sur }  \Gamma  \\[0.5mm]
	U_{p,\ell } & 0 & \text{ sur }  \borddroit
	\end{array}\right.
\end{equation}
\begin{equation}{\label{eq: eq vol champs de coin fr}}
	\left\{\begin{array}{r@{\;=\;}ll}
	\operatorname{div}(\mu  \nabla  S_{p,\ell } ) & -\omega ^2 \rho  \,S_{p-2,\ell } & \text{ dans }  \Omega _{1} \\[0.5mm]
	S_{p,\ell } & 0 & \text{ sur }  \partial  \Omega _{1}
	\end{array}\right.
\end{equation}
Par convention, on pose ici $u_{p,\ell } = 0$, $U_{p,\ell } = 0$ et $S_{p,\ell } = 0$ pour tout $(p,\ell ) \in (\mathbb{R} \setminus \NN)\times \mathbb{N}$, et $\delta _{i,j} := 1$ si $i = j$ et 0 sinon.
\end{textesansboite}

\begin{paragraphesansboite}{Remarques :}
\begin{itemize}
	\item La condition $\mu _{1} \partial _{Y} U_{p,\ell } = \mu _{0} \partial _{y} u_{p-1,\ell }$ est incluse dans le problème \eqref{eq: eq vol champs de couche fr} vérifié par $U_{p,\ell }$, tandis que $u_{p,\ell } = U_{p,\ell }$ est incluse dans le problème \eqref{eq: eq vol champs lointains fr} satisfait par $u_{p,\ell }$, de sorte que la construction est récursive par rapport à $p$ : $u_{p-1,\ell }$ permet de construire $U_{p,\ell }$, ce qui permet de construire $u_{p,\ell }$.
	\item Le champ de couche $U_{p,\ell }$ satisfait un problème unidimensionnel par rapport à $Y$, la variable $x$ jouant le rôle d'un paramètre.\\
\end{itemize}\fbi
\end{paragraphesansboite}

\begin{textesansboite}
Les équations \eqref{eq: eq vol champs lointains fr} à~\eqref{eq: eq vol champs de coin fr} suffiraient à définir de manière unique les champs s'ils appartenaient à leurs espaces variationnels naturels (par exemple, $H^{1} (\Omega )$ pour $u_{p,\ell }$). Mais il faut prendre en compte une \emph{hypothèse de raccord} (\og matching assumption \fg{} en anglais) : les champs lointains et de coin doivent coïncider asymptotiquement à l'intersection des zones lointaines et de coin, et de même pour les champs de couche et de coin. Ces intersections constituent la zone de raccord (cf.\ figure~\ref{fig: zones du dvlpmt raccorde fr}). Étant donné que $\varepsilon \rightarrow 0$ et $\sqrt{\varepsilon } /\varepsilon  \rightarrow \infty $, cette zone tend vers $(0,0)$ pour les champs lointains et de couche, mais elle tend vers l'infini pour les champs de coin. Ainsi, l'hypothèse de raccord relie formellement le comportement asymptotique des champs lointains et de couche au niveau du coin à celui des champs de coin à l'infini :
\begin{equation}{\label{eq: hypothese de raccord fr}}
\left\{\begin{array}{r@{}ll}
	\displaystyle \sum _{p,\ell } \varepsilon ^{p} \ln ^{\ell } \!\varepsilon \;u_{p,\ell } (x,y)\. &{}\approx \.\displaystyle \sum _{p,\ell } \varepsilon ^{p} \ln ^{\ell } \!\varepsilon \;S_{p,\ell } \Big( \frac{x}{\varepsilon } , \frac{y}{\varepsilon } \Big) & \text{ dans }  \Omega  \text{ quand }  r\rightarrow 0 \text{ et } \frac{r}{\varepsilon } \rightarrow \infty  \vrule width 0pt height 3.5ex \vspace{2pt}\\
	\displaystyle \sum _{p,\ell } \varepsilon ^{p} \ln ^{\ell } \!\varepsilon \;U_{p,\ell } \Big(x, \frac{y}{\varepsilon } \Big) &{}\approx \.\displaystyle \sum _{p,\ell } \varepsilon ^{p} \ln ^{\ell } \!\varepsilon \;S_{p,\ell } \Big( \frac{x}{\varepsilon } , \frac{y}{\varepsilon } \Big) & \text{ dans }  \couche \text{ quand }  x\rightarrow 0 \text{ et } \frac{x}{\varepsilon } \rightarrow \infty .
\end{array}\right.
\end{equation}
La signification rigoureuse de ces ``$\approx $'' sera précisée plus tard (cf.\ définition~\ref{def: raccord algebrique}).\\

Nous verrons que chaque champ lointain $u_{p,\ell }$ possède, au niveau du coin, un développement asymptotique par rapport à $r$ qui est grossièrement une somme de puissances de $r$, dont certaines sont positives (comme dans un développement limité par exemple). L'hypothèse de raccord implique que ces puissances positives de $r$ apparaissent également dans les développements asymptotiques des champs de coin à l'infini, où elles sont divergentes. Nous les appelons singularités pour les champs de coin. Inversement, les développements asymptotiques des champs de coin à l'infini contiennent des puissances négatives de $r$ correspondant à la décroissance de la partie variationnelle, et ces puissances doivent se retrouver dans les champs lointains, ce qui correspond à des singularités au niveau du coin. Ainsi, on ne peut pas rechercher ces champs dans leurs espaces variationnels naturels. Pour formuler des problèmes qui définissent ces champs de manière unique, il faut spécifier leur partie singulière, c'est-à-dire leur partie non variationnelle.\\

Les parties singulières des champs sont choisies de manière à satisfaire l'hypothèse de raccord \eqref{eq: hypothese de raccord fr}. Les parties singulières de tout champ $u_{p,\ell }$ ou $U_{p,\ell }$ dépendent des comportements à l'infini des champs de coin précédents, c'est-à-dire des $S_{p',\ell '}$ avec $p'<p$. Inversement, la partie singulière de n'importe quel $S_{p,\ell }$ dépend des comportements au coin des champs lointains et de couche précédents. Comme dans des travaux antérieurs \cite{BenMakTor11,DelSchSem16,SemDelSch18}, le choix de ces parties singulières peut se réduire à fixer correctement certains coefficients apparaissant dans les comportements singuliers, selon certaines relations que nous appelons \emph{conditions de raccord} (\og matching conditions \fg). Obtenir ces conditions est assez simple pour des couches périodiques infinies, mais s'avère beaucoup plus difficile dans les coins et en particulier dans notre cas en raison de la complexité exceptionnelle des singularités de coin à tout ordre. Ceci nous a amenés à introduire de nouveaux outils plus efficaces pour gérer le raccord sans être contraints par la complexité des singularités. Ce travail couvre les chapitres \ref{chap: SF} et~\ref{chap: raccord algebrique} et est présenté dans la section suivante.
\end{textesansboite}

\subsection{Chapitres \ref{chap: SF} et \ref{chap: raccord algebrique}: les conditions de raccord}
\label{sec: resume raccord en intro fr}

\begin{textesansboite}
\phantomsection\label{page: definition de Pi fr}
Pour effectuer simultanément le raccordement des champs de coin avec les champs lointains et de couche, nous fusionnons ces deux derniers en un seul champ noté $\bfu_{p,\ell }$ et appelé \og champ lointain concaténé \fg{}. Il~est défini sur un nouveau domaine $\Pi $, défini de la manière suivante.
\begin{itemize}
	\item Si $\Theta  \leqslant  \frac{3\pi}{2}$, alors $\Pi  := \Omega  \cup  \Gamma  \cup  \couche$ (union disjointe) et c'est un ouvert de $\mathbb{R}^2 $.
	\item Si $\Theta  > \frac{3\pi}{2}$, alors $\Omega $ et $\couche$ s'intersectent en tant que sous-ensembles de $\mathbb{R}^2 $, donc la définition précédente n'est plus valable (voir figure \ref{fig:Pi fr}). Nous définissons donc $\Pi $ comme le recollement disjoint de $\Omega $ et de $\couche$ le long de $\Gamma $ ($\Pi $ est alors une variété riemannienne plate).
\end{itemize}\fbi
\end{textesansboite}

\begin{figure}[h]
\def\widthFig{2.3} 
\def\heightFig{1.9} 
\def\depthFig{-1.55} 
\def\epsFig{0.8} 
\def\ThetaFig{215} 
\begin{tikzpicture}
	\coordinate (finLambda) at (-\widthFig, {-\widthFig*tan(\ThetaFig-180)});
	\fill [pattern={Lines[angle=45,line width=.2pt]}, pattern color=white!70!gray] (0,0) -- (\widthFig,0) -- (\widthFig,\heightFig) -- (-\widthFig,\heightFig) -- (finLambda) -- cycle;
	\fill [color=white!70!gray] (0,0) rectangle (\widthFig,-\epsFig);
	\draw [thick] (finLambda) -- (0,0) -- (\widthFig,0);
	\draw [thick] (\widthFig,-\epsFig) -- (0,-\epsFig) -- (0,0);
	\draw [>=stealth,->] (-\widthFig-0.3,0) -- (\widthFig+0.3,0);
	\draw [>=stealth,->] (0,{\depthFig-0.2}) -- (0,{\heightFig+0.2});
	\draw (-\widthFig*2/5,\heightFig/2) node {$\Omega$};
	\draw [->,>=latex] (0.35,0) arc (0:\ThetaFig-5:0.35) node[near start,above] {$\Theta$} -- (\ThetaFig:0.36);
	\draw [<-,>=latex] (-\widthFig/2, {-\widthFig/2*tan(\ThetaFig-180)}) to[bend right=20] ++(0.3,-0.4) node[right]{$\bordgauche$};
	\draw (\widthFig*0.45,-\epsFig/2) node {couche $\couche$};
	\draw [<->,>=latex] (\widthFig*0.85,0) -- node[midway,right]{1} (\widthFig*0.85,-\epsFig);
	\draw [<-,>=latex] (\widthFig/2,0) to[bend left=20] ++(0.15,0.5) node[above]{$\ \ \Gamma$};
	\draw [<-,>=latex] (\widthFig/2,-\epsFig) to[bend right=20] ++(0.15,-0.5) node[below]{$\ \ \borddroit$};
\end{tikzpicture}%
\hspace{10mm}\nolinebreak%
\def\ThetaFig{310}%
\begin{tikzpicture}
	\coordinate (finLambda) at ({\depthFig*cot(\ThetaFig)}, \depthFig);
	\fill [color=white!70!gray] (0,0) rectangle (\widthFig,-\epsFig);
	\fill [pattern={Lines[angle=45,line width=.2pt]}, pattern color=white!70!gray] (0,0) -- (\widthFig,0) -- (\widthFig,\heightFig) -- (-\widthFig,\heightFig) -- (-\widthFig,\depthFig) -- (finLambda) -- cycle;
	\fill [pattern={Hatch[angle=45,line width=1pt]}, pattern color=red!70!white] ({-\epsFig*cot(\ThetaFig)}, -\epsFig)  -- (0,0) -- (0,-\epsFig) -- cycle;
	\draw [thick] (finLambda) -- (0,0) -- (\widthFig,0);
	\draw [thick] (\widthFig,-\epsFig) -- (0,-\epsFig) -- (0,0);
	\draw [>=stealth,->] (-\widthFig-0.3,0) -- (\widthFig+0.3,0);
	\draw [>=stealth,->] (0,{\depthFig-0.2}) -- (0,{\heightFig+0.2});
	\draw (-\widthFig*2/5,\heightFig/2) node {$\Omega$};
	\draw [->,>=latex] (0.35,0) arc (0:\ThetaFig-5:0.35) node[midway,left] {$\Theta$} -- (\ThetaFig:0.36);
	\draw ({\widthFig*0.6},-\epsFig/2) node {couche $\couche$};
	\draw [<-,>=latex] ({\epsFig/3},{-\epsFig*2/3}) to[bend left=20] ++(-0.5,-0.5) node[left]{$\Omega\cap\couche$};
\end{tikzpicture}\par
\longCaption{Le domaine $\Pi $}{. Il est égal à $\Omega \sqcup \Gamma \sqcup \couche$ quand $\Theta \leqslant  \frac{3\pi }{2}$ (ouvert de $\mathbb{R}^2 $) et c'est une variété riemannienne plate quand $\Theta > \frac{3\pi }{2}$ ($\Omega \cap \couche \neq  \varnothing $).}
\label{fig:Pi fr}
\end{figure}

\begin{textesansboite}
Pour tout $(p, \ell ) \in \NN \times \mathbb{N}$, on définit le \emph{champ lointain concaténé} (\og far-and-layer field \fg{} en anglais) $\bfu_{p,\ell }$ et la variable radiale généralisée $\bfr$ par
\begin{equation}{\label{eq: def bfu bfr}}
    \bfu_{p,\ell } := \left\{\begin{array}{@{\;}l@{\,}l}
        u_{p,\ell } & \text{ dans }  \Omega  \\[0.5mm]
        U_{p,\ell } & \text{ dans }  \couche
    \end{array}\right. \qquad  \text{ et } \qquad  \bfr := \left\{\begin{array}{@{\;}l@{\,}l}
        r & \text{ dans }  \Omega  \\[0.5mm]
        x & \text{ dans }  \couche
    \end{array}\right.
\end{equation}

Une reformulation directe des équations \eqref{eq: eq vol champs lointains fr} et \eqref{eq: eq vol champs de couche fr} donne pour tout $(p, \ell ) \in \NN \times \mathbb{N}$ :
\begin{equation}{\label{eq: eq vol de bfu p l fr}}
\left\{\begin{array}{r@{\;=\;}ll}
	\mu _{0} \Delta  \bfu_{p,\ell } +\omega ^2 \rho _{0} \bfu_{p,\ell } & f\,\delta _{p,0} \,\delta _{\ell ,0} & \text{ dans }  \Omega  \\[0.5mm]
	\mu _{1} \partial _{Y} ^2  \bfu_{p,\ell } & -(\mu _{1} \partial _{x} ^2 +\omega ^2 \rho _{1} )\bfu_{p-2,\ell } & \text{ dans }  \couche \\[0.5mm]
	\mu _{1} \partial _{Y} \bfu_{p,\ell |Y=0^{-}} & \mu _{0} \,\partial _{y} \bfu_{p-1,\ell |y=0^{+}}  & \text{ sur }  \Gamma  \\[0.5mm]
	\bfu_{p,\ell |y=0^{+}} - \bfu_{p,\ell |Y=0^{-}} & 0 & \text{ sur }  \Gamma \\[0.5mm]
	\bfu_{p,\ell } & 0 & \text{ sur }  \bordgauche \cup \borddroit
\end{array}\right.\end{equation}
où par convention $\bfu_{p,\ell } =0$ quand $p\in \mathbb{R} \setminus \NN$.\\

Commençons par regarder les premiers termes pour donner une idée des mécanismes et difficultés de cette partie. %
L'hypothèse de raccord relie les comportements asymptotiques de $\sum _{p,\ell } \varepsilon ^{p} \ln ^{\ell } \!\varepsilon \, \bfu_{p,\ell }$ quand $\bfr\rightarrow 0$ et de $\sum _{p,\ell } \varepsilon ^{p} \ln ^{\ell } \!\varepsilon \, S_{p,\ell }$ quand $\bfr\rightarrow \infty $. Nous pouvons donc commencer par regarder l'asymptotique de $\bfu_{p,\ell }$ lorsque $\bfr\rightarrow 0$. Concentrons-nous sur $\Omega $, car c'est la partie la plus intéressante. Tout d'abord, d'après \eqref{eq: eq vol de bfu p l fr}, $u_{0,0}$ satisfait $\mu _{0} \Delta u_{0,0} + \omega ^2  \rho _{0} u_{0,0} = 0$ au voisinage du coin dans $\Omega $ avec une condition de Dirichlet homogène sur $\partial \Omega $. Ainsi, par séparation de variables, il est facile de montrer qu'on a dans un voisinage du coin :
\begin{equation}{\label{eq: decomposition de u 00 en fcts de Bessel fr}}
	u_{0,0} (r,\theta ) = \sum _{d\in \TN} s_{d} (\bfu_{0,0} ) \,J_{d} (k_{0} r) \.\sin (d\theta ) = \sum _{d\in \TN} s_{d} (\bfu_{0,0} ) \,\Bigg(\sum_{n=0}^\infty  j_{d,n} \.r^{d+2n} \Bigg) \sin (d\theta )
\end{equation}
avec $s_{d} (\bfu_{0,0} )$ et $j_{d,n}$ des coefficients dans $\mathbb{C}$, $J_{d}$ les fonctions de Bessel de première espèce et $k_{0} := \omega  \sqrt{\rho _{0} /\mu _{0} }$. Par ailleurs, nous verrons que \eqref{eq: eq vol de bfu p l fr} implique que $u_{1,0}$ satisfait l'équation de Helmholtz avec la condition $\mu _{1} u_{1,0} = \mu _{0} \.\partial _{y} u_{0,0|y=0^{+}}$ sur $\Gamma $. En utilisant \eqref{eq: decomposition de u 00 en fcts de Bessel fr}, on peut montrer qu'il existe des fonctions $f_{q,i}$ et des coefficients $s_{d} (\bfu_{1,0} )$ tels qu'on a proche du coin :
\[u_{1,0} = \sum _{d\in \NN} r^{d-1} (f_{d,0} (\theta )+\ln (r)\.f_{d,1} (\theta ))+ \sum _{d\in \TN} s_{d} (\bfu_{1,0} ) \,J_{d} (k_{0} r) \.\sin (d\theta ).\]
La première somme ci-dessus est une solution particulière de l'équation de Helmholtz dont la trace est $\frac{\mu _{0} }{\mu _{1}} \,\partial _{y} u_{0|y=0^{+}}$ sur $\Gamma $, et la seconde est une solution homogène.\\

Plus généralement, pour tout $(p,\ell ) \in \NN\times \mathbb{N}$, le comportement de $\bfu_{p,\ell }$ quand $\bfr\rightarrow 0$ (resp.\ $S_{p,\ell }$ quand $\bfr\rightarrow \infty $) a la forme suivante :
\[\bfu_{p,\ell } \underset{\bfr\rightarrow 0}\approx \left\{\begin{array}{@{\;}l@{\;}l@{}}
	\displaystyle \sum_{n=0}^\infty  r^{d_{n}} \ln ^{j_{n}} \!r \cdot  g^{0} _{n,1} (\theta ) & \text{ dans }  \Omega  \\[2pt]
	\displaystyle \sum_{n=0}^\infty  r^{d_{n}} \ln ^{j_{n}} \!r \cdot  g^{0} _{n,2} (Y) & \text{ dans }  \couche,
\end{array}\right.  \qquad  \text{resp.}  \qquad  
S_{p,\ell } \underset{\bfr\rightarrow \infty }\approx \left\{\begin{array}{@{\;}l@{\;}l@{}}
	\displaystyle \sum_{n=0}^\infty  r^{d'_{n}} \ln ^{j'_{n}} \!r \cdot  g^{\infty } _{n,1} (\theta ) & \text{ dans }  \Omega  \\[2pt]
	\displaystyle \sum_{n=0}^\infty  r^{d'_{n}} \ln ^{j'_{n}} \!r \cdot  g^{\infty } _{n,2} (Y) & \text{ dans }  \couche
\end{array}\right.\]
où $(d_{n} )$ est une suite croissante de $\mathbb{Z}+\TZ$ (resp., $(d'_{n} )$ est une suite décroissante de $\mathbb{Z}+\TZ$) et $j_{n} , j'_{n} \in \mathbb{N}$. Pour l'instant, il n'est pas du tout évident que les champs lointains concaténés et les champs de coin peuvent se raccorder. Cela dépend des fonctions $g^{0} _{n,i}$ et $g^{\infty } _{n,i}$, car si ces fonctions étaient trop différentes, le raccordement ne serait pas possible. Il faut donc décrire les fonctions $g^{0} _{n,i}$ et $g^{\infty } _{n,i}$. Dans la littérature, elles sont généralement construites par récurrence, en résolvant des problèmes unidimensionnels en $\theta $ et $Y$ \cite{CalCosDauVia06, DelSchSem16, SemDelSch18}. Cependant, en raison de l'inhomogénéité de l'opérateur de Helmholtz $\Delta +k_{0} ^2 $, les $g_{n,i}$ dans notre cas ne dépendent pas les unes des autres uniquement via des suites linéaires, mais plutôt via un arbre de dépendances (cf. remarque~\ref{rq: utilite des operateurs R}).\\

Pour trouver un moyen de gérer cette complexité, nous sommes revenus au cas plus simple d'une couche périodique infinie \cite{AmmHe97, AmmHe98, AmmLat98}. Dans cette situation, la zone proche correspond à la couche et à son voisinage, et les comportements des champs lointains à proximité de la couche s'expriment simplement par des développements limités par rapport à $y$, la distance à la couche. De même, les comportements à l'infini des champs proches sont aussi polynomiaux par rapport à $y$. Ainsi, le comportement de tout champ peut être décomposé sur une base très simple, donnée par les puissances $y^{n}$ avec $n\in \mathbb{N}$, et l'hypothèse de raccord se réduit au raccord des composantes des comportements des champs sur chaque $y^{n}$.\\

Afin de tenter de reproduire la simplicité de ce cas dans notre situation bien plus complexe, la première idée originale de cette thèse est de décomposer les comportements des champs sur une \emph{base de fonctions} bien choisie. Ainsi, seuls les coefficients des décompositions sur cette base doivent être raccordés. De plus, nous donnons des expressions explicites de ces fonctions de base (définition~\ref{def: espaces A}), de la forme :
\begin{equation}{\label{eq: def A dans intro fr}}
	\left\{\begin{array}{ll}
	\Im[(\alpha z)^{q} \,\overline{\alpha z}^{k} \log (\alpha z)^{j} ] & \text{ dans }  \Omega ,\; \text{ avec }  z:=x+\mathrm{i}y \text{ et } \alpha :=e^{-\mathrm{i}\Theta } \\[2pt]
	x^{d} \,Q(\ln x,Y) & \text{ dans }  \couche,\; \text{ avec }  Q \in \mathbb{C}[T,Y].
	\end{array}\right.
\end{equation}
Cette expression est loin d'être évidente et sera justifiée dans la section~\ref{sec: singularites quasi homogenes}. Nous notons $\mathcal{A}(\Omega )$, resp. $\mathcal{A}(\couche)$, l'espace vectoriel engendré par les fonctions de base sur $\Omega $, resp. $\couche$. Ce que nous appelons \og singularités de coin \fg{} (ou simplement \og singularités \fg{}) sont les éléments de ces espaces. Par exemple, $\mathcal{A}(\Omega )$ contient les singularités de coin $r^{d+2n} \sin(d\theta )$ déjà mentionnées en \eqref{eq: decomposition de u 00 en fcts de Bessel fr} car, pour tout $(d,n)\in  \TZ* \times \mathbb{N}$, on a $r^{d+2n} \sin(d\theta ) = \pm  \Im[r^{d+n} e^{\mathrm{i}(d+n)(\theta -\Theta )} \cdot  r^{n} e^{-\mathrm{i}n(\theta -\Theta )} ] = \pm \Im[(\alpha z)^{d+n} (\alpha z)^{n} ]$ avec $z=re^{\mathrm{i}\theta }$. L'équation \eqref{eq: def A dans intro fr} généralise une formule \og $\Re[(-z)^{q} \log(-z)^{j} ]$ \fg{} qui est brièvement mentionnée dans \cite{CosDau96} pour un cas plus simple sans couche, avec $\Theta  = \pi $ et pour l'équation de Poisson.\\

L'équation \eqref{eq: def A dans intro fr} est une définition appropriée des fonctions de base, car elle permet de définir certains \emph{opérateurs} clés sur les espaces~$\mathcal{A}$, notés $\ROmegaPi$, $\RGammaPi$ et $\RcouchePi$ (proposition~\ref{operateurs R}). Ces opérateurs servent à décrire comment les singularités sont liées entre elles et constituent la deuxième idée originale de cette partie. Ils sont définis à partir de la solution d'un problème canonique de type Poisson dans $\Pi $, donné en \eqref{eq: pb resolu dans A}. Ils peuvent être composés, sommés et factorisés de toutes sortes de façons, ce qui, contrairement aux suites linéaires de fonctions, offre une flexibilité suffisante pour décrire toutes les singularités. Un~résultat important est que toutes les singularités qui apparaissent dans le comportement des champs sont engendrées, via ces opérateurs, par les singularités fondamentales $\phi_{d} := r^{d} \sin(d \theta )$, $d\in \TZ*$ (qui font partie de la base). Plus précisément, le théorème~\ref{champs totaux en fct de sigma} montre les décompositions suivantes (écrites entre guillemets pour l'instant, car nous n'avons pas encore défini certains des objets impliqués).

\begin{equation}{\label{eq: champs totaux en fct de sigma en intro fr}}
	\fboxsep=4pt
	\fbox{$ \begin{array}{@{\,}l@{\,}}
	\text{``\,$\displaystyle \sum _{p,\ell } \powelne \,\bfuDA_{p,\ell } = \sum_{n=0}^\infty  R_{\varepsilon } ^{n} \bigg(\sum _{p,\ell } \powelne\, \sum _{\mathclap{d\in \TZ*}} \sigma _{d} (\bfuDA_{p,\ell } )\,\phi_{d} \bigg) $\,''}  \vspace{\jot}\\
	\text{``\,$\displaystyle \sum _{p,\ell } \powelne \,\SDA_{p,\ell } = \sum_{n=0}^\infty  \widetilde{R}_{\varepsilon } ^{n} \bigg(\sum _{p,\ell } \powelne\, \sum _{\mathclap{d\in \TZ*}} \sigma _{d} (\SDA_{p,\ell } )\,\phi_{d} \bigg) $\,''}  
	\end{array} $}
\end{equation}
Ici, $\bfuDA_{p,\ell }$ (resp. $\SDA_{p,\ell }$) est un objet qui décrit le comportement de $\bfu_{p,\ell }$ quand $\bfr\rightarrow 0$ (resp. de $S_{p,\ell }$ quand $\bfr\rightarrow \infty $) et qui est présenté ci-dessous. Les nombres $\sigma _{d} (\bfuDA_{p,\ell } )$ et $\sigma _{d} (\SDA_{p,\ell } )$ désignent les coefficients des comportements des champs sur les fonctions de base $\phi_{d}$. De plus, $R_{\varepsilon }$ et $\widetilde{R} _{\varepsilon }$ sont des opérateurs agissant sur les singularités de coin, qui sont définis à partir des opérateurs $\ROmegaPi$, $\RGammaPi$ et $\RcouchePi$, et de puissances de $\varepsilon $. L'équation \eqref{eq: champs totaux en fct de sigma en intro fr} implique que, pour que les champs lointains concaténés se raccordent avec les champs de coin, il suffit de fixer correctement les coefficients $\sigma _{d} (\bfuDA_{p,\ell } )$ et $\sigma _{d} (\SDA_{p,\ell } )$. Plus précisément, seuls les coefficients apparaissant dans les comportements non variationnels des champs doivent être imposés, c'est-à-dire les $\sigma _{d} (\bfuDA_{p,\ell } )$ avec $d<0$ et les $\sigma _{d} (\SDA_{p,\ell } )$ avec $d>0$.\\

La troisième idée principale consiste à utiliser des \emph{séries formelles} pour se concentrer sur la nature calculatoire du raccordement, sans avoir à traiter des troncatures de sommes et des restes asymptotiques tels que $o(r^{d} )$. Les séries formelles sont généralement utilisées en analyse asymptotique de manière intuitive, comme nous l'avons fait dans la section~\ref{sec: ansatz du dvlpmt raccorde fr}, ou en tant que simple notation (voir \cite{Won01}). Cependant, nous fournissons ici des fondements rigoureux aux séries formelles utilisées afin de s'en servir comme outil de démonstration et de justifier des opérations non triviales sur celles-ci, telles que l'application d'opérateurs, voire de séries formelles d'opérateurs comme dans \eqref{eq: champs totaux en fct de sigma en intro fr}. Par exemple, pour tout $(p, \ell )$, l'objet $\bfuDA_{p,\ell }$ mentionné ci-dessus est une série formelle qui décrit le développement asymptotique de $\bfu_{p,\ell }$ lorsque $\bfr \rightarrow 0$. Elle est grossièrement définie comme la suite $(v_{n} )_{n\in \mathbb{N}}$ des termes de ce développement. C'est-à-dire que, pour tout $n\in \mathbb{N}$, $v_{n}$ se comporte comme $\bfr^{d_{n}}$ pour une certaine suite croissante $(d_{n} )$, et on a $\bfu_{p,\ell } = \sum_{k=0}^n v_{k} + o(\bfr^{d_{n}} )$. De même pour $\SDA_{p,\ell }$ et $S_{p,\ell }$ lorsque $\bfr \rightarrow \infty $. La définition exacte des séries formelles est un peu moins directe, mais plus souple, ce~qui permet de sommer des séries et d'y appliquer des opérateurs.\\

Nous définissons plusieurs espaces de séries formelles, dans lesquels nous démontrons l'existence et l'unicité de la solution de certains problèmes de Poisson et de Helmholtz, et donnons une expression de la solution (propositions \ref{pb de Helmholtz dans A Pi} à~\ref{pb de Poisson dans A Pi}, \ref{pb de Helmholtz dans Aemacro} et~\ref{pb de Helmholtz dans Aemicro}). Cela nous permet d'exprimer le comportement des champs par rapport à\ $\bfr$ et de prouver \eqref{eq: champs totaux en fct de sigma en intro fr}, comme conséquence des équations vérifiées par \og $\sum _{p,\ell } \powelne \,\bfu_{p,\ell }$ \fg{} et \og $\sum _{p,\ell } \powelne \,S_{p,\ell }$ \fg{}.\\

L'introduction de la section~\ref{sec: C[[X]]} dresse un parallèle entre les séries formelles et les distributions. Toutes deux constituent des généralisations d'objets classiques (les séries convergentes, resp. les fonctions $L^{1} _{\mathrm{loc}}$) qui offrent davantage de souplesse pour certaines opérations (par exemple la dérivation, resp. les sommes infinies), permettant ainsi de résoudre des problèmes qui n'ont pas de solution bien définie dans le cadre classique. Notons que \cite{Fao02} a également utilisé une approche par séries formelles pour étudier un problème impliquant une couche mince, mais avec des outils algébriques plus élémentaires.\\

Ces trois idées --- la base explicite des singularités, les opérateurs sur les singularités et les séries formelles --- forment un ensemble d'outils de calcul puissants qui permettent de faire face à la grande complexité du raccordement. Ceci nous permet de prouver le théorème suivant, qui fournit les conditions de raccord et constitue le résultat principal du chapitre~\ref{chap: raccord algebrique}. L'une des originalités de notre méthode réside dans le fait que ces conditions sont obtenues à tous les ordres en même temps, contrairement au processus étape par étape souvent utilisé dans les problèmes de coin \cite{CalCosDauVia06,BenMakTor11,SemDelSch18,CosDau96}.
\end{textesansboite}

\begin{boitecarree}{Théorème \ref{condition de raccord}:}
\ Il existe deux familles de constantes $(\cmacro_{d,d',p,\ell } )$ et $(\cmicro_{d,d',p,\ell } )$ telles que l'hypothèse de raccord de \eqref{eq: hypothese de raccord fr} est équivalente au système d'équations suivant appelé \emph{conditions de raccord} :\vspace{-1mm}
\[\!\left\{\begin{array}{ll}
	\! \forall (d,p,\ell )\in (-\TN)\times \NN\times \mathbb{N},
	&\displaystyle \sigma _{d} (\bfuDA_{p,\ell } ) = \! \sum _{p'\in \intvl{0,p+d}} \sum _{\substack{d'\in \TZ*\\ p'-d'\in \intvl{0,p}}} \! \sum _{\ell '=0} ^{\ell } \,\cmacro_{d,d',p-p',\ell -\ell '} \cdot  \sigma _{d'} (\SDA_{p',\ell '} ) \\[-1mm]
	\! \forall (d,p,\ell )\in \TN\times \NN\times \mathbb{N},
	&\displaystyle \sigma _{d} (\SDA_{p,\ell } ) = \! \sum _{p'\in \intvl{0,p-d}} \sum _{\substack{d'\in \TZ*\\ p'+d'\in \intvl{0,p}}} \! \sum _{\ell '=0} ^{\ell } \,\cmicro_{d,d',p-p',\ell -\ell '} \cdot  \sigma _{d'} (\bfuDA_{p',\ell '} )
\end{array}\right.\]
avec $\intvl{0,a}:=\{b\in [0,a] \mid b\in \NN \text{ et } a-b\in \mathbb{N}\}$ pour tout $a\in \mathbb{R}$.
\end{boitecarree}

\begin{textesansboite}
Nous appelons \emph{coefficients de raccord} (\og matching coefficients \fg) les nombres $\cmacro_{d,d',p,\ell }$ et $\cmicro_{d,d',p,\ell }$. Ils ne dépendent que de $\Theta $ et de $(\mu _{0} ,\mu _{1} ,\rho _{0} ,\rho _{1} , \omega )$ ; en particulier ils ne dépendent ni du terme source $\fsource$, ni de la perturbation locale de la géométrie au niveau du coin. De plus, ces coefficients peuvent être obtenus par des calculs sur les singularités de coin impliquant les opérateurs $\ROmegaPi$, $\RGammaPi$ et $\RcouchePi$. Grâce à la connaissance d'une base des singularités et en utilisant le fait que les opérateurs $\ROmegaPi$, $\RGammaPi$ et $\RcouchePi$ ont des expressions explicites sur cette base, nous avons développé un \emph{module de calcul formel} spécialement conçu pour traiter ces calculs (section~\ref{sec: implementation des espaces A}). Il permet de calculer les coefficients de raccord, ainsi que les singularités de coin apparaissant dans le comportement des champs, de manière exacte et très rapide (généralement en une fraction de seconde). Ce dernier point constitue une contribution originale par rapport à la littérature existante.\\

Notons que les conditions de raccord du théorème~\ref{condition de raccord} possèdent une structure convolutive par rapport à $p$ et $\ell $, ce qui sera important pour les modèles effectifs. De plus, la complexité des conditions de raccord est en partie due au fait que nous considérons l'équation de Helmholtz. Pour l'équation de Poisson, nous aurions des sommes doubles au lieu de triples (car l'indice $d'$ de la deuxième somme peut être exprimé en fonction de l'indice $p'$).\\

Le chapitre~\ref{chap: raccord algebrique} constitue l'une des parties les plus originales de cette thèse. Elle permet non seulement de déduire les conditions de raccord, mais elle met aussi en place des outils qui seront utilisés tout au long des chapitres suivants. De plus, notre méthode repose sur des calculs abstraits, non spécifiques à un cas particulier, ce qui la rend généralisable à d'autres situations. Un aperçu de ces extensions est donné dans la section~\ref{sec: resume du raccord} et de futurs travaux exploreront cela plus en détail.\\

Notre approche utilise un type de série formelle que nous n'avons pas trouvé dans la littérature et qui est présenté dans le chapitre~\ref{chap: SF}. Il contient un indice de sommation pouvant être n'importe quel nombre réel, afin de traiter les puissances non entières de $\varepsilon $ et de $r$ qui apparaissent dans le développement. Le chapitre~\ref{chap: SF} définit principalement des notations, telles que le symbole de sommation formelle $\fsum$ (un $\sum $ gras), et fournit des règles de calcul sur les sommes formelles infinies autorisées ou non, les propriétés de commutation, ainsi que des théorèmes pour l'inversion d'opérateurs (utilisés pour résoudre des EDP parmi les séries formelles). Il s'agit d'une introduction progressive destinée à familiariser le lecteur avec les manipulations largement utilisées dans les chapitres suivants.
\end{textesansboite}

\subsection{Chapitre \ref{chap: construction du DA}: construction du développement}

\begin{textesansboite}
Dans le chapitre~\ref{chap: construction du DA}, nous construisons les termes du développement asymptotique de $u_{\varepsilon }$ et justifions ce~développement par des estimations d'erreur.\\

Tout d'abord, nous exprimons les champs de couche $U_{p,\ell }$ en fonction des champs lointains (lemme~\ref{expression de Upl}). Ensuite, nous introduisons des problèmes modèles pour construire chaque champ lointain et chaque champ de coin, et nous montrons qu'ils possèdent une unique solution. Par exemple, pour les champs lointains, ceci est donné par le théorème ci-dessous. Ces problèmes sont non variationnels car les champs présentent des comportements singuliers (au coin pour les champs lointains et à l'infini pour les champs de coin).
\end{textesansboite}

\begin{boitecarree}{Théorème \ref{cadre fonctionnel des champs lointains}:}
\ Notons $\mathcal{A}(\Gamma ) := \{\varphi _{|\Gamma } \mid \varphi \in \mathcal{A}(\Omega )\}$, où nous rappelons que $\mathcal{A}(\Omega )$ est l'espace des singularités de coin sur $\Omega $ introduit dans le chapitre~\ref{chap: raccord algebrique}. Soit $\khiz$ une fonction de troncature régulière localisée autour du coin. Soient $f\in (H^{1} (\Omega ))'$, $g\in H^{1/2} _{00} (\Gamma )+\khiz \mathcal{A}(\Gamma )$ et $(s_{d} )\in \mathbb{C}^{-\TN}$ à support fini. Le système suivant possède une unique solution dans $H^{1} (\Omega )+\khiz \mathcal{A}(\Omega )$.
\begin{equation}{\label{eq: pb modele de CL dans intro fr}}
	\left\{\begin{array}{r@{\;=\;}ll}
	\mu _{0} \Delta  u +\omega ^2 \rho _{0} u & f & \text{ dans }  \Omega  \\[0.5mm]
	u & g & \text{ sur }  \Gamma  \\[0.5mm]
	u & 0 & \text{ sur }  \bordgauche\\[0.5mm]
	\sigma _{d} (u)& s_{d} & \ \forall d\in -\TN
	\end{array}\right.
\end{equation}\fbeq
\end{boitecarree}

\begin{textesansboite}
La composante de $u$ dans $\khiz \mathcal{A}(\Omega )$ contient son comportement singulier. Par ailleurs, le théorème~\ref{cadre fonctionnel pour les champs de coin} fournit un résultat similaire pour les champs de coin, en remplaçant l'équation de Helmholtz par celle de Poisson et $d\in -\TN$ par $d\in \TN$.\\

Ces théorèmes montrent que les coefficients $\sigma _{d} (\bfuDA_{p,\ell } )$ avec $d\in -\TN$ et $\sigma _{d} (\SDA_{p,\ell } )$ avec $d\in \TN$ peuvent être imposés (par exemple, par la dernière ligne de \eqref{eq: pb modele de CL dans intro fr}). Pour les autres valeurs de $d$, ces coefficients, qui apparaissent au membre de droite des conditions de raccord, appartiennent à la partie variationnelle des champs et ne peuvent pas être fixés. Voir le tableau~\ref{tab: role de sigma d pour CL et CP fr}. Les comportements variationnels sont étudiés en démontrant des développements asymptotiques par rapport à $r$ à tout ordre (théorèmes \ref{DA des champs lointains concatenes} et~\ref{DA des champs de coin}). On utilise pour cela la théorie de Kondrat'ev, qui fait appel à des espaces à poids et à la transformée de Laplace par rapport à $\ln r$ (i.e. la transformée de Mellin par rapport à $r$). Ceci est classique, à ceci près que l'on exprime le développement à l'aide des singularités de coin dans les espaces~$\mathcal{A}$.
\end{textesansboite}

\begin{table}[h]
\renewcommand{\arraystretch}{1.4}
\definecolor{tabvaria}{rgb}{0.83, 1, 0.83} 
\definecolor{tabsing}{rgb}{1, 0.85, 0.85} 
\caption{Résumé du caractère variationnel ou non variationnel de $\phi_{d}$ et du rôle des coefficients $\sigma _{d}$ pour les différents champs}
\label{tab: role de sigma d pour CL et CP fr}
\begin{tblr}{Q[c,m]Q[c,m]Q[c,m]}
\cline{2-3}
& Champs lointains ($\bfr\rightarrow 0$) & Champs de coin ($\bfr\rightarrow \infty $) \\
\hline
$d\in \TN$ 
& \SetCell{tabvaria} {$\phi_{d}$ variationnel,\\ $\sigma _{d} (\bfuDA_{p,\ell } )$ calculé à partir de $u_{p,\ell }$ \footnotesize (chap.\,\ref{chap: validation numerique du DA})}
& \SetCell{tabsing} {$\phi_{d}$ non variationnel,\\ $\sigma _{d} (\SDA_{p,\ell } )$ fixé par une condition de raccord} \\
\hline
$d\in -\TN$ 
& \SetCell{tabsing} {$\phi_{d}$ non variationnel,\\ $\sigma _{d} (\bfuDA_{p,\ell } )$  fixé par une condition de raccord}
& \SetCell{tabvaria} {$\phi_{d}$ variationnel,\\ \vphantom{$\bfuDA$}$\sigma _{d} (\SDA_{p,\ell } )$ calculé à partir de $S_{p,\ell }$ \footnotesize (chap.\,\ref{chap: validation numerique du DA})} \\
\hline
\end{tblr}
\end{table}

\begin{textesansboite}
En rassemblant toutes les équations que les champs doivent vérifier (c'est-à-dire \ \eqref{eq: eq vol champs lointains fr} à~\eqref{eq: eq vol champs de coin fr} et les conditions de raccord), et en utilisant les résultats de ce chapitre, nous prouvons l'existence et l'unicité des champs $u_{p,\ell }$, $U_{p,\ell }$ et $S_{p,\ell }$ (théorème~\ref{validation construction et DA des champs}). De plus, nous fournissons une méthode permettant de construire directement les champs lointains, en éliminant les champs de couche et de coin, grâce au théorème suivant.
\end{textesansboite}

\begin{boitecarree}{Théorème \ref{construction directe des upl}:}
\ Soit $(T_{n} )_{n\in \mathbb{N}}$ la suite des coefficients de Taylor de la tangente. $(u_{p,\ell } )_{p,\ell }$ est l'unique famille dans $H^{1} (\Omega )+\khiz \mathcal{A}(\Omega )$ telle que pour tout $(p,\ell )\in \NN\times \mathbb{N}$ nous avons:
\begin{equation}{\label{eq: construction directe des upl en intro fr}}
	\!\!\left\{\begin{array}{r@{\;=\;}ll@{}}
	\mu _{0} \Delta  u_{p,\ell } +\omega ^2 \rho _{0} u_{p,\ell } & \fsource \,\delta _{p,0} \,\delta _{\ell ,0} & \text{ dans }  \Omega  \\[0.5mm]
	u_{p,\ell } & 0 & \text{ sur }  \bordgauche\\[0.5mm]
	u_{p,\ell } &\displaystyle \sum_{n=0}^\infty  \frac{\mu _{0} }{\mu _{1}} \,T_{n} \cdot  (\partial _{x} ^2 +k_{1} ^2 )^{n} \partial _{y} u_{p-1-2n,\ell } & \text{ sur }  \Gamma  \\[-1.5mm]
	\sigma _{d} (u_{p,\ell } )&\displaystyle \!\! \sum _{\substack{p'\in \NN\\ p-p'\in \NN+ \frac{2\pi }{\Theta } }} \; \sum _{\substack{d'\in \TZ*\\ p'+d'\in \intvl{0,p+d}}} \! \sum _{\ell '=0} ^{\ell } \cfinal_{d,d',p-p',\ell -\ell '} \cdot \sigma _{d'} (\bfuDA_{p',\ell '} ) & \ \forall d\in -\TN
	\end{array}\right.\!\!\!
\end{equation}\fbeq
\end{boitecarree}

\begin{textesansboite}
Dans \eqref{eq: construction directe des upl en intro fr}, les coefficients $\cfinal$ sont définis explicitement à partir des coefficients de raccord $\cmacro$ et $\cmicro$ figurant dans le théorème~\ref{condition de raccord}, ainsi qu'à partir de nouveaux coefficients $\sigma _{q} (\ccDA_{d,n} )$ avec $q\in -\TN$, voir \eqref{eq: def cfinal}. Les fonctions $\mathcal{S}_{d,n}$ forment une base des champs de coin nommée \emph{correcteurs de coin} (\og corner profiles \fg{} en anglais, cf. définition~\ref{def: correcteurs de coin}) qui ne dépend pas des champs lointains, et les $\sigma _{q} (\ccDA_{d,n} )$ décrivent leur comportement variationnel à l'infini.\\

Enfin, le développement asymptotique est rigoureusement justifié par l'estimation d'erreur optimale suivante.
\end{textesansboite}

\begin{boitecarree}{Théorème \ref{DA de SE avec CL et delta fixe}:}
\ Pour tout $\ptronc \geqslant 0$ et $\delta >0$, nous avons:
\[\bigg\|u_{\varepsilon } -\sum _{p\in \NN\mkern-1mu,\,p\leqslant \ptronc} \sum _{\ell =0} ^{\ell _{p}} \varepsilon ^{p} \ln ^{\ell } \!\varepsilon \, u_{p,\ell } \bigg\|_{H^{1} (\Omega  \setminus \closedBall(0,\delta ))} = o(\varepsilon ^{\ptronc} ) \qquad  \text{quand }  \varepsilon \rightarrow 0.\]\fbeq
\end{boitecarree}

\begin{textesansboite}
Ceci conclut la construction du développement, qui couvre les chapitres \ref{chap: SF} à~\ref{chap: construction du DA}.
\end{textesansboite}

\subsection{Chapitres \ref{chap: complements sur le DA} et \ref{chap: decompo modale CL}: propriétés additionnelles du développement}

\begin{textesansboite}
Les chapitres \ref{chap: complements sur le DA} et \ref{chap: decompo modale CL} établissent diverses propriétés supplémentaires des champs, en particulier sur les champs lointains. La plus remarquable d'entre elles est une simplification substantielle des conditions de raccord qui se produit dans \og la plupart des cas \fg{} (théorèmes \ref{simplifications cas irrationnel} et~\ref{simplifications cas rationnel}). Nous prouvons que, si $\frac{\Theta }{\pi }$ est irrationnel ou s'il est rationnel et que $p\in \NN$ est inférieur au dénominateur de $\frac{\Theta }{\pi }$ sous forme irréductible, on a :
\begin{itemize}
	\item le développement de $u_{\varepsilon }$ jusqu'à l'ordre $p$ ne contient pas de $\ln \varepsilon $,
    \item les conditions de raccord se simplifient en $\sigma _{d} (\bfuDA_{p,0} )=\sigma _{d} (\SDA_{p+d,0} )$,
    \item la condition de coin du théorème~\ref{construction directe des upl} devient
	\begin{equation}{\label{eq: condition de coin CL simplifiee en intro fr}}
		\sigma _{d} (u_{p,0} )= \sum _{\substack{n\in \mathbb{N}\\ 2n\leqslant p+d}} \sum _{\substack{d'\in \TN\\ d'\leqslant p+d-2n}} \omega ^{2n} \sigma _{d} (\ccDA_{d',n} )\cdot \sigma _{d'} (\bfuDA_{p-d'+d-2n,0} ).
	\end{equation}\fbeq
\end{itemize}
Ceci est généralisé à quelques autres cas particuliers pertinents dans l'annexe~\ref{sec: annexe des complements sur le DA}. Ainsi, à un ordre donné, par exemple~3 ou 5 comme dans les modèles effectifs, seul un petit nombre d'angles n'est pas concerné par ces simplifications, à savoir $\frac{2\pi }{3}$, $\pi $, $\frac{4\pi }{3}$ et $\frac{5\pi }{3}$ à l'ordre~5. Nous soulignons que l'approche du chapitre~\ref{chap: raccord algebrique} ne repose pas sur ces simplifications, ce qui nous permet de traiter n'importe quel angle à tout ordre.\\

Pour prouver ces simplifications, nous exprimons certaines des singularités des champs lointains à l'aide des fonctions de Bessel. Tout d'abord, nous montrons la décomposition \emph{formelle} suivante du comportement des champs lointains au coin (proposition~\ref{expression uDA avec mcl})
\begin{equation}{\label{eq: expression uDA avec mcl en intro fr}}
    \fbox{$ \,\text{``$\.\displaystyle u_{p,\ell } \, \underset{r\rightarrow 0}\approx \, \sum _{n=0} ^{\lfloor {p}\rfloor } \.\sum _{d\in \TZ*\!,\,d\geqslant n-p} \sigma _{d} (\bfuDA_{p-n,\ell } ) \,\mcla_{d,n} \.$''} \, $}
\end{equation}
où les $\mcla_{d,n}$ sont des séries formelles de singularités de coin qui satisfont formellement, pour tout $d\in \TZ*$,
\begin{equation}{\label{eq: pb verifie par mcl qd n=0 en intro fr}}
    \left\{\begin{array}{r@{\;=\;}ll}
    (\Delta +k_{0} ^2 )\mcla_{d,0} & 0 & \text{ dans }  \Omega  \\[3pt]
    \mcla_{d,0} & 0 & \text{ sur }  \Gamma  \\[3pt]
    \sigma _{q} (\mcla_{d,0} ) & \delta _{d,q} &\ \forall q\in \TZ*
    \end{array}\right.
\end{equation}
et, pour tout $n\in \mathbb{N}^*$,
\begin{equation}{\label{eq: pb verifie par mcl qd n>0 en intro fr}}
    \left\{\begin{array}{r@{\;=\;}ll}
    (\Delta +k_{0} ^2 )\mcla_{d,n} & 0 & \text{ dans }  \Omega  \\[2pt]
    \mcla_{d,n} & \sum\limits_{m=0}^\infty  \frac{\mu _{0} }{\mu _{1}} \,T_{m} \cdot  (\partial _{x} ^2 +k_{1} ^2 )^{m} \partial _{y} \mcla_{d,n-1-2m} & \text{ sur }  \Gamma  \vspace{2pt}\\
    \sigma _{q} (\mcla_{d,n} ) & 0 &\ \forall q\in \TZ*
    \end{array}\right.
\end{equation}
Nous montrons que, dans de nombreux cas, la série $\mcla_{d,n}$ converge vers une combinaison linéaire finie de fonctions de Bessel (proposition~\ref{expression sing heterog CL}). Par exemple, lorsque $d\not\in -\mathbb{N}$, $\mcla_{d,0}$ converge vers $J_{d} (k_{0} r) \sin(d\theta )$ multiplié par une constante, où $J_{d}$ est une fonction de Bessel de première espèce. Ceci est cohérent avec \eqref{eq: decomposition de u 00 en fcts de Bessel fr}. Comme les fonctions de Bessel de première espèce n'impliquent pas de $\ln r$, il en va de même pour les champs lointains dans la plupart des cas, ce qui implique l'absence de $\ln \varepsilon $ ainsi que les autres simplifications \og $\sigma _{d} (\bfuDA_{p,0} )=\sigma _{d} (\SDA_{p+d,0} )$ \fg{} et \eqref{eq: condition de coin CL simplifiee en intro fr}. Notez que les simplifications lorsque $\pst \not\in \mathbb{Q}$ étaient déjà connues pour l'équation de Poisson \cite{CalCosDauVia06}, mais elles sont nettement moins évidentes à prouver pour Helmholtz.\\

Dans la section~\ref{sec: exemples de structures de DA}, nous écrivons explicitement les conditions de raccord pour quelques exemples d'angles, et nous représentons graphiquement le comportement, en fonction de l'ordre du développement, de quelques grandeurs pertinentes telles que le nombre de champs lointains. Ces données ont été générées par notre module de calcul formel mentionné précédemment. Elles illustrent les simplifications annoncées ainsi que la complexité des autres cas. Nous utilisons la notion de \emph{généricité} pour caractériser les champs ou les coefficients qui ne s'annulent pas sauf pour des valeurs exceptionnelles de $\omega $, $\mu _{0}$, $\mu _{1}$, etc. (définition~\ref{def: genericite}).\\

Ensuite, dans le chapitre~\ref{chap: decompo modale CL}, nous prouvons que la série formelle $\mcla_{d,n}$ converge toujours, même si ce n'est pas vers des fonctions de Bessel (théorème~\ref{convergence des mcl}). Notons $\mcll_{d,n}$ la limite. Nous montrons que la décomposition modale des champs lointains dans \eqref{eq: expression uDA avec mcl en intro fr} est en fait aussi valable au sens classique, c'est-à-dire qu'elle converge.
\end{textesansboite}

\begin{boitecarree}{Théorème \ref{decompo modale des CL}:}
\ Soit $\Rsource:=\operatorname{dist}(0,\operatorname{supp}\fsource)>0$. Soient $p\in \NN$ et $\ell \in \mathbb{N}$. Dans $\Omega  \cap  B(0,\Rsource)$, la série ci-dessous converge et l'égalité suivante est vérifiée:
\[u_{p,\ell } \,=\, \sum _{n=0} ^{\lfloor {p}\rfloor } \sum _{\substack{d\in \TZ*\\ d\geqslant -p+n}} \sigma _{d} (\bfuDA_{p-n,\ell } ) \,\mcll_{d,n} .\]
De plus, pour tout $R\in (0,\Rsource)$, la convergence est uniforme sur $\Omega  \cap  B(0,R)$ et il existe une constante $C_{R} >0$ telle que, pour tout $N\in \mathbb{N}$, $r\in {]0,R[}$ et $\theta \in [0,\Theta ]$, on a
\begin{equation}{\label{eq: vitesse de cvg decompo modale des CL en intro fr}}
	\Bigg| u_{p,\ell } - \sum _{n=0} ^{\lfloor {p}\rfloor } \sum _{\substack{d\in \TZ*\\ -p+n\leqslant d\leqslant N\pst}} \sigma _{d} (\uDA_{p-n,\ell } ) \,\mcll_{d,n} \Bigg| \leqslant  C_{R} \left( \frac{r}{R} \right)^{(N+1)\pst-\lfloor {p}\rfloor } (1+|{\ln r}|)^{2p+1} 
\end{equation}
et de manière similaire pour le gradient.
\end{boitecarree}

\begin{textesansboite}
De plus, nous démontrons que l'estimation d'erreur \eqref{eq: vitesse de cvg decompo modale des CL en intro fr} est préservée lorsque l'on remplace $\mcll_{d,n}$ par une somme partielle de la série $\mcla_{d,n}$ contenant suffisamment de termes (théorème~\ref{decompo en modes tronques des CL}). Cette décomposition modale permet de justifier le calcul numérique des champs lointains dans le chapitre~\ref{chap: validation numerique du DA}, et est essentielle pour les modèles effectifs du chapitre~\ref{chap: ME}.\\

Pour démontrer les résultats du chapitre~\ref{chap: decompo modale CL}, en particulier le théorème~\ref{decompo modale des CL} ci-dessus, nous introduisons une norme sur les singularités de coin qui est égale à la somme des coefficients sur la base des singularités (définition~\ref{def: norme sur A_d}). À l'aide de cette norme, nous définissons ensuite des normes à poids sur les séries formelles de singularités, où le poids dépend du degré d'homogénéité des singularités (c'est-à-dire de la puissance de $\bfr$), cf. définition~\ref{def: norme sur Amacro}. Elles sont telles qu'une série formelle a une norme finie si et seulement si elle converge, et le poids donne la vitesse de convergence. Nous montrons que les opérateurs sur les singularités introduits dans le chapitre~\ref{chap: raccord algebrique} sont continus pour ces normes sur des sous-espaces bien choisis de singularités (lemme~\ref{norme des op R sur series formelles}), ce qui nous permet de majorer la norme de $\mcla_{d,n}$ (proposition~\ref{norme SF de mcl cas rationnel}). Ceci~entraîne que les $\mcla_{d,n}$ convergent, donne une majoration de leur limite, et enfin implique le théorème~\ref{decompo modale des CL} par somme (en bref, car les détails sont assez techniques). À notre connaissance, la méthode et les outils utilisés ici n'existent pas dans la littérature.
\end{textesansboite}

\subsection{Chapitre \ref{chap: validation numerique du DA}: validation numérique}

\begin{textesansboite}
Dans ce chapitre, nous illustrons numériquement la convergence du développement asymptotique à différents ordres et pour plusieurs angles. Ceci nous amène à concevoir des méthodes numériques spécifiques pour calculer les champs, qui restent pertinentes pour les modèles effectifs du chapitre~\ref{chap: ME}. Le calcul du développement comporte généralement quatre parties.
\begin{enumerate}
	\item Le calcul des correcteurs de coin (ou \og corner profiles \fg) $\mathcal{S}_{d,k}$, qui sont définis dans $\Omega _{1}$ et satisfont :
	\[\left\{\begin{array}{r@{\;=\;}l@{\;}l}
		\operatorname{div}(\mu  \nabla  \mathcal{S}_{d,0} ) & 0 & \text{ dans }  \Omega _{1} \\[0.5mm]
		\mathcal{S}_{d,0} & 0 & \text{ sur }  \partial \Omega _{1} \\[0.5mm]
		\sigma _{q} (\mathcal{S}_{d,0} )& \delta _{q,d} & \ \forall q\in \TN
	\end{array}\right. \quad  \text{ et } \quad  \ \forall k\in \mathbb{N}^*,\; \left\{\begin{array}{r@{\;=\;}l@{\;}l}
		\operatorname{div}(\mu  \nabla  \mathcal{S}_{d,k} ) & -\rho \,\mathcal{S}_{d,k-1} & \text{ dans }  \Omega _{1} \\[0.5mm]
		\mathcal{S}_{d,k} & 0 & \text{ sur }  \partial \Omega _{1} \\[0.5mm]
		\sigma _{q} (\mathcal{S}_{d,k} )& 0 & \ \forall q\in \TN
	\end{array}\right.\]
	Des correcteurs de coin, nous ne conservons que les coefficients $\sigma _{q} (\ccDA_{d,k} )$ avec $q\in -\TN$, qui décrivent le comportement variationnel à l'infini de $\mathcal{S}_{d,k}$. Ces nombres contiennent des informations sur la perturbation locale de la géométrie au niveau de l'angle, tandis que nous rappelons que les coefficients de raccord ne la prennent pas en compte.
	\item Le calcul des coefficients de raccord $\cmicro_{d,d',p,\ell }$ et $\cmacro_{d,d',p,\ell }$. Il est indépendant de la partie 1 et effectué par notre module de calcul formel.
    \item Le calcul des coefficients $\cfinal_{d,d',p,\ell }$ à partir des nombres $\sigma _{q} (\ccDA_{d,k} )$ et des coefficients de raccord $\cmacro$ et $\cmicro$.
    \item Enfin, le calcul des champs lointains $u_{p,\ell }$ par récurrence par rapport à $p$, d'après les équations du théorème~\ref{construction directe des upl}. Il faut également calculer certains des coefficients $\sigma _{d} (\uDA_{p,\ell } )$ avec $d\in \TN$, car ils interviennent dans la construction d'autres champs lointains (lorsque $d\in -\TN$, $\sigma _{d} (\uDA_{p,\ell } )$ est déjà connu, cf. tableau~\ref{tab: role de sigma d pour CL et CP fr}).
\end{enumerate}
Seules les étapes 1 et 4 nécessitent des calculs complexes, car elles impliquent la résolution d'une EDP sans solution explicite. De plus, il faut extraire les coefficients asymptotiques des solutions de l'EDP, ce qui requiert des méthodes numériques spécifiques.
\end{textesansboite}

\begin{textesansboite}
Pour les correcteurs de coin, nous introduisons en section~\ref{sec: calcul numerique des CC} une méthode de Galerkine enrichie par ce que l'on appelle des \emph{compléments à l'infini}. Tout d'abord, nous décomposons le correcteur de coin en une partie singulière, constituée de singularités de coin explicitement connues, et une partie variationnelle inconnue $\mathcal{S}^{V} _{d,k}$, qui doit être calculée numériquement. Nous savons que, pour tout $\nbcs \in \mathbb{N}$ donné, on a
\[\mathcal{S}^{V} _{d,k} = \sum _{n=1} ^{\nbcs} \sigma _{-n\pst} (\ccDA_{d,k} ) \,\varphi ^{\nbcs} _{n} +O(r^{-(\nbcs+1)\pst} ) \quad \text{quand } r\rightarrow \infty .\]
Les $\varphi ^{\nbcs} _{n}$ sont des singularités de coin explicites, tandis que les $\sigma _{-n\pst} (\ccDA_{d,k} )$ sont les coefficients que nous recherchons. Nous définissons un sous-domaine borné $\Omega _{L} :=\Omega \cap \{\bfr<L\}$ et une fonction de troncature régulière $\chi _{L}$ localisée dans $\{\bfr>L-1\}$. Nous recherchons une approximation de $\mathcal{S}^{V} _{d,k}$ dans l'espace de Galerkine
\[\big\{u:\Omega _{1} \rightarrow \mathbb{C} \mid u_{|\Omega _{L}} \in H^{1} _{0} (\Omega _{L} ) \text{ et } u_{|\Omega _{1} \setminus \ounl} =0 \big\} \oplus  \operatorname{Vect}(\{ \khiL \varphi ^{\nbcs} _{n} \mid n\in [\![1,\nbcs]\!]\}).\]
Les fonctions $\khiL \varphi ^{\nbcs} _{n}$ sont appelées \emph{compléments à l'infini}. En pratique, la partie dans $H^{1} _{0}$ est traitée par la méthode des éléments finis.\\

Nous montrons que la solution approchée converge vers $\mathcal{S}^{V} _{d,k}$ en $O\big(L^{-(\nbcs+1) \pst + \frac{1}{2} }\big)$ lorsque $L\rightarrow \infty $ (théorème~\ref{estimate of uL-SV}). De plus, nous prouvons que, pour tout $n\in [\![1,\nbcs]\!]$, la coordonnée de la solution approchée sur $\varphi ^{\nbcs} _{n}$ converge vers $\sigma _{-n\pst} (\ccDA_{d,k} )$ en $O\big(L^{(n-\nbcs-1)\pst + \frac{1}{2} }\big)$ (théorème~\ref{convergence des CCC}). Il s'agit donc d'une méthode d'ordre élevé, où plus $N$ est grand, meilleur est le taux de convergence. De plus, $\sigma _{-n\pst} (\ccDA_{d,k} )$ est mieux approché pour les petits $n$ que pour les grands.\\

Des résultats numériques sont présentés dans la section~\ref{sec: validation numerique complements singuliers}. Ils montrent que l'erreur commise sur les coefficients peut être très faible (de l'ordre de $10^{-5}$ à $10^{-10}$ en erreur relative) même pour un domaine $\Omega _{L}$ de taille raisonnable (par exemple, $L=10$).\\

Cette méthode présente des similitudes avec la méthode des fonctions singulières \cite{StrFix73, Loh98, LubNic94, CiaHe03} utilisée pour pallier le manque de régularité des solutions d'EDP elliptiques au voisinage de coins. L'originalité réside ici dans l'utilisation de plusieurs compléments, ce qui nécessite une technique différente pour démontrer la convergence.
\end{textesansboite}

\begin{figure}[h]
\def\widthFig{3.5} 
\def\heightFig{2.9} 
\def\ThetaFig{205} 
\def\RMEFig{1.3} 
\begin{tikzpicture}[scale=0.9]
	\coordinate (finLambda) at (-\widthFig, {-\widthFig*tan(\ThetaFig-180)});
	\coordinate (dbtLambda) at ({\RMEFig*cos(\ThetaFig)}, {\RMEFig*sin(\ThetaFig)});
	\fill [pattern={Lines[angle=45,line width=.2pt]}, pattern color=white!70!gray] (dbtLambda) arc (\ThetaFig:0:\RMEFig) -- (\widthFig,0) -- (\widthFig,\heightFig) -- (-\widthFig,\heightFig) -- (finLambda) -- cycle;
	\draw [thick] (finLambda) -- (dbtLambda) arc (\ThetaFig:0:\RMEFig) -- (\widthFig,0);
	\draw [dashed] (dbtLambda) -- (0,0);
	\draw [>=stealth,->] (-\widthFig-0.4,0) -- (\widthFig+0.4,0) node [right] {$x$};
	\draw [>=stealth,->] (0,{-\widthFig*tan(\ThetaFig-180)}) -- (0,{\heightFig+0.3}) node [right] {$y$};
	\draw (-\widthFig*0.5,\heightFig*0.6) node {$\OME$};
	\draw [->,>=latex] (0.35,0) arc (0:\ThetaFig-5:0.35) node[near start,above] {$\Theta$} -- (\ThetaFig:0.36);
	\draw [<-,>=latex] (-\widthFig*0.6, {-\widthFig*0.6*tan(\ThetaFig-180)}) to[bend right=20] ++(0.4,-0.4) node[right]{$\!\bordgaucheME$};
	\draw [<-,>=latex] (\widthFig*0.65,0) to[bend right=20] ++(0.15,-0.55) node[below]{$\ \ \GME$};
	\draw (\RMEFig,0) node{$\bullet$};
	\draw [<-,>=latex] (\RMEFig+0.05,0.05) to[bend left=10] ++(0.45,0.45) node[right]{\raisebox{1.5ex}{$\!\PME$}};
	\draw [<-,>=latex] ({\RMEFig*cos(55)}, {\RMEFig*sin(55)}) to[bend left=20] ++(0.4,0.45) node[right]{\raisebox{1ex}{$\SME$}};
	\draw [<->,>=latex] (0,-0.2) -- node[midway, below] {$\RME$} (\RMEFig, -0.2);
\end{tikzpicture}\par
\caption{Les domaines utilisés pour calculer les champs lointains et pour les modèles effectifs}
\label{fig: domaine ME en intro fr}
\end{figure}
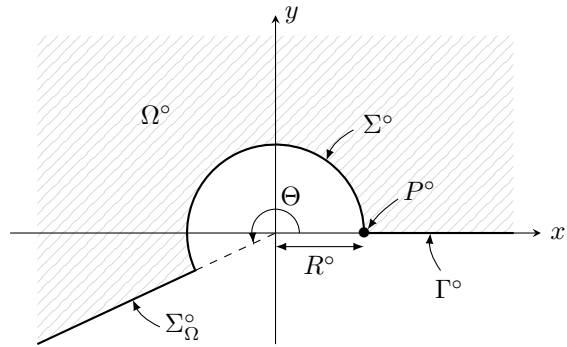

\begin{textesansboite}
Ensuite, nous calculons les champs lointains selon les équations \eqref{eq: construction directe des upl en intro fr} du théorème~\ref{construction directe des upl} (section \ref{sec: numerique des CL}). La première difficulté consiste à gérer le comportement singulier du champ lointain au niveau du coin. Nous y remédions en introduisant une frontière artificielle $\SME$ autour de l'angle (voir figure~\ref{fig: domaine ME en intro fr}) et en imposant une condition de Dirichlet-to-Neumann (DtN) non homogène de la forme
\[\partial _{r} u_{p,\ell } -\opF(u_{p,\ell } )=g \quad \  \text{ sur }  \SME.\]
Ici, $\opF$ est l'opérateur de DtN standard pour un problème de Helmholtz dans un coin de Dirichlet. Le terme de droite $g$ n'est pas standard ; il dépend des champs lointains précédents et tient compte de la condition sur $\sigma _{d} (u_{p,\ell } )$ pour $d\in -\TN$ du théorème~\ref{construction directe des upl}. Cette méthode est plus pratique que la méthode des fonctions singulières dans ce cas car elle résout trois problèmes à la fois : le fait que $u_{p,\ell }$ est non variationnel, le manque de régularité au coin de sa partie variationnelle, et le besoin de calculer $\sigma _{d} (\bfuDA_{p,\ell } )$ pour $d\in \TN$.\\

Une difficulté supplémentaire réside ici dans le traitement de la condition de Dirichlet $u_{p,\ell } = \sum_{n=0}^\infty  \frac{\mu _{0} }{\mu _{1}} \,T_{n} \cdot  (\partial _{x} ^2 +k_{1} ^2 )^{n} \partial _{y} u_{p-1-2n,\ell }$ sur $\Gamma $. Elle ne peut pas être calculée directement, car les dérivées d'éléments finis sont peu régulières et ne peuvent donc pas être une condition de Dirichlet pour une inconnue éléments finis. Nous résolvons ce problème en projetant les dérivées sur l'espace des éléments finis et en régularisant ces dérivées.\\

Nous fournissons des résultats numériques qui valident le développement asymptotique jusqu'à l'ordre~5 (section~\ref{sec: resultats numeriques CL}). À notre connaissance, c'est la première fois qu'un développement asymptotique est calculé numériquement au-delà de l'ordre 2 pour un problème comportant un angle et une couche mince. Cela est possible grâce à la généralité de notre méthode pour les correcteurs de coin et de la construction théorique présentée dans les chapitres~\ref{chap: raccord algebrique} et~\ref{chap: construction du DA}.\\

L'inconvénient d'approcher la solution exacte par une somme de champs lointains est que cela peut nécessiter de calculer un grand nombre de champs lointains (par exemple, 12 à l'ordre~5 pour $\Theta = \frac{3\pi }{2}$) et que la construction par récurrence des champs lointains induit une accumulation d'erreurs numériques. Les modèles effectifs du chapitre suivant ne présentent pas ce problème.
\end{textesansboite}

\subsection{Chapitre \ref{chap: ME}: modèles effectifs}

\begin{textesansboite}
Dans le chapitre~\ref{chap: ME}, nous proposons deux modèles effectifs d'ordre élevé, resp. d'ordre~3 et 5 par rapport à $\varepsilon $. Ils permettent de calculer une solution approchée sans avoir à mailler finement la couche, ni le voisinage du coin. En posant $\omega =0$, ce travail fournit également des modèles d'ordre élevé pour l'équation de Poisson (en remplaçant simplement $H^{1}$ par l'espace variationnel approprié), qui n'existaient pas non plus dans la littérature. Le modèle effectif d'ordre $m\in \{3,5\}$ est conçu de telle sorte que sa solution est proche de $\sum _{p<m,\ell } \powelne \,u_{p,\ell }$ (qui est une approximation de $u_{\varepsilon }$ d'ordre $m$) et tire parti de ce que l'on sait sur les champs lointains.\\

Le modèle d'ordre~3, resp.\ 5, est défini sur le domaine $\OME$ de la figure~\ref{fig: domaine ME en intro fr} selon le système d'équations suivant resp.
\begin{equation}{\label{eq: ME3 en intro fr}}
	\text{trouver\ \,}  u\in H^{1} (\OME) \quad  \text{tel que}  \quad  \left\{ \begin{array}{r@{\,=\,}ll}
		\mu _{0} \Delta  u+\omega ^2 \rho _{0} u & \fsource & \text{ dans }  \OME\\[0.5mm]
		\mu _{0} \partial _{n} u + \frac{\mu _{1} }{\varepsilon } u & 0 & \text{ sur }  \GME\\[0.5mm]
		u & 0 & \text{ sur }  \bordgaucheME\\[0.5mm]
		\partial _{r} u - \opS3 (u_{|\SME} ) & 0 & \text{ sur }  \SME 
	\end{array}\right. \qquad 
\end{equation}
\begin{equation}{\label{eq: ME5 en intro fr}}
	\text{trouver\ \,$u\in  H^{1} (\OME)\cap  H^{1} (\GME)$\ \;tel que\ \;}  \left\{ \begin{array}{r@{\,=\,}ll@{}}
		\mu _{0} \Delta  u+\omega ^2 \rho _{0} u & \fsource & \text{ dans }  \OME\\[0.5mm]
		\!\mu _{0} \,\partial _{n} u + \frac{\mu _{1} }{\varepsilon } u - \frac{\varepsilon }{3} (\mu _{1} \partial _{x} ^2 +\omega ^2 \rho _{1} )u & 0 & \text{ sur }  \GME\\[0.5mm]
		u & 0 & \text{ sur }  \bordgaucheME\\[0.5mm]
		\partial _{r} u - \opS5 (u_{|\SME} ) & 0 & \text{ sur }  \SME \\[0.5mm]
		\partial _{x} u - \opP5 (u_{|\SME} ) & 0 & \text{ en }  \PME
	\end{array}\right.
\end{equation}
Nous imposons sur $\GME$ les mêmes conditions au bord effectives que pour la couche infinie, cf. \eqref{eq: conditions effectives classiques intro fr}. De~plus, la nouveauté de nos modèles réside dans le remplacement du coin par une condition de DtN sur la frontière artificielle $\SME$. Dans le modèle d'ordre~5 de \eqref{eq: ME5 en intro fr}, nous imposons également une condition de DtN non classique au point $\PME$ (l'extrémité de $\GME$, cf. figure~\ref{fig: domaine ME en intro fr}), ce qui est nécessaire en raison de la présence d'un $\partial _{x} ^2 $ dans la condition au bord sur $\GME$.\\

Les opérateurs de DtN $\opS$ et $\opP$ prennent en compte le comportement singulier près du coin des champs lointains ainsi que les correcteurs de coin. Ils reposent sur une décomposition modale approchée de $\scl := \sum _{p<m,\ell } \powelne \,u_{p,\ell }$ à l'ordre $\varepsilon ^{m}$. Cette décomposition est tronquée, ne conservant que les $N$ premiers modes, où $N$ est un paramètre de la méthode. $\opS$ et $\opP$ sont des opérateurs non locaux et sont de rang fini $N$.\\

Une décomposition modale préliminaire de $\scl$ est obtenue en sommant les décompositions modales des champs lointains prouvées dans le chapitre~\ref{chap: decompo modale CL} (proposition~\ref{decompo modale SE 1}). Une version simplifiée est
\[\begin{array}{c}
	\displaystyle \scl = \sum _{d\in \TZ,\,d\geqslant -m} \sigma _{d} (\sclDA) \cdot \Phi ^{\varepsilon } _{d} \;+\, O(\varepsilon ^{m} \ln ^{\ell _{m}} \!\varepsilon ) \\
	\displaystyle \text{avec } \; \ell _{m} \in \mathbb{N} ,\quad \ \Phi ^{\varepsilon } _{d} := \sum _{n=0} ^{m-1} \varepsilon ^{n} \.\mcll_{d,n} 
	 \quad \  \text{ et } \quad \ \forall d\in \txtTZ*,\ \sigma _{d} (\sclDA) := \! \sum _{p\in \NN \cap [0,m)} \sum _{\ell \in \mathbb{N}} \powelne \cdot  \sigma _{d} (\bfuDA_{p,\ell } ) \in \mathbb{C} .
\end{array}\]
La condition au coin portant sur les coefficients $\sigma _{d} (\bfuDA_{p,\ell } )$ dans le théorème~\ref{construction directe des upl} induit des contraintes sur les nombres $\sigma _{d} (\sclDA)$ qu'il faut prendre en compte. Ces contraintes dépendent indirectement des coefficients à l'infini $\sigma _{q} (\ccDA_{d',k} )$ des correcteurs de coin, donc elles contiennent des informations sur la perturbation locale de la géométrie au niveau de l'angle. Nous en déduisons une seconde décomposition (proposition~\ref{decompo modale SE 2})
\begin{equation}{\label{eq: decompo modale SE 2 en intro fr}}
	\scl = \sum_{n=0}^\infty  \sigma _{n\pst} (\sclDA) \cdot \Psi ^{\varepsilon } _{n} \;+\, O(\varepsilon ^{m} \ln ^{\ell _{m}} \!\varepsilon )
\end{equation}
où les modes $\Psi ^{\varepsilon } _{n}$ sont des combinaisons linéaires des $\Phi ^{\varepsilon } _{d}$. Les modes $(\Psi ^{\varepsilon } _{1} ,\dots ,\Psi ^{\varepsilon } _{N} )$ sont ceux utilisés pour construire les opérateurs de DtN $\opS$ et $\opP5$, ou plus précisément leur orthonormalisation de Gram-Schmidt car ils ne sont pas orthogonaux en général.\\

Nous démontrons que pour $\varepsilon $ suffisamment petit les modèles effectifs sont bien posés et uniformément stables par rapport à $\varepsilon $ (théorème~\ref{caractere bien pose des ME}). Notre stratégie est de prouver la coercivité du problème en montrant que la partie principale de la formulation variationnelle (dépendant de l'équation de Helmholtz et de la condition sur $\GME$) est coercive, que $\opSref$ est positif pour $\varepsilon =0$, et que $\opS-\opSref$ et $\opP5$ ont des normes faibles lorsque $\varepsilon $ est assez petit. Nous utilisons le fait que, pour $\varepsilon =0$, les modes $\Psi ^{0} _{n}$ sont proportionnels à $J_{n\pst} (k_{0} r) \sin (n\pst \Theta )$, i.e. aux modes de l'équation de Helmholtz dans un coin de Dirichlet sans couche. De plus, $\|\opS-\opSref\|$ est majoré grâce à la continuité par rapport à $\varepsilon $ des fonctions $\Psi ^{\varepsilon } _{n}$. Ceci est rendu précis par des estimations fines des fonctions $\mcll_{d,n}$ dans le chapitre~\ref{chap: decompo modale CL}. Nous obtenons la condition de stabilité suivante, où $C$ est une constante indépendante de $(\varepsilon ,\RME,N)$ :
\[\varepsilon  \leqslant  C \. \frac{\RME}{N^{3} \.(1+|{\ln \RME}|)^{3} } \ \ \text{si } m=3 \qquad  \text{ et } \qquad  \varepsilon  \leqslant  C \. \frac{\RME^{4/3} }{N^{3} \.(1+|{\ln \RME}|)^{2} } \ \ \text{si } m=5.\]\\

De plus, en notant $\sme$ la solution du modèle effectif d'ordre $m$, nous montrons l'estimation d'erreur suivante (théorème~\ref{cvg des ME})
\begin{equation}{\label{eq: cvg des ME en intro fr}}
	\|\sme-u_{\varepsilon } \|_{H^{1} (\OME)} \leqslant  C \left[ \bigg( \frac{\varepsilon }{\RME} \bigg)^{m} (1+|{\ln \varepsilon }|)^{\ell _{m}} \.N^{3} + \left( \frac{\RME}{R} \right)^{N\pst} \right]
\end{equation}
où $C$ est indépendante de $(\varepsilon ,\RME,N)$ et $\ell _{m} \in [\![0,4]\!]$ dépend de $\Theta $. Le second terme de \eqref{eq: cvg des ME en intro fr} est dû à l'erreur de troncature de \eqref{eq: decompo modale SE 2 en intro fr}, que nous majorons grâce à \eqref{eq: vitesse de cvg decompo modale des CL en intro fr}. En particulier, en prenant par exemple $N:=\lceil {m \frac{\Theta }{\pi } \log _{2} \varepsilon }\rceil $ et $\RME<\Rsource/2$, on obtient 
\[\|\sme-u_{\varepsilon } \|_{H^{1} (\OME)} =O(\varepsilon ^{m} \ln ^{\ell _{m}} \!\varepsilon  ) \qquad  \text{quand }  \varepsilon \rightarrow 0\]
avec $\ell _{m} = 0$ pour la plupart des angles.\\

L'équation \eqref{eq: cvg des ME en intro fr} est démontrée en deux étapes. Tout d'abord, nous obtenons une erreur en $\varepsilon ^{m-1}$ à cause du remplacement de la condition de Dirichlet par $\mu _{0} \partial _{n} \sme3 = \red{\varepsilon ^{-1} } \mu _{1} \sme3$\footnote{Le rouge est utilisé dans ce document pour attirer l'attention sur des éléments importants ou des détails qui peuvent facilement être manqués.} (et de même pour $m=5$). Ensuite, comme pour les méthodes de pénalisation, l'ordre de convergence peut être amélioré à condition que $\partial _{n} \sme$ soit suffisamment régulier. Nous utilisons des arguments de régularité elliptique, qui sont rendus plus difficiles par le coin que possède le domaine $\OME$ en $\PME$.\\

Des résultats numériques sont fournis, montrant l'efficacité de nos modèles (section~\ref{sec: numerique ME}). Voir par exemple les figures \ref{fig: ME angle pi SE intro fr}, \ref{fig: ME angle pi err intro fr} et~\ref{fig: graphe err ME pi intro fr} pour le cas $\Theta =\pi$, un angle particulièrement pertinent et présentant l'un des développements asymptotiques les plus complexes. La section~\ref{sec: numerique ME} fournit également des résultats numériques pour d'autres angles pertinents.
\end{textesansboite}
\end{otherlanguage}

\begingroup
\def\nameslope{pente}
\pgfplotsset{/pgf/number format/use comma}

	\tikzset{external/export=true}%
	\tikzsetfigurename{fig-ME_pi_intro_fr-}%
	\begingroup
\def\largeurimages{53mm}
\def\colorbarheight{20mm}
\def\mycolorbar[#1]{\,\colorbar[#1]\hspace*{-7mm}}
\def\espaceapressubcaption{\hspace*{0mm}}
\def\boite#1{\raisebox{-.45\height}{\vbox{\vskip2mm\relax \hbox{#1}}}}%

\begin{figure}[h]
    \boite{%
    \dimen0=\largeurimages\relax
    \dimen0=1.375\dimen0\relax
    \edef\tempo{\the\dimen0}%
    \begin{tikzpicture}[x=\largeurimages, y=\largeurimages]
        \clip (-0.5,{-0.5*0.1/4}) -- (0.5,{-0.5*0.1/4}) -- (0.5,0) arc (0:180:0.5) -- cycle;
        \draw (0,0.335) node{\includegraphics[width=\tempo]{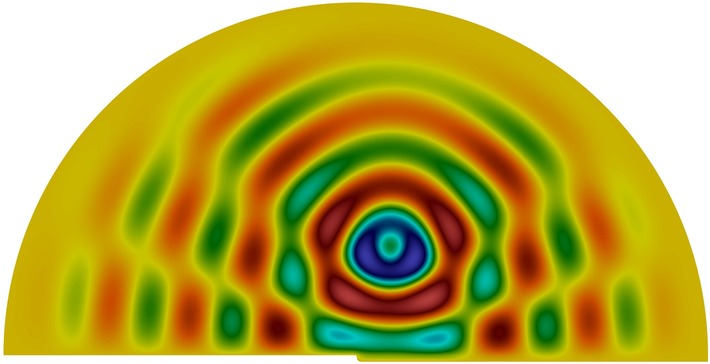}};
    \end{tikzpicture}}
    \mycolorbar[min=-2.1, max=1.4, 1e-2]

    \caption{Partie réelle de la solution exacte $u_\peps$ pour $\Theta=\pi$ et $\peps=0{,}1$}
    \label{fig: ME angle pi SE intro fr}
\end{figure}

\begin{figure}[h]
    \begin{subfigure}[c]{0.5\textwidth}
        \centering
        \boite{\includegraphics[width=\largeurimages]{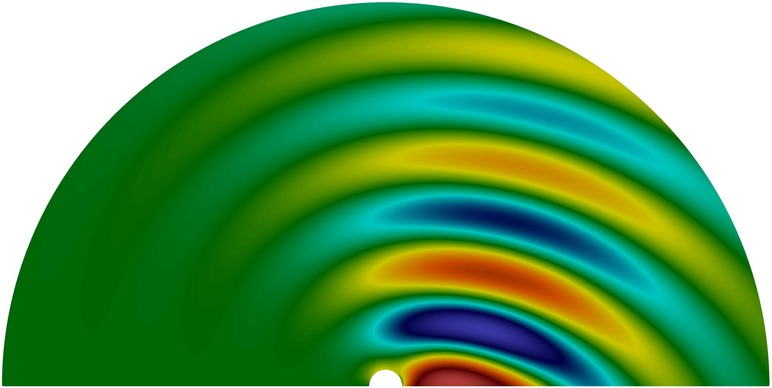}}
        \mycolorbar[min=-3.4, max=4.7, 1e-4]
        \subcaption{\normalsize $m=3$\espaceapressubcaption}
    \end{subfigure}%
    \begin{subfigure}[c]{0.5\textwidth}
        \centering
        \boite{\includegraphics[width=\largeurimages]{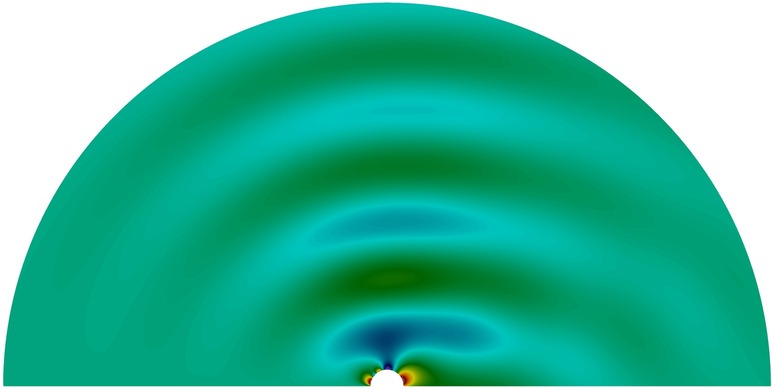}}
        \mycolorbar[min=-5, max=10, 1e-5]
        \subcaption{\normalsize $m=5$\espaceapressubcaption}
    \end{subfigure}

    \caption{Partie réelle de la différence $\protect\sme -u_\peps$ entre la solution des modèles effectifs et la solution exacte pour $\Theta=\pi$ et $\peps=0{,}1$}
    \label{fig: ME angle pi err intro fr}
\end{figure}

\endgroup%
	\tikzset{external/export=false}%
	\tikzsetfigurename{main-figure}
	\tikzset{external/export=true}%
	\tikzsetfigurename{fig-erreurs_ME_pi_intro_fr-}%
	\pgfplotsset{
    every axis/.append style={
        thick,
    },
    cycle list name=Dark2_avec_regressions,
    /pgfplots/table/x=epsilon,
    /pgfplots/table/y=nbCC,
    legend cell align=left,
}

\begin{figure}
\begin{tikzpicture}
\begin{axisME}[erreur relative $\dfrac{\|\sme-u_\peps\|_{H^1 (\OME)}}{\|u_{0,0}\|_{H^1 (\OME)}}$]{err_pi_tot.dat}{}
\courbeMEintro[y=errH1_ME_3, m=3, reg=0:5]
\courbeMEintro[y=errH1_ME_5, m=5, reg=1:5]
\end{axisME}
\end{tikzpicture}
\longCaption{Erreurs numériques des modèles effectifs pour $\Theta=\pi$}{ (les lignes continues représentent les régressions linéaires)}
\label{fig: graphe err ME pi intro fr}
\end{figure}%
	\tikzset{external/export=false}%
	\tikzsetfigurename{main-figure}

\endgroup

\let\varepsilon=\feps
\chapter{Algebraic formal series: a powerful tool for matching conditions}
\label{chap: SF}

\section{Preliminary: an introduction to the usual algebraic formal series}
\label{sec: C[[X]]}

\begin{textesansboite}
Before introducing the formal series that we will use in this document, that have been designed especially to study our problem, let us first present the usual formal series. Formal series are a concept that may seem complex at first glance, but it becomes a very powerful tool once mastered. The aim of this section is to let the reader become progressively familiar with this concept, by~starting with the simplest type of formal series. We will present these series not exactly as usually done, but with some adaptations, in order to better prepare the introduction of the formal series of corner singularities that we will actually use to derive the matching conditions. For other introductions on the subject, see \cite{Niv69, Sam23}.\\

Before giving the precise definitions, let us explain the idea. The usual formal series with coefficients in $\mathbb{C}$, whose set is denoted $\FCX$, are a generalization of convergent power series. Convergent power series designates functions of the type $x\mapsto \sum_{n=0}^\infty  a_{n} x^{n}$, with $(a_{n} )\in \mathbb{C}^{\mathbb{N}}$, that are well-defined on a neighborhood of 0. $\FCX$ contains those convergent power series, but also other objects  that cannot be associated with a function of $x$ but that behave and can be manipulated like power series. These objects represent ``non convergent power series''.  For example ``$\sum_{n=0}^\infty  n!\.x^{n}$'' is one of these series, as it does~not converge for~any $x\neq 0$. Although these non convergent series are~not well-defined in the classical sense (i.e.\ as functions), formal series give a meaning to them. Moreover, in $\FCX$, it will be possible to write calculations like the one below, as one would do on convergent power series:
\begin{equation}{\label{eq: eq hypothetique dans C[[X]]}}
    ``\; \sum_{n=0}^\infty  n!\.x^{n} + \sum_{n=0}^\infty  (-1)^{n} \.n! \.x^{n} = \sum _{n\in 2\mathbb{N}} 2\.n! \.x^{n} \;".
\end{equation}
This is what we meant by ``the formal series can be manipulated like convergent power series''. Note that for now, we put quotes as we have not defined those objects yet.\\

Non convergent power series appear naturally when one tries to describe the behavior of a function $f\in \mathcal{C}^{\infty } (\mathbb{R},\mathbb{C})$ near $x=0$. This behavior is given by the following family of approximations:
\[\forall N\in \mathbb{N}, \quad  f(x) \underset{x\rightarrow 0}{=} \sum _{n=0} ^{N} \frac{f^{(n)}(0)}{n!} x^{n} + o(x^{N} ).\]
When doing computations on asymptotic behaviors, like we have to do in the method of matched asymptotic expansions, one has to deal with the sum truncation at $N$ and the small $o$. This complicates the task, although these elements appear only for technical reasons and do~not contribute to the information of the asymptotic behavior, which is given by the terms $ \frac{f^{(n)}(0)}{n!} x^{n}$. If $f$ is analytic, this can be avoided by working directly on the power series $\sum_{n=0}^\infty  \frac{f^{(n)}(0)}{n!} x^{n}$, which is equal to $f$. However, for a general $\mathcal{C}^{\infty }$ function, this series may not converge. In fact, Borel's Lemma states that, for~any $(a_{n} )\in \mathbb{C}^{\mathbb{N}}$, there exists $f\in \mathcal{C}^{\infty } (\mathbb{R},\mathbb{C})$ s.t.\ $f^{(n)}(0)=a_{n}$ for~any $n$ \cite[p.\,22--25]{Erd56}. Taking $a_{n} := n!^2 $, we obtain a function $f$ s.t.\ $\sum_{n=0}^\infty  \frac{f^{(n)}(0)}{n!} x^{n}$ diverges for~any $x\neq 0$. It~is also possible that this series converges, but not towards $f$. Take for example $f:x\mapsto  e^{x} +e^{-1/|x| }$ (extended by continuity to $x=0$) and $\sum_{n=0}^\infty  \frac{1}{n!} x^{n}$. In these cases, the formal series ``$\sum_{n=0}^\infty  \frac{f^{(n)}(0)}{n!} x^{n}$'' offers an alternative.\\

In a way, formal series are similar to distributions. A distribution is an object that generalizes functions. It can be manipulated like a function on some aspects (it can be differentiated, restricted to a smaller domain, etc.)\ but it cannot be associated to a function of $x$. Distributions include objects that are~not difficult to understand intuitively, like the Dirac $\delta $, but that cannot be defined mathematically in the classical framework. Historically, distributions and formal series have both been used intuitively by mathematicians and physicists before their rigorous definition was introduced, but a more advance use of these notions requires solid foundations, that are given by a mathematical definition.\\

Even when working on objects defined in the classical sense, the generalized objects can be useful in an intermediate proof step. For instance, this happens often when solving PDEs: the PDE may be easier to solve in the distributions (e.g.\ using the Fourier transform), and one can show a~posteriori that the found solution(s) are actually functions in the classical sense. As~for formal series, in \cref{chap: raccord algebrique,chap: construction du DA} we build the asymptotic behavior of the far fields at the corner as formal series because we do~not know a~priori whether they are convergent series. Later, in \cref{chap: decompo modale CL} we prove that the obtained objects are actually convergent (which is highly non-trivial).\\

The mathematical definition of $\FCX$ is also comparable to the one of distributions. Any function $f\in L^{1} _{\mathrm{loc}} (\mathbb{R})$ is characterized by the $\langle f,\varphi \rangle :=\int _{\mathbb{R}} f\varphi $ with $\varphi \in \mathcal{D}(\mathbb{R})$, which is why an application $\mathcal{D}(\mathbb{R}) \rightarrow \mathbb{R}$ gives a generalization of a function. Similarly, any convergent power series $x\mapsto \sum_{n=0}^\infty  a_{n} x^{n}$ is characterized by its coefficients $(a_{n} )$, so a general sequence of coefficients provides a generalization of power series.
\end{textesansboite}

\begin{table}[h]
\renewcommand{\arraystretch}{1.2}
\caption{Comparison between the construction of distributions and formal series}
\begin{tabular}{>{\centering\arraybackslash}m{40mm}cc}
\cline{2-3}
& Distributions & Formal series \\
\hline
Set of classical objects & $L^{1} _{\mathrm{loc}} (\mathbb{R})$ & convergent power series \\
\hline
Set of generalized objects & $\mathcal{D}'(\mathbb{R})$ & $\FCX$ \\
\hline
Objects used to define the generalization & $\langle T,\varphi \rangle $ with $\varphi \in \mathcal{D}(\mathbb{R})$ & the coefficients of the series\\
\hline
\end{tabular}
\end{table}

\subsection{Definition of the formal series set $\FCX$}

\begin{definition}{}
\label{def: C[[X]]}%
We define $\FCX := \mathbb{C}^{\mathbb{N}}$.
\end{definition}

\begin{textesansboite}
For now $\FCX$ has no structure and is inert, but we will soon add one. Eventually we will be able to write formal series with a notation that looks like series. But in the beginning, let us write them simply as families in $\mathbb{C}^{\mathbb{N}}$.\\

The link between Definition~\cref{def: C[[X]]} and convergent power series is the following. For~any $(a_{n} )\in \FCX$, if the function $f_{(a_{n} )}$ written below is well-defined on a neighborhood of 0, we say that the formal series $(a_{n} )$ is \emph{convergent} and that its \emph{limit function} is
\[f_{(a_{n} )} : x \;\mapsto \; \sum_{n=0}^\infty  a_{n} x^{n} := \lim _{N\rightarrow \infty } \sum _{n=0} ^{N} a_{n} x^{n} .\]
The map $(a_{n} ) \mapsto f_{(a_{n} )}$ is a bijection from convergent formal series to convergent power series (defined as functions). On the contrary, a simple example of non convergent formal series is given by $(n!)_{n\in \mathbb{N}}$ as $\sum_{n=0}^\infty  n!\,x^{n}$ converges only for $x=0$.
\end{textesansboite}

\begin{exemplenumerote}
\label{ex: SF du DL d'une fct scalaire}%
The other important example of formal series, mentionned in the introduction, is the formal series associated with the asymptotic behavior at 0 of a function $f\in \mathcal{C}^{\infty } (\mathbb{R},\mathbb{C})$. This series is defined as $f^{0} := \big( \frac{f^{(n)}(0)}{n!} \big)_{n\in \mathbb{N}}$.
\end{exemplenumerote}

\begin{definition}{basic operations in $\FCX$}
\label{def: operations basiques dans C[[X]]}%
Let $f=(a_{n} )$ and $g=(b_{n} )$ in $\FCX$, and $\lambda \in \mathbb{C}$. We define in $\FCX$:
\[\lambda  f := (\lambda  a_{n} )_{n} 
\quad  \text{ and } \quad  f+g:= (a_{n} +b_{n} )_{n}
\quad  \text{ and } \quad  fg:= \bigg( \sum _{k=0} ^{n} a_{k} b_{n-k} \bigg)_{\!n} .\]\fbeq
\end{definition}

\begin{remarque}
This makes $\FCX$ a $\mathbb{C}$-vector space. The multiplication works like in Cauchy product of series, and one can verify that it~is associative, commutative and bilinear. Moreover, the operations defined above are natural because they are compatible with the operations on limit functions (if defined) and Taylor expansions:
\begin{itemize}
    \item If $f,g\in \FCX$ are convergent, then $\lambda f$, $f+g$ and $fg$ are convergent and their limit function is $\lambda $ times the limit function of $f$ (resp.\ the sum and product of the limit functions of $f$ and $g$).
    \item If $u,v$ belong to $\mathcal{C}^{\infty } (\mathbb{R},\mathbb{C})$ and $u^{0} ,v^{0}$ are the formal series associated to their asymptotic behavior at $x=0$ (like in \cref{ex: SF du DL d'une fct scalaire}), then $\lambda  u^{0}$ (resp.\ $u^{0} +v^{0}$ and $u^{0} \,v^{0}$) represents the asymptotic behavior of $\lambda  u$ (resp.\ $u+v$ and $uv$).
\end{itemize}
\end{remarque}

\refstepcounter{notation}
\begin{paragraphesansboite}{Notation \thenotation:}
\label{def: X dans C[[X]]}%
We denote $X := (\delta _{n,1} )_{n\in \mathbb{N}} \in \FCX$.
\end{paragraphesansboite}

\begin{lemmesansboite}{}
\label{X^k dans C[[X]]}%
For~any $k\in \mathbb{N}$, we have $X^{k} = (\delta _{n,k} )_{n\in \mathbb{N}} \in \FCX$ (where $X^{k}$ denotes the iterated multiplication of $X$ $k$ times).
\end{lemmesansboite}

\begin{demo}
It~is easy to prove by induction. First, one can verify that $(\delta _{n,0} )_{n\in \mathbb{N}}$ is the multiplicative neutral element of $\FCX$. Hence $X^{0} = (\delta _{n,0} )_{n}$. Then if we assume $X^{k} = (\delta _{n,k} )_{n}$ for a given $k$, we can deduce
\[X\cdot X^{k} = (\delta _{n,1} )_{n} \cdot  (\delta _{n,k} )_{n} = \bigg(\sum_{i=0}^n \delta _{i,1} \cdot \delta _{n-i,k} \bigg)_{\!n} = (\delta _{n,k+1} )_{n} .\]
\end{demo}

\begin{exemple}
Let $j\in \mathbb{N}$ and $(a_{0} ,\dots ,a_{j} )\in \mathbb{C}^{j+1}$. \Cref{X^k dans C[[X]],def: operations basiques dans C[[X]]} imply: \[\sum _{k=0} ^{j} a_{k} X^{k} = \sum _{k=0} ^{j} (a_{k} \delta _{n,k} )_{n} = \bigg( \sum _{k=0} ^{j} a_{k} \delta _{n,k} \bigg)_{\!n} =(a_{0} ,a_{1} ,\dots ,a_{j-1} ,a_{j} ,0,0,\dots ).\]
This type of formal series whose coefficients have a finite support is a \emph{formal polynomial}. The set of formal polynomials, usually denoted $\mathbb{C}[X]$, is thus included in $\FCX$. Formal series can be seen as a generalization of formal polynomial where an infinite number of coefficients may be non-zero.
\end{exemple}

\begin{textesansboite}
For now the sum defined in $\FCX$ (\cref{def: operations basiques dans C[[X]]}) allows us to sum only a finite number of elements (by iterating the sum of two elements). On the contrary, the next definition allows us to sum an infinite number of formal series in some situations.
\end{textesansboite}

\begin{definitionsecable}{}
\label{def: fsum dans C[[X]]}%
Let $I$ be a set (that may be infinite) and $(f_{i} )_{i\in I} \in \FCX^{I}$ be a family of formal series. For~any $i\in I$, let us denote $(a_{i,n} ):=f_{i}$ (each $a_{i,n}$ is a complex number).\\
We say that $(f_{i} )$ \emph{has a finite support on each coordinate} iff
\begin{center}
    for~any $n\in \mathbb{N}$, $\{i\in I \mid a_{i,n} \neq 0\}$ is finite.
\end{center}
If this property is verified, then we define
\begin{equation}{\label{eq: def fsum dans C[[X]]}}
    \fsum_{i\in I} f_{i} := \bigg( \sum _{i\in I} a_{i,n} \bigg)_{\!n\in \mathbb{N}} \in \FCX
\end{equation}
where each $\sum _{i\in I} a_{i,n}$ is a classical sum with only a finite number of non-zero terms. Note the boldness of the symbol $\fsum$, that differentiates it from a classical sum $\sum $.\\
If $I$ is an interval of integers, i.e.\ $I=\mathbb{Z}\cap [a,b]$ with $a,b\in \mathbb{Z}\cup \{\pm \infty \}$, then we also denote $\fsum\limits_{i=a} ^{b} f_{i} := \fsum\limits_{i\in I} f_{i}$.
\end{definitionsecable}

\begin{textesansboite}
The next proposition gives a way to finally write formal series in a way that looks like a series.
\end{textesansboite}

\begin{proposition}{}
\label{ecriture des elements de C[[X]] en puissances de X}%
For~any $(a_{n} )_{n\in \mathbb{N}} \in \FCX$, we have
\[(a_{n} )_{n} = \fsumn a_{n} X^{n} .\]\fbeq
\end{proposition}

\begin{demo}
By~\cref{X^k dans C[[X]],def: operations basiques dans C[[X]]}, we have $a_{k} X^{k} = (a_{k} \delta _{n,k} )_{n\in \mathbb{N}}$ for~any $k\in \mathbb{N}$. Next, for~any $n\in \mathbb{N}$, $\{n\in \mathbb{N} \mid a_{k} \delta _{n,k} \neq 0\}\subset \{n\}$ is finite, so $(a_{k} X^{k} )_{k\in \mathbb{N}}$ has finite support on each coordinate. Hence we can write
\begin{equation}{}
    \fsum_{k=0} ^{\infty } a_{k} X^{k} = \bigg(\sum_{k=0}^\infty  a_{k} \,\delta _{n,k} \bigg)_{\!n\in \mathbb{N}} = (a_{n} )_{n\in \mathbb{N}} .
\end{equation}
Note that we can change the index of $\fsum$ from $n$ to $k$ as it~is only a bound variable.
\end{demo}

\begin{textesansboite}
In~practice we will most of the time manipulate the notation $\fsumn a_{n} X^{n}$ instead of ``$(a_{n} )$'', as it looks more meaningful.
\end{textesansboite}

\begin{remarque}
The classical way to introduce formal series is to set ``$\sum_{n=0}^\infty  a_{n} X^{n} := (a_{n} )$'' by~definition and to use the usual symbol $\sum $. In this document, we chose another approach: to use a different symbol $\fsum$ to avoid confusion with classical sums, and to make this symbol denote generalized infinite sums, which will be very useful to express the asymptotic behaviors of the fields around the corner.
\end{remarque}

\begin{remarque}
More generally, for~any $I\subset \mathbb{N}$ and $(a_{n} )\in \mathbb{C}^{I}$, $(a_{n} X^{n} )_{n\in I}$ has finite support on each coordinate, so $\fsum_{n\in I} a_{n} X^{n}$ is well defined. Denoting $\tilde{a} _{n} :=a_{n}$ if $n\in I$ and $\tilde{a} _{n} :=0$ if $n\in \mathbb{N} \setminus I$, we have
\[\fsum_{n\in I} a_{n} X^{n} = ( \tilde{a} _{n} )_{n\in \mathbb{N}} .\]\fbeq
\end{remarque}

\begin{remarque}
Using \cref{ecriture des elements de C[[X]] en puissances de X}, \cref{eq: def fsum dans C[[X]]} can be rewritten as follows. For~any set $I$ and $(a_{i,n} )_{(i,n)} \in \mathbb{C}^{I\times \mathbb{N}}$ s.t.\ $((a_{i,n} )_{n} )_{i}$ has finite support on each coordinate, 
\[\fsum_{i\in I} \fsumn a_{i,n} X^{n} = \fsumn \bigg(\sum _{i\in I} a_{i,n} \bigg)X^{n}\]
where the last $\sum _{i\in I}$ is a classical sum with a finite number of non-zero terms.
\end{remarque}

\begin{textesansboite}
The next result states a few natural properties that one would expect for a summation symbol.
\end{textesansboite}

\begin{propositionsecable}{properties of $\fsum$}
\label{proprietes de fsum dans C[[X]]}%
Let $I,J$ be two sets, and let $(f_{i} )_{i\in I} \in \FCX^{I}$ and $(g_{j} )_{j\in J} \in \FCX^{J}$. We assume that $(f_{i} )$ and $(g_{j} )$ have a finite support on each coordinate (abbreviated here FSEC).
\begin{enumerate}
    \item If $I$ is finite, then $\displaystyle \fsum_{i\in I} f_{i} = \sum _{i\in I} f_{i}$.
    \item For~any $\lambda \in \mathbb{C}$, $(\lambda  f_{i} )_{i\in I}$ has FSEC and $\displaystyle \lambda \fsum_{i\in I} f_{i} = \fsum_{i\in I} \lambda  f_{i}$.
    \item If $I=J$, then $(f_{i} +g_{i} )_{i\in I}$ has FSEC and $\displaystyle \fsum_{i\in I} f_{i} +\fsum_{i\in I} g_{i} = \fsum_{i\in I} (f_{i} +g_{i} )$.
    \item If $I\cap J= \varnothing  $ and if we denote $h_{k} :=f_{k}$ if $k\in I$ and $h_{k} :=g_{k}$ if $k\in J$, then $(h_{k} )_{k\in I\cup J}$ has FSEC and $\displaystyle \fsum_{i\in I} f_{i} +\fsum_{j\in J} g_{j} = \fsum_{k\in I\cup J} h_{k}$.
    \item For~any $k\in \mathbb{N}$, $(X^{k} f_{i} )_{i\in I}$ has FSEC and $\displaystyle X^{k} \fsum_{i\in I} f_{i} = \fsum_{i\in I} X^{k} f_{i}$.
    \item If there exists a bijection $\varphi :I\rightarrow J$ s.t.\ $g_{\varphi (i)} =f_{i}$ for~any $i\in I$, then $\displaystyle \fsum_{i\in I} f_{i} =\fsum_{j\in J} g_{j}$. In other words, changes of index are possible.
\end{enumerate}\fbi
\end{propositionsecable}

\begin{demo}
For~any $i\in I$ and $j\in J$, let us denote $(a_{i,n} )_{n} :=f_{i}$ and $(b_{j,n} )_{n} :=g_{j}$.
\begin{enumerate}
    \item Direct by~definition of $\fsum$ and \cref{def: operations basiques dans C[[X]]}.
    \item \Cref{def: operations basiques dans C[[X]]} states that $\lambda f_{i} = (\lambda a_{i,n} )_{n}$ for~any $i\in I$. In addition, for~any $n\in \mathbb{N}$, $\{i\in I \mid \lambda  a_{i,n} \neq 0\}$ is included in $\{i\in I \mid a_{i,n} \neq 0\}$ which is finite because $(f_{i} )$ has FSEC. So $(\lambda  f_{i} )_{i\in I}$ has FSEC. Moreover, 
    \begin{align*}
        \lambda  \fsum_{i\in I} f_{i} &= \lambda  \cdot \bigg(\sum _{i\in I} a_{i,n} \bigg)_{\!n} && \text{by \cref{def: fsum dans C[[X]]}} \\
        &= \bigg(\lambda  \sum _{i\in I} a_{i,n} \bigg)_{\!n} && \text{by \cref{def: operations basiques dans C[[X]]}}  \\
        &= \bigg(\sum _{i\in I} \lambda  a_{i,n} \bigg)_{\!n} \\
        &= \fsum_{i\in I} \lambda  f_{i} && \text{by \cref{def: operations basiques dans C[[X]]}} 
    \end{align*}
    \item \Cref{def: operations basiques dans C[[X]]} states that $f_{i} +g_{i} = (a_{i,n} +b_{i,n} )_{n}$ for~any $i\in I$. For~any $n\in \mathbb{N}$, $\{i\in I \mid a_{i,n} +b_{i,n} \neq 0\}$ is included in $\{i\in I \mid a_{i,n} \neq 0\} \cup  \{i\in I \mid b_{i,n} \neq 0\}$ which is finite because $(f_{i} )$ and $(g_{i} )$ has FSEC. So $(f_{i} +g_{i} )_{i\in I}$ has FSEC. Moreover,
    \begin{align*}
        \fsum_{i\in I} f_{i} +\fsum_{i\in I} g_{i} &= \bigg(\sum _{i\in I} a_{i,n} \bigg)_{\!n} + \bigg(\sum _{i\in I} b_{i,n} \bigg)_{\!n} && \text{by \cref{def: fsum dans C[[X]]}} \\
        &= \bigg(\sum _{i\in I} (a_{i,n} +b_{i,n} ) \bigg)_{\!n} && \text{by \cref{def: operations basiques dans C[[X]]}}  \\
        &= \fsum_{i\in I} \,(a_{i,n} +b_{i,n} )_{n} && \text{by \cref{def: fsum dans C[[X]]}} \\
        &= \fsum_{i\in I} \,(f_{i} +g_{i} ) && \text{by \cref{def: operations basiques dans C[[X]]}} 
    \end{align*}
    \item One can verify that $(h_{k} )$ has FSEC similarly to point 2. Next, for~any $k\in I\cup J$ and $n\in \mathbb{N}$, let us define $c_{k,n} :=a_{k,n}$ if $k\in I$ and $c_{k,n} :=b_{k,n}$ if $k\in J$. We can write:
    \begin{align*}
        \fsum_{i\in I} f_{i} +\fsum_{j\in J} g_{j} &= \bigg(\sum _{i\in I} a_{i,n} \bigg)_{\!n} + \bigg(\sum _{j\in J} b_{j,n} \bigg)_{\!n} && \text{by \cref{def: fsum dans C[[X]]}} \\
        &= \bigg(\sum _{i\in I} a_{i,n} +\sum _{j\in J} b_{j,n} \bigg)_{\!n} && \text{by \cref{def: operations basiques dans C[[X]]}}  \\
        &= \bigg(\sum _{k\in I\cup J} c_{k,n} \bigg)_{\!n} \\
        &= \fsum_{k\in I\cup J} h_{k}  && \text{by \cref{def: operations basiques dans C[[X]]}} 
    \end{align*}
    \item It suffices to prove the property for $k=1$ and then to iterate on $k$. To prove the case $k=1$, let us first compute $X f$ for~any $f=(a_{n} )\in \FCX$. Using \cref{def: X dans C[[X]],def: operations basiques dans C[[X]]}, we obtain:
    \begin{equation}{\label{eq: preuve proprietes de fsum 5}}
        X f = (\delta _{n,1} )_{n} \cdot  (a_{n} )_{n} = \bigg(\sum_{m=0}^n \delta _{m,1} a_{n-m} \bigg)_{\!n} = (0, a_{0} , a_{1} , a_{2} ,\dots  ).
    \end{equation}
    For~any $n\in \mathbb{N}^*$, $\{i\in I \mid a_{i,n-1} \neq 0\}$ is finite, so $(Xf_{i} )_{i}$ has FSEC. Finally,
    \begin{align*}
        X\fsum_{i\in I} f_{i} &= X \cdot  \bigg(\sum _{i\in I} a_{i,n} \bigg)_{n} && \text{by \cref{def: fsum dans C[[X]]}} \\
        &= \bigg(0 , \sum _{i\in I} a_{i,0} , \dots ,\sum _{i\in I} a_{i,n-1} ,\dots \bigg)_{n} && \text{by \cref{eq: preuve proprietes de fsum 5}} \\
        &= \fsum_{i\in I} (0, a_{i,0} , \dots , a_{i,n-1} ,\dots )_{n} && \text{by \cref{def: fsum dans C[[X]]}} \\
        &= \fsum_{i\in I} Xf_{i} && \text{by \cref{eq: preuve proprietes de fsum 5}} 
    \end{align*}
    \item Left to the reader.
\end{enumerate}
\end{demo}

\begin{exemple}
Item 5 of \cref{proprietes de fsum dans C[[X]]} shows that, for~any $(a_{n} )\in \FCX$ and $k\in \mathbb{N}$, we have $X^{k} \fsumn a_{n} X^{n} = \fsumn a_{n} X^{n+k}$. In addition, Equation~\cref{eq: eq hypothetique dans C[[X]]}, that was written between quotes in this section's introduction, can now be rigorously written using the right notations:
\begin{align*}
    \fsumn n!\,X^{n} + \fsumn (-1)^{n} \,n! \,X^{n} &= \fsumn (n! + (-1)^{n} \,n!) \,X^{n} && \text{by \cref{proprietes de fsum dans C[[X]]}, item 3} \\
    &= \fsum_{n\in 2\mathbb{N}} 2\,n!\,X^{n} + \fsum_{n\in 2\mathbb{N}+1} 0 && \text{by \cref{proprietes de fsum dans C[[X]]}, item 4} \\
    &= \fsum_{n\in 2\mathbb{N}} 2\,n! \,X^{n} && \text{by \cref{def: fsum dans C[[X]]}} \\
\end{align*}\fbeq
\end{exemple}


\begin{definitionsansboite}{}
For~any $f=(a_{n} ) \in \FCX$, we define the order of $f$ as $\ord(f):=\inf \{n\in \mathbb{N} \mid a_{n} \neq 0\}$ with the convention $\inf (\varnothing ):=+\infty $ (i.e.\ $\ord(0)=+\infty $).
\end{definitionsansboite}

\begin{lemme}{}
\label{SF scalaires decalees}%
Let $(f_{k} )_{k\in \mathbb{N}}$ be a sequence of elements of $\FCX$. If $\ord(f_{k} )\rightarrow \infty $ when $k\rightarrow \infty $, then $(f_{k} )$ has a finite support on each coordinate (in the sense of \cref{def: fsum dans C[[X]]}).
\end{lemme}

\begin{demo}
For~any $k\in \mathbb{N}$, let us denote $(a_{k,n} )_{n} :=f_{k}$. Let $n\in \mathbb{N}$ and $k_{n}$ be s.t.\ $\ord(f_{k} )>n$ for~any $k>k_{n}$. Then $\{k\in \mathbb{N} \mid a_{k,n} \neq 0 \}$ is included in $[\![0,k_{n} ]\!]$, so it~is finite.
\end{demo}

\begin{remarque}
The implication of \cref{SF scalaires decalees} is actually an equivalence. However for the formal series of \cref{sec: SF 1D} (that are more general), we will not have equivalence.
\end{remarque}

\begin{definition}{composition of formal series}
\label{def: composition dans C[[X]]}%
Let $f =(a_{n} )$ and $g=(b_{n} )$ be two elements of $\FCX$. If $b_{0} =0$, we define $\displaystyle f\circ g :=f(g):=\fsum_{k=0} ^{\infty } \,a_{k} \,g^{k}$.
\end{definition}

\begin{demo}[Justification]
Since $b_{0} =0$, we have $g=Xh$ with $h:=\fsumn b_{n+1} X^{n}$. If one likes to break it into details, we can write:
\begin{align*}
    g &= \fsumn b_{n} X^{n} && \text{by \cref{ecriture des elements de C[[X]] en puissances de X}} \\
    &= \cancel{\fsum_{n\in \{0\}} b_{0} X^{0} } + \fsum_{n=1} ^{\infty } b_{n} X^{n} && \text{by \cref{proprietes de fsum dans C[[X]]}, item 3} \\
    &= \fsum_{n=0} ^{\infty } b_{n+1} X^{n+1} && \text{by \cref{proprietes de fsum dans C[[X]]}, item 6} \\
    &= X \fsumn b_{n+1} X^{n} =Xh && \text{by \cref{proprietes de fsum dans C[[X]]}, item 5} .
\end{align*}
Next, denoting $(c_{k,n} ):=h^{k} $, \cref{proprietes de fsum dans C[[X]]} item 5 gives $X^{k} h^{k} =\fsumn c_{k,n} X^{n+k}$. Therefore, $\ord(g^{k} )=\ord(X^{k} h^{k} )\geqslant k \rightarrow \infty $. Using \cref{SF scalaires decalees}, we deduce that the sequence $(a_{k} g^{k} )$ has a finite support on each coordinate. So $\fsum_{k=0} ^{\infty } \,a_{k} \,g^{k}$ is well-defined thanks to \cref{def: fsum dans C[[X]]}.
\end{demo}

\begin{exemple}
For~any $f\in \FCX$, we have $f(X)=f$. So~``$f(X)$'' may be used as an alternative writing of $f$.
\end{exemple}

\begin{remarque}
Similarly to the sum and the product, we have the following compatibility properties:
\begin{itemize}
    \item if the formal series $f$ and $g$ are convergent, then (at least in a neighborhood of 0) $f\circ g$ is convergent and its limit function is the composition of the limit functions of $f$ and $g$,
    \item and for~any $u,v\in \mathcal{C}^{\infty } (\mathbb{R},\mathbb{C})$ s.t.\ $v(0)=0$, $u^{0} \circ v^{0}$ represents the behavior at 0 of $u\circ v$ (in the sense of \cref{ex: SF du DL d'une fct scalaire}).
\end{itemize}
\end{remarque}

\begin{textesansboite}
To conclude this section, the important points to keep in mind are:
\begin{itemize}
    \item a formal series is the data of a sequence of coefficients,
    \item it can be written as $\fsumn a_{n} X^{n}$ which looks very similar to a usual power series and can be manipulated like one,
    \item it can be convergent (so associated to an analytic function) or non convergent (in this case it~is a generalization of classical power series),
    \item it can be used to represent efficiently the asymptotic behavior of a $\mathcal{C}^{\infty }$ function at 0,
    \item operations of sum, product and composition are allowed as long as they involve summing only a finite number of non-zero terms to compute each coefficient.
\end{itemize}
\end{textesansboite}

\subsection{Examples of problems solved in $\FCX$}

\begin{textesansboite}
In this section, we present two concrete examples in which formal series are especially useful compared to classical tools. They involve manipulations that are similar to the one we will do to derive the matching conditions. So these examples are also toy problems to prepare for \cref{chap: raccord algebrique}.
\end{textesansboite}

\begin{exemplenumerote}
\label{ex: 1e pb jouet dans C[[X]]}%
Let $f\in \mathcal{C}^{\infty } (\mathbb{R},\mathbb{C})$ and $g:x\mapsto  f(x(1+x))$, and let us compute the asymptotic behavior of $g$ at 0 in terms of the one of $f$. This is a simple example of the computation of an asymptotic behavior after a change of variable. We recall that, in the method of matched asymptotic expansion, the change of variable applied is $\bfr \leadsto  \peps^{\pm 1} \bfr$.\\

One strategy would be to express the derivatives $g^{(n)}(0)$ in terms of the $f^{(n)}(0)$. Denoting $h:x\mapsto x(1+x)$, one would need to use Faà di Bruno's formula:
\[g^{(n)}(0) =(f\circ h)^{(n)} (0)= \! \sum _{\substack{(m_{1} ,\dots ,m_{n} )\in \mathbb{N}^{n} \\ \mathclap{1\cdot m_{1} +2\cdot m_{2} +3\cdot m_{3} +\cdots +n\cdot m_{n} =n}}} \ \frac{n!}{m_{1} !\,1!^{m_{1}} \,m_{2} !\,2!^{m_{2}} \,\cdots \,m_{n} !\,n!^{m_{n}} }\cdot  f^{(m_{1} +\cdot \cdot +m_{n} )}(h(0))\cdot  \prod_{j=1}^n \big(h^{(j)}(0)\big)^{m_{j}} .\]
Needless to say that this formula is a nightmare. Another, better, possibility is to write the expansion as a sum and do calculations on it:
\[\forall N\in \mathbb{N}, \qquad  g(x) = \sum _{n=0} ^{N} \frac{f^{(n)}(0)}{n!} \big(x(1+x)\big)^{n} + o\big([x(1+x)]^{N} \big) = \mathrm{etc.}\]
This approach requires to deal with the sum truncation and the small $o$. In this simple example, it~is~not difficult, but in a more complex situation it would make calculations unnecessarily heavy.\\

Using formal series instead, one could do the following reasoning. Let $f^{0} :=\fsumn \frac{f^{(n)}(0)}{n!} X^{n}$ and $g^{0} :=\fsumn \frac{g^{(n)}(0)}{n!} X^{n}$ in $\FCX$. By~compatibility of the operations between $\mathcal{C}^{\infty }$ functions and the formal series that represent their behavior at 0, we have $g^{0} = f^{0} (X(1+X))$. So, by \cref{def: composition dans C[[X]]},
\[ g^{0} = \fsumn \frac{f^{(n)}(0)}{n!} \big(X(1+X)\big)^{n} .\]
Classical algebraic calculations in a ring (here $(\FCX,+,\times )$) give for~any $n\in \mathbb{N}$:
\[\big(X(1+X)\big)^{n} = X^{n} (1+X)^{n} = X^{n} \sum_{k=0}^n \binom{n}{k} X^{k} = \sum _{k=0} ^{n} \binom{n}{k} X^{n+k} = \sum _{k=n} ^{2n} \binom{n}{k-n} X^{k} .\]
Hence we deduce
\begin{align*}
    g^{0} &= \fsumn \sum _{k=n} ^{2n} \frac{f^{(n)}(0)}{n!} \binom{n}{k-n} X^{k} \\
    &= \fsumn \fsum_{k=n} ^{2n} \frac{f^{(n)}(0)}{n!} \binom{n}{k-n} X^{k} && \text{by \cref{proprietes de fsum dans C[[X]]}, item 1} \\
    &= \fsumn \fsum_{k=0} ^{\infty } \mathbf{1}_{n\leqslant k\leqslant 2n} \frac{f^{(n)}(0)}{n!} \binom{n}{k-n} X^{k} && \text{by \cref{proprietes de fsum dans C[[X]]}, item 4} \\
    &= \fsum_{k=0} ^{\infty } \bigg(\sum_{n=0}^\infty  \mathbf{1}_{n\leqslant k\leqslant 2n} \frac{f^{(n)}(0)}{n!} \binom{n}{k-n} \bigg) X^{k} && \text{by \cref{ecriture des elements de C[[X]] en puissances de X,def: fsum dans C[[X]]}} 
\end{align*}
where the last sum $\sum_{n=0}^\infty $ is a classical sum with a finite number of non-zero terms. Note that, to apply \cref{def: fsum dans C[[X]]} in the last line, it~is~not necessary to check that $\big(\fsum_{k=0} ^{\infty } \mathbf{1}_{n\leqslant k\leqslant 2n} \frac{f^{(n)}(0)}{n!} \binom{n}{k-n} X^{k} )_{n}$ has a finite support on each coordinate. Indeed, this family is equal to $\big( \frac{f^{(n)}(0)}{n!} \big(X(1+X)\big)^{n} )_{n}$ which already has a finite support on each coordinate according to the justification under \cref{def: composition dans C[[X]]}. Finally we obtain:
\[g^{0} = \fsum_{k=0} ^{\infty } \bigg(\sum _{n=\lceil {k/2}\rceil } ^{k} \frac{f^{(n)}(0)}{n!} \binom{n}{k-n} \bigg) X^{k} .\]
Here we have justified carefully each step, but with a bit more practice this computation can be done in just a few lines. We recall that if $f$ is analytic, convergent power series can be used instead of formal series, but in general it~is~not possible.
\end{exemplenumerote}

\begin{remarque}
In the method of matched asymptotic expansion the change of variable used is $\bfr \leadsto  \peps^{\pm 1} \bfr$. It replaces in the corner singularities the factors $\bfr^{d}$ by $\peps^{\pm d} \bfr^{d}$ and the $\ln ^{j} \!\bfr$ by $(\pm \ln \peps+\ln \bfr)^{j}$. In our example, the powers of $1+x$ act similarly to $(\pm \ln \peps+\ln \bfr)^{j}$.
\end{remarque}

\begin{exemplenumerote}
\label{ex: 2e pb jouet dans C[[X]]}%
Let us show how to compute the asymptotic behavior of the solution of a differential equation. We consider $f\in \mathcal{C}^{1} (\mathbb{R}_+,\mathbb{C})$ solution of
\begin{equation}{\label{eq: 2e pb jouet dans C[[X]]: EDO}}
    f'+f = g
\end{equation}
with $g\in \mathcal{C}^{0} (\mathbb{R}_+,\mathbb{C})$. We are interested in the behavior of $f$ at infinity in terms of $g$. Let us assume that we know that $f$ and $g$ have a behavior in powers of $ \frac{1}{x}$, i.e.\ there exist $(a_{n} )\in \mathbb{C}^{\mathbb{N}}$ and $(b_{n} )\in \mathbb{C}^{\mathbb{N}}$ s.t.
\[\forall N\in \mathbb{N}, \quad  f(x) = \sum _{n=0} ^{N} \frac{a_{n} }{x^{n}} + o\!\left( \frac{1}{x^{N}} \right)
\quad  \text{ and } \quad  g(x) = \sum _{n=0} ^{N} \frac{b_{n} }{x^{n}} + o\!\left( \frac{1}{x^{N}} \right)
\quad  \text{ when }  \quad  x\rightarrow +\infty  .\]
 We assume that $g$ and $(b_{n} )$ are a known data, whereas $f$ and $(a_{n} )$ are unknown. Let us assume also~that the expansion of $f$ at infinity may be differentiated:
\begin{equation}{\label{eq: 2e pb jouet dans C[[X]]: DA de f'}}
    \forall N\in \mathbb{N}, \qquad  f'(x) = \sum _{n=0} ^{N} a_{n} \cdot  \frac{\mathrm{d}}{\mathrm{d}x} \!\left[ \frac{1}{x^{n}} \right] + o\!\left( \frac{1}{x^{N+1}} \right).
\end{equation}
All these assumptions will be satisfied in practice when working on the far, layer and corner fields. Our goal is to compute the coefficients $a_{n}$ in terms of the $b_{n}$.\\

The standard way would be to replace each term in \cref{eq: 2e pb jouet dans C[[X]]: EDO} with its asymptotic expansion, which gives
\[\forall N\in \mathbb{N}, \qquad  \sum _{n=0} ^{N} a_{n} \frac{-n}{x^{n+1}} + o\!\left( \frac{1}{x^{N+1}} \right) +  \sum _{n=0} ^{N} \frac{a_{n} }{x^{n}} + o\!\left( \frac{1}{x^{N}} \right) = \sum _{n=0} ^{N} \frac{b_{n} }{x^{n}} + o\!\left( \frac{1}{x^{N}} \right).\]
Then, by identifying each power of $x$, we obtain
\[a_{0} =b_{0} \qquad  \text{ and } \qquad  \forall n\in \mathbb{N}^*,\ -(n-1)a_{n-1} +a_{n} =b_{n} .\]
This sure works, but it does~not provide a nice formula for $(a_{n} )$ that would make further calculations handy (for example a change of variable on $x$, cf.\ \cref{ex: 1e pb jouet dans C[[X]]}). In addition, the work needs to be repeated if the differential equation changes, because the solution does~not have a general expression.\\

Now let us show how to solve the problem using formal series. The behaviors at infinity of $f$ and $g$ can be represented by the formal series $f^{\infty } :=\fsumn a_{n} X^{n}$ and $g^{\infty } :=\fsumn b_{n} X^{n}$, where $X$ ``represents'' $ \frac{1}{x}$. Let us compute $f^{\infty }$ in terms of $g^{\infty }$. Let
\def\opA{A}
\[\opA:\left( \begin{array}{ccc} \FCX & \rightarrow  & \FCX \\ 
\displaystyle \fsumn c_{n} X^{n} & \mapsto  &\displaystyle \fsumn c_{n} nX^{n+1} \end{array} \right)\]
It~is a linear map. Since $a_{n} \frac{\mathrm{d}}{\mathrm{d}x} \!\left[ \frac{1}{x^{n}} \right] = -a_{n} \frac{n}{x^{n+1}}$ for~any $n$, and given \cref{eq: 2e pb jouet dans C[[X]]: DA de f'}, the behavior at infinity of $f'$ is described by $-A(f^{\infty } )$. Therefore, we deduce $-\opA(f^{\infty } )+f^{\infty } =g^{\infty }$. In other words,
\begin{equation}{}
    ( \mathrm{id}-\opA)(f^{\infty } ) = g^{\infty } .
\end{equation}
Now comes the magic of formal series : $ \mathrm{id}-\opA$ is invertible and $( \mathrm{id}-\opA)^{-1} = \fsumn \opA^{n}$, so
\begin{equation}{\label{eq: 2e pb jouet dans C[[X]]: solution}}
    f^{\infty } = \fsumn \opA^{n} (g^{\infty } ).
\end{equation}
Let us explain and justify this. First, note that $\ord(\opA(h))=\ord(h)+1$ for~any $h\in \FCX$. Hence, $\ord(\opA^{n} (h))\rightarrow \infty $. So, by \cref{SF scalaires decalees,def: fsum dans C[[X]]}, $\fsumn \opA^{n} (h )$ is well-defined for~any $h\in \FCX$. We denote $\fsumn \opA^{n} :\FCX\rightarrow \FCX$ the map
\[\fsumn \opA^{n} :h\mapsto  \fsumn \opA^{n} (h).\]
Next, we can write for~any $h\in \FCX$:
\[( \mathrm{id}-\opA) \circ \bigg(\fsumn \opA^{n} \bigg)(h) = \fsumn \opA^{n} (h) - \opA \bigg(\fsumn \opA^{n} (h)\bigg).\]
To continue, let us show that $A$ and $\fsum$ commute.
\begin{sousdemo}
    Let $(h_{k} )\in \FCX^{\mathbb{N}}$ that has finite support on each coordinate. For~any $k\in \mathbb{N}$, let $(c_{k,n} )_{n} :=h_{k}$. The following calculation uses the fact that, like the multiplication by $X$, $\opA$ shifts all coefficients by one step (and summing coefficients after shifting them all gives the same as shifting their sum):
    \begin{align*}
        A\bigg(\fsum_{k=0} ^{\infty } h_{k} \bigg)&= A\!\left( \bigg(\sum_{k=0}^\infty  c_{k,n} \bigg)_{\!n\in \mathbb{N}} \right) && \text{by \cref{def: fsum dans C[[X]]}} \\
        &= \fsumn \bigg(\sum_{k=0}^\infty  c_{k,n} \bigg) n X^{n+1} && \text{by~definition of } A \\
        &= X\cdot  \fsumn \bigg(\sum_{k=0}^\infty  c_{k,n} n\bigg) X^{n} && \text{by \cref{proprietes de fsum dans C[[X]]}, item 5}  \\
        &= X\cdot  \fsum_{k=0} ^{\infty } \fsumn c_{k,n} n X^{n} && \text{by \cref{ecriture des elements de C[[X]] en puissances de X,def: fsum dans C[[X]]}}  \\
        &= \fsum_{k=0} ^{\infty } \fsumn c_{k,n} nX^{n+1} && \text{by \cref{proprietes de fsum dans C[[X]]}, item 5} \\
        &= \fsum_{k=0} ^{\infty } A(h_{k} ) && \text{by~definition of } A
    \end{align*}
\end{sousdemo}
Therefore, we have
\begin{align*}
    ( \mathrm{id}-\opA) \circ \bigg(\fsumn \opA^{n} \bigg)(h) &= \fsumn \opA^{n} (h) - \fsumn \opA^{n+1} (h) \\
    &= \fsum_{n\in \{0\}} \opA^{n} (h) + \fsum_{n=1} ^{\infty } \opA^{n} (h) - \fsum_{n=1} ^{\infty } \opA^{n} (h) && \text{by \cref{proprietes de fsum dans C[[X]]}, items 4,6}  \\
    &= h && \text{by \cref{proprietes de fsum dans C[[X]]}, item 1} 
\end{align*}
Similarly,
\begin{align*}
    \bigg(\fsumn \opA^{n} \bigg)\circ ( \mathrm{id}-\opA)(h) &= \fsumn (\opA^{n} (h) -\opA ^{n+1} (h) ) && \text{by~definition of $\fsumn \opA^{n}$ and linearity of }  \opA\\
    &= \fsumn \opA^{n} (h) - \fsumn \opA^{n+1} (h) && \text{by \cref{proprietes de fsum dans C[[X]]}, item 3}  \\
    &= h && \text{by the previous computation} 
\end{align*}
We have shown that $\fsumn \opA^{n}$ is the inverse of $ \mathrm{id}-\opA$. This proves \cref{eq: 2e pb jouet dans C[[X]]: solution} and ends the resolution of the problem.\\

Note that this method does~not rely heavily on the form of the equation. It suffices to express the problem as $( \mathrm{id}-A)(f^{\infty } ) =g^{\infty }$ for some operator $\opA:\FCX \rightarrow \FCX$ s.t.\ $\fsumn A^{n}$ is well defined. For example the method could be applied to the equation $\pm f''+f=g$.\\

In \cref{inverse de id-f dans les SF 1D,inverse de id-f dans les SF 2D} we give a general condition to have $( \mathrm{id}-A)^{-1} =\fsumn A^{n}$, so~that we do~not have to repeat the calculations of this example each time we need to invert an operator $ \mathrm{id}-A$. This inversion is~not possible for~any $A$ (take e.g.\ $A:= \mathrm{id}$), but compared to the inversion of $ \mathrm{id}- \tilde{A}$ in a Banach space using a Neumann series when $\| \tilde{A} \|<1$, formal series are much more permissive. For example, if $ \mathrm{id}-A$ is invertible, then so is $ \mathrm{id}-10^{10} A$.\\

In the present example, there~is another reason why formal series are especially useful: even if $g^{\infty }$ is a convergent series and $g(x)=\sum_{n=0}^\infty  a_{n} \frac{1}{x^{n}} $ for~any $x$, $f^{\infty }$ is~not convergent in general. Take for example $g^{\infty } :=X$. In this case, it~is easy to verify by induction that $\opA^{n} (g^{\infty } )=n! X^{n+1}$ for~any $n$, so $f^{\infty } = \fsumn n! X^{n+1}$. This is the typical example of non convergent formal series. It seems to us that something similar happens for corner fields. We conjecture that the formal series associated with their behavior at infinity (defined in \cref{validation construction et DA des champs}) are~not convergent (we did not try to prove it as it was not useful).
\end{exemplenumerote}

\begin{remarque}
In \cref{ex: 1e pb jouet dans C[[X]],ex: 2e pb jouet dans C[[X]]}, we avoided as much as possible using the notation $(a_{n} )$ for an element of $\FCX$, and prefered using $\fsumn a_{n} X^{n}$. This shows how formal series can be manipulated like normal series, almost ignoring the formal construction of $\FCX$.
\end{remarque}


\section{A quick review of existing generalized algebraic formal series}
\label{sec: intro a mes SF}

\begin{textesansboite}
In the rest of this chapter, we introduce the specific kind of formal series that we will use to describe the behaviors of the fields $\bfu_{p,\ell }$ and $S_{p,\ell }$ around the corner, in order to study the matching assumption. These behaviors will be represented by formal series of corner singularities of the form $\sum_{i=0}^n r^{d} \ln ^{i} \!r\, f_{i} (\theta )$ in $\Omega $ and $\sum_{i=0}^n x^{d} \ln ^{i} \!x\, g_{i} (Y)$ in $\couche$. Analytic expressions of the corner singularities will be provided later in \cref{sec: singularites quasi homogenes}. In this chapter, we do~not try to describe the singularities ; we consider a set of singularities as a black box and we explain how to build formal series of these singularities.\\

The type of formal series introduced here is a novel generalization of usual formal series. To our knowledge, the most notable existing generalizations of $\FCX$ are the following:
\begin{itemize}
    \item $R[ [X] ]$ where $R$ may be any commutative ring.
    \item Formal series of several variables: for example $R[ [X,Y] ]$ is the set of the formal series ``$\sum _{n,m\in \mathbb{N}} a_{n,m} \* X^{n} Y^{m}$'' with $(a_{n,m} )\in R^{\mathbb{N}^2 }$.
    \item Laurent series: $R((X))$ is the set of the formal series ``$\sum _{n=N} ^{\infty } a_{n} X^{n}$'' with $N\in \mathbb{Z}$ and $a_{n} \in R$ for~any $n$. It~is a ring and if $K$ is a field, then $K((X))$ is a field and it~is the smallest field that contains $K[ [X] ]$.
    \item Puiseux series: for~any field $K$, $K\{\{X\}\} :=\bigcup _{k\in \mathbb{N}^*} K(( X^{ \frac{1}{k} } ))$ is the set of the series ``$\sum _{n=N} ^{\infty } a_{n} X^{ \frac{n}{k} }$'' with $k\in \mathbb{N}^*$, $N\in \mathbb{Z}$ and $a_{n} \in K$ for~any $n$. If $K$ is algebraically closed and has characteristic 0, then $K\{\{X\}\}$ is the algebraic closure of the field $K((X))$.
    \item Hahn series: for~any field $K$ and any totally ordered abelian group $G$, $K[ [X^{G} ] ]$ is the set of the formal series ``$\sum _{h\in H} a_{h} X^{h}$'' where $H$ is a well-ordered subset of $G$ (i.e.\ every non-empty subset of $H$ has a minimum) and $a_{h} \in K$ for~any $h$. Hahn series are a generalization of Puiseux series because we have $K\{\{X\}\} \subsetneq K[ [X^{\mathbb{Q}} ] ]$.
\end{itemize} 
There are also other generalizations in more distant areas of mathematics, like non commutative algebra.\\

Here we introduce two types of formal series.
\begin{itemize}
    \item Series with index in $\mathbb{R}$: These series describe the behavior around the corner of a single field $\bfu_{p,\ell }$ or $S_{p,\ell }$. They have the form
    \[\text{``$\displaystyle \sum _{d\in \mathbb{R}} \varphi _{d} $''} \]
    where each $\varphi _{d}$ is a corner singularity that behaves roughly like $r^{d}$. Moreover, we ask that any truncation of the form ``$\sum _{d\in [d_{1} ,d_{2} ]} \varphi _{d}$'' has only a finite number of non-zero terms.
    \item Series with index in $\mathbb{R}^2 $: These series describe the behaviors of $u_{\peps} \approx $ ``$\sum _{p,\ell } \peps^{p} \ln ^{\ell } \!\peps \,u_{p,\ell }$'' and ``$\sum _{p,\ell } \peps^{p} \ln ^{\ell } \!\peps \,S_{p,\ell }$'' around the corner and when $\peps\rightarrow 0$. They have the form
    \[\text{``$\displaystyle \sum _{(d,p)\in \mathbb{R}^2 } \peps^{p} \varphi _{d} $''} \]
    where each $\varphi _{d}$ is a corner singularity that behaves roughly like $r^{d}$ and that contains the dependence in $\ln \peps$. Again, we ask that any truncation on a bounded set of indexes $(d,p)$ has a finite support.\\
\end{itemize}

In comparison with the existing types of formal series, our series differ in several ways. First we adopt a point of view from linear algebra rather than commutative algebra: the terms of our series belong to vector spaces rather than a ring and we do~not define a multiplication operation on our formal series. Indeed, the spaces of corner singularities defined in \cref{sec: singularites quasi homogenes} are vector spaces that are~not stable by multiplication (and we never need to multiply singularities with each other anyway). This allows us to consider summation indexes that are more general in a way. 
Indeed, in the mentioned usual generalizations of $\FCX$, the summation index must belong to a set that is bounded from below, but we remove this restriction. In~practice the series associated to the $\bfu_{p,\ell }$ will involve indexes bounded from below and going to $+\infty $, because they describe a behavior when $\bfr\rightarrow 0$, but the series associated to the $S_{p,\ell }$ will involve indexes bounded from above and going to $-\infty $, because they describe a behavior when $\bfr\rightarrow \infty $. Moreover we consider formal series with an index in $\mathbb{R}^2 $, a set that is~not even ordered.\\

The other major difference with pre-existing formal series, is the way we build and use our series. We introduce a new operation symbol, the bold summation symbol $\fsum$, that provides a unified and rigorous way to both write formal series as the sum of their terms, like ``$\fsum_{d} \varphi _{d}$'', and to write sums of infinitely many formal series, which we will use a lot in \cref{chap: raccord algebrique,chap: complements sur le DA,chap: decompo modale CL}. There exist another, rather classical, way of summing infinitely many elements of $\FCX$, which is to introduce a special topology on $\FCX$. However, with the more general formal series that we are using, and given the calculational use we are aiming at, the approach with $\fsum$ is simpler. Moreover we introduce a specific system of calculational tools centered around the notion of operator with a degree, a kind of operator that includes differential operators on corner singularities and some important resolvant operators introduced in \cref{chap: raccord algebrique} (denoted $\ROmegaPi$, $\RGammaPi$ and $\RcouchePi$, see \cref{operateurs R}). We define some infinite series of these operators, prove useful commutation properties, and explain how to inverse these operators.\\

In \cref{sec: SF generales} we introduce general definitions that apply to formal series indexed by $\mathbb{R}$ or $\mathbb{R}^2 $. In \cref{sec: SF 1D}, resp.\ \ref{sec: SF 2D}, we present notions and properties that are specific to the formal series whose index lies in $\mathbb{R}$, resp.\ $\mathbb{R}^2 $. Before that, let us explain how the polynomial behavior in $\ln \varepsilon $ can be treated. The following definition defines the set of polynomials with coefficients in a vector space.
\end{textesansboite}

\begin{definition}{}
\label{def: E[T]}%
Let $E$ be a $\mathbb{C}$-vector space. Let $T$ denote a formal indeterminate. We define
\[E[T] :=\mathbb{C}[T] \otimes _{\mathbb{C}} E\]
which is naturally equipped with a structure of $\mathbb{C}[T]$-module.
\end{definition}

\begin{textesansboite}
In~practice $T$ will represent $\ln \peps$ and will even be denoted ``$\ln \varepsilon $'' (where ``$\varepsilon $'' is only a symbol and not a real number). If the sentence about $\mathbb{C}[T]$-modules makes sense to you, then this is all you need to know about $E[T]$. Otherwise, let us explain what it means. A module is the generalization of a vector space where we replace the field by a ring. Here the ring is $\mathbb{C}[T]$. Let us give a more basic construction of $E[T]$, but sufficient for our use. We define $E[T]$ as the set of the sequences $(\varphi _{n} )\in E^{\mathbb{N}}$ that have finite support: $\exists N\in \mathbb{N},\forall n>N,\;\varphi _{n} =0$. Then for~any $k\in \mathbb{N}$ and $(\varphi _{n} )\in E[T]$, we set $T^{k} \cdot (\varphi _{n} )_{n\in \mathbb{N}} :=( \widetilde{\varphi }_{n} )_{n\in \mathbb{N}}$ where $\widetilde{\varphi }_{n} :=\varphi _{n-k}$ if $n\geqslant k$ and $\widetilde{\varphi }_{n} :=0$ otherwise. We also consider the injection $\varphi \in E \mapsto  (\varphi ,0,0,\dots )\in E[T]$ as an inclusion. Then it~is easy to check that for~any $k,k'\in \mathbb{N}$ and $(\varphi _{n} )\in E[T]$, we have $T^{k} \cdot (T^{k'} \cdot (\varphi _{n} )_{n} )=T^{k+k'} \cdot (\varphi _{n} )_{n}$ and $(\varphi _{n} )_{n} =\sum_{n=0}^\infty  T^{n} \cdot \varphi _{n}$ where the last sum has a finite number of non-zero terms. Therefore we obtain:
\[E[T] = \bigg\{\sum _{n=0} ^{N} T^{n} \cdot \varphi _{n} \ \bigg|\ N\in \mathbb{N},\, (\varphi _{0} ,\dots ,\varphi _{N} )\in E^{N+1} \bigg\}.\]
Hence $E[T]$ has the intended behavior: it~is the set of polynomials in $T$ with coefficients in $E$.
\end{textesansboite}

\section{General definitions}
\label{sec: SF generales}

\newcommand{\dimsf}{m}
\newcommand{\Rsf}{\mathbb{R}^\dimsf}
\newcommand{\isf}{q}
\newcommand{\Esf}{E_{\rm lfs}}

\begin{textesansboite}
Let $\dimsf \in \{1,2\}$ and $(E_{\isf} )_{\isf\in \Rsf}$ be a family of $\mathbb{C}$-vector spaces. In~practice, $E_{d}$ (if $m=1$) and $E_{(d,p)}$ (if $m=2$) will be the space of corner singularities that are roughly $d$-homogeneous (up to some logarithms). If $m=2$, $p$ the second coordinate of $\isf$ will be used for powers of $\peps$ as explained later in \cref{sec: SF 2D}. The next paragraph defines what we consider to be a formal series.
\end{textesansboite}

\begin{definition}{}
\label{def: support locmt fini et SF}%
Let $\Phi =(\varphi _{\isf} )_{\isf\in \Rsf} \in  \prod _{\isf\in \Rsf} E_{\isf}$. We define the support of $(\varphi _{\isf} )$ as
\[\operatorname{supp}(\Phi ) := \{\isf \in  \Rsf \mid \varphi _{\isf} \neq 0\}.\]
We say that the support of $\Phi $ is \emph{locally finite} iff: 
\begin{center}
    for~any bounded set $B\subset \Rsf$, $\operatorname{supp}(\Phi ) \cap B$ is finite.
\end{center}
In this case, we say that $\Phi $ is a \emph{formal series} over $(E_{\isf} )_{\isf\in \Rsf}$. The $\varphi _{d}$ are called the \emph{terms} or the \emph{components} of $\Phi $.
\end{definition}
\nomenclature[28-30]{$\operatorname{supp}$}{Support of a formal series}

\begin{remarquenumerotee}
\label{rq: reformulation de LFS}%
Equivalently, $\operatorname{supp}(\Phi )$ is locally finite iff $\operatorname{supp}(\Phi ) \cap B_{\|.\|}(0,R)$ is finite for~any $R>0$, where $B_{\|.\|}(0,R)$ denotes the ball of radius $R$ centered at 0 for a given norm $\|.\|$ on $\Rsf$. Moreover, a locally finite support is at most countable.
\end{remarquenumerotee}

\begin{textesansboite}
Let us denote $\Esf$ the set of the $\Phi  \in \prod _{\isf\in \Rsf} E_{\isf}$ that have a locally finite support. This is a temporary notation, as \cref{ecriture avec fsum de l'ensemble des SF 1D,rq: ecriture avec fsum de l'ensemble des SF 2D} will give more natural ways to write the sets of formal series.
\end{textesansboite}

\begin{definition}{}
\label{def: operations basiques dans prod}%
Let $\Phi =(\varphi _{\isf} )$ and $\Psi =(\psi _{\isf} )$ in $\Esf$, and $\lambda \in \mathbb{C}$. We define in $\Esf$:
\[\lambda  \Phi  := (\lambda  \varphi _{\isf} )_{\isf} 
\quad  \text{ and } \quad  \Phi  + \Psi := (\varphi _{\isf} +\psi _{\isf} )_{\isf} .\]\fbeq[1.1]
\end{definition}

\begin{textesansboite}
With this operations, $\Esf$ is a $\mathbb{C}$-vector space. Contrary to $\FCX$, we have not defined a multiplication on $\Esf$ as it will not be useful in our case.\\

The next definition generalizes to our formal series the notion of having ``finite support on each coordinate'', that was introduced in \cref{def: fsum dans C[[X]]} for $\FCX$. Our aim is to describe when a $\fsum$-sum of formal series is well defined. Since we would like the result of the $\fsum$-sum to still be a formal series, we modify the condition to have a kind of locally uniform control on the supports of the summed family. 
\end{textesansboite}

\begin{definition}{family with collectively locally finite support}
\label{def: famille avec CLFS}%
Let $I$ be a set (that may be infinite) and $(\Phi _{i} )_{i\in I} \in \Esf^{I}$ be a family of formal series. For~any $i\in I$, let us denote $(\varphi _{i,\isf} ):=\Phi _{i}$ (so, for~any $\isf \in \Rsf$, $\varphi _{i,\isf} \in E_{\isf}$).\\
We say that $(\Phi _{i} )_{i\in I}$ has \emph{collectively locally finite support} iff
\begin{center}
    for~any bounded set $B\subset \Rsf$, $\{i\in I \mid \operatorname{supp}(\Phi _{i} )\cap B\neq  \varnothing  \}$ is finite.
\end{center}
\end{definition}

\begin{remarque}
Like in \cref{rq: reformulation de LFS}, we obtain an equivalent criterion if the bounded sets $B$ are replaced by the balls $B_{\|.\|}(0,R)$ for~any $R>0$ and a given norm $\|.\|$ on $\Rsf$. Moreover, if $(\Phi _{i} )_{i}$ has collectively locally finite support, then $\{i\in I \mid \Phi _{i} \neq 0\}$ is at most countable.
\end{remarque}

\SaveVerb{fsum}|\fsum|
\begin{SaveVerbatim}{deffsum}
    \usepackage{contour}
    \contourlength{0.05ex}
    \DeclareMathOperator*{\fsum}{\mathchoice
        {\contour*{.}{$\displaystyle\sum$}}%
        {\contour*{.}{$\textstyle\sum$}}%
        {\contour*{.}{$\scriptstyle\sum$}}%
        {\contour*{.}{$\scriptscriptstyle\sum$}}}
\end{SaveVerbatim}

\begin{definition}{symbol $\fsum$}
\label{def: fsum}%
Let $I$ be a set and $(\Phi _{i} )_{i\in I} \in \Esf^{I}$ be a family of formal series. For~any $i\in I$, let us denote $(\varphi _{i,\isf} ):=\Phi _{i}$.\\
If $(\Phi _{i} )_{i\in I}$ has collectively locally finite support, then we define
\begin{equation}{\label{eq: def fsum}}
    \fsum_{i\in I} \Phi _{i} := \bigg( \sum _{i\in I} \varphi _{i,\isf} \bigg)_{\!\isf\in \Rsf} \in \prod _{\isf\in \Rsf} E_{\isf}
\end{equation}
where each $\sum _{i\in I} \varphi _{i,\isf}$ is a classical sum with only a finite number of non-zero terms. Note the boldness of the symbol $\fsum$\footnote{To obtain this symbol, we use the following \LaTeX{} code that defines the macro \UseVerb{fsum}.\UseVerbatim{deffsum}}, that differentiates it from a classical sum $\sum $.\\
If $I$ is an interval of integers, i.e.\ $I=\mathbb{Z}\cap [a,b]$ with $a,b\in \mathbb{R}\cup \{\pm \infty \}$, then we also denote $\fsum\limits_{i=a} ^{b} \Phi _{i} := \fsum\limits_{i\in I} \Phi _{i}$.
\end{definition}
\nomenclature[53]{$\fsum$}{Generalized summation symbol for formal series}

\begin{lemme}{}
\label{les SF sont stables par fsum}%
Let $I$ be a set and $(\Phi _{i} )_{i\in I} \in \Esf^{I}$ be a family of formal series that has collectively locally finite support. Then $\operatorname{supp}\!\big(\fsum_{i\in I} \Phi _{i} \big)\subset  \bigcup _{i\in I} \operatorname{supp}(\Phi _{i} )$ and
\[\fsum_{i\in I} \Phi _{i} \in \Esf .\]\fbeq 
\end{lemme}

\begin{demo}
For~any $i\in I$, let us denote $(\varphi _{i,q} )_{q\in \Rsf} :=\Phi _{i}$. For~any $q\in \Rsf$, we have $\big(\sum _{i\in I} \varphi _{i,q} \neq 0 \Rightarrow  \exists i\in I,\, \varphi _{i,q} \neq 0)$. So \cref{def: fsum} implies $\operatorname{supp}\!\big(\fsum_{i\in I} \Phi _{i} \big)\subset  \bigcup _{i\in I} \operatorname{supp}(\Phi _{i} )$. Next, let $B\subset \Rsf$ be a bounded set and let $I_{B} :=\{i\in I \mid \operatorname{supp}(\Phi _{i} )\cap B\neq \varnothing \}$. We have
\[\operatorname{supp}\!\bigg(\fsum_{i\in I} \Phi _{i} \bigg)\cap B \;\subset \; \bigcup _{i\in I} \operatorname{supp}(\Phi _{i} )\cap B \;=\; \bigcup _{i\in \red{I_{B} }} \operatorname{supp}(\Phi _{i} )\cap B.\]
Since $I_{B}$ and the sets $\operatorname{supp}(\Phi _{i} )\cap B$ are all finite by assumption, $\operatorname{supp}\!\big(\fsum_{i\in I} \Phi _{i} \big)\cap B$ is finite and $\fsum_{i\in I} \Phi _{i}$ has locally finite support.
\end{demo}

\begin{textesansboite}
The next Proposition is very similar to \cref{proprietes de fsum dans C[[X]]} (stated for $\FCX$).
\end{textesansboite}

\begin{propositionsecable}{properties of $\fsum$}
\label{proprietes de fsum}%
Let $I,J$ be two sets, and let $(\Phi _{i} )_{i\in I} \in \Esf^{I}$ and $(\Psi _{j} )_{j\in J} \in \Esf^{J}$. We assume that $(\Phi _{i} )$ and $(\Psi _{j} )$ have collectively locally finite support (abbreviated here CLFS).
\begin{enumerate}
    \item If $I$ is finite, then $\displaystyle \fsum_{i\in I} \Phi _{i} = \sum _{i\in I} \Phi _{i}$.
    \item For~any $\lambda \in \mathbb{C}$, $(\lambda  \Phi _{i} )_{i\in I}$ has CLFS and $\displaystyle \lambda \fsum_{i\in I} \Phi _{i} = \fsum_{i\in I} \lambda  \Phi _{i}$.
    \item If $I=J$, then $(\Phi _{i} +\Psi _{i} )_{i\in I}$ has CLFS and $\displaystyle \fsum_{i\in I} \Phi _{i} +\fsum_{i\in I} \Psi _{i} = \fsum_{i\in I} (\Phi _{i} +\Psi _{i} )$.
    \item If $I\cap J= \varnothing  $ and if we denote $\Upsilon _{k} :=\Phi _{k}$ if $k\in I$ and $\Upsilon _{k} :=\Psi _{k}$ if $k\in J$, then $(\Upsilon _{k} )_{k\in I\cup J}$ has CLFS and $\displaystyle \fsum_{i\in I} \Phi _{i} +\fsum_{j\in J} \Psi _{j} = \fsum_{k\in I\cup J} \Upsilon _{k}$.
    \item If there exists a bijection $\varphi :I\rightarrow J$ s.t.\ $\Psi _{\varphi (i)} =\Phi _{i}$ for~any $i\in I$, then $\displaystyle \fsum_{i\in I} \Phi _{i} =\fsum_{j\in J} \Psi _{j}$. In other words, changes of index are possible.
\end{enumerate}\fbi[0.5]
\end{propositionsecable}

\begin{demo}
We only show that the involved families have CLFS, as the calculations to prove the formulas are almost identical to the proof of \cref{proprietes de fsum dans C[[X]]}.
\begin{itemize}
    \item[2.] For~any $i\in I$, we have $\operatorname{supp}(\lambda \Phi _{i} )\subset \operatorname{supp}(\Phi _{i} )$. So for~any bounded set $B\subset \Rsf$, $\{i\in I \mid \operatorname{supp}(\lambda \Phi _{i} )\cap B\neq \varnothing \}$ is included in $\{i\in I \mid \operatorname{supp}(\Phi _{i} )\cap B\neq \varnothing \}$ which is finite.
    \item[3.] For~any $i\in I$, we have $\operatorname{supp}(\Phi _{i} +\Psi _{i} )\subset \operatorname{supp}(\Phi _{i} )\cup \operatorname{supp}(\Psi _{i} )$. So for~any bounded set $B\subset \Rsf$, $\{i\in I \mid \operatorname{supp}(\Phi _{i} +\Psi _{i} )\cap B\neq \varnothing \}$ is included in $\{i\in I \mid \operatorname{supp}(\Phi _{i} )\cap B\neq \varnothing \}\cup \{i\in I \mid \operatorname{supp}(\Psi _{i} )\cap B\neq \varnothing \}$ which is finite.
    \item[4.] For~any bounded set $B\subset \Rsf$, $\{k\in I\cup J \mid \operatorname{supp}(\Upsilon _{k} )\cap B\neq \varnothing \}$ is included in $\{i\in I \mid \operatorname{supp}(\Phi _{i} )\cap B\neq \varnothing \}\cup \{j\in J \mid \operatorname{supp}(\Psi _{j} )\cap B\neq \varnothing \}$ which is finite.
\end{itemize}
\end{demo}

\begin{proposition}{Fubini property}
\label{Fubini pour fsum}%
Let $I,J$ be two sets and $(\Phi _{i,j} )_{(i,j)\in I\times J} \in \Esf^{I\times J}$ that has collectively locally finite support (CLFS). Then $(\Phi _{i,j'} )_{j'\in J}$ and $(\Phi _{i',j} )_{i'\in I}$ have CLFS for~any $i\in I$ and $j\in J$, $\big(\fsum_{j\in J} \Phi _{i,j} \big)_{i\in I}$ and $\big(\fsum_{i\in I} \Phi _{i,j} \big)_{j\in J}$ also have CLFS, and 
\[\fsum_{(i,j)\in I\times J} \Phi _{i,j} = \fsum_{i\in I} \fsum_{j\in J} \Phi _{i,j} = \fsum_{j\in J} \fsum_{i\in I} \Phi _{i,j} .\]\fbeq[0.7]
\end{proposition}

\begin{textesansboite}
Note that, like for the usual Fubini's theorem, we need a control over the terms on $I\times J$ entirely. It~is~not sufficient that $(\Phi _{i,j'} )_{j'\in J}$ and $(\Phi _{i',j} )_{i'\in I}$ have CLFS for~any $i\in I$ and $j\in J$.
\end{textesansboite}

\begin{demo}
For~any $(i,j)$, let us denote $(\varphi _{i,j,q} )_{q\in \Rsf} := \Phi _{i,j}$. Let $i\in I$. For~any bounded set $B\subset \Rsf$, $\{j'\in J \mid \operatorname{supp}(\Phi _{i,j'} )\cap B \neq  \varnothing \}$ is in bijection with $\{(i',j')\in I\times J \mid i'=i \text{ and } \operatorname{supp}(\Phi _{i,j'} )\cap B\neq \varnothing \}$, which is included in the finite set $\{(i',j')\in I\times J \mid \operatorname{supp}(\Phi _{i,j'} )\cap B\neq \varnothing \}$. Hence $(\Phi _{i,j'} )_{j'\in J}$ has CLFS. Similarly, one can verify that $(\Phi _{i',j} )_{i'\in I}$ has CLFS for~any $j\in J$.\\
Then, $\fsum_{j\in J} \Phi _{i,j}$ is well defined for~any $i\in I$ and equal to $\big(\sum _{j\in J} \varphi _{i,j,q} \big)_{q}$. In addition, using that $\operatorname{supp}\big(\fsum_{j\in J} \Phi _{i,j} \big) \subset  \bigcup _{j\in J} \operatorname{supp}(\Phi _{i,j} )$, we have for~any bounded set $B\subset \Rsf$:
\begin{align*}
    \bigg\{i\in I \;\bigg|\; \operatorname{supp}\bigg(\fsum_{j\in J} \Phi _{i,j} \bigg) \cap B\neq \varnothing  \bigg\}  
    &\subset  \bigg\{i\in I \;\bigg|\; \bigcup _{j\in J} \operatorname{supp}(\Phi _{i,j} ) \cap B\neq \varnothing  \bigg\} \\
    &= \{i\in I \mid \exists j\in J,\; \operatorname{supp}(\Phi _{i,j} )\cap B \neq \varnothing \} \\
    &= \pi _{1} \big(\{(i,j)\in I\times J \mid \operatorname{supp}(\Phi _{i,j} )\cap B\neq \varnothing  \}\big)
\end{align*}
with $\pi _{1} :(i,j)\mapsto  i$. Therefore $\big(\fsum_{j\in J} \Phi _{i,j} \big)_{i\in I}$ has CLFS. Similarly, $\big(\fsum_{i\in I} \Phi _{i,j} \big)_{j\in J}$ has CLFS. Now we can write:
\begin{align*}
    \fsum_{i\in I} \fsum_{j\in J} \Phi _{i,j} &= \fsum_{i\in I} \bigg(\sum _{j\in J} \varphi _{i,j,q} \bigg)_{\!q\in \Rsf} && \text{by \cref{def: fsum}}  \\
    &= \bigg( \sum _{i\in I} \sum _{j\in J} \varphi _{i,j,q} \bigg)_{\!q\in \Rsf} && \text{by \cref{def: fsum}}  \\
    &= \bigg( \sum _{(i,j)\in I\times J} \varphi _{i,j,q} \bigg)_{\!q\in \Rsf} && \text{because these sums have finite support} \\
    &= \fsum_{(i,j)\in I\times J} \Phi _{i,j} && \text{by \cref{def: fsum}} 
\end{align*}
Similarly, one can show $\fsum_{j\in J} \fsum_{i\in I} \Phi _{i,j} = \fsum_{(i,j)\in I\times J} \Phi _{i,j} $.
\end{demo}

\begin{definition}{}
\label{def: fsum d'ensembles}%
Let $I$ be a set and $(A_{i} )_{i\in I}$ be a family of subsets of $\Esf$. We denote
\[\fsum_{i\in I} A_{i} := \bigg\{ \fsum_{i\in I} \Phi _{i} \ \bigg|\ (\Phi _{i} )\in \prod _{i\in I} A_{i} \text{ that has collectively locally finite support} \bigg\}.\]\fbeq
\end{definition}

\begin{remarque}
In \cref{def: fsum d'ensembles}, $\fsum_{i\in I} A_{i}$ is included in $\Esf$ by \cref{les SF sont stables par fsum}.
\end{remarque}

\section{Formal series with index in $\mathbb{R}$}
\label{sec: SF 1D}

\begin{textesansboite}
In this section, we define the type of formal series that will represent the behaviors of one far-and-layer field or one corner field around the corner. These behaviors have the form ``$\sum _{d} \varphi _{d}$'' where each $\varphi _{d}$ is a corner singularity that is $d$-homogeneous up to some logarithms.\\

Let $E$ be a $\mathbb{C}$-vector space and $(E_{d} )_{d\in \mathbb{R}}$ be a family of vector subspaces of $E$. We assume that $E$ is generated by the spaces $E_{d}$, that is to say
\[E = \sum _{d\in \mathbb{R}} E_{d} \]
where $\sum _{d\in \mathbb{R}} E_{d}$ (with the usual symbol $\sum $) classically denotes all the \emph{finite} sums of elements of the spaces $E_{d}$
\[\sum _{d\in \mathbb{R}} E_{d} := \bigg\{ \sum _{d\in \mathbb{R}} \varphi _{d} \ \bigg|\ (\varphi _{d} )\in \prod _{d\in \mathbb{R}} E_{d} \text{ s.t.\ only a finite number of $\varphi _{d}$ are non zero}  \bigg\} .\]
Moreover, we assume that the $E_{d}$ are in direct sum, which is usually denoted as
\[E = \bigoplus_{d\in \mathbb{R}} E_{d} .\]
This means that, for~any $\varphi \in E$, there exists a unique family $(\varphi _{d} ) \in \prod _{d\in \mathbb{R}} E_{d}$ with finite support s.t.\ $\varphi =\sum _{d\in \mathbb{R}} \varphi _{d}$. We denote
\[\iota :E \rightarrow \prod _{d\in \mathbb{R}} E_{d}\]
the (linear) mapping of $\varphi $ to this family $(\varphi _{d} )$. Hence $\iota $ is an injection. In~practice, by abuse of notation, we will consider this injection as an inclusion and we will drop the ``$\iota $'' most of the time. This is similar to the indentification of functions to their associated distribution.
\end{textesansboite}

\begin{remarque}
It~is easy to check that for~any $\Phi \in \prod _{d\in \mathbb{R}} E_{d}$ with locally finite support, we have: $\Phi \in  \iota (E) \Leftrightarrow  \operatorname{supp}(\Phi )$ is finite $\Leftrightarrow  \operatorname{supp}(\Phi )$ is bounded.
\end{remarque}

\begin{textesansboite}
As hinted by the names, the notions of ``locally finite support'' and ``collectively locally finite support'' are closely related. This link is made clearer by the following lemma.
\end{textesansboite}

\begin{lemmesansboite}{}
\label{LFS implique CLFS}%
Let $(\varphi _{d} )\in \prod _{d\in \mathbb{R}} E_{d}$. Then $(\varphi _{d} )_{d\in \mathbb{R}}$ has locally finite support iff $(\iota (\varphi _{d} ))_{d\in \mathbb{R}}$ has collectively locally finite support.
\end{lemmesansboite}

\begin{demo}
For~any $d\in \mathbb{R}$, we have $\iota (\varphi _{d} )=( \widetilde{\varphi }_{d,d'} )_{d'\in \mathbb{R}}$ where $\widetilde{\varphi }_{d,d'} :=\varphi _{d}$ if $d'=d$ and $\widetilde{\varphi }_{d,d'} :=0$ otherwise. Hence $\operatorname{supp}(\iota (\varphi _{d} ))$ is equal to $\{d\}$ if $\varphi _{d} \neq 0$ and to $\varnothing $ otherwise. So, for~any bounded set $B\subset \mathbb{R}$, $\{d\in \mathbb{R} \mid \operatorname{supp}(\iota (\varphi _{d} ))\cap B\neq  \varnothing  \}$ is equal to $\{d\in \mathbb{R} \mid \varphi _{d} \neq 0 \text{ and } d\in B\} =\operatorname{supp}((\varphi _{d} )_{d} )\cap B$. If one of these sets is finite, then so is the other one.
\end{demo}

\begin{textesansboite}
\Cref{LFS implique CLFS} allows us to consider $\fsumd \iota (\varphi _{d} )$ in the next result. Points 1 and 2 of \cref{SF 1D de composantes pures} are two variants of the same idea.
\end{textesansboite}

\begin{lemme}{}
\label{SF 1D de composantes pures}%
\begin{enumerate}
    \item For~any $(\varphi _{d} )\in \prod _{d\in \mathbb{R}} E_{d}$ that has a locally finite support, $(\iota (\varphi _{d} ))_{d\in \mathbb{R}}$ has collectively locally finite support and $\fsumd \iota (\varphi _{d} )=(\varphi _{d} )_{d\in \mathbb{R}}$. Using the abuse of notation on $\iota $ mentioned earlier, this rewrites as:
    \[\fsumd \varphi _{d} =(\varphi _{d} )_{d\in \mathbb{R}} .\]
    \item Let $A\subset \mathbb{R}$ be a locally finite set (i.e.\ $A\cap B$ is finite for~any bounded set $B\subset \mathbb{R}$) and $(\varphi _{d} )_{d\in A} \in \prod _{d\in A} E_{d}$. Let us denote $\widetilde{\varphi }_{d} := \varphi _{d}$ for~any $d\in A$ and $\widetilde{\varphi }_{d} := 0$ for~any $d\in \mathbb{R} \setminus A$. Then $(\iota (\varphi _{d} ))_{d\in A}$ has collectively locally finite support, $( \widetilde{\varphi }_{d} )_{d\in \mathbb{R}}$ has locally finite support and $\fsum_{d\in A} \iota (\varphi _{d} )=( \widetilde{\varphi }_{d} )_{d\in \mathbb{R}}$. Using again the abuse of notation on $\iota $, this rewrites as:
    \[\fsum_{d\in A} \varphi _{d} =( \widetilde{\varphi }_{d} )_{d\in \mathbb{R}} .\]
\end{enumerate}\fbeq
\end{lemme}

\begin{demo}
\begin{enumerate}
    \item For~any $q,d\in \mathbb{R}$, let $\widetilde{\varphi }_{q,d} := \delta _{q,d} \varphi _{q}$. Then, for~any $q\in \mathbb{R}$, $( \widetilde{\varphi }_{q,d} )_{d\in \mathbb{R}}$ belongs to $\prod _{d\in \mathbb{R}} E_{d}$, has finite support and its sum in $E$ is $\varphi _{q}$, so $( \widetilde{\varphi }_{q,d} )_{d\in \mathbb{R}} =\iota (\varphi _{q} )$. Hence \cref{def: fsum} gives:
    \[\fsum_{q\in \mathbb{R}} \iota (\varphi _{q} ) = \fsum_{q\in \mathbb{R}} ( \widetilde{\varphi }_{q,d} )_{d\in \mathbb{R}} = \bigg(\sum _{q\in \mathbb{R}} \widetilde{\varphi }_{q,d} \bigg)_{\!d\in \mathbb{R}} = (\varphi _{d} )_{d\in \mathbb{R}} .\]
    \item By~\cref{proprietes de fsum} item 4, we have
    \[\fsum_{d\in A} \iota (\varphi _{d} ) = \fsum_{d\in A} \iota ( \widetilde{\varphi }_{d} ) + \fsum_{d\in \mathbb{R} \setminus A} \iota (   \widetilde{\varphi }_{d} ) = \fsumd \iota ( \widetilde{\varphi }_{d} ) .\]
    We can then apply Point 1.
\end{enumerate}
\end{demo}

\begin{remarque}
Like we did in \cref{sec: C[[X]]}, in practice we will write formal series as $\fsum_{d\in A} \varphi _{d}$ way more often than as $(\varphi _{d} )_{d\in \mathbb{R}}$, because the first notation is more natural for a series.
\end{remarque}

\begin{exemplenumerote}
\label{ex: trucs qui sont ou pas des series formelles}%
Let us give some examples of what is and what is~not a formal series. For~any $d\in \mathbb{R}$, let $E_{d} := \{x\in \mathbb{R}_+^* \mapsto x^{d} P(\ln x) \mid P\in \mathbb{C}[T]\}\subset  \mathcal{C}^{0} (\mathbb{R}_+^*,\mathbb{C})$ and let $E:=\sum _{d\in \mathbb{R}} E_{d}$ be the vector space generated by the spaces $E_{d}$. This example will correspond to the traces on $\Gamma $ of the corner singularities (where $\Gamma $ is identified with $\mathbb{R}_+^*$).
\begin{itemize}
    \item $\fsum_{n\in \mathbb{Z}} x^{n}$, $\fsum_{n\in \mathbb{Z}} x^{n\pst}$ are well defined formal series because their support are $\mathbb{Z}$ and $\TZ$ resp., which are locally finite. Similarly, $\fsum_{p\in \NN} x^{p}$ is a formal series. Although, when $\pst \not\in \mathbb{Q}$, the points of $\NN:=\mathbb{N}+\pst \mathbb{N}$ become more and more dense at infinity (i.e.\ $\operatorname{card}(\NN \cap [n,n+1]) \rightarrow \infty $ when $n\rightarrow \infty $), $\NN$ is still locally finite.
    \item However, ``$\fsum_{n=1} ^{\infty } x^{ \frac{1}{n} }$'' is~not a formal series as its support is $\{ \frac{1}{n} \mid n\in \mathbb{N}^*\}$, which is~not locally finite. When $\pst \not\in \mathbb{Q}$, ``$\fsum_{d\in \mathbb{Z}+\TZ} x^{d}$'' is neither a formal series as $\mathbb{Z}+\TZ$ is dense in $\mathbb{R}$.
    \item $\fsumn n!\, x^{n\pst}$ is a formal series: the very quick growth of $n!$ is~not a problem as no convergence is involved in the definition of formal series. 
    \item ``$\fsumn x \ln ^{n} \!x$'' is~not a valid formal series because we are trying to sum an infinite number of terms that belong to the same $E_{d}$ (with $d=1$). Note that the problem is~not the unboundedness of the degree in $\ln x$, as $\fsumn x^{n\pst} \ln ^{n^2 } \!x$ is a well defined formal series. Similarly, $\fsumn x^{n\pst} \big(\sum_{k=1}^n \ln ^{k} \!x \big) = \fsum_{k=0} ^{\infty } \fsum_{n=k} ^{\infty } x^{n\pst} \ln ^{k} \!x$ is well defined. Here the number of terms that belong to a same space $E_{d}$ may be arbitrarily large, but what matters is that it~is finite for~any $d$.
    \item Small exercise: Is ``$\fsum_{n=1} ^{\infty } x^{ \frac{1}{n} } \ln ^{n} x$'' a well defined formal series ? What about ``$\fsum_{n=1} ^{\infty } x^{n} \ln ^{ \frac{1}{n} } x$''?\footnote{Answer: No and no because the support of the first one is $\{ \frac{1}{n} \mid n\in \mathbb{N}^*\}$ and the $x^{n} \ln ^{ \frac{1}{n} } \!x$ do~not even belong to the spaces $E_{d}$ (we only consider integer powers of $\ln x$).}
    \item To finish, let us give an extreme example. The following expression is a valid formal series, although the density of the support $\{\ln n \mid n\in \mathbb{N}^*\}$ explodes at infinity, the coefficients are huge, and the number of terms in $E_{\ln n}$ grows very quickly when $n\rightarrow \infty $.
    \[\fsum_{n=1} ^{\infty } x^{\ln n} \bigg(\sum _{k=0} ^{2^{n}} k!^{n} \ln ^{k} \!x\bigg)\]
\end{itemize}\fbi
\end{exemplenumerote}

\begin{lemme}{}
\label{fsum de sev des Ed}%
Let $A\subset \mathbb{R}$ and, for~any $d\in A$, let $\tilde{E}_{d}$ be a vector subspace of $E_{d}$. We have
\[\fsum_{d\in A} \tilde{E}_{d} =\fsum_{d\in A} \iota ( \tilde{E}_{d} ) = \left\{(\varphi _{d} )_{d\in \mathbb{R}} \ \bigg|\ \begin{array}{@{\,}l@{\,}} 
    (\forall d\in A,\, \varphi _{d} \in  \tilde{E}_{d} ) \text{ and } (\forall d\in \mathbb{R} \setminus A,\, \varphi _{d} =0)\\
    \text{and }  (\varphi _{d} ) \text{ has locally finite support} 
\end{array} \right\}\]
where in the left-hand side, we use the abuse of notation of dropping $\iota $.\\
Moreover, $\fsum_{d\in A} \tilde{E}_{d}$ is closed under $\fsum$, in the following sense: for~any set $I$ and any $(\Phi _{i} )_{i\in I} \in \big(\fsum_{d\in A} \tilde{E}_{d} )^{I}$ that has collectively locally finite support, $\fsum_{i\in I} \Phi _{i} \in \fsum_{d\in A} \tilde{E}_{d}$.
\end{lemme}

\begin{demo}
The first statement is a direct consequence of \cref{def: fsum d'ensembles,SF 1D de composantes pures} item 2. Next, let $I$ be a set and let $(\Phi _{i} )_{i\in I} \in \big(\fsum_{d\in A} \tilde{E}_{d} )^{I}$ that has collectively locally finite support. For~any $i\in I$, let us denote $(\varphi _{i,d} )_{d\in \mathbb{R}} :=\Phi _{i}$. By~\cref{def: fsum}, we have $\fsum_{i\in I} \Phi _{i} =\big(\sum _{i\in I} \varphi _{i,d} \big)_{d\in \mathbb{R}}$ and this family has locally finite support by \cref{les SF sont stables par fsum}.
\begin{itemize}
    \item If $d\in A$, then $\sum _{i\in I} \varphi _{i,d}$ has only a finite number of non-zero terms and these terms all belong to $\tilde{E}_{d}$. Since $\tilde{E}_{d}$ is a vector space, we deduce $\sum _{i\in I} \varphi _{i,d} \in  \tilde{E}_{d}$.
    \item If $d\in \mathbb{R} \setminus A$, then all the $\varphi _{i,d}$ vanish, so $\sum _{i\in I} \varphi _{i,d} =0$.
\end{itemize}
Using the first part of the present lemma, we deduce $\fsum_{i\in I} \Phi _{i} \in \fsum_{d\in A} \tilde{E}_{d}$.
\end{demo}

\begin{corollaire}{}
\label{ecriture avec fsum de l'ensemble des SF 1D}%
The set of formal series on $(E_{d} )_{d\in \mathbb{R}}$, i.e.\ the subset of $\prod _{d\in \mathbb{R}} E_{d}$ made of the families that have a locally finite support, is equal to
$\displaystyle \fsumd E_{d} $.
\end{corollaire}

\begin{textesansboite}
In~practice we will use $\fsumd E_{d}$ to designate the set of formal series on $(E_{d} )_{d\in \mathbb{R}}$. In summary, we have the following sequence of inclusions for~any $d'\in \mathbb{R}$:
\[\xymatrixrowsep{5.5ex}\xymatrixcolsep{5.5ex}\xymatrix{
	E_{d'} \ar[r] &\displaystyle \bigoplus_{d\in \mathbb{R}} E_{d} \ar[r]^{\iota } &\displaystyle \fsumd E_{d} \ar[r] &\displaystyle \prod _{d\in \mathbb{R}} E_{d}
}\]
In general all these inclusions are strict.
\end{textesansboite}

\subsection{Operators with a degree}

\begin{textesansboite}
Let $F$ be another $\mathbb{C}$-vector space and let $(F_{d} )_{d\in \mathbb{R}}$ be a family of subspaces of $F$ s.t.\ $F= \bigoplus_{d\in \mathbb{R}} F_{d}$. 
\end{textesansboite}

\begin{definition}{operators that have a degree}
\label{def: degre d'un operateur}%
Let $f:E \rightarrow  F$ and $d_{0} \in \mathbb{R}$. We say that $f$ has degree $d_{0}$ w.r.t.\ $(E_{d} )$ and $(F_{d} )$ iff $f$ is linear and:
\[\forall d\in \mathbb{R},\  \forall \varphi \in E_{d} ,\quad  f(\varphi )\in F_{d+d_{0}} .\]
In this case we denote $\deg f:=d_{0}$. If there is no ambiguity, ``w.r.t.\ $(E_{d} )$ and $(F_{d} )$'' can be omitted (which will be the case in practice).
\end{definition}
\nomenclature[13-14-16]{$\deg $}{Degree of an operator}

\begin{textesansboite}
\cref{fig: degre d'un operateur} gives a schematic illustration of the action of an operator that has a degree.
\end{textesansboite}

\begin{figure}[h]
\def\widthFig{4}
\def\heightFig{0.9}
\def\degfFig{1}
\def\icourbe{(exp(\i/4)-1)*4}
\def\diFig{(\icourbe/3.2-1)*\widthFig}
\begin{tikzpicture}
	\draw [>=stealth,->] (-\widthFig,\heightFig) -- (\widthFig,\heightFig) node[above]{$\;\ d$} node[right]{$\;\ \ \Big\}\,E=\bigoplus\limits_{d\in\mathbb{R}}^{\rule{0ex}{1ex}} E_d$};
	\draw [>=stealth,->] (-\widthFig,0) -- (\widthFig,0) node[above]{$\;\ d$} node[right]{$\;\ \ \Big\}\,F=\bigoplus\limits_{d\in\mathbb{R}}^{\rule{0ex}{1ex}} F_d$};
	\draw  (-\widthFig-1,0) node{\ }; 
	\foreach \i in {1,2,3} {
		\draw ({\diFig},\heightFig) node{\footnotesize$\bullet$} node[above]{$E_{d_\i}$};
		\draw ({\diFig},0) node{\footnotesize$\bullet$} node[below]{$F_{d_\i}$};
		\draw ({\diFig+\degfFig},0) node{\footnotesize$\bullet$} node[below]{$F_{d_\i +d_0}$};
		\draw[arrows={->[length=6pt,width=6pt]}, dashed] ({\diFig},\heightFig) -- node[midway,right]{\rule[-6pt]{0pt}{1ex}$f$} ({\diFig+\degfFig},0);}
\end{tikzpicture}\par
\longCaption{Schematic illustration of an operator that has degree $d_{0}$}{ (here $d_{0} >0$)}
\label{fig: degre d'un operateur}
\end{figure}

\begin{exemplenumerote}
\label{ex: operateurs avec un degre}%
In practice, $E_{d}$ and $F_{d}$ will be sets of functions that are roughly $d$-homogeneous (up to some logarithms), for example $E_{d} =F_{d} = \{x\mapsto  x^{d} P(\ln x) \mid P\in \mathbb{C}[T]\}$. In this setting, homogeneous differential operators have a degree that is equal to minus their derivation order, e.g.\ $\deg (\partial _{x} )=-1$. In \cref{sec: operateurs de resolution dans A}, we also define operators on corner singularities that build solutions of a canonical Poisson problem in $\Pi $. These operators have a degree and are extremely important to build the behaviors of the fields around the corner.
\end{exemplenumerote}

\begin{remarquenumerotee}
\label{rq: deg d'une composition d'op}%
Let $G=\bigoplus_{d\in \mathbb{R}} G_{d}$ be another $\mathbb{C}$-vector space. If $f:E\rightarrow F$ and $g:F\rightarrow G$ have a degree, then $g\circ f$ has a degree and $\deg (g\circ f )=\deg f +\deg g$.
\end{remarquenumerotee}

\begin{textesansboite}
It~is easy to see that, for~any $(\varphi _{d} )\in \prod _{d\in \mathbb{R}} E_{d}$ that has a locally finite support, $\big(f(\varphi _{d-\deg f} )\big)_{d} \in \prod _{d\in \mathbb{R}} F_{d}$ also has a locally finite support. This allows to introduce the following definition.
\end{textesansboite}

\begin{definition}{extension of an operator that has a degree}
\label{def: extension d'un op a degre aux SF 1D}
Let $f:E \rightarrow  F$ be a map that has a degree (w.r.t.\ $(E_{d} )$ and $(F_{d} )$). We define
\[\tilde{f} : \left(\begin{array}{clc}
    \displaystyle \fsumd E_{d} &\rightarrow &\displaystyle \fsumd F_{d} \\[1ex]
    (\varphi _{d} )_{d\in \mathbb{R}} &\mapsto & \big(f(\varphi _{d-\deg f} )\big)_{d\in \mathbb{R}}
\end{array}\right) .\]
We call it the extension of $f$ to $\fsum_{d} E_{d}$. In~practice we will denote $\tilde{f}$ as $f$ by abuse of notation.
\end{definition}

\begin{textesansboite}
Let us explain why it makes sense to call $\tilde{f}$ the ``extension of $f$'' and to denote it like $f$. First, let us show that $\tilde{f} \circ \iota _{E} =\iota _{F} \circ f$, where $\iota _{E}$ is the $\iota $ of $E$ and $\iota _{F}$ is the $\iota $ of $F$. Let $\varphi \in E$ and $(\varphi _{d} )_{d\in \mathbb{R}} :=\iota _{E} (\varphi )$. By~definition of $\iota _{E}$, $(\varphi _{d} )$ has a finite support and $\varphi =\sum _{d\in \mathbb{R}} \varphi _{d}$. Hence $f(\varphi ) =\sum _{d} f(\varphi _{d} ) = \sum _{d} f(\varphi _{d-\deg f} )$ where we have $f(\varphi _{d-\deg f} ) \in F_{d}$ for~any $d$. Since $\iota _{F} (f(\varphi ))$ is, by~definition, the unique element of $\prod _{d\in \mathbb{R}} F_{d}$ with finite support whose sum is equal to $f(\varphi )$, we deduce $\iota _{F} (f(\varphi ))=\big(f(\varphi _{d-\deg f} )\big)_{d\in \mathbb{R}}$. Therefore:
\[\tilde{f} \circ \iota _{E} (\varphi ) = \tilde{f} \big((\varphi _{d} )_{d\in \mathbb{R}} \big) = \big(f(\varphi _{d-\deg f} ) \big)_{d\in \mathbb{R}} = \iota _{F} \circ f(\varphi ).\]
People in algebra like to write this kind of commutation property as a commutative diagram:
\[\xymatrixrowsep{5.5ex}\xymatrixcolsep{5.5ex}\xymatrix{
	E \ar[r]^{f} \ar[d]^-{\iota _{E} } & F \ar[d]^-{\iota _{F} } \\
    \displaystyle \fsumd E_{d} \ar[r]_{ \tilde{f} } & \displaystyle \fsumd F_{d}
}\]
Saying that a diagram is commutative means that, for~any two paths of arrows that share the same starting and ending nodes, the composition of the maps of their arrows coincide.
Now, if we consider $\iota _{E}$ as an inclusion $E\subset \fsum_{d} E_{d}$ by abuse of notation, and similarly for $\iota _{F}$, then $\tilde{f} \circ \iota _{E} =\iota _{F} \circ f$ simply becomes $\tilde{f} _{|E} = f$. This justifies the terminology ``extension'' and the non-ambiguity of the notation $f:= \tilde{f}$.\\

We can rewrite the definition of the extension of $f$ in a more natural manner. Let $(\varphi _{d} )\in \fsumd E_{d}$. By~\cref{SF 1D de composantes pures}, we have $(\varphi _{d} )=\fsumd \varphi _{d}$ and $\big(f(\varphi _{d-\deg f} )\big)_{d} =\fsumd f(\varphi _{d-\deg f} ) =\fsumd f(\varphi _{d} )$, where we used \cref{proprietes de fsum} item 5 in the last equality. Therefore (using the abuses of notations):
\[f \bigg(\fsumd \varphi _{d} \bigg) = \fsumd f(\varphi _{d} ).\]
As shown by the next result, this property generalizes to $\fsum$-sums of formal series. It can be interpreted as a commutation of $\fsum$ with operators that have a degree, as well as a kind of super-linearity of these operators.
\end{textesansboite}

\begin{proposition}{}
\label{fsum commute avec les op a degre dans les SF 1D}
Let $f:E\rightarrow F$ be a map that has a degree. Let $I$ be a set and $(\Phi _{i} )_{i\in I} \in \big(\fsumd E_{d} \big)^{I}$ be a family with collectively locally finite support. Then $(f(\Phi _{i} ))_{i\in I}$ has collectively locally finite support and
\[f\bigg(\fsum_{i\in I} \Phi _{i} \bigg) = \fsum_{i\in I} f(\Phi _{i} ).\]\fbeq
\end{proposition}

\begin{demo}
For~any $i\in I$, we have $\operatorname{supp}(f(\Phi _{i} )) \subset  \operatorname{supp}(\Phi _{i} )+\deg f$. Therefore, for~any bounded set $B\subset \mathbb{R}$, $\{i\in I \mid \operatorname{supp}(f(\Phi _{i} ))\cap B \neq  \varnothing \} \subset  \{i\in I \mid (\operatorname{supp}(\Phi _{i} )+\deg f)\cap B \neq  \varnothing \} = \{i\in I \mid \operatorname{supp}(\Phi _{i} )\cap (B-\deg f) \neq  \varnothing \}$ which is finite by assumption. So~$(f(\Phi _{i} ))$ has collectively locally finite support. Now, for~any $i\in I$, let $(\varphi _{i,d} )_{d\in \mathbb{R}} :=\Phi _{i}$. We have:
\begin{align*}
    \tilde{f}\bigg(\fsum_{i\in I} \Phi _{i} \bigg) &= \tilde{f}\left( \bigg(\sum _{i\in I} \varphi _{i,d} \bigg)_{\!d\in \mathbb{R}} \right) && \text{by~definition of $\fsum$ (\cref{def: fsum})}  \\
    &= \left(f\bigg(\sum _{i\in I} \varphi _{i,d-\deg f} \bigg)\right)_{\!d\in \mathbb{R}} && \text{by~definition of }  \tilde{f} \\
    &= \bigg(\sum _{i\in I} f(\varphi _{i,d-\deg f} )\bigg)_{\!d\in \mathbb{R}} && \text{by linearity of $f$ on E}  \\
    &= \fsum_{i\in I} \big(f(\varphi _{i,d-\deg f} )\big)_{d\in \mathbb{R}} && \text{by~definition of }  \fsum \\
    &= \fsum_{i\in I} \tilde{f} (\Phi _{i} ) && \text{by~definition of }  \tilde{f}
\end{align*}
Finally, with the abuse of notation $f:= \tilde{f}$, we obtain the aimed formula.
\end{demo}

\begin{remarque}
If $f$ is a differential operator like in \cref{ex: operateurs avec un degre}, we can continue the analogy with distributions. Let us take $f:=\partial _{x}$ e.g.. For most convergence modes (pointwise, uniform, etc.), the convergence of a function series $\sum_{n=0}^\infty  u_{n}$, with differentiable functions $u_{n}$, does~not imply the convergence of $\sum_{n=0}^\infty  \partial _{x} u_{n}$, nor the equality $\partial _{x} (\sum_{n=0}^\infty  u_{n} )=\sum_{n=0}^\infty  \partial _{x} u_{n}$. However in the generalized frameworks of distributions and formal series, this implication becomes true and straightforward. On the one hand, if $(T_{n} )\in \mathcal{D}'(\mathbb{R})^{\mathbb{N}}$ and $\sum_{n=0}^\infty  T_{n}$ converges in $\mathcal{D}'(\mathbb{R})$, then $\sum_{n=0}^\infty  \partial _{x} T_{n}$ converges in $\mathcal{D}'(\mathbb{R})$ and $\partial _{x} (\sum_{n=0}^\infty  T_{n} )=\sum_{n=0}^\infty  \partial _{x} T_{n}$. On the other hand, \cref{fsum commute avec les op a degre dans les SF 1D} states that, if $(\Phi _{n} )\in (\fsum E_{d} )^{\mathbb{N}}$ and $\fsumn \Phi _{n}$ is well defined as a formal series, then $\fsumn \partial _{x} \Phi _{n}$ is well defined too and $\partial _{x} (\fsumn \Phi _{n} )=\fsumn \partial _{x} \Phi _{n}$.
\end{remarque}

\begin{exemple}
The space $\FCX$ presented in \cref{sec: C[[X]]} can be seen as a particular example of the generalized formal series of this section (if we forget $\FCX$'s multiplicative structure). It suffices to take $E_{d} :=\{x\mapsto  ax^{d} \mid a\in \mathbb{C}\}$ when $d\in \mathbb{N}$ and $E_{d} :=0$ when $d\in \mathbb{R} \setminus \mathbb{N}$. The isomorphism between $\FCX$ and $\fsum_{d} E_{d}$ is given by
\[\beta  : \left(\begin{array}{clc}
    \FCX &\rightarrow &\displaystyle \fsumd E_{d} \\[0.5ex]
    \displaystyle \fsumn a_{n} X^{n} &\mapsto &\displaystyle \fsumn a_{n} \.x^{n} = \big(\mathbf{1}_{\mathbb{N}} (d) \,a_{d} \.x^{d} \big)_{d\in \mathbb{R}}
\end{array}\right) \]
(where $x^{n}$ is a shortcut for $x\mapsto  x^{n}$), 
because $\beta $ is a linear bijection that commutes with $\fsum$. Moreover, the multiplication by $X$ is a particular case of an operator that has a degree. Indeed, if we define $\xi :\fsumd E_{d} \rightarrow  \fsumd E_{d}$ by $\xi \big((\varphi _{d} )_{d\in \mathbb{R}} \big):=(x\cdot \varphi _{d-1} )_{d\in \mathbb{R}}$, then $\xi $ has obviously degree 1 and we have: $\forall f\in \FCX,\;\hskip0pt plus 1.5pt \xi \circ \beta (f)=\beta (Xf)$. This commutation property can be represented by the following commutative diagram:
\[\xymatrixrowsep{5.5ex}\xymatrixcolsep{5.5ex}\xymatrix{
	\FCX \ar[r]^{\times X} \ar[d]^-{\beta }_{\cong} & \FCX \ar[d]^-{\beta }_{\cong} \\
    \displaystyle \fsumd E_{d} \ar[r]_{\xi } & \displaystyle \fsumd E_{d}
}\]
Hence the property ``$X \fsum_{i\in I} f_{i} =\fsum_{i\in I} Xf_{i}$'' in $\FCX$ stated in \cref{proprietes de fsum dans C[[X]]} is a particular case of \cref{fsum commute avec les op a degre dans les SF 1D}. Likewise, the operator $A$ studied in \cref{ex: 2e pb jouet dans C[[X]]} is also an operator with a degree and we are generalizing here the work that was done in that particular example.
\end{exemple}

\begin{textesansboite}
We end this section with some useful truncation operators.
\end{textesansboite}

\begin{definition}{truncation operators}
\label{def: troncature de SF}%
For~any $q\in \mathbb{R}$ and $\Phi =(\varphi _{d} )\in \fsumd E_{d}$, we define
\[T_{<q} (\Phi ):=\fsum_{d<q} \varphi _{d} , \qquad  T_{\leqslant q} (\Phi ):=\fsum_{d\leqslant q} \varphi _{d} , \qquad  T_{>q} (\Phi ):=\fsum_{d>q} \varphi _{d} , \qquad  T_{\geqslant q} (\Phi ):=\fsum_{d\geqslant q} \varphi _{d} .\]\fbeq
\end{definition}
\nomenclature[29-60]{$T_{\leqslant d}$, $T_{\geqslant d}$, $T_{<d}$, $T_{>d}$}{Truncations operator over formal series with index in $\mathbb{R}$}

\begin{remarquenumerotee}
\label{rq: fsum commute avec les troncatures T}%
It~is easy to verify the following facts.
\begin{itemize}
    \item We have $T_{<q} = \mathrm{id}-T_{\geqslant q}$ and $T_{>q} =  \mathrm{id}-T_{\leqslant q}$.
    \item The four truncation operators have degree 0. For instance, $T_{<q}$ is the extension to formal series of $f:E\rightarrow E$ that maps each $\varphi _{d} \in E_{d}$ to $\mathbf{1}_{(-\infty ,q)}(d) \cdot \varphi _{d} \in E_{d}$. Therefore, by \cref{fsum commute avec les op a degre dans les SF 1D}, the truncation operators commute with $\fsum$.
\end{itemize}
\end{remarquenumerotee}

\begin{exemple}
The truncations operators will be used in \cref{chap: construction du DA} to link formal series with asymptotic expansions. For example, if we take $E_{d} :=\{x\mapsto  ax^{d} \mid a\in \mathbb{C}\}$ for~any $d\in \mathbb{R}$, then the rigorous formulation of ``{\it$\fsumn x^{n} /n!$ represents the behavior at 0 of $x\mapsto e^{x}$\/}'' is:
\[\forall d\in \mathbb{R},\quad \ e^{x} = T_{\leqslant d} \left(\fsumn \frac{x^{n} }{n!} \right) + o(x^{d} ) \quad  \text{when }  x\rightarrow 0.\]
Indeed, thanks to the identification $E=\iota (E)$, $T_{\leqslant d} \big(\fsum_{n} x^{n} /n! \big)$ is equal to the partial sum $\sum _{n=0} ^{\lceil {d}\rceil } \frac{x^{n} }{n!}$, which is a function $\mathbb{C} \rightarrow \mathbb{C}$ in the classical sense.
\end{exemple}

\subsection{Geometric series of operators}
\label{sec: series d'op de SF 1D}

\begin{textesansboite}
Let us denote $E^{+}$, resp.\ $E^{-}$, the subset of $\fsumd E_{d}$ made of elements whose support is bounded from below, resp.\ from above. For instance:
\[E^{+} = \bigg\{\Phi  \in \fsumd E_{d} \ \bigg|\ \exists d_{0} \in \mathbb{R},\; \operatorname{supp}(\Phi )\subset [d_{0} ,\infty )\bigg\}.\]
This can be rewritten as
\[E^{+} = \bigcup _{d_{0} \in \mathbb{R}} \fsum_{d\in [d_{0} ,\infty )} E_{d} \qquad  \text{ and } \qquad  E^{-} = \bigcup _{d_{0} \in \mathbb{R}} \fsum_{d\in (-\infty ,d_{0} ]} E_{d} .\]
Since an increasing union of vector subspaces is a vector subspace, $E^{+}$ and $E^{-}$ are vector subspaces of $\fsumd E_{d}$.
\end{textesansboite}

\begin{remarquenumerotee}
\label{rq: op de troncature et E+-}%
Let $\Phi \in \fsumd E_{d}$ and $q\in \mathbb{R}$.
\begin{itemize}
    \item We have $T_{<q} (\Phi ) \in E^{-}$ and $T_{>q} (\Phi ) \in E^{+}$, and likewise with large inequalities.
    \item If $\Phi  \in E^{+}$, then $T_{<q} (\Phi ) \in E$ and, if $\Phi  \in E^{-}$, then $T_{>q} (\Phi ) \in E$ (with the abuse of notation $E:=\iota (E)$). The same applies with large inequalities. In these cases, the truncation is simply a partial sum of the series $\Phi $. 
\end{itemize}
\end{remarquenumerotee}

\begin{remarquenumerotee}
\label{rq: SF 1D unilaterales non stables par fsum}%
$E^{\pm }$ is~not closed under $\fsum$. For example, if $\varphi _{n} \in E_{n} \setminus \{0\}$ for~any $n\in \mathbb{Z}$, then $\varphi _{n} \in E^{+} \cap E^{-}$ for~any $n$ and $\fsum_{n\in \mathbb{Z}} \varphi _{n}$ is well defined, but $\fsum_{n\in \mathbb{Z}} \varphi _{n}$ does~not belong to $E^{+}$, neither to $E^{-}$. However, \cref{existence des serie geometrique d'op 1D} gives a very useful situation in which a $\fsum$-sum of elements of $E^{\pm }$ stays in $E^{\pm }$.
\end{remarquenumerotee}

\begin{proposition}{}
\label{existence des serie geometrique d'op 1D}%
Let $f:\fsumd E_{d} \rightarrow \fsumd E_{d}$ be a sum of a finite number of operators that have a positive (resp.\ negative) degree. Let $f^{n}$ denote the iterated compositions of $f$.\\
For~any $\Phi \in E^{+}$ (resp.\ $\Phi \in E^{-}$), we have $f(\Phi )\in E^{+}$ (resp.\ $f(\Phi )\in E^{-}$), $\big(f^{n} (\Phi )\big)_{n\in \mathbb{N}}$ has collectively locally finite support and
\[\fsumn f^{n} (\Phi ) \in E^{+} , \qquad  \text{resp.}  \quad   \fsumn f^{n} (\Phi ) \in E^{-} .\]\fbeq[0.8]
\end{proposition}

\begin{demo}
We treat only the positive case, as the negative case is entirely similar. Let $k\in \mathbb{N}$ and let $f_{1} ,\dots ,f_{k}$ be linear maps $E\rightarrow E$ that have a positive degree s.t.\ $f=\sum_{i=1}^k f_{i}$. Let $d_{f} :=\min \{\deg  f_{i} \mid i\in [\![1,k]\!]\}$. Since $\Phi \in E^{+}$, there exists $d_{\Phi } \in \mathbb{R}$ s.t.\ $\operatorname{supp}(\Phi )\subset [d_{\Phi } ,\infty )$. So, for~any $i\in [\![1,k]\!]$, we have $\operatorname{supp}(f_{i} (\Phi ))\subset \operatorname{supp}(\Phi )+\deg  f_{i} \subset [d_{\Phi } +d_{f} ,\infty )$. This implies $\operatorname{supp}(f(\Phi ))\subset [d_{\Phi } +d_{f} ,\infty )$ and $f(\Phi )\in E^{+}$. Moreover, iterating this result gives:
\[\forall n\in \mathbb{N}, \qquad  \operatorname{supp}(f^{n} (\Phi ))\subset [d_{\Phi } +nd_{f} ,\infty ).\]
Hence for~any bounded set $B\subset \mathbb{R}$, $\{n\in \mathbb{N} \mid \operatorname{supp}(f^{n} (\Phi ))\cap B \neq  \varnothing \}$ is included in $\{n\in \mathbb{N} \mid d_{\Phi } +nd_{f} \leqslant \sup B\}$ which is finite because $d_{f} >0$. Thus $(f^{n} (\Phi ))_{n}$ has collectively locally finite support. Finally it~is easy to see that $\operatorname{supp}\big(\fsumn f^{n} (\Phi )\big)\subset  \bigcup _{n\in \mathbb{N}} \operatorname{supp}(f^{n} (\Phi )) \subset [d_{\Phi } ,\infty )$, so $\fsum_{n} f^{n} (\Phi ) \in E^{+}$.
\end{demo}

\begin{figure}[h]
\def\widthFig{4.3}
\def\dminFig{-\widthFig+1} 
\def\pasVFig{0.55} 
\def\pasHFig{0.25} 
\def\NFig{3} 
\def\NmoinsUnFig{2}
\def\degfFig{2}
\def\ij{(\i*\degfFig+\j)}
\def\Nj{(\NFig*\degfFig+\j)}
\def\basFig{-(\NFig+1.55)*\pasVFig} 
\begin{tikzpicture}
	\foreach \i in {0,...,\NFig} {
		\draw[>=stealth,->] (-\widthFig,-\i*\pasVFig) -- (\widthFig,-\i*\pasVFig) node[above]{$\;\ d$};}
	\draw[>=stealth,->] (-\widthFig,{\basFig}) -- (\widthFig,{\basFig}) node[above]{$\;\ d$};
	\foreach \i in {1,...,\NFig} {
		\draw (-\widthFig-0.2,{-(\i-0.5)*\pasVFig}) node[left]{$+$};}
	\draw (-\widthFig-0.2,{-(\NFig+0.5)*\pasVFig}) node[left]{\scalebox{0.75}{$\vdots$}\;};
	\draw (-\widthFig-0.2,{\basFig}) node[left]{$=$};
	\draw  (-\widthFig-1.2,0) node{\ }; 
	\draw (\widthFig,0) node[right]{$\;\ \ \big\}\; \fsum\varphi_d $};
	\draw (\widthFig,-\pasVFig) node[right]{$\;\ \ \big\}\; f\big(\fsum\varphi_d \big)$};
	\foreach \i in {2,...,\NFig} {
		\draw  (\widthFig,-\i*\pasVFig) node[right]{$\;\ \ \big\}\; f^\i\big(\fsum\varphi_d \big)$};}
	\draw (\widthFig,{(-\NFig*\pasVFig+\basFig)/2+0.05}) node[right]{$\qquad \quad$\scalebox{0.75}{$\vdots$}};
	\draw (\widthFig,{\basFig}) node[right]{$\;\ \ \big\}\; \fsum f^n \big(\fsum\varphi_d \big)$};
	\begin{scope}
		\clip (-\widthFig,0.1) rectangle (\widthFig-0.15,{\basFig-0.1});
		\foreach \j in {0,4,9,12,18,20,27} {
			\foreach \i in {0,...,\NFig} {
				\draw[dotted, gray] ({\dminFig+\ij*\pasHFig},{\basFig}) -- ({\dminFig+\ij*\pasHFig},{-\i*\pasVFig});
				\draw ({\dminFig+\ij*\pasHFig},{-\i*\pasVFig}) node{\footnotesize$\bullet$};}
			\foreach \i in {0,...,40} {
				\draw ({\dminFig+\ij*\pasHFig},{\basFig}) node{\footnotesize$\bullet$};}
			\foreach \i in {0,...,\NmoinsUnFig} {
				\draw[arrows={->[length=5pt,width=4pt]}, densely dashed] ({\dminFig+\ij*\pasHFig},{-\i*\pasVFig}) -- ({\dminFig+(\ij+\degfFig)*\pasHFig},{-(\i+1)*\pasVFig});}
			\draw[arrows={->[length=5pt,width=4pt]}, densely dashed] ({\dminFig+\Nj*\pasHFig},{-\NFig*\pasVFig}) -- ({\dminFig+(\Nj+\degfFig*0.6)*\pasHFig},{-(\NFig+0.6)*\pasVFig});} 
	\end{scope}
	\draw[<-,>=latex] ({\dminFig},0.05) -- ++(0,0.25) node[above]{$d_{\mathrm{inf}}$};
	\foreach \i in {0} {
		\foreach \j in {4} {
			\draw[<-,>=latex] ({\dminFig+\ij*\pasHFig+0.03},{-\i*\pasVFig+0.05}) to[bend left=25] ++(0.5,0.4) node[right]{$\varphi_d \neq0$};}
		\foreach \j in {15} {
			\draw[<-,>=latex] ({\dminFig+\ij*\pasHFig+0.03},{-\i*\pasVFig+0.05}) to[bend left=25] ++(0.5,0.4) node[right]{$\varphi_d =0$};}}
	\draw[<-,>=latex] ({\dminFig+4*\pasHFig+0.03},{\basFig-0.05}) to[bend right=25] ++(0.5,-0.35) node[right]{finite sum};
\end{tikzpicture}\par
\caption{Schematic illustration of $\fsumn f^{n}$ for $f$ a linear map that has a positive degree}
\label{fig: serie geometrique d'operateur}
\end{figure}

\begin{definition}{geometric series of operators}
\label{def: serie geometrique d'op 1D}%
Let $f:\fsumd E_{d} \rightarrow \fsumd E_{d}$ be a sum of a finite number of operators that have a positive (resp.\ negative) degree. We define $\fsumn f^{n} :E^{+} \rightarrow E^{+}$ (resp.\ $E^{-} \rightarrow E^{-}$) by:
\[\fsumn f^{n} : \ \Phi  \longmapsto  \fsumn f^{n} (\Phi ).\]
Moreover, if $k\in \mathbb{N}$, $f_{1} ,\dots ,f_{k}$ are distinct linear maps $E\rightarrow E$ with a positive degree s.t.\ $f=\sum_{i=1}^k f_{i}$ and $\mathcal{F}:=\{f_{1} ,\dots ,f_{k} \}$, we also denote
\[\langle \mathcal{F}\rangle  :=\langle f_{1} ,\dots ,f_{k} \rangle  :=\fsumn f^{n} .\]\fbeq
\end{definition}
\nomenclature[60-11a]{$\langle {\dots }\rangle $}{Geometric series of operators on formal series}

\begin{remarque}
Let $\Phi =(\varphi _{d} )_{d} \in \fsumd E_{d}$. One can verify that, with the notations of \cref{def: serie geometrique d'op 1D}, we have
\begin{align*}
    \langle \mathcal{F}\rangle (\Phi ) &= \fsumn \bigg(\sum _{g\in \mathcal{F}} g\bigg)^{n} (\Phi )\\
    &= \fsumn \sum _{(g_{1} ,\dots , g_{n} )\in \mathcal{F}^{n} } g_{1} \circ \cdots \circ g_{n} (\Phi ) \\
    &= \bigg( \sum _{n\in \mathbb{N}} \sum _{(g_{1} ,\dots , g_{n} )\in \mathcal{F}^{n} } g_{1} \circ \cdots \circ g_{n} (\varphi _{d-\Sigma _{1\leqslant i\leqslant n} \deg  g_{i} }) \bigg)_{\!d\in \mathbb{R}} .
\end{align*}
In the last line we used $\deg (g_{1} \circ \cdots \circ g_{n} )=\sum_{i=1}^n \deg  g_{i}$, see \cref{rq: deg d'une composition d'op}. As~you can see, $\langle \mathcal{F}\rangle (\Phi )$ can be quite a complex object. But thanks to the calculus tools that we introduce in this chapter, it will not be necessary to specify the terms of $\langle \mathcal{F}\rangle (\Phi )$. This is one of the main points of these tools.
\end{remarque}

\begin{textesansboite}
The next theorem is a key result that will allow us to solve equations among formal series in order to compute the behaviors of the fields around the corner. Indeed, these behaviors satisfy the same PDE as the fields themselves, but formulated among the formal series (\cref{ex: operateurs avec un degre} explains how differential operators can be defined on formal series). In \cref{ex: 2e pb jouet dans C[[X]]}, we proved a similar result on a particular case in $\FCX$. \Cref{inverse de id-f dans les SF 1D} gives the general setting in which the reasoning of this example applies.
\end{textesansboite}

\begin{theoreme}{operator inversion by Neumann series}
\label{inverse de id-f dans les SF 1D}%
Let $f:\fsumd E_{d} \rightarrow \fsumd E_{d}$ be a sum of a finite number of operators that have a positive (resp.\ negative) degree. The operator $ \mathrm{id}-f$ maps $E^{+}$ to $E^{+}$ (resp.\ $E^{-}$ to $E^{-}$), its restriction to $E^{+}$ (resp.\ $E^{-}$) is invertible and the inverse is
\[( \mathrm{id}-f)^{-1} = \fsumn f^{n} .\]\fbeq
\end{theoreme}

\begin{demo}
Let us treat only the positive case, the negative case being very similar. $( \mathrm{id}-f)(E^{+} )\subset E^{+}$ is a direct consequence of \cref{existence des serie geometrique d'op 1D}. The rest of the proof is very similar to \cref{ex: 2e pb jouet dans C[[X]]}. For~any $\Phi \in E^{+}$, we have:
\begin{align*}
    ( \mathrm{id}-f) \circ \bigg(\fsumn f^{n} \bigg)(\Phi ) &= \fsumn f^{n} (\Phi ) - \fsumn f^{n+1} (\Phi ) && \text{by \cref{fsum commute avec les op a degre dans les SF 1D}}  \\
    &= \fsum_{n\in \{0\}} f^{n} (\Phi ) + \fsum_{n=1} ^{\infty } f^{n} (\Phi ) - \fsum_{n=1} ^{\infty } f^{n} (\Phi ) && \text{by \cref{proprietes de fsum}, items 4,5}  \\
    &= \Phi  && \text{by \cref{proprietes de fsum}, item 1} 
\end{align*}
Similarly,
\begin{align*}
    \bigg(\fsumn f^{n} \bigg)\circ ( \mathrm{id}-f)(\Phi ) &= \fsumn (f^{n} (\Phi ) -f ^{n+1} (\Phi ) ) && \text{by~definition of $\fsumn f^{n}$ and linearity of }  f\\
    &= \fsumn f^{n} (\Phi ) - \fsumn f^{n+1} (\Phi ) && \text{by \cref{proprietes de fsum}, item 3}  \\
    &= \Phi  && \text{by the previous computation} 
\end{align*}
\end{demo}

\begin{proposition}{commutation of $\langle .\rangle $ and $\fsum$}
\label{fsum et le crochet commutent en 1D}%
Let $\mathcal{F}$ be a finite set of operators $E\rightarrow E$ that have a positive (resp.\ negative) degree. Let $I$ be a set and $(\Phi _{i} )_{i\in I}$ be a family of elements of $E^{+}$ (resp.\ $E^{-}$) that has collectively locally finite support. We assume that 
$\bigcup _{i\in I} \operatorname{supp}(\Phi _{i} )$ is bounded from below (resp.\ above).\\
Then $\fsum_{i\in I} \Phi _{i}$ belongs to $E^{+}$ (resp.\ $E^{-}$), $\big(\langle \mathcal{F}\rangle (\Phi _{i} )\big)_{i}$ has collectively locally finite support and
\[\langle \mathcal{F}\rangle \bigg(\fsum_{i\in I} \Phi _{i} \bigg) = \fsum_{i\in I} \,\langle \mathcal{F}\rangle (\Phi _{i} ).\]\fbeq[0.7]
\end{proposition}

\begin{remarque}
The proposition would be wrong without the assumption on $\bigcup _{i\in I} \operatorname{supp}(\Phi _{i} )$. In the positive case for instance, \cref{rq: SF 1D unilaterales non stables par fsum} gives an example where $\Phi _{i} \in E^{+}$ for~any $i\in I$, but $\fsum_{i\in I} \Phi _{i} \not\in E^{+}$. In this case, $\langle \mathcal{F}\rangle \big(\fsum_{i} \Phi _{i} \big)$ is~not even defined.\\
Moreover, one can verify that, if $I=\mathbb{N}$, the statement ``$(\Phi _{n} )$ has collectively locally finite support and $\bigcup _{n\in \mathbb{N}} \operatorname{supp}(\Phi _{n} )$ is bounded from below'' is equivalent to $\inf (\operatorname{supp}(\Phi _{n} ))\rightarrow \infty $ when $n\rightarrow \infty $.
\end{remarque}

\begin{demo}
Let us treat only the positive case. Let $d_{\mathcal{F}} :=\min \{\deg f \mid f\in \mathcal{F}\}$, let $d_{I} \in \mathbb{R}$ be s.t.\ $\bigcup _{i\in I} \operatorname{supp}(\Phi _{i} ) \subset [d_{I} ,\infty )$ and, for~any $i\in I$, let $d_{i} := \inf (\operatorname{supp}(\Phi _{i} ))$. We also denote $g:=\sum _{f\in \mathcal{F}} f$. Using \cref{les SF sont stables par fsum}, we get $\operatorname{supp}\!\big(\fsum_{i\in I} \Phi _{i} \big) \subset  \bigcup _{i\in I} \operatorname{supp}(\Phi _{i} )\subset [d_{I} ,\infty )$ which proves $\fsum_{i\in I} \Phi _{i} \in E^{+}$.\\
Moreover, like in the proof of \cref{existence des serie geometrique d'op 1D}, we have $\operatorname{supp}(g^{n} (\Phi _{i} )) \subset [d_{i} +nd_{\mathcal{F}} ,\infty )$ for~any $i\in I$ and $n\in \mathbb{N}$. So, for~any bounded set $B\subset \mathbb{R}$, $\{(i,n)\in I\times \mathbb{N} \mid \operatorname{supp}(g^{n} (\Phi _{i} ))\cap B \neq \varnothing  \}$ is included in $\{(i,n)\in I\times \mathbb{N} \mid d_{i} +nd_{\mathcal{F}} \leqslant \sup  B \}$. This last set is finite for the following reasons.
\begin{itemize}
    \item $d_{i} >-\infty $ for~any $i$ because $\Phi _{i} \in E^{+}$.
    \item $d_{\mathcal{F}} >0$.
    \item For~any $M\in \mathbb{R}$, $\{i\in I \mid d_{i} <M\}$ is finite. To prove this statement, let us take $M\in \mathbb{R}$. The set $B_{M} :=[d_{I} ,M]$ is bounded, so $I_{M} :=\{i\in I \mid \operatorname{supp}(\Phi _{i} )\cap B_{M} \neq  \varnothing \}$ is finite. Therefore $\{i\in I \mid d_{i} <M\}\subset \{i\in I \mid \operatorname{supp}(\Phi _{i} )\cap (-\infty ,M] \neq  \varnothing \}\subset I_{M}$ is finite.
\end{itemize}
Hence, $\big(g^{n} (\Phi _{i} )\big)_{(i,n)\in I\times \mathbb{N}}$ has collectively locally finite support. Therefore \cref{Fubini pour fsum} (Fubini) states
\begin{equation}{\label{eq: fsum et le crochet commutent en 1D preuve 1}}
    \fsumn \fsum_{i\in I} g^{n} (\Phi _{i} ) = \fsum_{i\in I} \fsumn g^{n} (\Phi _{i} ).
\end{equation}
\Cref{Fubini pour fsum} also states that the $\fsum$-sums in \cref{eq: fsum et le crochet commutent en 1D preuve 1} are well-defined, and in particular that the family $\big(\fsumn g^{n} (\Phi _{i} )\big)_{i\in I} = \big(\langle \mathcal{F}\rangle (\Phi _{i} ))_{i\in I}$ has collectively finite support. In addition, we have
\[g^{n} = \sum _{(f_{1} ,\dots ,f_{n} )\in \mathcal{F}^{n} } f_{1} \circ \cdots \circ f_{n} .\]
It~is easy to see that the maps $f_{1} \circ \cdots \circ f_{n}$ have a degree, so \cref{fsum commute avec les op a degre dans les SF 1D} allows us to swap $\fsum_{i\in I}$ and $g^{n}$ in the left hand side of \cref{eq: fsum et le crochet commutent en 1D preuve 1}, which gives
\[\fsumn g^{n} \bigg(\fsum_{i\in I} \Phi _{i} \bigg) = \fsum_{i\in I} \fsumn g^{n} (\Phi _{i} ).\]
\end{demo}

\section{Formal series with index in $\mathbb{R}^2 $ (involving $\feps$)}
\label{sec: SF 2D}

\begin{textesansboite}
In this section, we define the type of formal series that will represent the behavior of ``$\sum _{p,\ell } \peps^{p} \ln ^{\ell } \!\peps\, \bfu_{p,\ell }$'' and ``$\sum _{p,\ell } \peps^{p} \ln ^{\ell } \!\peps\, S_{p,\ell }$'' w.r.t.\ $\peps$ and $\bfr$. To describe this double asymptotic behavior, the index of these new formal series belongs to $\mathbb{R}^2 $.\\

We still consider a $\mathbb{C}$-vector space $E$ and a family $(E_{d} )_{d\in \mathbb{R}}$ of vector subspaces of $E$ s.t.\ $E=\bigoplus_{d\in \mathbb{R}} E_{d}$. Each $E_{d}$ is meant to be the set of polynomials in the formal indeterminate ``$\ln \varepsilon $'', whose coefficients are corner singularities that are $d$-homogeneous w.r.t.\ $\bfr$ (up to some $\ln \bfr$)\footnote{That is to say, of the form $r^{d} \sum _{j=0} ^{J} \ln ^{j} \!r \cdot  f_{j} (\theta )$ in $\Omega $ and $x^{d} \sum _{j=0} ^{J} \ln ^{j} \!x \cdot  g_{j} (\theta )$ in $\couche$, with $J\in \mathbb{N}$ and $f_{j} ,g_{j}$ some functions.}. In other words, the $E_{d}$ of this section should be the $E_{d} [\ln \varepsilon ]$ of \cref{sec: SF 1D} (see \cref{def: E[T]} for $E_{d} [\ln \varepsilon ]$). However, in this section, we will not look at the behavior in powers of $\ln \varepsilon $, which we know can be handled by \cref{def: E[T]}. We will focus on powers of $\varepsilon $. The powers of $\peps$ are more complicated than the powers of $\ln \peps$ as they may not be integers.\\

The formal series that we study in this section are the elements of
\[\prod _{(d,p)\in \mathbb{R}^2 } E_{d}\]
that have a locally finite support. $p$ will be the power of $\varepsilon $. The link with the formal series of \cref{sec: SF 1D} is made by the following linear injection
\[\isfud : (\varphi _{d} )_{d\in \mathbb{R}} \mapsto  \big(\delta _{p,0} \,\varphi _{d} \big)_{(d,p)\in \mathbb{R}^2 } \]
that maps $\fsumd E_{d}$ to $\prod _{(d,p)\in \mathbb{R}^2 } E_{d}$. It~is easy to see that the elements of $\iota\big(\fsumd E_{d} \big)$ have locally finite support. Like for $\iota :E\rightarrow \fsumd E_{d}$, we will consider $\isfud$ as a inclusion by abuse of notation. This abuse is legitimated by the following lemma.
\end{textesansboite}

\begin{lemme}{}
Let $I$ be a set and $(\Phi _{i} )_{i\in I} \in \big(\fsumd E_{d} \big)^{I}$ be a family of formal series. Then $(\Phi _{i} )_{i\in I}$ has collectively locally finite support (CLFS) iff $\big(\isfud(\Phi _{i} )\big)_{i\in I}$ has CLFS, and in this case we have
\[\fsum_{i\in I} \isfud(\Phi _{i} ) = \isfud\bigg(\fsum_{i\in I} \Phi _{i} \bigg).\]\fbeq
\end{lemme}

\begin{demo}
First, let us assume that $(\Phi _{i} )_{i}$ has CLFS. For~any $i\in I$, we have $\operatorname{supp}(\isfud(\Phi _{i} ))=\{0\}\times \operatorname{supp}(\Phi _{i} )$. So, for~any bounded set $B\subset \mathbb{R}^2 $, $I_{B} :=\{i\in I \mid \operatorname{supp}(\isfud(\Phi _{i} ))\cap B\neq \varnothing \}$ is equal to $\{i\in I \mid \operatorname{supp}(\Phi _{i} )\cap B' \neq \varnothing \}$ with $B':=\{d\in \mathbb{R} \mid (d,0)\in B\}$. $B'$ is bounded, so $I_{B}$ is finite and $\big(\isfud(\Phi _{i} )\big)_{i}$ has CLFS.\\
Conversely, let us assume that $\big(\isfud(\Phi _{i} )\big)_{i}$ has CLFS. Then, for~any bounded set $B\subset \mathbb{R}$, the set $\{i\in I \mid \operatorname{supp}(\Phi _{i} )\cap B\neq \varnothing \}$ is equal to $\{i\in I \mid \operatorname{supp}(\isfud(\Phi _{i} ))\cap (B\times \{0\}) \neq \varnothing \}$ which is finite, so $(\Phi _{i} )_{i}$ also has CLFS.\\
Next, let us denote $(\varphi _{i,d} ):=\Phi _{i}$ for~any $i\in I$. We have
\begin{align*}
    \fsum_{i\in I} \isfud(\Phi _{i} ) &= \fsum_{i\in I} \big(\delta _{p,0} \,\varphi _{i,d} \big)_{(d,p)\in \mathbb{R}^2 } && \text{by~definition of } \isfud \\
    &= \bigg(\sum _{i\in I} \delta _{p,0} \,\varphi _{i,d} \bigg)_{(d,p)\in \mathbb{R}^2 } && \text{by \cref{def: fsum}}  \\
    &= \isfud\left( \bigg(\sum _{i\in I} \varphi _{i,d} \bigg)_{(d,p)\in \mathbb{R}^2 } \right) && \text{by~definition of } \isfud \\
    &= \isfud\bigg(\fsum_{i\in I} \Phi _{i} \bigg) && \text{by \cref{def: fsum}} 
\end{align*}
\end{demo}

\subsection{Multiplication by powers of $\feps$}
\label{sec: feps dans les SF 2D}

\begin{definition}{powers of $\varepsilon $}
\label{def: feps}%
For~any $p\in \mathbb{R}$ and $(\varphi _{d,p'} )_{(d,p')} \in \prod _{(d,p')\in \mathbb{R}^2 } E_{d}$ that has locally finite support, we define
\[\varepsilon ^{p} \cdot  (\varphi _{d,p'} )_{(d,p')} := (\varphi _{d,p'-p} )_{(d,p')}.\]
Here $\varepsilon $ does~not denote a real number, it~is just a notation symbol. If $\Phi :=(\varphi _{d,p'} )_{d,p'}$, it can also be written more compactly as $\varepsilon ^{p} \Phi $.
\end{definition}
\nomenclature[40b]{$\feps$}{Formal indeterminate representing the layer's thickness}

\begin{textesansboite}
The notation ``$\varepsilon $'' is chosen to avoid the confusion with the real number $\peps>0$, although the physical interpretation of $\varepsilon ^{p}$ is that it represents a multiplication by $\peps^{p}$.
\end{textesansboite}

\begin{remarque}
If $\Phi  \in \prod _{(d,p')\in \mathbb{R}^2 } E_{d}$ has locally finite support, $\varepsilon ^{p} \Phi $ also has locally finite support, because $\operatorname{supp}(\varepsilon ^{p} \Phi )=\tau _{p} (\operatorname{supp}(\Phi ))$ where $\tau _{p}$ is the translation $(d,p') \mapsto  (d,p'+p)$.
\end{remarque}

\begin{lemmesansboite}{}
\label{action de groupe de eps dans les SF}%
``$\varepsilon ^{p} \,\cdot $'' is a group action of $(\mathbb{R},+)$. That is to say, for~any $p,p'\in \mathbb{R}$ and any $\Phi  \in \prod _{(d,p')\in \mathbb{R}^2 } E_{d}$ that has locally finite support, we have:
\begin{itemize}
    \item $\varepsilon ^{0} \cdot  \Phi  = \Phi $,
    \item $\varepsilon ^{p} \cdot (\varepsilon ^{p'} \cdot \Phi ) = \varepsilon ^{p+p'} \cdot \Phi $.
\end{itemize}
\end{lemmesansboite}

\begin{lemme}{commutation of $\varepsilon $ and $\fsum$}
\label{fsum et eps commutent}%
Let $p\in \mathbb{R}$, $I$ be a set and $(\Phi _{i} )_{i\in I} \in \big(\prod _{(d,p)\in \mathbb{R}^2 } E_{d} )$ be a family of formal series (i.e.\ each $\Phi _{i}$ has locally finite support) that has collectively locally finite support. Then $(\varepsilon ^{p} \Phi _{i} )_{i\in I}$ has collectively locally finite support and
\[\fsum_{i\in I} \varepsilon ^{p} \Phi _{i} = \varepsilon ^{p} \fsum_{i\in I} \Phi _{i} .\]\fbeq 
\end{lemme}

\begin{demo}
For~any $i\in I$, we have $\operatorname{supp}(\varepsilon ^{p} \Phi _{i} )=\tau _{p} (\operatorname{supp}(\Phi _{i} ))$ with $\tau _{p} :(d,p') \mapsto  (d,p'+p)$, so it~is easy to see that $(\varepsilon ^{p} \Phi _{i} )_{i\in I}$ has collectively locally finite support. Moreover, denoting $(\varphi _{i,d,p'} )_{d,p'} :=\Phi _{i}$ for~any $i\in I$, we have
\[\fsum_{i\in I} \varepsilon ^{p} \Phi _{i} =\fsum_{i\in I} (\varphi _{i,d,p'-p} )_{d,p'} =\bigg(\sum _{i\in I} \varphi _{i,d,p'-p} \bigg)_{\!d,p'} =\varepsilon ^{p} \cdot \bigg(\sum _{i\in I} \varphi _{i,d,p'} \bigg)_{\!d,p'} = \varepsilon ^{p} \fsum_{i\in I} \Phi _{i} .\]
\end{demo}

\begin{textesansboite}
The next result allows us to write formal series with index in $\mathbb{R}^2 $ as series in powers of $\varepsilon $. In this result, ``$\varepsilon ^{p} \varphi _{d,p}$'' stands for $\varepsilon ^{p} \cdot \isfud\circ \iota (\varphi _{d,p} )$ by abuse of notation, where $\isfud\circ \iota $ is the canonical injection $E\rightarrow \prod _{d,p} E_{d}$.
\end{textesansboite}

\newcommand{\secAd}[1][d]{P_{#1}}
\newcommand{\secAp}[1][p]{D_{#1}}
\newcommand{\projAp}{P}
\newcommand{\projAd}{D}

\begin{lemmesecable}{}
\label{SF 2D de composantes pures}%
\begin{enumerate}
    \item Let $(\varphi _{d,p} )\in \prod _{(d,p)\in \mathbb{R}^2 } E_{d}$ that has a locally finite support. The $\fsum$-sums below are well defined (i.e.\ the summed families have collectively locally finite support) and the following equalities hold:
    \[\fsumdp \varepsilon ^{p} \varphi _{d,p} =\fsumd \fsump \varepsilon ^{p} \varphi _{d,p} =\fsump \fsumd \varepsilon ^{p} \varphi _{d,p} =(\varphi _{d,p} )_{(d,p)\in \mathbb{R}^2 } .\]
    \item Let $A\subset \mathbb{R}^2 $ be a locally finite set (i.e.\ $A\cap B$ is finite for~any bounded set $B\subset \mathbb{R}^2 $) and $(\varphi _{d,p} )_{(d,p)\in A} \in \prod _{(d,p)\in A} E_{d}$. Let us denote $\widetilde{\varphi }_{d,p} := \varphi _{d,p}$ for~any $(d,p)\in A$ and $\widetilde{\varphi }_{d,p} := 0$ for~any $(d,p)\in \mathbb{R}^2  \setminus A$. We also denote:
    \begin{itemize}
        \item $\secAd \,:=\{p\in \mathbb{R} \mid (d,p)\in A\}$ for~any $d\in \mathbb{R}$, and $\projAd:=\{d\in \mathbb{R} \mid \secAd \,\neq \varnothing \}$,
        \item $\secAp :=\{d\in \mathbb{R} \mid (d,p)\in A\}$ for~any $p\in \mathbb{R}$, and $\projAp\,:=\{p\in \mathbb{R} \mid \secAp \neq \varnothing \}$.
    \end{itemize}Then the $\fsum$-sums below are well defined and the following equalities hold:
    \[\fsum_{(d,p)\in A} \varepsilon ^{p} \varphi _{d,p} =\fsum_{d\in \projAd} \fsum_{p\in \secAd} \varepsilon ^{p} \varphi _{d,p} =\fsum_{p\in \projAp} \fsum_{d\in \secAp} \varepsilon ^{p} \varphi _{d,p} =( \widetilde{\varphi }_{d,p} )_{(d,p)\in \mathbb{R}^2 } .\]
\end{enumerate}\fbeq
\end{lemmesecable}

\begin{demo}
\begin{enumerate}
    \item For~any $(d,p)\in \mathbb{R}^2 $, we have $\varepsilon ^{p} \varphi _{d,p} =(\delta _{d',d} \,\delta _{p',p} \,\varphi _{d,p} )_{(d',p')\in \mathbb{R}^2 }$, so $\operatorname{supp}(\varepsilon ^{p} \varphi _{d,p} )$ is equal to $\{(d,p)\}$ if $\varphi _{d,p} \neq 0$ and to $\varnothing $ otherwise. Like in the proof of \cref{SF 1D de composantes pures}, it~is easy to deduce that $(\varepsilon ^{p} \varphi _{d,p} )_{(d,p)\in \mathbb{R}^2 }$ has collectively locally finite support and that $\fsumdp \varepsilon ^{p} \varphi _{d,p} =(\varphi _{d,p} )_{d,p}$. Moreover the equality with $\fsumd \fsump \varepsilon ^{p} \varphi _{d,p}$ and $\fsump \fsumd \varepsilon ^{p} \varphi _{d,p}$ (and the well-definition of these $\fsum$-sums) is a consequence of Fubini's property (\cref{Fubini pour fsum}).
    \item Point 1 yields:
    \[\fsumdp \varepsilon ^{p} \widetilde{\varphi }_{d,p} =\fsumd \fsump \varepsilon ^{p} \widetilde{\varphi }_{d,p} =\fsump \fsumd \varepsilon ^{p} \widetilde{\varphi }_{d,p} =(\varphi _{d,p} )_{(d,p)\in \mathbb{R}^2 } .\]
    Then \cref{proprietes de fsum} item 4 gives:
    \[\fsumdp \varepsilon ^{p} \widetilde{\varphi }_{d,p} = \fsum_{(d,p)\in A} \varepsilon ^{p} \widetilde{\varphi }_{d,p} + \fsum_{(d,p)\in \mathbb{R}^2  \setminus A} \varepsilon ^{p} \widetilde{\varphi }_{d,p} = \fsum_{(d,p)\in A} \varepsilon ^{p} \varphi _{d,p}\]
    and similarly \[\fsumd \fsump \varepsilon ^{p} \widetilde{\varphi }_{d,p} = \fsumd \fsum_{p\in \secAd} \varepsilon ^{p} \varphi _{d,p} =\fsum_{d\in \projAd} \fsum_{p\in \secAd} \varepsilon ^{p} \varphi _{d,p}\]
    and $\fsump \fsumd \varepsilon ^{p} \widetilde{\varphi }_{d,p} = \fsum_{p\in \projAp} \fsum_{d\in \secAp} \varepsilon ^{p} \varphi _{d,p}$.
\end{enumerate}
\end{demo}

\begin{textesansboite}
In the next lemma, $\varepsilon ^{p} E_{d}$ denotes by abuse of notation $\varepsilon ^{p} \cdot \isfud\circ \iota (E_{d} ) := \{\varepsilon ^{p} \cdot \isfud\circ \iota (\varphi ) \mid \varphi \in E_{d} \}$.
\end{textesansboite}

\begin{lemme}{}
\label{fsum des Ed dans les SF 2D}%
Let $A\subset \mathbb{R}^2 $. For~any $d,p\in \mathbb{R}$, we define $\projAd$, $\secAd$, $\projAp$ and $\secAp$ like in \cref{SF 2D de composantes pures}. Then we have
\begin{align*}
    \fsum_{(d,p)\in A} \varepsilon ^{p} E_{d} &= \fsum_{d\in \projAd} \fsum_{p\in \secAd} \varepsilon ^{p} E_{d} =\fsum_{p\in \projAp} \fsum_{d\in \secAp} \varepsilon ^{p} E_{d} \\
    &= \bigg\{\Phi  \in  \!\prod _{(d,p)\in \mathbb{R}^2 } E_{d} \ \bigg|\ \Phi  \text{ has locally finite support and} \; \operatorname{supp}(\Phi ) \subset A\bigg\}.
\end{align*}
Moreover, this set is closed under $\fsum$: for~any set $I$ and any $(\Phi _{i} )_{i\in I} \in \big(\fsum_{(d,p)\in A} \varepsilon ^{p} E_{d} )^{I}$ that has collectively locally finite support, $\fsum_{i\in I} \Phi _{i} \in \fsum_{(d,p)\in A} \varepsilon ^{p} E_{d}$.
\end{lemme}

\begin{demo}
The set equalities are a consequence of \cref{def: fsum d'ensembles,SF 2D de composantes pures} item 2. Moreover the closure under $\fsum$ comes from the property $\operatorname{supp}\!\big(\fsum_{i\in I} \Phi _{i} \big) \subset \bigcup _{i\in I} \operatorname{supp}(\Phi _{i} )$ stated by \cref{les SF sont stables par fsum}.
\end{demo}

\begin{remarquenumerotee}
\label{rq: ecriture avec fsum de l'ensemble des SF 2D}%
In~particular, the set of the formal series of $\prod _{(d,p)\in \mathbb{R}^2 } E_{d}$ is equal to
\[\fsumdp \varepsilon ^{p} E_{d} = \fsumd \fsump \varepsilon ^{p} E_{d} .\]
This is how we will denote it in practice. In summary, we have the following diagram for~any $d',p'\in \mathbb{R}$.
\[\xymatrixrowsep{6ex}\xymatrixcolsep{6ex}\xymatrix{
	E_{d'} \ar[r] &\displaystyle \bigoplus_{d\in \mathbb{R}} E_{d} \ar[r]^-{\iota } &\displaystyle \fsumd E_{d} \ar[r]^-{\isfud} \ar[d]_{\varepsilon ^{p'\!} \isfud({\dots })}^{\cong} &\displaystyle \fsumdp \varepsilon ^{p} E_{d} \ar[r]  &\displaystyle \prod _{(d,p)\in \mathbb{R}^2 } E_{d} \\
    &&\displaystyle \fsumd \varepsilon ^{p'\!} E_{d} \ar[ru] &&
}\]
All the arrows are inclusions, except the vertical one which is a bijection. 
This bijection commutes with $\fsum$ (\cref{fsum et eps commutent}), so $\fsumd \varepsilon ^{p'\!} E_{d}$ behaves like a copy of $\fsumd E_{d}$.
\end{remarquenumerotee}

\subsection{Operators with a degree}

\begin{textesansboite}
Let $F=\bigoplus_{d\in \mathbb{R}} F_{d}$ be another $\mathbb{C}$-vector space.
\end{textesansboite}

\begin{definition}{operators that have a degree}
\label{def: op à degree dans les SF 2D}%
\begin{enumerate}
    \item Let $f:E\rightarrow F$ be a map that has a degree. We define its extension to $\fsumdp \varepsilon ^{p} E_{d}$ by
    \[\tilde{f} : \left(\begin{array}{clc}
       \displaystyle \fsumdp \varepsilon ^{p} E_{d} &\rightarrow &\displaystyle \fsumdp \varepsilon ^{p} F_{d} \\[3ex]
        (\varphi _{d,p} )_{(d,p)\in \mathbb{R}^2 } &\mapsto & \big(f(\varphi _{d-\deg f,p} )\big)_{(d,p)\in \mathbb{R}^2 }
    \end{array}\right) .\]
    By~abuse of notation, we also denote it $f$.
    \item Let $g:\fsumdp \varepsilon ^{p} E_{d} \rightarrow  \fsumdp \varepsilon ^{p} F_{d}$. We say that $g$ has a degree if there exist $p\in \mathbb{R}$ and a map $f:E\rightarrow F$ that has a degree s.t.\ 
    \[\forall \Phi \in \!\fsumdp \varepsilon ^{p} E_{d} \,,\qquad  g(\Phi ) = \varepsilon ^{p} \cdot  \tilde{f}(\Phi ) .\]
    In this case, we denote $\deg  g:=(\deg  f,p)$.
\end{enumerate}\fbi
\end{definition}

\begin{textesansboite}
It~is easy to see that this new extension coincides on $\fsumd E_{d}$ with the extension defined in \cref{def: extension d'un op a degre aux SF 1D}. So the terminology and notations are consistent. Moreover, Point 1 of \cref{def: op à degree dans les SF 2D} can be rewritten as
\[\tilde{f}\bigg(\fsumdp \varepsilon ^{p} \varphi _{d,p} \bigg) = \fsumdp \varepsilon ^{p} f(\varphi _{d,p} )\]
for~any $(\varphi _{d,p} )\in \fsumdp \varepsilon ^{p} E_{d}$.
\end{textesansboite}

\begin{exemple}
If $f:E\rightarrow F$ that has a degree, then $\tilde{f}$ has a degree and $\deg  \tilde{f} =(\deg  f,0)$.
\end{exemple}

\begin{proposition}{}
\label{les op a degre commutent avec fsum et eps dans les SF 2D}%
Let $f:\fsumdp \varepsilon ^{p} E_{d} \rightarrow  \fsumdp \varepsilon ^{p} F_{d}$ be an operator that has a degree.
\begin{enumerate}
    \item For~any $p\in \mathbb{R}$ and $\Phi \in \fsum_{(d,p')\in \mathbb{R}^2 } \varepsilon ^{p'} E_{d}$, we have
    \[f(\varepsilon ^{p} \Phi ) = \varepsilon ^{p} f(\Phi ).\]
    \item For~any set $I$ and any $(\Phi _{i} )_{i} \in \big(\fsumdp \varepsilon ^{p} E_{d} \big)^{I}$ that has collectively locally finite support, $\big(f(\Phi _{i} )\big)_{i\in I}$ has collectively locally finite support and 
    \[f\bigg(\fsum_{i\in I} \Phi _{i} \bigg) = \fsum_{i\in I} f(\Phi _{i} ).\]\fbeq
\end{enumerate}
\end{proposition}

\begin{textesansboite}
The proof is similar to lots of other proofs of \cref{sec: SF 1D,sec: SF 2D}, especially \cref{fsum commute avec les op a degre dans les SF 1D,fsum et eps commutent}, and is left to the reader.
\end{textesansboite}

\subsection{Geometric series of operators}
\label{sec: series d'op de SF 2D}

\begin{textesansboite}
Let $v\in \mathbb{R}^2  \setminus \{0\}$ and
\[E^{v} := \bigg\{\Phi  \in \!\fsumdp \varepsilon ^{p} E_{d} \ \bigg|\ \exists M\in \mathbb{R},\; \forall (d,p)\in \operatorname{supp}(\Phi ),\; v\cdot (d,p) \geqslant M\bigg\} = \bigcup _{M\in \mathbb{R}} \fsum_{\substack{(d,p)\in \mathbb{R}^2 \\ v\cdot (d,p)\geqslant M}} \varepsilon ^{p} E_{d} .\]
where $v\cdot (d,p)$ denotes the scalar product of $\mathbb{R}^2 $. $E^{v}$ is a vector subspace of $\fsumdp \varepsilon ^{p} E_{d}$.
\end{textesansboite}

\begin{proposition}{}
\label{existence des serie geometrique d'op 2D}%
Let $\mathcal{F}$ be a finite set of operators $\fsumdp \varepsilon ^{p} E_{d} \rightarrow \fsumdp \varepsilon ^{p} E_{d}$ that have a degree and $f:=\sum _{g\in \mathcal{F}} g$. We assume: $\forall g\in \mathcal{F},\; v\cdot \deg  g >0$.\\
For~any $\Phi \in E^{v}$, we have $f(\Phi )\in E^{v}$, $\big(f^{n} (\Phi )\big)_{n\in \mathbb{N}}$ has collectively locally finite support and
\[\fsumn f^{n} (\Phi ) \in E^{v} .\]\fbeq[0.8]
\end{proposition}

\begin{textesansboite}
The proof is similar to the one of \cref{existence des serie geometrique d'op 1D}. 
To reduce the problem to something similar to formal series with a one-dimensional index, we study the projection on $v$ of the supports $\operatorname{supp}(f^{n} (\Phi ))$.
\end{textesansboite}

\begin{demo}
For~any $g\in \mathcal{F}$ and $\Psi  \in \prod _{(d,p)\in \mathbb{R}^2 } E_{d}$, let us denote 
\[\deg _{v} g:=v\cdot \deg g \qquad \, \text{ and } \qquad  \operatorname{supp}_{v} (\Psi ) = \{v\cdot (d,p) \mid (d,p)\in \operatorname{supp}(\Psi )\}.\]
Let $M_{f} :=\min \{\deg _{v} g \mid g\in \mathcal{F}\}$. Since $\Phi \in E^{v}$, there exists $M_{\Phi } \in \mathbb{R}$ s.t.\ $\operatorname{supp}_{v} (\Phi ) \subset [M_{\Phi } ,\infty )$. For~any $g\in \mathcal{F}$, we have $\operatorname{supp}(g(\Phi ))\subset \operatorname{supp}(\Phi )+\deg  g$ so $\operatorname{supp}_{v} (g(\Phi ))\subset  [M_{\Phi } +\deg _{v} g,\infty )\subset [M_{\Phi } +M_{f} ,\infty )$. This implies $\operatorname{supp}_{v} (f(\Phi )) \subset [M_{\Phi } +M_{f} ,\infty )$, thus $f(\Phi )\in E^{v}$. Moreover, iterating this result gives: 
\[\forall n\in \mathbb{N}, \qquad  \operatorname{supp}_{v} (f^{n} (\Phi ))\subset  [M_{\Phi } +n M_{f} ,\infty ) .\]
Hence for~any bounded set $B\subset \mathbb{R}^2 $, $\{n\in \mathbb{N} \mid \operatorname{supp}(f^{n} (\Phi ))\cap B \neq  \varnothing \}$ is included in $\{n\in \mathbb{N} \mid M_{\Phi } +n M_{f} \leqslant \sup _{b\in B} v\cdot b \}$ which is finite because $M_{f} >0$. Thus $(f^{n} (\Phi ))_{n}$ has collectively locally finite support. Finally, $E^{v} _{M} :=\fsum_{(d,p)\in \mathbb{R}^2 ,\,v\cdot (d,p)\geqslant M} \varepsilon ^{p} E_{d}$ is closed under $\fsum$ by \cref{fsum des Ed dans les SF 2D}, so $\fsum_{n} f^{n} (\Phi ) \in E^{v} _{M} \subset E^{v}$.
\end{demo}

\begin{remarque}
In \cref{existence des serie geometrique d'op 2D}, the assumption $\Phi \in  E^{v}$ is~not optimal for $\fsumn f^{n} (\Phi )$ to be well defined. But this criterion is easy to check and will be sufficient for us.
\end{remarque}

\begin{definition}{geometric series of operators}
\label{def: serie geometrique d'op 2D}%
With the notations of \cref{existence des serie geometrique d'op 2D}, we define $\fsumn f^{n} :E^{v} \rightarrow  E^{v}$ by:
\[\fsumn f^{n} : \ \Phi  \longmapsto  \fsumn f^{n} (\Phi )\]
Moreover, if $k\in \mathbb{N}$, $f_{1} ,\dots ,f_{k}$ are distinct linear maps $E\rightarrow E$ s.t.\ $\mathcal{F}:=\{f_{1} ,\dots ,f_{k} \}$, we also denote $\langle \mathcal{F}\rangle  :=\langle f_{1} ,\dots ,f_{k} \rangle  :=\fsumn f^{n}$.
\end{definition}

\begin{theoreme}{operator inversion by Neumann series}
\label{inverse de id-f dans les SF 2D}%
With the notations of \cref{existence des serie geometrique d'op 2D}, the operator $ \mathrm{id}-f$ maps $E^{v}$ to $E^{v}$, its restriction to $E^{v}$ is invertible and the inverse is
\[( \mathrm{id}-f)^{-1} = \fsumn f^{n} .\]\fbeq
\end{theoreme}

\begin{demo}
The proof is almost identical to the one of \cref{inverse de id-f dans les SF 1D}. The steps that were justified by \cref{fsum commute avec les op a degre dans les SF 1D,existence des serie geometrique d'op 1D} are now justified by \cref{existence des serie geometrique d'op 2D,les op a degre commutent avec fsum et eps dans les SF 2D} resp., which are the analogous results for formal series with index in $\mathbb{R}^2 $.
\end{demo}

\begin{proposition}{commutation of $\langle .\rangle $ and $\fsum$}
\label{fsum et le crochet commutent en 2D}%
We keep the notations of \cref{existence des serie geometrique d'op 2D}. Let $I$ be a set and $(\Phi _{i} )_{i\in I}$ be a family of elements of $E^{v}$ that has collectively locally finite support. We assume there exists $M\in \mathbb{R}$ s.t.
\[\forall i\in I, \qquad  \Phi _{i} \in  \fsum_{(d,p)\in \mathbb{R}^2 ,\, v\cdot (d,p)\geqslant M} \varepsilon ^{p} E_{d} .\]
Then $\fsum_{i\in I} \Phi _{i}$ belongs to $E^{v}$, $\big(\langle \mathcal{F}\rangle (\Phi _{i} )\big)_{i}$ has collectively locally finite support and
\[\langle \mathcal{F}\rangle \bigg(\fsum_{i\in I} \Phi _{i} \bigg) = \fsum_{i\in I} \,\langle \mathcal{F}\rangle (\Phi _{i} ).\]\fbeq[0.7]
\end{proposition}

\begin{remarquenumerotee}
\label{rq: eps et le crochet commutent}%
\Cref{fsum et eps commutent,les op a degre commutent avec fsum et eps dans les SF 2D} show that $\langle \mathcal{F}\rangle $ and $\varepsilon ^{p}$ commute for~any $p\in \mathbb{R}$.
\end{remarquenumerotee}

\begin{demo}
\def\sfjn#1{\mathcal{I}_{#1}}
\def\sfan#1{\mathcal{I}'_{#1}}
\underline{Step 1:} Let $f:=\sum _{g\in \mathcal{F}} g$. Let us show that, for~any $n\in \mathbb{N}$, $\sfjn{n} :=\big(f^{k} (\Phi _{i} )\big)_{(i,k)\in I\times [\![0,n]\!]}$ has collectively locally finite support (CLFS). We start by proving that $\sfan{n}:=\big(f^{n} (\Phi _{i} )\big)_{i\in I}$ has CLFS for~any $n\in \mathbb{N}$, by induction on $n$.
\begin{itemize} 
    \item For $n=0$, $\sfan{0} =(\Phi _{i} )_{i\in I}$ has CLFS by assumption.
    \item Let $n>0$. Let us assume that $\sfan{n-1}$ has CLFS and prove that $\sfan{n}$ has CLFS. \Cref{les op a degre commutent avec fsum et eps dans les SF 2D} shows that $\big(g\circ f^{n-1} (\Phi _{i} )\big)_{i\in I}$ has CLFS for~any $g\in \mathcal{F}$. By~sum (i.e.using \cref{proprietes de fsum} item 3), we deduce that $\sfan{n}$ has CLFS.
\end{itemize}
Now, $\big(f^{k} (\Phi _{i} )\big)_{(i,k)\in I\times \{n\}}$ is in bijection with $\sfan{n}$, so it has CLFS by \cref{proprietes de fsum} item 5. Finally, by finite disjoint union of families that have CLFS (see \cref{proprietes de fsum} item 4), we get that $\sfjn{n}$ has CLFS.\medskip

\underline{Step 2:} Let us show that $\mathcal{I}:=\big(f^{k} (\Phi _{i} )\big)_{(i,k)\in I\times \mathbb{N}}$ has CLFS. Like in the proof of \cref{existence des serie geometrique d'op 2D}, we define for~any $g\in \mathcal{F}$ and $\Psi  \in \prod _{(d,p)\in \mathbb{R}^2 } E_{d}$:
\[\deg _{v} g:=v\cdot \deg g \qquad \, \text{ and } \qquad  \operatorname{supp}_{v} (\Psi ) = \{v\cdot (d,p) \mid (d,p)\in \operatorname{supp}(\Psi )\}.\]
Let $M_{\mathcal{F}} :=\min \{\deg _{v} g \mid g\in \mathcal{F}\}$ and let $M_{I} \in \mathbb{R}$ be s.t.\ $\bigcup _{i\in I} \operatorname{supp}_{v} (\Phi _{i} ) \subset [M_{I} ,\infty )$. Using \cref{les SF sont stables par fsum}, we get $\operatorname{supp}_{v} \!\big(\fsum_{i\in I} \Phi _{i} \big) \subset  \bigcup _{i\in I} \operatorname{supp}_{v} (\Phi _{i} )\subset [M_{I} ,\infty )$ which proves $\fsum_{i\in I} \Phi _{i} \in E^{v}$.\\
Let $B\subset \mathbb{R}^2 $ be a bounded set and $A :=\{(i,n)\in I\times \mathbb{N} \mid \operatorname{supp}(f^{n} (\Phi _{i} ))\cap B \neq \varnothing  \}$. We need to show that $A$ is finite. Like in the proof of \cref{existence des serie geometrique d'op 2D}, we have $\operatorname{supp}_{v} (f^{n} (\Phi _{i} )) \subset [M_{I} +n M_{\mathcal{F}} ,\infty )$ for~any $i\in I$ and $n\in \mathbb{N}$. Therefore
\begin{align*}
    A &= \{(i,n)\in I\times \mathbb{N} \,\mid\, \operatorname{supp}(f^{n} (\Phi _{i} ))\cap B \neq \varnothing  \text{ and } M_{I} +n M_{\mathcal{F}} \leqslant  \sup _{b\in B} v\cdot b\}\\
    & = \{(i,n)\in I\times [\![0,n_{B} ]\!] \,\mid\,  \operatorname{supp}(f^{n} (\Phi _{i} ))\cap B \neq \varnothing \} \qquad  \text{with }  n_{B} := \left\lfloor { \frac{1}{M_{\mathcal{F}}} \Big(\sup _{b\in B} v\cdot b -M_{I} \Big)}\right\rfloor 
\end{align*}
because $M_{\mathcal{F}} >0$. Since $\sfjn{n_{B} }$ has CLFS by step 1, $A$ is finite. Therefore, $\mathcal{I}$ has CLFS.\medskip

\underline{Step 3:} One can conclude using \cref{Fubini pour fsum} (Fubini) and \cref{les op a degre commutent avec fsum et eps dans les SF 2D} item 2. This part is entirely similar to the end of the proof of \cref{fsum et le crochet commutent en 1D} and we do~not repeat it.
\end{demo}

\begin{textesansboite}
\Cref{fig: relations de commutation dans les SF 2D} gives a summary of the commutation properties in $\fsumdp \varepsilon ^{p} E_{d}$. Basically everything commutes between $\fsum$, $\varepsilon ^{p}$ and operators (with a degree or of the form $\langle \mathcal{F}\rangle $), except that operators do~not commute with other operators in general. 
\end{textesansboite}

\begin{figure}[h]
\newcommand{\drawCycle}[6][]{
    \draw[#1] (#2.center) ++ ({#4/2*tan(#5/2)*cos(#3+90-#5/2)},{#4/2*tan(#5/2)*cos(#3-90+#5/2)}) arc ({#3+(180-#5/2)} : {#3-(180-#5/2)} : #4/2) #6;}
\begin{tikzpicture}[scale=0.4]
    \node (S) at (0,10) {$\fsum$};
    \node (E) at (10,10) {$\feps^p$};
    \node (F) at (0,0) {$f$ with a degree};
    \node (C) at (10,0) {$\langle\mathcal{F}\rangle$};
	\draw[<->] (F) -- node[midway,above,rotate=90]{\small Prop.~\ref{les op a degre commutent avec fsum et eps dans les SF 2D}} (S);
	\draw[<->] (S) -- node[midway,above]{\small \cref{fsum et eps commutent}} (E);
	\draw[<->] (E) -- node[midway,above,rotate=-90]{\cref{rq: eps et le crochet commutent}} (C);
	\draw[<->] (F) -- node[near end,below,rotate=45]{\small Prop.~\ref{les op a degre commutent avec fsum et eps dans les SF 2D}} (E);
	\draw[<->] (S) -- node[near start,above,rotate=-45]{\small Prop.~\ref{fsum et le crochet commutent en 2D}} (C);
    \drawCycle[<->]{S}{135}{2.1}{80}{node[midway,left]{Prop.~\ref{Fubini pour fsum} (Fubini)}}
    \drawCycle[<->]{E}{45}{2.1}{80}{node[midway,right]{\cref{action de groupe de eps dans les SF}}}
    \draw (-12,0) node{\ } (22,0) node{\ }; 
\end{tikzpicture}\par
\longCaption{Commutation relations among formal series with index in $\mathbb{R}^2 $}{ (each arrow represents a commutation property)}
\label{fig: relations de commutation dans les SF 2D}
\end{figure}

\chapter{Derivation of the matching conditions}
\label{chap: raccord algebrique}
\begin{textesansboite}
This chapter establishes the matching conditions linking corner fields to far and layer fields. This is by far the most difficult relation to derive, while all the others have been easily stated in \cref{eq: eq vol champs lointains}--\cref{eq: eq vol champs de coin}. In this chapter we assume that the various fields exist and that they satisfy \cref{eq: eq vol champs lointains}--\cref{eq: eq vol champs de coin} and we give a necessary and sufficient condition for the matching assumption \cref{eq: hypothese de raccord} to be satisfied. Our approach is based on an algebraic formulation of the problem, that reveals the structure of the matching assumption by a rigorous algebraic calculus. This relies on three original ideas:
\begin{itemize}
	\item decomposing the corner singularities on a basis of functions that have explicit expressions,
	\item building operators that allow to express all the singularities that appear in the fields behavior around the corner in terms of some fundamental singularities equal to $r^{d} \sin (d\theta )$ with $d\in \TZ*$,
	\item and using the formal series introduced in \cref{chap: SF} to efficiently manipulate the fields behaviors and the matching assumption.
\end{itemize}
See \cref{sec: resume raccord en intro} for a more detailed and motivated introduction to these tools.\\

The outline of this chapter is the following.
\begin{itemize}
	\item In \cref{sec: singularites quasi homogenes} we define spaces, denoted $\mathcal{A}(\Omega )$, $\mathcal{A}(\couche)$, etc., that contain the corner singularities appearing in the fields behaviors around the corner. These functions are quasi-homogeneous in the sense that they behave like $\bfr^{d} \ln ^{j} \!\bfr$ for some $d\in \mathbb{R}$ and $j\in \mathbb{N}$. We provide basis of these spaces made of explicit functions. In addition, we build operators that allow to solve a Poisson-like problem in $\Pi $ in the spaces $\mathcal{A}$.
	\item In \cref{sec: SF de singularites}, we use the tools of \cref{chap: SF} to define spaces of formal series of the previous singularities. We study the resolution of formal Helmholtz-like and Poisson-like problems in these spaces, using the operators built in \cref{sec: singularites quasi homogenes}. This allows us to express the fields behaviors around the corner in terms of their coefficients on the basis functions $r^{d} \sin (d\theta )$.
	\item In \cref{sec: raccordement algebrique}, we rigorously re-express the matching assumption \cref{eq: hypothese de raccord} in those spaces of formal series, and we prove that it~is equivalent to some matching conditions that link the coefficients of the fields on the basis functions $r^{d} \sin (d\theta )$.
	\item In \cref{sec: resume du raccord}, we summarize the main ingredients of our method, in order to extend it to other cases in future works.
	\item Finally, in \cref{sec: implementation des espaces A}, we explain how the tools of this chapter allow a numerical computation of the fields behaviors and the matching coefficients that is fast, memory-thrifty and exact.
\end{itemize}
\end{textesansboite}

\section{Quasi-homogeneous corner singularities}
\label{sec: singularites quasi homogenes}

\subsection{Definition of the spaces $\mathcal{A}$}
\label{sec: def de A}

\begin{textesansboite}
In \cite[p.10]{CosDau96}, Costabel and Dauge build a similar asymptotic expansion for the Poisson equation in the half plane with mixed boundary condition: Neumann in a part of the boundary and Robin $u+\peps \partial _{n} u=0$ in another. They quickly mention that their singularities can be written as $\Re[(-z)^{q} P(\ln (-z))]$ with $z=x+\mathrm{i}y$, $q\in \mathbb{R}$ and $P$ a real polynomial. To define the spaces $\mathcal{A}$, we adapted this idea to take into account the layer, the angle $\Theta $ and the Helmholtz equation. These simple expressions give both powerful algebraic tools for the theory and fast precise algorithms for the numerical resolution (see \cref{sec: operateurs de resolution dans A,sec: implementation des espaces A}).
\end{textesansboite}

\begin{definition}{the spaces $\mathcal{A}$}
\label{def: espaces A}%
Let $\alpha  := e^{-\mathrm{i}\Theta }$. We define in $\Omega $ the complex variable $z:=x+\mathrm{i}y=re^{\mathrm{i}\theta }$. For~any $q\in \mathbb{R}$, we take the following conventions: $(\alpha z)^{q} := r^{q} e^{q\mathrm{i}(\theta -\Theta )}$, $\overline{\alpha z}^{q} := r^{q} e^{-q\mathrm{i}(\theta -\Theta )}$ and $\log (\alpha z) := \ln r+\mathrm{i}(\theta -\Theta )$. Let $d\in \mathbb{R}$. We denote:
\begin{itemize}
	\item $\mathcal{A}_{d} (\Omega )$ the vector subspace of $\mathcal{C}^{0} (\Omega ,\mathbb{C})$ generated by the functions $z\mapsto \Im[(\alpha z)^{q} \,\overline{\alpha z}^{k} P(\log (\alpha z))]$ with $q\in \mathbb{R}$, $k\in \mathbb{N}$, $q+k=d$ and $P\in \mathbb{R}[T]$,
	\item $\mathcal{A}_{d} (\couche) := \{(x,Y)\mapsto  x^{d} \,Q(\ln x,Y) \mid Q\in \mathbb{C}[T,Y] \text{ and } Q(T,-1)=0\}$,
	\item $\mathcal{A}_{d} (\Gamma ) := \{x\mapsto  x^{d} \,Q(\ln x) \mid Q\in \mathbb{C}[T]\}$,
	\item $\mathcal{A}_{d} (\Pi ) := \{\varphi \in \mathcal{C}^{0} (\Pi ,\mathbb{C}) \mid \varphi _{|\Omega } \in \mathcal{A}_{d} (\Omega ) \text{ and } \varphi _{|\couche} \in \mathcal{A}_{d} (\couche)\}$,
	\item and for~any $D\in \{\Pi ,\Omega ,\couche,\Gamma \}$, $\mathcal{A}(D) :=\displaystyle \sum _{d\in \mathbb{R}} \mathcal{A}_{d} (D)$.
\end{itemize}\fbi
\end{definition}
\nomenclature[10]{$\mathcal{A}(\dots )$}{Space of corner singularities}
\nomenclature[10-13]{$\mathcal{A}_{d} (\dots )$}{Space of quasi-homogeneous corner singularities of degree $d$}

\begin{textesansboite}
We recall that $\sum _{d\in \mathbb{R}} \mathcal{A}_{d} (D)$ denotes the set of the \emph{finite} (classical) sums of elements of the spaces $\mathcal{A}_{d} (D)$, as opposed to $\fsumd \mathcal{A}_{d} (D)$ (with a bold ``$\fsum$'') which contains \emph{infinite} formal series of elements of the $\mathcal{A}_{d} (D)$ (cf.\ the introduction of \cref{sec: SF 1D}). Formal series of corner singularities appear later in \cref{sec: SF de singularites,sec: raccordement algebrique}, but for now we focus on individual singularities and their finite sums.\\

Note also that the elements of $\mathcal{A}(\Omega )$ and $\mathcal{A}(\Pi )$ vanish on $\bordgauche$, and the elements of $\mathcal{A}(\couche)$ and $\mathcal{A}(\Pi )$ vanish on $\borddroit$. So the Dirichlet condition on the boundary is included in the space definition.
\end{textesansboite}

\begin{remarquenumerotee}
\label{rq: espaces A pour 2 couches}%
Our definition of $\mathcal{A}(\Omega )$ differs from the formula ``$\Re[(-z)^{q} P(\ln (-z))]$'' given by Costabel and Dauge in three ways:
\begin{itemize}
	\item $-z$ is replaced by $\alpha z$ to impose a boundary condition on the halfline $\theta =\Theta $ instead of $\theta =\pi $,
	\item the real part is replaced by an imaginary part to impose a Dirichlet condition instead of Neumann,
	\item moreover, $\Re[(-z)^{q} P(\ln (-z))]$, as the real part of an analytic function, is harmonic, which is appropriate to study the Laplace equation but not Helmholtz, so we added a non-analytic factor (namely $\overline{\alpha z}^{k}$).
\end{itemize}
\Cref{def: espaces A} can be extended to the case where $\bordgauche$ is covered by another layer, by defining $\mathcal{A}_{d} (\Omega )$ as the vector space generated by the functions $z\mapsto \Im[(\alpha z)^{q} \,\overline{\alpha z}^{k} P(\log (\alpha z))]$ and $z\mapsto \Im[z^{q} \bar{z}^{k} P(\log  z)]$ with $q\in \mathbb{R}$, $k\in \mathbb{N}$, $q+k=d$ and $P\in \mathbb{R}[T]$.
\end{remarquenumerotee}

\begin{remarque}
The elements of $\mathcal{A}(\Omega )$ are naturally functions depending on the polar coordinates. For~instance, we have:
\begin{itemize}
	\item $\Im[(\alpha z)^{q} \,\overline{\alpha z}^{k} ] = r^{q+k} \sin \!\big((q-k)(\theta -\Theta )\big)$
	\item $\Im[(\alpha z)^{q} \,\overline{\alpha z}^{k} \log (\alpha z)] = r^{q+k} \left[ \ln (r)\,\cos \!\big((q-k)(\theta -\Theta )\big)-(\theta -\Theta )\sin \!\big((q-k)(\theta -\Theta )\big) \right]$
\end{itemize}
Note also~that in \cref{def: espaces A} we used the variables $x,y$, which are relevant for far fields, but all the tools developed in this section can also be used for corner fields, replacing $(x,y)$ by $(X,Y)$.
\end{remarque}

\begin{textesansboite}
In order to define operators in the spaces $\mathcal{A}$, and especially some operators that give particular solutions of PDEs (see \cref{sec: operateurs de resolution dans A}), we will need the three following lemmas. In particular, we provide the explicit basis of the spaces $\mathcal{A}$ that we announced in this chapter's introduction.
\end{textesansboite}

\begin{lemme}{}
\label{decomposition de A en sev}%
For~any $D\in \{\Pi ,\Omega ,\Gamma ,\couche\}$, we have:
\[\mathcal{A}(D) =\bigoplus_{d\in \mathbb{R}} \mathcal{A}_{d} (D).\]\fbeq
\end{lemme}

\begin{demo}
\def\polt#1{P_{\theta\mkern-1mu,\mkern1mu #1}}%
By~\cref{def: espaces A}, we have $\mathcal{A}(D)=\sum _{d\in \mathbb{R}} \mathcal{A}_{d} (D)$, so it suffices to show that the sum is direct. We will only treat the case $D=\Omega $ because the other three are similar.\\
Let $n\in \mathbb{N}$ and, for~any $j\in [\![1,n]\!]$, $d_{j} \in \mathbb{R}$ and $\varphi _{j} \in \mathcal{A}_{d_{j}} (\Omega )$. We assume that $d_{1} <d_{2} <\cdots <d_{n}$ and $\sum_{j=1}^n \varphi _{j} =0$. Let us show that $\forall j\in [\![1,n]\!],\; \varphi _{j} =0$. Let us fix $\theta \in (0,\Theta )$. By~definition of $\mathcal{A}_{d_{j}} (\Omega )$, we have for~any $j\in [\![1,n]\!]$:
\[\exists \polt{j} \in \mathbb{C}[T] ,\; \forall r>0,\qquad  \varphi _{j} (r,\theta )=r^{d_{j}} \polt{j} (\ln  r).\]
Let us assume by contradiction that there exists $j\in [\![1,n]\!]$ s.t.\ $\polt{j} \neq 0$. Let $j_{0} :=\max \{j\in [\![1,n]\!] \mid \polt{j} \neq 0 \}$. Then we have
\[0 = \sum_{i=1}^n \varphi _{i} (r,\theta ) \sim  r^{d_{j_{0}}} \polt{j_{0} } (\ln r) \quad \ \mathrm{when} \ r\rightarrow \infty ,\]
so $\polt{j_{0} } =0$. This is absurd, so we have: $\forall j\in [\![1,n]\!],\forall r>0,\; \varphi _{j} (r,\theta )=0$. Since this holds for~any $\theta $, we deduce $\varphi _{j} =0$ for~any $j$.
\end{demo}

\begin{remarquenumerotee}
\label{rq: base de A couche et Gamma}%
For~any $d\in \mathbb{R}$, the maps $Q\in \mathbb{C}[T] \mapsto  x^{d} \.Q(\ln x)\in \mathcal{A}(\Gamma )$ and $Q\in \mathbb{C}[T,Y] \mapsto  x^{d} \.Q(\ln x,Y)\in \mathcal{A}(\couche)$ are clearly injective. In other words:
\begin{itemize}
	\item $\{x\mapsto  x^{d} \ln ^{j} \!x \mid j\in \mathbb{N}\}$ is a basis of $\mathcal{A}_{d} (\Gamma )$,
	\item and $\{(x,Y) \mapsto  x^{d} \ln ^{j} \!x \,(Y-1)^{k} \mid (j,k)\in \mathbb{N}^2 \}$ is a basis of $\mathcal{A}_{d} (\couche)$.
\end{itemize}
However, given an explicit basis of $\mathcal{A}_{d} (\Omega )$ is less obvious and is done in the next result. 
\end{remarquenumerotee}

\begin{lemme}{basis functions of $\mathcal{A}_{d} (\Omega )$}
\label{base de A Omega}%
For~any $d\in \mathbb{R}$, the set
\[\Big\{z\mapsto  \Im[(\alpha z)^{q} \,\overline{\alpha z}^{k} \log (\alpha z)^{j} ] \ \Big|\ (q,k,j)\in \mathbb{R}\times \mathbb{N}^2  \text{ s.t.\ }  q+k=d \text{ and } \big[(q\in \mathbb{N} \text{ and } q\leqslant k) \Rightarrow  j>0\big] \Big\}\]
is a basis of $\mathcal{A}_{d} (\Omega )$.
\end{lemme}

\begin{remarque}
The condition ``$(q\in \mathbb{N} \text{ and } q\leqslant k) \Rightarrow  j>0$'' excludes the functions $z\mapsto  \Im[(\alpha z)^{q} \.\overline{\alpha z}^{k} ]$ with $q\in \mathbb{N}$ and $q\leqslant k$, because these functions vanish when $q=k$ and, when $q\neq k$, they are equal to $z\mapsto  -\Im[(\alpha z)^{k} \.\overline{\alpha z}^{q} ]$ which is already in the basis.
\end{remarque}

\begin{demo}
Let $d\in \mathbb{R}$ and $\mathcal{B} :=\{z\mapsto  \Im[(\alpha z)^{q} \,\overline{\alpha z}^{k} \log (\alpha z)^{j} ] \ \big|\ (q,k,j)\in \mathbb{R}\times \mathbb{N}^2  \text{ s.t.\ }  q+k=d \text{ and } [(q\in \mathbb{N} \text{ and } q\leqslant k) \Rightarrow  j>0]\}$. Separating real and imaginary parts, it~is enough to show that $\mathcal{B}$ is a basis of $\mathcal{A}_{d} (\Omega )\cap \mathcal{C}^{0} (\Omega ,\mathbb{R})$. The proof uses the fact that $\Im[(\alpha z)^{q} \,\overline{\alpha z}^{k} \log (\alpha z)^{j} ] = r^{q+k} \,\Im\big[e^{\mathrm{i}(q-k)(\theta -\Theta )} (\ln r+\mathrm{i}(\theta -\Theta ))^{j} \big]$ for~any $(q,k,j)$ and has two steps:
\begin{itemize}
	\item First, we show that the subspaces of $\mathcal{A}_{d} (\Omega )$ in which $|q-k|$ is equal to a given number are in direct sum. 
	Since $q+k=d$, for a given value of $|q-k|$, $(q,k)$ can take at most two values.
	\item Then we separate the different powers of $\ln r$, by looking at the behavior when $r\rightarrow \infty $, similarly to \cref{decomposition de A en sev}.
\end{itemize}\medskip\smallskip

\underline{Step 1:} For~any $a\in \mathbb{R}_+$, let us denote
\[E_{a} := \operatorname{Span}\!\big(\big\{z\mapsto  \Im[(\alpha z)^{q} \,\overline{\alpha z}^{k} P(\log (\alpha z))] \ \big|\ (q,k,P)\in \mathbb{R}\times \mathbb{N}\times \mathbb{R}[T],\; q+k=d \text{ and } |q-k|=a\big\}\big).\]
Let us first show
\[\mathcal{A}_{d} = \bigoplus_{a\geqslant 0} E_{a} .\]
Let $n\in \mathbb{N}$, $0\leqslant a_{1} <\cdots <a_{n}$ and, for~any $j\in [\![1,n]\!]$, $\varphi _{j} \in E_{a_{j}}$. Let us assume $\sum_{j=1}^n \varphi _{j} =0$ and prove: $\forall j\in [\![1,n]\!],\; \varphi _{j} =0$. For~any $j\in \mathbb{N}$, $\varphi _{j}$ is a sum of terms of the form $\psi :z\mapsto  \Im[(\alpha z)^{q} \,\overline{\alpha z}^{k} P(\log (\alpha z))]$ with $(q,k,P)\in \mathbb{R}\times \mathbb{N}\times \mathbb{R}[T]$, $q+k=d$ and $|q-k|=a_{j}$. Such a function can be written in polar coordinates as
\[\psi (r,\theta ) =  r^{d} \,\Im\big[e^{\pm \mathrm{i}a_{j} (\theta -\Theta )} P(\ln r+\mathrm{i}(\theta -\Theta ))\big].\]
It is well known from the theory of ODE that, for~any $a\in \mathbb{R}$ and $Q\in \mathbb{C}[T]$, $f:\theta \mapsto e^{\pm \mathrm{i}a\theta } Q(\theta )$ is a solution of $(\partial _{\theta } ^2 +a^2 )^{m_{Q}} f=0$ for some $m_{Q} \in \mathbb{N}$. Hence there exists an integer $m$ s.t.\ we have:
\[\forall j\in [\![1,n]\!],\qquad  (\partial _{\theta } ^2 +a_{j} ^2 )^{m} \varphi _{j} =0.\]
Thus, for~any $r_{0} >0$, the functions $\varphi _{j|\{r=r_{0} \}}$ with $j\in [\![1,n]\!]$ belong to distinct generalized eigenspaces of the operator $\partial _{\theta } ^2 $ in $\mathcal{C}^{\infty } (0,\Theta )$. Since these spaces are in direct sum, we deduce $\varphi _{j|\{r=r_{0} \}} =0$ for~any $j$, and finally $\varphi _{j} =0$ on $\Omega $.\medskip\smallskip

\underline{Step 2:} Let $a\in \mathbb{R}_+$ and let us show that $\mathcal{B}\cap E_{a}$ is a basis of $E_{a}$. By~construction $\mathcal{B}\cap E_{a}$ generates $E_{a}$, so it suffices to prove that $\mathcal{B}\cap E_{a}$ is linearly independent. There exists $(q,k)\in \mathbb{R}\times \mathbb{N}$ s.t.\ $q+k=d$ and $|q-k|=a$. There are two cases.
\begin{itemize}
	\item \underline{Case 1: $q\not\in \mathbb{N}$.} In this case, we have $\mathcal{B}\cap E_{a} =\{z\mapsto  \Im[(\alpha z)^{q} \,\overline{\alpha z}^{k} \log (\alpha z)^{j} ] \mid j\in \mathbb{N}\}$. Let $n\in \mathbb{N}$ and $(\lambda _{j} ) \in \mathbb{R}^{n+1}$ be s.t.\ $\varphi := \sum_{j=0}^n \lambda _{j} \,\Im[(\alpha z)^{q} \,\overline{\alpha z}^{k} \log (\alpha z)^{j} ] =0$ and let us prove $(\lambda _{j} )=0$. We assume by contradiction that $(\lambda _{j} ) \neq 0$. Without loss of generality, we can assume $\lambda _{n} \neq 0$. We have for~any fixed $\theta \in (0,\Theta )$:
	\begin{multline*}
		0=r^{-d} \varphi (r,\theta ) = \smash{\sum_{j=0}^n \lambda _{j} \Im\big[e^{\mathrm{i}(q-k)(\theta -\Theta )} (\ln r+\mathrm{i}(\theta -\Theta ))^{j} \big]} \\
		\underset{r\rightarrow \infty }{\sim} \lambda _{n} \Im\big[e^{\mathrm{i}(q-k)(\theta -\Theta )} \ln ^{j} r] = \lambda _{n} \sin ((q-k)(\theta -\Theta )) \ln ^{n} r.
	\end{multline*}
	Taking $\theta $ s.t.\ $\sin ((q-k)(\theta -\Theta )) \neq 0$ we get $\lambda _{n} =0$. This is absurd and concludes this part of the proof.\smallbreak

	\item \underline{Case 2: $q\in \mathbb{N}$ and $q=k$.} The reasoning is similar to the previous case, but $\Im[(\alpha z)^{q} \. \overline{\alpha z}^{k} ]$ does~not belong to $\mathcal{B}$, so we ignore $\lambda _{0}$, we have $n\geqslant 1$ and we have:
	\[0=r^{-d} \varphi (r,\theta ) = \sum_{j=1}^n \lambda _{j} \Im\big[(\ln r+\mathrm{i}(\theta -\Theta ))^{j} \big] \underset{r\rightarrow \infty }{\sim} \lambda _{n} \Im\!\left[ \binom{n}{1} \ln ^{n-1} \!r \cdot  \mathrm{i}(\theta -\Theta )\right] = \lambda _{n} n \ln ^{n-1} \!r\cdot  (\theta -\Theta ).\]
	Again this implies $\lambda _{n} =0$.\smallbreak

	\item \underline{Case 3: $q\in \mathbb{N}$ and $q\neq k$.} Up to switching $q$ and $k$, we can assume $q>k$. Then we have 
	\[\mathcal{B}\cap E_{a} =\{z\mapsto  \Im[(\alpha z)^{q} \,\overline{\alpha z}^{k} \log (\alpha z)^{j} ] \mid j\in \mathbb{N}\} \cup  \{z\mapsto  \Im[(\alpha z)^{\red{k}} \,\overline{\alpha z}^{\red{q}} \log (\alpha z)^{j} ] \mid j\in \red{\mathbb{N}^*}\}.\]
	Let $(\lambda _{j} )_{j\in [\![0,n]\!]}$ and $(\nu _{j} )_{j\in [\![1,n]\!]}$ be s.t.
	\[\varphi := \sum_{j=0}^n \lambda _{j} \.\Im[(\alpha z)^{q} \,\overline{\alpha z}^{k} \log (\alpha z)^{j} ] + \sum_{j=1}^n \nu _{j} \.\Im[(\alpha z)^{\red{k}} \,\overline{\alpha z}^{\red{q}} \log (\alpha z)^{j} ] =0\]
	and let us prove that $(\lambda _{j} )$ and $(\nu _{j} )$ vanish. If $n=0$, then it~is clear that we have $\lambda _{0} =0$, so $(\lambda _{j} )$ and $(\nu _{j} )$ vanish. Now let us assume $n>1$. We assume by contradiction $\lambda _{n} \neq 0$ or $\nu _{n} \neq 0$. We have
	\[0=r^{-d} \varphi (r,\theta ) =\sum_{j=0}^n \lambda _{j} \Im\big[e^{\mathrm{i}(q-k)(\theta -\Theta )} (\ln r+\mathrm{i}(\theta -\Theta ))^{j} \big] \red{ {}-{} }  \sum_{j=1}^n \nu _{j} \Im\big[e^{\mathrm{i}(\red{q-k})(\theta -\Theta )} (\ln r \red{ {}-{} } \mathrm{i}(\theta -\Theta ))^{j} \big]\]
	Applying $(r\partial _{r} )^{n-1}$ (which is like differentiating $n-1$ times w.r.t.\ $\ln r$), we obtain:
	\begin{align*}
		0 &= \lambda _{n} \Im\big[e^{\mathrm{i}(q-k)(\theta -\Theta )} \,n!\. (\ln r+\mathrm{i}(\theta -\Theta )) \big] + \lambda _{n-1} \Im\big[e^{\mathrm{i}(q-k)(\theta -\Theta )} \.(n-1)! \big]  \\
		&\quad \ - \nu _{n} \Im\big[e^{\mathrm{i}(q-k)(\theta -\Theta )} \,n!\. (\ln r-\mathrm{i}(\theta -\Theta ))\big] - \nu _{n-1} \Im\big[e^{\mathrm{i}(q-k)(\theta -\Theta )} \.(n-1)! \big] 
	\end{align*}
	Thus we reduced the problem to an expression similar to the case $n=1$. Dividing by $(n-1)!$, this yields
	\[0 = \big[n (\lambda _{n} -\nu _{n} )\ln r + (\lambda _{n-1} -\nu _{n-1} )\big] \,\Im[e^{\mathrm{i}(q-k)(\theta -\Theta )}] + n(\lambda _{n} +\nu _{n} ) \,\Im\big[e^{\mathrm{i}(q-k)(\theta -\Theta )}\,\mathrm{i}(\theta -\Theta )\big].\]
	Since $\Im[e^{\mathrm{i}(q-k)(\theta -\Theta )}]$ and $\Im\big[e^{\mathrm{i}(q-k)(\theta -\Theta )}\,\mathrm{i}(\theta -\Theta )\big]$ are linearly independent functions of $\theta $, and $\ln r$ and 1 are linearly independent functions of $r$, we deduce $\lambda _{n} +\nu _{n} =\lambda _{n} -\nu _{n} =0$, so $\lambda _{n} =\nu _{n} =0$. This is absurd and concludes the proof.
\end{itemize}
\end{demo}

\begin{remarque}
Giving a basis of $\mathcal{A}_{d} (\Pi )$ would be more complicated because the functions of $\mathcal{A}(\Pi )$ are continuous on $\Gamma $ whereas the basis functions that we gave for $\mathcal{A}_{d} (\Omega )$ and $\mathcal{A}_{d} (\couche)$ do~not coincide on $\Gamma $. Fortunately, it will not be necessary to know a basis of $\mathcal{A}_{d} (\Pi )$ in the sequel.
\end{remarque}

\begin{textesansboite}
Let $\varphi $ be a function of $\mathcal{A}(\Omega )$ of the form $\Im[(\alpha z)^{q} \,\overline{\alpha z}^{k} P(\log (\alpha z))]$. Note that on $\Gamma $, $\varphi $ is equal to $x^{q+k} \Im[\alpha ^{q -k } P(\ln x-\mathrm{i}\Theta )]$. Let us define
\begin{equation}{\label{eq: def Im[alphad P(T-iTheta)]}}
	\Im[\alpha ^{q-k} P(T-\mathrm{i}\Theta )] :=\sum _{i=1} ^{\deg P} \Im(a_{i} )\,T^{i}  \in  \mathbb{R}[T], \quad \ \text{where}  \quad \ \sum _{i=1} ^{\deg P} a_{i} \.T^{i} := \alpha ^{q-k} P(T-\mathrm{i}\Theta ) \in \mathbb{C}[T].
\end{equation}
Then $\varphi _{|\Gamma } (x) = x^{q+k} Q(\ln x)$ for some $Q:=\Im[\alpha ^{q-k} P(T-\mathrm{i}\Theta )]$, which implies $\varphi _{|\Gamma } \in \mathcal{A}(\Gamma )$. Conversely, for future constructions, it will be important to solve the equation:
\begin{equation}{\label{eq: resolution de Im[alphad P(T-iTheta)]=Q}}
	\text{given $Q\in \mathbb{R}[T]$, \;find $P\in \mathbb{R}[T]$ \;s.t.\ \ $\Im[\alpha ^{d} P(T-\mathrm{i}\Theta )]=Q(T)$.} \end{equation}\fbeq
\end{textesansboite}

\begin{lemme}{resolution of \cref{eq: resolution de Im[alphad P(T-iTheta)]=Q}}
\label{resolution de Im[alphad P(T-iTheta)]=Q}%
Let $\alpha :=e^{-\mathrm{i}\Theta }$, $d\in \mathbb{R}$ and $Q\in \mathbb{R}[T]$.
\begin{enumerate}
	\item If $d\in \mathbb{R} \setminus \TZ$, then there exists a unique solution $P\in \mathbb{R}[T]$ of \cref{eq: resolution de Im[alphad P(T-iTheta)]=Q}. Moreover $\deg P=\deg Q$. We denote the solution $\RIma{d}(Q)$.
	\item If $d\in \TZ$, then the set of solutions of \cref{eq: resolution de Im[alphad P(T-iTheta)]=Q} is of the form $\{P_{0} +c \mid c\in \mathbb{R}\}$ with $P_{0} \in \mathbb{R}[T]$ and $\deg (P_{0} )=\deg (Q)+1$. We denote $\RIma{d}(Q)$ the unique solution that vanishes at 0.
\end{enumerate}
In both cases, $\RIma{d}$ is a linear map from $\mathbb{R}[T]$ into itself.
\end{lemme}
\nomenclature[21-13]{$\RIma{d}$}{Maps any $Q\in \mathbb{R}[T]$ to a polynomial $P\in \mathbb{R}[T]$ s.t.\ $\Im[\alpha ^{d} P(T-\mathrm{i}\Theta )]=Q(T)$}

\begin{demo}
There are two cases whether the coefficient of degree $m$ of $\Im[\alpha ^{d} (T-\mathrm{i}\Theta )^{m} ]$ vanishes or not.
\begin{enumerate}
	\item If $d\in \mathbb{R} \setminus \TZ$, then $\alpha ^{d} \in \mathbb{C} \setminus \mathbb{R}$, so: $\forall m\in \mathbb{N},\; \deg \Im[\alpha ^{d} (T-\mathrm{i}\Theta )^{m} ] =m$. Therefore, $(\Im[\alpha ^{d} (T-\mathrm{i}\Theta )^{m} ])_{m\in \mathbb{N}}$ is a basis of $\mathbb{R}[T]$. So~writing $Q$ in this basis gives a unique solution of \cref{eq: resolution de Im[alphad P(T-iTheta)]=Q}.
	\item If $d\in \TZ$, since $\alpha ^{d} \in \mathbb{R}$, we have: $\forall m\in \mathbb{N},\; \deg \Im[\alpha ^{d} (T-\mathrm{i}\Theta )^{m} ] =m-1$. So in this case $(\Im[\alpha ^{d} (T-\mathrm{i}\Theta )^{m} ])_{m\in \mathbb{N}^*}$ is a basis of $\mathbb{R}[T]$. Thus, $F:P\mapsto \Im[\alpha ^{d} P(T-\mathrm{i}\Theta )]$ is surjective and its kernel is the set of constant polynomials. Its restriction to $E:=\{P\in \mathbb{R}[T] \mid P(0)=0\}$ is therefore an isomorphism and we set $\RIma{d}:= (F_{|E} )^{-1}$. Finally, for~any $Q\in \mathbb{R}[T]$ we have $F^{-1} (\{Q\})= \{\RIma{d}(Q)\}+\operatorname{Ker} F$.
\end{enumerate}\fbi
\end{demo}

\subsection{Solving the Poisson and Helmholtz equations in the spaces $\mathcal{A}$}
\label{sec: operateurs de resolution dans A}

\newcommand{\inc}{\varphi} 
\newcommand{\sm}{\psi} 
\begin{textesansboite}
In this section, we show how to solve canonical problems set in $\Pi $ in the spaces $\mathcal{A}$. More precisely, let $(\sm_{\Omega } ,\sm_{\couche} ,\sm_{\Gamma } ) \in \mathcal{A}(\Omega )\times \mathcal{A}(\couche) \times \mathcal{A}(\Gamma )$, we look for the solutions $\inc \in \mathcal{A}(\Pi )$ of problems of the form
\begin{equation}{\label{eq: pb resolu dans A}}
\left\{\begin{array}{r@{\;=\;}ll}
	\Delta  \inc & \sm_{\Omega } &\text{ in }\Omega  \\[0.5mm]
	\partial _{Y} ^2 \inc & \sm_{\couche} &\text{ in }\couche \\[0.5mm]
	\partial _{Y} \inc_{|Y=0^{-}} & \sm_{\Gamma } &\text{ on }\Gamma  
\end{array}\right.\end{equation}
Note that by~definition of $\mathcal{A}(\Pi )$, $\inc$ also satisfies $\inc_{|y=0^{+}} = \inc_{|Y=0^{-}}$ on $\Gamma $ and $\inc_{|\bordgauche \cup \borddroit} =0$. Solving this system will enable us to build in \cref{sec: SF 1D de singularites de coin} the asymptotic expansion of $\bfu_{p,\ell }$ and $S_{p,\ell }$. Indeed, note for instance that this system is identical to \cref{eq: eq vol de bfu p l} except for the first line. We first describe the homogeneous solutions of \cref{eq: pb resolu dans A}, then build explicitly some particular solutions. Since we know bases of the spaces $\mathcal{A}$ (see \cref{base de A Omega,rq: base de A couche et Gamma}), we are able to code an \emph{exact, very fast and memory-thrifty} solver of \cref{eq: pb resolu dans A}, see \cref{sec: implementation des espaces A}. This is one of the key advantages of the $\mathcal{A}$ framework.
\end{textesansboite}

\begin{definition}{}
\label{def: phi d}%
For~any $d\in \TZ*$, we define on $\Omega $ the function
\[\phiOmega_{d} := (-1)^{d\Theta /\pi } \.r^{d} \sin (d\theta ) = r^{d} \sin \!\big(d(\theta -\Theta )\big) =\Im[(\alpha z)^{d} ] \in \mathcal{A}_{d} (\Omega )\]
and $\phi_{d} \in \mathcal{A}_{d} (\Pi )$ its extension by 0 in $\couche$.
\end{definition}
\nomenclature[56-13a]{$\phi_{d}$}{Fundamental singularity, equal to $\pm r^{d} \sin (d\theta )$ in $\Omega $ and vanishing in $\couche$}
\nomenclature[56-13b]{$\phiOmega_{d}$}{Fundamental singularity on $\Omega $, equal to $\pm r^{d} \sin (d\theta )$}

\begin{textesansboite}
These functions play an important role in the sequel because they solve the homogeneous Laplace equation in $\Omega $, resp.\ $\Pi $. In particular, Equation~\cref{eq: pb de Laplace homogene dans A(Pi)} below is the homogeneous version of \cref{eq: pb resolu dans A}.
\end{textesansboite}

\begin{propositionsecable}{the Laplace problem in $\mathcal{A}$}
\label{pb de Laplace homogene dans A}\vspace{-2pt}%
\begin{enumerate}
	\item $\operatorname{Span}(\{\phiOmega_{d} \mid d\in \TZ*\})$ is the set of solutions in $\mathcal{A}(\Omega )$ of 
	\begin{equation}{\label{eq: pb de Laplace homogene dans A(Omega)}}
	\left\{\begin{array}{r@{\;=\;}ll}
		\Delta  \varphi  & 0 &\text{ in }\Omega \\[0.5mm]
		\varphi  & 0 &\text{ on }\Gamma \cup \bordgauche
	\end{array}\right.\end{equation}
	\item $\operatorname{Span}(\{\phi_{d} \mid d\in \TZ*\})$ is the set of solutions in $\mathcal{A}(\Pi )$ of 
	\begin{equation}{\label{eq: pb de Laplace homogene dans A(Pi)}}
	\left\{\begin{array}{r@{\;=\;}ll}
		\Delta  \varphi  & 0 &\text{ in }\Omega \\[0.5mm]
		\partial _{Y} ^2 \varphi  & 0 &\text{ in }\couche \\[0.5mm]
		\partial _{Y} \varphi _{|Y=0^{-}} & 0 &\text{ on }\Gamma 
	\end{array}\right.\end{equation}
\end{enumerate}\fbeq
\end{propositionsecable}

\begin{demo}
\begin{enumerate}
	\item $\phiOmega_{d}$ is clearly solution of \cref{eq: pb de Laplace homogene dans A(Omega)} for~any $d\in \TZ*$. Conversely, let $\varphi $ be a solution.	Let us denote 
	\[\forall d\in {\textstyle \TN},\;\forall r\in \mathbb{R}_+^*, \quad  c_{d} (r) := \frac{2}{\Theta } \int _{0} ^{\Theta } \varphi (r,\theta )\,\sin (d\theta )\,\mathrm{d}\theta .\]
	Using separation of variables and $\varphi _{|\partial \Omega } =0$, it~is easy to show that we have
	\[\varphi =\sum _{d\in \TN} c_{d} (r) \sin (d\theta )\]
	with convergence in $H^2 (\Omega \cap \{r_{1} <r<r_{2} \})$ for~any $0<r_{1} <r_{2} <\infty $. Since $\Delta \varphi =0$, we have $(r \frac{\mathrm{d}}{\mathrm{d}r} )^2 c_{d} = d^2 c_{d}$ for~any $d$. Hence we get: $\forall d\in \TN, \exists a_{d} ,a_{-d} \in \mathbb{C}, \forall r\in \mathbb{R}_+^*,\; c_{d} (r) = a_{d} r^{d} +a_{-d} r^{-d}$.\\
	Moreover, by~definition of $\mathcal{A}(\Omega )$, there exists $q\in \mathbb{R}_+^*$ s.t.\ $\varphi \in \sum _{d\in [-q,q]} \mathcal{A}_{d} (\Omega )$. Therefore
	\[\forall d\in {\textstyle \TZ*} \cap (q,\infty ),\quad  a_{d} = \lim _{r\rightarrow \infty } r^{-d} c_{d} (r) = \lim _{r\rightarrow \infty } \frac{2}{\Theta } \int _{0} ^{\Theta } \underbrace{r^{-d} \varphi (r,\theta )}_{\longrightarrow \,0}\,\sin (d\theta )\,\mathrm{d}\theta  =0.\]
	Similarly, looking at $r\rightarrow 0$ one gets: $\forall d\in \TZ* \cap (-\infty ,-q),\; a_{d} =0$. This yields
	\[\varphi  = \sum _{d\in \TZ*\cap [-q,q]} a_{d} \.\operatorname{sgn}(d) \,\phiOmega_{d}\]
	(where $\operatorname{sgn} (d):=d/|d|$), which is a finite sum. Therefore, $\varphi $ is in the desired span.
	\item Any solution of the system vanishes in $\couche$, so point 2 easily follows from point 1.
\end{enumerate}
\end{demo}

\begin{textesansboite}
Let us now define the following linear forms $\sigma _{d}$ which map each $\varphi \in \mathcal{A}(\Omega )$ to its coordinate on the basis function $\phiOmega_{d}$. 
These linear forms appear later as key singularity coefficients in the matching conditions.
\end{textesansboite}

\begin{definition}{linear forms $\sigma _{d}$}
\label{def: sigma}%
Let $d\in \TZ*$. For~any $(q,k,j)\in \mathbb{R}\times \mathbb{N}^2 $ s.t.\ we have $(q\in \mathbb{N} \text{ and } q\leqslant k)\Rightarrow j>0$, we define:
\[\sigma _{d} \big( \Im[(\alpha z)^{q} \. \overline{\alpha z}^{k} \log (\alpha z)^{j} ] \big) := \delta _{d,q} \,\delta _{k,0} \,\delta _{j,0} .\]
By~\cref{decomposition de A en sev,base de A Omega}, the functions above form a basis of $\mathcal{A}(\Omega )$, so $\sigma _{d}$ can be extended into a linear form $\mathcal{A}(\Omega ) \rightarrow \mathbb{C}$. Finally, for~any $\varphi \in \mathcal{A}(\Pi )$, we define $\sigma _{d} (\varphi ):=\sigma _{d} (\varphi _{|\Omega } )$.
\end{definition}
\nomenclature[53-13a]{$\sigma _{d}$}{Linear form on $\mathcal{A}(\dots )$ and $\Amaicro(\dots )$ that gives the coordinate on the basis function $\phi_{d}$}

\begin{textesansboite}
Let us now build particular solutions of \cref{eq: pb resolu dans A}. By~rearranging the equations, the resolution of \cref{eq: pb resolu dans A} can be decomposed into two steps:
\[\text{solve }  \left\{\begin{array}{r@{\;=\;}ll}
	\partial _{Y} ^2 \inc & \sm_{\couche} &\text{ in }\couche \\[0.5mm]
	\partial _{Y} \inc_{|Y=0^{-}} & \sm_{\Gamma } &\text{ on }\Gamma  \\[0.5mm]
	\inc & 0 &\text{ on }\borddroit 
\end{array}\right. \qquad  \text{then solve }  \left\{\begin{array}{r@{\;=\;}ll}
	\Delta  \inc & \sm_{\Omega } &\text{ in }\Omega  \\[0.5mm]
	\inc_{|y=0^{+}} & \inc_{|Y=0^{-}} &\text{ on }\Gamma  \\[0.5mm]
	\inc & 0 &\text{ on }\bordgauche 
\end{array}\right.\]
The first step is a 1D Poisson problem w.r.t.\ $Y$ for~any $x$, so it~is easy to solve. Moreover, the proposition below builds particular solutions to the problem of the second step. By~linearity, it suffices to build particular solutions of two sub-problems. We recall that, according to \cref{def: degre d'un operateur}, for~any $D_{1} ,D_{2} \in \{\Pi ,\Omega ,\Gamma ,\couche\}$, we say that a linear map $F:\mathcal{A}(D_{1} )\rightarrow \mathcal{A}(D_{2} )$ has degree $d\in \mathbb{R}$ iff: $\forall q\in \mathbb{R},\, \forall \varphi \in \mathcal{A}_{q} (D_{1} ),\; F(\varphi )\in \mathcal{A}_{q+d} (D_{2} )$.
\end{textesansboite}

\newcommand{\smOmega}{\sm_\Omega}
\newcommand{\smGamma}{\sm_\Gamma}
\newcommand{\incOmega}{\opensubscript{\inc_{\Delta}}}
\newcommand{\incGamma}{\opensubscript{\inc_{\mathrm{D}}}}
\begin{propositionsecable}{}
\label{operateurs R dans Omega}%
Let us denote $\Aperp(\Omega ) = \{\varphi \in \mathcal{A}(\Omega ) \mid \forall d\in \TZ*,\; \sigma _{d} (\varphi )=0\}$, which is a complementary subspace of $\operatorname{Span}(\{\phiOmega_{d} \mid d\in \TZ*\})$ in $\mathcal{A}(\Omega )$.
\begin{enumerate}
	\item For~any $\smOmega \in \mathcal{A}(\Omega )$ there exists a unique solution $\incOmega \in \Aperp(\Omega )$ of
	\begin{equation}{\label{eq: operateur R Omega Omega}}
	\left\{\begin{array}{r@{\;=\;}ll}
		\Delta  \inc & \smOmega &\text{ in }\Omega \\[0.5mm]
		\inc & 0 &\text{ on }\Gamma 
	\end{array}\right. \end{equation}
	The associated map $\ROmegaOmega:\smOmega \in \mathcal{A}(\Omega ) \mapsto \incOmega \in \Aperp(\Omega )$ is linear and has degree 2.
	\item For~any $\smGamma \in \mathcal{A}(\Gamma )$ there exists a unique solution $\incGamma \in \Aperp(\Omega )$ of
	\begin{equation}{\label{eq: operateur R Gamma Omega}}
	\left\{\begin{array}{r@{\;=\;}ll}
		\Delta  \inc & 0 &\text{ in }\Omega \\[0.5mm]
		\inc & \smGamma &\text{ on }\Gamma 
	\end{array}\right. \end{equation}
	The associated map $\RGammaOmega:\smGamma \in \mathcal{A}(\Gamma ) \mapsto \incGamma \in \Aperp(\Omega )$ is linear and has degree 0.
\end{enumerate}\fbi
\end{propositionsecable}
\nomenclature[27-23]{$\ROmegaOmega$}{Operator $\mathcal{A}(\Omega )\rightarrow \mathcal{A}(\Omega )$ inverting $\Delta $}
\nomenclature[27-13]{$\RGammaOmega$}{Operator $\mathcal{A}(\Gamma )\rightarrow \mathcal{A}(\Omega )$ inverting the trace operator.}

\begin{textesansboite}
Note that the proof below is entirely constructive as it provides explicit expressions of the solutions in terms of the data.
\end{textesansboite}

\begin{demo}
\Cref{pb de Laplace homogene dans A} shows the uniqueness of the solutions $\incOmega$ and $\incGamma$, so only their existence remains to prove. We will build $\incOmega$ (resp.\ $\incGamma$) when $\smOmega$ (resp.\ $\smGamma$) is part of the basis of $\mathcal{A}(\Omega )$ (resp.\ $\mathcal{A}(\Gamma )$) given by \cref{decomposition de A en sev,base de A Omega} (resp.\ \cref{rq: base de A couche et Gamma}). This induces a unique way to extend $\smOmega \mapsto \incOmega$ and $\smGamma \mapsto \incGamma$ into linear maps, and the images of these extensions will remain solutions of \cref{eq: operateur R Omega Omega,eq: operateur R Gamma Omega} by linearity of the problems.
\begin{itemize}
	\item First, let us solve \cref{eq: operateur R Gamma Omega} for the basis function $\smGamma:x\mapsto  x^{d} \ln ^{j} \!\,x$, with $d\in \mathbb{R}$ and $j\in \mathbb{N}$. Let us look for a solution $\incGamma$ of the form $z\mapsto \Im[(\alpha z)^{d} P(\log (\alpha z))]$ with $P\in \mathbb{R}[T]$. As~the real part of an analytic function, this function is harmonic. Moreover, we have $\incGamma_{|\Gamma } = x^{d} \.Q(\ln x)$ with $Q:=\Im[\alpha ^{d} P(T-i\Theta )]$ using the notation of \cref{eq: def Im[alphad P(T-iTheta)]}. Thus, by \cref{resolution de Im[alphad P(T-iTheta)]=Q}, if we take $P:=\RIma{d}(T^{j} )$, then we have $\Im[\alpha ^{d} P(T-i\Theta )]=T^{j}$ and $\incGamma_{|\Gamma } =\smGamma$. Hence $\incGamma$ solves \cref{eq: operateur R Gamma Omega}. In addition, it~is clear that $\smGamma \mapsto \incGamma$ has degree 0.
	\item Now let us solve \cref{eq: operateur R Omega Omega} for the basis function $\smOmega:z\mapsto  \Im[(\alpha z)^{q} \,\overline{\alpha z}^{k} \log (\alpha z)^{j} ]$ with $(q,k,j)\in \mathbb{R}\times \mathbb{N}^2 $ s.t.\ we have $(q\in \mathbb{N} \text{ and } q\leqslant k)\Rightarrow  j>0$ (see \cref{base de A Omega}). Given that $\Delta  =4\partial _{z} \partial _{ \bar{z}}$, we have for~any $q_{1} \in \mathbb{R}$, $k_{1} \in \mathbb{N}$, $P_{1} \in \mathbb{R}[T]$ and $\inc_{1} :z\mapsto  \Im[(\alpha z)^{q_{1}} \,\overline{\alpha z}^{k_{1}} P_{1} (\log (\alpha z))]$,
	\[\Delta \inc_{1} = 4\,\Im[(\alpha z)^{q_{1} -1} \, k_{1} \overline{\alpha z}^{k_{1} -1} (q_{1} P_{1} +P_{1} ')(\log (\alpha z))].\]
	So~taking $q_{1} :=q+1$, $k_{1} :=k+1$ and $P_{1} \in \mathbb{R}[T]$ a solution of $4k_{1} (q_{1} P_{1} +P_{1} ') = P_{\sm}$, we have $\Delta \inc_{1} = \smOmega$ in $\Omega $.
	\begin{itemize} 
		\item If $q_{1} =0$, $P_{1}$ is unique up to a constant, Moreover, writing $\sum_{i=0}^n a_{i} \.T^{i} :=P$, we have
		\[\varphi _{1} = -a_{0} \,\Im[(\alpha z)^{k_{1}} ] + \sum_{i=1}^m a_{i} \,\Im\big[\overline{\alpha z}^{k_{1}} \log (\alpha z)^{i} \big] \]
		where each term belongs to the basis of $\mathcal{A}(\Omega )$. So for~any $d\in \TZ*$ different from $k_{1}$, we have $\sigma _{d} (\varphi _{1} )=0$, while $\sigma _{k_{1}} (\varphi _{1} )=-a_{0} =-P_{1} (0)$ if $k_{1} \in \mathbb{N}\cap \TZ*$. Taking $P_{1} (0):=0$ thus gives: $\forall d\in \TZ*,\; \sigma _{d} (\inc_{1} )=0$.
		\item If $q_{1} \neq 0$, there exists a unique solution $P_{1}$. Given that $q_{1} \neq 0$ and $k_{1} \neq 0$, \cref{def: sigma} implies that: $\forall d\in \TZ*,\; \sigma _{d} (\inc_{1} )=0$.
	\end{itemize}
	Since $\inc_{1|\Gamma }$ may not vanish, we finally define $\incOmega :=\inc_{1} -\RGammaOmega(\inc_{1|\Gamma } )$. It is easy to verify that it solves \cref{eq: operateur R Omega Omega} and $\smOmega \mapsto \incOmega$ has degree 2.
\end{itemize}
\end{demo}

\begin{textesansboite}
The proposition below builds particular solutions of \cref{eq: pb resolu dans A}. Similarly to \cref{operateurs R dans Omega}, we decompose by linearity the problem into three sub-problems.
\end{textesansboite}

\newcommand{\smCouche}{\sm_\couche}
\newcommand{\incCouche}{\opensubscript{\inc_{\partial_Y^2}}}
\renewcommand{\incGamma}{\opensubscript{\inc_{\mathrm{N}}}}
\begin{propositionsecable}{particular solutions of \cref{eq: pb resolu dans A}}
\label{operateurs R}%
Let us denote $\Aperp(\Pi ) := \{\varphi \in \mathcal{A}(\Pi ) \mid \forall d\in \TZ*,\; \sigma _{d} (\varphi )=0\}$, which is a complementary subspace of $\operatorname{Span}(\{\phi_{d} \mid d\in \TZ*\})$ in $\mathcal{A}(\Pi )$.
\begin{enumerate}
	\item For~any $\smOmega \in \mathcal{A}(\Omega )$ there exists a unique solution $\incOmega \in \Aperp(\Pi )$ of
	\begin{equation}{\label{eq: operateur R Omega Pi}}
	\left\{\begin{array}{r@{\;=\;}ll}
		\Delta  \inc & \smOmega &\text{ in }\Omega \\[0.5mm]
		\partial _{Y} ^2 \inc & 0 &\text{ in }\couche \\[0.5mm]
		\partial _{Y} \inc_{|Y=0^{-}} & 0 &\text{ on }\Gamma 
	\end{array}\right. \end{equation}
	The associated map $\ROmegaPi:\smOmega \in \mathcal{A}(\Omega ) \mapsto \incOmega \in \Aperp(\Pi )$ is linear and has degree 2.
	\item For~any $\smCouche \in \mathcal{A}(\couche)$ there exists a unique solution $\incCouche \in \Aperp(\Pi )$ of
	\begin{equation}{\label{eq: operateur R couche Pi}}
	\left\{\begin{array}{r@{\;=\;}ll}
		\Delta  \inc & 0 &\text{ in }\Omega \\[0.5mm]
		\partial _{Y} ^2 \inc & \smCouche &\text{ in }\couche \\[0.5mm]
		\partial _{Y} \inc_{|Y=0^{-}} & 0 &\text{ on }\Gamma 
	\end{array}\right. \end{equation}
	The associated map $\RcouchePi:\smCouche \in \mathcal{A}(\couche) \mapsto \incCouche \in \Aperp(\Pi )$ is linear and has degree 0.
	\item For~any $\smGamma \in \mathcal{A}(\Gamma )$ there exists a unique solution $\incGamma \in \Aperp(\Pi )$ of
	\begin{equation}{\label{eq: operateur R Gamma Pi}}
	\left\{\begin{array}{r@{\;=\;}ll}
		\Delta  \inc & 0 &\text{ in }\Omega \\[0.5mm]
		\partial _{Y} ^2 \inc & 0 &\text{ in }\couche \\[0.5mm]
		\partial _{Y} \inc_{|Y=0^{-}} & \smGamma &\text{ on }\Gamma 
	\end{array}\right. \end{equation}
	The associated map $\RGammaPi:\smGamma \in \mathcal{A}(\Gamma ) \mapsto \incGamma \in \Aperp(\Pi )$ is linear and has degree 0.
\end{enumerate}\fbi
\end{propositionsecable}
\nomenclature[27-39]{$\ROmegaPi$}{Operator $\mathcal{A}(\Omega )\rightarrow \mathcal{A}(\Pi )$ inverting $\Delta $}
\nomenclature[27-13-33]{$\RcouchePi$}{Operator $\mathcal{A}(\couche)\rightarrow \mathcal{A}(\Pi )$ inverting $\partial _{Y\mid\couche} ^2 $} 
\nomenclature[27-23]{$\RGammaPi$}{Operator $\mathcal{A}(\Gamma )\rightarrow \mathcal{A}(\Pi )$ inverting $\partial _{Y\mid Y=0^{-}}$} 

\begin{remarque}
Using \cref{pb de Laplace homogene dans A,operateurs R}, it~is then easy to see that the set of solutions of \cref{eq: pb resolu dans A} is \[\ROmegaPi(\varphi _{\Omega } )+\RcouchePi(\varphi _{\couche} )+\RGammaPi(\varphi _{\Gamma } )+\operatorname{Span}(\{\phi_{d} \mid d\in \txtTZ*\}).\]
Moreover, the functions $\incOmega$, $\incCouche$ and $\incGamma$ in \cref{operateurs R} have explicit expressions (see the proof below), which allows to compute them easily in practice.
\end{remarque}

\begin{demo}
\Cref{pb de Laplace homogene dans A} shows the uniqueness of $\incOmega$, $\incCouche$ and $\incGamma$, so only their existence remains to prove.
\begin{enumerate}
	\item Let us define $\ROmegaPi:\mathcal{A}(\Omega ) \rightarrow \Aperp (\Pi )$ by:
	\[\forall \.\smOmega \in \mathcal{A}(\Omega ),\qquad  \ROmegaPi(\smOmega) := \left\{\begin{array}{ll}
		\ROmegaOmega(\smOmega) &\text{ in }\Omega  \\[0.5mm]
		0 &\text{ in }\couche
	\end{array}\right.\]
	Then it~is easy to verify that $\ROmegaPi(\smOmega)$ is a solution of \cref{eq: operateur R Omega Pi}. Moreover, $\ROmegaOmega$ is linear and has degree 2, so $\ROmegaPi$ has the same properties.

	\item By~\cref{decomposition de A en sev}, it suffices to build $\incCouche$ when $\smCouche= x^{d} \,Q_{\sm} (\ln x,Y)$ with $d\in \mathbb{R}$ and $Q_{\sm} \in \mathbb{C}[T,Y]$. In addition, by~\cref{rq: base de A couche et Gamma}, $Q_{\psi}$ is uniquely defined from $\smCouche$, so $\incCouche$ can be built in terms of $Q_{\sm}$. We look for $\incCouche$ of the form $x^{d} \,Q_{\inc} (\ln x,Y)$ in $\couche$ with $Q_{\inc} \in \mathbb{C}[T,Y]$. Then necessarily we have \[\left\{\begin{array}{r@{\;=\;}l}
		\partial _{Y} ^2  Q_{\inc} & Q_{\sm} \\[0.5mm]
		\partial _{Y} Q_{\inc} (\cdot ,0)&0 \\[0.5mm]
		Q_{\inc} (\cdot ,-1) & 0
	\end{array}\right.\]
	This uniquely determines $Q_{\inc}$. Now we define $\incCouche$ in $\Omega $ as $\RGammaOmega(x^{d} \.Q(\ln x,0))$. This makes $\incCouche$ continuous on $\Gamma $, so we have $\incCouche \in \mathcal{A}(\Pi )$. Moreover, $\incCouche_{|\Omega } \in \Aperp(\Omega )$, so $\incCouche \in \Aperp(\Pi )$. Finally, it~is easy to see that $\incCouche$ is a solution of \cref{eq: operateur R couche Pi}.
	
	\item Similarly, it suffices to build $\incGamma$ when $\smGamma= x^{d} \,P_{\sm} (\ln x)$ with $d\in \mathbb{R}$ and $P_{\sm} \in \mathbb{C}[T]$. Like in point 2, we take $\incGamma_{|\couche} :=x^{d} \,Q_{\inc} (\ln x,Y)$ and $\incGamma_{|\Omega } := \RGammaOmega(x^{d} \.Q(\ln x,0))$ for some $Q_{\inc} \in \mathbb{C}[T,Y]$. Then we have $\partial _{Y} {\incGamma}_{|Y=0^{-}} = x^{d} \,\partial _{Y} Q_{\inc} (\ln x,0)$ so it suffices to take $Q_{\inc}$ as a solution of: \[\left\{\begin{array}{r@{\;=\;}l}
		\partial _{Y} ^2  Q_{\inc} & 0 \\[0.5mm]
		\partial _{Y} Q_{\inc} (\cdot ,0) & P_{\sm} \\[0.5mm]
		Q_{\inc} (\cdot ,-1) & 0
	\end{array}\right.\]
\end{enumerate}
\end{demo}

\begin{textesansboite}
To end this section, let us show that the spaces $\mathcal{A}$ are stable under some differential operators.
\end{textesansboite}

\begin{lemme}{}
\label{operateurs de derivation dans A}%
\begin{enumerate}
	\item $\partial _{x|\couche} ^2  :\varphi  \mapsto  \partial _{x} ^2  \varphi _{|\couche}$ maps $\mathcal{A}(\Pi )$ to $\mathcal{A}(\couche)$ and has degree $-2$ (see \cref{def: degre d'un operateur}).
	\item $\partial _{y|\Gamma ,y=0^{+}} : \varphi  \mapsto  \partial _{y} \varphi _{|\Gamma ,y=0^{+}}$ maps $\mathcal{A}(\Pi )$ to $\mathcal{A}(\Gamma )$ and has degree $-1$.
	\item $\Delta :\varphi \mapsto  \Delta \varphi _{|\Omega }$ maps $\mathcal{A}(\Omega )$ (and $\mathcal{A}(\Pi )$) to $\mathcal{A}(\Omega )$ and has degree $-2$.
\end{enumerate}\fbi
\end{lemme}

\begin{textesansboite}
Note that, when used on corner fields (which depend on $(X,Y)$), the first two operators will rather be denoted $\partial _{X|\couche} ^2 $ and $\partial _{Y|\Gamma ,Y=0^{+}}$.
\end{textesansboite}

\begin{demo}
It suffices to verify it when $\varphi $ has the form $\Im[(\alpha z)^{q} \,\overline{\alpha z}^{k} P(\log (\alpha z)) ]$ in $\Omega $ and $x^{d} \,Q(\ln x,Y)$ in $\couche$, with $q\in \mathbb{R}$, $k\in \mathbb{N}$, $d=q+k$, $P\in \mathbb{C}[T]$ and $Q\in \mathbb{C}[T,Y]$.
\begin{enumerate}
	\item We have $\partial _{x} ^2 \varphi (x,Y) = x^{d-2} (d(d-1)+(2d-1)\.\partial _{T} +\partial _{T} ^2 )\.Q(\ln x,Y)$ which is in $\mathcal{A}(\couche)$.
	\item Let $\psi : z\mapsto  (\alpha z)^{q} \,\overline{\alpha z}^{k} P(\log (\alpha z))$. From $\partial _{z} \psi  = \frac{1}{2} (\partial _{x} \psi -\mathrm{i}\partial _{y} \psi )$ and $\partial _{ \bar{z}} \psi  = \frac{1}{2} (\partial _{x} \psi +\mathrm{i}\partial _{y} \psi )$ it follows that $\partial _{y} \varphi  = \Im[\partial _{y} \psi ] = \Re[\partial _{z} \psi  -\partial _{ \bar{z}} \psi ]$. Therefore, \begin{align*}
		\partial _{y} \varphi _{|\Gamma ,y=0^{+}} &= \Re[\alpha  (\alpha x)^{q-1} \,\overline{\alpha x}^{k} (qP+P')(\log (\alpha x)) -(\alpha x)^{q} \, \overline{\alpha }\,k\,\overline{\alpha x}^{k-1} P(\log (\alpha x))] \\
		&= x^{q+k-1} \Re[\alpha ^{q-k} ((q-k)P+P')(\ln x-\mathrm{i}\Theta )] \ \ \in \ \mathcal{A}(\Gamma )
	\end{align*}
	\item Since $\Delta  =4\partial _{z} \partial _{ \bar{z}}$, we have $\Delta \varphi  = 4\,\Im[(\alpha z)^{q-1} \,k\.\overline{\alpha z}^{k-1} (qP+P')(\log (\alpha z))]$ which is in $\mathcal{A}(\Omega )$.
\end{enumerate}\fbi
\end{demo}
\section{Formal series of corner singularities}
\label{sec: SF de singularites}

\begin{textesansboite}
In this section, we apply the tools of \cref{chap: SF} to the spaces of corner singularities $\mathcal{A}_{d} ({\dots })$ introduced in \cref{sec: singularites quasi homogenes}. We define subspaces of $\fsumd \mathcal{A}_{d} ({\dots })$ and $\fsum_{d,p\in \mathbb{R}} \varepsilon ^{p} \mathcal{A}({\dots })[\ln \varepsilon ]$ that are relevant for the far, layer and corner fields, and we prove the well-posedness of Helmholtz and Poisson-like problems in these spaces of formal series. These tools will be put into application in \cref{sec: raccordement algebrique} to derive the matching conditions. \Cref{sec: SF 1D de singularites de coin} focuses on formal series over $d$, i.e.\ without $\varepsilon $, while \cref{sec: SF 2D de singularites de coin} treats formal series of singularities that involve $\varepsilon $.
\end{textesansboite}

\subsection{Asymptotic behaviors w.r.t.\ $\bfr$ of solutions of model problems}
\label{sec: SF 1D de singularites de coin}

\begin{textesansboite}
In this section we give tools that will be used to represent and compute the asymptotic behaviors of $\bfu_{p,\ell }$ when $\bfr\rightarrow 0$ and $S_{p,\ell }$ when $\bfr\rightarrow \infty $. We apply the tools of \cref{sec: SF 1D} with $E:=\mathcal{A}(D)$ and $E_{d} :=\mathcal{A}_{d} (D)$ for~any $d\in \mathbb{R}$ and for a given $D\in \domsA$. This is possible because the assumption ``$E=\bigoplus_{d\in \mathbb{R}} E_{d}$'' made by \cref{sec: SF 1D} is satisfied here according to \cref{decomposition de A en sev}. We recall that the $\fsum$-sum of vector spaces is defined by \cref{def: fsum d'ensembles} and is studied in \cref{fsum de sev des Ed}.
\end{textesansboite}

\begin{definition}{the spaces $\Amaicro$}
\label{def: A barre}%
Let $D\in \domsA$. We denote $\Amacro(D)$, resp.\ $\Amicro(D)$, the set of elements of $\fsum_{d\in \mathbb{R}} \mathcal{A}_{d} (D)$ whose support is bounded from below, resp.\ above:
\[\Amacro(D) := \bigcup _{d_{0} \in \mathbb{R}} \.\fsum_{d\in [d_{0} ,\infty )} \mathcal{A}_{d} (D) \qquad  \text{ and } \qquad  \Amicro(D) := \bigcup _{d_{0} \in \mathbb{R}} \.\fsum_{d\in (-\infty ,d_{0} ]} \mathcal{A}_{d} (D).\]\fbeq
\end{definition}
\nomenclature[10-60]{$\Amaicro(\dots )$}{Space of formal series of corner singularities (not involving $\feps$)}

\begin{textesansboite}
By~construction, formal series have locally finite support (see \cref{def: support locmt fini et SF}), so any element of $\Amacro(D)$, resp.\ $\Amicro(D)$, can be written as $\fsumn \varphi _{n}$ where each $\varphi _{n}$ belongs to $\mathcal{A}_{d_{n}} (D)$ for some sequence $(d_{n} )\in \mathbb{R}^{\mathbb{N}}$ that is increasing and tends to $\infty $, resp.\ decreasing and tends to $-\infty $.
Since the asymptotic behavior of $\bfu_{p,\ell }$ when $\bfr\rightarrow 0$ involves increasing powers of $r$, it will be expressed in $\Amacro(D)$. Similarly, the behavior $S_{p,\ell }$ when $\bfr\rightarrow \infty $ involves decreasing powers of $r$, so it will be expressed in $\Amicro(D)$.
\end{textesansboite}

\begin{remarque}
If we take $E:=\mathcal{A}(D)$ and $E_{d} :=\mathcal{A}_{d} (D)$ for~any $d\in \mathbb{R}$, then $\Amaicro(D)$ coincides with the set $E^{\pm }$ studied in \cref{sec: series d'op de SF 1D}, so the results of these section apply. In particular, $\Amaicro(D)$ is a vector space. Moreover, we have the following diagram of inclusions for~any $d\in \mathbb{R}$:
\[\xymatrixrowsep{-0.5ex}\xymatrixcolsep{5ex}\xymatrix{
	& & \Amacro(D) \ar[rd] & \\
	\mathcal{A}_{d} (D) \ar[r] & \mathcal{A}(D) =\displaystyle \bigoplus_{q\in \mathbb{R}} \mathcal{A}_{q} (D) \ar[ru] \ar[rd] & &\displaystyle \raggedleft\fsum_{q\in \mathbb{R}} \mathcal{A}_{q} (D) \\
	& & \Amicro(D) \ar[ru] & \
}\]
The inclusion $\mathcal{A}(D)\subset \Amaicro(D)$ is strictly speaking an abuse of notation explained and justified in \cref{sec: SF 1D}. All the inclusions of the diagram are strict. For example, the elements of $\Amaicro(D)$ are~not in general functions from $D$ to $\mathbb{C}$, as the formal series may be divergent.
\end{remarque}

\begin{exemple}
\def\q{q}%
\Cref{ex: trucs qui sont ou pas des series formelles} already gives examples of objects that are or are~not elements of $\fsumd \mathcal{A}(\Gamma )$. Let us give similar examples for $\Amaicro$. Let $q\in \mathbb{R}$.
\begin{itemize}
	\item ``$\fsum_{d>0} \Im[(\alpha z)^{d} ]$'' is not well-defined because formal series must have locally finite support (see \cref{def: support locmt fini et SF}).
	\item The objects $\fsumn \Im[(\alpha z)^{\q+n} ]$, $\fsumn \Im[(\alpha z)^{\q} \.\overline{\alpha z}^{n} ]$ and $\fsumn \Im[(\alpha z)^{\q+n} \.\overline{\alpha z}^{n} ]$ are well-defined formal series and belong to $\Amacro(\Omega )$. However, ``$\fsumn \Im[(\alpha z)^{\q-n} \.\overline{\alpha z}^{n} ]$'' is not well-defined because we are trying to sum an infinite number of terms that belong to the same space $\mathcal{A}_{q} (\Omega )$.
	\item $\fsum_{d\in \TZ*\!,\.d\geqslant \q} \phi_{d}$ is well-defined and belongs to $\Amacro(\Pi )$. Similarly $\fsum_{d\in \TZ*\!,\.d\leqslant \q} \phi_{d}$ belongs to $\Amicro(\Pi )$, but not to $\Amacro(\Pi )$.
	\item $\fsumn \fsum_{d\in \TZ*\!,\.d\geqslant \q} \Im[(\alpha z)^{d+n} ]$ is well-defined and belongs to $\Amacro(\Omega )$ because its support is $\mathbb{N}+{(\TZ*\cap [q,\infty ))}$. \Cref{pb de Helmholtz dans A Pi} below involves a term of the form $\big\langle {-k_{0} ^2  \ROmegaPi}\big\rangle  \big( \fsum_{d\in \TZ*\!,\.d\geqslant \q} s_{d} \.\phi_{d} \big)$, with $s_{d} \in \mathbb{C}$ that has a similar support.
\end{itemize}
\end{exemple}

\begin{textesansboite}
According to \cref{operateurs R dans Omega,operateurs R,operateurs de derivation dans A}, the operators $\Delta _{|\Omega }$, $\partial _{x|\couche} ^2 $, $\partial _{y|\Gamma ,y=0^{+}}$, $\ROmegaOmega$, $\RGammaOmega$ $\RcouchePi$, $\RGammaPi$ and $\ROmegaPi$ have a degree (see \cref{def: degre d'un operateur}). So \cref{def: extension d'un op a degre aux SF 1D}, allows to naturally extend them to the spaces $\Amaicro$ (their extension is still denoted in the same way). Likewise, the operators $\partial _{Y|\couche} ^2 $ and $\partial _{Y|\Gamma ,Y=0^{-}}$ and the restriction operators $\mathcal{A}(\Pi )\rightarrow \mathcal{A}(\Omega )$, $\mathcal{A}(\Omega )\rightarrow \mathcal{A}(\Gamma )$, etc.\ have degree 0, so they can be extended to the spaces $\Amaicro$. In a similar way, the definition below extends the linear forms $\sigma _{d}$ to the spaces $\Amaicro$.
\end{textesansboite}

\begin{definition}{}
\label{def: sigma sur Amaicro}%
Let $D\in \{\Omega ,\Pi \}$ and $d\in \TZ*$. For~any $(\varphi _{q} )_{q\in \mathbb{R}} \in \fsum_{q\in \mathbb{R}} \mathcal{A}_{q} (D)$ (which implies $\varphi _{q} \in \mathcal{A}_{q} (D)$ for~any $q$), we define
\[\sigma _{d} \bigg(\fsum_{q\in \mathbb{R}} \varphi _{q} \bigg):=\sigma _{d} (\varphi _{d} ).\]\fbeq
\end{definition}

\begin{textesansboite}
This definition coincides with the original definition of $\sigma _{d}$ on $\mathcal{A}(D)$ given by \cref{def: sigma}. Indeed, if $(\varphi _{q} )_{q}$ has finite support, then we have $\varphi :=\sum _{q} \varphi _{q} =\fsum_{q} \varphi _{q}$ (by \cref{proprietes de fsum} item 1) and $\sigma _{d} (\varphi _{d} )=\sigma _{d} (\varphi _{d} )+\sum _{q\neq d} \sigma _{d} (\varphi _{q} )=\sigma _{d} (\varphi )$.
\end{textesansboite}

\begin{lemmesansboite}{}
\label{def: sigma et fsum commutent dans Aemaicro}%
Let $D\in \{\Omega ,\Pi \}$ and $d\in \TZ*$. Let $I$ be a set and $(\Phi _{i} )_{i\in I} \in \big(\fsum_{q\in \mathbb{R}} \mathcal{A}_{q} (D)\big)^{I}$ be a family of formal series that has collectively locally finite support. We have
\[\sigma _{d} \bigg(\fsum_{i\in I} \Phi _{i} \bigg) = \sum _{i\in I} \sigma _{d} (\Phi _{i} )\]
where the right-hand side sum has only a finite number of non-zero terms.
\end{lemmesansboite}

\begin{demo}
Let us denote $(\varphi _{i,d} )_{d\in \mathbb{R}} :=\Phi _{i}$ for~any $i\in I$. By~definition of $\fsum$ in \cref{def: fsum}, the component of $\fsum_{i} \Phi _{i}$ in $\mathcal{A}_{d} (D)$ is $\sum _{i} \varphi _{i,d}$ (which is a finite sum). Thus we have $\sigma _{d} \big(\fsum_{i} \Phi _{i} \big) =\sigma _{d} \big(\sum _{i} \varphi _{i,d} \big)=\sum _{i} \sigma _{d} (\varphi _{i,d} )=\sum _{i} \sigma _{d} (\Phi _{i} )$.
\end{demo}

\begin{textesansboite}
The next three result show the existence and uniqueness of the solution of some Helmholtz and Poisson-like problems in the spaces $\Amaicro$. These problems have a form that is similar to the equations satisfied by resp.\ $\bfu_{p,\ell }$, $u_{p,\ell }$ and $S_{p,\ell }$, which are resp.\ \cref{eq: eq vol de bfu p l,eq: eq vol champs lointains,eq: eq vol champs de coin}. These results will be used in \cref{sec: DA autour du coin} to build formal series that represent the behavior of the fields w.r.t.\ $\bfr$ (in the sense that any truncation of the formal series at a given order gives an asymptotic expansion of the associated field at this order). For instance, each far field $u_{p,\ell }$ satisfies $\Delta u_{p,\ell } =0$ in $\Omega $, so it~is natural that its behavior, denoted $\uDA_{p,\ell } \in \Amacro(\Omega )$, satisfies $\Delta \uDA_{p,\ell } =0$ in $\Amacro(\Omega )$. \Cref{pb de Helmholtz dans A Omega} proves the existence and uniqueness of a harmonic element of $\Amacro(\Omega )$ that has a given Dirichlet boundary condition on $\Gamma $, and gives an expression of this solution.\\

Note that, contrary to proving the existence of the fields $u_{p,\ell }$ and $S_{p,\ell }$, which is done in \cref{sec: pb modeles bien poses} and requires non elementary tools of functional analysis like the Lax-Milgram theorem, solving PDEs in $\Amaicro$ is essentially a question of calculus. For example, in \cref{pb de Helmholtz dans A Pi,pb de Helmholtz dans A Omega}, the Helmholtz equation is solved basically by writing the problem as
\[( \mathrm{id}-F)(\Phi ^{0} )=\Psi ^{0} ,\]
where $\Phi ^{0}$ is the unknown, $\Psi ^{0}$ is a formal series that can be computed from the data and $F$ is an operator that has a degree (namely $F=-k_{0} ^2 \ROmegaPi$, resp.\ $F=-k_{0} ^2 \ROmegaOmega$). Then we directly have
\[\Phi ^{0} =\fsumn F^{n} (\Psi ^{0} )=:\langle F\rangle (\Psi ^{0} )\]
(see \cref{def: serie geometrique d'op 1D} for the notation $\langle .\rangle $). This formal calculation is justified by the work done in \cref{chap: SF}, especially by \cref{inverse de id-f dans les SF 1D}. The easiness to invert operators among formal series in one of the key benefits of our algebraic approach.
\end{textesansboite}

\begin{proposition}{resolution of a Helmholtz-like problem in $\Amacro(\Pi )$}
\label{pb de Helmholtz dans A Pi}%
Let $g^{0} \in \Amacro(\couche)$ and $h^{0} \in \Amacro(\Gamma )$. Let $(s_{d} )_{d\in \TZ*} \in \mathbb{C}^{\TZ*}$ be a family whose support is bounded from below. The following problem has a unique solution in $\Amacro(\Pi )$: 
\begin{equation}{\label{eq: equations de bfuDA}}
	\left\{\begin{array}{r@{\;=\;}ll}
	\mu _{0} \Delta  \bfuDA +\omega ^2 \rho _{0} \bfuDA & 0 &\text{ in }\Omega \\[0.5mm]
	\mu _{1} \partial _{Y} ^2  \bfuDA & g^{0} &\text{ in }\couche \\[0.5mm]
	\mu _{1} \partial _{Y|Y=0^{-}} \bfuDA & h^{0} &\text{ on }\Gamma  \\[0.5mm]
	\sigma _{d} (\bfuDA) & s_{d} & \ \forall d\in \TZ*
	\end{array}\right.
\end{equation}
Moreover, the solution is equal to
\begin{equation}{\label{eq: bfuDA}}
	\bfuDA = \big\langle  {-}k_{0} ^2 \.\ROmegaPi \big\rangle  \bigg( \frac{1}{\mu _{1}} \RcouchePi(g^{0} ) + \frac{1}{\mu _{1}} \RGammaPi(h^{0} ) + \fsum_{d\in \TZ*} s_{d} \,\phi_{d} \bigg). 
\end{equation}\fbeq
\end{proposition}

\begin{textesansboite}
In \cref{eq: equations de bfuDA}, the equalities ``in'' or ``on'' a subdomain should be understood as the equality of the restrictions. For instance, ``$\mu _{1} \partial _{Y} ^2  \bfuDA = g^{0}$ in $\couche$'' means that we have $\mu _{1} \partial _{Y} ^2  \bfuDA_{|\couche} = g^{0}$ in $\Amacro(\couche)$. Besides, the right-hand side of \cref{eq: bfuDA} is well-defined thanks to \cref{existence des serie geometrique d'op 1D}.
\end{textesansboite}

\begin{demo}
Let $\bfv^{0} := ( \mathrm{id}+k_{0} ^2  \ROmegaPi)(\bfuDA)$. By~\cref{inverse de id-f dans les SF 1D}, $\langle -k_{0} ^2  \ROmegaPi \rangle $ maps $\Amacro(\Pi )$ to itself, and it~is the inverse of $ \mathrm{id}+k_{0} ^2  \ROmegaPi$ on $\Amacro(\Pi )$, so we have $\bfuDA = \langle -k_{0} ^2  \ROmegaPi \rangle (\bfv^{0} )$. In addition, \cref{operateurs R} states that we have $\Delta \circ \ROmegaPi= \mathrm{id}$, $\partial _{Y|\couche} ^2  \circ \ROmegaPi =0$, $\partial _{Y|\Gamma ,Y=0^{-}} \circ \ROmegaPi=0$ and $\sigma _{d} \circ \ROmegaPi =0$ for~any $d\in \TZ*$, as operators on $\mathcal{A}(\Pi )$. These relations then generalize to the extensions of these operators to $\Amacro(\Pi )$. Hence we deduce:
\begin{itemize}
	\item $(\mu _{0} \Delta  +\omega ^2 \rho _{0} )\circ \langle -k_{0} ^2  \ROmegaPi \rangle  = \mu _{0} \Delta  \circ ( \mathrm{id} +k_{0} ^2  \ROmegaPi)\circ \langle -k_{0} ^2  \ROmegaPi \rangle =\mu _{0} \Delta $,
	\item $\partial _{Y|\couche} ^2  \circ \langle -k_{0} ^2  \ROmegaPi \rangle  = \partial _{Y|\couche} ^2 $,
	\item $\partial _{Y|\Gamma ,Y=0^{-}} \circ \langle -k_{0} ^2  \ROmegaPi \rangle  = \partial _{Y|\Gamma ,Y=0^{-}}$,
	\item and $\sigma _{d} \circ \langle -k_{0} ^2  \ROmegaPi \rangle  = \sigma _{d}$ for~any $d\in \TZ*$.
\end{itemize}
Therefore, \cref{eq: equations de bfuDA} is equivalent to
\begin{equation}{\label{eq: pb de Helmholtz dans A Pi preuve 1}}
\left\{\begin{array}{r@{\;=\;}ll}
	\mu _{0} \Delta  \bfv^{0} & 0&\text{ in }\Omega  \\[0.5mm]
	\mu _{1} \partial _{Y} ^2  \bfv^{0} & g^{0} &\text{ in }\couche \\[0.5mm]
	\mu _{1} \partial _{Y} \bfv^{0} _{|Y=0^{-}} & h^{0} &\text{ on }\Gamma  \\[0.5mm]
	\sigma _{d} (\bfv^{0} ) & s_{d} & \ \forall d\in \TZ*
\end{array}\right.\end{equation}
Let us denote $(g^{0} _{d} )_{d\in \mathbb{R}} :=g^{0}$, $(h^{0} _{d} )_{d} :=h^{0}$ and $(\bfv^{0} _{d} )_{d} :=\bfv^{0}$ the coordinates of the formal series $g^{0}$, $h^{0}$ and $\bfv^{0}$. Since $\Delta _{|\Omega }$ has a degree, and $\partial _{Y|\couche} ^2 $ and $\partial _{Y|Y=0^{-}}$ have degree 0, \cref{eq: pb de Helmholtz dans A Pi preuve 1} is equivalent to:
\begin{equation}{\label{eq: pb de Helmholtz dans A Pi preuve 2}}
\forall d\in \mathbb{R},\qquad  \left\{\begin{array}{r@{\;=\;}ll}
	\mu _{0} \Delta  \bfv^{0} _{d} & 0&\text{ in }\Omega  \\[0.5mm]
	\mu _{1} \partial _{Y} ^2  \bfv^{0} _{d} & g^{0} _{d} &\text{ in }\couche \\[0.5mm]
	\mu _{1} \partial _{Y} \bfv^{0} _{d|Y=0^{-}} & h^{0} _{d} &\text{ on }\Gamma  \\[0.5mm]
	\sigma _{d} (\bfv^{0} _{d} ) & s_{d} & \ \text{if } d\in \TZ*
\end{array}\right.\end{equation}

Then \cref{pb de Laplace homogene dans A,operateurs R} show that, for~any $d\in \mathbb{R}$, the set of solutions of \cref{eq: pb de Helmholtz dans A Pi preuve 2} in $\mathcal{A}_{d} (\Pi )$ is
\begin{itemize}
	\item $\big\{\varphi = \frac{1}{\mu _{1}} \RcouchePi(g^{0} ) + \frac{1}{\mu _{1}} \RGammaPi(h^{0} ) + c\. \phi_{d} \ \big|\ c\in \mathbb{C},\, \sigma _{d} (\varphi )=s_{d} \big\}$ if $d\in \TZ*$,
	\item and $\big\{ \frac{1}{\mu _{1}} \RcouchePi(g^{0} ) + \frac{1}{\mu _{1}} \RGammaPi(h^{0} )\big\}$ otherwise.
\end{itemize}
Moreover, in the first case, we have $\sigma _{d} \circ \ROmegaPi =\sigma _{d} \circ \RcouchePi =\sigma _{d} \circ \RGammaPi=0$ by \cref{operateurs R}, so $c=s_{d}$. Therefore, \cref{eq: pb de Helmholtz dans A Pi preuve 2} is equivalent to $\bfv^{0} = \frac{1}{\mu _{1}} \RcouchePi(g^{0} ) + \frac{1}{\mu _{1}} \RGammaPi(h^{0} ) + \fsum_{d\in \TZ*} s_{d} \,\phi_{d}$. Thanks to the assumptions on $g^{0}$, $h^{0}$ and $(s_{d} )$, this $\bfv^{0}$ indeed belongs to $\Amacro(\Pi )$. This concludes the proof.
\end{demo}

\begin{proposition}{resolution of a Helmholtz problem in $\Amacro(\Omega )$}
\label{pb de Helmholtz dans A Omega}%
Let $g^{0} \in \Amacro(\Gamma )$, and $(s_{d} )_{d\in \TZ*} \in \mathbb{C}^{\TZ*}$ a family whose support is bounded from below. The following problem has a unique solution in $\Amacro(\Omega )$:
\begin{equation}{\label{eq: equations de uDA}}
	\left\{\begin{array}{r@{\;=\;}ll}
	\mu _{0} \Delta  \uDA +\omega ^2 \rho _{0} \uDA & 0 &\text{ in }\Omega \\[0.5mm]
	\uDA & g^{0} &\text{ on }\Gamma  \\[0.5mm]
	\sigma _{d} (\uDA) & s_{d} & \ \forall d\in \TZ*
	\end{array}\right.
\end{equation}
Moreover, the solution is equal to
\begin{equation}{\label{eq: uDA}}
	\uDA = \big\langle  {-}k_{0} ^2 \,\ROmegaOmega \big\rangle  \bigg( \RGammaOmega(g^{0} ) + \fsum_{d\in \TZ*} s_{d} \,\phi_{d} \bigg).
\end{equation}\fbeq
\end{proposition}

\begin{textesansboite}
The proof of \cref{pb de Helmholtz dans A Omega} is very similar to \cref{pb de Helmholtz dans A Pi}, so we do~not write it.
\end{textesansboite}

\begin{proposition}{resolution of a Poisson problem in $\Amicro(\Pi )$}
\label{pb de Poisson dans A Pi}%
Let $F^{\infty } _{\Omega } \in \Amicro(\Omega )$, $F^{\infty } _{\couche} \in \Amicro(\couche)$ and $(s_{d} )_{d\in \TZ*} \in \mathbb{C}^{\TZ*}$ a family whose support is bounded from \emph{above}. The following problem has a unique solution in $\Amicro(\Pi )$:
\begin{equation}{\label{eq: equations de SDA}}
	\left\{\begin{array}{r@{\;=\;}ll}
	\mu _{0} \Delta  \SDA & F^{\infty } _{\Omega } &\text{ in }\Omega  \\[0.5mm]
	\mu _{1} \Delta  \SDA & F^{\infty } _{\couche} &\text{ in }\couche \\[0.5mm]
	\mu _{0} \partial _{Y|Y=0^{+}} \SDA - \mu _{1} \partial _{Y|Y=0^{-}} \SDA & 0 &\text{ on }\Gamma \\[0.5mm]
	\sigma _{d} (\SDA) & s_{d} & \ \forall d\in \TZ*
	\end{array}\right.
\end{equation}
Moreover, the solution is equal to
\begin{equation}{\label{eq: SDA}}
	\SDA = \left\langle {-\RcouchePi\circ \partial _{X|\couche} ^2 ,\; \frac{\mu _{0} }{\mu _{1}} \RGammaPi\circ \partial _{Y|\Gamma ,Y=0^{+}} }\right\rangle  \bigg( \frac{1}{\mu _{0}} \ROmegaPi (F^{\infty } _{\Omega } )+ \frac{1}{\mu _{1}} \RcouchePi (F^{\infty } _{\couche} ) + \fsum_{d\in \TZ*} s_{d} \,\phi_{d} \bigg) . 
\end{equation}\fbeq
\end{proposition}

\begin{remarque}
$F^{\infty } _{\Omega }$ and $F^{\infty } _{\couche}$ cannot simply be represented by an element of $\mathcal{A}(\Pi )$ because their concatenation (as a function $\Pi \rightarrow \mathbb{C}$) may~not be continuous on $\Gamma $.
\end{remarque}

\begin{demo}
This is similar to \cref{pb de Helmholtz dans A Pi}. Let $R_{1} :=-\RcouchePi\circ \partial _{X|\couche} ^2 $, $R_{2} :=\frac{\mu _{0} }{\mu _{1}} \RGammaPi\circ \partial _{Y|\Gamma ,Y=0^{+}} $ and $\bfv^{\infty } := ( \mathrm{id}-R_{1} -R_{2} )(\SDA)$. By~\cref{inverse de id-f dans les SF 1D}, $\langle R_{1} ,R_{2} \rangle $ send $\Amicro(\Pi )$ to itself and is the inverse of $ \mathrm{id}-R_{1} -R_{2}$ on $\Amicro(\Pi )$, so we have $\SDA = \langle R_{1} ,R_{2} \rangle (\bfv^{\infty } )$. Moreover, \cref{operateurs R} implies that we have:
\begin{itemize}
	\item $\Delta \circ R_{1} =\Delta \circ R_{2} =0$,
	\item $\partial _{Y|\couche} ^2  \circ R_{1} =-\partial _{X|\couche} ^2 $ and $\partial _{Y|\couche} ^2  \circ R_{2} =0$,
	\item $\partial _{Y|\Gamma ,Y=0^{-}} \circ R_{1} =0$ and $\partial _{Y|\Gamma ,Y=0^{-}} \circ R_{2} =\frac{\mu _{0} }{\mu _{1}} \partial _{Y|\Gamma ,Y=0^{+}}$,
	\item and $\sigma _{d} \circ R_{1} =\sigma _{d} \circ R_{2} =0$ for~any $d\in \TZ*$.
\end{itemize}
These relations first hold for the original operators on the spaces $\mathcal{A}$, and are then generalized to their extension to $\Amicro$. This implies:
\begin{itemize}
	\item $\Delta \circ \langle R_{1} ,R_{2} \rangle =\Delta $,
	\item $(\partial _{Y|\couche} ^2  + \partial _{X|\couche} ^2 )\circ \langle R_{1} ,R_{2} \rangle  = \partial _{Y|\couche} ^2  \circ ( \mathrm{id}-R_{1} )\circ \langle R_{1} ,R_{2} \rangle  = \partial _{Y|\couche} ^2  \circ ( \mathrm{id}-R_{1} -R_{2} )\circ \langle R_{1} ,R_{2} \rangle  =\partial _{Y|\couche} ^2 $,
	\item similarly $(\mu _{0} \partial _{Y|\Gamma ,Y=0^{+}} - \mu _{1} \partial _{Y|\Gamma ,Y=0^{-}} )\circ \langle R_{1} ,R_{2} \rangle =\mu _{1} \partial _{Y|\Gamma ,Y=0^{-}}$.
	\item and $\sigma _{d} \circ \langle R_{1} ,R_{2} \rangle =\sigma _{d}$ for~any $d\in \TZ*$.
\end{itemize}
Therefore, \cref{eq: equations de SDA} is equivalent to
\begin{equation}{\label{eq: pb de Poisson dans A Pi preuve}}
\left\{\begin{array}{r@{\;=\;}ll}
	\mu _{0} \Delta  \bfv^{\infty } & F^{\infty } _{\Omega } &\text{ in }\Omega  \\[0.5mm]
	\mu _{1} \partial _{Y} ^2  \bfv^{\infty } & F^{\infty } _{\couche} &\text{ in }\couche \\[0.5mm]
	\mu _{1} \partial _{Y} \bfv^{\infty } _{|Y=0^{-}} & 0 &\text{ on }\Gamma  \\[0.5mm]
	\sigma _{d} (\bfv^{\infty } ) & s_{d} & \ \forall d\in \TZ*
\end{array}\right.\end{equation}
Similarly to the proof of \cref{pb de Helmholtz dans A Pi}, the only solution of \cref{eq: pb de Poisson dans A Pi preuve} in $\Amicro(\Pi )$ is $ \frac{1}{\mu _{0}} \ROmegaPi (F^{\infty } _{\Omega } )+ \frac{1}{\mu _{1}} \RcouchePi (F^{\infty } _{\couche} ) + \fsum_{d\in \TZ*} \! s_{d} \.\phi_{d}$, which ends the proof.
\end{demo}

\refstepcounter{ansatz}
\begin{paragraphesansboite}{Ansatz \theansatz:}
\label{ansatz: existence de bfuDApl et SDApl}%
In the rest of this chapter, we consider a family $(\bfuDA_{p,\ell } )_{p\in \NN,\ell \in \mathbb{N}}$ of elements of $\Amacro(\Pi )$ and a family $(\SDA_{p,\ell } )_{p\in \NN,\ell \in \mathbb{N}}$ of elements of $\Amicro(\Pi )$ s.t., for any $p\in \NN$ and $\ell \in \mathbb{N}$, we have:
\begin{equation}{\label{eq: eq vol bfuDApl}}
	\left\{\begin{array}{r@{\;=\;}ll}
	\mu _{0} \Delta  \bfuDA_{p,\ell } +\omega ^2 \rho _{0} \bfuDA_{p,\ell } & 0 &\text{ in }\Omega  \\[0.5mm]
	\mu _{1} \partial _{Y} ^2  \bfuDA_{p,\ell } & -(\mu _{1} \partial _{x} ^2 +\omega ^2 \rho _{1} )\bfuDA_{p-2,\ell } &\text{ in }\couche \\[0.5mm]
	\mu _{1} \partial _{Y} \bfuDA_{p,\ell |Y=0^{-}} & \mu _{0} \.\partial _{y} \bfuDA_{p-1,\ell |y=0^{+}}  &\text{ on }\Gamma  
	\end{array}\right.
\end{equation}
and
\begin{equation}{\label{eq: eq vol SDApl}}
	\left\{\begin{array}{r@{\;=\;}ll}
	\mu _{0} \Delta  \SDA_{p,\ell } & -\omega ^2 \rho _{0} \.\SDA_{p-2,\ell } &\text{ in }\Omega  \\[0.5mm]
	\mu _{1} \Delta  \SDA_{p,\ell } & -\omega ^2 \rho _{1} \.\SDA_{p-2,\ell } &\text{ in }\couche \\[0.5mm]
	\mu _{0} \.\partial _{Y|Y=0^{+}} \SDA_{p,\ell } - \mu _{1} \partial _{Y|Y=0^{-}} \SDA_{p,\ell } & 0 &\text{ on }\Gamma 
	\end{array}\right.
\end{equation}
By~convention, we are taking $\bfuDA_{p,\ell } =0$ and $\SDA_{p,\ell } =0$ for~any $p\in \mathbb{R} \setminus \NN$ and $\ell \in \mathbb{N}$.
\end{paragraphesansboite}
\nomenclature[30b-00]{$\bfuDA_{p,\ell }$}{Formal series representing the asymptotic behavior of $\bfu_{p,\ell }$ at the corner}
\nomenclature[28-60]{$\SDA_{p,\ell }$}{Formal series representing the asymptotic behavior of $S_{p,\ell }$ at infinity}

\begin{textesansboite}
Note that \cref{eq: eq vol bfuDApl,eq: eq vol SDApl} have the same structure as the equations satisfied by $\bfu_{p,\ell }$ and $S_{p,\ell }$, which are resp.\ \cref{eq: eq vol de bfu p l,eq: eq vol champs de coin}. We will prove in \cref{validation construction et DA des champs} the existence of such families $(\bfuDA_{p,\ell } )_{p,\ell }$ and $(\SDA_{p,\ell } )_{p,\ell }$ that satisfy also
\[\forall d\in \mathbb{R}, \quad  \left\{\begin{array}{r@{\;=\;}ll} 
	\bfu_{p,\ell } & T_{\leqslant d} (\bfuDA_{p,\ell } ) + o(\bfr^{d} ) & \mathrm{when} \ \bfr\rightarrow 0\\[3pt]
	S_{p,\ell } & T_{\geqslant d} (\SDA_{p,\ell } ) + o(\bfr^{d} )& \mathrm{when} \ \bfr\rightarrow \infty 
\end{array} \right.\]
and similar relations for the derivatives of the fields. However, in this chapter we assume that these formal series exist and we will establish matching conditions that guarantee the satisfaction of the matching assumption \cref{eq: hypothese de raccord} (which links the behavior at the corner of the far-and-layer fields to the behavior at infinity of the corner fields).\\

By~\cref{pb de Helmholtz dans A Pi,pb de Poisson dans A Pi}, we have for any $p\in \NN$ and $\ell \in \mathbb{N}$:
\begin{align}
	\label{eq: bfuDApl}
	\!\bfuDA_{p,\ell } &= \big\langle {-k_{0} ^2 \.\ROmegaPi}\big\rangle  \bigg( {-\RcouchePi\circ (\partial _{x|\couche} ^2 +k_{1} ^2 )(\bfuDA_{p-2,\ell } )} + \frac{\mu _{0} }{\mu _{1}} \RGammaPi\circ \partial _{y|\Gamma ,y=0^{+}} (\bfuDA_{p-1,\ell } ) + \!\fsum_{d\in \TZ*} \! \sigma _{d} (\bfuDA_{p,\ell } ) \,\phi_{d} \bigg) \!\\
	\label{eq: SDApl}
	\!\SDA_{p,\ell } &= \left\langle \! -\RcouchePi\circ \partial _{X|\couche} ^2 , \frac{\mu _{0} }{\mu _{1}} \RGammaPi\circ \partial _{Y|\Gamma ,Y=0^{+}} \!\right\rangle  \!\bigg( \!{- k_{0} ^2  \ROmegaPi (\SDA_{p-2,\ell |\Omega } )}- k_{1} ^2  \RcouchePi (\SDA_{p-2,\ell |\couche} ) + \!\!\fsum_{d\in \TZ*} \! \sigma _{d} (\SDA_{p,\ell } ) \,\phi_{d} \bigg) 
\end{align}
An important remark is that $(\bfuDA_{p,\ell } )_{p,\ell }$ and $(\SDA_{p,\ell } )_{p,\ell }$ are uniquely determined by $(\sigma _{d} (\bfuDA_{p,\ell } ))_{d,p,\ell }$ and $(\sigma _{d} (\SDA_{p,\ell } ))_{d,p,\ell }$. Thus, it suffices that the coefficients $\sigma _{d} (\bfuDA_{p,\ell } )$ and $\sigma _{d} (\SDA_{p,\ell } )$ have appropriate values in order to satisfy the matching assumption \cref{eq: hypothese de raccord}.\\

When $d<0$ (resp.\ $d>0$), $\phi_{d}$ is non-variational for the far-and-layer fields (resp.\ corner fields), and \cref{cadre fonctionnel des champs lointains} (resp.~\ref{cadre fonctionnel pour les champs de coin}) will show that $\sigma _{d} (\bfuDA_{p,\ell } )$ (resp.\ $\sigma _{d} (\SDA_{p,\ell } )$) can be fixed arbitrarily. \Cref{sec: raccordement algebrique} is devoted to finding how to fix them in order to satisfy the matching assumption~\cref{eq: hypothese de raccord}. On the contrary, when $d>0$ (resp.\ $d<0$), $\phi_{d}$ is variational for the far-and-layer fields (resp.\ corner fields) and the values of $\sigma _{d} (\bfuDA_{p,\ell } )$ (resp.\ $\sigma _{d} (\SDA_{p,\ell } )$) are uniquely defined once $(\sigma _{d'} (\bfuDA_{p,\ell } ))_{d'<0}$ (resp.\ $(\sigma _{d'} (\SDA_{p,\ell } ))_{d'>0}$) has been fixed. Ways to compute these values numerically are investigated in \cref{chap: validation numerique du DA}.
\end{textesansboite}

\subsection{Formal series of corner singularities involving $\feps$}
\label{sec: SF 2D de singularites de coin}

\begin{textesansboite}
In this section we give tools that will be used in \cref{sec: raccordement algebrique} to study the formal series $\fsum_{p,\ell } \powelne \cdot \bfuDA_{p,\ell }$ and $\fsum_{p,\ell } \varepsilon ^{p} \ln ^{\ell } \! \varepsilon \cdot \SDA_{p,\ell }$ and to deduce the matching conditions. These two formal series are rigorously defined objects that give a meaning to the ill-defined series ``$\sum _{p,\ell } \peps^{p} \ln ^{\ell } \!\peps \cdot \bfu_{p,\ell }$'' and ``$\sum _{p,\ell } \peps^{p} \ln ^{\ell } \! \peps\cdot S_{p,\ell }$'' that were used in \cref{eq: hypothese de raccord} to intuitively introduce the matching assumption.\\

\nomenclature[40-21-23]{$\ln \feps$}{Formal indeterminate used to describe logarithmic behaviors w.r.t.\ $\ln \peps$ with formal series}
Here, the symbols ``$\varepsilon $'' and ``$\ln \varepsilon $'' denote two algebraic indeterminates independent of each other (so they are~not numbers, contrary to $\peps$ and $\ln \peps$).
More precisely, for a given $D\in \domsA$, we consider the vector spaces $E:=\mathcal{A}(D)[\ln \varepsilon ]$ and $E_{d} :=\mathcal{A}_{d} (D)[\ln \varepsilon ]$ for~any $d\in \mathbb{R}$, which are defined by \cref{def: E[T]} (where we denoted the indeterminate ``$\ln \varepsilon $'' instead of ``$T$''). We recall that these spaces are the sets of formal polynomials with coefficients in $\mathcal{A}(D)$, resp.\ $\mathcal{A}_{d} (D)$. Then, we apply the tools of \cref{sec: SF 2D} with these spaces $E$ and $E_{d}$, which defines the formal series of the form 
\[\fsumd \.\fsump \varepsilon ^{p} \varphi _{d,p} \qquad  \text{with }  \varphi _{d,p} \in \mathcal{A}_{d} (D)[\ln \varepsilon ] \text{ for~any }  (d,p)\in \mathbb{R}^2 .\]
\Cref{sec: feps dans les SF 2D} shows that the powers of $\varepsilon $ behave like powers of a formal indeterminate would do, so $\varepsilon $ can be manipulated like a formal indeterminate. Since each $\varphi _{d,p}$ has the form $\sum _{\ell \in \mathbb{N}} \ln ^{\ell } \!\varepsilon \,\varphi _{d,p,\ell }$ (where this sum has only a finite number of non-zero terms), this kind of formal series can also be written as
\[\fsumd \.\fsump \.\sum _{\ell \in \mathbb{N}} \powelne \cdot  \varphi _{d,p,\ell } \quad  \text{ or } \quad  \fsumd \.\fsump \.\fsum_{\ell \in \mathbb{N}} \powelne \cdot  \varphi _{d,p,\ell } \quad  \text{ e.g.}  .\]
Note that it~is possible to use the work of \cref{sec: SF 2D} on these spaces because the assumption ``$E=\bigoplus_{d\in \mathbb{R}} E_{d}$'' is satisfied thanks to the property $\mathcal{A}(D)=\bigoplus_{d\in \mathbb{R}} \mathcal{A}_{d} (D)$ established in \cref{decomposition de A en sev}.\\

The next definition introduces the spaces of formal series that are relevant to the different fields. We recall that the $\fsum$-sum of vector spaces in defined by \cref{def: fsum d'ensembles} and is studied in \cref{fsum des Ed dans les SF 2D}.
\end{textesansboite}

\begin{definition}{spaces $\Aeps$ and $\Aemaicro$}
\label{def: Aeps et Aemaicro}%
Let $D\in \{\Omega ,\Gamma ,\layer,\Pi \}$. We define the spaces $\Aeps(D)$, $\Aemacro(D)$ and $\Aemicro(D)$ by:
\[\Aeps(D) := \fsum_{(d,p)\in \mathbb{R}^2 } \varepsilon ^{p} \mathcal{A}_{d} (D)[\ln \varepsilon ] \qquad \! \text{ and } \!\qquad 
\Aemaicro(D) := \bigcup _{(d_{0} ,p_{0} )\in \mathbb{R}^2 } \.\fsum_{(d,p) \mkern1mu \in  \mkern1mu \conemaicro+(d_{0} ,p_{0} )} \!\varepsilon ^{p} \mathcal{A}_{d} (D)[\ln \varepsilon ] ,\]
where $\conemaicro :=\{(d,p)\in \mathbb{R}^2  \mid p\geqslant 0 \text{ and } p\pm d\geqslant 0\}$ (see \cref{fig: conemaicro}).
\end{definition}
\nomenclature[10-60-40]{$\Aeps$, $\Aemaicro$}{Spaces of formal series of corner singularities with powers of $\feps$ and $\ln \feps$}
\nomenclature[10-60-40]{$\conemaicro$}{Cone that contains the degrees of the operators in $\Rmaicro$}

	\tikzset{external/export=true}%
	\tikzsetfigurename{fig-conemaicro-}%
	\begin{figure}[h]
\def\ymin{-1.8}
\def\ymax{1.9}
\def\xmax{2.5}
\def\pzero{-0.85}
\def\dzero{-0.6}
\def\traceraxes{
	\draw [->] (-\xmax-0.3,0) -- (\xmax+0.3,0) node [right] {$d$};
	\draw [->] (0,\ymin-0.1) -- (0,\ymax+0.3) node [right] {$p$};}
\begin{tikzpicture}[scale=0.95]
	\begin{scope}
		\clip (-\xmax,\ymin) rectangle (\xmax,\ymax); 
		\fill [green!70!black, opacity=0.2] (\dzero,\pzero) -- (10,\pzero) -- (10,10) -- (-10,10) -- (\dzero-10,\pzero+10) -- cycle;
	\end{scope}
	\traceraxes
	\begin{scope}
		\clip (-\xmax,\ymin) rectangle (\xmax,\ymax); 
		\draw [very thick, green!55!black] (\dzero,\pzero) -- (10,\pzero);
		\draw [very thick, green!55!black] (\dzero,\pzero) -- (\dzero-10,\pzero+10);
	\end{scope}
	\draw [very thick, green!30!black] (\xmax*0.5,\ymax*0.4) node {$\conemacro+(d_0,p_0)$};
	\draw (\dzero,\pzero) node{$\bullet$};
	\draw [<-,>=latex] (\dzero-0.02,\pzero-0.03) to[bend right=20] ++(-0.2,-0.4) node[below]{$(d_0,p_0)$};
\end{tikzpicture}%
\hspace{10mm}\nolinebreak%
\begin{tikzpicture}[scale=0.95]
	\begin{scope}
		\clip (-\xmax,\ymin) rectangle (\xmax,\ymax); 
		\fill [green!70!black, opacity=0.2] (\dzero,\pzero) -- (-10,\pzero) -- (-10,10) -- (10,10) -- (\dzero+10,\pzero+10) -- cycle;
	\end{scope}
	\traceraxes
	\begin{scope}
		\clip (-\xmax,\ymin) rectangle (\xmax,\ymax); 
		\draw [very thick, green!55!black] (\dzero,\pzero) -- (-10,\pzero);
		\draw [very thick, green!55!black] (\dzero,\pzero) -- (\dzero+10,\pzero+10);
	\end{scope}
	\draw [very thick, green!30!black] (-\xmax*0.5,\ymax*0.4) node {$\conemicro+(d_0,p_0)$};
	\draw (\dzero,\pzero) node{$\bullet$};
	\draw [<-,>=latex] (\dzero+0.01,\pzero-0.03) to[bend left=20] ++(-0.05,-0.4) node[below]{$(d_0,p_0)$};
\end{tikzpicture}

\caption{The sets $\conemaicro+(d_0,p_0)$ for a given $(d_0,p_0)\in\mathbb{R}^2$}
\label{fig: conemaicro}
\end{figure}%
	\tikzset{external/export=false}%
	\tikzsetfigurename{main-figure}

\begin{textesansboite}
The multiplication by powers of $\varepsilon $ in $\Aeps(D)$ is defined in \cref{def: feps}. As~for powers of $\ln \varepsilon $, for~any $\ell \in \mathbb{N}$, we define the map $\Phi \in \Aeps(D) \mapsto \ln ^{\ell } \!\varepsilon  \cdot \Phi  \in \Aeps(D)$ as the extension (in the sense of \cref{def: op à degree dans les SF 2D}) of the map $\varphi \in \mathcal{A}(D)[\ln \varepsilon ] \mapsto  \ln ^{\ell } \!\varepsilon  \cdot \varphi \in \mathcal{A}(D)[\ln \varepsilon ]$ which has degree 0 w.r.t.\ the decomposition $\mathcal{A}(D)[\ln \varepsilon ]=\bigoplus_{d\in \mathbb{R}} \mathcal{A}_{d} (D)[\ln \varepsilon ]$. By~\cref{les op a degre commutent avec fsum et eps dans les SF 2D}, the multiplication by $\ln ^{\ell } \!\varepsilon $ in $\Aeps(D)$ commutes with $\fsum$ and powers of $\varepsilon $.
\\

$\Aeps(D)$ contains all the formal series of corner singularities that involves powers of $\varepsilon $ and $\ln \varepsilon $, while the sub-spaces $\Aemacro(D)$ are relevant for the far-and-layer fields and the sub-spaces $\Aemicro(D)$ are relevant for the corner fields. We will actually have
\[\fsum_{p\in \NN} \fsum_{\ell \in \mathbb{N}} \powelne \cdot \bfuDA_{p,\ell } \in \Aemacro(\Pi ) \qquad  \text{ and } \qquad  \fsum_{p\in \NN} \fsum_{\ell \in \mathbb{N}} \varepsilon ^{p} \ln ^{\ell } \! \varepsilon \cdot \SDA_{p,\ell } \in \Aemicro(\Pi ).\]
See \cref{ansatz sur Aeps}. 
The exact definition of $\Aemaicro$ is chosen so that these spaces are the natural setting to solve Helmholtz problems that will be important to study $\fsum_{p,\ell } \powelne\,\bfuDA_{p,\ell }$ and $\fsum_{p,\ell } \smash{\powelne} \,\SDA_{p,\ell }$ (see \cref{pb de Helmholtz dans Aemacro,pb de Helmholtz dans Aemicro}).
\end{textesansboite}

\begin{remarque}
We have the following diagram of inclusions for~any $D\in \{\Omega ,\Gamma ,\layer,\Pi \}$ and $d\in \mathbb{R}$:
\[\xymatrixrowsep{0.5ex}\xymatrixcolsep{5ex}\xymatrix{
	& & \Amacro(D) \ar[r] & \Aemacro(D) \ar[rd] & \\
	\mathcal{A}_{d} (D) \ar[r] & \mathcal{A}(D) \ar[ru] \ar[rd] & & & \Aeps(D) \\
	& & \Amicro(D) \ar[r] & \Aemicro(D) \ar[ru] & \
}\]
Conversely, the components of any element of $\Aemacro(\Pi )$, resp.\ $\Aemicro(D)$, on a given power of $\varepsilon $ and $\ln \varepsilon $ belong to $\Amacro(\Pi )$, resp.\ $\Amicro(D)$. We mean by this that, for any $\fsum_{d,p,\ell } \powelne\,\varphi _{d,p,\ell } \in \Aemaicro(D)$ with $\varphi _{d,p,\ell } \in \mathcal{A}_{d} (D)$ for any $(d,p,\ell )$, we have $\fsum_{d} \varphi _{d,p,\ell } \in \Amaicro(D)$ for any $(p,\ell )$.
\end{remarque}

\begin{exemple}
$\fsum_{p\in \NN} \fsum_{d\in \NN-p} \sum _{\ell =0} ^{\lfloor {p}\rfloor } \powelne \,\Im[(\alpha z)^{d} ]$ is well-defined and belongs to $\Aemacro(\Pi )$ because its support is locally finite and included in $\conemacro$, and the powers of $\ln \varepsilon $ are bounded for each $(d,p)$. Similarly, $\fsum_{p\in \NN} \fsum_{d\in \NN-p} \sum _{\ell =0} ^{\lfloor {p\red{+|d|}}\rfloor } \powelne \,\Im[(\alpha z)^{d} ]$ belongs to $\Aemacro(\Pi )$, although in practice we will only consider series whose powers of $\ln \varepsilon $ are bounded at each power of $\varepsilon $. However ``$\fsum_{p\in \NN} \fsum_{d\in \NN-p} \fsum_{\ell =0} ^{\red{\infty }} \powelne \,\Im[(\alpha z)^{d} ]$'' is not well-defined.
\end{exemple}

\begin{textesansboite}
We generalize the linear forms $\sigma _{d}$ to $\Aeps$ by defining the linear forms $\sigma _{d,p,\ell }$, which evaluate $\sigma _{d}$ on the component of a formal series that is multiplied by $\varepsilon ^{p} \ln ^{\ell } \varepsilon $.
\end{textesansboite}

\begin{definition}{linear form $\sigma _{d,p,\ell }$}
\label{def: sigma sur Aeps}%
Let $D\in \{\Omega ,\Pi \}$, $d\in \TZ*$, $p\in \mathbb{R}$ and $\ell \in \mathbb{N}$. We define $\sigma _{d,p,\ell } :\Aeps(D) \rightarrow \mathbb{C}$ as follows.\\
Let $\Phi  \in \Aeps(D)$. There exists a unique family $(\varphi _{d',p',\ell '} )_{(d',p')\in \mathbb{R}^2 ,\ell \in \mathbb{N}}$ s.t.\ $\Phi  =\fsum_{(d',p')\in \mathbb{R}^2 } \sum _{\ell '\in \mathbb{N}} \varepsilon ^{p'} \ln ^{\ell '\!} \!\varepsilon \cdot  \varphi _{d',p',\ell '}$ and $\varphi _{d',p',\ell '} \in \mathcal{A}_{d'} (D)$ for any $(d',p',\ell ')$. We then define
\[\sigma _{d,p,\ell } (\Phi  ) = \sigma _{d,p,\ell } \bigg(\fsum_{d'\in \mathbb{R}} \fsum_{p'\in \mathbb{R}} \sum _{\ell '\in \mathbb{N}} \varepsilon ^{p'} \ln ^{\ell '\!} \!\varepsilon \cdot  \varphi _{d',p',\ell '} \bigg):=\sigma _{d} (\varphi _{d,p,\ell } ).\]\fbeq
\end{definition}
\nomenclature[53-13z-25-21]{$\sigma _{d,p,\ell }$}{Image by $\sigma _{d}$ of some term of a formal series}

\begin{exemple}
If $\Phi  \in \fsum_{q\in \mathbb{R}} \mathcal{A}_{q} (\Omega )$, then $\sigma _{d,p,\ell } (\powelne\,\Phi )= \sigma _{d} (\Phi )$ (see \cref{def: sigma sur Amaicro}) and $\sigma _{d,p',\ell '} (\powelne\,\Phi )= 0$ for any $(p',\ell ') \neq (p,\ell )$.
\end{exemple}

\begin{lemmesansboite}{}
\label{def: sigma et fsum commutent dans Aeps}%
Let $D\in \{\Omega ,\Pi \}$, $d\in \TZ*$, $p\in \mathbb{R}$ and $\ell \in \mathbb{N}$. Let $I$ be a set and $(\Phi _{i} )_{i\in I} \in \Aeps(D)^{I}$ be a family of formal series that has collectively locally finite support. We have
\[\sigma _{d,p,\ell } \bigg(\fsum_{i\in I} \Phi _{i} \bigg) = \sum _{i\in I} \sigma _{d,p,\ell } (\Phi _{i} )\]
where the right-hand side sum has only a finite number of non-zero terms.
\end{lemmesansboite}

\begin{demo}
See the proof of \cref{def: sigma sur Amaicro} which is very similar.
\end{demo}

\begin{textesansboite}
Let us explain how we extend to $\Aeps$ the operators that have a degree. Let $D_{1} ,D_{2} \in \{\Omega ,\Gamma ,\layer,\Pi \}$ and let $f:\mathcal{A}(D_{1} )\rightarrow \mathcal{A}(D_{2} )$ be a linear map that thas a degree $d$. We extend $f$ into a linear map $\tilde{f} :\mathcal{A}(D_{1} )[\ln \varepsilon ] \rightarrow \mathcal{A}(D_{2} )[\ln \varepsilon ]$ defined by:
\[\forall n\in \mathbb{N}, \;\forall (\varphi _{\ell } )\in \mathcal{A}(D_{1} )^{n+1} ,\qquad  \tilde{f} \bigg(\sum _{\ell =0} ^{n} \ln ^{\ell } \!\varepsilon  \, \varphi _{\ell } \bigg) :=\sum _{\ell =0} ^{n} \ln ^{\ell } \!\varepsilon  \, f(\varphi _{\ell } ).\]
From now on, this extension will also be denoted $f$. It is clear that it commutes with powers of $\ln \varepsilon $ and one can verify that it has degree $d$ w.r.t.\ the decomposition $\mathcal{A}(D_{i} )[\ln \varepsilon ] = \bigoplus_{q\in \mathbb{R}} \mathcal{A}_{q} (D_{i} )[\ln \varepsilon ]$. Then, by \cref{def: op à degree dans les SF 2D} $f$ can be again extended into an operator $\Aeps(D_{1} )\rightarrow \Aeps(D_{2} )$.\\

With this process, we extend the operators $\Delta _{|\Omega }$, $\partial _{x|\couche} ^2 $, $\partial _{y|\Gamma ,y=0^{+}}$, $\partial _{Y|\couche} ^2 $, $\partial _{Y|\Gamma ,Y=0^{-}}$, $\ROmegaOmega$, $\RGammaOmega$ $\RcouchePi$, $\RGammaPi$ and $\ROmegaPi$ to the spaces $\Aeps$ (these operators have a degree by \cref{operateurs R dans Omega,operateurs R,operateurs de derivation dans A}).
Moreover, for any operator $f:\Aeps(D_{1} )\rightarrow \Aeps(D_{2} )$ that has a degree (see \cref{def: op à degree dans les SF 2D}), we denote its degree as:
\[(\deg _{\mathcal{A}} f,\,\deg _{\varepsilon } f):=\deg  f \in \mathbb{R}^2 .\]
\nomenclature[13-14-16-10-40]{$\deg _{\mathcal{A}} ,\deg _{\feps}$}{Components of the degree of an operator acting the spaces $\Aeps$}

Next, let us introduce:
\begin{equation}{\label{eq: def Rmaicro}}
\begin{array}{c}
	\Rmacro :=\displaystyle \left\{-k_{0} ^2 \ROmegaPi,\ -\varepsilon ^2  \RcouchePi\circ \partial _{x|\couche} ^2 ,\ -\varepsilon ^2 k_{1} ^2 \RcouchePi,\ \varepsilon \. \frac{\mu _{0} }{\mu _{1}} \RGammaPi\circ \partial _{y|\Gamma ,y=0^{+}} \right\} \\[3ex]
	\Rmicro :=\displaystyle \left\{-\varepsilon ^2  k_{0} ^2 \ROmegaPi,\ -\RcouchePi\circ \partial _{X|\couche} ^2 ,\ -\varepsilon ^2  k_{1} ^2 \RcouchePi,\ \frac{\mu _{0} }{\mu _{1}} \RGammaPi\circ \partial _{Y|\Gamma ,Y=0^{+}} \right\}
\end{array} \end{equation}
$\Rmacro$, resp.\ $\Rmicro$, will be used to express the solutions of a Helmholtz problem in $\Aemacro(\Pi )$, resp.\ a rescaled Helmholtz problem in $\Aemicro(\Pi )$. As~you can see, the operators in $\Rmacro$ and $\Rmicro$ differ only on the powers of~$\varepsilon $. Their degrees are represented in \cref{fig: degres des elements de Rmacro et Rmicro}. Note that for~any $R\in \Rmaicro$, we have $\deg R \in \conemaicro$. This allows us to prove the following lemma.
\end{textesansboite}
\nomenclature[27-40-60]{$\Rmaicro$}{Particular sets of operators on formal series}

\begin{figure}[h]
\def\xmax{2.4}
\def\ymin{-0.3}
\def\ymax{2.3}
\begin{tikzpicture}[scale=1]
	\draw [thin, gray, dashed] (-\xmax,\ymin) grid (\xmax,\ymax);
	\draw [->] (-\xmax,0) -- (\xmax,0) node [right] {$\deg_{\mathcal{A}}$};
	\draw [->] (0,\ymin) -- (0,\ymax+0.2) node [right] {$\deg_\varepsilon$};
	\draw (1,0) node{$+$} node[below]{\ \ 1};
	\draw (0,1) node{$+$} node[right]{1};
	\draw [->,>=latex,very thick,blue] (0,0) --  node[midway, above] {$-k_{0} ^2 \,\ROmegaPi$} (2,0);
	\draw [->,>=latex,very thick,magenta] (-0.9,0.9) -- (-2,2) node[above] {$-\varepsilon ^2  \RcouchePi\circ \partial _{x|\couche} ^2 $};
	\draw [->,>=latex,very thick,red] (0,0) -- node[near end,right] {$-\varepsilon ^2 k_{1} ^2 \RcouchePi$} (0,2);
	\draw [->,>=latex,very thick,green!55!black] (0,0) -- node[midway,left] {$\varepsilon  \frac{\mu _{0} }{\mu _{1}} \RGammaPi\circ \partial _{y|\Gamma ,y=0^{+}}$} (-1,1);
\end{tikzpicture}%
\hspace{10mm}\nolinebreak%
\begin{tikzpicture}[scale=1]
	\draw [thin, gray, dashed] (-\xmax,\ymin) grid (\xmax,\ymax);
	\draw [->] (-\xmax,0) -- (\xmax,0) node [right] {$\deg_{\mathcal{A}}$};
	\draw [->] (0,\ymin) -- (0,\ymax+0.2) node [right] {$\deg_\varepsilon$};
	\draw (1,0) node{$+$} node[above]{\ \ 1};
	\draw (0,1) node{$+$} node[right]{1};
	\draw [->,>=latex,very thick,blue] (0,0) -- node[midway, right] {\vspace{-4pt}$-\varepsilon ^2  k_{0} ^2 \,\ROmegaPi$} (2,2);
	\draw [->,>=latex,very thick,magenta] (-0.9,0) -- (-2,0) node[above] {$-\RcouchePi\circ \partial _{x|\couche} ^2 $};
	\draw [->,>=latex,very thick,red] (0,0) -- (0,2) node[below left] {$-\varepsilon ^2  k_{1} ^2 \RcouchePi$};
	\draw [->,>=latex,very thick,green!55!black] (0,0) node[below left] {$\frac{\mu _{0} }{\mu _{1}} \RGammaPi\circ \partial _{y|\Gamma ,y=0^{+}} \!\!$} -- (-1,0);
\end{tikzpicture}\par
\caption[Degrees of the elements of $\protect\Rmacro$ and $\protect\Rmicro$]{Degrees of the elements of $\protect\Rmacro$ (on the left) and $\protect\Rmicro$ (on the right)}
\label{fig: degres des elements de Rmacro et Rmicro}
\end{figure}

\begin{lemme}{}
\label{Rmaicro defini et inversible dans Aemaicro}%
The operator $\langle \Rmaicro\rangle $ (see \cref{def: serie geometrique d'op 2D}) is well-defined on $\Aemaicro(\Pi )$. Moreover, $\langle \Rmaicro\rangle $ and $ \mathrm{id}-\sum _{R\in \Rmaicro} R$ both map $\Aemaicro(\Pi )$ into itself, are invertible on $\Aemaicro(\Pi )$ and satisfy:
\[\langle \Rmaicro\rangle  = \bigg( \mathrm{id} - \sum _{R\in \Rmaicro} R\bigg)^{-1} .\]\fbeq
\end{lemme}

\begin{demo}
For~any $d\in \mathbb{R}$, let us denote $E_{d} :=\mathcal{A}_{d} (\Pi )[\ln \varepsilon ]$, so that we have $\Aeps(\Pi ) =\fsum_{d,p\in \mathbb{R}} \varepsilon ^{p} E_{d}$. \Cref{sec: series d'op de SF 2D} considers the space
\[E^{v} := \bigcup _{M\in \mathbb{R}} \fsum_{\substack{(d,p)\in \mathbb{R}^2 \\ v\cdot (d,p)\geqslant M}} \varepsilon ^{p} E_{d}\]
where $v\cdot (d,p)$ denotes the scalar product of $\mathbb{R}^2 $. Moreover, \cref{existence des serie geometrique d'op 2D} shows that, if we have $v\cdot \deg R>0$ for~any $R\in \Rmaicro$, then $\langle \Rmaicro\rangle $ is well-defined on $E^{v}$. Thus, if we take $v:=(\pm 1,2)$ we have $v\cdot w>0$ for any $w\in \conemaicro \setminus \{0\}$, which implies $\Aemaicro(\Pi )\subset E^{v}$ and $(\forall R\in \Rmaicro,\, v\cdot \deg R>0)$. Hence $\langle \Rmaicro\rangle $ is well-defined on $\Aemaicro(\Pi )$.\\
\def\q{\mathbf{q}}%
Next, let $S:=\sum _{R\in \Rmaicro} R$ and $\Phi \in \Aemaicro(\Pi )$ and let us show that $( \mathrm{id}-S)(\Phi )$ and $\langle \Rmaicro\rangle (\Phi )$ belong to $\Aemaicro(\Pi )$. Let $\q \in \mathbb{R}^2 $ be s.t.\ $\operatorname{supp}(\Phi )\subset  \conemaicro+\q$. For~any $R\in \Rmaicro$, we have
\[\operatorname{supp}(R(\Phi )) \subset  \operatorname{supp}(\Phi ) +\deg R \subset \conemaicro+\q +\conemaicro \subset  \conemaicro+\q\]
because $\conemaicro$ is a cone, so it~is stable by addition. Hence $\operatorname{supp}(S(\Phi ))\subset \conemaicro+\q$, which implies in particular $( \mathrm{id}-S)(\Phi )\in \Aemaicro(\Pi )$. By~induction, we also deduce $\operatorname{supp}(S^{n} (\Phi ))\subset \conemaicro+\q$ for~any $n\in \mathbb{N}$. Hence, by \cref{fsum des Ed dans les SF 2D}, we obtain:
\[\langle \Rmaicro\rangle (\Phi ) \.:=\. \fsumn S^{n} (\Phi ) \.\in  \!\fsum_{(d,p)\in \conemaicro+\q} \!\varepsilon ^{p} E_{d} \,\subset \, \Aemaicro(\Pi ).\]
Finally, by \cref{inverse de id-f dans les SF 2D}, $ \mathrm{id}-S$ is invertible on $E^{v}$ and we have $( \mathrm{id}-S)^{-1} =\langle \Rmaicro\rangle $. Since $\Aemaicro(\Pi )$ is stable by these operators, we also have $( \mathrm{id}-S_{|\Aemaicro(\Pi )} )^{-1} =\langle \Rmaicro\rangle _{|\Aemaicro(\Pi )}$.
\end{demo}

\begin{textesansboite}
The next two results describe the solutions of two model problems in $\Aemaicro(\Pi )$. The first one is a Helmholtz problem that (except from singularity conditions involving $\sigma _{d,p,\ell }$) has the same form as the equations satisfied by $u_{\peps}$ when the layer is rescaled so that it has thickness 1. Hence it~is natural that $\fsum_{p,\ell } \powelne\,\bfuDA_{p,\ell }$ is a solution of this problem, which is used in \cref{sec: raccordement algebrique} to provide an expression of this series. Similarly, the second model problem has the same form as the equation satisfied by the rescaled exact solution $u_{\peps} (\peps x,\peps y)$. It is natural that $\fsum_{p,\ell } \powelne\,\SDA_{p,\ell }$ is a solution of this problem.
\end{textesansboite}

\begingroup
\def\inc{\Phi}%
\def\incBis{\Psi}%
\begin{proposition}{resolution of a Helmholtz-like problem in $\Aemacro(\Pi )$}
\label{pb de Helmholtz dans Aemacro}%
Let $(s_{d,p,\ell } )_{d,p,\ell } \in \mathbb{C}^{\TZ*\times \mathbb{R}\times \mathbb{N}}$ be a family of complex numbers s.t.\ $\fsum_{d\in \TZ*\!,\.p\in \mathbb{R},\ell \in \mathbb{N}} \powelne \cdot s_{d,p,\ell } \,\phi_{d} \in \Aemacro(\Pi )$. The following problem has a unique solution in $\Aemacro(\Pi )$:
\begin{equation}{\label{eq: pb de Helmholtz dans Aemacro}}
    \left\{\begin{array}{r@{\;=\;}ll}
        \mu _{0} \Delta  \inc +\omega ^2 \rho _{0} \inc & 0 &\text{ in }\Omega  \\[0.5mm]
        \mu _{1} \partial _{Y} ^2  \inc + \varepsilon ^2 (\mu _{1} \partial _{x} ^2 +\omega ^2 \rho _{1} )\inc &0&\text{ in }\couche \\[0.5mm]
        \mu _{1} \partial _{Y|Y=0^{-}} \inc - \varepsilon \,\mu _{0} \,\partial _{y|y=0^{+}} \inc &0 &\text{ on }\Gamma  \\[0.5mm]
		\sigma _{d,p,\ell } (\inc) & s_{d,p,\ell } &\ \forall (d,p,\ell ) \in \TZ*\times \mathbb{R}\times \mathbb{N}
    \end{array}\right.
\end{equation}
Moreover, the solution is equal to
\begin{equation}{\label{eq: solution pb de Helmholtz dans Aemacro}}
	\inc = \langle \Rmacro\rangle  \bigg( \fsum_{p\in \mathbb{R}} \fsum_{d\in \TZ*} \sum _{\ell \in \mathbb{N}} \powelne \cdot  s_{d,p,\ell } \,\phi_{d} \bigg)
\end{equation}\fbeq
\end{proposition}

\begin{proposition}{resolution of a rescaled Helmholtz-like problem in $\Aemicro(\Pi )$}
\label{pb de Helmholtz dans Aemicro}%
Let $(s_{d,p,\ell } )_{d,p,\ell } \in \mathbb{C}^{\TZ*\times \mathbb{R}\times \mathbb{N}}$ be a family of complex numbers s.t.\ $\fsum_{d\in \TZ*\!,\.p\in \mathbb{R},\ell \in \mathbb{N}} \powelne \cdot s_{d,p,\ell } \,\phi_{d} \in \Aemicro(\Pi )$. The following problem has a unique solution in $\Aemicro(\Pi )$:
\begin{equation}{\label{eq: pb de Helmholtz dans Aemicro}}
    \left\{\begin{array}{r@{\;=\;}ll}
        \mu _{0} \Delta  \inc +\varepsilon ^2 \omega ^2 \rho _{0} \inc & 0 &\text{ in }\Omega  \\[0.5mm]
        \mu _{1} \Delta  \inc +\varepsilon ^2 \omega ^2 \rho _{1} \inc & 0 &\text{ in }\couche \\[0.5mm]
        \mu _{1} \partial _{Y|Y=0^{-}} \inc - \mu _{0} \,\partial _{y|y=0^{+}} \inc &0 &\text{ on }\Gamma  \\[0.5mm]
		\sigma _{d,p,\ell } (\inc) & s_{d,p,\ell } &\ \forall (d,p,\ell ) \in \TZ*\times \mathbb{R}\times \mathbb{N}
    \end{array}\right.
\end{equation}
Moreover, the solution is equal to
\begin{equation}{\label{eq: solution pb de Helmholtz dans Aemicro}}
    \inc = \langle \Rmicro\rangle  \bigg( \fsum_{p\in \mathbb{R}} \fsum_{d\in \TZ*} \sum _{\ell \in \mathbb{N}} \powelne \cdot  s_{d,p,\ell } \,\phi_{d} \bigg)
\end{equation}\fbeq
\end{proposition}

\begin{remarque}
The assumptions ``$\fsum_{d\in \TZ*\!,\.p\in \mathbb{R},\ell \in \mathbb{N}} \powelne \cdot s_{d,p,\ell } \,\phi_{d} \in \Aemaicro(\Pi )$'' are quick ways to say that 
the support of $((s_{d,p,\ell } )_{\ell } )_{d,p}$ is locally finite in $(d,p)$ (especially in $p$ as we already have $d\in \TZ*$) and included in a translation of $\conemaicro$, and that the support of $(s_{d,p,\ell } )_{\ell }$ is finite for~any $(d,p)$.
\end{remarque}

\begin{textesansboite}
The proof method is very similar to \crefrange{pb de Helmholtz dans A Pi}{pb de Poisson dans A Pi}. We only give the proof of \cref{pb de Helmholtz dans Aemacro}.
\end{textesansboite}

\begin{demo}[Proof of \cref{pb de Helmholtz dans Aemacro}]
Let $S:=\sum _{R\in \Rmacro} R$ and $\incBis := ( \mathrm{id}-S)(\inc)$. By~\cref{Rmaicro defini et inversible dans Aemaicro}, we have $\incBis \in \Aemacro(\Pi )$ and $\inc = \langle \Rmacro\rangle (\incBis)$. Let us denote:
\[R_{1} := -k_{0} ^2 \ROmegaPi, \quad  R_{2} := -\varepsilon ^2  \RcouchePi\circ \partial _{x|\couche} ^2 ,\quad  R_{3} := -\varepsilon ^2 k_{1} ^2 \RcouchePi, \quad  R_{4} :=\varepsilon \. \frac{\mu _{0} }{\mu _{1}} \RGammaPi\circ \partial _{y|\Gamma ,y=0^{+}} .\]
Thus we have $\Rmacro =\{R_{1} ,\dots ,R_{4} \}$. \Cref{operateurs R} implies that we have:
\begin{itemize}
	\item $\Delta \circ R_{1} =-k_{0} ^2 \, \mathrm{id}$ and $\Delta \circ R_{i} =0$ for~any $i\neq 1$,
	\item $\partial _{Y|\couche} ^2 \circ R_{2} =-\varepsilon ^2  \partial _{X|\couche} ^2 $ and $\partial _{Y|\couche} ^2 \circ R_{3} =-\varepsilon ^2 k_{1} ^2 \, \mathrm{id}$ and $\partial _{Y|\couche} ^2  \circ R_{i} =0$ for~any $i\neq 2,3$,
	\item $\partial _{Y|\Gamma ,Y=0^{-}} \circ R_{4} =\varepsilon \. \frac{\mu _{0} }{\mu _{1}} \partial _{y|\Gamma ,y=0^{+}}$ and $\partial _{Y|\Gamma ,Y=0^{-}} \circ R_{i} =0$ for~any $i\neq 4$,
	\item and $\sigma _{d} \circ R_{i} =0$ for~any $i$ and $d\in \TZ*$.
\end{itemize}
These relations first hold for the original operators on the spaces $\mathcal{A}$, and are then generalized to their extension to $\Amicro$. This implies:
\begin{itemize}
	\item $(\mu _{0} \Delta +\omega ^2 \rho _{0} )\circ \langle \Rmacro\rangle =\mu _{0} \Delta \circ ( \mathrm{id}-R_{1} )\circ \langle \Rmacro\rangle  =\mu _{0} \Delta \circ ( \mathrm{id}-S)\circ \langle \Rmacro\rangle  =\mu _{0} \Delta $,
	\item $\big(\mu _{1} \partial _{Y|\couche} ^2  + \varepsilon ^2 (\mu _{1} \partial _{x|\couche} ^2 +\omega ^2 \rho _{1} )\big)\circ \langle \Rmacro\rangle  =\mu _{1} \partial _{Y|\couche} ^2  \circ ( \mathrm{id} -R_{2} -R_{3} )\circ \langle \Rmacro\rangle =\mu _{1} \partial _{Y|\couche} ^2  \circ {( \mathrm{id} -S)}\circ \langle \Rmacro\rangle =\mu _{1} \partial _{Y|\couche} ^2 $,
	\item similarly $(\mu _{1} \partial _{Y|\Gamma ,Y=0^{-}} - \varepsilon \.\mu _{0} \.\partial _{Y|\Gamma ,Y=0^{+}} )\circ \langle \Rmacro\rangle =\mu _{1} \partial _{Y|\Gamma ,Y=0^{-}}$.
	\item and $\sigma _{d,p,\ell } \circ \langle \Rmacro\rangle =\sigma _{d,p,\ell }$ for any $d\in \TZ*$, $p\in \mathbb{R}$ and $\ell \in \mathbb{N}$.
\end{itemize}
Therefore, \cref{eq: equations de SDA} is equivalent to
\[\left\{\begin{array}{r@{\;=\;}ll}
	\mu _{0} \Delta  \incBis & 0 &\text{ in }\Omega  \\[0.5mm]
	\mu _{1} \partial _{Y} ^2  \incBis & 0 &\text{ in }\couche \\[0.5mm]
	\mu _{1} \partial _{Y} \incBis_{|Y=0^{-}} & 0 &\text{ on }\Gamma  \\[0.5mm]
	\sigma _{d,p,\ell } (\incBis) & s_{d,p,\ell } & \ \forall (d,p,\ell )\in \TZ*\times \mathbb{R}\times \mathbb{N}
\end{array}\right.\]
Similarly to the proof of \cref{pb de Helmholtz dans A Pi}, the only solution of this system in $\Aemacro(\Pi )$ is $\fsum_{d,p,\ell } \powelne\cdot  s_{d,p,\ell } \,\phi_{d}$, which ends the proof.
\end{demo}
\endgroup 

\subsection{Scaling by $\feps$ of formal series of corner singularites}
\label{sec: Eps}

\begin{textesansboite}
The last tool we need to express and study the matching assumption \cref{eq: hypothese de raccord} in the framework of the spaces $\Aemaicro$ is a scaling operator $\Eps:\Aemacro(\Pi )\rightarrow \Aemicro(\Pi )$ s.t.:
\begin{equation}{\label{eq: E}}
\forall \varphi \in \Aeps(\Pi ), \qquad  \Eps(\varphi ) = \left\{\begin{array}{ll}
	(X,Y) \mapsto  \varphi (\varepsilon X,\varepsilon Y) &\text{ in }\Omega  \\[0.5mm]
	(X,Y) \mapsto  \varphi (\varepsilon X,Y) &\text{ in }\couche
\end{array}\right.\end{equation}
However, $\varepsilon $ is a formal indeterminate, so the above formula is ill-defined. Hence we introduce the rigorous definition below. It is basically obtained by expanding each ``$\ln (\varepsilon x)$'' into ``$\ln \varepsilon +\ln x$'' in the singularities formulas, and similarly for ``$\log (\varepsilon \alpha z)$''. The definition is a bit complicated but this section will demonstrate that $\Eps$ satisfies all the desired properties for a scaling operator.
\end{textesansboite}
\nomenclature[17-40]{$\Eps$}{Operator of scaling by $\varepsilon $ in the spaces $\Aeps$}

\begin{definitionsecable}{scaling operator $\Eps$}
\label{def: Eps}%
Let $D\in \{\Omega ,\Gamma ,\layer,\Pi \}$. We define $\Eps:\Aeps(D)\rightarrow \Aeps(D)$ by:
\begin{equation}{\label{eq: def Eps avec preEps}}
	\forall (\varphi _{d,p} )_{(d,p)\in \mathbb{R}^2 } \in \Aeps(D), \qquad  \Eps\big((\varphi _{d,p} )_{d,p} \big) := \big(\preEps(\varphi _{d,p-d} )\big)_{d,p}
\end{equation}
where $\preEps: \mathcal{A}_{d} (D)[\ln \varepsilon ]\rightarrow \mathcal{A}_{d} (D)[\ln \varepsilon ]$ is defined for~any $d\in \mathbb{R}$ as follows.
\begin{itemize}
	\item If $D=\Omega $, then for any $(q,k,j,\ell )\in \mathbb{R}\times \mathbb{N}^{3}$ s.t.\ $q+k=d$ and $[(q\in \mathbb{N} \text{ and } q\leqslant k)\Rightarrow j>0]$, we define 
	\begin{equation}{\label{eq: def preEps dans Omega}}
		\preEps\big({\ln ^{\ell } \!\varepsilon \cdot \Im[(\alpha z)^{q} \,\overline{\alpha z}^{k} \log (\alpha z)^{j} ]}\big) := \sum _{j'=0} ^{j} \binom{j}{j'} \ln ^{\ell +j'\!} \!\varepsilon  \cdot  \Im[(\alpha z)^{q} \,\overline{\alpha z}^{k} \log (\alpha z)^{j-j'} ].
	\end{equation}
	By~\cref{base de A Omega}, this defines $\preEps$ on a basis of $\mathcal{A}_{d} (\Omega )[\ln \varepsilon ]$. Then we extend it into a linear map.
	\item If $D=\Gamma $, then for any $Q\in \mathbb{C}[T]$ and $\ell \in \mathbb{N}$, we define
	\begin{equation}{\label{eq: def preEps dans Gamma}}
		\preEps(\ln ^{\ell } \!\varepsilon  \cdot x^{d} \.Q(\ln x)) := \sum _{j=0} ^{\deg Q} \frac{\ln ^{\ell +j} \!\varepsilon }{j!} \cdot x^{d} \,\partial _{T} ^{j} Q(\ln x) .
	\end{equation}
	\item Similarly, if $D=\layer$, then for any $Q\in \mathbb{C}[T,Y]$ and $\ell \in \mathbb{N}$, we define
	\begin{equation}{\label{eq: def preEps dans la couche}}
		\preEps(\ln ^{\ell } \!\varepsilon  \cdot x^{d} \.Q(\ln x,Y)) := \sum _{j=0} ^{\deg Q} \frac{\ln ^{\ell +j} \!\varepsilon }{j!} \cdot x^{d} \,\partial _{T} ^{j} Q(\ln x,Y) .
	\end{equation}
	\item if $D=\Pi $, then for any $\varphi \in \mathcal{A}_{d} (\Pi )[\ln \varepsilon ]$ we take $\preEps(\varphi )_{|\Omega } :=\preEps(\varphi _{|\Omega } )$ and $\preEps(\varphi )_{|\layer} :=\preEps(\varphi _{|\layer} )$.
\end{itemize}
\end{definitionsecable}

\begin{remarquenumerotee}
\label{rq: def Eps}%
\begin{itemize}
	\item \cref{eq: def Eps avec preEps} can be rewritten as
	\[\forall (\varphi _{d,p} )_{(d,p)\in \mathbb{R}^2 } \in \Aeps(D),\qquad  \Eps\bigg(\fsum_{(d,p)\in \mathbb{R}^2 } \varepsilon ^{p} \varphi _{d,p} \bigg) := \fsum_{(d,p)\in \mathbb{R}^2 } \varepsilon ^{p+d} \cdot \preEps(\varphi _{d,p} ).\]
	\item \Cref{Eps inversible} shows that $\Eps$ is invertible. Its inverse $\Eps^{-1}$ is a formalization of the scaling:
	\[\varphi  \mapsto  \left\{\begin{array}{l@{\,}l}
		(x,y) \mapsto  \varphi ( \frac{x}{\varepsilon } , \frac{y}{\varepsilon } ) &\text{ in }\Omega  \\[0.5ex]
		(x,Y) \mapsto  \varphi ( \frac{x}{\varepsilon } ,Y) &\text{ in }\couche
	\end{array}\right.\]
	\item The image of $\conemacro$ by $\tau :(d,p)\mapsto (d,p+d)$ is $\conemicro$. So $\Eps$ sends $\Aemacro(D)$ to $\Aemicro(D)$ and $\Eps^{-1}$ send $\Aemicro(D)$ to $\Aemacro(D)$. As~expected, $\Eps$ links the macroscopic world of $\Aemacro$ and the far fields to the microscopic world of $\Aemicro$ and the corner fields.
	\item For the first time, powers of $\ln \varepsilon $ naturally appear because of the powers of $\ln r$. This explains why the presence of these powers of $\varepsilon $ in the ansatz is necessary from the beginning.
\end{itemize}
\end{remarquenumerotee}

\begin{textesansboite}
Since the elements of $\mathcal{A}(\Pi )$ are continuous on $\Gamma $, the well-definition of $\preEps$ on $\mathcal{A}_{d} (\Pi )[\ln \varepsilon ]$ requires to have $\preEps(\varphi _{|\Omega } )_{|\Gamma } =\preEps(\varphi _{|\layer} )_{|\Gamma }$ for any $\varphi \in \mathcal{A}_{d} (\Pi )[\ln \varepsilon ]$. This is verified in the next lemma.
\end{textesansboite}

\begin{lemmesansboite}{}
$\Eps$ commutes with the restriction operators $\Aeps(\Omega )\rightarrow \Aeps(\Gamma )$ and $\Aeps(\layer)\rightarrow \Aeps(\Gamma )$ and is well defined on $\Aeps(\Pi )$.
\end{lemmesansboite}

\begin{demo}
The commutation with the restriction $\Aeps(\layer)\rightarrow \Aeps(\Gamma )$ is obvious given the definition. As~for $\Aeps(\Omega )\rightarrow \Aeps(\Gamma )$, we have for any basis function
\begin{align*}
	\preEps\big({\ln ^{\ell } \!\varepsilon \cdot \Im[(\alpha z)^{q} \,\overline{\alpha z}^{k} \log (\alpha z)^{j} ]}\big)_{|\Gamma } &= \sum _{j'=0} ^{j} \binom{j}{j'} \ln ^{\ell +j'\!} \!\varepsilon  \cdot  r^{q+k} \,\Im[\alpha ^{q-k} (\ln x-\mathrm{i}\Theta )^{j-j'} ]_{|\Gamma } \\
	&= \sum _{j'=0} ^{j} \frac{\ln ^{\ell +j'\!} \!\varepsilon }{j'!} \cdot  r^{q+k} \,\Im\!\left[ \alpha ^{q-k} \frac{j!}{(j-j') !} (\ln x-\mathrm{i}\Theta )^{j-j'} \right]_{|\Gamma } \\
	&= \sum _{j'=0} ^{j} \frac{\ln ^{\ell +j'\!} \!\varepsilon }{j'!} \cdot  r^{q+k} \.\partial _{T} Q(\ln x) \quad \ \text{with }  Q(T) :=\Im[\alpha ^{q-k} (T-\mathrm{i}\Theta )^{j} ] \\
	&= \preEps\big({\ln ^{\ell } \!\varepsilon \cdot \Im[(\alpha z)^{q} \,\overline{\alpha z}^{k} \log (\alpha z)^{j} ]}_{|\Gamma } \big)
\end{align*}
Therefore, for any $\varphi \in \mathcal{A}_{d} (\Pi )[\ln \varepsilon ]$, we have $\preEps(\varphi _{|\Omega } )_{|\Gamma } =\preEps\big((\varphi _{|\Omega } )_{|\Gamma } \big) =\preEps(\varphi _{|\Gamma } )= \preEps(\varphi _{|\layer} )_{|\Gamma }$. This shows that $\Eps$ is well defined on $\Aeps(\Pi )$.
\end{demo}

\phantomsection\label{page: def pEps}%
\begin{textesansboite}
Let us compare $\Eps$ to a classical scaling to show that $\Eps$ really behaves like a scaling operator. Let $\rho >0$. For~any $D\in \domsA$, we define the physical scaling operator $\pEps_{\rho } :\mathcal{A}(D)\rightarrow \mathcal{A}(D)$ (in opposition to the formal scaling $\Eps$) by:
\begin{itemize}
	\item $\forall \varphi \in \mathcal{A}(\Omega ),\ \pEps_{\rho } (\varphi ):(r,\theta ) \mapsto  \varphi (\rho \.r,\theta )$ in polar coordinates,
	\item $\forall \varphi \in \mathcal{A}(\Gamma ),\ \pEps_{\rho } (\varphi ):x \mapsto  \varphi (\rho  x)$,
	\item $\forall \varphi \in \mathcal{A}(\couche),\ \pEps_{\rho } (\varphi ):(x,Y) \mapsto  \varphi (\rho  x,Y)$,
	\item for any $\varphi \in \mathcal{A}(\Pi )$, $\pEps_{\rho } (\varphi )$ is equal to $\pEps_{\rho } (\varphi _{|\Omega } )$ in $\Omega $ and to $\pEps_{\rho } (\varphi _{|\couche} )$ in $\couche$.
\end{itemize}
This map is well-defined because $\mathcal{A}(D)$ is stable by scaling. Moreover, let $\Aefini(D)$ be the subset of $\Aeps(D)$ that contains the formal series that have a finite support. For~any $\peps>0$ (a real number, not an indeterminate), we define the evaluation map $\evaleps:\Aefini(D)\rightarrow \mathcal{A}(D)$ by:
\begin{equation}{\label{eq: def evaleps}}
	\evaleps\bigg(\fsum_{d\in \mathbb{R},p\in \mathbb{R}} \sum _{\ell \in \mathbb{N}} \powelne\, \varphi _{d,p,\ell } \bigg) := \sum _{d\in \mathbb{R},p\in \mathbb{R}} \sum _{\ell \in \mathbb{N}} \peps^{p} \ln ^{\ell } \!\peps\, \varphi _{d,p,\ell } 
\end{equation}
for any family $(\varphi _{d,p,\ell } )$ s.t.\ $\varphi _{d,p,\ell } \in \mathcal{A}_{d} (D)$ for any $(d,p,\ell )$ (where the above sums are finite). Then it~is easy to see that we have
\begin{equation}{\label{eq: relation Eps pEps}}
	\evaleps\circ \Eps = \pEps_{\peps} \circ \evaleps,
\end{equation}
which can be represented as a commutative diagram:
\[\xymatrixrowsep{8ex}\xymatrixcolsep{5ex}\xymatrix{
	\Aefini(D) \ar[r]^{\Eps} \ar[d]^-{\evaleps} & \Aefini(D) \ar[d]^-{\evaleps} \\
    \mathcal{A}(D) \ar[r]_{\pEps_{\peps}} & \mathcal{A}(D)
}\]
This commutation property shows that $\Eps$ behaves like a scaling in the world of formal series.\\

The next lemma will allow us to transfer to $\Eps$ some properties that are easy to verify on $\pEps_{\peps}$.
\end{textesansboite}

\begin{lemmesansboite}{}
\label{rigidite des SF 2D en eps}
Let $D\in \domsA$ and $\Phi ,\Psi  \in \Aefini(D)$. We have $\Phi  =\Psi $ iff $\evaleps(\Phi )=\evaleps(\Psi )$ for~any $\peps>0$.
\end{lemmesansboite}

\begin{demo}
By~subtraction, it suffices to treat the case $\Psi =0$. Let $(\varphi _{d,p,\ell } )_{d,p,\ell }$ be a family s.t.
\[\Phi  = \fsum_{(d,p)\in \mathbb{R}^2 } \sum _{\ell \in \mathbb{N}} \powelne\,\varphi _{d,p,\ell } \]
and $\varphi _{d,p,\ell } \in \mathcal{A}_{d} (D)$ for any $(d,p,\ell )$. For~any $(p,\ell )$, let $\varphi _{p,\ell } := \fsum_{d\in \mathbb{R}} \varphi _{d,p,\ell }$, so that we have $\Phi  =\fsum_{p} \sum _{\ell } \powelne\,\varphi _{p,\ell }$. Since $\Phi \in \Aefini(D)$, we have $\varphi _{p,\ell } \in \mathcal{A}(D)$. Let $\bfx \in D$. The map
\[f_{\bfx} :\peps \,\mapsto \, \evaleps(\Phi )(\bfx) = \sum _{p,\ell } \peps^{p} \ln ^{\ell } \!\peps\, \varphi _{p,\ell } (\bfx)\]
is well-defined because the above sum is finite and we have $f_{\bfx} \in \mathcal{A}(\Gamma )$ (for the variable $\peps$ instead of $x$). By~\cref{decomposition de A en sev,rq: base de A couche et Gamma} (which give a basis of $\mathcal{A}(\Gamma )$), we have the equivalence: $f_{\bfx} =0 \Leftrightarrow  \forall (p,\ell ),\, \varphi _{p,\ell } (\bfx)=0$. Therefore, we have: $(\forall \peps>0,\, \evaleps(\Phi )=0) \Leftrightarrow  (\forall \bfx \in D,f_{\bfx} =0) \Leftrightarrow  (\forall (p,\ell ),\, \varphi _{p,\ell } =0 )\Leftrightarrow  \Phi =0$.
\end{demo}

\begin{lemme}{commutation properties of $\Eps$}
\label{Eps commute avec fsum et eps et lneps}%
Let $D\in \{\Omega ,\Gamma ,\layer,\Pi \}$.
\begin{itemize}
	\item Let $I$ be a set and $(\Phi _{i} )_{i\in I} \in \Aeps(D)$ be a family of formal series that has collectively locally finite support. Then $(\Eps(\Phi _{i} ))_{i\in I}$ has collectively locally finite support and we have
	\[\fsum_{i\in I} \Eps(\Phi _{i} )=\Eps\bigg(\fsum_{i\in I} \Phi _{i} \bigg).\]
	\item For~any $\Phi \in \Aeps(D)$, $p\in \mathbb{R}$ and $\ell \in \mathbb{N}$, we have
	\[\powelne \cdot \Eps(\Phi ) = \Eps\big(\powelne \cdot \Phi \big).\]
\end{itemize}\fbeq
\end{lemme}

\begin{demo}
Let us extend $\preEps$ into an application $\mathcal{A}(D)[\ln \varepsilon ] \rightarrow \mathcal{A}(D)[\ln \varepsilon ]$ that has degree 0. It is then extended to $\Aeps(D)$ by \cref{def: op à degree dans les SF 2D}. Moreover, let $\tau :\Aeps(D)\rightarrow \Aeps(D)$ be defined by $\tau :(\varphi _{d,p} )_{d,p} \mapsto (\varphi _{d,p-d} )_{d,p}$. By~construction, we have $\Eps=\tau \circ \preEps$. Let us study $\preEps$ and $\tau $ separately.
\begin{itemize}
	\item $\preEps$ has a degree, so it commutes with $\fsum$ and powers of $\varepsilon $ by \cref{les op a degre commutent avec fsum et eps dans les SF 2D}. In addition, it commutes with powers of $\ln \varepsilon $ by construction.
	\item It is easy to see that $\tau $ commutes with $\fsum$, $\varepsilon $ and operators that have a degree (in particular with the multiplication by $\ln ^{\ell } \varepsilon $, which has degree 0).
\end{itemize}
\end{demo}

\begin{lemme}{}
\label{Eps inversible}%
For~any $D\in \{\Omega ,\Gamma ,\layer,\Pi \}$, \;$\Eps:\Aeps(D)\rightarrow \Aeps(D)$ is invertible
\end{lemme}


\begin{demo}
\def\invEps{\mathcal{H}^{\rm inv}_\feps}%
\def\preinvEps{\widetilde{\mathcal{H}}^{\rm inv}_\feps}%
Let us define $(\preinvEps,\invEps)$ in the same way as $(\preEps,\Eps)$ in \cref{def: Eps}, but replacing \crefrange{eq: def Eps avec preEps}{eq: def preEps dans la couche} by resp.:\vspace{-2pt}
\begin{align*}
	\Eps\big((\varphi _{d,p} )_{d,p} \big) &:= (\preEps(\varphi _{d,p\red+d} ))_{d,p}  \tag{\ref*{eq: def Eps avec preEps}'} \\[0.5ex]
	\preinvEps\big({\ln ^{\ell } \!\varepsilon \cdot \Im[(\alpha z)^{q} \,\overline{\alpha z}^{k} \log (\alpha z)^{j} ]}\big) &:= \sum _{j'=0} ^{j} \binom{j}{j'} \red{(-1)^{j'\!} } \ln ^{\ell +j'\!} \!\varepsilon  \cdot  \Im[(\alpha z)^{q} \,\overline{\alpha z}^{k} \log (\alpha z)^{j-j'} ], \tag{\ref*{eq: def preEps dans Omega}'} \\
	\preinvEps(\ln ^{\ell } \!\varepsilon  \cdot x^{d} Q(\ln x)) &:= \sum _{j=0} ^{\deg Q} \red{(-1)^{j} }\. \frac{\ln ^{\ell +j} \!\varepsilon }{j!} \cdot x^{d} \,\partial _{T} ^{j} Q(\ln x), \tag{\ref*{eq: def preEps dans Gamma}'} \\
	\preinvEps(\ln ^{\ell } \!\varepsilon  \cdot x^{d} Q(\ln x,Y)) &:= \sum _{j=0} ^{\deg Q} \red{(-1)^{j} }\. \frac{\ln ^{\ell +j} \!\varepsilon }{j!} \cdot x^{d} \,\partial _{T} ^{j} Q(\ln x,Y). \tag{\ref*{eq: def preEps dans la couche}'}
\end{align*}
Let us prove that $\invEps$ is the inverse of $\Eps$. We do it in two steps.
\begin{enumerate}
	\item Let $d\in \mathbb{R}$ and $\varphi \in \mathcal{A}_{d} (D)[\ln \varepsilon ]$. It is easy to see that, for~any $\peps>0$, we have $\evaleps\circ \,\invEps (\varphi ) = \pEps_{1/\peps} \circ \evaleps(\varphi )$. By~\cref{eq: relation Eps pEps}, this implies
	\[\evaleps\circ \Eps\circ \invEps (\varphi ) =  \pEps_{\peps} \circ  \pEps_{1/\peps} \circ \evaleps(\varphi ) =\evaleps(\varphi )\]
	which yields $\Eps\circ \invEps (\varphi ) = \varphi $ by \cref{rigidite des SF 2D en eps}.
	\item Now let $(\varphi _{d,p} )_{(d,p)\in \mathbb{R}^2 } \in \Aeps(D)$. Using the commutation of $\Eps$ (and similarly $\invEps$) with $\fsum$ and powers of $\varepsilon $ (see \cref{Eps commute avec fsum et eps et lneps}), and then using step 1, we obtain
	\[\Eps\circ \invEps\bigg(\fsum_{d,p} \varepsilon ^{p} \varphi _{d,p} \bigg) = \fsum_{d,p} \varepsilon ^{p} \.\Eps\circ \invEps(\varphi _{d,p} ) = \fsum_{d,p} \varphi _{d,p} .\]
	Similarly, one can show $\invEps\circ \Eps = \mathrm{id}_{\Aeps(D)}$.
\end{enumerate}
\end{demo}

\begin{textesansboite}
The last commutative property that we need links $\Eps$ to the operators that have a degree. There is a kind of commutativity, but only for a subclass of operators that is given by the following definition.
\end{textesansboite}

\begin{definition}{homogeneous operator}
\label{def: op homogene dans A}%
Let $D_{1} ,D_{2} \in \domsA$, $f:\mathcal{A}(D_{1} )\rightarrow \mathcal{A}(D_{2} )$ a linear map, and $d\in \mathbb{R}$. We say that $f$ is \emph{homogeneous} of order $d$ (or simply $d$-homogeneous) if we have:
\[\forall \varphi \in \mathcal{A}(D_{1} ), \;\forall \rho >0,\qquad   \pEps_{\rho } (f(\varphi )) = \rho ^{d} \.f(\pEps_{\rho } (\varphi )).\]\fbeq
\end{definition}

\begin{textesansboite}
The next result interestingly shows that there~is equivalence between the ``physical'' homogeneity of \cref{def: op homogene dans A} and the ``formal'' homogeneity w.r.t.\ $\Eps$ (although, in practice, we will only use the direct implication).
\end{textesansboite}

\begin{proposition}{}
\label{commutation Eps avec op homogene}%
Let $D_{1} ,D_{2} \in \domsA$, $f:\mathcal{A}(D_{1} )\rightarrow \mathcal{A}(D_{2} )$ a linear map, and $d\in \mathbb{R}$. We have equivalence between:
\begin{enumerate}
	\item $f$ is $d$-homogeneous,
	\item $f$ has degree $d$ (in the sense of \cref{def: degre d'un operateur}) and its natural extension $\Aeps(D_{1} )\rightarrow \Aeps(D_{2} )$ (see \cref{def: op à degree dans les SF 2D}) satisfies:
	\[\forall \Phi \in \Aeps(D_{1} ),\qquad   \Eps(f(\Phi )) = \varepsilon ^{d} \.f(\Eps(\Phi )).\]\fbeq\fbi
\end{enumerate}
\end{proposition}

\begin{exemplenumerote}
\label{ex: op homogenes ou pas dans A}%
The differential operators $\Delta _{|\Omega }$, $\partial _{x|\couche} ^2 $, $\partial _{y|\Gamma ,y=0^{+}}$, $\partial _{Y|\couche} ^2 $ and $\partial _{Y|\Gamma ,Y=0^{-}}$ (which send a space $\mathcal{A}$ into another one by \cref{operateurs de derivation dans A}) are homogeneous. However there exist operators that have a degree but are not homogeneous. Take for instance $f:\varphi \in \mathcal{A}(\Omega )\mapsto  \sigma _{d} (\varphi )\.\phiOmega_{d}$ with $d\in \mathbb{R}$. This operator keeps only the component of $\varphi $ on the basis function $\phiOmega_{d}$ and has degree 0. It satisfies
\begin{multline*}
	2^{0} \cdot f\circ \pEps_{2} \big(\Im[(\alpha z)^{d} \log (\alpha z)]\big) = f\big(2^{d} \Im[(\alpha z)^{d} \log (\alpha z)] + 2^{d} \ln 2 \,\Im[(\alpha z)^{d} ] \big) =2^{d} \ln 2 \cdot \phiOmega_{d} \\
	\neq  0 = \pEps_{2} \circ f\big(\Im[(\alpha z)^{d} \log (\alpha z)]\big)
\end{multline*}
so it is not homogeneous. The fundamental reason is that, although the basis of $\mathcal{A}(\Omega )$ used to define $\sigma _{d}$ seems natural given the explicit formulas of the elements of $\mathcal{A}(\Omega )$ (see \cref{def: espaces A}), it still relies on a rather arbitrary choice. Similarly, the operators $\ROmegaOmega$, $\RGammaOmega$, $\RcouchePi$, $\RGammaPi$ and $\ROmegaPi$ are not homogeneous because their definition depends on $\sigma _{d}$ (see \cref{operateurs R dans Omega,operateurs R}). This fact is responsible for the complexity of the matching conditions derived in \cref{sec: raccordement algebrique}, see \cref{rq: conjugaison de Rmaicro par Eps}.
\end{exemplenumerote}

\begin{demo}[Proof of \cref{commutation Eps avec op homogene}]
The implication $\Leftarrow$ can be proven by introducing $\evaleps[\rho ]$ and then commuting everything we can until we find what we want: $\forall \varphi \in \mathcal{A}(D_{1} ), \forall \rho >0,$
\begin{alignjustif}
	\rho ^{d} \.f\circ \pEps_{\rho } (\varphi ) &= \rho ^{d} \.f\circ \pEps_{\rho } \circ \evaleps[\rho ] (\varphi ) & \text{because } \varphi \in \mathcal{A}(D_{1} )\Rightarrow  \evaleps[\rho ](\varphi )=\varphi \\
	&= \rho ^{d} \.f\circ \evaleps[\rho ]\circ \Eps(\varphi ) & \text{by \cref{eq: relation Eps pEps}}  \\
	&= \rho ^{d} \,\evaleps[\rho ]\circ f\circ \Eps(\varphi ) & \text{by \cref{def: op à degree dans les SF 2D}} \\
	&= \evaleps[\rho ](\varepsilon ^{d} \.f\circ \Eps(\varphi )) \\
	&= \evaleps[\rho ]\circ \Eps\circ f(\varphi ) & \text{by assumption} \\
	&= \pEps_{\rho } \circ \evaleps[\rho ]\circ f(\varphi ) & \text{by \cref{eq: relation Eps pEps}}  \\
	&= \pEps_{\rho } \circ f(\varphi ) & \text{because } f(\varphi )\in \mathcal{A}(D_{2} )\Rightarrow  \evaleps[\rho ]\circ f(\varphi )=f(\varphi )
\end{alignjustif}
Incidentally, this calculation can be illustrated by the following diagram. The idea is that we already know that the back, sides, top and bottom faces of this cube are commutative diagrams, which implies that the front face is also commutative.
\def\diago#1{\begingroup
	\setbox0=\hbox to 0pt{\hss$\scriptstyle #1$\hss}%
	\dp0=-1pt%
	\box0\endgroup}
\def\plusbas{\vrule width0pt height8ex}
\begingroup
\catcode`\;=12\relax
\catcode`\!=12\relax 
\[\xymatrixrowsep{6ex}\xymatrixcolsep{3.5ex}\xymatrix {
    \Aeps(D_{1} ) \ar[rr]^{f} \ar[dd]_{\varepsilon ^{d} \Eps\!} \ar[dr]^*[@]{\diago{\evaleps[\rho ]}} && \Aeps(D_{2} ) \ar[dr]^*[@]{\diago{\evaleps[\rho ]}} \ar[dd]^{\Eps \plusbas} |!{[dl];[dr]}\hole \\ 
    & \mathcal{A}(D_{1} ) \ar[rr]^{f\qquad  \quad } \ar[dd]^{\rho ^{d} \pEps_{\rho } \plusbas} && \mathcal{A}(D_{2} ) \ar[dd]^{\pEps_{\rho } } \\
    \Aeps(D_{1} ) \ar[rr]^{f\qquad  \quad } |!{[ur];[dr]}\hole \ar[dr]_*[@]{\diago{\evaleps[\rho ]}} && \Aeps(D_{2} ) \ar[rd]^*[@]{\diago{\evaleps[\rho ]}} \\
    & \mathcal{A}(D_{1} ) \ar[rr]^{f\ } && \mathcal{A}(D_{2} ) \\
}\vspace{1pt}\]
\endgroup
Now let us prove the implication $\Rightarrow $. To begin, let us prove that, for~any $D\in \{D_{1} ,D_{2} \}$ and $q\in \mathbb{R}$, we have
\begin{equation}{\label{eq: Ad et Ker(pEps-truc)}}
	\mathcal{A}_{q} (D) = \bigcup _{j\in \mathbb{N}} \operatorname{Ker}\!\big[(\pEps_{2|\mathcal{A}(D)} -2^{q} )^{j} \big] =: E_{q} (D) . \vspace{-\jot}
\end{equation}
\begin{sousdemo}
	Let $q\in \mathbb{R}$ and $\varphi \in \mathcal{A}_{q} (D)$. By~construction, $\varphi $ is a linear combination of factors $\bfr^{q} \ln ^{j} \!\bfr$ (with $j\in \mathbb{N}$) multiplied by a function of $\theta $ for $D_{1} =\Omega $, or $Y$ for $D_{1} =\couche$. In addition, for any $P\in \mathbb{C}[T] \setminus \{0\}$, there exists $Q\in \mathbb{C}[T]$ s.t.\ that we have $(2\bfr)^{q} P(\ln (2\bfr))-2^{q} \bfr^{q} P(\ln \bfr) =\bfr^{q} Q(\ln \bfr)$ and $\deg (Q)=\deg (P)-1$. Hence we have $\varphi \in E_{d} (D)$.\\
	Conversely, let $\varphi \in E_{q} (D)$. By~definition of $\mathcal{A}(D)$, there exists $n\in \mathbb{N}$, $q_{1} <\cdots <q_{n}$ and $\varphi _{i} \in \mathcal{A}_{q_{i}} (D)$ for each $i\in [\![1,n]\!]$ s.t.\ $\varphi  = \sum_{i=1}^n \varphi _{i}$. We already know that $\varphi _{i} \in E_{q_{i}} (D)$ for~any $i$. Moreover, the $E_{q_{i}} (D)$ are in direct sum, as generalized eigenspaces of $\pEps_{2}$. Thus we have $\varphi _{i} =0$ for~any $i$ s.t.\ $q_{i} \neq i$, which yields $\varphi  \in \mathcal{A}_{q} (D)$.
\end{sousdemo}
Now, let $q\in \mathbb{R}$ and $\varphi \in \mathcal{A}_{q} (D_{1} )$. Let $j\in \mathbb{N}$ be s.t.\ $(\pEps_{2} -2^{q} )^{j} (\varphi )=0$. By~$d$-homogeneity of $f$, we have
\[(\pEps_{2} -2^{q+d} )^{j} \circ f(\varphi )=f\circ (2^{d} \pEps_{2} -2^{q+d} )^{j} (\varphi ) =2^{d\mkern0.5mu j} f\circ (\pEps_{2} -2^{q} )^{j} (\varphi )=0.\]
By~\cref{eq: Ad et Ker(pEps-truc)}, this implies $f(\varphi ) \in \mathcal{A}_{q+d} (D_{2} )$ and that $f$ as degree $d$ (in the sense of \cref{def: degre d'un operateur}).\\
Next, let $q\in \mathbb{R}$, $\varphi \in \mathcal{A}_{q} (D)$ and $\rho >0$. By~the same arguments as the implication $\Leftarrow$, we have
\[\evaleps[\rho ](\varepsilon ^{d} f\circ \Eps(\varphi )) =\rho ^{d} \.f\circ \pEps_{\rho } (\varphi ) \qquad  \text{ and } \qquad  \pEps_{\rho } \circ f(\varphi ) =\evaleps[\rho ]\circ \Eps\circ f(\varphi ).\]
By~$d$-homogeneity of $f$, we deduce $\evaleps[\rho ] (\varepsilon ^{d} f\circ \Eps(\varphi )) = \evaleps[\rho ](\Eps\circ f(\varphi ))$ for~any $\rho >0$, which implies $\varepsilon ^{d} f\circ \Eps(\varphi ) = \Eps\circ f(\varphi )$ by \cref{rigidite des SF 2D en eps}. Finally, this property can be extended to any element of $\Aeps(D)$ similarly to the proof of \cref{Eps inversible}, using the commutation of $f$, $\Eps$ and $\varepsilon ^{d}$ with $\fsum$ and the powers of $\varepsilon $ and $\ln \varepsilon $.
\end{demo}

\begingroup
\def\incA{\Phi}%
\def\incB{\Psi}%
\begin{corollaire}{}
\label{action de Eps sur Helmholtz dans Aemaicro}%
Let $\incA \in \Aemacro(\Pi )$. $\incA$ is a solution of 
\begin{equation}{\label{eq: pb de Helmholtz dans Aemacro bis}}
    \left\{\begin{array}{r@{\;=\;}ll}
        \mu _{0} \Delta  \incA +\omega ^2 \rho _{0} \incA & 0 &\text{ in }\Omega  \\[0.5mm]
        \mu _{1} \partial _{Y} ^2  \incA + \varepsilon ^2 (\mu _{1} \partial _{x} ^2 +\omega ^2 \rho _{1} )\incA &0&\text{ in }\couche \\[0.5mm]
        \mu _{1} \partial _{Y|Y=0^{-}} \incA - \varepsilon \,\mu _{0} \,\partial _{y|y=0^{+}} \incA &0 &\text{ on }\Gamma 
    \end{array}\right.
\end{equation}
iff $\incB:=\Eps(\incA)$ (which belongs to $\Aemicro(\Pi )$) is a solution of
\begin{equation}{\label{eq: pb de Helmholtz dans Aemicro bis}}
    \left\{\begin{array}{r@{\;=\;}ll}
        \mu _{0} \Delta  \incB +\varepsilon ^2 \omega ^2 \rho _{0} \incB & 0 &\text{ in }\Omega  \\[0.5mm]
        \mu _{1} \Delta  \incB +\varepsilon ^2 \omega ^2 \rho _{1} \incB & 0 &\text{ in }\couche \\[0.5mm]
        \mu _{1} \partial _{Y|Y=0^{-}} \incB - \mu _{0} \,\partial _{y|y=0^{+}} \incB &0 &\text{ on }\Gamma 
    \end{array}\right.
\end{equation}\fbeq
\end{corollaire}

\begin{remarque}
\cref{pb de Helmholtz dans Aemacro,pb de Helmholtz dans Aemicro}, which study the same equations, also considered the data of the coefficients $\sigma _{d,p,\ell } (\incA)$ and $\sigma _{d,p,\ell } (\incB)$. However we do~not include them here because there~is no simple relation between $(\sigma _{d,p,\ell } (\incA))_{d,p,\ell }$ and $(\sigma _{d,p,\ell } (\incB))_{d,p,\ell }$ when $\incB=\Eps(\incA)$. Indeed, there~is no commutative property between $\sigma _{d,p,\ell }$ and $\Eps$, see \cref{ex: op homogenes ou pas dans A}.
\end{remarque}

\begin{demo}
Let us assume \cref{eq: pb de Helmholtz dans Aemacro bis} and prove \cref{eq: pb de Helmholtz dans Aemicro bis}. The operator $\Delta _{|\Omega }$ sends $\Aeps(\Omega )$ to itself by \cref{operateurs de derivation dans A} and is homogeneous of order $-2$. So by \cref{commutation Eps avec op homogene} we have $\mu _{0} \Delta  \incB +\varepsilon ^2 \omega ^2 \rho _{0} \incB = \varepsilon ^2 \.\Eps(\mu _{0} \Delta  \incB +\omega ^2 \rho _{0} \incB)=0$. Using the homogeneity of the other involved differential operators, one can prove the rest of \cref{eq: pb de Helmholtz dans Aemicro bis}. Likewise for \cref{eq: pb de Helmholtz dans Aemicro bis}${} \Rightarrow  {}$\cref{eq: pb de Helmholtz dans Aemacro bis}.
\end{demo}
\endgroup 

\begin{textesansboite}
In conclusion, the technical details of the construction of $\Eps$ are not essential, but what should be remembered is:
\begin{enumerate}
	\item $\Eps$ behaves like we expect the scaling by $\varepsilon $ to behave : it~is invertible, it commutes with $\fsum$ and powers of $\varepsilon $ and $\ln \varepsilon $, it transforms PDEs like expected (\cref{action de Eps sur Helmholtz dans Aemaicro}), replacing $\varepsilon $ by $\peps\in \mathbb{R}_+^*$ in a series with finite support gives the real scaling by $\peps$, etc.,
	\item $\Eps$ has, in general, no commutation property with anything that depends on the linear forms $\sigma _{d}$.
\end{enumerate}
\end{textesansboite}


\section{Equivalence between the matching assumption and the matching conditions}
\label{sec: raccordement algebrique}

\begin{definition}{}
\label{def: raccord algebrique}%
We rigorously rewrite the matching assumption introduced in  \cref{eq: hypothese de raccord} as:
\begin{equation}{\label{eq: raccord algebrique}}
	\fsum_{p\in \NN,\ell \in \mathbb{N}} \powelne \cdot \bfuDA_{p,\ell } = \Eps^{-1} \bigg(\fsum_{p\in \NN,\ell \in \mathbb{N}} \powelne \cdot \SDA_{p,\ell } \bigg) .
\end{equation}\fbeq
\end{definition}

\begin{textesansboite}
In this section, we rigorously prove that the matching assumption \cref{eq: raccord algebrique} is equivalent to a set of scalar equations on the coefficients $\sigma _{d} (\bfuDA_{p,\ell } )$ and $\sigma _{d} (\SDA_{p,\ell } )$ that we call \emph{matching conditions} (see \cref{condition de raccord}). These matching conditions will be used in \cref{chap: construction du DA} to build the fields $\bfu_{p,\ell }$ and $S_{p,\ell }$ by induction.
\end{textesansboite}

\refstepcounter{ansatz}
\begin{paragraphesansboite}{Ansatz \theansatz:}
\label{ansatz sur Aeps}%
In this section, we assume
\[\fsum_{p\in \NN,\ell \in \mathbb{N}} \!\powelne \cdot \bfuDA_{p,\ell } \.\in \. \fsum_{p\in \NN} \fsum_{d\in \NN-p} \varepsilon ^{p} \mathcal{A}_{d} (\Pi )[\ln \varepsilon ] \quad \ \text{ and } \quad  \fsum_{p\in \NN,\ell \in \mathbb{N}} \!\powelne \cdot \SDA_{p,\ell } \.\in \. \fsum_{p\in \NN} \fsum_{d\in p-\NN} \varepsilon ^{p} \mathcal{A}_{d} (\Pi )[\ln \varepsilon ].\]
\end{paragraphesansboite}

\begin{figure}[h]
\def\ymin{-0.65}
\def\ymax{3.4}
\def\xmax{\ymax}
\foreach \piSurTheta in {1.5,3.14159265359/2} {
\begin{tikzpicture}[scale=0.75]
	\fill [green!70!black, opacity=0.2] (0,0) -- (-\xmax,\ymax) -- (\xmax,\ymax) -- (\xmax,0) -- cycle;
	\draw [very thick, green!55!black] (-\xmax,0) -- (\xmax,0) node [near start, below] {$p=0$};
	\draw [very thick, green!55!black] (-\xmax,\ymax) -- ({-\xmax/\ymax*\ymin},\ymin);
	\draw [very thick, green!55!black] (-2.55,2.55) node [left] {$d+p=0\;$};
	\draw [thin, gray, dashed] (-\xmax,\ymin) grid (\xmax,\ymax);
	\draw [->] (-\xmax-0.3,0) -- (\xmax+0.3,0) node [right] {$d$};
	\draw [->] (0,\ymin-0.1) -- (0,\ymax+0.3) node [right] {$p$};
	\draw (1,0) node{$+$} node[below]{\ \ 1};
	\draw (\piSurTheta,0) node{$+$} node[below]{$\,\frac\pi\Theta$};
	\draw (0,1) node{$+$} node[right]{1};
	\begin{scope}
		\clip (-\xmax,\ymin) rectangle (\xmax,\ymax); 
		\foreach \ip in {0,1,...,4} {
			\foreach \jp in {0,1,...,4} {
				\foreach \id in {0,1,...,4} {
					\foreach \jd in {0,1,...,4} {
						\draw ({\id+\piSurTheta*\jd-\ip-\piSurTheta*\jp},{\ip+\piSurTheta*\jp}) node{\footnotesize $\bullet$};}}}}
	\end{scope}
\end{tikzpicture}}\par
\caption[Points of $\{(d,p) \mid p\in \NN \text{ and } d\in \NN-p\}$ for $\Theta = \frac{2\pi }{3}$ and $\Theta =2$]{Points of $\{(d,p) \mid p\in \NN \text{ and } d\in \NN-p\}$ for $\Theta = \frac{2\pi }{3}$ (on the left) and $\Theta =2$ (on the right)}
\label{fig: cone de NN}
\end{figure}
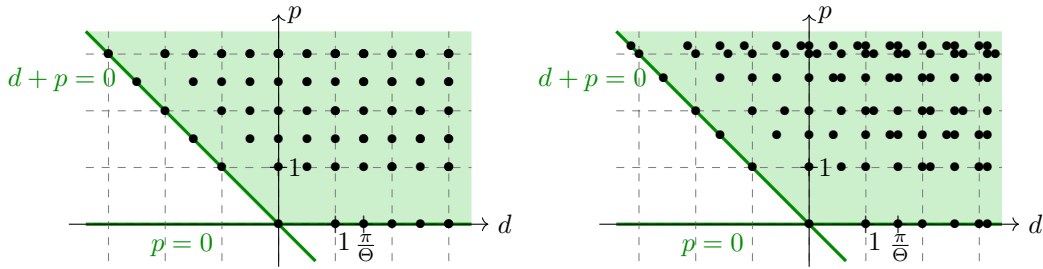

\begin{textesansboite}
See \cref{fig: cone de NN} for a representation of $\{(d,p) \mid p\in \NN, d\in \NN-p\}$ in the two cases $\pst \in \mathbb{Q}$ and $\pst \not\in \mathbb{Q}$. In particular, \cref{ansatz sur Aeps} implies $\fsum_{p,\ell } \powelne\,\bfuDA_{p,\ell } \in \Aemacro(\Pi )$ and $\fsum_{p,\ell } \powelne\,\SDA_{p,\ell } \in \Aemicro(\Pi )$. It also implies that, for~any $(p,\ell )\in \NN\times \mathbb{N}$, $\bfu_{p,\ell }$ has components only in the $\mathcal{A}_{d} (\Pi )$ s.t.\ $d\in \NN-p$ and $S_{p,\ell }$ has only in the $\mathcal{A}_{d} (\Pi )$ s.t.\ $d\in p-\NN$. 
Once the fields $\bfu_{p,\ell }$ and $S_{p,\ell }$ are built, \cref{ansatz sur Aeps} will be rigorously verified in \cref{confimation de l'hyp de raccord}.
\end{textesansboite}

\begin{remarque}
\Cref{ansatz sur Aeps} is natural given the matching assumption \cref{eq: raccord algebrique} and the initial ansatz that $u_{\peps}$ has an asymptotic expansion with powers of $\peps$ that belong to $\NN$ (see \cref{sec: ansatz du dvlpmt raccorde}). Indeed, if we assume \cref{eq: raccord algebrique} and that $\bfuDA_{p,\ell }$ and $\SDA_{p,\ell }$ vanish when $p\not\in \mathbb{R} \setminus \NN$, and if we denote $\bfuDA_{\varepsilon } := \fsum_{p,\ell } \powelne\,\bfuDA_{p,\ell }$, $\SDA_{\varepsilon } :=\fsum_{p,\ell } \powelne\,\SDA_{p,\ell }$ and $\tau :(d,p)\mapsto (d,p-d)$, then we have $\operatorname{supp}(\bfuDA_{\varepsilon } ) \subset  \mathbb{R}\times \NN$ and $\operatorname{supp}(\Eps^{-1} (\SDA_{\varepsilon } )) = \tau (\operatorname{supp}(\SDA_{\varepsilon } )) \subset  \tau (\mathbb{R}\times \NN)$, which imply $\operatorname{supp}(\bfuDA_{\varepsilon } ) \subset  (\mathbb{R}\times \NN) \cap  \tau (\mathbb{R}\times \NN) = \{(d,p) \mid p\in \NN, p+d\in \NN\}$. Likewise for the support of $\SDA_{\varepsilon }$.\\

\Cref{ansatz sur Aeps} and the ansatz of \cref{sec: ansatz du dvlpmt raccorde} simplify and clarify calculations. In particular, they directly fix the sum indexes of the matching conditions in \cref{condition de raccord}. However, we point out that these ansatze can be avoided. Indeed, one can manipulate formal series with unknown support (but still locally finite by \cref{def: support locmt fini et SF}), like e.g.\ ``$\fsum_{\red{p\in \mathbb{R}},\ell \in \mathbb{N}} \powelne\,\bfuDA_{p,\ell }$'', and find out which field vanish after the fields' construction. The only necessary a priori knowledge is the form of the corner singularities, given by the definition of the spaces $\mathcal{A}$ and $\Aeps$ (in particular the presence of logarithms of $r$ and $\varepsilon $, which is needed to define $\ROmegaPi$, $\RcouchePi$, $\RGammaPi$ and $\Eps$).
\end{remarque}

\begin{textesansboite}
The following theorem is an easy consequence of \cref{pb de Helmholtz dans Aemacro,pb de Helmholtz dans Aemicro}.
\end{textesansboite}

\begin{theoreme}{}
\label{champs totaux en fct de sigma}%
We have the following equalities in $\Aemacro(\Pi )$ and $\Aemicro(\Pi )$ respectively: \begin{align}
	\fsum_{p\in \NN,\ell \in \mathbb{N}} \powelne \cdot \bfuDA_{p,\ell } &= \langle \Rmacro\rangle \bigg(\fsum_{p\in \NN,\ell \in \mathbb{N}} \powelne\, \sumdmacro \sigma _{d} (\bfuDA_{p,\ell } )\,\phi_{d} \bigg)
	\label{eq: champ lointain total}\\
	\fsum_{p\in \NN,\ell \in \mathbb{N}} \powelne\cdot \SDA_{p,\ell } &= \langle \Rmicro\rangle \bigg(\fsum_{p\in \NN,\ell \in \mathbb{N}} \powelne\, \sumdmicro \sigma _{d} (\SDA_{p,\ell } )\,\phi_{d} \bigg)
	\label{eq: champ de coin total}
\end{align}\fbeq
\end{theoreme}

\begin{demo}
We only prove \cref{eq: champ lointain total}, as \cref{eq: champ de coin total} is similar. By~\cref{ansatz: existence de bfuDApl et SDApl}, we have for~any $(p,\ell )\in \NN\times \mathbb{N}$ :
\[\left\{\begin{array}{r@{\;=\;}ll}
	\mu _{0} \Delta  \bfuDA_{p,\ell } +\omega ^2 \rho _{0} \bfuDA_{p,\ell } & 0 &\text{ in }\Omega  \\[0.5mm]
	\mu _{1} \partial _{Y} ^2  \bfuDA_{p,\ell } & -(\mu _{1} \partial _{x} ^2 +\omega ^2 \rho _{1} )\bfuDA_{p-2,\ell } &\text{ in }\couche \\[0.5mm]
	\mu _{1} \partial _{Y|Y=0^{-}} \bfuDA_{p,\ell } & \mu _{0} \,\partial _{y|y=0^{+}} \bfuDA_{p-1,\ell } &\text{ on }\Gamma 
\end{array}\right.\]
Let $\bfuDA_{\varepsilon } :=\fsum_{p,\ell } \powelne\,\bfuDA_{p,\ell }$. By~multiplying the above equations by $\powelne$, summing over $(p,\ell )$ and using the commutation of $\fsum$ and $\powelne$ with operators that have a degree (see \cref{les op a degre commutent avec fsum et eps dans les SF 2D}), one easily obtains:
\[\left\{\begin{array}{r@{\;=\;}ll}
	\mu _{0} \Delta  \bfuDA_{\varepsilon } +\omega ^2 \rho _{0} \bfuDA_{\varepsilon } & 0 &\text{ in }\Omega  \\[0.5mm]
	\mu _{1} \partial _{Y} ^2  \bfuDA_{\varepsilon } + \varepsilon ^2 (\mu _{1} \partial _{x} ^2 +\omega ^2 \rho _{1} )\bfuDA_{\varepsilon } &0&\text{ in }\couche \\[0.5mm]
	\mu _{1} \partial _{Y|Y=0^{-}} \bfuDA_{\varepsilon } - \varepsilon \,\mu _{0} \,\partial _{y|y=0^{+}} \bfuDA_{\varepsilon } &0 &\text{ on }\Gamma 
\end{array}\right.\]
Finally, we have $\bfuDA_{\varepsilon } \in \Aemacro(\Pi )$ thanks to \cref{ansatz sur Aeps}. So we can apply \cref{pb de Helmholtz dans Aemacro}, which concludes the proof.
\end{demo}

\begin{remarquenumerotee}
\label{rq: utilite des operateurs R}%
\Cref{champs totaux en fct de sigma} is very important because it factorizes the behaviors around the corner of the fields into two parts:
\begin{itemize}
	\item a geometric series of operators, that encapsulates all the complexity of the singularities' expressions,
	\item and the coefficients $\sigma _{d} (\dots )$, that need to be fixed for the different fields to match.
\end{itemize}
This separation allows us to reduce the matching assumption to only fixing the coefficients $\sigma _{d} (\dots )$. This result is one of the main reasons why we introduced the operators $\ROmegaPi$, $\RcouchePi$ and $\RGammaPi$.\\

The alternative description of the fields' singularities that was used in the literature for Poisson's equation uses recursive sequences, see \cite{CalCosDauVia06,DelSchSem16}. Let us explain why this was possible for Poisson but is~not for Helmholtz's equation. When $\omega =0$, $\Rmaicro$ contains only two non-zero elements (instead of four) and their degrees are aligned in $\mathbb{R}^2 $ (see \cref{fig: degres des elements de Rmacro et Rmicro}). One can prove in this case that, for~any $d\in \TZ*$, we have $\langle \Rmicro\rangle (\phi_{d} )=\fsumn \phi_{d,n}$ where $\phi_{d,n} \in \mathcal{A}_{d-n} (\Pi )$ for~any $n$ and the sequence $(\phi_{d,n} )_{n}$ is defined by
\[\textstyle \phi_{d,-1} := 0, \quad  \phi_{d,0} :=\phi_{d} \quad  \text{ and } \quad  \forall n\in \mathbb{N}^*,\ \phi_{d,n} = \frac{\mu _{0} }{\mu _{1}} \RGammaPi\circ \partial _{y|\Gamma ,y=0^{+}} (\phi_{d,n-1} ) -\RcouchePi\circ \partial _{x|\couche} ^2 (\phi_{d,n-2} ).\]
Note that $\frac{\mu _{0} }{\mu _{1}} \RGammaPi\circ \partial _{y|\Gamma ,y=0^{+}}$ and $-\RcouchePi\circ \partial _{x|\couche} ^2 $ are the two non-zero elements of $\Rmicro$. In the literature, the inductive relation on $\phi_{d,n}$ is written in a way closer to:
\[\left\{\begin{array}{l@{\;=\;}ll}
	\Delta  \phi_{d,n} & 0 &\text{ in }\Omega  \\[0.5mm]
	\partial _{Y} ^2  \phi_{d,n} & - \partial _{x} ^2  \phi_{d,n-2} &\text{ in }\couche \\[0.5mm]
	\mu _{1} \partial _{Y} \phi_{d,n|Y=0^{-}} & \mu _{0} \.\partial _{y} \phi_{d,n-1|\Gamma ,y=0^{+}} &\text{ on }\Gamma  
\end{array}\right.\]
Similarly for $\langle \Rmacro\rangle (\phi_{d} )$. However, this description of $\langle \Rmaicro\rangle (\phi_{d} )$ is impossible when $\omega \neq 0$ because 
the elements of $\Rmaicro$ point towards different directions (see \cref{fig: degres des elements de Rmacro et Rmicro}). One could use sequences that are defined recursively w.r.t.\ two indexes $n,m$ at the same time, but this would amplify the maze of calculations that is typically obtained with this kind of asymptotic expansion. On the contrary, the operators $\ROmegaPi$, $\RcouchePi$ and $\RGammaPi$ offer a neat way to describe all those corner singularities, with a simple geometric series. See also \cref{ex: 2e pb jouet dans C[[X]]} for a toy example of this idea.
\end{remarquenumerotee}

\begin{remarque}
\Cref{champs totaux en fct de sigma} could be proven by inductively composing Equations~\cref{eq: bfuDApl} and \cref{eq: SDApl} that express each $\bfuDA_{p,\ell }$, resp.\ $\SDA_{p,\ell }$, in terms of the previous fields' behaviors. However, the proof would involve heavy sum manipulations. Instead, the results on the resolution of Helmholtz and Poisson problems in the spaces $\Amaicro$ and $\Aemaicro$ (which rely on the easiness to invert operators of the form $ \mathrm{id}-f$ among formal series, see \cref{inverse de id-f dans les SF 1D,inverse de id-f dans les SF 2D}) provide a much more elegant way to prove \cref{champs totaux en fct de sigma}.
\end{remarque}

\begin{textesansboite}
The next step in order to isolate the coefficients $\sigma _{d} (\bfuDA_{p,\ell } )$ and $\sigma _{d} (\SDA_{p,\ell } )$ from $\langle \Rmaicro\rangle $ is to project the matching assumption of \cref{def: raccord algebrique} on the subspace of $\Aeps(\Pi )$ generated by the functions $\phi_{d}$. Thus we introduce the projector $\projphi:\Aeps(\Pi )\rightarrow \Aeps(\Pi )$, defined by:
\[\forall \Phi \in \Aeps(\Pi ), \qquad  \projphi(\Phi ) := \fsum_{d\in \TZ*} \fsum_{p\in \mathbb{R}} \sum _{\ell \in \mathbb{N}}  \powelne\cdot  \sigma _{d,p,\ell } (\Phi )\,\phi_{d}\]
\end{textesansboite}

\begin{lemme}{}
\label{hyp raccord projetee sur les phid}%
The matching assumption of \cref{def: raccord algebrique} is equivalent to each of these equalities:
\begin{align}
	\fsum_{p\in \NN,\ell \in \mathbb{N}} \powelne \sumdmacro \sigma _{d} (\bfuDA_{p,\ell } )\,\phi_{d} &= \projphi \circ  \Eps^{-1} \circ  \langle \Rmicro\rangle  \bigg(\fsum_{p\in \NN,\ell \in \mathbb{N}} \powelne \sumdmicro \sigma _{d} (\SDA_{p,\ell } )\,\phi_{d} \bigg)
	\label{eq: hyp raccord projetee sur les phid pt de vue macro}\\
	\fsum_{p\in \NN,\ell \in \mathbb{N}} \powelne \sumdmicro \sigma _{d} (\SDA_{p,\ell } )\,\phi_{d} &= \projphi \circ  \Eps \circ  \langle \Rmacro\rangle  \bigg(\fsum_{p\in \NN,\ell \in \mathbb{N}} \powelne \sumdmacro \sigma _{d} (\bfuDA_{p,\ell } )\,\phi_{d} \bigg)
	\label{eq: hyp raccord projetee sur les phid pt de vue micro}
\end{align}\fbeq
\end{lemme}

\begin{remarque}
A priori, when applying the projector $\projphi$, one may lose some information (namely the components of the series $\bfuDA_{p,\ell }$ and $\SDA_{p,\ell }$ in $\operatorname{Ker}\projphi$). However, $\projphi$ keeps the data of the coefficients $\sigma _{d} (\bfuDA_{p,\ell } )$ and $\sigma _{d} (\SDA_{p,\ell } )$, which entirely determine the series $\bfuDA_{p,\ell }$ and $\SDA_{p,\ell }$ by \cref{champs totaux en fct de sigma}. Thus we are able to show that \cref{eq: hyp raccord projetee sur les phid pt de vue macro,eq: hyp raccord projetee sur les phid pt de vue micro} are not~only necessary conditions for the matching assumption but also sufficient ones.
\end{remarque}

\begin{demo}
Let us prove only \cref{eq: hyp raccord projetee sur les phid pt de vue macro}. By~\cref{champs totaux en fct de sigma}, the matching assumption of \cref{def: raccord algebrique} is equivalent to:
\begin{equation}{\label{eq: raccord algebrique presque prouve}}
\langle \Rmacro\rangle \bigg(\fsum_{p,\ell } \powelne\, \fsum_{d} \sigma _{d} (\bfuDA_{p,\ell } )\phi_{d} \bigg) = \Eps^{-1} \circ  \langle \Rmicro\rangle  \bigg(\fsum_{p,\ell } \powelne\, \fsum_{d} \sigma _{d} (\SDA_{p,\ell } )\phi_{d} \bigg)\end{equation}
Moreover, \cref{operateurs R} implies: $\forall R\in \Rmacro,\; \projphi\circ R=0$. Since $\langle \Rmacro\rangle $ is a sum of the identity and non-trivial products of elements of $\Rmacro$, we have $\projphi\circ \langle \Rmacro\rangle =\projphi\circ  \mathrm{id}=\projphi$. Thus applying $\projphi$ to \cref{eq: raccord algebrique presque prouve} gives \cref{eq: hyp raccord projetee sur les phid pt de vue macro}.\\
It remains to show conversely that \cref{eq: hyp raccord projetee sur les phid pt de vue macro} implies \cref{eq: raccord algebrique presque prouve}. Let us denote $\Phi _{1}$ and $\Phi _{2}$ the two sides of \cref{eq: raccord algebrique presque prouve}. Let us show that they satisfy 
\begin{equation}{\label{eq: hyp raccord projetee sur les phid preuve}}
	\forall i\in \{1,2\}, \quad \ \left\{\begin{array}{r@{\;=\;}ll}
		\mu _{0} \Delta  \Phi _{i} +\omega ^2 \rho _{0} \Phi _{i} & 0 &\text{ in }\Omega  \\[0.5mm]
		\mu _{1} \partial _{Y} ^2  \Phi _{i} + \varepsilon ^2 (\mu _{1} \partial _{x} ^2 +\omega ^2 \rho _{1} )\Phi _{i} &0&\text{ in }\couche \\[0.5mm]
		\mu _{1} \partial _{Y|Y=0^{-}} \Phi _{i} - \varepsilon \,\mu _{0} \,\partial _{y|y=0^{+}} \Phi _{i} &0 &\text{ on }\Gamma 
	\end{array}\right.
\end{equation}
For $\Phi _{1}$ this is a direct consequence of \cref{pb de Helmholtz dans Aemacro} (which can be applied thanks to \cref{ansatz sur Aeps}). As~for $\Phi _{2}$, \cref{pb de Helmholtz dans Aemicro} shows that $\widetilde{\Phi }_{2} := \langle \Rmicro\rangle  \big(\fsum_{p,\ell } \powelne\, \fsum_{d} \sigma _{d} (\SDA_{p,\ell } )\phi_{d} \big) \in \Aemicro(\Pi )$ is a solution of the rescaled Helmholtz problem \cref{eq: pb de Helmholtz dans Aemicro bis}, so $\Phi _{2} =\Eps^{-1} ( \widetilde{\Phi }_{2} )$ satisfies \cref{eq: hyp raccord projetee sur les phid preuve} by \cref{action de Eps sur Helmholtz dans Aemaicro}. Moreover, we have shown that \cref{eq: hyp raccord projetee sur les phid pt de vue macro} is equivalent to $\projphi(\Phi _{1} )=\projphi(\Phi _{2} )$. Hence for any $(d,p,\ell )\in \TZ*\times \mathbb{R}\times \mathbb{N}$, we have $\sigma _{d,p,\ell } (\Phi _{1} )= \sigma _{d,p,\ell } (\Phi _{2} )$ ; let us denote $s_{d,p,\ell }$ this common quantity. \Cref{pb de Helmholtz dans Aemacro} states the problem \cref{eq: hyp raccord projetee sur les phid preuve} with the additional condition $\sigma _{d,p,\ell } (\Phi _{i} )=s_{d,p,\ell }$ for any $(d,p,\ell )$ has a unique solution in $\Aemacro(\Pi )$. This proves $\Phi _{1} =\Phi _{2}$, i.e.\ \cref{eq: hyp raccord projetee sur les phid pt de vue macro}.
\end{demo}

\nomenclature[12-30-28]{$\cmacro$}{Matching coefficients that give the $\sigma _{d} (\bfuDA_{p,\ell } )$ in terms of the $\sigma _{d'} (\SDA_{p',\ell '} )$}
\nomenclature[12-28-30]{$\cmicro$}{Matching coefficients that give the $\sigma _{d} (\SDA_{p,\ell } )$ in terms of the $\sigma _{d'} (\bfuDA_{p',\ell '} )$}
\nomenclature[60]{$\intvl{a,b}$}{Finite subset of $[a,b]$}
\begin{textesansboite}
We define the \emph{matching coefficients} as, for~any $(d,d',p,\ell )\in (\TZ*)^2 \times \NN\times \mathbb{N}$, the following complex numbers:\vspace{-1mm}
\begin{equation}{\label{eq: def cmaicro}}
\cmacro_{d,d',p,\ell } := \sigma _{d,p,\ell } \circ  \Eps^{-1} \circ  \langle \Rmicro\rangle  (\phi_{d'} ) \qquad  \text{ and } \qquad 
\cmicro_{d,d',p,\ell } := \sigma _{d,p,\ell } \circ  \Eps \circ  \langle \Rmacro\rangle  (\phi_{d'} )
\end{equation}
\Cref{sec: implementation des espaces A} explains how to compute numerically these coefficients. In addition, for~any $a,b\in \mathbb{R}$, we denote
\begin{equation}{\label{eq: def intvl}}
	\intvl{a,b}:=\{c\in \mathbb{R} \mid c-a\in \NN \text{ and } b-c\in \mathbb{N}\}.
\end{equation}
It~is a finite subset of $[a,b]$.\\

By~identifying the components of \cref{eq: hyp raccord projetee sur les phid pt de vue macro,eq: hyp raccord projetee sur les phid pt de vue micro} on each basis function $\phi_{d}$, we are able to prove the following result, which re-expresses the matching assumption as a set of scalar relations on the coefficients $\sigma _{d} (\bfuDA_{p,\ell } )$ and $\sigma _{d} (\SDA_{p,\ell } )$. This is~not quite yet the matching conditions that will be used to build the fields, because only a part of these relations will be used, as explained after \cref{pre conditions de raccord}.
\end{textesansboite}

\begin{lemme}{relations between the $\sigma _{d} (\bfuDA_{p,\ell } )$ and the $\sigma _{d} (S_{p,\ell } )$}
\label{pre conditions de raccord}%
The matching assumption of \cref{def: raccord algebrique} is equivalent to each of the following statements:
\begin{enumerate}
	\item $\forall (d,p,\ell )\in \TZ*\times \NN\times \mathbb{N}, \quad \displaystyle \sigma _{d} (\bfuDA_{p,\ell } ) = \sum _{p'\in \intvl{0,p+d}} \sum _{\substack{d'\in \TZ*\\ p'-d'\in \intvl{0,p}}} \!\sum _{\ell '=0} ^{\ell } \,\cmacro_{d,d',p-p',\ell -\ell '} \cdot  \sigma _{d'} (\SDA_{p',\ell '} )$
	\item $\forall (d,p,\ell )\in \TZ*\times \NN\times \mathbb{N}, \quad \displaystyle \sigma _{d} (\SDA_{p,\ell } ) = \sum _{p'\in \intvl{0,p-d}} \sum _{\substack{d'\in \TZ*\\ p'+d'\in \intvl{0,p}}} \!\sum _{\ell '=0} ^{\ell } \,\cmicro_{d,d',p-p',\ell -\ell '} \cdot  \sigma _{d'} (\bfuDA_{p',\ell '} )$.
\end{enumerate}\fbi
\end{lemme}

\begin{demo}
We will only prove point 1, as point 2 is similar. We can rewrite \cref{eq: hyp raccord projetee sur les phid pt de vue macro} as
\[\fsum_{p\in \NN,\ell \in \mathbb{N}} \,\fsum_{d\in \TZ* \cap (\NN-p)} \!\! \powelne \cdot \sigma _{d} (\bfuDA_{p,\ell } )\,\phi_{d} = \fsum_{p\in \NN,\ell \in \mathbb{N}} \,\fsum_{d\in \TZ* \cap (\NN-p)} \!\! \powelne \cdot \sigma _{d} (\SDA_{p,\ell } ) \cdot  \projphi \circ  \Eps^{-1} \circ  \langle \Rmicro\rangle  (\phi_{d} )\]
This is possible thanks to the many commutation properties of \cref{chap: SF,chap: raccord algebrique}:
\begin{itemize}
	\item $\langle \Rmacro\rangle $ commutes with $\fsum$ on $\Amacro(\Pi )$ by \cref{fsum et le crochet commutent en 2D} and with $\powelne$ by \cref{rq: eps et le crochet commutent,sec: SF 2D de singularites de coin},
	\item $\Eps^{-1}$ commutes with $\fsum$ and $\powelne$ by \cref{Eps commute avec fsum et eps et lneps},
	\item it~is easy to see that $\projphi$ has degree 0, so it commutes with $\fsum$ and $\powelne$ by \cref{les op a degre commutent avec fsum et eps dans les SF 2D,sec: SF 2D de singularites de coin},
	\item $\powelne$ commutes with $\fsum$ by \cref{fsum et eps commutent,sec: SF 2D de singularites de coin}.
\end{itemize}
Next, let $(d,p,\ell )\in \TZ*\times \NN\times \mathbb{N}$. Applying $\sigma _{d,p,\ell }$ and using \cref{def: sigma et fsum commutent dans Aemaicro} gives:
\begin{equation}{\label{eq: pre conditions de raccord preuve 2}}
	\sigma _{d} (\bfuDA_{p,\ell } ) = \sum _{p'\in \NN,\ell '\in \mathbb{N}} \,\sum _{d'\in \TZ* \cap (\NN-p')} \!\! \sigma _{d,p,\ell } \big(\varepsilon ^{p'} \ln ^{\ell '\!} \!\varepsilon  \cdot  \projphi \circ  \Eps^{-1} \circ  \langle \Rmicro\rangle  (\phi_{d'} )\big) \cdot \sigma _{d'} (\SDA_{p',\ell '} )
\end{equation}
where the above sums are finite. Conversely, since both sides of \cref{eq: hyp raccord projetee sur les phid pt de vue macro} belong to the range of $\projphi$, having \cref{eq: pre conditions de raccord preuve 2} for any $(d,p,\ell )$ is equivalent to \cref{eq: hyp raccord projetee sur les phid pt de vue macro}, and thus to the matching assumption by \cref{hyp raccord projetee sur les phid}.\\
Now, we clearly have, for any $\Phi \in \Aeps(\Pi )$, $\sigma _{d,p,\ell } \big(\varepsilon ^{p'} \ln ^{\ell '\!} \!\varepsilon  \,\Phi \big) = \sigma _{d,p-p',\ell -\ell '} (\Phi )$ if $\ell \geqslant \ell '$ and $\sigma _{d,p,\ell } \big(\varepsilon ^{p'} \ln ^{\ell '\!} \!\varepsilon  \,\Phi \big) =0$ otherwise. Thus, \cref{eq: pre conditions de raccord preuve 2} is equivalent to
\begin{equation}{\label{eq: pre conditions de raccord preuve 3}}
	\sigma _{d} (\bfuDA_{p,\ell } ) = \sum _{p'\in \NN} \sum _{d'\in \TZ* \cap (\NN-p')} \sum _{\ell '=0} ^{\ell } \cmacro_{d,d',p-p',\ell -\ell '} \cdot  \sigma _{d'} (\SDA_{p',\ell '} )
\end{equation}
because $\sigma _{d,p-p',\ell -\ell '} \circ \projphi = \sigma _{d,p-p',\ell -\ell '}$. Moreover, for~any $R\in \Rmacro$, we have $\deg _{\varepsilon } \!R \in \mathbb{N}$ and $\deg _{\varepsilon } \!R+\deg _{\mathcal{A}} \!R\in \mathbb{N}$ (see \cref{fig: degres des elements de Rmacro et Rmicro}). So for $\cmacro_{d,d',p-p',\ell -\ell '}$ to be non-zero, we need $p-p'+d'\in \mathbb{N}$ and $(p-p'+d')+d-d'\in \mathbb{N}$. Therefore, \cref{eq: pre conditions de raccord preuve 3} is equivalent to the desired formula of $\sigma _{d} (\bfuDA_{p,\ell } )$.
\end{demo}

\begin{textesansboite}
Actually, we will use only half of the relations of \cref{pre conditions de raccord}: those giving $\sigma _{d} (\bfuDA_{p,\ell } )$ when $d<0$ and those giving $\sigma _{d} (\SDA_{p,\ell } )$ when $d>0$. First because these are the ones that depend only on the fields with a smaller index $p$, which allows to build the fields inductively. Secondly because these coefficients $\sigma _{d} ({\dots })$ are the ones associated with functions $\phi_{d}$ that are singular for these fields. Indeed, $\phi_{d}$ is singular at the corner only when $d<0$ and at infinity only when $d>0$. \Cref{sec: cadre fonctionnel des champs lointains,sec: cadre fonctionnel des champs de coin} will show that the singular behavior of the fields can be imposed, but not their variational part.
\end{textesansboite}

\begin{theoreme}{matching conditions}
\label{condition de raccord}%
We recall that we assume in this section that \cref{ansatz: existence de bfuDApl et SDApl,ansatz sur Aeps} are satisfied. Then the matching assumption of \cref{def: raccord algebrique} is equivalent to the following set of equations called the \emph{matching conditions}:\vspace{-1mm}
\begin{subequations}
\label{eq: condition de raccord}
\begin{empheq}[left={\hspace{-5mm}\empheqlbrace\,}]{alignat=3}
	& \textstyle \forall (d,p,\ell )\in (-\TN){\times }\NN{\times }\mathbb{N}, \ \
	& \sigma _{d} (\bfuDA_{p,\ell } ) = \!\!\! \sum _{p'\in \intvl{0,p+d}} \! \sum _{\substack{d'\in \TZ*\\ p'-d'\in \intvl{0,p}}} \!\! \sum _{\ell '=0} ^{\ell } \,\cmacro_{d,d',p-p',\ell -\ell '} \cdot  \sigma _{d'} (\SDA_{p',\ell '} ) 
	\label{eq: condition de raccord 1} \\[-1mm]
	& \textstyle \forall (d,p,\ell )\in \TN\times \NN\times \mathbb{N}, \ \ 
	& \sigma _{d} (\SDA_{p,\ell } ) = \!\!\! \sum _{p'\in \intvl{0,p-d}} \! \sum _{\substack{d'\in \TZ*\\ p'+d'\in \intvl{0,p}}} \!\! \sum _{\ell '=0} ^{\ell } \,\cmicro_{d,d',p-p',\ell -\ell '} \cdot  \sigma _{d'} (\bfuDA_{p',\ell '} )
	\label{eq: condition de raccord 2}
\end{empheq}\end{subequations}\fbeq
\end{theoreme}

\begin{textesansboite}
By~\cref{pre conditions de raccord} these conditions are necessary, and we must show that they are sufficient. This is done by re-expressing the problem as a simple linear algebra problem with $2\times 2$ block matrices (where the blocks have infinite size and separate the indexes $d>0$ and $d<0$) and doing some simple manipulations on those matrices.
\end{textesansboite}

\begin{demo}
Let $\Vmacro$ be the set of families of complex numbers $(\smacro_{p,d,\ell } )_{p\in \NN,d\in \TZ*\cap (\NN-p),\ell \in \mathbb{N}}$ s.t.\ for~any $(p,d)$ only a finite number of the $\smacro_{p,d,\ell }$ are non-zero. We define $\Vmicro$ similarly by replacing ``$d\in \NN-p$'' with ``$d\in p-\NN$''. We define the following linear maps:\vspace{-1mm}
\[\Pmacro : \Vmicro \rightarrow  \Vmacro,\ (\smicro_{p,d,\ell } )_{p,d,\ell } \ \mapsto \ \bigg( \sum _{p'\in \intvl{0,p+d}} \sum _{\substack{d'\in \TZ*\\ p'-d'\in \intvl{0,p}}} \!\sum _{\ell '=0} ^{\ell } \cmacro_{d,d',p-p',\ell -\ell '} \cdot  \smicro_{p',d',\ell '} \bigg)_{\!p,d,\ell } \vspace{-4mm}\]
\[\Pmicro : \Vmacro \rightarrow  \Vmicro,\ (\smacro_{p,d,\ell } )_{p,d,\ell } \ \mapsto \ \bigg( \sum _{p'\in \intvl{0,p-d}} \sum _{\substack{d'\in \TZ*\\ p'+d'\in \intvl{0,p}}} \!\sum _{\ell '=0} ^{\ell } \cmicro_{d,d',p-p',\ell -\ell '} \cdot  \smacro_{p',d',\ell '} \bigg)_{\!p,d,\ell }\]
These sums have a finite number of terms, so they are well-defined. Next, we define the subspaces $\Vmacro_{\pm} = \{\smacro \in \Vmacro \mid \forall (p,d,\ell ),\; \smacro_{p,d,\ell } \neq 0 \Rightarrow  \pm d>0\}$ and likewise $\Vmicro_{\pm}$. So we have $\Vmacro = \Vmacro_{+} \oplus  \Vmacro_{-}$ and similarly for $\Vmicro$. In addition, we decompose $\Pmacro$ and $\Pmicro$ on those subspaces, which gives in block matrix notation:\vspace{-1mm}
\begingroup
\def\reglerHauteurCase{\rule[-1.5ex]{0pt}{4ex}}
\def\reglerHauteurAcco{\rule{0pt}{2.2ex}}
\[\Pmacro \,=\, \raisebox{-2.2ex}{$\begin{array}{@{}c@{}c}
	\left(\begin{array}{c}
		\reglerHauteurCase \Pmacroo{}+ \\ \hline
		\reglerHauteurCase \Pmacroo{}-
	\end{array}\right)
	\! & \begin{array}{@{}c@{}}
		\left. \reglerHauteurAcco \right\} \Vmacro_+ \\[1ex]
		\left. \reglerHauteurAcco \right\} \Vmacro_-
	\end{array}\\[-1.9ex]
	\underbrace{\rule{5.5ex}{0pt}}_{\textstyle\Vmicro} & 
\end{array}$}
\qquad \text{ and } \qquad 
\Pmicro \,=\, \raisebox{-2.2ex}{$\begin{array}{@{}c@{}c@{}c}
	\left(\begin{array}{c}
		\reglerHauteurCase \Pmicroo{}+ \\ \hline
		\reglerHauteurCase \Pmicroo{}-
	\end{array}\right) 
	\;=\;& \left(\begin{array}{c|c}
		\reglerHauteurCase \Pmicroo++ & \Pmicroo-+ \\ \hline
		\reglerHauteurCase \Pmicroo+- & \Pmicroo--
	\end{array}\right)
	\! & \begin{array}{@{}c@{}}
		\left. \reglerHauteurAcco \right\} \Vmicro_+ \\[1ex]
		\left. \reglerHauteurAcco \right\} \Vmicro_-
	\end{array}\\[-1.9ex]
	& \underbrace{\rule{5ex}{0pt}}_{\textstyle\Vmacro_+} \underbrace{\rule{5ex}{0pt}}_{\textstyle\Vmacro_-} & 
\end{array}$} \vspace{-1mm}\]
\endgroup
Finally, we denote $\sMacro := (\sigma _{d} (\bfu_{p,\ell } ))_{p,d,\ell }$, $\sMicro := (\sigma _{d} (S_{p,\ell } ))_{p,d,\ell }$ and $\sMacro_{\pm}$, $\sMicro_{\pm}$ their components on $\Vmacro_{\pm}$ and $\Vmicro_{\pm}$. With these notations, the two conditions of \cref{pre conditions de raccord} rewrite as $\sMacro =\Pmacro\,\sMicro$ and $\sMicro =\Pmicro\,\sMacro$, while the condition of the present theorem is:
\begin{equation}{\label{eq: condition de raccord preuve}}
\left\{\begin{array}{r@{\;=\;}l}
	\sMacro_{-} & \Pmacroo{}- \,\sMicro \\[0.5mm]
	\sMicro_{+} & \Pmicroo{}+ \,\sMacro
\end{array}\right.\end{equation}
Let us show that \cref{eq: condition de raccord preuve} implies one of the two conditions of \cref{pre conditions de raccord}, e.g.\ $\sMacro =\Pmacro\,\sMicro$.\\
Note that by \cref{pre conditions de raccord} we have for~any $\smacro \in \Vmacro$ and $\smicro \in \Vmicro$: $\smacro=\Pmacro\smicro \Leftrightarrow  \smicro=\Pmicro\smacro$. So $\Pmacro$ and $\Pmicro$ are inverses of each other. Hence we have $\Pmicro(\sMacro-\Pmacro\,\sMicro) = \Pmicro\,\sMacro-\sMicro$. Projecting this onto $\Vmicro_{+}$, we deduce:
\[\Pmicroo++(\sMacro_{+} -\Pmacroo{}+\,\sMicro) + \Pmicroo-+(\underbrace{\sMacro_{-} -\Pmacroo{}-\,\sMicro}_{=0} ) = \underbrace{\Pmicroo{}+\,\sMacro-\sMicro_{+} }_{=0} .\vspace{-5mm}\]
Let us show that $\Pmicroo++$ is injective.
\begin{sousdemo}
	Let $\smacro \in \Vmacro_{+} \setminus \{0\}$ and let us show that $\Pmicroo++ \smacro \neq 0$. Let $(p,d,\ell )$ be the smallest triplet for lexicographic order s.t.\ $\smacro_{p,d,\ell } \neq 0$. The term of $\Pmicroo++ \smacro$ of index $(p+d,d,\ell )$ is :\vspace{-2mm}
	\[\sum _{p'\in \intvl{0,p}} \sum _{\substack{d'\in \TN\\ p'+d'\in \intvl{0,p+d}}} \!\sum _{\ell '=0} ^{\ell } \cmicro_{d,d',p+d-p',\ell -\ell '} \cdot  \smacro_{p',d',\ell '} \,=\, \cmicro_{d,d,d,0} \cdot  \smacro_{p,d,\ell } .\vspace{-1mm}\]
	We claim that it~is non-zero because $\cmicro_{d,d,d,0} =1$. Indeed, in the sum \vspace{-2mm}
	\[\langle \Rmicro\rangle  := \fsumn \sum _{(R_{1} ,\dots ,R_{n} )\in \Rmicro^{n} } R_{1} \circ R_{2} \circ \dots \circ R_{n} \]
	the only term of degree $(0,0)$ is the identity, which appears for $n=0$. So the component of $\Eps\circ \langle \Rmicro\rangle  (\phi_{d} )$ in $\varepsilon ^{d} \mathcal{A}_{d} (\Pi )[\ln \varepsilon ]$ is $\varepsilon ^{d} \phi_{d}$. Thus we have $\cmicro_{d,d,d,0} = \sigma _{d,d,0} \circ  \Eps\circ \langle \Rmicro\rangle  (\phi_{d} ) = \sigma _{d} (\phi_{d} ) =1$. So~$\Pmicroo++ \smacro \neq 0$.
\end{sousdemo}
We have proven that $\sMacro_{+} =\Pmacroo{}+\.\sMicro$. Given that $\sMacro_{-} = \Pmacroo{}- \.\sMicro$, we deduce $\sMacro =\Pmacro\,\sMicro$.
\end{demo}

\begin{remarquenumerotee}
\label{rq: conditions de raccord}%
\begin{itemize}
	\item The coefficients $\cmacro_{d,d',p,\ell }$ and $\cmicro_{d,d',p,\ell }$ can be computed exactly and very quickly, as explained in \cref{sec: implementation des espaces A}. Moreover, they depend only on $\Theta $, $\omega $ and $(\mu _{0} ,\mu _{1} ,\rho _{0} ,\rho _{1} )$, but not on the shape of $\Omega _{1}$ at the corner, the values taken by the functions $\mu $ and $\rho $ at the corner, nor the source term $\fsource$.
	\item The matching conditions have a convolutive structure w.r.t.\ $p$ and $\ell $ that will be crucial to design effective models, see \cref{chap: ME}.
	\item In the sums of \cref{condition de raccord}, the indexes $d$ and $d'$ satisfy $d-d'\in \mathbb{Z}\cap \TZ$. Indeed, on the one hand we have $d,d'\in \TZ*$. On the other hand, we have in \cref{eq: condition de raccord 1} $p'\in \intvl{0,p+d}$ and $p'-d'\in \intvl{0,p}$, so $d-d' = (p+d-p')-(p-(p'-d'))\in \mathbb{N}-\mathbb{N}\subset \mathbb{Z}$ (and likewise in \cref{eq: condition de raccord 2}). Note that the set $\mathbb{Z}\cap \TZ$ can be very sparse. If $\Theta  \in \pi \mathbb{Q}$, then $\mathbb{Z}\cap \TZ =b\mathbb{Z}$ where $\Theta =\pi  \frac{a}{b}$ with $(a,b)\in \mathbb{N}\times \mathbb{N}^*$ and $\pgcd{a,b}=1$. Otherwise, we have $\mathbb{Z}\cap \TZ=\{0\}$.
\end{itemize}\fbi
\end{remarquenumerotee}

\begin{remarquenumerotee}
\label{rq: conjugaison de Rmaicro par Eps}%
The complexity of the matching conditions in \cref{condition de raccord} is due to the lack of a commutation property between $\langle \Rmaicro\rangle $ and $\Eps$. Let us explain what kind of commutation property one could hope for, how it would simplify everything and why it fails to be satisfied. We indicate hypothetical equalities with quote marks.\\

It is easy to verify that, for any $R\in \Rmacro$ and $\Phi \in \Aeps(\Pi )$, $R(\Phi )$ and $\varepsilon ^{\deg _{\mathcal{A}} \!R\,} \Eps^{-1} \circ R\circ \Eps(\Phi )$ are solutions of the same problem of the type \cref{eq: pb resolu dans A} (the proof is similar to \cref{action de Eps sur Helmholtz dans Aemaicro}). For instance, if $R = -\varepsilon ^2 \RcouchePi\circ \partial _{x|\couche} ^2 $, then they are solutions of
\begin{equation}{\label{eq: rq conjugaison de Rmaicro par Eps systeme d'eq}}
	\left\{\begin{array}{r@{\;=\;}ll}
	\Delta  \Psi  & 0 &\text{ in }\Omega  \\[0.5mm]
	\partial _{Y} ^2 \Psi  & -\varepsilon ^2  \partial _{x} ^2  \Phi  &\text{ in }\couche \\[0.5mm]
	\partial _{Y} \Psi _{|Y=0^{-}} & 0 &\text{ on }\Gamma  
	\end{array}\right.
\end{equation}
Hence one could hope to have ``$R=\varepsilon ^{\deg _{\mathcal{A}} \!R\,} \Eps^{-1} \circ R\circ \Eps$'' for any $R\in \Rmacro$ (which means that $R$ would be homogeneous, see \cref{def: op homogene dans A}). Since we have $\langle \{\varepsilon ^{\deg _{\mathcal{A}} \!R\,} \Eps^{-1} \circ R\circ \Eps \mid R\in \Rmacro\}\rangle  = \Eps^{-1} \circ {\langle \{\varepsilon ^{\deg _{\mathcal{A}} \!R\,} R \mid R\in \Rmacro\}\rangle }\circ \Eps$ and $\{\varepsilon ^{\deg _{\mathcal{A}} \!R\,} R \mid R\in \Rmacro\}=\Rmicro$ by construction of $\Rmaicro$ in \cref{eq: def Rmaicro}, this would imply ``$\langle \Rmacro\rangle  = \Eps^{-1} \circ \langle \Rmicro\rangle \circ \Eps$'', i.e.
\begin{equation}{\label{eq: hypothetique commutation entre Rmaicro et Eps}}
	\text{``}  \.\Eps\circ \langle \Rmacro\rangle  = \langle \Rmicro\rangle \circ \Eps \text{''} .
\end{equation}
Since we have $\projphi\circ \langle \Rmaicro\rangle =\projphi$ according to the proof of \cref{hyp raccord projetee sur les phid}, \cref{hyp raccord projetee sur les phid} could be simplified into $\fsum_{p,\ell } \powelne \fsum_{d} \sigma _{d} (\bfuDA_{p,\ell } )\,\phi_{d} = \Eps^{-1} \big(\fsum_{p,\ell } \powelne \fsum_{d} \sigma _{d} (\SDA_{p,\ell } )\,\phi_{d} \big)$, which is equivalent to simple one-to-one relations between the fields' coefficients:
\begin{equation}{\label{eq: hypothetique relation de raccord super simple}}
	\text{``} \. \forall (d,p,\ell ), \quad \ \sigma _{d} (\bfuDA_{p,\ell } )=\sigma _{d} (\SDA_{p+d,\ell } ) \.\text{''} .
\end{equation}
Unfortunately, the elements of $\Rmaicro$ are not homogeneous. Indeed, $R:=-\varepsilon ^2 \RcouchePi\circ \partial _{x|\couche} ^2 \in \Rmacro$ for instance picks the particular solution of \cref{eq: rq conjugaison de Rmaicro par Eps systeme d'eq} whose image by the linear forms $\sigma _{d,p,\ell }$ vanish, i.e.\ which belongs to $\operatorname{Ker}\projphi$. It~is easy to see that $\operatorname{Ker}\projphi$ is not stable under $\Eps^{\pm 1}$ (e.g.\ consider an element of $\Aeps(\Pi )$ equal to $\Im[(\alpha z)^{ \frac{\pi }{\Theta } } \log (\alpha z)]$ on $\Omega $). Therefore $\varepsilon ^{\deg _{\mathcal{A}} \!R\,} \Eps^{-1} \circ R\circ \Eps$ selects other solutions than $R$, and \cref{eq: hypothetique commutation entre Rmaicro et Eps} does~not holds.\\

Nevertheless, \cref{sec: simplifications DA} shows that \cref{eq: hypothetique relation de raccord super simple} holds in some particular simple cases, especially when $\frac{\Theta }{\pi }$ is irrational or $\frac{\Theta }{\pi }$ is rational and $p$ is smaller than its denominator.
\end{remarquenumerotee}

\section{Summary of the method}
\label{sec: resume du raccord}

\newcounter{etapesRaccord}
\setcounter{etapesRaccord}{0}
\def\nvEtapeRaccord[#1]{\refstepcounter{etapesRaccord}\item[Step \theetapesRaccord: #1.]}

\begin{textesansboite}
\def\highlight#1{#1}%
Let us summarize the main steps of our method for deriving the matching conditions, and explain how this approach can be extended to other settings. Quick examples are given later in this section and we aim to give more examples in future works.
\begin{description}
	\nvEtapeRaccord[the canonical problem]\label{etape raccord: pb canonique}
	Write the systems of equations satisfied the near fields (called here ``corner fields'') and the far fields (actually the ``far-and-layer fields'' in our case). In our case, they are \crefrange{eq: eq vol champs lointains}{eq: eq vol champs de coin}. Derive from these equations a \emph{canonical problem} that involves only differential operators that are homogeneous w.r.t.\ $\bfr$ (the variable that is scaled in the matched asymptotic expansion method). This problem is obtained by keeping only the differential operators that are in the left-hand side for both the far and near fields. For instance, we have in $\couche$
	\[\left\{\begin{array}{r@{\;=\;}ll}
		\mu _{1} \partial _{Y} ^2  U_{p,\ell } & -(\mu _{1} \partial _{x} ^2 +\omega ^2 \rho _{1} )U_{p-2,\ell } & \text{ for the far fields} \\[2pt]
		\mu _{1} \partial _{X} ^2  S_{p,\ell } + \mu _{1} \partial _{Y} ^2  S_{p,\ell } & -\omega ^2 \rho  \,S_{p-2,\ell } & \text{ for the near fields} 
	\end{array}\right.\]
	so we keep only $\partial _{Y} ^2 $ in $\layer$. This operator is homogeneous (of degree 0) w.r.t.\ $\bfr$, contrary to $\partial _{X} ^2  + \partial _{Y} ^2 $ for instance. In our case, we obtain the problem \cref{eq: pb resolu dans A}, that we recall:
	\let\inc\varphi
	\let\sm\psi
	\[\left\{\begin{array}{r@{\;=\;}ll}
		\Delta  \inc & \sm_{\Omega } &\text{ in }\Omega  \\[0.5mm]
		\partial _{Y} ^2 \inc & \sm_{\couche} &\text{ in }\couche \\[0.5mm]
		\partial _{Y} \inc_{|Y=0^{-}} & \sm_{\Gamma } &\text{ on }\Gamma  \\[0.5mm]
		\inc_{|y=0^{+}} - \inc_{|Y=0^{-}} & 0 &\text{ on }\Gamma  \\[0.5mm]
		\inc & 0 &\text{ on }\bordgauche \cup \borddroit 
	\end{array}\right.\]
	Note that, depending on the problem satisfied by $u_{\peps}$, it may not be possible to derive this canonical problem that is common to the far and near fields (e.g.\ in some special high contrast settings). In such cases, the method must be adapted (which we leave for future work).

	\nvEtapeRaccord[the functions $\phi_{d}$]\label{etape raccord: phid}
	Let $\Pi $ be the scaling-invariant version of the geometry. Find the \highlight{homogeneous solutions} of the canonical problem in $\Pi $. Decompose these solutions on a basis made of quasi-homogenous functions (i.e.\ homogeneous up to some $\ln r$). In our case, they are the functions $\phi_{d}$, see \cref{pb de Laplace homogene dans A}.

	\nvEtapeRaccord[the space $\mathcal{A}$]\label{etape raccord: A}
	Build a \highlight{space of singularities} defined on $\Pi $, denoted $\mathcal{A}$, s.t.\ the following holds.
	\begin{enumerate}[label={(\roman*)},topsep=1pt]
		\item The elements of $\mathcal{A}$ are finite linear combinations of quasi-homogeneous functions.
		\item $\mathcal{A}$ contains the functions $\phi_{d}$.
		\item When the right-hand sides of the canonical problem belong to $\mathcal{A}$, $\mathcal{A}$ contains a particular solution of the problem. 
	\end{enumerate}
	Define operators that map these right-hand sides to these particular solutions, similarly to $\ROmegaPi$, $\RcouchePi$ and $\RGammaPi$ (see \cref{operateurs R}). In addition, for~any $d$, define on $\mathcal{A}$ a linear form $\sigma _{d}$ s.t.\ we have $\sigma _{d} (\phi_{q} )= \delta _{d,q}$ for~any $d$ and $q$. Note that, to make the computation of the matching coefficients easier, it~is better that $\mathcal{A}$ has a basis made of explicit functions (see \cref{sec: implementation des espaces A}).

	\nvEtapeRaccord[expression of the fields' behaviors]\label{etape raccord: champs tot}
	Using \cref{inverse de id-f dans les SF 2D}, express the formal series of the far and near fields behaviors w.r.t.\ $\bfr$, like $\fsum_{p,\ell } \powelne\,\bfuDA_{p,\ell }$, in the form
	\[\langle \mathcal{R}\rangle \bigg(\fsum_{p,\ell } \powelne \fsum_{d} \sigma _{d} (\bfuDA_{p,\ell } ) \,\phi_{d} \bigg)\]
	where $\mathcal{R}$ is a set of operators built from the resolvant operators of step \ref*{etape raccord: A} and \smash{$\langle \mathcal{R}\rangle :=\fsum\limits_{n=0} ^{\infty } \Big(\sum\limits _{R\in \mathcal{R}} R\Big)^{n}$}. Here this is done in \cref{champs totaux en fct de sigma}.

	\nvEtapeRaccord[conclusion]\label{etape raccord: conclusion}
	Project the matching assumption on the functions $\phi_{d}$ (like in \cref{hyp raccord projetee sur les phid}). This should not cause a loss of information as the fields' behaviors are entirely determined by their coefficients $\sigma _{d}$. Then identify each power of $\varepsilon $ and $\ln \varepsilon $ and apply $\sigma _{d}$ to obtain scalar relations on the coefficients $\sigma _{d}$ (the matching conditions), see \cref{pre conditions de raccord}. Finally, keep only the relations that are useful to build the fields (which should neither cause a loss of information, like in \cref{condition de raccord}).
\end{description}
\end{textesansboite}

\begin{remarque}
All the specificity of the studied problem (whether Helmholtz, Poisson, Maxwell, etc.)\ is localized in steps~1 to \ref*{etape raccord: A}. Once these first steps are completed, the final reasoning of steps~\ref*{etape raccord: champs tot} and \ref*{etape raccord: conclusion} is free of cumbersome specific calculations and depend little on the studied setting. This has two main advantages.
\begin{itemize}
	\item First, adapting the method to other situations is of course easier. Since step \ref*{etape raccord: pb canonique} is rather quick, most of the work needed to apply the method to another setting lies in steps~\ref*{etape raccord: phid} and \ref*{etape raccord: A}.
	\item Secondly, the calculatory effort is reduced, while being able to make the matching coefficients more explicit. Indeed, instead of having to deal from start to end with the specific construction of the singularities (the equations they satisfy, their inductive construction, etc.), the amount of case-dependent calculations is reduced to its minimum (mainly in the construction of the resolvant operators like $\ROmegaPi$, $\RcouchePi$ and $\RGammaPi$). Then each manipulation of a high-level object like $\langle \Rmaicro\rangle $ implies complex operations on singularities behind the scenes, that we do~not have to worry about (note that the implementation of the matching coefficients computation has a similar architecture, see \cref{sec: implementation des espaces A}). This allows us to consider singularities that have analytic expressions, leading to matching coefficients computed exactly, without causing an explosion of the quantity of calculations.
\end{itemize}
\end{remarque}

\begin{exemplenumerote}
\label{eq: generalisation a une double couche}%
It is easy to extend our work to the case of a corner covered by a layer on both sides, see \cref{fig: double couche}. Steps~\ref*{etape raccord: pb canonique} and \ref*{etape raccord: phid} are straightforward to adapt, while \cref{rq: espaces A pour 2 couches} explains how to modify the spaces $\mathcal{A}$ in this case. Then, in addition to $\ROmegaPi$, $\RcouchePi$ and $\RGammaPi$, similar operators $\RcouchePi'$ and $\RGammaPi'$ should be defined for the second layer. As~for step~\ref*{etape raccord: champs tot}, it involves sets of operators $\smash{ \widetilde{\mathcal{R}} ^{\pm } _{\varepsilon } }$ that are defined like $\Rmaicro$, but with duplicating the operators that involve $\RcouchePi$ and $\RGammaPi$ using $\RcouchePi'$ and $\RGammaPi'$. Finally, step~\ref*{etape raccord: conclusion} stays identical.
\end{exemplenumerote}

\begin{figure}[h]
\def\epsFigA{0.7} 
\def\epsFigB{0.75} 
\def\widthFig{3.5} 
\def\heightFig{2.3} 
\def\LFig{5} 
\def\ThetaFig{215} 
\def\xI{cot(\ThetaFig)*(-\epsFigA-\epsFigB*cos(\ThetaFig)) - \epsFigB*sin(\ThetaFig)}
\def\rLegEpsFigB{2.5} 
\begin{tikzpicture}[scale=0.9]
	\begin{scope}
		\clip (-\widthFig,-\heightFig) rectangle (\widthFig,\heightFig);
		\coordinate (finLambda) at ({\LFig*cos(\ThetaFig)}, {\LFig*sin(\ThetaFig)});
		\coordinate (finLambdaCoucheB) at ({\LFig*cos(\ThetaFig) - \epsFigB*sin(\ThetaFig)}, {\LFig*sin(\ThetaFig) + \epsFigB*cos(\ThetaFig)});
		\coordinate (coinCoucheAB) at ({\xI},-\epsFigA);
		\coordinate (finCoucheA) at ({\xI},0);
		\coordinate (finCoucheB) at ({\xI+\epsFigB*sin(\ThetaFig)},{-\epsFigA-\epsFigB*cos(\ThetaFig)});
		\fill [pattern={Lines[angle=45,line width=.2pt]}, pattern color=white!70!gray] (0,0) -- (\LFig,0) -- (\LFig,\LFig) -- (-\LFig,\LFig) -- (finLambda) -- cycle;
		\fill [color=white!70!gray] (finCoucheA) -- (coinCoucheAB) -- (\widthFig,-\epsFigA) -- (\widthFig,0) -- cycle;
		\fill [color=white!53!gray] (finCoucheB) -- (coinCoucheAB) -- (finLambdaCoucheB) -- (finLambda) -- cycle;
		\fill [color=white!35!gray] (0,0) -- (finCoucheA) -- (coinCoucheAB) -- (finCoucheB) -- cycle;
		\draw [thick] (finLambda) -- (0,0) -- (\widthFig,0);
		\draw [thick] (finLambdaCoucheB) -- (coinCoucheAB) -- (finCoucheA);
		\draw [thick] (finCoucheB) -- (coinCoucheAB) -- (\widthFig,-\epsFigA);
	\end{scope}
	\draw (-\widthFig-1.6,0) node {\ }; 
	\draw [>=stealth,->] (-\widthFig-0.4,0) -- (\widthFig+0.4,0) node [right] {$x$};
	\draw [>=stealth,->] (0,{-\heightFig-0.3}) -- (0,{\heightFig+0.3}) node [right] {$y$};
	\draw (\widthFig*0.4,{-\heightFig+0.1}) node {obstacle};
	\draw [->,>=latex] (0.35,0) arc (0:\ThetaFig-5:0.35) node[near start,above] {$\Theta$} -- (\ThetaFig:0.36);
	\draw [<->,>=latex] (\widthFig*0.85,0) -- node[midway,right]{$\peps$} (\widthFig*0.85,-\epsFigA);
	\draw [<->,>=latex] ({\rLegEpsFigB*cos(\ThetaFig)}, {\rLegEpsFigB*sin(\ThetaFig)}) -- ({\rLegEpsFigB*cos(\ThetaFig) - \epsFigB*sin(\ThetaFig)}, {\rLegEpsFigB*sin(\ThetaFig) + \epsFigB*cos(\ThetaFig)});
	\draw ({(\rLegEpsFigB+0.3)*cos(\ThetaFig) - \epsFigB/2*sin(\ThetaFig)}, {(\rLegEpsFigB+0.3)*sin(\ThetaFig) + \epsFigB/2*cos(\ThetaFig)}) node{$c\.\peps$};
	\draw (-\widthFig*0.55,\heightFig*0.6) node {$(\mu_0,\rho_0)$};
	\draw [<-,>=latex] (\widthFig*2/3,-\epsFigA*0.6) to[bend right=20] ++(0.2,-0.6) node[below]{$\ \ (\mu_1,\rho_1)$};
	\draw [<-,>=latex] ({\rLegEpsFigB/2*cos(\ThetaFig) - 0.4*\epsFigB*sin(\ThetaFig)}, {\rLegEpsFigB/2*sin(\ThetaFig) + 0.4*\epsFigB*cos(\ThetaFig)}) to[bend right=20] ++(-0.5,0.4) node[left]{$(\mu_2,\rho_2)$};
\end{tikzpicture}\par
\longCaption{An example of configuration with a layer on both sides of $\Omega$}{ ($\mu_2$, $\rho_2$ and $c$ are constants independent of $\peps$)}
\label{fig: double couche}
\end{figure}

	\tikzset{external/export=true}%
	\tikzsetfigurename{fig-couche_infinie_periodique-}%
	\begin{figure}[h]
\def\epsFig{0.8} 
\def\TFig{1} 
\def\widthFig{3} 
\def\heightFig{2.7} 
\newcommand{\cellulePeriodicite}[1]{
	\begin{scope}[shift={({#1*\epsFig*\TFig},0)}, xscale={\epsFig*\TFig}, yscale={-\epsFig}]
		\fill [color=white!60!gray] plot [smooth cycle] coordinates {(0.1,0.2) (0.3,0.2) (0.5,0.1) (0.7,0.3) (0.8,0.5) (0.7,0.9) (0.4,0.8) (0.3,0.9) (0.2,0.7) (0.1,0.5)} ;
		\fill [color=white!40!gray] plot [smooth cycle] coordinates {(0.3,0.25) (0.5,0.2) (0.6,0.3) (0.7,0.5) (0.7,0.8) (0.4,0.7) (0.3,0.7) (0.2,0.5)} ;
		\fill [color=white!20!gray] plot [smooth cycle] coordinates {(0.3,0.4) (0.5,0.3) (0.65,0.65) (0.4,0.6) (0.3,0.6) (0.3,0.5)} ;
	\end{scope}}
\begin{tikzpicture}
	\begin{scope}
		\clip (-\widthFig,0) rectangle (\widthFig,-\epsFig) ; 
		\fill [color=white!75!gray] (-\widthFig,0) rectangle (\widthFig,-\epsFig) ;
		\foreach \i in {-5,-4,...,4} {\cellulePeriodicite{\i}}
	\end{scope}
	\fill [pattern={Lines[angle=45,line width=.2pt]}, pattern color=white!70!gray] (-\widthFig,0) rectangle (\widthFig,\heightFig) ;
	\draw [thick] (-\widthFig,0) -- (\widthFig,0) ;
	\draw [thick] (-\widthFig,-\epsFig) -- (\widthFig,-\epsFig) ;
	\draw [>=stealth,->] (-\widthFig-0.4,0) -- (\widthFig+0.5,0) node [right] {$x$} ;
	\draw [>=stealth,->] (0,{-\epsFig*3/2-0.4}) -- (0,\heightFig+0.3) node [above] {$y$} ;
	\draw (-\widthFig/2,-\epsFig*3/2) node {obstacle\ } ;
	\draw  [<-,>=latex] (-\widthFig*0.54,-\epsFig/2) to[bend left=10] ++(0.1,{\epsFig/2+0.5}) node [above] {\ layer} ;
	\draw [decorate, decoration={brace,amplitude=10pt}, xshift=-4pt] (-\widthFig-0.3,-\epsFig) -- (-\widthFig-0.3,\heightFig) node [midway,xshift=-0.6cm] {$\widetilde\Omega_\peps$} ;
	\draw [<-,>=latex] (\widthFig/3,0) to[bend left=20] ++(0.2,0.6) node[above]{$\ \ \ \widetilde\Gamma$} ;
	\draw [<-,>=latex] (\widthFig/4,-\epsFig) to[bend right=30] ++(0.3,-0.6) node[right]{Dirichlet condition} ;
	\draw [<->,>=latex] ({\epsFig*\TFig*2.9},0) -- node[midway,right]{$\peps$} ({\epsFig*\TFig*2.9},-\epsFig) ;
	\draw [<->,>=latex] ({\epsFig*\TFig*1.9},0.1) -- node[midway,above]{$\peps$} ({\epsFig*\TFig*2.9},0.1) ;
\end{tikzpicture}

\caption{The infinite periodic layer}
\label{fig: couche infinie periodique}
\end{figure}
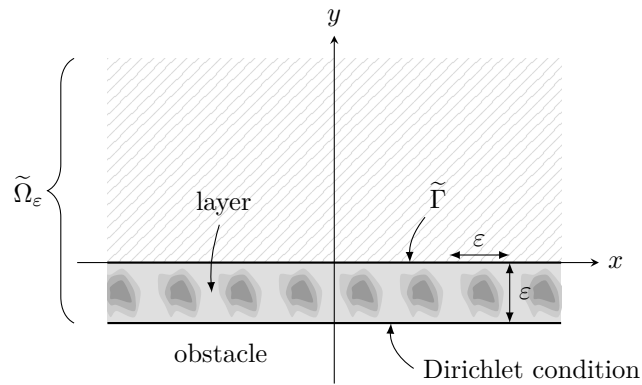%
	\tikzset{external/export=false}%
	\tikzsetfigurename{main-figure}

\begin{exemple}
\def\mci{\tilde{\mu}}%
\def\rci{\tilde{\rho}}%
\def\OCI{\widetilde{\Omega}}%
\def\GCI{\widetilde{\Gamma}}%
\def\fsourceCI{\tilde{f}_{\rm s}}%
\def\bfX{\mathrm{X}}%
\def\RCI{R_{\partial_y^2}}%
\def\RmaicroCI{\widetilde{\mathcal R}^\pm_{\feps}}%
\def\RmacroCI{\widetilde{\mathcal R}^+_{\feps}}%
\def\RmicroCI{\widetilde{\mathcal R}^-_{\feps}}%
Let us re-express the case of an infinite periodic layer using our method. This case is well-known in the literature, see \cite{AmmLat98, DelHadJol12, Bou23}. This example illustrates the generality our method and allows a comparison between the tools we used and what already existed in the literature for a simpler case.\\

\let\varepsilon=\peps
We consider the domains $\widetilde{\Omega }_{\varepsilon } :=\mathbb{R}\times (-\varepsilon ,\infty )$ and $\GCI := \mathbb{R}\times \{0\}$, and two functions $\mci,\rci : \widetilde{\Omega }_{1} \rightarrow  \mathbb{R}$ that are 1-periodic w.r.t.\ the first variable and equal to $\mu _{0} ,\rho _{0}$ outside of the layer $(-\varepsilon ,0)\times \mathbb{R}$. We denote $\mci_{\varepsilon } :(x,y) \mapsto  \mci( \frac{x}{\varepsilon } , \frac{y}{\varepsilon } )$ and similarly for $\rci_{\varepsilon }$. See \cref{fig: couche infinie periodique}. We consider a function $\tilde{u}_{\varepsilon } \in H^{1} _{\mathrm{loc}} ( \widetilde{\Omega }_{\varepsilon } )$ that satisfies
\[\operatorname{div}( \mci_{\varepsilon } \nabla  \tilde{u}_{\varepsilon } ) + \omega ^2  \rci_{\varepsilon } \tilde{u}_{\varepsilon } = \fsourceCI \quad \text{ in }\widetilde{\Omega }_{\varepsilon } \]
with $\fsourceCI \in H^{-1}$ s.t.\ $\operatorname{dist}(\operatorname{supp} \fsourceCI,\GCI)>0$, and e.g.\ $\tilde{u}_{\varepsilon } =0$ on $\partial  \widetilde{\Omega }_{1}$. We make the ansatz that the asymptotic expansion of $\tilde{u}_{\varepsilon }$ has the following form:
\begin{itemize}
	\item ``$\sum_{n=0}^\infty  \varepsilon ^{n} \.u_{n} (x,y)$'' far form the layer ($y\geqslant \sqrt{\varepsilon }$)
	\item ``$\sum_{n=0}^\infty  \varepsilon ^{n} \.U_{n} (x, \frac{x}{\varepsilon } , \frac{y}{\varepsilon } )$'' near the layer ($y\leqslant 2\sqrt{\varepsilon }$), where each $U_{n}$ is periodic w.r.t.\ its second variable.
\end{itemize}
Let us apply the steps presented above to this case.
\begin{description}
	\item[Step \ref*{etape raccord: pb canonique}:] Let us denote the microscopic variables $\frac{x}{\varepsilon } , \frac{y}{\varepsilon }$ as $X,Y$. In the matching zone ($\sqrt{\varepsilon } \leqslant y\leqslant 2\sqrt{\varepsilon }$), the far and near fields satisfy the equations:
	\begin{itemize}[topsep=0pt]
		\item $\mu _{0} \Delta  u_{n} +\omega ^2 \rho _{0} u_{n} = \delta _{n,0} \,\fsourceCI$
		\item and $- \mu _{0} \partial _{Y} ^2 U_{n} = (\mu _{0} \partial _{x} ^2 +\omega ^2  \rho _{0} )U_{n-2} $ using that $\mci$ and $\rci$ are constant in this zone, and neglecting the dependence of $U_{n}$ w.r.t.\ $X$ when $Y\rightarrow \infty $ (justified by the matching assumption and the fact that the far fields do~not depend on $X$).
	\end{itemize}
	The differential operator that always stays in the left-hand side is $\partial _{y} ^2 $ (denoted $\partial _{Y} ^2 $ for the near fields), so the canonical problem is
	\begin{equation}{\label{eq: pb canonique couche infinie}}
		\partial _{y} ^2  \varphi  = \psi  \qquad \text{ on }\mathbb{R}\times (0,\infty ).
	\end{equation}

	\item[Step \ref*{etape raccord: phid}:] The homogeneous solutions of \cref{eq: pb canonique couche infinie} are the functions of the form $(x,y) \mapsto  v_{0} (x)+v_{1} (x)\cdot y$ with $v_{0} ,v_{1} :\mathbb{R} \rightarrow \mathbb{C}$. They are generated by the two $y$-homogeneous functions $\phi_{0} =1$ and $\phi_{1} =y$ (which, in this case, are tensorized by functions of $x$).
	
	\item[Step \ref*{etape raccord: A}:] One can take e.g.
	\[\mathcal{A}:=\bigg\{(x,y) \mapsto  \sum_{j=0}^n v_{j} (x)\, y^{j} \ \bigg|\ n\in \mathbb{N},\; (v_{j} )_{j} \in \mathcal{C}^{\infty } (\mathbb{R})^{n+1} \bigg\}\]
	The $\mathcal{C}^{\infty }$ regularity of the functions $v_{j}$ is~not a problem as one can verify that the far and near fields are smooth in the matching zone. It can be decomposed as $\mathcal{A}=\bigoplus_{d\in \mathbb{R}} \mathcal{A}_{d}$ where $\mathcal{A}_{d} := \{v(x)\.y^{d} \mid v\in \mathcal{C}^{\infty } (\mathbb{R})\}$ if $d\in \mathbb{N}$ and $\mathcal{A}_{d} :=\{0\}$ otherwise, and each $\mathcal{A}_{d}$ only contains functions that are homogeneous w.r.t.\ $y$. Then, for~any $d\in \{0,1\}$, we define $\sigma _{d} :\mathcal{A} \rightarrow \mathcal{C}^{\infty } (\mathbb{R})$ by $\sigma _{d} \big(\sum _{j} v_{j} (x)\.y^{j} \big):=v_{d}$. In addition, we define $\RCI:\mathcal{A}\rightarrow \mathcal{A}$ by $\RCI\big(\sum _{j} v_{j} (x)\.y^{j} \big):=\sum _{j} v_{j} (x)\. \frac{y^{j+2} }{(j+1)(j+2)}$. It maps any $\psi \in \mathcal{A}$ to a particular solution $\varphi $ of \cref{eq: pb canonique couche infinie} s.t.\ $\sigma _{0} (\varphi )=\sigma _{1} (\varphi )=0$.
	
	\let\varepsilon=\feps
	\item[Step \ref*{etape raccord: champs tot}:] Let $u^{0} _{n}$ denote the behavior of $u_{n}$ when $y\rightarrow 0$, and $U^{\infty } _{n}$ the behavior of $U_{n}$ when $Y\rightarrow \infty $. Let $\RmacroCI :=\{-\RCI\circ (\partial _{x} ^2 +k_{0} ^2 )\}$ and $\RmicroCI :=\{-\varepsilon ^2  \RCI\circ (\partial _{x} ^2 +k_{0} ^2 )\}$. One can show
	\[\!\fsumn \varepsilon ^{n} u^{0} _{n} = \langle \RmacroCI\rangle \bigg(\fsumn \varepsilon ^{n} \!\! \sum _{d\in \{0,1\}} \! \sigma _{d} (u^{0} _{n} )\,\phi_{d} \bigg) 
	\quad \  \text{ and } \quad \ \fsumn \varepsilon ^{n} U^{\infty } _{n} = \langle \RmicroCI\rangle \bigg(\fsumn \varepsilon ^{n} \!\! \sum _{d\in \{0,1\}} \! \sigma _{d} (U^{\infty } _{n} )\,\phi_{d} \bigg).\]

	\item[Step \ref*{etape raccord: conclusion}:] One can check that, in this case, $\RCI$ is a homogeneous operator because there~is neither $\ln y$ nor $\ln \varepsilon $. Hence, by the reasoning of \cref{rq: conjugaison de Rmaicro par Eps}, we obtain:
	\[\forall n\in \mathbb{N},\ \forall d\in \{0,1\}, \quad \ \sigma _{d} (u^{0} _{n} ) = \sigma _{d} (U^{\infty } _{n+d} ).\]
	In practice, $\sigma _{1} (U^{\infty } _{n} )=\sigma _{1} (u^{0} _{n-1} )=\partial _{y} u_{n-1|\Gamma }$ determines the behavior at infinity of the near fields and $u_{n|\Gamma } =\sigma _{0} (u^{0} _{n} )= \sigma _{0} (U^{\infty } _{n} )$ provides a Dirichlet condition for the far fields.
\end{description}
\end{exemple}

\begin{textesansboite}
We will investigate other generalizations in future works, especially the Neumann boundary condition, the case of an edge in three dimensions and Maxwell's equations.
\end{textesansboite}

\section{Numerical computation of the fields behaviors and the matching coefficients}
\label{sec: implementation des espaces A}

\begin{textesansboite}
The matching coefficients $\cmacro_{d,d',p,\ell }$ and $\cmicro_{d,d',p,\ell }$ depend only on $(\Theta ,\omega ,\mu _{0} ,\mu _{1} ,\rho _{0} ,\rho _{1} )$ and this dependence is explicit. However this dependence is very indirect: the matching coefficients depend on $\langle \Rmaicro\rangle $, which is a formal series of iterated compositions of the elements of $\Rmaicro$, which depend on $(k_{0} ,k_{1} ,\mu _{0} ,\mu _{1} )$ and the operators $\ROmegaPi$, $\RcouchePi$ and $\RGammaPi$, which themselves involve $\Theta $ and corner singularities with analytic expressions. Doing all these calculations step by step by hand would be extremely laborious, as well as providing formulas of the matching coefficients in terms of $(\Theta ,\omega ,\mu _{0} ,\mu _{1} ,\rho _{0} ,\rho _{1} )$ (except in some particular simple cases that are studied in \cref{chap: complements sur le DA}). Thus we delegate this task to a computer. To this end, we developed a \emph{symbolic computation module} specifically designed to handle the objects introduced in this chapter. It can be implemented using any coding language that allows object-oriented programming. We used \Cplusplus{}, which is particularly appropriate for this type of programming paradigm.\\

The principle of object-oriented programming is to define custom variable types, called \emph{classes}. The data of a class' object is given by its \emph{attributes} (that are variables), while its behavior is determined by the \emph{methods} of the class (that are functions). Usually the attributes are hidden from outside the class so that only the class' behavior matters for the user, not its internal implementation. We used classes to build numerical representations of mathematical objects: each class is associated to a mathematical set like $\mathcal{A}(\Omega )$, while each object of the class represents an element of this set. One of the key advantages of \Cplusplus{} is that it allows to define the operators \texttt{+}, \texttt{-}, \texttt{*} and \texttt{/} for any class, like:
\begin{lstlisting}
	MyClass operator+(MyClass a1, MyClass a2) {/* implementation... */};
	// Now I can write "a1+a2" in my code
\end{lstlisting}
This allows to then almost copy-paste in the code mathematical formulas that involve e.g.\ elements of the spaces $\mathcal{A}$.\\

To begin, any implementation of this chapter's tools should take special care of the sensitivity w.r.t.\ $d$ of the definition of $\ROmegaOmega$, $\RGammaOmega$, $\ROmegaPi$, $\RcouchePi$ and $\RGammaPi$ on $\mathcal{A}_{d}$. For instance, the expression of $\RGammaOmega(x^{d} )$ is not continuous w.r.t.\ $d$, as we have
\[\forall d\in \mathbb{R} \setminus \txtTZ, \qquad  \RGammaOmega(x^{d} ) = \frac{r^{d} \sin (d(\theta -\Theta ))}{\sin (-d\Theta )} = \frac{\Im[(\alpha z)^{d} ]}{\sin (-d\Theta )}\]
while $\RGammaOmega(x^{d} )$ involves $\log (\alpha z)$ when $d\in \TZ$. The definition of $\ROmegaOmega$ and $\ROmegaPi$ are also sensitive to whether $d$ is in $\mathbb{Z}$ or~not. Handling the powers $d$ with floating-point numbers can lead to situations where $d$ theoretically belongs to $\TZ$ or $\mathbb{Z}$, but numerically does~not because of round-off errors. This problem can be solved by testing whether $d\in \TZ\cup \mathbb{Z}$ with a small tolerance error. However, there~is no guarantee of correctness as a too large, resp.\ too small, tolerance can cause false positives, resp.\ false negatives.\\

To address this issue neatly, we introduced two classes denoted \texttt{Ratio} and \texttt{DValue}. \texttt{Ratio}'s objects represent rational numbers using two integer attributes: the numerator and the denominator. We defined on \texttt{Ratio} the four operators \texttt{+}, \texttt{-}, \texttt{*} and \texttt{/}. Then \texttt{DValue}'s objects represent real numbers in the form
\[\qquad  \quad  r_{1} +r_{2} \cdot z^{s} \qquad  \text{where }  (r_{1} ,r_{2} ,z,s) \in \mathbb{Q}^2 \times (\mathbb{R} \setminus \mathbb{Q})\times \{-1,1\}.\]
Thus \texttt{DValue} has four attributes: two \texttt{Ratio} objects \texttt{r1}, \texttt{r2}, a floating-point number \texttt{z} and a boolean \texttt{inv} that describes $s$. We defined the operators \texttt{+}, \texttt{-}, \texttt{*} and \texttt{/} so that they do only exact operations on $(r_{1} ,r_{2} ,s)$ when possible and, otherwise, evaluate everything as a floating-point number and store the result in $z$. 
For instance, the code below is treated as explained in \cref{tab: calculs dans NbAlg}. Note that the handling of the two cases $\pst \in \mathbb{Q}$ and $\pst \not\in \mathbb{Q}$ is entirely transparent to the user.
\begin{lstlisting}
	DValue pi = 3.141592653589793;
	DValue Theta = pi*Ratio(2,3); // pi*2/3
	// other possibility: "DValue Theta = 1.5;" gives an angle of 1.5 rad
	DValue t = pi/Theta;
	DValue d1 = 4 + 3*t;
	DValue d2 = -2*t;
	bool b1 = (d1/t).isInt(); // false (isInt is a method of the class DValue)
	bool b2 = (d2/t).isInt(); // true
\end{lstlisting}

\begin{table}[h]
\renewcommand{\arraystretch}{1.35}%
\def\r{r}%
\longCaption{Values of $(r_{1} ,r_{2} ,z,s)$ used to represent a few relevant real numbers as $\texttt{DValue}$ objects}{. ``$*$'' indicates an attribute that can take any value.}
\label{tab: calculs dans NbAlg}
\begin{tabular}{>{\kern2pt}cc>{\;}c<{\;}}
    \cline{2-3}
    & case $\Theta  =\r \pi $, $\r \in \mathbb{Q}$ & case $\pst \not\in \mathbb{Q}$ \\
	\hline
	$\pi $ & $(0,1,\pi ,1)$ & $(0,1,\pi ,1)$ \\
	\hline
	$\Theta $ & $(0,\r,\pi ,1)$ & $(0,1,\Theta ,1)$\footnote{Or $(\Theta ,0,*,*)$ if $\Theta $ is an integer. In this case, $\pst$ is represented as $(0,\Theta ^{-1} ,\pi ,1)$. We do~not go into details but this case is correctly handled as well.}\!\!\\
	\hline
	$\pst$ & $(\r^{-1} ,0,*,*)$ & $(0,1,\pst,1)$\\
	\hline
	$d_{1} :=3+4\pst$ & $(d_{1} ,0,*,*)$ & $(3,4,\pst,1)$\\
	\hline
	$d_{1} /(\pst)$ & $(d_{1} \r,0,*,*)$ & $(4,3,\pst,-1)$\\
	\hline
	$(-2\pst) /(\pst)$ & $(-2,0,*,*)$ & $(-2,0,*,*)$\\
	\hline
\end{tabular}\smallskip
\end{table}
\\

Next, we defined four classes \texttt{AOmega}, \texttt{AGamma}, \texttt{ALayer} and \texttt{APi} that represent resp.\ $\mathcal{A}(\Omega )$, $\mathcal{A}(\Gamma )$, $\mathcal{A}(\couche)$ and $\mathcal{A}(\Pi )$. \texttt{APi} has two attributes of type \texttt{AOmega} and \texttt{ALayer} resp.\ that describe the restriction of an element of $\mathcal{A}(\Pi )$ to $\Omega $ and $\couche$. As~for the first three classes, they possess a single attribute that contains the coefficients of an element of $\mathcal{A}({\dots })$ on the basis of $\mathcal{A}$ given by \cref{base de A Omega,rq: base de A couche et Gamma}. For instance, the attribute of \texttt{AGamma} is:
\begin{lstlisting}
	map<DValue,vector<lcomplex>> coefs;
	// for any (d,k), coefs[d][k] is the coefficient in front of x^d ln(x)^k
\end{lstlisting}
Here, \texttt{lcomplex} is an alias of \texttt{complex<long double>}, which is a type of complex number whose real and imaginary parts are floating-point numbers with (on most computers) twice as much digits as usual ones (of type \texttt{double}). This improves the accuracy of computations, as round-off errors stay below the precision level of usual \texttt{double} numbers. Actually the attribute \texttt{z} of \texttt{DValue} also has type \texttt{long double}. Computations on those classes are so fast and memory-thrifty that the increased computational cost of using \texttt{long double} is~not a problem.\\

We defined on these classes the operators \texttt{+}, \texttt{-} and \texttt{*} (multiplication by a scalar), as well as the relevant differential operators and $\ROmegaOmega$, $\RGammaOmega$, $\ROmegaPi$, $\RcouchePi$ and $\RGammaPi$. These operators are implemented using their explicit expressions on the basis of the spaces $\mathcal{A}$ given by the proofs of \cref{operateurs R dans Omega,operateurs R,operateurs de derivation dans A}. They involve only simple operations on small collections of numbers, so these computations are very fast. Here is a quick example on how to use those new classes.
\begin{lstlisting}
	AOmega phi = AOmega::basisFunction(2*pi/Theta, 0, 1);
	//     phi = Im[(alpha z)^{2pi/Theta} log(alpha z)]
	AGamma dyPhi = phi.dyGamma(); // Neumann trace on Gamma
	APi psi = 3 * phi.Rlap(); // 3 R_Delta(phi)
	lcomplex s = psi.sigma(pi/Theta); // sigma_{pi/Theta}(psi) = 0
\end{lstlisting}\\

The next step is to compute formal series of functions of $\mathcal{A}$. Since a computer cannot store data of infinite size, we can only compute truncations of formal series. This is~not a problem in general because the asymptotic expansion of a field (e.g.\ $\bfu_{p,\ell }$) at a given order $q$ depends only on the truncation at order $q$ of the formal series associated to its asymptotic behavior (e.g.\ $T_{\leqslant q} (\bfuDA_{p,\ell } )$). In addition, although the definition of the matching coefficients $\cmacro$ and $\cmicro$ involve the formal series of operators $\langle \Rmaicro\rangle $, each matching coefficient depends only on a truncation of this series (more details in \cref{calcul numerique cmaicro} below).\\

Applying an operator $f$ (that has a degree or is of the form $\langle \dots \rangle $) to a formal series $\Phi $ truncated at order $q$ gives knowledge about a trunction of $f(\Phi )$ at some order $q'$. The rules to determine $q'$ from $q$ (or inversely to choose $q$ depending on a desired truncation order $q'$) are given by \cref{calcul op a degre sur SF tronquee,calcul op <F> sur SF tronquee} below. The first lemma is a straightforward consequence of \cref{def: extension d'un op a degre aux SF 1D,def: troncature de SF}.
\end{textesansboite}

\begin{lemmesansboite}{}
\label{calcul op a degre sur SF tronquee}%
Let $D_{1} ,D_{2} \in \domsA$ and let $f:\mathcal{A}(D_{1} )\rightarrow \mathcal{A}(D_{2} )$ be a linear map that has a degree. Let $d:=\deg f$ and $q\in \mathbb{R}$. We have:
\begin{itemize}
	\item $T_{\leqslant q} \circ f(\Phi ) =f\circ T_{\leqslant q-d} (\Phi )$ for~any $\Phi  \in \Amacro(D)$,
	\item $T_{\geqslant q} \circ f(\Phi ) =f\circ T_{\geqslant q-d} (\Phi )$ for~any $\Phi  \in \Amicro(D)$,
\end{itemize}
and similarly with strict inequalities.
\end{lemmesansboite}

\begin{textesansboite}
The next lemma studies the computation of $\langle \mathcal{F}\rangle (\Phi ):=\fsumn \big(\sum _{f\in \mathcal{F}} f\big)^{n} (\Phi )$ where $\mathcal{F}$ is a set of operators that have degrees of the same sign. By~definition of $\fsum$, for~any $q\in \mathbb{R}$, the component in $\mathcal{A}_{q}$ of $\fsum_{n=0} ^{N} \big(\sum _{f\in \mathcal{F}} f\big)^{n} (\Phi )$ becomes constant for sufficiently large $N$. Hence, to compute a truncation of this series, one needs to find a value of $N$ after which all the truncation's components are stationary. 
\end{textesansboite}

\begin{lemmesecable}{}
\label{calcul op <F> sur SF tronquee}%
Let $D\in \domsA$ and let $\mathcal{F}$ be a finite set of operators $\mathcal{A}(D)\rightarrow \mathcal{A}(D)$ that have a degree.
\begin{itemize}
	\item If $d:=\min\limits _{f\in \mathcal{F}} \deg f>0$, then we have for any $\Phi \in \Amacro(D)$ and $q\in \mathbb{R}$:
	\[T_{\leqslant q} \big( \langle \mathcal{F}\rangle (\Phi )\big) = T_{\leqslant q} \circ  \sum _{n=0} ^{N} \bigg(\sum _{f\in \mathcal{F}} f\bigg)^{n} \!\circ T_{\leqslant q} (\Phi ) \qquad  \text{with } \ N:= \left\lfloor { \frac{q-\min  \operatorname{supp}\Phi }{d} }\right\rfloor .\]
	\item If $d:=\max\limits _{f\in \mathcal{F}} \deg f<0$, then we have for any $\Phi \in \Amicro(D)$ and $q\in \mathbb{R}$:
	\[T_{\geqslant q} \big( \langle \mathcal{F}\rangle (\Phi )\big) = T_{\geqslant q} \circ \sum _{n=0} ^{N} \bigg(\sum _{f\in \mathcal{F}} f\bigg)^{n} \!\circ T_{\geqslant q} (\Phi ) \qquad  \text{with } \ N:= \left\lfloor { \frac{\max  \operatorname{supp}\Phi -q}{|d|} }\right\rfloor .\]
\end{itemize}
Similar results hold for strict inequalities.
\end{lemmesecable}

\begin{remarque}
The formulas in \cref{calcul op <F> sur SF tronquee} are optimal in the sense that taking $N$ smaller would give wrong formulas in general, and taking $N$ larger would give the same result.
\end{remarque}

\begin{demo}
We prove only the case $d>0$. The difference of the left and right-hand sides is equal to:
\begin{align*}
	T_{\leqslant q} \big( \langle \mathcal{F}\rangle (\Phi )\big) - T_{\leqslant q} \circ  \sum _{n=0} ^{N} \bigg(\sum _{f\in \mathcal{F}} f\bigg)^{n} \!\circ T_{\leqslant q} (\Phi )
	=  T_{\leqslant q} \circ \fsumn \bigg(\sum _{f\in \mathcal{F}} f\bigg)^{n} (\Phi ) -  T_{\leqslant q} \circ \sum _{n=0} ^{N} \bigg(\sum _{f\in \mathcal{F}} f\bigg)^{n} \!\circ T_{\leqslant q} (\Phi ) &\\
	= T_{\leqslant q} \circ  \underbrace{\fsum_{n=N+1} ^{\infty } \bigg(\sum _{f\in \mathcal{F}} f\bigg)^{n} (\Phi )}_{=:A} {}+ T_{\leqslant q} \circ  \underbrace{\sum _{n=0} ^{N} \bigg(\sum _{f\in \mathcal{F}} f\bigg)^{n} \!\circ T_{>q} (\Phi ) }_{=:B} &
\end{align*}
In addition, for any $\Psi  \in \Amacro(D)$, we have $\operatorname{supp}\big(\sum _{f\in \mathcal{F}} f(\Psi )\big) \subset \operatorname{supp}(\Psi )+d$. Hence, denoting $m:=\min  \operatorname{supp}\Phi $, we have $\operatorname{supp}(A) \subset \bigcup _{n>N} [m+nd,\infty )=[m+(N+1)d,\infty )$ by \cref{les SF sont stables par fsum}. By~construction of $N$, we have $m+(N+1)d>q$, so $T_{\leqslant q} (A)=0$. Moreover, we have $\operatorname{supp}(B)\subset \bigcup _{n=0} ^{N} (q+nd,\infty )=(q,\infty )$ so $T_{\leqslant q} (B)=0$.
\end{demo}

\begin{remarque}
In order to reduce the number of evaluations of the operators $f\in \mathcal{F}$, and thus the computational cost, one can use \cref{algo: calcul op <F> sur SF tronquee}. It is inspired by Horner's method to evaluate a polynomial.
\end{remarque}

\SetKwComment{Comment}{/* }{ */}
\begin{algorithm}[h]
\caption{Computation of $T_{\leqslant q} \big(\langle \mathcal{F}\rangle (\Phi )\big)$ in the first case of \cref{calcul op <F> sur SF tronquee}}
\label{algo: calcul op <F> sur SF tronquee}
\KwData{$\varphi  = T_{\leqslant q} (\Phi )$, $q$}
\KwResult{$\psi  = T_{\leqslant q} \big(\langle \mathcal{F}\rangle (\Phi )\big)$}
$\psi  \gets \varphi $ \;
$\chi  \gets \varphi $ \;
$N \gets \left\lfloor { (q-\min  \operatorname{supp}\varphi )/\min _{f\in \mathcal{F}} \deg f }\right\rfloor $ \;
\For{$n\gets 1$ \KwTo $N$}{
	$\chi  \gets \sum _{f\in \mathcal{F}} f(\chi )$ \Comment*[r]{Here we have $\chi =\big(\sum _{f\in \mathcal{F}} f \big)^{n} (\varphi )$}
	$\psi  \gets \psi  + \chi $ \; }
\end{algorithm}

\begin{exemplenumerote}
In order to compute a truncation of a field's behavior using truncations of the previous fields' behaviors, one can use the following formulas. They are obtained by applying \cref{calcul op a degre sur SF tronquee,calcul op <F> sur SF tronquee} to \cref{eq: bfuDApl,eq: SDApl}.
\setalignbox{\kern\hsize \kern-25mm\relax}
\begin{align*}
	T_{\leqslant q} (\bfuDA_{p,\ell } ) &\relax\.=\. T_{\leqslant q} \circ \sum _{n=0} ^{\lfloor {(p+q)/2}\rfloor } (-k_{0} ^2 \.\ROmegaPi)^{n} \bigg( {-\RcouchePi\circ (\partial _{x|\couche} ^2 +k_{1} ^2 )\circ T_{\leqslant q+2} (\bfuDA_{p-2,\ell } )} \\
	& \ralignbox{\displaystyle {} + \frac{\mu _{0} }{\mu _{1}} \RGammaPi\circ \partial _{y|\Gamma ,y=0^{+}} \circ T_{\leqslant q+1} (\bfuDA_{p-1,\ell } ) + \!\fsum_{d\in \TZ*\!,\.-p\leqslant d\leqslant q} \! \sigma _{d} (\bfuDA_{p,\ell } ) \,\phi_{d} \bigg) } \\
	T_{\geqslant q} (\SDA_{p,\ell } ) &\relax\.=\. T_{\geqslant q} \circ \sum _{n=0} ^{\lfloor {p-q}\rfloor }  \big( -\RcouchePi\circ \partial _{X|\couche} ^2  + \frac{\mu _{0} }{\mu _{1}} \RGammaPi\circ \partial _{Y|\Gamma ,Y=0^{+}} \big)^{n} \bigg( {- k_{0} ^2  \ROmegaPi \circ T_{\leqslant q-2} (\SDA_{p-2,\ell |\Omega } )}\\
	& \ralignbox{\displaystyle {} - k_{1} ^2  \RcouchePi \circ T_{\leqslant q-2} (\SDA_{p-2,\ell |\couche} ) + \!\!\fsum_{d\in \TZ*\!,\.q\leqslant d\leqslant p} \! \sigma _{d} (\SDA_{p,\ell } ) \,\phi_{d} \bigg). }
\end{align*}
\end{exemplenumerote}

\begin{textesansboite}
In order to compute the matching coefficients $\cmacro$ and $\cmicro$, we need one last class to represent $\Aeps(\Pi )$, or more precisely $\Aefini(\Pi )$, the set of elements of $\Aeps(\Pi )$ with finite support, as computers can store only finitely many objects. Let us denote this class \texttt{AEps}. It has only one attribute:
\begin{lstlisting}
	map<DValue,vector<APi>> terms; // for any (p,l), terms[p][l] is the term of the series that is multiplied by eps^p ln(eps)^l
\end{lstlisting}
As usual, we defined on this class \texttt{+}, \texttt{-} and \texttt{*}, as well as the natural extensions of $\ROmegaPi$, $\RcouchePi$, etc. We also defined the operators $\Eps$ and $\Eps^{-1}$ using their very descriptive construction (see \cref{def: Eps}). Here is a small example:
\begin{lstlisting}
	APi phi = /*...*/;
	AEps phi2 = phi; // casting to AEps
	AEps phi3 = phi2.timesEps(pi/Theta).timesLogEps(2); // eps^{pi/Theta} ln(eps)^2 phi
	AEps phi4 = phi2.hEps(); // H_eps(phi)
	AEps phi5 = phi2.hEpsInv(); // H_eps^{-1}(phi)
	lcomplex s = phi5.sigma(pi/Theta,3,1); // sigma_{pi/Theta,3,1}(phi_5)
\end{lstlisting}\\

Now, by \cref{eq: def cmaicro}, we have for any $(d,d',p,\ell )$
\[\cmacro_{d,d',p,\ell } := \sigma _{d,p,\ell } \circ \Eps^{-1} \circ \fsumn \bigg(\sum _{R\in \Rmicro} R \bigg)(\phi_{d'} )\]
and likewise for $\cmicro_{d,d',p,\ell }$. Similarly to \cref{calcul op <F> sur SF tronquee}, the sequence $N \mapsto  \sigma _{d,p,\ell } \circ \Eps^{-1} \circ \fsum_{n=0} ^{N} \big(\sum _{R\in \Rmicro} R \big)(\phi_{d'} )$ is stationary because the series over $n$ is formal and is~not a converging series. The next result gives an appropriate value of $N$ to truncate the series.
\end{textesansboite}

\begin{proposition}{}
\label{calcul numerique cmaicro}%
For~any $(d,d',p,\ell ) \in (\TZ*)^2 \times \NN\times \mathbb{N}$, we have:
\[\cmacro_{d,d',p,\ell } = \sigma _{d,p,\ell } \circ \Eps^{-1} \circ \sum _{n=0} ^{\lfloor {(3p+d+2d')/2}\rfloor } \bigg(\sum _{R\in \Rmicro} R\bigg)^{n} (\phi_{d'} ) \]
\vskip-0.5\baselineskip and\vskip-0.5\baselineskip
\[\cmicro_{d,d',p,\ell } = \sigma _{d,p,\ell } \circ \Eps \circ \sum _{n=0} ^{\lfloor {(3p-2d-d')/2}\rfloor } \bigg(\sum _{R\in \Rmacro} R\bigg)^{n} (\phi_{d'} )\]\fbeq[0.5]
\end{proposition}

\begin{demo}
We only prove the first case. Let $(d,d',p,\ell ) \in (\TZ*)^2 \times \NN\times \mathbb{N}$. For~any $N\in \mathbb{N}$, we have
\[A_{N} \.:=\, \cmacro_{d,d',p,\ell } - \sigma _{d,p,\ell } \circ \Eps^{-1} \circ \sum _{n=0} ^{N} \bigg(\sum _{R\in \Rmicro} R\bigg)^{n} (\phi_{d'} ) \,=\, \sigma _{d,p,\ell } \circ \Eps^{-1} \circ \underbrace{\fsum_{n=N+1} ^{\infty } \bigg(\sum _{R\in \Rmicro} R\bigg)^{n} (\phi_{d'} )}_{:=\Phi _{N}} .\]
$A_{N}$ depends only on the component of $\Phi _{N}$ in $\varepsilon ^{p+d} \mathcal{A}_{d} (\Pi )[\ln \varepsilon ]$. Let us find a value of $N$ s.t.\ this component vanishes.\\
Let $v:=(-2,3)$. For~any $\Phi \in \Aemicro(\Pi )$, we denote $\operatorname{supp}_{v} (\Phi ):=\{v\cdot (d,p) \mid (d,p)\in \operatorname{supp} \Phi \}$. We have $\min _{R\in \Rmicro} v\cdot \deg R = 2$ (see \cref{fig: degres des elements de Rmacro et Rmicro}), so, for any $\Phi \in \Aemicro(\Pi )$ and $M\in \mathbb{R}$ s.t.\ $\operatorname{supp}_{v} (\Phi )\subset [M,\infty )$, we have $\operatorname{supp}_{v} \big(\sum _{R\in \Rmicro} R(\Phi )\big) \subset  [M+2,\infty )$. This implies:
\[\operatorname{supp}_{v} (\Phi _{N} ) \subset  \bigcup _{N+1} ^{\infty } \operatorname{supp}_{v} \!\Bigg[ \bigg(\sum _{R\in \Rmicro} R\bigg)^{n} (\phi_{d'} ) \Bigg] \subset  [-2d'+2(N+1),\infty ).\]
Thus, to have $A_{N} =0$, it suffices that $-2d'+2(N+1)$ is greater than $v\cdot (d,p+d)=-2d+3(p+d)$, i.e.\ $N+1>(3p+d+2d')/2$. Hence, for $N:=\lfloor {(3p+d+2d')/2}\rfloor $, we have $A_{N} =0$.
\end{demo}

\begin{remarque}
One could take another vector $v$ in the proof of \cref{calcul numerique cmaicro} (provided that $\min _{R\in \Rmicro} v\cdot \deg R >0$), which would give another truncation order $N$. But this choice of $v$ seems to minimize $N$, which reduces the computations. Moreover, the series $\big(\sum _{R\in \Rmaicro} R\big)^{n} (\phi_{d'} )$ that appear in \cref{calcul numerique cmaicro} can be computed using an algorithm similar to \cref{algo: calcul op <F> sur SF tronquee}.
\end{remarque}

\begin{textesansboite}
In conclusion, although computing the matching coefficients involves several layers of classes, each layer only does simple operations on the layer beneath it, the bottom layer being integers and floating-point numbers, see \cref{fig: mes classes C++}. Hence, overall, these computations are just simple manipulations on small collections of numbers (typically, at most tens or hundreds of numbers). This is both much faster and memory-thriftier than computing singularities by solving ODES w.r.t.\ $\theta $, for example with finite elements. In addition, all our computations are exact (up to round-off errors), while finite elements require a fine mesh to be very precise. Finally, note that the algorithms and optimal truncation orders given in \cref{calcul op <F> sur SF tronquee,calcul numerique cmaicro,algo: calcul op <F> sur SF tronquee} are only refinements of our method and are not critical to its efficiency.
\end{textesansboite}

\begin{figure}[h]
\def\htEtage{0.88} 
\def\wdEtage{0.56} 
\def\topMargin{0.2}
\def\etage#1=#2{
	\def\x{(#1+\topMargin)}
	\draw ({-(\x-1)*\wdEtage},{-(\x-1)*\htEtage}) --  ({-\x*\wdEtage},{-\x*\htEtage}) -- ({\x*\wdEtage},{-\x*\htEtage}) -- ({(\x-1)*\wdEtage},{-(\x-1)*\htEtage});
	\draw (0,{-(\x-0.5 \if1#1 +0.1\fi)*\htEtage}) node {#2}; }
\let\code\lstinline
\begin{tikzpicture}
	\draw ({-\topMargin*\wdEtage},{-\topMargin*\htEtage}) -- (0,0) -- ({\topMargin*\wdEtage},{-\topMargin*\htEtage});
	\etage1={\code{AEps}}
	\etage2={\code{APi}}
	\etage3={\parbox{4cm}{\centering\strut\code{AOmega}, \code{AGamma},\\\code{ALayer}}}
	\etage4={\code{DValue}}
	\etage5={\code{Ratio}}
	\etage6={\code{int}, \code{complex<long double>}}
\end{tikzpicture}\par
\longCaption{Schematic illustration of the classes used to compute the matching coefficients}{. Each layer is built using the layers beneath it.}
\label{fig: mes classes C++}
\end{figure}

\let\varepsilon=\peps
\chapter{Contruction of the asymptotic expansion and error estimates}
\label{chap: construction du DA}

\begin{textesansboite}
This chapter is dedicated to the construction of the terms of the asymptotic expansion of $u_{\varepsilon }$, and to the mathematical justification of the expansion by error estimates. The outline is the following.
\begin{enumerate}
	\item In \cref{sec: pb modeles bien poses}, we introduce model problems that will be used to build each far field and corner field. We prove that they possess a unique solution. These problems are non-variational as the fields have singular behaviors (at the corner for far fields and at infinity for corner fields). However, the solution can be decomposed into a singular part and a variational part. The layer fields $U_{p,\ell }$ are treated separately by showing that they have an explicit expression in terms of the far fields $u_{p,\ell }$.
	\item In \cref{sec: DA autour du coin}, we prove asymptotic expansions w.r.t.\ $\bfr$ at any order for the solutions of those model problems using the Kondrat'ev theory. These expansions hold when $\bfr\rightarrow 0$ for far-and-layer fields and when $\bfr\rightarrow \infty $ for corner fields. The asymptotic behavior are represented by formal series in $\Amaicro(\Pi )$, which allows to use the results of \cref{chap: SF,chap: raccord algebrique}.
	\item In \cref{sec: construction des champs} we gather all the equations that the fields $u_{p,\ell }$, $U_{p,\ell }$ and $S_{p,\ell }$ must satisfy (i.e.\ \crefrange{eq: eq vol champs lointains}{eq: eq vol champs de coin} and \cref{condition de raccord}), and we prove that there exists a unique set of fields that satisfy them. Moreover we provide a way to build directly the far fields, without computing the layer and corner fields. This second construction uses some corner profiles.
	\item Finally, \cref{sec: estimations d'erreur} rigorously proves \cref{DA de SE avec CL et delta fixe}, which justifies the asymptotic expansion.
\end{enumerate}
\end{textesansboite}

\section{Existence and uniqueness for model problems}
\label{sec: pb modeles bien poses}

\begin{textesansboite}
First of all, let us express the layer fields with the far fields, so~that only two types of fields remain to build. Let $(\mathcal{U}_{n} )\in \mathbb{R}[Y]^{\mathbb{N}}$ be the unique sequence of polynomials s.t.\ for~any $n\in \mathbb{N}^*$:
\begin{equation}{\label{eq: correcteurs de couche}}
\left\{\begin{array}{r@{\;=\;}ll}
	\mathcal{U}_{0} '' & 0 \\[0.5mm]
	\mathcal{U}_{0} '(0) & \frac{\mu _{0} }{\mu _{1}} \\[0.5mm]
	\mathcal{U}_{0} (-1) & 0 
\end{array}\right. \qquad  \text{ and } \qquad  \left\{\begin{array}{r@{\;=\;}ll}
	\mathcal{U}_{n} '' & -\mathcal{U}_{n-1} \\[0.5mm]
	\mathcal{U}_{n} '(0) & 0 \\[0.5mm]
	\mathcal{U}_{n} (-1) & 0 
\end{array}\right.\end{equation}
\end{textesansboite}

\begin{lemme}{expression of the layer fields}
\label{expression de Upl}%
Let us assume that the fields $u_{p,\ell }$ and $U_{p,\ell }$ are regular enough (we will check later that they are). \cref{eq: eq vol champs de couche} implies for~any $(p,\ell )\in \NN \times \mathbb{N}$ and $(x,Y)\in \couche$:\vspace{-1mm}
\begin{equation}{\label{eq: expression de Upl}}
U_{p,\ell } (x,Y) = \sum_{n=0}^\infty  (\partial _{x} ^2 +k_{1} ^2 )^{n} \partial _{y} u_{p-1-2n,\ell } (x,0) \cdot  \mathcal{U}_{n} (Y).\vspace{-1mm}\end{equation}
where this sum has a finite number of non-zero terms by the convention: $\forall p\in \mathbb{R} \setminus \NN,\forall \ell \in \mathbb{N},\; u_{p,\ell } :=0$.
\end{lemme}

\begin{demo}
There exists an increasing sequence $(p_{m} )$ s.t.\ $\NN=\{p_{m} \mid m\in \mathbb{N}\}$. So we can prove the result by induction on $p\in \NN$. For $p=0$, \cref{eq: eq vol champs de couche} states: \[\left\{\begin{array}{r@{\;=\;}ll}
	\partial _{Y} ^2  U_{0,\ell } & 0 &\text{ in }\couche \\[0.5mm]
	\partial _{Y} U_{0,\ell } & 0 &\text{ on }\Gamma  \\[0.5mm]
	U_{0,\ell } & 0 &\text{ on }\borddroit
\end{array}\right.\] so $U_{0,\ell } =0$ for~any $\ell $. It~is coherent with \cref{eq: expression de Upl} (which is a sum of zeros in this case).\\
Next, for the inductive step, we assume that \cref{eq: expression de Upl} holds for ranks smaller that $p$. \cref{eq: eq vol champs de couche} gives: \[\left\{\begin{array}{r@{}l@{}l@{\,}l}
	\partial _{Y} ^2  U_{p,\ell } &\;=\; -(\partial _{x} ^2 +k_{1} ^2 )U_{p-2,\ell }
	&\;=\; -\sum_{n=0}^\infty  (\partial _{x} ^2 +k_{1} ^2 )^{n+1} \partial _{y} u_{(p-2)-1-2n,\ell } (x,0) \cdot  \mathcal{U}_{n} (Y) & \\[0.5mm]
	&&\;=\; -\sum_{n=1}^\infty  (\partial _{x} ^2 +k_{1} ^2 )^{\red{n}} \partial _{y} u_{p-1-2n,\ell } (x,0) \cdot  \mathcal{U}_{\red{n-1}} (Y) &\text{ in }\couche \\
	\partial _{Y} U_{p,\ell } &\;=\; \frac{\mu _{0} }{\mu _{1}} \,\partial _{y} u_{p-1,\ell } &&\text{ on }\Gamma  \\[0.5mm]
	U_{p,\ell } &\;=\; 0 &&\text{ on }\borddroit
\end{array}\right.\] It~is easy to see that $\sum_{n=0}^\infty  (\partial _{x} ^2 +k_{1} ^2 )^{n} \partial _{y} u_{p-1-2n,\ell } (x,0) \cdot  \mathcal{U}_{n} (Y)$ is the only solution of this.
\end{demo}

\begin{textesansboite}
We saw in \cref{chap: raccord algebrique} that far and corner fields possess singularities when $r\rightarrow 0$, resp.\ $\bfr\rightarrow \infty $. So the usual variational frameworks are~not sufficient to build these fields and we need to design new frameworks. This is done in \cref{sec: cadre fonctionnel des champs lointains,sec: cadre fonctionnel des champs de coin}. In both sections we start by introducing the natural space $H_{\mathrm{var}}$ in which an ad hoc variational problem is well-posed. Then we define a larger space $H := H_{\mathrm{var}} +\chi  \,\mathcal{A}(D) = \{u+\chi \varphi  \mid u\in H_{\mathrm{var}} , \varphi \in \mathcal{A}(D)\}$ that contains the singularities, where $\chi $ is a $\mathcal{C}^{\infty }$ truncation function in the vicinity of 0 (for $u_{p,\ell }$) or infinity (for $S_{p,\ell }$), and $D\in \{\Omega ,\Pi \}$. Next we determine the elements of $H_{\mathrm{var}} \cap  \chi  \mathcal{A}(D)$, which allows us to define on $H$ the linear forms $\sigma _{d}$ associated to the singularities. Finally we show that some model problems are well-posed in $H$.\\

Before we start, the following lemma will be a useful tool to estimate the behavior at 0 and $\infty $ of functions of the spaces $\mathcal{A}$.
\end{textesansboite}

\begin{lemme}{}
\label{croissance des fcts de A}%
Let $a<b$ in $\mathbb{R}$, $n\in \mathbb{N}^*$, $(d_{i} ,\ell _{i} )_{i\in [\![1,n]\!]}$ be $n$ distinct elements of $\mathbb{R}\times \mathbb{N}$, $(f_{i} )_{i\in [\![1,n]\!]} \in (\mathcal{C}([a,b],\mathbb{C}) \setminus \{0\})^{n}$ and:\vspace{-2mm}
\[\varphi :(r,\theta )\in \mathbb{R}_+^*\times [a,b] \;\mapsto \;\sum_{i=1}^n r^{d_{i}} \ln ^{\ell _{i}} \! r\cdot  f_{i} (\theta ).\]
Then there exist an interval $I\subset [a,b]$ with non-empty interior, $c\in \mathbb{R}_+^*$ and $r_{1} ,r_{2} \in \mathbb{R}$ such that:\vspace{-1mm}
\[\forall r\in (0,r_{1} ),\; \forall \theta \in I,\ \ |\varphi (r,\theta )|>c\, r^{\min _{j} d_{j} } \quad  \text{ and } \quad  \forall r\in (r_{2} ,\infty ),\; \forall \theta \in I,\ \ |\varphi (r,\theta )|>c\, r^{\max _{j} d_{j} } .\]\fbeq[1.8]
\end{lemme}

\begin{demo}
Let $j\in [\![1,n]\!]$ be s.t.\ $(d_{j} ,\ell _{j} )$ is maximal for the lexicographic order. Let $I\subset [a,b]$ be a non-trivial interval on which $|f_{j} |$ is greater than a positive constant. Since $\sum _{i\neq j} r^{d_{i}} \ln ^{\ell _{i}} \! r\cdot  f_{i} (\theta ) =o(r^{d_{j}} \ln ^{\ell _{j}} r)$, we have when $r\rightarrow \infty $ and $\theta \in I$: $|\varphi (r,\theta )| \gtrsim  r^{d_{j}} \ln ^{\ell _{j}} r \gtrsim  r^{d_{j}}$. And we can similarly treat the vicinity of 0.
\end{demo}

\begin{textesansboite}
We can apply \cref{croissance des fcts de A} to any $\varphi \in \mathcal{A}_{d} (\Omega )$ with $(a,b):=(0,\Theta )$, or to any $\varphi \in \mathcal{A}_{d} (\couche)$ with $(a,b):=(-1,0)$ (replacing the variables $(r,\theta )$ by $(x,Y)$). We can also apply it to $\partial _{r} \varphi $ and $\partial _{\theta } \varphi $ when $\varphi \in \mathcal{A}_{d} (\Omega )$ and to $\partial _{x} \varphi $ and $\partial _{Y} \varphi $ when $\varphi \in \mathcal{A}_{d} (\couche)$.
\end{textesansboite}

\begin{paragraphesansboite}{Definition:}
For~any $D\in \{\Pi ,\Omega ,\Gamma ,\couche\}$ and $\varphi \in \mathcal{A}(D)$, we denote $\mindeg \varphi :=\sup \{d\in \mathbb{R} \mid \varphi \in \sum _{q\geqslant d} \mathcal{A}_{q} (D)\}$ and $\maxdeg \varphi :=\inf \{d\in \mathbb{R} \mid \varphi \in \sum _{q\leqslant d} \mathcal{A}_{q} (D)\}$.
\end{paragraphesansboite}
\nomenclature[22-18-23]{$\mindeg,\maxdeg$}{Minimal or maximal homogeneity degree of an element of a space $\mathcal{A}$}

\subsection{Well-posedness of far fields-like problems}
\label{sec: cadre fonctionnel des champs lointains}

\begin{textesansboite}
Let us denote $H^{1/2} _{00} (\Gamma )$ the set of functions of $H^{1/2} (\Gamma )$ whose extension by 0 to $\partial \Omega =\Gamma \cup \{(0,0)\}\cup \bordgauche$ is in $H^{1/2} (\partial \Omega )$. Using Lax-Milgram theorem, it~is easy to prove the following lemma.
\end{textesansboite}

\begin{lemme}{the Helmholtz problem in $H^{1} (\Omega )$}
\label{pb de Helmholtz variationnel}%
Let $f\in (H^{1} (\Omega ))'$ and $g\in H^{1/2} _{00} (\Gamma )$. The following system has a unique solution in $H^{1} (\Omega )$. \[\left\{\begin{array}{r@{\;=\;}ll}
	\mu _{0} \Delta  u +\omega ^2 \rho _{0} u & f &\text{ in }\Omega  \\[0.5mm]
	u & g &\text{ on }\Gamma  \\[0.5mm]
	u & 0 &\text{ on }\bordgauche
\end{array}\right.\]\fbeq
\end{lemme}

\begin{textesansboite}
Let $\khiz$ be a radial function of $\mathcal{C}^{\infty } (\mathbb{R}^2 )$ equal to 1 in the vicinity of 0 and to 0 in the vicinity of infinity. The appropriate space to build the far fields is $H^{1} (\Omega )+\khiz \mathcal{A}(\Omega )$. One can check that it does~not depend on the choice of $\khiz$.
\end{textesansboite}
\nomenclature[57-00]{$\khiz$}{Truncation function equal to 1 in the vicinity of 0 and to 0 in the vicinity of infinity}

\begin{lemmesansboite}{}
\label{H1 inter A}%
$\displaystyle H^{1} (\Omega ) \cap  \khiz \mathcal{A}(\Omega ) = \khiz \sum _{d>0} \mathcal{A}_{d} (\Omega ).$\fbeq
\end{lemmesansboite}

\begin{demo}
The inclusion $\supset$ is easy to check, so we focus on $\subset $. Let $\varphi \in \mathcal{A}(\Omega ) \setminus \{0\}$ be s.t.\ $\khiz \varphi \in H^{1} (\Omega )$. There exist $d\in \mathbb{R}$, $\varphi _{1} \in \mathcal{A}_{d} (\Omega ) \setminus \{0\}$ and $\varphi _{2} \in \sum _{q>d} \mathcal{A}_{q} (\Omega )$ s.t.\ $\varphi =\varphi _{1} +\varphi _{2}$. If $\partial _{\theta } \varphi _{1}$ were null everywhere, then would so too $\varphi _{1}$ because $\varphi _{1|\bordgauche} =0$. But we assumed the contrary, so $\partial _{\theta } \varphi _{1} \neq 0$. Thus \cref{croissance des fcts de A} implies that there exists a non-trivial interval $I\subset [0,\Theta ]$ s.t.\ $|\partial _{\theta } \varphi _{1} (r,\theta )| \gtrsim  r^{d}$ when $r\rightarrow 0$ and $\theta \in I$. Finally $r^{-1} \partial _{\theta } (\khiz \varphi ) \in L^2 (\Omega )$ implies that $d>0$.
\end{demo}

\begin{definitionsansboite}{}
\label{def: sigma dans H1 + khi A}%
For~any $u\in H^{1} (\Omega )+\khiz \mathcal{A}(\Omega )$ and $d\in -\TN$, we denote $\sigma _{d} (u):=\sigma _{d} (\varphi )$ where $\varphi \in \mathcal{A}(\Omega )$ is s.t.\ $u-\khiz \varphi \in H^{1} (\Omega )$ (see \cref{def: sigma} for $\sigma _{d} (\varphi )$). It does~not depend on the choice of $\varphi $ thanks to \cref{H1 inter A}.
\end{definitionsansboite}
\nomenclature[53-13b]{$\sigma _{d}$}{Linear form on $H^{1} (\Omega )+\khiz \mathcal{A}(\Omega )$}

\begin{theoreme}{existence and uniqueness for a far fields-like model problem}
\label{cadre fonctionnel des champs lointains}%
Let $f\in (H^{1} (\Omega ))'$, $g\in H^{1/2} _{00} (\Gamma )+\khiz \mathcal{A}(\Gamma )$ and $(s_{d} )\in \mathbb{C}^{-\TN}$ with finite support. The following system has a unique solution in $H^{1} (\Omega )+\khiz \mathcal{A}(\Omega )$.
\[\left\{\begin{array}{r@{\;=\;}ll}
	\mu _{0} \Delta  u +\omega ^2 \rho _{0} u & f &\text{ in }\Omega  \\[0.5mm]
	u & g &\text{ on }\Gamma  \\[0.5mm]
	u & 0 &\text{ on }\bordgauche\\[0.5mm]
	\sigma _{d} (u)& s_{d} & \ \forall d\in -\TN
\end{array}\right.\]\fbeq
\end{theoreme}

\begin{remarque}
The coefficients $\sigma _{d} (u)$ with $d<0$ can be imposed because they describe a non-variational part of the behavior of the solution at the corner, while the coefficients $\sigma _{d} (u)$ with $d>0$ cannot be fixed because they belong to the variational part. This differs from \cref{pb de Helmholtz dans A Omega}, which solves a similar problem in $\Amacro(\Omega )$, where all the coefficients $\sigma _{d}$ can be fixed because this algebraic setting does~not depend on the variational space.
\end{remarque}

\begin{demo}
Let us show the existence first, and then the uniqueness.\\
\underline{Existence:} Let $\tilde{g} \in H^{1/2} _{00} (\Gamma )$ and $\varphi \in \mathcal{A}(\Gamma )$ be s.t.\ $g= \tilde{g} +\khiz \varphi $. We look for the solution in the form $u= \tilde{u}+\khiz \psi $ with $\tilde{u} \in H^{1} (\Omega )$ and $\psi \in \mathcal{A}(\Omega )$. Let
\[\psi ^{+} := \fsumn (-k_{0} ^2 \ROmegaOmega)^{n} \bigg(\RGammaOmega(\varphi )+\fsum_{d\in -\TN} s_{d} \, \phiOmega_{d} \bigg) \;\in \;\Amacro(\Omega )\]
where we recall that $\phiOmega_{d} :=\phi_{d|\Omega }$, and let $\psi :=T_{\leqslant 2} (\psi ^{+} )$. \Cref{pb de Helmholtz dans A Omega} shows that we have
\[\left\{\begin{array}{r@{\;=\;}ll}
	(\mu _{0} \Delta  +\omega ^2 \rho _{0} )\psi ^{+} & 0 &\text{ in }\Omega  \\[0.4mm]
	\psi ^{+} & \varphi  &\text{ on }\Gamma  \\[0.4mm]
	\psi ^{+} & 0 &\text{ on }\bordgauche\\[0.4mm]
	\sigma _{d} (\psi ^{+} )& s_{d} & \ \forall d\in -\TN
\end{array}\right.\]
Using $\psi \in \mathcal{A}(\Omega )$, $\psi  = \psi ^{+} -(T_{>2} )(\psi ^{+} )$ and $\deg \Delta =-2$ (see \cref{operateurs de derivation dans A}), this implies
\[\left\{\begin{array}{r@{}l}
	(\mu _{0} \Delta  +\omega ^2 \rho _{0} )\psi  &{}\in  \mathcal{A}(\Omega )\cap  \fsum\limits_{d>0} \mathcal{A}_{d} (\Omega ) = \sum\limits _{d>0} \mathcal{A}_{d} (\Omega ) \\[0.5mm]
	\psi _{|\Gamma } -\varphi  &{}\in  \mathcal{A}(\Gamma )\cap  \fsum\limits_{d>2} \mathcal{A}_{d} (\Gamma ) = \sum\limits _{d>2} \mathcal{A}_{d} (\Gamma ) \\
	\psi _{|\bordgauche} &{}= 0 \\[0.5mm]
	\sigma _{d} (\psi  )&{}= s_{d} \;\ \forall d\in -\TN
\end{array}\right.\]
Therefore $(\mu _{0} \Delta  +\omega ^2 \rho _{0} )(\khiz \psi )\in L^2 (\Omega )$ and $\khiz \,(\psi _{|\Gamma } -\varphi )\in H^{1/2} _{00} (\Gamma )$. Finally, thanks to \cref{pb de Helmholtz variationnel}, we can take $\tilde{u}$ as the unique solution in $H^{1} (\Omega )$ of: \[\left\{\begin{array}{r@{\;=\;}ll}
	\mu _{0} \Delta  \tilde{u} +\omega ^2 \rho _{0} \tilde{u} & f - (\mu _{0} \Delta +\omega ^2 \rho _{0} )(\khiz \psi ) &\text{ in }\Omega  \\[0.5mm]
	\tilde{u} & \tilde{g} -\khiz \,(\psi -\varphi )  &\text{ on }\Gamma  \\[0.5mm]
	\tilde{u} & 0 &\text{ on }\bordgauche
\end{array}\right.\]

\underline{Uniqueness:} Let $u$ be a homogeneous solution. There exist $\tilde{u} \in H^{1} (\Omega )$ and $\psi \in \mathcal{A}(\Omega )$ s.t.\ $u= \tilde{u}+\khiz \psi $. Let $\varphi :=\Delta \psi +k_{0} ^2 \psi $. By~\cref{operateurs de derivation dans A}, $\varphi \in \mathcal{A}(\Omega )$. Let us show that $\mindeg \varphi >-2$.
\begin{sousdemo}
	We denote $d:=\mindeg \varphi $ and $\widetilde{\varphi } :=(\Delta +k_{0} ^2 )(\khiz \psi )$. By~\cref{croissance des fcts de A}, there exists a non-trivial interval $I\subset [0,\Theta ]$ s.t.\ when $r\rightarrow 0$ and $\theta \in I$: $\widetilde{\varphi }(r,\theta ) =\varphi (r,\theta )\gtrsim  r^{d}$. Besides, $\widetilde{\varphi } = -(\Delta +k_{0} ^2 ) \tilde{u} \,\in  (H^{1} _{0} (\Omega ))'$. Let us test it with $\zeta _{q} :(r,\theta ) \mapsto  \widetilde{\varphi } (r,\theta )\cdot r^{q-2d-2} \,\big(1-\khiz(2^{1/q} r)\big)\,\chi (\theta )$ where $q>0$ and $\chi \in \mathcal{D}(0,\Theta ) \setminus \{0\}$ is everywhere non negative. Since $\zeta _{q} \in \mathcal{D}(\Omega )$:
	\[\langle  \widetilde{\varphi }, \zeta _{q} \rangle  := \int _{0} ^{\infty } \int _{0} ^{\Theta } \widetilde{\varphi } \cdot \zeta _{q} \cdot r\mathrm{d}\theta \,\mathrm{d}r \,\lesssim \, \|\zeta _{q} \|_{H^{1} (\Omega )} \qquad  \forall q>0\]
	Let us assume by contradiction that $d\leqslant -2$. It~is~not difficult to check that when $q\rightarrow 0^{+}$:
	\[\langle  \widetilde{\varphi }, \zeta _{q} \rangle  \gtrsim  \int _{c_{1} 2^{- \frac{1}{q} }} ^{c_{2}} r^{d} \,r^{q-d-2} \,r\mathrm{d}r \gtrsim  \frac{1}{q} \quad  \text{ and } \quad  \|\zeta _{q} \|_{H^{1}} \lesssim  \|\nabla  \zeta _{q} \|_{L^2 }  \lesssim  \bigg(\int _{c_{3} 2^{- \frac{1}{q} }} ^{c_{4}} r^{2(q-d-3)} \,r\mathrm{d}r \!\bigg)^{ \frac{1}{2} } \!\lesssim  \frac{1}{\sqrt{q}}\]
	with some constants $c_{i}$ and Poincaré's inequality. But it contradicts $\langle  \widetilde{\varphi }, \zeta _{q} \rangle  \lesssim  \|\zeta _{q} \|_{H^{1}}$. So $d>-2$.		
\end{sousdemo}
Now let us show that $\mindeg \psi >0$.
\begin{sousdemo}
	Let us assume the contrary. Then there exist $d\leqslant 0$, $\psi _{1} \in \mathcal{A}_{d} (\Omega ) \setminus \{0\}$ and $\psi _{2} \in \sum _{q>d} \mathcal{A}_{q} (\Omega )$ s.t.\ $\psi =\psi _{1} +\psi _{2}$. \Cref{operateurs de derivation dans A} states that $\deg  \Delta =-2$, so $\Delta \psi _{1} \in \mathcal{A}_{d-2} (\Omega )$ and $\Delta \psi _{2} \in \sum _{q>d-2} \mathcal{A}_{q} (\Omega )$. In addition, $\varphi  \in \sum _{q>-2} \mathcal{A}_{q} (\Omega )$ and $\mathcal{A}(\Omega ) = \bigoplus_{q\in \mathbb{R}} \mathcal{A}_{q} (\Omega )$ according to \cref{decomposition de A en sev}. So~identifying the coordinate of $\Delta \psi +k_{0} ^2 \psi  =\varphi $ in $\mathcal{A}_{d-2} (\Omega )$ gives $\Delta \psi _{1} =0$.\\
	Moreover, $\psi _{|\Gamma } = - \tilde{u}_{|\Gamma } \in H^{1/2} (\Gamma )$, so $\psi _{|\Gamma } \in \sum _{q>0} \mathcal{A}_{q} (\Gamma )$. Hence: \[\left\{\begin{array}{r@{\;=\;}ll}
		\Delta  \psi _{1} & 0 &\text{ in }\Omega \\[0.5mm]
		\psi _{1} & 0 &\text{ on }\bordgauche \cup \Gamma 
	\end{array}\right.\]
	By~\cref{pb de Laplace homogene dans A}, it implies $d\in -\TN$ and $\psi _{1} =\sigma _{d} (\psi _{1} )\,\phiOmega_{d}$. However, $\sigma _{d} (\psi _{1} )=\sigma _{d} (\psi )=\sigma _{d} (u)=0$, which contradicts $\psi _{1} \neq 0$.
\end{sousdemo}
Finally, \cref{H1 inter A} implies $\khiz \psi \in H^{1} (\Omega )$, so $u\in H^{1} (\Omega )$. This means that $u=0$ by \cref{pb de Helmholtz variationnel}.
\end{demo}

\subsection{Well-posedness of corner fields-like problems}
\label{sec: cadre fonctionnel des champs de coin}

\begin{definition}{the variational space $\Vcoin$}
\label{def: V}%
We define $\Vcoin := \{v\in H^{1} _{\mathrm{loc}} (\Omega _{1} ) \mid \nabla  v \in L^2 (\Omega _{1} ) \text{ and } v_{|\partial \Omega _{1}} = 0\}$ and the norm $\|v\|_{\Vcoin} := \|\nabla v\|_{L^2 (\Omega _{1} )}$.\\
Here ``$v_{|\partial \Omega _{1}} = 0$'' means that $\khiz( \frac{\cdot }{R} ) \,v\in H^{1} _{0} (\Omega _{1} )$ for~any $R>0$, where $\khiz( \frac{\cdot }{R} ): (X,Y)\mapsto  \khiz\big( \frac{X}{R} , \frac{Y}{R} \big)$. 
\end{definition}
\nomenclature[31]{$\Vcoin$}{Space of functions $\Omega _{1} \rightarrow \mathbb{C}$ whose gradient is $L^2 $}

\begin{lemmesansboite}{}
\label{inegalite de Hardy dans Omega1}%
Any $v\in \Vcoin$ satisfies $\|v\|_{H^{1} (\couche\cap \{X>\Rcoin\})} \lesssim  \|v\|_{\Vcoin}$ \,and \;$\displaystyle \Big\| \frac{v}{1+r} \Big\|_{L^2 (\Omega _{1} )} \lesssim  \|v\|_{\Vcoin}$.\fbeq[0.5]
\end{lemmesansboite}

\newcommand{\coucheR}{A}
\begin{demo}
We denote $\coucheR :=\couche \cap \{X>\Rcoin\}$ and $B:=B(0,\Rcoin)$, and we recall that $\Rcoin>0$ is s.t.\ $\Omega _{1} \setminus B=\Pi  \setminus B$.\\
Poincaré's inequality gives for a.e. $X>\Rcoin$: $\|v(X,\cdot )\|_{L^2 (-1,0)} \; \lesssim  \; \|\partial _{Y} v(X,\cdot )\|_{L^2 (-1,0)}$. Integrating w.r.t.\ $X$ the square of this then gives: $\|v\|_{L^2 (\coucheR)} \lesssim  \|\partial _{Y} v\|_{L^2 (\coucheR)} \leqslant  \|v\|_{\Vcoin}$. So~$\|v\|_{H^{1} (\coucheR)} \lesssim  \|v\|_{\Vcoin}$.\\
Next, a generalized Poincaré inequality on circular slices of $\Omega  \setminus B$ gives for a.e. $r>\Rcoin$: $\|v(r,\cdot )\|_{L^2 (0,\pi )} \lesssim  \|\partial _{\theta } v(r,\cdot )\|_{L^2 (0,\pi )} + |v(r,\theta =0)|$. Then integrating w.r.t.\ $r$ the square of this times $ \frac{1}{1+r}$ gives:\vspace{-2mm}
\[\int _{\Omega  \setminus B} \frac{v^2 }{(1+r)r} \,r\mathrm{d}\theta \,\mathrm{d}r \;\lesssim \; \int _{\Omega  \setminus B} \frac{(\partial _{\theta } v)^2 }{(1+r)r} \,r\mathrm{d}\theta \,\mathrm{d}r + \int _{\Gamma  \setminus B} \frac{v^2 }{1+r} \mathrm{d}r.\]
Therefore $\| \frac{v}{1+r} \|_{L^2 (\Omega  \setminus B)} ^2  \lesssim  \|\nabla  v\|_{L^2 (\Omega  \setminus B)} ^2 + \|v\|_{L^2 (\Gamma  \setminus B)} ^2 $. But we know that $\|v\|_{L^2 (\Gamma  \setminus B)} \lesssim  \|v\|_{H^{1} (A)} \lesssim  \|v\|_{\Vcoin}$, so $\| \frac{v}{1+r} \|_{L^2 (\Pi  \setminus B)} \lesssim \|v\|_{\Vcoin}$. Finally Poincaré's inequality in $B$ also gives $\| \frac{v}{1+r} \|_{L^2 (\Omega _{1} \cap  B)} \lesssim \|v\|_{\Vcoin}$.
\end{demo}

\begin{lemme}{the Poisson problem in $\Omega _{1}$}
\label{pb variationnel de Poisson dans Omega1}%
Let $F:\Omega _{1} \rightarrow \mathbb{C}$ be s.t.\ $(1+r)F\in L^2 (\Omega _{1} )$, and $g\in L^2 (\Gamma \cap \{X>\Rcoin\})$. The following system is well-posed in $\Vcoin$.\vspace{-2mm}
\[\left\{\begin{array}{r@{\;=\;}ll}
	\operatorname{div}(\mu  \nabla  S) &F &\text{ in }\Omega _{1} \setminus (\Gamma \cap \{X>\Rcoin\})  \\[0.5mm]
	S &0 &\text{ on }\partial \Omega _{1} \\[0.5mm]
	S_{|Y=0^{+}} - S_{|Y=0^{-}} & 0 &\text{ on }\Gamma \cap \{X>\Rcoin\} \\[0.5mm]
	\mu _{0} \partial _{Y} S_{|Y=0^{+}} - \mu _{1} \partial _{y} S_{|Y=0^{-}} & g &\text{ on }\Gamma \cap \{X>\Rcoin\} 
\end{array}\right.\]\fbeq
\end{lemme}

\begin{demo}
The variational formulation of this problem is: \[ \forall v\in \Vcoin,\qquad  \int _{\Omega _{1}} \mu  \nabla S\cdot \nabla v = -\int _{\Omega _{1}} Fv+ \int _{\Gamma \cap \{X>\Rcoin\}} gv.\]
The left-hand side is coercive by~definition of $\Vcoin$. \Cref{inegalite de Hardy dans Omega1} ensures that the right-hand side is continuous. Moreover it~is easy to see that $\Vcoin$ is complete. So we can conclude using the Lax-Milgram theorem.
\end{demo}

\begin{textesansboite}
Let $\khii$ be a function of $\mathcal{C}^{\infty } (\mathbb{R}^2 )$ equal to 0 on $B(0,\Rcoin)$ and 1 in a vicinity of infinity. The proper space to build the corner fields is $\Vcoin +\khii \mathcal{A}(\Pi )$. One can check that it does~not depend on the choice of $\khii$.
\end{textesansboite}
\nomenclature[57-60]{$\khii$}{Truncation function equal to 0 in the vicinity of 0 and to 1 in the vicinity of infinity}

\begin{lemmesansboite}{}
\label{V inter A}%
$\Vcoin\cap  \khii \mathcal{A}(\Pi ) = \big\{ \khii \varphi  \;\big|\; \varphi \in \mathcal{A}(\Pi ),\; \maxdeg(\varphi _{|\Omega } )<0\; \text{ and } \,\maxdeg(\varphi _{|\couche} )<- \frac{1}{2} \big\}.$
\end{lemmesansboite}

\begin{demo}
The inclusion $\supset$ is easy to check, so we focus on $\subset $. Let $\varphi \in \mathcal{A}(\Pi )$ be s.t.\ $\khii \varphi \in \Vcoin$. The same method as \cref{H1 inter A} shows that $\maxdeg \varphi _{|\Omega } <0$. And using that $\partial _{Y} (\khii \varphi ) \in L^2 (\couche)$ we likewise get $\maxdeg(\varphi _{|\couche} )<- \frac{1}{2}$.
\end{demo}

\begin{definitionsansboite}{}
\label{def: sigma dans V + khi A}%
For~any $S\in \Vcoin +\khii \mathcal{A}(\Pi )$ and $d\in \TN$, we denote $\sigma _{d} (S):=\sigma _{d} (\varphi )$ where $\varphi \in \mathcal{A}(\Pi )$ is s.t.\ $S-\khii \varphi \in \Vcoin$. It does~not depend on the choice of $\varphi $ thanks to \cref{V inter A}.
\end{definitionsansboite}
\nomenclature[53-13c]{$\sigma _{d}$}{Linear form on $V+\khii \mathcal{A}(\Pi )$}

\begin{theoreme}{existence and uniqueness for a corner fields-like model problem}
\label{cadre fonctionnel pour les champs de coin}%
Let $F:\Omega _{1} \rightarrow \mathbb{C}$ be s.t.\ $(1+r)F\in L^2 (\Omega _{1} )$, $\varphi :\Pi \rightarrow \mathbb{C}$ be s.t.\ $\varphi _{|\Omega } \in \mathcal{A}(\Omega )$ and $\varphi _{|\couche} \in \mathcal{A}(\couche)$, and $(s_{d} )\in \mathbb{C}^{\TN}$ with finite support. The following system has a unique solution in $V+\khii \mathcal{A}(\Pi )$.
\[\left\{\begin{array}{r@{\;=\;}ll}
	\operatorname{div}(\mu  \nabla  S) & F +\khii \varphi  &\text{ in }\Omega _{1} \\[0.5mm]
	S & 0 &\text{ on }\partial \Omega _{1} \\[0.5mm]
	\sigma _{d} (S)& s_{d} & \ \forall d\in \TN
\end{array}\right.\]\fbeq
\end{theoreme}

\begin{remarque}
$\varphi $ is~not just an element of $\mathcal{A}(\Pi )$ because it may be discontinuous on $\Gamma $.
\end{remarque}

\begin{demo}
This proof is similar to \cref{cadre fonctionnel des champs lointains}. Let us show the existence first, and then the uniqueness.\\
\underline{Existence:} We look for the solution in the form $S= \tilde{S}+\khii \psi $ with $\tilde{S} \in \Vcoin$ and $\psi \in \mathcal{A}(\Pi )$. More precisely we denote (using notation $\langle .\rangle $ from \cref{def: serie geometrique d'op 1D})
\[\varphi ^{-} := \left\langle -\RcouchePi\circ \partial _{X|\couche} ^2 ,\; \frac{\mu _{0} }{\mu _{1}} \RGammaPi\circ \partial _{Y|\Gamma ,Y=0^{+}} \right\rangle  \bigg( \frac{1}{\mu _{0}} \ROmegaPi (\varphi _{|\Omega } )+ \frac{1}{\mu _{1}} \RcouchePi (\varphi _{|\couche} ) + \fsum_{d\in \TN} s_{d} \,\phi_{d} \bigg) \,\in \, \Amicro(\Pi )\]
and $\psi :=T_{\geqslant -2} (\varphi ^{-} )$. By~\cref{pb de Poisson dans A Pi}, we have
\[\left\{\begin{array}{r@{\;=\;}ll}
	\mu _{0} \Delta  \varphi ^{-} & \varphi  &\text{ in }\Omega  \\[0.5mm]
	\mu _{1} \Delta  \varphi ^{-} & \varphi  &\text{ in }\couche \\[0.5mm]
	[\varphi ^{-} ]_{\Gamma } & 0 &\text{ on }\Gamma \\[0.5mm]
	[\mu \partial _{Y} \varphi ^{-} ]_{\Gamma } & 0 &\text{ on }\Gamma  \\[0.5mm]
	\varphi ^{-} & 0 &\text{ on }\bordgauche \cup \borddroit\\[0.5mm]
	\sigma _{d} (\varphi ^{-} ) & s_{d} & \ \forall d\in \TN
\end{array}\right.\]
Similarly to the proof of \cref{cadre fonctionnel des champs lointains}, one can check that it implies
\[\left\{\begin{array}{r@{}l}
	(\mu _{0} \Delta  \psi  -\varphi )_{|\Omega } &{}\in  \mathcal{A}(\Omega )\cap  \fsum\limits_{d<-4} \mathcal{A}_{d} (\Omega ) = \sum\limits _{d<-4} \mathcal{A}_{d} (\Omega ) \\[0.5mm]
	(\mu _{1} \Delta  \psi  -\varphi )_{|\couche} &{}\in  \mathcal{A}(\couche)\cap  \fsum\limits_{d<-2} \mathcal{A}_{d} (\couche) = \sum\limits _{d<-2} \mathcal{A}_{d} (\couche) \\[-0.5mm]
	[\psi ]_{\Gamma } &{}= 0\\[0.5mm]
	[\mu \partial _{Y} \psi ]_{\Gamma } &{}\in  \mathcal{A}(\Gamma )\cap  \fsum\limits_{d<-2} \mathcal{A}_{d} (\Gamma ) = \sum\limits _{d<-2} \mathcal{A}_{d} (\Gamma ) \\[-0.5mm]
	\psi _{|\bordgauche \cup \borddroit} &{}= 0 \\[0.5mm]
	\sigma _{d} (\psi )&{}= s_{d} \;\ \forall d\in \TN
\end{array}\right.\]
where $[{\dots }]_{\Gamma }$ stands for the jump on $\Gamma $. Therefore we have $(1+r)\,\big[\!\operatorname{div}\!\big(\mu  \nabla  (\khii \psi )\big)-\khii \varphi \big] \in L^2 (\Omega _{1} )$ and $[\mu \partial _{Y} (\khii \psi )]_{\Gamma } \in L^2 (\Gamma \cap \{X>\Rcoin\})$. Finally, thanks to \cref{pb variationnel de Poisson dans Omega1}, we take $\tilde{S}$ as the unique solution in $V$ of: \[\left\{\begin{array}{r@{\;=\;}ll}
	\operatorname{div}(\mu  \nabla  \tilde{S}) & F+\khii \varphi - \operatorname{div}\!\big(\mu  \nabla  (\khii \psi )\big)&\text{ in }\Omega _{1} \setminus (\Gamma \cap \{X>\Rcoin\}) \\[0.5mm]
	\tilde{S} &0 &\text{ on }\partial \Omega _{1} \\[0.5mm]
	[\mu \partial _{Y} \tilde{S}]_{\Gamma } & -[\mu \partial _{Y} (\khii \psi )]_{\Gamma } &\text{ on }\Gamma \cap \{X>\Rcoin\} 
\end{array}\right.\]

\underline{Uniqueness:} Let $S$ be a homogeneous solution and let us show that $S=0$. There exist $\tilde{S} \in \Vcoin$ and $\psi \in \mathcal{A}(\Omega )$ s.t.\ $S= \tilde{S}+\khii \psi $. The proof decomposes into the following steps, all proven by contradiction.
\begin{enumerate}
	\item \underline{$\maxdeg(\Delta \psi _{|\Omega } )<-2$\,:} Like in step 1 of the uniqueness proof of \cref{cadre fonctionnel des champs lointains}, we test the inequality:
	\[\int _{\Omega } \Delta (\khii \psi )\cdot \zeta  =-\int _{\Omega } \Delta  \tilde{S} \cdot \zeta  =\int _{\Omega } \nabla  \tilde{S} \cdot \nabla  \zeta  \lesssim  \|\zeta \|_{V} \qquad  \forall \zeta  \in \mathcal{D}(\Omega )\]
	with $\zeta _{q} :(r,\theta ) \mapsto  \Delta (\khii \psi )(r,\theta )\cdot r^{-q-2d-2} \,\khiz(2^{-1/q} r)\,\chi (\theta )$ when $q\rightarrow 0^{+}$, where $d:=\maxdeg(\Delta \psi _{|\Omega } )$ and $\chi \in \mathcal{D}(0,\Theta ) \setminus \{0\}$ is everywhere non negative.
	\item \underline{$\maxdeg(\Delta \psi _{|\couche} )<- \frac{1}{2}$\,:} This step works like the previous one.
	\item \underline{$\maxdeg [\mu \partial _{Y} \psi ]_{\Gamma } <- \frac{1}{2}$\,:} Let $d:= \maxdeg  [\mu \partial _{Y} \psi ]_{\Gamma }$ and $\varphi := [\mu \partial _{Y} (\khii \psi )]_{\Gamma }$.	We have $\nabla  \tilde{S} \in L^2 (\Omega _{1} )$ and, by steps 1 and 2, $\operatorname{div}(\mu  \nabla  \tilde{S}) =-\operatorname{div}\!\big(\mu  \nabla  (\khii \psi )\big) \in L^2 (\Omega _{1} )$. Thus $\varphi  =-[\mu \partial _{Y} \tilde{S}]_{\Gamma } \in H^{-1/2} (\Gamma )$. Like previously, we test it with $\zeta _{q} :x \mapsto  \varphi (x) \cdot x^{-q-2d-1} \,\khiz(2^{-1/q} x)$ when $q\rightarrow 0^{+}$. If $d\geqslant  - \frac{1}{2}$, we get $\langle  \varphi , \zeta _{q} \rangle _{\Gamma } \gtrsim  \frac{1}{q}$ and $\|\zeta _{q} \|_{H^{1/2} (\Gamma )} \lesssim  \|\zeta _{q} \|_{H^{1} (\Gamma )} \lesssim  \frac{1}{\sqrt{q}}$, which contradicts $\langle  \varphi , \zeta _{q} \rangle _{\Gamma } \lesssim  \|\zeta _{q} \|_{H^{1/2} (\Gamma )}$.
	\item \underline{$\maxdeg(\psi _{|\couche} )<- \frac{1}{2}$\,:} Let $d:=\maxdeg(\psi _{|\couche} )$. There exist $\psi _{1} \in \mathcal{A}_{d} (\couche) \setminus \{0\}$ and $\psi _{2} \in \sum _{q<d} \mathcal{A}_{q} (\couche)$ s.t.\ $\psi _{|\couche} =\psi _{1} +\psi _{2}$. By~\cref{operateurs de derivation dans A}, $\deg (\partial _{x|\couche} ^2 )=-2$, $\deg (\partial _{Y|\Gamma ,Y=0^{+}} )=-1$ and $\deg (\partial _{Y|\couche} ^2 )=\deg (\partial _{Y|\Gamma ,Y=0^{-}} )=0$. In addition, \cref{decomposition de A en sev} states that $\mathcal{A}(D) = \bigoplus_{q\in \mathbb{R}} \mathcal{A}_{q} (D)$ for~any $D\in \{\couche,\Gamma \}$. So by~taking the coordinates of $\Delta \psi _{|\couche}$ in $\mathcal{A}_{d} (\couche)$ and of $[\mu \partial _{Y} \psi ]_{\Gamma }$ in $\mathcal{A}_{d} (\Gamma )$, we get if $d\geqslant  - \frac{1}{2}$: \[\left\{\begin{array}{r@{\;=\;}ll}
		\partial _{Y} ^2  \psi _{1} & 0 &\text{ in }\couche \\[0.5mm]
		\partial _{Y|\Gamma ,Y=0^{-}} \psi _{1} & 0 &\text{ on }\Gamma 
	\end{array}\right.\]
	But since $\psi _{|\borddroit} =0$, it implies $\psi _{1} =0$, which is contradictory.
	\item \underline{$\maxdeg(\psi _{|\Omega } )<0$\,:} We denote $d:=\maxdeg(\psi _{|\Omega } )$. There exist $\psi _{1} \in \mathcal{A}_{d} (\Omega ) \setminus \{0\}$ and $\psi _{2} \in \sum _{q<d} \mathcal{A}_{q} (\Omega )$ s.t.\ $\psi _{|\Omega } =\psi _{1} +\psi _{2}$. Since $\deg (\Delta _{|\Omega } )=-2$, taking the coordinate of $\Delta \psi _{|\Omega }$ in $\mathcal{A}_{d-2} (\Omega )$ and of $\psi _{|\Gamma }$ in $\mathcal{A}_{d} (\couche)$ gives if $d\geqslant 0$: \[\left\{\begin{array}{r@{\;=\;}ll}
		\Delta  \psi _{1} & 0 &\text{ in }\Omega  \\[0.5mm]
		\psi _{1} & 0 &\text{ on }\Gamma 
	\end{array}\right.\]
	By~\cref{pb de Laplace homogene dans A}, it implies $d\in \TN$ and $\psi _{1} =\sigma _{d} (\psi _{1} )\phi_{d}$. But $\sigma _{d} (\psi _{1} )=\sigma _{d} (\psi )=\sigma _{d} (S)=0$.
\end{enumerate}
Finally, \cref{V inter A} implies $\khii \psi \in V$, so $S\in V$. So the uniqueness in \cref{pb variationnel de Poisson dans Omega1} implies $S=0$.
\end{demo}


\section{Asymptotic behaviors w.r.t.\ $\bfr$}
\label{sec: DA autour du coin}

\subsection{Statement of the results}

\begin{textesansboite}
In this section, we prove the existence of an asymptotic expansion w.r.t.\ $\bfr$ in $\Amaicro(\Pi )$ for solutions of the model problems of \cref{sec: pb modeles bien poses}. This new result will be applied in \cref{sec: construction des champs} to inductively build the fields $u_{p,\ell }$, $U_{p,\ell }$ and $S_{p,\ell }$ and describe their asymptotic behavior w.r.t.\ $\bfr$. We recall that we are interested in the behavior of the far-and-layer fields when $\bfr\rightarrow 0$, and of the corner fields when $\bfr\rightarrow \infty $.\\

\Crefrange{eq: eq vol champs lointains}{eq: eq vol champs de coin} show that the fields $u_{p,\ell }$, $U_{p,\ell }$ and $S_{p,\ell }$ depend inductively on the previous fields and some derivatives of the previous fields. Hence it will be useful to establish asymptotic expansions of the fields that still hold when differentiated. To this end, we introduce the notation $o_{\partial }$, which is a kind of differentiable small $o$.
\end{textesansboite}

\begin{definition}{the symbol $o_{\partial }$}
\label{def: petit o derivable}%
\def\varphi{f}%
Let $d\in \mathbb{R}$ and $a\in \{0,\infty \}$. 
\begin{itemize}
	\item For~any $\varphi :\Omega  \rightarrow \mathbb{C}$, we say that $\varphi  = o_{\partial } (r^{d} )$ when $r\rightarrow a$ if $\varphi $ is $\mathcal{C}^{\infty }$ in a vicinity of $r=a$ and $\forall (j,k)\in \mathbb{N}^2 ,\; \partial _{r} ^{j} \partial _{\theta } ^{k} \varphi  = o(r^{d-j} )$ uniformly w.r.t.\ $\theta $ when $r\rightarrow a$.
	\item For~any $\varphi :\couche \rightarrow \mathbb{C}$, we say that $\varphi  = o_{\partial } (x^{d} )$ when $x\rightarrow a$ if $\varphi $ is $\mathcal{C}^{\infty }$ in a vicinity of $x=a$ and $\forall (j,k)\in \mathbb{N}^2 ,\; \partial _{x} ^{j} \partial _{Y} ^{k} \varphi  = o(x^{d-j} )$ uniformly w.r.t.\ $Y$ when $x\rightarrow a$.
	\item For~any $\varphi :\Gamma \rightarrow \mathbb{C}$, we say that $\varphi  = o_{\partial } (x^{d} )$ when $x\rightarrow a$ if $\varphi $ is $\mathcal{C}^{\infty }$ in a vicinity of $x=a$ and $\forall j\in \mathbb{N},\; \partial _{r} ^{j} \varphi  = o(x^{d-j} )$ when $x\rightarrow a$.
	\item For~any $\varphi :\Pi  \rightarrow \mathbb{C}$, we say that $\varphi  = o_{\partial } (\bfr^{d} )$ when $\bfr\rightarrow a$ if $\varphi _{|\Omega } = o_{\partial } (r^{d} )$ and $\varphi _{|\couche} = o_{\partial } (x^{d} )$.
\end{itemize}\fbi
\end{definition}
\nomenclature[24-13]{$o_{\partial }$}{Asymptotic small $o$ over functions and their derivatives at any order}

\begin{textesansboite}
The rigorous formulation of the sentence ``\emph{$f^{0} \in \Amacro(\Pi )$ describes the behavior at the corner of the function $f$}'' that will be used from now on is 
\[\text{``$\.$for any $q\in \mathbb{R},\quad \ f=T_{\leqslant q} (f^{0} ) +o_{\partial } (\bfr^{q} )$ \ when $\bfr\rightarrow 0\.$''} .\]
Similarly for behaviors at infinity. The next two lemmas give some useful basic properties of this type of asymptotic behavior.
\end{textesansboite}

\begin{lemmesansboite}{}
\label{derivation des DA dans Amaicro}%
Let $\delta $ be a differential operator among $\Delta _{|\Omega }$, $\partial _{x|\Gamma } ^2 $, $\partial _{x|\couche} ^2 $, $\partial _{y|\Gamma ,y=0^{+}}$, $\partial _{Y|\couche} ^2 $ and $\partial _{Y|\Gamma ,Y=0^{-}}$. Let $D$ be the domain of definition of $\delta $ and let $f:D\rightarrow \mathbb{C}$.
\begin{enumerate}
	\item If $f^{0} \in \Amacro(D)$ satisfies $f \! \underset{\bfr\rightarrow 0}{=} \!T_{\leqslant q} (f^{0} )+o_{\partial } (\bfr^{q} )$ for~any $q\in \mathbb{R}$, then we have $\delta (f) \!\underset{\bfr\rightarrow 0}{=} \!T_{\leqslant q} (\delta (f^{0} ))+o_{\partial } (\bfr^{q} )$ for~any $q\in \mathbb{R}$.
	\item If $f^{\infty } \in \Amicro(D)$ satisfies $f \!\underset{\bfr\rightarrow \infty }{=} \!T_{\geqslant q} (f^{\infty } )+o_{\partial } (\bfr^{q} )$ for~any $q\in \mathbb{R}$, then we have $\delta (f) \!\underset{\bfr\rightarrow \infty }{=} \!T_{\geqslant q} (\delta (f^{\infty } ))+o_{\partial } (\bfr^{q} )$ for~any $q\in \mathbb{R}$.
\end{enumerate}
\end{lemmesansboite}

\begin{textesansboite}
In other words, if $\Phi \in \Amaicro$ gives the behavior of $f$, then $\delta (\Phi )$ gives the behavior of $\delta (f)$ (with $\Phi :=f^{0}$, resp.\ $\Phi :=f^{\infty ^{}}$).
\end{textesansboite}

\begin{demo}
Let us focus on the first point. By~\cref{operateurs de derivation dans A}, $\delta $ maps $\mathcal{A}(D)$ to $\mathcal{A}(D')$ for some $D'\in \{\Omega ,\Gamma ,\couche\}$ and has a degree (see \cref{def: degre d'un operateur}). Hence $\delta (f^{0} )$ is well-defined in $\Amacro(D')$ and the claimed statement makes sense. Let $d:=\deg  \delta $. It is easy to see that $\delta $ is a linear combination of terms of the form $r^{i} \partial _{r} ^{j} \partial _{\theta } ^{k}$ if $D=\Omega $, and $x^{i} \partial _{x} ^{j} \partial _{Y} ^{k}$ if $D\in \{\Gamma ,\couche\}$, with $(i,j,k)\in \mathbb{Z}\times \mathbb{N}^2 $ and $i-j=d$. Thus for~any $q\in \mathbb{R}$, we have $\delta (f)=\delta (T_{\leqslant q} (f))+o_{\partial } (\bfr^{q-d} ) = T_{\leqslant q-d} (\delta (f))+o_{\partial } (\bfr^{q-d} )$.
\end{demo}

\begin{textesansboite}
The next lemma shows that, when an asymptotic behavior exists in a space $\Amaicro$, it~is unique.
\end{textesansboite}

\begingroup
\def\fun{f}%
\begin{lemmesansboite}{}
\label{unicité du DA}%
Let $D\in \domsA$ and $\fun:D\rightarrow \mathbb{C}$.
\begin{itemize}
	\item Let $\fun^{0} _{1} ,\fun^{0} _{2} \in \Amacro(D)$ be two asymptotic behaviors of $\fun$ when $\bfr\rightarrow 0$, i.e.\ we have: $\forall d\in \mathbb{R},\; \fun =T_{\leqslant d} (\fun^{0} _{1} )+o_{\partial } (\bfr^{d} ) \text{ and } \fun =T_{\leqslant d} (\fun^{0} _{2} )+o_{\partial } (\bfr^{d} )$. Then $\fun^{0} _{1} =\fun^{0} _{2}$.
	\item Let $\fun^{\infty } _{1} ,\fun^{\infty } _{2} \in \Amicro(D)$ be two asymptotic behaviors of $\fun$ when $\bfr\rightarrow \infty $, i.e.\ we have: $\forall d\in \mathbb{R},\; \fun =T_{\geqslant d} (\fun^{\infty } _{1} )+o_{\partial } (\bfr^{d} ) \text{ and } \fun =T_{\geqslant d} (\fun^{\infty } _{2} )+o_{\partial } (\bfr^{d} )$. Then $\fun^{\infty } _{1} =\fun^{\infty } _{2}$.
\end{itemize} 
\end{lemmesansboite}

\begin{demo}
Let us prove only the first point. We assume by contradiction $f^{0} _{1} \neq f^{0} _{2} $. Then $d:=\min  \operatorname{supp}(f^{0} _{1} -f^{0} _{2} )$ is well-defined by construction of $\Amacro(D)$ and we have $\varphi := T_{\leqslant d} (f^{0} _{1} -f^{0} _{2} )\in \mathcal{A}_{d} (D) \setminus \{0\}$. Subtracting $\fun =T_{\leqslant d} (\fun^{0} _{i} )+o (\bfr^{d} )$ for $i=1,2$, we get $\varphi =o(\bfr^{d} )$. In addition, by \cref{croissance des fcts de A}, there exist $r_{0} >0$ and $C>0$ s.t.\ (if $D=\Omega $ e.g., but it~is similar in the other cases): $\forall r>r_{0} ,\,\forall \theta ,\, |\varphi (r,\theta )| \geqslant C\,r^{d}$. This contradicts $\varphi  = o(\bfr^{d} )$ and ends the proof.
\end{demo}
\endgroup

\begin{textesansboite}
The two following theorems are the main results of this section. They will be applied to the far-and-layer fields, resp.\ the corner fields.
\end{textesansboite}

\begin{theoremesecable}{asymptotic expansion of the solution of a far-and-layer field-like problem}
\label{DA des champs lointains concatenes}%
Let $\bfu\in H^{1} _{\mathrm{loc}} (\Pi )$, $f\in \mathcal{D}'(\Omega )$, $g\in L^2 _{\mathrm{loc}} (\couche)$ and $h\in L^2  _{\mathrm{loc}} (\Gamma )$ be s.t.
\begin{equation}{\label{eq: DA des CL concatenes : hyp}}
\left\{\begin{array}{r@{\;=\;}ll}
	\mu _{0} \Delta  \bfu +\omega ^2 \rho _{0} \bfu & f &\text{ in }\Omega  \\[0.5mm]
	\mu _{1} \partial _{Y} ^2  \bfu & g &\text{ in }\couche \\[0.5mm]
	\mu _{1} \partial _{Y|Y=0^{-}} \bfu & h &\text{ on }\Gamma  \\[0.5mm]
	\bfu_{|y=0^{+}} - \bfu_{|Y=0^{-}} & 0 &\text{ on }\Gamma \\[0.5mm]
	\bfu & 0 &\text{ on }\bordgauche \cup \borddroit
\end{array}\right.\end{equation}
We assume that:
\begin{itemize}
	\item $f$ vanishes in the vicinity of the corner $(0,0)$,
	\item there exists $g^{0} \in \Amacro(\couche)$ s.t.: $\forall d\in \mathbb{R},\; g(x,Y) = T_{\leqslant d} (g^{0} )(x,Y) + o_{\partial } (x^{d} )$ when $x\rightarrow 0$,
	\item there exists $h^{0} \in \Amacro(\Gamma )$ s.t.: $\forall d\in \mathbb{R},\; h(x) = T_{\leqslant d} (h^{0} )(x) + o_{\partial } (x^{d} )$ when $x\rightarrow 0$,
	\item there exist $\eta >0$, $u^{\mathrm{v}} \in H^{1} (\Omega \cap B(0,\eta ))$ and $\varphi \in \mathcal{A}(\Omega )$ s.t.\ $\bfu_{|\Omega \cap B(0,\eta )} =u^{\mathrm{v}} + \varphi $.
\end{itemize}
Then there exists $\bfuDA \in \Amacro(\Pi )$ that has the form
\[\bfuDA = \big\langle  {-}k_{0} ^2 \.\ROmegaPi \big\rangle  \bigg( \frac{1}{\mu _{1}} \RcouchePi(g^{0} ) + \frac{1}{\mu _{1}} \RGammaPi(h^{0} ) + \fsum_{d\in \TZ*} \sigma _{d} (\bfuDA) \,\phi_{d} \bigg)\]
s.t.\ we~have:
\[\forall d\in \mathbb{R}, \qquad  \bfu = T_{\leqslant d} (\bfuDA ) + o_{\partial } (\bfr^{d} ) \quad  \text{when }  \bfr\rightarrow 0.\]\fbeq[1.1]
\end{theoremesecable}

\begin{textesansboite}
Note that, since $g^{0} \in \Amacro(\couche)$, we have $T_{\leqslant d} (g^{0} )\in \mathcal{A}(\couche) \subset \mathcal{C}^{0} (\couche,\mathbb{C})$ for~any $d\in \mathbb{R}$, so the formula ``$g = T_{\leqslant d} (g^{0} ) +o_{\partial } (x^{d} )$'' makes sense. The same applies to the truncations of $h^{0}$ and $\bfu^{0}$, and of $F^{\infty } _{\Omega }$, $F^{\infty } _{\couche}$ and $\SDA$ in the next result.
\end{textesansboite}

\begin{theoremesecable}{asymptotic expansion of the solution of a corner field-like problem}
\label{DA des champs de coin}%
Let $S \in H^{1} _{\mathrm{loc}} (\Omega _{1} )$ and $F\in L^2 _{\mathrm{loc}} (\overline{\Omega _{1} })$ (i.e.\ $L^2 $ on any bounded subset of $\Omega _{1}$) s.t. 
\begin{equation}{\label{eq: DA des champs de coin, hyp}}
\left\{\begin{array}{r@{\;=\;}ll}
	\operatorname{div}(\mu  \nabla  S) & F &\text{ in }\Omega _{1} \\[0.5mm]
	S & 0 &\text{ on }\partial \Omega _{1}
\end{array}\right.\end{equation}
We assume that:
\begin{itemize}
	\item there exists $F^{\infty } _{\Omega } \in \Amicro(\Omega )$ such that $\forall d\in \mathbb{R},\; F_{|\Omega } = T_{\geqslant d} (F^{\infty } _{\Omega } ) + o_{\partial } (r^{d} )$ when $r\rightarrow \infty $,
	\item there exists $F^{\infty } _{\couche} \in \Amicro(\couche)$ such that $\forall d\in \mathbb{R},\; F_{|\couche} = T_{\geqslant d} (F^{\infty } _{\couche} ) + o_{\partial } (x^{d} )$ when $x\rightarrow \infty $,
	\item $S$ belongs to $\Vcoin+\khii \mathcal{A}(\Pi )$ (cf.\ \cref{def: V}).
\end{itemize}
Then there exists $\SDA \in \Amicro(\Pi )$ that has the form
\[\SDA = \left\langle {-\RcouchePi\circ \partial _{X|\couche} ^2 ,\; \frac{\mu _{0} }{\mu _{1}} \RGammaPi\circ \partial _{Y|\Gamma ,Y=0^{+}} }\right\rangle  \bigg( \frac{1}{\mu _{0}} \ROmegaPi (F^{\infty } _{\Omega } )+ \frac{1}{\mu _{1}} \RcouchePi (F^{\infty } _{\couche} ) + \fsum_{d\in \TZ*} \sigma _{d} (\SDA)\,\phi_{d} \bigg) .\]
s.t.\ we~have:
\[\forall d\in \mathbb{R}, \qquad  S = T_{\geqslant d} (\SDA ) +o_{\partial } (\bfr^{d} ) \quad  \text{when }  \bfr\rightarrow \infty .\]\fbeq[1.1]
\end{theoremesecable}

\begin{textesansboite}
We prove \cref{DA des champs de coin} in \cref{sec: DA des champs de coin} and \cref{DA des champs lointains concatenes} in \cref{sec: DA des champs lointains concatenes}. The second part is less detailed as it has similarities with the first part.
\end{textesansboite}

\subsection{Proof of \texorpdfstring{\cref{DA des champs de coin}}{Theorem 5.2.5}: asymptotic behavior for corner fields-like problems}
\label{sec: DA des champs de coin}

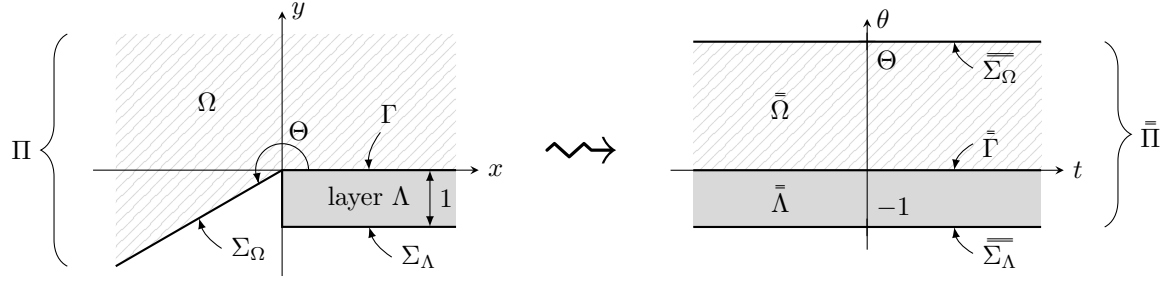
\begin{figure}[h]
\def\LwidthFig{2.2} 
\def\RwidthFig{2.3} 
\def\heightFig{1.8} 
\def\depthFig{-1.4} 
\def\epsFig{0.75} 
\def\ThetaFig{210} 
\begin{tikzpicture}
	\coordinate (finLambda) at (-\LwidthFig, {-\LwidthFig*tan(\ThetaFig-180)});
	\fill [pattern={Lines[angle=45,line width=.2pt]}, pattern color=white!70!gray] (0,0) -- (\RwidthFig,0) -- (\RwidthFig,\heightFig) -- (-\LwidthFig,\heightFig) -- (finLambda) -- cycle;
	\fill [color=white!70!gray] (0,0) rectangle (\RwidthFig,-\epsFig);
	\draw [thick] (finLambda) -- (0,0) -- (\RwidthFig,0);
	\draw [thick] (\RwidthFig,-\epsFig) -- (0,-\epsFig) -- (0,0);
	\draw [>=stealth,->] (-\LwidthFig-0.3,0) -- (\RwidthFig+0.3,0) node[right]{$x$};
	\draw [>=stealth,->] (0,\depthFig) -- (0,{\heightFig+0.3}) node[right]{$y$};
	\draw [decorate, decoration={brace,amplitude=10pt}, xshift=-4pt] (-\LwidthFig-0.5,{-\LwidthFig*tan(\ThetaFig-180)}) -- (-\LwidthFig-0.5,\heightFig) node [midway,left,xshift=-10pt] {$\Pi$};
	\draw (-\LwidthFig*0.45,\heightFig/2) node {$\Omega$};
	\draw [->,>=latex] (0.35,0) arc (0:\ThetaFig-5:0.35) node[near start,above] {$\Theta$} -- (\ThetaFig:0.36);
	\draw [<-,>=latex] (-\LwidthFig/2, {-\LwidthFig/2*tan(\ThetaFig-180)}) to[bend right=20] ++(0.3,-0.4) node[right]{$\bordgauche$};
	\draw (\RwidthFig/2,-\epsFig/2) node {layer $\couche$};
	\draw [<->,>=latex] (\RwidthFig*0.85,0) -- node[midway,right]{1} (\RwidthFig*0.85,-\epsFig);
	\draw [<-,>=latex] (\RwidthFig/2,0) to[bend left=20] ++(0.15,0.5) node[above]{$\ \ \Gamma$};
	\draw [<-,>=latex] (\RwidthFig/2,-\epsFig) to[bend right=20] ++(0.3,-0.4) node[right]{$\borddroit$};
	\draw (\RwidthFig+1,0.2) node[right] {\scalebox{3}{$\rightsquigarrow$}};
\end{tikzpicture}%
\hspace{3mm}\nolinebreak%
\def\widthFig{2.3} 
\def\heightFig{1.7} 
\begin{tikzpicture}
	\fill [color=white!70!gray] (-\widthFig,0) rectangle (\widthFig,-\epsFig);
	\fill [pattern={Lines[angle=45,line width=.2pt]}, pattern color=white!70!gray] (-\widthFig,0) rectangle (\widthFig,\heightFig);
	\draw [thick] (-\widthFig,\heightFig) -- (\widthFig,\heightFig);
	\draw [thick] (-\widthFig,0) -- (\widthFig,0);
	\draw [thick] (-\widthFig,-\epsFig) -- (\widthFig,-\epsFig);
	\draw [>=stealth,->] (-\widthFig-0.3,0) -- (\widthFig+0.3,0) node[right]{$t$};
	\draw [>=stealth,->] (0,{-\epsFig-0.3}) -- (0,{\heightFig+0.3}) node[right]{$\theta$};
	\draw (0,\heightFig) node{+} node[below right]{$\Theta$};
	\draw (0,-\epsFig) node{+} node[above right]{$-1$};
	\draw [decorate, decoration={brace,amplitude=10pt}, xshift=4pt] (\widthFig+0.7,\heightFig) -- (\widthFig+0.7,-\epsFig) node [midway,right,xshift=10pt] {$\sbd\Pi$};
	\draw ({-\widthFig/2},\heightFig/2) node {$\sbd\Omega$};
	\draw ({-\widthFig/2},-\epsFig/2) node {$\sbd\couche$};
	\draw [<-,>=latex] (\widthFig/2,\heightFig) to[bend right=25] ++(0.25,-0.35) node[right]{$\sbd\bordgauche$};
	\draw [<-,>=latex] (\widthFig/2,0) to[bend left=25] ++(0.25,0.35) node[right]{$\sbd\Gamma$};
	\draw [<-,>=latex] (\widthFig/2,-\epsFig) to[bend right=25] ++(0.25,-0.35) node[right]{$\sbd\borddroit$};
\end{tikzpicture}\par
\caption{The change of variables $(x,y) \rightsquigarrow (t,\theta )$ and the associated domains.}
\label{fig: sbd Pi}
\end{figure}

\begin{textesansboite}
The proof relies on the Kondrat'ev theory, usually used to analyse singularities of solutions of elliptic equations, see \cite{Kon67, KozMazRos97, MazNazPla00, CalCosDauVia06, Che12, Mak08}. We use it in a way that gives an expansion in $\Amacro(\Pi )$.\\

First we introduce the variables $(t,\theta )$, defined as $(\ln r,\theta )$ in $\Omega $ and $(\ln x,Y)$ in $\couche$. The pair $(t,\theta )$ lies in $\sbd \Pi  := \mathbb{R}\times (-1,\Theta )$. Moreover, we denote $\sbd \Omega $, $\sbd \couche$, $\sbd\bordgauche$, $\sbd \Gamma $ and $\sbd\borddroit$ the images of $\Omega $, $\couche$, $\bordgauche$, $\Gamma $ and $\borddroit$ by the change of variable $(x,y) \rightsquigarrow (t,\theta )$, see \cref{fig: sbd Pi}. The notation $\sbd{\raisebox{0pt}[1.5ex]{\dots }}$ is intended to remind the strip shape of $\sbd \Pi $. Finally, for~any $u:\Pi \rightarrow \mathbb{C}$, we denote \[\sbd u:\, (t,\theta )\in \sbd \Pi  \;\mapsto \; \left\{\begin{array}{ll}
u(r=e^{t} ,\theta ) &\text{ in }\sbd \Omega \\[0.5mm]
u(x=e^{t} ,Y=\theta ) &\text{ in }\sbd \couche
\end{array}\right.\]\fbeq
\end{textesansboite}

\begin{definition}{Kondrat'ev spaces}
\label{def: espaces de Kondrat'ev}%
Let $H$ be a Hilbert space, $m\in \mathbb{N}$ and $\beta  \in \mathbb{R}$. We define
\[K^{m} _{\beta } (\mathbb{R},H) := \{t\mapsto  e^{\beta t} u(t) \mid u\in H^{m} (\mathbb{R},H)\}\]
where $H^{m} (\mathbb{R},H)$ denotes the Sobolev space $\{v:\mathbb{R} \rightarrow H \mid \forall n\in [\![0,m]\!], \partial _{t} ^{n} v\in L^2 \}$. We equip $K^{m} _{\beta } (\mathbb{R},H)$ with the norm $\|u\|_{K^{m} _{\beta } (\mathbb{R},H)} := \|t\mapsto  e^{-\beta t} u(t)\|_{H^{m} (\mathbb{R},H)}$ where we classically have $\|v\|_{H^{m} (\mathbb{R},H)} ^2 :=\sum _{n=0} ^{m} \|\partial _{t} ^{n} v\|_{L^2 } ^2 $ for any $v$. We also denote
\[K^{\infty } _{\beta } (\mathbb{R},H) := \bigcap _{m\in \mathbb{N}} K^{m} _{\beta } (\mathbb{R},H) \qquad  \text{ and } \qquad  K^{m} _{\beta } (\sbd \Gamma ):=K^{m} _{\beta } (\mathbb{R},\mathbb{C})\]
(identifying the line $\sbd \Gamma $ with $\mathbb{R}$).
\end{definition}
\nomenclature[20-28-11]{$K^{m} _{\beta } (\mathbb{R},H)$}{In \cref{chap: construction du DA}, Kondrat'ev space of functions $\mathbb{R} \rightarrow H$ that have $m$ derivatives behaving roughly like $e^{\beta t}$}

\begin{remarques}
\begin{itemize}
	\item Kondrat'ev spaces are decreasing w.r.t.\ $m$ (i.e.\ $m\leqslant m' \Rightarrow  K^{m'} _{\beta } (\mathbb{R},H) \subset  K^{m} _{\beta } (\mathbb{R},H)$) but are not comparable w.r.t.\ $\beta $ (i.e.\ there~is no inclusion between $K^{m} _{\beta } (\mathbb{R},H)$ and $K^{m} _{\beta '} (\mathbb{R},H)$).
	\item $K^{m} _{\beta } (\mathbb{R},H) = \{u\in H^{m} _{\mathrm{loc}} (\mathbb{R},H) \mid \forall k\in [\![0,m]\!],\; e^{-\beta t} \partial _{t} ^{k} u(t) \in L^2 (\mathbb{R},H)\}$.
	\item Note that if $u\in K^{m} _{\beta } (\mathbb{R},H)$ and $u$ is zero in a vicinity of $-\infty $, then: $\forall \beta '\geqslant \beta ,\; u\in K^{m} _{\beta '} (\mathbb{R},H)$.
	\item Kondrat'ev's spaces are linked to $o_{\partial }$ (see \cref{def: petit o derivable}) in the following way. Let $\chi  \in \mathcal{C}^{\infty } (\mathbb{R})$ be equal to 0 in a vicinity of $-\infty $ and 1 in a vicinity of $+\infty $. Then for~any $u:\Omega \rightarrow \mathbb{C}$:
	\[\left\{\begin{array}{@{\;}l@{\quad\;}r@{\;\;\Rightarrow \;\;}l}
		\forall \beta \in \mathbb{R},&\chi \sbd{u}\in  \!\smash{\bigcap\limits _{m\in \mathbb{N}} } K^{\infty } _{\beta } (\mathbb{R},H^{m} (0,\Theta )) & \forall d>\beta ,\; u \! \underset{r\rightarrow \infty }{=} \! o_{\partial } (r^{d} )\\[0.5mm]
		\forall d\in \mathbb{R},&  u \! \underset{r\rightarrow \infty }{=} \! o_{\partial } (r^{d} ) & \exists t_{0} \in \mathbb{R},\,\forall \beta >d,\ \; \chi (\cdot +t_{0} )\,\sbd{u}\in  \!\smash{\bigcap\limits _{m\in \mathbb{N}} } K^{\infty } _{\beta } (\mathbb{R},H^{m} (0,\Theta ))
	\end{array}\right.\]
	And there are similar implications in $\couche$ and $\Gamma $.
\end{itemize}\fbi
\end{remarques}

\begin{textesansboite}
Let $H$ be a Hilbert space, $\beta \in \mathbb{R}$ and $\varphi \in K^{0} _{\beta } (\mathbb{R},H)$. For~any $\lambda \in \mathbb{C}$ s.t.\ $\Re(\lambda )=\beta $ we define the (bilateral) Laplace transform of $\varphi $ at $\lambda $ as\vspace{-1.5mm}
\begin{equation}{\label{eq: transformee de Fourier-Laplace}}
\widehat{\varphi }(\lambda ) := \int _{-\infty } ^{+\infty } e^{-\lambda t} \varphi (t)\,\mathrm{d}t = \mathcal{F}\big[t\mapsto  e^{-\beta t} \varphi (t)\big](\Im(\lambda ))\end{equation}
where $\mathcal{F}$ denotes the Fourier transform. By~properties of $\mathcal{F}$, we have $\widehat{\varphi } \in L^2 (\{\lambda \in \mathbb{C} \mid \Re(\lambda )=\beta \})$. If $\varphi $ depends on $t$ and $\theta $, $\widehat{\varphi }$ implies that we see $\varphi $ as a function from $\mathbb{R}$ to a space of functions of $\theta $.\\

To introduce the method, let us use the Laplace transform on $\sbd{S\,}$\!. Let $\khii$ be the truncation function introduced in \cref{sec: cadre fonctionnel des champs de coin} and $s:=\sbd{\khii S}$. Since $S\in \Vcoin+\khii \mathcal{A}(\Pi )$, there exists $\beta \in \mathbb{R}$ s.t.\ $s\in K^{0} _{\beta } (\mathbb{R},L^2 (-1,\Theta ))$, so $\hat{s}(\lambda ,\theta )$ is well-defined for~any $\lambda \in \mathbb{C}$ s.t.\ $\Re(\lambda )>\beta $. In addition we have for~any $f\in L^2 (\Omega )$, $g\in L^2 (\couche )$, $h\in H^{-1/2} (\Gamma )$ and $u\in H^{1} _{\mathrm{loc}} (\Pi )$:
\begin{equation}{\label{eq: transfert du probleme de Poisson de Pi a sbdPi}}
\left\{\begin{array}{@{\,}r@{\;\;}l@{}}
	-\mu _{0} \Delta  u =f &\text{ in }\Omega  \\[0.5mm]
	-\mu _{1} \Delta  u =g &\text{ in }\couche \\[0.5mm]
	u_{|y=0^{+}} - u_{|y=0^{-}} = 0 &\text{ on }\Gamma  \\[0.5mm]
	\mu _{0} \partial _{y} u_{|y=0^{+}} - \mu _{1} \partial _{y} u_{|y=0^{-}} = h &\text{ on }\Gamma  \\[0.5mm]
	u = 0 &\text{ on }\bordgauche \cup \borddroit
\end{array}\right. \;\Leftrightarrow \; \left\{\begin{array}{@{\,}r@{\;\;}l@{}}
	-e^{-2t} \mu _{0} \Delta  \sbd u = \sbd f &\text{ in }\sbd \Omega  \\[0.5mm]
	-\mu _{1} \,(e^{-2t} (\partial _{t} ^2 -\partial _{t} )+\partial _{\theta } ^2 ) \sbd u =\sbd g &\text{ in }\sbd\couche \\[0.5mm]
	\sbd u_{|\theta =0^{+}} - \sbd u_{|\theta =0^{-}} = 0 &\text{ on }\sbd \Gamma  \\[0.5mm]
	e^{-t} \mu _{0} \partial _{\theta } \sbd u_{|\theta =0^{+}} - \mu _{1} \partial _{\theta } \sbd u_{|\theta =0^{-}} = \sbd h &\text{ on }\sbd \Gamma  \\[0.5mm]
	\sbd u = 0 &\text{ on }\sbd\bordgauche \cup \sbd\borddroit
\end{array}\right.\end{equation}
Taking $u:=\khii S$ in \cref{eq: transfert du probleme de Poisson de Pi a sbdPi} and applying the Laplace transform yield that $\hat{s}$ satisfies:
\begin{equation}{\label{eq: pb verifie par la transformee de Laplace}}
\left\{\begin{array}{r@{\;=\;}ll}
	(\partial _{\theta } ^2 +\lambda ^2 ) \hat{s}(\lambda ,\theta ) & f & \text{if } \theta \in (0,\Theta )\\[0.5mm]
	\partial _{\theta } ^2  \hat{s}(\lambda ,\theta ) & g-\big[(\lambda +2)^2 -(\lambda +2)\big] \, \hat{s}(\lambda +2,\theta ) & \text{if } \theta \in (-1,0)\\[0.5mm]
	\hat{s}(\lambda ,0^{+} ) - \hat{s}(\lambda ,0^{-}) &0& \\[0.5mm]
	\mu _{1} \partial _{\theta } \hat{s}(\lambda ,0^{-}) & \mu _{0} \partial _{\theta } \hat{s}(\lambda +1,0^{+} ) & \\[0.5mm]
	\hat{s}(\lambda ,\Theta ) = \hat{s}(\lambda ,-1) & 0 & 
\end{array}\right.\end{equation}
for some functions $f$ and $g$ depending on $F$ and $\khii$. Solving this system w.r.t.\ $\theta $ allows us to extend $\hat{s}$ w.r.t.\ $\lambda $ further to the left in the complex plane, except at the $\lambda $ for which \cref{eq: pb verifie par la transformee de Laplace} is ill-posed. These $\lambda $ are poles of $\hat{s}$ and they will be used in \cref{etape de base du DA des champs de coin} to identify terms of the asymptotic expansion of $S$.\\

For~any $m\in \mathbb{N}^*$, let us define the Hilbert spaces $\mathcal{H}^{m} := \{u\in H^{1} _{0} (-1,\Theta ) \mid u_{|(0,\Theta )} \in H^{m} (0,\Theta ) \text{ and } u_{|(-1,0)} \in H^{m} (-1,0)\}$ with the norm $\|u\|_{\mathcal{H}^{m}} ^2  := \|u\|_{H^{m} (-1,0)} ^2 + \|u\|_{H^{m} (0,\Theta )} ^2 $, and $\prodH{m}:=H^{m} (0,\Theta )\times H^{m} (-1,0)\times \mathbb{C}$. We need to solve problems of the following form with $u\in \mathcal{H}^{m+2}$ and $(f,g,a)\in \prodH{m}$, $m\in \mathbb{N}$:
\begin{equation}{\label{eq: EDO en theta a resoudre}}
\left\{\begin{array}{r@{\;=\;}ll}
	u''+\lambda ^2 u & f &\text{ on }(0,\Theta ) \\[0.5mm]
	u'' & g &\text{ on }(-1,0) \\[0.5mm]
	u(0^{+} ) - u(0^{-}) &0& \\[0.5mm]
	u'(0^{-}) & a & \\[0.5mm]
	u(\Theta ) = u(-1) & 0 & 
\end{array}\right. \end{equation}
Denoting $A(\lambda ):u \mapsto  \big((u''+\lambda ^2 u)_{|(0,\Theta )} ,\, u''_{|(-1,0)} ,\, u'(0^{-}) \big)$, \cref{eq: EDO en theta a resoudre} is equivalent to $A(\lambda )(u) = (f,g,a)$. Note that for~any $m\in \mathbb{N}$, $A(\lambda ) \in  \mathcal{L}(\mathcal{H}^{m+2} ,\prodH{m})$ where $\mathcal{L}({\dots })$ denotes the space of continuous linear maps between two normed vector spaces.
\end{textesansboite}

\newcommand{\pole}{q}
\begin{lemmesecable}{properties of $A(\lambda )$}
\label{resolution de l'EDO en theta}%
Let $m\in \mathbb{N}$.
\begin{enumerate}
	\item $A(\lambda ):\mathcal{H}^{m+2} \rightarrow \prodH{m}$ is invertible iff $\lambda \in  \mathbb{C} \setminus \TZ*$.
	\item $\lambda \mapsto A(\lambda )^{-1}$ is meromorphic from $\mathbb{C}$ to $\mathcal{L}(\prodH{m},\mathcal{H}^{m+2} )$. Its poles are in $\TZ*$ and are simple. Moreover, for~any $\pole \in \TZ*$ and $v\in \prodH{m}$, the residue $(\Res_{\lambda =\pole} A(\lambda )^{-1} )(v)$ is proportional to the function $\theta  \in [-1,\Theta ] \mapsto  \mathbf{1}_{[0,\Theta ]} (\theta ) \cdot \sin (q\theta )$.
	\item Let $\beta _{1} < \beta _{2}$ be some reals. There exists $C>0$ depending only on $(m,\beta _{1} ,\beta _{2} )$ s.t., for~any $\lambda \in \mathbb{C}$ satisfying $\beta _{1} <\Re(\lambda )<\beta _{2}$ and $|\Im(\lambda )|>1$, we have:
	\begin{equation}{\label{eq: normes des operateurs resolvants R1 R2 et R3}}
		\|A(\lambda )^{-1} \|_{\mathcal{L}(\prodH{m},\mathcal{H}^{m+2} )} \leqslant  C\, |\Im(\lambda )|^{m+2} .
	\end{equation}
\end{enumerate}
\end{lemmesecable}


\begin{demo}
\begin{enumerate}
	\item Let $\lambda \in \mathbb{C}$. An easy calculation gives that any element of $\operatorname{Ker} A(\lambda )$ must be proportional to $\theta  \in [-1,\Theta ] \mapsto  \mathbf{1}_{[0,\Theta ]} (\theta ) \cdot \sin (\lambda \theta )$. This function belongs to $\mathcal{H}^{m} \setminus \{0\}$ iff $\lambda \in \TZ*$. Therefore $\operatorname{Ker} A(\lambda ) \neq 0 \Leftrightarrow  \lambda \in \TZ*$.\\
	Moreover, it~is easy to see that, for~any $(f,g,a)\in \prodH{m}$, \cref{eq: EDO en theta a resoudre} with $\lambda :=0$ has a unique solution in $\mathcal{H}^{m+2}$. Thus $A(0):\mathcal{H}^{m+2} \rightarrow \prodH{m}$ is invertible. Now, for~any $\lambda \in \mathbb{C}$, we have $A(\lambda )=A(0)+\lambda ^2 B$ with $B:u\mapsto (u_{|(0,\Theta )} ,0,0)$. $B$ is a compact operator from $\mathcal{H}^{m+2}$ to $\prodH{m}$, so the Fredholm alternative holds for $A(\lambda )$. Therefore $A(\lambda ):\mathcal{H}^{m+2} \rightarrow \prodH{m}$ is invertible iff $\lambda \in \mathbb{C} \setminus \TZ*$.

	\item Since $\lambda \mapsto A(\lambda )$ is holomorphic on $\mathbb{C} \setminus \TZ*$, so too is $\lambda \mapsto A(\lambda )^{-1}$. Let us describe its behavior near the points $q\in \TZ*$ using \cite[Theorem~5.1.1, p.\,147]{KozMazRos97}. It depends on the ``Jordan chains'' of $A(q)$, which are the sequences $(u_{0} ,\dots ,u_{n} )\in (\mathcal{H}^{m} )^{n+1}$, $n\in \mathbb{N}$, s.t.\ 
	\begin{equation}{\label{eq: chaine de Jordan}}
	\forall k\in [\![0,n]\!],\qquad  \sum _{j=0} ^{k} \frac{1}{j!} \frac{\mathrm{d}^{j} A}{\mathrm{d}\lambda ^{j}} \Big|_{\lambda =q} (u_{k-j} )= 0 \end{equation}
	Let us compute these chains. Taking $k:=0$ in \cref{eq: chaine de Jordan} gives $A(q)(u_{0} )=0$, so $u_{0}$ is proportional to $\theta \mapsto  \mathbf{1}_{[0,\Theta ]} (\theta ) \cdot \sin (q\theta )$ by step 1. In addition, if $n\geqslant 1$, taking $k:=1$ gives $A(q)(u_{1} )+ \frac{\mathrm{d}A}{\mathrm{d}\lambda } (q) (u_{0} )=0$. This implies on one hand $u_{1} '' +q^2 u_{1} +2qu_{0} =0$ in $(0,\Theta )$. On the other we get $u_{1|(-1,0)} '' =0$ and $u_{1} '(0^{-})=0$, so $u_{1|(-1,0)} =0$, which gives $u_{1} (0)=u_{1} (\Theta )=0$. Therefore
	\[0 \neq \int _{0} ^{\Theta } 2q|u_{0} |^2  = \int _{0} ^{\Theta } -(u_{1} ''+q^2 u_{1} ) \cdot  \overline{u_{0} }= \int _{0} ^{\Theta } -u_{1} \cdot (\underbrace{\overline{u_{0} ''+q^2 u_{0} }}_{=0} ) =0. \vspace{-2mm}\]
	This is absurd, so we must have $n=0$ for~any Jordan chain of $A(q)$. Therefore, \cite[Theorem~5.1.1, p.\,147]{KozMazRos97} states that $A(\lambda )^{-1}$ has a simple pole at $q$, and $\operatorname{Im}(\Res_{\lambda =q} A(\lambda )^{-1} )=\mathbb{C}u_{0}$.
	
	\item Let $(f,g,a)\in \prodH{m}$ and $u:=A(\lambda )^{-1} (f,g,a)\in \mathcal{H}^{m+2}$. We will write $\lesssim $ for inequalities valid up to a constant that depends on $(m,\beta _{1} ,\beta _{2} )$ but not on $(\lambda ,f,g,a)$. By~Poincaré's inequality and integration by parts, we have\vspace{-2mm}
	\[\|u\|_{H^{1} (-1,0)} ^2  \lesssim  \int _{-1} ^{0} |u'|^2  = - \int _{0} ^{-1} g \bar{u}+ a\, \bar{u}(0)\lesssim  (\|g\|_{L^2 (-1,0)} +|a|) \cdot  \|u\|_{H^{1} (-1,0)}\]
	thus $\|u\|_{H^{1} (-1,0)} \lesssim  \|g\|_{L^2 (-1,0)} +|a|$. Since $u''_{|(-1,0)} =g$, we deduce $\|u\|_{H^{m+2} (-1,0)} \lesssim  \|g\|_{H^{m}} +|a|$.\\
	Then, let $u_{1} :\theta  \in [0,\Theta ] \mapsto u(0) \big(1- \frac{\theta }{\Theta } \big)$, $v:=u-u_{1}$ and $f_{\lambda } :=f-\lambda ^2 u_{1}$. We have $v(0)=v(\Theta )=0$ and $v''+\lambda ^2 v=f_{\lambda }$, so
	\[\int _{0} ^{\Theta } (-|v'|^2 +\lambda ^2 |v|^2 )=\int _{0} ^{\Theta } f_{\lambda } \bar{v}\]
	Dividing by $\lambda $ and taking the absolute value of the imaginary part, we get:
	\[\frac{|\Im(\lambda )|}{|\lambda |^2 } \|v'\|_{L^2 } ^2  + |\Im(\lambda )|\, \|v\|_{L^2 } ^2  \leqslant  \frac{1}{|\lambda |} \|f_{\lambda } \|_{L^2 } \|v\|_{L^2 } .\]
	Now we assume $|\Im(\lambda )|>1$ and $\beta _{1} <\Re(\lambda )<\beta _{2}$, so $|\lambda |\lesssim |\Im(\lambda )|$. Thus $\|v'\|_{L^2 } ^2  + |\lambda |^2  \|v\|_{L^2 } ^2  \leqslant  \|f_{\lambda } \|_{L^2 } \|v\|_{L^2 }$, which yields $\|v\|_{L^2 } \lesssim  |\lambda |^{-2} \|f_{\lambda } \|_{L^2 }$ and then $\|v'\|_{L^2 } \lesssim  |\lambda |^{-1} \|f_{\lambda } \|_{L^2 }$. Now from $v''=f_{\lambda } -\lambda ^2 v$, one can easily derive by induction on $m$ that $\|v\|_{H^{m+2}}  \lesssim  |\lambda |^{m} \|f_{\lambda } \|_{H^{m}}$. This implies $\|u\|_{H^{m+2} (0,\Theta )} \lesssim  |\lambda |^{m} \|f\|_{H^{m}} +|\lambda |^{m+2} |u(0)|$, with $|u(0)|\lesssim  \|u\|_{H^{m+2} (-1,0)}$.\\
	So finally $\|u\|_{\mathcal{H}^{m+2}} \lesssim  |\Im(\lambda )|^{m+2} \|(f,g,a)\|_{\prodH{m}}$.
\end{enumerate}\fbi
\end{demo}

\begin{textesansboite}
\Cref{etape de base du DA des champs de coin} is the base step of the asymptotic expansion. It will be applied to $s$ the rest of the asymptotic expansion at a given order, and the function $s_{\mathrm{r}}$ below will be the rest at the next order. Iterating this process provides an asymptotic expansion of $S$ at any order. 
Since \cref{DA des champs de coin} implies infinite regularity on $S$, we work in the spaces $K^{\infty } _{\beta }$ and $\mathcal{H}^{m}$ for~any $m$.
\end{textesansboite}

\begin{proposition}{}
\label{etape de base du DA des champs de coin}%
Let $d\in \mathbb{R}$, $m\in \mathbb{N} \setminus \{0,1\}$, $s\in  \bigcap\limits _{\beta >d} K^{\infty } _{\beta } (\mathbb{R},\mathcal{H}^{m} )$, $f\in  \!\bigcap\limits _{\beta >d-1} \! K^{\infty } _{\beta } (\mathbb{R},H^{m-2} (0,\Theta ))$, $g\in  \!\bigcap\limits _{\beta >d-1} \! K^{\infty } _{\beta } (\mathbb{R},H^{m-2} (-1,0))$, $h\in  \!\bigcap\limits _{\beta >d-1} \! K^{\infty } _{\beta } (\sbd \Gamma )$ and $d_{\mathrm{\min }} := \min (\TZ*\cap (d-1,d])$. We assume:\vspace{-2mm}
\[\left\{\begin{array}{rl}
	\mu _{0} \Delta  s = f  &\text{ in }\sbd \Omega \\[0.5mm]
	\mu _{1} \,(e^{-2t} (\partial _{t} ^2 -\partial _{t} )+\partial _{\theta } ^2 )s = g  &\text{ in }\sbd\couche \\[0.5mm]
	s_{|\theta =0^{+}} - s_{|\theta =0^{-}} =0  &\text{ on }\sbd \Gamma  \\[0.5mm]
	e^{-t} \mu _{0} \partial _{\theta } s_{|\theta =0^{+}} - \mu _{1} \partial _{\theta } s_{|\theta =0^{-}} = h &\text{ on }\sbd \Gamma  \\[0.5mm]
	s =0 &\text{ on }\sbd\bordgauche \cup \sbd\borddroit
\end{array}\right.\]
Then there exist $s_{\mathrm{r}} \in  \!\bigcap\limits _{d-1<\beta <d_{\mathrm{\min }}}  \! K^{\infty } _{\beta } (\mathbb{R},\mathcal{H}^{m} )$ and $(c_{q} )\in \mathbb{C}^{\TZ*\cap  (d-1,d]}$ s.t.\ $\displaystyle s = \sum _{q\in \TZ*\cap  (d-1,d]} c_{q} \,\sbd{\phi_{q} } + s_{\mathrm{r}}$.\fbeq[0.5]
\end{proposition}

\newcommand{\poles}{\mathfrak{P}}
\begin{demo}
For~any $\beta \in \mathbb{R}$, let $\mathbb{C}_{\Re>\beta } :=\{\lambda \in \mathbb{C} \mid \Re(\lambda )>\beta \}$. The Laplace transform of $s$ is well-defined and holomorphic from $\mathbb{C}_{\Re>d}$ to $\mathcal{H}^{m}$. We will show that $\hat{s}$ has a meromorphic extension on $\mathbb{C}_{\Re>d-1}$, with poles belonging in $\poles:=\TZ*\cap (d-1,d]$, then we will apply the residue theorem on a rectangle surrounding these poles. The functions $\phi_{q}$ will appear in the residues.\medskip

\underline{Step 1 :} $\hat{s}$ satisfies in $\mathbb{C}_{\Re>d}$:
\[\left\{\begin{array}{r@{}ll}
	\mu _{0} (\partial _{\theta } ^2 +\lambda ^2 ) \hat{s}(\lambda ,\theta ) &\;=\; \hat{f}(\lambda ,\theta ) & \text{if}  \ \theta \in (0,\Theta )\\[0.5mm]
	\mu _{1} \partial _{\theta } ^2  \hat{s}(\lambda ,\theta ) &\;=\; \hat{g}(\lambda ,\theta ) -\mu _{1} \,\big[(\lambda +2)^2 -(\lambda +2)\big] \, \hat{s}(\lambda +2,\theta ) := \tilde{g}(\lambda ,\theta ) & \text{if}  \ \theta \in (-1,0)\\[0.5mm]
	\hat{s}(\lambda ,0^{+} ) - \hat{s}(\lambda ,0^{-}) &\;=\; 0& \\[0.5mm]
	\mu _{1} \partial _{\theta } \hat{s}(\lambda ,0^{-}) &\;=\; \hat{h}(\lambda ) +\mu _{0} \partial _{\theta } \hat{s}(\lambda +1,0^{+} ) := \tilde{h}(\lambda ) & \\[0.5mm]
	\hat{s}(\lambda ,\Theta ) = \hat{s}(\lambda ,-1) &\;=\; 0 &
\end{array}\right.\] For~any $\lambda \in  \mathbb{C}_{\Re>d-1} \setminus \TZ*$, let $\tilde{s} (\lambda ) := A(\lambda )^{-1} \big( \frac{1}{\mu _{0}} \hat{f}(\lambda ), \frac{1}{\mu _{1}} \tilde{g}(\lambda ), \frac{1}{\mu _{1}} \tilde{h}(\lambda )\big)$. Then $\tilde{s}$ is also solution of the above system, so $\tilde{s}$ and $\hat{s}$ coincide on $\mathbb{C}_{\Re>d} \setminus \TZ*$ by \cref{resolution de l'EDO en theta}. Hence $\tilde{s}$ is an extension of $\hat{s}$ on $\mathbb{C}_{\Re>d-1} \setminus \poles$, that we will still denote $\hat{s}$. \Cref{resolution de l'EDO en theta} implies that this extension is meromorphic with simple poles and: $\forall q\in \poles,\; \Res_{q} \hat{s} =(\Res_{\lambda =q} A(\lambda )^{-1} ) \big( \frac{1}{\mu _{0}} \hat{f}(q), \frac{1}{\mu _{1}} \tilde{g}(q), \frac{1}{\mu _{1}} \tilde{h}(q)\big)$.\medskip\smallskip

\underline{Step 2 :} Let $t\in \mathbb{R}$, $b\in (d-1,d_{\mathrm{\min }} )$ and $k\in \mathbb{R}_+^*$. The residue theorem applied to $\lambda  \mapsto  e^{\lambda t} \hat{s}(\lambda )$ on the rectangle $[b,d+1]\times [-k,k] \subset \mathbb{C}$ gives the following equality in $\mathcal{H}^{m}$. \begin{eqnarray*}
	2\mathrm{i}\pi  \sum _{q\in \poles} \Res_{q} (\lambda \mapsto e^{\lambda t} \hat{s}(\lambda )) &=& \int _{-k} ^{k} e^{(d+1+\mathrm{i}\gamma )t} \hat{s}(d+1+\mathrm{i}\gamma )\,\mathrm{d}\gamma  - \int _{-k} ^{k} e^{(b+\mathrm{i}\gamma )t} \hat{s}(b+\mathrm{i}\gamma )\,\mathrm{d}\gamma  \\[-1mm]
	&& + \int _{b} ^{d+1} e^{(\beta -\mathrm{i}k)t} \hat{s}(\beta -\mathrm{i}k)\,\mathrm{d}\beta  - \int _{b} ^{d+1} e^{(\beta +\mathrm{i}k)t} \hat{s}(\beta +\mathrm{i}k)\,\mathrm{d}\beta 
\end{eqnarray*}
Note that when $k\rightarrow \infty $ the first integral tends to the inverse Laplace transform of $s$ (up to a constant). We will show that the last two integrals tend to 0, by proving the following: 
\begin{equation}{\label{eq: etape de base du DA des champs de coin preuve}}
\forall n\in \mathbb{N}, \qquad  \sup _{\beta \in [b,d+1]} \int _{\mathbb{R} \setminus [-1,1]} \gamma ^{n} \| \hat{s}(\beta +\mathrm{i}\gamma )\|_{\mathcal{H}^{m}} \,\mathrm{d}\gamma  < \infty  \vspace{-2mm}\end{equation}
\begin{sousdemo}
	Let $n\in \mathbb{N}$ and $\beta \in [b,d+1]$. \cref{eq: normes des operateurs resolvants R1 R2 et R3} and the expressions of $\tilde{g}$ and $\tilde{h}$ yield
	\begin{align*}
		\int_{\mathbb{R} \setminus [-1,1]} \gamma ^{n} \,\| \hat{s}(\beta +\mathrm{i}\gamma )\|_{\mathcal{H}^{m}} \,\mathrm{d}\gamma 
		&\leqslant  \int_{\mathbb{R} \setminus [-1,1]} \gamma ^{n} \,\|A(\lambda )^{-1} \|_{\mathcal{L}(\prodH{m-2},\mathcal{H}^{m} )} \\
		& \qquad  \qquad  \quad  {}\times \Big\|\big({\textstyle  \frac{1}{\mu _{0}} \hat{f}(\beta +\mathrm{i}\gamma ), \frac{1}{\mu _{1}} \tilde{g}(\beta +\mathrm{i}\gamma ), \frac{1}{\mu _{1}} \tilde{h}(\beta +\mathrm{i}\gamma )}\big)\Big\|_{\prodH{m-2}} \,\mathrm{d}\gamma  \hspace{-20mm} \\
		&\lesssim  \int_{\mathbb{R} \setminus [-1,1]} \! \gamma ^{n+m+2} \,\Big( \| \hat{f} (\beta +\mathrm{i}\gamma )\|_{H^{m-2}} +\| \hat{g} (\beta +\mathrm{i}\gamma )\|_{H^{m-2}} +| \hat{h}(\beta +\mathrm{i}\gamma )| \\[-2mm]
		& \qquad  \qquad  \quad  \ {} + \| \hat{s}(\beta +\mathrm{i}\gamma +1)\|_{\mathcal{H}^{m}} + \gamma ^2 \| \hat{s}(\beta +\mathrm{i}\gamma +2)\|_{\mathcal{H}^{m}} \Big) \,\mathrm{d}\gamma  
	\end{align*}
	Let us treat the term with $f$, given that the others are similar.
	\begin{eqnarray*}
		\int _{\mathbb{R} \setminus [-1,1]} \gamma ^{n+m+2} \| \hat{f} (\beta +\mathrm{i}\gamma )\|_{H^{m-2}} \,\mathrm{d}\gamma  
		&\underset{\mathcal{\text{C.S.}}}{\leqslant }& \|\gamma ^{-1} \|_{L^2 (\mathbb{R} \setminus [-1,1])} \cdot  \Big\| \gamma ^{n+m+3} \| \hat{f}(\beta +\mathrm{i}\gamma )\|_{H^{m-2}} \Big\|_{L^2 (\mathbb{R} \setminus [-1,1])} \\[-0.5mm]
		&\lesssim & \|e^{-\beta t} \partial _{t} ^{n+m+3} f\|_{L^2 (\mathbb{R},H^{m-2} (0,\Theta ))} \\
		&\lesssim & \max _{\beta '\in \{b,d+1\}} \|e^{-\beta 't} \partial _{t} ^{n+m+3} f\|_{L^2 (\mathbb{R},H^{m-2} (0,\Theta ))} \\
		&=& \max _{\beta '\in \{b,d+1\}} \|f\|_{K^{n+m+3} _{\beta '} (\mathbb{R},H^{m-2} (0,\Theta ))}
	\end{eqnarray*}
	thanks to an interpolation between $b$ and $d+1$. We have majorized by a finite constant independent of $\beta $, so the supremum on $\beta $ is finite. Thus \cref{eq: etape de base du DA des champs de coin preuve} is proven.
\end{sousdemo}
Therefore: $\displaystyle \int _{1} ^{\infty } \bigg\| \int _{b} ^{d+1} e^{(\beta \pm \mathrm{i}k)t} \hat{s}(\beta \pm \mathrm{i}k) \,\mathrm{d}\beta \,\bigg\|_{\mathcal{H}^{m}} \mathrm{d}k \lesssim  \int _{1} ^{\infty } \! \int _{b} ^{d+1} \| \hat{s}(\beta \pm \mathrm{i}k)\|_{\mathcal{H}^{m}} \,\mathrm{d}\beta \,\mathrm{d}k < +\infty $.\\
So by~taking arbitrarily large $k$, we can make the terms $\int _{b} ^{d+1} e^{(\beta \pm \mathrm{i}\gamma )t} \hat{s}(\beta \pm \mathrm{i}k) \,\mathrm{d}\beta $ tend to 0 in $\mathcal{H}^{m}$. Hence:
\[\sum _{q\in \poles} \Res_{q} (\lambda \mapsto e^{\lambda t} \hat{s}(\lambda )) = \frac{1}{2\mathrm{i}\pi } \int _{\mathbb{R}} e^{(d+1+\mathrm{i}\gamma )t} \hat{s}(d+1+\mathrm{i}\gamma ) \,\mathrm{d}\gamma   - \frac{1}{2\mathrm{i}\pi } \int _{\mathbb{R}} e^{(b+\mathrm{i}\gamma )t} \hat{s}(b+\mathrm{i}\gamma ) \,\mathrm{d}\gamma .\]
Now, since $\hat{s}$ has only simple poles, the terms of the sum are equal to 
\[e^{qt} \Res_{q} \hat{s} = e^{qt} (\Res_{\lambda =q} A(\lambda )^{-1} )\big({\textstyle \frac{1}{\mu _{0}} \hat{f}(q), \frac{1}{\mu _{1}} \tilde{g}(q), \frac{1}{\mu _{1}} \tilde{h}(q)}\big) = e^{qt} c_{q} \mathbf{1}_{[0,\Theta ]} (\theta ) \sin (q\theta ) = c_{q} \,\sbd{\phi_{q} }(t,\theta )\]
for some constant $c_{q} \in \mathbb{C}$ by point 2 of \cref{resolution de l'EDO en theta}. The first integral is equal to $s(t)$ by inverse Laplace transform (because $s\in K_{d+1} ^{\infty } (\mathbb{R},\mathcal{H}^{m} )$). We define $s_{\mathrm{r}} (t)$ to be equal to the last integral. Thus we get the desired formula, and \cref{eq: etape de base du DA des champs de coin preuve} shows that $s_{\mathrm{r}} \in K^{\infty } _{b} (\mathbb{R},\mathcal{H}^{m} )$.
\end{demo}

\begin{textesansboite}
Let $S$ be the function set in \cref{DA des champs de coin} and $(\sigma _{d} (S))_{d\in \TN}$ the coefficients of \cref{def: sigma dans V + khi A}, which vanish for $d$ large enough. With the $\sigma _{d} (S)$, the non-variational part of $S$ (denoted $S^{\mathcal{A}}$ below) can be explicitly computed, so it remains to obtain an asymptotic expansion of its variational part (denoted $S^{\Vcoin}$ below). In order to apply \cref{etape de base du DA des champs de coin}, we need that $\sbd{\khii S^{\Vcoin} }$ belongs to some space $K^{\infty } _{\beta }$ (since $\Omega _{1} \setminus B(0,\Rcoin)=\Pi  \setminus B(0,\Rcoin)$, $\khii S^{\Vcoin}$ can be seen as a function defined on $\Pi $, which allows us to consider $\sbd{\khii S^{\Vcoin} }$). We will use again the notations $T_{\geqslant d}$ and $\langle .\rangle $ of \cref{def: serie geometrique d'op 1D,def: troncature de SF}.
\end{textesansboite}

\begin{proposition}{}
\label{regularite de S variationnel}%
Let $d\in \mathbb{R}_-$,\vspace{-1mm}
\[S^{\mathcal{A}} := T_{\geqslant d} \Bigg[ \left\langle  -\RcouchePi\circ \partial _{X|\couche} ^2  ,\; \frac{\mu _{0} }{\mu _{1}} \RGammaPi\circ \partial _{Y|\Gamma ,Y=0^{+}} \right\rangle  \bigg( \frac{1}{\mu _{0}} \ROmegaPi (F^{\infty } _{\Omega } )+ \frac{1}{\mu _{1}} \RcouchePi (F^{\infty } _{\couche} ) + \fsum_{d\in \TN} \sigma _{d} (S) \,\phi_{d} \bigg) \Bigg] \vspace{-2.5mm}\]
in $\mathcal{A}(\Pi )$ and $S^{\Vcoin} := S-\khii S^{\mathcal{A}}$. There exists a value of $d\in \mathbb{R}_-$ s.t.\ we have 
\[S^{\Vcoin} \in  \Vcoin \quad \  \text{ and } \quad \  \forall m\in \mathbb{N},\ \ \sbd{\khii S^{\Vcoin} } \in K^{\infty } _{1/2} (\mathbb{R},\mathcal{H}^{m} ).\]\fbeq[1.2]
\end{proposition}

\begin{demo}
We consider any $d\in \mathbb{R}_-$ and we will fix it later. The proof has three steps: showing that $S^{\Vcoin} \in \Vcoin$, showing that $S^{\Vcoin}$ is regular w.r.t $t$, and deducing that it~is regular w.r.t.\ $\theta $. For~any $r_{1} ,r_{2} \in \mathbb{R}$ s.t.\ $\Rcoin<r_{1} <r_{2}$, we will denote $\chi \in \mathcal{C}^{\infty } (\mathbb{R}^2 )$ a radial function equal to 0 on $B(0,r_{1} )$ and 1 on $\mathbb{R}^2  \setminus B(0,r_{2} )$ ($r_{1} ,r_{2} $ are implicit in this~notation).\medskip\smallskip

\underline{Step 1:} Let $\varphi :\Pi \rightarrow \mathbb{C}$ be equal to $T_{\leqslant d} (F^{\infty } _{D} )$ on $D$ for~any $D\in \{\Omega ,\couche\}$, and $F_{\mathrm{v}} = F-\khii \varphi $. If $d<-2$ then $(1+r)F_{\mathrm{v}} \in L^2 (\Omega _{1} )$. So the same construction than in \cref{cadre fonctionnel pour les champs de coin} shows that for $d$ small enough, the problem 
\begin{equation}{\label{eq: pb verifie par S V}}
\left\{\begin{array}{r@{\;=\;}ll}
\operatorname{div}(\mu  \nabla  S^{\Vcoin} ) & f_{1} := F_{\mathrm{v}} + \khii \varphi  - \operatorname{div}\!\big(\mu  \nabla  (\khii \cdot S^{\mathcal{A}} )\big)&\text{ in }\Omega _{1} \setminus (\Gamma \cap \{X>\Rcoin\})  \\[0.5mm]
[\mu  \partial _{Y} S^{\Vcoin} ]_{\Gamma } & g_{1} := -[\mu \partial _{Y} (\khii S^{\mathcal{A}} )]_{\Gamma } &\text{ on }\Gamma \cap \{X>\Rcoin\} \\[0.5mm]
S^{\Vcoin} &0 &\text{ on }\partial \Omega _{1}
\end{array}\right.\end{equation}
has a unique solution in $V$ and that $S = S^{V} + \khii S^{\mathcal{A}}$. In addition, it also gives that $(\mu  \Delta  S^{\mathcal{A}} - \varphi )_{|D} \in  \sum _{d'<d} \mathcal{A}_{d'} (D)$ for~any $D\in \{\Omega ,\couche\}$, and $[\mu \partial _{Y} S^{\mathcal{A}} ]_{\Gamma } \in  \sum _{d'<d} \mathcal{A}_{d'} (\Gamma )$. And we have $F_{\mathrm{v}} =o_{\partial } (\bfr^{d} )$ by hypothesis, so for $d$ small enough and up to increasing $\Rcoin$ (without loss of generality):
\begin{equation}{\label{eq: regularite de S V: snd membres et Kondrat'ev}}
	\forall m\in \mathbb{N}, \quad  \left\{\begin{array}{@{\;\,}r@{}l@{}}
		(\sbd{\chi  f_{1} })_{|\sbd \Omega } &{}\in  K^{\infty } _{-2} (\mathbb{R},H^{m} (0,\Theta )) \\[2pt]
		(\sbd{\chi  f_{1} })_{|\sbd \couche} &{}\in  K^{\infty } _{-3/2} (\mathbb{R},H^{m} (-1,0)) 
	\end{array}\right.
	\qquad  \text{ and } \qquad  \sbd{\chi  g_{1} } \in  K^{\infty } _{-1/2} (\sbd \Gamma )
\end{equation}\medskip

\underline{Step 2:} For~any $\beta \in \mathbb{R}$, we denote:
\begin{itemize}
	\item $K^{0} _{\beta } (\sbd \Omega ):=K^{0} _{\beta } (\mathbb{R},L^2 (0,\Theta )) = \{(t,\theta )\mapsto  e^{\beta t} u(t,\theta ) \mid u\in L^2 (\sbd \Omega )\}$,
	\item $K^{0} _{\beta } (\sbd\couche) :=K^{0} _{\beta } (\mathbb{R},L^2 (-1,0))$,
	\item and $\sbd \Vcoin := \{w\in H^{1} _{\mathrm{loc}} (\sbd \Pi ) \mid \nabla  w_{|\sbd \Omega } \in L^2 (\sbd \Omega ),\; \partial _{t} w_{|\sbd\couche } \in K^{0} _{1/2} (\sbd\couche ),\; \partial _{\theta } w_{|\sbd\couche } \in K^{0} _{-1/2} (\sbd\couche )\, \text{ and } \,w_{|\sbd\bordgauche \cup \sbd\borddroit }=0 \}$.
\end{itemize}
Note that if $\Omega _{1}$ were equal to $\Pi $, $\sbd \Vcoin$ would simply be $\{\sbd u \mid u\in \Vcoin \}$. So~it~is the natural variational space for Poisson's problem transferred into $\sbd \Pi $. Looking at $\chi  S^{\Vcoin}$ as a function defined on $\Pi $, we will show by induction on $n$ that:
\[\forall n\in \mathbb{N},\; \forall \Rcoin<r_{1} <r_{2} ,\qquad  \partial _{t} ^{n} \sbd{\chi  S^{\Vcoin} }\in \sbd \Vcoin.\]
Since $S^{\Vcoin} \in \Vcoin$, the initial case is trivial, so only the inductive step remains to prove. Let us assume the property at rank $n$ and show it at rank $n+1$. We will use the method of finite differences. Letting $\mu :=\mu _{0}$ on $\Omega $ and $\mu :=\mu _{1}$ on $\couche$, \cref{eq: pb verifie par S V} implies:
\[\left\{\begin{array}{r@{\;=\;}ll}
	\mu  \Delta (\chi  S^{\Vcoin} ) & f_{2} := \chi  f_{1} - 2\mu  \nabla  \chi  \cdot \nabla  S^{\Vcoin} -\mu  \Delta  \chi  \cdot S^{\Vcoin} &\text{ in }D,\ \forall D\in \{\Omega ,\couche\} \\[0.5mm]
	[\mu \partial _{Y} (\chi  S^{\Vcoin} )]_{\Gamma } & g_{2} := \chi  g_{1} &\text{ on }\Gamma  \\[0.5mm]
	\chi  S^{\Vcoin} &0 &\text{ on }\bordgauche \cup \borddroit
\end{array}\right.\]
By~induction hypothesis applied to $(r_{1} ',r_{2} '):= (\frac{\Rcoin+r_{1} }{2} ,r_{1} )$, $\partial _{t} ^{k} \nabla  S^{\Vcoin}$ and $\partial _{t} ^{k} S^{\Vcoin}$ are $L^2 $ on $\{r_{1} <r<r_{2} \}$ for~any $k\in [\![0,n]\!]$. So \cref{eq: regularite de S V: snd membres et Kondrat'ev} implies $\sbd{f_{2} }_{|\sbd \Omega } \in K^{n} _{-2} (\mathbb{R},L^2 (0,\Theta ))$ and $\sbd{f_{2} }_{|\sbd \couche} \in K^{n} _{-3/2} (\mathbb{R},L^2 (-1,0))$. Similarly $\sbd{g_{2} }\in K^{\infty } _{-1/2} (\sbd \Gamma )$. Next, we apply the change variables $(x,y) \rightsquigarrow (t,\theta )$ using \cref{eq: transfert du probleme de Poisson de Pi a sbdPi}, and then $\partial _{t} ^{n}$. We get that $s:=\partial _{t} ^{n} (\sbd{\chi  S^{\Vcoin} })$ satisfies: 
\begin{equation}{\label{eq: regularite de S: equations de s}}
\hspace{-5mm}\left\{\begin{array}{@{\;}r@{\;=\;}l@{\;\,}l@{}}
	-e^{-2t} \mu _{0} \Delta s & f_{\Omega } :=e^{-2t} \partial _{t} ^{n} (e^{2t} \sbd{f_{2} }_{|\Omega } ) &\text{ in }\sbd \Omega  \\[0.5mm]
	-\mu _{1} \,(e^{-2t} (\partial _{t} ^2 -\partial _{t} )+\partial _{\theta } ^2 )s & f_{\couche} :=\partial _{t} ^{n} \sbd{f_{2} }_{|\couche} + \sum\limits _{k=0} ^{n-1} (-2)^{n-k} e^{-2t} \mu _{1} (\partial _{t} ^2 -\partial _{t} ) \partial _{t} ^{k} (\sbd{\chi  S^{V} }) &\text{ in }\sbd\couche \\[0.5mm]
	s_{|\theta =0^{+}} - s_{|\theta =0^{-}} & 0 &\text{ on }\sbd \Gamma  \\[0.5mm]
	e^{-t} \mu _{0} \partial _{\theta } s_{|\theta =0^{+}} - \mu _{1} \partial _{\theta } s_{|\theta =0^{-}} & g :=\partial _{t} ^{n} \sbd{g_{2} } &\text{ on }\sbd \Gamma  \\[0.5mm]
	s & 0 &\text{ on }\sbd\bordgauche \cup \sbd\borddroit \hspace{-7mm}
\end{array}\right. \end{equation}
with $f_{\Omega } \in K^{0} _{-2} (\Omega )$, $f_{\couche} \in K^{0} _{-3/2} (\sbd\couche )$ and $g\in K^{\infty } _{-1/2} (\sbd \Gamma )$ by induction hypothesis. The variational formulation of \cref{eq: regularite de S: equations de s} is: $\forall \varphi \in \sbd \Vcoin ,$
\[ \int _{\sbd \Omega } \mu _{0} \nabla  s \cdot \nabla  \varphi +\int _{\sbd \couche} \mu _{1} \,(e^{-t} \,\partial _{t} s \,\partial _{t} \varphi + e^{t} \,\partial _{\theta } s \,\partial _{\theta } \varphi ) = \int _{\sbd \Omega } e^{2t} f_{\Omega } \varphi + \int _{\sbd\couche} e^{t} f_{\couche} \varphi + \int _{\sbd \Gamma } e^{t} g \varphi .\]
Let us denote $D_{\eta } \varphi (t,\theta ) := \frac{\varphi (t+\eta ,\theta )-\varphi (t,\theta )}{\eta }$ for~any $\eta \in \mathbb{R}^*$ and any function $\varphi $. Taking $\varphi  := D_{-\eta } D_{\eta } \bar{s}$ and discretely integrating by parts $D_{-\eta }$, we have:
\begin{multline*}
	\int _{\sbd \Omega } \mu _{0} |D_{\eta } \nabla  s|^2  +\int _{\sbd\couche} \mu _{1} \,(D_{\eta } (e^{-t} \partial _{t} s) \cdot D_{\eta } \partial _{t} \bar{s}+ e^{t} \,|D_{\eta } \partial _{\theta } s|^2 ) \\
	= \int _{\sbd \Omega } e^{2t} f_{\Omega } \cdot D_{-\eta } D_{\eta } \bar{s}
	+ \int _{\sbd\couche} e^{t} f_{\couche} \cdot D_{-\eta } D_{\eta } \bar{s}
	+ \int _{\sbd \Gamma } D_{-\eta } D_{\eta } (e^{t} g) \cdot  \bar{s}
\end{multline*}
But for~any functions $\varphi ,\psi $, we have $D_{\eta } (\varphi  \psi )=\varphi  \cdot  D_{\eta } \psi  + D_{\eta } \varphi  \cdot \psi $. Therefore:
\begin{itemize}
	\item $D_{\eta } (e^{-t} \partial _{t} s) \cdot D_{\eta } \partial _{t} \bar{s}=e^{-t} |D_{\eta } \partial _{t} s|^2  + \frac{e^{-\eta } -1}{\eta } e^{-t} \partial _{t} s\cdot D_{\eta } \partial _{t} \bar{s}$
	\item $D_{-\eta } D_{\eta } (e^{t} g) =e^{t} D_{-\eta } D_{\eta } g + 2 \frac{e^{\eta } -1}{\eta } e^{t} D_{\eta } g + \big(\frac{e^{\eta } -1}{\eta } \big)^2  e^{t} g$
\end{itemize}
Let us assume that $\eta $ is small enough so~that $\frac{|e^{\eta } -1|}{\eta } <2$ and $\frac{|e^{-\eta } -1|}{\eta } <2$. Then:\vspace{-2mm}
\begin{multline*}
	\|D_{\eta } \nabla  s\|_{L^2 (\sbd \Omega )}^2  + \|D_{\eta } \partial _{t} s\|_{K^{0} _{1/2} (\sbd\couche )}^2  +\|D_{\eta } \partial _{\theta } s\|_{K^{0} _{-1/2} (\sbd\couche )}^2  + \int _{\sbd\couche} \mu _{1} {\textstyle \frac{e^{-\eta } -1}{\eta } }e^{-t} \partial _{t} s\cdot D_{\eta } \partial _{t} \bar{s}\\
	\lesssim  \|f_{\Omega } \|_{K^{0} _{-2} (\sbd \Omega )} \|D_{-\eta } D_{\eta } s\|_{L^2 (\sbd \Omega )} 
	+ \|f_{\couche} \|_{K^{0} _{-3/2} (\sbd\couche )} \|D_{-\eta } D_{\eta } s\|_{K^{0} _{1/2} (\sbd\couche )}
	+ \|g\|_{K^{2} _{-1/2} (\sbd \Gamma )} \|s\|_{K^{0} _{-1/2} (\sbd \Gamma )}
\end{multline*}
Then, moving the integral to the right-hand side and majorizing some $\|D_{\eta } \cdot \|$ by $\|\partial _{t} \cdot \|$ or $\|\nabla  \cdot \|$, we get:
\begin{multline*}
	\|D_{\eta } \nabla  s\|_{L^2 (\sbd \Omega )}^2  + \|D_{\eta } \partial _{t} s\|_{K^{0} _{1/2} (\sbd\couche )}^2  +\|D_{\eta } \partial _{\theta } s\|_{K^{0} _{-1/2} (\sbd\couche )}^2  
	\lesssim  \|f_{\Omega } \|_{K^{0} _{-2} (\sbd \Omega )} \|D_{\eta } \nabla s\|_{L^2 (\sbd \Omega )} \\
	+ \|f_{\couche} \|_{K^{0} _{-3/2} (\sbd\couche )} \|D_{\eta } \partial _{t} s\|_{K^{0} _{1/2} (\sbd\couche )} 
	+ \|g\|_{K^{2} _{-1/2} (\sbd \Gamma )} \|s\|_{K^{0} _{-1/2} (\sbd \Gamma )}
	+ \|\partial _{t} s\|_{K^{0} _{1/2} (\sbd\couche )} \|D_{\eta } \partial _{t} s\|_{K^{0} _{1/2} (\sbd\couche )}
\end{multline*}
But, since $s_{|\sbd\borddroit} =0$, a Poincaré-type inequality gives $\|s\|_{K^{0} _{-1/2} (\sbd \Gamma )}^2  \lesssim  \|\partial _{\theta } s\|_{K^{0} _{-1/2} (\sbd\couche )} ^2 $. Finally, using Young's inequality $ab\lesssim  \frac{1}{\delta } a +\delta  b$ on norm products with $\delta $ small enough and moving the $\|D_{\eta } \partial _{t} s\|_{K^{0} _{1/2} (\sbd\couche )}^2 $ and $\|D_{\eta } \partial _{\theta } s\|_{K^{0} _{-1/2} (\sbd\couche )}^2 $ from the right-hand side to the left-hand one, we get:
\begin{align*}
	\|\nabla  \partial _{t} s\|_{L^2 (\sbd \Omega )}^2  &+ \|\partial _{t} ^2  s\|_{K^{0} _{1/2} (\sbd\couche )}^2  +\|\partial _{\theta } \partial _{t} s\|_{K^{0} _{-1/2} (\sbd\couche )}^2  \\
	&\leqslant  \limsup _{\eta \rightarrow 0} \|D_{\eta } \nabla  s\|_{L^2 (\sbd \Omega )}^2  + \|D_{\eta } \partial _{t} s\|_{K^{0} _{1/2} (\sbd\couche )}^2  +\|D_{\eta } \partial _{\theta } s\|_{K^{0} _{-1/2} (\sbd\couche )}^2  \\
	&\lesssim  \|f_{\Omega } \|_{K^{0} _{-2} (\sbd \Omega )}^2  + \|f_{\couche} \|_{K^{0} _{-3/2} (\sbd\couche )}^2  + \|g\|_{K^{1} _{-1/2} (\sbd \Gamma )}^2  + \|\partial _{\theta } s\|_{K^{0} _{-1/2} (\sbd\couche )}^2  + \|\partial _{t} s\|_{K^{0} _{1/2} (\sbd\couche )}^2  \\
	&<\infty 
\end{align*}
By~definition of $\sbd \Vcoin$, it implies that $\partial _{t} s\in \sbd \Vcoin$ and completes the induction.\medskip\smallskip

\underline{Step 3:} Let $K^{0} _{1/2} (\sbd \Pi ) := K^{0} _{1/2} (\mathbb{R},L^2 (-1,\Theta ))$. For now we have proven that, for~any $(n,i)\in (\mathbb{N}\times \{0,1\}) \setminus (\{(0,0)\})$ and $\Rcoin<r_{1} <r_{2}$, $\partial _{t} ^{n} \partial _{\theta } ^{i} \sbd{\chi  S^{\Vcoin} } \in  K^{0} _{1/2} (\sbd \Pi )$. It generalizes to the case $(n,i)=(0,0)$ thanks to a Poincaré inequality:
\[\|\sbd{\chi  S^{\Vcoin} }\|_{K^{0} _{1/2} (\sbd \Pi )} \lesssim  \|\partial _{\theta } \sbd{\chi  S^{\Vcoin} }\|_{K^{0} _{1/2} (\sbd \Pi )} \lesssim  \|\partial _{\theta } \sbd{\chi  S^{\Vcoin} }\|_{L^2 (\sbd \Omega )} + \|\partial _{\theta } \sbd{\chi  S^{\Vcoin} }\|_{K^{0} _{-1/2} (\sbd\couche)} < \infty .\]
To treat higher-order $\theta $-derivatives, we start from the equality $\mu  \Delta (\chi  S^{\Vcoin} ) = \chi  f_{1} - 2\mu  \nabla  \chi  \cdot \nabla  S^{\Vcoin} -\mu  \Delta  \chi  \cdot S^{\Vcoin}$ in $\Omega \cup \couche$ proven in step 2. Applying the change of variables $(x,y) \rightsquigarrow (t,\theta )$ gives by \cref{eq: transfert du probleme de Poisson de Pi a sbdPi}: 
\begin{equation}{\label{eq: regularite de S etape 3}}
\left\{\begin{array}{r@{\;=\;}ll}
	-e^{-2t} \mu _{0} \Delta \sbd{\chi  S^{\Vcoin} } & 
	\sbd{\chi  f_{1} } - 2\mu _{0} \sbd{\nabla  \chi  \cdot \nabla  S^{\Vcoin} } -\mu _{0} \sbd{\Delta  \chi  \cdot S^{\Vcoin} } &\text{ in }\sbd \Omega  \\[0.5mm]
	-\mu _{1} \,(e^{-2t} (\partial _{t} ^2 -\partial _{t} )+\partial _{\theta } ^2 )(\sbd{\chi  S^{\Vcoin} }) & 
	\sbd{\chi  f_{1} } - 2\mu _{1} \sbd{\nabla  \chi  \cdot \nabla  S^{\Vcoin} } -\mu _{1} \sbd{\Delta  \chi  \cdot S^{\Vcoin} } &\text{ in }\sbd\couche 
\end{array}\right.\end{equation}
Moreover, \cref{eq: regularite de S V: snd membres et Kondrat'ev} implies that, for~any $m\in \mathbb{N}$, $\sbd{\chi  f_{1} }_{|\sbd \Omega } \in K^{\infty } _{-3/2} (\mathbb{R},H^{m} (0,\Theta ))$ and $\sbd{\chi  f_{1} }_{|\sbd \couche} \in K^{\infty } _{1/2} (\mathbb{R},H^{m} (-1,0))$. So~deriving \cref{eq: regularite de S etape 3} w.r.t.\ $t$ and $\theta $ enough times gives by induction: $\forall (n,i)\in \mathbb{N}^2 , \forall \Rcoin<r_{1} <r_{2} ,\forall D\in \{\sbd \Omega ,\sbd\couche\},\; \partial _{t} ^{n} \partial _{\theta } ^{i} \sbd{\chi  S^{\Vcoin} }_{|D} \in  K^{0} _{1/2} (D)$. Finally, applying this to $\chi :=\khii$ completes the proof.
\end{demo}

\newcommand{\dTronc}[1]{d_{#1}}
\newcommand{\Sreste}[1]{S_{\mathrm{r},#1}}
\newcommand{\SAbarre}[1]{S^\infty_{#1}}
\newcommand{\SresteTilde}{\tilde{S}_{\mathrm{r},n+1}}
\newcommand{\nvCoefs}{\mathfrak{D}}
\phantomsection\label{preuve: DA des champs de coin 2}
\begin{demo}[Proof of \cref{DA des champs de coin}]
For~any $q\in \TN$, let $\sigma _{q} (\SDA):=\sigma _{q} (S)$. We denote for~any $d\in \mathbb{R}$
\[\SAbarre{d} := \left\langle  -\RcouchePi\circ \partial _{x|\couche} ^2 ,\; \frac{\mu _{0} }{\mu _{1}} \RGammaPi\circ \partial _{y|\Gamma ,y=0^{+}} \right\rangle  \bigg( \frac{1}{\mu _{0}} \ROmegaPi (F^{\infty } _{\Omega } )+ \frac{1}{\mu _{1}} \RcouchePi (F^{\infty } _{\couche} ) +\! \fsum_{q\in \TZ*\cap (d,\infty )} \! \sigma _{q} (\SDA) \,\phi_{q} \bigg) \]
and, for~any $n\in \mathbb{N}$, $\dTronc{n}:= \frac{1}{2} -n$ and $\Sreste{n} := \khii S- \khii T_{>\dTronc{n}} (\SAbarre{\dTronc{n}})$. Let $m\in \mathbb{N}$. We will show by induction on $n$ that there are coefficients $(\sigma _{q} (\SDA))_{q\in -\TN}$ s.t.: \[\forall n\in \mathbb{N},\; \forall \beta >\dTronc{n},\qquad  \sbd{\Sreste{n} }\in  K^{\infty } _{\beta } (\mathbb{R},\mathcal{H}^{m} ).\]
For now we only know $\sigma _{q} (\SDA)$ for $q\geqslant  \frac{\pi }{\Theta } >d_{0}$, and $\Sreste{n}$ involves $\sigma _{q} (\SDA)$ only when $q>\dTronc{n}$, so the $n$-th inductive step involves fixing $\sigma _{q} (\SDA)$ for all $q\in \TZ*\cap (\dTronc{n},\dTronc{n-1}]$.

\begin{description}
\item[Initial case:] By~\cref{regularite de S variationnel}, there exist $d\in \mathbb{R}_-$ and $S^{\Vcoin}$ s.t.\ $S=S^{\Vcoin} +\khii T_{\geqslant d} (\SAbarre{\dTronc{0}})$ and $\sbd{\khii S^{\Vcoin} } \in K^{\infty } _{1/2} (\mathbb{R},\mathcal{H}^{m} )$. Hence
\begin{align*}
	\Sreste{0} &= \khii S- \khii T_{>\dTronc{0}} (\SAbarre{\dTronc{0}})\\
	&= \khii \big(S^{\Vcoin} + \khii T_{\geqslant d} (\SAbarre{\dTronc{0}}) \big) - \khii T_{>\dTronc{0}} (\SAbarre{\dTronc{0}})\\
	&= \khii S^{\Vcoin} + \khii \cdot (\khii -1)\cdot T_{\geqslant d} (\SAbarre{\dTronc{0}}) +\,\khii \cdot \big(T_{\geqslant d} (\SAbarre{\dTronc{0}}) -T_{>0} (\SAbarre{\dTronc{0}}) \big)
\end{align*}
So $\sbd{\Sreste{0}}\in K^{\infty } _{1/2} (\mathbb{R},\mathcal{H}^{m} )$. But $\sbd{\Sreste{0}}$ is null in a vicinity of $-\infty $, so: $\forall \beta >d_{0} = \frac{1}{2} ,\; \sbd{\Sreste{0}}\in K^{\infty } _{\beta } (\mathbb{R},\mathcal{H}^{m} )$.

\item[Inductive step:] We assume the property at rank $n$ and will show it at rank $n+1$. Let $\SresteTilde:= \khii S- \khii T_{>\dTronc{n+1}} (\SAbarre{\dTronc{n}})$ and $\nvCoefs:=\TZ*\cap (\dTronc{n+1},\dTronc{n}]$. $\SresteTilde$ is a variant of $\Sreste{n+1}$ that does~not involve $\sigma _{q} (\SDA)$ for $q\in \nvCoefs$. We will apply \cref{etape de base du DA des champs de coin} to $d:=\dTronc{n}$ and $s:=\sbd{\SresteTilde}$. To do so, we must check that: \[\left\{\begin{array}{r@{\;\,\in\;}l}
\mu _{0} \Delta  \sbd{\SresteTilde}_{|\sbd \Omega } & \bigcap\limits _{\beta >\dTronc{n+1}} \! K^{\infty } _{\beta } (\mathbb{R},H^{m-2} (0,\Theta )) \\
\mu _{1} \,(e^{-2t} (\partial _{t} ^2 -\partial _{t} )+\partial _{\theta } ^2 )\sbd{\SresteTilde}_{|\sbd\couche} & \bigcap\limits _{\beta >\dTronc{n+1}} \! K^{\infty } _{\beta } (\mathbb{R},H^{m-2} (-1,0)) \\
(e^{-t} \mu _{0} \partial _{\theta |\theta =0^{+}} - \mu _{1} \partial _{\theta |\theta =0^{-}} ) \sbd{\SresteTilde} & \bigcap\limits _{\beta >\dTronc{n+1}} \! K^{\infty } _{\beta } (\sbd \Gamma )
\end{array}\right.\]
Let us show only the first line, the others being similar.
\begin{sousdemo}
	Since $\sbd{\Delta \SresteTilde} = e^{-2t} \Delta \sbd{\SresteTilde}$, it suffices to have: $\forall \beta >\dTronc{n+1}-2,\;\sbd{\Delta \SresteTilde} \in K^{\infty } _{\beta } (\mathbb{R},H^{m-2} (0,\Theta ))$. And this is true because we have in $\Omega $ in a vicinity of $r\rightarrow \infty $:\vspace{-1mm} \begin{align*}
		\mu _{0} \Delta \SresteTilde &= F - \mu _{0} \Delta  \big[T_{>\dTronc{n+1}} (\SAbarre{\dTronc{n}})\big] \tag*{by~definition of $\SresteTilde$, and $\mu _{0} \Delta S=F\text{ in }\Omega $}\\
		&= T_{\geqslant \dTronc{n+1}-2} (F^{\infty } _{\Omega } )+o_{\partial } (r^{\dTronc{n+1}-2} ) \tag*{by hypothesis on $F$}\\
		&\ \ \ -\mu _{0} T_{>\dTronc{n+1}-2} (\Delta \SAbarre{\dTronc{n}}) \tag*{by $\deg \Delta =-2$ (\cref{operateurs de derivation dans A})}\\
		&= (T_{\leqslant \dTronc{n+1}-2} -T_{<\dTronc{n+1}-2} )(F^{\infty } _{\Omega } ) + o_{\partial } (r^{\dTronc{n+1}-2} ) \tag*{by $\mu _{0} \Delta \SAbarre{\dTronc{n}}=F^{\infty } _{\Omega }$ (\cref{pb de Poisson dans A Pi})}
	\end{align*}\fbeq
\end{sousdemo}
Therefore, \cref{etape de base du DA des champs de coin} states that there are coefficients $(c_{q} )_{q\in \nvCoefs}$, $s_{\mathrm{r}}$ and $d_{\mathrm{\min }} :=\min (\nvCoefs)$ s.t.: $\forall \beta \in (\dTronc{n+1},d_{\mathrm{\min }} ), \; s_{\mathrm{r}} \in K^{\infty } _{\beta } (\mathbb{R},\mathcal{H}^{m} )$ and $\sbd{\SresteTilde} = \sum _{q\in \nvCoefs} c_{q} \,\sbd{\phi_{q} } + s_{\mathrm{r}}$. Let $\sigma _{q} (\SDA) :=-c_{q}$ for~any $q\in \TZ*\cap (\dTronc{n+1},\dTronc{n}]$. Then:
\[\sbd{\Sreste{n+1}} 
\;=\; \sbd{\SresteTilde} + \sum _{q\in \nvCoefs} \sigma _{q} (\SDA) \,\sbd{\khii \phi_{q} } 
\;=\; s_{\mathrm{r}} + \sum _{q\in \nvCoefs} c_{q} \,\sbd{(1-\khii )\phi_{q} } 
\;\in  \! \bigcap _{\beta \in (\dTronc{n+1},d_{\mathrm{\min }} )} \!\! K^{\infty } _{\beta } (\mathbb{R},\mathcal{H}^{m} )\]
Moreover, $\sbd{\Sreste{n+1}}$ is null in a vicinity of $-\infty $, so: $\forall \beta >\dTronc{n+1},\; \sbd{\Sreste{n+1}}\in K^{\infty } _{\beta } (\mathbb{R},\mathcal{H}^{m} )$. This concludes the induction.
\end{description}

To complete the proof of the theorem, we must show that: $\forall d\in \mathbb{R},\; S = T_{\geqslant d} (\SDA ) + o_{\partial } (\bfr^{d} )$ when $\bfr\rightarrow \infty $. We will do it in $\Omega $, but it works the same in $\couche$. Let $d\in \mathbb{R}$ and $n\in \mathbb{N}$ s.t.\ $\dTronc{n}<d$. We have in $\Omega $:
\[s:= \khii S -\khii T_{\geqslant d} (\SDA ) = \Sreste{n} + \khii \cdot T_{<d} \circ  T_{>\dTronc{n}} (\SDA ) \quad \text{ with }\quad  T_{<d} \circ  T_{>\dTronc{n}} (\SDA ) \in \sum\limits _{d'<d} \mathcal{A}_{d'} (\Omega ).\]
So there exists $d'<d$ s.t.: $\forall m\in \mathbb{N},\; \sbd{s\,}\in K^{\infty } _{d'} (\mathbb{R},H^{m} (0,\Theta ))$. Thus for~any $(i,j)\in \mathbb{N}^2 $:
\[\sbd{r^{-d'} (r\partial _{r} )^{i} \partial _{\theta } ^{j} s\,} = e^{-d't} \partial _{t} ^{i} \partial _{\theta } ^{j} \sbd{s\,} \in  H^2 (\mathbb{R},H^2 (0,\Theta )) \subset L^{\infty } (\sbd \Omega ).\]
So $r^{-d'+i} \partial _{r} ^{i} \partial _{\theta } ^{j} s$ is also bounded. By~definition of $o_{\partial }$ (see \cref{def: petit o derivable}), this concludes the proof.
\end{demo}

\subsection{Proof of \texorpdfstring{\cref{DA des champs lointains concatenes}}{Theorem 5.2.4}: asymptotic behavior for far-and-layer fields-like problems}
\label{sec: DA des champs lointains concatenes}

\begin{textesansboite}
This proof is very similar to \cref{sec: DA des champs de coin}. The main difference is that we look at the asymptotic expansion when $\bfr\rightarrow 0$ (i.e.\ $t\rightarrow -\infty $) instead of $\bfr\rightarrow \infty $ (i.e.\ $t\rightarrow +\infty $). So this time we use a truncation function $\khif \in \mathcal{C}^{\infty } (\Pi )$ that is equal to 1 in a vicinity of $\bfr=0$ and to 0 in a vicinity of infinity. In addition, we can assume that $\khif f=0$. Moreover the Laplace transform of $\sbd{\khif\bfu}$ is first defined in a left half-plane of the complex plane (instead of a right one), and then extended to the right.\\

Applying the change of variables $(x,y) \rightsquigarrow (t,\theta )$ on the equations satisfied by $\khif\bfu$ (that one can easily deduce from \cref{eq: DA des CL concatenes : hyp}) and then the Laplace transform yields that $s:=\sbd{\khif\bfu}$ satisfies:
\begin{equation}{\label{eq: pb verifie par la transformee de Laplace du champ lointain}}
\left\{\begin{array}{r@{\;=\;}ll}
	\mu _{0} (\partial _{\theta } ^2 +\lambda ^2 ) \hat{s}(\lambda ,\theta ) & \tilde{f}-\omega ^2 \rho _{0} \hat{s}(\lambda -2,\theta ) & \text{if } \ \theta \in (0,\Theta )\\[0.5mm]
	\partial _{\theta } ^2  \hat{s}(\lambda ,\theta ) & \tilde{g} & \text{if } \ \theta \in (-1,0)\\[0.5mm]
	\hat{s}(\lambda ,0^{+} ) - \hat{s}(\lambda ,0^{-}) &0& \\[0.5mm]
	\partial _{\theta } \hat{s}(\lambda ,0^{-}) & \tilde{h} & \\[0.5mm]
	\hat{s}(\lambda ,\Theta ) = \hat{s}(\lambda ,-1) & 0 & 
\end{array}\right.\end{equation}
for some functions $\tilde{f}, \tilde{g}, \tilde{h}$ depending on $g,h,\khif$. This system has the form of \cref{eq: EDO en theta a resoudre}, so \cref{resolution de l'EDO en theta} gives the tools to solve it and to extend $\hat{s}$ to the right (by steps of 2 here). This is stated in \cref{etape de base du DA des champs lointains}, whose proof is very similar to \cref{etape de base du DA des champs de coin}.
\end{textesansboite}

\begin{proposition}{}
\label{etape de base du DA des champs lointains}%
Let $d\in \mathbb{R}$, $m\in \mathbb{N} \setminus \{0,1\}$, $s\in  \bigcap\limits _{\beta <d} K^{\infty } _{\beta } (\mathbb{R},\mathcal{H}^{m} )$, $\tilde{f}\in  \!\bigcap\limits _{\beta <d+2} \! K^{\infty } _{\beta } (\mathbb{R},H^{m-2} (0,\Theta ))$, $\tilde{g}\in  \!\bigcap\limits _{\beta <d+2} \! K^{\infty } _{\beta } (\mathbb{R},H^{m-2} (-1,0))$, $\tilde{h}\in  \!\bigcap\limits _{\beta <d+2} \! K^{\infty } _{\beta } (\sbd \Gamma )$ and $d_{\mathrm{\max }} = \max (\TZ*\cap [d,d+2))$. We assume: \vspace{-2mm}
\[\left\{\begin{array}{rl}
	\mu _{0} \Delta  s+e^{2t} \omega ^2 \rho _{0} s = \tilde{f} &\text{ in }\sbd \Omega \\[0.5mm]
	\partial _{\theta } ^2 s = \tilde{g}  &\text{ in }\sbd\couche \\[0.5mm]
	s_{|\theta =0^{+}} - s_{|\theta =0^{-}} =0  &\text{ on }\sbd \Gamma  \\[0.5mm]
	\mu _{1} \partial _{\theta } s_{|\theta =0^{-}} = \tilde{h} &\text{ on }\sbd \Gamma  \\[0.5mm]
	s =0 &\text{ on }\sbd\bordgauche \cup \sbd\borddroit
\end{array}\right.\] Then there exist $s_{\mathrm{r}} \in  \!\bigcap\limits _{d_{\mathrm{\max }} <\beta <d+2}  \! K^{\infty } _{\beta } (\mathbb{R},\mathcal{H}^{m} )$ and $(c_{q} )\in \mathbb{C}^{\TZ*\cap [d,d+2)}$ s.t.\ $\displaystyle s = \sum _{q\in \TZ*\cap [d,d+2)} c_{q} \,\sbd{\phi_{q} } + s_{\mathrm{r}}$.\fbeq[0.5]
\end{proposition}

\refstepcounter{hypothese}
\begin{paragraphesansboite}{Assumption \thehypothese:}
\label{hyp: u a support compact}%
Replacing $\bfu$ by $\khii \bfu$, we assume without loss of generality that $\bfu$ has a compact support.
\end{paragraphesansboite}

\begin{textesansboite}
For~any $q\in -\TN$, let $\sigma _{q} (\bfu) := \sigma _{q} (\bfu_{|\Omega } )$, where $\sigma _{q} (\bfu_{|\Omega } )$ is set in \cref{def: sigma dans H1 + khi A}. This quantity vanishes when $q$ is small enough. \Cref{regularite de u variationnel} is the analogue of \cref{regularite de S variationnel}.
\end{textesansboite}

\begin{proposition}{}
\label{regularite de u variationnel}%
Let $d\in \mathbb{R}_+$,\vspace{-1mm}
\[\bfu^{\mathcal{A}} := T_{\leqslant d} \Bigg[ \big\langle  {-}k_{0} ^2 \,\ROmegaPi \big\rangle  \bigg( \frac{1}{\mu _{1}} \RcouchePi(g^{0} ) + \frac{1}{\mu _{1}} \RGammaPi(h^{0} ) + \fsum_{q\in -\TN} \sigma _{q} (\bfu) \,\phi_{q} \bigg)\Bigg] \vspace{-2mm}\]
in $\mathcal{A}(\Pi )$ and $\bfu^{\mathrm{v}} := \bfu-\khif \bfu^{\mathcal{A}}$. There exists a value of $d\in \mathbb{R}_+$ s.t.\ we have
\[\bfu^{\mathrm{v}} _{|\Omega } \in  H^{1} (\Omega ) \quad \  \text{ and } \quad \  \forall m\in \mathbb{N},\ \ \sbd{\khif \bfu^{\mathrm{v}} } \in K^{\infty } _{-1} (\mathbb{R},\mathcal{H}^{m} ).\]\fbeq[1.1]
\end{proposition}

\newcommand{\uAbarreRef}{\bfu^{\Abarre}_-} 
\newcommand{\dregu}{0} 
\begin{demo}
We consider any $d\in \mathbb{R}_+$ and we will fix it later. The proof has five steps: writing the equations satisfied by $\bfu^{\mathrm{v}}$, and then showing that $\bfu^{\mathrm{v}}$ is regular in the layer, that it~is $H^{1}$ in $\Omega $, that it~is regular in $\Omega $ w.r.t $t$, and that it~is regular in $\Omega $ w.r.t.\ $\theta $. For~any $r_{1} ,r_{2} \in \mathbb{R}$ s.t.\ $r_{1} <r_{2} <r_{f}$, we denote $\chi \in \mathcal{C}^{\infty } (\mathbb{R}^2 )$ a radial function equal to 1 on $B(0,r_{1} )$ and 0 on $\mathbb{R}^2  \setminus B(0,r_{2} )$ ($r_{1} ,r_{2} $ are implicit in this~notation).\\

\underline{Step 1:} Similarly to the proofs of \cref{cadre fonctionnel des champs lointains,cadre fonctionnel pour les champs de coin,regularite de S variationnel}, one can show that:
\[\left\{\begin{array}{@{\,}r@{\;=\;}l@{\,}l@{}}
	(\mu _{0} \Delta  +\omega ^2 \rho _{0} )\bfu^{\mathrm{v}} & f_{1} := f-(\mu _{0} \Delta  +\omega ^2 \rho _{0} )(\khif \bfu^{\mathcal{A}} ) &\text{ in }\Omega  \\[0.5mm]
	\mu _{1} \partial _{Y} ^2 \bfu^{\mathrm{v}} & g_{1} :=g-\mu _{1} \partial _{Y} ^2 (\khif \bfu^{\mathcal{A}} ) &\text{ in }\couche \\[0.5mm]
	\bfu^{\mathrm{v}} _{|y=0^{+}} - \bfu^{\mathrm{v}} _{|y=0^{-}} & 0 &\text{ on }\Gamma  \\[1mm]
	\mu _{1} \partial _{Y} \bfu^{\mathrm{v}} _{|Y=0^{-}} & h_{1} :=h-\mu _{1} \partial _{Y} (\khif \bfu^{\mathcal{A}} ) &\text{ on }\Gamma  \\[1mm]
	\bfu^{\mathrm{v}} & 0 &\text{ on }\bordgauche \cup \borddroit
\end{array}\right.\]
with
\[\left\{\begin{array}{@{}r@{\;\in \;}l@{}}
	(\mu _{0} \Delta  +\omega ^2 \rho _{0} )\bfu^{\mathcal{A}} & \sum\limits _{d'>d} \! \mathcal{A}_{d'} (\Omega ) \\[0.5mm]
	\partial _{Y} ^2 \bfu^{\mathcal{A}} _{|\couche} - T_{\leqslant d} (g^{0} ) & \sum\limits _{d'>d} \! \mathcal{A}_{d'} (\couche)  \\[0.5mm]
	\partial _{Y} \bfu^{\mathcal{A}} _{|Y=0^{-}} - T_{\leqslant d} (h^{0} ) & \sum\limits _{d'>d} \! \mathcal{A}_{d'} (\Gamma )
\end{array}\right.\]
Moreover, given that $\chi  f=0$ for~any $r_{1} <r_{2} <r_{f}$, one can check that for $d$ large enough and up to decreasing $r_{f}$ (without loss of generality), we have $\sbd{\chi  f_{1}} \in K^{\infty } _{\dregu} (\mathbb{R},H^{m} (0,\Theta ))$, $\sbd{\chi  g_{1}} \in K^{\infty } _{\dregu} (\mathbb{R},H^{m} (-1,0))$ and $\sbd{\chi  h_{1}} \in K^{\infty } _{\dregu} (\sbd \Gamma )$ for~any $m\in \mathbb{N}$.\\

\underline{Step 2:} Let us show that: $\forall r_{1} <r_{2} <r_{f} , \forall m\in \mathbb{N},\; \sbd{\chi  \bfu^{\mathrm{v}} _{|\couche} } \in K^{\infty } _{0} (\mathbb{R},H^{m} (-1,0))$. By~assumption, for~any $x>0$, $\bfu^{\mathrm{v}}$ in $H^{1}$ in a vicinity of $\{x\}\times [-1,0] \subset \couche$, so $\bfu^{\mathrm{v}} _{|\{x\}\times [-1,0]}$ is well-defined in $L^2 (\{x\}\times [-1,0])$. Moreover, $g_{1}$ and $h_{1}$ are $o_{\partial } (x^{d} )$, so they are differentiable on $(0,r_{f} )\times [-1,0]$ (up to decreasing $r_{f}$). Hence, $Y\in (-1,0)\mapsto  \bfu^{\mathrm{v}} (x,Y)$ is $\mathcal{C}^2 $ and the ODE it satisfies (see Step 1) is explicitly solvable:
\[\forall (x,Y)\in (0,r_{f} )\times [-1,0], \quad  \bfu^{\mathrm{v}} (x,Y) = \int _{-1} ^{0} \big( (Y-Y')^{+} -Y-1 \big)\, g_{1} (x,Y')\,\mathrm{d}Y' +h_{1} (x)\cdot (Y+1).\]
Since $\sbd{\chi  g_{1}} \in K^{\infty } _{\dregu} (\mathbb{R},H^{m} (-1,0))$ for~any $m$ and $\sbd{\chi  h_{1}} \in K^{\infty } _{\dregu} (\sbd \Gamma )$, we deduce $\sbd{\chi  \bfu^{\mathrm{v}} _{|\couche} } \in K^{\infty } _{\dregu} (\mathbb{R},H^{m} (-1,0))$.\\

\underline{Step 3:} Let us show that $\bfu^{\mathrm{v}} _{|\Omega } \in H^{1} (\Omega )$. First, step 2 implies $\sbd{\chi  \bfu^{\mathrm{v}} _{|\Gamma } } \in K^{\infty } _{\dregu} (\sbd \Gamma )$ so $\chi  \bfu^{\mathrm{v}} _{|\Gamma } \in H^{1/2} (\Gamma )$. Since $\bfu \in H^{1} _{\mathrm{loc}}$ and $\bfu$ has a compact support by \cref{hyp: u a support compact}, we get $\bfu^{\mathrm{v}} _{|\Gamma } \in H^{1/2} (\Gamma )$. So the following system has a solution in $H^{1} (\Omega )$ by \cref{pb de Helmholtz variationnel}: \[\left\{\begin{array}{r@{\;=\;}ll}
	\mu _{0} \Delta  u^{\mathrm{v}} +\omega ^2 \rho _{0} u^{\mathrm{v}} & f_{1} &\text{ in }\Omega  \\[0.5mm]
	u^{\mathrm{v}} _{|y=0^{+}} & \bfu^{\mathrm{v}} _{|y=0^{-}} &\text{ on }\Gamma  \\[1mm]
	u^{\mathrm{v}} & 0 &\text{ on }\bordgauche
\end{array}\right.\]
But $\bfu^{\mathrm{v}}$ also satisfies these equations, so $u^{\mathrm{v}} +\khif \bfu^{\mathcal{A}} _{|\Omega }$ satisfies the same problem as $\bfu_{|\Omega }$, and they both belong to $H^{1} (\Omega )+\khiz \mathcal{A}(\Omega )$. Then by uniqueness in \cref{cadre fonctionnel des champs lointains}, $u^{\mathrm{v}} +\khif \bfu^{\mathcal{A}} _{|\Omega } =\bfu_{|\Omega }$, i.e.\ $u^{\mathrm{v}} = \bfu^{\mathrm{v}} _{|\Omega }$. Thus, $\bfu^{\mathrm{v}} _{|\Omega } \in H^{1} (\Omega )$.\\

\underline{Step 4:} Let $v:\Omega \rightarrow \mathbb{C}$ defined by $v := \bfu^{\mathrm{v}} _{|\Omega } - v_{\mathrm{lift}}$ with $v_{\mathrm{lift}} (r,\theta ) := \big(1- \frac{\theta }{\Theta } \big) \bfu^{\mathrm{v}} _{|\Gamma } (r)$ in polar coordinates. Let $K^{0} _{\beta } (\sbd \Omega ):=K^{0} _{\beta } (\mathbb{R},L^2 (0,\Theta ))$ for~any $\beta \in \mathbb{R}$. We will show by induction on $n$ that: \[\forall n\in \mathbb{N}^*,\; \forall r_{1} <r_{2} <r_{f} ,\qquad  \partial _{t} ^{n} \sbd{\chi \, v}\in K^{0} _{-1} (\sbd \Omega ) \quad  \text{ and } \quad  \nabla  \partial _{t} ^{n} \sbd{\chi \, v}\in L^2 (\sbd \Omega ).\]
For the initial case, we have on the one hand $\chi \bfu^{\mathrm{v}} _{|\Omega } \in H^{1} (\Omega )$ so $\sbd{\chi \bfu^{\mathrm{v}} _{|\Omega } }\in K^{0} _{-1} (\sbd \Omega )$ and $\nabla  \sbd{\chi \bfu^{\mathrm{v}} _{|\Omega } }\in L^2 (\sbd \Omega )$, and on the other $\sbd{\chi  \bfu^{\mathrm{v}} _{|\Gamma } }\in K^{1} _{\dregu} (\sbd \Gamma )$ so $\sbd{\chi  v_{\mathrm{lift}} }\in K^{0} _{-1} (\sbd \Omega )$ and $\nabla  \sbd{\chi  v_{\mathrm{lift}} }\in L^2 (\sbd \Omega )$. Thus the initial case is proven and only the inductive step remains to prove. We assume the property at rank $n$ and we will show it at rank $n+1$, using the method of finite differences. Step 1 implies that: \[\left\{\begin{array}{r@{\;=\;}ll}
	(\mu _{0} \Delta +\omega ^2 \rho _{0} )(\chi  v) & f_{2} := \chi  f_{1} -\chi \,\mu _{0} \Delta  v_{\mathrm{lift}} - 2\mu _{0} \nabla  \chi  \cdot \nabla  \bfu^{\mathrm{v}} -\mu _{0} \Delta  \chi  \cdot \bfu^{\mathrm{v}} &\text{ in }\Omega  \\[0.5mm]
	\chi  v &0 &\text{ on }\bordgauche \cup \Gamma 
\end{array}\right.\]
In addition, we have $\sbd{f_{2} } \in K^{n} _{-2} (\mathbb{R},L^2 (0,\Theta ))$ since:
\begin{itemize}
	\item by step 1, $\sbd{\chi  f_{1}} \in K^{\infty } _{\dregu} (\mathbb{R},L^2 (0,\Theta ))$, so $\sbd{\chi  f_{1}} \in K^{\infty } _{-2} (\mathbb{R},L^2 (0,\Theta ))$,
	\item by step 2, $\sbd{\chi  \bfu^{\mathrm{v}} _{|\Gamma } }\in K^{\infty } _{\dregu} (\sbd \Gamma )$, so $\sbd{\chi  v_{\mathrm{lift}} } \in K^{\infty } _{\dregu} (\mathbb{R},H^2 (0,\Theta ))$, and applying it at $(r_{1} ',r_{2} '):= (r_{2} ,\frac{r_{2} +r_{f} }{2} )$ gives $\sbd{\chi  \Delta  v_{\mathrm{lift}} } \in K^{\infty } _{-2} (\mathbb{R},L^2 (0,\Theta ))$,
	\item and by induction hypothesis applied to $(r_{1} ',r_{2} '):= (r_{2} ,\frac{r_{2} +r_{f} }{2} )$, $\partial _{t} ^{k} \nabla  v$ and $\partial _{t} ^{k} v$ are $L^2 $ on $\{r_{1} <r<r_{2} \}$ for~any $k\in [\![0,n]\!]$.
\end{itemize}
Let $s:=\partial _{t} ^{n} (\sbd{\chi  v})$. By~changing variables $(x,y) \rightsquigarrow (t,\theta )$ in the previous system and applying $\partial _{t} ^{n}$, we get: \[\left\{\begin{array}{r@{\;=\;}ll}
	e^{-2t} \mu _{0} \Delta s+\omega ^2 \rho _{0} s & f_{3} := e^{-2t} \partial _{t} ^{n} (e^{2t} \sbd{f_{2} }) - \sum\limits _{k=0} ^{n-1} \omega ^2 \rho _{0} 2^{n-k} \partial _{t} ^{k} (\sbd{\chi  v}) &\text{ in }\sbd \Omega  \\[-1mm]
	s & 0 &\text{ on }\sbd\bordgauche \cup \sbd \Gamma 
\end{array}\right.\] with $f_{3} \in K^{0} _{-2} (\sbd \Omega )$ by induction hypothesis. And the variational formulation of this is: \[\forall \varphi \in K^{0} _{-1} (\sbd \Omega ), \ \nabla  \varphi \in L^2 (\sbd \Omega ) \text{ and } \varphi _{|\sbd\bordgauche \cup \sbd \Gamma }=0 \ \Rightarrow  \quad  \int _{\sbd \Omega } (\mu _{0} \nabla  s \cdot \nabla  \varphi +\omega ^2 \rho _{0} e^{2t} s \varphi ) = \int _{\sbd \Omega } e^{2t} f_{3} \varphi .\]
Let us denote $D_{\eta } \varphi (t,\theta ) := \frac{\varphi (t+\eta ,\theta )-\varphi (t,\theta )}{\eta }$ for~any $\eta \in \mathbb{R}^*$ and any function $\varphi $. Taking $\varphi  := D_{-\eta } D_{\eta } \bar{s}$ and discretely integrating by parts $D_{-\eta }$ gives:
\[\int _{\sbd \Omega } (\mu _{0} |D_{\eta } \nabla  s|^2  -\omega ^2 \rho _{0} e^{2t} |s|^2 ) = \int _{\sbd \Omega } e^{2t} f_{3} \cdot  D_{-\eta } D_{\eta } \bar{s}\]
Then by coercivity (since $\operatorname{Im}(\omega ) \neq 0$), we have for~any $\delta >0$:
\[\|D_{\eta } \nabla  s\|_{L^2 (\sbd \Omega )}^2  + \|D_{\eta } s\|_{K^{0} _{-1} (\sbd \Omega )}^2  \lesssim  \|f_{3} \|_{K^{0} _{-2} (\sbd \Omega )} \|D_{-\eta } D_{\eta } s\|_{K^{0} _{-1} (\sbd \Omega )} 
\lesssim  {\textstyle \frac{1}{\delta } } \|f_{3} \|_{K^{0} _{-2} (\sbd \Omega )}^2  + \delta  \|D_{\eta } \partial _{t} s\|_{L^2 (\sbd \Omega )}^2 \]
Taking $\delta $ small enough and moving $\|D_{\eta } \partial _{t} s\|_{L^2 (\sbd \Omega )} ^2 $ from the right-hand side to the left-hand one, we get:
\[\|\partial _{t} \nabla  s\|_{L^2 (\sbd \Omega )}^2  + \|\partial _{t} s\|_{K^{0} _{-1} (\sbd \Omega )}^2  
\leqslant  \limsup _{\eta \rightarrow 0} \|D_{\eta } \nabla  s\|_{L^2 (\sbd \Omega )}^2  + \|D_{\eta } s\|_{K^{0} _{-1} (\sbd \Omega )}^2 
\lesssim  \|f_{3} \|_{K^{0} _{-2} (\sbd \Omega )}^2  <\infty .\]
This completes the induction.\\

\underline{Step 5:} We have shown that for~any $(n,i)\in \mathbb{N}\times \{0,1\}$ and $r_{1} <r_{2} <r_{f}$, $\partial _{t} ^{n} \partial _{\theta } ^{i} \sbd{\chi  v} \in  K^{0} _{-1} (\sbd \Omega )$. Since $\sbd{\chi  v_{\mathrm{lift}} } \in  K^{\infty } _{-1} (\mathbb{R},H^{1} (0,\Theta ))$, we also have $\partial _{t} ^{n} \partial _{\theta } ^{i} \sbd{\chi  \bfu^{\mathrm{v}} _{|\Omega } } \in  K^{0} _{-1} (\sbd \Omega )$ for these $(n,i)$.\\
To treat higher-order $\theta $-derivatives, we begin with the equality $(\mu _{0} \Delta +\omega ^2 \rho _{0} )(\chi  \bfu^{\mathrm{v}} ) = \chi  f_{1} - 2\mu _{0} \nabla  \chi  \cdot \nabla  \bfu^{\mathrm{v}} -\mu _{0} \Delta  \chi  \cdot \bfu^{\mathrm{v}}$ in $\Omega $, which follows from step 1. Applying the change of variables $(x,y) \rightsquigarrow (t,\theta )$, we get:
\begin{equation}{\label{eq: regularite de u variationnel etape 5}}
(e^{-2t} \mu _{0} \Delta +\omega ^2 \rho _{0} )(\chi  \bfu^{\mathrm{v}} ) = \sbd{\chi  f_{1} }- 2\mu _{0} \sbd{\nabla  \chi  \cdot \nabla  \bfu^{\mathrm{v}} } -\mu _{0} \sbd{\Delta  \chi  \cdot \bfu^{\mathrm{v}} }.\end{equation}
Now $\sbd{\chi  f_{1}} \in K^{\infty } _{\dregu} (\mathbb{R},H^{m} (0,\Theta ))$ for~any $m\in \mathbb{N}$ by step 1. So~deriving \cref{eq: regularite de u variationnel etape 5} w.r.t.\ $t$ and $\theta $ enough times gives by induction: $\forall (n,i)\in \mathbb{N}^2 , \forall r_{1} <r_{2} <r_{f} ,\; \partial _{t} ^{n} \partial _{\theta } ^{i} \sbd{\chi  \bfu^{\mathrm{v}} _{|\Omega } } \in  K^{0} _{-1} (\sbd \Omega )$. Combining it with step 2, we get $\sbd{\chi  \bfu^{\mathrm{v}} } \in K^{\infty } _{-1} (\mathbb{R},\mathcal{H}^{m} )$. Finally, we apply it to $\chi :=\khii$ to conclude.
\end{demo}

\begin{textesansboite}
Finally, \cref{etape de base du DA des champs lointains,regularite de u variationnel} give all the ingredients to prove \cref{DA des champs lointains concatenes}. The proof is very similar to the one of \cref{DA des champs de coin} on page~\pageref{preuve: DA des champs de coin 2}, so we do~not go into details again.
\end{textesansboite}

\section{Construction of the asymptotic expansion}
\label{sec: construction des champs}

\subsection{Inductive construction of the fields}
\label{sec: construction des termes}

\begin{textesansboite}
The next definition summarizes how to build all the fields $u_{p,\ell }$, $U_{p,\ell }$ and $S_{p,\ell }$. We recall that the spaces $\mathcal{A}(\Omega )$ and $\mathcal{A}(\Pi )$ are defined in \cref{def: espaces A}, the coefficients $\cmacro$ and $\cmicro$ are defined by \cref{eq: def cmaicro}, $\intvl{a,b}$ is equal to $\{c\in \mathbb{R} \mid c-a\in \NN \text{ and } b-c\in \mathbb{N}\}$ and $(\mathcal{U}_{n} )$ is defined in \cref{eq: correcteurs de couche}. The truncation function $\khiz$ is introduced in \cref{sec: cadre fonctionnel des champs lointains} and localized near 0, while $\khii$ is defined in \cref{sec: cadre fonctionnel des champs de coin} and localized near infinity.
\end{textesansboite}
	
\begin{definitionsecable}{construction of the fields}
\label{def: construction des termes}%
For~any $(p,\ell )\in (\mathbb{R} \setminus \NN)\times \mathbb{N}$ we denote by convention $u_{p,\ell } =0$, $U_{p,\ell } =0$ and $S_{p,\ell } =0$. We define by induction on $p\in \NN$ that for~any $\ell \in \mathbb{N}$:
\begin{itemize}
	\item $\displaystyle U_{p,\ell } :(x,Y)\in \couche \mapsto  \sum_{n=0}^\infty  (\partial _{x} ^2 +k_{1} ^2 )^{n} \partial _{y} u_{p-1-2n,\ell } (x,0) \cdot  \mathcal{U}_{n} (Y)$. \leavevmode\hfill\inlineEqNb\label{eq: definition des champs de couche}
	\item $u_{p,\ell }$ is the unique solution in $H^{1} (\Omega )+\khiz \mathcal{A}(\Omega )$ of :
	\begin{equation}{\label{eq: definition des champs lointains}}
		\hspace*{-5mm}\left\{\begin{array}{r@{\;=\;}ll@{}}
		\mu _{0} \Delta  u_{p,\ell } +\omega ^2 \rho _{0} u_{p,\ell } & \fsource\,\delta _{p,0} \,\delta _{\ell ,0} &\text{ in }\Omega  \\[0.5mm]
		u_{p,\ell } & 0 &\text{ on }\bordgauche\\[0.5mm]
		u_{p,\ell } & U_{p,\ell } &\text{ on }\Gamma \\[-2mm]
		\sigma _{d} (u_{p,\ell } ) &\displaystyle \!\!\sum _{p'\in \intvl{0,p+d}} \sum _{\substack{d'\in \TZ*\\ p'-d'\in \intvl{0,p}}} \!\sum _{\ell '=0} ^{\ell } \,\cmacro_{d,d',p-p',\ell -\ell '} \cdot  \sigma _{d'} (\SDA_{p',\ell '} ) &\ \forall d\in -\TN
	\end{array}\right.\vspace{-2mm}\end{equation}
	\item $S_{p,\ell }$ is the unique solution in $\Vcoin+\khii \mathcal{A}(\Pi )$ of :
	\begin{equation}{\label{eq: definition des champs de coin}}
		\left\{\begin{array}{r@{\;=\;}ll}
		\operatorname{div}(\mu  \nabla  S_{p,\ell } ) & -\omega ^2 \rho  \,S_{p-2,\ell } &\text{ in }\Omega _{1} \\[0.5mm]
		S_{p,\ell } & 0 &\text{ on }\partial  \Omega _{1} \\[-2mm]
		\sigma _{d} (S_{p,\ell } ) &\displaystyle \!\!\sum _{p'\in \intvl{0,p-d}} \sum _{\substack{d'\in \TZ*\\ p'+d'\in \intvl{0,p}}} \!\sum _{\ell '=0} ^{\ell } \,\cmicro_{d,d',p-p',\ell -\ell '} \cdot  \sigma _{d'} (\bfuDA_{p',\ell '} ) &\ \forall d\in \TN
	\end{array}\right.\end{equation}\fbeq
\end{itemize}\fbi
\end{definitionsecable}

\begin{textesansboite}
In \cref{def: construction des termes}, $\sigma _{d} (u_{p,\ell } )$ and $\sigma _{d} (S_{p,\ell } )$ are to be understood as in \cref{def: sigma dans H1 + khi A,def: sigma dans V + khi A} resp., whereas $\sigma _{d} (\bfuDA_{p,\ell } )$ and $\sigma _{d} (\SDA_{p,\ell } )$ follow \cref{def: sigma sur Amaicro}. The next result proves the well-definition of the fields $u_{p,\ell }$, $U_{p,\ell }$ and $S_{p,\ell }$ and their asymptotic behaviors $\bfuDA_{p,\ell }$ and $\SDA_{p,\ell }$. It also verifies that the $\bfuDA_{p,\ell }$ and $\SDA_{p,\ell }$ satisfies the inductive relations \cref{eq: bfuDA SDA} that were assumed in \cref{ansatz: existence de bfuDApl et SDApl}. Moreover we prove an additional property of regularity on $\Gamma $ that is needed to prove the rest, and will be useful later in \cref{chap: ME}.
\end{textesansboite}

\begin{theoremesecable}{}
\label{validation construction et DA des champs}%
The families $(u_{p,\ell } )_{p,\ell }$, $(U_{p,\ell } )_{p,\ell }$ and $(S_{p,\ell } )_{p,\ell }$ as in \cref{def: construction des termes} are well-defined. Moreover, for~any $(p,\ell )\in \NN\times \mathbb{N}$, they satisfy the following.
\begin{enumerate}
	\item There exists $\bfuDA_{p,\ell } \in \Amacro(\Pi )$ and $\SDA_{p,\ell } \in \Amicro(\Pi )$ s.t., for~any $d\in \mathbb{R}$, $\bfu_{p,\ell } = T_{\leqslant d} (\bfuDA_{p,\ell } ) +o_{\partial } (\bfr^{d} )$ when $\bfr\rightarrow 0$ and $S_{p,\ell } = T_{\geqslant d} (\SDA_{p,\ell } ) + o_{\partial } (\bfr^{d} )$ when $\bfr\rightarrow \infty $ (where $\bfu_{p,\ell }$ is defined by \cref{eq: def bfu bfr}). These formal series have the following inductive expressions
	\makeatletter
	\def\pushright#1{\ifmeasuring@#1\else\omit\hfill$\displaystyle #1$\fi\ignorespaces} 
	\makeatother
	\begin{equation}{\label{eq: bfuDA SDA}}
		\left\{\begin{array}{@{\,}r@{}l@{}}
		\bfuDA_{p,\ell } &{}=\displaystyle \big\langle  {-k_{0} ^2 \,\ROmegaPi} \big\rangle  \bigg(\displaystyle {-\RcouchePi}\circ (\partial _{x|\couche} ^2 +k_{1} ^2 )(\bfuDA_{p-2,\ell } ) + \frac{\mu _{0} }{\mu _{1}} \RGammaPi\circ \partial _{y|\Gamma ,y=0^{+}} (\bfuDA_{p-1,\ell } ) \qquad  \qquad  \quad  \\
		& \pushright{ {}+ \fsum_{d\in \TZ*} \sigma _{d} (\bfuDA_{p,\ell } ) \,\phi_{d} \bigg) }\\
		\SDA_{p,\ell } &{}=\displaystyle \left\langle  {-\RcouchePi}\circ \partial _{X|\couche} ^2 , \frac{\mu _{0} }{\mu _{1}} \RGammaPi\circ \partial _{Y|\Gamma ,Y=0^{+}} \right\rangle  \bigg( {- k_{0} ^2  \ROmegaPi} (\SDA_{p-2,\ell |\Omega } ) \\
		& \pushright{ {}- k_{1} ^2  \RcouchePi (\SDA_{p-2,\ell |\couche} ) + \fsum_{d\in \TZ*} \sigma _{d} (\SDA_{p,\ell } ) \,\phi_{d} \bigg) }
		\end{array}\right.
	\end{equation}
	with the convention $\bfuDA_{p',\ell '} :=\SDA_{p',\ell '} :=0$ when $p'\in \mathbb{R} \setminus \NN$.
	\item For~any $m\in \mathbb{N}$, $\partial _{x} ^{m} u_{p,\ell |\Gamma }$ and $\partial _{x} ^{m} \partial _{y} u_{p,\ell |\Gamma }$ belong to $H^{1} _{0} (\Gamma )+\khiz \mathcal{A}(\Gamma )$.
\end{enumerate}
\end{theoremesecable}

\begin{demo}
Let us show by induction on $p$ that for~any $\ell $:
\begin{itemize}[itemindent=5mm]
	\item[(H1)] $u_{p,\ell }$, $U_{p,\ell }$ and $S_{p,\ell }$ exist and are unique,
	\item[(H2)] $\bfuDA_{p,\ell }$ and $\SDA_{p,\ell }$ exist, are the asymptotic expansions w.r.t.\ $\bfr$ of $\bfu_{p,\ell }$ and $S_{p,\ell }$ and satisfy \cref{eq: bfuDA SDA},
	\item[(H3)] $\forall j\in \{0,1\}, \forall m\in \mathbb{N},\; \partial _{x} ^{m} \partial _{y} ^{j} u_{p,\ell |\Gamma } \in  H^{1} _{0} (\Gamma )+\khiz \mathcal{A}(\Gamma )$.
\end{itemize}
Since all fields vanish for $p<0$, the initial case is trivial and only the inductive step remains to prove. Let $p\in \NN$. Let us assume (H1)--(H3) at any rank $p'<p$ and prove it at rank $p$.\medskip

\underline{Existence and uniqueness of the fields:}
\begin{itemize}
	\item \underline{$U_{p,\ell }$\,:} By~(H3), the $\partial _{y} u_{p',\ell |\Gamma }$ with $p' <p$ are all in $H^{m} _{\mathrm{loc}} (\Gamma )$ for~any $m\in \mathbb{N}$. So the functions $(\partial _{x} ^2 +k_{1} ^2 )^{n} \partial _{y} u_{p-1-2n,\ell |\Gamma }$ are continuous and Formula \cref{eq: definition des champs de couche} is well-defined at any point of $\couche$.
	\item \underline{$u_{p,\ell }$\,:} By~(H2), $\SDA_{p',\ell '}$ exists for~any $p'<p$ and $\ell '$, so the coefficients $\sigma _{d'} (\SDA_{p',\ell '} )$ are well-defined. Thus, \cref{eq: definition des champs lointains} defines well $\sigma _{d} (u_{p,\ell } )$. And we have: $\sigma _{d} (u_{p,\ell } ) \neq 0 \Rightarrow  \intvl{0,p+d} \neq  \varnothing  \Rightarrow  d\geqslant -p$, so $(\sigma _{d} (u_{p,\ell } ))_{d\in -\TN}$ has finite support. Next, we apply \cref{cadre fonctionnel des champs lointains}. Its hypotheses are satisfied, because (H3) implies:\vspace{-1mm}
	\[U_{p,\ell |\Gamma } = \sum_{n=0}^\infty  (\partial _{x} ^2 +k_{1} ^2 )^{n} \partial _{y} u_{p-1-2n,\ell } (x,0) \cdot  \mathcal{U}_{n} (0) \,\in \,H^{1} _{0} (\Gamma )+\khiz \mathcal{A}(\Gamma ) \,\subset \, H^{1/2} _{00} (\Gamma )+\khiz \mathcal{A}(\Gamma ) \vspace{-1mm}\] 
	\item \underline{$S_{p,\ell }$\,:} Similarly, (H2) implies that $(\sigma _{d} (S_{p,\ell } ))_{d\in \TN}$ is well defined and has finite support. Next, we apply \cref{DA des champs de coin}. We need to check its hypotheses, i.e.\ there exist $f:\Omega _{1} \rightarrow \mathbb{C}$ and $\varphi :\Pi \rightarrow \mathbb{C}$ s.t.\ $\omega ^2 \rho \,S_{p-2,\ell } = F+\khii \varphi $, $(1+r)F\in L^2 (\Omega _{1} )$, $\varphi _{|\Omega } \in \mathcal{A}(\Omega )$ and $\varphi _{|\couche} \in \mathcal{A}(\couche)$. By~(H2), there exists $\psi \in \mathcal{A}(\Pi )$ s.t.\ $S_{p-2,\ell } =\psi  +o_{\partial } (\bfr^{-4} )$. So it suffices to set $\varphi :=\omega ^2 \rho  \psi $ and $F:=\omega ^2 \rho  \,S_{p-2,\ell } -\khii \varphi $.
\end{itemize}\medskip

\underline{Asymptotic expansions:}
\begin{itemize}
	\item \underline{$\bfu_{p,\ell }$\,:} We use \cref{DA des champs lointains concatenes}. To do so, we must check that there exist $g_{\Gamma } \in \Amacro(\Gamma )$ and $g_{\couche} \in \Amacro(\couche)$ s.t.:\vspace{-2mm}
	\[\forall d\in \mathbb{R}, \qquad  \left\{\begin{array}{r@{\;=\;}ll}
		\partial _{Y} ^2  U_{p,\ell } & T_{\leqslant d} (g_{\couche} )+o_{\partial } (x^{d} ) &\text{ in }\couche\\[0.5mm]
		\partial _{Y} U_{p,\ell |Y=0^{-}} & T_{\leqslant d} (g_{\Gamma } )+o_{\partial } (x^{d} ) &\text{ on }\Gamma 
	\end{array}\right.\]
	Given the definition of $U_{p,\ell }$ in \cref{eq: definition des champs de couche} and using \cref{derivation des DA dans Amaicro}, it suffices to show that, for~any $p'<p$ there exists $h\in \Amacro(\Gamma )$ s.t.: $\forall d\in \mathbb{R},\; \partial _{n} u_{p'\!,\ell |\Gamma } = T_{\leqslant d} (h)+o_{\partial } (x^{d} )$. But it derives from (H2).
	\item \underline{$S_{p,\ell }$\,:} Similarly, we use \cref{DA des champs de coin} thanks to (H2).
\end{itemize}\medskip

\underline{Regularity on $\Gamma $:} Let $j\in \{0,1\}$ and $m\in \mathbb{N}$. The asymptotic expansion of $\bfu_{p,\ell }$ implies: $\exists h\in \mathcal{A}(\Gamma ),\allowbreak\;\partial _{x} ^{m} \partial _{y} ^{j} u_{p,\ell |\Gamma } = h+o_{\partial } (x^{1} )$. So there exists $x_{0} \in \mathbb{R}_+^*$ s.t.\ $\partial _{x} ^{m} \partial _{y} ^{j} u_{p,\ell |\Gamma } -h$ is $H^{1}$ on $\Gamma \cap \{x<x_{0} \}=(0,x_{0} )\times \{0\}$ and it vanishes at 0. To prove (H3) at rank $p$, it remains to show that $\partial _{x} ^{m} \partial _{y} ^{j} u_{p,\ell } \in H^{1} ((\frac{x_{0} }{2} ,\infty )\times \{0\})$. To do so, it suffices to get $u \in H^{m+3} ((\frac{x_{0} }{2} ,\infty )\times (0, \frac{\delta }{2} ))$ with $\delta  :=\operatorname{dist}(\operatorname{supp}(f),\Gamma )$. But it follows from classical elliptic regularity because on one side $\mu _{0} \Delta u_{p,\ell } +\omega ^2 \rho _{0} u_{p,\ell } =0$ on $(\frac{x_{0} }{4} ,\infty )\times (0,\delta )$, and on the other $u_{p,\ell |\Gamma } =U_{p,\ell |\Gamma } \in H^{m+3} ((\frac{x_{0} }{4} ,\infty )\times \{0\})$ by \cref{eq: definition des champs de couche} and (H3).
\end{demo}

\begin{textesansboite}
Now, let us show some basic properties of the fields that are useful to prove the error estimates of \cref{sec: estimations d'erreur}. Further properties will be given in \cref{chap: complements sur le DA}.
\end{textesansboite}

\begin{lemmesansboite}{}
\label{coincidence des sigma singuliers}%
For~any $(p,\ell )\in \NN\times \mathbb{N}$, we have:
\begin{itemize}
	\item $\forall d\in -\TN,\; \sigma _{d} (u_{p,\ell } )=\sigma _{d} (\bfuDA_{p,\ell } )$,
	\item $\forall d\in \TN,\hphantom{-} \; \sigma _{d} (S_{p,\ell } )=\sigma _{d} (\SDA_{p,\ell } )$.
\end{itemize}
\end{lemmesansboite}

\begin{remarque}
We recall that, by \cref{def: sigma dans H1 + khi A,def: sigma dans V + khi A}, $\sigma _{d} (u_{p,\ell } )$ and $\sigma _{d} (S_{p,\ell } )$ are defined only for these values of $d$, that are the one associated to a diverging behavior (when $\bfr\rightarrow 0$, resp.\ $\bfr\rightarrow \infty $). Note that, in \cref{def: construction des termes}, we used the formulas for $\sigma _{d} (\bfuDA_{p,\ell } )$ and $\sigma _{d} (\SDA_{p,\ell } )$ given by \cref{condition de raccord} to define instead $\sigma _{d} (u_{p,\ell } )$ and $\sigma _{d} (S_{p,\ell } )$. \Cref{coincidence des sigma singuliers} shows that this substitution is~not a problem.
\end{remarque}

\begin{demo}
We have $u_{p,\ell } -T_{\leqslant 1} (\bfuDA_{p,\ell } )_{|\Omega } = o_{\partial } (r^{1} )$, so $\khiz \cdot (u_{p,\ell } -T_{\leqslant 1} (\bfuDA_{p,\ell } )_{|\Omega } )\in H^{1} (\Omega )$. By~\cref{def: sigma dans H1 + khi A}, this implies $\sigma _{d} (u_{p,\ell } )=\sigma _{d} \big(T_{\leqslant 1} (\bfuDA_{p,\ell } )_{|\Omega } \big) =\sigma _{d} (\bfuDA_{p,\ell } )$. Likewise for $S_{p,\ell }$.
\end{demo}

\begin{propositionsansboite}{}
\label{borne sur l}%
Let $p\in \NN$. There exists an integer $\ell _{p}$ that depends only on $(\Theta ,p)$ (and not on $\mu ,\rho ,\omega ,\Omega _{1}$,$\fsource$...) s.t., for~any $\ell >\ell _{p}$, 
the fields $u_{p,\ell }$, $U_{p,\ell }$ and $S_{p,\ell }$ all vanish.
\end{propositionsansboite}

\begin{remarque}
$\ell _{p}$ can be numerically computed thanks to the tools of \cref{sec: implementation des espaces A}. See examples in \cref{sec: exemples de structures de DA}, especially \cref{fig: graphe deg ln eps,tab: ordre du premier ln eps}.
\end{remarque}

\begin{demo}
Let $(d,d',p)\in (\TZ*)^2 \times \NN$. Let us show that there exists $n\in \mathbb{N}$ that depends only on $(d,d',p,\Theta )$ s.t.: $\forall \ell >n,\; \cmacro_{d,d',p,\ell } =0$.
\begin{sousdemo}
	\let\n=k%
	By~definition, we have $\cmacro_{d,d',p,\ell } := \sigma _{d,p,\ell } \circ \Eps^{-1} \circ  \langle \Rmicro\rangle  (\phi_{d'} )$ with $\Eps^{-1} \circ  \langle \Rmicro\rangle  (\phi_{d'} )\in \Aemacro(\Pi )$.	So, by construction of $\Aemacro(\Pi )$, there exists $n$ (that may depend on $\mu _{0} ,\mu _{1} ,\rho _{0} ,\rho _{1} ,\omega $) s.t.: $\forall \ell >n,\; \cmacro_{d,d',p,\ell } =0$. It remains to show that $n$ can be taken independent of $(\mu _{0} ,\mu _{1} ,\rho _{0} ,\rho _{1} ,\omega )$.\\
	By~construction, $\cmacro_{d,d',p,\ell }$ depends only on the coordinate of $\Eps^{-1} \circ  \langle \Rmicro\rangle  (\phi_{d'} )$ in $\feps^{p} \mathcal{A}_{d} (\Pi )[\ln \feps]$
	, which in turn depends on the coordinate of $\langle \Rmicro\rangle  (\phi_{d'} )$ in $\feps^{p+d} \mathcal{A}_{d} (\Pi )[\ln \feps]$. %
	Moreover, we have 
	\[\langle \Rmicro\rangle  (\phi_{d'} ) :=\fsum_{\n=0} ^{\infty } \bigg(\sum _{R\in \smash\Rmicro} \bigg)^{\!\n} (\phi_{d'} ) =\fsum_{\n=0} ^{\infty } \sum _{(R_{1} ,\dots ,R_{k} )\in \smash{\Rmicro^{\n} }} R_{1} \circ \cdots \circ R_{\n} (\phi_{d'} )\]
	with $R_{1} \circ \cdots \circ R_{\n} (\phi_{d'} ) \in \feps^{\Sigma _{i} \deg _{\feps} \!R_{i} } \mathcal{A}_{d'+\Sigma _{i} \deg _{\mathcal{A}} \!R_{i} }(\Pi )$ for any term in the above sum. Let $f:(x,y)\in \mathbb{R}^2 \mapsto  2x-y$. Since, we~have $f(\deg _{\mathcal{A}} R,\deg _{\feps} R )\geqslant 1$ for any $R\in \Rmicro$, there exists $\n_{0} \in \mathbb{N}$ that depends only on $(d,d',p)$ s.t.:
	\[\forall \n>\n_{0} ,\; \forall (R_{i} )\in \Rmicro^{\n} ,\qquad  f\bigg(d'+\sum _{i=1} ^{\n} \deg _{\mathcal{A}} R_{i} ,\sum _{i=1} ^{\n} \deg _{\feps} R_{i} \bigg) > f(d,p+d).\]
	Hence we have
	\[\langle \Rmicro\rangle  (\phi_{d'} ) :=\sum _{\n=0} ^{\red{\n_{0} }} \sum _{(R_{1} ,\dots ,R_{n} )\in \smash{\Rmicro^{\n} }} R_{1} \circ \cdots \circ R_{\n} (\phi_{d'} ).\]
	Now, the construction of the operators $\ROmegaPi$, $\RcouchePi$ and $\RGammaPi$ in \cref{operateurs R dans Omega,operateurs R} shows that applying any of these operators increases the degree in $\log (\alpha z)$ and $\ln r$ by at most one, so the same holds for the elements of $\Rmacro$. Moreover, the construction of $\Eps$ in \cref{def: Eps} shows that, for~any $\varphi \in \mathcal{A}(\Pi )$, the degree in $\ln \feps$ of $\Eps^{\pm 1} (\varphi )$ is less or equal to the degree in $\log (\alpha z)$ and $\ln r$ of $\varphi $. Thus, taking $n:=\n_{0}$, we have: $\forall \ell >n,\; \cmacro_{d,d',p,\ell } =0$.
\end{sousdemo}
The same property is true for the coefficients $\cmicro_{d,d',p,\ell }$. Finally the result follows by induction from \cref{def: construction des termes}.
\end{demo}

\begin{propositionsansboite}{}
\label{confimation de l'hyp de raccord}%
\Cref{ansatz sur Aeps} and the matching assumption of \cref{def: raccord algebrique} are satisfied.
\end{propositionsansboite}

\begin{demo}
\Cref{borne sur l} shows that the degree in $\ln \feps$ is finite for~any $p$. So, to prove \cref{ansatz sur Aeps} it suffices to check that, for~any $(p,\ell )$, $\bfuDA_{p,\ell } \in \fsum_{d\in \NN-p} \mathcal{A}_{d} (\Pi )$ and $\SDA_{p,\ell } \in \fsum_{d\in p-\NN} \mathcal{A}_{d} (\Pi )$. Let us show it for $\bfuDA_{p,\ell }$ only, as it~is similar for $\SDA_{p,\ell }$. We proceed by induction.\vspace{1mm}
\begin{sousdemo}
	Since all fields vanish for $p<0$, only the inductive step is non-trivial. 
	\cref{eq: bfuDApl} states that
	\[\bfuDA_{p,\ell } =\sum_{n=0}^\infty  (-k_{0} ^2 \,\ROmegaPi)^{n} \bigg(\underbrace{\RcouchePi\circ (\partial _{x|\couche} ^2 +k_{1} ^2 )(\bfuDA_{p-2,\ell } ) + {\textstyle \frac{\mu _{0} }{\mu _{1}} }\RGammaPi\circ \partial _{y|\Gamma ,y=0^{+}} (\bfuDA_{p-1,\ell } )}_{:=A} + \!\! \fsum_{d\in \TZ*} \!\!\! \sigma _{d} (\bfuDA_{p,\ell } ) \,\phi_{d} \bigg). \vspace{-1mm}\]
	Since $\deg  \ROmegaPi \in \mathbb{N}$, it suffices to show that the big brackets belong to $\fsum_{d\in \NN-p} \mathcal{A}_{d} (\Pi )$. This is true for $A$ using the induction hypothesis and $\deg \RcouchePi=0$, $\deg (\partial _{x|\couche} ^2 )=-2$ and $\deg (\RGammaPi\circ \partial _{y|\Gamma ,y=0^{+}} )=-1$ (see \cref{operateurs R,operateurs de derivation dans A}). It remains to show that $\forall d\in \TZ*\!,\; \sigma _{d} (\bfuDA_{p,\ell } ) \neq 0 \Rightarrow d\in \NN-p$. This last assertion holds because on one hand $\TN\subset \NN-p$, and on the other hand \cref{coincidence des sigma singuliers} and \cref{eq: definition des champs lointains} imply: $\forall d\in -\TN,\; \sigma _{d} (\bfuDA_{p,\ell } ) \neq 0 \Leftrightarrow  \sigma _{d} (u_{p,\ell } ) \neq 0\Rightarrow  \intvl{0,p+d} \neq  \varnothing  \Rightarrow  d\in \NN-p$.
\end{sousdemo}
Hence the assumptions of \cref{sec: raccordement algebrique}, which were \cref{ansatz: existence de bfuDApl et SDApl,ansatz sur Aeps}, are satisfied. So we can apply \cref{condition de raccord}. \Cref{eq: definition des champs lointains,eq: definition des champs de coin,coincidence des sigma singuliers} show that the matching conditions \cref{eq: condition de raccord 1,eq: condition de raccord 2} are satisfied. Thus \cref{condition de raccord} states that the matching assumption of \cref{def: raccord algebrique} is satisfied.
\end{demo}

\subsection{Practical way to build the far fields}
\label{sec: construction directe des upl}

\begin{textesansboite}
This section shows how to build directly the far fields without computing the layer and corner fields. Thanks to the explicit expression of the layer fields in \cref{eq: definition des champs de couche}, the layer is replaced by boundary conditions on $\Gamma $, that are identical to the case of an infinite layer. Moreover, the corner fields are replaced by corner conditions that depend on corner profiles (``correcteurs de coin'' in French).
\end{textesansboite}

\begin{definition}{corner profiles}
\label{def: correcteurs de coin}%
Let $d\in \TN$. We denote $(\mathcal{S}_{d,n} )_{n\in \mathbb{N}}$ the unique sequence of $\Vcoin+\khii \mathcal{A}(\Pi )$ s.t.\ for~any $n\in \mathbb{N}^*$:
\[\left\{\begin{array}{r@{\;=\;}l@{\;}l}
	\operatorname{div}(\mu  \nabla  \mathcal{S}_{d,0} ) & 0 &\text{ in }\Omega _{1} \\[0.5mm]
	\mathcal{S}_{d,0} & 0 &\text{ on }\partial \Omega _{1} \\[0.5mm]
	\sigma _{q} (\mathcal{S}_{d,0} )& \delta _{d,q} & \ \forall q\in \TN
\end{array}\right. \quad  \text{ and } \quad  \ \left\{\begin{array}{r@{\;=\;}l@{\;}l}
	\operatorname{div}(\mu  \nabla  \mathcal{S}_{d,n} ) & -\rho \,\mathcal{S}_{d,n-1} &\text{ in }\Omega _{1} \\[0.5mm]
	\mathcal{S}_{d,n} & 0 &\text{ on }\partial \Omega _{1} \\[0.5mm]
	\sigma _{q} (\mathcal{S}_{d,n} )& 0 & \ \forall q\in \TN
\end{array}\right.\]
And for~any $(d,n)$ we denote $\ccDA_{d,n}$ the element of $\Amicro(\Pi )$ s.t.\ $\forall d\in \mathbb{R},\; \mathcal{S}_{d,n} = T_{\geqslant d} (\ccDA_{d,n} )+o_{\partial } (\bfr^{d} )$.\\
These objects are well-defined thanks to \cref{cadre fonctionnel pour les champs de coin,DA des champs de coin}. The proof is the same as for \cref{def: construction des termes,validation construction et DA des champs}.
\end{definition}
\nomenclature[28b]{$\mathcal{S}_{d,n}$}{Corner profile}
\nomenclature[28b-60]{$\ccDA_{d,n}$}{Formal series representing the asymptotic behavior of $\mathcal{S}_{d,n}$ at infinity}

\begin{textesansboite}
Using \cref{eq: definition des champs de coin} and the uniqueness in \cref{cadre fonctionnel pour les champs de coin}, one can easily show by induction on $p$:
\begin{equation}{\label{eq: expression de Spl avec les correcteurs}}
\forall (p,\ell )\in \NN\times \mathbb{N},\quad  S_{p,\ell } = \sum_{n=0}^\infty  \;\sum _{d\in \TN \cap (p-2n-\NN)} \omega ^{2n} \sigma _{d} (S_{p-2n,\ell } )\,\mathcal{S}_{d,n} \end{equation}
Note that the sum above has a finite number of non-zero terms. By~\cref{unicité du DA}, the same holds replacing $S_{p,\ell }$ and $\mathcal{S}_{d,n}$ by resp.\ $\SDA_{p,\ell }$ and $\ccDA_{d,n}$.\\

We could have set $\operatorname{div}(\mu  \nabla  \mathcal{S}_{d,n} ) = -\omega ^2 \rho \,\mathcal{S}_{d,n-1}$ to remove $\omega ^{2n}$ from \cref{eq: expression de Spl avec les correcteurs}, but the chosen definition makes the corner profiles independent of $\fsource$ and $\omega $. Thus they can often be computed once in advance in numerical applications. In addition, corner profiles only take real values.
\end{textesansboite}

\nomenclature[12-30-30]{$\cfinal$}{Matching coefficients that give the $\sigma _{d} (\uDA_{p,\ell } )$ in terms of the $\sigma _{d'} (\uDA_{p',\ell '} )$ with $p'<p$}
\begin{textesansboite}
We recall that the matching conditions express:
\begin{itemize}
	\item the coefficients $\sigma _{d} (\bfuDA_{p,\ell } )$ with $d<0$ in terms of the $\sigma _{d'} (\SDA_{p',\ell '} )$ (with $p'<p$ and $d'$ positive or negative),
	\item and the $\sigma _{d} (\SDA_{p,\ell } )$ with $d>0$ in terms of the $\sigma _{d'} (\bfuDA_{p',\ell '} )$ (also with $p'<p$).
\end{itemize}
We will show in \cref{construction directe des upl} that, by composing the two formulas, we can express the $\sigma _{d} (\bfuDA_{p,\ell } )$ (with $d<0$) in terms of the coefficients $\sigma _{d'} (\bfu_{p',\ell '} )$ of the previous fields (i.e.\ with $p'<p$). When $d'<0$, $\sigma _{d'} (\SDA_{p',\ell '} )$ cannot be substitued by the matching condition, so we will use \cref{eq: expression de Spl avec les correcteurs} to express $\sigma _{d'} (\SDA_{p',\ell '} )$ in terms of the $\sigma _{d} (\SDA_{p,\ell } )$ with $d>0$ and the corner profiles. In the resulting formula, some coefficients denoted $\cfinal$ will appear. They are defined for~any $(d,d',p,\ell )\in (-\TN)\times \TZ*\times \NN\times \mathbb{N}$ by:
\begin{equation}{\label{eq: def cfinal}}
	\cfinal_{d,d'\!,p,\ell } := \! \sum _{\substack{(p_{1} ,p_{2} )\in (\mathbb{N}-d)\times (\mathbb{N}+d')\\ n\in \mathbb{N},\, p_{1} +p_{2} +2n=p}}
	\ \sum _{\substack{\red{d_{1} \in \TZ*}\\ p_{1} +d_{1} \in \mathbb{N}}} \ \sum _{\substack{d_{2} \in \TN\\ p_{2} -d_{2} \in \mathbb{N}}}
	\ \sum _{\substack{(\ell _{1} ,\ell _{2} )\in \mathbb{N}^2 \\ \ell _{1} +\ell _{2} =\ell }}
	\cmacro_{d,d_{1} ,p_{1} ,\ell _{1} }\cdot \omega ^{2n} \sigma _{d_{1}} (\mathcal{S} ^{\infty } _{d_{2} ,n} )\cdot \cmicro_{d_{2} ,d'\!,p_{2} ,\ell _{2} } .
\end{equation}
Separating the positive and negative values of $d_{1}$ in the above sum, and using that, when $d_{1} >0$, we have $\sigma _{d_{1}} (\ccDA_{d_{2} ,n} )=\sigma _{d_{1}} (\mathcal{S}_{d_{2} ,n} ) =\delta _{d_{1} ,d_{2} } \,\delta _{n,0}$, the formula can be written differently:
\begin{subequations}
\begin{align}
	\cfinal_{d,d'\!,p,\ell } &= \! \sum _{\substack{(p_{1} ,p_{2} )\in (\mathbb{N}-d)\times (\mathbb{N}+d')\\ n\in \mathbb{N},\, p_{1} +p_{2} +2n=p}}
	\ \sum _{\substack{\red{d_{1} \in -\TN}\\ p_{1} +d_{1} \in \mathbb{N}}} \ \sum _{\substack{d_{2} \in \TN\\ p_{2} -d_{2} \in \mathbb{N}}}
	\ \sum _{\substack{(\ell _{1} ,\ell _{2} )\in \mathbb{N}^2 \\ \ell _{1} +\ell _{2} =\ell }}
	\cmacro_{d,d_{1} ,p_{1} ,\ell _{1} }\cdot \omega ^{2n} \sigma _{d_{1}} (\mathcal{S} ^{\infty } _{d_{2} ,n} )\cdot \cmicro_{d_{2} ,d'\!,p_{2} ,\ell _{2} }\\
	&\quad  + \!\! \sum _{\substack{(p_{1} ,p_{2} )\in (\mathbb{N}-d)\times (\mathbb{N}+d')\\ p_{1} +p_{2} =p}}
	\ \sum _{\substack{\red{d_{1} \in \TN}\\ p_{1} +d_{1} \in \mathbb{N},\, p_{2} -d_{1} \in \mathbb{N}}}
	\ \sum _{\substack{(\ell _{1} ,\ell _{2} )\in \mathbb{N}^2 \\ \ell _{1} +\ell _{2} =\ell }}
	\cmacro_{d,d_{1} ,p_{1} ,\ell _{1} }\cdot \cmicro_{d_{1} ,d'\!,p_{2} ,\ell _{2} }
	\label{eq: cfinal 3}
\end{align}
\end{subequations}
Like in \cref{rq: conditions de raccord}, we have $d-d_{1} ,d_{2} -d' \in \mathbb{Z}\cap \TZ$ in \cref{eq: def cfinal}. Thus, if $d-d'\not\in \mathbb{Z}\cap \TZ$, line~\cref{eq: cfinal 3} vanishes.
\end{textesansboite}

\begin{remarque}
In~practice, the coefficients $\cfinal$ involved in the first far fields have formulas that are much more simple and contain only a few terms. See \cref{sec: tableaux exemples DA} for examples. In addition, when $\pst$ is irrational, \cref{eq: def cfinal} has at most one non-zero term, so the coefficients $\cfinal$ have a very simple expression (see \cref{valeur de cmaicro cas irrationnel}).
\end{remarque}

\begin{lemmesansboite}{}
\label{coefs cfinal non nuls}%
For~any $(d,d',p,\ell )\in (-\TN)\times \TZ*\times \NN\times \mathbb{N}$, we have: $\cfinal_{d,d',p,\ell } \neq 0 \Rightarrow  p\in \NN+2\pst$ and $p+d-d'\in \mathbb{N}$.
\end{lemmesansboite}

\begin{demo}
If $\cfinal_{d,d',p,\ell } \neq 0$, then the sum \cref{eq: def cfinal} contains at least one term, of indexes denoted $(p_{1,} p_{2} ,d_{1} ,d_{2} ,\ell _{1} ,\ell _{2} )$. The conditions on the sum indexes imply
\begin{itemize}
	\item $p=(p_{1} +d)+2n+(p_{2} -d_{2} )-d+d_{2} \in \mathbb{N}+\TN+\TN \subset \NN+ 2\pst $
	\item and $p+d-d'=(p_{1} +d)+2n+(p_{2} -d')\in \mathbb{N}$.
\end{itemize}
\end{demo}

\nomenclature[30-00]{$\uDA_{p,\ell }$}{Formal series representing the asymptotic behavior of $u_{p,\ell }$ at the corner}
\begin{textesansboite}
For~any  $(p,\ell )$, we denote:
\[\uDA _{p,\ell } := \bfuDA_{p,\ell |\Omega } \in \Amacro(\Omega ).\]
It satisfies: $\forall d\in \mathbb{R},\; u_{p,\ell } \smash{\underset{r\rightarrow 0}{=} } T_{\leqslant d} (\uDA_{p,\ell } )+o_{\partial } (r^{d} )$.
\end{textesansboite}

\begin{theoreme}{direct construction of the far fields}
\label{construction directe des upl}%
Let $(T_{n} )_{n\in \mathbb{N}}$ be the sequence of Taylor coefficients of the tangent: $\forall t\in (- \frac{\pi }{2} , \frac{\pi }{2} ),\; \tan  t = \sum_{n=0}^\infty  T_{n} t^{2n+1}$.\\
$(u_{p,\ell } )_{p\in \NN,\ell \in \mathbb{N}}$ is the unique family of $H^{1} (\Omega )+\khiz \mathcal{A}(\Omega )$ s.t.\ for~any $(p,\ell )\in \NN\times \mathbb{N}$:
\[\left\{\begin{array}{r@{\;=\;}ll@{}}
	\mu _{0} \Delta  u_{p,\ell } +\omega ^2 \rho _{0} u_{p,\ell } & \fsource \,\delta _{p,0} \,\delta _{\ell ,0} &\text{ in }\Omega  \\[0.5mm]
	u_{p,\ell } & 0 &\text{ on }\bordgauche\\[0.5mm]
	u_{p,\ell } &\displaystyle \sum_{n=0}^\infty  \frac{\mu _{0} }{\mu _{1}} \,T_{n} \cdot  (\partial _{x} ^2 +k_{1} ^2 )^{n} \partial _{y} u_{p-1-2n,\ell } &\text{ on }\Gamma  \\[-1.5mm]
	\sigma _{d} (u_{p,\ell } )&\displaystyle \!\! \sum _{\substack{p'\in \NN\\ p-p'\in \NN+ \frac{2\pi }{\Theta } }} \; \sum _{\substack{d'\in \TZ*\\ p'+d'\in \intvl{0,p+d}}} \! \sum _{\ell '=0} ^{\ell } \cfinal_{d,d',p-p',\ell -\ell '} \cdot \sigma _{d'} (\uDA_{p',\ell '} ) & \ \forall d\in -\TN
\end{array}\right.\]\fbeq
\end{theoreme}
\nomenclature[29-23]{$T_{n}$}{Taylor coefficient of the tangent}

\begin{remarque}
In these relations, the far fields depend on the shape of $\Omega _{1}$ at the corner through the corner profiles coefficients $\sigma _{q} (\ccDA_{d,n} )$ (with $q<0$) that appear in the definition of $\cfinal$.
\end{remarque}

\newcommand{\bull}{\text{\tiny $\bullet$},\text{\tiny $\bullet$}}
\begin{demo}
Uniqueness follows from the uniqueness in \cref{cadre fonctionnel des champs lointains}. So it suffices to prove that the far fields satisfy the equations above.\medskip

\underline{Boundary condition:} By~\cref{eq: definition des champs de couche} and \cref{eq: definition des champs lointains}, we have
\[u_{p,\ell } (x,0) =U_{p,\ell } (x,0) = \sum_{n=0}^\infty  (\partial _{x} ^2 +k_{1} ^2 )^{n} \partial _{y} u_{p-1-2n,\ell } (x,0) \cdot  \mathcal{U}_{n} (0).\]
So we need to calculate $\mathcal{U}_{n} (0)$. Let $\mathrm{U} :(Y,t) \mapsto  \sum_{n=0}^\infty  \mathcal{U}_{n} (Y)\, t^{2n+1}$. The definition of $(\mathcal{U}_{n} )$ (reminded in \cref{eq: construction directe des upl eq1} below) formally implies a differential equation on $\mathrm{U}$ given in \cref{eq: construction directe des upl eq2}.\par
\vspace{2.5mm}
\begin{tabular}{@{}c@{}c@{}}
\begin{minipage}{0.5\linewidth}
\begin{equation}{\label{eq: construction directe des upl eq1}}
\forall n\in \mathbb{N},\ \left\{\begin{array}{r@{\;=\;}ll}
	\mathcal{U}_{n} '' & -\mathcal{U}_{n-1} \\[0.5mm]
	\mathcal{U}_{n} ' (0) & \frac{\mu _{0} }{\mu _{1}} \delta _{n,0} \\[0.5mm]
	\mathcal{U}_{n} (-1) & 0 
\end{array}\right.\end{equation}
\end{minipage}
&
\begin{minipage}{0.5\linewidth}
\begin{equation}{\label{eq: construction directe des upl eq2}}
\left\{\begin{array}{r@{\;=\;}ll}
\partial _{Y} ^2 \text{U}  & -t^2 \text{U}  \\[0.5mm]
\partial _{Y} \text{U}  _{|Y=0} & \mathcal{U}_{0} '(0)\cdot t = \frac{\mu _{0} }{\mu _{1}} t \\[0.5mm]
\text{U}  _{|Y=-1} & 0 
\end{array}\right.\end{equation}
\end{minipage}
\end{tabular}\par
\vspace{2.5mm}
Thus $\mathrm{U} (Y,t) = \frac{\mu _{0} }{\mu _{1}} \frac{\sin (t(Y+1))}{\cos (t)} =:\varphi (Y,t)$. This is formal, as we do~not know whether the series $\mathrm{U}$ converges. However, there exists a sequence of polynomial functions $(\Phi _{n} )$ s.t., for~any $(Y,t)\in (-1,0)\times (- \frac{\pi }{2} , \frac{\pi }{2} )$, $\varphi (Y,t) = \sum_{n=0}^\infty  \Phi _{n} (Y)\, t^{2n+1}$ (because $\varphi $ is odd w.r.t.\ $t$). Since $\varphi $ satisfies \cref{eq: construction directe des upl eq2}, $(\Phi _{n} )$ is a solution of \cref{eq: construction directe des upl eq1}. But this solution is unique, so $(\Phi _{n} )=(\mathcal{U}_{n} )$ and $\mathrm{U} =\varphi $. Hence $\forall t\in (- \frac{\pi }{2} , \frac{\pi }{2} ),\; \sum_{n=0}^\infty  \mathcal{U}_{n} (0)\, t^{2n+1} = \frac{\mu _{0} }{\mu _{1}} \tan (t)$. That is to say $\forall n\in \mathbb{N},\; \mathcal{U}_{n} (0) = \frac{\mu _{0} }{\mu _{1}} \,T_{n}$.\medskip\smallskip

\underline{Corner condition:} For~any $(d,p,\ell )\in \TN\times \mathbb{R}\times \mathbb{N}$, let $\widetilde{\mathcal{S}} ^{\infty } _{d,p,\ell } := \omega ^{p} \ccDA_{d,p/2}$ if $\frac{p}{2} \in \mathbb{N}$ and $\ell =0$ and $\widetilde{\mathcal{S}} ^{\infty } _{d,p,\ell } :=0$ otherwise. For~any family $(x_{p,\ell } )$, we denote $x_{\bull} :=(x_{p,\ell } )$ the family itself. We also denote $*$ the convolution product w.r.t.\ $(p,\ell )$. Let $d\in -\TN$. We have
\begin{align*}
	\sigma _{d} (u_{\bull} ) &= \sum _{d_{1} \in \TZ*} \cmacro_{d,d_{1} ,\bull} * \sigma _{d_{1}} (\SDA_{\bull} ) 
	\tag*{by \cref{eq: definition des champs lointains}}\\
	&= \sum _{d_{1} \in \TZ*} \cmacro_{d,d_{1} ,\bull} * \sum _{d_{2} \in \TN} \sigma _{d_{1}} ( \widetilde{\mathcal{S}} ^{\infty } _{d_{2} ,\bull} )*\sigma _{d_{2}} (S_{\bull} ) \tag*{by \cref{eq: expression de Spl avec les correcteurs}} \\
	&= \sum _{d_{1} \in \TZ*} \cmacro_{d,d_{1} ,\bull} * \sum _{d_{2} \in \TN} \sigma _{d_{1}} ( \widetilde{\mathcal{S}} ^{\infty } _{d_{2} ,\bull} )*\sum _{d'\in \TZ*} \cmicro_{d_{2} ,d'\!,\bull} * \sigma _{d'} (\bfuDA_{\bull} )
	\tag*{by \cref{eq: definition des champs de coin}}\\
	&= \sum _{d'\in \TZ*} \bigg(\sum _{d_{1} \in \TZ*} \sum _{d_{2} \in \TN} \cmacro_{d,d_{1} ,\bull} * \sigma _{d_{1}} ( \widetilde{\mathcal{S}}_{d_{2} ,\bull} )* \cmicro_{d_{2} ,d'\!,\bull} \bigg)* \sigma _{d'} (\uDA_{\bull} )\\
	&= \sum _{d'\in \TZ*} \cfinalbis_{d,d'\!,\bull} * \sigma _{d'} (\uDA_{\bull} )
\end{align*}
with $\displaystyle \cfinalbis_{d,d',p,\ell } := \sum _{\substack{(p_{1} ,p_{2} ),n\in \mathbb{N}\\ p_{1} +p_{2} +2n=p}} \sum _{d_{1} \in \TZ*} \sum _{d_{2} \in \TN} \sum _{\substack{(\ell _{1} ,\ell _{2} )\in \mathbb{N}^2 \\ \ell _{1} +\ell _{2} =\ell }} \cmacro_{d,d_{1} ,p_{1} ,\ell _{1} }\cdot \omega ^{2n} \sigma _{d_{1}} (\mathcal{S} ^{\infty } _{d_{2} ,n} )\cdot \cmicro_{d_{2} ,d'\!,p_{2} ,\ell _{2} }$.\\[1mm]
Moreover, for~any $R\in \Rmaicro$, we have $\deg _{\varepsilon } \!R \in \mathbb{N}$ and $\deg _{\varepsilon } \!R \pm  \deg _{\mathcal{A}} \!R\in \mathbb{N}$ (see \cref{fig: degres des elements de Rmacro et Rmicro}). So by~definition of the coefficients $\cmacro_{\dots }$ and $\cmicro_{\dots }$:
\begin{itemize}
	\item $\cmacro_{d,d_{1} ,p_{1} ,\ell _{1} } \neq 0 \Rightarrow  p_{1} +d\in \mathbb{N} \text{ and } p_{1} +d_{1} \in \mathbb{N}$,
	\item $\cmicro_{d_{2} ,d',p_{2} ,\ell _{2} } \neq 0 \Rightarrow  p_{2} -d_{2} \in \mathbb{N} \text{ and } p_{2} -d'\in \mathbb{N}$.
\end{itemize}
This implies that $\cfinalbis_{d,d',p,\ell } =\cfinal_{d,d',p,\ell }$ for~any $(d,d',p,\ell )$. Thus by \cref{coefs cfinal non nuls}, if $\cfinal_{d,d',p,\ell } \neq 0$ then $p\in \NN+ 2\pst$ and $p+d-d'\in \mathbb{N}$. This and the property $(\sigma _{d'} (u_{p',\ell '} ) \neq 0 \Rightarrow  p'\in \NN \text{ and } p'+d'\in \NN)$ explain the sum indexes of the formula given for $\sigma _{d} (u_{p,\ell } )$ in \cref{construction directe des upl}.
\end{demo}

\begin{exemplenumerote}
\label{ex: eq u00 et u10}%
Using \cref{construction directe des upl}, one can check that, if $\pst \not\in \mathbb{N}$, then $u_{\pi /\Theta ,\ell }$ vanishes for~any $\ell $. In addition $\frac{2\pi }{\Theta } > 1$, so the first non-zero far fields are $u_{0,0}$ and $u_{1,0}$ and they satisfy
\[\left\{\begin{array}{@{\;}r@{\;=\;}l@{\;\;}l@{}}
	\mu _{0} \Delta  u_{0,0} +\omega ^2 \rho _{0} u_{0,0} & \fsource &\text{ in }\Omega  \\[0.5mm]
	u_{0,0} & 0 &\text{ on }\bordgauche \\[0.5mm]
	u_{0,0} & 0 &\text{ on }\Gamma  \\[0.5mm]
	\sigma _{d} (u_{0,0} )& 0 & \ \forall d\in -\TN
\end{array}\right. \qquad  \text{ and } \quad  \left\{\begin{array}{@{\;}r@{\;=\;}l@{\;\;}l@{}}
	\mu _{0} \Delta  u_{1,0} +\omega ^2 \rho _{0} u_{1,0} & 0 &\text{ in }\Omega  \\[0.5mm]
	u_{1,0} & 0 &\text{ on }\bordgauche\\[0.5mm]
	u_{1,0} & \frac{\mu _{0} }{\mu _{1}} \partial _{y} u_{0,0} &\text{ on }\Gamma  \\[0.5mm]
	\sigma _{d} (u_{1,0} )& 0 & \ \forall d\in -\TN
\end{array}\right.\]
Much more detailed examples are given in \cref{sec: simplifications DA,sec: exemples de structures de DA}.
\end{exemplenumerote}

\begin{textesansboite}
The next result is a small improvement of \cref{ansatz sur Aeps} that  will be useful for effective models (see \cref{chap: ME}). We recall that \cref{ansatz sur Aeps} implies $\bfuDA_{p,\ell } \in  \fsum_{d\in \NN-p} \mathcal{A}_{d} (\Pi )$ for~any $p,\ell $. A similar improvement does~not exist for the corner profiles.
\end{textesansboite}

\begin{corollairesansboite}{}
\label{precision de l'ansatz sur uDA}%
For~any $p\in \NN$ and $\ell \in \mathbb{N}$, we have $\bfuDA_{p,\ell } \in  \fsum_{d\in \NN-p\red{+\pst}} \mathcal{A}_{d} (\Pi )$.
\end{corollairesansboite}

\begin{demo}
Let us prove it by induction on $p$.
\begin{itemize}
	\item First, let $p=0$. For~any $\ell \in \mathbb{N}$ and $d\in -\TN$, $\sigma _{d} (\bfuDA_{p,\ell } )$ is equal to $\sigma _{d} (u_{p,\ell } )$ by \cref{coincidence des sigma singuliers}, and $\sigma _{d} (u_{p,\ell } )$ vanishes by \cref{construction directe des upl}. Thus \cref{validation construction et DA des champs} shows
	\[\bfuDA_{0,\ell } =\langle -k_{0} ^2 \ROmegaPi\rangle \bigg(\fsum_{d\in \TN} \sigma _{d} (\bfuDA_{0,\ell } ) \,\phi_{d} \bigg).\]
	Hence we have $\bfuDA_{0,\ell } \in  \fsum_{q\in \mathbb{N}+\TN} \mathcal{A}_{q} (\Pi ) = \fsum_{q\in \NN+\pst} \mathcal{A}_{q} (\Pi )$.
	\item Now we assume $p>0$. Let $\ell \in \mathbb{N}$ and $d\in -\TN$ be s.t.\ $\sigma _{d} (\bfu_{p,\ell } ) \neq 0$. Then the expression of $\sigma _{d} (u_{p,\ell } )$ given by \cref{construction directe des upl} has at least one non-zero term $\cfinal_{d,d',p-p',\ell -\ell '} \cdot \sigma _{d'} (\uDA_{p',\ell '} )$. This term satisfies $p+d-(p'+d')\in \mathbb{N}$ by the conditions on the sum indexes, and $d'\in \NN-p'+\pst$ by induction hypothesis. Hence we have $d+p=[p+d-(p'+d')] + (p'+d') \in \mathbb{N} + (\NN+\pst) \subset  \NN+\pst$. Thus, by \cref{validation construction et DA des champs}, we have
	\[\bfuDA_{p,\ell } = \big\langle  {-k_{0} ^2 \,\ROmegaPi} \big\rangle  \bigg(\displaystyle \underbrace{ {-\RcouchePi}\circ (\partial _{x|\couche} ^2 +k_{1} ^2 )}_{=:R_{2} } (\bfuDA_{p-2,\ell } ) + \underbrace{\frac{\mu _{0} }{\mu _{1}} \RGammaPi\circ \partial _{y|\Gamma ,y=0^{+}} }_{=:R_{1} } (\bfuDA_{p-1,\ell } )+ \fsum_{\substack{d\in \TZ*\\ \mathclap{d\in \NN-p+\pst}}} \sigma _{d} (\bfuDA_{p,\ell } ) \,\phi_{d} \bigg).\]
	For~any $i\in \{1,2\}$, $R_{i}$ has degree $i$, so $R_{i} (\bfuDA_{p-i,\ell } )\in  \fsum_{q\in \NN-i+\pst+i} \mathcal{A}_{q} (\Pi )$ by induction hypotheses. Thus $\bfuDA_{p,\ell } \in \fsum_{q\in \NN-p+\pst} \mathcal{A}_{q} (\Pi )$.
\end{itemize}
\end{demo}

\section{Error estimates}
\label{sec: estimations d'erreur}

\begin{textesansboite}
Let $\chi  \in \mathcal{C}^{\infty } (\mathbb{R}^2 )$ be equal to 1 on $B(0, \frac{1}{2} )$ and 0 outside $B(0,1)$, and, for~any $\eta >0$, $\chi _{\eta } : x \mapsto  \chi ( \frac{x}{\eta } )$. We denote $\couche_{\varepsilon } := (-\varepsilon ,1)\times \mathbb{R}_+^*$ and $\Pi _{\varepsilon } :=\Omega \sqcup \Gamma \sqcup \couche_{\varepsilon }$ (defined similarly as $\Pi $ at page~\pageref{page: definition de Pi}). We define on $\Pi _{\varepsilon }$ the following variant of the far-and-layer fields $\bfu_{p,\ell }$ \[\bfue_{p,\ell } (x,y) := \left\{\begin{array}{ll}
	u_{p,\ell } (x,y) &\text{ in }\Omega  \\[0.5mm]
	U_{p,\ell } (x, \frac{y}{\varepsilon } )&\text{ in }\couche_{\varepsilon }
\end{array}\right.\]
We also denote, for~any $p\in \NN$, $n_{p} :=\max \{\ell \in \mathbb{N} \mid u_{p,\ell } \neq 0 \text{ or } U_{p,\ell } \neq 0 \text{ or } S_{p,\ell } \neq 0\}$.\\

Let us define the approximate global field at order $\ptronc\in \mathbb{R}_+$ as follows for~any $\varepsilon $ small enough: $\forall (x,y)\in \Omega _{\varepsilon } ,$
\begin{equation}{\label{eq: approx SE sur Omega eps}}
	u_{\varepsilon ,\ptronc} (x,y) := (1-\chi _{\eta } (x,y)) \! \sum _{\substack{p\in \NN \cap [0,\ptronc]\\ \ell \in [\![0,n_{p} ]\!]}} \varepsilon ^{p} \ln ^{\ell } \!\varepsilon \, \bfue_{p,\ell } (x,y) \;+\; \chi _{\eta } (x,y) \! \sum _{\substack{p\in \NN \cap [0,\ptronc]\\ \ell \in [\![0,n_{p} ]\!]}} \varepsilon ^{p} \ln ^{\ell } \!\varepsilon \, S_{p,\ell } ({\textstyle \frac{x}{\varepsilon } , \frac{y}{\varepsilon } })
\end{equation}
where $\eta :=\sqrt{\varepsilon }$. Note that $(1-\chi _{\eta } (x,y))\bfue_{p,\ell } (x,y)$ is well-defined on $\Omega _{\varepsilon }$ when $\eta /2>\varepsilon \Rcoin$ since $\Pi _{\varepsilon } \setminus B(0,\varepsilon \Rcoin) = \Omega _{\varepsilon } \setminus B(0,\varepsilon \Rcoin)$. We will see as a consequence of \cref{estimation d'erreur globale} that $u_{\varepsilon ,\ptronc} \in H^{1} (\Omega _{\varepsilon } )$.
\end{textesansboite}

\subsection{Estimation of the matching error}

\begin{textesansboite}
The matching zone is $\Omega _{\varepsilon } \cap \crown$ where $\crown$ is the annulus $\crown :=B(0,\eta ) \setminus B(0,\eta /2)$. Letting $\eta =\sqrt{\varepsilon }$ makes the matching zone tend to 0 w.r.t.\ the far fields (because $\eta  \rightarrow 0$ when $\varepsilon \rightarrow 0$) and to infinity w.r.t.\ the corner fields (because $\frac{\eta }{\varepsilon } \rightarrow \infty $). Thanks to the matching assumption, we can state the following first error estimate concerning the error in the matching zone. We use the symbol $\lesssim $ for majorations valid up to a constant independent of $\varepsilon $.
\end{textesansboite}

\begin{lemme}{matching error}
\label{erreur de raccord}%
Let $\ptronc \in \mathbb{R}_+$. For $\varepsilon $ small enough, we have
\[\bigg\| \sum _{p\in \NN\cap [0,\ptronc]} \sum _{\ell \in \mathbb{N}} \varepsilon ^{p} \ln ^{\ell } \!\varepsilon \,\big(\bfue_{p,\ell } (x,y)-S_{p,\ell } ({\textstyle \frac{x}{\varepsilon } , \frac{y}{\varepsilon } })\big)\bigg\|_{H^{1} (\Omega _{\varepsilon } \cap \crown)} \lesssim  \varepsilon ^{ \frac{\ptronc}{2} -1}.\]\fbeq[1.3]
\end{lemme}

\begin{demo}
We will compare $\bfue_{p,\ell }$ and $S_{p,\ell }$ in $\crown$ to their asymptotic expansions at 0, resp $\infty $. Let us denote:
\begin{itemize}
	\item $\bfu_{\varepsilon ,\ptronc} := \sum _{p\in \NN\cap [0,\ptronc]} \sum _{\ell \in \mathbb{N}} \varepsilon ^{p} \ln ^{\ell } \!\varepsilon \, \bfue_{p,\ell } $ in $\Pi _{\varepsilon }$
	\item $S_{\varepsilon ,\ptronc} := \sum _{p\in \NN\cap [0,\ptronc]} \sum _{\ell \in \mathbb{N}} \varepsilon ^{p} \ln ^{\ell } \!\varepsilon \, S_{p,\ell } ({\textstyle \frac{x}{\varepsilon } , \frac{y}{\varepsilon } })$ in $\Omega _{\varepsilon }$
	\item $\bfuDA_{\varepsilon ,\ptronc} := \sum _{p\in \NN\cap [0,\ptronc]} \sum _{\ell \in \mathbb{N}} \varepsilon ^{p} \ln ^{\ell } \!\varepsilon \cdot \left\{\begin{array}{@{}l@{}l}
		T_{\leqslant \ptronc-p} (\bfuDA_{p,\ell } )(x,y) &\text{ in }\Omega  \\[0.5mm]
		T_{\leqslant \ptronc-p} (\bfuDA_{p,\ell } )(x,{\textstyle \frac{y}{\varepsilon } }) &\text{ in }\couche_{\varepsilon } 
	\end{array}\right.$
	\item $\SDA_{\varepsilon ,\ptronc} := \sum _{p\in \NN\cap [0,\ptronc]} \sum _{\ell \in \mathbb{N}} \varepsilon ^{p} \ln ^{\ell } \!\varepsilon \cdot T_{\geqslant p-\ptronc} (\SDA_{p,\ell } )({\textstyle \frac{x}{\varepsilon } , \frac{y}{\varepsilon } })$ in $\Pi _{\varepsilon }$.
\end{itemize}
(see \cref{def: troncature de SF} for $T_{\bullet}$). We split the estimate into three parts that we will bound separately:
\[\|\bfu_{\varepsilon ,\ptronc} -S_{\varepsilon ,\ptronc} \|_{H^{1} (\Omega _{\varepsilon } \cap \crown)} \leqslant  \|\bfu_{\varepsilon ,\ptronc} -\bfuDA_{\varepsilon ,\ptronc} \|_{H^{1} (\Omega _{\varepsilon } \cap \crown)} + \|\bfuDA_{\varepsilon ,\ptronc} -\SDA_{\varepsilon ,\ptronc} \|_{H^{1} (\Omega _{\varepsilon } \cap \crown)} + \|\SDA_{\varepsilon ,\ptronc} -S_{\varepsilon ,\ptronc} \|_{H^{1} (\Omega _{\varepsilon } \cap \crown)} .\]
In addition, we will split some of the norms $\|{\dots }\|_{H^{1} (\Omega _{\varepsilon } \cap \crown)}$ into $\|{\dots }\|_{H^{1} (\Omega \cap \crown)} +\|{\dots }\|_{H^{1} (\couche_{\varepsilon } \cap \crown)}$.
\begin{enumerate}
	\item \underline{$\|\bfu_{\varepsilon ,\ptronc} -\bfuDA_{\varepsilon ,\ptronc} \|_{H^{1} (\Omega \cap \crown)}$ :} By~\cref{validation construction et DA des champs}, for~any $(p,\ell )$ we have $u_{p,\ell } -T_{\leqslant \ptronc-p} (\bfuDA _{p,\ell } )=o_{\partial } (r^{\ptronc-p} )$ in $\Omega $ when $r\rightarrow 0$. By~\cref{def: petit o derivable} of $o_{\partial }$, it implies $u_{p,\ell } -T_{\leqslant \ptronc-p} (\bfuDA_{p,\ell } )=O(r^{\ptronc-p} )=O(r^{\ptronc-p-1} )$ and $\nabla \big[u_{p,\ell } -T_{\leqslant \ptronc-p} (\bfuDA_{p,\ell } )\big]=O(r^{\ptronc-p-1} )$ uniformly in $\theta $. Thus:\vspace{-2mm} \begin{align*}
		\|\bfu_{\varepsilon ,\ptronc} -\bfuDA_{\varepsilon ,\ptronc} \|_{H^{1} (\Omega \cap \crown)}
		&\leqslant  \sum _{p\in \NN\cap [0,\ptronc]} \sum _{\ell =0} ^{n_{p}}  \varepsilon ^{p} \ln ^{\ell } \!\varepsilon \cdot \|u_{p,\ell } -T_{\leqslant \ptronc-p} (\bfuDA_{p,\ell } )\|_{H^{1} (\Omega \cap \crown)} \\
		&\lesssim  \sum _{p\in \NN\cap [0,\ptronc]} \sum _{\ell =0} ^{n_{p}}  \varepsilon ^{p} \ln ^{\ell } \!\varepsilon \cdot \left( \int _{B(0,\eta )} r^{2(\ptronc-p-1)} \,\mathrm{d}r\,r\mathrm{d}\theta  \right)^{1/2} \\
		&\lesssim  \sum _{p\in \NN\cap [0,\ptronc]} \sum _{\ell =0} ^{n_{p}}  \varepsilon ^{p} \ln ^{\ell } \!\varepsilon  \cdot \eta ^{\ptronc-p} \\
		&\lesssim \; \eta ^{\ptronc-2} \qquad  \quad  \text{because }  \varepsilon ^{p} \lesssim  \eta ^{p} \text{ and } \ln ^{\ell } \!\varepsilon \lesssim \eta ^{-2}
	\end{align*}
	
	\item \underline{$\|\bfu_{\varepsilon ,\ptronc} -\bfuDA_{\varepsilon ,\ptronc} \|_{H^{1} (\couche_{\varepsilon } \cap \crown)}$ :} Similarly, $U_{p,\ell } -T_{\leqslant \ptronc-p} (\bfuDA_{p,\ell } )=o_{\partial } (x^{\ptronc-p} )$ in $\couche$ when $x\rightarrow \infty $. Thus
	\newlength{\widthGauche}%
	\settowidth{\widthGauche}{$\partial _{x} \big[U_{p,\ell } -T_{\leqslant \ptronc-p} (\bfuDA_{p,\ell } )\big]=O(x^{\ptronc-p-1} )$\,}%
	\begin{itemize}
		\item \parbox{\widthGauche}{$U_{p,\ell } -T_{\leqslant \ptronc-p} (\bfuDA_{p,\ell } )=O(x^{\ptronc-p} )$} so $U_{p,\ell } (x, \frac{y}{\varepsilon } ) -T_{\leqslant \ptronc-p} (\bfuDA_{p,\ell } )(x, \frac{y}{\varepsilon } )=O(x^{\ptronc-p} )$
		\item \parbox{\widthGauche}{$\partial _{x} \big[U_{p,\ell } -T_{\leqslant \ptronc-p} (\bfuDA_{p,\ell } )\big]=O(x^{\ptronc-p-1} )$} so $\partial _{x} \big[U_{p,\ell } (x, \frac{y}{\varepsilon } ) -T_{\leqslant \ptronc-p} (\bfuDA_{p,\ell } )(x, \frac{y}{\varepsilon } )\big] =O(\varepsilon ^{-1} x^{\ptronc-p-1} )$
		\item \parbox{\widthGauche}{$\partial _{Y} \big[U_{p,\ell } -T_{\leqslant \ptronc-p} (\bfuDA_{p,\ell } )\big]=O(x^{\ptronc-p} )$} so $\partial _{Y} \big[U_{p,\ell } (x, \frac{y}{\varepsilon } ) -T_{\leqslant \ptronc-p} (\bfuDA_{p,\ell } )(x, \frac{y}{\varepsilon } )\big]=O(\varepsilon ^{-1} x^{\ptronc-p} )$ 
	\end{itemize}
	which are all $O(\varepsilon ^{-1} x^{\ptronc-p-1} )$ (uniformly in $Y$). Hence:\vspace{-2mm} \begin{align*}
		\|\bfu_{\varepsilon ,\ptronc} -\bfuDA_{\varepsilon ,\ptronc} \|_{H^{1} (\couche_{\varepsilon } \cap \crown)}
		&\lesssim  \sum _{p\in \NN\cap [0,\ptronc]} \sum _{\ell =0} ^{n_{p}}  \varepsilon ^{p} \ln ^{\ell } \!\varepsilon \cdot \varepsilon ^{-1} \,\left( \int _{[0,\eta ]\times [-\varepsilon ,0]} x^{2(\ptronc-p-1)} \,\mathrm{d}x\,\mathrm{d}y \right)^{1/2} \\
		&\lesssim  \sum _{p\in \NN\cap [0,\ptronc]} \sum _{\ell =0} ^{n_{p}}  \varepsilon ^{p} \ln ^{\ell } \!\varepsilon \cdot \varepsilon ^{-1} \cdot \varepsilon ^{ \frac{1}{2} } \eta ^{\ptronc-p- \frac{1}{2} } \\
		&\lesssim \; \eta ^{\ptronc-2} \qquad  \quad  \text{because }  \varepsilon =\eta ^2  \text{ and } \ln ^{\ell } \!\varepsilon \lesssim \eta ^{- \frac{1}{2} }
	\end{align*}
	
	\item \underline{$\|S_{\varepsilon ,\ptronc} -\SDA_{\varepsilon ,\ptronc} \|_{H^{1} (\Omega \cap \crown)}$ :} Similarly $S_{p,\ell } -T_{\geqslant p-\ptronc} (\SDA _{p,\ell } )=o_{\partial } (r^{p-\ptronc} )$ in $\Omega $ when $r\rightarrow \infty $, so $S_{p,\ell } -T_{\geqslant p-\ptronc} (\SDA_{p,\ell } )=O(r^{p-\ptronc} ) \text{ and } \nabla \!_{(X,Y)}\big[S_{p,\ell } -T_{\geqslant p-\ptronc} (\SDA_{p,\ell } )\big]=O(r^{p-\ptronc-1} )=O(r^{p-\ptronc} )$ uniformly in $\theta $. Since $\nabla \!_{(x,y)} =\varepsilon ^{-1} \nabla \!_{(X,Y)}$, we deduce: \vspace{-2mm} \begin{align*}
		\|S_{\varepsilon ,\ptronc} -\SDA_{\varepsilon ,\ptronc} \|_{H^{1} (\Omega \cap \crown)}
		&\leqslant  \sum _{p\in \NN\cap [0,\ptronc]} \sum _{\ell =0} ^{n_{p}}  \varepsilon ^{p} \ln ^{\ell } \!\varepsilon \cdot \Big\|\big[S_{p,\ell } -T_{\geqslant p-\ptronc} (\SDA_{p,\ell } )\big]({\textstyle \frac{x}{\varepsilon } , \frac{y}{\varepsilon } })\Big\|_{H^{1} (\Omega \cap \crown)} \\
		&\lesssim  \sum _{p\in \NN\cap [0,\ptronc]} \sum _{\ell =0} ^{n_{p}}  \varepsilon ^{p} \ln ^{\ell } \!\varepsilon \cdot  \varepsilon ^{-1} \left( \int _{\mathbb{R}^2  \setminus B(0,\eta /2)} \left( \frac{r}{\varepsilon } \right)^{2(p-\ptronc)} \,\mathrm{d}r\,r\mathrm{d}\theta  \right)^{1/2} \\
		&\lesssim  \sum _{p\in \NN\cap [0,\ptronc]} \sum _{\ell =0} ^{n_{p}}  \varepsilon ^{p} \ln ^{\ell } \!\varepsilon  \cdot \varepsilon ^{\ptronc-p-1} \cdot \eta ^{p-\ptronc+1} \\
		&\lesssim \; \eta ^{\ptronc-2} \qquad  \quad  \text{because }  \varepsilon =\eta ^2  \text{ and } \ln ^{\ell } \!\varepsilon \lesssim \eta ^{-1}
	\end{align*}
	
	\item \underline{$\|S_{\varepsilon ,\ptronc} -\SDA_{\varepsilon ,\ptronc} \|_{H^{1} (\couche_{\varepsilon } \cap \crown)}$ :} $S_{p,\ell } -T_{\geqslant p-\ptronc} (\SDA _{p,\ell } )=o_{\partial } (X^{p-\ptronc} )$ in $\couche$ when $x\rightarrow 0$, so:
	\begin{itemize}
		\item $S_{p,\ell } -T_{\geqslant p-\ptronc} (\bfuDA_{p,\ell } )=O(X^{p-\ptronc} )$ 
		\item $\partial _{X} \big[S_{p,\ell } -T_{\geqslant p-\ptronc} (\SDA_{p,\ell } )\big]=O(X^{p-\ptronc-1} )=O(X^{p-\ptronc} )$
		\item $\partial _{Y} \big[S_{p,\ell } -T_{\geqslant p-\ptronc} (\SDA_{p,\ell } )\big]=O(X^{p-\ptronc} )$ 
	\end{itemize}
	uniformly in $Y$. Hence:\vspace{-2mm} \begin{align*}
		\|S_{\varepsilon ,\ptronc} -\SDA_{\varepsilon ,\ptronc} \|_{H^{1} (\couche_{\varepsilon } \cap \crown)}
		&\leqslant  \sum _{p\in \NN\cap [0,\ptronc]} \sum _{\ell =0} ^{n_{p}}  \varepsilon ^{p} \ln ^{\ell } \!\varepsilon \, \Big\|\big[S_{p,\ell } -T_{\geqslant p-\ptronc} (\SDA_{p,\ell } )\big]({\textstyle \frac{x}{\varepsilon } , \frac{y}{\varepsilon } })\Big\|_{H^{1} (\couche_{\varepsilon } \cap \crown)} \\
		&\lesssim  \sum _{p\in \NN\cap [0,\ptronc]} \sum _{\ell =0} ^{n_{p}}  \varepsilon ^{p} \ln ^{\ell } \!\varepsilon \cdot \varepsilon ^{-1} \,\left( \int _{[\eta /2,\infty [\times [-\varepsilon ,0]} \left( \frac{x}{\varepsilon } \right)^{2(p-\ptronc)} \,\mathrm{d}x\,\mathrm{d}y \right)^{1/2} \\
		&\lesssim  \sum _{p\in \NN\cap [0,\ptronc]} \sum _{\ell =0} ^{n_{p}}  \varepsilon ^{p} \ln ^{\ell } \!\varepsilon \cdot  \varepsilon ^{\ptronc-p-1} \cdot \varepsilon ^{ \frac{1}{2} } \eta ^{p-\ptronc+ \frac{1}{2} } \\
		&\lesssim \; \eta ^{\ptronc-2} \qquad  \quad  \text{because }  \varepsilon =\eta ^2  \text{ and } \ln ^{\ell } \!\varepsilon \lesssim \eta ^{- \frac{3}{2} }
	\end{align*}
	
	\item \underline{$\|\bfuDA_{\varepsilon ,\ptronc} -\SDA_{\varepsilon ,\ptronc} \|_{H^{1} (\Omega _{\varepsilon } \cap \crown)}$ :} %
	\def\termeu{\tilde{\bfu}}%
	\def\termeS{\tilde{S}}%
	Let us show that this norm vanishes. We recall that the physical scaling by $\varepsilon $, denoted $\pEps_{\varepsilon }$, is defined at page~\pageref{page: def pEps} while the formal scaling by $\feps$, denoted $\Eps$ in introduced in \cref{def: Eps}. Here $\feps$ is the formal indeterminate of $\Aeps(\Pi )$, see \cref{sec: SF 2D de singularites de coin}. In addition, $\evaleps$ maps any element of $\Aeps(\Pi )$ with finite support to the element of $\mathcal{A}(\Pi )$ obtained by replacing $\feps$ with $\varepsilon $, see \cref{eq: def evaleps}. By~\cref{eq: relation Eps pEps}, we have
	\begin{align*}
		\sum _{p\in \NN\cap [0,\ptronc]} \sum _{\ell \in \mathbb{N}} \powelne \cdot T_{\leqslant \ptronc-p} (\bfuDA_{p,\ell } )&= \evaleps\bigg(\sum _{p\in \NN\cap [0,\ptronc]} \sum _{\ell \in \mathbb{N}} \feps^{p} \ln ^{\ell } \!\feps \cdot T_{\leqslant \ptronc-p} (\bfuDA_{p,\ell } ) \bigg) \\
		\text{and }  \sum _{p\in \NN\cap [0,\ptronc]} \sum _{\ell \in \mathbb{N}} \powelne \cdot \pEps_{\varepsilon } ^{-1} \big(T_{\geqslant p-\ptronc} (\SDA_{p,\ell } )\big) &= \evaleps\circ \Eps^{-1} \bigg(\sum _{p\in \NN\cap [0,\ptronc]} \sum _{\ell \in \mathbb{N}} \feps^{p} \ln ^{\ell } \!\feps \cdot T_{\geqslant p-\ptronc} (\SDA_{p,\ell } ) \bigg). \quad 
	\end{align*}
	The above functions are equal to $\bfuDA_{\varepsilon ,\ptronc}$, resp.\ $\SDA_{\varepsilon ,\ptronc}$, up to a scaling w.r.t.\ $y$ in the layer. Hence, to prove $\bfuDA_{\varepsilon ,\ptronc} =\SDA_{\varepsilon ,\ptronc}$, it suffices to show the following equality in $\Aemacro(\Pi )$:
	\begin{equation}{\label{eq: exactitude du raccord algebrique tronque}}
	\sum _{p\in \NN\cap [0,\ptronc]} \sum _{\ell \in \mathbb{N}} \feps^{p} \ln ^{\ell } \!\feps \cdot T_{\leqslant \ptronc-p} (\bfuDA_{p,\ell } ) = \Eps^{-1} \bigg(\sum _{p\in \NN\cap [0,\ptronc]} \sum _{\ell \in \mathbb{N}} \feps^{p} \ln ^{\ell } \!\feps \cdot T_{\geqslant p-\ptronc} (\SDA_{p,\ell } ) \bigg) \end{equation}
	For~any $(p,d)$ let $\termeu_{p,d}$ be the coordinate of $\fsum_{p,\ell } \feps^{p} \ln ^{\ell } \!\feps\, \bfuDA_{p,\ell }$ in $\mathcal{A}_{d} (\Pi )[\ln \feps]$ and $\termeS_{p,d}$ be the one of $\fsum_{p,\ell } \feps^{p} \ln ^{\ell } \!\feps\, \SDA_{p,\ell }$. Moreover, for~any $d$ and $\varphi \in \mathcal{A}_{d} (\Pi )[\ln \feps]$, let $\preEps^{-1} (\varphi ):= \feps^{d} \.\Eps^{-1} (\varphi )$. By~\cref{def: Eps}, $\preEps$ maps $\mathcal{A}_{d} (\Pi )[\ln \feps]$ to itself. Then we have:
	\def\iff{\relax\,\Leftrightarrow\,}
	\begin{align*}
		\text{\cref{eq: exactitude du raccord algebrique tronque}}  &\iff \sum _{p\leqslant \ptronc} \sum _{d\leqslant \ptronc-p} \feps^{p} \termeu_{p,d} = \Eps^{-1} \bigg(\sum _{p\leqslant \ptronc} \sum _{d\geqslant p-\ptronc} \feps^{p} \termeS_{p,d} \bigg)\\
		&\iff \sum _{p\leqslant \ptronc} \sum _{p+d\leqslant \ptronc} \feps^{p} \termeu_{p,d} = \sum _{p\leqslant \ptronc} \sum _{p-d\leqslant \ptronc} \feps^{p-d} \.\preEps^{-1} ( \termeS_{p,d} )\\
		&\iff \sum _{p\leqslant \ptronc} \sum _{p+d\leqslant \ptronc} \feps^{p} \termeu_{p,d} = \sum _{p+d\leqslant \ptronc} \sum _{p\leqslant \ptronc} \feps^{p} \.\preEps^{-1} ( \termeS_{p+d,d} )\\
		&\iff \forall (p,d) \text{ s.t.\ }  p\leqslant \ptronc \text{ and } p+d\leqslant \ptronc, \quad  \termeu_{p,d} = \preEps^{-1} ( \termeS_{p+d,d} ) \stepcounter{equation}\tag{\theequation}\label{eq: exactitude du raccord algebrique tronque bis}
	\end{align*}
	because two formal series coincide iff their coordinates coincide one by one.\\
	But \cref{confimation de l'hyp de raccord} shows the matching assumption
	\[ \fsum\limits_{p,d} \feps^{p} \termeu_{p,d} = \Eps^{-1} \bigg(\fsum\limits_{p,d} \feps^{p} \termeS_{p,d} \bigg) =\fsum\limits_{p,d} \feps^{p} \.\preEps^{-1} ( \termeS_{p+d,d} ),\]
	which is equivalent to: $\forall (p,d),\;\termeu_{p,d} = \preEps^{-1} ( \termeS_{p+d,d} )$. Thus we get \cref{eq: exactitude du raccord algebrique tronque bis}, and then \cref{eq: exactitude du raccord algebrique tronque}.
\end{enumerate}\fbi
\end{demo}

\subsection{Global and far fields error estimates}

\begin{theoreme}{global error estimate}
\label{estimation d'erreur globale}%
For~any $\ptronc \in \mathbb{R}_+$ we have\ $\big\|u_{\varepsilon } -u_{\varepsilon ,\ptronc} \big\|_{H^{1} (\Omega _{\varepsilon } )} = O(\varepsilon ^{ \frac{\ptronc}{2} -2} )$.\vspace{-1mm}
\end{theoreme}

\begin{remarque}
\Cref{estimation d'erreur globale} can be improved to $\|u_{\varepsilon } -u_{\varepsilon ,\ptronc} \|_{H^{1} (\Omega _{\varepsilon } )} = o(\varepsilon ^{\ptronc/2} )$. Indeed with the same notations as the above proof, one can show that, for~any $(p,\ell )\in \NN\times \mathbb{N}$, $\|(1-\chi _{\eta } )\bfue_{p,\ell } \|_{H^{1} (\Omega _{\varepsilon } )} \lesssim  \eta ^{p}$ and $\|\chi _{\eta } S_{p,\ell } \|_{H^{1} (\Omega _{\varepsilon } )} \lesssim  \eta ^{-p}$. Therefore:
\[\|r_{\varepsilon ,\ptronc} \|_{H^{1} (\Omega _{\varepsilon } )} \lesssim  \|r_{\varepsilon ,\ptronc+4} \|_{H^{1} (\Omega _{\varepsilon } )} + \!\! \sum _{p\in \NN\cap (\ptronc,\ptronc+4]} \sum _{\ell \in \mathbb{N}} \varepsilon ^{p} \ln ^{\ell } \!\varepsilon  \cdot \big(\|(1-\chi _{\eta } )\bfue_{p,\ell } \|_{H^{1} (\Omega _{\varepsilon } )} +\|\chi _{\eta } S_{p,\ell } \|_{H^{1} (\Omega _{\varepsilon } )}\big) = o(\varepsilon ^{ \frac{\ptronc}{2} }).\]
Even with this improvement, the obtained order of convergence is only half of what we could hope. The same limitation in present in similar works, see \cite{DelSchSem16,SemDelSch18}. To prove an estimate of order $\varepsilon ^{\ptronc}$ on the far fields, we will have to use the asymptotic expansion at order $2\ptronc$. 
\end{remarque}

\begin{demo}
Let $r_{\varepsilon ,\ptronc} := u_{\varepsilon ,\ptronc} -u_{\varepsilon }$. It satisfies, for some functions $f_{\varepsilon }$ and $g_{\varepsilon }$,
\[\left\{\begin{array}{r@{\;=\;}ll}
	\operatorname{div}(\mu _{\varepsilon } \nabla  r_{\varepsilon ,\ptronc} ) +\omega ^2 \rho _{\varepsilon } r_{\varepsilon ,\ptronc} & f_{\varepsilon } &\text{ in }\Omega _{\varepsilon } \\[0.5mm]
	r_{\varepsilon ,\ptronc} & 0 &\text{ on }\partial \Omega _{\varepsilon } \\[0.5mm]
	r_{\varepsilon ,\ptronc|y=0^{+}} - r_{\varepsilon ,\ptronc|y=0^{-}} & 0 &\text{ on }\Gamma \cap \{x>\eta /2\}\\[0.5mm]
	\mu _{0} \partial _{y} r_{\varepsilon ,\ptronc|y=0^{+}} - \mu _{1} \partial _{y} r_{\varepsilon ,\ptronc|y=0^{-}} & g_{\varepsilon } &\text{ on }\Gamma \cap \{x>\eta /2\}
\end{array}\right.\]
As for \cref{eq: uniforme coercivite}, this problem is well-posed with a stability constant independent of $\varepsilon $:
\[\|r_{\varepsilon ,\ptronc} \|_{H^{1} (\Omega _{\varepsilon } )} \lesssim  \|f_{\varepsilon } \|_{L^2 (\Omega _{\varepsilon } )} + \|g_{\varepsilon } \|_{L^2 (\Gamma \cap \{x>\eta /2\})} .\]
To get our error estimate, it suffices to show that $\|f_{\varepsilon } \|_{L^2 } +\|g_{\varepsilon } \|_{L^2 } \lesssim  \varepsilon ^{ \frac{\ptronc}{2} -2} =\eta ^{\ptronc-4}$.

\medskip\underline{Estimate of $\|f_{\varepsilon } \|_{L^2 }$:} We denote:
\begin{itemize}
	\item $\bfu_{\varepsilon ,\ptronc} := \sum _{p\in \NN\cap [0,\ptronc]} \sum _{\ell \in \mathbb{N}} \varepsilon ^{p} \ln ^{\ell } \!\varepsilon \, \bfue_{p,\ell } $ in $\Pi _{\varepsilon }$
	\item $S_{\varepsilon ,\ptronc} := \sum _{p\in \NN\cap [0,\ptronc]} \sum _{\ell \in \mathbb{N}} \varepsilon ^{p} \ln ^{\ell } \!\varepsilon \, S_{p,\ell } ({\textstyle \frac{x}{\varepsilon } , \frac{y}{\varepsilon } })$ in $\Omega _{\varepsilon }$
	\item and $\mathcal{D}_{\varepsilon } :u\mapsto  \operatorname{div}(\mu _{\varepsilon } \nabla  u) +\omega ^2 \rho _{\varepsilon } u$ the differential operator of Helmholtz's equation.
\end{itemize}
Since $\mathcal{D}_{\varepsilon } u_{\varepsilon } =f$ by~definition of $u_{\varepsilon }$, we have \begin{align*}
	f_{\varepsilon } &= \mathcal{D}_{\varepsilon } r_{\varepsilon ,\ptronc} \\
	&= \mathcal{D}_{\varepsilon } \big((1-\chi _{\eta } )\bfu_{\varepsilon ,\ptronc} + \chi _{\eta } S_{\varepsilon ,\ptronc} -u_{\varepsilon } \big) \\
	&= \big((1-\chi _{\eta } ) \mathcal{D}_{\varepsilon } \bfu_{\varepsilon ,\ptronc} -f\big)+ \chi _{\eta } \mathcal{D}_{\varepsilon } S_{\varepsilon ,\ptronc} +[\mathcal{D}_{\varepsilon } ,\chi _{\eta } ](S_{\varepsilon ,\ptronc} -\bfu_{\varepsilon ,\ptronc} )
\end{align*}
where $[.,.]$ is the commutator. Let us estimate these terms one by one.
\begin{enumerate}
	\item For $\varepsilon $ small enough, $1-\chi _{\eta }$ is equal to 0 in $B(0,\varepsilon \Rcoin)$ and 1 in $\operatorname{supp}(f)$, so $(1-\chi _{\eta } ) \mathcal{D}_{\varepsilon } \bfu_{\varepsilon ,\ptronc|\Omega } = (1-\chi _{\eta } )(\mu _{0} \Delta +\omega ^2 \rho _{0} ) \bfu_{\varepsilon ,\ptronc|\Omega } = (1-\chi _{\eta } )f=f$.\\
	Moreover, using that $\mu _{1} \partial _{Y} ^2  U_{p,\ell } = -(\mu _{1} \partial _{x} ^2 +\omega ^2 \rho _{1} )U_{p-2,\ell }$ for~any $(p,\ell )$, we get in $\couche_{\varepsilon } \setminus B(0,\varepsilon \Rcoin)$:
	\begin{align*}
		\mathcal{D}_{\varepsilon } \bfu_{\varepsilon ,\ptronc} &= (\mu _{1} \Delta +\omega ^2 \rho _{1} )\bfu_{\varepsilon ,\ptronc} \\
		&= \sum _{p\in \NN\cap (\ptronc-2,\ptronc]} \sum _{\ell \in \mathbb{N}} \varepsilon ^{p} \ln ^{\ell } \!\varepsilon \,(\mu _{1} \partial _{x} ^2 +\omega ^2 \rho _{1} ) U_{p,\ell } \\
		&= \sum _{p\in \NN\cap (\ptronc-2,\ptronc]} \sum _{\ell \in \mathbb{N}} \varepsilon ^{p+ \frac{1}{2} } \ln ^{\ell } \!\varepsilon  \sum _{n\in \mathbb{N}} \mu _{1} \,(\partial _{x} ^2 +k_{1} ^2 )^{n+1} \partial _{y} u_{p-1-2n,\ell } (x,0) \cdot  \mathcal{U}_{n} ({\textstyle \frac{y}{\varepsilon } }) \tag*{by \cref{eq: definition des champs de couche}}
	\end{align*}
	Since $\mathcal{U}_{n} \in L^{\infty } (-1,0)$ for~any $n$ (see \cref{eq: correcteurs de couche}), we deduce:
	\[\|(1-\chi _{\eta } ) \mathcal{D}_{\varepsilon } \bfu_{\varepsilon ,\ptronc} \|_{L^2 (\couche_{\varepsilon } )} \lesssim  \!\! \sum _{p\in \NN\cap (\ptronc-2,\ptronc]} \sum _{\ell =0} ^{n_{p}}  \sum _{n=0} ^{\lfloor {(p-1)/2}\rfloor }  \varepsilon ^{p+ \frac{1}{2} } \ln ^{\ell } \!\varepsilon  \cdot \big\|(\partial _{x} ^2 +k_{1} ^2 )^{n+1} \partial _{y} u_{p-1-2n,\ell } \big\|_{L^2 (\Gamma \cap \{x>\eta /2\})} .\]
	Let us show that, for~any $p,\ell ,m$, $\|\partial _{x} ^{m} \partial _{y} u_{p,\ell |\Gamma } \|_{L^2 (\Gamma \cap \{x>\eta /2\})} \lesssim  \eta ^{-p-m-1}$ \inlineEqNb\label{eq: estimation de derivees de upl sur Gamma}.
	\begin{itemize}
		\item \Cref{validation construction et DA des champs} states that $\partial _{x} ^{m} \partial _{y} u_{p,\ell } \in H^{1} _{0} (\Gamma )+\khiz \mathcal{A}(\Gamma )$ for~any $p,\ell ,m$. So $\|\partial _{x} ^{m} \partial _{y} u_{p,\ell |\Gamma } \|_{L^2 (\Gamma \cap \{x>1\})} <\infty $.
		\item By~\cref{validation construction et DA des champs}, $\partial _{y} u_{p,\ell |\Gamma } = \partial _{y} T_{\leqslant -p- \frac{1}{2} } (\bfuDA_{p,\ell } ) +o_{\partial } (x^{-p- \frac{3}{2} })$ when $x\rightarrow 0$, and \cref{ansatz sur Aeps} is satisfied, thus $T_{\leqslant -p- \frac{1}{2} } (\bfuDA_{p,\ell } )=0$. So: $\forall m\in \mathbb{N},\; \partial _{x} ^{m} \partial _{y} u_{p,\ell |\Gamma } =O(x^{-p- \frac{3}{2} -m})$ when $x\rightarrow 0$. It implies $\|\partial _{x} ^{m} \partial _{y} u_{p,\ell |\Gamma } \|_{L^2 (\Gamma \cap \{\eta /2<x<1\})} \lesssim  \eta ^{-p-m-1}$.
	\end{itemize}
	Hence:\vspace{-2mm} \begin{align*}
		\|(1-\chi _{\eta } ) \mathcal{D}_{\varepsilon } \bfu_{\varepsilon ,\ptronc} \|_{L^2 (\couche_{\varepsilon } )}
		&\lesssim  \sum _{p\in \NN\cap (\ptronc-2,\ptronc]} \sum _{\ell =0} ^{n_{p}}  \varepsilon ^{p+ \frac{1}{2} } \ln ^{\ell } \!\varepsilon  \sum _{n=0} ^{\lfloor {(p-1)/2}\rfloor }  \eta ^{-(p-1-2n)-2(n+1)-1} \\
		&\lesssim  \sum _{p\in \NN\cap (\ptronc-2,\ptronc]} \sum _{\ell =0} ^{n_{p}}  \varepsilon ^{p+ \frac{1}{2} } \ln ^{\ell } \!\varepsilon  \cdot \eta ^{-p-2} \\
		&\lesssim  \eta ^{\ptronc-4} \qquad  \quad  \text{because }  \ln ^{\ell } \!\varepsilon  \lesssim \eta ^{-1}
	\end{align*}
	So $\|(1-\chi _{\eta } )\mathcal{D}_{\varepsilon } \bfu_{\varepsilon ,\ptronc} -f\|_{L^2 (\Omega _{\varepsilon } )} \lesssim  \eta ^{\ptronc-4}$.
	
	\item For the second term of $f_{\varepsilon }$, \cref{eq: definition des champs de coin} implies that $\mathcal{D}_{\varepsilon } S_{\varepsilon ,\ptronc} = \sum _{p\in \NN\cap (\ptronc-2,\ptronc]} \sum _{\ell \in \mathbb{N}} \varepsilon ^{p} \ln ^{\ell } \!\varepsilon \,\omega ^2 \rho  S_{p,\ell }$. But by \cref{validation construction et DA des champs,ansatz sur Aeps}: $\forall (p,\ell ),\; S_{p,\ell } =O(r^{p} )$. Thus, using that $\rho $ is bounded and $\|S_{p,\ell } ( \frac{x}{\varepsilon } , \frac{y}{\varepsilon } ) \|_{L^2 (\Omega _{\varepsilon } \cap B(0,\eta ))} = \varepsilon \, \|S_{p,\ell } \|_{L^2 (\Omega _{1} \cap B(0,\eta /\varepsilon ))}$, we get \begin{align*}
		\|\chi _{\eta } \mathcal{D}_{\varepsilon } S_{\varepsilon ,\ptronc} \|_{L^2 (\Omega _{\varepsilon } )} &\lesssim  \sum _{p\in \NN\cap (\ptronc-2,\ptronc]} \sum _{\ell =0} ^{n_{p}}  \varepsilon ^{p} \ln ^{\ell } \!\varepsilon  \cdot \varepsilon \, \|S_{p,\ell } \|_{L^2 (\Omega _{1} \cap B(0,\eta /\varepsilon ))} \\
		&\lesssim  \sum _{p\in \NN\cap (\ptronc-2,\ptronc]} \sum _{\ell =0} ^{n_{p}}  \varepsilon ^{p+1} \ln ^{\ell } \!\varepsilon \cdot \Big( \frac{\eta }{\varepsilon } \Big)^{p+1} \\
		&\lesssim  \eta ^{\ptronc-4} \qquad  \qquad  \text{because }  \ln ^{\ell } \!\varepsilon  \lesssim \eta ^{-3}
	\end{align*}
	
	\item Using that $\|\Delta  \chi _{\eta } \|_{L^2 } +\|\nabla  \chi _{\eta } \|_{L^2 } =O(\eta ^{-2} )$, the last term of $f_{\varepsilon }$ satisfies \begin{align*}
		\|[\mathcal{D}_{\varepsilon } ,\chi _{\eta } ](S_{\varepsilon ,\ptronc} -\bfu_{\varepsilon ,\ptronc} )\|_{L^2 (\Omega _{\varepsilon } )}
		&\lesssim  \|\Delta  \chi _{\eta } \cdot (S_{\varepsilon ,\ptronc} -\bfu_{\varepsilon ,\ptronc} ) + 2 \nabla  \chi _{\eta } \cdot \nabla (S_{\varepsilon ,\ptronc} -\bfu_{\varepsilon ,\ptronc} )\|_{L^2 (\Omega _{\varepsilon } )} \\
		&\lesssim  \eta ^{-2} \|S_{\varepsilon ,\ptronc} -\bfu_{\varepsilon ,\ptronc} \|_{H^{1} (\Omega _{\varepsilon } \cap \crown)} \\
		&\lesssim  \eta ^{\ptronc-4} \qquad  \qquad  \text{by \cref{erreur de raccord}.} 
	\end{align*}\fbeq
\end{enumerate}

\underline{Estimate of $\|g_{\varepsilon } \|_{L^2 }$:} Without loss of generality, we assume $\partial _{y} \chi _{\eta |\Gamma } =0$. \cref{eq: definition des champs de coin} implies that $\mu _{0} \partial _{Y} S_{p,\ell |Y=0^{+}} = \mu _{1} \partial _{Y} S_{p,\ell |Y=0^{-}}$ on $\Gamma \cap \{X>\Rcoin\}$ for~any $(p,\ell )$. Thus we have on $\Gamma \cap \{x>\eta /2\}$: \begin{align*}
		g_{\varepsilon } &= \mu _{0} \partial _{y} r_{\varepsilon ,\ptronc|y=0^{+}} - \mu _{1} \partial _{y} r_{\varepsilon ,\ptronc|y=0^{-}} \\
		&= (1-\chi _{\eta } ) \cdot \big(\mu _{0} \partial _{y} \bfu_{\varepsilon ,\ptronc|y=0^{+}} - \mu _{1} \partial _{y} \bfu_{\varepsilon ,\ptronc|y=0^{-}} \big) \\
		&= (1-\chi _{\eta } ) \sum _{p\in \NN\cap (\ptronc-1,\ptronc]} \sum _{\ell \in \mathbb{N}} \varepsilon ^{p} \ln ^{\ell } \!\varepsilon \cdot \mu _{0} \partial _{y} u_{p,\ell |y=0^{+}} 
	\end{align*}
	By~\cref{eq: estimation de derivees de upl sur Gamma}, we have that $\|\partial _{y} u_{p,\ell |\Gamma } \|_{L^2 (\Gamma \cap \{x>\eta /2\})} \lesssim  \eta ^{-p-1}$. So $\|g_{\varepsilon } \|_{L^2 (\Gamma \cap \{x>\eta /2\})} \lesssim  \eta ^{\ptronc-4}$ since $\ln ^{\ell } \varepsilon \lesssim \eta ^{-2}$.
\end{demo}

\begin{theoreme}{error estimate for the asymptotic expansion of $u_{\varepsilon }$}
\label{DA de SE avec CL et delta variable}%
Let $\ptronc \in \mathbb{R}_+$ and $\ell _{\ptronc} :=\max \{\ell \in \mathbb{N} \mid u_{\ptronc,\ell } \neq 0\}$ (with the convention $\ell _{\ptronc} :=0$ if this set is empty). Let $0<\delta \leqslant 1$. There exists a constant $C_{\ptronc}$ that depends on $(\Omega _{1} ,\mu ,\rho ,\omega ,\ptronc,\fsource)$ but not on $(\varepsilon ,\delta )$ s.t.\ we have for~any $\varepsilon  \in (0,\delta ^2 )$:
\begin{equation}{\label{eq: DA de SE avec CL et delta variable}}
	\bigg\|u_{\varepsilon } -\sum _{p\in \NN\mkern-1mu,\,p<\ptronc} \sum _{\ell \in \mathbb{N}} \varepsilon ^{p} \ln ^{\ell } \!\varepsilon \, u_{p,\ell } \bigg\|_{H^{1} (\Omega  \setminus \closedBall(0,\delta ))} \leqslant  C_{\ptronc} \,\. \delta ^{-\ptronc} \varepsilon ^{\ptronc} (1+|{\ln \varepsilon }|)^{\ell _{\ptronc}} 
\end{equation}\fbeq
\end{theoreme}

\begin{remarquenumerotee}
\label{rq: condition sur delta dans DA de SE avec CL}%
The condition ``$\varepsilon \leqslant \delta ^2 $'' is due to the fact that $u_{\varepsilon ,\ptronc}$ is a sum of far fields only in the far zone $\Omega  \cap \closedBall(0,\eta )$. Since we took $\eta :=\sqrt{\varepsilon }$, this zone contains $\Omega  \setminus \closedBall(0,\delta )$ only when $\delta ^2 \geqslant \varepsilon $. However, that condition can be relaxed to ``$\varepsilon \leqslant M\delta ^{a}$'' for any fixed $M>0$ and $a>1$ (but the constant in \cref{eq: DA de SE avec CL et delta variable} then depends on $(M,a)$). Indeed, \cref{erreur de raccord,estimation d'erreur globale} can be adapted to $\eta :=M^{-1/a} \varepsilon ^{1/a}$ (which modifies the matching zone, but not the far fields). The error estimate of \cref{estimation d'erreur globale} may be deteriorated ($\varepsilon ^{\ptronc/2}$ is replaced by $\varepsilon ^{b\ptronc}$ for some $b$ linked to $a$), but the proof technique used below can be used anyway. 
\end{remarquenumerotee}

\begin{demo}
\def\OD{\Omega^\delta}%
Let $\OD:=\Omega  \setminus \closedBall(0,\delta )$. We will use the symbol $\lesssim $ for inequalities that hold up to a constant independent of $\varepsilon $ and $\delta $. To begin, let us show that, for any $(p,\ell )\in \NN\times \mathbb{N}$, we have $\|u_{p,\ell } \|_{H^{1} (\OD)} \lesssim  \delta ^{-p}$.
\begin{sousdemo}
	\Cref{validation construction et DA des champs} shows $u_{p,\ell } = T_{\leqslant 1} (\uDA_{p,\ell } )+o_{\partial } (r)$, so there exists $r_{0} >0$ s.t.\ we have on $\Omega \cap  B(0,r_{0} )$:
	\[|u_{p,\ell } -T_{\leqslant 1} (\uDA_{p,\ell } )|\lesssim  r \text{ and } |\nabla  (u_{p,\ell } -T_{\leqslant 1} (\uDA_{p,\ell } ))|\lesssim  1.\]
	Moreover \cref{precision de l'ansatz sur uDA} yields $T_{\leqslant 1} (\uDA_{p,\ell } ) \in \sum _{d\geqslant -p+\pst} \mathcal{A}_{d} (\Omega ) \subset \sum _{d>-p} \mathcal{A}_{d} (\Omega )$. Thus, we have on $\Omega \cap  B(0,r_{0} )$: $|u_{p,\ell } |\lesssim  r^{-p}$ and $|\nabla  u_{p,\ell } |\lesssim  r^{-p-1}$. Finally, we have $u_{p,\ell } \in H^{1} (\Omega )+\khiz \mathcal{A}(\Omega )$ by construction, so $u_{p,\ell }$ is also $H^{1}$ on $\Omega  \setminus \closedBall(0,r_{0} )$. The bound on $\|u_{p,\ell } \|_{H^{1} (\OD)}$ follows by integration.
\end{sousdemo}
Next, let us use \cref{estimation d'erreur globale} at order $\ptronc' := 2\ptronc+4$. It gives $\|u_{\varepsilon } -u_{\varepsilon ,\ptronc'} \|_{H^{1} (\Omega _{\varepsilon } )} =o(\varepsilon ^{\ptronc} )$, thus:
\begin{align*}
\bigg\|u_{\varepsilon } - \sum _{p\in \NN \cap [0,\ptronc)} \sum _{\ell \in \mathbb{N}} \varepsilon ^{p} \ln ^{\ell } \!\varepsilon \, u_{p,\ell } \bigg\|_{H^{1} (\OD)} &\lesssim  \|u_{\varepsilon } -u_{\varepsilon ,\ptronc'} \|_{H^{1} (\Omega _{\varepsilon } )} +  \sum _{p\in \NN\cap [\ptronc,\ptronc']} \sum _{\ell \in \mathbb{N}} \varepsilon ^{p} \ln ^{\ell } \!\varepsilon  \,\|u_{p,\ell } \|_{H^{1} (\OD)} \\
&\lesssim  \varepsilon ^{\ptronc} + \sum _{p\in \NN\cap [\ptronc,\ptronc']} \sum _{\ell =0} ^{\ell _{p}} \varepsilon ^{p} \.|{\ln \varepsilon }|^{\ell _{p}} \cdot  \delta ^{-p} 
\end{align*}
where we denoted $\ell _{p} :=\max \{\ell \in \mathbb{N} \mid u_{p,\ell } \neq 0\}$, similarly to $\ell _{\ptronc}$. For~any $a>0$ and $\ell \in \mathbb{N}$, we have $\varepsilon ^{a} |{\ln \varepsilon }|^{\ell } \.\delta ^{-a} \lesssim  \varepsilon ^{a/2} \.\delta ^{-a} \leqslant  1$ by assumption. Applying this to $a:=p-\ptronc$ when $p>\ptronc$, we obtain
\[\bigg\|u_{\varepsilon } - \sum _{p\in \NN \cap [0,\ptronc)} \sum _{\ell \in \mathbb{N}} \varepsilon ^{p} \ln ^{\ell } \!\varepsilon \, u_{p,\ell } \bigg\|_{H^{1} (\OD)} \lesssim \; \varepsilon ^{\ptronc} + \sum _{\ell =0} ^{\ell _{\ptronc}} \varepsilon ^{\ptronc} \.|{\ln \varepsilon }|^{\ell _{\ptronc}} \delta ^{-\ptronc} + \varepsilon ^{\ptronc} \delta ^{-\ptronc} 
\;\lesssim \; \delta ^{-\ptronc} \varepsilon ^{\ptronc} (1+|{\ln \varepsilon }|)^{\ell _{\ptronc}} .\]
\end{demo}

\begin{textesansboite}
In particular, \cref{DA de SE avec CL et delta variable} implies the following weaker result that was announced in introduction.
\end{textesansboite}

\begin{theoreme}{}
\label{DA de SE avec CL et delta fixe}%
For~any $\ptronc \geqslant 0$ and $\delta >0$ fixed, we have:
\[\bigg\|u_{\varepsilon } -\sum _{p\in \NN\mkern-1mu,\,\red{p\leqslant \ptronc}} \sum _{\ell \in \mathbb{N}} \varepsilon ^{p} \ln ^{\ell } \!\varepsilon \, u_{p,\ell } \bigg\|_{H^{1} (\Omega  \setminus \closedBall(0,\delta ))} = o(\varepsilon ^{\ptronc} ) \qquad  \text{when }  \varepsilon \rightarrow 0.\]\fbeq
\end{theoreme}

\begin{demo}
It suffices to apply \cref{DA de SE avec CL et delta variable} with a truncation order $\ptronc'>\ptronc$ s.t.\ $\NN \cap [0,\ptronc')=\NN \cap [0,\ptronc]$.
\end{demo}

\begin{remarque}
The stronger error estimate of \cref{DA de SE avec CL et delta variable} is useful to obtain precise error estimates on effectives models, see \cref{cvg des ME}.
\end{remarque}

\chapter{Further properties of the asymptotic expansion}
\label{chap: complements sur le DA}

\section{Expression of far fields singularities using Bessel functions}
\label{sec: expression des mcl}

\begin{textesansboite}
In this section, we express some of the corner singularities of the far fields, like $\langle -k_{0} ^2 \ROmegaPi\rangle (\phi_{d} )_{|\Omega }$, using the Bessel functions. \cref{eq: bfuDApl} shows that for~any $(p,\ell )$, $\uDA_{p,\ell } :=\bfuDA_{p,\ell |\Omega }$ depends on the $\sigma _{d} (\bfuDA_{p',\ell '} )$ with $p'<p$. By~singularity of a far field, we mean here all the terms of $\uDA_{p,\ell }$ that are multiplied by a given $\sigma _{d} (\bfuDA_{p',\ell '} )$. Expressions of these singularities not only may be used to compute numerically the singularities, but they will also give us more information on the structure of the asymptotic expansion of $u_{\peps}$ (see \cref{sec: simplifications DA}). In addition, it~is an important part of the proof of the far fields' modal decomposition in \cref{chap: decompo modale CL}, which will help designing numerical methods to compute the far fields and effective models (see \cref{sec: numerique des CL,sec: decompo modale de SE}).\\

\nomenclature[19-10]{$J_{a}$}{Bessel function of the first kind}
We recall that, for~any $a\in \mathbb{R}$, the Bessel function of the first kind $J_{a} $ is defined by the following converging series 
\begin{equation}{\label{eq: def BesselJ}}
    \forall z\in \mathbb{C}^*, \qquad  J_{a} (z) := \sum_{k=0}^\infty  \frac{(-1)^{k} }{k!\, \Gamma (a+k+1)} \left( \frac{z}{2} \right)^{a+2k}
\end{equation}
where $\Gamma $ denotes the Gamma function. $\Gamma $ has a pole at every $n\in -\mathbb{N}$, so when $a+k+1 \in  -\mathbb{N}$ we set by convention $\frac{1}{\Gamma (a+k+1)} :=0$. See \cite{Leb72} for an introduction to Bessel functions.
\end{textesansboite}
\nomenclature[38b]{$\Gamma $}{The Gamma function (depending on the context)}

\begin{definitionsansboite}{}
\label{def: convergence de SF dans Amacro}%
Let $D\in \{\Omega ,\Gamma \}$, $\varphi  \in  \Amacro(D)$ and $f:D\rightarrow \mathbb{C}$. We say that $\varphi $ converges to $f$, which we denote $\varphi  \longrightarrow  f$, iff for~any $R>0$, $T_{<q} (\varphi )$ converges to $f$ uniformly on $D\cap B(0,R)$ when $q\rightarrow \infty $. In this case we denote $\lim  \varphi  :=f$.
\end{definitionsansboite}

\nomenclature[19b-10]{$\BesselGamma_{a} ,\BesselOmega_{a}$}{Formal series associated to Bessel functions of the first kind}
\begin{textesansboite}
For~any $a\in \mathbb{R}$, let us define the following formal series
\begin{align}
    \label{eq: def BesselGamma}
    \BesselGamma_{a} &:= \fsumn \frac{(-1)^{n} }{n!\, \Gamma (a+n+1)} \left( \frac{k_{0} }{2} \right)^{a+2n} x^{a+2n} \ \in \; \Amacro(\Gamma ) \\
    \label{eq: def BesselOmega}
    \BesselOmega_{a} &:= \fsumn \frac{(-1)^{n} }{n!\, \Gamma (a+n+1)} \left( \frac{k_{0} }{2} \right)^{a+2n} \Im[(\alpha z)^{a+n} \,\overline{\alpha z}^{n} ] \ \in \; \Amacro(\Omega ).
\end{align}
Be cautious that, in \cref{eq: def BesselGamma}, the letter $\Gamma $ is used twice with two different meanings : the Gamma function and the boundary $\Gamma \subset \partial \Omega $. By~definition of $J_{a}$ and using $\Im[(\alpha z)^{a+n} \,\overline{\alpha z}^{n} ] =r^{a+2n} \sin (a(\theta -\Theta ))$, we have 
\[\BesselGamma_{a} \longrightarrow  J_{a} (k_{0} \.x) \qquad  \text{ and }  \qquad  \BesselOmega_{a} \longrightarrow  J_{a} (k_{0} \.r) \.\sin (a(\theta -\Theta )).\]

The following lemma is the transcription to formal series of well known properties of the Bessel functions, that will be useful to express the far fields singularities using Bessel functions. Note that, compared to properties expressed on the Bessel functions that only give information on the limit of the series, equalities between formal series are term by term equalities.
\end{textesansboite}

\begin{lemmesansboite}{}
\label{relations sur BesselGamma/Omega}%
Let $a\in \mathbb{R}$. We have the following equalities in $\Amacro(\Omega )$ and $\Amacro(\Gamma )$:
\begin{enumerate}
    \item $(\Delta +k_{0} ^2 )\BesselOmega_{a} = 0$.
    \item $\BesselOmega_{a|\Gamma } = -\sin (a\Theta ) \BesselGamma_{a}$.
    \item $\partial _{y} \BesselOmega_{a|\Gamma } = \frac{k_{0} }{2} \cos (a\Theta ) (\BesselGamma_{a-1} + \BesselGamma_{a+1} )$.
    \item $\partial _{x} \BesselGamma_{a} = \frac{k_{0} }{2} (\BesselGamma_{a-1} -\BesselGamma_{a+1} )$
    \item If $a\in \mathbb{Z}$, then $\BesselGamma_{-a} =(-1)^{a} \BesselGamma_{a}$ and $\BesselOmega_{-a} = -(-1)^{a} \BesselOmega_{a}$.
\end{enumerate}\fbi
\end{lemmesansboite}

\begin{demo}
First, let us note that $\lim \BesselGamma_{a} =(r\mapsto  J_{a} (k_{0} r))$ and $\lim\BesselOmega_{a} = \big[(r,\theta )\mapsto  J_{a} (k_{0} r)\sin (a(\theta -\Theta ))\big]$ satisfy the above equalities.
\begin{enumerate}[label=(\arabic*')]
    \item It~is well known that $J_{a}$ satisfies $x^2 J_{a} ''(x)+xJ_{a} '(x)+(x^2 -a^2 )J_{a} (x) =0$ (see \cite{Leb72}), which yields $(\Delta +k_{0} ^2 )[J_{a} (k_{0} r)\sin (a(\theta -\Theta ))] = 0$ thanks to the expression of the laplacian in polar coordinates.
    \item $J_{a} (k_{0} r)\sin (a(\theta -\Theta ))|_{\theta =0} = -\sin (a\Theta ) J_{a} (k_{0} r)$
    \item We have, for~any $z$, $ \frac{2a}{z} J_{a} (z) = J_{a-1} (z)+J_{a+1} (z)$ \cite{Leb72}. Therefore $\partial _{y} [J_{a} (k_{0} r)\sin (a(\theta -\Theta ))]_{|\Gamma } = \frac{1}{r} J_{a} (k_{0} r) \cdot  a \cos (a\Theta )= \frac{k_{0} }{2} \cos (a\Theta ) (J_{a-1} (k_{0} r)+J_{a+1} (k_{0} r))$.
    \item We have $2 J_{a} '(z) = J_{a-1} (z)-J_{a+1} (z)$ \cite{Leb72}. Therefore $\partial _{x} [J_{a} (k_{0} x)] = \frac{k_{0} }{2} (J_{a-1} (k_{0} x)-J_{a+1} (k_{0} x))$.
    \item When $a\in \mathbb{Z}$, \cite{Leb72} states $J_{-a} =(-1)^{a} J_{a}$. Hence $J_{-a} (k_{0} r)= (-1)^{a} J_{a} (k_{0} r)$ and $J_{-a} (k_{0} r)\sin (-a(\theta -\Theta ))= -(-1)^{a} J_{a} (k_{0} r)\sin (a(\theta -\Theta ))$ for~any $(r,\theta )$.
\end{enumerate}
Now $r^{-a+1} J_{a} (k_{0} r)$, $r^{-a+1} J_{a-1} (k_{0} r)$ and $r^{-a+1} J_{a+1} (k_{0} r)$ are power series of $r$ (with infinite convergence radius), so we can identify the powers of $r$ in equalities (1')--(4'). This gives exactly the desired properties.
\end{demo}

\begin{lemmesansboite}{}
\label{sigma de BesselOmega}%
For~any $a\in \mathbb{R} \setminus (-\mathbb{N}^*)$ and $d\in \TZ*$, we have $\sigma _{d} (\BesselOmega_{a} )= \delta _{a,d} \, \frac{1}{\Gamma (a+1)} (\frac{k_{0} }{2} )^{a}$.
\end{lemmesansboite}

\begin{remarque}
$a\in -\mathbb{N}^*$ is a special case as this implies that $a+1$ is a pole of the Gamma function.
\end{remarque}

\begin{demo}
By~definition of $\BesselOmega_{a}$ and commutation of $\fsum$ with $\sigma _{d}$ (see \cref{def: sigma et fsum commutent dans Aemaicro}), it suffices to prove $\sigma _{d} \big(\Im[(\alpha z)^{a+k} \.\overline{\alpha z}^{k} ]\big) = \delta _{d,a} \,\delta _{k,0}$ for~any $k\in \mathbb{N}$. Since $a\not\in -\mathbb{N}^*$, there are three possible cases:
\begin{itemize}
    \item either $a=0$, in which case both $\Im[(\alpha z)^{a+k} \.\overline{\alpha z}^{k} ]$ and $\delta _{d,a}$ vanish,
    \item or $a\in \mathbb{N}^*$, which implies $a+k>k$,
    \item or $a\not\in \mathbb{Z}$, which implies $a+k\not\in \mathbb{N}$.
\end{itemize}
In the two last cases, $\Im[(\alpha z)^{a+k} \.\overline{\alpha z}^{k} ]$ belongs to the basis of $\mathcal{A}(\Omega )$ used to define $\sigma _{d}$, so we have $\sigma _{d} \big(\Im[(\alpha z)^{a+k} \.\overline{\alpha z}^{k} ]\big) = \delta _{d,a+2k} \,\delta _{k,0} = \delta _{d,a} \,\delta _{k,0}$  by \cref{def: sigma}.
\end{demo}

\begin{textesansboite}
First, let us look at the singularities of $u_{0,0}$, which are the $\langle -k_{0} ^2 \ROmegaPi\rangle (\phi_{d} )_{|\Omega }$. Indeed, by \cref{eq: bfuDApl}, we have 
\[\uDA_{0,0} := \bfuDA_{0,0|\Omega } = \langle -k_{0} ^2 \ROmegaPi\rangle  \bigg(\fsum_{d\in \TN} \sigma _{d} (\bfuDA_{0,0} ) \,\phi_{d} \bigg)_{|\Omega } .\]
Meanwhile, it~is classical that a homogeneous solution of Helmholtz equation in a corner with Dirichlet condition is a series of Bessel functions times $\sin (d\theta )$. The following proposition links the two results by showing that the $\langle -k_{0} ^2  \ROmegaPi\rangle (\phi_{d} )_{|\Omega }$ can be expressed using Bessel functions. Most of the time $J_{d} (k_{0} r)$ is used because it behaves like $r^{d}$ when $r\rightarrow 0$, while the leading term of $\langle -k_{0} ^2  \ROmegaPi\rangle (\phi_{d} )_{|\Omega }$ is $\phi_{d|\Omega } = \pm  r^{d} \sin (d\theta )$. However, when $d\in -\mathbb{N}^*$, we have $J_{d} =(-1)^{d} J_{-d}$ \cite{Leb72}, which does~not behave like the right power of $r$. In this case, a Bessel function of the second kind shall be used, see Appendix~\ref{sec: annexe des complements sur le DA}.
\end{textesansboite}

\begin{propositionsecable}{}
\label{expression sing homog CL}%
Let $d\in \TZ*$. We have $\langle -k_{0} ^2  \ROmegaPi\rangle (\phi_{d} )_{|\Omega } = \langle -k_{0} ^2  \ROmegaOmega\rangle (\phiOmega_{d} )$.\\
If moreover $d \not\in -\mathbb{N}^*$, then
\[\langle -k_{0} ^2  \ROmegaOmega\rangle (\phiOmega_{d} ) = \Gamma (d+1) \left( \frac{2}{k_{0}} \right)^{d} \BesselOmega_{d} .\]
Therefore, when $d \not\in -\mathbb{N}^*$, \[\langle -k_{0} ^2  \ROmegaOmega\rangle (\phiOmega_{d} ) \longrightarrow  (-1)^{d\Theta /\pi } \.\Gamma (d+1) \left( \frac{2}{k_{0}} \right)^{d} J_{d} (k_{0} r) \sin (d\theta ).\]\fbeq
\end{propositionsecable}

\begin{remarque}
\Cref{expression sing homog CL avec BesselY} gives a similar (but more complicated) expression of $\langle -k_{0} ^2  \ROmegaOmega\rangle (\phiOmega_{d} )$ when $d\in \TZ \cap (-\mathbb{N}^*)$. It involves the Bessel function of the second kind $Y_{d}$. Note however that $d\in \TZ* \cap  (-\mathbb{N}^*)$ may happen only if $\pst \in \mathbb{Q}$. In this case, there exist $a,b\in \mathbb{N}^*$ s.t.\ $\Theta = \frac{a}{b} \pi $ and $\pgcd{a,b}=1$. It~is then easy to check that $\TZ \cap  \mathbb{Z}=b\mathbb{Z}$, which implies $\TZ* \cap  (-\mathbb{N}^*) = -b\mathbb{N}^*$.
\end{remarque}

\begin{textesansboite}
To prove \cref{expression sing homog CL}, we could calculate in details of the terms $(-k_{0} ^2 \ROmegaPi)^{n} (\phi_{d} )$ of the corner singularity $\langle -k_{0} ^2  \ROmegaPi\rangle (\phi_{d} )$ and compare them with the terms of the series of the Bessel functions \cref{eq: def BesselJ} and \cref{eq: expression de BesselY_n}. But to avoid too many calculations, we rather use the Helmholtz problem in $\Amacro(\Omega )$ that both $\langle -k_{0} ^2  \ROmegaPi\rangle (\phi_{d} )_{|\Omega }$ and the formal series associated to the Bessel functions satisfy, because it has a unique solution.
\end{textesansboite}

\begin{demo}
\underline{Step 1:} Let $d\in \TZ*$. First, note that $\langle -k_{0} ^2 \ROmegaOmega\rangle (\phiOmega_{d} )$ is a solution in $\Amacro(\Omega )$ of 
\begin{equation}{\label{eq: preuve expression sing homog CL}}
    \left\{\begin{array}{r@{\;}ll}
    \Delta u + k_{0} ^2 u &= 0 &\text{ in }\Omega \\
    u &= 0 &\text{ on }\Gamma \\
    \sigma _{q} (u) &= \delta _{q,d} & \ \forall q\in \TZ*
    \end{array}\right.
\end{equation}
while $\langle -k_{0} ^2 \ROmegaPi\rangle (\phi_{d} ) \in \Amacro(\Pi )$ solves
\[\left\{\begin{array}{r@{\;}ll}
\Delta u + k_{0} ^2 u &= 0 &\text{ in }\Omega \\
\partial _{Y} ^2 u &= 0 &\text{ in }\couche\\
\partial _{Y} u_{|Y=0^{-}} &= 0 &\text{ on }\Gamma \\
\sigma _{q} (u) &= \delta _{q,d} & \ \forall q\in \TZ*
\end{array}\right.\]
In~particular, $\langle -k_{0} ^2 \ROmegaPi\rangle (\phi_{d} )_{|\Omega } \in \Amacro(\Omega )$ is also solution of \cref{eq: preuve expression sing homog CL}. Since \cref{eq: preuve expression sing homog CL} has a unique solution in $\Amacro(\Omega )$ by \cref{pb de Helmholtz dans A Omega}, we have $\langle -k_{0} ^2  \ROmegaPi\rangle (\phi_{d} )_{|\Omega } = \langle -k_{0} ^2  \ROmegaOmega\rangle (\phiOmega_{d} )$.\\
Now let us assume $d\not\in -\mathbb{N}^*$. By~\cref{relations sur BesselGamma/Omega,sigma de BesselOmega}, $\Gamma (d+1) ( \frac{2}{k_{0}} )^{d} \BesselOmega_{d}$ is also a solution of \cref{eq: preuve expression sing homog CL}. Therefore it~is equal to $\langle -k_{0} ^2 \ROmegaOmega\rangle (\phiOmega_{d} )$. Moreover we have $\BesselOmega_{d} \longrightarrow  J_{a} (k_{0} r) \sin (d(\theta -\Theta )) = (-1)^{d\Theta /\pi } J_{a} (k_{0} r) \sin (d\theta )$, which explains the presence of $(-1)^{d\Theta /\pi }$.
\end{demo}

\begin{textesansboite}
Now let us define the singularities of the higher order far fields.
\end{textesansboite}

\begin{definition}{far field singularities $\mcla_{d,n}$}
\label{def: mcl}%
Let $(T_{n} )_{n\in \mathbb{N}}$ be the sequence of Taylor coefficients of the tangent: $\forall t\in (- \frac{\pi }{2} , \frac{\pi }{2} ),\; \tan  t = \sum_{n=0}^\infty  T_{n} \.t^{2n+1}$. For~any $d\in \TZ*$ and $n\in \mathbb{N}^*$, we define, in $\Amaicro(\Omega )$, $\mcla_{d,-n} :=0$ and then inductively:
\[\mcla_{d,0} := \langle -k_{0} ^2 \ROmegaOmega\rangle (\phiOmega_{d} ) \quad  \text{ and } \quad  
\mcla_{d,n} := \langle -k_{0} ^2 \ROmegaOmega\rangle \circ \RGammaOmega\bigg( \sum_{m=0}^\infty  \frac{\mu _{0} }{\mu _{1}} \,T_{m} \cdot  (\partial _{x} ^2 +k_{1} ^2 )^{m} \partial _{y} \mcla_{d,n-1-2m|\Gamma } \bigg).\]
Note that the sum over $m$ has only a finite number of non-zero terms.
\end{definition}
\nomenclature[56-13-23]{$\mcla_{d,n}$}{Far field singularity (as a formal series)}

\begin{remarque}
By~\cref{pb de Helmholtz dans A Omega}, $\mcla_{d,0}$ is the unique solution in $\Amacro(\Omega )$ of
\begin{equation}{\label{eq: pb verifie par mcl qd n=0}}
    \left\{\begin{array}{r@{\;=\;}ll}
    (\Delta +k_{0} ^2 )\mcla_{d,0} & 0 &\text{ in }\Omega  \\[3pt]
    \mcla_{d,0} & 0 &\text{ on }\Gamma  \\[3pt]
    \sigma _{q} (\mcla_{d,0} ) & \delta _{d,q} &\ \forall q\in \TZ*
    \end{array}\right.
\end{equation}
and, for~any $n\in \mathbb{N}^*$, $\mcla_{d,n}$ is the unique solution in $\Amacro(\Omega )$ of
\begin{equation}{\label{eq: pb verifie par mcl qd n>0}}
    \left\{\begin{array}{r@{\;=\;}ll}
    (\Delta +k_{0} ^2 )\mcla_{d,n} & 0 &\text{ in }\Omega  \\[2pt]
    \mcla_{d,n} & \sum\limits_{m=0}^\infty  \frac{\mu _{0} }{\mu _{1}} \,T_{m} \cdot  (\partial _{x} ^2 +k_{1} ^2 )^{m} \partial _{y} \mcla_{d,n-1-2m} &\text{ on }\Gamma  \vspace{2pt}\\
    \sigma _{q} (\mcla_{d,n} ) & 0 &\ \forall q\in \TZ*
    \end{array}\right.
\end{equation}
\end{remarque}

\begin{proposition}{}
\label{expression uDA avec mcl}%
For~any $p\in \NN$ and $\ell \in \mathbb{N}$, we have:\vspace{-0.5ex}
\begin{equation}{\label{eq: expression uDA avec mcl}}
    \uDA_{p,\ell } := \bfuDA_{p,\ell |\Omega } = \sum _{n=0} ^{\lfloor {p}\rfloor } \fsum_{\substack{d\in \TZ*\\ d\in n-p+\NN}} \sigma _{d} (\uDA_{p-n,\ell } ) \,\mcla_{d,n} .
\end{equation}\fbeq
\end{proposition}

\begin{remarque}
In its structure, this relation is similar to \cref{eq: expression de Spl avec les correcteurs}, the expression of the corner fields in terms of the corner profiles.
\end{remarque}

\begin{demo}
\underline{Step 1:} Let us show that, for any $p\in \NN$ and $\ell \in \mathbb{N}$, we have:
\begin{equation}{\label{eq: expression uDA avec mcl preuve 1}}
    \left\{\begin{array}{r@{\;=\;}ll}
    (\Delta +k_{0} ^2 )\uDA_{p,\ell } & 0 &\text{ in }\Omega  \\
    \uDA_{p,\ell } & \sum\limits_{n=0}^\infty  \frac{\mu _{0} }{\mu _{1}} \,T_{n} \cdot  (\partial _{x} ^2 +k_{1} ^2 )^{n} \partial _{y} \uDA_{p-1-2n,\ell } &\text{ on }\Gamma 
    \end{array}\right.
\end{equation}
On the one hand, by~\cref{validation construction et DA des champs,pb de Helmholtz dans A Pi}, we have $(\Delta +k_{0} ^2 )\uDA_{p,\ell } =(\Delta +k_{0} ^2 )\bfuDA_{p,\ell |\Omega } =0$. On the other hand, \cref{construction directe des upl} states that we have
\[u_{p,\ell } = \sum_{n=0}^\infty  \frac{\mu _{0} }{\mu _{1}} \,T_{n} \cdot  (\partial _{x} ^2 +k_{1} ^2 )^{n} \partial _{y} u_{p-1-2n,\ell } \quad \text{ on }\Gamma .\]
In addition, by \cref{validation construction et DA des champs}, we have $\bfu_{p',\ell } = T_{\leqslant q} (\bfuDA_{p',\ell } )+o_{\partial } (\bfr^{q} )$ for~any $p'\in \NN$ and $q\in \mathbb{R}$ ($\bfr$ is defined in \cref{eq: def bfu bfr} and the truncation operator $T_{\leqslant q}$ in introduced in \cref{def: troncature de SF}). Restricting this to $\Omega $ gives $u_{p',\ell } = T_{\leqslant q} (\uDA_{p',\ell } )+o_{\partial } (r^{q} )$. Thus by \cref{derivation des DA dans Amaicro}, we deduce
\[\forall q\in \mathbb{R}, \qquad  u_{p,\ell } = T_{\leqslant q} \bigg(\sum_{n=0}^\infty  \frac{\mu _{0} }{\mu _{1}} \,T_{n} \cdot  (\partial _{x} ^2 +k_{1} ^2 )^{n} \partial _{y|\Gamma } \uDA_{p-1-2n,\ell } \bigg) +o_{\partial } (r^{q} ).\]
By~uniqueness of the asymptotic expansion in $\Amacro(\Omega )$, see \cref{unicité du DA}, we deduce the second line of \cref{eq: expression uDA avec mcl preuve 1}.\smallbreak

\underline{Step 2:} Let us prove \cref{expression uDA avec mcl} by induction on $p$. We can initialize for $p<0$, which is trivial because by convention $u_{p,\ell }$ vanishes for~any $p<0$. Now let us assume that the property holds for~any $p'<p$ and let us show it at rank $p$. By~\cref{pb de Helmholtz dans A Omega}, \cref{eq: expression uDA avec mcl preuve 1} implies
\begin{align*}
\uDA_{p,\ell } &= \langle -k_{0} ^2 \ROmegaOmega\rangle  \Bigg( \RGammaOmega \bigg( \sum_{m=0}^\infty  \frac{\mu _{0} }{\mu _{1}} \,T_{m} \cdot  (\partial _{x} ^2 +k_{1} ^2 )^{m} \partial _{y} \uDA_{p-1-2m,\ell |\Gamma } \bigg)+ \fsum_{d\in \TZ*} \sigma _{d} (\uDA_{p,\ell } )\,\phiOmega_{d} \Bigg) \\
&= \langle -k_{0} ^2 \ROmegaOmega\rangle  \circ  \RGammaOmega \bigg( \sum_{m=0}^\infty  \frac{\mu _{0} }{\mu _{1}} \,T_{m} \cdot  (\partial _{x} ^2 +k_{1} ^2 )^{m} \partial _{y} \uDA_{p-1-2m,\ell |\Gamma } \bigg) + \fsum_{d\in \TZ*} \sigma _{d} (\uDA_{p,\ell } )\,\mcla_{d,0}
\end{align*}
Then the induction hypothesis yields:
\begin{align*}
\uDA_{p,\ell } &= \langle -k_{0} ^2 \ROmegaOmega\rangle  \circ  \RGammaOmega \bigg(\sum_{m=0}^\infty  \frac{\mu _{0} }{\mu _{1}} \,T_{m} \cdot  (\partial _{x} ^2 +k_{1} ^2 )^{m} \partial _{y} \sum_{n=0}^\infty  \fsum_{d\in \TZ*} \sigma _{d} (\uDA_{p-1-2m-n,\ell } ) \,\mcla_{d,n|\Gamma } \bigg) + \fsum_{d} \sigma _{d} (\uDA_{p,\ell } )\,\mcla_{d,0} \\
&= \langle -k_{0} ^2 \ROmegaOmega\rangle  \circ  \RGammaOmega \bigg(\sum_{m=0}^\infty  \frac{\mu _{0} }{\mu _{1}} \,T_{m} \cdot  (\partial _{x} ^2 +k_{1} ^2 )^{m} \partial _{y} \sum _{N=1} ^{\infty } \fsum_{d\in \TZ*} \sigma _{d} (\uDA_{p-N,\ell } ) \,\mcla_{d,N-1-2m|\Gamma } \bigg) + \fsum_{d} \sigma _{d} (\uDA_{p,\ell } )\,\mcla_{d,0} \\
& \qquad  \text{with }  N:=1+2m+n\\
&= \sum _{N=1} ^{\infty } \fsum_{d\in \TZ*}  \sigma _{d} (\uDA_{p-N,\ell } ) \cdot \langle -k_{0} ^2 \ROmegaOmega\rangle  \circ  \RGammaOmega \bigg( \sum_{m=0}^\infty  \frac{\mu _{0} }{\mu _{1}} \,T_{m} \cdot  (\partial _{x} ^2 +k_{1} ^2 )^{m} \partial _{y} \mcla_{d,N-1-2m|\Gamma } \bigg) + \fsum_{d} \sigma _{d} (\uDA_{p,\ell } )\,\mcla_{d,0} \\
&= \sum _{N=1} ^{\infty } \fsum_{d\in \TZ*}  \sigma _{d} (\uDA_{p-N,\ell } )\,\mcla_{d,N} + \fsum_{d\in \TZ*} \sigma _{d} (\uDA_{p,\ell } )\,\mcla_{d,0} \\
\end{align*}
\end{demo}

\begin{proposition}{expression of $\mcla_{d,n}$ using Bessel functions}
\label{expression sing heterog CL}%
\begin{enumerate}
    \item Let us assume that $\pst$ is irrational. For~any $n\in \mathbb{N}$, there exists a sequence $(a_{n,m} )_{m} \in \mathbb{C}^{[\![-n,n]\!]}$ s.t.\ for~any $d\in \TZ*$, we have
    \[\mcla_{d,n} = \Gamma (d+1) \left( \frac{2}{k_{0}} \right)^{d} \sum _{m=-n} ^{n} a_{n,m} \, \BesselOmega_{d+m} .\]
    Moreover we have $(n>0 \Rightarrow  a_{n,0} =0)$ and, for~any $m$, $a_{n,-m} = (-1)^{n} a_{n,m}$ and $(a_{n,m} \neq 0 \Rightarrow  m \equiv  n\,[2])$.
    \item Let us assume that $\pst$ is rational. Let $a,b\in \mathbb{N}^*$ be s.t.\ $\Theta  = \frac{a}{b} \pi $ and $\pgcd{a,b}=1$. Then the statement of point 1 still holds as long as $n\in [\![0,b-1]\!]$ and $d\not\in -\mathbb{N}$.
\end{enumerate}
\end{proposition}
\nomenclature[11]{$b$}{Denominator of $\frac{\Theta }{\pi }$ when $\Theta  \in \pi \mathbb{Q}$ (in general, depends on context)}

\begin{demo}
Let us prove it by induction on $n$. The case $n=0$ follows from \cref{expression sing homog CL}. Let $n\geqslant 1$ and let us assume that the property is true at the ranks $n'<n$.\\
For~any $(a_{m} )\in \mathbb{C}^{\mathbb{Z}}$ and $s\in \{-1,1\}$, we say that $(a_{m} )$ is $s$-symmetric iff: $\forall m\in \mathbb{Z},\; a_{-m} = s\,a_{m}$. Let $\iota : (a_{m} )\in \mathbb{C}^{\mathbb{Z}} \mapsto  (a_{-m} )$ and $\tau : (a_{m} )\in \mathbb{C}^{\mathbb{Z}} \mapsto  (a_{m+1} )$. We have $\iota  \tau =\tau ^{-1} \iota $ and $\iota  \tau ^{-1} =\tau  \iota $, so $\iota (\tau +\tau ^{-1} )=(\tau +\tau ^{-1} )\iota $ and $\iota (\tau ^2 -2+\tau ^{-2} )=(\tau ^2 -2+\tau ^{-2} )\iota $. Therefore, for~any $s\in \{-1,1\}$ and $a\in \mathbb{C}^{\mathbb{Z}}$, if $a$ is $s$-symmetric, then so are $(\tau +\tau ^{-1} )(a)$ and $(\tau ^2 -2+\tau ^{-2} )(a)$.\\
Now, \cref{relations sur BesselGamma/Omega} shows that for~any $(a_{m} )\in \mathbb{C}^{\mathbb{Z}}$ with finite support, we have:
\begin{itemize}
    \item $\partial _{y} \big( \sum _{m} a_{m} \BesselOmega_{d+m} \big)_{|\Gamma } = \sum _{m} b_{m} \BesselGamma_{d+m}$ with $(b_{m} ) := (-1)^{d\Theta /\pi } \frac{k_{0} }{2} (\tau +\tau ^{-1} )\big((\cos (m\Theta )a_{m} )_{m} \big)$
    \item $\partial _{x} ^2  \big( \sum _{m} a_{m} \BesselGamma_{d+m} \big) = \sum _{m} b'_{m} \BesselGamma_{d+m}$ with $(b'_{m} ) := \big( \frac{k_{0} }{2} \big)^2  (\tau -\tau ^{-1} )^2 \big((a_{m} )_{m} \big) = \big( \frac{k_{0} }{2} \big)^2  (\tau ^2 -2+\tau ^{-2} )\big((a_{m} ) \big)$
\end{itemize}
If $(a_{m} )$ is $s$-symmetric, then so are $(b_{m} )$ and $(b'_{m} )$. Thus, the induction hypothesis shows that there exists a $(-1)^{n-1}$-symmetric sequence $(c_{n,m} )\in \mathbb{C}^{\mathbb{Z}}$, independent of $d$, s.t.\ \[\mcla_{d,n|\Gamma } := \sum_{j=0}^\infty  \frac{\mu _{0} }{\mu _{1}} \,T_{j} \cdot  (\partial _{x} ^2 +k_{1} ^2 )^{j} \partial _{y} \mcla_{d,n-1-2j|\Gamma } = \Gamma (d+1) \left( \frac{2}{k_{0}} \right)^{d} (-1)^{d\Theta /\pi } \sum _{m} c_{n,m} \,\BesselGamma_{d+m} .\]
Moreover, 
the support of $(c_{n,m} )$ is included in $[\![-n,n]\!]$ (because applying $\partial _{y|\Gamma }$, resp.\ $\partial _{x} ^2 $, makes the support 1 step larger, resp.\ 2) and we have: $\forall m,\; c_{n,m} \neq 0 \Rightarrow  m\equiv n\,[2]$ (because applying $\partial _{y|\Gamma }$ moves the support from even integers to odd and conversely, while $\partial _{x} ^2 $ does~not change the parity of the support). Now let
\[a_{n,0} :=0 \quad  \text{ and } \quad  \forall m\in \mathbb{Z}^*,\,\ a_{n,m} :=- \frac{c_{n,m} }{\sin (m\Theta )} \qquad  \text{ and } \qquad  \varphi := \Gamma (d+1) \left( \frac{2}{k_{0}} \right)^{d} \sum _{m=-n} ^{n} a_{n,m} \,\BesselOmega_{d+m} .\]
If $\pst \in \mathbb{R} \setminus \mathbb{Q}$ or $(\pst \in \mathbb{Q} \text{ and } n\in [\![0,b-1]\!])$, then none of the above sinus vanishes, so $(a_{n,m} )$ and $\varphi $ are well defined. Let us show that $\varphi =\mcla_{d,n}$. First, note that $c_{n,0} =0$ (by $(-1)$-symmetry of $(c_{m} )$ if $n$ is even, and by oddness of the support of $(c_{n,m} )$ if $n$ is odd). In addition, $(a_{n,m} )$ is independent of $d$ and $(-1)^{n}$-symmetric and $(\forall m,\; a_{n,m} \neq 0 \Rightarrow  m\equiv n\,[2])$ by properties of $(c_{n,m} )$. Finally, let us show that $\varphi $ is a solution of \cref{eq: pb verifie par mcl qd n>0}, the problem solved by $\mcla_{d,n}$. Since this problem has a unique solution (see \cref{pb de Helmholtz dans A Omega}), it will imply that $\mcla_{d,n} =\varphi $.
\begin{itemize}
    \item $(\Delta +k_{0} ^2 )\varphi  = 0$ by \cref{relations sur BesselGamma/Omega}.
    \item \Cref{relations sur BesselGamma/Omega} shows $\varphi _{|\Gamma } = \Gamma (d+1) ( \frac{2}{k_{0}} )^{d} \sum\limits _{m=-n} ^{n} a_{n,m} \,\sin ((d+m)\Theta )\,\BesselGamma_{d+m}$. So, using $c_{n,0} =0$, we deduce
    \[\varphi _{|\Gamma } = \Gamma (d+1) \left( \frac{2}{k_{0}} \right)^{d} (-1)^{d\Theta /\pi } \sum\limits _{m\in \mathbb{Z}} c_{n,m} \,\BesselGamma_{d+m} =\mcla_{d,n|\Gamma } .\]
    \item Let $q\in \TZ*$ and let us show that $\sigma _{q} (\varphi ) = 0$. $\varphi $ is a linear combination of the $\BesselOmega_{d+m}$ with $m\in [\![-n,n]\!] \setminus \{0\}$ (because $a_{n,0} =0$), so it suffices to prove $\sigma _{q} (\BesselOmega_{d+m} )=0$ for~any such $m$. Let $m\in [\![-n,n]\!] \setminus \{0\}$. 
    \begin{itemize}
        \item If $\pst \not\in \mathbb{Q}$, then $\mathbb{Z}\cap \TZ=\{0\}$, so $m\not\in \TZ$.
        \item Otherwise, $\pst \in \mathbb{Q}$ and it~is easy to check that $\mathbb{Z}\cap \TZ =b\mathbb{Z}$. By~assumption $n<b$, so $0<|m|<b$, which yields $m\not\in \TZ$.
    \end{itemize}
    In both cases, we obtain $m\not\in \TZ$, which yields $d+m\not\in \TZ$. Now if $d+m \not\in -\mathbb{N}^*$, we can apply \cref{sigma de BesselOmega}, which gives $\sigma _{q} (\BesselGamma_{d+m} )=0$ because $d+m\neq q$. Otherwise, $d+m\in -\mathbb{N}^*$, so $\BesselGamma_{d+m} =(-1)^{d+m} \BesselGamma_{-d-m}$ by \cref{relations sur BesselGamma/Omega}, item 5. In this case, \cref{sigma de BesselOmega} gives again $\sigma _{q} (\BesselGamma_{-d-m} )=0$.
\end{itemize}
\end{demo}

\begin{remarque}
The additional properties of the sequence $(a_{n,m} )$ stated in \cref{expression sing heterog CL} are a key ingredient of the proof. Indeed, the above proof shows by induction that we have $\mcla_{d,n|\Gamma } =\sum _{m=-n} ^{n} \tilde{a}_{d,n,m} \BesselGamma_{d+m}$ for some coefficients $\tilde{a}_{d,n,m} \in \mathbb{C}$. Then \cref{relations sur BesselGamma/Omega} states that we have $\BesselOmega_{a|\Gamma } = -\sin (a\Theta ) \BesselGamma_{a}$ for~any $a$, so each term $\BesselGamma_{d+m}$ in $\mcla_{d,n|\Gamma }$ can be lifted into $\BesselOmega_{d+m}$ in $\Omega $ provided that $\sin ((d+m)\Theta ) \neq 0$. In particular, one needs $m\neq 0$ as $d\in \TZ*$. The additional properties of $(a_{n,m} )$ allow us to prove by induction on $n$ that we have $\tilde{a}_{d,n,0} =0$, so $\BesselGamma_{d}$ does~not appear in the expression of $\mcla_{d,n|\Gamma }$ and we are able to express $\mcla_{d,n}$ with the $\BesselOmega_{d+m}$, $m\neq 0$.
\end{remarque}

\begin{exemplenumerote}
\label{ex: expression sing heterog CL}%
\newcommand{\rappmu}{\mu_{\rm r}}%
\newcommand{\rappk}{k_{\rm r}}%
\newcommand{\cotm}[1]{t_{#1}}%
Let $\rappmu:= \frac{\mu _{0} }{\mu _{1}}$, $\rappk := \big( \frac{2k_{1} }{k_{0}} \big)^2 $ and $\cotm{m} :=\operatorname{cotan}(m\Theta )$ for~any $m\in \mathbb{N}$. In the cases where \cref{expression sing heterog CL} holds (i.e.\ $\pst \in \mathbb{R} \setminus \mathbb{Q}$ or $(n<b \text{ and } d\not\in -\mathbb{N})$), we have:
\begin{itemize}
    \item $\displaystyle \mcla_{d,0} =  \Gamma (d+1)\. \left( \frac{k_{0} }{2} \right)^{-d} \,\BesselOmega_{d}$
    \item $\displaystyle \mcla_{d,1} =  \Gamma (d+1)\. \left( \frac{k_{0} }{2} \right)^{1-d} \frac{\rappmu}{\sin \Theta } \,(\BesselOmega_{d-1} - \BesselOmega_{d+1} )$
    \item $\displaystyle \mcla_{d,2} = \Gamma (d+1)\. \left( \frac{k_{0} }{2} \right)^{2-d} \frac{\rappmu^2  \cotm{1}}{\sin (2\Theta )} \,(\BesselOmega_{d-2} + \BesselOmega_{d+2} )$
    \item $\displaystyle \mcla_{d,3} = \Gamma (d+1)\. \left( \frac{k_{0} }{2} \right)^{3-d} \left[\frac{3\rappmu^{3} \cotm{1}\cotm{2} + \rappmu}{3\sin (3\Theta )} \,(\BesselOmega_{d-3} - \BesselOmega_{d+3} )
    \,+\, \frac{3\rappmu^{3} \cotm{1}\cotm{2} +\rappmu(\rappk-1)}{3\sin \Theta } \,(\BesselOmega_{d-1} - \BesselOmega_{d+1} )\right]$
    \item $\displaystyle \mcla_{d,4} = \Gamma (d+1)\. \left( \frac{k_{0} }{2} \right)^{4-d} \bigg[ \frac{3\rappmu^{4} \cotm{1}\cotm{2}\cotm{3} + \rappmu^2 (\cotm{1}+\cotm{3})}{3\sin (4\Theta )} \,(\BesselOmega_{d-4} + \BesselOmega_{d+4} )$\\
    \hbox to \linewidth{\hfill $\displaystyle +\; \frac{3\rappmu^{4} \cotm{1}\cotm{2}(\cotm{1}+\cotm{3}) + \rappmu^2  \cotm{1}(2\rappk-3)+\rappmu^2 \cotm{3} }{3\sin (2\Theta )} \,(\BesselOmega_{d-2} + \BesselOmega_{d+2} ) \bigg]$.}
\end{itemize}
To protect against calculation errors, these formulas have been verified numerically.
\end{exemplenumerote}

\begin{remarque}
\Cref{expression sing heterog CL avec BesselY} shows how to extend \cref{expression sing heterog CL} to the case $\pst \in \mathbb{Q}$ and $d\in -\mathbb{N}$, using Bessel functions of the second kind. As~for the case $\pst \in \mathbb{Q}$ and $n\geqslant b$, in general $\mcla_{d,n}$ involves powers of $\log (\alpha z)$ and thus cannot be expressed using only the series $\BesselOmega_{a}$. However, in some special cases like $b\in \{2,4\}$, the constraint $n<b$ is not optimal, see \cref{expression sing heterog CL b=2 ou 4}.
\end{remarque}

\section{Simplifications depending on the angle}
\label{sec: simplifications DA}

\let\varepsilon=\feps

\begin{textesansboite}
\Cref{sec: construction des champs} gives a general construction of the fields $u_{p,\ell }$, $U_{p,\ell }$ and $S_{p,\ell }$ for~any angle $\Theta $. It involves heavy formulas, like the definition of $\cfinal_{d,d',p,\ell }$ in \cref{eq: def cfinal}. However, when looking at specific angles or specific orders, we can simplify those results. In this section, we explain how. We first study angles s.t.\ $\pst\not\in \mathbb{Q}$, then angles s.t.\ $\pst \in \mathbb{Q}$. We show that in the first case no $\ln \peps$ appears, which simplifies a lot (\cref{sec: simplifications DA cas irrationnel}), and that the second case behaves like the first one for small orders (\cref{sec: simplifications DA cas rationnel}).\\

The following preliminary lemma will allow us to use the properties of the $\mcla_{d,n}$ from \cref{sec: expression des mcl} to study the matching coefficients $\cmacro_{d,d',p,\ell }$ and $\cmicro_{d,d',p,\ell }$ (whose definition involve $\langle \Rmaicro\rangle (\phi_{d'} )$, see \cref{eq: def cmaicro}).
\end{textesansboite}

\begin{remarque}
This section uses a lot the tools introduced in \cref{chap: SF,chap: raccord algebrique}: formal series and their properties (\cref{sec: SF 1D,sec: SF 2D}), the spaces $\Aemaicro$ (\cref{def: Aeps et Aemaicro}), the sets of operators $\Rmaicro$ (\cref{eq: def Rmaicro}), the resolution of the Helmholtz equation in $\Aemaicro$ (\cref{pb de Helmholtz dans Aemacro,pb de Helmholtz dans Aemicro}), the scaling operator $\Eps$ (\cref{def: Eps}), etc. We will stop justifying the usage of the most basic properties, like commutation properties for instance, as we expect that by now the reader has become more familiar with them. Besides, we recall that $\varepsilon $ and $\ln \varepsilon $ denote the formal indeterminates of the space $\Aeps$ and are~not real numbers.
\end{remarque}

\begin{lemmesansboite}{}
\label{expression de Rmacro phi avec les mcl}%
For~any $d\in \TZ*$, we have $\displaystyle \langle \Rmacro\rangle (\phi_{d} )_{|\Omega } = \fsum_{n=0} ^{\infty } \varepsilon ^{n} \.\mcla_{d,n}$.
\end{lemmesansboite}

\begin{demo}
Let $d\in \TZ*$ and $\Phi _{\varepsilon } :=\langle \Rmacro\rangle (\phi_{d} )$. For~any $R\in \Rmacro$, we have $\deg _{\varepsilon } R \in \mathbb{N}$, so there exists $(\Phi _{n} )_{n} \in \Amacro(\Pi )^{\mathbb{N}}$ s.t.\ $\Phi _{\varepsilon } = \fsum_{n=0} ^{\infty } \varepsilon ^{n} \.\Phi _{n}$. By~\cref{pb de Helmholtz dans Aemacro}, we have 
\[\left\{\begin{array}{r@{\;=\;}ll}
	\mu _{0} \Delta  \Phi _{\varepsilon } +\omega ^2 \rho _{0} \Phi _{\varepsilon } & 0 &\text{ in }\Omega  \\[0.5mm]
	\mu _{1} \partial _{Y} ^2  \Phi _{\varepsilon } + \varepsilon ^2 (\mu _{1} \partial _{x} ^2 +\omega ^2 \rho _{1} )\Phi _{\varepsilon } &0&\text{ in }\couche \\[0.5mm]
	\mu _{1} \partial _{Y|Y=0^{-}} \Phi _{\varepsilon } - \varepsilon \,\mu _{0} \,\partial _{y|y=0^{+}} \Phi _{\varepsilon } &0 &\text{ on }\Gamma 
\end{array}\right.\]
Identifying each power of $\varepsilon $ gives for~any $n\in \mathbb{N}$ (denoting $\Phi _{-1} :=\Phi _{-2} :=0$):
\[\left\{\begin{array}{r@{\;=\;}ll}
	\mu _{0} \Delta  \Phi _{n} +\omega ^2 \rho _{0} \Phi _{n} & 0 &\text{ in }\Omega  \\[0.5mm]
	\mu _{1} \partial _{Y} ^2  \Phi _{n} & -(\mu _{1} \partial _{x} ^2 +\omega ^2 \rho _{1} )\Phi _{n-2} &\text{ in }\couche \\[0.5mm]
	\mu _{1} \partial _{Y|Y=0^{-}} \Phi _{n} & \mu _{0} \,\partial _{y|y=0^{+}} \Phi _{n-1} &\text{ on }\Gamma 
\end{array}\right.\]
This system has the same structure as the equations satisfied by $\bfu_{n,0}$. Thus, similarly to \cref{expression de Upl} one can express $\Phi _{n|\couche}$ in terms of the $\partial _{y} \Phi _{n'|\Gamma ,y=0^{+}}$ with $n'<n$ and deduce:
\[\forall n\in \mathbb{N}, \qquad  \Phi _{n} = \sum_{m=0}^\infty  \frac{\mu _{0} }{\mu _{1}} \,T_{m} \cdot  (\partial _{x} ^2 +k_{1} ^2 )^{m} \partial _{y} \mcla_{d,n-1-2m|y=0^{+}} \quad \text{ on }\Gamma .\]
Moreover, by \cref{pb de Helmholtz dans Aemacro}, we have $\sigma _{q} (\Phi _{n} )=\sigma _{q,n,0} (\Phi _{\varepsilon } ) =\delta _{d,q} \,\delta _{n,0}$ for~any $q\in \TZ*$. Hence, $(\Phi _{n|\Omega } )$ satisfies the same equations \cref{eq: pb verifie par mcl qd n=0}--\cref{eq: pb verifie par mcl qd n>0} as $(\mcla_{d,n} )$. Since these equations have unique solution according to \cref{pb de Helmholtz dans A Omega}, we deduce: $\forall n\in \mathbb{N},\; \Phi _{n} =\mcla_{d,n}$.
\end{demo}

\subsection{The irrational case}
\label{sec: simplifications DA cas irrationnel}

\begin{textesansboite}
In this section, we assume that $\pst$ is irrational. In this case, the behaviors of the fields $u_{p,\ell }$, $U_{p,\ell }$ and $S_{p,\ell }$ around the corner do~not involve logarithms of $r$ or $x$. Indeed, \cref{expression sing heterog CL} item 1 and the definition of the $\BesselOmega_{a}$ ($a\in \mathbb{R}$) yield that, for~any $d\in \TZ*$ and $n\in \mathbb{N}$, $\mcla_{d,n}$ contains no $\log (\alpha z)$. By~\cref{expression uDA avec mcl}, the same holds for $\bfuDA_{p,\ell |\Omega }$. Using the expression of $U_{p,\ell }$ in terms of the $u_{p'\!,\ell }$ given by \cref{expression de Upl}, we can extend this fact to $\bfuDA_{p,\ell |\couche}$: it does~not involve $\ln x$. Then by the matching assumption $\fsum_{p,\ell } \varepsilon ^{p} \ln ^{\ell } \!\varepsilon \cdot \SDA_{p,\ell } = \Eps (\fsum_{p,\ell } \varepsilon ^{p} \ln ^{\ell } \!\varepsilon \cdot \bfuDA_{p,\ell } )$, the $\SDA_{p,\ell }$ also contain no $\log (\alpha z)$ or $\ln x$.\\

\Cref{valeur de cmaicro cas irrationnel,simplifications cas irrationnel} below show that this implies that the expansion of $u_{\peps}$ does~not involve $\ln \peps$. First, the lemma calculates the matching coefficients, then the theorem concludes using the matching conditions that involve these coefficients.
\end{textesansboite}

\begin{lemme}{}
\label{valeur de cmaicro cas irrationnel}%
Let us assume that $\pst$ is irrational. For~any $d,d'\in \TZ*$, $p\in \NN$ and $\ell \in \mathbb{N}$, we have $\cmacro_{d,d',p,\ell } = \delta _{d,d'} \,\delta _{p,-d} \,\delta _{\ell ,0}$ and $\cmicro_{d,d',p,\ell } = \delta _{d,d'} \,\delta _{p,d} \,\delta _{\ell ,0}$.\\
Moreover if $d<0$, $d'>0$ and $n\in \mathbb{N}$, we have $\cfinal_{d,d',d'-d+2n,0} = \omega ^{2n} \sigma _{d} (\ccDA_{d',n} )$. In all other cases, $\cfinal_{d,d',p,\ell }$ vanishes.
\end{lemme}

\begin{textesansboite}
The proof of \cref{valeur de cmaicro cas irrationnel} uses the expression of the $\mcla_{d,n}$ in terms of the Bessel functions (\cref{expression sing heterog CL}) and \cref{expression de Rmacro phi avec les mcl} to show that $\langle \Rmacro\rangle (\phi_{d'} )$ contains no $\log (\alpha z)$ or $\ln x$. Thus $\Eps\circ \langle \Rmacro\rangle (\phi_{d'} )$ contains no $\ln \varepsilon $, which allows us to compute the coefficients $\cmicro_{d,d',p,\ell } := \sigma _{d,p,\ell } \circ  \Eps \circ  \langle \Rmacro\rangle  (\phi_{d'} )$. As~for the coefficients $\cmacro_{d,d',p,\ell }$, they depend on $\Eps^{-1} \circ  \langle \Rmicro\rangle  (\phi_{d'} )$. 
Thanks to the tools introduced in \cref{chap: raccord algebrique}, the property of $\Eps\circ \langle \Rmacro\rangle (\phi_{d'} )$ of having no $\ln \varepsilon $ is easily transfered to its twin $\Eps^{-1} \circ  \langle \Rmicro\rangle  (\phi_{d'} )$. We especially use the resolution of scaled and non-scaled Helmholtz problems in $\Aemaicro(\Pi )$ (\cref{pb de Helmholtz dans Aemacro,pb de Helmholtz dans Aemicro}) and the action of the scaling operator $\Eps$ on these problems (\cref{action de Eps sur Helmholtz dans Aemaicro}). Finally, $\cfinal$ is computed thanks to its formula.
\end{textesansboite}

\begin{demo}
\underline{Step 1:} Let us first treat $\cmicro_{d,d',p,\ell }$. We have $\cmicro_{d,d',p,\ell } := \sigma _{d,p,\ell } \circ  \Eps \circ  \langle \Rmacro\rangle  (\phi_{d'} )$, where $ \langle \Rmacro\rangle  (\phi_{d'} ) =\fsum_{n=0} ^{\infty } \varepsilon ^{n} \mcla_{d',n} $ by~\cref{expression de Rmacro phi avec les mcl}. Thanks to \cref{expression sing heterog CL} and the fact that $\pst$ is irrational, we know that, for~any $n\in \mathbb{N}$, $\mcla_{d',n}$ is a linear combination of the $\BesselOmega_{a}$, $a\in \mathbb{R}$. So, by \cref{eq: def BesselOmega}, $\mcla_{d',n}$ is a formal series of terms of the form $c\,\Im[(\alpha z)^{q} \overline{\alpha z}^{k} ]$ with $c\in \mathbb{C}$. In other words, $\mcla_{d',n}$ does~not involve $\log (\alpha z)$. Now, by~definition of $\Eps$ in \cref{def: Eps}, we have for~any $(q,k,n)\in \mathbb{R}\times \mathbb{N}^2 $:
\[\Eps\big(\varepsilon ^{n} \Im[(\alpha z)^{q} \,\overline{\alpha z}^{k} ]\big) = \varepsilon ^{n+q+k} \,\Im[(\alpha z)^{q} \,\overline{\alpha z}^{k} ].\]
By~definition of $\sigma _{d,p,\ell }$ in \cref{def: sigma sur Aeps}, this implies:
\[\sigma _{d,p,\ell } \circ \Eps\big(\varepsilon ^{n} \Im[(\alpha z)^{q} \,\overline{\alpha z}^{k} ]\big) = \delta _{d,q+k} \,\delta _{p,n+q+k} \,\delta _{\ell ,0} \cdot \sigma _{d} \big(\Im[(\alpha z)^{q} \,\overline{\alpha z}^{k} ]\big) = \sigma _{d,p-d,\ell } \big(\varepsilon ^{n} \Im[(\alpha z)^{q} \,\overline{\alpha z}^{k} ]\big) .\]
As we have seen, $\langle \Rmacro\rangle  (\phi_{d'} )$ is a formal series of terms of the form $c\,\varepsilon ^{n} \.\Im[(\alpha z)^{q} \,\overline{\alpha z}^{k} ]$ with $c\in \mathbb{C}$. So we deduce
\begin{equation}{\label{eq: valeur de cfinal cas irrationnel preuve EpsRmacro}}
    \cmicro_{d,d',p,\ell } := \sigma _{d,p,\ell } \circ  \Eps \circ  \langle \Rmacro\rangle  (\phi_{d'} ) =\sigma _{d,p-d,\ell } \circ  \langle \Rmacro\rangle  (\phi_{d'} ) = \sigma _{d,p-d,\ell } (\phi_{d'} ) =\delta _{d,d'} \,\delta _{p,d} \,\delta _{\ell ,0}
\end{equation}
because we have $\sigma _{d} \circ R=0$ for any $R\in \{\ROmegaPi,\RGammaPi,\RcouchePi\}$, so $\sigma _{d,p-d,\ell } \circ R' =0$ for~any $R'\in \Rmacro$, see \cref{operateurs R,eq: def Rmaicro}.\medskip\smallskip

\underline{Step 2:} Let us treat $\cmacro_{d,d',p,\ell } :=\sigma _{d,p,\ell } \circ  \Eps^{-1} \circ  \langle \Rmicro\rangle  (\phi_{d'} )$. By~\cref{eq: valeur de cfinal cas irrationnel preuve EpsRmacro}, we have $\sigma _{d,p,\ell } \big(\varepsilon ^{-d'} \langle \Rmacro\rangle  (\phi_{d'} )\big)= \sigma _{d,p+d',\ell } \circ  \langle \Rmacro\rangle  (\phi_{d'} ) = \delta _{d,d'} \,\delta _{p+d',0} \,\delta _{\ell ,0} = \delta _{d,d'} \,\delta _{p,-d} \,\delta _{\ell ,0}$. Hence, if we show
\begin{equation}{\label{eq: valeur de cfinal cas irrationnel eq Aemacro}}
    \Eps^{-1} \circ  \langle \Rmicro\rangle  (\phi_{d'} ) = \varepsilon ^{-d'} \langle \Rmacro\rangle  (\phi_{d'} ),
\end{equation}
then this implies $\cmacro_{d,d',p,\ell } = \delta _{d,d'} \,\delta _{p,-d} \,\delta _{\ell ,0}$ and concludes this step. Proving \cref{eq: valeur de cfinal cas irrationnel eq Aemacro} is equivalent to proving
\begin{equation}{\label{eq: valeur de cfinal cas irrationnel eq Aemicro}}
    \langle \Rmicro\rangle  (\phi_{d'} ) = \varepsilon ^{-d'} \Eps \circ \langle \Rmacro\rangle  (\phi_{d'} ).
\end{equation}
Let us apply \cref{pb de Helmholtz dans Aemicro} (which shows the uniqueness of the solution of a scaled Helmholtz problem in $\Aemicro(\Pi )$).
\begin{itemize}
    \item First, both sides of \cref{eq: valeur de cfinal cas irrationnel eq Aemicro} belong to $\Aemicro(\Pi )$ by \cref{Rmaicro defini et inversible dans Aemaicro,rq: def Eps}.
    \item Secondly, $\langle \Rmicro\rangle  (\phi_{d'} )$ is a solution of the scaled Helmholtz problem \cref{eq: pb de Helmholtz dans Aemicro bis} by \cref{pb de Helmholtz dans Aemicro}. In addition, $\langle \Rmacro\rangle  (\phi_{d'} )$ is solution of the non-scaled Helmholtz problem \cref{eq: pb de Helmholtz dans Aemacro bis} by \cref{pb de Helmholtz dans Aemacro}, so $\varepsilon ^{-d'} \Eps \circ \langle \Rmacro\rangle  (\phi_{d'} )$ is also solution of \cref{eq: pb de Helmholtz dans Aemicro bis} by \cref{action de Eps sur Helmholtz dans Aemaicro}.
    \item Thirdly, for any $(d,p,\ell ) \in \TZ*\times \mathbb{R}\times \mathbb{N}$, we have $\sigma _{d,p,\ell } \big(\langle \Rmicro\rangle  (\phi_{d'} )\big) =\sigma _{d,p,\ell } (\phi_{d'} ) =\delta _{d,d'} \,\delta _{p,0} \, \delta _{\ell ,0}$ (similarly to \cref{eq: valeur de cfinal cas irrationnel preuve EpsRmacro}) and $\sigma _{d,p,\ell } \big(\varepsilon ^{-d'} \Eps \circ \langle \Rmacro\rangle  (\phi_{d'} )\big) =\delta _{d,d'} \,\delta _{p,0} \, \delta _{\ell ,0}$ by \cref{eq: valeur de cfinal cas irrationnel preuve EpsRmacro}. So $\sigma _{d,p,\ell } \big(\langle \Rmicro\rangle  (\phi_{d'} )\big)$ and $\sigma _{d,p,\ell } \big(\varepsilon ^{-d'} \Eps \circ \langle \Rmacro\rangle  (\phi_{d'} )\big)$ coincide.
\end{itemize}
By~\cref{pb de Helmholtz dans Aemicro}, we deduce that \cref{eq: valeur de cfinal cas irrationnel eq Aemicro} holds.\medskip\smallskip

\underline{Step 3:} Let us treat $\cfinal_{d,d',p,\ell }$. We recall that it~is defined only for $d\in -\TN$ and that:
\[\cfinal_{d,d',p,\ell } := \! \sum _{\substack{(p_{1} ,p_{2} )\in (\mathbb{N}-d)\times (\mathbb{N}+d')\\ n\in \mathbb{N},\, p_{1} +p_{2} +2n=p}}
\ \sum _{\substack{d_{1} \in \TZ*\\ p_{1} +d_{1} \in \mathbb{N}}} \ \sum _{\substack{d_{2} \in \TN\\ p_{2} -d_{2} \in \mathbb{N}}}
\ \sum _{\substack{(\ell _{1} ,\ell _{2} )\in \mathbb{N}^2 \\ \ell _{1} +\ell _{2} =\ell }}
\cmacro_{d,d_{1} ,p_{1} ,\ell _{1} }\cdot \omega ^{2n} \sigma _{d_{1}} (\mathcal{S} ^{\infty } _{d_{2} ,n} )\cdot \cmicro_{d_{2} ,d'\!,p_{2} ,\ell _{2} } .\]
By~steps~1 and 2, the only non-zero terms are obtained for $d_{1} =d$, $d_{2} =d'$, $p_{1} =-d$, $p_{2} =d'$, $\ell _{1} =\ell _{2} =0$ and $\cmacro_{d,d_{1} ,p_{1} ,\ell _{1} }=\cmicro_{d_{2} ,d'\!,p_{2} ,\ell _{2} }=1$. One can deduce that if $d'<0$ then $\cfinal_{d,d',p,\ell } =0$, and otherwise the formula reduces to:
\[\cfinal_{d,d',p,\ell } := \delta _{\ell ,0} \! \sum _{\substack{n\in \mathbb{N},\, d'-d+2n=p}}\omega ^{2n} \sigma _{d} (\mathcal{S} ^{\infty } _{d',n} ).\]
Finally, this allows to prove the last result of the lemma.
\end{demo}

\begin{theoreme}{}
\label{simplifications cas irrationnel}%
Let us assume that $\pst$ is irrational. Then the asymptotic expansion of $u_{\peps}$ involves no power $\ln \peps$: $\forall p\in \NN,\forall \ell \in \mathbb{N}^*,\; \bfu_{p,\ell } =0$ and $S_{p,\ell } =0$.\\
Moreover, for~any $p\in \NN$, let us denote $u_{p} := u_{p,0}$, $\bfuDA_{p} :=\bfuDA_{p,0} $ and $\SDA_{p} :=\SDA_{p,0}$. The matching conditions simplify as follows:
\begin{itemize}
    \item $\forall p\in \NN,\; \forall d\in \TZ*,\ \sigma _{d} (\bfuDA_{p} )=\sigma _{d} (\SDA_{p+d} )$
    \item $\displaystyle \forall p\in \NN,\; {\textstyle \forall d\in -\TN},\ \sigma _{d} (u_{p} )= \sum _{\substack{n\in \mathbb{N}\\ 2n\leqslant p+d}} \sum _{\substack{d'\in \TN\\ d'\leqslant p+d-2n}} \omega ^{2n} \sigma _{d} (\ccDA_{d',n} )\cdot \sigma _{d'} (\uDA_{p-d'+d-2n} )$ \hfill\inlineEqNb\label{eq: simplification des sigma de CL cas irrationnel}
\end{itemize}\fbi
\end{theoreme}

\begin{demo}
Applying \cref{valeur de cmaicro cas irrationnel} to the matching conditions in \cref{def: construction des termes}, one can deduce $\sigma _{d} (\bfuDA_{p,\ell } )=\sigma _{d} (\SDA_{p+d,\ell } )$ for~any $(p,\ell ,d)$. Moreover, \cref{construction directe des upl} gives
\begin{align*}
    \sigma _{d} (u_{p,\ell } ) &= \!\! \sum _{\substack{p'\in \NN\\ p-p'\in \NN+ 2\pst}} \; \sum _{\substack{d'\in \TZ*\\ p'+d'\in \intvl{0,p+d}}} \! \sum _{\ell '=0} ^{\ell } \cfinal_{d,d',p-p',\ell -\ell '} \cdot \sigma _{d'} (\uDA_{p',\ell '} ) \\
    &= \!\! \sum _{\substack{p'\in \NN\\ p-p'\in \NN+ 2\pst}} \; \sum _{\substack{d'\in \TZ*\\ p'+d'\in \intvl{0,p+d}}} \! \cfinal_{d,d',p-p',0} \cdot \sigma _{d'} (\uDA_{p',0} ) \\
    &= \!\! \sum _{\substack{p'\in \NN\\ p-p'\in \NN+ 2\pst}} \; \sum _{\substack{d'\in \TN\\ p'+d'\in \intvl{0,p+d}}} \sum _{\substack{n\in \mathbb{N}\\ p-p'=d'-d+2n}} \! \cfinal_{d,d',d'-d+2n,0} \cdot \sigma _{d'} (\uDA_{p',\ell } )
\end{align*}
where we have used \cref{valeur de cmaicro cas irrationnel} to establish the last two equalities. Let us rewrite the sum conditions on $p'$ w.r.t.\ $n$. $p'\in \NN$ means $p-d'+d-2n\in \NN$, $p-p'\in \NN+ 2\pst$ means $d'-d-2\pst+2n \in \NN$ (automatically satisfied), and $p'+d'\in \intvl{0,p+d}$ means $p+d-2n\in \intvl{0,p+d}$, which is equivalent to $p+d-2n\in \NN$. So, potentially adding a few vanishing terms, we deduce
\[\sigma _{d} (u_{p,\ell } )= \sum _{\substack{n\in \mathbb{N}\\ 2n\leqslant p+d}} \sum _{\substack{d'\in \TN\\ d'\leqslant p+d-2n}} \omega ^{2n} \sigma _{d} (\ccDA_{d',n} )\cdot \sigma _{d'} (\uDA_{p-d'+d-2n,\ell } ).\]
In particular, taking $\ell =0$, we obtain \cref{eq: simplification des sigma de CL cas irrationnel}. Moreover, by construction of the fields in \cref{def: construction des termes}, the established relations show that the families $((\bfu_{p,\ell } ,S_{p,\ell } ))_{p\in \NN}$ are independent for different values of $\ell $. Only the family with $\ell =0$ has a non trivial initialization, namely $(\Delta +k_{0} ^2 )u_{0,0} =\fsource$. So it~is easy to deduce by induction that $\bfu_{p,\ell }$ and $S_{p,\ell }$ vanish for $p\in \NN$ and $\ell >0$. 
\end{demo}

\begin{remarque}
The absence of powers of $\ln r$ and $\ln \peps$ when $\pst\not\in \mathbb{Q}$ is well known for the Poisson equation, see \cite{CalCosDauVia06}. But it~is significantly less obvious to prove for the Helmholtz equation. Let us show a way to prove that there~is no $\ln r$ for Poisson and how it fails for Helmholtz (in both cases, the absence of $\ln \peps$ follows from the absence of $\ln r$).\\

\Cref{champs totaux en fct de sigma} states that the $\bfuDA_{p,\ell }$ can be expressed using the $\langle \Rmaicro\rangle (\phi_{d'} )$. So it suffices to prove that the $\langle \Rmaicro\rangle (\phi_{d'} )$ contain neither $\log (\alpha z)$ nor $\ln x$. They are defined in terms of $\ROmegaPi$, $\RGammaPi$ and $\RcouchePi$. The proof of \cref{operateurs R} shows when these operators may introduce some logarithm: for~any $R\in \{\ROmegaPi, \RGammaPi, \RcouchePi\}$ and $D\in \{\Omega ,\Gamma ,\couche\}$ s.t.\ $R:\mathcal{A}(D)\rightarrow \mathcal{A}(\Pi )$, and for any $\varphi \in \mathcal{A}_{d-\deg R} (D)$ that contains no logarithm, $R(\varphi )\in \mathcal{A}_{d} (\Pi )$ may contain logarithms only if $d\in \TZ \cup \mathbb{Z}$.\\

For the Poisson equation, we have $\omega =0$ so $k_{0} =k_{1} =0$ (the work of \cref{chap: raccord algebrique} remains valid in this case). Hence, we have $R=0$ or $\deg_{\mathcal{A}} R\in -\mathbb{N}^*$ for~any $R\in \Rmaicro$ (see \cref{fig: Rmaicro pour Poisson vs Helmholtz}). So~$\langle \Rmaicro\rangle (\phi_{d'} )$ only involves terms that belong to the spaces $\mathcal{A}_{d} (\Pi )$ with $d\in d'-\mathbb{N}^*$ (except the term $\phi_{d'}$ itself). Since $(d'-\mathbb{N}^*)\cap (\TZ \cup \mathbb{Z}) = \varnothing $ when $\pst\not\in \mathbb{Q}$, we deduce that $\langle \Rmaicro\rangle (\phi_{d'} )$ contains no $\log (\alpha z)$ or $\ln x$ for Poisson problems.\\

However, for Helmholtz problems, $\omega \neq 0$ and $\Rmaicro$ also contains operators $R$ s.t.\ $\deg _{\mathcal{A}} R\in \mathbb{N}$, for example $-k_{0} ^2 \ROmegaPi$. In this case, $\langle \Rmaicro\rangle (\phi_{d'} )$ involves terms in the spaces $\mathcal{A}_{d} (\Pi )$ with $d\in d'+\mathbb{Z}$. Since $(d'+\mathbb{Z})\cap (\TZ \cup \mathbb{Z})=\{d'\}$, some logarithms could appear in $\langle \Rmaicro\rangle (\phi_{d'} )$ at degree $d'$. See \cref{fig: Rmaicro pour Poisson vs Helmholtz}. In a sense, they do appear, but they compensate each other. For example, the term of $\langle \Rmacro\rangle (\phi_{d'} )$ in $\varepsilon ^2  \mathcal{A}_{d'} (\Pi )[\ln \varepsilon ]$ is
\newcommand{\pttFigVecs}[1]{
    \raisebox{-3.5mm}{\begin{tikzpicture}[scale=0.23]
	    \draw [thin, gray, dotted] (-2,0) grid (2,2);
        #1
    \end{tikzpicture}}}
\newcommand{\pttDrawVecs}[1]{\draw [->,>=latex, thick] #1;}
\[\underbrace{-\varepsilon ^2 k_{1} ^2 \RcouchePi}_{\pttFigVecs{
    \pttDrawVecs{(0,0) -- (0,2)}
}} (\phi_{d'} ) 
+ \underbrace{(-k_{0} ^2 \ROmegaPi)\circ (-\varepsilon ^2  \RcouchePi\circ \partial _{x|\couche} ^2 )}_{\pttFigVecs{
    \pttDrawVecs{(0,0) -- (-2,2)}
    \pttDrawVecs{(-2,2) -- (0,2)}
}} (\phi_{d'} ) 
+ \underbrace{(-\varepsilon ^2  \RcouchePi\circ \partial _{x|\couche} ^2 )\circ (-k_{0} ^2 \ROmegaPi)}_{\pttFigVecs{
    \pttDrawVecs{(0,0) -- (2,0)}
    \pttDrawVecs{(2,0) -- (0,2)}
}} (\phi_{d'} ) + \mathrm{etc} 
\]
Under each term we drew the path taken by the applied operators in coordinates $(\deg _{\mathcal{A}} ,\deg _{\varepsilon } )$, as in \cref{fig: Rmaicro pour Poisson vs Helmholtz}. Some terms of this sum do contain $\log (\alpha z)$ or $\ln r$ in general, but these logarithmic parts sum up to 0. Hence there~is no logarithm in the end. But this is~not obvious to prove at all. We observed these cancellations of logarithmic terms using the symbolic computation module of \cref{sec: implementation des espaces A}, but we were not able to explain it for around a year. The key observation is the fact that any $\mcla_{d,n}$ is a linear combination of Bessel functions of the first kind (which contain no $\ln r$), see \cref{expression sing heterog CL}. It relies on the calculatory properties of Bessel functions in \cref{relations sur BesselGamma/Omega}. This seems quite specific to the Helmholtz equation and we do~not know is there~is a more general and deeper reason for asymptotic expansions in corner not to involve $\ln \peps$ when $\pst\not\in \mathbb{Q}$.
\end{remarque}

	\tikzset{external/export=true}%
	\tikzsetfigurename{fig-Rmaicro_Poisson_vs_Helmholtz-}%
	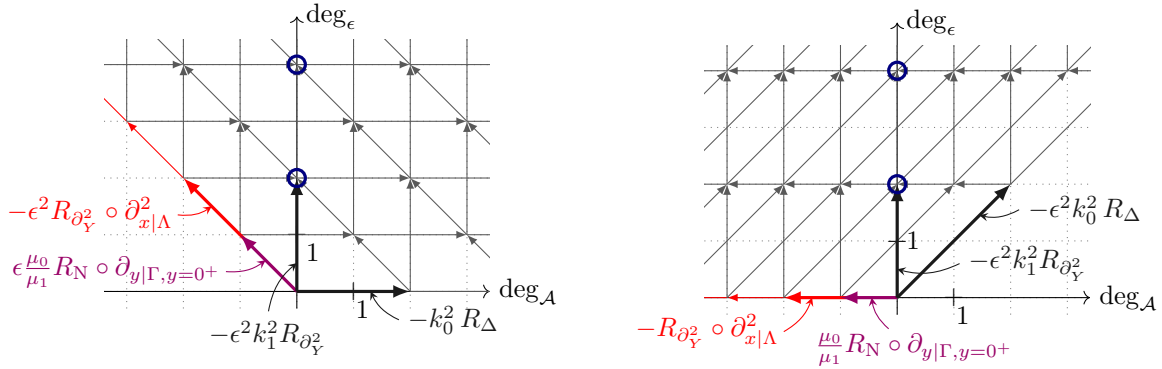
\begin{figure}[h]
\def\xmin{-3.4}
\def\xmax{3.4}
\def\ymin{-0.3}
\def\ymax{4.45}
\def\couleurPoisson{red}
\def\couleurPoissonA{red}
\def\couleurPoissonB{purple!80!blue}
\def\couleurHelmholtz{black!30!gray}
\def\couleurHelmholtzA{black!80!gray}
\newcommand{\cercleFig}[1]{\draw[very thick, blue!50!black] (0,#1) circle (1.5mm);}

\begin{tikzpicture}[scale=0.75]
	\draw [thin, gray, dotted] (\xmin,\ymin) grid (\xmax,\ymax);
	\draw [->] (\xmin,0) -- (\xmax,0) node [right] {$\deg_{\mathcal{A}}$};
	\draw [->] (0,\ymin) -- (0,\ymax+0.4) node [right] {$\deg_\varepsilon$};
	\draw (1,0) node{$+$} node[below]{\ \ 1};
	\draw (0,1) node{$+$} node[below right]{1};

	\begin{scope}
		\clip (\xmin,\ymin) rectangle (\xmax,\ymax);

		\draw [->,>=latex,\couleurPoisson] (-2,2) -- (-3,3);
		\draw [->,>=latex,\couleurPoisson] (-3,3) -- (-4,4);
		\draw [->,>=latex,thin,\couleurHelmholtz] (2,0) -- (4,0);
		\foreach \i in {0,1,...,4} {
			\draw [->,>=latex,thin,\couleurHelmholtz] (-1-\i,1+\i) -- (-1-\i,3+\i); 
			\foreach \j in {0,2,4,6} {
				\draw [->,>=latex,thin,\couleurHelmholtz] (2+\j-\i,\i) -- (1+\j-\i,\i+1); 
				\draw [->,>=latex,thin,\couleurHelmholtz] (-1+\j-\i,1+\i) -- (1+\j-\i,1+\i); 
				\draw [->,>=latex,thin,\couleurHelmholtz] (\j+2-\i,\i) -- (\j+2-\i,2+\i); 
			}
		}

		\draw [->,>=latex,very thick,\couleurPoissonA] (-0.9,0.9) -- (-2,2);
		\draw [->,>=latex,very thick,\couleurPoissonB] (0,0) --  (-1,1);
		\draw [->,>=latex,very thick,\couleurHelmholtzA] (0,0) -- (0,2);
		\draw [->,>=latex,very thick,\couleurHelmholtzA] (0,0) -- (2,0);

		\cercleFig{2}
		\cercleFig{4}
	\end{scope}

	\draw[<-,>=stealth,\couleurPoissonA] (-1.5,1.5) to[bend left=20] ++(-0.6,-0.15) node[left] {$-\varepsilon^2 \RcouchePi\circ\partial_{x|\couche}^2 \!$};
	\draw[<-,>=stealth,\couleurPoissonB] (-0.55,0.55) to[bend left=20] ++(-0.7,-0.2) node[left]  {$\varepsilon \frac{\mu_0}{\mu_1} \RGammaPi\circ\partial_{y|\Gamma,y=0^+} \!$};
	\draw[<-,>=stealth,\couleurHelmholtzA] (0,0.5) to[bend right=20] ++(-0.5,-0.85) node[below] {$-\varepsilon^2 k_1^2 \RcouchePi$};
	\draw[<-,>=stealth,\couleurHelmholtzA] (1.3,0) to[bend right=20] ++(0.5,-0.5) node[right] {$-k_0^2\,\ROmegaPi$};

	\draw (0,-1.3) node {\ }; 
\end{tikzpicture}%
\hspace{8mm}\nolinebreak%
\begin{tikzpicture}[scale=0.75]
	\draw [thin, gray, dotted] (\xmin,\ymin) grid (\xmax,\ymax);
	\draw [->] (\xmin,0) -- (\xmax,0) node [right] {$\deg_{\mathcal{A}}$};
	\draw [->] (0,\ymin) -- (0,\ymax+0.4) node [right] {$\deg_\varepsilon$};
	\draw (1,0) node{$+$} node[below]{\ \ 1};
	\draw (0,1) node{$+$} node[right]{1};

	\begin{scope}
		\clip (\xmin,\ymin) rectangle (\xmax,\ymax);
			
		\draw [->,>=latex,\couleurPoisson] (-2,0) -- (-3,0);
		\draw [->,>=latex,\couleurPoisson] (-3,0) -- (-4,0);
		\foreach \i in {0,1,...,8} {
			\draw [->,>=latex,thin,\couleurHelmholtz] (-1-\i,0) -- (1-\i,2); 
			\draw [->,>=latex,thin,\couleurHelmholtz] (-\i,0) -- (-\i,2); 
			\foreach \j in {0,2} {
				\draw [->,>=latex,thin,\couleurHelmholtz] (2-\i+\j,2+\j) -- (4-\i+\j,4+\j); 
				\draw [->,>=latex,thin,\couleurHelmholtz] (2-\i+\j,2+\j) -- (2-\i+\j,4+\j); 
				\draw [->,>=latex,thin,\couleurHelmholtz] (2-\i+\j,2+\j) -- (1-\i+\j,2+\j); 
			}
		}

		\draw [->,>=latex,very thick,\couleurPoissonA] (-0.9,0) -- (-2,0);
		\draw [->,>=latex,very thick,\couleurPoissonB] (0,0) -- (-1,0);
		\draw [->,>=latex,very thick,\couleurHelmholtzA] (0,0) -- (2,2);
		\draw [->,>=latex,very thick,\couleurHelmholtzA] (0,0) -- (0,2);

		\cercleFig{2}
		\cercleFig{4}
	\end{scope}

	\draw[<-,>=stealth,\couleurPoissonA] (-1.5,0) to[bend left=20] ++(-0.5,-0.6) node[left] {$-\RcouchePi\circ\partial_{x|\couche}^2 \!$};
	\draw[<-,>=stealth,\couleurPoissonB] (-0.45,0) to[bend right=20] ++(0.3,-0.55);
	\draw[\couleurPoissonB] (0.3,-0.45) node[below] {$\frac{\mu_0}{\mu_1} \RGammaPi\circ\partial_{y|\Gamma,y=0^+}$};
	\draw[<-,>=stealth,\couleurHelmholtzA] (0,0.5) to[bend right=20] ++(1.1,0.1) node[right] {$-\varepsilon^2 k_1^2 \RcouchePi$};
	\draw[<-,>=stealth,\couleurHelmholtzA] (1.45,1.45) to[bend right=20] ++(0.7,0.1) node[right] {$-\varepsilon^2 k_0^2\,\ROmegaPi$};
\end{tikzpicture}
\vspace{-2mm}

\caption[Degrees of the elements of $\Rmaicro$ and their iterated composition for Poisson and Helmholtz equations]{Degrees of the elements of $\Rmacro$ (on the left) and $\Rmicro$ (on the right) in bold and their iterated composition. When $\omega=0$, only the colored arrows appear. The circled points are the degrees where some $\ln r$ could naively appear given the properties of the elements of $\Rmaicro$.}
\label{fig: Rmaicro pour Poisson vs Helmholtz}
\end{figure}%
	\tikzset{external/export=false}%
	\tikzsetfigurename{main-figure}

\subsection{The rational case}
\label{sec: simplifications DA cas rationnel}

\begin{textesansboite}
In this section, we assume that $\pst$ is rational and we show that the asymptotic expansion of $u_{\peps}$ behaves like in the irrational case at the first orders.
\end{textesansboite}

\begin{lemme}{}
\label{valeur de cmaicro cas rationnel}%
Let us assume $\pst \in \mathbb{Q}$. Let $(a,b)\in \mathbb{N}^*$ be s.t.\ $\Theta = \frac{a}{b} \pi $ and $\pgcd{a,b}=1$.\\
Let $d,d'\in \TZ*$, $p\in \NN$ and $\ell \in \mathbb{N}$. Let us assume $d'>-b$.
\begin{itemize}
    \item If $p+d'<b$, we have $\cmacro_{d,d',p,\ell } = \delta _{d,d'} \,\delta _{p,-d} \,\delta _{\ell ,0}$. If $p-d<b$, we have $\cmicro_{d,d',p,\ell } = \delta _{d,d'} \,\delta _{p,d} \,\delta _{\ell ,0}$.
    \item Let us assume $d<0$ and $p<b$. If there exists $n\in \mathbb{N}$ s.t.\ $p=d'-d+2n$, then $\cfinal_{d,d',p,0} = \omega ^{2n} \sigma _{d} (\ccDA_{d',n} )$. Otherwise, $\cfinal_{d,d',p,\ell } =0$.
\end{itemize}\fbi
\end{lemme}

\begin{textesansboite}
The proof is similar to the one of \cref{valeur de cmaicro cas irrationnel}, but more subtle as we study only certain values of $(d,d',p)$.
\end{textesansboite}

\begin{demo}
\underline{Step 1:} Let us first treat $\cmicro_{d,d',p,\ell }$ when $p-d<b$. By~\cref{expression de Rmacro phi avec les mcl}, we have
\begin{align*}
    \cmicro_{d,d',p,\ell } &:= \sigma _{d,p,\ell } \circ  \Eps \circ  \langle \Rmacro\rangle  (\phi_{d'} ) \\
    &= \sigma _{d,p,\ell } \circ  \Eps \circ  \bigg( \fsum_{n=0} ^{b-1} \varepsilon ^{n} \.\mcla_{d',n} \bigg) + \sigma _{d,p,\ell } \circ  \Eps \circ  \bigg( \fsum_{n=b} ^{\infty } \varepsilon ^{n} \.\mcla_{d',n} \bigg) \\
    &= \sigma _{d,p,\ell } \circ  \Eps \circ  \bigg( \fsum_{n=0} ^{b-1} \varepsilon ^{n} \.\mcla_{d',n} \bigg) \quad\qquad  \mathrm{for} \ p-d<b
\end{align*}
since $\Eps$ send $\varepsilon ^{n} \mathcal{A}_{d} (\Omega )[\ln \varepsilon ]$ to $\varepsilon ^{n+d} \mathcal{A}_{d} (\Omega )[\ln \varepsilon ]$ for~any $n$. Moreover, the assumption $d'>-b$ yields $d'\not\in -b\mathbb{N}^*$. So \cref{expression sing heterog CL} states that, for~any $n<b$, $\mcla_{d',n}$ is a linear combination of the $\BesselOmega_{a}$ ($a\in \mathbb{R}$). Hence $\mcla_{d',n}$ does~not involve $\log (\alpha z)$. Using the same reasoning as in \cref{valeur de cmaicro cas irrationnel} we deduce $\cmicro_{d,d',p,\ell } = \delta _{d,d'} \,\delta _{p,d} \,\delta _{\ell ,0}$.
\medskip\smallskip

\underline{Step 2:} Let us treat $\cmacro_{d,d',p,\ell }$ when $p+d'<b$. The proof of \cref{valeur de cmaicro cas irrationnel} used that $\varepsilon ^{d'} \langle \Rmicro\rangle  (\phi_{d'} ) = \Eps \circ \langle \Rmacro\rangle  (\phi_{d'} )$. In the case $\pst \in \mathbb{Q}$, this equality does~not hold anymore but there~is still equality for some terms of these formal sums. Indeed, by \cref{Rmaicro defini et inversible dans Aemaicro,rq: def Eps}, $\varepsilon ^{-d'} \Eps \circ \langle \Rmacro\rangle  (\phi_{d'} )$ belongs to $\Aemicro(\Pi )$ and, by \cref{pb de Helmholtz dans Aemacro,action de Eps sur Helmholtz dans Aemaicro}, it~is solution of the Helmholtz problem \cref{eq: pb de Helmholtz dans Aemicro bis}. So \cref{pb de Helmholtz dans Aemicro} shows that we have
\[\varepsilon ^{-d'} \Eps \circ \langle \Rmacro\rangle  (\phi_{d'} ) = \langle \Rmicro\rangle  \bigg( \fsum_{p'\in \NN} \,\fsum_{d''\in \TZ*\cap (p'+d'-\NN)}\, \sum _{\ell '\in \mathbb{N}} \sigma _{d'',p',\ell '} \big(\varepsilon ^{-d'}  \Eps \circ \langle \Rmacro\rangle  (\phi_{d'} )\big) \cdot  \varepsilon ^{p'\!} \,\ln ^{\ell '\!} \!\varepsilon  \,\phi_{d''} \bigg)\]
so
\begin{subequations}
\label{eq: preuve valeur de cmaicro cas rationnel 1}
\begin{align}
    \varepsilon ^{-d'} \langle \Rmacro\rangle  (\phi_{d'} ) &= \Eps^{-1} \circ \langle \Rmicro\rangle \bigg(\fsum_{p'\in \NN} \,\fsum_{d''\in \TZ*\cap (p'+d'-\NN)}\, \sum _{\ell '\in \mathbb{N}} \,\cmicro_{d'',d',p'+d',\ell '} \cdot  \varepsilon ^{p'\!} \,\ln ^{\ell '\!} \!\varepsilon  \,\phi_{d''} \bigg) \\
    &= \fsum_{p'\in \NN} \,\fsum_{d''\in \TZ*\cap (p'+d'-\NN)}\, \sum _{\ell '\in \mathbb{N}} \,\cmicro_{d'',d',p'+d',\ell '} \cdot  \varepsilon ^{p'\!} \,\ln ^{\ell '\!} \!\varepsilon  \cdot \Eps^{-1} \circ \langle \Rmicro\rangle (\phi_{d''} )
    \label{eq: preuve valeur de cmaicro cas rationnel 1b} \\
    &= \Eps^{-1} \circ \langle \Rmicro\rangle (\phi_{d'} ) + \text{other terms} . 
    \label{eq: preuve valeur de cmaicro cas rationnel 1c}
\end{align}
\end{subequations}
As one can see in \cref{fig: preuve valeur de cmaicro cas rationnel}, there~is a set of terms (the hatched zone) on which $\varepsilon^{-d'}\langle\Rmacro\rangle(\phi_{d'})$ and $\Eps^{-1}\circ\langle\Rmicro\rangle(\phi_{d'})$ coincide, although they are~not globally equal (because of the other terms mentioned in \cref{eq: preuve valeur de cmaicro cas rationnel 1c}). We will see that, because of the assumption $p+d'<b$, $\cmacro_{d,d',p,\ell }$ only ``sees'' the terms on which $\varepsilon^{-d'}\langle\Rmacro\rangle(\phi_{d'})$ and $\Eps^{-1}\circ\langle\Rmicro\rangle(\phi_{d'})$ coincide. This will allow us to obtain for $\cmacro_{d,d'p,\ell }$ the same formula as in \cref{valeur de cmaicro cas irrationnel}.

	\tikzset{external/export=true}%
	\tikzsetfigurename{fig-preuve_valeur_cmacro_cas_rationnel-}%

\begin{figure}[h]
\def\couleurGdConeFig{blue}
\def\couleurPttConeFig{red}
\def\couleurGdTexteFig{blue!70!black}
\def\couleurPttTexteFig{red!70!black}
\def\vertFig{green!60!black}

\newcommand{\fillGrandCone}[2]{
	\begin{scope}[xshift={#1cm}, yshift={#2cm}]
		\clip (-10,10) -- (0,0) -- (10,0) -- (10,10) -- cycle;
		\shade[inner color=\couleurGdConeFig, outer color=white, opacity=0.3] ({-2.5*cos(67.5)},{-2*sin(67.5)}) circle (5); 
		\shade[inner color=\couleurGdConeFig, outer color=white, opacity=0.35] ({-2.5*cos(67.5)},{-2*sin(67.5)}) circle (9);
	\end{scope}
}
\newcommand{\drawGrandCone}[2]{
	\begin{scope}[xshift={#1cm}, yshift={#2cm}]
		\draw[very thick, \couleurGdConeFig] (-10,10) -- (0,0) -- (10,0);
	\end{scope}
}

\tikzfading[name=fade out, inner color=transparent!67, outer color=transparent!100]
\newcommand{\fillPetitCone}[3][1]{
	\begin{scope}[xshift={#2cm}, yshift={#3cm}]
		\clip (-10,10) -- (0,0) -- (10,0) -- (10,10) -- cycle;
		\fill[\couleurPttConeFig, path fading=fade out] ({-1.7*cos(67.5)*#1},{-1.7*sin(67.5)*#1}) circle (2.5*#1);
		\fill[\couleurPttConeFig, path fading=fade out] ({-1.7*cos(67.5)*#1},{-1.7*sin(67.5)*#1}) circle (4*#1);
	\end{scope}
}
\newcommand{\drawPetitCone}[3][1]{
	\begin{scope}[xshift={#2cm}, yshift={#3cm}]
		\clip (-10,10) -- (0,0) -- (10,0) -- (10,10) -- cycle;
		\draw[very thick, \couleurPttConeFig, path fading=east] (0,0) -- (3.5*#1,0) ++(0,1pt);
		\draw[very thick, \couleurPttConeFig, path fading=west] (0,0) -- ({-3.5/sqrt(2)*#1},{3.5/sqrt(2)*#1}) ++(1pt,1pt);
	\end{scope}
	\begin{scope}[xshift={#2cm}, yshift={#3cm}]
		\clip (-10,10) -- (0,0) -- (10,0) -- (10,-10) -- (-10,-10) -- cycle;
		\draw[very thick, \couleurPttConeFig, path fading=east] (0,0) -- (3.5*#1,0) ++(0,-1pt);
		\draw[very thick, \couleurPttConeFig, path fading=west] (0,0) -- ({-3.5/sqrt(2)*#1},{3.5/sqrt(2)*#1}) ++(-1pt,-1pt);
	\end{scope}
}
\newcommand{\filldrawPetitCone}[3][1]{
	\fillPetitCone[#1]{#2}{#3}
	\drawPetitCone[#1]{#2}{#3}
}

\def\croixFig{\raisebox{1pt}{\small\bf+}}

\def\xmin{-3.5}
\def\xmax{4.5}
\def\ymin{-1.6}
\def\ymax{4.2}
\def\dprimeFig{0.85}
\def\dsecondeFig{2.9}
\def\pprimeFig{3}
\def\bFig{1.2}


\raggedleft 
\begin{tikzpicture}[scale=0.9]
	\begin{scope}
		\clip (\xmin,\ymin) rectangle (\xmax,\ymax);

		\fillGrandCone{\dprimeFig}{-\dprimeFig}
		\fillPetitCone[0.65]{\dprimeFig}{-\dprimeFig}
		\fillPetitCone[0.65]{\dprimeFig}{-\dprimeFig} 
		\filldrawPetitCone[0.5]{-1.8}{2.6}
		\filldrawPetitCone[0.4]{-0.25}{3.4}
		\filldrawPetitCone[0.55]{1.6}{1.7}
		\filldrawPetitCone[0.5]{\dsecondeFig}{\pprimeFig}

		\fill [pattern={Lines[angle=45,line width=.5pt,distance=7pt]}, pattern color=\vertFig] (-\bFig,\bFig) -- (\dprimeFig,-\dprimeFig) -- (10,-\dprimeFig) -- (10,\bFig) -- cycle;
		\draw[thick, green!65!black] (-\bFig,\bFig) -- (10,\bFig);
		\draw (0,\bFig) node{\croixFig} node[above left]{$b-d'\!$};
		\draw [<->,>=latex] (-\bFig,\bFig-0.1) -- node[near start,left]{\raisebox{-5ex}{$b$}} (-\bFig,-\dprimeFig);
		\draw[dotted] (-\bFig,-\dprimeFig) -- (\dprimeFig,-\dprimeFig);

		\drawGrandCone{\dprimeFig}{-\dprimeFig}
		\drawPetitCone[0.65]{\dprimeFig}{-\dprimeFig}
		\drawPetitCone[0.65]{\dprimeFig}{-\dprimeFig} 

		\draw[\couleurGdTexteFig] (\dprimeFig,-\dprimeFig) 
			node {\croixFig}
			node[below] {$(d',-d')$};
		\draw[\couleurPttTexteFig] (\dsecondeFig,\pprimeFig) 
			node {\raisebox{1pt}{\croixFig}}
			node [below] {$(d'',p'-d'')$};
	\end{scope}

	\draw[<-,>=latex,\couleurGdTexteFig] ({-\ymax*0.55},{\ymax*0.55}) to[bend right=20] ++(-0.5,-0.5) node[below]{\parbox{28mm}{
		support of\\ $\varepsilon^{-d'}\langle\Rmacro\rangle(\phi_{d'})$}}; 
	\draw[<-,>=latex,\couleurPttTexteFig] (\dprimeFig+0.8,-\dprimeFig) to[bend right=20] ++(0.6,-0.4) node[right]{\parbox{42mm}{
		support of $\Eps^{-1}\circ\langle\Rmicro\rangle(\phi_{d'})$}}; 
	\draw[<-,>=latex,\couleurPttTexteFig] (\dsecondeFig+1.2,\pprimeFig) to[bend right=20] ++(0.7,-0.3) node[right]{\parbox{45mm}{
		support of a term $\cmicro_{d''\!,d'\!,p'+d'\!,\ell'}$\\ $\cdot\,\varepsilon^{p'}\ln^{\ell'}\!\varepsilon\cdot\Eps^{-1}\circ\langle\Rmicro\rangle(\phi_{d''})$}}; 
	\draw[<-,>=latex,\vertFig] (\xmax,{\bFig*0.55}) to[bend right=20] ++(0.5,0.4) node[right]{\parbox{30mm}{
		only one red cone\\ in this zone}}; 

	\draw [->] (\xmin,0) -- ({\xmax+0.25},0) node [right] {degree in the spaces $\mathcal{A}_q(\Pi)$};
	\draw [->] (0,\ymin) -- (0,{\ymax+0.3}) node [right]{power of $\varepsilon$};
\end{tikzpicture}

\longCaption{Schematic illustration of \cref{eq: preuve valeur de cmaicro cas rationnel 1}}{. Each cone represents the convex closure of the support of a formal series in the equation. They are infinite cones, but they are truncated for visual clarity. In the hatched zone, we have $\cmicro_{d''\!,d'\!,p'+d'\!,\ell'} = \delta_{d'\!,d''}\,\delta_{p'\!,-d'}\,\delta_{\ell'\!,0}$, so $\varepsilon^{-d'}\langle\Rmacro\rangle(\phi_{d'})$ coincides with $\Eps^{-1}\circ\langle\Rmicro\rangle(\phi_{d'})$ in this zone.}
\label{fig: preuve valeur de cmaicro cas rationnel}
\end{figure}%
	\tikzset{external/export=false}%
	\tikzsetfigurename{main-figure}

Applying $\sigma _{d,p,\ell }$ on both sides of \cref{eq: preuve valeur de cmaicro cas rationnel 1b} gives
\begin{equation}{\label{eq: preuve valeur de cmaicro cas rationnel 2}}
    \delta _{d,d'} \,\delta _{p,-d'} \,\delta _{\ell ,0} = \sum _{p'\in \NN} \,\sum _{d''\in \TZ*\cap (p'+d'-\NN)}\, \sum _{\ell '\in \mathbb{N}} \,\cmicro_{d'',d',p'+d',\ell '} \cdot  \cmacro_{d,d'',p-p',\ell -\ell '} .
\end{equation}
For the left-hand side we used that $\sigma _{d} \circ R=0$ for~any $R\in \{\ROmegaPi,\RGammaPi,\RcouchePi\}$, which yields $\sigma _{d,p,\ell } \circ \langle \Rmacro\rangle  (\varepsilon ^{-d'} \phi_{d'} ) =\sigma _{d,p,\ell } (\varepsilon ^{-d'} \phi_{d'} )$. As~for the right-hand side, it contains a finite number of non-zero terms by \cref{def: sigma et fsum commutent dans Aeps}. Let us precise what are these non-zero terms. Our final goal is to show that the right-hand side is equal to $\cmacro_{d,d',p,\ell }$.\\
Any $R\in \Rmicro$ satisfies $\deg _{\varepsilon } R\geqslant 0$ and $\deg _{\mathcal{A}} R\leqslant \deg _{\varepsilon } R$, so we have
\[\langle \Rmicro\rangle (\phi_{d''} ) \in \fsum_{p'\geqslant 0} \fsum_{q\leqslant p'+d''} \varepsilon ^{p'} \mathcal{A}_{q} (\Pi )[\ln \varepsilon ] \qquad  \text{so}  \qquad  
\Eps^{-1} \circ  \langle \Rmicro\rangle (\phi_{d''} ) \in \fsum_{\substack{(p',q)\\ p'\geqslant -q, 0\leqslant p'+d''}} \varepsilon ^{p'} \mathcal{A}_{q} (\Pi )[\ln \varepsilon ].\]
Thus $\cmacro_{d,d'',p-p',\ell -\ell '} \neq 0$ only if $0\leqslant p-p'+d''$, which implies $p'+d'-d''<p+d'<b$ by assumption. Now using step~1, $p'+d'-d''<b$ yields $\cmicro_{d'',d',p'+d',\ell } =\delta _{d'',d'} \,\delta _{p'+d',d'} \,\delta _{\ell ',0}$. Therefore \cref{eq: preuve valeur de cmaicro cas rationnel 2} simplifies and becomes $\delta _{d,d'} \,\delta _{p,-d'} \,\delta _{\ell ,0} =\cmacro_{d,d',p,\ell }$.\medskip\smallskip

\underline{Step 3:} Let us treat $\cfinal_{d,d',p,\ell }$ when $d<0$ and $p<b$. We recall its definition (for $d<0$):
\begin{equation}{\label{eq: preuve valeur de cmaicro cas rationnel 3}}
    \cfinal_{d,d',p,\ell } := \! \sum _{\substack{(p_{1} ,p_{2} )\in (\mathbb{N}-d)\times (\mathbb{N}+d')\\ n\in \mathbb{N},\, p_{1} +p_{2} +2n=p}}
    \ \sum _{\substack{d_{1} \in \TZ*\\ p_{1} +d_{1} \in \mathbb{N}}} \ \sum _{\substack{d_{2} \in \TN\\ p_{2} -d_{2} \in \mathbb{N}}}
    \ \sum _{\substack{(\ell _{1} ,\ell _{2} )\in \mathbb{N}^2 \\ \ell _{1} +\ell _{2} =\ell }}
    \cmacro_{d,d_{1} ,p_{1} ,\ell _{1} }\cdot \omega ^{2n} \sigma _{d_{1}} (\mathcal{S} ^{\infty } _{d_{2} ,n} )\cdot \cmicro_{d_{2} ,d'\!,p_{2} ,\ell _{2} } .
\end{equation}
Let us show that any non-zero term satisfies $d_{1} >-b$, $p_{1} +d_{1} <b$ and $p_{2} -d_{2} <b$ in order to apply steps~1 and 2. We do it by contraposition, by showing that if a non-zero term verifies $d_{1} \leqslant -b$ or $p_{1} +d_{1} \geqslant b$ or $p_{2} -d_{2} \geqslant b$ then $p\geqslant b$.
\begin{itemize}
    \item If $d_{1} \leqslant -b$, then $p=p_{1} +p_{2} +2n=(p_{1} +d_{1} )-d_{1} +(p_{2} -d_{2} )+d_{2} +2n\geqslant 0+b+0+0+0=b$.
    \item Let us assume $p_{1} +d_{1} \geqslant b$. If $d_{1} >0$, then $\sigma _{d_{1}} (\mathcal{S} ^{\infty } _{d_{2} ,n} ) =\delta _{d_{1} ,d_{2}} \delta _{n,0}$. Since we are looking at a non-zero term, we have $d_{1} =d_{2}$ and $n=0$. Thus $p=(p_{1} +d_{1} ) +(p_{2} -d_{2} ) \geqslant b+0$. Otherwise, we have $d_{1} <0$, which yields $p=(p_{1} +d_{1} )-d_{1} +(p_{2} -d_{2} )+d_{2} +2n \geqslant b+0+0+0+0=b$.
    \item If $p_{2} -d_{2} \geqslant b$, then $p=(p_{1} +d)-d +(p_{2} -d_{2} )+d_{2} +2n\geqslant 0+0+b+0+0=b$.
\end{itemize}
By~steps~1 and 2 we deduce that any non-zero term of \cref{eq: preuve valeur de cmaicro cas rationnel 3} satisfies $d=d_{1}$, $d_{2} =d'$, $p_{1} =-d$, $p_{2} =d$, $\ell _{1} =\ell _{2} =0$ and $\cmacro_{d,d_{1} ,p_{1} ,\ell _{1} }=\cmicro_{d_{2} ,d'\!,p_{2} ,\ell _{2} }=1$. Similarly to \cref{valeur de cmaicro cas irrationnel}, this implies the desired formula.
\end{demo}

\begin{theoremesecable}{}
\label{simplifications cas rationnel}%
Let us assume that $\pst \in \mathbb{Q}$. Let $a,b\in \mathbb{N}^*$ be s.t.\ $\Theta = \frac{a}{b} \pi $ and $\pgcd{a,b}=1$ and let
\[p_{\min } := \left\{\begin{array}{ll} 2b & \text{if } a=1\\
b & \text{otherwise} \end{array}\right.\]
We have: \[\forall p\in \NN,\quad  p<p_{\min }  \Rightarrow (\forall \ell  \neq 0, u_{p,\ell } =0).\]
In other words,  there~is no $\ln \peps$ in the asymptotic expansion of $u_{\peps}$ before the order $p_{\min } $.\\
Moreover, for~any $p\in \NN$ s.t.\ $p<p_{\min }$, we have
\begin{equation}{\label{eq: simplifications cas rationnel}}
    {\textstyle \forall d\in -\TN}, \qquad  \sigma _{d} (u_{p,0} )= \sum _{\substack{n\in \mathbb{N}\\ 2n\leqslant p+d}} \sum _{\substack{d'\in \TN\\ d'\leqslant p+d-2n}} \omega ^{2n} \sigma _{d} (\ccDA_{d',n} )\cdot \sigma _{d'} (\uDA_{p-d'+d-2n,0} ).
\end{equation}\fbeq
\end{theoremesecable}

\begin{remarque}
In~practice, the asymptotic expansion of $u_{\peps}$ is~not used at very high orders. For example, in \cref{chap: ME}, we build effective models of order~3 and 5. \Cref{simplifications cas rationnel} shows that, except for a few angles $\Theta $, it~is as if $\pst$ were irrational for~any $p\leqslant 5$. Moreover, $p_{\min }$ is~not optimal in some cases. \Cref{simplifications cas rationnel version forte} improves the value of $p_{\min }$ when $b\in \{2,4\}$ (with theoretical justification) and \cref{tab: ordre du premier ln eps} gives the optimal values of $p_{\min }$ for the angles with $b\leqslant 7$ (computed numerically). They tell us that except for $\Theta  \in \{\pi , \frac{2\pi }{3} , \frac{4\pi }{3} , \frac{5\pi }{3} \}$ the asymptotic expansion of $u_{\peps}$ contains no $\ln \peps$ at orders $p<5$ (which gives an approximation of $u_{\peps}$ in $O(\peps^{5} \ln ^{\ell } \!\peps)$ for some $\ell $).
\end{remarque}

\begin{demo}
Let $p<b$, $\ell \in \mathbb{N}$ and $d\in -\TN$. \Cref{construction directe des upl} states:
\begin{equation}{\label{eq: preuve simplifications cas rationnel 1}}
    \sigma _{d} (u_{p,\ell } ) = \!\! \sum _{\substack{p'\in \NN\\ p-p'\in \NN+ 2\pst}} \; \sum _{\substack{d'\in \TZ*\\ p'+d'\in \intvl{0,p+d}}} \! \sum _{\ell '=0} ^{\ell } \cfinal_{d,d',p-p',\ell -\ell '} \cdot \sigma _{d'} (\uDA_{p',\ell '} ) .
\end{equation}
The terms of this sum satisfy $p-p'\leqslant p<b$ and $d'\geqslant -p'>-p>-b$, so \cref{valeur de cmaicro cas rationnel} applies. Hence, similarly to \cref{simplifications cas irrationnel}, one can deduce 
\[{\textstyle \forall d\in -\TN}, \qquad  \sigma _{d} (u_{p,\ell } )= \sum _{\substack{n\in \mathbb{N}\\ 2n\leqslant p+d}} \sum _{\substack{d'\in \TN\\ d'\leqslant p+d-2n}} \omega ^{2n} \sigma _{d} (\ccDA_{d',n} )\cdot \sigma _{d'} (\uDA_{p-d'+d-2n,\ell } ).\]
Taking $\ell =0$, it gives \cref{eq: simplifications cas rationnel}. In addition, like in \cref{simplifications cas rationnel}, it~is easy to deduce by induction that $u_{p,\ell } =0$ for~any $p<b$ and $\ell >0$. Moreover, for~any $p<2\pst$, the sums in the right-hand sides of \cref{eq: preuve simplifications cas rationnel 1,eq: simplifications cas rationnel} are empty, so \cref{eq: simplifications cas rationnel} and the property $(\ell >0 \Rightarrow  u_{p,\ell } =0)$ can be extended to any $p<\max (b,2\pst)$. Finally, since $\pst = \frac{b}{a}$, we have $\max (b,2\pst)=p_{\min }$.
\end{demo}

\begin{textesansboite}
Let us now state a different but complementary result. While \cref{valeur de cmaicro cas rationnel} gives the precise value of \emph{some} coefficients $\cmacro$ and $\cmicro$, \cref{classification coefs cmaicro} gives a partial description of \emph{any} coefficient $\cmacro$ and $\cmicro$.
\end{textesansboite}

\begin{proposition}{classification of the non-zero coefficients $\cmacro$ and $\cmicro$}
\label{classification coefs cmaicro}%
Let $(d,d',p,\ell )\in (\TZ*)^2 \times \NN\times \mathbb{N}$ be s.t.\ $\cmacro_{d,d',p,\ell } \neq 0$, resp.\ $\cmicro_{d,d',p,\ell } \neq 0$. There are two possibilities:
\begin{enumerate}
    \item either $(d',p)=(d,-d)$, resp.\ $(d',p)=(d,d)$, and $\ell =0$ and the coefficient is equal to 1.
    \item or $(d',p) \neq (d,-d)$, resp.\ $(d',p) \neq (d,d)$, and $\ell >0$ and the coefficient depends on at least one variable among $\omega ,\mu _{0} ,\mu _{1} ,\rho _{0} ,\rho _{1}$.
\end{enumerate}
\end{proposition}

\begin{remarque}
Each case is associated with a different type of matching relation.
\begin{enumerate}
    \item The first case describes the simplest matching relations : if $\bfuDA_{p,\ell }$ contains a term $c\,\phi_{d}$ with $d\in \TZ*$ and $c\in \mathbb{C}$, then $\fsum_{p',\ell '} \varepsilon ^{p'} \ln ^{\ell '\!} \!\varepsilon \, \SDA_{p'\!,\ell '}$ must contain $c\cdot \varepsilon ^{p} \ln ^{\ell } \!\varepsilon \,\Eps(\phi_{d} )=c\,\varepsilon ^{p+d} \ln ^{\ell } \!\varepsilon \cdot \phi_{d}$, so $\SDA_{p+d,\ell }$ must contain $c\,\phi_{d}$. In terms of matching coefficients, it means $\cmicro_{d,d,d,0} =1$, because $\sigma _{d} (\SDA_{p+d,\ell } ) =\cmicro_{d,d,d,0} \cdot \sigma _{d} (\uDA_{p,\ell } )$ + other terms (see \cref{condition de raccord}). Similarly for $\cmacro_{d,d,-d,0}$.
    \item The second case happens when the behaviors of the fields around the corner $\bfuDA_{p,\ell }$ and $\SDA_{p,\ell }$ contain some $\log (\alpha z)$. These logarithms appear when applying $\ROmegaPi$, $\RGammaOmega$ and $\RcouchePi$ to generate the field behaviors by induction, see \cref{eq: bfuDApl}. Due to the lack of uniqueness of the problems~\cref{eq: operateur R Omega Pi}--\cref{eq: operateur R Gamma Pi} that define $\ROmegaPi$, $\RGammaOmega$ and $\RcouchePi$, the logarithmic terms do~not match automatically between far and corner fields, so specific matching relations are needed. This is also linked to the fact that $\Eps^{-1} \circ \langle \Rmicro\rangle  \neq \langle \Rmacro\rangle \circ \Eps^{-1}$ in general, see \cref{rq: conjugaison de Rmaicro par Eps} for an explanation with other words. This type of matching relation depends heavily on the angle $\Theta $ and is tough to compute by hand.
\end{enumerate}
\end{remarque}

\begin{demo}
Let us prove points 1 and 2 only for $\cmicro$ as it~is similar for $\cmacro$. For~any $\varphi \in \mathcal{A}(\Pi )$ s.t.\ $\varphi _{|\Omega } = \Im[(\alpha z)^{q} \,\overline{\alpha z}^{k} P(\log (\alpha z))]$ with $(q,k,P)\in \mathbb{R}\times \mathbb{N}\times \mathbb{R}[T]$, we have by~definition of $\Eps$
\[\Eps(\varphi )_{|\Omega } = \sum _{j=0} ^{\deg P} \varepsilon ^{d} \,\frac{\ln ^{j} \!\varepsilon }{j!} \,\Im[(\alpha z)^{q} \,\overline{\alpha z}^{k} P^{(j)} (\log (\alpha z))]\]
So, for~any $d\in \TZ*$, $\sigma _{d,p,0} \big(\Eps(\varphi )\big) = \delta _{p,d} \,\sigma _{d} (\varphi )=\sigma _{d,p-d,0} (\varphi )$. Therefore, $\cmicro_{d,d',p,0} := \sigma _{d,p,0} \circ \Eps\circ  \langle \Rmacro\rangle (\phi_{d'} ) = \sigma _{d,p-d,0} \circ  \langle \Rmacro\rangle (\phi_{d'} )=\delta _{d,d'} \,\delta _{p,d}$, using that $\sigma _{d} \circ R=0$ for~any $R\in \{\ROmegaPi,\RGammaPi,\RcouchePi\}$.\\
Now let $\ell >0$ and $(d',p):=(d,d)$. The term of $\Eps \circ  \langle \Rmacro\rangle (\phi_{d} )$ in $\varepsilon ^{d} \mathcal{A}_{d} (\Pi )[\ln \varepsilon ]$ is $\Eps(\phi_{d} ) =\varepsilon ^{d} \phi_{d}$. Hence, $\cmicro_{d,d',p,\ell } := \sigma _{d,d,\ell } \circ \Eps \circ  \langle \Rmacro\rangle (\phi_{d} ) =\sigma _{d,d,\ell } (\varepsilon ^{d} \phi_{d} )=0$. Finally, if $(d',p) \neq (d,d)$ and $\cmicro_{d,d',p,\ell } \neq 0$, then this coefficient depends on at least one variable among $k_{0}$, $k_{1}$ and $\frac{\mu _{0} }{\mu _{1}}$ by~definition of $\Rmacro$.
\end{demo}

\let\varepsilon=\peps

\section{Examples of the structure of the asymptotic expansion for different angles}
\label{sec: exemples de structures de DA}

\begin{textesansboite}
In this section, we give and comment explicit examples of the asymptotic expansion at small orders for different values of $\Theta $. We provide the list and number of far fields, their matching conditions, the list and values of the matching coefficients, the corner profiles, etc.\\

Given the complexity of the asymptotic expansion, it would have been laborious to produce by hand the data given in this section. Instead, it was entirely generated by the symbolic computation module described in \cref{sec: implementation des espaces A}. This module allows extremely fast computations ; for example, all the tables given in \cref{sec: tableaux exemples DA} took only 0.1 second to be generated on a laptop. This does~not include finite element resolutions, which are much more costly computations.
\end{textesansboite}

\subsection{Genericity}
\label{sec: genericite}

\begin{textesansboite}
The structure of the asymptotic expansion of $u_{\peps}$ depends mostly on the angle $\Theta $. It does~not vary w.r.t the other parameters of the problem, except for exceptional values. To describe the expansion's structure regarless of these exceptions, we introduce the notion of \emph{genericity}. We say that a coefficient or a far field is \emph{generically} non-zero at a fixed angle $\Theta $, if it~is non-zero except for exceptional values of the parameters $(\mu _{0} ,\mu _{1} ,\rho _{0} ,\rho _{1} ,\omega )$, of the variational corner profiles coefficients $\sigma _{q} (\ccDA_{d,n} )$ (due to a particular geometry at the corner) or for very particular sources $\fsource$. The notion of genericity is common in many mathematical fields (for a diversified list of recent examples, see \cite{EleGehKanKatKel21, FelHenPas25, CalMarMerSte24, KabTen17, GolElaLod21}). For instance:
\begin{itemize}
    \item a generic triangle is~not flat, nor right-angled or isosceles,
    \item a generic polynomial has no multiple root, etc.
\end{itemize}
The definition below details rigorously what we consider generic in our context.
\end{textesansboite}

\begin{definitionsansboite}{}
\label{def: genericite}%
Let us take a fixed $\Theta $.
\begin{enumerate}
    \item Let $c$ be a coefficient of the type $\cmacro_{d,d',p,\ell }$ or $\cmicro_{d,d',p,\ell }$. It depends on $(\mu _{0} ,\mu _{1} ,\rho _{0} ,\rho _{1} ,\omega )$. We say that $c$ is \emph{generically} non-zero iff the set of values of $(\mu _{0} ,\mu _{1} ,\rho _{0} ,\rho _{1} ,\omega )$ s.t.\ $c=0$ has Lebesgue measure 0 in $(\mathbb{R}_+^*)^{4} \times \mathbb{C}$.
    \item Let $c$ be a coefficient of the type $\cfinal_{d,d',p,\ell }$. It depends on some coefficients of the type $\cmacro$ or $\cmicro$ and on some corner profile coefficients $(\sigma _{q_{i}} (\ccDA_{d_{i} ,n_{i} }))_{i\in [\![1,I]\!]}$ with $I\in \mathbb{N}$ and $q_{i} <0$ for~any $i$. So $c$ depends indirectly on $v:=(\mu _{0} ,\mu _{1} ,\rho _{0} ,\rho _{1} ,\omega ,$ $\sigma _{q_{1}} (\ccDA_{d_{1} ,n_{1} }),\dots ,\sigma _{q_{I}} (\ccDA_{d_{I} ,n_{I} }))$. We say that $c$ is \emph{generically} non-zero iff the set of values of $v$ s.t.\ $c=0$ has Lebesgue measure 0 in $(\mathbb{R}_+^*)^{4} \times \mathbb{C}\times \mathbb{R}^{I}$.
    \item We define the generic non-nullity of $u_{p,\ell }$ and $\sigma _{d} (\uDA_{p,\ell } )$ inductively on $p$.\vspace{-1mm}
    \begin{itemize}
        \item When $p\in \mathbb{R} \setminus \NN$, $u_{p,\ell } :=0$ by convention, so we say that the $u_{p,\ell }$ and $\sigma _{d} (\uDA_{p,\ell } )$ are~not generically non-zero.
        \item $u_{0,\ell }$ is generically non-zero iff $\ell =0$, and $\sigma _{d} (\uDA_{0,\ell } )$ is generically non-zero iff $\ell =0$ and $d>0$.
        \item Let $p\in \NN \setminus \{0\}$ and $\ell \in \mathbb{N}$. \Cref{construction directe des upl} states that, for~any $d\in -\TN$,
        \begin{equation}{\label{eq: def genericite : formule des CCL}}
            \sigma _{d} (\uDA_{p,\ell } ) =  \sum _{\substack{p'\in \NN\\ p-p'\in \NN+ 2\pst}} \; \sum _{\substack{d'\in \TZ*\\ p'+d'\in \intvl{0,p+d}}} \! \sum _{\ell '=0} ^{\ell } \cfinal_{d,d',p-p',\ell -\ell '} \cdot \sigma _{d'} (\uDA_{p',\ell '} ).
        \end{equation}
        We say that $\sigma _{d} (\uDA_{p,\ell } )$ is generically non-zero iff at least one term of this sum is s.t.\ $\cfinal_{d,d',p-p',\ell -\ell '}$ and $\sigma _{d'} (\uDA_{p',\ell '} )$ are both generically non-zero.\par 
        Then by \cref{construction directe des upl}, $u_{p,\ell }$ depends on the $\sigma _{d} (\uDA_{p,\ell } )$ with $d\in -\TN$ and on the $u_{p-1-2n,\ell }$ with $n\in \mathbb{N}$ (through the Dirichlet condition on $\Gamma $). We say that $u_{p,\ell } \neq 0$ generically iff at least one of these objects is generically non-zero.\par 
        Finally, for~any $d\in \TN$, $\sigma _{d} (\uDA_{p,\ell } ) \neq 0$ generically iff $u_{p,\ell } \neq 0$ generically.
    \end{itemize}
\end{enumerate}
\end{definitionsansboite}

\begin{remarque}
Unlike points 1 and 2, point 3 does~not rely on a measure property. 
Indeed the far fields depend on the source $\fsource \in H^{-1} (\Omega )$, and defining what is generic in an infinite dimensional vector space like $H^{-1} (\Omega )$ is non trivial, see \cite{Hun92}.
\end{remarque}

\begin{exemple}
\def\pst{\psttab}
Let $\pst:= \frac{\pi }{\Theta }$. It~is easy to check that $\cmacro_{-\pst,-\pst,\pst,0}=\cmicro_{\pst,\pst,\pst,0}=1$, so $\cfinal_{-\pst,\pst,0,0}:=\cmacro_{-\pst,-\pst,\pst,0} \cdot \sigma _{-\pst}(\ccDA_{\pst,0})\cdot \cmicro_{\pst,\pst,\pst,0}=\sigma _{-\pst}(\ccDA_{\pst,0})$. This is generically non-zero. Therefore $\sigma _{-\pst}(u_{2\pst,0})=\cfinal_{-\pst,\pst,0,0} \cdot \sigma _{\pst}(u_{0,0} )$ is generically non-zero, and so is $u_{2\pst,0}$.
\end{exemple}

\begin{remarque}
There is (at least) one natural setting that satisfies non-generic properties: the rectangular layer with $\Theta  = \frac{\pi }{2}$ (i.e.\ $\Omega _{1} := \mathbb{R}_+^*\times (-1,\infty )$). By~the method of images, solving this problem comes down to solving a Helmholtz problem in a half plane covered by an infinite thin layer. It~is well known that, in this situation without corner, the asymptotic expansion contains no $\ln \varepsilon $, see \cite{BenLem96,EngNed93}. Thus so do the expansion for the rectangular layer and $\Theta  = \frac{\pi }{2}$. However, as we will see in \cref{sec: graphes exemples DA} (\cref{tab: ordre du premier ln eps}), there~is generically some $\ln \varepsilon $ for this angle. In other words, the property ``$\exists p\in \NN, \exists \ell >0,\; u_{p,\ell } \neq 0$'' is generically true (for $\Theta = \frac{\pi }{2}$), but it does not hold in the very specific, yet natural, case of a rectangular layer.
\end{remarque}

\begin{remarque}
Although the matching conditions derived in \cref{chap: raccord algebrique} remain valid for the Poisson equation ($\omega :=0$), generic properties of the Helmholtz equation say nothing on what happens for the Poisson equation, as $\omega =0$ is non generic.
\end{remarque}

\begin{propositionsansboite}{}
\label{generiqmt non nul ou tjs nul}%
Let $w$ be a coefficient of the type $\cmacro_{d,d',p,\ell }$, $\cmicro_{d,d',p,\ell }$, $\cfinal_{d,d',p,\ell }$ or $\sigma _{d} (\uDA_{p,\ell } )$, or a far~field $u_{p,\ell }$. 
Then either $w$ is generically non-zero or $w$ is always zero (i.e.\ $w$ vanishes for any value of $(\mu _{0} ,\mu _{1} ,\rho _{0} ,\rho _{1} ,\omega ,\fsource)$ and of the corner profile coefficients $\sigma _{q} (\ccDA_{d,n} )$).
\end{propositionsansboite}

\begin{demo}
Let us denote $\mathcal{P} :=(\mu _{0} ,\mu _{1} ,\rho _{0} ,\rho _{1} ,\omega )$ for simplicity.
\begin{enumerate}
    \item Let us assume that $w =\cmicro_{d,d',p,\ell }$ for some $(d,d',p,\ell )$. Then, by~definition of $\cmicro_{d,d',p,\ell }$ and $\langle {\dots }\rangle $, we have\vspace{-1ex}
    \[w = \sigma _{d,p,\ell } \circ  \Eps \circ  \langle \Rmacro\rangle  (\phi_{d'} ) = \sigma _{d,p,\ell } \circ \Eps \Bigg( \fsumn \sum _{(R_{1} ,\dots ,R_{n} )\in \Rmacro^{n}} R_{1} \circ \cdots \circ R_{n} (\phi_{d'} ) \Bigg).\]
    Moreover, by construction of $\Rmacro$ in \cref{eq: def Rmaicro}, we have $\deg _{\mathcal{A}} R+2\deg _{\varepsilon } R>0$ for any $R\in \Rmacro$. Hence, there exists $n_{0}$ that depends only on $(d,d',p)$ and not on $\mathcal{P}$ s.t. 
    \[w = \sigma _{d,p,\ell } \circ \Eps \Bigg( \sum _{n=0} ^{n_{0}} \sum _{(R_{1} ,\dots ,R_{n} )\in \Rmacro^{n}} R_{1} \circ \cdots \circ R_{n} (\phi_{d'} ) \Bigg) .\]
    Now, given how $\Rmacro$ depends on $\mathcal{P}$, it~is easy to see that $F:\mathcal{P} \mapsto w$ is analytic. So, by properties of analytic functions, $F\neq 0$ iff $F$ is almost everywhere non-zero. This means that $w$ is not always zero iff it~is generically non-zero. Similarly, one can show this property if $w =\cmacro_{d,d',p,\ell }$.
    \item If $w =\cfinal_{d,d',p,\ell }$, then $w$ is an analytic function of some coefficients $\cmacro$, $\cmicro$ and some corner profile coefficients $(\sigma _{q_{1}} (\ccDA_{d_{1} ,n_{1} }),\dots ,\sigma _{q_{I}} (\ccDA_{d_{I} ,n_{I} })) =: \mathcal{Q}$. By~point 1, we deduce that $w$ is analytic in $(\mathcal{P},\mathcal{Q})$ and the same property of analytic functions can be used to conclude.
    \item Finally, let us assume that $w$ is equal to some $\sigma _{d} (\uDA_{p,\ell } )$ or $\uDA_{p,\ell }$. Then, using \cref{def: genericite}, it~is easy to check by induction on $p$ that, if $w$ is~not generically non-zero, then $w$ is always zero.
\end{enumerate}
\end{demo}

\begin{textesansboite}
The next result characterizes completely the powers of $\varepsilon $ that appear generically in the expansion of $u_{\peps}$.
\end{textesansboite}

\begin{proposition}{}
\label{ordre des CL generiquement non nuls}%
For~any $p\in \mathbb{N}+(\mathbb{N} \setminus \{1\})\pst$, $u_{p,0}$ is generically non-zero. For~any other $p$ and any $\ell \in \mathbb{N}$, $u_{p,\ell } =0$.
\end{proposition}

\begin{remarque}
This proposition says nothing about the fields $u_{p,\ell }$ with $p\in \mathbb{N}+(\mathbb{N} \setminus \{1\})\pst$ and $\ell \geqslant 1$. These fields determine which powers of $\ln \varepsilon $ appear in the expansion. Thanks to \cref{simplifications cas irrationnel}, we already know that those fields vanish if $\pst \not\in \mathbb{Q}$. However, when $\pst \in \mathbb{Q}$, the indexes $(p,\ell )\in (\mathbb{N}+(\mathbb{N} \setminus \{1\})\pst)\times \mathbb{N}^*$ s.t.\ $u_{p,\ell }$ vanishes depend a lot on the angle $\Theta $. \Cref{sec: tableaux exemples DA,sec: graphes exemples DA} will provide examples.
\end{remarque}

\begin{demo}
\def\pst{\tau}
It suffices to consider $p\in \NN$, as otherwise $u_{p,\ell } :=0$ by convention. Let us prove the property by induction on $p$. The case $p=0$ is obvious by \cref{def: genericite}. Let us prove the property for $p\in \NN \setminus \{0\}$. Let $\pst:= \frac{\pi }{\Theta }$ and $E:=\mathbb{N}+(\mathbb{N} \setminus \{1\})\pst$.\\
First let us assume $p\in E$. There exist $n\in \mathbb{N}$ and $m\in \mathbb{N} \setminus \{1\}$ s.t.\ $p=n+m\pst$.
\begin{itemize}
    \item If $n\geqslant 1$, then the Dirichlet condition on $\Gamma $ of $u_{p,0}$ involves $u_{p-1,0}$, which is generically non-zero by $p-1\in E$ and induction hypothesis.
    \item Otherwise, $m\geqslant 2$. In this case, \cref{eq: def genericite : formule des CCL} with $d:=-\pst$ contains a term $\cfinal_{-\pst,\pst,2\pst,0} \cdot \sigma _{\pst} (u_{p-2\pst,} 0)$. We know that $\sigma _{\pst} (u_{p-2\pst,0} )$ is generically non-zero by induction hypothesis because $p-2\pst \in E$. Hence it suffices to prove that $\cfinal_{-\pst,\pst,2\pst,0}$ is generically non-zero to deduce that $\sigma _{-\pst} (u_{p,0} )$ and $u_{p,0}$ are generically non-zero. Let us compute this coefficient:
    \begin{align*}
        \cfinal_{-\pst,\pst,2\pst,0} &:= \! \sum _{\substack{(p_{1} ,p_{2} )\in (\mathbb{N}+\pst)\times (\mathbb{N}+\pst)\\ n\in \mathbb{N},\, p_{1} +p_{2} +2n=2\pst}}\ \sum _{\substack{d_{1} \in \TZ*\\ p_{1} +d_{1} \in \mathbb{N}}} \ \sum _{\substack{d_{2} \in \TN\\ p_{2} -d_{2} \in \mathbb{N}}}	\cmacro_{-\pst,d_{1} ,p_{1} ,0}\cdot \omega ^{2n} \sigma _{d_{1}} (\mathcal{S} ^{\infty } _{d_{2} ,n} )\cdot \cmicro_{d_{2} ,\pst,p_{2} ,0} \\
        &= \sum _{\substack{d_{1} \in \TZ*\\ \pst+d_{1} \in \mathbb{N}}} \ \sum _{\substack{d_{2} \in \TN\\ \pst-d_{2} \in \mathbb{N}}}	\cmacro_{-\pst,d_{1} ,\pst,0}\cdot \sigma _{d_{1}} (\mathcal{S} ^{\infty } _{d_{2} ,0} )\cdot \cmicro_{d_{2} ,\pst,\pst,0} 
    \end{align*}
    because the condition in the first sum over $(p_{1} ,p_{2} ,n)$ implies $p_{1} =p_{2} =\pst$ and $n=0$. Then:
    \begin{alignjustif}[\quad]
        \cfinal_{-\pst,\pst,2\pst,0} &= \sum _{\substack{d_{1} \in \TZ*\\ \pst+d_{1} \in \mathbb{N}}} \cmacro_{-\pst,d_{1} ,\pst,0}\cdot \sigma _{d_{1}} (\mathcal{S} ^{\infty } _{\pst,0} )\cdot \cmicro_{\pst,\pst,\pst,0} & \text{because }  d_{2} =\pst \\[5ex]
        &= \cmacro_{-\pst,-\pst,\pst,0}\cdot \sigma _{-\pst} (\mathcal{S} ^{\infty } _{\pst,0} )\cdot \cmicro_{\pst,\pst,\pst,0} + \cmacro_{-\pst,\pst,\pst,0}\cdot \cmicro_{\pst,\pst,\pst,0} & \text{because }  \sigma _{d_{1}} (\mathcal{S} ^{\infty } _{\pst,0} ) =\delta _{d_{1} ,\pst} \text{ for }  d_{1} >0\\[1ex]
        &= \sigma _{-\pst} (\mathcal{S} ^{\infty } _{\pst,0} ) & \text{by \cref{classification coefs cmaicro}} 
    \end{alignjustif}
    So $\cfinal_{-\pst,\pst,2\pst,0}$ is generically non-zero, which proves that $u_{p,0}$ is generically non-zero.
\end{itemize}
Conversely, let us show that, if there exists $\ell \in \mathbb{N}$ s.t.\ $u_{p,\ell } \neq 0$, then $p\in E$. Without loss of generality, we can assume $p>0$. If $u_{p,\ell } \neq 0$, then either $u_{p,\ell |\Gamma } \neq 0$ or there exists $d\in -\TN$ s.t.\ $\sigma _{d} (u_{p,\ell } ) \neq 0$.
\begin{itemize}
    \item In the first case, since $u_{p,\ell |\Gamma }$ depends on the $u_{p-1-2n,\ell }$, there exists $n\in \mathbb{N}$ s.t.\ $u_{p-1-2n,\ell } \neq 0$. By~induction hypothesis, we have $p-1-2n\in E$, so $p\in E$.
    \item In the other case, there exist a term ``$\cfinal_{d,d',p-p',\ell -\ell '} \cdot \sigma _{d'} (\uDA_{p',\ell '} )$'' in \cref{eq: def genericite : formule des CCL} that does~not vanish. It satisfies $p-p'\in \NN+2\pst$ and $u_{p',\ell '} \neq 0$, so $p'\in E$. Thus, $p=(p-p')+p'\in E$.
\end{itemize}
\end{demo}

\begin{textesansboite}
\Cref{fig: ordres des CL} represents the orders $p$ s.t.\ at least one $u_{p,\ell }$ is generically non-zero, in terms of $\Theta $. These orders are given by \cref{ordre des CL generiquement non nuls}. One can see that, for small angles, they are only integer orders in the first terms of the expansion, whereas large angles involve many non integer orders.
\end{textesansboite}

	\tikzset{external/export=true}%
	\tikzsetfigurename{fig-ordres_des_CL-}%
	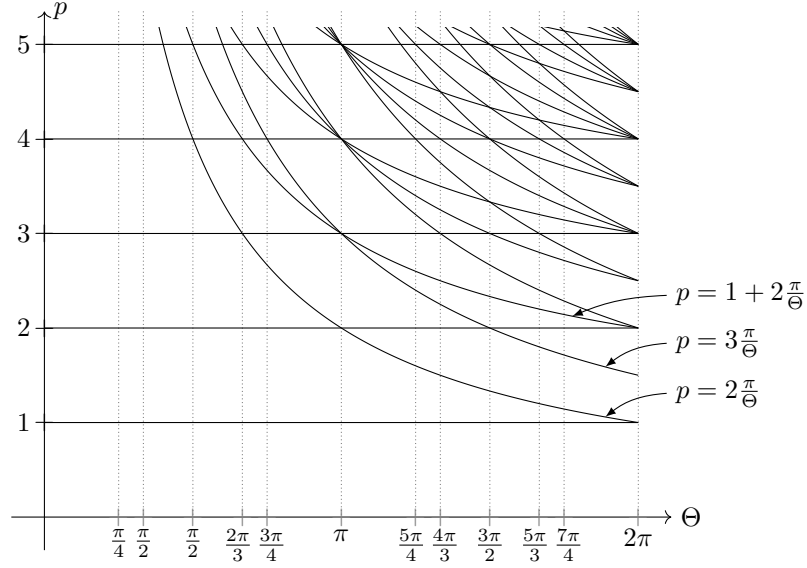
\begin{figure}[h]
\def\ymax{5}
\def\margeAxeFig{0.35}
\def\margeLegFig{\ }
\newcommand{\axevert}[2]{
	\draw [thin, gray, densely dotted] (#1,\ymax+\margeAxeFig) -- (#1,0) node {$+$} ++(0,-0.03) node [below] {\color{black}$#2$};}
\newcommand{\ptintersec}[1]{\draw (#1) node {\tiny$\bullet$};}

\begin{tikzpicture}[scale=1.25]
	\draw [->] (-\margeAxeFig,0) -- (2*pi+\margeAxeFig,0) node [right] {$\Theta$};
	\draw [->] (0,-\margeAxeFig) -- (0,\ymax+\margeAxeFig) node [right] {$p$};

    \foreach \i in {1,2,...,5} {
        \draw (0,\i) node{$+$} node[left]{\i\,}; }
    \axevert{1/4*pi}{\frac\pi4}
    \axevert{1/3*pi}{\frac\pi2}
    \axevert{1/2*pi}{\frac\pi2}
    \axevert{2/3*pi}{\frac{2\pi}3\margeLegFig}
    \axevert{3/4*pi}{\margeLegFig\frac{3\pi}4}
    \axevert{pi}{\pi}
    \axevert{5/4*pi}{\frac{5\pi}4\margeLegFig}
    \axevert{4/3*pi}{\margeLegFig\frac{4\pi}3}
    \axevert{3/2*pi}{\frac{3\pi}2}
    \axevert{5/3*pi}{\frac{5\pi}3\margeLegFig}
    \axevert{7/4*pi}{\margeLegFig\frac{7\pi}4}
    \axevert{2*pi}{2\pi}

    \foreach \i in {1,2,...,5} {
	    \draw (0,\i) -- (2*pi,\i); }
	\begin{scope}
		\clip (0,0) rectangle (2*pi,\ymax+\margeAxeFig/2);
        \foreach \i in {2,3,...,10} {
            \foreach \j in {0,1,...,4} {
                \draw [domain=0.1:{2*pi},samples=200] plot(\x,{\i*pi/\x+\j});}}
    \end{scope}

    \def\ThetaLegFig{(2*pi-\dxLegFig)}
    \def\dxLegFig{0.35}
    \draw[<-,>=latex] ({\ThetaLegFig},{2*pi/\ThetaLegFig}) to[bend left=17] ++(0.3+\dxLegFig,0.25) node[right]{$p=2\pst$};
    \draw[<-,>=latex] ({\ThetaLegFig},{3*pi/\ThetaLegFig}) to[bend left=17] ++(0.3+\dxLegFig,0.25) node[right]{$p=3\pst$};
    \def\dxLegFig{0.7}
    \draw[<-,>=latex] ({\ThetaLegFig},{1+2*pi/\ThetaLegFig}) to[bend left=10] ++(0.3+\dxLegFig,0.22) node[right]{$p=1+2\pst$};

    \draw (-2,0) node{\ }; 
\end{tikzpicture}%
\caption{
    Points $(\Theta,p)$ s.t.\ at least one $u_{p,\ell}$ is generically non-zero for the angle $\Theta$
}
\label{fig: ordres des CL}
\end{figure}%
	\tikzset{external/export=false}%
	\tikzsetfigurename{main-figure}

\begin{remarque}
Things are different of the corner fields: one can show that $S_{0,\ell }$ always vanishes whereas $S_{\pst,0}$ is generically non-zero.
\end{remarque}

\subsection{Detailed examples for a few angles}
\label{sec: tableaux exemples DA}

\begin{textesansboite}
In \crefrange{tab: liste CL pour 2.5}{tab: liste cfinal pour pi} we detail the structure of the asymptotic expansion of $u_{\peps}$ up to order $p=4$ for three different angles: $\Theta =2.5$, $\Theta =4$ and $\Theta =\pi $. For each angle, we provide 4 tables. The first one enumerates the generically non-zero far fields $u_{p,\ell }$ with $p\leqslant 4$ and their corner conditions, i.e.\ the generically non-zero coefficients $\sigma _{d} (u_{p,\ell } )$ with $d\in -\TN$ and their expression. The second table enumerates the corner profiles and their coefficients $\sigma _{q} (\ccDA_{d,n} )$ that need to be computed to build the $u_{p,\ell }$ for $p\leqslant 4$. The third table enumerates the generically non-zero coefficients $\cmacro_{d,d',p,\ell }$ and $\cmicro_{d,d',p,\ell }$ involved in the asymptotic expansion up to order 4. We give an example of their values for a random choice of the parameters $(\mu _{0} ,\mu _{1} ,\rho _{0} ,\rho _{1} ,\omega )$. The last table enumerates the generically non-zero coefficients $\cfinal_{d,d',p,\ell }$ involved in the asymptotic expansion up to order 4 and provide their expression.
\end{textesansboite}

\begin{paragraphesansboite}{Notation:}
To avoid too heavy formulas in the tables, we denote $\psttab:=\pst$ in this section.
\end{paragraphesansboite}
\nomenclature[54]{$\tau $}{In \cref{sec: tableaux exemples DA}, $\pi /\Theta $}

\begin{textesansboite}
The expressions of the $\sigma _{d} (u_{p,\ell } )$ and $\cfinal_{d,d',p,\ell }$ change little when only $\ell $ is modified. So, to save space, we write them for all values of $\ell $ at once in \crefrange{tab: liste CL pour 2.5}{tab: liste cfinal pour pi}. The values of $\ell $ for which the involved coefficient are generically non-zero are indicated in index, like this: $\sigma _{\psttab} (u_{2,\ell \leqslant 1} )$, $\cmacro_{-\psttab,\psttab,2,\ell =1}$ and $\cmicro_{-\psttab,\psttab,2,0}$ (if only $\ell =0$ gives a non-zero value). When two coefficients that may be non-zero for $\ell >0$ are multiplied with the symbol ``$*$'', it should be understood as a convolution product w.r.t.\ $\ell $. For example, 
\[``\cfinal_{-2\psttab,\psttab,2,\ell \leqslant 1} * \sigma _{\psttab} (u_{2,\ell \leqslant 1} )\text{''}  \qquad  \mathrm{means} \qquad  \sum _{(\ell _{1} ,\ell _{2} )\in \{0,1\}^2 , \ell _{1} +\ell _{2} =\ell }
\cfinal_{-2\psttab,\psttab,2,\ell _{1}} \cdot  \sigma _{\psttab} (u_{2,\ell _{2}} ).\]

\Crefrange{tab: liste CL pour 2.5}{tab: liste cfinal pour pi} are large, so let us point out the interesting observations. First of all, the reader will probably be relieved to see that the expressions of the coefficients $\cfinal$ and $\sigma _{d} (u_{p,\ell } )$ with $d<0$, which are very heavy sums over up to 7 variables (see \cref{sec: construction directe des upl}), only contain a few terms in practice. Next, the cases $\Theta =2.5$ and $\Theta =4$ (\crefrange{tab: liste CL pour 2.5}{tab: liste cfinal pour 4}) show all the simplifications predicted by \cref{simplifications cas irrationnel} when $\pst\not\in \mathbb{Q}$: there~is no $\ln \varepsilon $, the coefficients $\cmacro$ and $\cmicro$ are all equal to 0 or 1 and the $\cfinal$ have simple expressions. Moreover, the greater the angle, the heavier the expansion. For $\Theta =4$, the first corner condition for a far field appears earlier than for $\Theta =2.5$ ($p\approx 1.57$ versus $2.51$) and the expansion involves much more far fields, corner profiles and matching coefficients.\\

In contrast, $\Theta =\pi $ is a good example of the complexity of the expansion in the rational case. The denominator of $\frac{\Theta }{\pi }$ is minimal, so $\ln \varepsilon $ appears in the expansion very soon (at $p=2$) and the expressions of the coefficients $\cfinal$ are quite complex even at small orders. This is actually the angle with the most intricate expansion at small orders, as angles with a larger denominator behave like the irrational case up to greater orders according to \cref{simplifications cas rationnel}. But the flat angle is very relevant for applications. However, with a small denominator, the orders $p$ s.t.\ there~is a generically non-zero far field $u_{p,\ell }$ are simpler: instead of having irrational orders like for $\Theta \in \{2.5,4\}$, the far field orders are all rational, and even integers in our example. Hence, for~any angle $\Theta $, there~is always some simplifications compared to the general results stated in \cref{chap: construction du DA}.\\

Finally \cref{tab: liste cmaicro pour pi} is coherent with the classification of non-zero coefficients $\cmacro$ and $\cmicro$ given in \cref{classification coefs cmaicro}. In addition, as expected, the coefficients of the second category (i.e.\ other than the $\cmacro_{d,d,-d,0}$ and $\cmicro_{d,d,d,0}$) are generically different from 1. Some of them are real, other ones are complex, depending on whether they depend on $\omega $ or not.
\end{textesansboite}

%
	\tikzset{external/export=true}%
	\tikzsetfigurename{tab-structure_DA_2,5-}%
	\renewcommand{\arraystretch}{1.25}

\begin{table}[h]
\caption{List of far fields $u_{p,\ell}$ with $p\leqslant4$ that are generically non-zero for $\Theta=2.5$}
\label{tab: liste CL pour 2.5}
\begin{tabular}{|c|c|>{\centering\arraybackslash}m{19ex}|l|}
\hline
$p$ & $\ell$ & cumulated number of far fields & \multicolumn{1}{c|}{corner conditions} \\
\hline
$0$ & 0 & 1 & \multicolumn{1}{c|}{none} \\
\hline
$1$ & 0 & 2 & \multicolumn{1}{c|}{none} \\
\hline
$2$ & 0 & 3 & \multicolumn{1}{c|}{none} \\
\hline
$2\psttab \approx 2.51$ & 0 & 4 & $\begin{array}{@{}l@{\;}l@{\;}l@{\;}l@{}l@{\;}l@{}} \sigma_{-\psttab}(u_{2\psttab,0})  & =\; \cfinal_{-\psttab,\psttab,2\psttab,0}  & \cdot\; \sigma_{\psttab}(u_{0,0}) \\
[0.5ex] \end{array}$ \\
\hline
$3$ & 0 & 5 & \multicolumn{1}{c|}{none} \\
\hline
$1 + 2\psttab \approx 3.51$ & 0 & 6 & $\begin{array}{@{}l@{\;}l@{\;}l@{\;}l@{}l@{\;}l@{}} \sigma_{-\psttab}(u_{1 + 2\psttab,0})  & =\; \cfinal_{-\psttab,\psttab,2\psttab,0}  & \cdot\; \sigma_{\psttab}(u_{1,0}) \\
[0.5ex] \end{array}$ \\
\hline
$3\psttab \approx 3.77$ & 0 & 7 & $\begin{array}{@{}l@{\;}l@{\;}l@{\;}l@{}l@{\;}l@{}} \sigma_{-\psttab}(u_{3\psttab,0})  & =\; \cfinal_{-\psttab,2\psttab,3\psttab,0}  & \cdot\; \sigma_{2\psttab}(u_{0,0}) \\
\sigma_{-2\psttab}(u_{3\psttab,0})  & =\; \cfinal_{-2\psttab,\psttab,3\psttab,0}  & \cdot\; \sigma_{\psttab}(u_{0,0}) \\
[0.5ex] \end{array}$ \\
\hline
$4$ & 0 & 8 & \multicolumn{1}{c|}{none} \\
\hline
\end{tabular}
\end{table}

\begin{table}[h]
\caption{List of the coefficients $\sigma_q(\mathcal{S}^\infty_{d,n})$ that are needed to compute the $u_{p,\ell}$ with $p\leqslant4$ and $\Theta=2.5$}
\label{tab: liste CCC pour 2.5}
\begin{tabular}{|c|c|l|}
\hline
$d$ & $n$ & \multicolumn{1}{c|}{$q$} \\
\hline
$\psttab$ & 0 & $-\psttab$, $-2\psttab$ \\
\hline
$2\psttab$ & 0 & $-\psttab$ \\
\hline
\end{tabular}
\end{table}

\renewcommand{\arraystretch}{1.2}

\begin{table}[h]
\caption{List of the 4 non-zero coefficients $\cmacro_{d,d',p,\ell}$ and $\cmicro_{d,d',p,\ell}$ needed to compute the $u_{p,\ell}$ with $p\leqslant4$ and $\Theta=2.5$}
\label{tab: liste cmaicro pour 2.5}
\begin{tabular}[t]{|c|c|c|c|c|}
\hline
$d$ & $d'$ & $p$ & $\ell$ & $\cmacro_{d,d',p,\ell}$\rule[-1.4ex]{0pt}{4ex} \\
\hline
$-2\psttab$ & $-2\psttab$ & $2\psttab$ & 0 & $1$ \\
\hline
$-\psttab$ & $-\psttab$ & $\psttab$ & 0 & $1$ \\
\hline
\end{tabular}
\hspace{10mm}
\begin{tabular}[t]{|c|c|c|c|c|}
\hline
$d$ & $d'$ & $p$ & $\ell$ & $\cmicro_{d,d',p,\ell}$\rule[-1.4ex]{0pt}{4ex} \\
\hline
$\psttab$ & $\psttab$ & $\psttab$ & 0 & $1$ \\
\hline
$2\psttab$ & $2\psttab$ & $2\psttab$ & 0 & $1$ \\
\hline
\end{tabular}
\end{table}
\renewcommand{\arraystretch}{1.25}

\begin{table}[h]
\caption{List of the 3 generically non-zero coefficients $\cfinal_{d,d',p,\ell}$ needed to compute the $u_{p,\ell}$ with $p\leqslant4$ and $\Theta=2.5$}
\label{tab: liste cfinal pour 2.5}
\begin{tabular}{|c|c|c|c|l|}
\hline
$d$ & $d'$ & $p$ & $\ell$ & \multicolumn{1}{c|}{expression of $\cfinal_{d,d',p,\ell}$} \\
\hline
$-2\psttab$ & $\psttab$ & $3\psttab$ & 0 & $\begin{array}{@{}l@{\;}l@{\;}l@{\;}l@{\;}l@{\;}l@{}}
\rule{0pt}{2.8ex} \!\! =\; \cmacro_{-2\psttab,-2\psttab,2\psttab,0}  & \cdot\; \sigma_{-2\psttab}(\mathcal{S}^\infty_{\psttab,0})  & \cdot\; \cmicro_{\psttab,\psttab,\psttab,0} \\
[0.8ex] \end{array}$ \\
\hline
\multirow{2}{*}{$-\psttab$} & $\psttab$ & $2\psttab$ & 0 & $\begin{array}{@{}l@{\;}l@{\;}l@{\;}l@{\;}l@{\;}l@{}}
\rule{0pt}{2.8ex} \!\! =\; \cmacro_{-\psttab,-\psttab,\psttab,0}  & \cdot\; \sigma_{-\psttab}(\mathcal{S}^\infty_{\psttab,0})  & \cdot\; \cmicro_{\psttab,\psttab,\psttab,0} \\
[0.8ex] \end{array}$ \\
\cline{2-5}
 & $2\psttab$ & $3\psttab$ & 0 & $\begin{array}{@{}l@{\;}l@{\;}l@{\;}l@{\;}l@{\;}l@{}}
\rule{0pt}{2.8ex} \!\! =\; \cmacro_{-\psttab,-\psttab,\psttab,0}  & \cdot\; \sigma_{-\psttab}(\mathcal{S}^\infty_{2\psttab,0})  & \cdot\; \cmicro_{2\psttab,2\psttab,2\psttab,0} \\
[0.8ex] \end{array}$ \\
\hline
\end{tabular}
\end{table}

	\tikzset{external/export=false}%
	\tikzsetfigurename{main-figure}
	\tikzset{external/export=true}%
	\tikzsetfigurename{tab-structure_DA_4-}%
	\renewcommand{\arraystretch}{1.25}

\begin{table}[h]
\caption{List of far fields $u_{p,\ell}$ with $p\leqslant4$ that are generically non-zero for $\Theta=4$}
\label{tab: liste CL pour 4}
\begin{tabular}{|c|c|>{\centering\arraybackslash}m{11ex}|l|}
\hline
$p$ & $\ell$ & cumulated number of far fields & \multicolumn{1}{c|}{corner conditions} \\
\hline
$0$ & 0 & 1 & \multicolumn{1}{c|}{none} \\
\hline
$1$ & 0 & 2 & \multicolumn{1}{c|}{none} \\
\hline
$2\psttab \approx 1.57$ & 0 & 3 & $\begin{array}{@{}l@{\;}l@{\;}l@{\;}l@{}l@{\;}l@{}} \sigma_{-\psttab}(u_{2\psttab,0})  & =\; \cfinal_{-\psttab,\psttab,2\psttab,0}  & \cdot\; \sigma_{\psttab}(u_{0,0}) \\
[0.5ex] \end{array}$ \\
\hline
$2$ & 0 & 4 & \multicolumn{1}{c|}{none} \\
\hline
$3\psttab \approx 2.36$ & 0 & 5 & $\begin{array}{@{}l@{\;}l@{\;}l@{\;}l@{}l@{\;}l@{}} \sigma_{-\psttab}(u_{3\psttab,0})  & =\; \cfinal_{-\psttab,2\psttab,3\psttab,0}  & \cdot\; \sigma_{2\psttab}(u_{0,0}) \\
\sigma_{-2\psttab}(u_{3\psttab,0})  & =\; \cfinal_{-2\psttab,\psttab,3\psttab,0}  & \cdot\; \sigma_{\psttab}(u_{0,0}) \\
[0.5ex] \end{array}$ \\
\hline
$1 + 2\psttab \approx 2.57$ & 0 & 6 & $\begin{array}{@{}l@{\;}l@{\;}l@{\;}l@{}l@{\;}l@{}} \sigma_{-\psttab}(u_{1 + 2\psttab,0})  & =\; \cfinal_{-\psttab,\psttab,2\psttab,0}  & \cdot\; \sigma_{\psttab}(u_{1,0}) \\
[0.5ex] \end{array}$ \\
\hline
$3$ & 0 & 7 & \multicolumn{1}{c|}{none} \\
\hline
$4\psttab \approx 3.14$ & 0 & 8 & $\begin{array}{@{}l@{\;}l@{\;}l@{\;}l@{}l@{\;}l@{}} \sigma_{-\psttab}(u_{4\psttab,0})  & =\; \cfinal_{-\psttab,3\psttab,4\psttab,0}  & \cdot\; \sigma_{3\psttab}(u_{0,0}) &  & +\; \cfinal_{-\psttab,\psttab,2\psttab,0}  & \cdot\; \sigma_{\psttab}(u_{2\psttab,0}) \\
\sigma_{-2\psttab}(u_{4\psttab,0})  & =\; \cfinal_{-2\psttab,2\psttab,4\psttab,0}  & \cdot\; \sigma_{2\psttab}(u_{0,0}) \\
\sigma_{-3\psttab}(u_{4\psttab,0})  & =\; \cfinal_{-3\psttab,\psttab,4\psttab,0}  & \cdot\; \sigma_{\psttab}(u_{0,0}) \\
[0.5ex] \end{array}$ \\
\hline
$1 + 3\psttab \approx 3.36$ & 0 & 9 & $\begin{array}{@{}l@{\;}l@{\;}l@{\;}l@{}l@{\;}l@{}} \sigma_{-\psttab}(u_{1 + 3\psttab,0})  & =\; \cfinal_{-\psttab,2\psttab,3\psttab,0}  & \cdot\; \sigma_{2\psttab}(u_{1,0}) \\
\sigma_{-2\psttab}(u_{1 + 3\psttab,0})  & =\; \cfinal_{-2\psttab,\psttab,3\psttab,0}  & \cdot\; \sigma_{\psttab}(u_{1,0}) \\
[0.5ex] \end{array}$ \\
\hline
$2 + 2\psttab \approx 3.57$ & 0 & 10 & $\begin{array}{@{}l@{\;}l@{\;}l@{\;}l@{}l@{\;}l@{}} \sigma_{-\psttab}(u_{2 + 2\psttab,0})  & =\; \cfinal_{-\psttab,\psttab,2 + 2\psttab,0}  & \cdot\; \sigma_{\psttab}(u_{0,0}) &  & +\; \cfinal_{-\psttab,\psttab,2\psttab,0}  & \cdot\; \sigma_{\psttab}(u_{2,0}) \\
[0.5ex] \end{array}$ \\
\hline
$5\psttab \approx 3.93$ & 0 & 11 & $\begin{array}{@{}l@{\;}l@{\;}l@{\;}l@{}l@{\;}l@{}} \sigma_{-\psttab}(u_{5\psttab,0})  & =\; \cfinal_{-\psttab,4\psttab,5\psttab,0}  & \cdot\; \sigma_{4\psttab}(u_{0,0}) &  & +\; \cfinal_{-\psttab,2\psttab,3\psttab,0}  & \cdot\; \sigma_{2\psttab}(u_{2\psttab,0}) \\
 & +\; \cfinal_{-\psttab,\psttab,2\psttab,0}  & \cdot\; \sigma_{\psttab}(u_{3\psttab,0}) &  &  &  \\
\sigma_{-2\psttab}(u_{5\psttab,0})  & =\; \cfinal_{-2\psttab,3\psttab,5\psttab,0}  & \cdot\; \sigma_{3\psttab}(u_{0,0}) &  & +\; \cfinal_{-2\psttab,\psttab,3\psttab,0}  & \cdot\; \sigma_{\psttab}(u_{2\psttab,0}) \\
\sigma_{-3\psttab}(u_{5\psttab,0})  & =\; \cfinal_{-3\psttab,2\psttab,5\psttab,0}  & \cdot\; \sigma_{2\psttab}(u_{0,0}) \\
\sigma_{-4\psttab}(u_{5\psttab,0})  & =\; \cfinal_{-4\psttab,\psttab,5\psttab,0}  & \cdot\; \sigma_{\psttab}(u_{0,0}) \\
[0.5ex] \end{array}$ \\
\hline
$4$ & 0 & 12 & \multicolumn{1}{c|}{none} \\
\hline
\end{tabular}
\end{table}

\begin{table}[h]
\caption{List of the coefficients $\sigma_q(\mathcal{S}^\infty_{d,n})$ that are needed to compute the $u_{p,\ell}$ with $p\leqslant4$ and $\Theta=4$}
\label{tab: liste CCC pour 4}
\begin{tabular}{|c|c|l|}
\hline
$d$ & $n$ & \multicolumn{1}{c|}{$q$} \\
\hline
\multirow{2}{*}{$\psttab$} & 0 & $-\psttab$, $-2\psttab$, $-3\psttab$, $-4\psttab$ \\
\cline{2-3}
 & 1 & $-\psttab$ \\
\hline
$2\psttab$ & 0 & $-\psttab$, $-2\psttab$, $-3\psttab$ \\
\hline
$3\psttab$ & 0 & $-\psttab$, $-2\psttab$ \\
\hline
$4\psttab$ & 0 & $-\psttab$ \\
\hline
\end{tabular}
\end{table}

\renewcommand{\arraystretch}{1.2}

\begin{table}[h]
\caption{List of the 8 non-zero coefficients $\cmacro_{d,d',p,\ell}$ and $\cmicro_{d,d',p,\ell}$ needed to compute the $u_{p,\ell}$ with $p\leqslant4$ and $\Theta=4$}
\label{tab: liste cmaicro pour 4}
\begin{tabular}[t]{|c|c|c|c|c|}
\hline
$d$ & $d'$ & $p$ & $\ell$ & $\cmacro_{d,d',p,\ell}$\rule[-1.4ex]{0pt}{4ex} \\
\hline
$-4\psttab$ & $-4\psttab$ & $4\psttab$ & 0 & $1$ \\
\hline
$-3\psttab$ & $-3\psttab$ & $3\psttab$ & 0 & $1$ \\
\hline
$-2\psttab$ & $-2\psttab$ & $2\psttab$ & 0 & $1$ \\
\hline
$-\psttab$ & $-\psttab$ & $\psttab$ & 0 & $1$ \\
\hline
\end{tabular}
\hspace{10mm}
\begin{tabular}[t]{|c|c|c|c|c|}
\hline
$d$ & $d'$ & $p$ & $\ell$ & $\cmicro_{d,d',p,\ell}$\rule[-1.4ex]{0pt}{4ex} \\
\hline
$\psttab$ & $\psttab$ & $\psttab$ & 0 & $1$ \\
\hline
$2\psttab$ & $2\psttab$ & $2\psttab$ & 0 & $1$ \\
\hline
$3\psttab$ & $3\psttab$ & $3\psttab$ & 0 & $1$ \\
\hline
$4\psttab$ & $4\psttab$ & $4\psttab$ & 0 & $1$ \\
\hline
\end{tabular}
\end{table}
\renewcommand{\arraystretch}{1.25}

\begin{table}[h]
\caption{List of the 11 generically non-zero coefficients $\cfinal_{d,d',p,\ell}$ needed to compute the $u_{p,\ell}$ with $p\leqslant4$ and $\Theta=4$}
\label{tab: liste cfinal pour 4}
\begin{tabular}{|c|c|c|c|l|}
\hline
$d$ & $d'$ & $p$ & $\ell$ & \multicolumn{1}{c|}{expression of $\cfinal_{d,d',p,\ell}$} \\
\hline
$-4\psttab$ & $\psttab$ & $5\psttab$ & 0 & $\begin{array}{@{}l@{\;}l@{\;}l@{\;}l@{\;}l@{\;}l@{}}
\rule{0pt}{2.8ex} \!\! =\; \cmacro_{-4\psttab,-4\psttab,4\psttab,0}  & \cdot\; \sigma_{-4\psttab}(\mathcal{S}^\infty_{\psttab,0})  & \cdot\; \cmicro_{\psttab,\psttab,\psttab,0} \\
[0.8ex] \end{array}$ \\
\hline
\multirow{2}{*}[-0.5ex]{$-3\psttab$} & $\psttab$ & $4\psttab$ & 0 & $\begin{array}{@{}l@{\;}l@{\;}l@{\;}l@{\;}l@{\;}l@{}}
\rule{0pt}{2.8ex} \!\! =\; \cmacro_{-3\psttab,-3\psttab,3\psttab,0}  & \cdot\; \sigma_{-3\psttab}(\mathcal{S}^\infty_{\psttab,0})  & \cdot\; \cmicro_{\psttab,\psttab,\psttab,0} \\
[0.8ex] \end{array}$ \\
\cline{2-5}
 & $2\psttab$ & $5\psttab$ & 0 & $\begin{array}{@{}l@{\;}l@{\;}l@{\;}l@{\;}l@{\;}l@{}}
\rule{0pt}{2.8ex} \!\! =\; \cmacro_{-3\psttab,-3\psttab,3\psttab,0}  & \cdot\; \sigma_{-3\psttab}(\mathcal{S}^\infty_{2\psttab,0})  & \cdot\; \cmicro_{2\psttab,2\psttab,2\psttab,0} \\
[0.8ex] \end{array}$ \\
\hline
\multirow{3}{*}[-1ex]{$-2\psttab$} & $\psttab$ & $3\psttab$ & 0 & $\begin{array}{@{}l@{\;}l@{\;}l@{\;}l@{\;}l@{\;}l@{}}
\rule{0pt}{2.8ex} \!\! =\; \cmacro_{-2\psttab,-2\psttab,2\psttab,0}  & \cdot\; \sigma_{-2\psttab}(\mathcal{S}^\infty_{\psttab,0})  & \cdot\; \cmicro_{\psttab,\psttab,\psttab,0} \\
[0.8ex] \end{array}$ \\
\cline{2-5}
 & $2\psttab$ & $4\psttab$ & 0 & $\begin{array}{@{}l@{\;}l@{\;}l@{\;}l@{\;}l@{\;}l@{}}
\rule{0pt}{2.8ex} \!\! =\; \cmacro_{-2\psttab,-2\psttab,2\psttab,0}  & \cdot\; \sigma_{-2\psttab}(\mathcal{S}^\infty_{2\psttab,0})  & \cdot\; \cmicro_{2\psttab,2\psttab,2\psttab,0} \\
[0.8ex] \end{array}$ \\
\cline{2-5}
 & $3\psttab$ & $5\psttab$ & 0 & $\begin{array}{@{}l@{\;}l@{\;}l@{\;}l@{\;}l@{\;}l@{}}
\rule{0pt}{2.8ex} \!\! =\; \cmacro_{-2\psttab,-2\psttab,2\psttab,0}  & \cdot\; \sigma_{-2\psttab}(\mathcal{S}^\infty_{3\psttab,0})  & \cdot\; \cmicro_{3\psttab,3\psttab,3\psttab,0} \\
[0.8ex] \end{array}$ \\
\hline
\multirow{5}{*}[-2.5ex]{$-\psttab$} & \multirow{2}{*}{$\psttab$} & $2\psttab$ & 0 & $\begin{array}{@{}l@{\;}l@{\;}l@{\;}l@{\;}l@{\;}l@{}}
\rule{0pt}{2.8ex} \!\! =\; \cmacro_{-\psttab,-\psttab,\psttab,0}  & \cdot\; \sigma_{-\psttab}(\mathcal{S}^\infty_{\psttab,0})  & \cdot\; \cmicro_{\psttab,\psttab,\psttab,0} \\
[0.8ex] \end{array}$ \\
\cline{3-5}
 &  & $2 + 2\psttab$ & 0 & $\begin{array}{@{}l@{\;}l@{\;}l@{\;}l@{\;}l@{\;}l@{}}
\rule{0pt}{2.8ex} \!\! =\; \cmacro_{-\psttab,-\psttab,\psttab,0}  & \cdot\; \omega^2 \sigma_{-\psttab}(\mathcal{S}^\infty_{\psttab,1})  & \cdot\; \cmicro_{\psttab,\psttab,\psttab,0} \\
[0.8ex] \end{array}$ \\
\cline{2-5}
 & $2\psttab$ & $3\psttab$ & 0 & $\begin{array}{@{}l@{\;}l@{\;}l@{\;}l@{\;}l@{\;}l@{}}
\rule{0pt}{2.8ex} \!\! =\; \cmacro_{-\psttab,-\psttab,\psttab,0}  & \cdot\; \sigma_{-\psttab}(\mathcal{S}^\infty_{2\psttab,0})  & \cdot\; \cmicro_{2\psttab,2\psttab,2\psttab,0} \\
[0.8ex] \end{array}$ \\
\cline{2-5}
 & $3\psttab$ & $4\psttab$ & 0 & $\begin{array}{@{}l@{\;}l@{\;}l@{\;}l@{\;}l@{\;}l@{}}
\rule{0pt}{2.8ex} \!\! =\; \cmacro_{-\psttab,-\psttab,\psttab,0}  & \cdot\; \sigma_{-\psttab}(\mathcal{S}^\infty_{3\psttab,0})  & \cdot\; \cmicro_{3\psttab,3\psttab,3\psttab,0} \\
[0.8ex] \end{array}$ \\
\cline{2-5}
 & $4\psttab$ & $5\psttab$ & 0 & $\begin{array}{@{}l@{\;}l@{\;}l@{\;}l@{\;}l@{\;}l@{}}
\rule{0pt}{2.8ex} \!\! =\; \cmacro_{-\psttab,-\psttab,\psttab,0}  & \cdot\; \sigma_{-\psttab}(\mathcal{S}^\infty_{4\psttab,0})  & \cdot\; \cmicro_{4\psttab,4\psttab,4\psttab,0} \\
[0.8ex] \end{array}$ \\
\hline
\end{tabular}
\end{table}

	\tikzset{external/export=false}%
	\tikzsetfigurename{main-figure}
%
	\tikzset{external/export=true}%
	\tikzsetfigurename{tab-structure_DA_pi-}%
	\renewcommand{\arraystretch}{1.25}

\begin{table}[h]
\caption{List of far fields $u_{p,\ell}$ with $p\leqslant4$ that are generically non-zero for $\Theta=\pi$}
\label{tab: liste CL pour pi}
\begin{tabular}{|c|c|>{\centering\arraybackslash}m{11ex}|l|}
\hline
$p$ & $\ell$ & cumulated number of far fields & \multicolumn{1}{c|}{corner conditions} \\
\hline
$0$ & 0 & 1 & \multicolumn{1}{c|}{none} \\
\hline
$1$ & 0 & 2 & \multicolumn{1}{c|}{none} \\
\hline
$2$ & 0, 1 & 4 & $\begin{array}{@{}l@{\;}l@{\;}l@{\;}l@{}l@{\;}l@{}} \sigma_{-\psttab}(u_{2,\ell})  & =\; \cfinal_{-\psttab,\psttab,2,\ell\leqslant1}  & \cdot\; \sigma_{\psttab}(u_{0,0}) \\
[0.5ex] \end{array}$ \\
\hline
$3$ & 0, 1, 2 & 7 & $\begin{array}{@{}l@{\;}l@{\;}l@{\;}l@{}l@{\;}l@{}} \sigma_{-\psttab}(u_{3,\ell})  & =\; \cfinal_{-\psttab,2\psttab,3,\ell\leqslant2}  & \cdot\; \sigma_{2\psttab}(u_{0,0}) &  & +\; \cfinal_{-\psttab,\psttab,2,\ell\leqslant1}  & \cdot\; \sigma_{\psttab}(u_{1,0}) \\
\sigma_{-2\psttab}(u_{3,\ell})  & =\; \cfinal_{-2\psttab,\psttab,3,\ell\leqslant2}  & \cdot\; \sigma_{\psttab}(u_{0,0}) \\
[0.5ex] \end{array}$ \\
\hline
$4$ & 0, 1, 2, 3 & 11 & $\begin{array}{@{}l@{\;}l@{\;}l@{\;}l@{}l@{\;}l@{}} \sigma_{-\psttab}(u_{4,\ell})  & =\; \cfinal_{-\psttab,\psttab,4,\ell\leqslant3}  & \cdot\; \sigma_{\psttab}(u_{0,0}) &  & +\; \cfinal_{-\psttab,3\psttab,4,\ell\leqslant3}  & \cdot\; \sigma_{3\psttab}(u_{0,0}) \\
 & +\; \cfinal_{-\psttab,2\psttab,3,\ell\leqslant2}  & \cdot\; \sigma_{2\psttab}(u_{1,0}) &  & +\; \cfinal_{-\psttab,-\psttab,2,\ell\leqslant2}  & *\; \sigma_{-\psttab}(u_{2,\ell\leqslant1}) \\
 & +\; \cfinal_{-\psttab,\psttab,2,\ell\leqslant1}  & *\; \sigma_{\psttab}(u_{2,\ell\leqslant1}) &  &  &  \\
\sigma_{-2\psttab}(u_{4,\ell})  & =\; \cfinal_{-2\psttab,2\psttab,4,\ell\leqslant3}  & \cdot\; \sigma_{2\psttab}(u_{0,0}) &  & +\; \cfinal_{-2\psttab,\psttab,3,\ell\leqslant2}  & \cdot\; \sigma_{\psttab}(u_{1,0}) \\
\sigma_{-3\psttab}(u_{4,\ell})  & =\; \cfinal_{-3\psttab,\psttab,4,\ell\leqslant3}  & \cdot\; \sigma_{\psttab}(u_{0,0}) \\
[0.5ex] \end{array}$ \\
\hline
\end{tabular}
\end{table}

\begin{table}[h]
\caption{List of the coefficients $\sigma_q(\mathcal{S}^\infty_{d,n})$ that are needed to compute the $u_{p,\ell}$ with $p\leqslant4$ and $\Theta=\pi$}
\label{tab: liste CCC pour pi}
\begin{tabular}{|c|c|l|}
\hline
$d$ & $n$ & \multicolumn{1}{c|}{$q$} \\
\hline
\multirow{2}{*}{$\psttab$} & 0 & $-\psttab$, $-2\psttab$, $-3\psttab$ \\
\cline{2-3}
 & 1 & $-\psttab$ \\
\hline
$2\psttab$ & 0 & $-\psttab$, $-2\psttab$ \\
\hline
$3\psttab$ & 0 & $-\psttab$ \\
\hline
\end{tabular}
\end{table}

\renewcommand{\arraystretch}{1.2}

\begin{table}[h]
\caption{List of the 37 non-zero coefficients $\cmacro_{d,d',p,\ell}$ and $\cmicro_{d,d',p,\ell}$ needed to compute the $u_{p,\ell}$ with $p\leqslant4$, $\Theta=\pi$ and arbitrary values of $(\mu_0,\mu_1,\rho_0,\rho_1,\omega)$}
\label{tab: liste cmaicro pour pi}
\begin{tabular}[t]{|c|c|c|c|c|}
\hline
$d$ & $d'$ & $p$ & $\ell$ & $\cmacro_{d,d',p,\ell}$\rule[-1.4ex]{0pt}{4ex} \\
\hline
\multirow{7}{*}{$-3\psttab$} & $-3\psttab$ & $3$ & 0 & $1$ \\
\cline{2-5}
 & $-2\psttab$ & $3$ & 1 & $0.6852$ \\
\cline{2-5}
 & \multirow{2}{*}{$-\psttab$} & \multirow{2}{*}{$3$} & 1 & $0.1174$ \\
 &  &  & 2 & $0.1174$ \\
\cline{2-5}
 & \multirow{3}{*}{$\psttab$} & \multirow{3}{*}{$3$} & 1 & $-0.001376$ \\
 &  &  & 2 & $-0.02067$ \\
 &  &  & 3 & $-0.004593$ \\
\hline
\multirow{6}{*}{$-2\psttab$} & $-2\psttab$ & $2$ & 0 & $1$ \\
\cline{2-5}
 & $-\psttab$ & $2$ & 1 & $0.3426$ \\
\cline{2-5}
 & \multirow{2}{*}{$\psttab$} & \multirow{2}{*}{$2$} & 1 & $-0.04021$ \\
 &  &  & 2 & $-0.02011$ \\
\cline{2-5}
 & \multirow{2}{*}{$2\psttab$} & \multirow{2}{*}{$2$} & 1 & $-0.03996$ \\
 &  &  & 3 & $0.004593$ \\
\hline
\multirow{13}{*}{$-\psttab$} & $-2\psttab$ & $3$ & 1 & $2.013 + 1.911\,\mathrm{i}$ \\
\cline{2-5}
 & \multirow{3}{*}{$-\psttab$} & $1$ & 0 & $1$ \\
\cline{3-5}
 &  & \multirow{2}{*}{$3$} & 1 & $-1.379 - 1.31\,\mathrm{i}$ \\
 &  &  & 2 & $1.035 + 0.9823\,\mathrm{i}$ \\
\cline{2-5}
 & \multirow{4}{*}{$\psttab$} & $1$ & 1 & $-0.1174$ \\
\cline{3-5}
 &  & \multirow{3}{*}{$3$} & 1 & $0.1798 + 0.1708\,\mathrm{i}$ \\
 &  &  & 2 & $0.06071 + 0.05765\,\mathrm{i}$ \\
 &  &  & 3 & $-0.05397 - 0.05124\,\mathrm{i}$ \\
\cline{2-5}
 & \multirow{2}{*}{$2\psttab$} & \multirow{2}{*}{$1$} & 1 & $-0.08043$ \\
 &  &  & 2 & $0.04021$ \\
\cline{2-5}
 & \multirow{3}{*}{$3\psttab$} & \multirow{3}{*}{$1$} & 1 & $0.03308$ \\
 &  &  & 2 & $0.062$ \\
 &  &  & 3 & $-0.01378$ \\
\hline
\end{tabular}
\hspace{10mm}
\begin{tabular}[t]{|c|c|c|c|c|}
\hline
$d$ & $d'$ & $p$ & $\ell$ & $\cmicro_{d,d',p,\ell}$\rule[-1.4ex]{0pt}{4ex} \\
\hline
\multirow{7}{*}{$\psttab$} & $-\psttab$ & $1$ & 1 & $11.75 + 11.16\,\mathrm{i}$ \\
\cline{2-5}
 & \multirow{3}{*}{$\psttab$} & $1$ & 0 & $1$ \\
\cline{3-5}
 &  & \multirow{2}{*}{$3$} & 1 & $-0.3448 - 0.3274\,\mathrm{i}$ \\
 &  &  & 2 & $1.035 + 0.9823\,\mathrm{i}$ \\
\cline{2-5}
 & $2\psttab$ & $2$ & 1 & $0.6852$ \\
\cline{2-5}
 & \multirow{2}{*}{$3\psttab$} & \multirow{2}{*}{$3$} & 1 & $0.3521$ \\
 &  &  & 2 & $0.3521$ \\
\hline
\multirow{3}{*}{$2\psttab$} & $\psttab$ & $3$ & 1 & $1.007 + 0.9557\,\mathrm{i}$ \\
\cline{2-5}
 & $2\psttab$ & $2$ & 0 & $1$ \\
\cline{2-5}
 & $3\psttab$ & $3$ & 1 & $1.028$ \\
\hline
$3\psttab$ & $3\psttab$ & $3$ & 0 & $1$ \\
\hline
\end{tabular}
\end{table}
\renewcommand{\arraystretch}{1.25}

\begin{table}[h]
\caption{List of the 29 generically non-zero coefficients $\cfinal_{d,d',p,\ell}$ needed to compute the $u_{p,\ell}$ with $p\leqslant4$ and $\Theta=\pi$}
\label{tab: liste cfinal pour pi}
\begin{tabular}{|c|c|c|c|l|}
\hline
$d$ & $d'$ & $p$ & $\ell$ & \multicolumn{1}{c|}{expression of $\cfinal_{d,d',p,\ell}$} \\
\hline
$-3\psttab$ & $\psttab$ & $4$ & 0, 1, 2, 3 & $\begin{array}{@{}l@{\;}l@{\;}l@{\;}l@{\;}l@{\;}l@{}}
\rule{0pt}{2.8ex} \!\! =\; \cmacro_{-3\psttab,-\psttab,3,\ell\leqslant2}  & \cdot\; \sigma_{-\psttab}(\mathcal{S}^\infty_{\psttab,0})  & \cdot\; \cmicro_{\psttab,\psttab,1,0} & +\; \cmacro_{-3\psttab,-2\psttab,3,\ell=1}  & \cdot\; \sigma_{-2\psttab}(\mathcal{S}^\infty_{\psttab,0})  & \cdot\; \cmicro_{\psttab,\psttab,1,0} \\
+\; \cmacro_{-3\psttab,-3\psttab,3,0}  & \cdot\; \sigma_{-3\psttab}(\mathcal{S}^\infty_{\psttab,0})  & \cdot\; \cmicro_{\psttab,\psttab,1,0} & +\; \cmacro_{-3\psttab,\psttab,3,\ell\leqslant3}  & \cdot\; \cmicro_{\psttab,\psttab,1,0} &  \\
[0.8ex] \end{array}$ \\
\hline
\multirow{3}{*}[-9mm]{$-2\psttab$} & $-\psttab$ & $3$ & 1, 2, 3 & $\begin{array}{@{}l@{\;}l@{\;}l@{\;}l@{\;}l@{\;}l@{}}
\rule{0pt}{2.8ex} \!\! =\; \cmacro_{-2\psttab,-\psttab,2,\ell=1}  & *\; \sigma_{-\psttab}(\mathcal{S}^\infty_{\psttab,0})  & *\; \cmicro_{\psttab,-\psttab,1,\ell=1} & +\; \cmacro_{-2\psttab,-2\psttab,2,0}  & \cdot\; \sigma_{-2\psttab}(\mathcal{S}^\infty_{\psttab,0})  & \cdot\; \cmicro_{\psttab,-\psttab,1,\ell=1} \\
+\; \cmacro_{-2\psttab,\psttab,2,\ell\leqslant2}  & *\; \cmicro_{\psttab,-\psttab,1,\ell=1} &  &  &  &  \\
[0.8ex] \end{array}$ \\
\cline{2-5}
 & $\psttab$ & $3$ & 0, 1, 2 & $\begin{array}{@{}l@{\;}l@{\;}l@{\;}l@{\;}l@{\;}l@{}}
\rule{0pt}{2.8ex} \!\! =\; \cmacro_{-2\psttab,-\psttab,2,\ell=1}  & \cdot\; \sigma_{-\psttab}(\mathcal{S}^\infty_{\psttab,0})  & \cdot\; \cmicro_{\psttab,\psttab,1,0} & +\; \cmacro_{-2\psttab,-2\psttab,2,0}  & \cdot\; \sigma_{-2\psttab}(\mathcal{S}^\infty_{\psttab,0})  & \cdot\; \cmicro_{\psttab,\psttab,1,0} \\
+\; \cmacro_{-2\psttab,\psttab,2,\ell\leqslant2}  & \cdot\; \cmicro_{\psttab,\psttab,1,0} &  &  &  &  \\
[0.8ex] \end{array}$ \\
\cline{2-5}
 & $2\psttab$ & $4$ & 0, 1, 2, 3 & $\begin{array}{@{}l@{\;}l@{\;}l@{\;}l@{\;}l@{\;}l@{}}
\rule{0pt}{2.8ex} \!\! =\; \cmacro_{-2\psttab,-\psttab,2,\ell=1}  & *\; \sigma_{-\psttab}(\mathcal{S}^\infty_{\psttab,0})  & *\; \cmicro_{\psttab,2\psttab,2,\ell=1} & +\; \cmacro_{-2\psttab,-\psttab,2,\ell=1}  & \cdot\; \sigma_{-\psttab}(\mathcal{S}^\infty_{2\psttab,0})  & \cdot\; \cmicro_{2\psttab,2\psttab,2,0} \\
+\; \cmacro_{-2\psttab,-2\psttab,2,0}  & \cdot\; \sigma_{-2\psttab}(\mathcal{S}^\infty_{\psttab,0})  & \cdot\; \cmicro_{\psttab,2\psttab,2,\ell=1} & +\; \cmacro_{-2\psttab,-2\psttab,2,0}  & \cdot\; \sigma_{-2\psttab}(\mathcal{S}^\infty_{2\psttab,0})  & \cdot\; \cmicro_{2\psttab,2\psttab,2,0} \\
+\; \cmacro_{-2\psttab,\psttab,2,\ell\leqslant2}  & *\; \cmicro_{\psttab,2\psttab,2,\ell=1} &  & +\; \cmacro_{-2\psttab,2\psttab,2,\ell\leqslant3}  & \cdot\; \cmicro_{2\psttab,2\psttab,2,0} &  \\
[0.8ex] \end{array}$ \\
\hline
\multirow{5}{*}[-20mm]{$-\psttab$} & $-\psttab$ & $2$ & 1, 2 & $\begin{array}{@{}l@{\;}l@{\;}l@{\;}l@{\;}l@{\;}l@{}}
\rule{0pt}{2.8ex} \!\! =\; \cmacro_{-\psttab,-\psttab,1,0}  & \cdot\; \sigma_{-\psttab}(\mathcal{S}^\infty_{\psttab,0})  & \cdot\; \cmicro_{\psttab,-\psttab,1,\ell=1} & +\; \cmacro_{-\psttab,\psttab,1,\ell=1}  & *\; \cmicro_{\psttab,-\psttab,1,\ell=1} &  \\
[0.8ex] \end{array}$ \\
\cline{2-5}
 & \multirow{2}{*}[-9mm]{$\psttab$} & $2$ & 0, 1 & $\begin{array}{@{}l@{\;}l@{\;}l@{\;}l@{\;}l@{\;}l@{}}
\rule{0pt}{2.8ex} \!\! =\; \cmacro_{-\psttab,-\psttab,1,0}  & \cdot\; \sigma_{-\psttab}(\mathcal{S}^\infty_{\psttab,0})  & \cdot\; \cmicro_{\psttab,\psttab,1,0} & +\; \cmacro_{-\psttab,\psttab,1,\ell=1}  & \cdot\; \cmicro_{\psttab,\psttab,1,0} &  \\
[0.8ex] \end{array}$ \\
\cline{3-5}
 &  & $4$ & 0, 1, 2, 3 & $\begin{array}{@{}l@{\;}l@{\;}l@{\;}l@{\;}l@{\;}l@{}}
\rule{0pt}{2.8ex} \!\! =\; \cmacro_{-\psttab,-\psttab,1,0}  & \cdot\; \sigma_{-\psttab}(\mathcal{S}^\infty_{\psttab,0})  & \cdot\; \cmicro_{\psttab,\psttab,3,\ell\leqslant2} & +\; \cmacro_{-\psttab,-\psttab,1,0}  & \cdot\; \sigma_{-\psttab}(\mathcal{S}^\infty_{2\psttab,0})  & \cdot\; \cmicro_{2\psttab,\psttab,3,\ell=1} \\
+\; \cmacro_{-\psttab,-\psttab,1,0} \cdot\omega^2 & \cdot\; \sigma_{-\psttab}(\mathcal{S}^\infty_{\psttab,1})  & \cdot\; \cmicro_{\psttab,\psttab,1,0} & +\; \cmacro_{-\psttab,-\psttab,3,\ell\leqslant2}  & \cdot\; \sigma_{-\psttab}(\mathcal{S}^\infty_{\psttab,0})  & \cdot\; \cmicro_{\psttab,\psttab,1,0} \\
+\; \cmacro_{-\psttab,-2\psttab,3,\ell=1}  & \cdot\; \sigma_{-2\psttab}(\mathcal{S}^\infty_{\psttab,0})  & \cdot\; \cmicro_{\psttab,\psttab,1,0} & +\; \cmacro_{-\psttab,\psttab,1,\ell=1}  & *\; \cmicro_{\psttab,\psttab,3,\ell\leqslant2} &  \\
+\; \cmacro_{-\psttab,2\psttab,1,\ell\leqslant2}  & *\; \cmicro_{2\psttab,\psttab,3,\ell=1} &  & +\; \cmacro_{-\psttab,\psttab,3,\ell\leqslant3}  & \cdot\; \cmicro_{\psttab,\psttab,1,0} &  \\
[0.8ex] \end{array}$ \\
\cline{2-5}
 & $2\psttab$ & $3$ & 0, 1, 2 & $\begin{array}{@{}l@{\;}l@{\;}l@{\;}l@{\;}l@{\;}l@{}}
\rule{0pt}{2.8ex} \!\! =\; \cmacro_{-\psttab,-\psttab,1,0}  & \cdot\; \sigma_{-\psttab}(\mathcal{S}^\infty_{\psttab,0})  & \cdot\; \cmicro_{\psttab,2\psttab,2,\ell=1} & +\; \cmacro_{-\psttab,-\psttab,1,0}  & \cdot\; \sigma_{-\psttab}(\mathcal{S}^\infty_{2\psttab,0})  & \cdot\; \cmicro_{2\psttab,2\psttab,2,0} \\
+\; \cmacro_{-\psttab,\psttab,1,\ell=1}  & *\; \cmicro_{\psttab,2\psttab,2,\ell=1} &  & +\; \cmacro_{-\psttab,2\psttab,1,\ell\leqslant2}  & \cdot\; \cmicro_{2\psttab,2\psttab,2,0} &  \\
[0.8ex] \end{array}$ \\
\cline{2-5}
 & $3\psttab$ & $4$ & 0, 1, 2, 3 & $\begin{array}{@{}l@{\;}l@{\;}l@{\;}l@{\;}l@{\;}l@{}}
\rule{0pt}{2.8ex} \!\! =\; \cmacro_{-\psttab,-\psttab,1,0}  & \cdot\; \sigma_{-\psttab}(\mathcal{S}^\infty_{\psttab,0})  & \cdot\; \cmicro_{\psttab,3\psttab,3,\ell\leqslant2} & +\; \cmacro_{-\psttab,-\psttab,1,0}  & \cdot\; \sigma_{-\psttab}(\mathcal{S}^\infty_{2\psttab,0})  & \cdot\; \cmicro_{2\psttab,3\psttab,3,\ell=1} \\
+\; \cmacro_{-\psttab,-\psttab,1,0}  & \cdot\; \sigma_{-\psttab}(\mathcal{S}^\infty_{3\psttab,0})  & \cdot\; \cmicro_{3\psttab,3\psttab,3,0} & +\; \cmacro_{-\psttab,\psttab,1,\ell=1}  & *\; \cmicro_{\psttab,3\psttab,3,\ell\leqslant2} &  \\
+\; \cmacro_{-\psttab,2\psttab,1,\ell\leqslant2}  & *\; \cmicro_{2\psttab,3\psttab,3,\ell=1} &  & +\; \cmacro_{-\psttab,3\psttab,1,\ell\leqslant3}  & \cdot\; \cmicro_{3\psttab,3\psttab,3,0} &  \\
[0.8ex] \end{array}$ \\
\hline
\end{tabular}
\end{table}

	\tikzset{external/export=false}%
	\tikzsetfigurename{main-figure}

\FloatBarrier 
\subsection{Synthetic representation of properties of the expansion for various orders and angles}
\label{sec: graphes exemples DA}

\begin{textesansboite}
\Crefrange{fig: graphe nb de CL}{fig: graphe nb de CCC} plot three quantities linked to computational costs: the numbers of far fields, of corner profiles and of corner profiles coefficients $\sigma _{q} (\ccDA_{d,n} )$ that need to be computed to build an asymptotic expansion of a given order. The abscissa of each point of the graphs is a $p\in \NN$ s.t.\ there~is at least one generically non-zero $u_{p,\ell }$, while the ordinate is the number of objects that need to be computed to build the $u_{p',\ell }$ with $p'<p$. Since $u_{\peps} -\sum _{p'<p} \sum _{\ell } \varepsilon ^{p'} \ln ^{\ell } \!\varepsilon \, u_{p',\ell } =O(\varepsilon ^{p} \ln ^{\ell _{p}} \varepsilon )$ for some $\ell _{p} \in \mathbb{N}$ (see \cref{DA de SE avec CL et delta variable}), $p$ is the convergence order of the obtained expansion.\\

Unsurprisingly, all curves in \crefrange{fig: graphe nb de CL}{fig: graphe nb de CCC} are increasing. In addition, curves associated to large angles tend to be above curves with a smaller angle: often, the greater the angle, the greater the computational cost. There are two notable exceptions in the represented examples.
\begin{enumerate}
    \item $\Theta =4$ requires more far fields than $\Theta = \frac{3\pi }{2} \approx  4.7$. Indeed, since $\frac{3\pi }{2}$ belongs to $\pi \mathbb{Q}$ but not 4, some far field orders $p=n_{1} +n_{2} \pst$ with distinct couples of integers $(n_{1} ,n_{2} )$ end up being equal for $\Theta = \frac{3\pi }{2}$ while there are different for $\Theta =4$. These coincidences correspond to intersection points of the curves in \cref{fig: ordres des CL}. Moreover, the expansion contains no $\ln \varepsilon $ for $\Theta  \in \{4 , \frac{3\pi }{2} \}$ at the values of $p$ considered in \cref{fig: graphe nb de CL}. So, for both angles, there~is only one far field to be computed for~any far field order $p$.
    \item Similarly, the curves associated to $\Theta =3$ are sometimes above the curves of $\Theta =\pi $. However, for the number of far fields for example, these curves meet regularly. Indeed the inferior number of far field orders for $\Theta =\pi $ is compensated by the presence of some $\ln \varepsilon $. Several far fields usually need to be computed for $\Theta =\pi $ at a given far field order.
\end{enumerate}

Note that the number of far fields to be computed can be very high. For example, for $\Theta = \frac{3\pi }{2}$ (a case that our industrial contacts where interested in), 12 far fields are needed to obtain an approximation of $u_{\peps}$ of order~5. When $\Theta  \rightarrow  2\pi ^{-}$, this number grows as high as 31 (although these are extreme angles). This drawback is a strong motivation for effective models, as they replace the computation of all these far fields, by the computational of one approximate solution. However, note that the number of corner profiles to be computed cannot be shrunk, even with effective models. But these functions are independent of the frequency $\omega $ and the source $\fsource$, so they can often be computed once in advance.
\end{textesansboite}

	\tikzset{external/export=true}%
	\tikzsetfigurename{fig-nb_CL-}%
	\pgfplotsset{
    every axis/.append style={
        thick,
    },
    cycle list/Dark2,
    cycle multiindex* list={
        mark list*\nextlist
        Dark2\nextlist
    },
    /pgfplots/table/x=ordrecvg,
    /pgfplots/table/y=nbCL,
    legend cell align=left,
}

\begin{figure}
\begin{tikzpicture}
\begin{axis}[
    width=10cm,
    xlabel={convergence order},
    ylabel={number of far fields},
    enlargelimits=0,
    restrict x to domain=0:5,
    xtick={1,2,...,5},
    ymin=0,
    xmajorgrids=true,
    ymajorgrids=true,
    legend pos=north west,
    legend style={legend columns=1},
]

\addplot
    table {../figures/stats_DA/stats_DA_pi2.dat};
    \addlegendentry{$\Theta=\frac\pi2 \approx 1.57$}
\addplot
    table {../figures/stats_DA/stats_DA_2pi3.dat};
    \addlegendentry{$\Theta=\frac{2\pi}3 \approx 2.09$}
\addplot
    table {../figures/stats_DA/stats_DA_3.dat};
    \addlegendentry{$\Theta=3$}
\addplot
    table {../figures/stats_DA/stats_DA_pi.dat};
    \addlegendentry{$\Theta=\pi \approx 3.14$}
\addplot
    table {../figures/stats_DA/stats_DA_4.dat};
    \addlegendentry{$\Theta=4$}
\addplot
    table {../figures/stats_DA/stats_DA_3pi2.dat};
    \addlegendentry{$\Theta=\frac{3\pi}2 \approx 4.71 $}
\addplot
    table {../figures/stats_DA/stats_DA_5.dat};
    \addlegendentry{$\Theta=5$}
    
\end{axis}
\end{tikzpicture}
\caption{
    Generic number of far fields needed to be computed for a desired rate of convergence.
}
\label{fig: graphe nb de CL}
\end{figure}
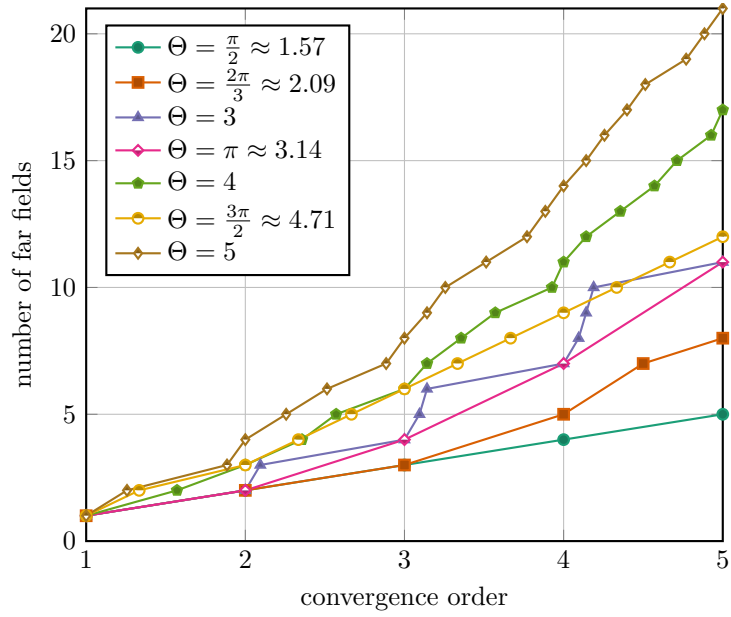%
	\tikzset{external/export=false}%
	\tikzsetfigurename{main-figure}
	\tikzset{external/export=true}%
	\tikzsetfigurename{fig-nb_CC-}%
	\pgfplotsset{
    every axis/.append style={
        thick,
    },
    cycle list/Dark2,
    cycle multiindex* list={
        mark list*\nextlist
        Dark2\nextlist
    },
    /pgfplots/table/x=ordrecvg,
    /pgfplots/table/y=nbCC,
    legend cell align=left,
}

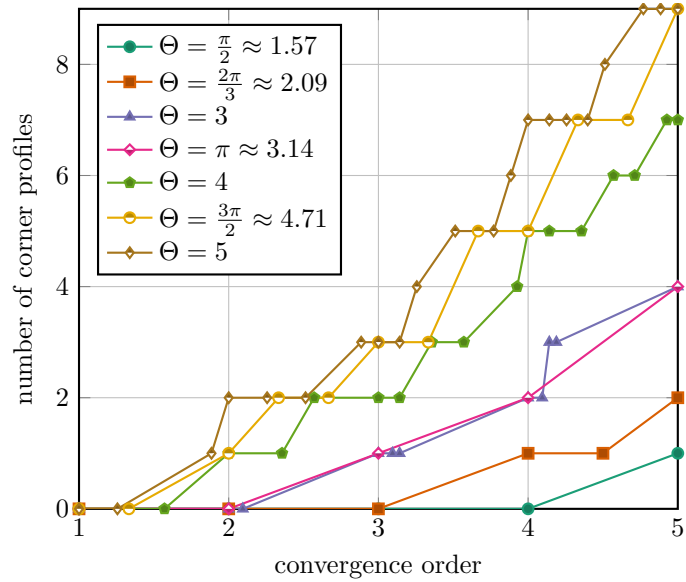
\begin{figure}
\begin{tikzpicture}
\begin{axis}[
    width=95mm,
    xlabel={convergence order},
    ylabel={number of corner profiles},
    enlargelimits=0,
    restrict x to domain=0:5,
    xtick={1,2,...,5},
    ymin=0,
    xmajorgrids=true,
    ymajorgrids=true,
    legend pos=north west,
    legend style={legend columns=1},
]

\addplot
    table {../figures/stats_DA/stats_DA_pi2.dat};
    \addlegendentry{$\Theta=\frac\pi2 \approx 1.57$}
\addplot
    table {../figures/stats_DA/stats_DA_2pi3.dat};
    \addlegendentry{$\Theta=\frac{2\pi}3 \approx 2.09$}
\addplot
    table {../figures/stats_DA/stats_DA_3.dat};
    \addlegendentry{$\Theta=3$}
\addplot
    table {../figures/stats_DA/stats_DA_pi.dat};
    \addlegendentry{$\Theta=\pi \approx 3.14$}
\addplot
    table {../figures/stats_DA/stats_DA_4.dat};
    \addlegendentry{$\Theta=4$}
\addplot
    table {../figures/stats_DA/stats_DA_3pi2.dat};
    \addlegendentry{$\Theta=\frac{3\pi}2 \approx 4.71 $}
\addplot
    table {../figures/stats_DA/stats_DA_5.dat};
    \addlegendentry{$\Theta=5$}
    
\end{axis}
\end{tikzpicture}
\caption{
    Generic number of corner profiles needed to be computed for a desired rate of convergence of the far fields
}
\label{fig: graphe nb de CC}
\end{figure}%
	\tikzset{external/export=false}%
	\tikzsetfigurename{main-figure}
	\tikzset{external/export=true}%
	\tikzsetfigurename{fig-nb_CCC-}%
	\pgfplotsset{
    every axis/.append style={
        thick,
    },
    cycle list/Dark2,
    cycle multiindex* list={
        mark list*\nextlist
        Dark2\nextlist
    },
    /pgfplots/table/x=ordrecvg,
    /pgfplots/table/y=nbCCC,
    legend cell align=left,
}

\begin{figure}
\begin{tikzpicture}
\begin{axis}[
    width=95mm,
    xlabel={convergence order},
    ylabel={number of corner profiles coefficients},
    enlargelimits=0,
    restrict x to domain=0:5,
    xtick={1,2,...,5},
    ymin=0,
    xmajorgrids=true,
    ymajorgrids=true,
    legend pos=north west,
    legend style={legend columns=1},
]

\addplot
    table {../figures/stats_DA/stats_DA_pi2.dat};
    \addlegendentry{$\Theta=\frac\pi2 \approx 1.57$}
\addplot
    table {../figures/stats_DA/stats_DA_2pi3.dat};
    \addlegendentry{$\Theta=\frac{2\pi}3 \approx 2.09$}
\addplot
    table {../figures/stats_DA/stats_DA_3.dat};
    \addlegendentry{$\Theta=3$}
\addplot
    table {../figures/stats_DA/stats_DA_pi.dat};
    \addlegendentry{$\Theta=\pi \approx 3.14$}
\addplot
    table {../figures/stats_DA/stats_DA_4.dat};
    \addlegendentry{$\Theta=4$}
\addplot
    table {../figures/stats_DA/stats_DA_3pi2.dat};
    \addlegendentry{$\Theta=\frac{3\pi}2 \approx 4.71 $}
\addplot
    table {../figures/stats_DA/stats_DA_5.dat};
    \addlegendentry{$\Theta=5$}
    
\end{axis}
\end{tikzpicture}
\caption{
    Generic number of coefficients corner profiles $\sigma_q(\mathcal{S}_{d,n})$ with $q<0$ needed for a desired rate of convergence of the far fields
}
\label{fig: graphe nb de CCC}
\end{figure}
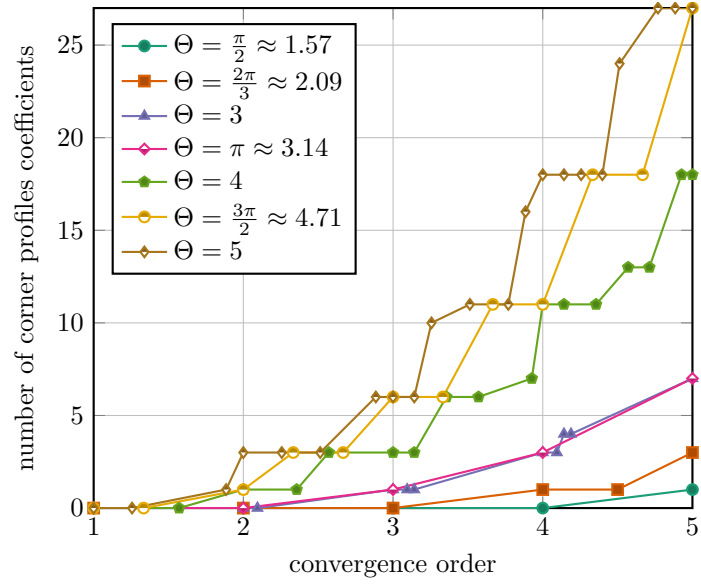%
	\tikzset{external/export=false}%
	\tikzsetfigurename{main-figure}

\begin{textesansboite}
\Cref{tab: ordre du premier ln eps} shows the first order where a $\ln \varepsilon $ appears for angles $\Theta $ in $\pi \mathbb{Q}$ with a denominator smaller or equal to 6. For larger denominators, \cref{simplifications cas rationnel} states that $\ln \varepsilon $ do~not appear before order $p=7$ at least, which is enough for most applications. Comparing \cref{tab: ordre du premier ln eps} to \cref{simplifications cas rationnel}, it seems that the bound of the theorem is optimal when the denominator is odd, but not when it~is even. \Cref{simplifications cas rationnel version forte} justifies this for $b\in \{2,4\}$ and we conjecture that similar techniques could explain the non-optimality of \cref{simplifications cas rationnel} for the other even denominators.\\

\Cref{fig: graphe deg ln eps} shows the evolution of the maximal power of $\ln \varepsilon $ w.r.t.\ the convergence order. Note that in this figure $p\in [0,10]$, whereas $p\in [0,5]$ for \crefrange{fig: graphe nb de CL}{fig: graphe nb de CCC}. One can see that each curve starts with a plateau at 0, as there~is no $\ln \varepsilon $ in the first terms of the expansion.
\end{textesansboite}

	\tikzset{external/export=true}%
	\tikzsetfigurename{tab-premier_ln_eps-}%
	
\renewcommand{\arraystretch}{1.3}

\begin{table}[h]
\caption{
    Values of $\min\{p\in\NN \mid \exists\ell>0, u_{p,\ell}\neq0 \text{ generically}\}$ for $\Theta = \frac{a}{b}\pi$ with $a,b\in\mathbb{N}^*$ and $\pgcd{a,b}=1$
}
\label{tab: ordre du premier ln eps}
\def\rien{$-$}
\begin{tabular}{|c||*{10}{c|}}
\hline
\backslashbox{$b$}{$a$} & 1 & 2 & 3 & 4 & 5 & 6 & 7 & 8 & 9 & 11 \\ \hline\hline
1 & 2 & \rien & \rien & \rien & \rien & \rien & \rien & \rien & \rien & \rien \\ \hline
2 & 8 & \rien & $\frac{16}{3} \approx 5.33$ & \rien & \rien & \rien & \rien & \rien & \rien & \rien \\ \hline
3 & 6 & 3 & \rien & 3 & 3 & \rien & \rien & \rien & \rien & \rien \\ \hline
4 & 12 & \rien & 10 & \rien & 9.2 & \rien & $\frac{58}{7} \approx 8.29$ & \rien & \rien & \rien \\ \hline
5 & 10 & 5 & 5 & 5 & \rien & 5 & 5 & 5 & 5 & \rien \\ \hline
6 & 14 & \rien & \rien & \rien & 8.4 & \rien & $\frac{54}{7} \approx 7.71$ & \rien & \rien & $\frac{78}{11} \approx 7.09$ \\ \hline
\end{tabular}
\end{table}
	\tikzset{external/export=false}%
	\tikzsetfigurename{main-figure}
	\tikzset{external/export=true}%
	\tikzsetfigurename{fig-deglogmax-}%
	\pgfplotsset{
    every axis/.append style={
        thick,
    },
    cycle list/Dark2,
    cycle multiindex* list={
        mark list*\nextlist
        Dark2\nextlist
    },
    /pgfplots/table/x=ordrecvg,
    /pgfplots/table/y=deglogmax,
    legend cell align=left,
}

\begin{figure}
\begin{tikzpicture}
\begin{axis}[
    width=80mm,
    xlabel={convergence order},
    ylabel={maximal power of $\ln\varepsilon$},
    enlargelimits=0,
    restrict x to domain=0:10,
    ymin=0,
    xmajorgrids=true,
    ymajorgrids=true,
    legend pos=north west,
    legend style={legend columns=1},
]

\addplot
    table {../figures/stats_DA/stats_DA_pi2.dat};
    \addlegendentry{$\Theta=\frac\pi2 \approx 1.57$}
\addplot
    table {../figures/stats_DA/stats_DA_2pi3.dat};
    \addlegendentry{$\Theta=\frac{2\pi}3 \approx 2.09$}
\addplot
    table {../figures/stats_DA/stats_DA_pi.dat};
    \addlegendentry{$\Theta=\pi \approx 3.14$}
\addplot
    table {../figures/stats_DA/stats_DA_3pi2.dat};
    \addlegendentry{$\Theta=\frac{3\pi}2 \approx 4.71 $}
    
\end{axis}
\end{tikzpicture}
\longCaption{Values of $\max\{\ell\in\mathbb{N} \mid \exists p'<p,\, u_{p,\ell}\neq0 \text{ generically}\}$ in terms of $p$}{ (the convergence order of $\sum_{p'<p,\ell} u_{p',\ell}$ towards $u_\varepsilon$)}
\label{fig: graphe deg ln eps}
\end{figure}%
	\tikzset{external/export=false}%
	\tikzsetfigurename{main-figure}

\begin{subappendices}
\section{Appendix: extension of the results of \cref{sec: expression des mcl,sec: simplifications DA cas rationnel} to some particular cases}
\label{sec: annexe des complements sur le DA}

\let\varepsilon=\feps

\nomenclature[34-10]{$Y_{a}$}{Bessel function of the second kind}
\begin{textesansboite}
In this section we assume $\pst \in \mathbb{Q}$. We extend the results of \cref{sec: expression des mcl,sec: simplifications DA cas rationnel} to some particular cases that were not covered due to the assumptions of \cref{expression sing heterog CL}. In particular, we study the series $\mcla_{d,n}$ that can be expressed with Bessel function of the \emph{second} kind, and improve \cref{simplifications cas rationnel} for some angles at some orders that are relevant for the effective models of \cref{chap: ME}. We do~not give all the proofs' details as the reasoning used in this section are very similar to \cref{sec: expression des mcl,sec: simplifications DA}.
\end{textesansboite}

\subsection{Extension of \texorpdfstring{\cref{sec: expression des mcl}}{Section 6.1}}

\begin{textesansboite}
For~any $a\in \mathbb{R} \setminus \mathbb{Z}$, the Bessel function of the second kind $Y_{a}$ is defined by
\begin{equation}{\label{eq: def BesselY_a}}
    \forall z\in \mathbb{C}^*, \qquad  Y_{a} (z) := \frac{J_{a} (z) \cos (\pi a) - J_{-a} (z) }{\sin (\pi a).}
\end{equation}
It~is extended to the case of an integer parameter $n\in \mathbb{Z}$ by $Y_{n} (z):=\lim _{a\rightarrow n} Y_{a} (z)$. It satisfies for~any $n\in \mathbb{N}$:
\begin{equation}{\label{eq: expression de BesselY_n}}
    Y_{n} (z) = - \sum _{k=0} ^{n-1} \frac{(n-k-1) !}{\pi \.k!} \left( \frac{z}{2} \right)^{-n+2k} + \frac{2}{\pi } J_{n} (z) \ln \!\left( \frac{z}{2} \right) - \sum_{k=0}^\infty  (-1)^{k} \frac{\psi (k+1)+\psi (n+k+1)}{\pi \. k!\, (n+k) !} \left( \frac{z}{2} \right)^{n+2k} 
\end{equation}
where $\psi (z):= \frac{\Gamma '(z)}{\Gamma (z)}$ denotes the Digamma function, and
\begin{equation}{\label{eq: BesselY pour n et -n}}
    Y_{-n} (z)=(-1)^{n} \.Y_{n} (z).
\end{equation}
See \cite{Leb72}. We have seen in \cref{sec: expression des mcl} that $J_{a}$ behaves like $r^{a}$ for~any $a$, except when $a\in -\mathbb{N}^*$. The Bessel functions of the second kind fill this gap as $Y_{a}$ behaves like $r^{a}$ when $a\in -\mathbb{N}^*$ according to \cref{eq: expression de BesselY_n,eq: BesselY pour n et -n}.\\

Similarly to $\BesselOmega_{a}$ and $\BesselGamma_{a}$, we introduce formal series $\BesselYGamma_{a} \in \Amacro(\Gamma )$ and $\BesselYOmega_{a} \in \Amacro(\Omega )$ that converge to the Bessel functions of the second kind (which is verified in \cref{BesselYGamma/Omega cvg vers Y} below):
\begin{itemize}
    \item For~any $a\in \mathbb{R} \setminus \mathbb{Z}$, we define
    \begin{equation}{\label{eq: def BesselYGamma/Omega_a}}
        \BesselYGamma_{a} := \frac{\BesselGamma_{a} \cos (\pi a) - \BesselGamma_{-a} }{\sin (\pi a)} \qquad  \text{ and } \qquad  \BesselYOmega_{a} := \frac{\BesselOmega_{a} \cos (\pi a) + \BesselOmega_{-a} }{\sin (\pi a).}
    \end{equation}
    \item For~any $n\in \mathbb{N}$, we define
    \begin{align}
        \BesselYGamma_{n} &:= -\fsum_{k=0} ^{n-1} \frac{(n-k-1) !}{\pi \.k!} \left( \frac{k_{0} x}{2} \right)^{-n+2k} +\, \fsum_{k=0} ^{\infty } \frac{2\.(-1)^{k} }{\pi \.k!\, \Gamma (n+k+1)} \left( \frac{k_{0} x}{2} \right)^{n+2k} \ln \!\left( \frac{k_{0} x}{2} \right) \nonumber\\
        &\quad {}- \fsum_{k=0} ^{\infty } (-1)^{k} \frac{\psi (k+1)+\psi (n+k+1)}{\pi \. k!\, (n+k) !} \left( \frac{k_{0} x}{2} \right)^{n+2k} \label{eq: def BesselYGamma_n}
    \end{align}
    and
    \begin{align}
        \BesselYOmega_{n} &:= \fsum_{k=0} ^{n-1} \frac{(n-k-1) !}{\pi \.k!} \left( \frac{k_{0} }{2} \right)^{-n+2k} \Im[(\alpha z)^{-n+k} \,\overline{\alpha z}^{k} ] \nonumber\\
        &\quad {}+ \fsum_{k=0} ^{\infty } \frac{2\.(-1)^{k} }{\pi \.k!\, \Gamma (n+k+1)} \left( \frac{k_{0} }{2} \right)^{n+2k} \bigg( \begin{array}[t]{@{}l}
           \displaystyle \frac{1}{2} \Im[(\alpha z)^{n+k} \,\overline{\alpha z}^{k} \log (\alpha z)] - \frac{1}{2} \Im[(\alpha z)^{k} \,\overline{\alpha z}^{n+k} \log (\alpha z)] \\[1.5ex]
           \displaystyle + \ln \!\Big( \frac{k_{0} }{2} \Big)\, \Im[(\alpha z)^{n+k} \,\overline{\alpha z}^{k} ] \bigg)
        \end{array} \nonumber\\
        &\displaystyle - \fsum_{k=0} ^{\infty } (-1)^{k} \frac{\psi (k+1)+\psi (n+k+1)}{\pi \. k!\, (n+k) !} \left( \frac{k_{0} }{2} \right)^{n+2k} \Im[(\alpha z)^{n+k} \,\overline{\alpha z}^{k} ] \label{eq: def BesselYOmega_n}
    \end{align}
    and $\BesselYGamma_{-n} := (-1)^{n} \BesselYGamma_{n}$ and $\BesselYOmega_{-n} := -(-1)^{n} \BesselYOmega_{n}$.
\end{itemize}
\end{textesansboite}

\begin{lemmesansboite}{}
\label{BesselYGamma/Omega cvg vers Y}%
For~any $a\in \mathbb{R}$, we have $\BesselYGamma_{a} \longrightarrow  Y_{a} (k_{0} x)$ \;and\; $\BesselYOmega_{a} \longrightarrow  Y_{a} (k_{0} r) \sin (a(\theta -\Theta ))$ (see \cref{def: convergence de SF dans Amacro}).
\end{lemmesansboite}

\begin{demo}
By~\crefrange{eq: def BesselY_a}{eq: BesselY pour n et -n} and the definition of $J_{a} ,\BesselGamma_{a}$, $\BesselOmega_{a}$ in \crefrange{eq: def BesselJ}{eq: def BesselOmega}, the only non obvious case is $\BesselYOmega_{n}$ when $n\in \mathbb{N}$. We need the following calculation: for~any $z\in \mathbb{C}$ with polar decomposition $re^{\mathrm{i}\theta }$, we have
\begin{multline*}
    \Im[(\alpha z)^{n+k} \,\overline{\alpha z}^{k} \log (\alpha z)] - \Im[(\alpha z)^{k} \,\overline{\alpha z}^{n+k} \log (\alpha z)] \\
    \begin{array}{l}
    = \Im[r^{n+2k} e^{\mathrm{i}n(\theta -\Theta )} (\ln r+\mathrm{i}(\theta -\Theta ))]- \Im[r^{n+2k} e^{-\mathrm{i}n(\theta -\Theta )} (\ln r+\mathrm{i}(\theta -\Theta ))] \\[\jot]
    = \Im[r^{n+2k} e^{\mathrm{i}n(\theta -\Theta )} (\ln r+\cancel{\mathrm{i}(\theta -\Theta )})]+ \Im[r^{n+2k} e^{\mathrm{i}n(\theta -\Theta )} (\ln r-\cancel{\mathrm{i}(\theta -\Theta )})] \\[\jot]
    = 2\, r^{n+2k} \sin (n(\theta -\Theta )) \ln r 
    \end{array}
\end{multline*}
This shows the second sum in \cref{eq: def BesselYOmega_n} converges to
$\frac{2}{\pi } J_{n} (k_{0} r) \ln (k_{0} r/2) \sin (n(\theta -\Theta ))$. Therefore we have $\BesselYOmega_{n} \longrightarrow  Y_{n} (k_{0} r) \sin (n(\theta -\Theta ))$.
\end{demo}

\begin{lemmesansboite}{}
\label{relations sur BesselYGamma/Omega}%
The relations of \cref{relations sur BesselGamma/Omega} also hold when replacing $\BesselGamma_{a}$ by $\BesselYGamma_{a}$ and $\BesselOmega_{a}$ by $\BesselYOmega_{a}$. Namely, we have for~any $a\in \mathbb{R}$:
\begin{enumerate}
    \item $(\Delta +k_{0} ^2 )\BesselYOmega_{a} = 0$.
    \item $\BesselYOmega_{a|\Gamma } = -\sin (a\Theta ) \BesselYGamma_{a}$.
    \item $\partial _{y} \BesselYOmega_{a|\Gamma } = \frac{k_{0} }{2} \cos (a\Theta ) (\BesselYGamma_{a-1} + \BesselYGamma_{a+1} )$.
    \item $\partial _{x} \BesselYGamma_{a} = \frac{k_{0} }{2} (\BesselYGamma_{a-1} -\BesselYGamma_{a+1} )$.
\end{enumerate}\fbi
\end{lemmesansboite}

\begin{demo}
The properties of $J_{a}$ that are used in the proof of \cref{relations sur BesselGamma/Omega}, e.g.\ $2 J_{a} '(z) = J_{a-1} (z)-J_{a+1} (z)$, also hold for $Y_{a}$ by \cite{Leb72}. Hence points 1 to 4 hold for $x\mapsto Y_{a} (k_{0} x)$ and $(r,\theta ) \mapsto  Y_{a} (k_{0} r) \sin (a(\theta -\Theta ))$. It remains to show that this implies the same relations on the formal series $\BesselYGamma_{a}$ and $\BesselYOmega_{a}$, i.e.\ term by term.\smallskip\par
For~any $a\in \mathbb{R}$, there exist $(d_{n} )\in \mathbb{R}^{\mathbb{N}}$ and $(c_{n} ),(c'_{n} ) \in \mathbb{C}^{\mathbb{N}}$ s.t.\ $d_{n} \rightarrow \infty $ and
\[\forall x>0, \qquad  Y_{a} (k_{0} x) = \sum_{n=0}^\infty  x^{d_{n}} (c_{n} + c'_{n} \ln x).\]
In addition we have $c_{n} =O(\rho ^{d_{n}} )$ and $c'_{n} =O(\rho ^{d_{n}} )$ for~any $\rho >0$ (using $\psi (k)=O(\ln k)$ when $k\rightarrow \infty $). One can verify that this implies
\[\forall j\in \mathbb{N},\forall q\in \mathbb{R}, \qquad  \partial _{x} ^{j} \big[Y_{a} (k_{0} x)\big] = \sum _{n\in \mathbb{N},\. d_{n} \leqslant q} \partial _{x} ^{j} \big(x^{d_{n}} (c_{n} + c'_{n} \ln x)\big) + o(x^{q-j} ).\]
Hence we deduce:
\[\forall q\in \mathbb{R}, \qquad  Y_{a} (k_{0} x) = T_{\leqslant q} (\BesselYGamma_{a} ) + o_{\partial } (x^{q} ) \quad  \text{ and } \quad  Y_{a} (k_{0} r) \sin (a(\theta -\Theta )) = T_{\leqslant q} (\BesselYOmega_{a} ) + o_{\partial } (r^{q} ).\]
Using \cref{derivation des DA dans Amaicro}, we obtain
\[\forall q\in \mathbb{R}, \qquad  0 = (\Delta +k_{0} ^2 )\big[Y_{a} (k_{0} r) \sin (a(\theta -\Theta )) \big] = T_{\leqslant q} \circ (\Delta +k_{0} ^2 )(\BesselYOmega_{a} ) + o_{\partial } (r^{q} ).\]
By~uniqueness of the asymptotic expansion in $\Amacro$ (see \cref{unicité du DA}), this yields $(\Delta +k_{0} ^2 )(\BesselYOmega_{a} ) = 0$. Similarly for the other relations.
\end{demo}

\begin{lemmesansboite}{}
\label{sigma de BesselYOmega}%
Let $n\in \mathbb{N}$. We have $\sigma _{d} (\BesselYOmega_{n} )=0$ for~any $d\in \TZ* \setminus \{-n,n\}$ and, if $n\in \TZ*$, then
\[\sigma _{-n} (\BesselYOmega_{n} ) = \frac{(n-1)!}{\pi } \left( \frac{k_{0} }{2} \right)^{-n} \qquad  \text{ and } \qquad  \sigma _{n} (\BesselYOmega_{n} ) = \frac{2\ln (k_{0} /2)-H_{n} +2\gamma }{\pi \.n!} \left( \frac{k_{0} }{2} \right)^{n} \]
where $H_{n} := \sum _{k=1} ^{n} \frac{1}{k}$ denotes the $n$-th harmonic number and $\gamma \approx 0.577$ is the Euler-Mascheroni constant.
\end{lemmesansboite}

\begin{remarque}
When $a\in \mathbb{R} \setminus \mathbb{Z}$, the values of $\sigma _{d} (\BesselYOmega_{a} )$ are given by \cref{sigma de BesselOmega,eq: def BesselYGamma/Omega_a}.
\end{remarque}

\begin{demo}
Similarly to \cref{sigma de BesselOmega}, one can verify that each function of the form $\Im[(\alpha z)^{q} \overline{\alpha z}^{k} \log (\alpha z)^{j} ]$ that appears in \cref{eq: def BesselYOmega_n} either vanishes or belongs to the basis of $\mathcal{A}(\Omega )$ used to define $\sigma _{d}$ in \cref{def: sigma}. Thus the values of $\sigma _{d} (\BesselYOmega_{n} )$ (for $d\in \TZ* \setminus \{-n,n\}$) and $\sigma _{-n} (\BesselYOmega_{n} )$ (when $n\in \TZ*$) follow easily, and we have (when $n\in \TZ*$)
\[\sigma _{n} (\BesselYOmega_{n} ) = \frac{2}{\pi \.\Gamma (n+1)} \left( \frac{k_{0} }{2} \right)^{n} \ln \!\left( \frac{k_{0} }{2} \right) - \frac{\psi (1)+\psi (n+1)}{\pi \.n!} \left( \frac{k_{0} }{2} \right)^{n} = \frac{2\ln (k_{0} /2)-H_{n} +2\gamma }{\pi \.n!} \left( \frac{k_{0} }{2} \right)^{n} \]
because $\Gamma (n+1)=n!$ and $\psi (m)=H_{m-1} -\gamma $ for~any $m\in \mathbb{N}^*$ (see \cite{Leb72}).
\end{demo}

\begin{textesansboite}
Using \cref{BesselYGamma/Omega cvg vers Y,relations sur BesselYGamma/Omega,sigma de BesselYOmega} and the same reasoning as \cref{expression sing homog CL}, one can prove the following result. Together with \cref{expression sing homog CL}, this gives an explicit expression of $\mcla_{d,0}$ for~any $d\in \TZ*$.
\end{textesansboite}

\begin{proposition}{}
\label{expression sing homog CL avec BesselY}%
For~any $n\in \mathbb{N} \cap  \TZ*$, we have
\[\mcla_{-n,0} = \langle -k_{0} ^2  \ROmegaOmega\rangle (\phiOmega_{-n} ) = \frac{1}{(n-1) !} \left( \frac{k_{0} }{2} \right)^{n} \left( \pi  \BesselYOmega_{n} +\Big(2\ln \!\Big( \frac{2}{k_{0}} \Big)+H_{n} -2\gamma \Big) \BesselOmega_{n} \right) \]
and
\[\mcla_{-n,0} \longrightarrow  \frac{(-1)^{n\Theta /\pi } }{(n-1) !} \left( \frac{k_{0} }{2} \right)^{n} \left( \pi  Y_{n} (k_{0} r) +\Big(2\ln \!\Big( \frac{2}{k_{0}} \Big)+H_{n} -2\gamma \Big) J_{n} (k_{0} r) \right) \sin (n\theta ).\]\fbeq[0.8]
\end{proposition}

\begin{textesansboite}
In the rest of this section, we consider $a,b\in \mathbb{N}^*$ s.t.\ $\Theta  = \frac{a}{b} \pi $ and $\pgcd{a,b}=1$. The next result generalizes \cref{expression sing heterog CL} to the case $d\in -\mathbb{N}$.
\end{textesansboite}

\begin{proposition}{}
\label{expression sing heterog CL avec BesselY}%
For~any $d\in \TZ* \cap (-\mathbb{N})$ and $n\in [\![0,b-1]\!]$ we have
\[\mcla_{d,n} = \frac{1}{(|d|-1) !} \left( \frac{2}{k_{0}} \right)^{d} \sum _{m=-n} ^{n} a_{n,m} \left( \pi  \BesselYOmega_{|d|+m} +\Big(2\ln \!\Big( \frac{2}{k_{0}} \Big)+H_{|d|} -2\gamma \Big) \BesselOmega_{|d|+m} \right) .\]
where the numbers $a_{n,m}$ were introduced in \cref{expression sing heterog CL}.
\end{proposition}

\begin{demo}
\def\ZG{\mathcal{Z}}%
\def\ZO{\mathcal{Z}^\Omega}%
\def\Zrec{\widetilde{\mathcal{Z}}}%
The reasoning in the proof of \cref{expression sing homog CL} shows the following result.
\begin{quote}
    ``Let $(\ZG_{a} )\in \Amacro(\Gamma )^{\mathbb{R}}$ and $(\ZO_{a} )\in \Amacro(\Omega )^{\mathbb{R}}$ be two families. We assume that they satisfy the relations of \cref{relations sur BesselGamma/Omega,relations sur BesselYGamma/Omega} (which are the same) and that we have: $\forall d\in \mathbb{R} \setminus \TZ*,\,\allowbreak \forall q\in \TZ*,\; \sigma _{q} (\ZO_{d} )=0$. For~any $d\in \TZ*$, let us denote
    \[\forall n\in -\mathbb{N}^*,\ \ \Zrec_{d,n} :=0 \qquad  \text{ and } \qquad  \forall n\in \mathbb{N},\quad  \Zrec_{d,n} := \sum _{m=-n} ^{n} a_{n,m} \ZO_{d+m} .\]
    Then for any $d\in \TZ*$ and $n\in [\![1,b-1]\!]$, $\Zrec_{d,n}$ is the unique solution in $\Amacro(\Omega )$ of 
    \[ \left\{\begin{array}{r@{\;=\;}ll}
        (\Delta +k_{0} ^2 )\Zrec_{d,n} & 0 &\text{ in }\Omega  \\
        \Zrec_{d,n} & \sum\limits_{m=0}^\infty  \frac{\mu _{0} }{\mu _{1}} \,T_{m} \cdot  (\partial _{x} ^2 +k_{1} ^2 )^{m} \partial _{y} \Zrec_{d,n-1-2m} &\text{ on }\Gamma  \\
        \sigma _{q} (\Zrec_{d,n} ) & 0 &\ \forall q\in \TZ*\, \text{''} 
        \end{array}\right.\]
\end{quote}
\Cref{expression sing homog CL} applies this principle to $(\ZG_{a} ):=(\BesselGamma_{a} )$ and $(\ZO_{a} ):=(\BesselOmega_{a} )$. Thanks to \cref{relations sur BesselYGamma/Omega}, \ref{sigma de BesselYOmega} and \ref{sigma de BesselOmega}, it can also be applied to $(\ZG_{a} ):=(\BesselYGamma_{a} )$ and $(\ZO_{a} ):=(\BesselYOmega_{a} )$. Given \cref{expression sing homog CL avec BesselY}, this gives the desired result.
\end{demo}

\begin{textesansboite}
The next result in another (slight) generalization of \cref{expression sing homog CL} for certain angles, this time w.r.t.\ $n$.
\end{textesansboite}

\begin{proposition}{}
\label{expression sing heterog CL b=2 ou 4}%
Let us assume $b\in \{2,4\}$. For~any $d\in \TZ*$ and $n\in [\![0,4]\!]$, we have:
\begin{itemize}
    \item if $d\not\in -\mathbb{N}$, then $\mcla_{d,n}$ is a linear combination of the $\BesselOmega_{d+m}$ with $m\in [\![-n,n]\!]$,
    \item if $d\in -\mathbb{N}$, then $\mcla_{d,n}$ is a linear combination of the $\BesselOmega_{d+m}$ and $\BesselYOmega_{d+m}$ with $m\in [\![-n,n]\!]$.
\end{itemize}\fbi
\end{proposition}

\begin{textesansboite}
We recall that \cref{expression sing heterog CL} shows this result only for $n\in [\![0,b-1]\!]$ (with $b-1=1,3$ here). 
\Cref{expression sing heterog CL b=2 ou 4} allows us to show that, when $b\in \{2,4\}$, the construction of our effective models of order~3 and 5 (see \cref{chap: ME}) is as simple as the case $\pst \not\in \mathbb{Q}$. Thus only a few angles s.t.\ $\pst \in \mathbb{Q}$ require more complicated effective models (namely those s.t.\ $b\in \{1,3\}$).\\

The proof of \cref{expression sing heterog CL b=2 ou 4} is purely calculatory. We take the expressions of $\mcla_{d,n}$ given by \cref{ex: expression sing heterog CL} for $n\leqslant 4$ (which hold only when $n\leqslant b-1$ due to some divisions by vanishing factors when $n\geqslant b$) and we extend them to $n\in [\![b,4]\!]$ by showing that the ill-defined terms of these expression actually do~not appear when $b\in \{2,4\}$. This proof is very specific to these values of $b$ and $n$, but it seems to us that there~is something more general for even values of $b$. However, we did not investigate further as studying $\mcla_{d,n}$ for $b\in \{2,4\}$ and $n\leqslant 4$ is sufficient for effective models of order at most 5.
\end{textesansboite}

\begin{demo}
\def\rappmu{\mu_{\rm r}}%
\def\rappk{k_{\rm r}}%
\def\cotm#1{t_{#1}}%
We only prove the case $d\not\in -\mathbb{N}$, but the result can be extended to the case $d\in -\mathbb{N}$ using the same method as \cref{expression sing heterog CL avec BesselY}. We take the notations of \cref{ex: expression sing heterog CL}: $\rappmu:= \frac{\mu _{0} }{\mu _{1}}$, $\rappk := \big( \frac{2k_{1} }{k_{0}} \big)^2 $ and $\cotm{m} :=\operatorname{cotan}(m\Theta )$ for~any $m\in \mathbb{N}$. We also denote $\eta  :=(-1)^{d\Theta /\pi }$. Similarly to \cref{ex: expression sing heterog CL}, one can show that the following formulas of $\mcla_{d,n|\Gamma }$ hold when $n\leqslant b$ (this is actually an intermediate step to obtain the formulas of \cref{ex: expression sing heterog CL}):
\begin{itemize}
    \item $\displaystyle \mcla_{d,1|\Gamma } = \eta \, \Gamma (d+1)\. \left( \frac{k_{0} }{2} \right)^{1-d} \rappmu \.(\BesselGamma_{d-1} + \BesselGamma_{d+1} )$
    \item $\displaystyle \mcla_{d,2|\Gamma } = \eta \, \Gamma (d+1)\. \left( \frac{k_{0} }{2} \right)^{2-d} \rappmu^2  \cotm{1} \.(\BesselGamma_{d-2} - \BesselGamma_{d+2} )$
    \item $\displaystyle \mcla_{d,3|\Gamma } = \eta \, \Gamma (d+1)\. \left( \frac{k_{0} }{2} \right)^{3-d} \bigg[ \frac{3\rappmu^{3} \cotm{1}\cotm{2} + \rappmu}{3} (\BesselGamma_{d-3} + \BesselGamma_{d+3} )
    + \frac{3\rappmu^{3} \cotm{1}\cotm{2} +\rappmu(\rappk-1)}{3}  (\BesselGamma_{d-1} + \BesselGamma_{d+1} )\bigg]$
    \item $\displaystyle \mcla_{d,4|\Gamma } = \eta \, \Gamma (d+1)\. \left( \frac{k_{0} }{2} \right)^{4-d} \bigg[ \frac{3\rappmu^{4} \cotm{1}\cotm{2}\cotm{3} + \rappmu^2 (\cotm{1}+\cotm{3})}{3} (\BesselGamma_{d-4} - \BesselGamma_{d+4} )$\\[-1ex]
    \hbox to \linewidth{\hfill $\displaystyle {}+ \frac{3\rappmu^{4} \cotm{1}\cotm{2}(\cotm{1}+\cotm{3}) + \rappmu^2  \cotm{1}(2\rappk-3)+\rappmu^2 \cotm{3} }{3} \,(\BesselGamma_{d-2} - \BesselGamma_{d+2} ) \bigg]$.}
\end{itemize}

\underline{Case $b=2$:} In this case, $a$ is odd (because $\gcd(a,b)=1$) so we have $\cotm1=\operatorname{cotan}(a \frac{\pi }{2} )=0$ and $\cotm3=\operatorname{cotan}(3a \frac{\pi }{2} )=0$. Thus $\mcla_{d,2|\Gamma } =0$, which implies $\mcla_{d,2} =0$ by \cref{eq: pb verifie par mcl qd n>0,pb de Helmholtz dans A Omega}. If we compare this to the formula ``$\mcla_{d,2} = \Gamma (d+1) \left( \frac{k_{0} }{2} \right)^{2-d} \frac{\rappmu^2  \cotm{1}}{\sin (2\Theta )} \,(\BesselOmega_{d-2} + \BesselOmega_{d+2} )$'' given by \cref{ex: expression sing heterog CL} when $2\leqslant b-1$, we notice that the problem of the division by $\sin (2\Theta )=0$ is solved here by the fact that $\cotm1=0$.\\
Next, one can verify that we have \begin{align*}
    \mcla_{d,3|\Gamma } &= \cancel{\rappmu \partial _{y} \mcla_{d,2|\Gamma } } + \frac{\rappmu}{3} (\partial _{x} ^2 +k_{1} ^2 )\partial _{y} \mcla_{d,0|\Gamma } \\
    & =\eta \, \Gamma (d+1)\. \left( \frac{k_{0} }{2} \right)^{3-d} \Big[ \frac{\rappmu}{3} (\BesselGamma_{d-3} + \BesselGamma_{d+3} ) + \frac{\rappmu}{3} (\rappk-1) \. (\BesselGamma_{d-1} + \BesselGamma_{d+1} )\Big].
\end{align*}
This is similar to the formula of $\mcla_{d,3|\Gamma }$ when $3\leqslant b$, but with a few terms less because here $\mcla_{d,2} =0$. Continuing the calculations, one similarly obtains the reduced formula:
\[\mcla_{d,4|\Gamma } = \eta \,\Gamma (d+1)\. \left( \frac{k_{0} }{2} \right)^{4-d} \bigg[ \frac{\rappmu^2 (\cotm{1}+\cotm{3})}{3} (\BesselGamma_{d-4} - \BesselGamma_{d+4} ) + \frac{\rappmu^2  \cotm{1}(2\rappk-3)+\rappmu^2 \cotm{3} }{3} (\BesselGamma_{d-2} - \BesselGamma_{d+2} ) \bigg].\]
Since $\cotm1=\cotm3=0$, this yields $\mcla_{d,4|\Gamma } =0$ and $\mcla_{d,4} =0$.\medskip\smallskip

\underline{Case $b=4$:} Given \cref{expression sing heterog CL}, we already know that, when $n\leqslant 3$, $\mcla_{d,n}$ is a linear combination of the $\BesselOmega_{d+m}$ with $|m|\leqslant n$. As~for $\mcla_{d,4}$, since $\gcd(a,b)=1$, we have $2a\equiv 2\,[{4}]$ and $\{a,3a\}\equiv \{a,-a\}\equiv \{1,-1\}$ modulo 4. Hence $\cotm2=0$ and $\cotm1 +\cotm3 =\operatorname{cotan}( \frac{\pi }{4} )+\operatorname{cotan}(- \frac{\pi }{4} )=0$. Thus the expression of $\mcla_{d,4|\Gamma }$ given at the beginning of this proof reduces to
\[\mcla_{d,4|\Gamma } = \eta \, \Gamma (d+1)\, \left( \frac{k_{0} }{2} \right)^{4-d} \frac{\rappmu^2  \cotm{1}(2\rappk-3)+\rappmu^2 \cotm{3} }{3} (\BesselGamma_{d-2} - \BesselGamma_{d+2} ).\]
So $\mcla_{d,4} = \Gamma (d+1)\. \left( \frac{k_{0} }{2} \right)^{4-d} \frac{\rappmu^2  \cotm{1}(2\rappk-3)+\rappmu^2 \cotm{3}}{3\sin (2\Theta )} (\BesselOmega_{d-2} + \BesselOmega_{d+2} )$, which is well-defined because we do~not have to divide by $\sin (4\Theta )=0$ thanks to the absence of terms involving $\BesselOmega_{d\pm 4}$.
\end{demo}

\subsection{Extension of \texorpdfstring{\cref{sec: simplifications DA cas rationnel}}{Section 6.2.2}}
\label{sec: annexe simplifications DA cas rationnel}

\begin{textesansboite}
Let $\beta $ be the greatest element of $\mathbb{N}\cup \{\infty \}$ s.t., for any $d\in \TZ*$ and $n\in [\![0,\beta -1]\!]$, $\mcla_{d,n}$ is a linear combination of the $\BesselOmega_{d+m}$ with $m\in [\![-n,n]\!]$ and, if $d\in -\mathbb{N}$, of the $\BesselYOmega_{d+m}$ with $m\in [\![-n,n]\!]$. We have $\beta \geqslant b$ by \cref{expression sing heterog CL,expression sing heterog CL avec BesselY} and, if $b\in \{2,4\}$, then $\beta \geqslant 5$ by \cref{expression sing heterog CL b=2 ou 4}. Applying the reasoning of \cref{valeur de cmaicro cas rationnel} with $\beta $ instead of $b$ gives the following straightforward generalization of the first point of \cref{valeur de cmaicro cas rationnel}.
\end{textesansboite}

\begin{lemmesansboite}{}
\label{valeur de cmaicro cas rationnel generalise 1}%
Let $d,d'\in \TZ*$, $p\in \NN$ and $\ell \in \mathbb{N}$. Let us assume $d'\not\in -\mathbb{N}$.
\begin{itemize}
    \item If $p+d'<\beta $, then $\cmacro_{d,d',p,\ell } = \delta _{d,d'} \,\delta _{p,-d} \,\delta _{\ell ,0}$.
    \item If $p-d<\beta $, then $\cmicro_{d,d',p,\ell } = \delta _{d,d'} \,\delta _{p,d} \,\delta _{\ell ,0}$.
\end{itemize}
\end{lemmesansboite}

\begin{textesansboite}
Using the expression of $\mcla_{d,n}$ given by \cref{expression sing homog CL avec BesselY,expression sing heterog CL avec BesselY} when $d\in -\mathbb{N}$, we obtain the following (less obvious) generalization of \cref{valeur de cmaicro cas rationnel}, point 1.
\end{textesansboite}

\begin{lemmesansboite}{}
\label{valeur de cmaicro cas rationnel generalise 2}%
Let $d,d'\in \TZ*$, $p\in \NN$ and $\ell \in \mathbb{N}$. Let us assume $\red{d'\in -\mathbb{N}}$.
\begin{itemize}
    \item If $p+d'<\beta $ and \red{$p+d<-2d'$}, then $\cmacro_{d,d',p,\ell } = \delta _{d,d'} \,\delta _{p,-d} \,\delta _{\ell ,0}$.
    \item If $p-d<\beta $ and \red{$p+d'<0$}, then $\cmicro_{d,d',p,\ell } = \delta _{d,d'} \,\delta _{p,d} \,\delta _{\ell ,0}$.
\end{itemize}
\end{lemmesansboite}

\begin{demo}
The proof is similar to \cref{valeur de cmaicro cas rationnel}. We only give the main ideas, as the details are quite cumbersome, especially for $\cmacro$, like in \cref{valeur de cmaicro cas rationnel}.\medskip

\underline{Step 1:} Let us first treat $\cmicro_{d,d',p,\ell }$. By~\cref{expression de Rmacro phi avec les mcl}, we have when $p-d<\beta $:
\[\cmicro_{d,d',p,\ell } := \sigma _{d,p,\ell } \circ  \Eps \circ  \langle \Rmacro\rangle  (\phi_{d'} ) = \sigma _{d,p,\ell } \circ  \Eps \circ  \bigg( \fsumn \varepsilon ^{n} \.\mcla_{d',n} \bigg) = \sigma _{d,p,\ell } \circ  \Eps \circ  \bigg( \fsum_{n=0} ^{\red{\beta -1}} \varepsilon ^{n} \.\mcla_{d',n} \bigg).\]
By~definition of $\beta $, there exist complex numbers $a_{d',n,m}$ and $b_{d',n,m}$ s.t.\ we have for~any $n\in [\![0,\beta -1]\!]$:
\[\mcla_{d',n} = \sum _{m=-n} ^{n} \big(a_{d',n,m} \,\BesselYOmega_{d'+m} + b_{d',n,m} \,\BesselOmega_{d'+m} \big).\]
By~definition of $\BesselYOmega_{j}$ for $j\in \mathbb{Z}$ in \cref{eq: def BesselYOmega_n}, this implies that, for~any $q<-d'-n$, the component of $\mcla_{d',n}$ in $\mathcal{A}_{q} (\Omega )$ does~not involve $\log (\alpha z)$. Hence, under the assumptions $p-d<\beta $ and $p+d'<0$, the component of $\fsum_{n=0} ^{\beta -1} \varepsilon ^{n} \.\mcla_{d',n}$ in $\varepsilon ^{p-d} \mathcal{A}_{d} (\Omega )$ does~not involve $\log (\alpha z)$ (because $d<-d'-(p-d)$), which yields that the component of $\Eps\big(\fsum_{n=0} ^{\beta -1} \varepsilon ^{n} \.\mcla_{d',n} \big)$ in $\varepsilon ^{p} \mathcal{A}_{d} (\Omega )[\ln \varepsilon ]$ does~not involve $\ln \varepsilon $. The rest of this step is similar to \cref{valeur de cmaicro cas irrationnel,valeur de cmaicro cas rationnel}.\medskip\smallskip

\underline{Step 2:} This is essentially a technical adaptation of the second step of the proof of \cref{valeur de cmaicro cas rationnel} (that is already technical). The main idea of that previous proof is explained in \cref{fig: preuve valeur de cmaicro cas rationnel}. Here we only illustrate in \cref{fig: preuve valeur de cmaicro cas rationnel generalise} how this idea is modified. The main difference between \cref{fig: preuve valeur de cmaicro cas rationnel generalise,fig: preuve valeur de cmaicro cas rationnel} is that the green zone is an infinite band in the first one whereas it~is a parallelogram in the latter one. This is linked to the new equalities of \cref{valeur de cmaicro cas rationnel generalise 2} (written in red) that do~not appear in  \cref{valeur de cmaicro cas rationnel}.
\end{demo}

	\tikzset{external/export=true}%
	\tikzsetfigurename{fig-preuve_valeur_cmacro_cas_rationnel_bis-}%

\begin{figure}[h]
\def\couleurGdConeFig{blue}
\def\couleurPttConeFig{red}
\def\couleurGdTexteFig{blue!70!black}
\def\couleurPttTexteFig{red!70!black}
\def\vertFig{green!60!black}

\newcommand{\fillGrandCone}[2]{
	\begin{scope}[xshift={#1cm}, yshift={#2cm}]
		\clip (-10,10) -- (0,0) -- (10,0) -- (10,10) -- cycle;
		\shade[inner color=\couleurGdConeFig, outer color=white, opacity=0.3] ({-2.5*cos(67.5)},{-2*sin(67.5)}) circle (5); 
		\shade[inner color=\couleurGdConeFig, outer color=white, opacity=0.35] ({-2.5*cos(67.5)},{-2*sin(67.5)}) circle (9);
	\end{scope}
}
\newcommand{\drawGrandCone}[2]{
	\begin{scope}[xshift={#1cm}, yshift={#2cm}]
		\draw[very thick, \couleurGdConeFig] (-10,10) -- (0,0) -- (10,0);
	\end{scope}
}

\tikzfading[name=fade out, inner color=transparent!67, outer color=transparent!100]
\newcommand{\fillPetitCone}[3][1]{
	\begin{scope}[xshift={#2cm}, yshift={#3cm}]
		\clip (-10,10) -- (0,0) -- (10,0) -- (10,10) -- cycle;
		\fill[\couleurPttConeFig, path fading=fade out] ({-1.7*cos(67.5)*#1},{-1.7*sin(67.5)*#1}) circle (2.5*#1);
		\fill[\couleurPttConeFig, path fading=fade out] ({-1.7*cos(67.5)*#1},{-1.7*sin(67.5)*#1}) circle (4*#1);
	\end{scope}
}
\newcommand{\drawPetitCone}[3][1]{
	\begin{scope}[xshift={#2cm}, yshift={#3cm}]
		\clip (-10,10) -- (0,0) -- (10,0) -- (10,10) -- cycle;
		\draw[very thick, \couleurPttConeFig, path fading=east] (0,0) -- (3.5*#1,0) ++(0,1pt);
		\draw[very thick, \couleurPttConeFig, path fading=west] (0,0) -- ({-3.5/sqrt(2)*#1},{3.5/sqrt(2)*#1}) ++(1pt,1pt);
	\end{scope}
	\begin{scope}[xshift={#2cm}, yshift={#3cm}]
		\clip (-10,10) -- (0,0) -- (10,0) -- (10,-10) -- (-10,-10) -- cycle;
		\draw[very thick, \couleurPttConeFig, path fading=east] (0,0) -- (3.5*#1,0) ++(0,-1pt);
		\draw[very thick, \couleurPttConeFig, path fading=west] (0,0) -- ({-3.5/sqrt(2)*#1},{3.5/sqrt(2)*#1}) ++(-1pt,-1pt);
	\end{scope}
}
\newcommand{\filldrawPetitCone}[3][1]{
	\fillPetitCone[#1]{#2}{#3}
	\drawPetitCone[#1]{#2}{#3}
}

\def\croixFig{\raisebox{1pt}{\small\bf+}}

\def\xmin{-4.4}
\def\xmax{4.4}
\def\ymin{-0.4}
\def\ymax{6.2}
\def\dprimeFig{-1.87}
\def\dsecondeFig{2.65}
\def\pprimeFig{4.45}
\def\bFig{3.1}


\raggedleft 
\begin{tikzpicture}[scale=0.9]
	\begin{scope}
		\clip (\xmin,\ymin) rectangle (\xmax,\ymax);

		\fillGrandCone{\dprimeFig}{-\dprimeFig}
		\fillPetitCone[0.65]{\dprimeFig}{-\dprimeFig}
		\fillPetitCone[0.65]{\dprimeFig}{-\dprimeFig} 
		\filldrawPetitCone[0.5]{-2.8}{4.05}
		\filldrawPetitCone[0.55]{-0.65}{5}
		\filldrawPetitCone[0.45]{2.65}{2.35}
		\filldrawPetitCone[0.5]{\dsecondeFig}{\pprimeFig}

		\fill [pattern={Lines[angle=45,line width=.5pt,distance=7pt]}, pattern color=\vertFig] (-\bFig,\bFig) -- (\dprimeFig,-\dprimeFig) -- (-\dprimeFig,-\dprimeFig) -- (-2*\dprimeFig-\bFig,\bFig) -- cycle;
		\draw[thick, green!65!black] (-\bFig,\bFig) -- (-2*\dprimeFig-\bFig,\bFig) -- (-\dprimeFig,-\dprimeFig);
		\draw[thick, dashed, green!65!black] (-2*\dprimeFig-\bFig,\bFig) -- (0,-2*\dprimeFig);
		\draw (0,-2*\dprimeFig) node{\croixFig} node[right]{$\.-2d'\mathstrut$};
		\draw (0,\bFig) node{\croixFig} ++(0,-0.05) node[above left]{$b-d'\!$};
		\draw [<->,>=latex] (-\bFig,\bFig-0.1) -- node[near start,left]{\raisebox{-5ex}{$b$}} (-\bFig,-\dprimeFig);
		\draw[dotted] (-\bFig,-\dprimeFig) -- (\dprimeFig,-\dprimeFig) -- (0,0);

		\drawGrandCone{\dprimeFig}{-\dprimeFig}
		\drawPetitCone[0.65]{\dprimeFig}{-\dprimeFig}
		\drawPetitCone[0.65]{\dprimeFig}{-\dprimeFig} 

		\draw[\couleurGdTexteFig] (\dprimeFig,-\dprimeFig) 
			node {\croixFig}
			node[below] {$(d',-d')$};
		\draw[\couleurPttTexteFig] (\dsecondeFig,\pprimeFig) 
			node {\raisebox{1pt}{\croixFig}}
			node [below] {$(d'',p'-d'')$};
	\end{scope}

	\draw[<-,>=latex,\couleurGdTexteFig] ({-\ymax*0.625},{\ymax*0.625}) to[bend right=20] ++(-0.5,-0.5) node[below left]{\parbox{28mm}{\raggedleft support of\\ $\varepsilon^{-d'}\langle\Rmacro\rangle(\phi_{d'})$}\!\!\!}; 
	\draw[<-,>=latex,\couleurPttTexteFig] (\dprimeFig+0.9,-\dprimeFig) to[bend left=25] ++(-0.1,-1.15) node[left]{\parbox{42mm}{
		support of $\Eps^{-1}\circ\langle\Rmicro\rangle(\phi_{d'})$}}; 
	\draw[<-,>=latex,\couleurPttTexteFig] (\dsecondeFig+1.2,\pprimeFig) to[bend right=20] ++(0.7,-0.3) node[right]{\parbox{45mm}{
		support of a term $\cmicro_{d''\!,d'\!,p'+d'\!,\ell'}$\\ $\cdot\,\varepsilon^{p'}\ln^{\ell'}\!\varepsilon\cdot\Eps^{-1}\circ\langle\Rmicro\rangle(\phi_{d''})$}}; 
	\draw[<-,>=latex,\vertFig] (\dprimeFig+2.7,-\dprimeFig) to[bend right=20] ++(0.3,-0.8) node[right]{\parbox{30mm}{
		only one red cone\\ in the green zone}}; 

	\draw [->] (\xmin,0) -- ({\xmax+0.25},0) node [right] {degree in the spaces $\mathcal{A}_q(\Pi)$};
	\draw [->] (0,\ymin) -- (0,{\ymax+0.3}) node [right]{power of $\varepsilon$};
\end{tikzpicture}

\longCaption{Schematic illustration of \cref{eq: preuve valeur de cmaicro cas rationnel 1} in the case of \cref{valeur de cmaicro cas rationnel generalise 2}}{. Each cone represents the convex closure of the support of a formal series in the equation. They are infinite cones, but they are truncated for visual clarity. In the hatched zone, we have $\cmicro_{d''\!,d'\!,p'+d'\!,\ell'} = \delta_{d'\!,d''}\,\delta_{p'\!,-d'}\,\delta_{\ell'\!,0}$, so $\varepsilon^{-d'}\langle\Rmacro\rangle(\phi_{d'})$ coincides with $\Eps^{-1}\circ\langle\Rmicro\rangle(\phi_{d'})$ in this zone.}
\label{fig: preuve valeur de cmaicro cas rationnel generalise}
\end{figure}%
	\tikzset{external/export=false}%
	\tikzsetfigurename{main-figure}

\begin{textesansboite}
\Cref{valeur de cmaicro cas rationnel generalise 1,valeur de cmaicro cas rationnel generalise 2} imply the following generalization of the second point of \cref{valeur de cmaicro cas rationnel}. Compared to \cref{valeur de cmaicro cas rationnel}, the assumptions ``$d'>-b$'' and ``$p<b$'' are relaxed.
\end{textesansboite}

\begin{lemmesansboite}{}
\label{valeur de cfinal cas rationnel generalise}%
Let $d\in -\TN$, $d'\in \TZ*$, $p\in \NN$ and $\ell \in \mathbb{N}$. Let us assume that we have $p<\min (\beta ,2b+2\pst)$ and, if $d'\in -\mathbb{N}$, then $p<b+\pst$. Then we have:
\[\cfinal_{d,d',p,\ell } = \left\{\begin{array}{ll}
	\omega ^{2n} \sigma _{d} (\ccDA_{d',n} ) & \text{if $d'\in \TN$, $\ell =0$ and there exists $n\in \mathbb{N}$ s.t.\ }  p=d'-d+2n\\
	0 & \text{otherwise}.
\end{array}\right.\]
\end{lemmesansboite}

\begin{demo}
We recall that we have
\begin{equation}{\label{eq: valeur de cfinal cas rationnel generalisation preuve 1}}
    \cfinal_{d,d',p,\ell } := \! \sum _{\substack{(p_{1} ,p_{2} )\in (\mathbb{N}-d)\times (\mathbb{N}+d')\\ n\in \mathbb{N},\, p_{1} +p_{2} +2n=p}}
    \ \sum _{\substack{d_{1} \in \TZ*\\ p_{1} +d_{1} \in \mathbb{N}}} \ \sum _{\substack{d_{2} \in \TN\\ p_{2} -d_{2} \in \mathbb{N}}}
    \ \sum _{\substack{(\ell _{1} ,\ell _{2} )\in \mathbb{N}^2 \\ \ell _{1} +\ell _{2} =\ell }}
    \cmacro_{d,d_{1} ,p_{1} ,\ell _{1} }\cdot \omega ^{2n} \sigma _{d_{1}} (\mathcal{S} ^{\infty } _{d_{2} ,n} )\cdot \cmicro_{d_{2} ,d'\!,p_{2} ,\ell _{2} } .
\end{equation}
Let us show that any non-zero term satisfies
\[p_{1} +d_{1} <\beta  \text{ and } p_{2} -d_{2} <\beta  \text{ and } (d_{1} \not\in -\mathbb{N} \text{ or } p_{1} +d<-2d_{1} ) \text{ and } (d'\not\in -\mathbb{N} \text{ or } p_{2} +d'<0)\]
in order to apply \cref{valeur de cmaicro cas rationnel generalise 1,valeur de cmaicro cas rationnel generalise 2}. We do it by contraposition, by showing that if a non-zero term verifies
\begin{equation}{\label{eq: valeur de cfinal cas rationnel generalisation preuve 2}}
    p_{1} +d_{1} \geqslant \beta  \text{ or } p_{2} -d_{2} \geqslant \beta  \text{ or } (d_{1} \in -\mathbb{N} \text{ and } p_{1} +d\geqslant -2d_{1} ) \text{ or } (d'\in -\mathbb{N} \text{ and } p_{2} +d'\geqslant 0)
\end{equation}
then we do~not have $\big[p<\min (\beta ,2b+2\pst) \text{ and } (d'\in -\mathbb{N} \Rightarrow  p<b+\pst)\big]$, which means that we have 
\[\textstyle p\geqslant \beta  \text{ or } p\geqslant 2b+2\pst \text{ or } (d'\in -\mathbb{N} \text{ and } p\geqslant b+\pst).\]
Let us treat separately each subcase of \cref{eq: valeur de cfinal cas rationnel generalisation preuve 2}.
\begin{itemize}
    \item Let us assume $p_{1} +d_{1} \geqslant \beta $. If $d_{1} >0$, then $\sigma _{d_{1}} (\mathcal{S} ^{\infty } _{d_{2} ,n} ) =\delta _{d_{1} ,d_{2}} \delta _{n,0}$. Since we are looking at a non-zero term, we have $d_{1} =d_{2}$ and $n=0$. Thus $p=(p_{1} +d_{1} ) +(p_{2} -d_{2} ) \geqslant \beta +0$. Otherwise, we have $d_{1} <0$, which yields $p=(p_{1} +d_{1} )-d_{1} +(p_{2} -d_{2} )+d_{2} +2n \geqslant \beta +0+0+0+0=\beta $.
    \item If $p_{2} -d_{2} \geqslant \beta $, then $p=(p_{1} +d)-d +(p_{2} -d_{2} )+d_{2} +2n\geqslant 0+0+\beta +0+0=\beta $.
    \item If $d_{1} \in -\mathbb{N}$ and $p_{1} +d\geqslant -2d_{1}$, then $d_{1} \leqslant -b$ (because $\mathbb{Z}\cap \TZ=b\mathbb{Z}$), so $p=(p_{1} +d)-d +(p_{2} -d_{2} )+d_{2} +2n\geqslant 2b+\pst+0+\pst+0=2b+2\pst$.
    \item If $d' \in -\mathbb{N}$ and $p_{2} +d'\geqslant 0$, then $d' \leqslant -b$, so $p=(p_{1} +d)-d +(p_{2} +d')-d' +2n\geqslant 0+\pst+0+b+0=b+\pst$.
\end{itemize}
By~\cref{valeur de cmaicro cas rationnel generalise 1,valeur de cmaicro cas rationnel generalise 2} we deduce that any non-zero term of \cref{eq: valeur de cfinal cas rationnel generalisation preuve 1} satisfies $d=d_{1}$, $d_{2} =d'$, $p_{1} =-d$, $p_{2} =d$, $\ell _{1} =\ell _{2} =0$ and $\cmacro_{d,d_{1} ,p_{1} ,\ell _{1} }=\cmicro_{d_{2} ,d'\!,p_{2} ,\ell _{2} }=1$. This implies the desired formula.
\end{demo}

\begin{theoreme}{}
\label{simplifications cas rationnel version forte}%
Let 
\[p_{\min } := \max \!\left( \frac{2b}{a} ,\, \min \!\left( \beta  ,\, 2b\. \frac{a+1}{a} \right)\right) = \min \!\left( \max \!\left(\beta , \frac{2b}{a} \right)\!,\, 2b\. \frac{a+1}{a} \right).\]
We have: \[\forall p\in \NN,\quad  p<p_{\min }  \Rightarrow (\forall \ell  \neq 0, u_{p,\ell } =0).\]
Moreover, \cref{eq: simplifications cas rationnel} (i.e.\ the explicit expression of $\sigma _{d} (u_{p,\ell } )$ in terms of the coefficients $\sigma _{d} (\ccDA_{q,n} )$) holds for~any $p\in \NN$ s.t.\ $p<p_{\min }$.\\
In particular, when $b\not\in \{1,3\}$, we have $p_{\min } \geqslant 5$.
\end{theoreme}

\begin{textesansboite}
This extends \cref{simplifications cas rationnel}, which states the same result with $p_{\min } =\max ( \frac{2b}{a} ,b)$. The value of $p_{\min }$ is especially improved when $b\in \{2,4\}$. This is coherent with \cref{tab: ordre du premier ln eps}, which shows that the value of $p_{\min }$ given by \cref{simplifications cas rationnel} was not optimal.
\end{textesansboite}

\begin{demo}
The proof is similar to \cref{simplifications cas rationnel}. Let $p<p_{\min }$, $\ell \in \mathbb{N}$ and $d\in -\TN$. \Cref{construction directe des upl} states:
\begin{equation}{\label{eq: simplifications cas rationnel version forte preuve 1}}
    \sigma _{d} (u_{p,\ell } ) = \!\! \sum _{\substack{p'\in \NN\\ p-p'\in \NN+ 2\pst}} \; \sum _{\substack{d'\in \TZ*\\ p'+d'\in \intvl{0,p+d}}} \! \sum _{\ell '=0} ^{\ell } \cfinal_{d,d',p-p',\ell -\ell '} \cdot \sigma _{d'} (\uDA_{p',\ell '} ) .
\end{equation}
Let us show that the right-hand side of \cref{eq: simplifications cas rationnel version forte preuve 1} is equal to
\begin{equation}{\label{eq: simplifications cas rationnel version forte preuve 2}}
    \sum _{\substack{n\in \mathbb{N}\\ 2n\leqslant p+d}} \sum _{\substack{d'\in \TN\\ d'\leqslant p+d-2n}} \omega ^{2n} \sigma _{d} (\ccDA_{d',n} )\cdot \sigma _{d'} (\uDA_{p-d'+d-2n,0} )
\end{equation}
(which, when $\ell =0$, is the right-hand side of \cref{eq: simplifications cas rationnel}). If $p< \frac{2b}{a} =2\pst$, then the sums in \cref{eq: simplifications cas rationnel version forte preuve 1,eq: simplifications cas rationnel version forte preuve 2} are both empty, so they vanish. Now, let us assume $p<\min (\beta ,2b \frac{a+1}{a} )$. Any non-zero term in the right-hand side of \cref{eq: simplifications cas rationnel version forte preuve 1} satisfies $p-p'\leqslant p<\min (\beta ,2b+2\pst)$ because $2b \frac{a+1}{a} = 2b+2\pst$. In addition, if $d'\in -\mathbb{N}$, we also have $-b\geqslant d'\geqslant p'-\pst$ by \cref{precision de l'ansatz sur uDA}, so $p'\geqslant b+\pst$, which yields $p-p'\leqslant 2b+2\pst-(b+\pst)=b+\pst$. Therefore we can apply \cref{valeur de cfinal cas rationnel generalise} to any non-zero term of \cref{eq: simplifications cas rationnel version forte preuve 1}, which shows that it~is equal to \cref{eq: simplifications cas rationnel version forte preuve 2}. The rest is identical to the proof of \cref{simplifications cas rationnel}, except the last statement about $b\not\in \{1,3\}$.\\
Finally, let us prove the last statement.
\begin{itemize}
    \item If $b\geqslant 5$, then we have $p_{\min } \geqslant  \min \!\left( \beta  ,\, 2b\. \frac{a+1}{a} \right) \geqslant  \min (b,2b)\geqslant 5$ because $\beta \geqslant b$ by \cref{expression sing heterog CL}.
    \item If $b\in \{2,4\}$, then $p_{\min } \geqslant  \min \!\left( \beta  ,\, 2b\. \frac{a+1}{a} \right) =\min (\beta ,2b+2\pst)\geqslant  \min (5,2\cdot 2+2\cdot  \frac{1}{2} )=5$ by \cref{expression sing heterog CL b=2 ou 4}.
\end{itemize}
\end{demo}

\end{subappendices}

\chapter{Modal decomposition of the far fields}
\label{chap: decompo modale CL}

\section{Statement of the results}
\label{sec: intro convergence des mcl et uDA}

\begin{textesansboite}
In this section, we establish a modal decomposition of the far fields $u_{p,\ell }$ in the vicinity of the layer. This is an important step towards effective models. \Cref{chap: ME} will explain how summing the decompositions of the first far fields provides an approximate modal decomposition of the exact solution $u_{\varepsilon }$, which can then be used to approximate the behavior of $u_{\varepsilon }$ around the vicinity of the layer with a Dirichlet-to-Neumann operator. The modal decomposition of $u_{p,\ell }$ will also be useful to justify the numerical computation of the far fields by induction in \cref{sec: numerique des CL}.\\

The modes of the decomposition of $u_{p,\ell }$ are the $\mcll_{d,n}$, which are the limit functions of the series $\mcla_{d,n} \in \Amacro(\Omega )$:
\[\forall d\in \txtTZ*,\; \forall n\in \mathbb{N},\; \forall r>0,\;\forall \theta \in (0,\Theta ), \qquad  \mcll_{d,n} (r,\theta ) := \lim _{q\rightarrow \infty } T_{<q} (\mcla_{d,n} )(r,\theta ).\]
\Cref{convergence des mcl} proves that the $\mcla_{d,n}$ are (locally uniformly) convergent series, so $\mcll_{d,n}$ is well defined. The coefficients of the modal decomposition are the $\sigma _{d} (u_{p,\ell } )$, which are involved in the matching conditions (see \cref{condition de raccord,def: construction des termes}). The main result of this section is the following theorem.
\end{textesansboite}

\begin{paragraphesansboite}{Notations:}
Let $\Rsource := \operatorname{dist}(0,\operatorname{supp}\fsource)$. Since we assumed $\operatorname{dist}(\Gamma ,\operatorname{supp}\fsource)>0$, we have $\Rsource>0$. In addition, all over this section, $C_{\dots }$ will denote any positive constant that depends \emph{only} on the variables in ``...'' and on $(\Theta ,\mu _{0} ,\mu _{1} ,\rho _{0} ,\rho _{1} ,\omega )$. Its value may change from line to line.
\end{paragraphesansboite}

\begin{theoremesecable}{modal decomposition of the far fields around the corner}
\label{decompo modale des CL}%
Let $p\in \NN$ and $\ell \in \mathbb{N}$.
\begin{enumerate}
    \item On $\Omega  \cap  B(0,\Rsource)$, the series below converges and the following equality holds:
    \begin{equation}{\label{eq: decompo modale des CL}}
        u_{p,\ell } \,=\, \sum _{n=0} ^{\lfloor {p}\rfloor } \sum _{\substack{d\in \TZ*\\ d\geqslant -p+n}} \sigma _{d} (\uDA_{p-n,\ell } ) \,\mcll_{d,n} .
    \end{equation}
    In addition, the convergence is uniform on $\Omega  \cap  B(0,R)$ for~any $R\in (0,\Rsource)$.
    \item For~any $R\in (0,\Rsource)$, there exists $C_{p,\fsource,R} >0$ s.t.\ we have: $\forall d\in \TZ*,\, \sigma _{d} (u_{p,\ell } ) \leqslant  C_{p,\fsource,R} \,R^{-d} $.
    \item For~any $N\in \mathbb{N}$, let us denote $\rho _{N}$ the remainder at order $N\pst$ of the series \cref{eq: decompo modale des CL}:
    \[\rho _{N} := u_{p,\ell } - \sum _{n=0} ^{\lfloor {p}\rfloor } \sum _{\substack{d\in \TZ*\\ -p+n\leqslant d\leqslant N\pst}} \sigma _{d} (\uDA_{p-n,\ell } ) \,\mcll_{d,n} .\]
    Let $R\in (0,\Rsource)$. There exists $C_{p,\fsource,R} >0$ s.t., for~any $N\in \mathbb{N}$, $r\in (0,R)$ and $\theta \in [0,\Theta ]$, we have
    \begin{equation}{\label{eq: vitesse de cvg decompo modale des CL}}
        |\rho _{N} (r,\theta )| + r|\nabla  \rho _{N} (r,\theta )| \leqslant  C_{p,\fsource,R} \left( \frac{r}{R} \right)^{(N+1)\pst-\lfloor {p}\rfloor } (1+|{\ln r}|)^{2p+1} .
    \end{equation}
\end{enumerate}
\end{theoremesecable}

\begin{remarquenumerotee}
\label{rq: decompo modale des CL}%
\begin{itemize}
    \item \cref{eq: decompo modale des CL} has the same structure as the formula
    \begin{equation*}
         \uDA_{p,\ell } = \sum _{n=0} ^{\lfloor {p}\rfloor } \fsum_{\substack{d\in \TZ*\\ d\in p-n+\NN}} \sigma _{d} (\uDA_{p-n,\ell } ) \,\mcla_{d,n}  \tag{\ref{eq: expression uDA avec mcl}}
    \end{equation*}
    stated by \cref{expression uDA avec mcl}. Note that \cref{eq: expression uDA avec mcl} is an equality in $\Amacro(\Omega )$, so %
    it means that the formal series on both sides are equal term by term%
    , with no guarantee of convergence. 
    Unlike \cref{decompo modale des CL}, it was proven with a simple algebraic calculation.
    \item In the particular case $(p,\ell )=(0,0)$, \cref{eq: decompo modale des CL} is the usual modal decomposition into Bessel functions of a solution of $\Delta u+k_{0} ^2 u=0$ in a corner with homogeneous Dirichlet condition. Indeed, \cref{expression sing homog CL} states that $\mcll_{d,0}$ is proportional to $J_{d} (k_{0} r)\sin (d\theta )$ for~any $d\in \TN$.
    \item $\Rsource$ is an optimal radius of convergence, as the decomposition of $u_{0,0}$ using Bessel functions is valid only where $\Delta u_{0,0} +k_{0} ^2 u_{0,0} =0$, i.e.\ in $B(0,\Rsource)$.
    \item The error estimate \cref{eq: vitesse de cvg decompo modale des CL} gives a geometric convergence of the remainder w.r.t.\ $N$ when $N\rightarrow \infty $, as~well~as an algebraic convergence w.r.t.\ $r$ when $r\rightarrow \infty $, because the constant $C_{p,\fsource,R}$ is independent of $N$ and $r$. Moreover, the factor $r^{(N+1)\pst-\lfloor {p}\rfloor }$ can be roughly interpreted as the order of magnitude of the greatest neglected term, i.e.\ $\mcll_{(N+1)\smallpst,\lfloor {p}\rfloor }$. Indeed, this function is a sum of a non zero element of $\mathcal{A}_{(N+1)\smallpst-\lfloor {p}\rfloor } (\Omega )$ and of elements of the $\mathcal{A}_{d} (\Omega )$ with $d>(N+1)\pst-\lfloor {p}\rfloor $.
    \item The factor $(1+|{\ln r}|)^{2p+1}$ is relevant only if $\pst \in \mathbb{Q}$, as otherwise the far fields involve no $\ln r$ (see \cref{sec: simplifications DA cas irrationnel}). Moreover, the power $2p+1$ is~not optimal in most cases, even if $\pst \in \mathbb{Q}$. But we do~not seek optimality here as the logarithmic behavior is secondary compared to the power or $r$.
    \item Since $u_{p,\ell }$ belongs to $H^{1} _{\mathrm{loc}} (\Omega )$, \cref{eq: vitesse de cvg decompo modale des CL} is to be understood almost everywhere in $(r,\theta )$, or for a particular representative of $u_{p,\ell }$.
\end{itemize}
\end{remarquenumerotee}

\begin{textesansboite}
If $\pst$ is irrational, then \cref{expression sing heterog CL} shows that the $\mcll_{d,n}$ are finite linear combinations of Bessel functions (times some sinus of $\theta $). Since Bessel functions are well implemented in most coding languages, the $\mcll_{d,n}$ can be precisely and efficiently computed at any point $(r,\theta )$.\\

\def\dmax{N\pst}\def\dmaxplst{(N+1)\pst}\def\qmax{q_0}
However, when $\pst$ is rational, the $\smash{\mcla_{d,n} }$ have in general more complex expressions involving $\ln r$. In this case, we do~not know any expression of $\smash{\mcll_{d,n} }$ using usual special functions like the Bessel functions, so the only way to compute them is by summing a large number of terms of the series $\smash{\mcla_{d,n} }$. Hence one could choose an order $\dmax$ to truncate the modal decomposition and then a truncation order $\qmax$ to compute each $\mcll_{d,n}$:
\begin{equation}{\label{eq: troncation naive decompo modale CL}}
    u_{p,\ell } \;\approx \; \sum _{n=0} ^{\lfloor {p}\rfloor } \sum _{\substack{d\in \TZ*\\ -p+n\leqslant d\leqslant \dmax}} \!\sigma _{d} (\uDA_{p-n,\ell } ) \,\mcll_{d,n}
    \;\approx \; \sum _{n=0} ^{\lfloor {p}\rfloor } \sum _{\substack{d\in \TZ*\\ -p+n\leqslant d\leqslant \dmax}} \!\sigma _{d} (\uDA_{p-n,\ell } ) \cdot T_{<\qmax} (\mcla_{d,n} )
\end{equation}
The order $\qmax$ should be chosen in function of $N$ so that it does~not deteriorate the error of the truncated modal decomposition given by \cref{eq: vitesse de cvg decompo modale des CL}, while minimizing the number of terms of the $\mcla_{d,n}$ that are computed. The next theorem shows that one can take $\qmax:=\dmaxplst-\lfloor {p}\rfloor $. Any higher value of $\qmax$ does~not improve the error made by \cref{eq: troncation naive decompo modale CL}, or at least its order of magnitude. The idea is that any term of order $q>\dmaxplst-\lfloor {p}\rfloor $ is negligible compared to the error estimate \cref{eq: vitesse de cvg decompo modale des CL}, which is roughly $r^{\dmaxplst-\lfloor {p}\rfloor }$. Note however that this theorem is~not obvious to prove as it approximates the last considered modes with very few terms.
\end{textesansboite}

\begin{theoreme}{}
\label{decompo en modes tronques des CL}%
Let $p\in \NN$ and $\ell \in \mathbb{N}$. For~any $N\in \mathbb{N}$ and $q\in \mathbb{R}$, let us denote
\begin{equation}{\label{eq: reste decompo en modes tronques des CL}}
    \rho _{N,q} := u_{p,\ell } - \sum _{n=0} ^{\lfloor {p}\rfloor } \sum _{\substack{d\in \TZ*\\ -p+n\leqslant d\leqslant N\pst}} \sigma _{d} (\uDA_{p-n,\ell } ) \cdot T_{<q} (\mcla_{d,n} ).
\end{equation}
Let $R\in (0,\Rsource)$. There exists $C_{p,\fsource,R} >0$ s.t., for~any $N\in \mathbb{N}$, $q\in \mathbb{R}$, $r\in (0,R)$ and $\theta \in [0,\Theta ]$, we have
\begin{equation}{\label{eq: vitesse de cvg decompo en modes tronques des CL}}
    |\rho _{N,q} (r,\theta )|+r|\nabla  \rho _{N,q} (r,\theta )| \leqslant  C_{p,\fsource,R} \left( \frac{r}{R} \right)^{\min \left((N+1)\pst -\lfloor {p}\rfloor ,\, q\right)} (1+|{\ln r}|)^{2p+1} .
\end{equation}\fbeq[1.4]
\end{theoreme}

\begin{textesansboite}
One can show that, with an appropriate choice of $N$, \cref{eq: reste decompo en modes tronques des CL} is equal to $u_{p,\ell } -T_{<q} (\uDA_{p,\ell } )$. Hence \cref{decompo en modes tronques des CL} implies an interesting result about the behavior of $u_{p,\ell }$ near the corner. Let us recall that \cref{validation construction et DA des champs} states that $u_{p,\ell } = T_{\leqslant q} (\uDA_{p,\ell } )+o_{\partial } (r^{q} )$ for~any $q\in \mathbb{R}$. By~\cref{def: petit o derivable}, this implies that \emph{there exists} $R_{p,\fsource,q} >0$ (which can be very small) and $n\in \mathbb{N}$ s.t.\ for~any $r\in (0,R_{p,\fsource,q} )$ and $\theta  \in [0,\Theta ]$:
\[|\rho _{q} (r,\theta )|+r|\nabla  \rho _{q} (r,\theta )|  \leqslant  C_{p,\fsource,q} \,r^{q} \,(1+|{\ln r}|)^{n} \qquad \text{ with }\rho _{q} := u_{p,\ell } -T_{<q} (\uDA_{p,\ell } ).\]
Since the constant depends on $q$, no convergence is guaranteed when $q\rightarrow \infty $. In comparison, the next result proves the convergence w.r.t.\ $q$ and it gives an error estimate that is explicit w.r.t.\ $(q,r)$ and holds in $B(0,\Rsource)$ entirely.
\end{textesansboite}

\begin{corollaire}{}
\label{convergence de uDA}%
Let $p\in \NN$ and $\ell \in \mathbb{N}$. We have
\[\qquad  T_{<q} (\uDA_{p,\ell } ) \xrightarrow[{q\rightarrow \infty }]{} u_{p,\ell } \qquad \text{ on }\Omega  \cap B(0,\Rsource)\]
with uniform convergence on $\Omega  \cap  B(0,R)$ for~any $R\in (0,\Rsource)$. Moreover, let us denote the remainder $\rho _{q} := u_{p,\ell } -T_{<q} (\uDA_{p,\ell } )$ for~any $q$. Let $R\in (0,\Rsource)$. There exists $C_{p,\fsource,R} >0$ s.t., for~any $q\in \mathbb{R}$, $r\in (0,R)$ and $\theta \in [0,\Theta ]$, we have
\begin{equation}{\label{eq: vitesse de cvg de uDA}}
    |\rho _{q} (r,\theta )|+r|\nabla  \rho _{q} (r,\theta )| \leqslant  C_{p,\fsource,R} \left( \frac{r}{R} \right)^{q} (1+|{\ln r}|)^{2p+1} .
\end{equation}\fbeq[1.4]
\end{corollaire}

\begin{remarque}
Like in \cref{decompo modale des CL}, the error estimates \cref{eq: vitesse de cvg decompo en modes tronques des CL} and \cref{eq: vitesse de cvg de uDA} behaves like the first neglected term of $\uDA_{p,\ell }$. The other points of \cref{rq: decompo modale des CL} about the optimality of the radius of convergence $\Rsource$ and the logarithmic factor also apply here.
\end{remarque}

\begin{textesansboite}
The overall outline of this section is as follows.
\begin{itemize}
    \item In~\cref{sec: convergence des mcl}, we prove the convergence of the formal series $\mcla_{d,n}$, so~that $\mcll_{d,n}$ is well defined. We also provide a precise estimate of $\mcll_{d,n} -T_{<q} (\mcla_{d,n} )$ in order to deduce the error estimates \cref{eq: vitesse de cvg decompo modale des CL} and \cref{eq: vitesse de cvg decompo en modes tronques des CL}. This result is stated in \cref{convergence des mcl} and its proof is the longest part.
    \item Then, in \cref{sec: preuve decompo modale CL}, we deduce \cref{decompo modale des CL,decompo en modes tronques des CL,convergence de uDA} (in this order, see resp.\ p.~\pageref{preuve: decompo modale des CL}, \pageref{preuve: decompo en modes tronques des CL}, \pageref{preuve: convergence de uDA}).\\
\end{itemize}

Let us detail the first point. If $\pst$ is irrational, we already know that the $\mcla_{d,n}$ are convergent because \cref{expression sing heterog CL} states that $\mcla_{d,n}$ is a linear combination of the $\BesselOmega_{a}$. Obtaining an error estimate for this convergence requires to look more closely at the series defining the Bessel functions. This is technical but reasonable.\\

The main difficulty lies in the case $\pst \in \mathbb{Q}$. No nice expression of $\mcla_{d,n}$ is known in this situation, but only its definition by induction using the operators $\ROmegaOmega$, $\RGammaOmega$, $\partial _{y|\Gamma }$ and $\partial _{x|\Gamma }$ (see \cref{def: mcl}). Our strategy is to estimate the terms of $\mcla_{d,n}$, in order to deduce that this series converges absolutely and estimate its remainder. It~is divided in three steps.
\begin{itemize}
    \item In \cref{sec: convergence des mcl cas rationnel 1}, we define norms on the spaces $\mathcal{A}_{d} (\Omega )$ and $\mathcal{A}_{d} (\Gamma )$ and we compute continuity constants for the operators $\RGammaOmega$, $\ROmegaOmega$, $\partial _{y|\Gamma }$ and $\partial _{x|\Gamma }$ on some subspaces of $\mathcal{A}_{d} (\Omega )$ and $\mathcal{A}_{d} (\Gamma )$ that are denoted $\AOprecis{d,n,\delta }$ and $\AGprecis{d,n}$.
    \item In \cref{sec: convergence des mcl cas rationnel 2}, we do a similar work on formal series. We define discrete weighted norms on $\Amacro(\Omega )$ and $\Amacro(\Gamma )$, related to the norms on $\mathcal{A}_{d}$. They are s.t.\ any element with finite norm is an~absolutely converging series, and they involve a parameter $\eAMP$ that measures the convergence speed of the series. Then we deduce, from the estimates on $\AOprecis{d,n,\delta }$ and $\AGprecis{d,n}$, continuity constants for $\RGammaOmega$, $\ROmegaOmega$, $\partial _{y|\Gamma }$ and $\partial _{x|\Gamma }$ for these new norms on $\Amacro$. Finally, we show in \cref{norme SF de mcl cas rationnel} that any $\mcla_{d,n} \in \Amacro(\Omega )$ has a finite norm and we provide a precise bound of its norm.
    \item In \cref{sec: convergence des mcl cas general}, we explain how to control pointwise evaluations of functions of $\mathcal{A}_{d} (\Omega )$ with their norm in $\mathcal{A}_{d}$. We deduce \cref{convergence des mcl} that proves the convergence of $\mcla_{d,n}$ and estimates its speed of convergence.
\end{itemize}

This outline was designed to handle the difficult case $\pst \in \mathbb{Q}$. To avoid repeating similar works for $\pst \not\in \mathbb{Q}$, we unify the two cases as soon as possible, by fitting the case $\pst \not\in \mathbb{Q}$ into the framework designed for $\pst \in \mathbb{Q}$. Since \cref{sec: convergence des mcl cas rationnel 1,sec: convergence des mcl cas rationnel 2} need the assumption $\pst \in \mathbb{Q}$, the unification happens at the beginning of \cref{sec: convergence des mcl cas general}.
\end{textesansboite}

\section{Convergence of the far field modes}
\label{sec: convergence des mcl}

\subsection{The rational case: estimates on individual terms}
\label{sec: convergence des mcl cas rationnel 1}

\begin{textesansboite}
In this section and the next one, we assume that $\pst \in \mathbb{Q}$. We introduce tools and give preliminary results to prove \cref{convergence des mcl} in this case. Note however that those tools are still well defined when $\pst \not\in \mathbb{Q}$. They will be used again in \cref{sec: convergence des mcl cas general} for the case $\pst \not\in \mathbb{Q}$.
\end{textesansboite}

\begin{definition}{norm on $\mathcal{A}_{d} (\Omega )$ and $\mathcal{A}_{d} (\Gamma )$}
\label{def: norme sur A_d}%
Let $d\in \mathbb{R}$. By~\cref{base de A Omega,rq: base de A couche et Gamma}, $\{\Im[(\alpha z)^{q} \,\overline{\alpha z}^{k} \log ^{j} \!(\alpha z)] \mid (q,k,j)\in \mathbb{R}\times \mathbb{N}^2  \text{ s.t.\ }  q+k=d \text{ and } [(q\in \mathbb{N} \text{ and } q\leqslant k) \Rightarrow  j>0] \}$ is a basis of $\mathcal{A}_{d} (\Omega )$ and $\{x^{d} \ln ^{j} x \mid j\in \mathbb{N}\}$ is a basis of $\mathcal{A}_{d} (\Gamma )$. We define $\|.\|$ on $\mathcal{A}_{d} (\Omega )$, resp.\ $\mathcal{A}_{d} (\Gamma )$, as the sum of the absolute values of the coordinates on this basis of $\mathcal{A}_{d} (\Omega )$, resp.\ $\mathcal{A}_{d} (\Gamma )$. 
\end{definition}

\begin{exemple}
We have $\|\phiOmega_{d} \|=1$ for~any $d\in \TZ*$, and $\big\|\Im[(\alpha z)^{q} \overline{\alpha z}^{k} P(\log (\alpha z))]\big\| = \sum_{k=0}^n |a_{k} |$ for~any $(q,k)\in \mathbb{R}\times \mathbb{N}$ and $P=\sum_{k=0}^n a_{k} T^{k} \in \mathbb{R}[T]$.
\end{exemple}

\begin{remarque}
\Cref{majoration des fonctions de A par leur norme} will show how to control pointwise evaluations of functions in $\mathcal{A}_{d} (\Omega )$ using the introduced norm.
\end{remarque}

\begin{textesansboite}
To estimate the terms of $\mcla_{d,n}$, it will be useful to control their powers of $\alpha z$, $\overline{\alpha z}$ and $\log (\alpha z)$. To this aim, we define the following spaces.
\end{textesansboite}

\begin{definition}{$\AOprecis{d,n,\delta }$ and $\AGprecis{d,n}$}
\label{def: Aprecis}%
For~any $d,\delta \in \mathbb{R}$ and $n\in \mathbb{N}$, we define:
\begin{enumerate}
    \item $\AOprecis{d,n,\delta }$ the vector subspace of $\mathcal{A}_{d} (\Omega )$ generated by the functions $z\mapsto \Im[(\alpha z)^{q} \,\overline{\alpha z}^{k} P(\log (\alpha z))]$ with $q\in \mathbb{R}$, $k\in \mathbb{N}$, $q+k=d$, $P\in \mathbb{R}[T]$, $q-k\geqslant \delta $ and $\deg P\leqslant 2n-1+ \mathbf{1}_{\TZ} (q-k) +\mathbf{1}_{\mathbb{N}} (q)$,
	\item and $\AGprecis{d,n} := \{x\mapsto  x^{d} \,Q(\ln x) \mid Q\in \mathbb{C}[T] \text{ and } \deg Q\leqslant 2n-1\}$.
\end{enumerate}\fbi
\end{definition}
\nomenclature[10-13-39-23]{$\AOprecis{d,n,\delta },\AGprecis{d,n}$}{In \cref{chap: decompo modale CL}, subspaces of $\mathcal{A}_{d} (\Omega )$ and $\mathcal{A}_{d} (\Gamma )$}

\begin{textesansboite}
Compared to the definitions of $\mathcal{A}_{d} (\Omega )$ and $\mathcal{A}_{d} (\Gamma )$, the inequalities on $q-k$, $\deg P$ and $\deg Q$ are new. These inequalities are necessary to obtain continuity constants for $\RGammaOmega$, $\ROmegaOmega$, $\partial _{y|\Gamma }$ and $\partial _{x|\Gamma }$. For example, the expression of $\partial _{y|\Gamma }$ involves $q-k$, which is unbounded in a given $\mathcal{A}_{d} (\Omega )$ (i.e.\ the set $\{q-k \mid {q\in \mathbb{R}},{k\in \mathbb{N}},{q+k=d} \}$ is unbounded) but it~is~not on $\AOprecis{d,n,\delta }$ because $\delta \leqslant q-k\leqslant q+k=d$. \Cref{rq: lien entre d q k delta} shows that the condition $q-k\geqslant \delta $ also provides control over $q$, $k$ and $d$, which is useful to bound the other operators. As~for the bounds on $\deg P$ and $\deg Q$, they allow us to treat all operations on polynomials of $\log (\alpha z)$ and $\ln x$ as linear maps in a finite dimensional vector space, so~that they have a continuity constant that depends on $n$.\\

The complicated condition on $\deg P$ in \cref{def: Aprecis}, instead of simply ``$\deg P\leqslant n$'', was chosen so~that it~is preserved when applying $\ROmegaOmega$. This is crucial because $\ROmegaOmega$ is iterated in the definition of $\mcla_{d,n}$. Thus, we are able to show that the terms of $\mcla_{d,n}$ belong to spaces $\AOprecis{d,n,\delta }$ that share a common value of $n$. This implies that the degrees in $\log (\alpha z)$ are bounded in a given $\mcla_{d,n}$, and that the operations that build the polynomials in $\log (\alpha z)$ of $\mcla_{d,n}$ are controlled by a common constant. Without this condition on $\deg P$, the continuity constants that bound the terms of $\mcla_{d,n}$ would explode at high order, which would not allow the proof of convergence of $\mcla_{d,n}$.
\end{textesansboite}

\begin{remarquenumerotee}
\label{rq: inclusions entre les Aprecis}%
Let $d\in \mathbb{R}$, $\delta \leqslant \delta '$ in $\mathbb{R}$, and $n\leqslant n'$ in $\mathbb{N}$. We have:
\begin{itemize}
    \item $\AOprecis{d,n,\delta } \subset  \AOprecis{d,n',\delta }$ and $\AGprecis{d,n} \subset  \AGprecis{d,n'}$
    \item and $\AOprecis{d,n,\delta } \supset \AOprecis{d,n,\delta '}$.
\end{itemize}
\end{remarquenumerotee}

\begin{remarquenumerotee}
\label{rq: lien entre d q k delta}%
The relations $q+k=d$ and $q-k\geqslant \delta $ in \cref{def: norme sur A_d} imply $2q = (q+k)+(q-k)\geqslant d+\delta $ and $2k = q+k-(q-k) \leqslant  d-\delta $ and $d\geqslant  2k+\delta  \geqslant \delta $ (because $k\in \mathbb{N}$). So the inequality $q-k\geqslant \delta $ is sufficient to have control over $q$, $k$ and $q\pm k$. In~particular, we deduce the implication \[\AOprecis{d,n,\delta }\neq \{0\} \Rightarrow  d\geqslant \delta .\]
One can show that this is in fact an equivalence.
\end{remarquenumerotee}

\begin{lemmesansboite}{}
\label{action des op R sur A precis}
Let $d,\delta \in \mathbb{R}$ and $n\in \mathbb{N}$. We have:
\begin{enumerate}
    \item $\RGammaOmega(\AGprecis{d,n}) \subset  \AOprecis{d,n,\delta }$ provided that $\delta \leqslant d$,
    \item $\ROmegaOmega(\AOprecis{d,n,\delta }) \subset  \AOprecis{d+2,n,\delta }$ for~any $\delta $,
    \item $\partial _{y|\Gamma } (\AOprecis{d,n,\delta }) \subset  \AGprecis{d-1,n+1}$ for~any $\delta $,
    \item $\partial _{x} (\AGprecis{d,n}) \subset  \AGprecis{d-1,n}$.
\end{enumerate}
\end{lemmesansboite}

\begin{demo}
Let $\varphi $ be in $\AOprecis{d,n,\delta }$, resp.\ $\AGprecis{d,n}$. It~is easy to see that the decomposition of $\varphi $ on the basis of $\mathcal{A}_{d} (\Omega )$, resp.\ $\mathcal{A}_{d} (\Gamma )$, mentioned in \cref{def: norme sur A_d} only involves elements of the basis that are also in $\AOprecis{d,n,\delta }$, resp.\ $\AGprecis{d,n}$. Therefore, it suffices to verify that $\RGammaOmega$, $\ROmegaOmega$, $\partial _{y|\Gamma }$ and $\partial _{x|\Gamma }$ send the elements of these basis into the right spaces.
\begin{enumerate}
    \item Let $j\in [\![0,2n]\!]$ and $\varphi :x\mapsto x^{d} \ln ^{j} x$. The proof of \cref{operateurs R dans Omega} shows that $\RGammaOmega(\varphi )=\Im[(\alpha z)^{d} P(\log (\alpha z))]$ with $P:=\RIma{d}(T^{j} )$. By~\cref{resolution de Im[alphad P(T-iTheta)]=Q}, we have $\deg P=j+1$ if $d\in \TZ$ and $\deg P=j$ otherwise. Therefore, it~is easy to check that $\RGammaOmega(\varphi )\in \AOprecis{d,n,\delta }$.
    \item Let $\varphi :=\Im[(\alpha z)^{q} \,\overline{\alpha z}^{k} \log ^{j} \!(\alpha z)]$ be an element of the basis of $\mathcal{A}_{d} (\Omega )$ s.t.\ $\varphi \in \AOprecis{d,n,\delta }$. The proof of \cref{operateurs R dans Omega} shows that $\ROmegaOmega(\varphi )=\psi -\RGammaOmega(\psi _{|\Gamma } )$ with $\psi :=\Im[(\alpha z)^{q+1} \,\overline{\alpha z}^{k+1} P(\log (\alpha z))]$ and $4(k+1)\big((q+1)P+P'\big)=T^{j}$. Using point 1, it suffices to show that $\psi \in \AOprecis{d+2,n,\delta }$. But this is easy by writing $q-k\geqslant \delta  \Rightarrow  (q+1)-(k+1)\geqslant \delta $ and:
    \begin{alignat*}{2}
        \deg P &\leqslant  j + \delta _{q,-1} && \text{ because }  \deg ((q+1)P+P')=j \\
        &\leqslant  2n-1+ \mathbf{1}_{\TZ} (q-k) +\mathbf{1}_{\mathbb{N}} (q) + \delta _{q,-1} && \text{ because }  \varphi \in \AOprecis{d,n,\delta }\\
        &\leqslant  2n-1+ \mathbf{1}_{\TZ} \big((q+1)-(k+1)\big) +\mathbf{1}_{\mathbb{N}} (q+1)
    \end{alignat*}
    \item Let $\varphi :=\Im[(\alpha z)^{q} \,\overline{\alpha z}^{k} \log ^{j} \!(\alpha z)]$ be an element of the basis of $\mathcal{A}_{d} (\Omega )$ s.t.\ $\varphi \in \AOprecis{d,n,\delta }$. The proof of \cref{operateurs de derivation dans A} shows that $\partial _{y} \varphi _{|\Gamma } = x^{q+k-1} \Re[\alpha ^{q-k} P(\ln x-\mathrm{i}\Theta )]$ with $P:=((q-k)T^{j} +jT^{j-1} )$. Since $\deg P\leqslant j\leqslant 2n+1$, we have $\partial _{y} \varphi _{|\Gamma } \in \AGprecis{d-1,n+1}$.
    \item Let $j\in [\![0,2n-1]\!]$ and $\varphi :x\mapsto x^{d} \ln ^{j} x$. We have $\partial _{x} \varphi  = x^{d-1} (d\ln ^{j} x+ j\ln ^{j-1} x) \in \AGprecis{d,n}$.
\end{enumerate}
\end{demo}

\begin{textesansboite}
Let $(a,b)\in (\mathbb{N}^*)^2 $ be s.t.\ $\Theta = \frac{a}{b} \pi $ and $\pgcd{a,b}=1$. Note that $\mathbb{Z}+\TZ=\ZTZ$. 
\end{textesansboite}

\begin{lemme}{}
\label{norme naive des op R}%
Let $d\in \ZTZ$, $\delta \in \mathbb{R}$ and $n\in \mathbb{N}$. We have:
\begin{enumerate}
    \item $\forall \varphi \in \AGprecis{d,n}, \ \ \|\RGammaOmega(\varphi )\| \leqslant  C_{n} \,\|\varphi \|$,
    \item $\forall \varphi \in \AOprecis{d,n,\delta }, \ \ \|\ROmegaOmega(\varphi )\| \leqslant  \left\{\begin{array}{@{\;}ll}
        \frac{C_{n} }{d+1} \|\varphi \| & \text{if } d\geqslant |\delta | \\
        C_{\delta ,n} \,\|\varphi \| & \text{if } d<|\delta |
    \end{array}\right.$,
    \item $\forall \varphi \in \AOprecis{d,n,\delta }, \ \ \|\partial _{y} \varphi _{|\Gamma } \| \leqslant  C_{n} \,(1+\max (|\delta |,d)) \,\|\varphi \|$,
    \item $\forall \varphi \in \AGprecis{d,n}, \ \ \|\partial _{x} \varphi \| \leqslant  C_{n} \,(1+\max (|\delta |,d)) \,\|\varphi \|$ provided that $d\geqslant \delta $.
\end{enumerate}\fbi
\end{lemme}

\begin{remarquenumerotee}
\label{rq: norme naive des op R}
\Cref{norme naive des op R} crucially uses the assumptions $\pst \in \mathbb{Q}$ and $d\in \ZTZ$ in the first estimate. For instance, if $\varphi  \in \AGprecis{d,n}$ is the function $x^{d}$ and $d\not\in \TZ$, then $\RGammaOmega(\varphi ) = \pm  \sin (d\Theta )^{-1} \Im[(\alpha z)^{d} ]$. For a general $d$, $\sin (d\Theta )^{-1}$ may be arbitrarily large, which happens when $d$ is very close to $\TZ$. However, when $\Theta = \frac{a}{b} \pi $ and $d\in \ZTZ$, $\sin (d\Theta )^{-1}$ can take only a finite number of values, which let us bound it by a constant. 
If $d\in  \frac{1}{a} \mathbb{Z}$, it~is also possible to have $d\in \TZ$ and $\sin (d\Theta )=0$, but in this case $\RGammaOmega(\varphi )$ has a different expression involving $\log (\alpha z)$ and its coefficients can still be controlled. 
\end{remarquenumerotee}

\begin{demo}
Similarly to the proof of \cref{action des op R sur A precis}, it suffices to show these inequalities when $\varphi $ is an element of the basis mentioned in \cref{def: norme sur A_d}. For~any $P=\sum_{j=0}^m a_{j} T^{j} \in \mathbb{C}[T]$, we denote $\|P\|:=\sum_{j=0}^m |a_{j} |$.
\begin{enumerate}
    \item Let $j\in [\![0,2n-1]\!]$ and $\varphi :x\mapsto x^{d} \ln ^{j} x$. We have $\RGammaOmega(\varphi )=\Im[(\alpha z)^{d} P(\log (\alpha z))]$ with $P(T):=\RIma{d}(T^{j} )$. In addition, $\RIma{d}$ (which maps any $Q\in \mathbb{R}[T]$ to a solution of $\Im[\alpha ^{d} P(T-\mathrm{i}\Theta )]=Q$, see \cref{resolution de Im[alphad P(T-iTheta)]=Q}) only depends on $\alpha ^{d} = e^{-\mathrm{i}\Theta d}$, so it only depends on $d \,\mathrm{mod}\, \pst=d \,\mathrm{mod}\, \frac{b}{a}$, which can take only a finite number of values since $d\in \ZTZ$. Thus
    \[\sup _{d\in \ZTZ} \sup _{j\in [\![0,2n-1]\!]} \|\RIma{d}(T^{j} )\| = C_{n} < +\infty \]
    Therefore $\|\RGammaOmega(\varphi )\| = \|P\|\leqslant C_{n} = C_{n} \|\varphi \|$.
    \item Let $\varphi :=\Im[(\alpha z)^{q} \,\overline{\alpha z}^{k} \log ^{j} \!(\alpha z)]$ be an element of the basis of $\mathcal{A}_{d} (\Omega )$ s.t.\ $\varphi \in \AOprecis{d,n,\delta }$. The proof of \cref{operateurs R dans Omega} shows that $\ROmegaOmega(\varphi )=\psi -\RGammaOmega(\psi _{|\Gamma } )$ with $\psi :=\Im[(\alpha z)^{q+1} \,\overline{\alpha z}^{k+1} P(\log (\alpha z))]$ and $4(k+1)\big((q+1)P+P'\big)=T^{j}$. Using point 1, it suffices to show that
    \[\|\psi \|\leqslant  \left\{\begin{array}{@{\;}ll}
        \frac{C_{n} }{d+1}  & \text{if } d>|\delta |\\[1ex]
        C_{\delta ,n}  & \text{if } d<|\delta |
    \end{array}\right.\]
    First let us asume that $d\geqslant |\delta |$. \Cref{rq: lien entre d q k delta} implies $2q\geqslant d+\delta \geqslant 0$. Let $D$ be the derivation operator of $\mathbb{R}_{2n+1} [T] :=\{Q\in \mathbb{R}[T] \mid \deg Q\leqslant 2n+1\}$ and $C^{D} _{n}$ its continuity constant for $\|.\|$. Using $j\leqslant 2n+1$, $D^{2n+3} =0$ and $q\geqslant 0$, we get:
    \[P = \frac{1}{4(k+1)(q+1)} \left( \mathrm{id}+ \frac{1}{q+1} D\right)^{-1} (T^{j} ) = \frac{1}{4(k+1)} \sum _{m=0} ^{2n+1} \frac{(-1)^{m} }{(q+1)^{m+1} } D^{m} (T^{j} )\]
    so 
    \[\|P\| \leqslant  \frac{1}{4(k+1)} \sum _{m=0} ^{2n+1} \frac{\|D^{m} (T^{j} )\|}{(q+1)^{m+1} } \leqslant  \frac{1}{4(k+1)} \sum _{m=0} ^{2n+1} \frac{(2n+1) !}{q+1} \leqslant  \frac{C_{n} }{(q+1)(k+1)}.\]
    Moreover, $(q+1)(k+1)=qk+q+k+1\geqslant d+1$, as $d=q+k$ and $q,k\geqslant 0$. Therefore $\|\psi \| = \|P\|\leqslant  C_{n} /(d+1)$.\\
    Now, let us assume $d<|\delta |$. By~\cref{rq: lien entre d q k delta}, we also have $d\geqslant \delta $, so $\delta <0$ and $d\in [\delta ,-\delta ]\cap \ZTZ$. In addition, $k\leqslant  \frac{d-\delta }{2} \leqslant |\delta |$ and $j\leqslant 2n+1$. Thus $q=d-k$ and $k$ and $j$ can take only a finite number of values that depends on $(\delta ,n)$. Hence $\|\psi \|\leqslant C_{\delta ,n}$.
    \item Let $\varphi :=\Im[(\alpha z)^{q} \,\overline{\alpha z}^{k} \log ^{j} \!(\alpha z)]$ be an element of the basis of $\mathcal{A}_{d} (\Omega )$ s.t.\ $\varphi \in \AOprecis{d,n,\delta }$. We know that $\partial _{y} \varphi _{|\Gamma } = x^{q+k-1} \Re[\alpha ^{q-k} P(\ln x-\mathrm{i}\Theta )]$ with $P(T):=((q-k)T^{j} +jT^{j-1} )$. Hence $\|\partial _{y} \varphi _{|\Gamma } \|\leqslant \|P(\ln x-\mathrm{i}\Theta )\| \leqslant  C_{n} \|P\|\leqslant C_{n} \,(|q-k|+1)$ because $j\leqslant 2n+1$. Finally, we have $q-k\geqslant \delta $ by \cref{def: norme sur A_d} and $q-k\leqslant q+k=d$, so $|q-k|\leqslant \max (-\delta ,d)\leqslant \max (|\delta |,d)$.
    \item Let $j\in [\![0,2n-1]\!]$ and $\varphi :x\mapsto x^{d} \ln ^{j} x$. We have $\|\partial _{x} \varphi \| = \|x^{d-1} (d\ln ^{j} x+ j\ln ^{j-1} x)\|\leqslant C_{n} (1+|d|)\leqslant C_{n} (1+\max (|\delta |,d))$ as $d\geqslant \delta $.
\end{enumerate}
\end{demo}

\begin{textesansboite}
\Cref{norme naive des op R} gives relatively direct continuity estimates of $\RGammaOmega$, $\ROmegaOmega$, $\partial _{y|\Gamma }$ and $\partial _{x|\Gamma }$ for a natural norm on $\mathcal{A}_{d}$. Now let us reformulate this using a more suitable norm in order to deduce easily estimates on formal series.\\

Among $\RGammaOmega$, $\ROmegaOmega$, $\partial _{y|\Gamma }$ and $\partial _{x|\Gamma }$, the most critical operator is $\ROmegaOmega$ because it~is iterated in the definition of $\mcla_{d,n}$. Fortunately, \cref{norme naive des op R} shows that it acts like a multiplication by $\frac{C_{n} }{d+1}$, while increasing $d$ by 2 (it sends $\mathcal{A}_{d} (\Omega )$ to $\mathcal{A}_{d+2} (\Omega )$). Intuitively this should make the series $\mcla_{d,n}$ converge at a factorial speed, similarly to $\sum_{n=1}^\infty  \frac{x^{n} }{n!}$. To make this idea rigourous, let us introduce the following function
\nomenclature[38]{$\gamma $}{In \cref{chap: decompo modale CL}, $x\mapsto  \sqrt{\Gamma (\max (x,1))}$ where $\Gamma $ is the Gamma function}
\begin{equation}{\label{eq: def fct petit gamma}}
    \gamma (d):= \left\{\begin{array}{ll}
        \sqrt{\Gamma (d)} & \text{if } d\geqslant 1\\[1ex]
        1 & \text{if } d\leqslant 1
    \end{array}\right.
\end{equation}
\nomenclature[60-38]{$\|.\|_{\gamma }$}{In \cref{chap: decompo modale CL}, norm on $\mathcal{A}(\Omega )$ normalized using the function $\gamma $}%
where $\Gamma $ is the usual Gamma function and, for~any $D\in \{\Omega ,\Gamma \}$ and $\varphi \in \mathcal{A}_{d} (D)$, $\normg{\varphi } :=\gamma (d) \|\varphi \|$. We~use the Gamma function rather than factorials, as non integer degrees appear in our problem. But unlike the Gamma function, $\gamma $ has no pole (see \cref{fig: fonction petit gamma}). Moreover, the square root \cref{eq: def fct petit gamma} is related to the fact that $\deg \ROmegaOmega =2$. \Cref{norme-gamma des op R} gives estimates for the norm $\normg{.}$ that have a simpler form than \cref{norme naive des op R}.
\end{textesansboite}

\pgfplotsset{
    every axis/.append style={
        thick,
    },
    cycle list/Dark2,
    cycle multiindex* list={
        mark list*\nextlist
        Dark2\nextlist
    },
    /pgfplots/table/x=ordrecvg,
    /pgfplots/table/y=nbCCC,
    legend cell align=left,
}

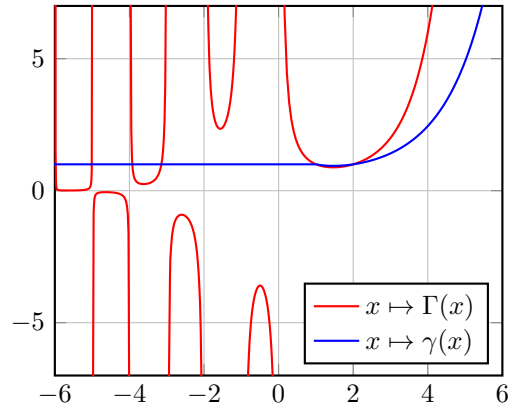
\begin{figure}[h]
\newcommand{\approxFctGamma}[1]{sqrt(2*pi)*(#1)^(#1-0.5)*exp(-(#1))*exp(1/(12*(#1)))} 
\begin{tikzpicture}
\begin{axis}[
    width=75mm,
    enlargelimits=0,
    xmajorgrids=true,
    ymajorgrids=true,
    legend pos=south east,
    legend style={legend columns=1},
    xmin=-6,
    xmax=6,
    ymin=-7,
    ymax=7,
]

\addplot[domain=0.04:6, color=red, samples=100] {\approxFctGamma{x}};
\addlegendentry{$x\mapsto\Gamma(x)$}
\addplot[domain=1:6, color=blue, smooth] {sqrt(\approxFctGamma{x})};
\addlegendentry{$x\mapsto\gamma(x)$}

\addplot[domain=-0.95:-0.05, color=red, smooth] {\approxFctGamma{x+1}/x};
\addplot[domain=-1.95:-1.05, color=red, smooth] {\approxFctGamma{x+2}/x/(x+1)};
\addplot[domain=-2.95:-2.05, color=red, smooth] {\approxFctGamma{x+3}/x/(x+1)/(x+2)};
\addplot[domain=-3.97:-3.02, color=red, samples=50] {\approxFctGamma{x+4}/x/(x+1)/(x+2)/(x+3)};
\addplot[domain=-4.98:-4.005, color=red, samples=50] {\approxFctGamma{x+5}/x/(x+1)/(x+2)/(x+3)/(x+4)};
\addplot[domain=-5.99:-5.001, color=red, samples=100] {\approxFctGamma{x+6}/x/(x+1)/(x+2)/(x+3)/(x+4)/(x+5)};

\addplot[domain=-10:1, color=blue] {1};

\end{axis}
\end{tikzpicture}
\caption{The usual Gamma function and the function $\gamma$}
\label{fig: fonction petit gamma}
\end{figure}

\begin{corollairesansboite}{}
\label{norme-gamma des op R}%
Let $d\in \ZTZ$, $\delta \in \mathbb{R}$ and $n\in \mathbb{N}$. We have:
\setalignbox{$\forall \varphi \in \AOprecis{d,n,\delta },$}\disablealignbox
\begin{enumerate}
    \item $\alignbox{\forall \varphi \in \AGprecis{d,n},} \ \ \normg{\RGammaOmega(\varphi )} \leqslant  C_{n} \normg{\varphi }$,
    \item $\alignbox{\forall \varphi \in \AOprecis{d,n,\delta },} \ \ \normg{\ROmegaOmega(\varphi )} \leqslant  C_{\delta ,n} \normg{\varphi }$,
    \setalignbox{$\forall \eAMP>1,\; \forall \varphi \in \AOprecis{d,n,\delta },$}\disablealignbox
    \item $\alignbox{\forall \eAMP>1,\; \forall \varphi \in \AOprecis{d,n,\delta },} \ \ \normg{\partial _{y} \varphi _{|\Gamma } } \leqslant  C_{\delta ,n,\eAMP} \,\eAMP^{d} \,\normg{\varphi }$,
    \item $\alignbox{\forall \eAMP>1,\; \forall \varphi \in \AGprecis{d,n},} \ \ \normg{\partial _{x} \varphi } \leqslant  C_{\delta ,n,\eAMP} \,\eAMP^{d} \,\normg{\varphi }$ provided that $d\geqslant \delta $.
\end{enumerate}
\end{corollairesansboite}

\begin{remarque}
$s$ is related to the parameter of the norm defined on $\Amacro(\Omega )$ in \cref{def: norme sur Amacro}. Its role will become clearer in \cref{sec: convergence des mcl cas rationnel 2}.
\end{remarque}

\begin{demo}
The proof uses the results of \cref{action des op R sur A precis,norme naive des op R}.
\begin{enumerate}
    \item It suffices to multiply $\|\RGammaOmega(\varphi )\| \leqslant  C_{n} \,\|\varphi \|$ by $\gamma (d)$ as $\RGammaOmega(\varphi )\in \AOprecis{d,n,\delta }$.
    \item If $d\geqslant \max (|\delta |,1)$, then $\normg{\ROmegaOmega(\varphi )}=\sqrt{\Gamma (d+2)}\,\|\ROmegaOmega(\varphi )\|\leqslant  \sqrt{\Gamma (d+2)} \frac{C_{n} }{d+1} \|\varphi \| = \sqrt{\frac{\Gamma (d)d}{d+1}} C_{n} \|\varphi \|$ $\leqslant $ $\sqrt{\Gamma (d)} \,C_{n} \|\varphi \|=C_{n} \normg{\varphi }$.\\
    If $|\delta |\leqslant d<\max (|\delta |,1)$, then $0\leqslant d<1$, so $\normg{\ROmegaOmega(\varphi )}\leqslant  C \|\ROmegaOmega(\varphi )\|\leqslant  C_{n} \|\varphi \| \leqslant  C_{n} \normg{\varphi }$.\\
    Otherwise, $d<|\delta |$ so $\normg{\ROmegaOmega(\varphi )}\leqslant  C_{\delta } \|\ROmegaOmega(\varphi )\|\leqslant  C_{n,\delta } \|\varphi \| \leqslant  C_{n,\delta } \normg{\varphi }$.
    \item We have $\normg{\partial _{y} \varphi _{|\Gamma } } \leqslant  C_{n} \,f_{\delta } (d) \normg{\varphi }$ with $f_{\delta } (d):=(1+\max (d,-\delta )) \frac{\gamma (d-1)}{\gamma (d)}$. $f_{\delta }$ satisfies $f_{\delta } (d)=O(\sqrt{d} )=O(s^{d} )$ when $d\rightarrow \infty $. Moreover, \cref{rq: lien entre d q k delta} states that $\AOprecis{d,n,\delta } \neq \{0\}\Rightarrow d\geqslant \delta $, so we can assume that $d\geqslant \delta $. Hence we have $\normg{\partial _{y} \varphi _{|\Gamma } } \leqslant  C_{\delta ,n,s} s^{d} \normg{\varphi }$.
    \item Similar to 3.
\end{enumerate}
\end{demo}

\begin{remarque}
All constants in \cref{norme-gamma des op R} may depend on $(\delta ,n,\eAMP)$, but they are independent of $d$, which will be useful to prove \cref{decompo modale des CL,decompo en modes tronques des CL}.
\end{remarque}

\subsection{The rational case: estimates on formal series}
\label{sec: convergence des mcl cas rationnel 2}

\begin{textesansboite}
Now that we have tools to estimate elements of the spaces $\mathcal{A}_{d} (\Omega )$ and $\mathcal{A}_{d} (\Gamma )$, let us build tools to estimate formal series. We still assume $\pst \in \mathbb{Q}$. 
\end{textesansboite}

\begin{definition}{norms $\|.\|_{\eAMP}$ on formal series}
\label{def: norme sur Amacro}%
Let $D\in \{\Omega ,\Gamma \}$ and $\eAMP\in \mathbb{R}_+^*$. For~any $\Phi =(\varphi _{d} )_{d\in \mathbb{R}} \in \Amacro(D)$, we define:
\[\|\Phi \|_{\eAMP} := \sum _{d\in \mathbb{R}} \eAMP^{-d} \normg{\varphi _{d} } \ \in \ [0,+\infty ].\]
Note that, by~definition of $\Amacro(D)$, the above sum has a countable support.
\end{definition}
\nomenclature[60-28]{$\|.\|_{\eAMP}$}{In \cref{chap: decompo modale CL}, norm on $\Amacro(D)$}

\begin{remarque}
If a series $\Phi =(\varphi _{d} )_{d} \in \Amacro(D)$ satisfies $\|\Phi \|_{\eAMP} <\infty $, then $\|\varphi _{d} \|=O\!\left(\displaystyle \frac{\eAMP^{d} }{\gamma (d)} \right) =o(1)$ as $d\rightarrow \infty $. The norm $\|.\|_{\eAMP}$ controls the speed at which $\|\varphi _{d} \|$ tends to 0.
\end{remarque}

\begin{definition}{}
\label{def: AMprecis}%
For~any $\delta \in \mathbb{R}$ and $n\in \mathbb{N}$, we denote
\[\AMprecis{n,\delta }(\Omega ) := \Amacro(\Omega ) \cap  \fsum_{d\in (\mathbb{Z}+\TZ)\cap [\delta ,\infty )} \AOprecis{d,n,\delta } \quad \;\  \text{ and } \quad \;\ \AMprecis{n,\delta }(\Gamma ) := \Amacro(\Gamma ) \cap  \fsum_{d\in (\mathbb{Z}+\TZ)\cap [\delta ,\infty )} \AGprecis{d,n}.\]\fbeq
\end{definition}

\begin{remarque}
We recall that here we have $\mathbb{Z}+\TZ=\ZTZ$ because we assumed $\Theta = \frac{a}{b} \pi $ with $a,b\in \mathbb{N}^*$ and $\pgcd{a,b}=1$. But we write \cref{def: AMprecis} using $\mathbb{Z}+\TZ$ as we will use it again in the case $\pst \not\in  \mathbb{Q}$ in \cref{sec: convergence des mcl cas general}.
\end{remarque}

\begin{remarquenumerotee}
\label{rq: inclusions entre les AMprecis}%
Let $\delta \leqslant \delta '$ in $\mathbb{R}$, $n\leqslant n'$ in $\mathbb{N}$ and $D\in \{\Omega ,\Gamma \}$. \Cref{rq: inclusions entre les Aprecis} implies $\AMprecis{n,\delta }(D) \subset  \AMprecis{n',\delta }(D)$ and $\AMprecis{n,\delta }(D) \supset \AMprecis{n,\delta '}(D)$.
\end{remarquenumerotee}

\begin{textesansboite}
\Cref{norme SF de mcl cas rationnel} will show that $\mcla_{d,n} \in \AMprecis{n,d-n}(\Omega )$ for~any $d\in \TZ*$ and $n\in \mathbb{N}$. The first step to prove this is the following lemma.
\end{textesansboite}

\begin{lemmesansboite}{}
\label{action des op R sur AMprecis}
Let $\delta \in \mathbb{R}$ and $n\in \mathbb{N}$. We have:
\begin{enumerate}
    \item $\RGammaOmega(\AMprecis{n,\delta }(\Gamma )) \subset  \AMprecis{n,\delta }(\Omega )$,
    \item $\ROmegaOmega(\AMprecis{n,\delta }(\Omega )) \subset  \AMprecis{n,\delta +2}(\Omega )$ and 
    $\langle -k_{0} ^2 \ROmegaOmega\rangle (\AMprecis{n,\delta }(\Omega )) \subset  \AMprecis{n,\delta }(\Omega )$,
    \item $\partial _{y|\Gamma } (\AMprecis{n,\delta }(\Omega )) \subset  \AMprecis{n+1,\delta -1}(\Gamma )$,
    \item $\partial _{x} (\AMprecis{n,\delta }(\Gamma )) \subset  \AMprecis{n,\delta -1}(\Gamma )$.
\end{enumerate}
\end{lemmesansboite}

\begin{demo}
All the inclusions are straightforward consequences of \cref{action des op R sur A precis}, except $\langle -k_{0} ^2 \ROmegaOmega\rangle (\AMprecis{n,\delta }(\Omega )) \subset  \AMprecis{n,\delta }(\Omega )$. Let us prove it. Let $\Phi  \in \AMprecis{n,\delta }(\Omega )$. By~$\ROmegaOmega(\AMprecis{n,\delta }(\Omega )) \subset  \AMprecis{n,\delta +2}(\Omega )$ and \cref{rq: inclusions entre les AMprecis}, $\AMprecis{n,\delta }(\Omega )$ is stable under $\ROmegaOmega$. Thus, $(-k_{0} ^2 \ROmegaOmega)^{m} \Phi  \in  \AMprecis{n,\delta }(\Omega )$ for~any $m\in \mathbb{N}$. Now, $\AMprecis{n,\delta }(\Omega )$ is closed under $\fsum$ according to \cref{fsum de sev des Ed}, so $\langle -k_{0} ^2 \ROmegaOmega\rangle (\Phi ) = \fsum_{m=0} ^{\infty } (-k_{0} ^2 \ROmegaOmega)^{m} (\Phi ) \in \AMprecis{n,\delta }(\Omega )$.
\end{demo}

\begin{lemme}{}
\label{norme des op R sur series formelles}%
Let $\delta \in \mathbb{R}$, $n\in \mathbb{N}$, $\eAMP>0$ and $\eAMP'>\eAMP$. We have:
\begin{enumerate}
    \item $\forall \Phi \in \AMprecis{n,\delta }(\Gamma ), \ \ \|\RGammaOmega(\Phi )\|_{\eAMP} \leqslant  C_{n} \|\Phi \|_{\eAMP}$,
    \item $\forall \Phi \in \AMprecis{n,\delta }(\Omega ), \ \ \|\ROmegaOmega(\Phi )\|_{\eAMP} \leqslant  C_{\delta ,n} \,\eAMP^{-2} \, \|\Phi \|_{\eAMP}$,
    \item $\forall \Phi \in \AMprecis{n,\delta }(\Omega ), \ \ \|\partial _{y} \Phi _{|\Gamma } \|_{\eAMP'} \leqslant  C_{\delta ,n,\eAMP,\eAMP'} \,\|\Phi \|_{\eAMP}$,
    \item $\forall \Phi \in \AMprecis{n,\delta }(\Gamma ), \ \ \|\partial _{x} \Phi \|_{\eAMP'} \leqslant  C_{\delta ,n,\eAMP,\eAMP'} \,\|\Phi \|_{\eAMP}$.
\end{enumerate}\fbi
\end{lemme}

\begin{demo}
The proof uses the results of \cref{action des op R sur A precis,norme-gamma des op R}.
\begin{enumerate}
    \item Writing $\Phi $ as $\fsum_{d} \varphi _{d}$, it suffices to multiply $\normg{\RGammaOmega(\varphi _{d} )} \leqslant  C_{n} \normg{\varphi _{d} }$ by $\eAMP^{-d}$ and to sum over $d$.
    \item Let $(\varphi _{d} )$ be s.t.\ $\varphi _{d} \in \AOprecis{d,n,\delta }$ for~any $d$ and $\Phi =\fsum_{d} \varphi _{d}$. Then
    \[\|\ROmegaOmega(\Phi )\|_{\eAMP} =\Big\| \fsum_{d} \underbrace{\ROmegaOmega(\varphi _{d-2} )}_{\in \,\mathcal{A}_{d} (\Omega )} \Big\|_{\eAMP} = \sum _{d} \eAMP^{-d} \,\normg{\ROmegaOmega(\varphi _{d-2} )} \leqslant  C_{\delta ,n} \sum _{d} \eAMP^{-d} \normg{\varphi _{d-2} } = C_{\delta ,n} \,\eAMP^{-2} \,\|\Phi \|_{\eAMP} .\]
    \item Let $(\varphi _{d} )$ be s.t.\ $\varphi _{d} \in \AGprecis{d,n}$ for~any $d$ and $\Phi =\fsum_{d} \varphi _{d}$. Then
    \begin{multline*}
        \|\partial _{y} \Phi _{|\Gamma } \|_{\eAMP'} = \Big\|\fsum_{d} \underbrace{\partial _{y} \varphi _{d+1|\Gamma } }_{\in \,\mathcal{A}_{d} (\Gamma )} \|_{\eAMP'} = \sum _{d} \eAMP'^{-d} \,\normg{\partial _{y} \varphi _{d+1|\Gamma } } = \eAMP' \sum _{d} \eAMP'^{-d} \,\normg{\partial _{y} \varphi _{d|\Gamma } } \\[-3ex]
        \leqslant  \eAMP'\, C_{\delta ,n, \eAMP'/\eAMP } \sum _{d} \eAMP'^{-d} \, \Big( \frac{\eAMP'}{\eAMP} \Big)^{d} \normg{\varphi _{d} } = C_{\delta ,n,\eAMP,\eAMP'} \,\|\Phi \|_{\eAMP} .
    \end{multline*}\vspace{-3.5ex}
    \item Very similar to 3.
\end{enumerate}
\end{demo}

\begin{textesansboite}
As the definition of $\mcla_{d,n}$ involves $\langle -k_{0} ^2 \ROmegaOmega\rangle $ (see \cref{def: mcl}), we also need the following result about the continuity of $\langle -k_{0} ^2 \ROmegaOmega\rangle $.
\end{textesansboite}

\begin{proposition}{}
\label{norme de <Rlap>}
Let $\delta \in \mathbb{R}$ and $n\in \mathbb{N}$. There exists $\eAMPmin{\delta ,n}\in \mathbb{R}$ s.t.\ for~any $\eAMP \geqslant \eAMPmin{\delta ,n}$ and $\Phi \in \AMprecis{n,\delta }(\Omega )$, we have $\big\|\langle -k_{0} ^2 \ROmegaOmega\rangle (\Phi )\big\|_{\eAMP} \leqslant  C_{\delta ,n} \, \|\Phi \|_{\eAMP}$.
\end{proposition}

\begin{demobis}
\underline{Step 1:} Let us start by showing that, for~any $\eAMP>0$ and any sequence $(\Phi _{m} )_{m\in \mathbb{N}}$ of elements of $\Amacro(\Omega )$ that has collectively locally finite support (see \cref{def: famille avec CLFS}), we have
\[\bigg\|\fsum_{m=0} ^{\infty } \Phi _{m} \bigg\|_{\eAMP} \leqslant  \sum_{m=0}^\infty  \|\Phi _{m} \|_{\eAMP} \qquad \text{ in }[0,+\infty ].\]
For~any $m\in \mathbb{N}$, let us denote $(\varphi _{m,d} )_{d\in \mathbb{R}} :=\Phi _{m}$ the family of the terms of the formal series $\Phi _{m}$. We have:
\begin{alignjustif}
    \bigg\|\fsum_{m=0} ^{\infty } \Phi _{m} \bigg\|_{\eAMP} &= \bigg\|\bigg(\sum_{m=0}^\infty  \varphi _{m,d} \bigg)_{\!d\in \mathbb{R}} \bigg\|_{\eAMP} & \text{by \cref{def: fsum}} \\
    &= \sum _{d\in \mathbb{R}} \eAMP^{-d} \,\bigg\|\sum_{m=0}^\infty  \varphi _{m,d} \bigg\|_{\gamma } & \text{by~definition of }  \|.\|_{\eAMP} \\
    &= \sum _{d\in \mathbb{R}} \eAMP^{-d} \sum_{m=0}^\infty  \|\varphi _{m,d} \|_{\gamma } & \begin{array}{@{}l@{}}
        \text{because, for~any $d$, the sum over $m$} \\
        \text{has a finite number of non-zero terms} 
    \end{array} \\
    &= \sum_{m=0}^\infty  \sum _{d\in \mathbb{R}} \eAMP^{-d} \|\varphi _{m,d} \|_{\gamma }  \\
    &= \sum_{m=0}^\infty  \|\Phi _{m} \|_{\eAMP} & \text{by~definition of }  \|.\|_{\eAMP}
\end{alignjustif}

\underline{Step 2:} Let $\Phi \in \AMprecis{n,\delta }(\Omega )$. \Cref{action des op R sur AMprecis,rq: inclusions entre les AMprecis} show that $\AMprecis{n,\delta }(\Omega )$ is stable under $\ROmegaOmega$, thus: $\forall m\in \mathbb{N},\; (-k_{0} ^2 \ROmegaOmega)^{m} (\Phi ) \in  \AMprecis{n,\delta }(\Omega )$. Moreover, according to \cref{norme-gamma des op R,norme des op R sur series formelles}, there exists $\Clap{\delta ,n}>0$ s.t.: $\forall \Psi \in \AMprecis{n,\delta }(\Omega ), \; \|\ROmegaOmega(\Psi )\|_{\eAMP} \leqslant  \Clap{\delta ,n} \,\eAMP^{-2} \, \|\Psi \|_{\eAMP}$. Let us fix the value of this constant : unlike $C_{\delta ,n}$, $\Clap{\delta ,n}$ will not change from line to line. We deduce that, for~any $m\in \mathbb{N}$, we have $\big\|(-k_{0} ^2 \ROmegaOmega)^{m} (\Phi )\big\|_{\eAMP} \leqslant \big|k_{0} ^2 \,\Clap{\delta ,n} \,\eAMP^{-2} \big|^{m} \|\Phi \|_{\eAMP}$. Therefore, we can write:\\
\begin{alignjustif}
    \bigg\| \fsum_{m=0} ^{\infty } (-k_{0} ^2 \ROmegaOmega)^{m} (\Phi )\bigg\|_{\eAMP}
    &\leqslant  \sum_{m=0}^\infty  \big\|(-k_{0} ^2 \ROmegaOmega)^{m} (\Phi )\big\|_{\eAMP} & \text{by step 1} \\
    &\leqslant  \sum_{m=0}^\infty  \big|k_{0} ^2 \,\Clap{\delta ,n} \,\eAMP^{-2} \big|^{m} \|\Phi \|_{\eAMP}  &\\
    &= \frac{1}{1-|k_{0} ^2 | \, \Clap{\delta ,n} \,\eAMP^{-2} } \, \|\Phi \|_{\eAMP} & \text{provided that }  \eAMP>|k_{0} | \sqrt{\Clap{\delta ,n}} \\
    &\leqslant  C_{\delta ,n} \,\|\Phi \|_{\eAMP} & \text{provided that }  \eAMP\geqslant \eAMPmin{\delta ,n}:=|k_{0} | \sqrt{\Clap{\delta ,n}} +1.
\end{alignjustif}
\end{demobis}

\begin{textesansboite}
By~considering $\max \{\eAMPmin{\delta ,m} \mid m\in [\![0,n]\!]\}$ instead of $\eAMPmin{\delta ,n}$, it~is easy to see that $(\eAMPmin{\delta ,n})$ can be chosen as an~increasing sequence w.r.t.\ $n$ for~any $\delta $. This will be important for the proof by induction of the following proposition.
\end{textesansboite}

\begin{proposition}{}
\label{norme SF de mcl cas rationnel}
Let $\delta \in \mathbb{R}$, $d\in \TZ*\cap [\delta ,\infty )$, $n\in \mathbb{N}$ and $\eAMP'>\eAMP>\eAMPmin{\delta ,n}$.\\
We have $\mcla_{d,n} \in \AMprecis{n,d-n}(\Omega )$ and $\|\mcla_{d,n} \|_{\eAMP'} \leqslant  C_{\delta ,n,\eAMP,\eAMP'} \,\eAMP^{-d} \,\gamma (d)$.
\end{proposition}

\begin{remarque}
The explicit behavior w.r.t.\ $d$ of the estimate in \cref{norme SF de mcl cas rationnel} will be useful to prove \cref{decompo modale des CL,decompo en modes tronques des CL}.
\end{remarque}

\begin{demo}
Let us prove it by induction on $n$. First, let $n=0$. We have $\|\phi_{d} \|=1$ so, seeing $\phi_{d}$ as an element of $\Amacro(\Omega )$, we get $\|\phi_{d} \|_{\eAMP} = \eAMP^{-d} \,\gamma (d)$. In addition, it~is easy to check that $\phi_{d} \in \AMprecis{0,d}(\Omega )$, which yields $\mcla_{d,n} :=\langle -k_{0} ^2 \ROmegaOmega\rangle (\phi_{d} ) \in \AMprecis{n,d}(\Omega )$ by \cref{action des op R sur AMprecis}. Then, since $\AMprecis{0,d}(\Omega ) \subset  \AMprecis{0,\delta }(\Omega )$ and $\eAMP>\eAMPmin{\delta ,n}$, \cref{norme de <Rlap>} shows that we have $\|\mcla_{d,0} \|_{\eAMP} \leqslant  C_{\delta } \;\eAMP^{-d} \;\gamma (d)$. Besides, note that 
\begin{equation}{\label{eq: comparaison normes sur series formelles}}
    \forall \delta '\in \mathbb{R},\ \forall n'\in \mathbb{N},\ \forall \eAMP_{1} \leqslant \eAMP_{2} ,\ \forall \Phi \in \AMprecis{n',\delta '}(\Omega ),\quad  \|\Phi \|_{\eAMP_{2}} \leqslant  \left( \frac{\eAMP_{1} }{\eAMP_{2}} \right)^{\delta '} \|\Phi \|_{\eAMP_{1}} = C_{\delta ',\eAMP_{1} ,\eAMP_{2} } \,\|\Phi \|_{\eAMP_{1}} .
\end{equation}
This implies that $\|\mcla_{d,0} \|_{\eAMP'} \leqslant  C_{\delta ,\eAMP,\eAMP'} \;\eAMP^{-d} \;\gamma (d)$.\medskip\smallskip

Secondly, let $n>0$. Let us assume that the property holds for~any $n'<n$ and let us show it at rank $n$. By~induction hypothesis, we have $\mcla_{d,n-1-2m} \in \AMprecis{n-1-2m,d-n+1+2m}(\Omega )$ for~any $m\in \mathbb{N}$. Hence, using \cref{action des op R sur AMprecis}, we get successively that $\partial _{y} \mcla_{d,n-1-2m|\Gamma } \in \AMprecis{n-2m,d-n+2m}(\Omega ) \subset  \AMprecis{n,d-n+2m}(\Omega )$ and $(\partial _{x} ^2 +k_{1} ^2 )^{m} \partial _{y} \mcla_{d,n-1-2m|\Gamma } \in  \AMprecis{n,d-n}(\Omega )$, and finally by \cref{def: mcl}:
\[\mcla_{d,n} := \langle -k_{0} ^2 \ROmegaOmega\rangle \circ \RGammaOmega\bigg( \sum_{m=0}^\infty  \frac{\mu _{0} }{\mu _{1}} \,T_{m} \cdot  (\partial _{x} ^2 +k_{1} ^2 )^{m} \partial _{y} \mcla_{d,n-1-2m|\Gamma } \bigg) \in  \AMprecis{n,d-n}(\Omega ).\]
Now it only remains to show that $\|\mcla_{d,n} \|_{\eAMP'} \leqslant  C_{\delta ,n,\eAMP,\eAMP'} \;\eAMP^{-d} \;\gamma (d)$. Let $M_{n} :=\max \{m\in \mathbb{N} \mid n-1-2m\geqslant 0\}$. Using first \cref{norme de <Rlap>} and $\mcla_{d,n} \in \AMprecis{n,d-n}(\Omega ) \subset  \AMprecis{n,\delta -n}(\Omega )$, and then \cref{norme des op R sur series formelles}, we obtain
\begin{align*}
    \|\mcla_{d,n} \|_{\eAMP'} &\leqslant  C_{\delta -n,n,\eAMP'} \,\bigg\|\RGammaOmega\bigg( \sum _{m=0} ^{M_{n}} \frac{\mu _{0} }{\mu _{1}} \,T_{m} \cdot  (\partial _{x} ^2 +k_{1} ^2 )^{m} \partial _{y} \mcla_{d,n-1-2m|\Gamma } \bigg) \bigg\|_{\eAMP'} \\
    &\leqslant  C_{\delta ,n,\eAMP'} \sum _{m=0} ^{M_{n}} \big\| (\partial _{x} ^2 +k_{1} ^2 )^{m} \partial _{y} \mcla_{d,n-1-2m|\Gamma } \big\|_{\eAMP'} 
\end{align*}
\newcommand{\eAMPmed}[1]{\eAMP_{#1}}%
For~any $k\in [\![0,n]\!]$, let $\eAMPmed{k} := \frac{(k+1)\eAMP' + (n-k)\eAMP}{n+1}$. We have $\eAMP<\eAMPmed{0} <\eAMPmed{1} <\cdots <\eAMPmed{n} =\eAMP'$. \Cref{norme des op R sur series formelles} and \cref{eq: comparaison normes sur series formelles} imply, for~any $k\in [\![2,n]\!]$ and $\Phi \in \AMprecis{n,\delta -n}(\Gamma )$,
\begin{align*}
    \big\|(\partial _{x} ^2 +k_{1} ^2 ) \Phi  \big\|_{\eAMPmed{k}} 
    &\leqslant  \big\|\partial _{x} ^2  \Phi  \big\|_{\eAMPmed{k}} + |k_{1} ^2 |\,\|\Phi \|_{\eAMPmed{k}} \\
    &\leqslant  C_{\delta ,n,\eAMPmed{k},\eAMPmed{k-1}} \,\|\partial _{x} \Phi \|_{\eAMPmed{k-1}} + |k_{1} ^2 |\,\|\Phi \|_{\eAMPmed{k}} \\[1ex]
    &\leqslant  C_{\delta ,n,\eAMPmed{k},\eAMPmed{k-1}} \,C_{\delta ,n,\eAMPmed{k-1},\eAMPmed{k-2}}\,\|\Phi \|_{\eAMPmed{k-2}} +C_{\delta ,n,\eAMPmed{k},\eAMPmed{k-2}} \,\|\Phi \|_{\eAMPmed{k-2}} \\[1ex]
    &\leqslant  C_{\delta ,n,\eAMP,\eAMP'} \,\|\Phi \|_{\eAMPmed{k-2}} .
\end{align*}
Therefore\vspace{-1ex}
\begin{align*}
    \|\mcla_{d,n} \|_{\eAMP'} &\leqslant  C_{\delta ,n,\eAMP,\eAMP'} \sum _{m=0} ^{M_{n}} \big\|\partial _{y} \mcla_{d,n-1-2m|\Gamma } \big\|_{\eAMPmed{n-2m}} \\
    &\leqslant  C_{\delta ,n,\eAMP,\eAMP'} \sum _{m=0} ^{M_{n}} \big\|\mcla_{d,n-1-2m} \big\|_{\eAMPmed{n-1-2m}} \qquad  \text{ by \cref{norme des op R sur series formelles}}  \\
    &\leqslant  C_{\delta ,n,\eAMP,\eAMP'} \;\eAMP^{-d} \;\gamma (d)
\end{align*}
where we used in the last line the induction hypothesis, which applies because $\eAMPmed{n-1-2m} \geqslant \eAMPmed{0} >\eAMP$ and $\eAMP > \eAMPmin{\delta ,n} \geqslant  \eAMPmin{\delta ,n-1-2m}$.
\end{demo}

\subsection{The general case}
\label{sec: convergence des mcl cas general}

\begin{textesansboite}
In this section, we do~not assume $\pst \in \mathbb{Q}$ anymore -- $\Theta $ can be any real in $(0,2\pi )$ -- and we prove that $\mcla_{d,n}$ is a converging series for~any $(d,n)$ (\cref{convergence des mcl}). We will continue to use the following objects introduced in the two previous sections, which are still well-defined for~any $\Theta \in (0,2\pi )$: $\|.\|$ on $\mathcal{A}_{d} (\Omega )$, $\AOprecis{d,n,\delta }$, $\gamma $, $\normg{.}$, $\AMprecis{n,\delta }(\Omega )$ and $\|.\|_{\eAMP}$ on $\Amacro(\Omega )$.\\

The final result of the previous section is \cref{norme SF de mcl cas rationnel}. We start by extending it to any $\Theta $ in \cref{norme SF de mcl cas general}. We could not treat the case $\pst \not\in  \mathbb{Q}$ together with $\pst \in \mathbb{Q}$ in the previous sections, because $\RGammaOmega$ and $\ROmegaOmega$ do~not have continuity constants that depend nicely on $d$ when $\pst \not\in \mathbb{Q}$. See \cref{rq: norme naive des op R} for more details. Instead we rely in this case on the explicit expression of $\mcla_{d,n}$ given by \cref{expression sing heterog CL} when $\pst\not\in \mathbb{Q}$. Hence, the proof of the next result is mainly calculatory manipulations of the series expansions of the Bessel functions.
\end{textesansboite}

\begin{proposition}{}
\label{norme SF de mcl cas general}%
For~any $\Theta  \in (0,2\pi )$, there exists a family $(\eAMPmin{\delta ,n})_{d\in \mathbb{R},n\in \mathbb{N}} \in (\mathbb{R}_+^*)^{\mathbb{R}\times \mathbb{N}}$, that is increasing w.r.t.\ $n$, s.t.\ the results of \cref{norme SF de mcl cas rationnel} hold. Namely, we have for~any $\delta \in \mathbb{R}$, $d\in \TZ*\cap [\delta ,\infty )$, $n\in \mathbb{N}$ and $\eAMP'>\eAMP>\eAMPmin{\delta ,n}$:
\[\mcla_{d,n} \in \AMprecis{n,d-n}(\Omega ) \qquad   \text{ and } \qquad  \|\mcla_{d,n} \|_{\eAMP'} \leqslant  C_{\delta ,n,\eAMP,\eAMP'} \cdot \eAMP^{-d} \cdot \gamma (d).\]\fbeq
\end{proposition}

\begin{demo}
It suffices to prove the case $\pst\not\in \mathbb{Q}$. According to \cref{expression sing heterog CL}, in this case there exists $(a_{n,m} )_{m} \in \mathbb{C}^{[\![-n,n]\!]}$, that does~not depend on $d$, s.t.\ for~any $d\in \TZ*$, we have
\[\mcla_{d,n} = \Gamma (d+1) \left( \frac{2}{k_{0}} \right)^{d} \sum _{m=-n} ^{n} a_{n,m} \, \BesselOmega_{d+m} .\]
It~is easy to check that $\BesselOmega_{q'} \in \AMprecis{0,q}(\Omega )$ for~any reals $q<q'$ in $\mathbb{Z}+\TZ$, which yields $\mcla_{d,n} \in \AMprecis{n,d-n}(\Omega )$. Thus, it only remains to show the following:
\[\forall m\in [\![-n,n]\!], \qquad  A_{m} := \bigg\|\Gamma (d+1) \left( \frac{2}{k_{0}} \right)^{d} \BesselOmega_{d+m} \bigg\|_{\eAMP'} \leqslant  C_{\delta ,n,\eAMP,\eAMP'} \cdot \eAMP^{-d} \cdot \gamma (d).\]
By~definition of $\|.\|_{s}$ and $\BesselOmega_{d+m}$, we have 
\begin{align*}
    A_{m} &= |\Gamma (d+1)| \left| \frac{2}{k_{0}} \right|^{d}  \sum_{k=0}^\infty  \eAMP'^{-d-m-2k} \gamma (d+m+2k) \bigg\| \frac{(-1)^{k} }{k!\, \Gamma (d+m+k+1)} \left( \frac{k_{0} }{2} \right)^{d+m+2k} \Im[(\alpha z)^{d+m+k} \overline{\alpha z}^{k} ] \bigg\| \\
    &= \eAMP'^{-d} \,\gamma (d) \sum_{k=0}^\infty  \underbrace{ \frac{\gamma (d+m+2k)}{\gamma (d)} \cdot  \left| \frac{\Gamma (d+1)}{k!\, \Gamma (d+m+k+1)} \right| }_{=:B_{m,k}} \,\cdot \; \left| \frac{k_{0} }{2\eAMP'} \right|^{m+2k}
    \stepcounter{equation}\tag{\theequation}\label{eq: preuve norme SF de mcl cas irrationnel 1}
\end{align*}
Let us estimate $B_{m,k}$. It involves two ratios of factorial-like expressions. Similarly to $j \mapsto  \frac{(j+\ell )!}{j!}$ with $j,\ell \in \mathbb{N}$ for instance, we will see that these ratios have a polynomial behavior. We will then bound this behavior by an exponential, as it~is sufficient for us. However that behavior holds only when the arguments of $\gamma $ and $\Gamma $ are greater than 1, as $\gamma $ is constant on $(-\infty ,1]$ and $\Gamma $ has poles in $(-\infty ,1]$. Hence we will need to treat separately the first terms.\\
To begin, note that, for~any $x\in [1,\infty )$ and $k\in \mathbb{N}$, we have
\begin{equation}{\label{eq: preuve norme SF de mcl cas irrationnel 2}}
    \frac{\gamma (x+2k)}{\gamma (x)} = \sqrt{ \frac{\Gamma (x+2k)}{\Gamma (x)} } = \prod _{j=0} ^{2k-1} \sqrt{x+j} = \prod _{j=0} ^{k-1} \sqrt{(x+2j)(x+2j+1)} \leqslant  \prod _{j=0} ^{k-1} (x+2j+1)
\end{equation}
Let us also show that, for~any $a,b\in \mathbb{Z}$ and $\eAMP''>1$,
\begin{equation}{\label{eq: majoration de quotients des fonctions gamma}}
    \frac{\gamma (d+b)}{\gamma (d+a)} \leqslant C_{\delta ,a,b,\eAMP''} \,\eAMP''^{d} \qquad  \text{ and } \qquad  \left| \frac{\Gamma (d+b)}{\Gamma (d+a)} \right| \leqslant C_{\delta ,a,b,\eAMP''} \,\eAMP''^{d} .
\end{equation}
\begin{itemize}
    \item When $d+a\geqslant 1$ and $d+b\geqslant 1$, \cref{eq: preuve norme SF de mcl cas irrationnel 2} and the properties of the Gamma function imply that there exists $P_{a,b} \in \mathbb{R}[X]$ s.t $\frac{\gamma (d+b)}{\gamma (d+a)}$ and $\frac{\Gamma (d+b)}{\Gamma (d+a)}$ are bounded by $P_{a,b} (d)$. Then, there exists $C_{\delta ,a,b,\eAMP''}$ s.t.\ $P_{a,b} (d) \leqslant  C_{\delta ,a,b,\eAMP''} \,\eAMP''^{d}$ for~any $d\geqslant \delta $.
    \item Otherwise, we have $d\in \TZ*\cap [\delta ,1-\min (a,b)]$ by assumption on $d$, so $d$ can take only a finite number of values.
\end{itemize}
Now, let $\ell :=\max (0,\lceil {-\delta +n+1}\rceil )\in \mathbb{N}$. With this choice of $\ell $, we have, for~any $k>\ell $, $d+m+\ell \geqslant 1$ so the arguments of $\gamma $ and $\Gamma $ in the expression of $B_{m,k}$ are greater than 1. Thus, for~any $k>\ell $, we can write:
\begin{alignat*}{2}
    B_{m,k} &\leqslant  C_{\delta ,n,\eAMP''} \,\eAMP''^{d} \,\frac{\gamma (d+m+2k)}{\gamma (d+m+2\ell )} \, \left| \frac{\Gamma (d+m+\ell +1)}{\Gamma (d+m+k+1)} \right| \, \frac{1}{k!} &\ & \text{ by \cref{eq: majoration de quotients des fonctions gamma}} \\
    &\leqslant  C_{\delta ,n,\eAMP''} \,\eAMP''^{d} \times \prod _{j=0} ^{k-\ell -1} (d+m+2\ell +2j+1) \times \prod _{j=\ell } ^{k-1}  \frac{1}{d+m+j} \times \frac{1}{\ell !} \prod _{j=\ell } ^{k-1}  \frac{1}{j+1} && \text{ by \cref{eq: preuve norme SF de mcl cas irrationnel 2}} \\
    &\leqslant  C_{\delta ,n,\eAMP''} \,\eAMP''^{d} \times \prod _{j=\ell } ^{k-1} (d+m+2j+1) \times \frac{1}{\ell !} \prod _{j=\ell } ^{k-1}  \frac{1}{d+m+2j} && \text{ using ``$ab\geqslant a+b-1$''}  \\
    &\leqslant  C_{\delta ,n,\eAMP''} \,\eAMP''^{d} \times \prod _{j=\ell } ^{k-1} \left( 1+ \frac{1}{d+m+2j} \right) && \\
    &\leqslant  C_{\delta ,n,\eAMP''} \,\eAMP''^{d} \, 2^{k} 
\end{alignat*}
In addition, \cref{eq: majoration de quotients des fonctions gamma} shows that, for~any $k\in [\![0,\ell ]\!]$, we have $B_{m,k} \leqslant  C_{\delta ,n,\eAMP''} \eAMP''^{d}$. Combining the estimates of $B_{m,k}$ with \cref{eq: preuve norme SF de mcl cas irrationnel 1}, we obtain:
\[A_{m} \leqslant  C_{\delta ,n,\eAMP''} \left( \frac{\eAMP''}{\eAMP'} \right)^{d} \gamma (d) \left( 1+ \left| \frac{k_{0} }{2\eAMP} \right|^{m} \sum _{k>\ell } \left| \frac{k_{0} }{\sqrt{2} \eAMP'} \right|^{2k} \right) 
\leqslant  C_{\delta ,n,\eAMP,\eAMP''} \left( \frac{\eAMP''}{\eAMP'} \right)^{d} \gamma (d) \]
provided that $\eAMP'> \frac{|k_{0} |}{\sqrt{2}} +1$. Finally, it suffices to take $\eAMP'' := \frac{\eAMP'}{\eAMP}$ and to define $\eAMPmin{\delta ,n} := \frac{|k_{0} |}{\sqrt{2}} +1 +n$ for this angle $\Theta $ (the ``$+\,n$'' just makes the sequence increasing).
\end{demo}

\begin{textesansboite}
The next two results link the pointwise evaluations of corner singularities to their norms $\|.\|$ and $\|.\|_{\eAMP}$ and show that any $\Phi \in \AMprecis{n,\delta }(\Omega )$ s.t.\ $\|\Phi \|_{\eAMP} <\infty $ is a converging series.
\end{textesansboite}

\begin{lemmesansboite}{}
\label{majoration des fonctions de A par leur norme}
Let $\delta ,d\in \mathbb{R}$, $n\in \mathbb{N}$ and $\varphi \in \AOprecis{d,n,\delta }$. For~any $r>0$ and $\theta \in [0,\Theta ]$ we have in polar coordinates
\[|\varphi (r,\theta )| + r|\nabla  \varphi (r,\theta )| \;\leqslant \; C_{n} \, (1+\max (|\delta |,d)) \,r^{d} \, (1+|{\ln r}|)^{2n+1} \|\varphi \|.\]\fbeq[1.5]
\end{lemmesansboite}

\begin{demo}
By~definition of $\|.\|$ on $\AOprecis{d,n,\delta }$, and similarly to the proof of \cref{action des op R sur A precis,norme naive des op R}, it suffices to consider the case where $\varphi =\Im[(\alpha z)^{q} \,\overline{\alpha z}^{k} \log ^{j} \!(\alpha z)]$ with $ (q,k,j)\in \mathbb{R}\times \mathbb{N}^2 $, $q+k=d$ and $q-k\geqslant \delta $. For~any $r>0$ and $\theta \in [0,\Theta ]$, we have
\begin{align*}
    |\varphi (r,\theta )| &= \big|r^{d} \Im\big[e^{\mathrm{i}(q-k)(\theta -\Theta )} (\ln r+\mathrm{i}(\theta -\Theta ))^{j} \big]\big|\\
    &\leqslant  C_{j} \,r^{d} \,(1+|{\ln r}|)^{j} \\
    &\leqslant  C_{n} \,r^{d} \,(1+|{\ln r}|)^{2n+1} \|\varphi \|
\end{align*}
because $j\leqslant 2n+1$ (by~definition of $\AOprecis{d,n,\delta }$) and $\|\varphi \|=1$. Similarly, we have
\begin{align*}
    |\partial _{r} \varphi (r,\theta )| &= \big|r^{d-1} \Im\big[e^{\mathrm{i}(q-k)(\theta -\Theta )} \big[d(\ln r+\mathrm{i}(\theta -\Theta ))^{j} +j(\ln r+\mathrm{i}(\theta -\Theta ))^{j-1} \big]\big]\big|\\
    &\leqslant  C_{n} \,r^{d-1} (1+|d|)\,(1+|{\ln r}|)^{2n+1} \|\varphi \|
\end{align*}
with $|d|\leqslant \max (|\delta |,d)$ because $d\geqslant \delta $ (see \cref{rq: lien entre d q k delta}). Likewise:
\begin{align*}
    |r^{-1} \partial _{\theta } \varphi (r,\theta )| &= \big|r^{d-1} \Im\big[\mathrm{i}e^{\mathrm{i}(q-k)(\theta -\Theta )} \big[(q-k)(\ln r+\mathrm{i}(\theta -\Theta ))^{j} +j(\ln r+\mathrm{i}(\theta -\Theta ))^{j-1} \big]\big]\big|\\
    &\leqslant  C_{n} \,r^{d-1} (1+|q-k|)\,(1+|{\ln r}|)^{2n+1} \|\varphi \|
\end{align*}
with $|q-k|\leqslant \max (|\delta |,d)$ because $q-k\leqslant q+k=d$ and $q-k\geqslant \delta $.
\end{demo}

\newcommand{\restPhi}[1][q]{\rho_{#1}}%
\begin{proposition}{}
\label{cvg des SF de norme finie}%
Let $\delta \in \mathbb{R}$, $n\in \mathbb{N}$, $\Phi \in \AMprecis{n,\delta }(\Omega )$ and $\eAMP\in \mathbb{R}$ s.t.\ $\|\Phi \|_{\eAMP} <\infty $. There exists $f\in \mathcal{C}^{1} (\overline{\Omega } \setminus \{0\})$ s.t.
\[T_{<q} (\Phi ) \longrightarrow  f \quad  \text{ and } \quad  \nabla \big[T_{<q} (\Phi )\big] \longrightarrow  \nabla  f\]
pointwise when $q\rightarrow \infty $. Moreover, for~any $q\in \mathbb{R}$, let us denote $\restPhi := f - T_{<q} (\Phi )$ the remainder of the series $\Phi $. Let $R>0$ and $\eAMP'>\eAMP>0$. For~any $q'\geqslant q\geqslant \delta $, $r\in (0,R)$ and $\theta \in [0,\Theta ]$ we have
\vspace{-1ex}%
\begin{equation}{\label{eq: vitesse de cvg des SF de norme finie}}
    |\restPhi[q'](r,\theta )| + r|\nabla  \restPhi[q'](r,\theta )| \leqslant  C_{\delta ,n,R,\eAMP,\eAMP'} \, \frac{(\eAMP' r)^{q} }{\gamma (q)} \, (1+|{\ln r}|)^{2n+1} \|\Phi \|_{\eAMP} .
\end{equation}\fbeq[1.3]
\end{proposition}

\begin{demo}
Let $\eAMP'>\eAMP>0$, $q'\geqslant q\geqslant \delta $, $0<r\leqslant R$ and $\theta \in [0,\Theta ]$. Let us denote $(\varphi _{d} )_{d\in \mathbb{R}} :=\Phi $ the family of the components of $\Phi $. Since $\Phi \in \AMprecis{n,\delta }(\Omega )$, we have $\varphi _{d} \in \AOprecis{d,n,\delta }$ for~any $d\in \mathbb{R}$. Using \cref{majoration des fonctions de A par leur norme}, we get
\begin{align*}
    |\varphi _{d} (r,\theta )| + r|\nabla  \varphi _{d} (r,\theta )| &\leqslant  C_{n} \,(1+\max (|\delta |,d)) \,r^{d} \, (1+|{\ln r}|)^{2n+1} \,\|\varphi _{d} \| \\
    &\leqslant  C_{\delta ,n,\eAMP'/\eAMP} \,\Big( \frac{\eAMP'r}{\eAMP} \Big)^{d} (1+|{\ln r}|)^{2n+1} \,\|\varphi _{d} \| \\[.5ex]
    &= C_{\delta ,n,\eAMP,\eAMP'} \, \frac{(\eAMP'r)^{d} }{\gamma (d)} \, (1+|{\ln r}|)^{2n+1} \cdot  s^{-d} \|\varphi _{d} \|_{\gamma } .
\end{align*}
Let us show that, for~any $d\geqslant q$ we have $\displaystyle \frac{(\eAMP'r)^{d} }{\gamma (d)} \leqslant  C_{\delta ,R,\eAMP'} \frac{(\eAMP'r)^{q} }{\gamma (q)}$. There are three possible cases.
\begin{enumerate}
    \item Let us assume $d\geqslant q\geqslant 1$. Let $I:=[1,\infty )$ and
    \[g_{r} :d\in I\longmapsto \ln \!\left( \frac{(\eAMP'r)^{d} }{\gamma (d)} \right) = d\ln (\eAMP'r)- \frac{1}{2} \ln  \Gamma (d).\]
    \def\dmax{m_r}%
    Since $d\mapsto \ln \Gamma (d)$ is convex on $I$, $g_{r}$ is concave. In addition, $g_{r} (d)\rightarrow -\infty $ as $d\rightarrow \infty $. Therefore, there exits $\dmax \in I$ (that may be equal to 1) s.t.\ $g_{r}$ increases on $[1,\dmax]$ and decreases on $[\dmax,\infty )$. If $q\geqslant \dmax$, then we directly have $(\eAMP'r)^{d} /\gamma (d) \leqslant  (\eAMP'r)^{q} /\gamma (q)$. Now, let us assume $q<\dmax$. In this case, we have $\dmax>1$, so
    \[0=g'(\dmax)=\ln (\eAMP'r)- \frac{1}{2} \psi (\dmax) \qquad  \text{which implies}  \qquad  \dmax = \psi ^{-1} ( 2\ln (\eAMP'r))\]
    because the Digamma function $\psi $ is an increasing bijection from $I$ to $[\psi (1),\infty )$. 
    In~particular, this yields $\ln (\eAMP'r)= \frac{1}{2} \psi (\dmax)\geqslant  \frac{1}{2} \psi (1)$. Hence
    \begin{align*}
        g_{r} (d)-g_{r} (q) &\leqslant  \max _{I} g_{r} - \min _{[1,\dmax]} g_{r} \\
        &= g_{r} (\dmax) - g_{r} (1)\\
        &= \psi ^{-1} ( 2\ln (\eAMP'r)) \ln (\eAMP'r)- \frac{1}{2} \ln  \Gamma (\dmax) - \ln (\eAMP'r) \\
        &\leqslant  \psi ^{-1} ( 2\ln (\eAMP'R)) \ln (\eAMP'R)- \frac{1}{2} \ln (\min \nolimits_{I} \Gamma  )- \frac{1}{2} \psi (1) \\
        &= C_{R,\eAMP'}
    \end{align*}
    This proves $\frac{(\eAMP'r)^{d} }{\gamma (d)} \leqslant  C_{\delta ,R,\eAMP'} \frac{(\eAMP'r)^{q} }{\gamma (q)}$ in this case.
    \item Let us assume $q\leqslant d\leqslant 1$. We have $(\eAMP'r)^{d} = (\eAMP'r)^{d-q} (\eAMP'r)^{q} \leqslant  C_{\delta ,R,\eAMP'} (\eAMP'r)^{q}$ with $C_{\delta ,R,\eAMP'} := 1$ if $\eAMP'r<1$ and $C_{\delta ,R,\eAMP'} := (\eAMP'R)^{1-\delta }$ otherwise (because $q\geqslant \delta $). In addition, we have $\gamma (q)=\gamma (d)=1$ by~definition of $\gamma $.
    \item Otherwise, we have $q\leqslant 1\leqslant d$. Using the two previous cases, we obtain $\frac{(\eAMP'r)^{d} }{\gamma (d)} \leqslant  C_{\delta ,R,\eAMP'} \frac{(\eAMP'r)^{1} }{\gamma (1)} \leqslant  C_{\delta ,R,\eAMP'} \frac{(\eAMP'r)^{q} }{\gamma (q)}$.
\end{enumerate}
Therefore we have shown that, for~any $d\geqslant q$, we have
\begin{equation}{\label{eq: preuve cvg des SF de norme finie 1}}
    |\varphi _{d} (r,\theta )| + r|\nabla  \varphi _{d} (r,\theta )| \leqslant  C_{\delta ,n,R,\eAMP,\eAMP'} \, \frac{(\eAMP'r)^{q} }{\gamma (q)} \, (1+|{\ln r}|)^{2n+1} \cdot  s^{-d} \|\varphi _{d} \|_{\gamma } .
\end{equation}
Let $r_{1} \in (0,R)$ and $D:=\{\mathbf{x}\in \Omega  \mid |\mathbf{x}|\in [r_{1} ,R]\}$. If $q\geqslant 0$, \cref{eq: preuve cvg des SF de norme finie 1} implies
\begin{align*}
    \|\varphi _{d} \|_{L^{\infty } (D)} + r_{1} \|\nabla  \varphi _{d} \|_{L^{\infty } (D)} 
    &\leqslant  C_{\delta ,n,R,\eAMP,\eAMP'} \, \frac{(\eAMP'R)^{q} }{\gamma (q)} \, (1+|{\ln r_{1} }|+|{\ln R}|)^{2n+1} \cdot  s^{-d} \|\varphi _{d} \|_{\gamma } 
\end{align*}
where we used $q\geqslant 0$ to bound $r^{q}$ by $R^{q}$. By~sum, we get
\begin{align*}
    \smash{\sum _{\substack{d\in \operatorname{supp}(\Phi )\\d\geqslant q}}} \|\varphi _{d} \|_{L^{\infty } (D)} + r_{1} \|\nabla  \varphi _{d} \|_{L^{\infty } (D)} &\leqslant  C_{\delta ,n,r_{1} ,R,\eAMP,\eAMP'} \, \frac{(\eAMP'R)^{q} }{\gamma (q)} \, \sum _{d\geqslant q} s^{-d} \|\varphi _{d} \|_{\gamma } \\
    &\leqslant  C_{\delta ,n,r_{1} ,R,\eAMP,\eAMP'} \, \frac{(\eAMP'R)^{q} }{\gamma (q)} \, \|\Phi \|_{\eAMP} \ \underset{q\rightarrow \infty }{\longrightarrow} \; 0 
\end{align*}
This proves that $T_{<q} (\Phi )$ and its~derivatives converge on $D$ as $q\rightarrow \infty $. This convergence is uniform, which ensures that the limit $f$ is $\mathcal{C}^{1}$ and $\nabla  [\lim _{q\rightarrow \infty } T_{<q} (\Phi )] = \lim _{q\rightarrow \infty } \nabla  T_{<q} (\Phi )$. Finally, using
\[\restPhi[q'] := f - T_{<q'} (\Phi ) = \sum _{d\in \operatorname{supp}(\Phi ),\,d\geqslant q'} \varphi _{d}\]
and summing \cref{eq: preuve cvg des SF de norme finie 1} over $d\geqslant q'$ proves \cref{eq: vitesse de cvg des SF de norme finie}.
\end{demo}

\begin{textesansboite}
The next theorem is the main result of this section. We recall that $\eAMPmin{\delta ,n}$ is introduced in \cref{norme SF de mcl cas general}.
\end{textesansboite}

\newcommand{\restmcl}{\rho_{d,n,q}}%
\begin{theoremesecable}{convergence of the formal series $\mcla_{d,n}$}
\label{convergence des mcl}%
\begin{enumerate}
    \item Let $d\in \TZ*$ and $n\in \mathbb{N}$. There exists a function in $\mathcal{C}^{1} (\overline{\Omega } \setminus \{0\},\mathbb{C})$, that we will denote $\mcll_{d,n}$, s.t.
    \vspace{-1ex}%
    \[T_{<q} (\mcla_{d,n} ) \longrightarrow  \mcll_{d,n} \quad  \text{ and } \quad  \nabla \big[T_{<q} (\mcla_{d,n} )\big] \longrightarrow  \nabla  \mcll_{d,n}\]
    when $q\rightarrow \infty $. The convergence is uniform on $(\overline{\Omega } \setminus \{0\}) \cap  B(0,R )$ for~any $R>0$.
    \item Let $d\in \TZ*$. For~any $n\in -\mathbb{N}^*$, let us denote $\mcll_{d,n} :=0$. With this convention, $\mcll_{d,n}$ satisfies for~any $n\in \mathbb{N}$:
    \begin{equation}{\label{eq: pb verifie par lim mcl}}
        \left\{\begin{array}{r@{\;=\;}ll}
            (\Delta +k_{0} ^2 )\mcll_{d,n} & 0 &\text{ in }\Omega  \\
            \mcll_{d,n} & \sum\limits_{m=0}^\infty  \frac{\mu _{0} }{\mu _{1}} \,T_{m} \cdot  (\partial _{x} ^2 +k_{1} ^2 )^{m} \partial _{y} \mcll_{d,n-1-2m} &\text{ on }\Gamma 
        \end{array}\right.
    \end{equation}
    in the sense of distributions.
    \item Let $\delta \in \mathbb{R}$, $n\in \mathbb{N}$, $R>0$ and $\eAMP'>\eAMP>\eAMPmin{\delta ,n}$. There exists $C_{\delta ,n,R,\eAMP,\eAMP'} >0$ s.t., for~any $d\in \TZ*\cap [\delta ,\infty )$, $q\geqslant \delta $, $0<r\leqslant R$ and $\forall \theta \in [0,\Theta ]$, we have in polar coordinates
    \begin{equation}{\label{eq: vitesse de cvg des mcl}}
        |\restmcl (r,\theta )| + r|\nabla  \restmcl (r,\theta )| \;\leqslant \; C_{\delta ,n,R,\eAMP,\eAMP'} \, \frac{\gamma (d)}{\gamma (q)} \,\eAMP^{-d} \,(\eAMP' r)^{q} \, (1+|{\ln r}|)^{2n+1} 
    \end{equation}
    where $\restmcl := \mcll_{d,n} - T_{<q} (\mcla_{d,n} )$ denotes the remainder of the series $\mcla_{d,n}$ at rank $q$.
\end{enumerate}
\end{theoremesecable}
\nomenclature[56-13-23-21]{$\mcll_{d,n}$}{Far field singularity (as a function on $\Omega $)}

\begin{remarques}
\begin{itemize}
    \item \cref{eq: pb verifie par lim mcl} is the limit version of \cref{eq: pb verifie par mcl qd n=0}--\cref{eq: pb verifie par mcl qd n>0}, a similar relation on $\mcla_{d,n}$. The laplacian and derivatives $\partial _{x}$ on $\Gamma $ are written in the weak sens (i.e.\ in $\mathcal{D}'(\Omega )$ and $\mathcal{D}'(\Gamma )$ resp.), but one can actually show that $\mcll_{d,n}$ is $\mathcal{C}^{\infty }$ on $\overline{\Omega } \setminus \{0\}$.
    \item Like in \cref{decompo modale des CL,decompo en modes tronques des CL,convergence de uDA}, the logarithmic factor of \cref{eq: vitesse de cvg des mcl} is~not optimal and can be ignored if $\pst \not\in \mathbb{Q}$. See \cref{rq: decompo modale des CL} for more details.
    \item \cref{eq: vitesse de cvg des mcl} implies: $\forall q\in \mathbb{R},\; \mcll_{d,n} = T_{\leqslant q} (\mcla_{d,n} )+o(r^{q} )$. In other words, the asymptotic expansion of $\mcll_{d,n}$ at the corner is described by $\mcla_{d,n}$. This explains the notation $\mcla_{d,n}$, that is similar to $\bfuDA_{p,\ell }$ and $\SDA_{p,\ell }$.
\end{itemize}
\end{remarques}

\begin{demo}
\def\eAMPmed{\eAMP^*}
Let $\eAMPmed := (\eAMP+\eAMP')/2$. Combining \cref{norme SF de mcl cas general,cvg des SF de norme finie} directly shows that there exists $\mcll_{d,n} \in \mathcal{C}^{1} (\overline{\Omega } \setminus \{0\},\mathbb{C})$ s.t.\ $T_{<q} (\mcla_{d,n} )\rightarrow \mcll_{d,n}$ and $\nabla \big[T_{<q} (\mcla_{d,n} )\big]\rightarrow \nabla \mcll_{d,n}$ pointwise. It also gives:
\begin{align*}
    |\restmcl(r,\theta )| + r|\nabla  \restmcl(r,\theta )| &\leqslant  C_{\delta ,n,\eAMPmed,\eAMP'} \, \frac{(\eAMP'r)^{q} }{\gamma (q)} \, (1+|{\ln r}|)^{2n+1} \, \|\mcla_{d,n} \|_{\eAMPmed} \\
    &\leqslant  C_{\delta ,n,\eAMPmed,\eAMP'} \, \frac{(\eAMP'r)^{q} }{\gamma (q)} \, (1+|{\ln r}|)^{2n+1} \times C_{\delta ,n,\eAMP,\eAMPmed} \, \eAMP^{-d} \,\gamma (d)
\end{align*}
which proves \cref{eq: vitesse de cvg des mcl}. Next, let us prove that the convergence is uniform on every ball $B(0,R )$. For~any $0<r\leqslant R$, \cref{eq: vitesse de cvg des mcl} yields
\begin{align*}
    |\restmcl(r,\theta )| + |\nabla  \restmcl(r,\theta )| 
    &\leqslant  \frac{C_{R} }{r} \cdot  \big( |\restmcl(r,\theta )| + r|\nabla  \restmcl(r,\theta )| \big) \\
    &\leqslant  \frac{C_{R} }{r} \cdot  C_{\red{d},n,\eAMP,\eAMP'} \, \frac{(\eAMP'r)^{q} }{\gamma (q)} (1+|{\ln r}|)^{2n+1} \\
    &\leqslant  C_{d,n,\eAMP,\eAMP',R} \, \frac{(\eAMP'r)^{q-2} }{\gamma (q)} \cdot  r \, (1+|{\ln r}|)^{2n+1} \\
    &\leqslant  C_{d,n,\eAMP,\eAMP',R} \, \frac{(\eAMP'R)^{q-2} }{\gamma (q)} \quad\qquad  \text{if }  q\geqslant 2
\end{align*}
because $r\mapsto  r\,(1+|{\ln r}|)^{2n+1} $ is bounded on $(0,R]$.  Therefore the convergence of $T_{<q} (\mcla_{d,n} )$ and its gradient is uniform on $(\overline{\Omega } \setminus \{0\}) \cap  B(0,R )$ for~any $R>0$. Thus it remains only to prove \cref{eq: pb verifie par lim mcl}.\smallskip\medskip

Since $T_{<q} (\mcla_{d,n} )\rightarrow \mcll_{d,n}$ locally uniformly, this convergence also holds in $\mathcal{D}'(\Omega )$. Moreover, we have $(\Delta +k_{0} ^2 )\mcla_{d,n} =0$ by \cref{eq: pb verifie par mcl qd n=0}--\cref{eq: pb verifie par mcl qd n>0}, and $\Delta $ has degrees $-2$ as an operator on $\mathcal{A}(\Omega )$. Therefore we have
\[0= T_{<q} \big[(\Delta +k_{0} ^2 )\mcla_{d,n} \big] = \Delta \big[T_{<q+2} (\mcla_{d,n} )\big]+k_{0} ^2  T_{<q} (\mcla_{d,n} ) \xrightarrow[{q\rightarrow \infty }]{} \Delta \mcll_{d,n} + k_{0} ^2  \mcll_{d,n}\]
with convergence in $\mathcal{D}'(\Omega )$. Hence, $(\Delta + k_{0} ^2 )\mcll_{d,n} =0$.\\
Similarly, for~any $m\in \mathbb{N}$, we have $\partial _{y} \big[T_{<q} (\mcla_{d,n-1-2m} )\big]_{|\Gamma } \rightarrow  \partial _{y} \mcll_{d,n-1-2m|\Gamma }$ locally uniformly, so also in $\mathcal{D}'(\Gamma )$. Thus, using that $\partial _{x}$ and $\partial _{y|\Gamma ,y=0^{+}} $ have degre $-1$ as operators in the spaces $\mathcal{A}$, we get for~any $j\in \mathbb{N}$:
\begin{equation}{\label{eq: preuve pb verifie par lim mcl}}
    T_{<q} \big[ \partial _{x} ^{j} \partial _{y} \mcla_{d,n-1-2m} \big]_{|\Gamma } 
    = \partial _{x} ^{j} \partial _{y} \big[T_{<q+j+1} (\mcla_{d,n-1-2m} )\big]_{|\Gamma } \xrightarrow[q\rightarrow \infty ]{\mathcal{D}'(\Gamma )} \partial _{x} ^{j} \partial _{y} \mcll_{d,n-1-2m|\Gamma }
\end{equation}
Finally, we deduce the following limit calculation in $\mathcal{D}'(\Gamma )$ (note that the sums in $m$ below have only a~finite number of non zero terms)
\begin{align*}
    \sum_{m=0}^\infty  \frac{\mu _{0} }{\mu _{1}} \,T_{m} \cdot  (\partial _{x} ^2 +k_{1} ^2 )^{m} \partial _{y} \mcll_{d,n-1-2m} 
    &= \sum_{m=0}^\infty  \frac{\mu _{0} }{\mu _{1}} \,T_{m} \cdot  \lim _{q\rightarrow \infty } T_{<q} \big[ (\partial _{x} ^2 +k_{1} ^2 )^{m} \partial _{y} \mcla_{d,n-1-2m} \big] \qquad  \text{by \cref{eq: preuve pb verifie par lim mcl}} \\[-0.5ex]
    &= \lim _{q\rightarrow \infty } T_{<q} \bigg( \sum_{m=0}^\infty  \frac{\mu _{0} }{\mu _{1}} \,T_{m} \cdot  (\partial _{x} ^2 +k_{1}  ^2 )^{m} \partial _{y} \mcla_{d,n-1-2m|\Gamma } \bigg) \\
    &= \lim _{q\rightarrow \infty } T_{<q} (\mcla_{d,n} ) 
    \setbox0=\hbox{by \cref{eq: pb verifie par mcl qd n=0}--\cref{eq: pb verifie par mcl qd n>0}}
    \hspace{73mm-\wd0} \box0  \\
    &= \mcll_{d,n|\Gamma }    
\end{align*}
Note that here $T_{<q}$ denotes a truncation operator $\Amacro(\Gamma )\rightarrow \mathcal{A}(\Gamma )$ whereas $T_{m}$ denotes the Taylor coefficients of the tangent.
\end{demo}

\begin{corollaire}{estimate on $\mcll_{d,n}$}
\label{estimation des mcl}%
Let $\delta \in \mathbb{R}$, $n\in \mathbb{N}$, $R>0$ and $\eAMP>1$. There exists $C_{\delta ,n,R,\eAMP} >0$ s.t., for~any $d\in \TZ*\cap [\delta ,\infty )$, $0<r\leqslant R$ and $\forall \theta \in [0,\Theta ]$, we have
\[|\mcll_{d,n} (r,\theta )| + r|\nabla  \mcll_{d,n} (r,\theta )| \;\leqslant \; C_{\delta ,n,R,\eAMP} \,(\eAMP r)^{d-n} \, (1+|{\ln r}|)^{2n+1} .\]\fbeq
\end{corollaire}

\begin{demo}
Let us apply \cref{eq: vitesse de cvg des mcl} with $q:=d-n$. In this case, we have $T_{<q} (\mcla_{d,n} )=0$ so the remainder of $\mcla_{d,n}$ at order $q$ is equal to $\mcll_{d,n}$. Hence, \cref{eq: vitesse de cvg des mcl} states
\[|\mcll_{d,n} (r,\theta )| + r|\nabla  \mcll_{d,n} (r,\theta )| \leqslant  C_{\delta ',n,R,\eAMP',\eAMP''} \, \frac{\gamma (d)}{\gamma (d-n)} \,\eAMP'^{-d} \,(\eAMP'' r)^{d-n} \, (1+|{\ln r}|)^{2n+1} \]
with $\delta ':=\delta -n$ (so that $q\geqslant \delta '$), $\eAMP':=\eAMPmin{\delta ,n}+1$ and $\eAMP'':=\eAMP' \sqrt{\eAMP}$. Thanks to \cref{eq: majoration de quotients des fonctions gamma} from the proof of \cref{norme SF de mcl cas general}, we have $\frac{\gamma (d)}{\gamma (d-n)} \leqslant  C_{\delta ,n,\eAMP} \sqrt{\eAMP}^{d}$. Thus we get
\[|\mcll_{d,n} (r,\theta )| + r|\nabla  \mcll_{d,n} (r,\theta )| \leqslant  C_{\delta ,n,R,\eAMP} \,\Big( \sqrt{\eAMP} \, \frac{\eAMP''r}{\eAMP'} \Big)^{d-n} \, (1+|{\ln r}|)^{2n+1} \leqslant  C_{\delta ,n,R,\eAMP} \,(\eAMP r)^{d-n} \, (1+|{\ln r}|)^{2n+1} .\]
\end{demo}

\section{Proof of the modal decomposition}
\label{sec: preuve decompo modale CL}

\begin{textesansboite}
In this section, we prove \cref{decompo modale des CL,decompo en modes tronques des CL}, which are the main results of \cref{chap: decompo modale CL}.
\end{textesansboite}

\subsection{Proof of \texorpdfstring{\cref{decompo modale des CL}}{Theorem 7.1.1}}
\label{sec: preuve decompo modale CL 1}

\begin{textesansboite}
Let us start with two technical lemmas.
\end{textesansboite}

\begin{lemmesansboite}{}
\label{majoration coefs decompo modale Helmholtz homog}
Let $u\in H^{1} (\Omega  \cap B(0,\Rsource))$ be s.t.
\[\left\{\begin{array}{r@{\;=\;}ll}
    \Delta u + k_{0} ^2 u & 0 &\text{ in }\Omega \cap B(0,\Rsource) \\
    u & 0 &\text{ on }\partial \Omega  \cap B(0,\Rsource)  
\end{array}\right.\]
There exists a family $(s_{d} )_{d\in \TN}$ of complex numbers s.t.\ we have:
\begin{enumerate}
    \item $u = \sum _{d\in \TN} s_{d} \,\mcll_{d,0} $ with locally uniform convergence,
    \item and $\forall R\in (0,\Rsource),\, \exists  C_{R,u} >0,\, \forall d\in \TN, \ |s_{d} |\leqslant  C_{R,u} \.R^{-d}$.
\end{enumerate}
\end{lemmesansboite}

\begin{textesansboite}
The modal decomposition of a solution of the Helmholtz equation in a corner is classical ; the important result of this lemma is the bound on $s_{d}$. 
\end{textesansboite}

\begin{demo}
By~elliptic regularity, $u$ is $H^2 $ in $\Omega \cap \{r_{1} <r<r_{2} \}$ for~any $0<r_{1} <r_{2} <R$. Using Fourier series w.r.t.\ $\theta $, it~is easy to deduce that 
\begin{equation}{\label{eq: preuve majoration coefs decompo modale Helmholtz homog 1}}
    u(r,\theta ) = \sum _{d\in \TN} c_{d} (r) \sin (d\theta )
\end{equation}
where $c_{d} (r) := \frac{2}{\Theta } \int _{0} ^{\Theta } u(r,\theta )\,\sin (d\theta )\,\mathrm{d}\theta $, with convergence in $L^2 (\Omega \cap \partial B(0,r))$ for~any $r\in (0,\Rsource)$. This implies
\begin{equation}{\label{eq: preuve majoration coefs decompo modale Helmholtz homog 2}}
    \forall r\in (0,\Rsource),\qquad  \|u(r,\cdot )\|_{L^2 (0,\Theta )} ^2  = \frac{\Theta }{2} \sum _{d\in \TN} |c_{d} (r)|^2 .
\end{equation}
Integrating this relation over $r$ shows that $c_{d} \in L^2 (0,\Rsource)$ for~any $d$. Hence \cref{eq: preuve majoration coefs decompo modale Helmholtz homog 1} holds also with convergence in $L^2 (\Omega \cap B(0,\Rsource))$, so in $\mathcal{D}'(\Omega \cap B(0,\Rsource))$ as well. Thus we can differentiate \cref{eq: preuve majoration coefs decompo modale Helmholtz homog 1}. Using the equation $\Delta u + k_{0} ^2 u = 0$, we get $(r \frac{\mathrm{d}}{\mathrm{d}r} )^2 c_{d} + (k_{0} ^2 r^2 -d^2 )c_{d} =0$. The properties of the Bessel functions then show that there exist coefficients $a_{d} \in \mathbb{C}$
\[\forall d\in {\textstyle\TN},\ \forall r\in (0,\Rsource), \quad  c_{d} (r) = a_{d} Z_{d} (k_{0} r) + a_{-d} Z_{-d} (k_{0} r) \quad \text{ with }\quad  Z_{d} := \left\{\begin{array}{ll}
    Y_{d} & \text{if } d\in -\mathbb{N}^* \\
    J_{d} & \text{otherwise}
\end{array}\right.\]
Now, according to \cite{Leb72}, for~any $q\in \mathbb{R}^*$, there exists $b_{q} \in \mathbb{C}^*$ s.t.\ $Z_{q} (r) \sim b_{q} r^{q}$ when $q\rightarrow \infty $. Since $d\geqslant \pst> \frac{1}{2}$ and $c_{d} \in L^2 $, it yields $a_{-d} =0$ for~any $d\in \TN$. 
Using \cref{expression sing homog CL}, we deduce
\[u(r,\theta ) = \sum _{d\in \TN} a_{d} \,J_{d} (k_{0} r) \sin (d\theta ) = \sum _{d\in \TN} s_{d} \mcll_{d,0} (r,\theta ) \qquad \text{ with }\qquad  s_{d} := a_{d} \frac{(-1)^{d\pi /\Theta } }{\Gamma (d+1)} \left( \frac{k_{0} }{2} \right)^{d} \]
Moreover, \cite[p.\,44]{Wat95} states:
\[\forall q\in \mathbb{R}_+,\ \forall z\in \mathbb{C}^*, \qquad  \left| \Gamma (q+1) \left( \frac{2}{z} \right)^{q} J_{q} (z) - 1 \right| \leqslant  \exp \!\left( \frac{|z|^2 }{4(q+1)} \right) -1.\]
Let $C>0$ be s.t.\ $\exp (t^2 /4)-1< \frac{1}{2}$ for~any $t\in [0,C]$. We deduce:
\[\forall q\in \mathbb{R}_+,\ \forall z\in \mathbb{C}^*, \quad  \ |z|< C \sqrt{q+1} \ \Rightarrow \ \left| \Gamma (q+1) \left( \frac{2}{z} \right)^{q} J_{q} (z) \right| \geqslant  \frac{1}{2} .\]
Combined with \cref{eq: preuve majoration coefs decompo modale Helmholtz homog 2}, it yields, for~any $r\in (0,\Rsource)$ and $d_{0} \in \TN$ s.t.\ $C\sqrt{d_{0} +1} >|k_{0} |\Rsource$,
\[\infty > \sum _{d\geqslant d_{0}} |c_{d} (r)|^2  = \sum _{d\geqslant d_{0}} |a_{d} J_{d} (k_{0} r)|^2  = \sum _{d\geqslant d_{0}} \bigg|s_{d} \,\Gamma (d+1) \left( \frac{2}{k_{0}} \right)^{d} J_{d} (k_{0} r)\bigg|^2  \geqslant  \sum _{d\geqslant d_{0}} \frac{1}{4} |s_{d} r^{d} |^2  \]
so $s_{d} = O(r^{-d} )$ when $d\rightarrow \infty $. Finally it implies that \cref{eq: preuve majoration coefs decompo modale Helmholtz homog 1} converges locally uniformly.
\end{demo}

\begin{textesansboite}
The following lemma will be used in the proof of \cref{decompo modale des CL} to show that some series of the $\sigma _{d} (u_{p,\ell } )\,\mcll_{d,n}$ converge, without specifying the limit. What will remain to do to prove \cref{decompo modale des CL} is mainly to link those limits to the far fields.
\end{textesansboite}

\begin{lemme}{}
\label{convergence des series de mcl}%
Let $\delta \in \mathbb{R}$ and $(s_{d} )\in \mathbb{C}^{\TZ*}$ be s.t.\ $s_{d} =0$ for~any $d<\delta $ and $s_{d} = O(R^{-d} )$ for~any $R\in (0,\Rsource)$. Let us denote $D :=(\overline{\Omega } \setminus \{0\}) \cap B(0,\Rsource)$. Let $n\in \mathbb{N}$ and, for~any $N\in \mathbb{N}$,
\[f_{N} = \sum _{d\in \TZ*,\, \delta \leqslant d\leqslant N\pst} s_{d} \,\mcll_{d,n} .\]
Then $(f_{N} )$ converges pointwise on $D$ when $N\rightarrow \infty $ and its limit $f$ belongs to $\mathcal{C}^{1} (D,\mathbb{C})$. Moreover, for~any $R\in (0,\Rsource)$, $N\in \mathbb{N}$, $r\in (0,R)$ and $\theta \in [0,\Theta ]$, we have:
\begin{equation}{\label{eq: vitesse de cvg des series de mcl}}
    |(f-f_{N} )(r,\theta )| + r|\nabla  (f-f_{N} )(r,\theta )| \leqslant  C_{\delta ,n,(s_{d} ),R,\Rsource} \left( \frac{r}{R} \right)^{(N+1)\pst-n} (1+|{\ln r}|)^{2n+1} .
\end{equation}\fbeq[1.3]
\end{lemme}

\begin{demo}
Let $0<r_{1} \leqslant r_{2} \leqslant R<\Rsource$ and $D':=\overline{\Omega }\cap  \{r_{1} \leqslant r\leqslant r_{2} \}$. Let $\eAMP>1$ (to be chosen later). \Cref{estimation des mcl} implies that, for~any $d\geqslant \delta $, we have
\begin{equation}{\label{eq: convergence des series de mcl preuve 1}}
    \|\mcll_{d,n} \|_{L^{\infty } (D')} + r_{1} \|\nabla  \mcll_{d,n} \|_{L^{\infty } (D')} \leqslant  C_{\delta ,n,R,\eAMP} \,(\eAMP r_{2} )^{d-n\,} (1+|{\ln r_{1} }|+|{\ln r_{2} }|)^{2n+1} 
\end{equation}
provide that $d-n\geqslant 0$ or $r_{1} =r_{2}$, so that $r^{d-n}$ can be bounded by $r_{2} ^{d-n}$ for~any $r\in [r_{1} ,r_{2} ]$. Let $R'=(R+\Rsource)/2$ and $N\in \mathbb{N}$. By~assumption on $(s_{d} )$, we have $|s_{d} |\leqslant C_{(s_{d} ),R'} \,R'^{-d}$ for~any $d$. Moreover, let us assume that either $(N+1)\pst-n\geqslant 0$ or $r_{1} =r_{2}$. Using \cref{eq: convergence des series de mcl preuve 1} and the bound on $|s_{d} |$, we get:
\begin{multline*}
    A:= \sum _{\substack{d\in \TZ*\\ \mathclap{d\geqslant (N+1)\pst}}} |s_{d} |\,\Big(\|\mcll_{d,n} \|_{L^{\infty } (D')} + r_{1} \|\nabla  \mcll_{d,n} \|_{L^{\infty } (D')} \Big) \\[-4ex]
    \leqslant  C_{\delta ,n,(s_{d} ),R,\Rsource,\eAMP} \sum _{\substack{d\in \TZ*\\ \mathclap{d\geqslant (N+1)\pst}}} \Big(\frac{\eAMP r_{2} }{R'} \Big)^{d-n} \, (1+|{\ln r_{1} }|+|{\ln r_{2} }|)^{2n+1} .
\end{multline*}
Let us take $\eAMP:= \frac{R'}{R''}$ with $R'' := \frac{R+R'}{2} \in (R,R')$. Thus $\eAMP>1$ and $\displaystyle \frac{\eAMP r_{2} }{R'} = \frac{r_{2} }{R''} \leqslant  \frac{R}{R''} <1$. Therefore, we~have:\vspace{-0.5ex}
\begin{multline*}
    A \leqslant  C_{\delta ,n,(s_{d} ),R,\Rsource} \sum _{\substack{d\in \TZ*\\ \mathclap{d\geqslant (N+1)\pst}}} \left( \frac{r_{2} }{R''} \right)^{d-n} \, (1+|{\ln r_{1} }|+|{\ln r_{2} }|)^{2n+1} \\[-4ex]
    \begin{array}{r}
        =\displaystyle C_{\delta ,n,(s_{d} ),R,\Rsource} \; \frac{1}{1-\left( \frac{r_{2} }{R''} \right)^{\smallpst} }  \left( \frac{r_{2} }{R''} \right)^{(N+1)\pst-n} \, (1+|{\ln r_{1} }|+|{\ln r_{2} }|)^{2n+1} \\[3ex]
        \leqslant \displaystyle C_{\delta ,n,(s_{d} ),R,\Rsource} \left( \frac{r_{2} }{R} \right)^{(N+1)\pst-n} \, (1+|{\ln r_{1} }|+|{\ln r_{2} }|)^{2n+1}  .
    \end{array}
\end{multline*}
Taking $r_{1} <r_{2}$ and $N\rightarrow \infty $, this proves that $f_{N}$ and $\nabla  f_{N}$ converge uniformly on $D'$. Hence there exists $f\in C^{1} (D)$ s.t.\ $f_{N} \rightarrow f$ and $\nabla  f_{N} \rightarrow \nabla  f$ pointwise on $D$. Finally, taking $r_{1} :=r_{2} :=r$ and any $N\in \mathbb{N}$, we obtain \cref{eq: vitesse de cvg des series de mcl}.
\end{demo}

\begin{textesansboite}
We are ready to prove the first major theorem of this section, which is stated in \cref{sec: intro convergence des mcl et uDA}.
\end{textesansboite}

\begin{demo}[Proof of \cref{decompo modale des CL}]
\phantomsection\label{preuve: decompo modale des CL}%
Let us prove the result by induction of $p\in \NN$. First, let $p=0$. The only $\ell $ s.t.\ $u_{0,\ell } \neq 0$ is $\ell =0$. Since $u_{0,0} \in H^{1} (\Omega )$, \cref{majoration coefs decompo modale Helmholtz homog} shows that there exists $(s_{d} )\in \mathbb{C}^{\TN}$ s.t.\ $s_{d} =O(R^{-d} )$ for~any $R\in (0,\Rsource)$ and $u_{0,0} =\sum _{d\in \TN} s_{d} \,\mcll_{d,0}$ in $B(0,\Rsource)$ with locally uniform convergence. \Cref{convergence des series de mcl} then states that
\[\widetilde{\rho }_{N} := u_{0,0} - \sum _{\substack{d\in \TN\\ d\leqslant N\pst}} s_{d} \,\mcll_{d,0} \quad  \mathrm{satisfies} \quad  
| \widetilde{\rho }_{N} (r,\theta )| + r|\nabla  \widetilde{\rho }_{N} (r,\theta )| \leqslant  C_{(s_{d} ),R} \left( \frac{r}{R} \right)^{(N+1)\pst} (1+|{\ln r}|)^{2} .\]
Let us show that $s_{d} =\sigma _{d} (\uDA_{0,0} )$ for~any $d$.
\begin{sousdemo}
    Let us assume by contradiction that there exists $d\in \TN$ s.t.\ $s_{d} \neq \sigma _{d} (\uDA_{0,0} )$. Let us assume that $d$ is minimal for this property. On the one hand, the estimate on $\widetilde{\rho }_{N}$ implies that 
    \[u_{0,0} = \sum _{d'\leqslant d} s_{d'} \,\mcll_{d',0} + o(r^{d} ) \quad  \text{when }  r\rightarrow 0.\]
    On the other hand, we have
    \begin{alignjustif}
        u_{0,0} &= T_{\leqslant d} (\uDA_{0,0} )+o(r^{d} )\\
        &= \sum _{d'\leqslant d} \sigma _{d'} (\uDA_{0,0} ) \cdot T_{\leqslant d} (\mcla_{d',0} ) + o(r^{d} ) & \text{by \cref{expression uDA avec mcl}} \\
        &= \sum _{d'\leqslant d} \sigma _{d'} (\uDA_{0,0} ) \cdot \mcll_{d',0} + o(r^{d} )
    \end{alignjustif}
    because $\mcll_{d',0} = T_{\leqslant d} (\mcla_{d',0} )+o(r^{d} )$ by \cref{eq: vitesse de cvg des mcl}. Subtracting the two gives $(s_{d} -\sigma _{d} (\uDA_{0,0} )) \,\mcll_{d,0} = o(r^{d} )$ when $r\rightarrow 0$. Now, \cref{expression sing homog CL} states that $\mcll_{d,0}$ is proportional to $J_{d} (k_{0} r) \sin (d\theta )$, which is~not a $o(r^{d} )$ by the properties of Bessel functions, see \cite{Leb72}. Hence $s_{d} -\sigma _{d} (\uDA_{0,0} ) =0$, which is absurd.
\end{sousdemo}
This concludes the initialization.\medskip\smallskip

Let $p\in \NN \setminus \{0\}$ and $\ell \in \mathbb{N}$. Let us assume that the property holds for~any $p'<p$ and let us show it at rank $p$. By~induction hypothesis, we have $\sigma _{d} (\uDA_{p-n,\ell } ) =O(R^{-d} )$ for~any $n\geqslant 1$ and $R\in (0,\Rsource)$. So~\cref{convergence des series de mcl} states that the following series, denoted $f$, converges locally uniformly in $\overline{\Omega }\cap B(0,\Rsource)$: 
\[f := \sum _{n=\red{1}} ^{\lfloor {p}\rfloor } \sum _{\substack{d\in \TZ*\\ d\geqslant -p+n}} \sigma _{d} (\uDA_{p-n,\ell } ) \,\mcll_{d,n} .\]
Note that $n\geqslant 1$ in the above sum, which means that $f$ takes into account the contributions of the previous far fields, but not the $\sigma _{d} (\uDA_{p,\ell } )$. Since $(\Delta +k_{0} ^2 )\mcll_{d,n} =0$ for~any $(d,n)$, we have $(\Delta +k_{0} ^2 )f=0$ in $\mathcal{D}'(\Omega )$. In addition, we have on $\Gamma $:
\begin{alignat*}{2}
    f &= \sum _{n=1} ^{\lfloor {p}\rfloor } \sum _{d\geqslant -p+n} \sigma _{d} (\uDA_{p-n,\ell } ) \,\sum _{m=0} ^{\lfloor {(n-1)/2}\rfloor } \frac{\mu _{0} }{\mu _{1}} T_{m} \cdot  (\partial _{x} ^2 +k_{1} ^2 )^{m} \partial _{y} \mcll_{d,n-1-2m} &\ \ & \text{by \cref{convergence des mcl}} \\ 
    &= \sum _{m=0} ^{\lfloor {p}\rfloor } \frac{\mu _{0} }{\mu _{1}} T_{m} \cdot  (\partial _{x} ^2 +k_{1} ^2 )^{m} \partial _{y} \Bigg( \sum _{n=2m+1} ^{\lfloor {p}\rfloor } \sum _{d\geqslant -p+n} \sigma _{d} (\uDA_{p-n,\ell } ) \, \mcll_{d,n-1-2m} \Bigg) && \text{with convergence in }  \mathcal{D}'(\Gamma )\\
    &= \sum _{m=0} ^{\lfloor {p}\rfloor } \frac{\mu _{0} }{\mu _{1}} T_{m} \cdot  (\partial _{x} ^2 +k_{1} ^2 )^{m} \partial _{y} \Bigg( \sum _{n=0} ^{\lfloor {p-1-2m}\rfloor } \!\!\! \sum _{d\geqslant -p+1+2m+n} \sigma _{d} (\uDA_{p-1-2m-n,\ell } ) \, \mcll_{d,n} \Bigg) && \text{with }  n \leadsto  n-1-2m \\
    &= \sum _{m=0} ^{\lfloor {p}\rfloor } \frac{\mu _{0} }{\mu _{1}} T_{m} \cdot  (\partial _{x} ^2 +k_{1} ^2 )^{m} \partial _{y} u_{p-1-2m,\ell } && \text{by induction hypothesis} 
\end{alignat*}
Hence $f$ satisfies the same Dirichlet condition on $\Gamma $ as $u_{p,\ell }$ (see \cref{construction directe des upl}). Now, let
\begin{equation}{\label{eq: decompo modale des CL preuve def g}}
    g:=u_{p,\ell } -f-\sum\limits _{d\in \red{-\TN},d\geqslant -p} \sigma _{d} (\uDA_{p,\ell } ) \,\mcll_{d,0} .
\end{equation}
We will show that $g$ is equal to $\sum _{d\in \red{\TN}} \sigma _{d} (\uDA_{p,\ell } )\,\mcll_{d,0}$. Using the problems satisfied by $u_{p,\ell }$ (see \cref{construction directe des upl}), $f$ (thanks to the previous steps) and $\mcll_{d,0}$ (see \cref{convergence des mcl}), we~obtain
\[\left\{\begin{array}{r@{\;=\;}ll}
    (\Delta  + k_{0} ^2 )g & 0 &\text{ in }\Omega \cap B(0,\Rsource) \\
    g & 0 &\text{ on }\partial \Omega  \cap B(0,\Rsource)  
\end{array}\right.\]
\def\rsg{R'}%
In order to apply \cref{majoration coefs decompo modale Helmholtz homog} to $g$, let us show that $g\in H^{1} (\Omega \cap B(0,\rsg))$ for~any $\rsg \in (0,\Rsource)$. Note that the sum over $d$ in \cref{eq: decompo modale des CL preuve def g} is the singular behavior at 0 of $u_{p,\ell } -f$ and subtracting it was necessary to make $g$ regular at the corner.
\begin{sousdemo}
    Let $\rsg \in (0,\Rsource)$. By~construction of $u_{p,\ell }$ in \cref{def: construction des termes}, we know that $u_{p,\ell }$ is in $H^{1} (\Omega  \setminus B(0,\eta ))$ for~any $\eta >0$. By~\cref{convergence des mcl}, the same holds for~any $\mcll_{d,n}$. \Cref{convergence des series de mcl} also shows that $f \in H^{1} (\Omega  \cap B(0,\rsg) \setminus B(0,\eta ))$ for~any $\eta >0$. Therefore it only remains to prove that $g$ is $H^{1}$ in the vicinity of the corner.\\
    Let $N\in \mathbb{N}$. \cref{eq: vitesse de cvg des series de mcl} gives 
    \begin{align}
        \Bigg|f(r,\theta ) - \sum _{n=1} ^{\lfloor {p}\rfloor } \sum _{\substack{d\in \TZ*\\ -p+n\leqslant d\leqslant N\pst}} \sigma _{d} (\uDA_{p-n,\ell } ) \,\mcll_{d,n} (r,\theta )\Bigg| 
        &\,\leqslant \, C_{\delta ,n,\fsource,R} \left( \frac{r}{R} \right)^{(N+1)\pst-\lfloor {p}\rfloor } (1+|{\ln r}|)^{2n+1} \nonumber\\[-3ex]
        & \underset{r\rightarrow 0}{\longrightarrow} 0 \label{eq: decompo modale des mcl preuve ineg f}
    \end{align}
    provided that $N$ is chosen s.t.\ $(N+1)\pst-\lfloor {p}\rfloor  >0$. Using $u_{p,\ell } =T_{\leqslant 0} (\uDA_{p,\ell } )+o(1)$ and $\mcll_{d,n} = T_{\leqslant 0} (\mcla_{d,n} )+o(1)$ when $r\rightarrow 0$, we deduce 
    \[g = T_{\leqslant 0} (\uDA_{p,\ell } ) - \sum _{n=1} ^{\lfloor {p}\rfloor } \sum _{\substack{d\in \TZ*\\ -p+n\leqslant d\leqslant N\pst}} \sigma _{d} (\uDA_{p-n,\ell } )\, T_{\leqslant 0} (\mcla_{d,n} ) \,-\! \sum\limits _{d\in -\TN,d\geqslant -p} \sigma _{d} (\uDA_{p,\ell } )\, T_{\leqslant 0} (\mcla_{d,0} ) +o(1).\]
    This is valid because we are summing a finite number of $o(1)$ (this is why approximating $f$ by a finite sum using \cref{eq: decompo modale des mcl preuve ineg f} was necessary). This implies
    \begin{align*}
        g &= T_{\leqslant 0} \bigg( \uDA_{p,\ell } - \sum _{n=1} ^{\lfloor {p}\rfloor } \sum _{-p+n\leqslant d\leqslant N\pst} \sigma _{d} (\uDA_{p-n,\ell } )\, \mcla_{d,n} \,-\! \sum\limits _{d\in -\TN,d\geqslant -p} \sigma _{d} (\uDA_{p,\ell } )\, \mcla_{d,0} \bigg) +o(1)\\
        &= T_{\leqslant 0} \bigg( \sum _{n=1} ^{\lfloor {p}\rfloor } \fsum_{d\geqslant (N+1)\pst} \sigma _{d} (\uDA_{p-n,\ell } )\, \mcla_{d,n} \,+\! \fsum_{d\in \TN} \sigma _{d} (\uDA_{p,\ell } )\, \mcla_{d,0} \bigg) +o(1)
    \end{align*}
    thanks to \cref{expression uDA avec mcl}. Next, for~any $d\geqslant (N+1)\pst$ and $n\geqslant 1$, we~have $d-n\geqslant (N+1)\pst-\lfloor {p}\rfloor  >0$, so $T_{\leqslant 0} (\mcla_{d,n} ) =0$. We also have $T_{\leqslant 0} (\mcla_{d,0} )=0$ for~any $d>0$. Therefore, we~get $g=o(1)$. Similarly, one can show $\nabla g =o(1)$. This proves that $g$ is $H^{1}$ in the vicinity of the corner.
\end{sousdemo}
Similarly to the case $p=0$ (replacing $\Rsource$ by $\rsg$), one can deduce using \cref{majoration coefs decompo modale Helmholtz homog} that $g=\sum _{d\in \TN} \sigma _{d} (\uDA_{p,\ell } )\,\mcll_{d,0}$ and $\sigma _{d} (\uDA_{p,\ell } )=O(R^{-d} )$ for~any $R\in (0,\Rsource)$. This proves \cref{eq: decompo modale des CL}. Then, \cref{eq: vitesse de cvg decompo modale des CL} follows from \cref{eq: vitesse de cvg des series de mcl}. Finally, the convergence is uniform on $\Omega  \cap B(0,R)$ for~any $R\in (0,\Rsource)$ because, with $R'':= (R+\Rsource)/2$, \cref{eq: vitesse de cvg decompo modale des CL} implies
\[|\rho _{N} (r,\theta )| \leqslant  C_{p,\fsource,R} \left( \frac{r}{R''} \right)^{(N+1)\pst-\lfloor {p}\rfloor } (1+|{\ln r}|)^{2p+1} 
\leqslant  C_{p,\fsource,R} \left( \frac{r}{R''} \right)^{N\pst-\lfloor {p}\rfloor } 
\leqslant  C_{p,\fsource,R} \left( \frac{R}{R''} \right)^{N\pst}.\]
\end{demo}

\begin{remarque}
\Cref{convergence des mcl} states that $\mcla_{d,n}$ converges on $\Omega $ entirely. So what restricts the domain of convergence of the modal decomposition in \cref{decompo modale des CL} is the behavior of the coefficients $\sigma _{d} (u_{p-n,\ell } )$ when $d\rightarrow \infty $.
\end{remarque}

\subsection{Proof of \texorpdfstring{\cref{decompo en modes tronques des CL,convergence de uDA}}{Theorem 7.1.3 and Corollary 7.1.4}}
\label{sec: preuve decompo modale CL 2}

\begin{demo}[Proof of \cref{decompo en modes tronques des CL}]
\phantomsection\label{preuve: decompo en modes tronques des CL}%
\def\minnq{m_{N,q}}%
\def\nq{M}%
\def\rnq{\rho_{N,q}}%
\def\termea{f_{\nq,q}}%
\def\termeb{g_{\nq,q}}%
\def\termebp{g_{\nq,-p}}%
\def\domd{ -p+n\leqslant d\leqslant \nq\pst }%
To begin, note that, if $N$ is very large compared to $q$, then a lot of the truncated modes $T_{<q} (\mcla_{d,n} )$ in \cref{eq: reste decompo en modes tronques des CL} vanish. Let us eliminate them by defining $\nq$ as the greatest integer in $(-\infty ,N]$ s.t.\ $\nq\pst-\lfloor {p}\rfloor \leqslant q$. Thus, for~any $n\in [\![0,\lfloor {p}\rfloor ]\!]$ and $d\in \TZ* \cap  (\nq\pst,N\pst]$ (which may be empty) we have $d-n\geqslant (\nq+1)\pst -\lfloor {p}\rfloor >q$, so $T_{<q} (\mcla_{d,n} ) = 0$. Therefore, we have
\[\rnq = \rho _{\nq,q} .\]
Next, we can decompose the remainder as
\[\rnq \;=\; \underbrace{u_{p,\ell } - \sum _{n=0} ^{\lfloor {p}\rfloor } \! \sum _{\substack{d\in \TZ*\\ \domd}} \!\sigma _{d} (\uDA_{p-n,\ell } ) \,\mcll_{d,n} }_{\textstyle :=\termea} \;-\, \underbrace{\sum _{n=0} ^{\lfloor {p}\rfloor } \! \sum _{\substack{d\in \TZ*\\ \domd}} \! \sigma _{d} (\uDA_{p-n,\ell } ) \,\restmcl }_{\textstyle :=\termeb} .\]
with $\restmcl := \mcll_{d,n} -T_{<q} (\mcla_{d,n} )$ for~any $(d,n,q)$. Let $\minnq:=\min ((N+1)\pst-\lfloor {p}\rfloor ,q)$. We will bound $\termea$ and $\termeb$ by $C_{p,\fsource,R} \left( \frac{r}{R} \right)^{\minnq} (1+|{\ln r}|)^{2p+1}$ to prove the theorem.\\
Let us start with $\termea$. Thanks to \cref{decompo modale des CL}, we know that
\[|\termea(r,\theta )|+r|\nabla  \termea(r,\theta )| \leqslant  C_{p,\fsource,R} \left( \frac{r}{R} \right)^{(\nq+1)\pst-\lfloor {p}\rfloor } (1+|{\ln r}|)^{2p+1} .\]
In addition, by~definition of $\nq$, we have either $\nq=N$ or $(\nq+1)\pst-\lfloor {p}\rfloor >q$. In both cases, $(\nq+1)\pst-\lfloor {p}\rfloor \geqslant  \minnq$, which gives the desired bound on $\termea$.\smallskip\medskip

Let us now treat $\termeb$ and prove more precisely
\begin{equation}{\label{eq: preuve decompo en modes tronques des CL} }
    |\termeb(r,\theta )|+r|\nabla  \termeb(r,\theta )| \leqslant  C_{p,\fsource,R} \left( \frac{r}{R} \right)^{q} (1+|{\ln r}|)^{2p+1} .
\end{equation}
Let $\delta :=-p$ and $\eAMP'>\eAMP>\eAMPmin{\delta ,n}$ (to be choosen later). Let $R' :=(R+\Rsource)/2$. By~\cref{decompo modale des CL,eq: vitesse de cvg des mcl}, we have for~any $q\geqslant \delta $
\begin{align*}
    |\termeb (r,\theta )|+r|\nabla  \termeb (r,\theta )|
    &\leqslant  \sum _{n=0} ^{\lfloor {p}\rfloor } \sum _{\domd} |\sigma _{d} (\uDA_{p-n,\ell } )| \cdot \big(|\restmcl(r,\theta )| + r|\nabla  \restmcl(r,\theta )|\big) \\
    &\leqslant  \sum _{n=0} ^{\lfloor {p}\rfloor } \sum _{\domd} C_{p,\fsource,R} \, R'^{-d} \cdot  C_{\delta ,n,\eAMP,\eAMP'} \, \frac{\gamma (d)}{\gamma (q)} \,\eAMP^{-d} \,(\eAMP' r)^{q} \, (1+|{\ln r}|)^{2n+1} \\
    &\leqslant  C_{p,\fsource,R,\eAMP,\eAMP'} \, \frac{(\eAMP' r)^{q} }{\gamma (q)} \, (1+|{\ln r}|)^{2p+1} \sum _{n=0} ^{\lfloor {p}\rfloor } \sum _{\domd} (\eAMP R')^{-d} \, \gamma (d) \\
    &\leqslant  C_{p,\fsource,R,\eAMP,\eAMP'} \, \frac{(\eAMP' r)^{q} }{\gamma (q)} \, (1+|{\ln r}|)^{2p+1} \underbrace{\sum _{-p\leqslant d\leqslant \nq\pst} (\eAMP R')^{-d} \, \gamma (d) }_{:=A}
\end{align*}
Let us show that $\displaystyle A\leqslant  C_{p,\fsource,R,\eAMP,\eAMP''} \left( \frac{\eAMP''}{\eAMP R'}\right)^{q} \gamma (q)$ for~any $q\geqslant -p$ and $\eAMP''>1$.
\begin{sousdemo}
    \def\qp{q'}%
    The proof of \cref{cvg des SF de norme finie} showed that we have  \[\frac{(\eAMP R')^{q'} }{\gamma (q')} \leqslant  C_{p,R',\eAMP} \frac{(\eAMP R')^{d} }{\gamma (d)}\]
    for~any $q'\geqslant d\geqslant -p$. Since $\nq\pst\leqslant q+\lfloor {p}\rfloor $ by~definition of $\nq$, we can applying this to $d\in [-p,\nq\pst]$ and $q':=q+\lfloor {p}\rfloor $, which gives:
    \[\textstyle A \leqslant  C_{p,\fsource,R,\eAMP} \cdot  \operatorname{card}(\TZ* \cap [-p,\nq\pst]) \cdot  (\eAMP R')^{-q'} \gamma (q')\] 
    Next, the above cardinal is bounded by $C\cdot (\nq\pst+p) \leqslant C\cdot (q+\lfloor {p}\rfloor +p)\leqslant  C_{p} \, q$ because $q\geqslant -p$. 
    Then, $q\cdot  (\eAMP R')^{-\qp}$ can be bounded by $C_{p,\fsource,R,\eAMP,\eAMP''} \cdot  \eAMP''^{q} \cdot  (\eAMP R')^{-q}$. Finally, using again \cref{eq: majoration de quotients des fonctions gamma} from the proof of \cref{norme SF de mcl cas general} to bound $\gamma (q+\lfloor {p}\rfloor ) / \gamma (q)$ by $C_{p,\eAMP''} \eAMP''^{q}$, we deduce the desired estimate of $A$.
\end{sousdemo} 
Therefore we have shown
\begin{align*}
    |\termeb (r,\theta )|+r|\nabla  \termeb (r,\theta )|
    &\leqslant  C_{p,\fsource,R,\eAMP,\eAMP',\eAMP''} \left( \frac{\eAMP' \eAMP'' r}{\eAMP R'} \right)^{q} \, (1+|{\ln r}|)^{2p+1}
\end{align*}
\def\qptt{q_0}%
Let us take $\eAMP:=\eAMPmin{\delta ,n}+1$ and $\eAMP':=\eAMP\cdot \eAMP''$, and choose $\eAMP''$ s.t.\ $\frac{\eAMP''^2 }{R'} = \frac{1}{R}$. This proves \cref{eq: preuve decompo en modes tronques des CL} when $q\geqslant -p$.\\
Finally, for~any $q<-p$, we have $\termeb = \termebp$ because $T_{<q} (\mcla_{d,n} )=0=T_{<-p} (\mcla_{d,n} )$ for~any $n\in [\![0,\lfloor {p}\rfloor ]\!]$ and $d\geqslant -p+n$. Using \cref{eq: preuve decompo en modes tronques des CL} on $\termebp$, we obtain
\[|\termeb (r,\theta )|+r|\nabla  \termeb (r,\theta )| \leqslant  C_{p,\fsource,R} \left( \frac{r}{R} \right)^{-p} (1+|{\ln r}|)^{2p+1} \leqslant  C_{p,\fsource,R} \left( \frac{r}{R} \right)^{q} (1+|{\ln r}|)^{2p+1} .\]
because $r/R \in (0,1)$. This concludes the proof.
\end{demo}

\begin{demo}[Proof of \cref{convergence de uDA}]
\phantomsection\label{preuve: convergence de uDA}%
Let us write $T_{<q} (\uDA_{p,\ell } )$ as a truncation of the modal decompostion of $u_{p,\ell }$ in order to use \cref{decompo en modes tronques des CL}. We have:
\begin{alignjustif}
T_{<q} (\uDA_{p,\ell } ) &= T_{<q} \bigg( \sum _{n=0} ^{\lfloor {p}\rfloor } \fsum_{\substack{d\in \TZ*\\ d\geqslant -p+n }} \!\sigma _{d} (\uDA_{p-n,\ell } ) \,\mcla_{d,n} \bigg) & \text{by \cref{expression uDA avec mcl}.} \\
&= \sum _{n=0} ^{\lfloor {p}\rfloor } \fsum_{\substack{d\in \TZ*\\ d\geqslant -p+n }} \!\sigma _{d} (\uDA_{p-n,\ell } ) \cdot  T_{<q} (\mcla_{d,n} ) & \text{by \cref{rq: fsum commute avec les troncatures T}} 
\end{alignjustif}
Let $N$ be the smallest integer s.t.\ $(N+1)\pst -\lfloor {p}\rfloor \geqslant q$. Then, for~any $n\leqslant \lfloor {p}\rfloor $ and $d\geqslant (N+1)\pst$, we have $d-n\geqslant q$, so $T_{<q} (\mcla_{d,n} )=0$. We deduce:
\[T_{<q} (\uDA_{p,\ell } ) 
= \sum _{n=0} ^{\lfloor {p}\rfloor } \sum _{\substack{d\in \TZ*\\ -p+n\leqslant d\leqslant N\pst }} \!\sigma _{d} (\uDA_{p-n,\ell } ) \cdot  T_{<q} (\mcla_{d,n} )\]
because any finite $\fsum$-sum can be written as a classical sum (see \cref{proprietes de fsum}). Hence, according to \cref{decompo en modes tronques des CL}, we have
\[|\rho _{q} (r,\theta )|+r|\nabla  \rho _{q} (r,\theta )| \leqslant  C_{p,\fsource,R} \left( \frac{r}{R} \right)^{\min \left((N+1)\pst -\lfloor {p}\rfloor ,\, q\right)} (1+|{\ln r}|)^{2p+1} = C_{p,\fsource,R} \left( \frac{r}{R} \right)^{q} (1+|{\ln r}|)^{2p+1}\]
which proves \cref{eq: vitesse de cvg de uDA}. Similarly to the end of the proof of \cref{decompo modale des CL}, this implies $T_{<q} (\uDA_{p,\ell } )\rightarrow  u_{p,\ell }$ uniformly on $\Omega \cap B(0,R)$ for~any $R\in (0,\Rsource)$.
\end{demo}

\chapter{Numerical validation of the asymptotic expansion}
\label{chap: validation numerique du DA}
\begin{textesansboite}
In this chapter, we provide numerical results regarding the asymptotic expansion of $u_{\varepsilon }$ given in \cref{chap: construction du DA}, in order to illustrate and validate our theoretical results. We detail the numerical methods chosen, in particular the one used to compute the corner profiles, as it~is novel.\\

Let us recall and summarise the construction algorithm. See \cref{sec: construction directe des upl} for more details. We can divide the construction into 4 parts.
\begin{enumerate}
	\item The first one is the computation of the corner profiles $\mathcal{S}_{d,n}$. We recall that we denote $V:=\{v\in H^{1} _{\mathrm{loc}} (\Omega _{1} ) \mid \allowbreak \nabla u\in L^2  \text{ and } u_{|\partial \Omega _{1}} =0\}$ (\cref{def: V}) and that, for~any fixed $d\in \TN$, the corner profiles are defined inductively w.r.t.\ $n$ as the unique solutions in $\Vcoin+\khii \mathcal{A}(\Pi )$  of:
	\[\left\{\begin{array}{r@{\;=\;}l@{\;}l}
		\operatorname{div}(\mu  \nabla  \mathcal{S}_{d,0} ) & 0 &\text{ in }\Omega _{1} \\[0.5mm]
		\mathcal{S}_{d,0} & 0 &\text{ on }\partial \Omega _{1} \\[0.5mm]
		\sigma _{q} (\mathcal{S}_{d,0} )& \delta _{q,d} & \ \forall q\in \TN
	\end{array}\right. \quad  \text{ and } \quad  \ \forall n\in \mathbb{N}^*,\; \left\{\begin{array}{r@{\;=\;}l@{\;}l}
		\operatorname{div}(\mu  \nabla  \mathcal{S}_{d,n} ) & -\rho \,\mathcal{S}_{d,n-1} &\text{ in }\Omega _{1} \\[0.5mm]
		\mathcal{S}_{d,n} & 0 &\text{ on }\partial \Omega _{1} \\[0.5mm]
		\sigma _{q} (\mathcal{S}_{d,n} )& 0 & \ \forall q\in \TN
	\end{array}\right.\]
	For different values of $d$, the sequences $(\mathcal{S}_{d,n} )_{n}$ are independent from one another. From the corner profiles, we only keep the coefficients $\sigma _{q} (\ccDA_{d,n} )$ with $q\in -\TN$ for the rest of the computations. Computing these coefficients is difficult as they describe the behavior at infinity of the profiles, while we need to truncate the domain for the numerical computation. Hence a specific numerical method is required.
	\item The second part is the computation of the matching coefficients $\cmicro_{d,d',p,\ell }$ and $\cmacro_{d,d',p,\ell }$. It~is independent from part~1. These coefficients can be computed exactly and very quickly through algebraic calculations because they are equal to $\sigma _{d,p,\ell } \circ \Eps^{\pm } \circ \langle \Rmaicro\rangle (\phi_{d'} )$. The implementation of this computation is explained in \cref{sec: implementation des espaces A}.
	\item Then, one must compute the coefficients $\cfinal_{d,d',p,\ell }$. They depend on the results obtained in parts~1 and 2 through the following formula
	\[\cfinal_{d,d'\!,p,\ell } := \! \sum _{\substack{(p_{1} ,p_{2} )\in (\mathbb{N}-d)\times (\mathbb{N}+d')\\ n\in \mathbb{N},\, p_{1} +p_{2} +2n=p}}
	\ \sum _{\substack{d_{1} \in \TZ*\\ p_{1} +d_{1} \in \mathbb{N}}} \ \sum _{\substack{d_{2} \in \TN\\ p_{2} -d_{2} \in \mathbb{N}}}
	\ \sum _{\substack{(\ell _{1} ,\ell _{2} )\in \mathbb{N}^2 \\ \ell _{1} +\ell _{2} =\ell }}
	\cmacro_{d,d_{1} ,p_{1} ,\ell _{1} }\cdot \omega ^{2n} \sigma _{d_{1}} (\mathcal{S} ^{\infty } _{d_{2} ,n} )\cdot \cmicro_{d_{2} ,d'\!,p_{2} ,\ell _{2} } .\]
	Note that the above sum is actually finite and can be simplified depending on the value of $\pst$ (see \cref{sec: simplifications DA}).
	\item The final step is the computation of the far fields $u_{p,\ell }$, inductively w.r.t.\ $p$. Each far field is solution of the following problem.
	\[\left\{\begin{array}{r@{\;=\;}ll@{}}
		\mu _{0} \Delta  u_{p,\ell } +\omega ^2 \rho _{0} u_{p,\ell } & \fsource \,\delta _{p,0} \,\delta _{\ell ,0} &\text{ in }\Omega  \\[0.5mm]
		u_{p,\ell } & 0 &\text{ on }\bordgauche\\[0.5mm]
		u_{p,\ell } &\displaystyle \sum_{n=0}^\infty  \frac{\mu _{0} }{\mu _{1}} \,T_{n} \cdot  (\partial _{x} ^2 +k_{1} ^2 )^{n} \partial _{y} u_{p-1-2n,\ell } &\text{ on }\Gamma  \\[-1.5mm]
		\sigma _{d} (u_{p,\ell } )&\displaystyle \!\! \sum _{\substack{p'\in \NN\\ p-p'\in \NN+ 2\pst }} \; \sum _{\substack{d'\in \TZ*\\ p'+d'\in \intvl{0,p+d}}} \! \sum _{\ell '=0} ^{\ell } \cfinal_{d,d',p-p',\ell -\ell '} \cdot \sigma _{d'} (\uDA_{p',\ell '} ) & \ \forall d\in -\TN
	\end{array}\right.\]
	Note that the above sums contain a finite number of non-zero terms. One needs also to compute some of the coefficients $\sigma _{d} (\uDA_{p,\ell } )$ with $d\in \TN$, as they are involved in the construction of other far fields (when $d\in -\TN$, $\sigma _{d} (\uDA_{p,\ell } )=\sigma _{d} (u_{p,\ell } )$ is given by the last equation in the system above). Thus, similarly to the corner profiles, an appropriate method is needed to compute the far fields. Additional difficulties arise from the singular behavior at the corner of $u_{p,\ell }$, and its Dirichlet condition on $\Gamma $ as it involves high order derivatives of the previous far fields. 
\end{enumerate}
Only steps 1 and 4 require heavy computations as they involve solving a PDE with no explicit solution. Moreover one needs to extract accurately enough asymptotic coefficients of the PDE's solutions, so special numerical methods are needed. \Cref{sec: calcul numerique des CC} introduces a numerical method to compute the corner profiles, based on the use of complements at infinity. \cref{sec: numerique des CL} explains the computation of the far fields and gives convergence plots for the asymptotic expansion of $u_{\varepsilon }$ that validate the theoretical error estimates. 
\end{textesansboite}

\begin{remarquenumerotee}
\label{eq: parallelisation des calculs}%
The computations can be parallelized. Indeed, the corner profiles $\mathcal{S}_{d,n}$ with different values of $d$ can be computed independently. Moreover, although in principle computing the far fields requires to know the coefficients $\sigma _{q} (\mathcal{S} ^{\infty } _{d,n} )$, in practice the $\sigma _{q} (\ccDA_{d,n} )$ with a large $n$ are only needed by the $u_{p,\ell }$ with a large $p$. One can check that the first far fields (with $p<2\pst$) do~not even depend on the $\sigma _{q} (\ccDA_{d,n} )$. Hence the computation of the corner profiles (which is inductive w.r.t.\ $n$) and the far fields (which is inductive w.r.t.\ $p$) can be run in parallel
\ (although the construction of the far fields may sometimes have to wait for some coefficients $\sigma _{q} (\ccDA_{d,n} )$). We implemented this parallelization (more details in \cref{sec: resultats numeriques CL}).
\end{remarquenumerotee}

\begin{remarque}
Steps 1 to 3 remain for the effective models of \cref{chap: ME}. In particular, the numerical method used for the corner profiles remains relevant for the effective models.
\end{remarque}

\section{Complements at infinity: a numerical method to compute the corner profiles}
\label{sec: calcul numerique des CC}

\begin{textesansboite}
In this section, we propose a numerical method to compute the corner profiles $\mathcal{S}_{d,n}$ and their asymptotic coefficients $\sigma _{q} (\ccDA_{d,n} )$. We prove convergence estimates and we provide numerical results to validate the theoretical results.
\end{textesansboite}

\subsection{Model problem and state of the art}
\label{sec: numerique CC pb modele et etat de l'art}

\begin{paragraphesansboite}{The non-variational problem:}
\def\douq{q}
All the corner profiles are solutions of problems of the form
\begin{equation}{\label{eq: model pb for corner fields}}
\left\{\begin{array}{r@{\;=\;}ll}
	\operatorname{div}(\mu  \nabla  S) & F &\text{ in }\Omega _{1} \\[0.5mm]
	S & 0 &\text{ on }\partial \Omega _{1} \\[0.5mm]
	\sigma _{\douq} (S)& s_{\douq} & \ \forall \douq \in \TN
\end{array}\right.\end{equation}
with $S\in \Vcoin+\khii \mathcal{A}(\Pi )$, $F\in L^2 _{\mathrm{loc}} (\overline{\Omega _{1} })$ (i.e.\ $L^2 $ on any bounded subset of $\Omega _{1}$) and $(s_{\douq} )_{\douq} \in \mathbb{C}^{\TN}$ with a finite support. Moreover, $F$ satisfies the following assumptions:
\begin{itemize}
	\item there exists $F^{\infty } _{\Omega } \in \Amicro(\Omega )$ s.t.\ $\forall \douq \in \mathbb{R},\; F_{|\Omega } = T_{\geqslant \douq} (F^{\infty } _{\Omega } ) + o_{\partial } (r^{\douq} )$ when $r\rightarrow \infty $,
	\item there exists $F^{\infty } _{\couche} \in \Amicro(\couche)$ s.t.\ $\forall \douq \in \mathbb{R},\; F_{|\couche} = T_{\geqslant \douq} (F^{\infty } _{\couche} ) + o_{\partial } (x^{\douq} )$ when $x\rightarrow \infty $.
\end{itemize}
Using \cref{cadre fonctionnel pour les champs de coin}, one can show that \cref{eq: model pb for corner fields} has a unique solution (this is similar to the proof of existence and uniqueness of the corner fields done in \cref{sec: construction des termes}). To compute the corner profiles, it suffices to be able to solve numerically the model problem \cref{eq: model pb for corner fields}.\\

Let us denote
\[\rmicro := \fsum \left( -\RcouchePi\circ \partial _{x|\couche} ^2  + \frac{\mu _{0} }{\mu _{1}} \RGammaPi\circ \partial _{Y|\Gamma ,Y=0^{+}} \right)^{n} .\]
This is an operator $\Amicro(\Pi )\rightarrow \Amicro(\Pi )$. \Cref{DA des champs de coin} states that there exists $\SDA \in \Amicro(\Pi )$ that has the following form 
\begin{equation}{\label{eq: SDA compl sing}}
	\SDA = \rmicro\bigg( \frac{1}{\mu _{0}} \ROmegaPi (F^{\infty } _{\Omega } )+ \frac{1}{\mu _{1}} \RcouchePi (F^{\infty } _{\couche} ) + \!\fsum_{\douq \in \TZ*} \! \sigma _{\douq} (\SDA) \,\phi_{\douq} \bigg)
\end{equation}
such that: \ $\forall \douq\in \mathbb{R},\ S = T_{\geqslant \douq} (\SDA ) + o_{\partial } (\bfr^{\douq} )$ when $\bfr\rightarrow \infty $. For~any $\douq\in \TN$, $\sigma _{\douq} (\SDA)$ is known because it~is equal to $\sigma _{\douq} (S)$, i.e.\ to $s_{\douq}$. It belongs to the non-variational behavior of $S$. But the $\sigma _{\douq} (\SDA)$ with $\douq\in -\TN$ are unknown and belong to the non-variational behavior of $S$. We are interested in a method to compute these coefficients.\\

We need to answer two questions:
\begin{itemize}
	\item Question 1: how to compute numerically a good approximation of $S$ ?
	\item Question 2: how to compute numerically a good approximation of $\sigma _{\douq} (\SDA)$ when $\douq\in -\TN$ ?
\end{itemize}
Although in the end we are only interested in the coefficients $\sigma _{q} (\ccDA_{d,n} )$, the inductive construction of the corner profiles through the relation $\operatorname{div}(\mu  \nabla  \mathcal{S}_{d,n} )=-\rho \,\mathcal{S}_{d,n-1}$ makes it also necessary to compute a good approximation of $\mathcal{S}_{d,n-1}$. So we must compute accurately both $S$ and the $\sigma _{\douq} (\SDA)$.
\end{paragraphesansboite}

\begin{paragraphesansboite}{The variational problem:}
The first step is to rewrite the problem as a variational one. Let $\chi \in \mathcal{C}^{\infty } (\mathbb{R}^2 )$ be a truncation function that is equal to 0 in $B(0,\Rcoin)$ and to 1 in a vicinity of infinity. Let $\beta \leqslant -1$ and 
\begin{equation}{\label{eq: def S survaria et S varia}}
	S^{\mathcal{A}} := T_{\geqslant \beta } \circ \rmicro\bigg( \frac{1}{\mu _{0}} \ROmegaPi (F^{\infty } _{\Omega } )+ \frac{1}{\mu _{1}} \RcouchePi (F^{\infty } _{\couche} ) + \!\fsum_{\red{d\in \TN}} \! s_{d} \,\phi_{d} \bigg) \in \mathcal{A}(\Pi )
\quad  \text{ and } \quad  \SV := S - \chi \,S^{\mathcal{A}} .
\end{equation}
$\beta $ will be chosen later (see \cref{estimate of uL-SV}). $\SV$ satisfies
\begin{equation}{\label{eq: eq of the variational part for corner fields}}
\left\{\begin{array}{r@{\;=\;}ll}
	\operatorname{div}(\mu  \nabla  \SV) &f := F-\operatorname{div}[\mu  \nabla (\chi  S^{\mathcal{A}} )]&\text{ in }\Omega _{1} \setminus (\Gamma \cap \{x>\Rcoin\})  \\[0.5mm]
	\SV &0 &\text{ on }\partial \Omega _{1} \\[0.5mm]
	\SV_{|Y=0^{+}} - \SV_{|Y=0^{-}} & 0 &\text{ on }\Gamma \cap \{x>\Rcoin\} \\[0.5mm]
	[\mu  \partial _{y} \SV] & g := -[\mu  \partial _{y} (\chi  S^{\mathcal{A}} )] &\text{ on }\Gamma \cap \{x>\Rcoin\} 
\end{array}\right.\end{equation}
with $[\mu \partial _{y} \dots ] :=\mu _{0} \partial _{y|y=0^{+}} - \mu _{1} \partial _{y|y=0^{-}}$. Similarly to the proof of \cref{cadre fonctionnel pour les champs de coin}, we have $\SV \in V$, $(1+r)f\in L^2 (\Omega _{1} )$ and $g\in L^2 (\Gamma \cap \{x>\Rcoin\})$, and $\SV$ is the unique solution of \cref{eq: eq of the variational part for corner fields} in $V$. The couple $(\SV,S^{\mathcal{A}} )$ is a decomposition of $S$ into a variational part and a diverging part.\\

In~practice, $F$ is $-\rho $ times some previously computed corner profile, say $\ccDA_{d,n-1}$. The computation of the first coefficients $\sigma _{q} (\ccDA_{d,n-1} )$, provides knowledge of $F^{\infty } _{\Omega }$ and $F^{\infty } _{\couche}$ up to a given order (i.e.\ we know $T_{\geqslant q} (F^{\infty } _{\Omega } )$ and $T_{\geqslant q} (F^{\infty } _{\couche} )$ for some $q$). If enough coefficients $\sigma _{q} (\ccDA_{d,n-1} )$ have been computed, then enough terms of $F^{\infty } _{\Omega }$ and $F^{\infty } _{\couche}$ are known so~that $S^{\mathcal{A}}$ is a known data. Let us assume that we are in this situation. The number of coefficients $\sigma _{q} (\ccDA_{d,n-1} )$ needed depends on $\beta $. Note that $S^{\mathcal{A}}$ has an explicit expression, thanks to the explicit expressions of $\ROmegaPi$, $\RGammaPi$, $\RcouchePi$ (see \cref{operateurs R}). Therefore, to compute $S$ it suffices to compute $\SV$.
\end{paragraphesansboite}

\begin{paragraphesansboite}{Existing truncation methods:}
Since $\SV$ is defined on an unbounded domain, it~is natural to introduce a truncated domain. Let us recall that the variable $\bfr$ is equal to $r$ in $\Omega $ and $x$ in $\Lambda $. For~any $L>\Rcoin+1$, we define the truncated domain $\ounl:= \Omega _{1} \cap \{\bfr<L\}$.\\

One must choose a boundary condition to apply on the truncated edge. A homogeneous Dirichlet condition is a possibility that provides convergence of the approximated solution towards $\SV$, because $\SV \rightarrow  0$ when $\bfr \rightarrow 0$. However, $\SV$ may converge to 0 slowly because its behavior at infinity involves $\phi_{-\smallpst}$, which behaves like $r^{-\smallpst}$, with $\pst \in ( \frac{1}{2} ,\infty )$. Thus, the approximated solution may converge slowly to $\SV$. A homogeneous Neumann condition leads to the same issue. It~is thus not easy to answer Question 1 of this section.\\

The truncation condition can be improved using the known expansion of $\SV$ at infinity (which is given by $\SDA-S^{\mathcal{A}}$). For instance, \cite{Auv18} achieved a better rate of convergence by using a Robin condition in a similar setting with $\Theta =\pi $. One can check that there~is a known function $\varphi  \in \mathcal{A}(\Pi )$ (depending on the data $F^{\infty } _{\Omega }$, $F^{\infty } _{\couche}$ and $(s_{d} )$) s.t.
\[\SV = \varphi  +\sigma _{-\pst} (\SDA)\,\phi_{-\pst} + o_{\partial } (\bfr^{-\pst} ).\]
However $\sigma _{-\pst} (\SDA)$ is unknown. Since we have $\partial _{n} \phi_{-\pst} + \frac{\pi }{\Theta L} \phi_{-\pst} = 0$ on $\partial \ounl \setminus \partial \Omega _{1}$, the condition
\[\textstyle \partial _{n} u + \frac{\pi }{\Theta L} u = \partial _{n} \varphi  + \frac{\pi }{\Theta L} \varphi  \quad \text{ on }\partial \ounl \setminus \partial \Omega _{1} \]
in an approximated transparent boundary condition with a known right hand side. One can show that the associated approximated solution converges to $\SV$ faster than for Dirichlet and Neumann conditions. Then \cite{Auv18} computes $\sigma _{-\pst} (\SDA)$ a~posteriori in the following way. Using the calculation below
\[\sigma _{-\pst} (\SDA) = 2r^{\pst} \int _{0} ^{\Theta } \sigma _{-\pst} (\SDA) \,\phi_{-\pst} (r,\theta ) \cdot \sin \!\left(\pst \theta \right) \mathrm{d}\theta  = \lim _{r\rightarrow \infty } 2r^{\pst} \int _{0} ^{\Theta } (\SV -\varphi )(r,\theta ) \cdot \sin \!\left(\pst \theta \right) \mathrm{d}\theta  ,\]
one deduces the following approximated formula for $\sigma _{-\pst} (\SDA)$ (where $\uL$ denotes the approximated solution)
\[\sigma _{-\pst} (\SDA) \approx  2L^{\pst} \int _{0} ^{\Theta } \big(\uL(L,\theta ) -\varphi (L,\theta )\big) \sin \!\left(\pst \theta \right) \mathrm{d}\theta  .\]\fbeq
\end{paragraphesansboite}

\begin{paragraphesansboite}{The singular function method:}
The approach of \cite{Auv18} works fine to obtain $\sigma _{-\pst} (\SDA)$. However to obtain more coefficients, that describe a finer behavior of $\SV$ at infinity, one needs a finer truncation boundary condition. So we chose a different approach, closer to the method of singular functions, see \cite[chap.\,8]{StrFix73} and \cite{Loh98, LubNic94, CiaHe03}. 
This method is used on elliptic problems in a domain with a corner. In a neighborhood of the corner, the solution is not very regular, e.g.\ it~is~not $H^{2}$, which deteriorates the convergence of the finite element method w.r.t.\ the mesh step. This difficulty is overcome by adding to the Galerkin space a ``singular function'' $\varphi $ s.t.\ the solution has the form $c\,\varphi +v$ with $c\in \mathbb{C}$ and $v$ a regular function (e.g.\ $H^2 $). In this method, $c$ is unknown, but not $\varphi $. For instance, if the problem is the Poisson equation with homogeneous Dirichlet conditions and the corner has angle $\Theta $, then $\varphi  =\chi \,\phi_{\pst}$, where $\chi $ is a smooth truncation function supported in a neighborhood of the corner. Using this method, the rate of convergence of the finite elements is limited only by the regularity of $v$, and not by the poor regularity of $\varphi $.\\

In our case, we are~not fighting against a lack of regularity at the corner but rather against the slow decay of $\SV$ at infinity. There is a kind of correspondence between the two as the regularity at 0 of $x\mapsto x^{a}$ and the decay at infinity of $x\mapsto  ( \frac{1}{x} )^{a}$ both get better for $a$ large and worse for $a$ small ($a\in \mathbb{R} \setminus \mathbb{Z}$). To solve our problem \cref{eq: eq of the variational part for corner fields}, we add to the Galerkin space functions located near infinity rather than near the corner. We call these functions ``complements at infinity''. An important advantage of this approach is that we can use as many complements at infinity as we want, providing a high order method (although it seems that usually only one singular function is used).
\end{paragraphesansboite}

\subsection{Presentation of the numerical method}

\begin{paragraphesansboite}{Space of approximation:}
Let $L>\Rcoin+1$. We define the truncated domains $\Omega ^{L} := \Omega \cap \{r<L\}$, $\Lambda ^{L} := \Lambda \cap \{x<L\}$, $\Gamma ^{L} := \Gamma \cap \{x<L\}$ and $\ounl:= \Omega _{1} \cap \{\bfr<L\} =\Omega ^{L} \cup \Gamma ^{L} \cup \layer^{L}$. Let $\nbcs \in \mathbb{N}^*$, $\gamma \in \mathbb{R}$ and
\begin{equation}{\label{eq: def complement singulier}}
	\fcs := T_{\geqslant \gamma } \circ \rmicro(\phi_{-n\smallpst} )\in \mathcal{A}(\Pi ) \qquad  \text{for~any }  n\in [\![1,\nbcs]\!].
\end{equation}
The parameter $\gamma $ will be chosen later in terms of $\nbcs$ (see \cref{estimate of uL-SV}). Let $\widetilde{\chi } \in \mathcal{C}^{\infty } (\mathbb{R})$ be a truncation function that is equal to 0 on $(-\infty ,-1]$ and 1 on $[0,\infty )$, and $\khiL \in \mathcal{C}^{\infty } (\Omega _{1} )$ be equal to $\widetilde{\chi }(\bfr-L)$. The approximation space is defined as
\begin{equation}{\label{eq: def evcs}}
	\evcs := \{u\in V \mid u_{|\Omega _{1} \setminus \ounl} =0\} \oplus  \operatorname{Span}(\{ \khiL \fcs \mid n\in [\![1,\nbcs]\!]\}).
\end{equation}
It~is included in $V$ (see \cref{def: V}). Note that the first term for the direct sum is equal to the set of $H^{1} _{0} (\ounl)$-functions extended by 0 in $\Omega _{1} \setminus \ounl$. In practice, it~is discretised with finite elements, but we do~not consider this discretisation in our theoretical analysis. Besides, the functions $\khiL \fcs$ are what we call the \emph{complements at infinity}.\\

Moreover, for~any $v\in \evcs$, we denote $(\pef(v),\pcs{1}(v),\dots ,\pcs{\nbcs}(v))$ the only element of $V \times \mathbb{R}^{\nbcs}$ s.t.\ $\pef(v)_{|\Omega _{1} \setminus \ounl} =0$ and
\begin{equation}{\label{eq: decomposition of elements of evcs}}
	v=\pef(v)+\sum _{n=1} ^{\nbcs} \pcs{n}(v) \cdot \khiL \fcs .
\end{equation}\fbeq
\end{paragraphesansboite}

\begin{paragraphesansboite}{Variational formulation:}
Given \cref{eq: eq of the variational part for corner fields}, $\SV$ is solution of the following variational formulation
\[\forall v\in V, \qquad  \int _{\Omega _{1}} \mu  \nabla  \SV\cdot \nabla  v = -\int _{\Omega _{1}} fv + \int _{\Gamma \cap \{x>\Rcoin\}} gv.\]
Hence we consider the following approximate problem: finding $\uL \in \evcs$ solution of 
\begin{equation}{\label{eq: approximated problem}}
\fbox{$\displaystyle \forall v\in \evcs, \qquad  \int _{\Omega _{1}} \mu  \nabla  \uL\cdot \nabla  v = -\int _{\Omega _{1}} fv + \int _{\Gamma \cap \{x>\Rcoin\}} gv. $} \end{equation}
Using Lax-Milgram theorem in $\evcs$ with $\|.\|_{V}$, it~is easy to see that \cref{eq: approximated problem} is well posed. The total approximate solution of \cref{eq: model pb for corner fields} is given by $\chi  S^{\mathcal{A}} + \uL$. 
Moreover, the coefficients $\sigma _{-n\pst} (\SDA)$ with $n\in [\![1,\nbcs]\!]$ are approximated by the numbers $\pcs{n}(\uL)$.\\

Let us explain how to compute numerically the integrals in \cref{eq: approximated problem}. Using \cref{eq: decomposition of elements of evcs}, \cref{eq: approximated problem} is equivalent to: $\forall v\in \evcs$,
\begin{multline}\label{eq: approximated problem with splitted integrals}
	\int _{\ounl} \mu  \nabla  \uL\cdot \nabla  v + \int _{\Omega _{1} \setminus \ounl} \mu \,\bigg(\sum _{n=1} ^{\nbcs} \pcs{n}(\uL)\,\nabla  \fcs \bigg)\cdot \bigg(\sum _{n=1} ^{\nbcs} \pcs{n}(v)\,\nabla  \fcs \bigg) \\
	= -\int _{\ounl} fv -\int _{\Omega _{1} \setminus \ounl} f\cdot \bigg(\sum _{n=1} ^{\nbcs} \pcs{n}(v)\,\nabla  \fcs \bigg) + \int _{\Gamma ^{L} \cap \{x>\Rcoin\}} gv + \int _{(\Gamma  \setminus \Gamma ^{L} )\cap \{x>\Rcoin\}} g\cdot \bigg(\sum _{n=1} ^{\nbcs} \pcs{n}(v)\,\nabla  \fcs \bigg)
\end{multline}
The integrals on the bounded domains $\ounl$ and $\Gamma ^{L}$ do~not raise any difficulty. As~for the other integrals, they have explicit integrands. Indeed, first, the $\fcs$ are explicit functions because the operators $\ROmegaPi$, $\RGammaPi$ and $\RcouchePi$ are explicit and so is $\rmicro$. Secondly, $f$ and $g$ depend on $S^{\mathcal{A}}$ (explicit) and on $F$. In~practice $F$ is proportional to another corner profile, computed using the same method, so it~is also explicitly known outside $\ounl$.\\

Theoretically, one could calculate analytically the integrals on $\Omega _{1} \setminus \ounl$ and $\Gamma  \setminus \Gamma ^{L}$. But it requires to express the $\fcs$, $f$ and $g$ in $\Omega  \setminus \Omega ^{L}$ with the variables $(r,\theta )$ rather than $z=x+\mathrm{i}y$. The obtained expressions usually involve tens of terms, so the integrands can have hundreds or thousands of terms (after expanding the products in the integrals). In the end it~is~not very efficient compared to a quadrature method, and it requires more work to implement.\\

Therefore we chose to compute the integrals on $\Omega _{1} \setminus \ounl$ and $\Gamma  \setminus \Gamma ^{L}$ in \cref{eq: approximated problem with splitted integrals} using a high order quadrature method. Before the quadrature, we apply the following changes of variables in order to integrate on bounded domains ($\varphi $ denotes one of the integrands that appear in \cref{eq: approximated problem with splitted integrals} and $a$ is a parameter that belongs to $\mathbb{R}_-^*$)
\[\int _{\Gamma  \setminus \Gamma ^{L} } \varphi  = \int _{L} ^{\infty } \varphi (x)\,\mathrm{d}x = \int _{0} ^{L^{1/a}} \varphi (t^{a} )\,at^{a-1} \,\mathrm{d}t, \qquad  \quad  \int _{\layer \setminus \layer^{L} } \varphi  = \int _{0} ^{L^{1/a}} \int _{-1} ^{0} \varphi (t^{a} ,y)\,at^{a-1} \,\mathrm{d}y\,\mathrm{d}t,\]
\[\int _{\Omega  \setminus \Omega ^{L} } \varphi  = \int _{L} ^{\infty } \int _{0} ^{\Theta } \varphi (r,\theta )\,r\mathrm{d}\theta \,\mathrm{d}r = \int _{0} ^{L^{1/a}} \int _{0} ^{\Theta } \varphi (t^{a} ,\theta )\,at^{2a-1} \,\mathrm{d}\theta \,\mathrm{d}t.\]
Note that the algebraic behavior of $\varphi $ at infinity implies that the integrand after the change of variable has a limited regularity at 0. 
This regularity must be sufficiently high for the quadrature method to converge at its the optimal rate.
By~choosing an appropriate value of $a$, we can achieve the needed regularity.
\end{paragraphesansboite}

\subsection{Proofs of convergence}

\begin{textesansboite}
In this section we prove the convergence of the method and provide error estimates rates, when $L\rightarrow \infty $ at a fixed $\nbcs$. We do~not know if there~is convergence when $\nbcs\rightarrow \infty $ at a fixed $L$. To prove it, one would need $\SDA$ to be a converging series (i.e.\ $T_{\geqslant q} (\SDA)$ converges pointwise when $q\rightarrow -\infty $). A heuristic reasoning and some numerical experiments suggest that it~is probably not the case, even far from the origin.
\\

One could also study the convergence of the method w.r.t.\ to the discretisation in space of $\{\pef(v) \mid v\in \evcs\}$, i.e.\ the first term of the direct sum in \cref{eq: def evcs}, by finite elements. Here we do~not consider it, because it~is classical and we want to focus on the truncation at infinity of the domain. An additional discretisation error is due to the quadrature of the integrals on $\Omega _{1} \setminus \ounl$ in the variational formulation. But numerical tests suggest that it~is close to the machine error, thanks to the high-order quadrature method we use.
\end{textesansboite}

\begin{theoreme}{}
\label{estimate of uL-SV}%
Let us assume $\beta ,\gamma  \leqslant  -(\nbcs+1) \pst$ (see \cref{eq: def S survaria et S varia,eq: def complement singulier}). There exists a constant $C$ that depends on $\nbcs$, the problem's data $(f,(s_{q} )_{q\in \TN} )$ and $(\beta ,\gamma ,\Theta ,\omega ,\mu _{0} ,\mu _{1} ,\rho _{0} ,\rho _{1} )$, but does~not depend on $L$, s.t.\ we have
\[\forall L>\Rcoin+1, \qquad  \|\uL - \SV\|_{V} \leqslant  C \.L^{-(\nbcs+1) \pst + \frac{1}{2} }.\]\fbeq
\end{theoreme}

\begin{remarque}
Note that the greater $\nbcs$, the faster the convergence. In~particular, for a large $\nbcs$, our method is a high order method. Moreover the condition $\beta \leqslant -1$, required when we defined $\beta $, is satisfied because $\Theta <2\pi $ implies that $(\nbcs+1) \pst \geqslant 2 \pst >1$. The optimal value of $\beta $ and $\gamma $ is $-(\nbcs+1)\pst$ as it minimizes the number of functions in $\mathcal{A}(\Pi )$ to compute, while ensuring that \cref{estimate of uL-SV} applies. Smaller values of $\beta $ and $\gamma $ would not improve the convergence speed. In addition, note that, if one takes $\gamma :=-(\nbcs+1)\pst$ e.g., then the complements at infinity $\fcs$ depend on $N$, although we did not write it in the notation.
\end{remarque}

\begin{demo}
To compare $\uL$ to $\SV$, let us first build an approximation of $\SV$ in $\evcs$, denoted $\SVL$, defined as follows
\begin{equation}{\label{eq: def of SVL}}
\SVL := \underbrace{(1-\khiL)\,\SV}_{\displaystyle =\pef(\SVL)} + \sum _{n=1} ^{\nbcs} \underbrace{\sigma _{-n\smallpst} (S^{\infty } )}_{\displaystyle =\pcs{n}(\SVL)} \cdot \khiL \fcs .\end{equation}
Céa's lemma yields \[\|\SV-\uL\|_{V} \lesssim  \|\SV-\SVL\|_{V}\] with a constant independent of $L$. Let us estimate $\|\SV-\SVL\|_{V}$. We have $\SV-\SVL = \khiL\Srest$ with
\begin{align*}
	\Srest &:= \SV - \sum _{n=1} ^{\nbcs} \sigma _{-n\smallpst} (S^{\infty } ) \,\fcs \\
	&= S - \chi  \cdot  T_{\geqslant \beta } \circ \rmicro\bigg( \frac{1}{\mu _{0}} \ROmegaPi (F^{\infty } _{\Omega } )+ \frac{1}{\mu _{1}} \RcouchePi (F^{\infty } _{\couche} ) + \!\fsum_{d\in \TN} \! s_{d} \,\phi_{d} \bigg) - \sum _{n=1} ^{\nbcs} \sigma _{-n\smallpst} (S^{\infty } ) \cdot T_{\geqslant \gamma } \circ \rmicro(\phi_{-n\smallpst} )
\end{align*}
by~definition of $\SV$ and $\fcs$, see \cref{eq: def S survaria et S varia,eq: def complement singulier}. Let $\delta  :=-(\nbcs+1) \pst$. Provided that $\beta \leqslant \delta $ and $\gamma \leqslant \delta $, we deduce
\[\Srest = S-T_{\geqslant \delta } (S^{\infty } ) + \sigma _{\delta } (S^{\infty } )\,\phi_{\delta } +o_{\partial } (\bfr^{\delta } )\]
because: $\forall d\in \TN,\; s_{d} = \sigma _{d} (S^{\infty } )$ and $T_{\geqslant \delta } \circ \rmicro(\phi_{\delta } )= \phi_{\delta }$ and $\forall d\in \TZ*\cap (-\infty ,\delta ),\; T_{\geqslant \delta } \circ \rmicro(\phi_{d} )=0$. Therefore $\Srest = \sigma _{\delta } (S^{\infty } )\,\phi_{\delta } +o_{\partial } (\bfr^{\delta } )$, which implies
\begin{equation}{\label{eq: estimates of Srest}}
	\forall (i,j)\in \mathbb{N}^2 , \qquad  \partial _{x} ^{i} \partial _{y} ^{j} \Srest = 
	\left\{\begin{array}{@{\;}l@{\;}l}
	O(\bfr^{\delta -i-j} ) &\text{ in }\Omega  \;\text{ when }  r\rightarrow \infty  \\
	O(\bfr^{\delta -i} ) &\text{ in }\layer \;\text{ when }  x\rightarrow \infty .
	\end{array}\right.\end{equation}
This yields
\begin{align}
	\|\SV-\SVL\|_{V} &= \|\nabla (\khiL\Srest)\|_{L^2 } = \|\nabla \khiL\|_{L^{\infty }} \|\Srest\|_{L^2 (\Omega _{1} \cap \{L-1<\bfr<L\})} + \|\khiL\|_{L^{\infty }} \|\nabla  \Srest\|_{L^2 (\Omega _{1} \cap \{\bfr>L-1\})} \nonumber\\
	&= O(L^{\delta + \frac{1}{2} } ) \label{eq: estimate of SV-SVL}
\end{align}
since $\|\khiL\|_{L^{\infty }}$ and $\|\nabla \khiL\|_{L^{\infty }}$ do~not depend on $L$.
\end{demo}

\begin{remarque}
Instead of \cref{eq: approximated problem}, another approximate variational formulation that would seem natural, although obviously less precise, is the following.
\def\evcsbis{\widetilde{V}_L}
\begin{equation}{\label{eq: alternative approximated problem}}
\forall v\in  \evcsbis, \qquad  \int _{\red{\ounl}} \mu  \nabla  \uL\cdot \nabla  v = -\int _{\red{\ounl}} fv + \int _{\red{\Gamma ^{L} }\cap \{x>\Rcoin\}} gv. \end{equation}
with $\evcsbis := \{u_{|\ounl} \mid u\in \evcs\}$. However this formulation would impose the following weak Neumann condition on $\Omega _{1} \cap \{\bfr=L\}$:
\[\forall n\in [\![1,\nbcs]\!],\qquad  \int _{\Omega _{1} \cap \{\bfr=L\}} \partial _{n} \uL \cdot \fcs =0.\]
It seems to us that the error introduced by this condition is roughly the same as a true Neumann condition $\partial _{n} \uL =0$. That error depends on $\partial _{n} {\SVL}_{|\Omega _{1} \cap \{\bfr=L\}}$, which is equivalent to $\sigma _{-\smallpst} (S^{\infty } )\,\partial _{n} \fcs[1]_{|\Omega _{1} \cap \{\bfr=L\}}$ when $L\rightarrow \infty $ (provided that $\sigma _{-\smallpst} (S^{\infty } ) \neq 0$). Therefore the error $\|\nabla (\SVL-\uL)\|_{L^2 (\ounl)}$ would not decrease faster when $\nbcs$ is larger. So \cref{eq: alternative approximated problem} seems to be a bad option.
\end{remarque}

\begin{theoreme}{}
\label{convergence des CCC}%
We still assume $\beta ,\gamma  \leqslant  -(\nbcs+1)\pst$. There exists a constant $C$ that depends on $\nbcs$, the problem's data $(f,(s_{q)_{q\in \TN}} )$ and $(\beta ,\gamma ,\Theta ,\omega ,\mu _{0} ,\mu _{1} ,\rho _{0} ,\rho _{1} )$, but does~not depend on $L$, s.t.\ we have
\[\forall n\in [\![1,\nbcs]\!],\ \forall L>\Rcoin+1, \qquad  |\pcs{n}(\uL) - \sigma _{-n\smallpst} (S^{\infty } )| \leqslant  C \.L^{(n-\nbcs-1)\pst + \frac{1}{2} }.\]\fbeq
\end{theoreme}

\begin{remarque}
For~any $n\in [\![1,\nbcs]\!]$, we have $(n-\nbcs-1)\pst + \frac{1}{2} \leqslant  -\pst+ \frac{1}{2} <0$ as $\Theta <2\pi $. So the given bound always tends to 0 when $L\rightarrow \infty $. Moreover the convergence is quick for $n$ small while it~is slow when $n$ approches $\nbcs$. Therefore, it~is a good idea to take $\nbcs$ a bit larger than the real number of coefficients $\sigma _{q} (\SDA)$ desired. The last computed coefficients will~not be used, but they permit a greater precision on the coefficients that are kept. This behavior is confirmed by numerical experiments, see \cref{sec: validation numerique complements singuliers}.
\end{remarque}

\begin{textesansboite}
The proof of \cref{convergence des CCC} relies on \cref{estimate of uL-SV}, but it uses its estimate only in $\Omega _{1} \setminus \ounl$. Indeed, in this region we have $\uL=\sum _{n=1} ^{\nbcs} \pcs{n}(\uL) \,\fcs$ ($\pef(\uL)$ does~not intervene because it vanishes) which makes it easier to obtain control over the coefficients $\pcs{n}(\uL)$. The error is this region is dominated by the global error thanks to the simple inequality $\|\nabla (\uL-\SV)\|_{L^2 (\Omega  \setminus \Omega ^{L} )} \leqslant  \|\uL-\SV\|_{V}$. However it does~not seem to be sharp. Since $\Omega  \setminus \Omega ^{L}$ tends to an empty domain when $L\rightarrow \infty $, one can expect $\|\nabla (\uL-\SV)\|_{L^2 (\Omega  \setminus \Omega ^{L} )}$ to decrease faster than $\|\uL-\SV\|_{V}$. Therefore, we think that the estimates of \cref{convergence des CCC} are~not optimal, even if we were not able to improve it.\\

The rest of the proof of \cref{convergence des CCC} isolates each coefficient $\pcs{n}(\uL)$. It~is difficult as the $\fcs$ are~not orthogonal. If they were, one could simply write something like
\begin{align*}
	\text{``} \ \big\|\nabla (\uL-\SV)\big\|_{L^2 (\Omega  \setminus \Omega ^{L} )} ^2  &\approx 
	\bigg\|\sum _{n=1} ^{\nbcs} \big(\pcs{n}(\uL)-\sigma _{-n\pst} (\SDA)\big)\,\nabla  \fcs \bigg\|_{L^2 (\Omega  \setminus \Omega ^{L} )} ^2  \\
	&= \sum _{n=1} ^{\nbcs} |\pcs{n}(\uL)-\sigma _{-n\pst} (\SDA)|^2  \cdot \|\nabla  \fcs \|_{L^2 (\Omega  \setminus \Omega ^{L} )} ^2  \ \text{''.} 
\end{align*}
As for the method of singular functions, error estimates on the coefficients were provided in the case of a unique singular function ($\nbcs=1$), see \cite{Loh98}. The proof uses orthogonality properties w.r.t.\ the bilinear form of the variational formulation. The difficulty here is that $\nbcs$ is arbitrary. To overcome this difficulty, we use an original proof by contradiction and the fact that the $\fcs$ are equivalent at infinity to orthogonal functions (the $\phi_{-n\pst}$).
\end{textesansboite}

\begin{demo}[Proof of \cref{convergence des CCC}]
For~any $n\in [\![1,\nbcs]\!]$, let us denote $\snL{n}{L} := L^{-(n-\nbcs-1) \pst - \frac{1}{2} } \big(\pcs{n}(\uL) - \sigma _{-n\smallpst} (S^{\infty } )\big) = L^{-(n-\nbcs-1) \pst - \frac{1}{2} } \pcs{n}(\uL-\SVL)$ and $\psinL{L} := L^{n \pst } \partial _{r} \fcs$. We have
\begin{align*}
	\bigg\|\sum _{n=1} ^{\nbcs} \snL{n}{L}\,\psinL{L} \bigg\|_{L^2 (\Omega  \setminus \Omega ^{L} )} &= L^{(\nbcs+1) \pst - \frac{1}{2} } \bigg\|\sum _{n=1} ^{\nbcs} \pcs{n}(\uL-\SVL)\,\partial _{r} \fcs \bigg\|_{L^2 (\Omega  \setminus \Omega ^{L} )} \\
	&\leqslant  L^{(\nbcs+1) \pst - \frac{1}{2} } \|\uL-\SVL\|_{V} \\
	&\leqslant  L^{(\nbcs+1) \pst - \frac{1}{2} } (\|\uL-\SV\|_{V} +\|\SV-\SVL\|_{V} ) =O(1). \addtocounter{equation}{1}\tag{\theequation}\label{eq: sum snL psinL is bounded}
\end{align*}
by \cref{estimate of uL-SV} and \cref{eq: estimate of SV-SVL}.\smallskip

Let us assume by contradiction that there exists $m\in [\![1,\nbcs]\!]$ s.t.\ $\pcs{m}(\uL) - \sigma _{-m\smallpst} (S^{\infty } ) \neq  O\big(L^{(m-\nbcs-1) \pst + \frac{1}{2} }\big)$, which means $\snL{m}{L} \neq O(1)$ when $L\rightarrow \infty $. Without loss of generality, we can assume that there exists a sequence $(L_{j} )_{j\in \mathbb{N}}$ s.t.\ $L_{j} \rightarrow \infty $, 
\begin{equation}{\label{eq: properties of m}}
	|\snL{m}{L_{j} }| \underset{j\rightarrow \infty }{\longrightarrow} \infty  \quad  \text{ and } \quad  \forall j\in \mathbb{N},\ \ |\snL{m}{L_{j} }| = \max \!\big(|\snL{1}{L_{j} }|,\dots , |\snL{\nbcs}{L_{j} }|\big).\end{equation}
Up to a subsequence, we can also assume that for every $n\in [\![1,\nbcs]\!]$, $\snL{n}{L_{j} } / \snL{m}{L_{j} }$ converges. Let $\ell _{n}$ be its limit. We have for~any $j\in \mathbb{N}$
\begin{align*}
	\bigg\| \frac{1}{\snL{m}{L_{j} }} \sum _{n=1} ^{\nbcs} \snL{n}{L_{j} }\,\psinL{L_{j} } \bigg\|_{L^2 (\Omega  \setminus \Omega ^{L_{j}} )}
	&= \int _{L_{j}} ^{\infty } \int _{0} ^{\Theta } \bigg(\sum _{n=1} ^{\nbcs} \frac{\snL{n}{L_{j} }}{\snL{m}{L_{j} }} \,\psinL{L_{j} }(r,\theta ) \bigg)^2  r\,\mathrm{d}\theta \,\mathrm{d}r \\
	&= \int _{1} ^{\infty } \int _{0} ^{\Theta } \bigg(\sum _{n=1} ^{\nbcs} \frac{\snL{n}{L_{j} }}{\snL{m}{L_{j} }} \,\psinL{L_{j} }(L_{j} t,\theta ) \,L_{j} \bigg)^2  t\,\mathrm{d}\theta \,\mathrm{d}t
\end{align*}
with $t := \frac{r}{L_{j}}$. Now we have in $\Omega $ when $r\rightarrow \infty $ that $\fcs =\phi_{-n\smallpst} +o_{\partial } (r^{-n\smallpst})$, so $\partial _{r} \fcs \sim \partial _{r} \phi_{-n\smallpst} = -n \pst \cdot r^{-n \pst -1} \sin \!\big(-n \pst \theta \big)$. Applying the dominated convergence theorem gives
\begin{equation}{\label{eq: limit of sum snL psinL over smL}}
	\bigg\| \frac{1}{\snL{m}{L_{j} }} \sum _{n=1} ^{\nbcs} \snL{n}{L_{j} }\,\psinL{L_{j} } \bigg\|_{L^2 (\Omega  \setminus \Omega ^{L_{j}} )} \underset{j\rightarrow \infty }{\longrightarrow} \int _{1} ^{\infty } \int _{0} ^{\Theta } \bigg(\sum _{n=1} ^{\nbcs} \ell _{n} \cdot  n \pst \cdot t^{-n \pst -1} \sin \!\Big({-n \pst \theta }\Big) \bigg)^2  t\,\mathrm{d}\theta \,\mathrm{d}t	\end{equation}
(the domination hypothesis is satisfied because $\exists r_{0} ,\forall r>r_{0} , \forall n\in [\![1,\nbcs]\!],\; |\partial _{r} \fcs |\leqslant  2n \pst \cdot r^{-n \pst -1}$). Since the functions $t^{-n \pst -1} \sin \!\big({-n \pst \theta }\big)$, $n\in [\![1,\nbcs]\!]$, are linearly independent and $\ell _{m} =1 \neq 0$, the limit in \cref{eq: limit of sum snL psinL over smL} is non-zero. However by \cref{eq: sum snL psinL is bounded} and \cref{eq: properties of m}, that limit must be zero. This is contradictory, and it concludes the proof.
\end{demo}

\begin{remarque}
\Cref{convergence des CCC} shows that the $\pcs{n}(\uL)$ are estimators of the $\sigma _{-n\pst} (\SDA)$. One could also consider a~posteriori estimators, i.e.\ computed after $\uL$. In the method of singular functions for instance, an a~posteriori estimator can be derived using a Green formula. It converges faster than the first obtained estimator. However, we were not able to extend this a~posteriori estimator method to our case because of the presence of the layer, the arbitrary number $\nbcs$ of complements and the non orthogonality of the $\fcs$.
\end{remarque}

\subsection{Numerical validation}
\label{sec: validation numerique complements singuliers}

\begin{textesansboite}
Let us illustrate the convergence of the method on a case where the corner profiles are known explicitly. We consider $\Theta \neq \pi $, $\xcoin:=-\operatorname{cotan}(\Theta )$ and
\[\Omega _{1} := \{(r \cos \theta  + \xcoin, r \sin \theta  - 1) \mid r\in \mathbb{R}_+^*, \theta \in (0,\Theta )\}.\]
See \cref{fig:Omega1 en forme de secteur}. We assume that the functions $\mu $ and $\rho $ are constant on $\Omega _{1}$. Hence the layer is virtual. For~any $d\in \TN$, we denote $\phicoin_{d} :(x,y) \mapsto  \phi_{d} (x-\xcoin,y+1)$. The couple $(\Omega _{1} ,\phicoin_{d} )$ is a translated version of $(\Omega ,\phi_{d} )$.

\begin{figure}[h]
\begin{center}
\def\widthFig{2.3} 
\def\heightFig{2} 
\def\ThetaFig{145} 
\def\epsFig{0.7} 
\begin{tikzpicture}[scale=0.9]
	\coordinate (finLambda) at (-\widthFig, {-\widthFig*tan(\ThetaFig)});
	\coordinate (coin) at ({-\epsFig/tan(\ThetaFig)},-\epsFig);
	\fill [pattern={Lines[angle=45,line width=.2pt]}, pattern color=white!70!gray] (coin) -- (\widthFig,-\epsFig) -- (\widthFig,\heightFig) -- (-\widthFig,\heightFig) -- (finLambda) -- cycle;
	\draw [thick] (finLambda) -- (coin) -- (\widthFig,-\epsFig);
	\draw [>=stealth,->] (-\widthFig-0.4,0) -- (\widthFig+0.4,0) node [right] {$x$};
	\draw [>=stealth,->] (0,-\epsFig-0.4) -- (0,{\heightFig+0.3}) node [right] {$y$};
	\draw (\widthFig/2,\heightFig*3/5) node {$\Omega$};
	\draw [->,>=latex] (0.35,0) arc (0:\ThetaFig-5:0.35) node[near start,above] {$\Theta$} -- (\ThetaFig:0.355);
	\draw [<->,>=latex] (\widthFig*0.85,0) -- node[midway,right]{1} (\widthFig*0.85,-\epsFig);
	\draw [<->,>=latex] (0,-\epsFig-0.1) -- node[midway,below]{$\xcoin$} ({-\epsFig/tan(\ThetaFig)},-\epsFig-0.1);
\end{tikzpicture}
\caption{The domain $\Omega_1$ used for the numerical validation}
\label{fig:Omega1 en forme de secteur}
\end{center}
\end{figure}
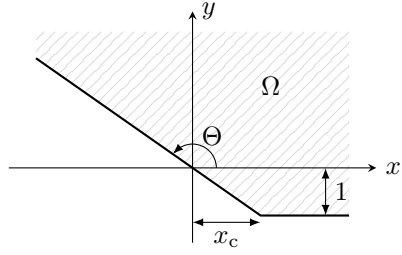

The following lemma shows that the corner profiles $\mathcal{S}_{q,0}$ have analytic expressions. We do~not write the explicit expressions, as they can be heavy, but they can be deduced from the proof below.
\end{textesansboite}

\begin{lemme}{}
For~any $q\in \TN$, the corner profile $\mathcal{S}_{q,0}$ is an explicit linear combination of the functions $\phicoin_{d}$, $d\in \TN$.
\end{lemme}

\begin{demo}
\def\zcoin{z_{\mathrm{c}}}%
We denote $z$ the complex variable $x+\mathrm{i}y$ and $\zcoin = \xcoin-\mathrm{i} = -\frac{\cos \Theta +\mathrm{i}\sin \Theta }{\sin \Theta } =- \frac{1}{\alpha \sin (\Theta )}$. For~any $(x,y)\in \Omega $ with $x^2 +y^2 $ large enough, we have:
\begin{align*}
	\phicoin_{d} :(x,y) 
	&= \Im\big[(\alpha  (z-\zcoin))^{d} \big]\\
	&= \Im\!\left[(\alpha z)^{d} \left(1- \frac{\zcoin}{z} \right)^{d} \right]\\
	&= \Im\!\left[(\alpha z)^{d} \sum_{n=0}^\infty  \binom{d}{n} \left(- \frac{\zcoin}{z} \right)^{n} \right]\\
	&= \sum_{n=0}^\infty  \binom{d}{n} (\sin \Theta )^{-n} \,\Im\!\left[(\alpha z)^{d-n} \right]
\end{align*}
where $\binom{d}{n} := \frac{1}{n!} \prod _{k=0} ^{n-1} (d-k)$. Similarly, when $(x,y)\in \couche$, we have
\begin{align*}
\phicoin_{d} :(x,y) &= \Im\left[(\alpha x)^{d} \left(1+ \frac{\mathrm{i}y-\zcoin}{x} \right)^{d} \right]\\
&= \Im\left[(\alpha x)^{d} \sum_{n=0}^\infty  \binom{d}{n} \left( \frac{\mathrm{i}y-\zcoin}{x} \right)^{n} \right]\\
&= \sum_{n=0}^\infty  \binom{d}{n} \,x^{d-n} \,\Im[\alpha ^{d} (\mathrm{i}y-\zcoin)^{n} ]
\end{align*}
Let \[\varphi (x,y) := \sum _{n=0} ^{\lceil {d}\rceil +1} \binom{d}{n} \times \left\{\begin{array}{ll}
	(\sin \Theta )^{-n} \,\Im\!\left[(\alpha z)^{d-n} \right] &\text{ in }\Omega \\
	x^{d-n} \,\Im[\alpha ^{d} (\mathrm{i}y-\zcoin)^{n} ] &\text{ in }\couche
\end{array}\right.\]
We have $\varphi \in \mathcal{A}(\Pi )$ and $\nabla (\phicoin_{d} -\varphi ) = O(r^{-3} )$ at infinity, so $\phicoin_{d} \in V+\khiz \mathcal{A}(\Pi )$. Therefore by the uniqueness in \cref{cadre fonctionnel pour les champs de coin}
\begin{equation}{\label{eq: expression de phicoin}}
	\phicoin_{d} = \sum_{q\in \TN ,q\leqslant d} \sigma _{q} (\phicoin_{d} )\,\mathcal{S}_{q,0} 
	\qquad  \text{with}  \qquad  \forall q\in {\textstyle\TN},\ \ \sigma _{q} (\phicoin_{d} ) = \mathbf{1}_{\mathbb{N}} (d-q) \cdot  \binom{d}{d-q} (\sin \Theta )^{q-d} .
\end{equation}
To conclude, it suffices to notice that the matrix $(\sigma _{q} (\phicoin_{d} ))_{q,d\in \TN}$ is triangular with coefficients 1 on the diagonal, hence it~is invertible.
\end{demo}

\begin{textesansboite}
\cref{eq: expression de phicoin} shows that $\sigma _{q} (\phicoin_{d} )$ can be non-zero only if $d-q$ is an integer. In order to analyse relative errors, we would like to have as much non-zero coefficients as possible. Hence, for the numerical tests, we chose $\Theta = \frac{2\pi }{3}$ because the small denominator of $\pst = \frac{3}{2}$ makes $d-q$ integer more often. Now let us choose the target function that we will compute numerically.
\begin{enumerate}
	\item $\mathcal{S}_{\pst,0}$ (which is equal to $\phicoin_{\pst}$ by \cref{eq: expression de phicoin}) is~not very regular at the corner $(\xcoin,-1)$ (it~is $H^{5/2-\delta } _{\mathrm{loc}}$ for~any $\delta >0$ but not $H^{5/2} _{\mathrm{loc}}$). Indeed, $\phicoin_{\pst}$ behaves like $r_{\mathrm{c}} ^{3/2}$, where $r_{\mathrm{c}}$ is the distance to $(\xcoin,-1)$. This deteriorates the numerical error when using Lagrange finite elements of high order, which appeared during our tests to be necessary in order to avoid numerical instability. Therefore, we do~not choose $\mathcal{S}_{\pst,0}$. This problem could be fixed by adding to $\evcs$ the singular function $\khiz \,\phicoin_{\pst}$ located at the corner (see \cite{Loh98}), but it~is~not in the spirit of this work.
	\item Next, $\mathcal{S}_{2\pst,0} =\phicoin_{2\pst}$ satisfies $\sigma _{q} (\phicoin_{2\pst} )=0$ for~any $q\in -\TN$, so we do~not choose it either. 
	\item Finally, we choose \[\phicoin_{3\pst} =\mathcal{S}_{3\pst,0} + \frac{35}{2\sqrt{3}} \,\mathcal{S}_{\pst,0} \]
	because it~is very regular (it~is $H^{11/2-\delta } _{\mathrm{loc}}$ for~any $\delta >0$) and we have
	\begin{equation}{\label{eq: valeurs exactes des CCC calcules numeriquement}}
		\begin{array}{c}
			\displaystyle \sigma _{-\pst} \big(\phicoin_{3\pst} \big) = \frac{-7}{144} \approx  -4.9 \times 10^{-2} , \qquad  \sigma _{-2\pst} \big(\phicoin_{3\pst} \big) =0, \\[3.5ex]
			\displaystyle \sigma _{-3\pst} \big(\phicoin_{3\pst} \big) = \frac{35}{10368\sqrt{3}} \approx  1.9 \times 10^{-3} , \qquad  \sigma _{-4\pst} \big(\phicoin_{3\pst} \big) =0, \quad  \mathrm{etc.}
		\end{array}
	\end{equation}
\end{enumerate}
We used P4 Lagrange finite elements, and relied on the FEM-BEM library XLiFE++ \cite{XLiFE}. \Cref{fig: paraview CC solution exacte} represents the exact solution, and the differences $\uL^{h} -\SV$ are represented when $\nbcs=0$ in \cref{fig: paraview CC erreur Dirichlet} (i.e.\ for a homogeneous Dirichlet condition on $\partial \Omega _{1} ^{L}$) and when $\nbcs=4$ in \cref{fig: paraview CC erreur N=4}. One can see that the error for $\nbcs=4$ is much smaller than for the Dirichlet condition. In addition, the relative error obtained for $\nbcs=4$ is very small (the maximal error is $8\cdot 10^{-5}$ whereas the exact solution is in $[-9.7,0]$ on $\Omega _{1} ^{L}$) although $L$ is not very large ($L=4$).\\

\begin{figure}[p]
	\includegraphics[width=70mm]{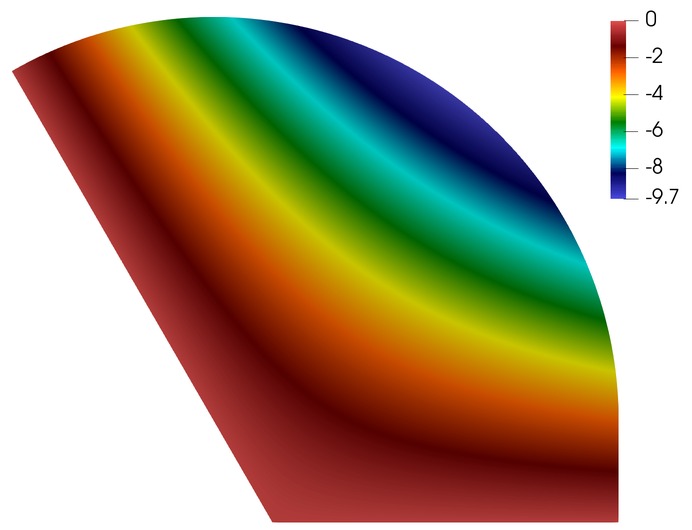}
	\caption{Exact solution $S=\phicoin_{3\pst}$ represented on $\Omega _{1} ^{2.5}$}
	\label{fig: paraview CC solution exacte}
\end{figure}

\begin{figure}[p]
	\includegraphics[width=70mm]{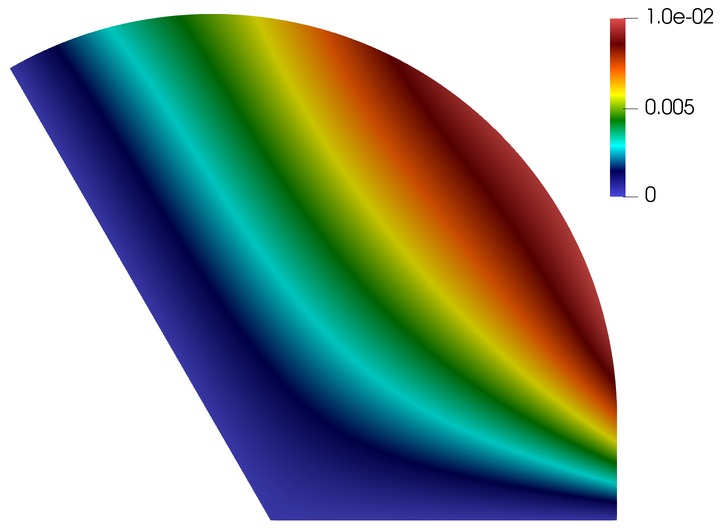}
	\caption{Difference $\uL^{h} -\SV$ for a homogeneous Dirichlet condition on $\partial \Omega _{1} ^{L}$ (with no complement at infinity)}
	\label{fig: paraview CC erreur Dirichlet}
\end{figure}

\begin{figure}[p]
	\includegraphics[width=70mm]{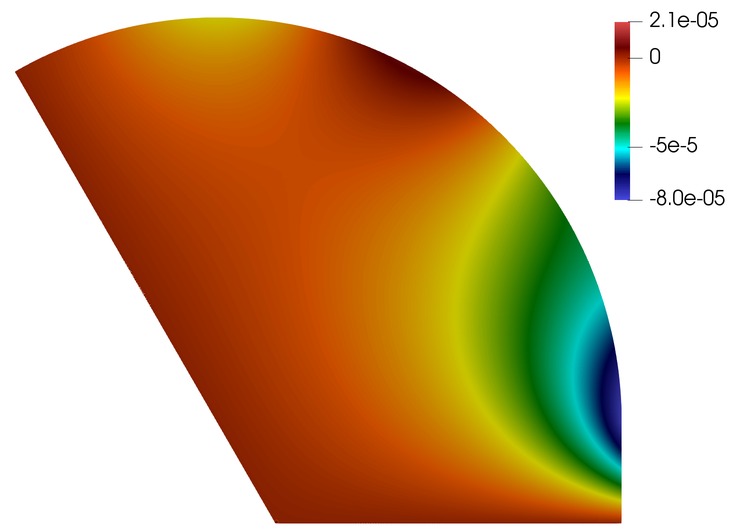}
	\caption{Difference $\uL^{h} -\SV$ on $\Omega _{1} ^{L}$ for $\nbcs=4$}
	\label{fig: paraview CC erreur N=4}
\end{figure}

Moreover, \cref{fig: erreurs CC norme,fig: erreurs CCC} plot the errors obtained for $\nbcs=2$ and $\nbcs=4$ (column 1, resp.\ 2). \Cref{fig: erreurs CC norme} plots the norm of $\uL^{h} -\SV$ while \cref{fig: erreurs CCC} shows the error on the asymptotic coefficients. In each graph, the errors are plotted w.r.t.\ $L$ and each curve is associated to a different mesh step (denoted $h$).\\

For $\nbcs=2$, the errors on the norm on the fixed bounded subdomain $\Omega _{1} ^{2.5}$ first becomes smaller when $L$ increases, before reaching a plateau. This plateau seems to be related to the mesh discretisation, as its value becomes smaller when the mesh step decreases. For the error on the entire domain $\Omega _{1}$, the behavior is similar, but the plateau is replaced by an increasing phase. This behavior also seems to be due to the discretisation, as it~is~not predicted by the theory on the continuous problem (see \cref{estimate of uL-SV}). We believe that it~is linked to the natural norm of the source term, which is $\|(1+r)f\|_{L^2 (\Omega _{1} )}$ (see \cref{pb variationnel de Poisson dans Omega1}). The weight $1+r$ diverges at infinity, so the discretisation error depends on $L$. Note however that the slope of the curve of $\|\SV-\uL\|_{V}$ in the first phase ($-3.61$) is coherent with the one predicted by \cref{estimate of uL-SV} (which is 4). The slope for the norm on $\Omega _{1} ^{2.5}$ (which is $-6.59$) indicates that there may be a higher rate of convergence on bounded fixed subdomains.\\

As for the case $\nbcs=4$, we cannot observe the first phase of convergence, but only plateaux and increasing phases. Indeed, the converge is much faster than for $\nbcs=2$, and it seems that the error we want to observe immediatly falls under the discretisation error. This issue cannot be solved by reducing $L$, as our minimal $L$ is already small (it~is equal to $3.2$) and some space needs to be left for the shape of the tip of the layer and for the truncation functions.\\

Therefore it~is difficult to observe numerically the convergence of $\uL$ in norm. However, what interests us the most are the coefficients $\sigma _{-n\pst} (\SDA)$. \Cref{fig: erreurs CCC} shows that they are computed with very high precision, even for values of $L$ that are~not large. For example, for $\nbcs=4$ and $L=5$, we obtain $|\pcs{1}(\uL^{h} )-\sigma _{-\pst} (\SDA)| \approx  10^{-10}$. Given that $\sigma _{-\pst} (\SDA) \approx -5\cdot 10^{-2}$ by \cref{eq: valeurs exactes des CCC calcules numeriquement}, this a relative error of $2\cdot 10^{-9}$. Taking $L=12.5$ and $h=0.035$, the relative error reaches $4\cdot 10^{-11}$.\\

The curves of $|\pcs{n}(\uL^{h} ) - \sigma _{-n\pst} (\SDA)|$ in \cref{fig: erreurs CCC} have roughly the shape of a line in log-log scale. Note that, surprisingly, the convergence of the coefficients $\pcs{n}(\uL^{h} )$ is~not affected by the plateaux of the norm of $\uL^{h}$. All the measured slopes (indicated in the graphs) are greater in absolute value than the slope announced by \cref{convergence des CCC} (written in a table in \cref{fig: erreurs CCC}). This is consistent with the theorem and it seems to confirm that the speed of convergence given by \cref{convergence des CCC} is~not optimal (for the reasons explained under the theorem).\\

Moreover, according to \cref{convergence des CCC}, the convergence should be faster when $\nbcs$ increases and slower when $n$ increases. This behavior can indeed be observed when comparing the measured slopes between each other (with just one exception between $(\nbcs,n)=(4,2)$ and $(4,3)$). In addition, the values of the errors have a similar behavior: they become smaller when $\nbcs$ increases and larger when $n$ increases. For example, for $(\nbcs,n)=(2,1)$, the obtained errors span between $10^{-5}$ and $10^{-8}$, while for $(\nbcs,n)=(4,1)$ they go from $10^{-8}$ to $10^{-12}$. This is an interesting point, not predicted by \cref{convergence des CCC}. Of course it depends on the values of $L$ used for the computations, the chosen values (between $3.2$ and $15.5$) seem reasonable.
\end{textesansboite}

	\tikzset{external/export=true}%
	\tikzsetfigurename{fig-erreurs_CC_norme-}%
	\pgfplotsset{
    every axis/.append style={
        thick,
    },
    cycle list/Dark2,
    cycle multiindex* list={
        mark list*\nextlist
        Dark2\nextlist
    },
    legend cell align=left,
    xmajorgrids=true,
    ymajorgrids=true,
}

\catcode`\;=12
\newcommand{\tracerAxisCC}[6][0,025]{ 
    \begin{tikzpicture}
    \begin{axis}[
        width=70mm,
        /pgfplots/table/x=L,
        /pgfplots/table/y=#3,
        xlabel={$L$},
        xlabel shift=-0.3ex,
        ylabel shift=1.3ex,
        title style={yshift=0.5ex},
        xmode=log,
        ymode=log,
        #6
    ]
    \if 2#2 
        \addplot table {../figures/erreurs_CC/erreurs_N_2_h_0,05.dat};
        \addplot table {../figures/erreurs_CC/erreurs_N_2_h_0,035.dat};
        \addplot table {../figures/erreurs_CC/erreurs_N_2_h_0,025.dat};
        \if -#4 \else
            {\pgfplotstablesave[
                columns={L,#3},
                skip rows between index={0}{#4},
                skip rows between index={#5}{8}
            ]{../figures/erreurs_CC/erreurs_N_2_h_#1.dat}{tempo.dat}}
        \fi
    \else
        \addplot table {../figures/erreurs_CC/erreurs_N_4_h_0,05.dat};
        \addplot table {../figures/erreurs_CC/erreurs_N_4_h_0,035.dat};
        \addplot table {../figures/erreurs_CC/erreurs_N_4_h_0,025.dat};
        \if -#4 \else
            {\pgfplotstablesave[
                columns={L,#3},
                skip rows between index={0}{#4},
                skip rows between index={#5}{8}
            ]{../figures/erreurs_CC/erreurs_N_4_h_#1.dat}{tempo.dat}}
        \fi
    \fi
    \if -#4 \else
        \addplot[mark=none, very thick] table[y={create col/linear regression={y=#3}}] {tempo.dat};
        \node [draw, fill=white, anchor=north east] at (rel axis cs: 0.95,0.95) {slope $= \pgfmathprintnumber{\pgfplotstableregressiona}$};
    \fi
    \end{axis}
    \end{tikzpicture}
}

\newcommand{\uLh}{u_L^h}
\def\eeltf{1.3ex} 

\begin{figure}[h]
\begin{tabular}{rr}
\tracerAxisCC{2}{errTot}{0}{5}{
    title={$N=2$},
    legend to name=fig: legende partagee,
    legend entries={$h=0.05$, $h=0.035$, $h=0.025$, linear regression},
    ylabel={$\|\SV-\uLh\|_V$},
} &
\tracerAxisCC{4}{errTot}{-}{-}{
    title={$N=4$},
    legend style={at={(1.1,1)},anchor=north west},
} \\[\eeltf]

\tracerAxisCC{2}{errLmin}{0}{3}{
    ylabel={\parbox{27ex}{\centering $\|\nabla(\SV-\uLh)\|_{L^2}$ on\\[.3ex] the fixed subdomain $\Omega_1^{2.5}$}},
} &
\tracerAxisCC{4}{errLmin}{-}{-}{
} \\[\eeltf]

\multicolumn{2}{c}{\hspace{15mm}\ref*{fig: legende partagee}}
\end{tabular}
\longCaption{Observed numerical error between $\uLh$ and $\SV$ in terms of $L$}{ ($h$ is the mesh step)}
\label{fig: erreurs CC norme}
\end{figure}%
	\tikzset{external/export=false}%
	\tikzsetfigurename{main-figure}
	\tikzset{external/export=true}%
	\tikzsetfigurename{fig-erreurs_CCC-}%

\begin{figure}
\begin{tabular}{rr}
\tracerAxisCC{2}{abs1}{0}{8}{
    title={$N=2$},
    ylabel={$|\pcs{1}(\uLh)-\sigma_{-\pst}(\SDA)|$},
} &
\tracerAxisCC[0,035]{4}{abs1}{0}{7}{
    title={$N=4$},
} \\[\eeltf]

\tracerAxisCC{2}{abs2}{0}{8}{
    ylabel={$|\pcs{2}(\uLh)-\sigma_{-1\pst}(\SDA)|$},
} &
\tracerAxisCC{4}{abs2}{0}{6}{
} \\[\eeltf]

\multicolumn{1}{c}{\raisebox{20mm}{\ref*{fig: legende partagee}}\hspace*{10ex}}
&
\hspace{-5ex}
\tracerAxisCC{4}{abs3}{0}{8}{
    ylabel={$|\pcs{3}(\uLh)-\sigma_{-3\pst}(\SDA)|$},
} \\[\eeltf]

\multicolumn{1}{c}{\raisebox{25mm}{\parbox[b]{35ex}{
\centering Order of convergence of $|\pcs{n}(\uLh)-\sigma_{-n\pst}(\SDA)|$ given by \cref{convergence des CCC} (for comparison)\\[2ex]
\renewcommand{\arraystretch}{1.15}
\begin{tabular}{ccc}
    \cline{2-3}
    & $N=2$ & $N=4$ \\
    \hline
    $n=1$ & $-2.5$ & $-5.5$ \\
    \hline
    $n=2$ & $-1$ & $-4$ \\
    \hline
    $n=3$ & $\diagup$ & $-2.5$ \\
    \hline
    $n=4$ & $\diagup$ & $-1$ \\
    \hline
\end{tabular}}}
\renewcommand{\arraystretch}{1}\hspace*{10ex}}
&
\hspace{-5ex}
\tracerAxisCC{4}{abs4}{0}{8}{
    ylabel={$|\pcs{4}(\uLh)-\sigma_{-4\pst}(\SDA)|$},
}
\end{tabular}

\longCaption{Observed numerical error between the numbers $\pcs{n}(\uLh)$ and the asymptotic coefficients $\sigma_{-n\pst}(\SDA)$ in terms of $L$}{ ($h$ is the mesh step)}
\label{fig: erreurs CCC}
\end{figure}%
	\tikzset{external/export=false}%
	\tikzsetfigurename{main-figure}

\begin{remarque}
In \cref{chap: ME}, we introduce effective models with a Dirichlet-to-Neumann condition (DtN) near the corner. The DtN's modes are non-orthogonal but almost orthogonal (in some sense) which allows proving the convergence of the method. A posteriori, it seems to us that a similar DtN might be used to truncate the corner profiles' domain at infinity, instead of the complements at infinity presented in this section. Chronologically, we designed the complements at infinity before searching for effective models, and we did not have time to try a DtN for corner profiles, but this is an interesting research perspective.
\end{remarque}

\section{Numerical computation of the far fields}
\label{sec: numerique des CL}

\begin{textesansboite}
In this section, we explain how we computed numerically the far fields for various angles $\Theta $, and we provide numerical results that confirm the error estimate of \cref{DA de SE avec CL et delta variable} and the validity of the asymptotic expansion of $u_{\varepsilon }$. We do~not perform the error analysis of the numerical method because the discretization in space and the DtN that we use are classical, while taking into account the boundary condition on $\Gamma $ is less obvious but is related to existing methods in the literature (see below). In addition, our final objective is to replace the computation of the far fields by effective models, which we justify in greater details in \cref{chap: ME}.
\end{textesansboite}

\subsection{Presentation of the numerical method}

\begin{textesansboite}
The far fields $u_{p,\ell }$ are built inductively w.r.t.\ $p\in \NN$. In addition, \cref{borne sur l} states that, for~any $p$, there exists an integer $\ell _{p}$ that depends only on $(\Theta ,p)$ s.t.\ we have $u_{p,\ell } =0$ for~any $\ell >\ell _{p}$. This number $\ell _{p}$ can be computed at low cost using the tools of \cref{sec: implementation des espaces A}, like we did in \cref{sec: exemples de structures de DA} for a few example angles $\Theta $. Hence, for~any $p\in \NN$, we build all at once the finitely many fields $u_{p,\ell }$ for $\ell  \in [\![0,\ell _{p} ]\!]$.\\

We recall that, by \cref{construction directe des upl}, for~any $p\in \NN$ and $\ell \in \mathbb{N}$, $u_{p,\ell }$ is the unique solution in $H^{1} (\Omega )+\khiz \mathcal{A}(\Omega )$ of:
\begin{equation}{\label{eq: numerique CL rappel eq CL}}
\left\{\begin{array}{r@{\;=\;}ll@{}}
    \mu _{0} \Delta  u_{p,\ell } +\omega ^2 \rho _{0} u_{p,\ell } & \fsource \,\delta _{p,0} \,\delta _{\ell ,0} &\text{ in }\Omega  \\[0.5mm]
    u_{p,\ell } & 0 &\text{ on }\bordgauche\\[0.5mm]
    u_{p,\ell } &\displaystyle \sum_{n=0}^\infty  \frac{\mu _{0} }{\mu _{1}} \,T_{n} \cdot  (\partial _{x} ^2 +k_{1} ^2 )^{n} \partial _{y} u_{p-1-2n,\ell } &\text{ on }\Gamma  \\[-1.5mm]
    \sigma _{d} (u_{p,\ell } )&\displaystyle \!\! \sum _{\substack{p'\in \NN\\ p-p'\in \NN+ 2\pst }} \; \sum _{\substack{d'\in \TZ*\\ p'+d'\in \intvl{0,p+d}}} \! \sum _{\ell '=0} ^{\ell } \cfinal_{d,d',p-p',\ell -\ell '} \cdot \sigma _{d'} (\uDA_{p',\ell '} ) & \ \forall d\in -\TN \!\!\!\!
\end{array}\right.
\end{equation}
where the above sums contain a finite number of non-zero terms. Computing the far fields presents two main challenges.
\begin{enumerate}
    \item First, we need to compute the coefficients $\sigma _{d'} (\uDA_{p',\ell '} )$ with $d'\in \TN$ that are involved in the formula of $\sigma _{d} (u_{p,\ell } )$ for some $d\in -\TN$ and $p>p'$ s.t.\ we are interested in computing $u_{p,\ell }$. These coefficients $\sigma _{d'} (\uDA_{p',\ell '} )$ belong to the variational part of $u_{p',\ell '}$, so computing them is a challenge similar to the computation of the coefficients $\sigma _{q} (\ccDA_{d,n} )$ with $q\in -\TN$ addressed in \cref{sec: calcul numerique des CC}. An additional difficulty is that the solutions of Helmholtz's equation in $H^{1} (\Omega )$ are not very regular at the corner, especially for a large $\Theta $, so they are poorly approximated by finite elements. This issue must be handled in order to approximate accurately $u_{p,\ell }$ and its behavior around the corner.
    \item Secondly, each far field has a heterogeneous Dirichlet condition on $\Gamma $ that depends on possibly high order derivatives of previous far fields. Mathematically, this is fine because the far fields are $\mathcal{C}^{\infty ^{}}$ in a neighborhood of $\Gamma $, see \cref{validation construction et DA des champs} and its proof. However, numerically, this is challenging because derivatives of finite elements functions (e.g.\ of type Lagrange) are not continuous in general, so they cannot be directly used as a Dirichlet condition for another finite elements computation.
\end{enumerate}
In addition, one has to deal with the fact that far fields may diverge at the corner, so they are not variational solutions of Helmholtz equations. This will be handled at the same time as point 1. Besides, the unbounded domain $\Omega $ is truncated numerically using a perfectly matched layer (PML), which is classical for the Helmholtz equation. Now let us explain how we address the two above challenges one by one.
\end{textesansboite}

\begin{paragraphesansboite}{Computation of the coefficients $\sigma _{d} (\uDA_{p,\ell } )$ with $d>0$:}
This issue is addressed using a Dirichlet-to-Neumann (DtN) map on an artificial boundary around the corner. Let us introduce the relevant domains. Let $\Rsource := \operatorname{dist}(0,\operatorname{supp} f)$ and $\RME \in (0,\Rsource)$. We have $\Rsource>0$ because we assumed $\operatorname{dist}(\operatorname{supp}(\fsource),\Gamma )>0$. We define $\OME:= \Omega  \setminus \overline{B}(0,\RME)$, $\GME:= \Gamma  \setminus \overline{B}(0,\RME)$ and $\SME := \{(\RME \cos \theta , \RME \sin \theta ) \mid \theta \in (0,\Theta )\}$. Note that $\partial \OME = \overline{\bordgaucheME \cup  \SME \cup  \GME}$. See \cref{fig: domaine numerique CL}. The ``$\circ $'' exponents are intended to recall that a ball is excluded around the corner.\\

\begin{figure}[h]
\def\widthFig{3.3} 
\def\heightFig{2.7} 
\def\ThetaFig{205} 
\def\RMEFig{1.3} 
\begin{tikzpicture}[scale=0.8]
	\coordinate (finLambda) at (-\widthFig, {-\widthFig*tan(\ThetaFig-180)});
	\coordinate (dbtLambda) at ({\RMEFig*cos(\ThetaFig)}, {\RMEFig*sin(\ThetaFig)});
	\fill [pattern={Lines[angle=45,line width=.2pt]}, pattern color=white!70!gray] (dbtLambda) arc (\ThetaFig:0:\RMEFig) -- (\widthFig,0) -- (\widthFig,\heightFig) -- (-\widthFig,\heightFig) -- (finLambda) -- cycle;
	\draw [thick] (finLambda) -- (dbtLambda) arc (\ThetaFig:0:\RMEFig) -- (\widthFig,0);
	\draw [dashed] (dbtLambda) -- (0,0);
	\draw [>=stealth,->] (-\widthFig-0.4,0) -- (\widthFig+0.4,0) node [right] {$x$};
	\draw [>=stealth,->] (0,{-\widthFig*tan(\ThetaFig-180)}) -- (0,{\heightFig+0.3}) node [right] {$y$};
	\draw (-\widthFig*0.5,\heightFig*0.6) node {$\OME$};
	\draw [->,>=latex] (0.35,0) arc (0:\ThetaFig-5:0.35) node[near start,above] {$\Theta$} -- (\ThetaFig:0.36);
	\draw [<-,>=latex] (\widthFig*0.65,0) to[bend right=20] ++(0.15,-0.55) node[below]{$\ \ \GME$};
	\draw [<-,>=latex] ({\RMEFig*cos(55)}, {\RMEFig*sin(55)}) to[bend left=20] ++(0.4,0.45) node[right]{\raisebox{1ex}{$\SME$}};
	\draw [<->,>=latex] (0,-0.2) -- node[midway, below] {$\RME$} (\RMEFig, -0.2);
\end{tikzpicture}\par
\caption{The domains used to compute the far fields}
\label{fig: domaine numerique CL}
\end{figure}
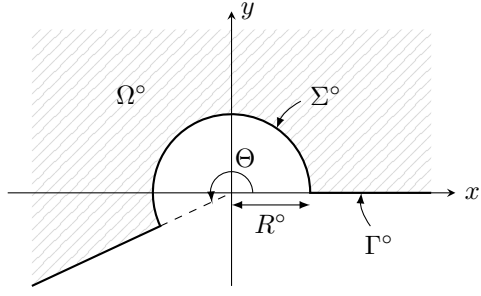

It is classical that any function $u\in H^{1} (\Omega )$ that satisfies
\begin{equation}{\label{eq: numerique CL pb Helmholtz homog}}
    \left\{\begin{array}{r@{\;=\;}ll@{}}
        \mu _{0} \Delta  u +\omega ^2 \rho _{0} u & f &\text{ in }\Omega  \\[0.5mm]
        u & 0 &\text{ on }\partial \Omega 
    \end{array}\right.
\end{equation}
with $f\in \{\fsource,0\}$ has the following modal decomposition in $\Omega  \cap B(0,\Rsource)$ for some coefficients $s_{n} \in \mathbb{C}$:
\begin{equation}{\label{eq: decompo modale solution H1 Helmholtz}}
    u(r,\theta ) = \sum_{n=0}^\infty  s_{n} \cdot J_{n\pst} (k_{0} r) \sin  (n\txtpst \theta ).
\end{equation}
Here $J_{n\pst}$ denotes a Bessel function of the first kind. In addition, the above series and its gradient converge uniformly on $B(0,R)$ for~any $R\in (0,\Rsource)$. So we also have
\begin{equation}{\label{eq: dr decompo modale solution H1 Helmholtz}}
    \partial _{r} u = \sum_{n=0}^\infty  s_{n} \cdot k_{0} J_{n\pst} '(k_{0} r) \sin  (n\txtpst \theta ).
\end{equation}
Moreover, we have
\begin{equation}{\label{eq: formule coef Fourier solution H1 Helmholtz}}
    \forall n\in \mathbb{N}^*, \qquad  s_{n} = \sqrt{ \frac{2}{\RME\Theta } } \cdot  \frac{\langle u ,\varphi _{n} \rangle _{L^2 (\SME)} }{J_{n\pst} (k_{0} \RME)} \qquad  \text{with } \ \varphi _{n} :\theta  \mapsto \sqrt{ \frac{2}{\RME\Theta } } \sin (n\theta )
\end{equation}
by orthonormality in $L^2 (\SME)$ of the functions $\varphi _{n}$. Combining \cref{eq: dr decompo modale solution H1 Helmholtz,eq: formule coef Fourier solution H1 Helmholtz}, we obtain
\begin{equation}{\label{eq: def opF}}
    \partial _{r} u_{|\SME} = \opF(u_{|\SME} ) \qquad  \text{with}  \qquad  \opF : v \mapsto  \sum_{n=0}^\infty  \frac{k_{0} \.J'_{n\pst} (k_{0} \RME)}{J_{n\pst} (k_{0} \RME)} \,\langle v,\varphi _{n} \rangle _{L^2 (\SME)} \,\varphi _{n.}
\end{equation}
One can show that $\opF$ is well-defined from $H^{1/2} _{00} (\SME)$, the set of elements of $H^{1/2} (\SME)$ whose extension by 0 to $\partial \OME$ remains $H^{1/2}$, to $(H^{1/2} _{00} (\SME))'$. $\opF$ is the DtN map that we use.\\

Let us use this on $u_{p,\ell }$ for a given $(p,\ell )\in \NN\times \mathbb{N}$. By~\cref{decompo modale des CL}, we have the following modal decomposition in $\Omega  \cap B(0,\Rsource)$:
\begin{equation}{\label{eq: rappel decompo modale CL pour numerique CL}}
    u_{p,\ell } \.=\. \sum _{n=0} ^{\lfloor {p}\rfloor } \sum _{\substack{d\in \TZ*\\ d\geqslant -p+n}} \sigma _{d} (\uDA_{p-n,\ell } ) \,\mcll_{d,n} .
\end{equation}
We recall $\mcll_{d,n}$ is introduced in \cref{convergence des mcl} as the limit function of the converging series $\mcla_{d,n}$, which is defined in \cref{def: mcl}. In addition, \cref{decompo modale des CL} states that \cref{eq: rappel decompo modale CL pour numerique CL} and its gradient converge uniformly on $B(0,R)$ for~any $R\in (0,\Rsource)$. Let us write
\[u_{p,\ell } := \uV_{p,\ell } + \uA_{p,\ell } \quad  \text{with }  \left\{\begin{array}{l}
    \uV_{p,\ell } :=\displaystyle \sum _{d\in \red\TN} \sigma _{d} (\uDA_{p,\ell } ) \,\mcll_{d,0} \vrule width0pt height 3ex\relax \\
    \uA_{p,\ell } :=\displaystyle \smash{\sum _{\substack{d\in \red{-\TN}\\ d\geqslant -p}} \sigma _{d} (\uDA_{p,\ell } ) \,\mcll_{d,0} } + \sum _{n=\red1} ^{\lfloor {p}\rfloor } \smash{\sum _{\substack{d\in \TZ*\\ d\geqslant -p+n}} \sigma _{d} (\uDA_{p-n,\ell } ) \,\mcll_{d,n} } \vrule width0pt depth 4ex\relax
\end{array}\right. \vspace{1ex} \]
When computing $u_{p,\ell }$, $\uA_{p,\ell }$ is a known piece of data (at least approximately) because
\begin{itemize}
    \item the functions $\mcll_{d,n}$ are described in \cref{convergence des mcl,def: mcl},
    \item the fields $u_{p',\ell '}$ are built inductively w.r.t.\ $p'$ so the coefficients $\sigma _{d} (\uDA_{p-n,\ell } )$ with $n>0$ are known,
    \item the coefficients $\sigma _{d} (\uDA_{p,\ell } )$ with $d<0$ are given by the last equation in \cref{eq: numerique CL rappel eq CL}.
\end{itemize}
On the contrary, $\uV_{p,\ell }$ is unknown and is a solution in $H^{1} (\Omega )$ of \cref{eq: numerique CL pb Helmholtz homog} for a certain $f$. Therefore, we have $\partial _{r} \uV_{p,\ell |\SME} =\opF(\uV_{p,\ell |\SME} )$, which yields:
\begin{equation}{\label{eq: condition DtN CL}}
    \partial _{r} u_{p,\ell |\SME} - \opF(u_{p,\ell |\SME} ) = \partial _{r} \uA_{p,\ell |\SME} - \opF(\uA_{p,\ell |\SME} ).
\end{equation}
This gives a non-homogeneous DtN condition for $u_{p,\ell }$ on the boundary $\SME$. We use it to compute $u_{p,\ell }$ in the domain $\OME$. Numerically, this condition can be approximated by truncating the series that defines $\opF$ in \cref{eq: def opF}.\\

Moreover, for~any $d\in \TN$, we have $\mcll_{d,0} = (-1)^{d\Theta /\pi } \,\Gamma (d+1) \big( \frac{2}{k_{0}} \big)^{d} J_{d} (k_{0} r) \sin (d \theta )$ by \cref{expression sing homog CL}. So, when applying \cref{eq: decompo modale solution H1 Helmholtz} to $u:=\uV_{p,\ell }$, we have
\begin{equation}{\label{eq: lien coef Fourier CL et sigma}}
    \forall n\in \mathbb{N}^*, \qquad  s_{n} = (-1)^{n} \,\Gamma \!\left(n\pst+1\right) \left( \frac{2}{k_{0}} \right)^{n\pst} \sigma _{n\pst} (\uDA_{p,\ell } ).
\end{equation}
\Cref{eq: lien coef Fourier CL et sigma,eq: formule coef Fourier solution H1 Helmholtz} allow us to compute the coefficient $\sigma _{d} (\uDA_{p,\ell } )$ for~any $d\in \TN$ once $u_{p,\ell }$ is computed.\\

Thus, the DtN condition \cref{eq: condition DtN CL} addresses at once and in a neat way three issues:
\begin{itemize}
    \item It handles the fact that $u_{p,\ell }$ is in general a non variational solution of Helmholtz's equation, due to its diverging behavior at the corner. Note that this does~not require introducing an auxiliary variational unknown (like we did for the corner profiles), as the condition \cref{eq: condition DtN CL} is directly formulated in terms of $u_{p,\ell }$, not in terms of $\uV_{p,\ell }$.
    \item It allows a good approximation of $u_{p,\ell }$ at the corner despite its lack of regularity.
    \item It provides a way to compute the coefficients $\sigma _{d} (u_{p,\ell } )$ with $d\in \TN$.
\end{itemize}
\end{paragraphesansboite}


\begin{paragraphesansboite}{Handling the Dirichlet condition of the far fields:}
\def\approxh#1[#2]{#1[#2#1]^h}
\def\msh{X^h}%
\def\ptmsh{\mathbf{x}}%
\def\ptmshbis{\mathbf{x'}}%
This Dirichlet condition on $\Gamma $ involves the functions $(\partial _{x} ^2 +k_{1} ^2 )^{n} \partial _{y} u_{p,\ell |\Gamma }$, see \cref{eq: numerique CL rappel eq CL}. Let us denote $h$ the mesh step, $\Vh$ the Galerkin space used for the far fields (e.g.\ the Lagrange finite elements P1 on $\OME$) and $\Vhg := \{v_{|\GME} \mid v\in \Vh\}$. For~any function $f$ that is defined in the continuous setting, we denote $\approxh[f]$ its discrete approximation, which will be properly defined  case by case. In particular, to shorten notations, we denote, for any $(p,\ell ,n)\in \NN\times \mathbb{N}^2 $, $\uh :=\approxh[u_{p,\ell } ]\in \Vh$ the computed far field and $\duh{n} :=\approxh\big[(\partial _{x} ^2 +k_{1} ^2 )^{n} \partial _{y} u_{p,\ell |\GME} ]$. We look for $\duh{n}$ in $\Vhg$, so that $\duh{n}$ can be directly used in the Dirichlet condition of a far field. Now, let us explain how $\duh{n}$ is computed from $\uh$, starting with the case $n=0$.\\

For~any $v\in H^{1/2} (\GME)$, we have:
\[\int _{\GME} \partial _{y} u_{p,\ell } \,v = -\int _{\OME} (\nabla  u_{p,\ell } \cdot  \nabla  w + \Delta  u_{p,\ell } \,w) + \int _{\partial \OME \setminus \GME} \partial _{n} u_{p,\ell } \,w \]
where $w\in H^{1} (\OME)$ is s.t.\ $w_{|\GME} =v$. Assuming that $w$ vanishes on $\partial \OME \setminus (\GME \cup \SME)$ and using that $u_{p,\ell }$ satisfies $\Delta u_{p,\ell } +k_{0} ^2  \.u_{p,\ell } =f_{p,\ell }$ with $f_{p,\ell } := \delta _{p,0} \.\delta _{\ell ,0} \.\mu _{0} ^{-1} \fsource$, we obtain:
\[\int _{\GME} \partial _{y} u_{p,\ell } \,v = \int _{\OME} \big({-\nabla  u_{p,\ell } }\cdot  \nabla  w + k_{0} ^2  \.u_{p,\ell } \.w + f_{p,\ell } \.w \big) - \int _{\SME} \partial _{r} u_{p,\ell } \,w.\]
This equation gives a way to compute an approximation of $\partial _{y} u_{p,\ell }$ by solving the following problem:
\begin{equation}{\label{eq: calcul dy upl Gamma methode 1}}
    \begin{multlined}
    \text{find $\duh0 \in \Vhg$ s.t.\ }  \forall v\in \Vhg, \qquad  \qquad \\
    \qquad  \qquad  \int _{\GME} \duh0 \,v = \int _{\Omega } \big({-\nabla  \uh}\cdot  \nabla  w + k_{0} ^2  \.\uh \.w + f_{p,\ell } \.w \big) - \int _{\SME} \approxh[\partial _{r} u_{p,\ell } ]\,w.
    \end{multlined}
\end{equation}
Here $w$ is taken in $\Vh$ and built in the following way, that is classical. Let us denote $\msh$ the set of vertices of the mesh and, for any $\ptmsh \in \msh$, $w_{\ptmsh}$ the canonical basis function of $\Vh$ associated to $\ptmsh$ (i.e.\ the unique element of $\Vh$ s.t.: $\forall \ptmshbis \in \msh,\, w_{\ptmsh} (\ptmshbis) = \delta _{\ptmsh,\ptmshbis}$). We define
\[w := \sum _{\ptmsh \in \msh \cap \GME} v(\ptmsh) \,w_{\ptmsh} .\]
In addition, in \cref{eq: calcul dy upl Gamma methode 1}, $\approxh[\partial _{r} u_{p,\ell } ]$ is computed thanks to the DtN condition \cref{eq: condition DtN CL}.\\

Actually, we use a variant of \cref{eq: calcul dy upl Gamma methode 1} to compute $\duh0$. The computation of high order derivatives of $\uh$ introduces numerical instabilities, as any small parasite oscillation at the scale of the mesh step gets amplificated when taking derivatives. To solve this issue, we introduced a regularization in the computation of $\duh{n}$, by solving instead the following problem
\begin{equation}{\label{eq: calcul dy upl Gamma methode 2}}
    \begin{array}{l}
    \text{find $\duh0 \in \Vhg$ s.t.\ }  \forall v\in \Vhg, \\[1.3\jot]
    \displaystyle \int _{\GME} \big(\red{h^2 }\.\partial _{x} \duh0 \,\partial _{x} v + \duh0 \,v \big)= \int _{\OME} \big({-\nabla  \uh}\cdot  \nabla  w + k_{0} ^2  \.\uh \.w + f_{p,\ell } \.w \big) - \int _{\SME} \approxh[\partial _{r} u_{p,\ell } ]\,w. \hspace{-1em}
\end{array} \end{equation}
This regularization reduces the high frequency oscillations of $\duh0$. The factor $h^2 $ in \cref{eq: calcul dy upl Gamma methode 2} is a relevant regularization parameter because, by dimensional analysis, the parameter in front of $\partial _{x} \duh0 \,\partial _{x} v$ must be the square of a length. See \cite{LazSig09,LazSig11} for similar filters used in topology optimization and \cite{Sag06} in the context of large eddy simulations.\\

Next, let $n\geqslant 1$. For~any $u\in H^2 (\GME)$ and $v\in H^{1} (\GME)$, we have
\[\int _{\GME} (\partial _{x} ^2 +k_{1} ^2 )u \cdot v = \int _{\GME} \big({-\partial _{x} u} \, \partial _{x} ^2 v +k_{1} ^2 \.u\.v\big) - \partial _{x} u(\PME) \cdot v(\PME)\]
with $\PME :=(\RME,0)$ the tip of $\GME$. Applying this to $u:=(\partial _{x} ^2 +k_{1} ^2 )^{n-1} \partial _{y} u_{p,\ell |\Gamma }$, and employing the same regularization technique as in \cref{eq: calcul dy upl Gamma methode 2}, we obtain the following problem which we use to compute inductively the functions $\duh{n}$ with $n\geqslant 1$
\begin{equation}{\label{eq: calcul (dx+k)n dy upl}}
    \begin{multlined}    
    \text{find $\duh{n} \in \Vhg$ s.t.\ }  \forall v\in \Vhg, \quad   \int _{\GME} \big(4h^2 \partial _{x} \duh{n} \,\partial _{x} v + \duh{n} \,v \big) = \qquad  \qquad  \qquad  \quad \ \\
    \qquad  \qquad  \int _{\GME} \big({-\partial _{x} \duh{n-1}} \, \partial _{x} ^2 v +k_{1} ^2 \.\duh{n-1}\.v\big) - \approxh\big[\partial _{x} (\partial _{x} ^2 +k_{1} ^2 )^{n-1} \partial _{y} u_{p,\ell } (\PME)] \cdot v(\PME) \hspace*{-1em}
    \end{multlined}
\end{equation}
The factor 4 in the right-hand side of \cref{eq: calcul (dx+k)n dy upl} seemed to give the best results in practice. Moreover, $\approxh\big[\partial _{x} (\partial _{x} ^2 +k_{1} ^2 )^{n-1} \partial _{y} u_{p,\ell } (\PME)]$ is calculated in the following way. As~explained before, \cref{eq: lien coef Fourier CL et sigma,eq: formule coef Fourier solution H1 Helmholtz} give a way to compute the values of the coefficients $\sigma _{d} (u_{p,\ell } )$. Thanks to a truncation of the modal decomposition \cref{eq: rappel decompo modale CL pour numerique CL}, this then allows us to reconstitute an approximation of $u_{p,\ell }$ in $\Omega  \cap B(0,\RME)$, i.e.\ outside the computation domain $\OME$. Let us denote $\approxh[\uDA_{p,\ell } ]$ this approximate reconstitution. It can be:
\begin{itemize}
    \item either a linear combination of Bessel functions, which is possible if and only if the functions $\mcll_{d,n}$ that appear can be expressed using Bessel functions (see \cref{sec: expression des mcl,sec: annexe des complements sur le DA}),
    \item or a linear combination of truncations of the form $T_{<q} (\mcla_{d,n} )$, which is always possible but can be less easy to implement. In this case, $\approxh[\uDA_{p,\ell } ]$ is an approximation of $T_{<q} (\uDA_{p,\ell } )$ (see \cref{expression uDA avec mcl}). In addition, it belongs to $\mathcal{A}(\Omega )$, so it can be numerically represented as explained in \cref{sec: implementation des espaces A}. We used this second option because \cref{chap: complements sur le DA} was chronologically posterior to these numerical experiments and this options allows to unified implementation for any angle $\Theta $.
\end{itemize}
In both cases, we define $\approxh\big[\partial _{x} (\partial _{x} ^2 +k_{1} ^2 )^{n-1} \partial _{y} u_{p,\ell } (\PME)]:= \partial _{x} (\partial _{x} ^2 +k_{1} ^2 )^{n-1} \partial _{y} \approxh[\uDA_{p,\ell } ](\PME)$. This derivative can be computed exactly from $\approxh[\uDA_{p,\ell } ]$ as $\approxh[\uDA_{p,\ell } ]$ has an explicit expression in both cases.
\end{paragraphesansboite}

\begin{remarquenumerotee}
\label{rq: cumul d'erreurs numerique sur CL}%
The main drawback in using the far fields to compute an approximation of the exact solution $u_{\varepsilon }$ is that their inductive construction generates accumulated numerical errors. This is made worse by the far fields' Dirichlet condition on $\Gamma $, as the way we handle it in \cref{eq: calcul dy upl Gamma methode 2,eq: calcul (dx+k)n dy upl} involves a projection on the finite element space $\Vhg$ and a numerical regularization parameter linked to the mesh step.\\

On the contrary, the effective models proposed in \cref{chap: ME} directly provide an approximate solution, which reduces the computational cost and numerical errors. In addition, the boundary condition of the effective model of order~3, resp.\ 5, is of Robin type, resp.\ Ventcel, which is much easier to handle numerically. Hence computing the far fields is interesting to validate the error estimate of \cref{DA de SE avec CL et delta variable}, but we will see that, in practice, effective models are more effective to compute an approximate solution.
\end{remarquenumerotee}

\subsection{Numerical results}
\label{sec: resultats numeriques CL}

\begin{textesansboite}
We verified numerically the convergence of the asymptotic expansion of $u_{\varepsilon }$ up to order~5. To our knowledge, this is the first time that an asymptotic expansion is computed further than order 2 for a problem in a corner with a layer. We consider three different angles: $\Theta =2.5$, $\Theta =\pi $ and $\Theta =4$. This is a representative sample as we have small and large angles (the latter ones having more far fields, corner profiles and far field corner conditions at a given order), as well as the two cases $\pst \in \mathbb{Q}$ and $\pst \not\in \mathbb{Q}$ which lead to different types of expansion (see \cref{sec: simplifications DA}). A detailed description of the expansion's structure is given in \cref{sec: tableaux exemples DA} for these angles (up to order 4).\\

	\tikzset{external/export=true}%
	\tikzsetfigurename{fig-couche_arrondie-}%
	\begin{figure}[h!]
\def\epsFig{3.2} 
\def\xMax{4.5}
\def\xMin{-2.2}
\def\yMax{1.6}
\def\yMin{-3.35}
\def\ThetaFig{120} 

\begin{tikzpicture}[scale=1]
    \begin{scope}
        \clip (\xMin,\yMin) rectangle (\xMax,\yMax);
        \coordinate (finLambda) at (\xMin, {\xMin*tan(\ThetaFig-180)});
        \fill [pattern={Lines[angle=45,line width=.2pt]}, pattern color=white!70!gray] (0,0) -- (\xMax,0) -- (\xMax,\yMax) -- (\xMin,\yMax) -- (finLambda) -- cycle;
        \draw [thick] (finLambda) -- (0,0) -- (\xMax,0);
    \end{scope}
    \fill [color=white!70!gray] (0,0) -- (\xMax,0) -- (\xMax,-\epsFig) -- (0.8*\epsFig,-\epsFig) arc (270:180:0.6*\epsFig) arc (0:53.13:0.5*\epsFig) -- cycle;
    \draw[gray, opacity=0.75, thin, step=0.1*\epsFig] (-0.5*\epsFig-0.2,-\epsFig-0.2) grid (\epsFig+0.2,0.1);
    \draw [thick] (\xMax,-\epsFig) -- (0.8*\epsFig,-\epsFig) arc (270:180:0.6*\epsFig) arc (0:53.13:0.5*\epsFig);
	\draw [>=stealth,->] (\xMin-0.3,0) -- (\xMax+0.3,0) node [right] {$X$};
	\draw [>=stealth,->] (0,{\yMin-0.3}) -- (0,{\yMax+0.3}) node [right] {$Y$};
	\draw (\xMax*0.5,\yMax*0.6) node {$(\mu_0,\rho_0)$};
	\draw [->,>=latex] (0.35,0) arc (0:\ThetaFig-5:0.35) node[midway,above right] {$\Theta$} -- (\ThetaFig:0.36);
	\draw (0.82*\xMax,-0.5*\epsFig) node {layer};
	\draw (0.82*\xMax,-0.63*\epsFig) node{$(\mu_1,\rho_1)$};
    \draw ( 0.5*\epsFig,0) node{\bf+} ++(0,0.05) node[above]{0.5};
    \draw (-0.5*\epsFig,0) node{\bf+} ++(0,0.05) node[above]{$-0.5$};
    \draw (     \epsFig,0) node{\bf+} ++(0,0.05) node[above]{1};
    \draw (0,-0.5*\epsFig) node{\bf+} ++(-0.05,0) node[left]{0.5};
    \draw (0,    -\epsFig) node{\bf+} ++(-0.05,0) node[left]{1};
    \draw (0,0) node{$\bullet$};
    \draw (0.8*\epsFig,-\epsFig) node{$\bullet$};
    \draw (0.2*\epsFig,-0.4*\epsFig) node{$\bullet$};
    \draw[blue] (0.8*\epsFig,-0.4*\epsFig) node{$\bullet$};
    \draw[blue] (-0.3*\epsFig,-0.4*\epsFig) node{$\bullet$};
\end{tikzpicture}

\longCaption{The shape of the layer used for the simulations}{. The layer's boundary contains two circular arcs centered at the blue points.}
\label{fig: forme couche simus}
\end{figure}%
	\tikzset{external/export=false}%
	\tikzsetfigurename{main-figure}

We consider a domain $\Omega _{1}$ with a rounded layer at the corner as in \cref{fig: forme couche simus}. We also consider a compactly supported source term $\fsource$ that is a small bump centered somewhere on the line $\theta =0.4\cdot \Theta $. This avoids any apparent symmetry. In addition, we take $\mu _{0} =1$, $\mu _{1} =2$, $\rho _{0} =1.2$, $\rho _{1} =3$ and $\omega =5+0.2\mathrm{i}$. The radius of the DtN boundary is $\RME=0.15\cdot \lambda _{0}$, with $\lambda _{0} :=2\pi /\Re(k_{0} )\approx 1.15$ the wave length in $\Omega $. We consider the layer thicknesses $\varepsilon  \in \{0.04, 0.065, 0.1, 0.15, 0.25\}$ (which are roughly equally spaced in log scale). We use P1 Lagrange finite elements, because elements of higher order seem to badly handle the computation of $\partial _{y} u_{p,\ell }$ by \cref{eq: calcul dy upl Gamma methode 2} for reasons that we do~not understand and did not investigate in depth. The reference solutions $u_{\varepsilon }$ are computed by finely meshing the layer.\\

We use the FEM-BEM library XLiFE++ \cite{XLiFE}. It is a \Cplusplus{} library, which is why we wrote the symbolic computation module of \cref{sec: implementation des espaces A} in \Cplusplus{} in the first place. Moreover, we parallelized the computations as much as possible to reduce computational time:
\begin{itemize}
    \item One process computes the far fields $u_{p,\ell }$ inductively w.r.t.\ $p$.
    \item Several processes compute the corner profiles $\mathcal{S}_{d,n}$ inductively w.r.t.\ $n$, one process per value of $d$. The list of the needed corner profiles is determined by the program according to the angle and desired order of the expansion.
    \item Several processes compute the exact solutions $u_{\varepsilon }$, one per value of $\varepsilon $.
    \item One process receives the far fields and exact solutions and gradually computes the errors.
\end{itemize}
The computation of the corner singularities and matching coefficients is so quick that it does~not need to be parallelized.\\

As explained in \cref{eq: parallelisation des calculs}, the first process may wait sometimes for a corner profile coefficients. In order to handle these situations, and make it go back to work as soon as the needed coefficient is available, I developed a \emph{task dispatching library} in \Cplusplus{} using MPI. It relies on a manager process that decides what each process does (e.g.\ it can import a previously computed corner profile if available, rather than using a process to recompute it), and centralizes the communications between processes. It can also indicate the remaining computational time and pause some processes if the RAM may become full. These features may not be essential, but they helped me save time in my simulations. This library has been published on git \cite{OPTEC} under the name OPTEC (for ``Optimized Parallelized Task Execution in \Cplusplus''). See \url{https://plmlab.math.cnrs.fr/cbaudet/optec}. A complete presentation of this library and a documentation are available at this url.\\

Moreover, our code can handle any angle and any expansion order (theoretically, as numerical errors accumulate at high orders). This is made possible especially by the generality of the symbolic computation module of \cref{sec: implementation des espaces A}.
\end{textesansboite}

\newcommand{\rda}[1][p]{R_{#1,\peps}}
\begin{textesansboite}
\Cref{tab: CL angle 4}, resp.\ \cref{tab: CL angle pi}, shows the real parts of the first far fields and the remainders
\[\rda := u_{\varepsilon } -\sum _{p'\red{<p}} \sum _{\ell } \varepsilon ^{p'} \ln ^{\ell } \!\varepsilon \,u_{p',\ell } \vspace{-1.5pt}\]
for $\Theta =4$, resp.\ $\Theta =\pi $ (without the PML). In both cases, we observe that the larger $p$, the smaller the remainder $\rda$.\\

Let us first comment $\Theta =4$ (\cref{tab: CL angle 4}). We recall that, since $\frac{\pi }{4} \not\in \mathbb{Q}$, there~is no $\ln \varepsilon $ in this case (\cref{simplifications cas irrationnel}). In addition, by \cref{DA de SE avec CL et delta fixe} we have $u_{\varepsilon } -\sum _{p'\leqslant p} \varepsilon ^{p'} u_{p',0} =o(\varepsilon ^{p} )$ for~any $p\in \NN$, which implies
\[\varepsilon ^{-p} \.\rda = \varepsilon ^{-p} \.\bigg(u_{\varepsilon } -\sum _{p'<p} \varepsilon ^{p'} u_{p',0} \bigg) \underset{\varepsilon \rightarrow 0}{\longrightarrow} u_{p,0} .\]
Thus each far field $u_{p,0}$ should look similar to the associated remainder $\rda$. This can be observed in \cref{tab: CL angle 4} and allows to validate qualitatively the computations. It is especially clear for $p\in \{0,1,2\}$. As~for $p=2\pst$, $\rda$ looks like a superposition of $u_{p,0}$ and the next far field $u_{2,0}$. Finally, $u_{3\pst,0}$ looks different from $\rda[3\pst]$, which seems linked to the plateaux in the convergence curves of \cref{fig: graphe err CL 2.5}.\\

Moreover, there~is a visible difference between $\{u_{1,0} , u_{2,0} \}$ and $\{u_{2\pst} , u_{3\pst} \}$. The former have a non-homogeneous Dirichlet condition on $\Gamma $ and no corner condition (i.e.\ $\sigma _{d} (u_{p,0} )=0$ for~any $d\in -\TN$), whereas the latter have a homogeneous condition on $\Gamma $ and a non-trivial corner condition. Note that some higher order far fields have at the same time a non-homogeneous condition on $\Gamma $ and corner conditions, e.g.\ $u_{1+2\pst,0}$. Besides, $u_{2\pst}$, and $u_{3\pst}$ take their largest values near the corner, which is coherent with the fact they theoretically diverge at the corner (generically). Finally, $u_{2\pst}$ has one bump near the corner and globally looks like $J_{-\pst} (k_{0} r)\sin (\pst \theta )$, because only $\sigma _{-\pst} (u_{2\pst} )$ is non-zero, while $u_{3\pst}$ has two bumps because $\sigma _{-\red{2}\pst} (u_{3\pst} ) \neq 0$.\\

When $\Theta =\pi $ (\cref{tab: CL angle pi}), the expansion contains powers of $\ln \varepsilon $. Thus, each remainder $\rda$ should be close to a superposition of the far fields $u_{p,\ell }$ because we have
\[\varepsilon ^{-p} \.\rda = \.\sum _{\ell \in \mathbb{N}} \.\ln ^{\ell } \!\varepsilon \,u_{p,\ell } +o(1) \qquad  \text{when }  \varepsilon \rightarrow 0.\]
This can be observed in \cref{tab: CL angle pi} (in practice, $\rda$ is especially close to $u_{p,0}$ here). Again, the far fields with a non-homogeneous condition on $\Gamma $ look very different from the ones with only a corner condition (which happens only for $\ell >0$ when $\Theta =\pi $).
\end{textesansboite}

	\tikzset{external/export=true}%
	\tikzsetfigurename{tab-CL_4-}%
	\begin{table}[h]
\def\largeurimages{44mm}
\def\boite#1{\raisebox{-.5\height}{\hskip1mm \vbox{\vskip2.5mm\relax \hbox{#1} \vskip2mm}\!}}%
\def\addgrapheCL#1{\boite{\includegraphics[width=\largeurimages]{images/angle4_CL#1.jpg}}}%
\def\addgrapheErr#1{\boite{\includegraphics[width=\largeurimages]{images/angle4_errSE-CL#1_eps0,065.jpg}}}%
\def\colorbarheight{20mm}
\renewcommand{\arraystretch}{1.8}%

\def\imageSE{\boite{%
    \dimen0=\largeurimages\relax
    \dimen0=1.375\dimen0\relax
    \edef\tempo{\the\dimen0}%
    \begin{tikzpicture}[x=\largeurimages, y=\largeurimages]
        \clip (0.5,-0.1) -- (0.5,0) arc (0:{4*180/3.14}:0.5) -- cycle;
        \draw (0,0.085) node{\includegraphics[width=\tempo]{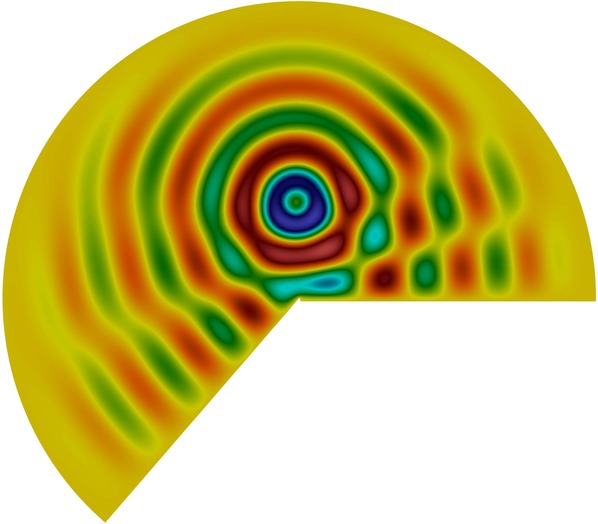}};
    \end{tikzpicture}%
}}

\caption{Real part of the first far fields and the errors w.r.t.\ the exact solution for $\Theta=4$ and $\peps=0.065$}
\label{tab: CL angle 4}

\makebox[\textwidth][c]{
\begin{tabular}{|c|cl<{\!}|cl<{\!}|}
    \hline
    $p$ &
    \multicolumn{2}{c|}{Far field $u_{p,0}$ (NB: $u_{p,\ell} =0$ for $\ell>0$)} & 
    \multicolumn{2}{c|}{Error $\displaystyle u_\peps - \sum_{p'<p\mathstrut} \peps^{p'} u_{p',0}$} \\
    \hline
    0 &
    \addgrapheCL{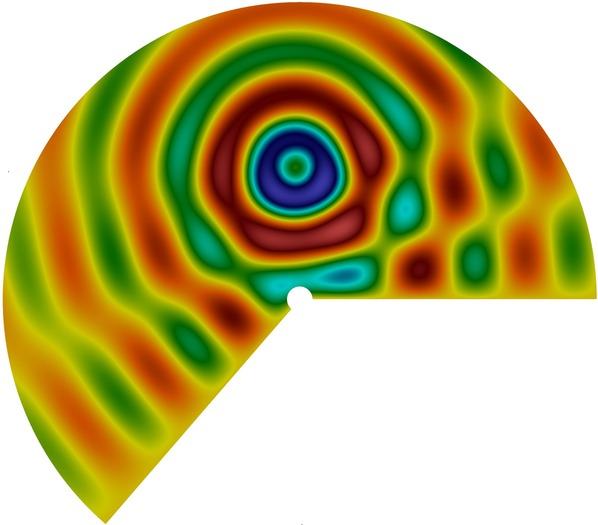} & \colorbar[min=-2, max=1.4, 1e-2] &
    \imageSE & \colorbar[min=-2, max=1.4, 1e-2] \\[-3ex]
    & & & \multicolumn{2}{c|}{(exact solution)} \\
    \hline
    1 &
    \addgrapheCL{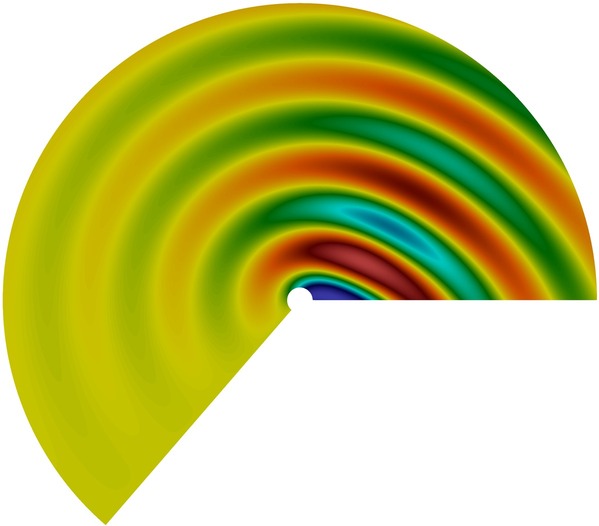} & \colorbar[min=-3.2, max=2.5, 1e-2] &
    \addgrapheErr{1_p0} & \colorbar[min=-2.2, max=1.7, 1e-3] \\
    \hline
    $2\pst \approx 1.6$ &
    \addgrapheCL{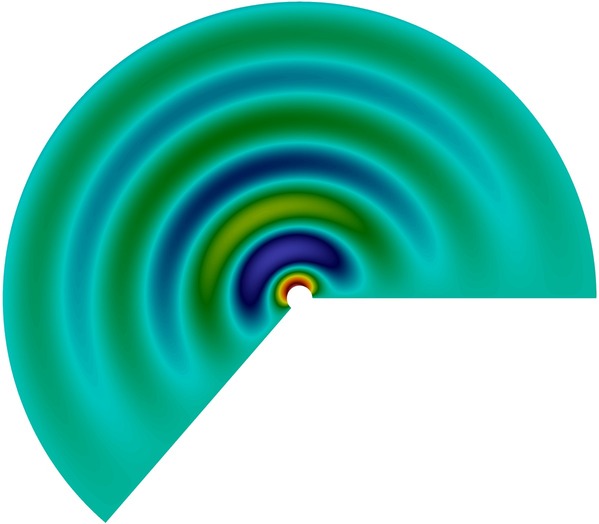} & \colorbar[min=-1.6, max=3.6, 1e-2] &
    \addgrapheErr{2_p1} & \colorbar[min=-2.5, max=5, 1e-4] \\
    \hline
    $2$ &
    \addgrapheCL{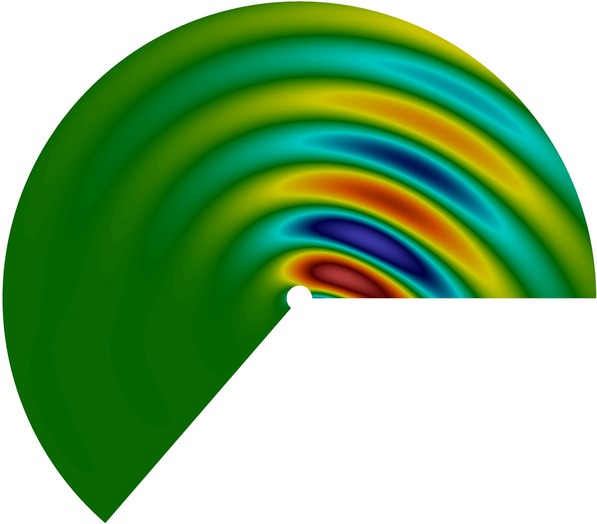} & \colorbar[min=-5.6, max=7.3, 1e-2] &
    \addgrapheErr{3_p1,57} & \colorbar[min=-2.5, max=5, 1e-4] \\
    \hline
    $3\pst \approx 2.4$ &
    \addgrapheCL{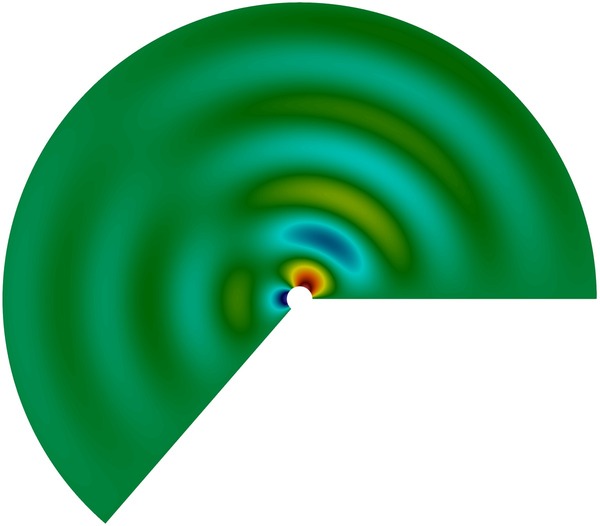} & \colorbar[min=-7.9, max=12, 1e-2] &
    \addgrapheErr{4_p2} & \colorbar[min=-6.1, max=10, 1e-5] \\
    \hline
\end{tabular}}
\end{table}%
	\tikzset{external/export=false}%
	\tikzsetfigurename{main-figure}
	\tikzset{external/export=true}%
	\tikzsetfigurename{tab-CL_pi-}%
	\begin{table}[h]
\def\largeurimages{52mm}
\def\boite#1{\raisebox{-.5\height}{\hskip1mm \vbox{\vskip3.5mm\relax \hbox{#1} \vskip4mm}\!}}%
\def\addgrapheCL#1{\boite{\includegraphics[width=\largeurimages]{images/anglepi_CL#1.jpg}}}%
\def\addgrapheErr#1{\boite{\includegraphics[width=\largeurimages]{images/anglepi_errSE-CL#1_eps0,01.jpg}}}%
\def\colorbarheight{20mm}
\renewcommand{\arraystretch}{1.8}%

\def\imageSE{\boite{%
    \dimen0=\largeurimages\relax
    \dimen0=1.375\dimen0\relax
    \edef\tempo{\the\dimen0}%
    \begin{tikzpicture}[x=\largeurimages, y=\largeurimages]
        \clip (-0.5,{-0.5*0.1/4}) -- (0.5,{-0.5*0.1/4}) -- (0.5,0) arc (0:180:0.5) -- cycle;
        \draw (0,0.335) node{\includegraphics[width=\tempo]{images/anglepi_SE_eps0,1.jpg}};
    \end{tikzpicture}%
}}

\caption{Real part of the first far fields and the errors w.r.t.\ the exact solution for $\Theta=\pi$ and $\peps=0.1$}
\label{tab: CL angle pi}

\makebox[\textwidth][c]{
\begin{tabular}{|c|c|cl<{\!}|cl<{\!}|}
    \hline
    \;$p$\; & \;$\ell$\; &
    \multicolumn{2}{c|}{Far field $u_{p,\ell}$} & 
    \multicolumn{2}{c|}{Error $\displaystyle u_\peps - \sum_{p'<p\mathstrut} \sum_{\ell'} \peps^{p'} \ln^{\ell'}\!\peps \,u_{p',0}$} \\
    \hline
    0 & 0 &
    \addgrapheCL{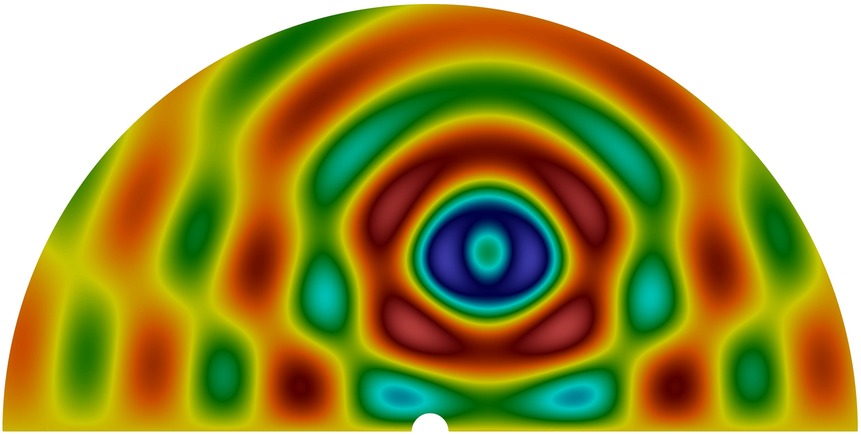} & \colorbar[min=-2.1, max=1.4, 1e-2] &
    \imageSE & \colorbar[min=-2.1, max=1.4, 1e-2] \\[-2.5ex]
    & & & & \multicolumn{2}{c|}{(exact solution)} \\
    \hline
    1 & 0 &
    \addgrapheCL{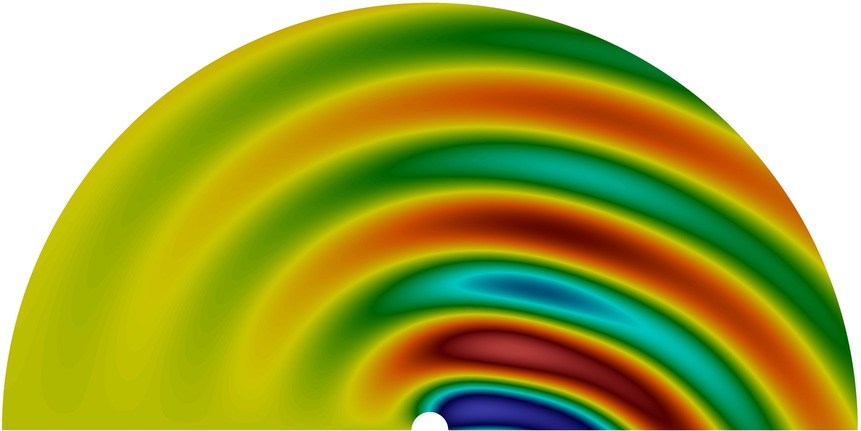} & \colorbar[min=-3.5, max=2.8, 1e-2] &
    \addgrapheErr{1_p0} & \colorbar[min=-3.9, max=3.1, 1e-3] \\
    \hline
    \multirow{2}{*}[-13mm]{2} & 0 &
    \addgrapheCL{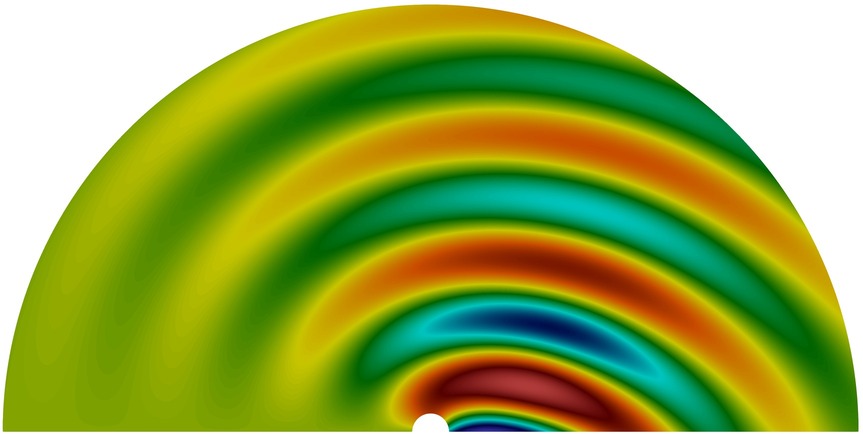} & \colorbar[min=-8.8, max=8.1, 1e-2] &
    \multirow{2}{*}[-5mm]{\addgrapheErr{2_p1}} & \multirow{2}{*}[-7mm]{\colorbar[min=-11, max=9.6, 1e-4]} \\
    \cline{2-4}
    & 1 &
    \addgrapheCL{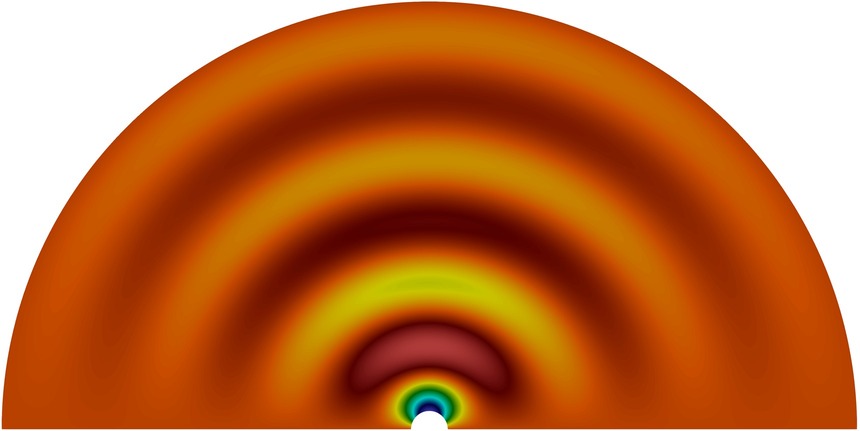} & \colorbar[min=-13, max=4.8, 1e-3] & & \\
    \hline
    \multirow{3}{*}[-24.5mm]{3} & 0 &
    \addgrapheCL{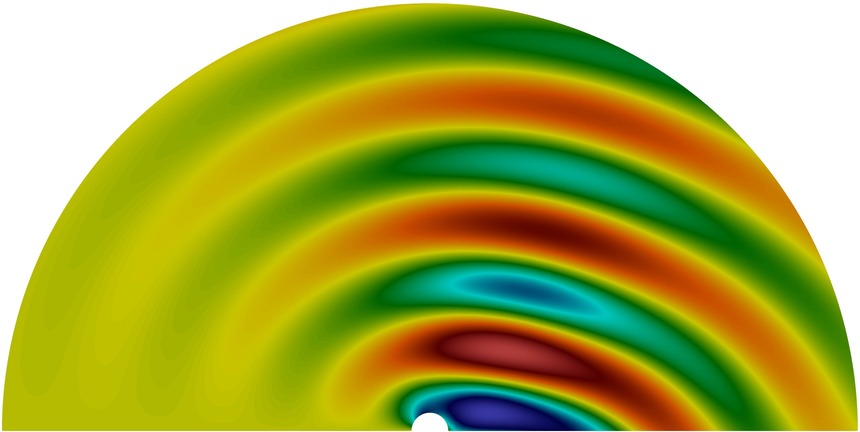} & \colorbar[min=-1.9, max=1.5, 1e-1] &
    \multirow{3}{*}[-22mm]{\addgrapheErr{4_p2}} & \multirow{3}{*}[-25mm]{\colorbar[min=-2.7, max=2.2, 1e-4]} \\
    \cline{2-4}
    & 1 &
    \addgrapheCL{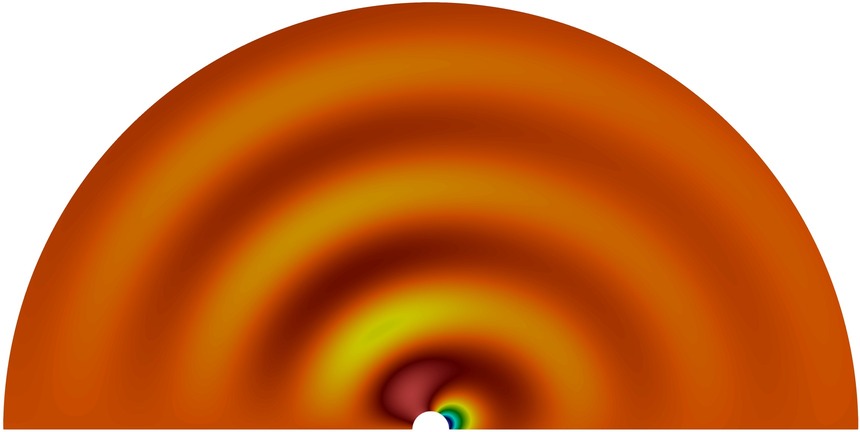} & \colorbar[min=-4.5, max=1.8, 1e-2] & & \\
    \cline{2-4}
    & 2 &
    \addgrapheCL{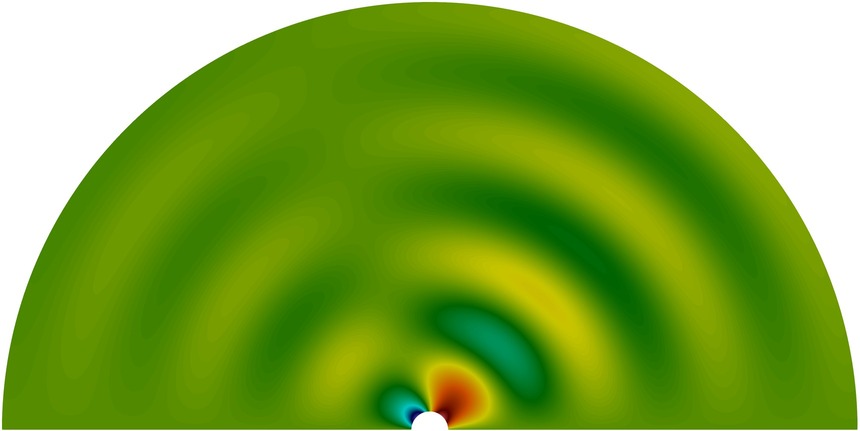} & \colorbar[min=-4.9, max=5.1, 1e-3] & & \\
    \hline
\end{tabular}}
\end{table}%
	\tikzset{external/export=false}%
	\tikzsetfigurename{main-figure}

\begin{textesansboite}
On a quantitative level, \crefrange{fig: graphe err CL 2.5}{fig: graphe err CL 4} show convergence curves. For each $p$, we plot the relative error $\|\rda\|_{H^{1}} / \|u_{0,0} \|_{H^{1}}$ w.r.t.\ $\varepsilon $, for the $H^{1}$~norm on the computational domain (except the PML). $u_{0,0}$ is used as a normalization as it~is independent from $\varepsilon $, unlike $u_{\varepsilon }$. We plot the linear regressions with solid lines of the same color as the associated points. The curves have plateaux that we believe are linked to the numerical errors mentioned in \cref{rq: cumul d'erreurs numerique sur CL}. We can still obtain the right order of convergence by making the linear regressions only on the first points. Some curves for consecutive values of $p$ almost coincide, because some far fields appear to be small and contribute little to the expansion for the parameters we use.
\end{textesansboite}

	\tikzset{external/export=true}%
	\tikzsetfigurename{fig-erreurs_CL_2,5-}%
	\pgfplotsset{
    every axis/.append style={
        thick,
    },
    cycle list name=Dark2_CL_2.5,
    /pgfplots/table/x=epsilon,
    legend cell align=left,
}

\begin{figure}[h]

\begin{tikzpicture}
\begin{axisCL}{err_2,5_tot.dat}
    \courbeCL[no=1, p=1, reg=0:5]
    \courbeCL[no=2, p=2, reg=0:5]
    \courbeCL[no=3, p=2\pst\approx2.51, reg=0:5]
    \courbeCL[no=4, p=3, reg=0:5]
    \courbeCL[no=5, p=1+2\pst\approx3.51, reg=1:5]
    \courbeCL[no=6, p=3\pst\approx3.77, reg=1:5]
    \courbeCL[no=7, p=4, reg=1:5]
    \courbeCL[no=8, p=2+2\pst\approx4.51, reg=2:5]
    \courbeCL[no=9, p=1+3\pst\approx4.77, reg=2:5]
    \courbeCL[no=10, p=5, reg=2:5]
\end{axisCL}
    \nodepenteCL[nodbt=1, nb=1]
    \nodepenteCL[nodbt=2, nb=1]
    \nodepenteCL[nodbt=3, nb=2]
    \nodepenteCL[nodbt=5, nb=3]
    \nodepenteCL[nodbt=8, nb=3]
\end{tikzpicture}
\caption{
    Numerical errors of the asymptotic expansion for $\Theta=2.5$ (for the $H^1$~norm on $\OME$) 
}
\label{fig: graphe err CL 2.5}
\end{figure}%
	\tikzset{external/export=false}%
	\tikzsetfigurename{main-figure}
	\tikzset{external/export=true}%
	\tikzsetfigurename{fig-erreurs_CL_pi-}%
	\pgfplotsset{
    every axis/.append style={
        thick,
    },
    cycle list name=Dark2_avec_regressions,
    /pgfplots/table/x=epsilon,
    legend cell align=left,
}

\begin{figure}[h]
\begin{tikzpicture} 
\begin{axisCL}{err_pi_tot.dat}
\grapheCLpcourt
    \courbeCL[no=1, p=1, reg=0:5]
    \courbeCL[no=2, p=2, reg=0:5]
    \courbeCL[no=3, p=3, reg=0:5]
    \courbeCL[no=4, p=4, reg=1:5]
    \courbeCL[no=5, p=5, reg=2:5]
\end{axisCL}
\grapheCLpcourt
    \nodepenteCL[nodbt=1, nb=1]
    \nodepenteCL[nodbt=2, nb=1]
    \nodepenteCL[nodbt=3, nb=1]
    \nodepenteCL[nodbt=4, nb=1]
    \nodepenteCL[nodbt=5, nb=1]
\end{tikzpicture}

\caption{Numerical errors of the asymptotic expansion for $\Theta=\pi$ (for the $H^1$~norm on $\OME$)}
\label{fig: graphe err CL pi}
\end{figure}%
	\tikzset{external/export=false}%
	\tikzsetfigurename{main-figure}
	\tikzset{external/export=true}%
	\tikzsetfigurename{fig-erreurs_CL_4-}%
	\pgfplotsset{
    every axis/.append style={
        thick,
    },
    cycle list name=Dark2_CL_4a,
    /pgfplots/table/x=epsilon,
    legend cell align=left,
}

\begin{figure}[h]

\begin{tikzpicture}
\begin{axisCL}{err_4_tot.dat}
    \courbeCL[no=1, p=1, reg=0:5]
    \courbeCL[no=2, p=2\pst\approx1.57, reg=0:5]
    \courbeCL[no=3, p=2, reg=0:5]
    \courbeCL[no=4, p=3\pst\approx2.36, reg=1:5]
    \courbeCL[no=5, p=1+2\pst\approx2.57, reg=1:5]
    \courbeCL[no=6, p=3, reg=1:5]
    \courbeCL[no=7, p=4\pst\approx3.14, reg=1:5]
    \courbeCL[no=8, p=1+3\pst\approx3.36, reg=1:5]
    \courbeCL[no=9, p=2+2\pst\approx3.57, reg=2:5]
    \courbeCL[no=10, p=5\pst\approx3.93, reg=2:5]
    \courbeCL[no=11, p=4, reg=2:5]
\end{axisCL}
    \nodepenteCL[nodbt=1, nb=1]
    \nodepenteCL[nodbt=2, nb=2]
    \nodepenteCL[nodbt=4, nb=1]
    \nodepenteCL[nodbt=5, nb=1]
    \nodepenteCL[nodbt=6, nb=1]
    \nodepenteCL[nodbt=7, nb=2]
    \nodepenteCL[nodbt=9, nb=3]
\end{tikzpicture}
\smallskip\smallskip\smallskip

\pgfplotsset{
    cycle list name=Dark2_CL_4b,
}

\begin{tikzpicture}
\begin{axisCL}{err_4_tot.dat}
    \courbeCL[no=11, p=4, reg=2:5]
    \courbeCL[no=12, p=1+4\pst\approx4.14, reg=2:5]
    \courbeCL[no=13, p=2+3\pst\approx4.36, reg=2:5]
    \courbeCL[no=14, p=3+2\pst\approx4.57, reg=2:5]
    \courbeCL[no=15, p=6\pst\approx4.71, reg=2:5]
    \courbeCL[no=16, p=1+5\pst\approx4.93, reg=2:5]
    \courbeCL[no=17, p=5, reg=2:5]
\end{axisCL}
    \nodepenteCL[nodbt=11, nb=1]
    \nodepenteCL[nodbt=12, nb=1]
    \nodepenteCL[nodbt=13, nb=4]
    \nodepenteCL[nodbt=17, nb=1]
\end{tikzpicture}

\longCaption{Numerical errors of the asymptotic expansion for $\Theta=4$ (for the $H^1$~norm on $\OME$)}{. The graph is splitted in two to improve lisibility.}
\label{fig: graphe err CL 4}
\end{figure}%
	\tikzset{external/export=false}%
	\tikzsetfigurename{main-figure}

\chapter{High order effective models}
\label{chap: ME}
\begin{textesansboite}
In this chapter, we propose two effective models of order~3, resp.\ 5, when $\varepsilon \rightarrow 0$. They allow computing directly an approximation of the exact solution $u_{\varepsilon }$ instead of summing lots of far fields, which reduces by a lot the computational cost. For example, for $\Theta = \frac{3\pi }{2}$, 12 far fields are needed to obtain an approximation of $u_{\varepsilon }$ of order~5. When $\Theta  \rightarrow  2\pi ^{-}$, this number becomes 31. See \cref{fig: graphe nb de CL} for more examples. Note however that the number of corner profiles to be computed is~not reduced by using effective models. But these functions are independent of the frequency $\omega $ and the source term $\fsource$, so they can be computed once in advance.\\

The proposed models replace, far from the corner, the layer by an effective Robin, resp.\ Ventcel, boundary condition. 
In addition, %
we introduce an artificial boundary near the corner and impose a Dirichlet-to-Neumann (DtN) condition. The DtN map takes into account the behavior of the exact solution near the corner that was established in \cref{chap: raccord algebrique,chap: construction du DA}.\\

In \cref{sec: intro des ME}, we describe the effective models and state the main results of well-posedness and convergence. The construction of the DtN map requires special care and is explained in \cref{sec: construction du DtN des ME}. It relies on an approximate modal decomposition of the exact solution, that is introduced in \cref{sec: decompo modale de SE}. The well-posedness of the effective models is proven in \cref{sec: ME bien poses}, while their convergence is shown in \cref{sec: cvg ME}. Numerical results are provided in \cref{sec: numerique ME}. Finally, some technical proofs can be found in Appendices~\ref{sec: preuves lemmes de regularite elliptique ME} and \ref{sec: annexe comportement en N des ME}.
\end{textesansboite}

\section{Introduction}
\label{sec: intro des ME}

\subsection{Description of the models}

\begin{textesansboite}
Let us define the geometry. It is similar to the computational domain of the far fields in \cref{sec: numerique des CL}. Let $\RME>0$ be s.t.\ $\overline{B}(0,\RME) \cap  \operatorname{supp}(\fsource) = \varnothing $. We know that such a number exists because we assumed $\operatorname{dist}(\operatorname{supp}(\fsource),\Gamma )>0$. We define $\OME:= \Omega  \setminus \overline{B}(0,\RME)$, $\GME:= \Gamma  \setminus \overline{B}(0,\RME)$,  $\bordgaucheME:= \bordgauche \setminus \overline{B}(0,\RME)$, $\SME := \{(\RME \cos \theta , \RME \sin \theta ) \mid \theta \in (0,\Theta )\}$ and $\PME:=(\RME,0)\in \mathbb{R}^2 $. Note that $\partial \OME = \overline{\bordgaucheME \cup  \SME \cup  \GME}$. See \cref{fig: domaine ME}. The ``$\circ $'' exponents are intended to recall that a ball is excluded around the corner.\\
\nomenclature[27-00]{$\RME$}{Radius of the DtN artificial boundary used for the far fields and effective models}
\nomenclature[59-00]{$\OME$}{Domain of the effective models, $\Omega  \setminus \closedBall(0,\RME)$}
\nomenclature[38-00]{$\GME$}{$\Gamma  \cap  \{x>\RME\}$}
\nomenclature[53-00]{$\SME$}{Artificial boundary around the corner associated to a DtN condition}
\nomenclature[25-00]{$\PME$}{The point $(\RME,0)$}

\begin{figure}[h]
\def\widthFig{3.5} 
\def\heightFig{2.9} 
\def\ThetaFig{205} 
\def\RMEFig{1.3} 
\begin{tikzpicture}[scale=0.9]
	\coordinate (finLambda) at (-\widthFig, {-\widthFig*tan(\ThetaFig-180)});
	\coordinate (dbtLambda) at ({\RMEFig*cos(\ThetaFig)}, {\RMEFig*sin(\ThetaFig)});
	\fill [pattern={Lines[angle=45,line width=.2pt]}, pattern color=white!70!gray] (dbtLambda) arc (\ThetaFig:0:\RMEFig) -- (\widthFig,0) -- (\widthFig,\heightFig) -- (-\widthFig,\heightFig) -- (finLambda) -- cycle;
	\draw [thick] (finLambda) -- (dbtLambda) arc (\ThetaFig:0:\RMEFig) -- (\widthFig,0);
	\draw [dashed] (dbtLambda) -- (0,0);
	\draw [>=stealth,->] (-\widthFig-0.4,0) -- (\widthFig+0.4,0) node [right] {$x$};
	\draw [>=stealth,->] (0,{-\widthFig*tan(\ThetaFig-180)}) -- (0,{\heightFig+0.3}) node [right] {$y$};
	\draw (-\widthFig*0.5,\heightFig*0.6) node {$\OME$};
	\draw [->,>=latex] (0.35,0) arc (0:\ThetaFig-5:0.35) node[near start,above] {$\Theta$} -- (\ThetaFig:0.36);
	\draw [<-,>=latex] (-\widthFig*0.6, {-\widthFig*0.6*tan(\ThetaFig-180)}) to[bend right=20] ++(0.4,-0.4) node[right]{$\!\bordgaucheME$};
	\draw [<-,>=latex] (\widthFig*0.65,0) to[bend right=20] ++(0.15,-0.55) node[below]{$\ \ \GME$};
	\draw (\RMEFig,0) node{$\bullet$};
	\draw [<-,>=latex] (\RMEFig+0.05,0.05) to[bend left=10] ++(0.45,0.45) node[right]{\raisebox{1.5ex}{$\!\PME$}};
	\draw [<-,>=latex] ({\RMEFig*cos(55)}, {\RMEFig*sin(55)}) to[bend left=20] ++(0.4,0.45) node[right]{\raisebox{1ex}{$\SME$}};
	\draw [<->,>=latex] (0,-0.2) -- node[midway, below] {$\RME$} (\RMEFig, -0.2);
\end{tikzpicture}\par
\caption{The domains of the effective models}
\label{fig: domaine ME}
\end{figure}
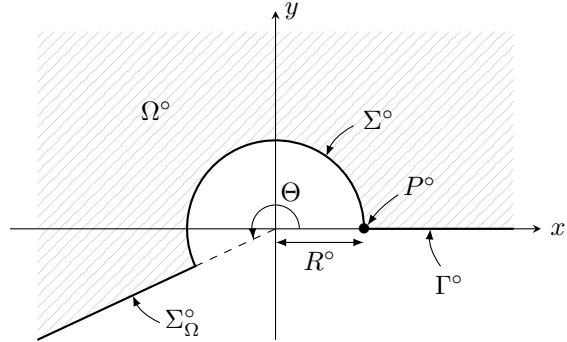

For~any $m\in \{3,5\}$, let us denote $\sme$ the solution of the effective model of order $m$. This function is defined on $\OME$. It satisfies the Helmholtz equation
\[\mu _{0} \Delta \sme + \omega ^2 \rho _{0} \sme =\fsource \quad \text{ in }\OME\]
with a Dirichlet condition $\sme =0$ on $\bordgaucheME$. Moreover the layer is replaced by the following effective boundary conditions:
\begin{equation}{\label{eq: condition effective sur GME}}
	\mu _{0} \partial _{n} \sme3 = -\frac{\mu _{1} }{\varepsilon } \sme3 , \quad \ \text{resp.\ }  \,
	\mu _{0} \partial _{n} \sme5 = -\frac{\mu _{1} }{\varepsilon } \sme5 + \frac{\varepsilon }{3} (\mu _{1} \partial _{x} ^2 +\omega ^2 \rho _{1} )\sme5, \quad \text{ on }\GME.
\end{equation}
These conditions are classical in presence of an infinite layer, see \cite{BenLem96}.\\

The difficulty here lies in the vicinity of the corner, where the exact solution $u_{\varepsilon }$ has an intricate behavior. We propose a Dirichlet-to-Neumann condition (DtN) on the artificial boundary $\SME$:
\begin{equation}{\label{eq: condition DtN ME}}
	\partial _{r} \sme = \opS(\sme_{|\SME} ) \quad \text{ on }\SME
\end{equation}
where $\smash\opS \in \mathcal{L}(L^2 (\SME))$ is a non-local operator of finite rank $N$, which is an additional parameter of the method. The construction of $\opS$ is described in \cref{sec: construction du DtN des ME}. Using the modal decomposition of the far fields in \cref{decompo modale des CL} and the approximation of $u_{\varepsilon }$ by sums of several far fields up to order $m$ (see \cref{DA de SE avec CL et delta variable}), we build in \cref{sec: decompo modale de SE} an approximate modal decomposition of $u_{\varepsilon }$ at order $\varepsilon ^{m}$. Then, $\opS$ is built s.t.\ it sends each mode of $u_{\varepsilon }$ to its derivative w.r.t.\ $r$ on $\SME$. More precisely, this is ensured only for the first $N$ modes because,  for technical reasons, we need to restrict ourselves to a finite number of modes to prove the well-posedness of the effective models. The modes of $u_{\varepsilon }$ depend on the coefficients $\cfinal_{d,d',p,\ell }$, which depend on the algebraic matching coefficients $\cmacro_{d,d',p,\ell }$ and $\cmicro_{d,d',p,\ell }$ and on the corner profile coefficients $\sigma _{q} (\ccDA_{d,k} )$, so the modes contain the information of the matching conditions. Intuitively, Condition \cref{eq: condition DtN ME} ensures that $\smash\sme \vrule height 2ex width 0pt$ has a behavior around the corner that is similar to $u_{\varepsilon }$.\\

Note that, contrary to the classical definition of DtN maps, the construction of $\opS(v)$ for $v:\SME\rightarrow \mathbb{C}$ does~not involve building a solution $v_{\varepsilon }$ of Helmholtz's equation behind the DtN (i.e.\ in $\Omega _{\varepsilon } \cap B(0,\RME)$) that has trace $v$ on $\SME$. Indeed having eliminated the layer implies that the boundaries $\Omega _{\varepsilon } \cap \partial B(0,\RME)$ and $\SME=\Omega  \cap \partial B(0,\RME)$ do~not match. So defining $v_{\varepsilon }$ would require extending $v$ to $(\varepsilon \layer) \cap \partial B(0,\RME)$ in the layer, but there~is no clear way on how to do so. Instead, $\opS$ is simply built s.t.\ $\opS(u_{\varepsilon |\SME} )$ is an approximation of $\partial _{r} u_{\varepsilon |\SME}$ and $\opS$ has good properties that allow to prove the well-posedness and convergence of the effective models.\\

Moreover, when $m=5$, \cref{eq: condition effective sur GME} involves the tangential derivative $\partial _{x} ^2 $ on $\GME$, that is treated by moving one $\partial _{x}$ to the test function. This 1D integration by parts requires a boundary condition at the point $\PME$, which is given by an additional condition of the form
\begin{equation}{\label{eq: condition DtN en PME}}
	\partial _{x} \sme5 (\PME)=\opP5(\sme5_{|\SME} )
\end{equation}
with $\opP5\in \mathcal{L}(L^2 (\SME),\mathbb{C})$. The construction of $\opP5$ is similar to $\opS$ and is also explained in \cref{sec: construction du DtN des ME}.\\

The tangential derivatives on $\GME$ in the model of order~5 also require to look for solutions $\sme5$ that are $H^{1}$ on $\GME$, in order to obtain a well-defined and well-posed variational formulation. Therefore we define the following spaces
\[\VME3 := \{v\in H^{1} (\OME) \mid v_{|\bordgaucheME} = 0\} \quad  \text{ and } \quad  \VME5 :=  \{v\in \VME3 \mid v_{|\GME} \in H^{1} (\GME)\}.\]
$\VME3$ is endowed with the norm $\|v\|_{H^{1} (\OME)}$, while $\VME5$ is equipped with the norm $(\|v\|_{H^{1} (\OME)} ^2 +\varepsilon  \mkern1mu \|\partial _{x} v\|_{L^2 (\GME)} ^2 )^{1/2}$ where the weight $\varepsilon $ is linked to the variational formulations given below.\\

In summary, the approximate solutions $\sme3$ and $\sme5$ solve the following problems resp.
\begin{equation}{\label{eq: ME3}}
	\text{find  }  u\in \VME3 \quad  \text{s.t.}  \qquad  \left\{ \begin{array}{r@{\,=\,}ll}
		\mu _{0} \Delta  u+\omega ^2 \rho _{0} u & \fsource &\text{ in }\OME\\[0.5mm]
		\mu _{0} \partial _{n} u + \frac{\mu _{1} }{\varepsilon } u & 0 &\text{ on }\GME\\[0.5mm]
		u & 0 &\text{ on }\bordgaucheME\\[0.5mm]
		\partial _{r} u - \opS3 (u_{|\SME} ) & 0 &\text{ on }\SME 
	\end{array}\right. \qquad 
\end{equation}
\begin{equation}{\label{eq: ME5}}
	\text{find  }  u\in \VME5 \quad  \text{s.t.}  \qquad  \left\{ \begin{array}{r@{\,=\,}ll}
		\mu _{0} \Delta  u+\omega ^2 \rho _{0} u & \fsource &\text{ in }\OME\\[0.5mm]
		\mu _{0} \,\partial _{n} u + \frac{\mu _{1} }{\varepsilon } u - \frac{\varepsilon }{3} (\mu _{1} \partial _{x} ^2 +\omega ^2 \rho _{1} )u & 0 &\text{ on }\GME\\
		u & 0 &\text{ on }\bordgaucheME\\[0.5mm]
		\partial _{r} u - \opS5 (u_{|\SME} ) & 0 &\text{ on }\SME \\[0.5mm]
		\partial _{x} u - \opP5 (u_{|\SME} ) & 0 & \text{ at }  \PME
	\end{array}\right. \qquad 
\end{equation}
It is easy to see that the variational formulations associated to \cref{eq: ME3,eq: ME5} are resp.:
\begin{equation}{\label{eq: FV ME3}}
	\forall v\in H_{3} , \qquad  \int _{\OME} (\mu _{0} \nabla  u \cdot  \nabla  \bar{v}- \omega ^2  \rho _{0} u \bar{v}) + \int _{\GME} \frac{\mu _{1} }{\varepsilon } u \bar{v}+ \int _{\SME} \mu _{0} \,\opS3 (u_{|\SME} )\, \bar{v}= 
	\langle \fsource, \bar{v}\rangle _{H^{-1} ,H^{1} },
\end{equation}
\begin{multline}\label{eq: FV ME5}
	\forall v\in H_{5} , \qquad  \int _{\OME} (\mu _{0} \nabla  u \cdot  \nabla  \bar{v}- \omega ^2 \rho _{0} u \bar{v}) + \int _{\GME} \left( \frac{\mu _{1} }{\varepsilon } u \bar{v} + \frac{\varepsilon }{3} (\mu _{1} \partial _{x} u \,\partial _{x} \bar{v}- \omega ^2 \rho _{1} u \bar{v}) \right) \\
	+ \int _{\SME} \mu _{0} \,\opS5 (u_{|\SME} )\, \bar{v}\,+\, \frac{\varepsilon  \mu _{1} }{3} \,\opP5 (u_{|\SME} ) \cdot  \overline{v(\PME)}= 
	\langle \fsource, \bar{v}\rangle _{H^{-1} ,H^{1} },
\end{multline}
The parameters of these models are:
\begin{itemize}
	\item $m\in \{3,5\}$ the order of the method,
	\item $\RME$ the radius of the DtN boundary,
	\item and $N$ the number of modes taken into account in the DtN.
\end{itemize}
Note that $\varepsilon $ appears in the effective model, but this does~not require a fine mesh.
\end{textesansboite}
\nomenclature[22]{$m$}{In \cref{chap: ME}, order of the effective model}
\nomenclature[23]{$N$}{In \cref{chap: ME}, number of modes taken into account in the DtN maps of the effective models}

\begin{remarque}
The previous approaches in the literature \cite{AuvVia18,AuvVia19,Auv18} tried to set the effective models directly in $\Omega $ (so without a transparent condition around the corner). Effective boundary conditions similar to \cref{eq: condition effective sur GME} were imposed on $\Gamma $ entirely, with some modification (in general at the corner or its vicinity) to try to approach more accurately the exact solution, see \cref{sec: etat de l'art}. Thus these previous attempts did~not directly involve the far fields' corner singularities and only include pointwise differential boundary conditions.\\

It seems to us difficult to obtain significantly better results with this kind of approach for a general angle $\Theta $ for two reasons. First, it seems reasonable that the effective solution has a behavior at the corner that is similar to the exact solution. As~we have seen, this behavior is intricate and involves corner singularities -- the $\mcll_{d,n}$ -- with complicated expressions (especially when $\pst$ is rational), and the coefficients of these singularities -- the $\sigma _{d} (\bfuDA_{p,\ell } )$ -- are linked to each other through involved matching conditions. Thus it seems difficult to capture this complex behavior with only a differential boundary condition on $\Gamma $. Moreover, this behavior is a local information in the vicinity of the corner, while a differential condition with constant coefficients is rather pointwise (local perturbations of the boundary condition on $\Gamma $ like in \cite{AuvVia19} avoid this problem, but they are not satisfying as they require finely meshing near the corner).
Secondly, an effective model of order $m$ should take into account the corner profile coefficients that are involved in the asymptotic expansion of $u_{\varepsilon }$ up to order $m$. When $m=5$ and $\Theta \rightarrow 2\pi ^{-}$, that amounts to 46 coefficients. There is not enough space in effective conditions with a few parameters to include all this complexity.\\

For these reasons, we chose to include the singularities $\smash{\mcll_{d,n} }$ in the model, which is done via the DtN operators $\opS$ and $\opP5$. Note that, unlike the computation of the corner profiles in \cref{sec: calcul numerique des CC}, we could not include these singularities in the variational space as singular complements (and define the effective models in $\Omega $), as some singularities diverge at 0 and are~not $H^{1}$ at the corner.
\end{remarque}

\subsection{Statement of the main results}

\refstepcounter{notation}
\begin{paragraphesansboite}{Notations \thenotation:}
\label{notations ME}%
\begin{itemize}
	\item We denote $\Rsource := \operatorname{dist}(0,\operatorname{supp}\fsource)$. Since we assumed $\operatorname{dist}(\Gamma ,\operatorname{supp}\fsource)>0$, we have $\Rsource>0$. We also assumed $\closedBall(0,\RME)\cap \operatorname{supp}(\fsource)= \varnothing $, which is equivalent to $\RME<\Rsource$.
	\item For~any $p\in \mathbb{N}$, we define 
	\begin{align*}
		\ell _{p} &:= \max _{p'\in \NN\cap [0,p]} \max \{\ell \in \mathbb{N} \mid u_{p',\ell } \text{ is generically non-zero} \}\\
		&\,= \, \min \{\ell \in \mathbb{N} \mid \forall p'\in \NN\cap [0,p],\; \forall \ell '>\ell ,\ u_{p',\ell '} \text{ is always zero} \}
	\end{align*}
	where ``always zero'' means that the far field vanishes for any value of $(\mu _{0} ,\mu _{1} ,\rho _{0} ,\rho _{1} ,\omega ,\fsource)$ and the corner profile coefficients $\sigma _{q} (\ccDA_{d,n} )$. The equality between the two expressions of $\ell _{p}$ is a consequence of \cref{generiqmt non nul ou tjs nul}. By~construction, $\ell _{p}$ depends only on $\Theta $ and $p$, and is independent of $(\mu _{0} ,\mu _{1} ,\rho _{0} ,\rho _{1} ,\omega ,\fsource)$ and the corner profiles. The existence of $\ell _{p}$ is guaranteed by \cref{borne sur l}.
	\item All over this section, $C_{\dots }$ will denote any positive constant that depends \emph{only} on the variables in ``...'' and $(\Theta ,\mu _{0} ,\mu _{1} ,\rho _{0} ,\rho _{1} ,\omega )$ and on $m$. Its value may change from one line to another.
	\item For~any $x\in \mathbb{R}$, we denote $\domi{x}:=1+|x|$.
	\item For~any domain $D$ and $u\in H^{1} (D)$, $\snh{u}(D) :=\|\nabla  u\|_{L^2 (D)}$ is the $H^{1}$ semi-norm.
\end{itemize}
\end{paragraphesansboite}
\nomenclature[60-11b]{$\langle {\dots }\rangle $}{In \cref{chap: ME}, $\langle x\rangle :=1+|x|$}
\nomenclature[17-60]{$\snh{u}(D)$}{The $H^{1}$ semi-norm: $\snh{u}(D) :=\|\nabla  u\|_{L^2 (D)}$ for~any domain $D$}

\begin{textesansboite}
The following two theorems are the main results of this chapter. They are proven resp.\ in \cref{sec: ME bien poses} and \ref{sec: cvg ME}.
\end{textesansboite}

\begin{theoreme}{well-posedness and uniform stability of the effective models}
\label{caractere bien pose des ME}%
Let $\RME \in (0,\Rsource)$, $N\in \mathbb{N}$ and $\varepsilon >0$. There exist two positive constants $\CstabME$ and $\CcontME$ s.t.\ we have
\def\monstrut{\vrule height0pt depth1.5pt width0pt}
\begin{itemize}
	\item \underline{for $m=3$:\monstrut} if $\displaystyle \varepsilon \leqslant \CstabME\, \frac{\RME}{\domi{\ln \RME}^{3}}$, then \cref{eq: FV ME3} has a unique solution $\sme3 \in \VME3$ and $\sme3$ satisfies\vspace{-.5ex}
	\[\|\sme3\|_{H^{1} (\OME)} \leqslant  \CcontME\, \|\fsource\|_{H^{-1} (\Omega )} .\vspace{-.5ex}\]
	\item \underline{for $m=5$:\monstrut} if $\displaystyle \varepsilon \leqslant \CstabME\, \frac{\RME^{4/3} }{\domi{\ln \RME}^2 }$, then \cref{eq: FV ME5} has a unique solution $\sme5 \in \VME5$ and $\sme5$ satisfies
	\[\|\sme5\|_{H^{1} (\OME)} + \sqrt{\varepsilon } \,\snhg{\sme5} \leqslant  \CcontME\, \|\fsource\|_{H^{-1} (\Omega )} .\]
\end{itemize}\fbi
\end{theoreme}

\begin{definitionsansboite}{}
\label{def: condition de stabilite}%
We define the \emph{stability condition} as ``$\varepsilon \leqslant \CstabME\, \RME \,\domi{\ln \RME}^{-3}$'' if $m=3$ and ``$\varepsilon \leqslant \CstabME\, \RME^{4/3} \,\domi{\ln \RME}^{-2}$'' if $m=5$.
\end{definitionsansboite}

\begin{remarque}
\label{rq: Cstab et N3}%
\begin{itemize}
	\item $\CstabME$ and $\CcontME$ follow the convention of \cref{notations ME}, so $\CstabME$ is independent of $(\RME,\varepsilon )$, while $\CcontME$ depends only on $(\Theta ,\mu _{0} ,\mu _{1} ,\rho _{0} ,\rho _{1} ,\omega ,m)$.
	\item The dependence of $\CstabME$ w.r.t.\ $N$ can be made explicit at the cost of more technical proofs. One can actually take
	\[\CstabME = \frac{C_{\Rsource} }{N^{3}} .\]
	Appendix~\ref{sec: annexe comportement en N des ME} explains how to modify the intermediate results and proofs to obtain this result. There is no~guarantee that $N^{3}$ is optimal.
	\item To prove \cref{caractere bien pose des ME}, we first show that the part of the variational formulation that does~not involve the DtN is uniformly coercive w.r.t.\ $(\varepsilon ,\RME,N)$. Then we prove that the DtN part is close at $O(\varepsilon )$ to a positive operator (namely the DtN operator for a Dirichlet corner without layer). Thus, when $\varepsilon $ is small enough, the overall variational formulation is uniformly coercive. This explains the presence of a stability condition.
	\item Using the analytic Fredholm alternative, one can prove that, when the stability condition is~not satisfied, the effective models are still well-posed except for a discrete set of values of $\varepsilon $ (that depends on $\RME$, $N$, etc.). So in~practice, even without the stability condition, it~is likely that there~is no numerical problem to compute $\sme$ (except for a discrete set of $\varepsilon $).
\end{itemize}
\end{remarque}

\begin{theoreme}{error estimate for the effective models}
\label{cvg des ME}%
Let $R>0$ be s.t.\ $\closedBall(0,R) \cap  \operatorname{supp}(\fsource)=\varnothing $, i.e.\ $R<\Rsource$. We assume that $\RME<R$, that the stability condition is satisfied and that $\varepsilon \leqslant \RME^2 $.
Then there exists $C_{\fsource,R} >0$ s.t.\ we have 
\begin{equation}{\label{eq: cvg des ME}}
	\|\sme-u_{\varepsilon } \|_{H^{1} (\OME)} \leqslant  C_{\fsource,R} \left[ \bigg( \frac{\varepsilon }{\RME} \bigg)^{m} \domi{\ln \varepsilon }^{\ell _{m}} N^{3} + \left( \frac{\RME}{R} \right)^{N\pst} \right].
\end{equation}\fbeq
\end{theoreme}

\begin{remarquenumerotee}
\label{rq: commentaire estimee de cvg des ME}%
The first term of \cref{eq: cvg des ME} is the $O(\varepsilon ^{m} )$ that was announced. It is weighted by $1/\RME^{m}$ because of the singular behavior of the exact solution and the far fields near the corner. It is also multiplied by $N^{3}$ for technical reasons. This factor can be improved into $N^{2.5+\eta }$ for~any $\eta >0$ (replacing $C_{\fsource,R}$ by $C_{\fsource,R,\eta }$ in \cref{eq: cvg des ME}), using more subtle arguments that are given in Appendix~\ref{sec: annexe comportement en N des ME}. The second term of \cref{eq: cvg des ME} is due to the truncation of the DtN operators $\opS$ and $\opP5$. It tends to 0 geometrically w.r.t.\ $N$, so it is easy to make it very small, which makes the $O(\varepsilon ^{m} )$ the dominant error. The constant $C_{\fsource,R}$ is likely to explode when $R\rightarrow \Rsource^{-}$. Moreover, \cref{eq: cvg des ME} is sharp w.r.t.\ $\varepsilon $, but probably not w.r.t.\ $\RME$ and $N$.
\end{remarquenumerotee}

\begin{remarque}
The condition ``$\varepsilon \leqslant \RME^2 $'' is necessary in the proof of \cref{cvg des ME} to use \cref{DA de SE avec CL et delta variable}, which quantifies the error $u_{\varepsilon } -\sum _{p<m,\ell } \varepsilon ^{p} \ln ^{\ell } \!\varepsilon \,u_{p,\ell }$ made by the asymptotic expansion.
\Cref{rq: condition sur delta dans DA de SE avec CL} explains that, by modifying the choice of the matching zone, this condition can be relaxed to ``$\varepsilon \leqslant M\.\RME^{1+\delta }$'' for~any fixed $M,\delta >0$ (if the constant $C_{\fsource,R}$ in \cref{eq: cvg des ME} is replaced by $C_{\fsource,R,M,\delta }$). Roughly speaking, this means that we must have
\begin{equation}{\label{eq: eps<<R dans ME}}
	\fboxsep=4pt
	\fbox{$ \mkern1mu \varepsilon  \ll \RME $}.
\end{equation}
The stability condition of \cref{def: condition de stabilite} can be seen as another rigorous elaboration of \cref{eq: eps<<R dans ME}. Moreover, the factor $(\varepsilon /\RME)^{m}$ in \cref{eq: cvg des ME} also encourages to take $\varepsilon \ll \RME$.
\end{remarque}

\begin{remarquenumerotee}
\label{rq: choix parametres ME}%
$\RME$ cannot be taken too large (because it deteriorates the second term of \cref{eq: cvg des ME}), nor too small (because it deteriorates the first term of \cref{eq: cvg des ME} and would require a fine mesh near the corner). The same holds for $N$: if it~is too large, then the stability condition is more restrictive and the first term of \cref{eq: cvg des ME} is deteriorated, and if it~is too small then the second term of \cref{eq: cvg des ME} increases.\\

In~practice, one can choose the model's parameters in the following way. First, we take $R$ rather close to $\Rsource$, e.g.\ $R :=0.9\,\Rsource$. Then we take $\RME := R/2$, so that $\frac{\RME}{R} = \frac{1}{2}$. Next, $\varepsilon $ should satisfy \cref{eq: eps<<R dans ME}, which is equivalent here to $\varepsilon \ll \Rsource$. This condition is not very restrictive as we would clearly not be in an asymptotic regime if we had $\varepsilon  \sim \Rsource$ or $\varepsilon \geqslant \Rsource$. Then, we take
\[N := \lceil {-m \frac{\Theta }{\pi } \log _{2} \varepsilon  }\rceil .\]
Thus we have $( \RME/R )^{N\smallpst} \leqslant  \varepsilon ^{m}$, so \cref{eq: cvg des ME} reduces to:
\[\|\sme-u_{\varepsilon } \|_{H^{1} (\OME)} \leqslant  C_{\fsource} \,\varepsilon ^{m} \domi{\ln \varepsilon }^{\ell _{m} +3} .\]
As for the stability condition of \cref{def: condition de stabilite}, by \cref{rq: Cstab et N3} it becomes here $\varepsilon  \ln ^{6} \varepsilon  \leqslant  C_{\Rsource} \RME \domi{\ln \RME}^{-3}$ or $\varepsilon  \ln ^{6} \varepsilon  \leqslant  C_{\Rsource} \RME^{4/3} \domi{\ln \RME}^{-2}$. This condition is satisfied because we assumed $\varepsilon \ll \RME$.
\end{remarquenumerotee}

\begin{remarquenumerotee}
\label{rq: degre en ln eps dans ME}%
Let $E_{3} :=\{ \frac{2\pi }{3} , \pi  , \frac{4\pi }{3} , \frac{5\pi }{3} \}$ and $E_{5} =E_{3} \cup \{ \frac{a\pi }{5} \mid a\in [\![2,9]\!]\}$.
\begin{itemize}
	\item If $\Theta  \not\in E_{m}$ (i.e.\ for most angles), then we have $\ell _{m} =0$, see \cref{simplifications cas irrationnel,simplifications cas rationnel version forte}. Thus the factor $\domi{\ln \varepsilon }^{\ell _{m}}$ disappears from \cref{cvg des ME}.
	\item Otherwise, we have $\ell _{m} =1$ if $\Theta \neq \pi $, and $\ell _{m} =m-1$ if $\Theta =\pi $. These values were computed in the same way as the data of \cref{sec: exemples de structures de DA}.
\end{itemize}
\end{remarquenumerotee}

\begin{remarque}
By~taking $\omega =0$, it~is easy to adapt our work to obtain effective models of order~3 and 5 for the Poisson equation, which is also novel. This requires replacing the variational spaces used for $u_{\varepsilon }$, $\sme3$ and $\sme5$ by resp.\ $V(\Omega _{\varepsilon } )$, $V(\OME)$ and $V(\OME) \cap  V(\GME)$, where we denote $V(D) := \{v\in H^{1} _{\mathrm{loc}} (D) \mid \nabla  v \in L^2 (D)\}$ for any domain $D$. In addition, the particular expressions of the functions $\mcll_{d,n}$ and of the modes of the DtN that involve the Bessel functions are not valid anymore. They are replaced by simpler ones, involving only individual powers of $r$, so the case of Poisson is simpler to study than Helmholtz.
\end{remarque}

\section{Approximate modal decomposition of the exact solution}
\label{sec: decompo modale de SE}

\subsection{Preliminary decomposition}
\label{sec: decompo modale de SE 1}

\begin{textesansboite}
The construction of DtN operators $\opS$ and $\opP5$ is based on a modal decomposition of
\begin{equation}{\label{eq: def notation somme CL}}
	\scl := \sum _{p\in \NN,\, p<m} \; \sum _{\ell \in \mathbb{N}} \, \varepsilon ^{p} \ln ^{\ell } \!\varepsilon \, u_{p,\ell } ,
\end{equation}
which is an approximation of $u_{\varepsilon }$ at order $\varepsilon ^{m}$ according to \cref{DA de SE avec CL et delta variable} (in the notation, ``f{}f'' stands for ``far fields''). This decomposition is introduced in this section, using the already known modal decomposition of the far fields $u_{p,\ell }$ stated in \cref{decompo modale des CL}. We recall that this previous result states that we have
\begin{equation}{\label{eq: decompo modale CL repetee pour ME}}
	u_{p,\ell } \,=\, \sum _{n=0} ^{\lfloor {p}\rfloor } \sum _{\substack{d\in \TZ*\\ d\geqslant -p+n}} \sigma _{d} (\uDA_{p-n,\ell } ) \,\mcll_{d,n} 
\end{equation}
in $\Omega \cap B(0,\Rsource)$, where the functions $\mcll_{d,n}$ are built in \cref{convergence des mcl}. Moreover, when truncating the above series at $d=N\pst$, one makes an error that is smaller than $C_{p,\fsource,R} \left( \frac{r}{R} \right)^{(N+1)\pst-\lfloor {p}\rfloor } \domi{\ln r}^{2p+1}$ for~any $R<\Rsource$.\\

In most cases (specified in \cref{rq: def mcle} below), the functions $\mcll_{d,n}$ are easy to compute because they are finite linear combinations of Bessel functions (times some sines of $\theta $). Otherwise, the series $\mcla_{d,n}$ is more complex as it involves $\ln r$, and there~is no known expression of $\mcll_{d,n}$. Hence, in that case, $\mcll_{d,n}$ can only be approximated by a series truncation $T_{<q} (\mcla_{d,n} )$ for a given $q$. The following definition introduces the case-dependent functions that are used numerically in place of the modes $\mcll_{d,n}$. We restrict ourselves to $n<m$ as we will see that the effective models do~not involve the functions $\mcll_{d,n}$ s.t.\ $n\geqslant m$.
\end{textesansboite}

\begin{definition}{}
\label{def: mcle}%
Let $E_{3} := \{\pi \}$ and $E_{5} := \left\{ \frac{\pi }{3} , \frac{2\pi }{3} , \pi  , \frac{4\pi }{3} , \frac{5\pi }{3} \right\}$. For~any $N\in \mathbb{N}$, $d\in \TZ*$ s.t.\ $d\leqslant N\pst$ and $n\in [\![0,m-1]\!]$, we define the function $\mcle_{d,n} : \overline{\Omega } \setminus \{0\} \rightarrow \mathbb{C}$ as follows.
\begin{itemize}
	\item If $\Theta  \not\in E_{m}$, then $\mcle_{d,n} :=\mcll_{d,n}$.
	\item If $\Theta  \in E_{m}$, then $\mcle_{d,n} := \left\{\begin{array}{@{\;}l@{\;\;}l}
		\mcll_{d,n} & \text{if } d>0 \text{ and } n=0\\[1pt]
		T_{<(N+1)\pst} (\mcla_{d,n} ) & \text{otherwise}.
	\end{array}\right.$
\end{itemize}\fbi[0.8]
\end{definition}

\begin{remarquenumerotee}
\label{rq: def mcle}%
The angles $\Theta  \not\in E_{m}$ are exactly the ones for which all the functions $\mcll_{d,n}$ with $n\in [\![0,m-1]\!]$ have explicit expressions in terms of Bessel functions and thus can be directly computed. Indeed,
\begin{itemize}
	\item if $\frac{\Theta }{\pi }$ is irrational, or if $\frac{\Theta }{\pi }$ is rational and its denominator $b$ (as an irreducible fraction) satisfies $m\leqslant b$, then \cref{expression sing heterog CL,ex: expression sing heterog CL} provide a formula for $\mcll_{d,n}$ using Bessel functions,
	\item and \cref{expression sing heterog CL b=2 ou 4} extends this result to the case where $\frac{\Theta }{\pi }$ is rational $b\in \{2,4\}$ and $n<5$\,\footnote{In this case, the Bessel functions of the second kind may be needed}.
\end{itemize}
The remaining case is ($b\in \{1,3\}$ and $b<m$), i.e.\ $\Theta  \in E_{m}$. Using the tools of \cref{sec: implementation des espaces A}, we verified numerically that in this case not all functions $\mcll_{d,n}$ are linear combinations of Bessel functions. However, even here, the functions $\mcll_{d,0}$ with $d>0$ can be expressed using a Bessel function according to \cref{eq: preuve expression sing homog CL}. This allows us to take $\mcle_{d,0} :=\mcll_{d,0}$ in \cref{def: mcle}, which will simplify a bit the stability analysis.\\

Nevertheless, note that it~is possible to take $\mcle_{d,n} := T_{<(N+1)\pst} (\mcla_{d,n} )$ for any angle $\Theta $ without changing the stability and error analysis. This could be interesting to have a unique implementation of the effective models for any angle $\Theta $, whereas taking $\mcle_{d,n}$ as in \cref{def: mcle} simplifies the implementation for the angles $\Theta  \not\in E_{m}$.\\

As for the truncation order of $\mcla_{d,n}$ used in \cref{def: mcle}, which is $(N+1)\pst$, it~is chosen so that the error made by the truncation of \cref{eq: decompo modale CL repetee pour ME} w.r.t.\ $d$ stays the same when replacing $\mcll_{d,n}$ by $\mcle_{d,n}$. This result is justified by \cref{decompo en modes tronques des CL}. This truncation order corresponds to the order of the first neglected mode $\mcll_{(N+1)\pst,0}$, as this function behaves at 0 like $r^{(N+1)\pst}$. Truncating further the series $\mcla_{d,n}$ would not improve the final error estimate, and truncating before would deteriorate it.
\end{remarquenumerotee}

\begin{textesansboite}
We will start by proving a preliminary modal decomposition of $\scl$, which is~not yet the decomposition used to build $\opS$ and $\opP5$. Its modes are defined for any $N\in \mathbb{N}$, $\varepsilon \geqslant 0$ and $d\in \TZ*$ s.t.\ $d\leqslant N\pst$ as
\begin{equation}{\label{eq: def mse}}
	\mse_{d} := \sum _{n=0} ^{m-1} \varepsilon ^{n} \mcle_{d,n} .
\end{equation}
\nomenclature[56-23-40-13]{$\mse_{d}$}{Functions used to build the modes $\mme_{n}$ for the effective models}%
Note that $\mse_{d}$ depends on $m$, but we do~not  include it in the notation for simplicity. Each function $\mcle_{d,n}$ satisfies the Helmholtz equation in $\Omega $ (exactly if $\mcle_{d,n} = \mcll_{d,n}$ and approximately otherwise), so $\mse_{d}$ has the same property. Moreover, $\mcle_{d,\red0}$ is a mode of Helmholtz's equation in a Dirichlet corner \emph{without} layer, and the terms $\varepsilon ^{n} \mcle_{d,n}$ with $n\geqslant 1$ act like corrections to $\mcle_{d,0}$ that take into account the presence of the layer. Using that the functions $\mcll_{d,n}$ satisfy the same inductive Dirichlet condition on $\Gamma $ (see \cref{convergence des mcl,construction directe des upl}), one can prove that $\mse_{d}$ approximately satisfy the effective condition \cref{eq: condition effective sur GME} on $\Gamma $. However, the modes $\mse_{d}$ do~not take into account the matching conditions around the corner ; this part will be considered later using other modes.\\

Next, for~any $\varepsilon >0$, let us denote:
\begin{equation}{\label{eq: def sclDA}}
	\sclDA := \sum _{p\in \NN \cap [0,m)} \sum _{\ell \in \mathbb{N}} \varepsilon ^{p} \ln ^{\ell } \!\varepsilon \, \uDA_{p,\ell } \in \Amacro(\Omega ) .
\end{equation}
where $\uDA_{p,\ell } :=\bfuDA_{p,\ell |\Omega }$ describes the behavior at the corner of $u_{p,\ell }$ (see \cref{validation construction et DA des champs}). Thus $\sclDA$ describes the behavior at the corner of $\scl$. $\sclDA$ will actually only be used through its coefficients $\sigma _{d}$ which satisfy:
\begin{equation}{\label{eq: sigma de sclDA}}
	\forall d\in \txtTZ*, \qquad  \sigma _{d} (\sclDA) = \! \sum _{p\in \NN \cap [0,m)} \sum _{\ell \in \mathbb{N}} \varepsilon ^{p} \ln ^{\ell } \!\varepsilon \cdot  \sigma _{d} (\uDA_{p,\ell } ) \in \mathbb{C} .
\end{equation}
The following proposition is based on the results of \cref{chap: decompo modale CL}.
\end{textesansboite}

\begin{proposition}{decomposition of $\scl$ on the modes $\mse_{d}$}
\label{decompo modale SE 1}%
For~any $N\in \mathbb{N}$ and $\varepsilon >0$, let us denote
\[\rho _{N,\varepsilon } := \scl - \sum _{\substack{d\in \TZ*\\ -m\leqslant d\leqslant N\pst}} \sigma _{d} (\sclDA) \cdot \mse_{d} .\]
Let $R\in (0,\Rsource)$ and $M>0$. There exists a positive constant $C_{\fsource,R,M}$ s.t., for~any $\varepsilon >0$, $N\in \mathbb{N}$, $r\in (0,R)$ and $\theta \in [0,\Theta ]$ that satisfy $\varepsilon \leqslant Mr$, we have\vspace{-1ex}
\[|\rho _{N,\varepsilon } (r,\theta )|+r|\nabla  \rho _{N,\varepsilon } (r,\theta )| \leqslant  C_{\fsource,R,M} \left[ \left( \frac{r}{R} \right)^{N\pst+ \frac{1}{2} }  + \frac{\varepsilon ^{m} }{r^{m-\pst}} \,\domi{\ln \varepsilon }^{\ell _{m}} \right] \domi{\ln r}^{2m+1} .\]\fbeq
\end{proposition}

\begin{remarque}
Taking $\mcle_{d,n} :=\mcll_{d,n}$ and $N\rightarrow \infty $, this implies in particular
\begin{equation}{\label{eq: rq decompo modale SE 1}}
	\scl = \sum _{d\in \TZ*\cap [-m,\infty )} \sigma _{d} (\sclDA) \cdot \Phi ^{\varepsilon } _{d} \;+\, O(\varepsilon ^{m} \ln ^{\ell _{m}} \!\varepsilon ) \qquad  \text{with}  \qquad  \Phi ^{\varepsilon } _{d} := \sum _{n=0} ^{m-1} \varepsilon ^{n} \mcll_{d,n}
\end{equation}
where the big $O$ depends on $(\fsource,R,M,r)$. \Cref{eq: rq decompo modale SE 1} gives an approximate modal decomposition of $\scl$ on the $\Phi ^{\varepsilon } _{d}$ at order $\varepsilon ^{m} \ln ^{\ell _{m}} \!\varepsilon $. This is also an approximate decomposition of the exact solution, since $u_{\varepsilon } =\scl+O(\varepsilon ^{m} \ln ^{\ell _{m}} \!\varepsilon )$. \Cref{decompo modale SE 1} is a refinement of this result, that provides a precise error estimate.\\

The condition $\varepsilon \leqslant Mr$ is~not essential to obtain an error estimate, but it improves a bit the estimate. In~practice, this condition will be satisfied for $r:=\RME$ and $M:=\CstabME$ because of the stability condition (see \cref{def: condition de stabilite}).
\end{remarque}

\begin{demo}
\def\pdm#1{D^N_{#1,\ell,n}}
\def\lmax{\ell_{m}}%
\underline{Step 1:} First, let us express $\sum _{d} \sigma _{d} (\sclDA) \,\mse_{d} $ in terms of the modal decompositions of the far fields. By~\cref{eq: def mse,eq: sigma de sclDA}, we have
\begin{align*}
	\sum _{\substack{d\in \TZ*\\ -m\leqslant d\leqslant N\pst}} \sigma _{d} (\sclDA) \cdot \mse_{d} 
	&= \sum _{\substack{d\in \TZ*\\ -m\leqslant d\leqslant N\pst}} \left(\sum _{p\in \NN \cap [0,m)} \sum _{\ell =0} ^{\lmax} \varepsilon ^{p} \ln ^{\ell } \!\varepsilon \cdot  \sigma _{d} (\uDA_{p,\ell } ) \right) \left( \sum _{n=0} ^{m-1} \varepsilon ^{n} \mcle_{d,n} \right) \\
	&= \sum _{p\in \NN \cap [0,m)} \sum _{n=0} ^{m-1} \sum _{\ell =0} ^{\lmax} \varepsilon ^{p+n} \ln ^{\ell } \!\varepsilon \cdot  \sum _{\substack{d\in \TZ*\\ -m\leqslant d\leqslant N\pst}} \sigma _{d} (\uDA_{p,\ell } )\cdot  \mcle_{d,n} \\
	&= \sum _{p'\in \NN \cap [0,2m)} \sum _{\substack{n\in [\![0,m-1]\!]\\ p'-n\in \NN \cap [0,m)}} \sum _{\ell =0} ^{\lmax} \varepsilon ^{p'} \ln ^{\ell } \!\varepsilon \cdot  \! \sum _{\substack{d\in \TZ*\\ -m\leqslant d\leqslant N\pst}} \sigma _{d} (\uDA_{p'-n,\ell } )\cdot  \mcle_{d,n}
	\refstepcounter{equation}\tag{\theequation}\label{eq: decompo modale SE 1 preuve 1}
\end{align*}
using the change of index $p':=p+n$. Since we have $u_{p,\ell } =0$ for any $p\in \mathbb{R} \setminus \NN$ by convention, the condition $p'-n \in \NN$ in \cref{eq: decompo modale SE 1 preuve 1} can be deleted. Similarly, \Cref{precision de l'ansatz sur uDA} implies that $\sigma _{d} (u_{p,\ell } )=0$ when $d<-p+\pst$. So the condition $-m\leqslant d$ in \cref{eq: decompo modale SE 1 preuve 1} can be replaced by $-p'+n+\pst\leqslant d$ (as $p'-n=p\leqslant m$). In addition, the condition $p'-n\in [0,m)$ is equivalent to $p'-m<n\leqslant p'$. To simplify the calculations, let us denote
\[\pdm{p} := \sum _{\substack{d\in \TZ*\\ \red{-p+n+\pst}\leqslant d\leqslant N\pst}} \sigma _{d} (\uDA_{p-n,\ell } )\cdot  \mcle_{d,n} .\]
Thus, we obtain:
\begin{align}
	\smash{\sum _{\substack{d\in \TZ*\\ -p+m\leqslant d\leqslant N\pst}}} \sigma _{d} (\sclDA) \cdot \mse_{d} 
	&= \sum _{p'\in \NN \cap [0,2m)} \, \sum _{n\in [\![0,m-1]\!]\cap (p'-m,p']}\, \sum _{\ell =0} ^{\lmax} \varepsilon ^{p'} \ln ^{\ell } \!\varepsilon \cdot  \pdm{p'} \nonumber \\
	&= \sum _{p\in \NN \cap [0,m)} \sum _{\ell =0} ^{\lmax} \varepsilon ^{p} \ln ^{\ell } \!\varepsilon  \sum _{n=0} ^{\lfloor {p}\rfloor } \pdm{p}
	+ \sum _{p\in \NN \cap [m,2m)} \sum _{\ell =0} ^{\lmax} \varepsilon ^{p} \ln ^{\ell } \!\varepsilon  \sum _{n=\lfloor {p-m}\rfloor +1} ^{m} \pdm{p} . \label{eq: decompo modale SE 1 preuve 2}
\end{align}
We notice that the first term of \cref{eq: decompo modale SE 1 preuve 2} is the sum of the total modal decompositions of the far fields of order $p\leqslant m$, while the second term is a kind of boundary term that contains only parts of the modal decompositions of the next far fields.\medskip\smallskip

\underline{Step 2:} For~any $f\in H^{1} _{\mathrm{loc}} (\Omega )$ and $(r,\theta )$, let us denote $|f|_{r,\theta } := |f(r,\theta )|+r|\nabla  f(r,\theta )|$. We claim that we have for~any $p\in \NN$ and $\ell \in \mathbb{N}$:
\begin{equation}{\label{eq: decompo modale SE 1 preuve decompo modale des CL generalisee}}
	\bigg|u_{p,\ell } - \sum _{n=0} ^{\lfloor {p}\rfloor } \pdm{p} \bigg|_{r,\theta } \leqslant C_{p,\fsource,R} \left( \frac{r}{R} \right)^{(N+1)\pst -\lfloor {p}\rfloor } \domi{\ln r}^{2p+1} .
\end{equation}
Indeed, \cref{decompo modale des CL} proves this result in the case where $\mcle_{d,n} =\mcll_{d,n}$ for~any $(d,n)$, while \cref{decompo en modes tronques des CL} treats the case where $\mcle_{d,n} =T_{<(N+1)\pst} (\mcla_{d,n} )$ for~any $(d,n)$. Here we are in a mixed situation where, for~any $n$, we have either $\mcle_{d,n} =\mcll_{d,n}$ for~any $d$, or $\mcle_{d,n} =T_{<(N+1)\pst} (\mcla_{d,n} )$ for~any $d$. Since the proofs of \cref{decompo modale des CL,decompo en modes tronques des CL} treat separately the different values of $n$, it~is easy to generalize these theorems to our situation, which proves \cref{eq: decompo modale SE 1 preuve decompo modale des CL generalisee}. Multiplying \cref{eq: decompo modale SE 1 preuve decompo modale des CL generalisee} by $\varepsilon ^{p} \ln ^{\ell } \!\varepsilon $ and summing over $(p,\ell )$, we get:
\[\bigg|\scl - \!\! \sum _{p\in \NN \cap [0,m)} \sum _{\ell =0} ^{\lmax} \varepsilon ^{p} \ln ^{\ell } \!\varepsilon  \sum _{n=0} ^{\lfloor {p}\rfloor } \pdm{p} \bigg|_{r,\theta } \leqslant  C_{\fsource,R} \sum _{p\in \NN \cap [0,m)} \sum _{\ell =0} ^{\lmax} \varepsilon ^{p} \ln ^{\ell } \!\varepsilon  \cdot  \left( \frac{r}{R} \right)^{(N+1)\pst -p} \domi{\ln r}^{2m+1} .\]
Let $\delta :=\pst- \frac{1}{2} >0$. For~any $p>0$, we have $\varepsilon ^{p} \ln ^{\ell } \!\varepsilon \cdot  r^{\pst-p} = \varepsilon ^{\delta } \ln ^{\ell } \!\varepsilon  \cdot  \big( \frac{\varepsilon }{r} \big)^{p-\delta } \cdot  r^{1/2} \leqslant  C_{\Rsource,\ell } \,M^{p-\delta } \,r^{1/2}$ because $\varepsilon  \leqslant Mr\leqslant M\Rsource$. This yields
\[\bigg|\scl - \!\! \sum _{p\in \NN \cap [0,m)} \sum _{\ell =0} ^{\lmax} \varepsilon ^{p} \ln ^{\ell } \!\varepsilon  \sum _{n=0} ^{\lfloor {p}\rfloor } \pdm{p} \bigg|_{r,\theta } \leqslant  C_{\fsource,R,M} \left( \frac{r}{R} \right)^{N\pst + \frac{1}{2} } \domi{\ln r}^{2m+1.}\]
Given \cref{eq: decompo modale SE 1 preuve 2}, it only remains to prove
\begin{equation}{\label{eq: decompo modale SE 1 preuve estimation terme 2}}
	\bigg|\sum _{p\in \NN \cap [m,2m)} \sum _{\ell =0} ^{\lmax} \varepsilon ^{p} \ln ^{\ell } \!\varepsilon  \sum _{n=\lfloor {p-m}\rfloor +1} ^{m} \pdm{p} \bigg|_{r,\theta } \leqslant C_{\fsource,R,M} \,\. \varepsilon ^{m} \.r^{\pst-m} \.\domi{\ln \varepsilon }^{\ell _{m}} \domi{\ln r}^{2m+1} .
\end{equation}
To begin, let us estimate each $\pdm{p}$. This part is similar to the proof of \cref{convergence des series de mcl}. By~\cref{estimation des mcl}, we know that, for~any $\eAMP>1$, $d\in \TZ*\cap [-m,\infty )$, $n\in \mathbb{N}$, $0<r\leqslant R$ and $\forall \theta \in [0,\Theta ]$, we have
\[|\mcle_{d,n} (r,\theta )|_{r,\theta } \;\leqslant \; C_{n,R,\eAMP} \,(\eAMP r)^{d-n} \, \domi{\ln r}^{2n+1} .\]
Let $R':=(2R+\Rsource)/3$ and $R'' :=(R+2\Rsource)/3$. We have $R<R'<R''<\Rsource$. \Cref{decompo modale des CL} states that $\sigma _{d} (u_{p,\ell } ) \leqslant  C_{p,\fsource,R} \,R''^{-d}$. Hence we have
\[|\pdm{p}|_{r,\theta } \leqslant  C_{p,n,\fsource,\eAMP,R} \sum _{\substack{d\in \TZ*\\ -p+n+\pst\leqslant d\leqslant N\pst}} R''^{-d} \cdot  (\eAMP r)^{d-n} \, \domi{\ln r}^{2n+1} \vspace{-1ex}\]
Taking $\eAMP := \frac{R''}{R'}$, we obtain:
\begin{align*}
	\varepsilon ^{p} \,|\pdm{p}|_{r,\theta } &\leqslant  \varepsilon ^{p} \cdot C_{p,n,\fsource,R} \sum _{\substack{d\in \TZ*\\ -p+n+\pst\leqslant d\leqslant N\pst}} \left( \frac{r}{R'} \right)^{d-n} \, \domi{\ln r}^{2n+1} \\
	&= \varepsilon ^{p} \cdot  C_{p,n,\fsource,R} \, \frac{1}{1-\left( \frac{r}{R'} \right)^{\smallpst} } \left( \frac{r}{R'} \right)^{-p+\pst} \, \domi{\ln r}^{2n+1} \\
	&\leqslant  C_{p,n,\fsource,R,M} \,\. \varepsilon ^{m} \.r^{\pst-m} \.\domi{\ln r}^{2m+1} 
\end{align*}
because $(\varepsilon /r)^{p-m} \leqslant  M^{p-m}$. Summing over $(p,\ell ,n)$, this yields \cref{eq: decompo modale SE 1 preuve estimation terme 2} and concludes the proof.
\end{demo}

\subsection{Reexpression of the far fields' corner conditions}
\label{sec: coefs cme}

\begin{textesansboite}
As explained under \cref{eq: def mse}, the functions $\mse_{d}$ take into account the layer of thickness $\varepsilon $ and the fact that $u_{\varepsilon }$ solves Helmholtz's equation. However, they do~not depend on the algebraic matching coefficients $\cmacro$ and $\cmicro$ and the corner profile coefficients $\sigma _{q} (\ccDA_{d,n} )$, so they lack the information of the geometry perturbation at the corner and the matching conditions. This information must be taken into account to approximate accurately $u_{\varepsilon }$ and is contained in the corner conditions of the far fields given by \cref{construction directe des upl}:
\begin{equation}{\label{eq: rappel ME de condition de coin CL}}
	\forall d\in -\txtTN, \qquad  \sigma _{d} (u_{p,\ell } ) = \!\! \sum _{\substack{p'\in \NN\\ p-p'\in \NN+ \frac{2\pi }{\Theta } }} \; \sum _{\substack{d'\in \TZ*\\ p'+d'\in \intvl{0,p+d}}} \! \sum _{\ell '=0} ^{\ell } \cfinal_{d,d',p-p',\ell -\ell '} \cdot \sigma _{d'} (\uDA_{p',\ell '} )
\end{equation}
(with $\sigma _{d} (u_{p,\ell } )=\sigma _{d} (\bfuDA_{p,\ell } )=\sigma _{d} (\uDA_{p,\ell } )$ in the right-hand side, by \cref{coincidence des sigma singuliers} and $\uDA_{p,\ell } :=\bfuDA_{p,\ell |\Omega }$).
Another point of view is that \cref{eq: rappel ME de condition de coin CL} implies linear constraints on the coefficients $\sigma _{d} (\scl)$ of the modal decomposition of \cref{decompo modale SE 1} that must be taken into account.\\

This issue is solved by introducing a second modal decomposition of $\scl$ (\cref{decompo modale SE 2}) involving new modes denoted $\mme_{n}$ ($n\in \mathbb{N}^*$). It has roughly speaking the form
\[\scl \approx  \sum _{n=1} ^{N} \sigma _{n\pst} (\sclDA) \cdot \mme_{n} .\]
For~any $n\in \mathbb{N}^*$, $\mme_{n}$ is equal to $\mse_{n\pst}$ plus a linear combination of the $\mse_{d}$ with $d<0$ that depends on the coefficients $\cfinal$ of \cref{eq: rappel ME de condition de coin CL}. The constraints on the coefficients $\smash{\sigma _{d} (\scl)}$ due to \cref{eq: rappel ME de condition de coin CL} are taken into account directly in the modes $\mme_{n}$.\\

To go from a decomposition involving $\sigma _{d} (\scl)$ for $d\in \TZ*$ to a decomposition involving $\sigma _{d} (\scl)$ only when $d>0$, we express the $\sigma _{d} (\uDA_{p,\ell } )$ with $d<0$ in terms of the $\sigma _{d'} (\uDA_{p',\ell '} )$ with $d'>0$. Notice that \cref{eq: rappel ME de condition de coin CL} expresses the $\sigma _{d} (\uDA_{p,\ell } )$ with $d'<0$ in terms of the $\sigma _{d} (\uDA_{p',\ell '} )$ with $d'\in \TZ*$ and $p'<p$. By~replacing, in the right-hand side of \cref{eq: rappel ME de condition de coin CL}, the $\sigma _{d'} (\uDA_{p',\ell '} )$ with $d'<0$ with their expression also given by \cref{eq: rappel ME de condition de coin CL}, and repeating this iteratively, we can indeed express the coefficients with $d<0$ in terms of the ones with $d'>0$.\\

To make this idea rigorous, we introduce the coefficients $\cme_{d,d',p,\ell }$. They are defined for~any $(d,d',p,\ell )\in (-\TN)\times (\TN)\times \NN\times \mathbb{N}$ inductively w.r.t.\ $p$ by:
\begin{itemize}
	\item $\cme_{d,d',0,\ell } := 0$ for any $(d,d',\ell )$,
	\item and then when $p>0$: $\forall (d,d',\ell ),$ 
	\begin{equation}{\label{eq: def cme}}
		\cme_{d,d',p,\ell } := \cfinal_{d,d',p,\ell } + \sum _{\substack{(p_{1} ,p_{2} )\in (\NN+2\pst)^2  \\p_{1} +p_{2} =p}} 
		\sum _{\substack{d''\in -\TN\\ d''\in (d+p_{1} -\mathbb{N})\cap (d'-p_{2} +\mathbb{N})}}
		\sum _{\substack{(\ell _{1} ,\ell _{2} )\in \mathbb{N}^2  \\\ell _{1} +\ell _{2} =\ell }} \cfinal_{d,d'',p_{1} ,\ell _{1} } \cdot \cme_{d'',d',p_{2} ,\ell _{2} }.
	\end{equation}
\end{itemize}
\nomenclature[12-30-30]{$\cme$}{Matching coefficients that give the $\sigma _{d} (\uDA_{p,\ell } )$ with $d<0$ in terms of the $\sigma _{d'} (\uDA_{p',\ell '} )$ with $d>0$ and $p'<p$}
Since $p_{2} =p-p_{1} \leqslant p-2\pst<p$ in the above sum, this is indeed an inductive definition. Moreover this sum has a finite number of terms. The claimed expression of the $\sigma _{d} (\uDA_{p,\ell } )$ ($d<0$) with the $\sigma _{d'} (\uDA_{p',\ell '} )$ ($d'>0$) is proven in \cref{condition de coin CL iteree}. Its proof uses the following small lemma to simplify sum manipulations.
\end{textesansboite}

\begin{lemmesansboite}{}
\label{coefs cme non nuls}%
For~any $(d,d',p,\ell )\in (-\TN)\times (\TN)\times \NN\times \mathbb{N}$, we have: $\cme_{d,d',p,\ell } \neq 0 \Rightarrow  p\in \NN+2\pst$ and $p+d-d'\in \mathbb{N}$.
\end{lemmesansboite}

\begin{demo}
Let us prove by induction on $p$. The property is trivial for $p=0$, so let us assume $p>0$. Let us assume that we have $\cme_{d,d',p,\ell } \neq 0$. Then we have
\begin{itemize}
	\item either $\cfinal_{d,d',p,\ell } \neq 0$, which directly implies $p\in \NN+2\pst$ and $p+d-d'\in \mathbb{N}$ by \cref{coefs cfinal non nuls},
	\item or the complicated sum in \cref{eq: def cme} is non empty. In this case, let $(p_{1} ,p_{2} ,\ell _{1} ,\ell _{2} ,d'')$ be the index of one of the terms of this sum. Then we have $p=p_{1} +p_{2} \in (\NN+2\pst)+(\NN+2\pst)\subset \NN+2\pst$ and $p+d-d' = (p_{1} +d-d'')+(p_{2} +d''-d') \in \mathbb{N}+\mathbb{N}=\mathbb{N}$.
\end{itemize}
\end{demo}

\begin{proposition}{reexpression of the far fields' corner conditions}
\label{condition de coin CL iteree}%
For~any $p\in \NN$, $\ell \in \mathbb{N}$ and $d\in -\TN$, we have\vspace{-.5ex}
\begin{equation}{\label{eq: condition de coin CL iteree}}
	\sigma _{d} (\uDA_{p,\ell } ) = \!\! \sum _{\substack{p'\in \NN\\ p-p'\in \NN+ \frac{2\pi }{\Theta } }} \; \sum _{\substack{\red{d'\in \TN}\\ p'+d'\in \intvl{0,p+d}}} \! \sum _{\ell '=0} ^{\ell } \cme_{d,d',p-p',\ell -\ell '} \cdot \sigma _{d'} (\uDA_{p',\ell '} ).
\end{equation}\fbeq
\end{proposition}

\begin{remarque}
\def\cpath{c^{\rm paths}}%
\def\flecheinv{\,\reflectbox{$\leadsto $}\,}%
Let us give another point of view on the coefficients $\cme$ to explain \cref{eq: def cme,eq: condition de coin CL iteree}, which is actually how I found these formulas. Let us consider the oriented graph $G=(V,E)$ defined as follows.
\begin{itemize}
	\item The vertex set is $V:=\TZ*\times \NN\times \mathbb{N}$. In addition, for any $(d,p,\ell )\in V$, we denote $s_{v} :=\sigma _{d} (\uDA_{p,\ell } )$.
	\item For~any $v=(d,p,\ell )$ and $v'=(d',p',\ell ')$ in $V$, we denote $c_{(v',v)} :=\cfinal_{d,d',p-p',\ell -\ell '}$ if it is well-defined (i.e.\ if $d\in -\TN$, $p-p'\in \NN$ and $\ell -\ell '\in \mathbb{N}$) and $c_{(v',v)} :=0$ otherwise. The set of edges is $E:=\{e\in V^2  \mid c_{e} \neq 0\}$.
\end{itemize}
Let $V_{\pm } :=\{(d,p,\ell )\in V \mid \pm d>0\}$. \Cref{eq: rappel ME de condition de coin CL} can be written as
\[\forall v\in V_{-} ,\qquad  s_{v} = \sum _{v'\in V,v'\leadsto v} c_{(v',v)} \cdot  s_{v'}\]
where $v'\leadsto v$ means $(v',v)\in E$. Note that, for~any edge $(v',v)\in E$, its end point $v$ belongs to $V_{-}$ because the $\cfinal_{d,d',p,\ell }$ are defined only when $d<0$. Under some assumptions that are satisfied here (namely $V= V_{+} \sqcup V_{-}$ and $G$ contains no loop nor semi-infinite path of the form $v_{1} \flecheinv v_{2} \flecheinv v_{3} \flecheinv\cdots $), graph theory shows that we have
\begin{equation}{\label{eq: point de vue sur cme avec graphes}}
	\forall v\in V_{-} , \qquad  s_{v} = \sum _{v'\in V_{+} } \cpath_{(v',v)} \cdot  s_{v'}
\end{equation}
\begin{equation}{\label{eq: def cpath}}
	\text{where}  \qquad  \cpath_{(v',v)} := \sum _{\substack{\text{path} \; v_{0} \leadsto v_{1} \leadsto \cdots \leadsto v_{n} \\ \text{with} \; n\in \mathbb{N}^*,\, v_{0} =v'\!,\, v_{n} =v,\\ v_{1} ,\dots ,v_{n-1} \in V_{-} }} \! c_{(v_{0} ,v_{1} )} \,c_{(v_{1} ,v_{2} )}\cdots  c_{(v_{n-1} ,v_{n} )}.
\end{equation}
The sum above is indexed by all the possible paths that go from $v'$ to $v$. Since any edge of $G$ ends in $V_{-}$, the intermediate points of those paths belong to $V_{-}$. See \cref{fig: graphe cfinal et cme}. We claim that $\cpath_{e} = \cme_{d,d',p-p',\ell -\ell '}$ for any $e=((d,p,\ell ),(d',p',\ell '))\in E$. To see it, let us take apart the terms of \cref{eq: def cpath} s.t.\ $n=1$ and denote $v'':=v_{n-1}$ for the other terms. This gives:
\begin{align*}
	\cpath_{(v',v)} &= c_{(v',v)} + \sum _{\substack{v''\in V_{-} \\v''\leadsto v}} \; \sum _{\substack{\text{path} \; v_{0} \leadsto v_{1} \leadsto \cdots \leadsto v_{n-1} \\ \text{with} \; n\in \mathbb{N}^*,\, n\geqslant 2 ,\, v_{0} =v'\!,\, v_{n-1} =v''\\ v_{1} ,\dots ,v_{n-1} \in V_{-} }} \! c_{(v_{0} ,v_{1} )} \,c_{(v_{1} ,v_{2} )}\cdots  c_{(v_{n-2} ,v_{n-1} )} \,c_{(v'',v)} \\
	&= c_{(v',v)} + \sum _{\substack{v''\in V_{-} \\v''\leadsto v}} \cpath_{(v',v'')} \cdot c_{(v'',v)}
\end{align*}
Thus $\cpath$ satisfies an inductive relation analoguous to the inductive definition of $\cme$ (to properly verify that $\cpath$ and $\cme$ coincide, one would also need to check the initialization). Therefore, \cref{eq: point de vue sur cme avec graphes} is a rewritting of \cref{eq: condition de coin CL iteree}. However, we chose to prove \cref{condition de coin CL iteree} using more elementary tools, mainly sum manipulations.
\end{remarque}

	\tikzset{external/export=true}%
	\tikzsetfigurename{fig-graphe_cfinal_cme-}%
	\begin{figure}[h]
\catcode`\;=12
\def\Lgraphe{4.5} 
\def\Hgraphe{4.5} 
\def\Laxes{4.4} 
\def\Haxes{4.6} 
\def\nanobox#1{\vbox to 0pt{\vss\hbox to 0pt{\hss#1\hss}\vss}} 
\def\sommet#1(#2,#3);{
	\node (#1) at ({(#2)*\Lgraphe/100}, {(#3)*\Hgraphe/100}) {\nanobox{$\bullet$}}; node[above right]{\tiny #1}; 
}
\def\arete[#1] #2 -> #3;{
	\if\relax#1\relax \def\couleurArete{black}
	\else \def\couleurArete{#1} \fi
	\draw [->,>=latex,thick,\couleurArete] #2 -- #3;}

\begin{tikzpicture}[scale=0.75]
	\draw [->] (-\Laxes,0) -- (\Laxes,0) node [right] {$d$};
	\draw [->] (0,-0.25) -- (0,\Haxes) node [right] {$p$};
	\draw [decorate, decoration={brace,amplitude=10pt}, yshift=-6pt] (\Laxes-0.2,0) -- (0.2,0) node [midway,below,yshift=-10pt] {$V_+$};
	\draw [decorate, decoration={brace,amplitude=10pt}, yshift=-6pt] (-0.2,0) -- (-\Laxes+0.2,0) node [midway,below,yshift=-10pt] {$V_-$};
	\sommet Ap(27,15);
	\sommet Bp(50,34);
	\sommet Cp(35,55);
	\sommet Dp(70,65);
	\sommet An(-22,22);
	\sommet Bn(-50,50);
	\sommet Cn(-55,75);
	\sommet Dn(-30,85);
	\sommet En(-85,90);
	\node (Bng) at (-50.4*\Lgraphe/100,49.9*\Hgraphe/100) {};
	\node (Bnd) at (-49.6*\Lgraphe/100,50.1*\Hgraphe/100) {};
	\node (Cng) at (-55.4*\Lgraphe/100,74.9*\Hgraphe/100) {};
	\node (Cnd) at (-54.6*\Lgraphe/100,75.1*\Hgraphe/100) {};
	\arete[red] (Ap) -> (An);
	\arete[green!70!black] (Ap) -> (Bn);
	\arete[] (Bp) -> (Cn);
	\arete[blue] (Ap) -> (Cn);
	\arete[] (Cp) -> (Dn);
	\arete[] (Dp) -> (Dn);
	\arete[] (Dp) -> ++(-0.5,0.4);
	\arete[] (Dp) -> ++(-0.2,0.5);
	\arete[red] (An) -> (Bn);
	\draw[-{Latex[left,length=4pt]},thick,red] (Bng) -- (Cng);
	\draw[-{Latex[right,length=4pt]},thick,green!70!black] (Bnd) -- (Cnd);
	\arete[] (Bn) -> (En);
	\arete[] (Cn) -> (Dn);
	\arete[] (Cn) -> (En);
	\arete[] (Dn) -> ++(-0.1,0.5);
	\arete[] (En) -> ++(-0.3,0.5);
	\arete[] (En) -> ++(0.3,0.5);
	\draw (Ap) node[right]{$v'$};
	\draw (Cn) ++(0,0.1) node[above]{$v$};
\end{tikzpicture}

\label{fig: graphe cfinal et cme}
\end{figure}%
	\tikzset{external/export=false}%
	\tikzsetfigurename{main-figure}

\begin{demo}[Proof of \cref{condition de coin CL iteree}]
Let us prove the property by induction on $p$. When $p=0$, the right-hand sides of \cref{eq: condition de coin CL iteree} and \cref{eq: rappel ME de condition de coin CL} both vanish thanks to the conditions on the index $p'$. This proves \cref{eq: condition de coin CL iteree} in this case.\\
Now, let us assume $p>0$. Let $\ell \in \mathbb{N}$ and $d\in -\TN$. For~any $(d',p',\ell ')$, let us denote $\cmebis_{d,d',p',\ell '} := \cme_{d,d',p',\ell '} - \cfinal_{d,d',p',\ell '}$. In addition, for any family $(x_{p',\ell '} )$, let us denote $x_{\bullet} :=(x_{p',\ell '} )$ the family itself. We also denote $*$ the convolution product w.r.t.\ $(p',\ell ')$.\\
\Cref{coefs cfinal non nuls,coefs cme non nuls} show that the conditions ``$(p_{1} ,p_{2} )\in (\NN+2\pst)^2 $'' and ``$d''\in (d+p_{1} -\mathbb{N})\cap (d'-p_{2} +\mathbb{N})$'' in \cref{eq: def cme} are not necessary because removing them only adds null terms. Therefore, we have
\begin{equation}{\label{eq: condition de coin CL iteree preuve cmebis}}
	\cmebis_{d,d',\bullet} = \sum _{d''\in -\TN} \cfinal_{d,d'',\bullet} *\cme_{d'',d',\bullet} .
\end{equation}
Similarly, using \cref{ansatz sur Aeps,coefs cfinal non nuls,coefs cme non nuls}, one can verify that \cref{eq: condition de coin CL iteree} rewrites as
\[\sigma _{d} (\uDA_{\bullet} ) = \!\sum _{d'\in \TN} \cme_{d,d',\bullet} *\sigma _{d'} (\uDA_{\bullet} )\]
while \cref{eq: rappel ME de condition de coin CL} is equivalent to
\[\sigma _{d} (\uDA_{\bullet} ) = \!\sum _{d'\in \TZ*} \cfinal_{d,d',\bullet} *\sigma _{d'} (\uDA_{\bullet} ) .\]
Now, by~definition of $\cmebis$, we have
\begin{equation}{\label{eq: condition de coin CL iteree preuve sep en 2 bouts}}
	\sum _{d'\in \TN} \cme_{d,d',\bullet} *\sigma _{d'} (\uDA_{\bullet} ) = \sum _{d'\in \TN} \cfinal_{d,d',\bullet} *\sigma _{d'} (\uDA_{\bullet} ) + \sum _{d'\in \TN} \cmebis_{d,d',\bullet} *\sigma _{d'} (\uDA_{\bullet} ).
\end{equation}
Let us focus on the second term. \Cref{eq: condition de coin CL iteree preuve cmebis} implies
\begin{align*}
	\sum _{d'\in \TN} \cmebis_{d,d',\bullet} *\sigma _{d'} (\uDA_{\bullet} )
	&= \sum _{d'\in \TN} \bigg(\sum _{d''\in -\TN} \cfinal_{d,d'',\bullet} *\cme_{d'',d',\bullet} \bigg) *\sigma _{d'} (\uDA_{\bullet} )\\
	&= \sum _{d''\in -\TN} \cfinal_{d,d'',\bullet} *\bigg( \sum _{d'\in \TN}  \cme_{d'',d',\bullet} *\sigma _{d'} (\uDA_{\bullet} ) \bigg)
\end{align*}
because $*$ is associative and the above infinite sums are actually finite for each value of $(p,\ell )$. In~particular, at rank $(p,\ell )$ we have
\[\sum _{d'\in \TN} \cmebis_{d,d',\bullet} *\sigma _{d'} (\uDA_{\bullet} ) \Big|_{p,\ell } = \sum _{p',\ell '} \sum _{d''\in -\TN} \cfinal_{d,d'',p-p',\ell -\ell '} \cdot \bigg( \sum _{d'\in \TN}  \cme_{d'',d',\bullet} *\sigma _{d'} (\uDA_{\bullet} ) \bigg)\bigg|_{p',\ell '}\]
where ``$|_{p',\ell '}$'' means that the result of the convolution is evaluated at a given $(p',\ell ')$. By~\cref{coefs cfinal non nuls}, we have $\cfinal_{d,d'',p-p',\ell -\ell '} \neq 0 \Rightarrow  p-p'\in \NN+2\pst$, so the non zero terms of the above sum satisfy $p'<p$. Hence we can apply the induction hypothesis on them, which yields
\[\sum _{d'\in \TN} \cmebis_{d,d',\bullet} *\sigma _{d'} (\uDA_{\bullet} ) \Big|_{p,\ell } = \sum _{p',\ell '} \sum _{d''\in -\TN} \cfinal_{d,d'',p-p',\ell -\ell '} \cdot \sigma _{d''} (\uDA_{p',\ell '} ).\]
Adding the first term of \cref{eq: condition de coin CL iteree preuve sep en 2 bouts}, we finally obtain
\[\sum _{d'\in \TN} \cme_{d,d',\bullet} *\sigma _{d'} (\uDA_{\bullet} )\Big|_{p,\ell } =\sum _{p',\ell '} \sum _{d''\in \red{\TZ*}} \cfinal_{d,d'',p-p',\ell -\ell '} \cdot \sigma _{d''} (\uDA_{p',\ell '} ) =\sigma _{d} (\uDA_{p,\ell } ) .\]
\end{demo}

\subsection{Final decomposition}
\label{sec: decompo modale de SE 2}

\begin{definition}{modes of the DtN}
\label{def: mme}%
For~any $N\in \mathbb{N}$, $\varepsilon \geqslant 0$ and $n\in [\![1,N]\!]$, we define
\[\mme_{n} := \mse_{n\smallpst} + \sum _{\substack{d\in \TZ*\\n\pst-m<d<0}} \cme_{d,n\pst,\varepsilon } \cdot  \mse_{d} .\vspace{-2ex}\]
where $\displaystyle \cme_{d,d',\varepsilon } := \sum _{p\in \NN \cap [0,m)} \, \sum _{\ell =0} ^{\ell _{p}} \varepsilon ^{p} \ln ^{\ell } \!\varepsilon  \cdot  \cme_{d,d',p,\ell } $ for~any $d\in -\TN$ and $d'\in \TN$.
\end{definition}
\nomenclature[58-23-40-23]{$\mme_{n}$}{Modes of the effective models' DtN maps}

\begin{remarquenumerotee}
\label{rq: mme}%
\begin{itemize}
	\item The modes $\mme_{n}$ and the coefficients $\cme_{d,d',\varepsilon }$ depend on $m$, although we do~not write it in the notation.
	\item Like the functions $\mse_{d}$, the $\mme_{n}$ approximately satisfy the Helmholtz equation in $\Omega $ and the effective condition \cref{eq: condition effective sur GME} on $\Gamma $ that takes into account the presence of the layer. What $\mme_{n}$ does that $\mse_{d}$ does~not is taking into account the matching conditions.
	\item Using \cref{coefs cme non nuls}, one can verify that the condition $n\pst-m<d$ in the definition of $\mme_{n}$ does~not exclude any non-zero term. In addition, this condition implies $\mme_{n} =\mse_{n\smallpst}$ when $n\pst\geqslant m$, so only the first modes $\smash{\mme_{n} }$ and $\smash{\mse_{n\smallpst} }$ differ. \Cref{sec: explicitation des mme} describes more explicitly the modes $\mme_{n}$ depending on the angle $\Theta $.
	\item \Cref{coefs cme non nuls} shows that $\cme_{d,d',p,\ell } \neq 0 \Rightarrow  p\geqslant 2\pst>1$, so $\cme_{d,d',\varepsilon } =O(\varepsilon )$ for~any $(d,d')$ and $\mme_{n} =\mse_{n\smallpst} +O(\varepsilon )$ with a big $O$ that depends on $n$. Hence $\mme_{n}$ is a small perturbation of $\mse_{n\smallpst}$ when $\varepsilon $ is small. Moreover, denoting $d:=n\pst$, we have when $\varepsilon =0$:
	\begin{equation}{\label{eq: mme qd eps=0}}
	\begin{array}{r@{}ll}
		\mme0_{n} &{}= \mse0_{d} & \text{by \cref{def: mme}}  \\[1.5\jot]
		&{}=\mcle_{d} & \text{by \cref{eq: def mse}}  \\[1.5\jot]
		&{}=\mcll_{d,0} & \text{by \cref{def: mcle}}  \\[1.5\jot]
		&{}=\displaystyle (-1)^{n} \,\Gamma (d+1) \left( \frac{2}{k_{0}} \right)^{\!d} J_{d} (k_{0} r) \sin (d\theta ) & \text{by \cref{expression sing homog CL}} 
	\end{array}
	\end{equation}
	where $J_{d} (k_{0} r) \sin (d\theta )$ is a mode of Helmholtz's equation in a Dirichlet corner without layer.
\end{itemize}
\end{remarquenumerotee}

\begin{textesansboite}
\Cref{decompo modale SE 2} will prove that $\scl$ can be decomposed on the $\mme_{n}$ and provide error estimates, very similarly to \cref{decompo modale SE 1}. To show this result, we first express the $\sigma _{d} (\sclDA)$ with $d<0$ in terms of the $\sigma _{d} (\sclDA)$ with $d>0$, using the coefficients $\cme_{d,d',\varepsilon }$. This is done in \cref{condition de coin sur scl et cme_eps}. The proof uses \cref{coefs cme non nuls wrt l} below to control the powers of $\ln \varepsilon $ in the error estimates. Note that \cref{coefs cme non nuls wrt l} shows that stopping the sum in $\ell $ at $\ell _{p}$ in the definition of $\cme_{d,d',\varepsilon }$ (\cref{def: mme}) does~not exclude any non-zero coefficient $\cme_{d,d'p,\ell }$.
\end{textesansboite}

\begin{lemmesansboite}{}
\label{coefs cme non nuls wrt l}%
For~any $(d,d',p,\ell )\in (-\TN)\times (\TN)\times \NN\times \mathbb{N}$, we have: $\cme_{d,d',p,\ell } \neq 0 \Rightarrow  \ell \leqslant \ell _{p}$.
\end{lemmesansboite}

\begin{demo}
Let us show the property $\mathcal{P}$ defined as ``for any $(d,d',p,\ell )\in (-\TN)\times (\TN)\times \NN\times \mathbb{N}$, if $\cme_{d,d',p,\ell } \neq 0$ then $\sigma _{d} (u_{p,\ell } )$ is generically non-zero''. Note that the statement ``$\cme_{d,d',p,\ell } \neq 0$'' depends on the values of $\mu _{0} ,\mu _{1} ,\rho _{0} ,\rho _{1} ,\omega $ and the corner profile coefficients, while ``$\sigma _{d} (u_{p,\ell } )$ is generically non-zero'' is independant of it. The property $\mathcal{P}$ implies the desired result because, by construction of $\ell _{p}$, if $\sigma _{d} (u_{p,\ell } )$ is generically non-zero then $\ell \leqslant \ell _{p}$.\\
The general idea is the following. The expression of $\sigma _{d} (\uDA_{p,\ell } )$ given by \cref{condition de coin CL iteree} contains a term $\cme_{d,d',p,\ell } \cdot \sigma _{d'} (\uDA_{0,0} )$ where $\sigma _{d'} (\uDA_{0,0} )$ is generically non-zero because $d'\in \TN$. So, intuitively, if $\cme_{d,d',p,\ell } \neq 0$ then $\sigma _{d} (\uDA_{p,\ell } )$ is very likely to be non-zero. To make this idea rigorous, we need to use the actual definition of genericity (\cref{def: genericite}), that relies on the coefficients $\cfinal$. Hence we will use the inductive definition of the $\cme$ in terms of the $\cfinal$.\\
Let us prove $\mathcal{P}$ by induction on $p\in \NN$. The case $p=0$ is trivial because we have $\cme_{d,d',p,\ell } =0$ by \cref{coefs cme non nuls}. Now we assume $p>0$. Let us assume $\cme_{d,d',p,\ell } \neq 0$ and prove that $\sigma _{d} (u_{p,\ell } )$ is generically non-zero. By~definition of $\cme_{d,d',p,\ell }$ in \cref{eq: def cme}, we are in one of the two following cases.
\begin{itemize}
	\item \underline{Case 1: $\cfinal_{d,d',p,\ell } \neq 0$.} By~\cref{coefs cfinal non nuls}, we have $p\in \NN+2\pst$ and $p+d-d'\in \mathbb{N}$. By~assumption, we also have $d'\in \TN\subset \NN$. Therefore, $\cfinal_{d,d',p,\ell }$ appears in the right-hand side of \cref{eq: rappel ME de condition de coin CL} (with $p':=\ell ':=0$), which we recall here:
	\begin{equation*} 
		\sigma _{d} (u_{p,\ell } ) = \!\! \sum _{\substack{p'\in \NN\\ p-p'\in \NN+ \frac{2\pi }{\Theta } }} \; \sum _{\substack{d'\in \TZ*\\ p'+d'\in \intvl{0,p+d}}} \! \sum _{\ell '=0} ^{\ell } \cfinal_{d,d',p-p',\ell -\ell '} \cdot \sigma _{d'} (\uDA_{p',\ell '} ). \tag*{\cref{eq: rappel ME de condition de coin CL}}
	\end{equation*}
	In this expression, $\cfinal_{d,d',p,\ell }$ is multiplied by $\sigma _{d'} (\uDA_{0,0} )$, which is generically non-zero by \cref{def: genericite} because $d'\in \TN$. The coefficient $\cfinal_{d,d',p,\ell }$ itself is also generically non-zero by \cref{generiqmt non nul ou tjs nul} because of it~is non-zero (for the current value of $\mu _{0} ,\mu _{1} ,\rho _{0} $, etc.). Hence, by \cref{def: genericite}, $\sigma _{d} (\uDA_{p,\ell } )$ is generically non-zero.
	\item \underline{Case 2.} Otherwise, the expression of $\cme_{d,d',p,\ell }$ written in \cref{eq: def cme} contains a non-zero term of the form $\cfinal_{d,d'',p_{1} ,\ell _{1} }\cdot \cme_{d'',d',p_{2} ,\ell _{2} }$. By~induction hypothesis, $\sigma _{d''} (\uDA_{p_{2} ,\ell _{2} })$ is generically non-zero. Using again \cref{eq: rappel ME de condition de coin CL} (this time with $p':=p_{2}$ and $\ell ':=\ell _{2}$) and \cref{def: genericite}, one can verify similarly that $\sigma _{d} (\uDA_{p,\ell } )$ is generically non-zero.
\end{itemize}
\end{demo}

\begin{textesansboite}
We recall that $\sclDA$ is defined in \cref{eq: def sclDA}.
\end{textesansboite}

\begin{lemme}{}
\label{condition de coin sur scl et cme_eps}%
For~any $d\in -\TN$, we have:
\[\sigma _{d} (\sclDA) = \sum _{d'\in \TN \cap (0,d+m)} \cme_{d,d',\varepsilon } \cdot \sigma _{d'} (\sclDA) + O(\varepsilon ^{m} \ln ^{\ell _{m}} \!\varepsilon )\]
where, more precisely, the $O(\varepsilon ^{m} \ln ^{\ell _{m}} \!\varepsilon )$ is smaller than $C_{\fsource} \,\varepsilon ^{m} \.\domi{\ln \varepsilon }^{\ell _{m}}$.
\end{lemme}

\begin{remarque}
If $d\leqslant -m$, then the above sum is empty. This is coherent with \cref{ansatz sur Aeps}, which implies $\sigma _{d} (\sclDA) =0$ when $d\leqslant -m$.
\end{remarque}

\begin{demo}
Let $d\in -\TN$ and $d'\in \TN \cap (0,d+m)$. By~\cref{def: mme,eq: sigma de sclDA}, and using \cref{coefs cme non nuls wrt l} to drop the condition $\ell \leqslant \ell _{p}$ in the definition of $\cme_{d,d',\varepsilon }$, we have
\begin{align*}
	\cme_{d,d',\varepsilon } \cdot \sigma _{d'} (\sclDA) &= \bigg(\sum _{p\in \NN \cap [0,m)} \, \sum _{\ell \in \mathbb{N}} \varepsilon ^{p} \ln ^{\ell } \!\varepsilon  \cdot  \cme_{d,d',p,\ell } \bigg) \bigg(\sum _{p\in \NN \cap [0,m)} \sum _{\ell \in \mathbb{N}} \varepsilon ^{p} \ln ^{\ell } \!\varepsilon \cdot  \sigma _{d'} (\uDA_{p,\ell } ) \bigg) \\
	&= \sum _{p\in \NN \cap [0,2m)} \sum _{\ell \in \mathbb{N}} \varepsilon ^{p} \ln ^{\ell } \!\varepsilon \cdot  \! \sum _{\substack{p'\in \NN \cap [0,m)\\ p-p'\in \NN \cap [0,m)}} \sum _{\ell '=0} ^{\ell }  \cme_{d,d',p-p',\ell -\ell '} \cdot  \sigma _{d'} (\uDA_{p',\ell '} )  \\
	&= \sum _{p\in \NN \cap [0,\red{m]}} \sum _{\ell \in \mathbb{N}} \varepsilon ^{p} \ln ^{\ell } \!\varepsilon \cdot  \! \sum _{\substack{p'\in \NN \cap [0,m)\\ p-p'\in \NN \cap [0,m)}} \sum _{\ell '=0} ^{\ell } \cme_{d,d',p-p',\ell -\ell '} \cdot  \sigma _{d'} (\uDA_{p',\ell '} ) +  O_{\fsource} (\varepsilon ^{m} )
\end{align*}
where $O_{\fsource} (x)$ denotes a quantity $y$ s.t.\ $|y|\leqslant  C_{\fsource} \,x$. Summing over $d'$, we obtain
\begin{multline*}
	\smash{\sum _{\substack{d'\in \TN\\ d'\in (0,d+m)}} \!\! \cme_{d,d',\varepsilon } \cdot \sigma _{d'} (\sclDA)} = \sum _{p\in \NN \cap [0,m]} \! A_{p} +  O_{\fsource} (\varepsilon ^{m} ) \\
	\text{with}  \quad  A_{p} :=\sum _{\ell \in \mathbb{N}} \varepsilon ^{p} \ln ^{\ell } \!\varepsilon \cdot  \!\! \sum _{\substack{p'\in \NN \cap [0,m)\\ p-p'\in \NN \cap [0,m)}} \; \sum _{\substack{d'\in \TN\\ d'\in (0,d+m)}} \! \sum _{\ell '=0} ^{\ell } \sigma _{d} (\uDA_{p',\ell '} ) \cdot  \cme_{d,d',p-p',\ell -\ell '} .
\end{multline*}
By~\cref{ansatz sur Aeps,coefs cme non nuls}, the non-zero terms of $A_{p}$ satisfy $p-p'\in \NN+2\pst$ and $p'+d'\in \intvl{0,p+d}$. Moreover, we have $\textstyle \{d'\in \TN \mid p'+d'\in \intvl{0,p+d}\} \subset  \TN \cap [-p',p-p'+d] \subset  (0,d+m)$. We deduce
\[\forall p\in [0,m], \qquad  A_{p} =\sum _{\ell \in \mathbb{N}} \varepsilon ^{p} \ln ^{\ell } \!\varepsilon \cdot  \!\! \sum _{\substack{p'\in \NN \cap [0,m)\\ p-p'\in (\NN+2\pst) \cap [0,m)}} \; \sum _{\substack{d'\in \TN\\ p'+d'\in \intvl{0,p+d}}} \! \sum _{\ell '=0} ^{\ell } \sigma _{d'} (\uDA_{p',\ell '} ) \cdot  \cme_{d,d',p-p',\ell -\ell '} .\]
In addition, when $p<m$, the conditions $p'\in [0,m)$ and $p-p'\in [0,m)$ are~not necessary. Using \cref{condition de coin CL iteree} and then \cref{eq: sigma de sclDA}, we deduce
\[\sum _{p\in \NN \cap [0,\red{m)}} A_{p} \, = \sum _{p\in \NN \cap [0,m)} \sum _{\ell \in \mathbb{N}} \varepsilon ^{p} \ln ^{\ell } \!\varepsilon \cdot  \sigma _{d} (\uDA_{p,\ell } ) \,=\, \sigma _{d} (\sclDA).\]
Hence, it only remains to prove that $A_{m} = O_{\fsource} (\varepsilon ^{m} \domi{\ln \varepsilon }^{\ell _{m}} )$. When $p=m$, the condition $(p'<m \text{ and } p-p'<m)$ is equivalent to $(p'>0 \text{ and } p-p'>0)$, so:
\begin{align*}
	A_{m} &= \sum _{\ell \in \mathbb{N}} \varepsilon ^{m} \ln ^{\ell } \!\varepsilon \cdot  \!\!\! \sum _{\substack{p'\in \NN \setminus \{0\}\\ m-p'\in \NN +2\pst}} \; \sum _{\substack{d'\in \TN\\ p'+d'\in \intvl{0,m+d}}} \sum _{\ell '=0} ^{\ell } \cme_{d,d',m-p',\ell -\ell '} \cdot  \sigma _{d'} (\uDA_{p',\ell '} ) \\
	&= \sum _{\ell \in \mathbb{N}} \varepsilon ^{m} \ln ^{\ell } \!\varepsilon \cdot  \!\!\! \sum _{\substack{p'\in \red\NN\\ m-p'\in \NN +2\pst}} \; \sum _{\substack{d'\in \TN\\ p'+d'\in \intvl{0,m+d}}} \sum _{\ell '=0} ^{\ell } \cme_{d,d',m-p',\ell -\ell '} \cdot  \sigma _{d'} (\uDA_{p',\ell '} )\\
	& \quad \  - \sum _{\ell \in \mathbb{N}} \varepsilon ^{m} \ln ^{\ell } \!\varepsilon \cdot  \!\! \sum _{\substack{d'\in \TN\\ d'\in \intvl{0,m+d}}} \sum _{\ell '=0} ^{\ell } \cme_{d,d',\red{m},\ell -\ell '} \cdot  \sigma _{d'} (\uDA_{\red0,\ell '} ) \\
	&= \sum _{\ell \in \mathbb{N}} \varepsilon ^{m} \ln ^{\ell } \varepsilon \cdot \sigma _{d} (\uDA_{m,\ell } ) \,-\, \sum _{\ell \in \mathbb{N}} \varepsilon ^{m} \ln ^{\ell } \!\varepsilon \cdot  \!\! \sum _{\substack{d'\in \TN\\ d'\in \intvl{0,m+d}}} \sum _{\ell '=0} ^{\ell } \cme_{d,d',m,\ell -\ell '} \cdot  \sigma _{d'} (\uDA_{0,\ell '} )
\end{align*}
where we added and substracted the term with $p'=0$, and then used \cref{condition de coin CL iteree}. According to \cref{ansatz sur Aeps,coefs cme non nuls wrt l}, only the terms with $\ell \leqslant \ell _{m}$ in the left sum, and with $\ell '=0$ and $\ell -\ell '\leqslant \ell _{m}$ in the right one, are non-zero. Thus
\[A_{m} \,=\, \sum _{\ell =0} ^{\red{\ell _{m} }} \varepsilon ^{m} \ln ^{\ell } \varepsilon \cdot \sigma _{d} (\uDA_{m,\ell } ) \,-\, \sum _{\ell =0} ^{\red{\ell _{m} }} \varepsilon ^{m} \ln ^{\ell } \!\varepsilon \cdot  \!\! \smash{\sum _{\substack{d'\in \TN\\ d'\in \intvl{0,m+d}}} \cme_{d,d',m,\ell } \cdot  \sigma _{d'} (\uDA_{0,0} )} \,=\, O_{\fsource} (\varepsilon ^{m} \domi{\ln \varepsilon }^{\ell _{m}} ).\]
This concludes the proof.
\end{demo}

\begin{proposition}{decomposition of $\scl$ on the modes $\mme_{n}$}
\label{decompo modale SE 2}%
For~any $N\in \mathbb{N}$ and $\varepsilon >0$, let us denote
\[\rho _{N,\varepsilon } := \scl - \sum _{n=1} ^{N} \sigma _{n\pst} (\sclDA) \cdot \mme_{n} .\]
Let $R\in (0,\Rsource)$ and $M>0$. There exists a positive constant $C_{\fsource,R,M}$ s.t., for~any $\varepsilon >0$, $N\in \mathbb{N}$, $r\in (0,R)$ and $\theta \in [0,\Theta ]$ that satisfy $\varepsilon \.\domi{\ln r}^{3} \leqslant Mr$, we have\vspace{-1ex}
\[|\rho _{N,\varepsilon } (r,\theta )|+r|\nabla  \rho _{N,\varepsilon } (r,\theta )| \leqslant  C_{\fsource,R,M} \left[ \left( \frac{r}{R} \right)^{N\pst+ \frac{1}{2} }  + \frac{\varepsilon ^{m} }{r^{m-\pst}} \,\domi{\ln \varepsilon }^{\ell _{m}} \right] \domi{\ln r}^{2m+1} .\]\fbeq
\end{proposition}

\begin{remarque}
If $\mcle_{d,n} :=\mcll_{d,n}$ for~any $(d,n)$, then it~is easy to see that $\mme_{n}$ does~not depend on $N$, so let us denote it $\Psi ^{\varepsilon } _{n}$. Taking $N\rightarrow \infty $ in \cref{decompo modale SE 2}, we obtain
\[\scl = \sum_{n=1}^\infty  \sigma _{n\pst} (\sclDA) \cdot \Psi ^{\varepsilon } _{n} \;+\, O(\varepsilon ^{m} \ln ^{\ell _{m}} \!\varepsilon  )\]
where the big $O$ depends on $(\fsource,R,M,r)$. This a new approximate modal decomposition of $\scl$ (and of $u_{\varepsilon }$) at order $\varepsilon ^{m} \ln ^{\ell _{m}} \!\varepsilon $, like \cref{eq: rq decompo modale SE 1}. \Cref{decompo modale SE 2} quantifies the error made by truncations of this decomposition. The given error estimate is identical to \cref{decompo modale SE 1}.
\end{remarque}

\begin{remarque}
We prove \cref{decompo modale SE 2} under the simplifying additional assumption $N\geqslant m \Theta /\pi $ (which is not very restrictive because $\Theta /\pi  <2$). The simplification of the proof in this case is related to the fact that we have $\mme_{n} =\mse_{n\pst}$ when $n\pst\geqslant m$ (see \cref{rq: mme}). This property implies that, when $N\pst\geqslant m$, the modes $\mse_{d}$ s.t.\ $d\leqslant N\pst$ are entirely taken into account in $(\mme_{1} ,\dots ,\mme_{n} )$, whereas otherwise there may exist $d\in -\TN$ and $n>N$ s.t.\ $\mme_{n}$ depends on $\mse_{d}$ (see \cref{def: mme}).\\

The assumption $N\geqslant m \Theta /\pi $ can be avoided by using an estimate on $\mse_{d}$ that appears later in \cref{eq: estimation mse} and implies
\begin{equation}{\label{eq: estimation mse pour decompo modale 2}}
	|\mse_{d} (r,\theta )|+r|\nabla  \mse_{d} (r,\theta )| \leqslant  C_{R,M} \,r^{d} \domi{\ln r}
\end{equation}
for~any $r>0$, $\theta \in (0,\Theta )$, $d\in -\TN$ and $M>0$ s.t.\ $\varepsilon  \.\domi{\ln r}^{3} \leqslant  Mr$. Using this estimate, \cref{coefs cme non nuls} and the assumption $\varepsilon \.\domi{\ln r}^{3} \leqslant Mr$ in \cref{decompo modale SE 2}, one can control the error made by the ``missing'' modes $\mse_{d}$ in the step \cref{eq: decompo modale SE 2 preuve 1} of the proof below. Since this is a technical detail, we do~not expand further on this. Note that \cref{eq: estimation mse pour decompo modale 2} does~not rely on \cref{decompo modale SE 2}, so the proof is~not circular. 
\end{remarque}

\begin{demo}[Proof of \cref{decompo modale SE 2}]
\def\diff{\tilde{\rho}_{\Nvirg\peps}}%
\def\pst{\tau}%
For~any $f\in H^{1} _{\mathrm{loc}} (\Omega )$ and $(r,\theta )$, let us denote $|f|_{r,\theta } := |f(r,\theta )|+r|\nabla  f(r,\theta )|$. We also denote $\pst:= \pi /\Theta $ for simplicity, and
\[\diff := \sum _{n=1} ^{N} \sigma _{n\pst} (\sclDA) \cdot \mme_{n} -  \sum _{\substack{d\in \TZ*\\ -m\leqslant d\leqslant N\pst}} \sigma _{d} (\sclDA) \cdot \mse_{d} .\]
Given \cref{decompo modale SE 1}, it suffices to prove
\[|\diff(r,\theta )|_{r,\theta } \leqslant  C_{\fsource,R,M} \, \frac{\varepsilon ^{m} }{r^{m-\pst}} \.\domi{\ln \varepsilon }^{\ell _{m}} \domi{\ln r}^{2m+1} .\]
By~\cref{def: mme}, we have
\begin{align*}
	\sum _{n=1} ^{N} \sigma _{n\pst} (\sclDA) \cdot \mme_{n} &= \sum _{n=1} ^{N} \sigma _{n\pst} (\sclDA) \cdot \mse_{n\pst} \,+\, \sum _{n=1} ^{N} \sigma _{n\pst} (\sclDA) \cdot \! \sum _{\substack{d\in \TZ*\\ n\pst-m<d<0}} \cme_{d,n\pst,\varepsilon } \cdot  \mse_{d} \\
	&= \sum _{d\in \TZ*\cap [0,N\pst]} \sigma _{d} (\sclDA) \cdot \mse_{d} + \sum _{\substack{d\in \TZ*\\ d<0}} \Bigg( \underbrace{\sum _{\substack{d'\in \TZ*\cap [0,N\pst]\\ d'<d+m}} \cme_{d,d',\varepsilon } \cdot \sigma _{d'} (\sclDA) }_{\displaystyle =:A_{d} } \Bigg) \cdot  \mse_{d}
\end{align*}
Let $d\in -\TN$. We have $(d'<d+m \Rightarrow  d'<m\leqslant N\pst)$ because we assumed $N\geqslant m/\pst$. So \cref{condition de coin sur scl et cme_eps} shows
\begin{equation}{\label{eq: decompo modale SE 2 preuve 1}}
	A_{d} = \sum _{d'\in \TZ*\cap (0,d+m)} \cme_{d,d',\varepsilon } \cdot \sigma _{d'} (\sclDA) \,=\, \sigma _{d} (\sclDA) +O_{\fsource} (\varepsilon ^{m} \domi{\ln \varepsilon }^{\ell _{m}} )
\end{equation}
where $O_{\fsource} (x)$ denotes a quantity $y$ s.t.\ $|y|\leqslant C_{\fsource} \, x$. In addition, if $A_{d} \neq 0$, then there exists $d'\in \TN$ s.t.\ $d'<d+m$, so we must have $d>-m+\pst$. We deduce
\begin{align*}
	|\diff(r,\theta )|_{r,\theta } &= \Bigg|\sum _{d\in \TZ*\cap (\red{-m+\pst},0)} \big(A_{d} - \sigma _{d} (\sclDA) \big) \cdot  \mse_{d} \Bigg|_{r,\theta } \\
	&\leqslant  C_{\fsource} \,\varepsilon ^{m} \domi{\ln \varepsilon }^{\ell _{m}} \sum _{d\in \TZ*\cap (-m+\pst,0)} |\mse_{d} |_{r,\theta } \\
	&\leqslant  C_{\fsource} \,\varepsilon ^{m} \domi{\ln \varepsilon }^{\ell _{m}} \sum _{d\in \TZ*\cap (-m+\pst,0)} \sum _{n=0} ^{\min (m-1,\lfloor {d+m}\rfloor )} \varepsilon ^{n} \,|\mcle_{d,n} |_{r,\theta } \\
	&\leqslant  C_{\fsource} \,\varepsilon ^{m} \domi{\ln \varepsilon }^{\ell _{m}} \sum _{d\in \TZ*\cap (-m+\pst,0)} \sum _{n=0} ^{\lfloor {d+m}\rfloor } \varepsilon ^{n} \cdot C_{d,n,R} \, r^{d-n} \domi{\ln r}^{2n+1} \ \ \text{by \cref{estimation des mcl} }  (r\leqslant R)\\
	&\leqslant  C_{\fsource,R,M} \cdot  \varepsilon ^{m} \domi{\ln \varepsilon }^{\ell _{m}} \cdot  r^{-m+\pst} \domi{\ln r}^{2m+1} 
\end{align*}
\end{demo}

\section{Construction of the transparent condition around the corner}
\label{sec: construction du DtN des ME}

\subsection{Definition of the DtN operators}
\label{sec: construction de opS et opP}

\begin{textesansboite}
In this section, we build the DtN operators $\opS$ and $\opP$ (the operator $\opP$ intervenes in the effective model only for $m=5$, but it will also be useful in the error analysis when $m=3$). These two operators are designed so that we have $\opS(\scl_{|\SME} ) \approx  \partial _{r} \scl_{|\SME}$ and $\opP(\scl_{|\SME} ) \approx  \partial _{x} \scl(\PME)$. We recall that $\scl$ is defined in \cref{eq: def notation somme CL} and that $u_{\varepsilon } -\scl=O(\varepsilon ^{m} \ln ^{\ell _{m}} \varepsilon )$ by \cref{DA de SE avec CL et delta variable}. \Cref{decompo modale SE 2} states that $\scl$ can be approximated by $\sum _{n=1} ^{N} \sigma _{n\smallpst} (\sclDA)\cdot \mme_{n}$ up to $O(\varepsilon ^{m} \ln ^{\ell _{m}} \varepsilon )$ and a truncation error that decays to 0 geometrically w.r.t.\ $N$. Hence we will build $\opS$ and $\opP$ s.t.\ we have
\[\forall v\in \operatorname{Span}(\mme_{1} ,\dots ,\mme_{N} ), \qquad  \opS(v_{|\SME} )=\partial _{r} v_{|\SME} \ \text{ and } \ \opP(v_{|\SME} )=\partial _{x} v(\PME).\]
\def\Vmme{V_N}%
Let us denote $\Vmme:=\operatorname{Span}(\mme_{1} ,\dots ,\mme_{N} )$. The construction has three steps:
\begin{itemize}
	\item first, defining linear forms $\lambda _{1} ,\dots ,\lambda _{N} :\Vmme\rightarrow \mathbb{C}$ s.t.\ $v=\smash{\sum\limits _{n=1} ^{N} \lambda _{n} (v_{|\SME} )\,\mme_{n} }$ for~any $v\in \Vmme$,
	\item then setting $\opS(v) :=\sum\limits _{n=1} ^{N} \lambda _{n} (v)\cdot \partial _{r} \mme_{n|\SME} $ and $\opP(v) :=\sum\limits _{n=1} ^{N} \lambda _{n} (v)\cdot \partial _{x} \mme_{n} (\PME)$ for~any $v\in \Vmme$,
	\item finally extending $\opS$ and $\opP$ to $L^2 (\SME)$.
\end{itemize}
The linear forms $\lambda _{n}$ exist iff the functions $\mme_{n|\SME}$ are linearly independent. We will see in \cref{rq: inversibilite de gram} below that this is the case when $\varepsilon $ is small enough. Assuming that we are in that situation, let us express the forms $\lambda _{n}$ using scalar products in $L^2 (\SME)$. We denote the scalar product in $L^2 (\SME)$ as:
\begin{equation}{\label{eq: def produit scalaire L2(SME)}}
	\forall v,w\in L^2 (\SME), \quad \  \langle v,w\rangle _{L^2 (\SME)} := \int _{\SME} v \, \bar{w}=\int _{0} ^{\Theta } v(\theta ) \, \overline{w(\theta )}\,\RME \mathrm{d}\theta 
\end{equation}
or simply $\langle v,w\rangle $ if there~is no ambiguity. Since the modes $\mme_{n}$ are not orthogonal in general, $\lambda _{n} (v)$ is~not simply equal to $\langle v,\mme_{n} \rangle _{L^2 (\SME)}$. However, for any $(\nu _{n} )\in \mathbb{C}^{N}$ and $v_{\nu } :=\sum _{n=1} ^{N} \nu _{n} \mme_{n}$, we have $\big[\langle v_{\nu } ,\mme_{n} \rangle _{L^2 (\SME)} \big]_{n\in [\![1,N]\!]} = \gram^{\mathsf{T}} \cdot \big[\nu _{n} \big]_{n\in [\![1,N]\!]}$, where the brackets denote column vectors and
\begin{equation}{\label{eq: def Gram}}
	\gram :=\Big[\langle \mme_{n} ,\mme_{n'} \rangle _{L^2 (\SME)} \Big]_{(n,n')\in [\![1,N]\!]^2 } \in  \mathcal{M}_{N} (\mathbb{C})
\end{equation}
is the Gram matrix of the family $(\mme_{n|\SME} )_{n}$ (and ${\dots }^{\mathsf{T}} $ denotes the transpose). If the $\mme_{n|\SME}$ are linearly independent, then $\gram$ is invertible and $\lambda _{1} ,\dots ,\lambda _{N}$ satisfy $\big[\lambda _{n} (v)\big]_{n} =(\gram^{\mathsf{T}} )^{-1} \big[\langle v ,\mme_{n} \rangle _{L^2 (\SME)} \big]_{n}$ for~any $v\in L^2 (\SME)$. This yields:
\begin{align}
	\forall v\in \Vmme,\;\forall \bfx \in \SME,\quad  \ \; \sum _{n=1} ^{N} \lambda _{n} (v)\cdot \partial _{r} \mme_{n} (\bfx) &= \Big[\lambda _{n} (v)\Big]_{n} ^{\mathsf{T}} \cdot  \Big[\partial _{r} \mme_{n} (\bfx)\Big]_{n} \nonumber \\
	&= \Big[\langle v ,\mme_{n} \rangle _{L^2 (\SME)} \Big]_{n} ^{\mathsf{T}} \cdot \gram^{-1} \cdot  \Big[\partial _{r} \mme_{n} (\bfx)\Big]_{n} \label{eq: presentation def opS}
\end{align}
\Cref{eq: presentation def opS} will be used as the definition of $\opS$ on $\Vmme$, and similarly for $\opP$. Moreover, \cref{eq: presentation def opS} is also well-defined for~any $v\in L^2 (\SME)$, which gives a natural way to extend $\opS$ to $L^2 (\SME)$. This approach motivates the following definition.
\end{textesansboite}

\begin{definition}{operators of the DtN}
\label{def: opS opP}%
Let $N\in \mathbb{N}$ and $\varepsilon \geqslant 0$ be s.t.\ the matrix $\gram$ defined in \cref{eq: def Gram} is invertible. For~any $v\in L^2 (\SME)$ and $\bfx \in \SME$, we define
\[\opS(v)(\bfx) := \Big[\langle v ,\mme_{n} \rangle _{L^2 (\SME)} \Big]_{n} ^{\mathsf{T}} \cdot \gram^{-1} \cdot  \Big[\partial _{r} \mme_{n} (\bfx)\Big]_{n} \qquad  \! \text{ and } \qquad  \opP(v) := \opS(v)(\PME).\]\fbeq
\end{definition}
\nomenclature[28-22-23-40]{$\opS$}{DtN map on $\SME$ in the effective models}
\nomenclature[25-22-23-40]{$\opP5$}{DtN map on $\PME$ in the effective model of order~5}

\begin{remarque}
It is clear that $\opS$ is continuous from $L^2 (\SME)$ to itself, and $\opP:L^2 (\SME)\rightarrow \mathbb{C}$ is also continuous. In addition, for~any $v,w\in L^2 (\SME)$, we have
\[\int _{\SME} \opS(v)\cdot  \bar{w}= \Big[\langle v ,\mme_{n} \rangle _{L^2 (\SME)} \Big]_{n} ^{\mathsf{T}} \cdot \gram^{-1} \cdot  \Big[\langle \partial _{r} \mme_{n} ,w\rangle \Big]_{n} .\]
$\opS$ appears in this form in the variational formulations \cref{eq: FV ME3,eq: FV ME5}.
\end{remarque}

\begin{remarquenumerotee}
\label{rq: inversibilite de gram}%
For~any $N$, the function $\varepsilon  \mapsto \det(\gram)$ is analytic on $\mathbb{R}_+^*$ because the modes $\mme_{n}$ are linear combinations of the functions $\mcle_{d,n'}$ multiplied by powers of $\varepsilon $ and $\ln \varepsilon $. In addition, \cref{eq: mme qd eps=0} shows that, for $\varepsilon =0$, each mode $\mme0_{n}$ is proportional to $J_{n\pst} (k_{0} r) \sin (n\pst \theta )$ and non-zero, which implies that $\gram0$ is invertible. Thus, the matrix $\gram$ is invertible outside of a discrete set of values of $\varepsilon $ (that depends on $N$) and also on a neighborhood of $\varepsilon =0$. For these $\varepsilon $, $\opS$ and $\opP$ are well-defined. \Cref{opS - opS0} will quantify the size of this neighborhood of $\varepsilon =0$ in terms of $\RME$.
\end{remarquenumerotee}

\begin{remarque}
Another equivalent approach to build the DtN operators is to apply the Gram-Schmidt algorithm to $(\mme_{n} )_{n}$ to obtain a family $\big(\mmeo{n}\big)_{n}$ s.t.\ $\big(\mmeo{n}_{|\SME} \big)_{n}$ is orthonormal in $L^2 (\SME)$. By~considering the $\mmeo{n}$ as the actual modes of the DtN, one could define $\opS$ as:\vspace{-.5ex}
\begin{equation}{\label{eq: def opS avec Gram-Schimdt}}
	\forall v\in L^2 (\SME),\qquad  \opS(v) := \sum _{n=1} ^{N} \langle v,\mmeo{n}\rangle  \cdot \partial _{r} \mmeo{n}.
\end{equation}
With this definition, $\opS$ is well-defined iff $(\mmeo{n})_{n}$ is well-defined, which is equivalent to the linear independence of $(\mme_{n} )_{n}$, which is also equivalent to the inversibility of $\gram$ by properties of the Gram matrix. Let $V:=\operatorname{Span}(\mme_{1} ,\dots ,\mme_{N} )\subset H^{1} _{\mathrm{loc}} (\Omega )$ and $V_{|\SME} :=\operatorname{Span}(\mme_{1|\SME} ,\dots ,\mme_{N|\SME} )\subset L^2 (\SME)$. Since $V=\operatorname{Span}(\mmeo{1},\dots ,\mmeo{N})$, \cref{def: opS opP,eq: def opS avec Gram-Schimdt} both imply $\opS(v_{|\SME} )=\partial _{r} v$ for~any $v\in V$. So these two definitions agree on $V_{|\SME}$. Moreover, they both imply $\opS(v)=0$ for any $v\in (V_{|\SME} )^{\perp}$, which proves that the two definitions of $\opS$ are equivalent. However, numerically, \cref{def: opS opP} is easier to implement (see \cref{sec: numerique ME} for the numerical inversion of $\gram$).
\end{remarque}

\subsection{Explicit description of the DtN modes except for particular angles}
\label{sec: explicitation des mme}

\begin{textesansboite}
We recall that \cref{def: mme} implies $\mme_{n} =\mse_{n\smallpst}$ for any $n\geqslant m \Theta /\pi $. This section describes the modes $\mme_{n}$ when $n<m \Theta /\pi $, except for a few particular angles $\Theta $.\\

According to \cref{simplifications cas irrationnel} (when $\pst \not\in \mathbb{Q}$) and \cref{simplifications cas rationnel,simplifications cas rationnel version forte} (when $\pst \in \mathbb{Q}$), for all angles $\Theta $ except a few ones (given by the sets $E_{m}$ in the next proposition), the matching conditions of the far fields simplify at orders $p<5$ into
\[\forall \ell \in \mathbb{N},\, \forall d\in -\txtTN, \qquad  \sigma _{d} (u_{p,\ell } )= \sum _{\substack{n\in \mathbb{N}\\ 2k\leqslant p+d}} \sum _{\substack{d'\in \TN\\ d'\leqslant p+d-2k}} \omega ^{2k} \sigma _{d} (\ccDA_{d',k} )\cdot \sigma _{d'} (\uDA_{p-d'+d-2k,\ell } ).\]
Using the same techniques as \cref{chap: complements sur le DA}, we are able to show similar simplifications on the expression of the modes $\mme_{n}$, as explained in the next result. All these simplifications come from the fact that, in those cases, the matching coefficients have very simple expressions.
\end{textesansboite}

\begin{proposition}{}
\label{expression mme sauf cas speciaux}%
Let us denote $E_{3} := \{\pi \}$ and $E_{5} := \left\{ \frac{2\pi }{3} , \pi  , \frac{4\pi }{3} , \frac{5\pi }{3} \right\}$. For~any $\Theta  \in (0,2\pi ) \setminus E_{m}$, $\varepsilon >0$, $N\in \mathbb{N}$ and $n\in [\![1,N]\!]$, we have:
\[\mme_{n} = \mse_{n\pst} + \sum _{\substack{d\in -\TN,\, k\in \mathbb{N}\\|d|+2k<m-n\pst}} \varepsilon ^{n\pst-d+2k} \, \omega ^{2k} \, \sigma _{d} (\ccDA_{n\pst,k} ) \cdot \mse_{d} .\]\fbeq
\end{proposition}

\begin{demo}
We recall that, by \cref{def: mme,eq: def cme}, the modes $\mme_{n}$ depends on the coefficients $\cme_{d,d',p,\ell }$, which are defined using the coefficients $\cfinal_{d,d',p,\ell }$. Let us express these objects one after another, beginning with $\cfinal_{d,d',p,\ell }$. We mainly use the results of \cref{sec: simplifications DA} and its extension Appendix~\ref{sec: annexe simplifications DA cas rationnel}.\smallskip

\underline{Step 1:} If $\pst \in \mathbb{Q}$, let us denote $b$ the denominator of $\frac{\Theta }{\pi }$ in irreducible form (i.e.\ $\Theta  = \frac{a\pi }{b}$ with $a\in \mathbb{N}$ and $\gcd(a,b)=1$). Otherwise, we define $b:=+\infty $. Let $(d,d',p,\ell )\in (-\TN)\times \TZ*\times \NN\times \mathbb{N}$ be s.t.
\begin{equation}{\label{eq: expression mme sauf cas speciaux preuve hyp sur p}}
	\textstyle  p<m \, \text{ and } \, (d'\in -\mathbb{N} \Rightarrow  p<b+\pst).
\end{equation}
Let us show:
\begin{equation}{\label{eq: expression mme sauf cas speciaux preuve cfinal}}
	\cfinal_{d,d',p,\ell } = \left\{\begin{array}{ll}
	\omega ^{2k} \sigma _{d} (\ccDA_{d',k} ) & \text{if $d'>0$, $\ell =0$ and there exists $k\in \mathbb{N}$ s.t.\ }  p=d'-d+2k\\
	0 & \text{otherwise}.
	\end{array}\right.
\end{equation}
There are several cases, given that we assumed $\Theta  \not\in E_{m}$.
\begin{enumerate}
	\item If $\pst \not\in \mathbb{Q}$, then \cref{valeur de cmaicro cas irrationnel} shows \cref{eq: expression mme sauf cas speciaux preuve cfinal} (unconditionally on $p$ in fact).
	\item Let us assume $\pst \in \mathbb{Q}$ and $b\not\in F_{m}$, with $F_{3} :=\{1\}$ and $F_{5} :=\{1,3\}$. Let us use the notation $\beta $ from \cref{sec: annexe simplifications DA cas rationnel}, which is an integer s.t.\ for any $d\in \TZ*$ and $n<\beta $, $\mcla_{d,n}$ can be expressed with Bessel functions (more details in \cref{sec: annexe simplifications DA cas rationnel}). \Cref{valeur de cfinal cas rationnel generalise} shows \cref{eq: expression mme sauf cas speciaux preuve cfinal} for~any $p$ s.t.
	\begin{equation}{\label{eq: expression mme sauf cas speciaux preuve cfinal bis}}
		\textstyle p<\min (\beta ,2b+2\pst) \, \text{ and } \, (d'\in -\mathbb{N} \Rightarrow  p<b+\pst),
	\end{equation}
	Let us check that \cref{eq: expression mme sauf cas speciaux preuve cfinal bis} is satisfied. Given the assumption \cref{eq: expression mme sauf cas speciaux preuve hyp sur p}, it suffices to prove $\min (\beta ,2b+2\pst)\geqslant m$, i.e. $\beta \geqslant m$ and $2b+2\pst\geqslant m$. Propositions~\ref{expression sing heterog CL} and \ref{expression sing heterog CL avec BesselY} show that we have $\beta \geqslant b$, which yields $\beta \geqslant m$ when $b\geqslant m$. Otherwise, we have $b\in [\![1,m-1]\!] \setminus F_{m} \subset  \{2,4\}$ and \cref{expression sing heterog CL b=2 ou 4} shows $\beta \geqslant 5\geqslant m$. In addition, we have $b\geqslant 2$ and $\pst\geqslant  \frac{1}{2}$, so $2b+2\pst\geqslant 5$. Therefore \cref{eq: expression mme sauf cas speciaux preuve cfinal bis} holds, which implies \cref{eq: expression mme sauf cas speciaux preuve cfinal}.
	\item The only remaining case is ($\Theta  = \frac{\pi }{3}$ and $m=5$). By~\cref{coefs cfinal non nuls}, $\cfinal_{d,d',p,\ell }$ vanishes when $p<2\pst=6$. The right-hand side of \cref{eq: expression mme sauf cas speciaux preuve cfinal} also vanishes when $p<2\pst$. This proves \cref{eq: expression mme sauf cas speciaux preuve cfinal} in this case.
\end{enumerate}
\smallskip

\underline{Step 2:} Let us show a similar expression for the coefficients $\cme_{d,d',p,\ell }$: for any $(d,d',p,\ell )\in (-\TN)\times \TN\times \NN\times \mathbb{N}$ s.t.\ $p<m$, we have
\begin{equation}{\label{eq: expression mme sauf cas speciaux preuve cme}}
	\cme_{d,d',p,\ell } = \left\{\begin{array}{ll}
	\omega ^{2k} \sigma _{d} (\ccDA_{d',k} ) & \text{if $\ell =0$ and there exists $k\in \mathbb{N}$ s.t.\ }  p=d'-d+2k\\
	0 & \text{otherwise}.
	\end{array}\right.
\end{equation}
By~\cref{eq: def cme}, $\cme_{d,d',p,\ell } -\cfinal_{d,d',p,\ell }$ is equal to a sum of terms of the form $\cfinal_{d,d'',p_{1} ,\ell _{1} } \cdot  \cme_{d'',d',p_{2} ,\ell _{2} }$ with $p_{1} +p_{2} =p$, $p_{1} <p$, $d''<0$ and $d'' \geqslant d'-p_{2}$. If $d''\not\in -\mathbb{N}$, then step~1 shows that we have $\cfinal_{d,d'',p_{1} ,\ell _{1} }=0$ because $p_{1} <m$. Otherwise, we have $p_{2} \geqslant d'-d''\geqslant \pst+b$, which implies $p_{1} =p-p_{2} <m-(b+\pst) \leqslant b+\pst$ (because $2(b+\pst)\geqslant 2(2+ \frac{1}{2} )\geqslant m$). So we can also apply step~1, which gives again  $\cfinal_{d,d'',p_{1} ,\ell _{1} }=0$. Hence we have $\cme_{d,d',p,\ell } =\cfinal_{d,d',p,\ell }$. We conclude using step~1 one more time.\medskip

\underline{Step 3:} According to \cref{def: mme}, step~2 implies for any $d\in -\TN$ and $d'\in \TN$:\vspace{-1ex}
\[\cme_{d,d',\varepsilon } := \sum _{p\in \NN \cap [0,m)} \, \sum _{\ell =0} ^{\ell _{p}} \varepsilon ^{p} \ln ^{\ell } \!\varepsilon  \cdot  \cme_{d,d',p,\ell } = \sum _{k\in \mathbb{N},\,n\pst-d+2k<m} \! \varepsilon ^{n\pst-d+2k} \, \omega ^{2k} \, \sigma _{d} (\ccDA_{n\pst,k} ),\]
which yields ${\displaystyle \mme_{n} := \mse_{n\pst} + \; \sum _{\substack{d\in \TZ*\\ \mathclap{n\pst-m<d<0}}} \; \cme_{d,n\pst,\varepsilon } \cdot  \mse_{d} =\mse_{n\pst} + \!\! \sum _{\substack{d\in -\TN,\, k\in \mathbb{N}\\|d|+2k<m-n\pst}} \! \varepsilon ^{n\pst-d+2k} \, \omega ^{2k} \, \sigma _{d} (\ccDA_{n\pst,k} ) \cdot \mse_{d} }$.
\end{demo}

\begin{exemple}
\def\pst{\tau}%
Let us denote $\pst:= \pi /\Theta $. For small angles, \cref{expression mme sauf cas speciaux} gives the following expressions of the modes $\mme_{n}$.
\begin{itemize}
	\item If $\Theta  \leqslant  \frac{2\pi }{m}$, then $\mme_{n} =\mse_{n\pst}$ for~any $n$.
	\item If $\frac{2\pi }{m} < \Theta  \leqslant  \frac{3\pi }{m} $ and $\Theta \neq \pi $, then $\mme_{1} =\mse_{\pst} + \varepsilon ^{2\pst} \sigma _{-\pst} (\ccDA_{\pst,0} ) \,\mse_{-\pst}$ and $\mme_{n} =\mse_{n\pst}$ for~any $n>1$.
	\item If $\frac{3\pi }{m} < \Theta $ and $\Theta  \leqslant  \frac{4\pi }{3}$ if $m=3$ (resp.\ $\Theta  < \frac{2\pi }{3}$ if $m=5$), then 
	\[\left\{\begin{array}{l@{}l}
		\mme_{1} =\mse_{\pst} + \varepsilon ^{2\pst} \sigma _{-\pst} (\ccDA_{\pst,0} ) \,\mse_{-\pst} + \varepsilon ^{3\pst} \sigma _{-2\pst} (\ccDA_{\pst,0} ) \,\mse_{-2\pst} \\[.5ex]
		\mme_{2} =\mse_{2\pst} + \varepsilon ^{3\pst} \sigma _{-\pst} (\ccDA_{2\pst,0} ) \,\mse_{-\pst} \\[.5ex]
		\mme_{n} =\mse_{n\pst} & \text{for~any }  n>2.
	\end{array}\right.\]
\end{itemize}
\end{exemple}


\section{Well-posedness of the effective models}
\label{sec: ME bien poses}

\makeatletter
\def\aME{\use@m{a^{(\@m)}_\peps}}%
\def\bME{\use@m{b^{(\@m)}_\peps}}%
\makeatother

\begin{textesansboite}
In this section, we prove \cref{caractere bien pose des ME}. The general strategy is to prove the coercivity of the effective models for $\varepsilon $ small enough in order to deduce their well-posedness and uniform stability w.r.t.\ $\varepsilon $. This is possible because we assumed $\Im(\omega ) \neq 0$. Let us define the following bilinear forms
\begin{itemize}
	\item $\forall u,v\in \VME3,\quad  \displaystyle \aME3(u,v) := \int _{\OME} (\mu _{0} \nabla  u \cdot  \nabla  \bar{v}- \omega ^2  \rho _{0} u \bar{v}) + \int _{\GME} \frac{\mu _{1} }{\varepsilon } u \bar{v}$,
	\item $\forall u,v\in \VME5,\quad  \displaystyle \aME5(u,v) := \int _{\OME} (\mu _{0} \nabla  u \cdot  \nabla  \bar{v}- \omega ^2 \rho _{0} u \bar{v}) + \int _{\GME} \left( \frac{\mu _{1} }{\varepsilon } u \bar{v} + \frac{\varepsilon }{3} (\mu _{1} \partial _{x} u \,\partial _{x} \bar{v}- \omega ^2 \rho _{1} u \bar{v}) \right)$,
	\item $\forall u,v\in \VME3,\quad  \displaystyle \bME3(u,v) := \int _{\SME} \mu _{0} \,\opS3 (u_{|\SME} )\, \bar{v}$,
	\item $\forall u,v\in \VME5,\quad  \displaystyle \bME5(u,v) := \int _{\SME} \mu _{0} \,\opS5 (u_{|\SME} )\, \bar{v}\,+\, \frac{\varepsilon  \mu _{1} }{3} \,\opP5 (u_{|\SME} ) \cdot  \overline{v(\PME)}$,
\end{itemize}
The effective model of order $m=3$ given in \cref{eq: FV ME3}, resp.\ $m=5$ given in \cref{eq: FV ME5}, can be written as:
\[\text{find $u\in \VME$ s.t.:}  \qquad  \forall v\in \VME,\qquad  \aME(u,v)+\bME(u,v) = \langle \fsource, \bar{v}\rangle _{H^{-1} ,H^{1} }.\]
The next lemma shows that the bilinear forms $\aME$ are uniformly coercive w.r.t.\ $\varepsilon $. We recall the notation $\snhg{u} :=\|\partial _{x} u\|_{L^2 (\GME)}$.
\end{textesansboite}

\begin{lemmesansboite}{}
\label{coercivite partie hors DtN des ME}%
There exists $C$ that is independent of $\varepsilon $ s.t.:
\begin{itemize}
	\item $\forall u\in \VME3, \quad  \Im( \overline{\omega }) \,\Im( \overline{\omega }\,\aME3(u,u))\geqslant  C\, \|u\|_{H^{1} (\OME)}^2 $.
	\item $\forall u\in \VME5, \quad  \Im( \overline{\omega }) \,\Im( \overline{\omega }\,\aME5(u,u))\geqslant  C\, \big(\|u\|_{H^{1} (\OME)}^2  + \varepsilon  \.\snhg{u}^2 \big)$.
\end{itemize}
\end{lemmesansboite}

\begin{demo}
In both cases, $\aME(u,u)$ has the form $x-\omega ^2 y$ with $x,y\in \mathbb{R}_+$ that depend on $u$. It suffices to notice that we have $\Im( \overline{\omega }(x-\omega ^2 y))=\Im( \overline{\omega }x-\omega |\omega |^2  y)=\Im( \overline{\omega })(x+|\omega |^2 y)$ with $|\omega |^2 >0$.
\end{demo}

\begin{remarque}
$\big(\|u\|_{H^{1} (\OME)}^2  + \varepsilon \.\snhg{u}^2 \big)^{1/2}$ is a natural norm on $\VME5$ because it makes it complete, and the factor $\varepsilon $ comes from the variational formulation \cref{eq: FV ME5}. Besides, the multiplication by $\Im( \overline{\omega })$ 
correctly flips the inequality if $\Im(\omega )$ does~not have the right sign.
\end{remarque}

\begin{textesansboite}
Hence, to prove that $\aME+\bME$ is coercive it suffices to prove that $\bME$ is positive (in the sense $\Im( \overline{\omega }) \,\Im( \overline{\omega }\,\bME(u,u))\geqslant 0$) up to some perturbation error that is controlled by the coercivity of $\aME$. This work is done in \cref{sec: positivite de opS}. \Cref{caractere bien pose des ME} is finally proven in \cref{sec: ME bien poses preuve}.
\end{textesansboite}

\subsection{Asymptotic behavior of the DtN operators for a small $\varepsilon $}
\label{sec: positivite de opS}

\begin{textesansboite}
Let us recall that, according to \cref{eq: mme qd eps=0}, for~any $n\in [\![1,N]\!]$, the $n$-th mode of the DtN for $\varepsilon =0$ is equal to
\[\mme0_{n} \.=\. \mse0_{d} \.=\. \mcle_{d} \.=\. \mcll_{d,0} \.=\. (-1)^{n} \,\Gamma (d+1) \left( \frac{2}{k_{0}} \right)^{\!d} J_{d} (k_{0} r) \sin (d\theta ) \qquad  \text{with }  d:=n\txtpst.\]
Note that the functions $J_{d} (k_{0} r)\sin (d \theta )$ are the modes of Helmholtz's problem in a Dirichlet corner \emph{without} layer. Hence $\opSref$ is a DtN operator for this Helmholtz's problem without layer (and it actually does~not depend on $m$). This is a key argument for the proof of stability. We prove in \cref{opS positif qd eps=0} that the operator $\opSref$ adds a positive contribution to the variational formulation. Then \cref{opS - opS0} shows that $\|\opS-\opSref\|_{\mathcal{L}(L^2 (\SME))} =O(\varepsilon )$ when $\varepsilon \rightarrow 0$. We start with a preliminary result on Bessel functions.
\end{textesansboite}

\begin{lemmesansboite}{}
\label{minoration majoration BesselJ}%
For~any $R>0$, there exists $C_{R} >0$ s.t.\ we have
\[\forall d\in \txtTN,\; r\in (0,R], \qquad  (C_{R} )^{-1} \,r^{d} \leqslant  \bigg|\Gamma (d+1) \left( \frac{2}{k_{0}} \right)^{d} J_{d} (k_{0} r)\bigg| \leqslant  C_{R} \,r^{d} .\]
\end{lemmesansboite}

\begin{demo}
To begin, let us show that, for~any $d\in \TN$ and $r>0$, we have $J_{d} (k_{0} r) \neq 0$. It is a known property of Bessel functions that all the zeros of $J_{a}$ are real, but it~is easy and interesting to prove in this context.
\begin{sousdemo}
	Let us assume by contradiction that there exist $d\in \TN$ and $r_{1} >0$ s.t.\ $J_{d} (k_{0} r_{1} ) =0$. Let $f:(r,\theta ) \mapsto  J_{d} (k_{0} r)\sin (d\theta )$ be defined on $D:=\Omega  \cap B(0,r_{1} )$. By~assumption, we have $f_{|\partial D} =0$. By~properties of the Bessel functions, we also have $\Delta f-k_{0} ^2 f=0$ (see \cite{Leb72} and the proof of \cref{relations sur BesselGamma/Omega}). Moreover, we have $f\in H^{1} (D)$ because $J_{d} (k_{0} r)=O(r^{d} )$ when $r\rightarrow 0$ by \cite{Leb72}. Now, Helmholtz's problem with Dirichlet condition is well-posed in $H^{1} _{0} (D)$ because $\Im(\omega ) \neq 0$ (see e.g.\ \cref{pb de Helmholtz variationnel}). Thus $f=0$ in $D$, which is absurd.
\end{sousdemo}
Next, for~any $d\in \TN$, let $f_{d} :r\in \mathbb{R}_+^*\mapsto \big|\Gamma (d+1) \big( \frac{2}{k_{0} \red{r}} \big)^{d} J_{d} (k_{0} r)\big|$. We need to prove $\smash{\inf\limits _{d\in \TN\!,\,r\in (0,R]} f_{d} }>0$. \cite[p.\,44]{Wat95} states:
\begin{equation}{\label{estimation de Wat95 pour ME}}
	\forall q\in \mathbb{R}_+,\ \forall z\in \mathbb{C}^*, \qquad  \left| \Gamma (q+1) \left( \frac{2}{z} \right)^{q} J_{q} (z) - 1 \right| \leqslant  \exp \!\left( \frac{|z|^2 }{4(q+1)} \right) -1.
\end{equation}
Let $d_{0} \in \TN$ be s.t.\ $\exp \!\big( \frac{R^2 }{4(d_{0} +1)}\big)-1\leqslant  \frac{1}{2}$. Then we have $ \frac{1}{2} \leqslant  f_{d} (r) \leqslant  \frac{3}{2}$ for~any $d\geqslant d_{0}$ and $r\in (0,R]$. Moreover, for any $d\in \TN$ s.t.\ $d<d_{0}$, we have $f_{d} (r)\rightarrow 1$ when $r\rightarrow 0^{+}$ by \cite{Leb72}, and $f_{d}$ has no zero in $(0,R]$ by the beginning of this proof. Hence $f_{d}$ can be extended into a positive continuous function on $[0,R]$. By~compacity, we deduce $\min _{d\in \TN\!,d<d_{0} } \inf _{(0,R]} f_{d} >0$ and $\max _{d\in \TN\!,d<d_{0} } \sup _{(0,R]} f_{d} <\infty $ which concludes the proof.
\end{demo}

\begin{proposition}{}
\label{opS positif qd eps=0}%
$\forall v\in L^2 (\SME),\ \Im( \overline{\omega })\,\Im\big( \overline{\omega }\,\langle \opSref(v),v\rangle _{L^2 (\SME)} \big) \geqslant  0.$
\end{proposition}

\begin{demo}
By~\cref{eq: mme qd eps=0}, each $\mme0_{n}$ is equal to $\mcll_{n\smallpst,0}$ and is proportional on $\SME$ to $J_{d} (k_{0} \RME) \sin (d\theta )$ with $d:=n\pst$. The functions $\theta \mapsto \sin (d\theta )$ are orthogonal between each other and the $J_{d} (k_{0} \RME)$ do~not vanish by \cref{minoration majoration BesselJ}, so $\gram0$ is a diagonal inversible matrix. Hence we have: $\forall v\in L^2 (\SME),\;\forall \bfx \in \SME,$
\begin{align}
	\opSref(v)(\bfx) 
	&= \sum _{n=1} ^{N} \langle v,\mcll_{n\pst,0} \rangle  \cdot  \|\mcll_{n\pst,0} \|_{L^2 (\SME)} ^{-2} \cdot  \partial _{r} \mcll_{n\pst,0} (\bfx) \nonumber\\
	&= \sum _{n=1} ^{N} \langle v,\varphi _{n} \rangle  \, \partial _{r} \varphi _{n} (\bfx) \quad\qquad  \text{with }  \varphi _{n} := \frac{\mcll_{n\smallpst,0} }{\|\mcll_{n\smallpst,0} \|_{L^2 (\SME)} }
	\label{eq: opS positif qd eps=0 preuve expression opSref}
\end{align}
Let $w:\Omega \rightarrow \mathbb{C}$ be defined by $w := \sum _{n=1} ^{N} \langle v,\varphi _{n} \rangle  \, \varphi _{n} $. We have shown that $\opSref(v) =\partial _{r} w_{|\SME}$. Let $V :=\operatorname{Span}(\varphi _{1|\SME} ,\dots ,\varphi _{N|\SME} ) \subset L^2 (\SME)$. Using again \cref{eq: mme qd eps=0}, we know that $\varphi _{n|\SME}$ and $\partial _{r} \varphi _{n|\SME}$ are both proportional to $\sin (n\pst \theta )$, so they are colinear in $L^2 (\SME)$. This yields $\partial _{r} w\in V$. Moreover, $(\varphi _{n|\SME} )_{n}$ is an orthonormal family, so $w_{|\SME}$ is the $L^2 $ orthogonal projection of $v$ on $V$. We deduce $\langle \partial _{r} w, v-w\rangle _{L^2 (\SME)} =0$. Therefore, we obtain:
\[\langle \opSref(v),v\rangle _{L^2 (\SME)} = \langle \partial _{r} w, v\rangle _{L^2 (\SME)} = \langle \partial _{r} w, w\rangle _{L^2 (\SME)} = \int _{\Omega \cap B(0,\RME)} \big(|\nabla  w|^2 -k_{0} ^2 |w|^2 \big)\]
because \cref{convergence des mcl} implies that $w$ satisfies $\Delta w+k_{0} ^2 w=0$ in $\Omega $ and $w=0$ on $\partial \Omega $. We conclude as in the proof of \cref{coercivite partie hors DtN des ME}, by noticing that $\langle \opSref(v),v\rangle _{L^2 (\SME)}$ has the form $x-\omega ^2 y$ with $x,y\in \mathbb{R}_+$ and that $\Im( \overline{\omega })\,\Im( \overline{\omega }(x-\omega ^2 y))=\Im( \overline{\omega })^2 (x+|\omega |^2 y)\geqslant 0$.
\end{demo}

\begin{textesansboite}
In order to compare $\opS$ to $\opSref$, we need some estimates on the modes $\mme_{n}$, that are given by the two next lemmas.
\end{textesansboite}

\begin{lemmesansboite}{}
\label{estimation mcld0}%
Let $R>0$, $n\in [\![1,N]\!]$ and $j\in \{0,1,2\}$. If $\RME\leqslant R$, then we have
\[\|\mcle_{n\pst,0} \|_{L^{\infty } (\SME)} \leqslant  C_{R} \,(\RME)^{n\pst} \quad \  \text{ and } \quad \ \|\partial _{\theta } ^{j} \partial _{r} \mcle_{n\pst,0} \|_{L^{\infty } (\SME)} \leqslant  C_{R} \,N^{j+1} \.(\RME)^{n\pst-1} .\]
\end{lemmesansboite}

\begin{demo}
Let $d:=n\pst$. \Cref{def: mcle,expression sing homog CL} shows that we have $\mcle_{d,0} = (-1)^{n} \,\Gamma (d+1) \.\big( \frac{2}{k_{0}} \big)^{d} J_{d} (k_{0} r) \sin (d\theta )$. With \cref{minoration majoration BesselJ}, this directly implies $\|\mcle_{d,0} \|_{L^2 (\SME)} \leqslant C_{R} \,(\RME)^{n\pst}$. Moreover, we have $2 J_{d} '(z) = J_{d-1} (z)-J_{d+1} (z)$ by \cite{Leb72}. Using again \cref{minoration majoration BesselJ}, we get $|\Gamma (d+1) \big( \frac{2}{k_{0}} \big)^{d} J_{d} '(k_{0} r)| \leqslant  C_{R} \,d\.r^{d-1}$. Since $\partial _{r} \mcle_{d,0} = (-1)^{n} \,\Gamma (d+1) \big( \frac{2}{k_{0}} \big)^{d} k_{0} J_{d} '(k_{0} r) \sin (d\theta )$, we deduce $\|\partial _{\theta } ^{j} \partial _{r} \mcle_{d,0} \|_{L^{\infty } (\SME)} \leqslant  C_{R} \,d^{j+1} \*\.(\RME)^{n\pst-1}$.
\end{demo}

\begin{lemmesansboite}{}
\label{estimation mme}%
Let $R>0$, $M>0$ and $j\in \{0,1,2\}$. If $\RME\leqslant R$ and $\varepsilon  \,\domi{\ln \RME}^{3} \leqslant M\RME$, then
\setalignbox{$\|\partial _{\theta } ^{j} \partial _{r} (\mme_{n} -\mme0_{n} )\|_{L^{\infty } (\SME)}$}
\begin{itemize}
	\item $\alignbox{\|\mme_{n} \|_{L^{\infty } (\SME)}} \leqslant  C_{R,N,M} \,(\RME)^{n\pst}$,
	\item $\alignbox{\|\partial _{\theta } ^{j} \partial _{r} \mme_{n} \|_{L^{\infty } (\SME)}} \leqslant  C_{R,N,M} \,(\RME)^{n\pst-1}$,
	\item $\alignbox{\|\mme_{n} -\mme0_{n} \|_{L^{\infty } (\SME)}} \leqslant  C_{R,N,M} \,\varepsilon \,(\RME)^{n\pst-1} \,\domi{\ln \RME}^{3}$,
	\item $\alignbox{\|\partial _{\theta } ^{j} \partial _{r} (\mme_{n} -\mme0_{n} ) \|_{L^{\infty } (\SME)}} \leqslant  C_{R,N,M} \,\varepsilon \,(\RME)^{n\pst-2} \,\domi{\ln \RME}^{3}$.
\end{itemize}
\end{lemmesansboite}

\begin{demo}
We recall that, according to \cref{estimation des mcl}, we have for any $\delta \in \mathbb{R}$, $d\in \TZ*$ s.t.\ $d\geqslant \delta $, $n\in \mathbb{N}$, $\eAMP>1$, $r\in (0,R]$ and $\theta \in [0,\Theta ]$:
\begin{equation}{\label{eq: estimation mme preuve 1}}
	|\mcle_{d,n} (r,\theta )| + r|\nabla  \mcle_{d,n} (r,\theta )| \leqslant  C_{\delta ,n,R,\eAMP} \,(\eAMP r)^{d-n} \, \domi{\ln r}^{2n+1} .
\end{equation}
It is not difficult to extend this result to higher order derivatives and show:
\begin{equation}{\label{eq: estimation mme preuve 1 generalise}}
	\forall (i,j)\in \mathbb{N}^2 ,\qquad  r^{i+j} |\partial _{x} ^{i} \partial _{y} ^{j} \mcle_{d,n} (r,\theta )| \leqslant  C_{\delta ,n,R,\eAMP,i,j} \,(\eAMP r)^{d-n} \, \domi{\ln r}^{2n+1} .
\end{equation}
This requires to modify \cref{majoration des fonctions de A par leur norme,cvg des SF de norme finie,convergence des mcl} to control the derivatives of high order in the same way we controlled the first derivative. Each derivation introduces an additional factor $1/r$ in the bound, which explains the factors $r$ in \cref{eq: estimation mme preuve 1} and $r^{i+j}$ in \cref{eq: estimation mme preuve 1 generalise}. Next, for any $f\in H^{3} _{\mathrm{loc}} (\Omega )$ and $(r,\theta )$, let us denote
\[|f|_{r,\theta } := |f(r,\theta )|+r|\partial _{r} f(r,\theta )|+r|\partial _{\theta } \partial _{r} f(r,\theta )|+r|\partial _{\theta } ^2  \partial _{r} f(r,\theta )|.\]
Taking $\delta :=-m$, $\eAMP :=2$, $d\leqslant N\pst$ and $n\in [\![0,m]\!]$, we obtain in particular $|\mcle_{d,n} |_{r,\theta } \leqslant  C_{R,N} \, r^{d-n} \, \domi{\ln r}^{2n+1}$. By~\cref{eq: def mse}, this yields for $r:=\RME$:
\begin{equation}{\label{eq: estimation mse}}
	\big|\mse_{d} \big|_{r,\theta } \leqslant  \sum _{n=0} ^{m-1} \varepsilon ^{n} \big|\mcle_{d,n} \big|_{r,\theta } 
	\leqslant  \sum _{n=0} ^{m-1} C_{R,N} \, r^{d} \left( \frac{\varepsilon \,\domi{\ln r}^2 }{r} \right)^{n} \domi{\ln r}
	\leqslant  C_{R,N,M} \, r^{d} \,\domi{\ln r}
\end{equation}
because $\varepsilon \,\domi{\ln r}^2  \leqslant  Mr$ by assumption. Similarly, we have 
\[\big|\mse_{d} -\mcle_{d,0} \big|_{r,\theta } \leqslant  \sum _{n=\red{1}} ^{m-1} C_{R,N} \, r^{d} \left( \frac{\varepsilon \,\domi{\ln r}^2 }{r} \right)^{n} \domi{\ln r} \leqslant  C_{R,N,M} \, \varepsilon \, r^{d-1} \,\domi{\ln r}^{3} .\]
Now, in order to bound the modes $\mme_{n}$, which involve the coefficients $\cme_{d,d',\varepsilon }$, let us show: ${\forall d\in -\TN},\allowbreak\, \forall d'\in \TN,\; |\cme_{d,d',\varepsilon } | \leqslant  C_{R,M,d,d'} \,\varepsilon ^{d'-d} \.\domi{\ln \varepsilon }^{3(d'-d-1)}$.
\begin{sousdemo}
	\def\q{q}%
	Let $\q:= d'-d$. By~\cref{def: mme,coefs cme non nuls}, we have
	\[\cme_{d,d',\varepsilon } := \sum _{p\in \NN \cap [0,m)} \, \sum _{\ell =0} ^{\ell _{p}} \varepsilon ^{p} \ln ^{\ell } \!\varepsilon  \cdot  \cme_{d,d',p,\ell } = \sum _{p\in \NN \cap \red{[\q,m)}} \, \sum _{\ell =0} ^{\ell _{p}} \varepsilon ^{p} \ln ^{\ell } \!\varepsilon  \cdot  \cme_{d,d',p,\ell } .\]
	If $\q\geqslant m$, this yields $\cme_{d,d',\varepsilon } =0$. Let us assume $\q<m$. Since $\varepsilon \leqslant Mr\leqslant MR$, we deduce $|\cme_{d,d',\varepsilon } |\leqslant  C_{R,M,d,d'} \,\varepsilon ^{\q} \.\domi{\ln \varepsilon }^{\ell _{\q}} $. To conclude, it suffices to prove $\ell _{\q} \leqslant 3(\q-1)$.
	\begin{itemize}
		\item If $\Theta =\pi $, then we have $\q \in \TZ \cap [0,5)=[\![0,4]\!]$ so \cref{tab: liste CL pour pi} shows $\ell _{\q} \leqslant  \q -1$. In addition, we have $\q\geqslant 2\pst\geqslant 1$, so we deduce $\ell _{\q} \leqslant  3(\q-1)$.
		\item Otherwise, we have $\ell _{\q} \leqslant \ell _{m} \leqslant 1$ by \cref{rq: degre en ln eps dans ME} and, if $\q<2$, then $\ell _{\q} =0$ by \cref{simplifications cas irrationnel,simplifications cas rationnel}. Since $q\geqslant 2\pst\geqslant 1$, this implies $\ell _{\q} \leqslant 3(\q-1)$.
	\end{itemize}
\end{sousdemo}
Let us also show that we have $\varepsilon \.\domi{\ln \varepsilon }^{3} \leqslant  C_{R,M} \.r$ (still with $r:=\RME$).
\begin{itemize}
	\item If $\varepsilon \leqslant r^2 $, then we have $\varepsilon \.\domi{\ln \varepsilon }^{3} \leqslant  C_{R,M} \sqrt{\varepsilon } \leqslant C_{R,M} \.r$ because $\varepsilon \leqslant Mr\leqslant M R$.
	\item Otherwise, we have $\varepsilon \.\domi{\ln \varepsilon }^{3} \leqslant  \varepsilon  \cdot C \domi{\ln r}^{3} \leqslant C_{M} \.r$. 
\end{itemize}
Recalling that $\mme0_{n} =\mcle_{n\pst,0}$ by~\cref{eq: mme qd eps=0}, and using \cref{def: mme}, we deduce
\begin{align*}
	\big|\mme_{n} -\mme0_{n} \big|_{r,\theta } &= \bigg|\mse_{n\pst} -\mcle_{n\pst,0} + \sum _{d\in \TZ*\cap (-m,0)} \cme_{d,n\pst,\varepsilon } \cdot  \mse_{d} \bigg|_{r,\theta } \\
	&\leqslant  C_{R,N,M} \bigg( \varepsilon \,r^{n\pst-1} \langle \ln r\rangle ^{3} + \sum _{d\in \TZ*\cap (-m,0)} \varepsilon ^{n\pst-d} \,\domi{\ln \varepsilon }^{3(n\pst-d-1)} \cdot  r^{d} \,\domi{\ln r} \bigg) \\
	&\leqslant  C_{R,N,M} \,\varepsilon \,r^{n\pst-1} \domi{\ln r}^{3} \bigg( 1 + \sum _{d\in \TZ*\cap (-m,0)} \left( \frac{\varepsilon \.\domi{\ln \varepsilon }^{3} }{r} \right)^{n\pst-d-1} \bigg) \\
	&\leqslant  C_{R,N,M} \,\varepsilon \, r^{n\pst-1} \domi{\ln r}^{3}
\end{align*}
using $n\pst-d-1\geqslant 2\pst-1\geqslant 0$ and $\varepsilon \.\domi{\ln \varepsilon }^{3} \leqslant  C_{R,M} \.r$. Moreover \cref{estimation mcld0} shows that $|\mcle_{n\smallpst,0} |_{r,\theta } \leqslant C_{R,N} \,r^{n\smallpst}$. Compared to the previous bound on $|\mcle_{d,n} |_{r,\theta }$, this one does~not involve $\ln r$, which will slightly improve the powers of $\ln \RME$ (although this is not essential). Thanks to the hypothesis $\varepsilon \,\domi{\ln r}^{3} \leqslant  Mr$, this yields $\big|\mme_{n} \big|_{r,\theta } \leqslant \big|\mme_{n} -\mme0_{n} \big|_{r,\theta } +\big|\mcle_{n\pst,0} \big|_{r,\theta } \leqslant  C_{R,N,M} \,r^{n\pst} $.
\end{demo}

\begin{textesansboite}
In order to describe the behavior w.r.t.\ $\RME$ of the stability condition, and later of the error estimates, we use scaling techniques. They involve some variants of the Sobolev spaces that are introduced in the definition below. We denote them $H_{\mathrm{sc}} ^{s}$, where ``sc'' stands for ``scaled''. We recall that the norm on $L^2 (\SME)$ is defined by \cref{eq: def produit scalaire L2(SME)}.
\end{textesansboite}

\begin{definitionsansboite}{}
\label{def: sobolev pondere pour ME}%
Let $D\in \{\GME,\SME\}$ and $s\in \mathbb{R}$. We define $\HME{s}(D):=H^{s} (D)$ and equip it with the following norm. For~any $f\in H^{s} (D)$, let us denote $D_{1} :=\RME^{-1} D$ (which is independent of $\RME$) and ${f_{1} :x\in D_{1} \mapsto {}} f(\RME x)$. We define $\|f\|_{\HME{s}(D)} := \sqrt{\RME} \,\|f_{1} \|_{H^{s} (D_{1} )}$.
\end{definitionsansboite}
\nomenclature[17-28-12]{$\HME{s}$}{In \cref{chap: ME}, scaled variant of the Sobolev space $H^{s}$}

\begin{exemplenumerote}
\label{ex: sobolev pondere pour ME}%
\begin{itemize}
	\item $\|f\|_{\HME0(D)} = \|f\|_{L^2 (D)}$.
	\item $\|f\|_{\HME1(\GME)}$ is equivalent to $\|f\|_{L^2 (\GME)} + \RME \|\partial _{x} f\|_{L^2 (\GME)}$.
	\item For~any $j\in \mathbb{N}$, we have $\|f\|_{\HME{j}(\SME)} = \big(\sum _{i=0} ^{j} \|\partial _{\theta } ^{i} f\|_{L^2 (\SME)} ^2  \big)^{1/2}$. For comparison, the classical $H^{1}$~norm is defined by $\|f\|_{H^{j} (\SME)} = \big( \sum _{i=0} ^{j} \|\RME^{-i} \partial _{\theta } ^{i} f\|_{L^2 (\SME)} ^2  \big)^{1/2}$.
	\item However, $\|f\|_{\HME{1/2}(D)}$ cannot be expressed directly with $\|f\|_{H^{1/2} (D)}$.
\end{itemize}
\end{exemplenumerote}

\begin{textesansboite}
The next result rigorously proves $\opS-\opSref =O(\varepsilon )$ when $\varepsilon \rightarrow 0$. $\opS-\opSref$ is estimated for the $\HME2(\SME)$ norm. For the stability analysis, only the $L^2 $ estimate will be used, but the bound in $\HME2(\SME)$ will be useful in the error estimates, see \cref{regularite elliptique pour ME 2}.
\end{textesansboite}

\begin{proposition}{estimates on $\opS-\opSref$ and $\opP$ for a small $\varepsilon $}
\label{opS - opS0}%
Let $R>0$ and $N\in \mathbb{N}$. There exist two positive constants $\CdtnME[R]1$ and $\CdtnME[R]2$ s.t., for~any $\RME\in (0,R]$ and $0\leqslant \varepsilon \leqslant \CdtnME[R]1 \,\RME \,\domi{\ln \RME}^{-3}$, we have:
\begin{enumerate}
	\item $\gram$ is invertible (which means that $\opS$ and $\opP$ are well defined),
	\item $\displaystyle \|\opS-\opSref\|_{\mathcal{L}(L^2 (\SME),\HME2(\SME))} + \sqrt{\RME} \,\|\opP\|_{\mathcal{L}(L^2 (\SME),\mathbb{C})} \leqslant  \CdtnME[R]2 \. \frac{\varepsilon }{\RME^{2}} \,\domi{\ln \RME}^{3}$.
\end{enumerate}\fbi[0.8]
\end{proposition}

\begin{demo}
Let us introduce the following normalized variables:
\begin{itemize}
	\item $\forall n\in [\![1,N]\!],\ a_{n} :=\sqrt{\RME\Theta /2} \cdot  {\Gamma (n\pst+1)}\* \big( \frac{2}{k_{0}} \big)^{n\smallpst} J_{n\smallpst} (k_{0} \RME)$, which is non-zero by \cref{minoration majoration BesselJ},
	\item $\forall n\in [\![1,N]\!],\ \mmen_{n} :=\mme_{n} /a_{n}$,
	\item $\gramn :=\big[\langle \mmen_{n} ,\mmen_{n'} \rangle _{L^2 (\SME)} \big]_{(n,n')\in [\![1,N]\!]^2 }$,
	\item $D = \mathrm{Diag} \big[a_{1} ,\dots ,a_{N} \big]\in \mathrm{GL}_{N} (\mathbb{C})$.
\end{itemize}
$\gram$ is invertible iff $\gramn = D^{-1} \gram D^{-1}$ is invertible, and in this case we have: $\forall v\in L^2 (\SME),\;\forall \bfx \in \SME,$
\begin{align}
	\opS(v)(\bfx) &= \Big[\langle v ,\mme_{n} \rangle _{L^2 (\SME)} \Big]_{n} ^{\mathsf{T}} D^{-1} \cdot (D^{-1} \gram D^{-1} )^{-1} \cdot  D^{-1} \Big[\partial _{r} \mme_{n} (\bfx)\Big]_{n} \nonumber\\[.5ex]
	&= \Big[\langle v ,\mmen_{n} \rangle _{L^2 (\SME)} \Big]_{n} ^{\mathsf{T}} \cdot  \gramn^{-1} \cdot \Big[\partial _{r} \mmen_{n} (\bfx)\Big]_{n} . \label{eq: opS - opS0 preuve expression opS normalisee}
\end{align}
Moreover \cref{eq: mme qd eps=0} shows that we have for~any $n\in [\![1,N]\!]$: $\mme0_{n|\SME} = \mcll_{n\pst,0|\SME} = (-1)^{n} \,a_{n} \sqrt{ \frac{2}{\RME\Theta } } \*\sin (n\pst \theta )$. Since the functions $\sqrt{ \frac{2}{\RME\Theta } } \sin (n\pst \theta )$ are orthonormal in $L^2 (\SME)$, we deduce
\begin{equation}{\label{eq: opS - opS0 preuve gram normalisee eps=0}}
	\gramn0 = \bigg[\Big\langle  \frac{1}{a_{n}} \mcll_{n\pst,0} ,\, \frac{1}{a_{n'}} \mcll_{n'\pst,0} \Big\rangle _{\!L^2 (\SME)} \bigg]_{n,n'} = I_{N} .
\end{equation}

\underline{Step 1:} We will prove that $\gramn$ is invertible for small values of $\varepsilon $ using the Neumann series ``$\gramn^{-1} =(I_{N} -(\gramn0-\gramn))^{-1} =\sum_{j=0}^\infty  (\gramn0-\gramn)^{j}$''. To begin, \cref{minoration majoration BesselJ} implies $|a_{n} |\geqslant C_{R} \sqrt{\RME} \,\RME^{n\pst}$. Therefore, \cref{estimation mme} used with $M:=1$ implies that, for any $j\in \{0,1,2\}$, if $\varepsilon \,\domi{\ln \RME}^{3} \leqslant \RME$ then:
\setalignbox{$\|\partial _{\theta ^{j}} \partial _{r} (\mmen_{n} -\mmen0_{n} )\|_{L^{\infty } (\SME)}$}
\begin{itemize}
	\item $\alignbox{\|\mmen_{n} \|_{L^{\infty } (\SME)}} \leqslant  C_{R,N} \,(\RME)^{-1/2}$,
	\item $\alignbox{\|\partial _{\theta } ^{j} \partial _{r} \mmen_{n} \|_{L^{\infty } (\SME)}} \leqslant  C_{R,N} \,(\RME)^{-3/2}$,
	\item $\alignbox{\|\mmen_{n} -\mmen0_{n} \|_{L^{\infty } (\SME)}} \leqslant  C_{R,N} \,\varepsilon \,(\RME)^{-3/2} \,\domi{\ln \RME}^{3}$,
	\item $\alignbox{\|\partial _{\theta } ^{j} \partial _{r} (\mmen_{n} -\mmen0_{n} )\|_{L^{\infty } (\SME)}} \leqslant  C_{R,N} \,\varepsilon \,(\RME)^{-5/2} \,\domi{\ln \RME}^{3}$.
\end{itemize}
Besides, we have
\begin{align*}
	\langle \mmen_{n} ,\mmen_{n'} \rangle _{L^2 (\SME)} &= \langle \mmen0_{n} ,\mmen_{n'} \rangle  + \langle \mmen_{n} -\mmen0_{n} ,\,\mmen_{n'} \rangle  \\
	&= \langle \mmen0_{n} ,\mmen0_{n'} \rangle  + \langle \mmen0_{n} ,\,\mmen_{n'} -\mmen0_{n'} \rangle  + \langle \mmen_{n} -\mmen0_{n} ,\,\mme_{n'} \rangle  \vspace{-1ex}
\end{align*}
Hence
\begin{align*}
	|\langle \mmen_{n} ,\mmen_{n'} \rangle _{L^2 } -\langle \mmen0_{n} ,\mmen0_{n'} \rangle _{L^2 } | &\leqslant  \RME \, \|\mmen0_{n} \|_{L^{\infty } (\SME)} \|\mmen_{n'} -\mmen0_{n'} \|_{L^{\infty } (\SME)} \\
	& \qquad  {}+ \RME \, \|\mmen_{n} -\mmen0_{n} \|_{L^{\infty } (\SME)} \|\mmen_{n'} \|_{L^{\infty } (\SME)} \\
	&\leqslant  C_{R,N,M} \,\varepsilon \, \RME^{-1} \,\domi{\ln \RME}^{3}
\end{align*}
Next, we introduce the following norms for any $(z_{n} )\in \mathbb{C}^{N}$ and $A\in \mathcal{M}_{N} (\mathbb{C})$:
\[\|(z_{n} )\|_{\infty } :=\max (|z_{1} |,\dots ,|z_{N} |) \qquad  \text{ and } \qquad  \|A\|_{\infty } := \sup _{\mathbf z\in \mathbb{C}^{N} \setminus \{0\}} \frac{\|A\bf z\|_{\infty } }{\|\bf z\|_{\infty } }.\]
For~any $A=\big[a_{i,j} \big]\in \mathcal{M}_{N} (\mathbb{C})$, we have $\|A\|_{\infty } =\max _{i\in [\![1,N]\!]} \sum _{j=1} ^{N} |a_{i,j} |$. Hence we have $\|\gramn-\gramn0\|_{\infty } \leqslant  C^{\mathrm{gram}} _{R,N} \,\varepsilon \, \RME^{-1} \,\domi{\ln \RME}^{3}$ for some constant $C^{\mathrm{gram}} _{R,N} >0$, provided that $\varepsilon \,\domi{\ln \RME}^{3} \leqslant \RME$. Now let us define $\CdtnME1 := \min (1, \frac{1}{2} (C^{\mathrm{gram}} _{R,N} )^{-1} )$. Then for any $\varepsilon  \leqslant  \CdtnME1 \,\RME \domi{\ln \RME}^{-3}$ we have $\|\gramn-\gramn0\|_{\infty } \leqslant  \frac{1}{2}$, which implies that $\gramn$ is invertible and
\begin{itemize}
	\item $\displaystyle \|\gramn^{-1} \|_{\infty } \leqslant  \sum\limits_{j=0}^\infty  \|\gramn-\gramn0\|_{\infty } ^{j} \leqslant  2$\vspace{-0.5ex}
	\item $\displaystyle \|\gramn^{-1} -\gramn0^{-1} \|_{\infty } \leqslant  \sum\limits_{j=1}^\infty  \|\gramn-\gramn0\|_{\infty } ^{j} \leqslant  2\,C^{\mathrm{gram}} _{R,N,M} \, \frac{\varepsilon }{\RME} \,\domi{\ln \RME}^{4} = C_{R,N} \, \frac{\varepsilon }{\RME} \,\domi{\ln \RME}^{3}$.
\end{itemize}
\medskip\smallskip

\underline{Step 2:} It remains to estimate the operator norms of $\opS$, $\opS-\opSref$ and $\opP$. Let us assume $\varepsilon  \leqslant  \CdtnME[R]1 \,\RME \*\domi{\ln \RME}^{-3}$. Let $v\in L^2 (\SME)$ and $\bfx \in \SME$. We define $V_{\varepsilon } (v):=\big[\langle v,\mmen_{n} \rangle \big]_{n} \in \mathbb{C}^{N}$ and $V'_{\varepsilon } (\bfx):=\big[\partial _{r} \mmen_{n} (\bfx)\big]_{n} \in \mathbb{C}^{N}$. Step 1 provides the following estimates (taking $M:=1$ and $j\in \{0,1,2\}$):
\setalignbox{$\|\partial _{\theta } ^{j} (V'_{\varepsilon } -V'_{0} )(\bfx)\|_{\infty }$}
\begin{itemize}
	\item $\alignbox{\|V_{\varepsilon } (v)\|_{\infty } } \leqslant  C_{R,N} \,\|v\|_{L^2 (\SME)}$
	\item $\alignbox{\|V_{\varepsilon } (v)-V_{0} (v)\|_{\infty } } \leqslant  C_{R,N} \,\varepsilon \,\RME^{-1} \,\domi{\ln \RME}^{3} \,\|v\|_{L^2 (\SME)}$
	\item $\alignbox{\|\partial _{\theta } ^{j} V'_{\varepsilon } (\bfx)\|_{\infty } } \leqslant  C_{R,N} \,\RME^{-3/2} $
	\item $\alignbox{\|\partial _{\theta } ^{j} (V'_{\varepsilon } -V'_{0} )(\bfx)\|_{\infty } } \leqslant  C_{R,N} \,\varepsilon \,\RME^{-5/2} \,\domi{\ln \RME}^{3}$
\end{itemize}
\Cref{eq: opS - opS0 preuve expression opS normalisee} shows that $\opS(v)(\bfx)=V_{\varepsilon } (v)^{\mathsf{T}} \cdot \gramn^{-1} \cdot V'_{\varepsilon } (\bfx)$. Hence we~have
\def\alignbox#1{\begingroup
	\setbox0=\hbox{$\big|(\opS-\opSref)(v)(\bfx)\big|$}
	\dimen0=\linewidth
	\advance\dimen0 by -\wd0
	\advance\dimen0 by -4.5em
	\hbox to \dimen0{${}#1$}\endgroup}
\begin{align*}
	\big|(\opS-\opSref)(v)(\bfx)\big| &\alignbox{= \big|(V_{\varepsilon } -V_{0} )(v)^{\mathsf{T}} \cdot \gramn^{-1} \cdot V'_{\varepsilon } (\bfx) \,+\, V_{0} (v)^{\mathsf{T}} \cdot (\gramn^{-1} -\gramn0^{-1} ) \cdot V'_{\varepsilon } (\bfx) \hfill}\\
	&\alignbox{\hfill {}+ V_{0} (v)^{\mathsf{T}} \cdot \gramn0^{-1} \cdot (V'_{\varepsilon } -V'_{0} )(\bfx)\big|}\\
	&\alignbox{\leqslant  N\, \|V_{\varepsilon } (v)-V_{0} (v)\|_{\infty } \,\|\gramn^{-1} \|_{\infty } \,\|V'_{\varepsilon } (\bfx)\|_{\infty } \hfill}\\
	&\alignbox{\hfill {} + N\, \|V_{0} (v)\|_{\infty } \,\|\gramn^{-1} -\gramn0^{-1} \|_{\infty } \,\|V'_{\varepsilon } (\bfx)\|_{\infty } }\\
	& \alignbox{\hfill {} + N\, \|V_{0} (v)\|_{\infty } \,\|\gramn0^{-1} \|_{\infty } \,\|V'_{\varepsilon } (\bfx)-V'_{0} (\bfx)\|_{\infty } }\\
	&\leqslant  C_{R,N} \,\varepsilon \,\RME^{-5/2} \,\domi{\ln \RME}^{3} \,\|v\|_{L^2 (\SME)}
\end{align*}
$\big|\partial _{\theta } ^{j} (\opS-\opSref)(v)(\bfx)\big|$ can be bounded in the same way. Thus $\big\|(\opS-\opSref)(v)\big\|_{\HME2(\SME)} \leqslant  C_{R,N} \*\varepsilon \,\RME^{-2} \,\domi{\ln \RME}^{3} \,\|v\|_{L^2 (\SME)}$.\\
As for $\opP$, we have $V'_{0} (\PME)=\big[\partial _{r} \mcll_{n\pst,0} (\PME)\big]_{n}$ by \cref{eq: mme qd eps=0}. Since $\mcll_{n\pst,0|\partial \Omega } =0$ by \cref{convergence des mcl}, we get $V'_{0} (\PME) =0$. Therefore:\vspace{-.5ex}
\begin{align*}
	|\opP(v)| &= \big|V_{\varepsilon } (v)^{\mathsf{T}} \cdot \gramn^{-1} \cdot (V'_{\varepsilon } -V'_{0} )(\bfx)\big| \\
	&\leqslant  N\, \|V_{\varepsilon } (v)\|_{\infty } \,\|\gramn^{-1} \|_{\infty } \,\|(V'_{\varepsilon } -V'_{0} )(\bfx)\|_{\infty } \\
	&\leqslant  C_{R,N} \,\varepsilon \,\RME^{-5/2} \,\domi{\ln \RME}^{3} \,\|v\|_{L^2 (\SME)} .
\end{align*}
\end{demo}

\begin{remarque}
The proof of the well-posedness and uniform stability of the effective models gives another motivation to introduce the modal decomposition of $\scl$ in terms of the $\mme_{n}$ (\cref{decompo modale SE 2}), rather than using the first decomposition in terms of the $\mse_{d}$ (\cref{decompo modale SE 1}).
Indeed, in this section we use among others the fact that, for~any $n$, $\mme_{n|\SME}$ is close to $\sin(n\pst \theta )$ times a constant when $\varepsilon $ is small. The orthogonality of these sines contributes to the good properties of $\opSref$ and of $\opS$ when $\varepsilon \rightarrow 0$. However, by \cref{eq: mme qd eps=0}, each $\sin (n\pst \theta )$ is close to \emph{two} modes $\mse_{d}$ at once when $\varepsilon $ is small, namely $\mse_{n\pst}$ and $\mse_{-n\pst}$. So at the limit $\varepsilon \rightarrow 0$, $\mse_{n\pst}$ and $\mse_{-n\pst}$ become indistinguishable. Hence, if we defined $\opS$ using the functions $\mse_{d}$ instead of the $\mme_{n}$, $\opS$ would likely not converge when $\varepsilon \rightarrow 0$ and we would not be able to prove the uniform stability w.r.t.\ $\varepsilon $ of the effective models.
\end{remarque}

\subsection{Proof of well-posedness and stability}
\label{sec: ME bien poses preuve}

\newcommand{\sourceOME}{\tilde{f}_{\rm s}}%
\newcommand{\sourceGME}{g_{\rm s}}%
\newcommand{\sourceSME}{h_{\rm s}}%
\newcommand{\sourcePME}{p_{\rm s}}%
\def\lME#1{\ell^{(#1)}}%
\def\lMEm{\lME{m}}%
\begin{textesansboite}
When proving the error estimate of \cref{cvg des ME}, we will use that the error $\sme-\scl$ satisfies a variant of the effective model with small source terms on $\OME$, $\GME$, $\SME$ and $\PME$. Let us introduce this variant. Let $\sourceOME \in H^{-1} (\OME)$, $\sourceGME \in L^2 (\GME)$, $\sourceSME \in L^2 (\SME)$ and $\sourcePME \in \mathbb{C}$. We define the linear forms
\begin{itemize}
	\item $\forall u\in \VME3, \quad  \lME3(v) :=\displaystyle \langle \sourceOME, v \rangle _{H^{-1} ,H^{1} } + \int _{\GME} \sourceGME v + \int _{\SME} \sourceSME v$
	\item $\forall u\in \VME5, \quad  \lME5(v) :=\displaystyle \langle \sourceOME, v \rangle _{H^{-1} ,H^{1} } + \int _{\GME} \sourceGME v + \int _{\SME} \sourceSME v + \sourcePME \,v(\PME)$.
\end{itemize}
We consider the following problem:
\begin{equation}{\label{eq: FV ME augmente}}
	\text{find $u\in \VME$ s.t.:}  \ \; \forall v\in \VME,\quad  \ \aME(u,v) + \bME(u,v) = \lMEm( \bar{v}) ,
\end{equation}
\phantomsection\label{preuve: caractere bien pose des ME}%
To prove \cref{caractere bien pose des ME}, it suffices to prove the following result that is strictly stronger.
\end{textesansboite}

\begin{theoreme}{}
\label{caractere bien pose des ME augmentes}%
Let $R \in (0,\Rsource)$, $N\in \mathbb{N}$ and $\varepsilon >0$. There exist two positive constants $\CstabME$ and $\CcontME$ s.t., under the stability condition of \cref{def: condition de stabilite}, \cref{eq: FV ME augmente} has a unique solution $v_{\varepsilon } \in \VME$ and $v_{\varepsilon }$ satisfies\vspace{-.5ex}
\[\left\{\begin{array}{@{\,}r@{\;\leqslant \;}l@{\;\,}l}
	\|v_{\varepsilon } \|_{H^{1} (\OME)} & \CcontME\, \big(\|\sourceOME\|_{H^{-1} (\OME)} + \|\sourceGME\|_{L^2 (\GME)} + \sqrt{\RME} \|\sourceSME\|_{L^2 (\SME)} \big) & \text{if } m=3,\\[1ex]
	\|v_{\varepsilon } \|_{H^{1} (\OME)} + \sqrt{\varepsilon } \,\snhg{v_{\varepsilon }} & \CcontME \left( \|\sourceOME\|_{H^{-1} (\OME)} + \|\sourceGME\|_{L^2 (\GME)} + \sqrt{\RME} \|\sourceSME\|_{L^2 (\SME)} +\displaystyle \frac{|\sourcePME|}{\sqrt{\varepsilon }} \right){} & \text{if } m=5.
\end{array}\right.\]\fbeq[0.5]
\end{theoreme}

\begin{demo}
\def\Ccoer{C_{\rm coer}}%
\def\pso{\ell_\omega}
\def\bMEref#1{b^{(#1)}_0}%
\def\cME#1{c^{(#1)}_\peps}%
\def\cMEm{\cME{m}}%
Let $\pso:z\in \mathbb{C}\mapsto \Im( \overline{\omega })\,\Im( \overline{\omega }z)$. We define the following norms
\begin{equation}{\label{eq: norme VME}}
	\left\{\begin{array}{l}
	\forall u\in \VME3, \quad  \|u\|_{\VME3,\varepsilon } := \|u\|_{H^{1} (\OME)} \\
	\forall u\in \VME5, \quad  \|u\|_{\VME5,\varepsilon } := \big( \|u\|_{H^{1} (\OME)}^2  + \varepsilon  \,\snhg{u}^2  \big)^{1/2} 
	\end{array}\right.
\end{equation}
\Cref{coercivite partie hors DtN des ME} states that there exists $\Ccoer>0$ s.t.: $\forall u\in \VME,\; \pso(\aME(u,u))\geqslant  \Ccoer \|u\|_{\VME,\varepsilon }$. We will prove that, for an appropriate choice of $\CstabME$, $\aME+\bME$ is also coercive for $\|.\|_{\VME,\varepsilon }$, with constant $\Ccoer/2$.\\
Let $\cMEm:=\bME -\bMEref{m}$. We recall that, by \cref{opS - opS0}, there exists $\CdtnME1 , \CdtnME2 >0$ s.t., for~any $\RME\in (0,\Rsource)$ and $\varepsilon \leqslant \CdtnME1 \,\RME \,\domi{\ln \RME}^{-3}$, $\opS$ and $\opP$ are well defined and 
\[\|\opS-\opSref\|_{\mathcal{L}(L^2 (\SME))} + \sqrt{\RME} \,\|\opP\|_{\mathcal{L}(L^2 (\SME),\mathbb{C})} \leqslant   \CdtnME2 \, \frac{\varepsilon }{\RME^{2}} \,\domi{\ln \RME}^{3} .\]

\underline{Case $m=3$:} For~any $u\in \VME3$ we have:
\[|\cME3(u,u)| = \big|\big\langle \mu _{0} (\opS3-\opSref3)(u_{|\SME} ),u \big\rangle _{L^2 (\SME)} \big| \leqslant  \mu _{0} \.\CdtnME2 \, \frac{\varepsilon }{\RME^{2}} \,\domi{\ln \RME}^{3} \,\|u\|_{L^2 (\SME)} ^2  .\]
Let us show that we have
\begin{equation}{\label{eq: norme L2 SME vs H1 OME}}
	\|u\|_{L^2 (\SME)} ^2  \leqslant C\RME \|u\|_{H^{1} (\OME)}^2  . \vspace{-1ex}
\end{equation}
\begin{sousdemo}
	Let $\widetilde{\Omega } :=\Omega  \setminus \closedBall(0,1)$ and $\widetilde{\Sigma } := \Omega  \cap \partial B(0,1)$. We have $\OME =\RME \widetilde{\Omega }$ and $\SME=\RME \widetilde{\Sigma }$. Let $\tilde{u} \in H^{1} ( \widetilde{\Omega })$ be defined by $\tilde{u}:x \mapsto u(\RME x)$. It is easy to see that we have $\|u\|_{L^2 (\SME)} ^2 =\RME \| \tilde{u}\|_{L^2 ( \widetilde{\Sigma })} ^2 $ and $\|\nabla  u\|_{L^2 (\OME)} ^2 =\|\nabla  \tilde{u}\|_{L^2 ( \widetilde{\Omega })}^2 $. Finally, since $u_{|\bordgaucheME} =0$, Poincaré's inequality gives $\| \tilde{u}\|_{L^2 ( \widetilde{\Sigma })} \leqslant C\|\nabla  \tilde{u}\|_{L^2 ( \widetilde{\Omega })}$.
\end{sousdemo}
Thus we obtain $|\cME3(u,u)| \leqslant  C_{\Rsource,N} \, \varepsilon  \.\RME^{-1} \.\domi{\ln \RME}^{3} \.\|u\|_{H^{1} (\OME)} ^2 $.\\
Besides, we have $\bMEref3(u,u)=\langle \opSref(u_{|\SME} ), u\rangle _{L^2 (\SME)}$, so \cref{opS positif qd eps=0} shows that $\pso(\bMEref3(u,u))\geqslant  0$. Therefore:
\def\Cneg{C^\star_{\Rsource,N}}
\begin{align*}
	|\omega |^2 \big|(\aME3+\bME3)(u,u)\big| &\geqslant  \pso\big((\aME3+\bME3)(u,u)\big)\\
	&\geqslant  \pso(\aME3(u,u)) + \pso(\bMEref3(u,u)) + \pso(\cME3(u,u)) \\
	&\geqslant  \Ccoer \|u\|_{H^{1} (\OME)} ^2  + 0 - \Cneg \, \frac{\varepsilon }{\RME} \,\domi{\ln \RME}^{3} \,\|u\|_{H^{1} (\OME)} ^2  
\end{align*}
for some constant $\Cneg$. Hence $|\omega |^2 \big|(\aME3+\bME3)(u,u)\big| \geqslant  \frac{\Ccoer}{2} \,\|u\|_{H^{1} (\OME)}$ provided that $\Cneg \, \varepsilon  \,\RME^{-1} \*\,\domi{\ln \RME}^{3} \leqslant  \Ccoer/2$. Let us define
\begin{equation}{\label{eq: caractere bien pose des ME augmentes preuve def Cstab m=3}}
	\CstabME := \min \left(\CdtnME1 , \frac{\Ccoer}{2\,\Cneg }\right).
\end{equation}
We have proven that, when $\varepsilon \leqslant  \CstabME \RME\domi{\ln \RME}^{-3}$ (i.e.\ the stability condition of \cref{def: condition de stabilite}), $\opS$ is well defined and $\aME3+\bME3$ is uniformly coercive w.r.t.\ $(\varepsilon ,\RME,N)$. Moreover, using again \cref{eq: norme L2 SME vs H1 OME}, we have:
\[\forall v\in H^{1} (\OME), \qquad  \lME3(v) \leqslant  C\,\big(\|\sourceOME\|_{H^{-1} (\OME)} + \|\sourceGME\|_{L^2 (\GME)} + \sqrt{\RME} \|\sourceSME\|_{L^2 (\SME)} \big) \,\|v\|_{H^{1} (\OME)} .\]
Thus Lax-Milgram's theorem implies the well-posedness and uniform stability of \cref{eq: FV ME augmente}.\medskip\smallskip

\underline{Case $m=5$:} This is very similar to the previous case, but this time we have
\begin{align}
	|\cME5(u,u)| &= \left|\big\langle \mu _{0} (\opS5-\opSref5)(u_{|\SME} ),u \big\rangle _{L^2 (\SME)} + \frac{\varepsilon \mu _{1} }{3} \,\opP5(u_{|\SME} ) \,\overline{u(\PME)}\right| \nonumber\\
	&\leqslant  C_{\Rsource,N} \,\domi{\ln \RME}^{3} \left( \frac{\varepsilon }{\RME^2 } \, \|u\|_{L^2 (\SME)} ^2  + \frac{\varepsilon ^2 }{\RME^{5/2}} \,\|u\|_{L^2 (\SME)} \|u\|_{H^{1} (\GME)} \right) \nonumber\\
	&\leqslant  C_{\Rsource,N} \,\domi{\ln \RME}^{3} \left( \frac{\varepsilon }{\RME} \, \|u\|_{H^{1} (\OME)} ^2  + \frac{\varepsilon ^2 }{\RME^2 } \,\|u\|_{H^{1} (\OME)} \|u\|_{H^{1} (\GME)} \right) \nonumber\\
	&\leqslant  C_{\Rsource,N} \,\domi{\ln \RME}^{3} \left( \frac{\varepsilon }{\RME} \, \|u\|_{H^{1} (\OME)} ^2  + \frac{\varepsilon ^{3/2} }{2\RME^2 } \,\|u\|_{H^{1} (\OME)}^2  + \frac{\varepsilon ^{3/2} }{2\RME^2 } \, \varepsilon \,\|u\|_{H^{1} (\GME)}^2  \right) \nonumber\\
	&\leqslant  C_{\Rsource,N} \,\domi{\ln \RME}^{3} \left( \frac{\varepsilon }{\RME} + \frac{\varepsilon ^{3/2} }{\RME^2 } \right) \|u\|_{\VME5,\varepsilon } ^2  \label{eq: caractere bien pose des ME augmentes preuve m=5 forme c}
\end{align}
Now note that, since $\RME\leqslant \Rsource$, we have for~any $M>0$:
\begin{align}
	\varepsilon \leqslant M \, \frac{\RME^{4/3} }{\domi{\ln \RME}^2 } \quad  &\Rightarrow  \quad  \varepsilon ^{3/2} \leqslant  M^{3/2} \, \frac{\RME^2 }{\domi{\ln \RME}^{3}} \; \text{ and } \; \varepsilon \leqslant M\,C_{\Rsource} \frac{\RME}{\domi{\ln \RME}^{3}}  \nonumber\\
	&\Rightarrow  \quad  \domi{\ln \RME}^{3} \left( \frac{\varepsilon }{\RME} + \frac{\varepsilon ^{3/2} }{\RME^2 } \right) \leqslant  C_{\Rsource} \,(M+M^{3/2} ). \label{eq: caractere bien pose des ME augmentes preuve m=5 etude condition stab}
\end{align}
Hence, for a sufficiently small value of $\CstabME$, the stability condition ``$\varepsilon \leqslant \CstabME \,\RME^{4/3} \,\domi{\ln \RME}^{-2}$'' implies:
\begin{enumerate}[label={(\roman*)}]
	\item $\varepsilon  \leqslant \CdtnME1 \,\RME \,\domi{\ln \RME}^{-3}$, so $\opS5$ and $\opP5$ are well-defined,
	\item and $|\cME5(u,u)| \leqslant  \frac{\Ccoer}{2} \|u\|_{\VME5,\varepsilon } ^2 $, which implies that $\aME5+\bME5$ is uniformly coercive (same proof as $m=3$). Moreover, it~is easy to see that we have $|\lME5(v)| \leqslant  C\,(\|\sourceOME\|_{H^{-1}} + \|\sourceGME\|_{L^2 } + \sqrt{\RME} \|\sourceSME\|_{L^2 } +\varepsilon ^{- \frac{1}{2} }|\sourcePME|) \|v\|_{\VME5,\varepsilon }$ for~any $v\in \VME5$. So \cref{eq: FV ME augmente} is again well-posed and uniformly stable.
\end{enumerate}
\end{demo}

\begin{remarquenumerotee}
\label{rq: condition de stabilite et CdtnME}%
It will be useful in the next section to notice the following.
\begin{enumerate}
	\item The above proof shows that the stability condition implies $\varepsilon  \leqslant \CdtnME1 \,\RME \,\domi{\ln \RME}^{-3}$, which is the hypothesis of \cref{opS - opS0}. See \cref{eq: caractere bien pose des ME augmentes preuve def Cstab m=3} for $m=3$ and point (ii) of the proof for $m=5$.
	\item By~reducing the value of $\CstabME$, we can assume that the stability condition implies $\CdtnME2 \.\varepsilon  \.\RME^{-1} \*\domi{\ln \RME}^{3} \leqslant  C_{\Rsource}$ (actually a careful examination of the above proof shows that the given construction of $\CstabME$ directly implies this inequality without having to reduce the value of $\CstabME$).
\end{enumerate}
\end{remarquenumerotee}

\refstepcounter{hypothese}
\begin{paragraphesansboite}{Assumption \thehypothese:}
\label{hyp: Cstab majore}%
Replacing $\CstabME$ by $\min (\CstabME,1)$, we assume without loss of generality that $\sup _{N\in \mathbb{N}} \CstabME$ is finite. Thus the stability condition implies $\varepsilon  \leqslant  C_{\Rsource} \.\RME$, and even $\varepsilon  \leqslant  C_{\Rsource} \RME^{4/3}$ if $m=5$, with constants $C_{\Rsource}$ that are independent of $N$.
\end{paragraphesansboite}

\begin{remarque}
It is easy to extend \cref{caractere bien pose des ME augmentes} to any continuous linear form $\ell $ on $\VME$. Denoting
\[\|\ell \|_{\VME3'} := \sup _{v\in \VME3 \setminus \{0\}} \frac{|\ell (v)| }{\|v\|_{H^{1}}} 
\qquad  \text{ and } \qquad  \|\ell \|_{\VME5'{},\varepsilon } := \sup _{v\in \VME5 \setminus \{0\}} \frac{|\ell (v)| }{\|v\|_{\VME5,\varepsilon }} ,\]
where $\|.\|_{\VME5,\varepsilon }$ is defined in \cref{eq: norme VME}, one would obtain
\[\left\{\begin{array}{@{\,}r@{\;\leqslant \;}l@{\;\,}l}
	\|v_{\varepsilon } \|_{H^{1} (\OME)} & C\. \|\ell \|_{\VME3'} & \text{if } m=3,\\[1ex]
	\|v_{\varepsilon } \|_{H^{1} (\OME)} + \sqrt{\varepsilon } \,\snhg{v_{\varepsilon }} & C\. \|\ell \|_{\VME5'{},\varepsilon } & \text{if } m=5.
\end{array}\right.\]
Note that, in the second case, $\|\ell \|_{\VME5'{},\varepsilon }$ depends on $\varepsilon $.
\end{remarque}

\section{Convergence of the effective models}
\label{sec: cvg ME}

\begin{textesansboite}
In this section, we prove \cref{cvg des ME}. We recall that $\sme$ denotes the solution of the effective model while $\scl:=\sum _{p<m} \sum _{\ell } \varepsilon ^{p} \ln ^{\ell } \!\varepsilon \, u_{p,\ell }$. The proof has three steps.
\begin{enumerate}
	\item First, in \cref{erreur ME v1}, we use the uniform stability of the effective models to prove that $\|\sme-\scl\|_{H^{1} (\OME)}$ is a $O(\varepsilon ^{m-1} )$ plus a truncation error that depends on the number of modes $N$. Thus, in the beginning, we obtain an estimate that is~not the target order $O(\varepsilon ^{m} )$.
	This is classical for effective models associated to Dirichlet problems where a condition of the type $u=\varepsilon  \partial _{n} u$ appears as $ \frac{1}{\varepsilon } \int _{\Gamma } u \bar{v}$ in the variational formulation. This has similarities with penalisation methods for Dirichlet conditions, see e.g.\ \cite{Mau09}.
	\item Then, we improve the error estimate of $\sme-\scl$. For $m=3$ e.g., we use that $w_{\varepsilon } :=\sme3-\scl3$ satisfies $\mu _{1} w_{\varepsilon |\GME} = -\varepsilon  \mu _{0} \partial _{n} w_{\varepsilon }$, so $w_{\varepsilon }$ is the solution of a problem with a Dirichlet data on $\GME$ that is a $O(\varepsilon ^{3} )$. This requires to prove $\|\partial _{n} w_{\varepsilon } \|_{H^{1/2} (\GME)} =O(\varepsilon ^2 )$, which is done by elliptic regularity in \cref{regularite elliptique pour ME}. Special care must be taken at $\PME$ (the tip of $\GME$) where $\OME$ is not regular. Note that an alternative method used in the literature for infinite layers is to build an asymptotic expansion w.r.t.\ $\varepsilon $ of the solution of the effective model whose first terms coincide with the expansion of $u_{\varepsilon }$ (i.e.\ the far fields). See \cite{Via03,Del10,DelHadJol12}. However, this would be complicated here due to the intricate way the DtN maps depend on $\varepsilon $.
	\item Finally we use the error estimate of $\scl-u_{\varepsilon }$ given by \cref{sec: estimations d'erreur} to bound $\|\sme-u_{\varepsilon } \|_{H^{1} (\OME)}$.
\end{enumerate}
\end{textesansboite}

\subsection{Preliminary error estimates of order 2 and 4}

\begin{textesansboite}
The error estimate given by \cref{cvg des ME} is rather clean (the powers in $(\varepsilon /\RME)^{m}$ and $(\RME/R)^{N\smallpst}$ are simple and there~is no $\ln \RME$), but the intermediate steps involve dirty powers of $\RME$. To simplify calculations, we introduce below the notations $\errmdl[s][m+k]$ and $\errdtn[s]$.
\end{textesansboite}

\refstepcounter{notation}
\begin{paragraphesansboite}{Notation \thenotation:}
\label{def: errmdl errdtn}%
We fix a number $\delta $ s.t.\ $0<\delta <\min (1, \pst - \frac{1}{2} )$, for instance $\delta := \frac{1}{2} \min (1, \pst - \frac{1}{2} )$. This is possible because we have $\Theta <2\pi  \Leftrightarrow  \pst > \frac{1}{2}$. Then, for~any $k\in \mathbb{N}$, $s\in \mathbb{R}$ and $R>0$, we denote
\[\errmdl[s][k] :=\displaystyle \varepsilon ^{k} \,\RME^{-m+\delta +s} \,\domi{\ln \varepsilon }^{\ell _{m}} \domi{\ln \RME}^{2m+1} \quad \ \text{ and } \quad \ \errdtn[s]:=\displaystyle \left( \frac{\RME}{R} \right)^{N\pst+s} \domi{\ln \RME}^{2m+1} .\]

In addition, for any $\varepsilon >0$, we introduce the notation $\opD$ to denote the differential operator associated to the effective boundary condition on $\Gamma $:
\[\left\{\begin{array}{ll}
	\opD3(v) := \mu _{0} \partial _{y} v -\frac{\mu _{1} }{\varepsilon } v & \text{for }  m=3\\[2pt]
	\opD5(v) := \mu _{0} \partial _{y} v -\frac{\mu _{1} }{\varepsilon } v + \frac{\varepsilon }{3} (\mu _{1} \partial _{x} ^2 +\omega ^2 \rho _{1} )v & \text{for }  m=5.
\end{array}\right.\]
\end{paragraphesansboite}

\begin{remarque}
$\delta $ should be seen as a small number s.t.\ $\smash{\pst> \frac{1}{2} +\delta }$.
This typically allows simplifications like $\RME^{\smallpst} \leqslant  \RME^{ \frac{1}{2} +\delta }$ and $\domi{\RME^{-p+\delta } }\leqslant \RME^{-p+\delta }$ when $p\geqslant 1$. Although we have $\inf _{\Theta \in (0,2\pi )} \pst = \frac{1}{2}$, we do~not simply replace $\delta $ by 0 (which would give weaker estimates), because we will need the boundedness of $\RME^{\delta } \.\domi{\ln \RME}^{c}$ ($c>0$) at 0 to make the $\ln \RME$ disappear in the final error estimate.
\end{remarque}


\begin{lemmesansboite}{}
\label{estimations CL}%
For~any $p\in \NN$, $\ell \in \mathbb{N}$ and $i,j\in \{0,1\}$, we have $\|\partial _{x} ^{i} \partial _{y} ^{j} u_{p,\ell } \|_{L^2 (\GME)} \leqslant C_{\fsource,p} \*\,\domi{\RME^{-p+\delta +1-i-j }}$.
\end{lemmesansboite}

\begin{demo}
\Cref{validation construction et DA des champs} shows $u_{p,\ell } = T_{\leqslant 2} (\uDA_{p,\ell } )+o_{\partial } (r^2 )$, so there exists $r_{0} >0$ s.t.\ we have on $\Gamma \cap  \{r<r_{0} \}$: $|\partial _{x} ^{i} \partial _{y} ^{j} u_{p,\ell } -\partial _{x} ^{i} \partial _{y} ^{j} T_{\leqslant 2} (\uDA_{p,\ell } )|\leqslant  C_{\fsource,p} \,r^{2-i-j}$. Moreover \cref{precision de l'ansatz sur uDA} yields 
\[T_{\leqslant 2} (\uDA_{p,\ell } ) \,\in \sum _{d\geqslant -p+\pst} \mathcal{A}_{d} (\Omega )\, \subset \! \sum _{d>-p+ \frac{1}{2} +\delta } \mathcal{A}_{d} (\Omega ).\]
Thus, we have on $\Gamma \cap  \{r<r_{0} \}$: $|\partial _{x} ^{i} \partial _{y} ^{j} u_{p,\ell } |\leqslant  C_{\fsource,p} \,r^{-p+ \frac{1}{2} +\delta -i-j}$. Finally, we have $\partial _{y} ^{j} u_{p,\ell } \in H^{1} (\Omega )+\khiz \mathcal{A}(\Omega )$ by \cref{validation construction et DA des champs}, so $\partial _{x} ^{i} \partial _{y} ^{j} u_{p,\ell }$ is also $L^2 $ on $\Gamma  \cap  \{r>r_{0} /2\}$. The bound on $\|\partial _{x} ^{i} \partial _{y} ^{j} u_{p,\ell } \|_{L^2 (\GME)}$ follows by integration.
\end{demo}

\begin{lemmesansboite}{}
\label{norme opS}%
Under the stability condition of \cref{def: condition de stabilite}, we have $\|\opS\|_{\mathcal{L}(L^2 (\SME),\HME{j}(\SME))} \leqslant  C_{\Rsource} \.N^{j+1} /\RME$ for~any $j\in \{0,1,2\}$.
\end{lemmesansboite}


\begin{demo}
By~\cref{rq: condition de stabilite et CdtnME}, the stability condition of \cref{def: condition de stabilite} implies the hypothesis ``$\varepsilon \leqslant \CdtnME1 \,\RME \,\*\domi{\ln \RME}^{-3}$'' of \cref{opS - opS0}. Thus \cref{opS - opS0,rq: condition de stabilite et CdtnME} show that we have
\[\|\opS-\opSref\|_{\mathcal{L}(L^2 (\SME),\HME2(\SME))} \leqslant  \CdtnME[R]2 \. \frac{\varepsilon }{\RME^{2}} \,\domi{\ln \RME}^{3} \leqslant  \frac{C_{\Rsource} }{\RME} .\]
Hence it suffices to prove $\|\opSref\|_{\mathcal{L}(L^2 (\SME),\HME{j}(\SME))} \leqslant  C_{\Rsource} N^{j+1} /\RME$. By~\cref{eq: opS positif qd eps=0 preuve expression opSref}, we have for any $v\in L^2 (\SME)$ and $\bfx \in \SME$:
\[\opSref(v)(\bfx) = \sum _{n=1} ^{N} \langle v,\varphi _{n} \rangle  \, \partial _{r} \varphi _{n} (\bfx) \qquad  \text{with}  \quad  \varphi _{n} := \frac{\mcll_{n\smallpst,0} }{\|\mcll_{n\smallpst,0} \|_{L^2 (\SME)} }.\]
By~\cref{expression sing homog CL}, we have $\mcll_{n\smallpst,0} (r,\theta )=(-1)^{n} \,\Gamma (d+1) \big( \frac{2}{k_{0}} \big)^{d} J_{d} (k_{0} r) \sin (d\theta )$ with $d:=n\pst$. This yields
\[\varphi _{n|\SME} :\theta  \mapsto  \sqrt{ \frac{2}{\RME\Theta } } \sin (d\theta ) \.\beta _{n} \qquad  \text{ and } \qquad  \partial _{r} \varphi _{n|\SME} = \frac{J_{d} '(k_{0} \RME)}{J_{d} (k_{0} \RME)} \varphi _{n|\SME} \]
for some $\beta _{n} \in \mathbb{C}$ with $|\beta _{n} |=1$ (actually $\beta _{n} =(-1)^{n} (k_{0} /|k_{0} |)^{-d}$). Thus $(\varphi _{n|\SME} )_{n}$ is an orthonormal family of $L^2 (\SME)$ and $(\partial _{\theta } ^{j} \partial _{r} \varphi _{n|\SME} )_{n}$ is orthogonal for any $j\in \{0,1,2\}$. We deduce:
\[\big\|\partial _{\theta } ^{j} \opSref(v) \big\|_{L^2 (\SME)} ^2  = \sum _{n=1} ^{N} |\langle v,\varphi _{n} \rangle |^2   \, \|\partial _{\theta } ^{j} \partial _{r} \varphi _{n} \|_{L^2 (\SME)} ^2  \leqslant  \|v\|_{L^2 (\SME)} ^2  \cdot  \max _{n\in [\![1,N]\!]} \|\partial _{\theta } ^{j} \partial _{r} \varphi _{n} \|_{L^2 (\SME)} ^2  .\]
Moreover, by \cite{Leb72} and \cref{minoration majoration BesselJ}, we have for~any $n\in [\![1,N]\!]$:
\[\|\partial _{\theta } ^{j} \partial _{r} \varphi _{n} \|_{L^2 (\SME)} = d^{j} \left| \frac{J_{d} '(k_{0} \RME)}{J_{d} (k_{0} \RME)} \right| = d^{j} \left| \frac{J_{d-1} (k_{0} \RME)-J_{d+1} (k_{0} \RME)}{2\.J_{d} (k_{0} \RME)} \right| \leqslant  \frac{C_{\Rsource} d^{j+1} }{\RME} \leqslant  \frac{C_{\Rsource} N^{j+1} }{\RME} .\]
See \cref{estimation mcld0} for a similar reasoning. To conclude, it suffices to use
\[\big\|\opSref(v) \big\|_{\HME{j}(\SME)} ^2 = \sum\limits _{i=0} ^{j} \big\|\partial _{\theta } ^{j} \opSref(v) \big\|_{L^2 (\SME)} ^2  ,\]
see \cref{ex: sobolev pondere pour ME}.
\end{demo}

\begin{lemme}{consistency errors}
\label{erreur consistence ME}%
Let us denote $w_{\varepsilon } :=\scl-\sme$. Let $R\in (\RME,\Rsource)$ and $j\in \{0,1\}$. Under the stability condition of \cref{def: condition de stabilite}, we have:
\begin{itemize}
	\item $\big\|\opD(w_{\varepsilon } )\big\|_{H^{j} (\GME)} \leqslant C_{\fsource} \,\errmdl[2-j][m-1]$, 
	\item $\displaystyle \big\|\partial _{r} w_{\varepsilon } - \opS(w_{\varepsilon |\SME} )\big\|_{L^2 (\SME)} \leqslant  C_{\fsource,R} \,N \Big( \errmdl + \errdtn \Big)$,
	\item if $m=5$: $\displaystyle \big|\partial _{x} w_{\varepsilon } (\PME) - \opP5(w_{\varepsilon |\SME} )\big| \leqslant  C_{\fsource,R} \.\Big( \errmdl[-1/2] + \errdtn[-1/2] \Big)$.
\end{itemize}
\end{lemme}

\begin{textesansboite}
The consistency error on $\opD(w_{\varepsilon } )$ has only order $m-1$ w.r.t.\ $\varepsilon $, so we are only able to prove $\|\sme-\scl\|_{H^{1} (\OME)} =O(\varepsilon ^{m-1} )$ at first (see \cref{erreur ME v1}). The above consistency errors will be used again in \cref{sec: amelioration erreur ME} to improve the convergence error into a $O(\varepsilon ^{m} )$.
\end{textesansboite}

\begin{demo}
Note that, by~definition, $\sme$ is a solution of \cref{eq: FV ME augmente} with homogeneous boundary conditions on $\GME$, $\SME$ and $\PME$. So $\opD(w_{\varepsilon } ) = \opD(\scl)$, $\partial _{r} w_{\varepsilon } - \opS(w_{\varepsilon |\SME} ) = \partial _{r} \scl - \opS(\scl_{|\SME} )$ and likewise for $\opP$.
\begin{itemize}
	\item \underline{Bound on $\opD3(w_{\varepsilon } )$:} By~\cref{construction directe des upl}, we have on $\Gamma $:
	\begin{align}
		\mu _{0} \partial _{y} \scl3 -\frac{\mu _{1} }{\varepsilon } \scl3 &= \sum _{p\in \NN,p<3} \sum _{\ell \in \mathbb{N}} \varepsilon ^{p} \ln ^{\ell } \!\varepsilon \, \mu _{0} \partial _{y} u_{p,\ell } - \sum _{p\in \NN,p<3} \sum _{\ell \in \mathbb{N}} \varepsilon ^{p-1} \ln ^{\ell } \!\varepsilon \, \mu _{0} \partial _{y} u_{p-1,\ell } \nonumber \\
		&= \sum _{p\in \NN\cap [2,3)} \sum _{\ell \in \mathbb{N}} \varepsilon ^{p} \ln ^{\ell } \!\varepsilon \, \mu _{0} \partial _{y} u_{p,\ell } . \label{eq: erreur ME v1 preuve 1}
	\end{align}
	Let $j\in \{0,1\}$. By~\cref{estimations CL}, we have for~any $p\geqslant 2$: $\|\partial _{y} u_{p,\ell } \|_{H^{j} (\GME)} \leqslant  C_{\fsource,p} \,\domi{\RME^{-p+\delta +1-j}} \leqslant  C_{\fsource,p} \,\RME^{-p+\delta +1-j }$ because $-p+\delta +1-j \leqslant -2+1+1\leqslant 0$. Thus we obtain:
	\begin{align*}
		\big\|\opD3(w_{\varepsilon } )\big\|_{H^{j} (\GME)} &\leqslant  \sum _{p\in \NN\cap [2,3)} \varepsilon ^{p} \domi{\ln \varepsilon }^{\ell _{3}} \cdot  C_{\fsource,p} \,\RME^{-p+\delta +1-j} \\
		&\leqslant  C_{\fsource} \,\varepsilon ^{2} \domi{\ln \varepsilon }^{\ell _{3}} \.\RME^{\delta -1-j}
	\end{align*}
	because the stability condition implies $\varepsilon /\RME \leqslant C_{\Rsource}$ by \cref{hyp: Cstab majore}.

	\item \underline{Bound on $\opD5(w_{\varepsilon } )$:} Similarly, we have $\opD5(w_{\varepsilon } ) = \mu _{0} \partial _{y} \scl5_{|\GME} -\frac{\mu _{1} }{\varepsilon } \scl5_{|\GME} + \frac{\varepsilon }{3} (\mu _{1} \partial _{x} ^2 +\omega ^2 \rho _{1} )\scl5_{|\GME}$. By~\cref{construction directe des upl}, the second term of this is equal to:
	\begin{align*}
		\frac{\mu _{1} }{\varepsilon } \scl5 &= \sum _{p\in \NN,p<5} \sum _{\ell \in \mathbb{N}} \varepsilon ^{p-1} \ln ^{\ell } \!\varepsilon \, \mu _{0} \partial _{y} u_{p-1,\ell } + \sum _{p\in \NN,p<5} \sum _{\ell \in \mathbb{N}} \varepsilon ^{p-1} \ln ^{\ell } \!\varepsilon \, \frac{\mu _{0} }{3} (\partial _{x} ^2 +k_{1} ^2 ) \partial _{y} u_{p-3,\ell } \\
		&= \mu _{0} \partial _{y} \scl4 + \frac{\varepsilon ^2  \mu _{0} }{3} (\partial _{x} ^2 +k_{1} ^2 ) \partial _{y} \scl2 \qquad \text{ on }\Gamma .
	\end{align*}
	Hence we have:
	\begin{align*}
		\|\opD5(w_{\varepsilon } )\|_{H^{j} (\GME)} &= \bigg\|\sum _{p\in \NN\cap [4,5)} \sum _{\ell \in \mathbb{N}} \varepsilon ^{p} \ln ^{\ell } \!\varepsilon \, \mu _{0} \partial _{y} u_{p,\ell } + \frac{\varepsilon  \mu _{1} }{3} (\partial _{x} ^2 +k_{1} ^2 )\scl5 -  \frac{\varepsilon ^2  \mu _{0} }{3} (\partial _{x} ^2 +k_{1} ^2 ) \partial _{y} \scl2 \bigg\|_{H^{j} (\GME)} \\
		&\leqslant  \bigg\|\sum _{p\in \NN\cap [4,5)} \sum _{\ell \in \mathbb{N}} \varepsilon ^{p} \ln ^{\ell } \!\varepsilon \, \mu _{0} \partial _{y} u_{p,\ell } \bigg\|_{H^{j} (\GME)} + \frac{\varepsilon  \mu _{1} }{3} \big\|(\partial _{x} ^2 +k_{1} ^2 )(\scl5-\scl3)\big\|_{H^{j} (\GME)} \\
		& \qquad  + \frac{\varepsilon ^2 }{3} \Big\|(\partial _{x} ^2 +k_{1} ^2 )\underbrace{\left(\frac{\mu _{1} }{\varepsilon } \scl3-\mu _{0} \partial _{y} \scl2\right)}_{\textstyle =0 \text{ by \cref{eq: erreur ME v1 preuve 1}} } \Big\|_{H^{j} (\GME)}
	\end{align*}
	The first term can be bounded in the same way as for $\opD3(w_{\varepsilon } )$. Moreover, by \cref{estimations CL} we have for~any $p\geqslant 3$: $\|\partial _{x} ^2  u_{p,\ell } \|_{H^{j} (\GME)} \leqslant  C_{\fsource,p} \,\RME^{-p+\delta -j}$ because $-p+\delta -j\leqslant -3+1\leqslant 0$. This allows us to bound the second term: 
	\[\frac{\varepsilon  \mu _{1} }{3} \.\big\|(\partial _{x} ^2 +k_{1} ^2 )(\scl5-\scl3)\big\|_{H^{j} (\GME)} \leqslant  \sum _{p\in \NN\cap [3,5)} \varepsilon ^{p+1} \domi{\ln \varepsilon }^{\ell _{5}} \cdot  C_{\fsource,p} \,\RME^{-p+\delta -j} .\]
	It is easy to deduce similarly $\|\opD5(w_{\varepsilon } )\|_{H^{j} (\GME)} \leqslant  C_{\fsource} \,\varepsilon ^{4} \.\domi{\ln \varepsilon }^{\ell _{5}} \.\RME^{\delta -3-j}$.

	\item \underline{Bound on $\opS(w_{\varepsilon |\SME} )$:} Let $h:=\partial _{r} w_{\varepsilon } - \opS(w_{\varepsilon |\SME} )$ and
	\[\tilde{u} := \sum _{n=1} ^{N} \sigma _{n\pst} (\sclDA) \cdot \mme_{n} \,\in \,\operatorname{Span}(\mme_{1} ,\dots ,\mme_{N} ).\]
	By~construction of $\opS$ (see \cref{sec: construction de opS et opP}) we have $\opS( \tilde{u}_{|\SME} )=\partial _{r} \tilde{u}_{|\SME}$. Thus
	\[h = \partial _{r} \scl_{|\SME} -\opS(\scl_{|\SME} ) = \partial _{r} (\scl- \tilde{u})_{|\SME} -\opS(\scl_{|\SME} - \tilde{u}_{|\SME} ) \]
	which implies
	\[\|h\|_{L^2 (\SME)} \leqslant  \|\partial _{r} (\scl- \tilde{u})\|_{L^2 (\SME)}  + \|\opS\|_{\mathcal{L}(L^2 (\SME))} \|\scl- \tilde{u}\|_{L^2 (\SME)} .\]
	Besides, thanks to \cref{hyp: Cstab majore}, the stability condition implies $\varepsilon \leqslant  C_{\Rsource} ^{\star} \.\RME\.\domi{\ln \RME}^{-3}$ for some constant $C_{\Rsource} ^{\star}$. Applying \cref{decompo modale SE 2} with $M:=C_{\Rsource} ^{\star}$ and \cref{norme opS}, we obtain
	\[\|h\|_{L^2 (\SME)} \leqslant  C_{\fsource,R} \, \frac{N}{\sqrt{\RME}} \left[ \left( \frac{\RME}{R} \right)^{N\pst+ \frac{1}{2} }  + \frac{\varepsilon ^{m} }{\RME^{m-\pst}} \.\domi{\ln \varepsilon }^{\ell _{m}} \right] \domi{\ln \RME}^{2m+1} .\]
	Finally, we bound the factor $1/\RME^{m-\pst}$ by $\RME^{-m+\delta + \frac{1}{2} }$.

	\item \underline{Bound on $\opP5(w_{\varepsilon |\SME} )$:} By~\cref{rq: condition de stabilite et CdtnME}, the assumptions of \cref{opS - opS0} are fulfilled, so we have $\sqrt{\RME} \|\opP5\|_{\mathcal{L}(L^2 (\SME),\mathbb{C})} \leqslant  \CdtnME2 \,\varepsilon \,\RME^{-2} \domi{\ln \RME}^{3}$. Using again \cref{rq: condition de stabilite et CdtnME}, this implies $\sqrt{\RME} \|\opP5\|_{\mathcal{L}(L^2 (\SME),\mathbb{C})} \leqslant  C_{\Rsource} /\RME$, where this bound is independent of $N$. Moreover, we have $\opP5( \tilde{u}_{|\SME} )=\partial _{x} \tilde{u}(\PME)$, so \cref{decompo modale SE 2} gives
	\begin{align*}
		\big|\partial _{x} w_{\varepsilon } (\PME) - \opP5(w_{\varepsilon |\SME} )\big| &= \frac{\mu _{1} }{3} \Big| \partial _{x} (\scl5- \tilde{u})(\PME) -\opP5(\scl5_{|\SME} - \tilde{u}_{|\SME} )\Big| \\
		&\leqslant  C\,| \partial _{x} (\scl5- \tilde{u})(\PME)| + C \.\|\opP5\|_{\mathcal{L}(L^2 (\SME),\mathbb{C})} \,\|\scl5- \tilde{u}\|_{L^2 (\SME)} \\
		&\leqslant  C_{\fsource,R} \, \frac{1}{\RME} \left[ \left( \frac{\RME}{R} \right)^{N\pst+ \frac{1}{2} }  + \frac{\varepsilon ^{m} }{\RME^{m-\pst}} \.\domi{\ln \varepsilon }^{\ell _{m}} \right] \domi{\ln \RME}^{2\cdot 5+1}
	\end{align*}
\end{itemize}
\end{demo}

\begin{textesansboite}
Let us denote $\eme:= \scl-\sme$. Let $\emeg$ be the unique solution in $\VME$ of:
\begin{equation}{\label{eq: def emeg}}
\def\emeg{w}
\left\{ \begin{array}{r@{\,=\,}ll}
	(\Delta  +k_{0} ^2 )\emeg & 0 &\text{ in }\OME\\
	\opD(\emeg) & \opD(\eme) &\text{ on }\GME\\
	\emeg & 0 &\text{ on }\bordgaucheME\\
	\partial _{r} \emeg - \opS (\emeg_{|\SME} ) & 0 &\text{ on }\SME \\
	\partial _{x} \emeg - \opP (\emeg_{|\SME} ) & 0 & \text{ at }  \PME\ \text{if } m=5
\end{array}\right.\end{equation}
and $\emes$ be the unique solution in $\VME$ of:
\begin{equation}{\label{eq: def emes}}
\def\emes{w}
\left\{ \begin{array}{r@{\,=\,}ll}
	(\Delta  +k_{0} ^2 )\emes & 0 &\text{ in }\OME\\
	\opD(\emes) & 0 &\text{ on }\GME\\
	\emes & 0 &\text{ on }\bordgaucheME\\
	\partial _{r} \emes - \opS (\emes_{|\SME} ) & \partial _{r} \eme - \opS (\eme_{|\SME} ) &\text{ on }\SME \\
	\partial _{x} \emes - \opP (\emes_{|\SME} ) & \partial _{x} \eme - \opP (\eme_{|\SME} ) & \text{ at }  \PME\ \text{if } m=5
\end{array}\right.\end{equation}
They are well defined thanks to \cref{caractere bien pose des ME augmentes}. Moreover, by linearity and uniqueness of the solution of \cref{eq: FV ME augmente}, we have $\eme = \emeg+\emes$.
\end{textesansboite}

\begin{lemme}{}
\label{erreur ME v1}%
Let $R\in (\RME,\Rsource)$. Under the stability condition of \cref{def: condition de stabilite}, we have if $m=3$:
\begin{itemize}
	\item $\|\emeg3\|_{H^{1} (\OME)} \leqslant  C_{\fsource} \errmdl[2][m-1]$,
	\item $\|\emes3\|_{H^{1} (\OME)} \leqslant  C_{\fsource,R} \,N \big(\errmdl[1/2] + \errdtn[1/2] \big)$,
\end{itemize}
and if $m=5$:
\begin{itemize}
	\item $\|\emeg5\|_{H^{1} (\OME)} + \sqrt{\varepsilon } \,\snhg{\emeg5} \leqslant  C_{\fsource} \errmdl[2][m-1]$,
	\item $\|\emes5\|_{H^{1} (\OME)} + \sqrt{\varepsilon } \,\snhg{\emes5} \leqslant  C_{\fsource,R} \,N \big(\errmdl[1/6] + \errdtn[1/6] \big)$.
\end{itemize}
\end{lemme}

\begin{demo}
It suffices to combine \cref{caractere bien pose des ME augmentes,erreur consistence ME}. For $\emeg$, we have if $m=3$
\[\|\emeg3\|_{H^{1} (\OME)} \leqslant  \CcontME \,\|\opD3(\eme3)\|_{L^2 (\GME)} \leqslant  C_{\fsource} \,\errmdl[2][m-1] ,\]
and likewise for $m=5$. As~for $\emes$, we have when $m=3$:
\[\|\emes3\|_{H^{1} (\OME)} \leqslant  \CcontME \.\sqrt{\RME} \.\big\|\partial _{r} w_{\varepsilon } - \opS3(w_{\varepsilon |\SME} )\big\|_{L^2 (\SME)} \leqslant  C_{\fsource,R} \,N \big( \errmdl[1/2] + \errdtn[1/2] \big).\]
Similarly, when $m=5$, we have
\begin{multline*}
	\|\emes5\|_{H^{1} (\OME)} + \sqrt{\varepsilon } \.\snhg{\emes5} \leqslant  \CcontME \.\Big( \sqrt{\RME} \.\big\|\partial _{r} w_{\varepsilon } - \opS5(w_{\varepsilon |\SME} )\big\|_{L^2 (\SME)} \\
	{}+\sqrt{\varepsilon } \.\big|\partial _{x} w_{\varepsilon } (\PME) - \opP5(w_{\varepsilon |\SME} )\big| \Big)
\end{multline*}
Then we use that the stability condition implies $\varepsilon \leqslant C_{\Rsource} \RME^{4/3}$ (see \cref{hyp: Cstab majore}), so $\sqrt{\varepsilon } \leqslant C_{\Rsource} \RME^{2/3}$:
\begin{align*}
	\|\emes5\|_{H^{1} (\OME)} + \sqrt{\varepsilon } \.\snhg{\emes5} &\leqslant  C_{\fsource,R} \,N \big( \errmdl[1/2] + \errdtn[1/2] \big) + C_{\fsource,R} \.\big( \errmdl[1/6] + \errdtn[1/6] \big)\\
	&\leqslant  C_{\fsource,R} \,N \big(\errmdl[1/6] + \errdtn[1/6] \big)
\end{align*}
\end{demo}

\subsection{Tools of elliptic regularity}
\label{sec: lemmes de regularite elliptique ME}

\begin{textesansboite}
The proofs of the results of this section can be found in Appendix~\ref{sec: preuves lemmes de regularite elliptique ME}.\\

The next result shows that, if a function $u$ satisfies the Poisson equation and a boundary condition similar to the effective conditions on $\GME$, then the trace of $u$ on the boundary is a $O(\varepsilon )$ with a constant that depends on $\|u\|_{H^{1}}$ and the data. We will apply this result to $\emeg$. We already know that it~is $O(\varepsilon ^{m-1} )$ by \cref{erreur ME v1}, and \cref{regularite elliptique pour ME} will help us prove that it~is $O(\varepsilon ^{m} )$. 
\end{textesansboite}

\begingroup
\def\recI{R}%
\def\recJ{R'}%
\def\J{J}%
\begin{lemmesecable}{}
\label{regularite elliptique pour ME}%
We consider $\eta >0$ and $I,\J$ two open intervals of $\mathbb{R}$ s.t.\ $\overline{\J}\subset I$. We denote $\recI:=I\times (0,\eta )$, $I_{0} :=I\times \{0\}$ and $\J_{0} :=\J\times \{0\}$. Let $0<\varepsilon \leqslant \varepsilon _{0}$ and $\nu _{1} ,\nu _{2} \in \mathbb{R}_+^*$. Let $A:\overline{\recI}\rightarrow  \mathcal{M}_{2} (\mathbb{R})$ be a field of symmetric matrices. We assume that $A$ is elliptic and Lipschitz, and that $A(x)$ is the identity matrix for~any $x\in I_{0}$. Let $\mathcal{P}:=(\eta ,I,\J,\varepsilon _{0} ,\nu _{1} ,\nu _{2} ,A)$ denote the tuple of those parameters. Let $f\in L^2 (\recI)$, $g\in H^{1} (I_{0} )$ and $u\in H^{1} (\recI)$.
\begin{enumerate}
	\item If $u$ satisfies 
	\[\left\{\begin{array}{rl}
		\operatorname{div}(A \nabla  u)=f &\text{ in }\recI\\[2pt]
		\partial _{n} u + \frac{\nu _{1} }{\varepsilon } u=g &\text{ on }I_{0} ,
	\end{array}\right.\]
	then we have
	\[\|u\|_{H^{1/2} (J_{0} )} \leqslant  C_{\mathcal{P}} \,\varepsilon \,\big(\|f\|_{L^2 (\recI)} +\|g\|_{H^{1} (I_{0} )} + \|u\|_{H^{1} (\recI)}\big) .\]
	\item If $u_{|I_{0}} \in H^{1} (I_{0} )$ and $u$ satisfies 
	\[\left\{\begin{array}{rl}
		\operatorname{div}(A \nabla  u) =f &\text{ in }\recI\\[2pt]
		\partial _{n} u + \frac{\nu _{1} }{\varepsilon } u- \varepsilon \,\nu _{2} \partial _{x} ^2 u=g &\text{ on }I_{0} ,
	\end{array}\right.\]
	then we have
	\[\|u\|_{H^{1/2} (\J_{0} )} \leqslant  C_{\mathcal{P}} \,\varepsilon \,\big(\|f\|_{L^2 (\recI)} + \|g\|_{H^{1} (I_{0} )} + \|u\|_{H^{1} (\recI)} + \sqrt{\varepsilon } \,\snh{u}(I_{0} ) \big) .\]\fbeq[1.1]
\end{enumerate}
\end{lemmesecable}

\begin{remarque}
Note that the constants $C_{I,\J,\eta ,A,(\nu _{i} ),\varepsilon _{0} }$ do~not depend on $\varepsilon $. The lemma can be applied for instance to $I=(a,d)$ and $\J=(b,c)$ with $a<b<c<d$ in $\mathbb{R}$, or to $I=(a,\infty )$ and $\J=(b,\infty )$ with $a<b$ in $\mathbb{R}$, or to $I=\J=\mathbb{R}$. The proof relies on classical tools of elliptic regularity.
\end{remarque}

\begin{demo}
See Appendix~\ref{sec: preuves lemmes de regularite elliptique ME}, page~\pageref{preuve: regularite elliptique pour ME}.
\end{demo}

\begin{textesansboite}
\Cref{regularite elliptique pour ME} can be directly applied to $u:=\emeg$ in a vicinity of $\GME$ excluding a vicinity of $\PME$. However, the domain $\OME$ has a corner at $\PME$, so \cref{regularite elliptique pour ME} cannot be directly used around $\PME$. This corner has an angle of $\pi /2$, so by applying a diffeomorphism around $\OME$, we can locally deform the geometry into a quarter plane. Then, $u$ will be extended by symmetry to a half plane, where we can finally use \cref{regularite elliptique pour ME}. More details will be given in the proof of \cref{regularite elliptique pour ME 2}. The use of that diffeomorphism explains why we introduced a matrix $A$ in \cref{regularite elliptique pour ME} instead of simply assuming $\Delta u=f$. But we kept the coefficients $\nu _{1}$ and $\nu _{2}$ constant as this is needed in the proof. Thus the diffeomorphism must be chosen wisely so that it does~not introduce variable coefficients in the boundary condition (this is done in \cref{diffeo ME}).
\end{textesansboite}
\endgroup

\begin{lemme}{}
\label{diffeo ME}%
There exist a neighborhood $V_{1}$ of $\PME$ in $\mathbb{R}^2 $, a neighborhood $V_{2}$ of $(0,0)$ in $\mathbb{R}^2 $ and a $\mathcal{C}^{1}$-diffeomorphism $\varphi  :V_{1} \rightarrow V_{2} $ s.t.:
\begin{itemize}
	\item $\varphi (\PME) = (0,0)$,
	\item $\varphi (\OME \cap V_{1} )=(\mathbb{R}_+^*)^2  \cap V_{2} $,
	\item $\varphi (\GME \cap V_{1} )=(\mathbb{R}_+^*\times \{0\}) \cap V_{2} $,
	\item $\varphi (\SME \cap V_{1} )=(\{0\}\times \mathbb{R}_+^*) \cap V_{2} $,
	\item the differential $\diff \varphi $ is Lipschitz,
	\item $\forall \bfx \in \GME \cap V_{1} , \; \diff \varphi (\bfx) =  \mathrm{id}_{\mathbb{R}^2 }$,
	\item $\forall \bfx \in \SME \cap V_{1} , \; \diff \varphi (\bfx) = R_{-\theta }$, where $\theta $ is the polar coordinate and $R_{-\theta }$ is the rotation of $\mathbb{R}^2 $ of angle $-\theta $.
\end{itemize}
\end{lemme}

\begin{textesansboite}
See \cref{fig: diffeo ME} for a representation of $\varphi $. The conditions on $D\varphi $ on $\partial \OME$ ensure that the boundary conditions on $\partial \OME$ (of type Neumann, Robin or Ventcel) are preserved when applying $\varphi $, which is needed to use \cref{regularite elliptique pour ME}. The Lipschitz property of $D\varphi $ is also useful to obtain a Lipschitz matrix in the divergence as in \cref{regularite elliptique pour ME}.\\

Note that the existence of $\varphi $ is not obvious because of these conditions on $D\varphi $. For example, the function $\varphi _{1}$ that maps $(x,y)$ to $(r,\RME \theta )$ (using the polar coordinates) satisfies all the desired properties except the last one because: $\forall (x,y)\in \GME \cap V_{1} ,\; \frac{\partial }{\partial y} (\RME\theta  )\big|_{x,y} = \RME/r \neq 1$. The function $\varphi _{2} :(x,y)\mapsto (r,r\theta )$ solves this problem, but adds another one: $\forall (x,y)\in \SME \cap V_{1} ,\; \frac{\partial }{\partial r} (r\theta )\big|_{x,y} = \theta  \neq 0$. The proof of \cref{diffeo ME} makes $\varphi $ equal to $\varphi _{1}$ near $\SME$, equal to $\varphi _{2}$ near $\GME$, and uses a partition of unity in the middle to make the transition. With this construction, the difficulty is to prove that $\varphi $ is differentiable at $\PME$ and that $D\varphi $ is Lipschitz.
\end{textesansboite}

	\tikzset{external/export=true}%
	\tikzsetfigurename{fig-diffeo_ME-}%

\begin{figure}[h]
\setlength{\fboxsep}{1.5pt}
\def\xMin{-0.05}
\def\xMax{3.2}
\def\yMin{-0.4}
\def\yMax{2.2}
\def\ThetaFig{105}
\def\RMEFig{1.3}
\def\sigmo(#1,#2){sin(min(45,max(-45,30*ln(#2/(#1)))) +45)^2}
\def\cos(#1){cos(180/pi*#1)}
\def\sin(#1){sin(180/pi*#1)}
\def\diffeoC#1#2{{\RMEFig+#1},{#2}} 
\def\diffeoP#1#2{{(\RMEFig+#1)*\cos(#2/\RMEFig)},{(\RMEFig+#1)*\sin(#2/\RMEFig)}} 
\def\diffeoME#1#2{{(\RMEFig+#1)*(1+(\cos(#2/\RMEFig)-1)*\sigmo(#1,#2))},{#2*(1-\sigmo(#1,#2)) + (\RMEFig+#1)*\sin(#2/\RMEFig)*\sigmo(#1,#2)}}

\def\patateAvt#1{plot [smooth cycle,tension=0.6] coordinates {(\diffeoC{1.5}{0}) (\diffeoME{1.5}{0.55}) (\diffeoME{1.6}{1}) (\diffeoME{1.3}{1.27}) (\diffeoME{0.6}{1.3}) (\diffeoP{0}{1.33}) (\diffeoP{-0.45}{0.7}) (\diffeoC{-0.2}{-0.3}) (\diffeoC{1}{-0.4})}}
\def\patateApr{plot [smooth cycle,tension=0.6] coordinates {(1.5,0) (1.5,0.55) (1.6,1) (1.3,1.27) (0.6,1.3) (0,1.33) (-0.45,0.7) (-0.2,-0.3) (1,-0.4)}}

\hbox to \linewidth{\hfill%
\begin{tikzpicture}[scale=1.2]
	\begin{scope}
	\clip \patateAvt{\RMEFig+};
	\foreach \i in {1,2,...,12} { 
		\draw[domain=0.01:2,samples=100,gray,thin] plot (\diffeoME{\x}{\i/6});
		\draw[domain=-1:0,gray,thin] plot (\diffeoP{\x}{\i/6});
	}
	\foreach \i in {1,2,...,10} { 
		\draw[domain=0.01:1.5,samples=100,gray,thin] plot (\diffeoME{\i/6}{\x});
		\draw[domain=-1:0,gray,thin] plot (\diffeoC{\i/6}{\x});
	} 
	\foreach \i in {-1,-2,...,-7} { 
		\draw[domain=-1:2,gray,thin] plot (\diffeoC{\x}{\i/6});
	}
	\foreach \i in {0,-1,...,-7} { 
		\draw[domain=0.01:1.5,gray,thin] plot (\diffeoP{\i/6}{\x});
		\draw[domain=-1:0,gray,thin] plot (\diffeoC{\i/6}{\x});
	}
	\end{scope}
	\draw [>=stealth,->] (\xMin-0.4,0) -- (\xMax+0.4,0) node [right] {$x$};
	\draw [>=stealth,->] (0,\yMin-0.3) -- (0,\yMax+0.3) node [right] {$y$};
	\draw \patateAvt{\RMEFig+};
	\draw [thick] ({\RMEFig*cos(\ThetaFig)}, {\RMEFig*sin(\ThetaFig)}) arc (\ThetaFig:0:\RMEFig) -- (\xMax,0);
	\draw (\xMax*0.87,\yMax*0.87) node {$\OME$};
	\draw (\xMax*0.58,\yMax*0.38) node {\colorbox{white}{$V_1$}};
	\draw [<-,>=latex] (\xMax*0.95,0) to[bend right=20] ++(0.15,-0.55) node[below]{$\ \ \GME$};
	\draw (\RMEFig,0) node{$\bullet$};
	\draw [<-,>=latex] (\RMEFig+0.05,-0.05) to[bend left=15] ++(0.3,-0.55) node[below]{$\ \PME$};
	\draw [<-,>=latex] ({\RMEFig*cos(70)}, {\RMEFig*sin(70)}) to[bend right=15] ++(0.1,0.45) node[above]{$\ \SME$};
	\draw [<->,>=latex] (0,-0.2) -- node[midway, below] {$\RME$} (\RMEFig, -0.2);
\end{tikzpicture}%
\quad
\raisebox{17mm}{\large
	$\overset{\varphi}{\relbar\joinrel\longrightarrow}$}
\ \quad
\def\xMin{-0.3}%
\def\xMax{2}%
\begin{tikzpicture}[scale=1.2]
	\begin{scope}
	\clip \patateApr;
	\foreach \i in {-2,-1,...,12} { 
		\draw[domain=-1:2,gray,thin] plot (\x,\i/6);
	}
	\foreach \i in {-2,-1,...,10} { 
		\draw[domain=-1:2,gray,thin] plot (\i/6,\x);
	}
	\end{scope}
	\draw [>=stealth,->] (\xMin-0.4,0) -- (\xMax+0.4,0) node [right] {$x$};
	\draw [>=stealth,->] (0,\yMin-0.3) -- (0,\yMax+0.3) node [right] {$y$};
	\draw \patateApr;
	\draw [thick] (0,\yMax) -- (0,0) -- (\xMax,0);
	\draw (\xMax*0.35,\yMax*0.28) node {\colorbox{white}{$V_2$}};
	\draw (0,0) node{$\bullet$};
	\draw [<-,>=latex,opacity=0] (\RMEFig+0.05,-0.05) to[bend left=15] ++(0.3,-0.55) node[below]{$\ \PME$};
\end{tikzpicture}%
\hfill}

\longCaption{The diffeomorphism $\varphi$ that sends $\OME$ to $(\mathbb{R}_+^*)^2$ in the vicinity of $\PME$}{. $\varphi$ looks like the cartesian coordinates near $\GME$ and like the polar coordinates near $\SME$.}
\label{fig: diffeo ME}
\end{figure}%
	\tikzset{external/export=false}%
	\tikzsetfigurename{main-figure}

\begin{demo}
See Appendix~\ref{sec: preuves lemmes de regularite elliptique ME}, page~\pageref{preuve: diffeo ME}.
\end{demo}

\begin{textesansboite}
We recall the scaled Sobolev spaces $\HME{s}$ are introduced in \cref{def: sobolev pondere pour ME}.
\end{textesansboite}

\begin{lemmesecable}{}
\label{regularite elliptique pour ME 2}%
Let $0<\RME\leqslant R$, $M>0$ and $0<\varepsilon \leqslant M\RME$. Let $g\in H^{1} (\GME)$, $h\in H^2 (\SME)$ and $u\in H^{1} (\OME)$.
\begin{enumerate}
	\item If $u$ satisfies
	\[\left\{\begin{array}{r@{}ll}
		\mu _{0} \Delta u+\omega ^2 \rho _{0} u &{}= 0 &\text{ in }\OME\\
		\mu _{0} \partial _{n} u + \frac{\mu _{1} }{\varepsilon } u &{}=g &\text{ on }\GME\\
		\partial _{n} u &{}= h &\text{ on }\SME
	\end{array}\right.\]
	\[\|u\|_{\HME{1/2}(\GME)} \leqslant  C_{R,M} \,\varepsilon \. \left( \|g\|_{\HME1(\GME)} + \|h\|_{\HME2(\SME)} +  \frac{\|u\|_{H^{1} (\OME)}}{\RME^{3/2}} \right) + C_{M} \,\RME^{3/2} \, |h(\PME)|.\]
	\item If $u_{|\GME} \in H^{1} (\GME)$ and $u$ satisfies
	\[\left\{\begin{array}{r@{}ll}
		\mu _{0} \Delta u+\omega ^2 \rho _{0} u &{}= 0 &\text{ in }\OME\\
		\mu _{0} \partial _{n} u + \frac{\mu _{1} }{\varepsilon } u - \varepsilon \,(\mu _{1} \partial _{x} ^2 +\omega ^2 \rho _{1} )u &{}=g &\text{ on }\GME\\
		\partial _{n} u &{}= h &\text{ on }\SME \text{ (a.e.)} \\
		\partial _{x} u &{}= h &\text{ at }  \PME
	\end{array}\right.\]
	then we have
	\begin{align*}
		\|u\|_{\HME{1/2}(\GME)} \leqslant {}& C_{R,M} \,\varepsilon \.\left( \|g\|_{\HME1(\GME)} + \|h\|_{\HME2(\SME)} +  \frac{\|u\|_{H^{1} (\OME)}}{\RME^{3/2}} + \sqrt{ \frac{\varepsilon }{\RME} } \.\snhg{u} \right) \\
		&+ C_{M} \,\RME^{3/2} \, |h(\PME)|.
	\end{align*}\fbeq\fbi
\end{enumerate}
\end{lemmesecable}

\begin{demo}
See Appendix~\ref{sec: preuves lemmes de regularite elliptique ME}, page~\pageref{preuve: regularite elliptique pour ME 2}.
\end{demo}

\begin{remarque}
The term involving $|h(\PME)|$ is~not multiplied by $\varepsilon $, but in practice \cref{regularite elliptique pour ME 2} will be applied to a function $h$ that satisfies $|h(\PME)| = O(\varepsilon )\. \|h\|_{L^2 (\SME)}$. See the proof of \cref{cvg des ME} page~\pageref{preuve: cvg des ME}.
\end{remarque}

\subsection{Final error estimates of order~3 and 5}
\label{sec: amelioration erreur ME}

\begin{lemme}{}
\label{ME avec Dirichlet sur Gamma}%
\def\inc{u}%
There exists $C^{\mathrm{dir}} _{\Rsource,N} >0$ s.t., for any $\varepsilon \leqslant  C^{\mathrm{dir}} _{\Rsource,N} \,\RME \,\domi{\ln \RME}^{-3}$ and $\sourceGME \in H^{1/2} (\GME)$, the following problem has a unique solution in $H^{1} (\OME)$
\begin{equation}{\label{eq: ME avec Dirichlet sur Gamma}}
\left\{ \begin{array}{r@{\,=\,}ll}
	(\Delta  +k_{0} ^2 )\inc & 0 &\text{ in }\OME\\
	\inc & \sourceGME &\text{ on }\GME\\
	\inc & 0 &\text{ on }\bordgaucheME\\
	\partial _{r} \inc - \opS (\inc_{|\SME} ) & 0 &\text{ on }\SME 
\end{array}\right.
\end{equation}
Moreover, the solution satisfies $\|\inc\|_{H^{1} (\OME)} \leqslant  C_{\Rsource} \. \RME^{-1/2} \,\|\sourceGME\|_{\HME{1/2}(\GME)}$
\end{lemme}

\begin{remarquenumerotee}
\label{rq: ME avec Dirichlet sur Gamma}%
By~reducing the value of $\CstabME$, we can assume that the stability condition of \cref{def: condition de stabilite} implies $\varepsilon \leqslant  C^{\mathrm{dir}} _{\Rsource,N} \,\RME \,\domi{\ln \RME}^{-3}$.
\end{remarquenumerotee}

\begin{demo}
\def\bbMEm{\tilde{b}^{(m)}_\peps}%
\def\lift{\opensubscript{u_{\rm lift}}}%
\def\usc{\opensubscript{u_{\rm sc}}}%
For~any $u,v\in H^{1} (\Omega )$ let us define
\[a(u,v) := \int _{\OME} (\mu _{0} \nabla  u \cdot  \nabla  \bar{v}- \omega ^2  \rho _{0} \.u \bar{v}) \qquad  \text{ and } \qquad  \bbMEm(u,v) := \int _{\SME} \mu _{0} \,\opS(u_{|\SME} )\, \bar{v}.\]
Similarly to \cref{coercivite partie hors DtN des ME}, $a$ satisfies the uniform coercivity property: $\forall u\in H^{1} (\Omega ),\; \Im( \overline{\omega }) \*\,\Im( \overline{\omega }\,a(u,u))\geqslant  C\.\|u\|_{H^{1} (\OME)}^2 $. In a very similar way to \cref{caractere bien pose des ME augmentes}, one can prove that there exists $C^{\mathrm{dir}} _{\Rsource,N} >0$ s.t., under the assumptions $\varepsilon \leqslant  C^{\mathrm{dir}} _{\Rsource,N} \,\RME \,\domi{\ln \RME}^{-3}$ and $R\leqslant \Rsource$, $a+\bbMEm$ is uniformly coercive and uniformly continuous:
\[\forall u\in H^{1} (\Omega ), \quad  C\.\|u\|_{H^{1} (\OME)}^2  \leqslant  \big|a(u,u)+\bbMEm(u,u)\big| \leqslant  C_{\Rsource} \. \|u\|_{H^{1} (\OME)}^2 .\]
Let us show that there exists $\lift \in \VME3$ s.t.\ $\lift_{|\GME} =\sourceGME$ and $\|\lift\|_{H^{1} (\OME)} \leqslant  C_{\Rsource} \RME^{-1/2} \.\|\sourceGME\|_{\HME{1/2}(\GME)}$. We recall that $\VME3:=\{u\in H^{1} (\OME) \mid u_{|\bordgaucheME} =0\}$.
\begin{sousdemo}
	Let $\usc \in H^{1} (\Omega  \setminus \closedBall(0,1))$ be s.t.\ $\usc_{|\Gamma } =\sourceGME(\RME\varlibre)$, $\usc_{|\bordgauche} =0$ and $\|\usc\|_{H^{1}} \leqslant  C\.\|\sourceGME(\RME\varlibre)\|_{H^{1/2} }$. Let $\lift:x\mapsto \usc(x/\RME)$. We have $\lift \in \VME3$ and, by~definition of $\HME1$ (see \cref{def: sobolev pondere pour ME}),
	\[\|\lift\|_{H^{1} (\OME)} = \left(\RME^2  \.\|\usc\|_{L^2 (\OME)} ^2  +\|\nabla \usc\|_{L^2 (\OME)} ^2 \right)^{ \frac{1}{2} } \leqslant  C_{\Rsource} \|\usc\|_{H^{1}} \leqslant  C_{\Rsource} \RME^{-1/2} \.\|\sourceGME\|_{\HME{1/2}(\GME)} .\vspace{-2ex}\]
\end{sousdemo}
Let $H:=\{u\in H^{1} (\OME) \mid u_{|\bordgaucheME\cup \GME} =0\}$. $u$ is a solution of \cref{eq: ME avec Dirichlet sur Gamma} iff $u_{0} :=u-\lift$ satisfies $u_{0} \in H$ and
\[\forall v\in H,\qquad  a(u_{0} ,v)+\bbMEm(u_{0} ,v) = -a(\lift,v)-\bbMEm(\lift,v).\]
By~uniform coercivity and continity of $a+\bbMEm$ and using Lax-Milgram's theorem, this problem has a unique solution and it satisfies $\|u_{0} \|_{H^{1} (\OME)} \leqslant  C\.\|\lift\|_{H^{1} (\OME)} \leqslant C_{\Rsource} \RME^{-1/2} \.\|\sourceGME\|_{\HME{1/2}(\GME)}$.
\end{demo}

\begin{textesansboite}
The next lemma is a refinement of the bound of $\|\opP\|_{\mathcal{L}(L^2 (\SME),\mathbb{C})}$ given by \cref{opS - opS0}. It makes the behavior of the estimate w.r.t.\ $N$ explicit. This requires to study finely the behavior w.r.t.\ $n$ of the functions $\mme_{n}$, which is done in Appendix~\ref{sec: annexe comportement en N des ME}, see page~\pageref{page: preuve norme opP}.
\end{textesansboite}

\begin{lemmesansboite}{}
\label{norme opP}%
Under the stability condition of \cref{def: condition de stabilite}, we have 
\[\|\opP\|_{\mathcal{L}(L^2 (\SME),\mathbb{C})} \leqslant  C_{\Rsource} \, \varepsilon  \left( \frac{N}{\RME} \right)^{5/2} \domi{\ln \RME}^{3} .\]
\end{lemmesansboite}


\begin{textesansboite}
We are finally ready to prove \cref{cvg des ME}. We recall that $\errmdl[s][m+k]$ and $\errdtn[s]$ are defined in \cref{def: errmdl errdtn}.
\end{textesansboite}

\begin{demo}[Proof of \cref{cvg des ME}]
\phantomsection\label{preuve: cvg des ME}%
By~\cref{hyp: Cstab majore}, the stability condition implies $\varepsilon \leqslant C_{\Rsource} \.\RME$, so we can apply \cref{regularite elliptique pour ME 2} with $M:=C_{\Rsource}$. By~construction of $\emeg$ in \cref{eq: def emeg}, this gives 
\begin{align*}
	\|\emeg\|_{\HME{1/2}(\GME)} &\leqslant  C_{\Rsource} \,\varepsilon \. \bigg( \|\opD(\eme)\|_{\HME1(\GME)} + \|\opS(\emeg_{|\SME} )\|_{\HME2(\SME)} +  \frac{\|\emeg\|_{H^{1} (\OME)}}{\RME^{3/2}} \\
	&\quad {}+ \delta _{m,5} \.\sqrt{ \frac{\varepsilon }{\RME} } \.\snhg{\emeg} \bigg) \,+\, C_{\Rsource} \,\RME^{3/2} \, \big|\opS(\emeg_{|\SME} )(\PME)\big| \\[2pt]
	&\leqslant  C_{\Rsource} \,\varepsilon \. \bigg( \|\opD(\eme)\|_{\HME1(\GME)} + \|\opS(\emeg_{|\SME} )\|_{\HME2(\SME)} \\
	&\quad {}+ \frac{\|\emeg\|_{H^{1} (\OME)} + \delta _{m,5} \.\sqrt{\varepsilon } \.\snhg{\emeg}}{\RME^{3/2}} \bigg) \,+\, C_{\Rsource} \,\RME^{3/2} \, \big|\opP(\emeg_{|\SME} )\big| 
\end{align*}
Here $\delta _{m,5}$ indicates that the associated term should be ignored when $m=3$. By~\cref{erreur consistence ME,erreur ME v1}, this yields\vspace{-0.5ex}
\begin{align*}
	\|\emeg\|_{\HME{1/2}(\GME)} &\leqslant  C_{\Rsource} \,\varepsilon \. \left( C_{\fsource} \,\errmdl[2][m-1] + \big\|\opS(\emeg_{|\SME} )\big\|_{\HME2(\SME)} + C_{\fsource} \errmdl[1/2][m-1] \right) \\
	&\quad {}+ C_{\Rsource} \,\RME^{3/2} \, \big|\opP(\emeg_{|\SME} )\big| .
\end{align*}
Besides, we have
\begin{alignjustif}[\quad]
	\big\|\opS(\emeg_{|\SME} )\big\|_{\HME2(\SME)} &\leqslant  \|\opS\|_{\mathcal{L}(L^2 (\SME),\HME2(\SME))} \,\|\emeg\|_{L^2 (\SME)} \\
	&\leqslant  C_{\Rsource} \frac{N^{3} }{\RME} \cdot  \sqrt{\RME} \.\|\emeg\|_{H^{1} (\OME)} & \text{by \cref{opS - opS0,eq: norme L2 SME vs H1 OME}} \\
	&\leqslant  C_{\fsource} \.\errmdl[3/2][m-1] N^{3} & \text{by \cref{erreur ME v1}} . 
\end{alignjustif}
Similarly,
\[\big|\opP(\emeg_{|\SME} )\big| \leqslant  C_{\fsource} \,\varepsilon  \left( \frac{N}{\RME} \right)^{5/2} \domi{\ln \RME}^{3}  \cdot \sqrt{\RME} \,\errmdl[2][m-1] \leqslant  C_{\fsource} \,\errmdl \,\domi{\ln \RME}^{3} \,N^{5/2} .\]
Hence
\begin{align*}
	\|\emeg\|_{\HME{1/2}(\GME)} &\leqslant  C_{\fsource} \,\varepsilon  \left( \errmdl[2][m-1] + \errmdl[3/2][m-1] N^{3} + \errmdl[1/2][m-1] \right) + C_{\fsource} \,\errmdl[3/2] \,N^{5/2} \\
	&\leqslant  C_{\fsource} \errmdl[1/2] \,\domi{\ln \RME}^{3} \,N^{3} .
\end{align*}
Then \cref{ME avec Dirichlet sur Gamma} implies $\|\emeg\|_{H^{1} (\OME)} \leqslant  C_{\fsource} \errmdl \,\domi{\ln \RME}^{3} \,N^2 $. Next, using one more time \cref{erreur ME v1}, we obtain
\begin{align*}
	\|\eme\|_{H^{1} (\OME)} &\leqslant  \|\emeg\|_{H^{1} (\OME)} + \|\emes\|_{H^{1} (\OME)} \\
	&\leqslant  C_{\fsource} \errmdl \,\domi{\ln \RME}^{3} \,N^{3} + C_{\fsource,R} \,N \left( \errmdl[1/6] + \errdtn[1/6] \right) \\
	&\leqslant  C_{\fsource,R} \left( \errmdl \,\domi{\ln \RME}^{3} \,N^{3} + C_{R} \,\errdtn[1/6] N \right)
\end{align*}
By~\cref{def: errmdl errdtn}, this means that we have
\begin{align}
	\|\eme\|_{H^{1} (\OME)} &\leqslant  C_{\fsource,R} \left[ \left( \frac{\varepsilon \mathstrut}{\RME} \right)^{m} \domi{\ln \varepsilon }^{\ell _{m}} \cdot  \RME^{\delta } \.\domi{\ln \RME}^{2m+4} \cdot  N^{3} + \left( \frac{\RME}{R} \right)^{N\pst+ \frac{1}{6} } \domi{\ln \RME}^{2m+1} \cdot  N \right] \nonumber\\
	&\leqslant  C_{\fsource,R} \left[ \bigg( \frac{\varepsilon \mathstrut}{\RME} \bigg)^{m} \domi{\ln \varepsilon }^{\ell _{m}} \,N^{3} + \left( \frac{\RME}{R} \right)^{N\pst} \!N \right] \label{eq: cvg des ME preuve}
\end{align}
Let $R' :=(R+\Rsource)/2 \in (R,\Rsource)$. Applying \cref{eq: cvg des ME preuve} to $R'$ instead of $R$, and using the following bound
\[\left( \frac{\RME}{R'} \right)^{N\pst} \!N = \left( \frac{\RME}{R} \right)^{N\pst} \left( \frac{R}{R'} \right)^{N\pst} \!N \leqslant  C_{\Rsource,R} \left( \frac{\RME}{R} \right)^{N\pst} ,\]
we can eliminate $N$ from the right term of \cref{eq: cvg des ME preuve}, which gives:
\[\|\eme\|_{H^{1} (\OME)} \leqslant  C_{\fsource,R} \left[ \bigg( \frac{\varepsilon \mathstrut}{\RME} \bigg)^{m} \domi{\ln \varepsilon }^{\ell _{m}} \,N^{3} + \left( \frac{\RME}{R} \right)^{N\pst} \right] .\]
Finally, by \cref{DA de SE avec CL et delta variable}, $\|u_{\varepsilon } -\scl\|_{H^{1} (\OME)}$ is controlled by the same bound. We conclude by writing $\|u_{\varepsilon } -\sme\|_{H^{1} (\OME)} \leqslant  \|u_{\varepsilon } -\scl\|_{H^{1} (\OME)} +\|\eme\|_{H^{1} (\OME)}$.
\end{demo}

\section{Numerical results}
\label{sec: numerique ME}

\begin{textesansboite}
We have implemented our effective models and obtained numerical results that confirm the theoretical results and that are presented in this section.
\end{textesansboite}

\begin{textesansboite}
To begin with, a notable advantage of the effective models over the far fields to compute an approximation of $u_{\varepsilon }$ is that they are easier to discretise using finite elements. Indeed, we recall that, by \cref{construction directe des upl}, we have for any $(p,\ell )\in \NN\times \mathbb{N}$
\[u_{p,\ell } = \sum_{n=0}^\infty  \frac{\mu _{0} }{\mu _{1}} \,T_{n} \cdot  (\partial _{x} ^2 +k_{1} ^2 )^{n} \partial _{y} u_{p-1-2n,\ell } \qquad \text{ on }\Gamma \]
with $T_{n} \in \mathbb{R}$. This Dirichlet condition is difficult to take into account because it involves high order derivatives of the previous far fields. We handled it in \cref{sec: numerique des CL} by projecting the derivatives on the finite element space, which introduces additional numerical errors. On the contrary, the Robin and Ventcel conditions of the effective models add simple contributions to the variational formulations \cref{eq: FV ME3,eq: FV ME5} and are straightforward to take into account numerically.\\

The DtN condition $\partial _{r} \sme =\opS(\sme_{|\SME} )$ on $\SME$ is not difficult to include in the code as DtN conditions are quite classical and are available in some FEM libraries. The specificity here is the computation of the modes $\mme_{n}$ of the DtN. By~\cref{def: mme} and \cref{eq: def cme}, these modes depend on the functions $\mcle_{d,n}$ and on the coefficients $\cfinal$, which involve the matching coefficients $\cmacro$ and $\cmicro$ and the corner profile coefficients $\sigma _{q} (\ccDA_{d,n} )$, see \cref{eq: def cfinal}. The computation of the $\mcle_{d,n}$, $\cmacro$ and $\cmicro$ is explained in \cref{sec: implementation des espaces A}, while a method to compute the $\sigma _{q} (\ccDA_{d,n} )$ is presented in \cref{sec: calcul numerique des CC}.\\

The effective model of order 5 also involves a less classical non-local condition which is $\partial _{x} \sme5(\PME) =\opP5(\sme5_{|\SME} )$. It appears in the variational formulation \cref{eq: FV ME5} as a term \smash{$\opP5 (u_{|\SME} ) \cdot  \overline{v(\PME)}$} (times a coefficient), where $u\in \VME5$ denotes the unknown and $v\in \VME5$ is a test function. Using \cref{def: opS opP}, this can be rewritten as
\def\fctp{p^{(5)}_{\Nvirg\peps}}%
\begin{align}
	\opP5 (u_{|\SME} ) \cdot  \overline{v(\PME)} &= \Big[\langle u ,\mme_{n} \rangle _{L^2 (\SME)} \Big]_{n} ^{\mathsf{T}} \cdot \gram^{-1} \cdot  \Big[\partial _{r} \mme_{n} (\PME)\Big]_{n} \cdot  \overline{v(\PME)}\nonumber\\
	&= \int _{\SME} u(\bfx) \cdot  \fctp(\bfx) \,\mathrm{d}\bfx \cdot  \overline{v(\PME)}\label{eq: terme de opP dans FV}\\
	&\quad \ \text{with } \ \fctp: \bfx \in \SME \mapsto  \Big[\,\overline{\mme_{n} (\bfx)}\,\Big]_{n} ^{\mathsf{T}} \cdot \gram^{-1} \cdot  \Big[\partial _{r} \mme_{n} (\PME)\Big]_{n} \in \mathbb{C} \nonumber
\end{align}
This term is the product of a 1D integral of $u$ on a boundary and the evaluation of $v$ at a point, which is non-classical.\\

Besides, special care should be taken in the inversion of the Gram matrix $\gram$ when building $\opS$ and \smash{$\opP5$}, because this matrix can be ill-conditioned. Indeed, the proof of \cref{opS - opS0} shows that, when $\varepsilon $ is small, $\gram$ is close to $\gram0$. By~\cref{eq: mme qd eps=0}, $\gram0$ is a diagonal matrix and its diagonal coefficients are equal to $\|\mme0_{n} \|_{L^2 (\SME)} ^2  = \|\mcll_{n\pst,0} \|_{L^2 (\SME)} ^2 $ ($n\in [\![1,N]\!]$). Keeping only the first power of $r$ in $\mcll_{n\pst,0} $, we have roughly $\mcll_{n\pst,0} (r,\theta ) \approx  \phi_{d} (r,\theta ) = \pm  r^{n\pst} \sin (n\pst \theta )$ near the corner, so the diagonal coefficients of $\gram0$ have a geometric behavior w.r.t.\ $n$. This explains the ill-conditioning of $\gram$, which we observed numerically. In some cases, the ratio between the largest and the smallest diagonal coefficient of $\gram$ was around $10^{15}$. We handled this by computing $\gram^{-1}$ using a preconditioner:
\[\gram^{-1} = D \, (D\, \gram\, D)^{-1} D \qquad  \text{with }  D:= \mathrm{Diag} \big(\|\mme_{1} \|_{L^2 (\SME)} ^{-1} ,\dots ,\|\mme_{N} \|_{L^2 (\SME)} ^{-1} \big).\]
The matrix that is numerically inverted is
\[D\, \gram\, D = \bigg[ \bigg\langle  \frac{\mme_{n} }{\|\mme_{n} \|_{L^2 (\SME)} }, \frac{\mme_{n'} }{\|\mme_{n'} \|_{L^2 (\SME)} } \bigg\rangle  \bigg]_{(n,n')\in [\![1,N]\!]^2 },\]
i.e.\ the Gram matrix of the normalized modes. This approach is equivalent to replacing the modes $\mme_{n}$ with $\mme_{n} /\|\mme_{n} \|_{L^2 (\SME)}$ in formula of $\opS$ in \cref{def: opS opP}.
\end{textesansboite}

\begin{textesansboite}
The numerical experiments are performed for the same parameters as the computation of the far fields in \cref{sec: numerique des CL}: we consider the same values of $\Theta $, $\Omega _{1}$ (described in \cref{fig: forme couche simus}), the source term $\fsource$, and $\mu _{0}$, $\mu _{1}$, $\rho _{0}$, $\rho _{1}$, $\omega $, $\RME$ and $\varepsilon $. In addition, the number of modes is taken as
\begin{equation}{\label{eq: N des ME numerique}}
	\textstyle N :=\lceil {-m \frac{\Theta }{\pi } \log _{2} (\varepsilon ) +c}\rceil 
\end{equation}
where $c\in \mathbb{R}$ is fixed, because this choice yields $\|\sme-u_{\varepsilon } \|_{H^{1} (\OME)} = O(\varepsilon ^{m} \domi{\ln \varepsilon }^{\ell _{m} +2} )$ as explained in \cref{rq: choix parametres ME}. The constant $c$ is chosen depending on $(\Theta ,m)$ so that we have $N=4$ for the largest considered value of $\varepsilon $. With this parameters, we have $4\leqslant N\leqslant 12$ in any of our numerical computations. This may seem small, but we noticed that the error was larger with greater values of $N$. This is coherent with \cref{cvg des ME} which states an error estimate that contains a factor $N^{3}$ that increases w.r.t.\ $N$ (although $N^{3}$ is not optimal, see \cref{rq: commentaire estimee de cvg des ME}). Thus $N$ cannot be taken too large, which a counterintuitive behavior for a DtN condition. We expect to be able to take $N$ larger for smaller $\varepsilon $ according to \cref{eq: N des ME numerique}, but for small $\varepsilon $ the computation of the reference solution $u_{\varepsilon }$ is more costly as it requires finely meshing the layer.\\

Moreover, we use P2 Lagrange finite elements, because we noticed that P1 elements lead to larger errors for the model of order 5 (but we observed no difference between P2, P3 and P4). We conjecture that this is related to the non-classical DtN condition $\partial _{x} \sme5(\PME) =\opP5(\sme5_{|\SME} )$, but we have not investigated this question.\\

We use again the FEM-BEM library XLiFE++ \cite{XLiFE}. I thank Eric Luneville for adding a feature to XLiFE++ that allowed me to handle the non-classical term \cref{eq: terme de opP dans FV} in the variational formulation of the effective model of order 5.
\end{textesansboite}

\begin{textesansboite}
\Crefrange{tab: ME angle 2.5}{tab: ME angle 4} show the real parts of the computed solutions and the associated difference with $u_{\varepsilon }$. We represent $\sme3$ and $\sme5$, but also the first far field $u_{0,0}$ for comparison because it~is the solution of the limit problem (see \cref{ex: eq u00 et u10}), which is an effective model of order 1. Besides, we recall that the considered angles, namely $\Theta  \in \{2.5,\pi ,4\}$, are a representative sample as this include small and large angles, as well as the two cases $\pst \in \mathbb{Q}$ and $\pst \not\in \mathbb{Q}$.\\

Qualitatively, the $\sme3$ and $\sme5$ are very similar to $u_{\varepsilon }$ and, denoting $\sme1 := u_{0,0}$, the order of magnitude of \smash{$\sme-u_{\varepsilon }$} decreases with $m$. In addition, one can see that the error of the effective model of order 5 is mainly localized near the DtN boundary $\SME$. We noticed that this error was greater for a larger $N$ and we could~not really take a smaller $N$ as we have $N=4$ for the largest considered value of $\varepsilon $, which is already very small.\\

Moreover, \crefrange{fig: graphe err ME 2.5}{fig: graphe err ME 4} plot the relative error made by $\sme$ for $m\in \{1,3,5\}$ w.r.t.\ $\varepsilon $.
We use the $H^{1}$~norm that we normalize by $\|u_{0,0} \|_{H^{1}}$, like we did for the far fields in \cref{sec: numerique des CL}. For each $\Theta $, we provide two graphs: one that plots the errors on $\OME$ and one for the errors on a subdomain $\OME \cap  \{r>1\}$ that excludes a vicinity of the DtN boundary $\SME$ (as here $\RME \leqslant  1$). This allows to ignore the part of the error that is localized near $\SME$.\\

We observe that $u_{0,0} -u_{\varepsilon }$ and $\sme3-u_{\varepsilon }$ converge quite precisely at the expected orders (resp.\ 1 and 3). In addition, the error plot of the effective model of order 5 for $\Theta =\pi $ (\cref{fig: graphe err ME pi}) is a clear line that reaches a plateau that seems linked to the discretisation errors. The slope of this line is 4.7 on $\OME$, so a little less than 5, and more than 5 on $\OME \cap  \{r>1\}$. For $\Theta =4$, the points are a bit less aligned, but we also obtain a slope of almost 5 (4.7 to 4.8), see \cref{fig: graphe err ME 4}. As~for $\Theta =2.5$, we are not able to observe a convergence at order 5, see \cref{fig: graphe err ME 2.5}. In this case the error is especially localized near $\SME$, so the issue may be linked to the choice of the DtN's parameters $(N,\RME)$. These parameters are not obvious to choose (see \cref{rq: choix parametres ME}) and we did~not have time to investigate this question in more details. Note however that the convergence of naive effective models (that consider the effective boundary conditions on $\Gamma $ entirely with no special treatement of the corner) is less deteriorated compared to the infinite layer for small angles, see \cite{CalCosDauVia06}, \cref{sec: etat de l'art,chap: conclusion}. So our effective models are more relevant for larger angles like $\Theta \in \{\pi ,4\}$.\\

Therefore, our effective models are validated numerically in most cases, except the model of order 5 for $\Theta =2.5$ which requires more investigations on the choice of the parameters. In addition, the smallest obtained relative errors are around $10^{-4}$, which is satisfying.\\
\end{textesansboite}

	\tikzset{external/export=true}%
	\tikzsetfigurename{tab-ME_2,5-}%
	\begin{table}[h]
\def\largeurimages{47mm}
\def\boite#1{\raisebox{-.5\height}{\hskip1mm \vbox{\vskip3.5mm\relax \hbox{#1} \vskip4mm}\!}}%
\def\addgraphe#1{\boite{\includegraphics[width=\largeurimages]{images/angle2,5_#1.jpg}}}%
\def\addgrapheE#1{\addgraphe{#1_eps0,01}}%
\def\colorbarSE{\colorbar[min=-2.3, max=1.8, 1e-2]}
\def\colorbarheight{20mm}
\renewcommand{\arraystretch}{1.8}%

\def\imageSE{\boite{%
    \dimen0=\largeurimages\relax
    \dimen0=1.39\dimen0\relax
    \edef\tempo{\the\dimen0}%
    \pgfmathsetmacro{\rapp}{1/(1-cos(2.5*180/3.14159))}
    \begin{tikzpicture}[x=\largeurimages, y=\largeurimages]
        \clip (-0.1,{-\rapp*0.1/4}) -- (\rapp,{-\rapp*0.1/4}) -- (\rapp,0) arc (0:{2.5*180/3.14}:\rapp) -- cycle;
        \draw (0.075,0.378) node{\includegraphics[width=\tempo]{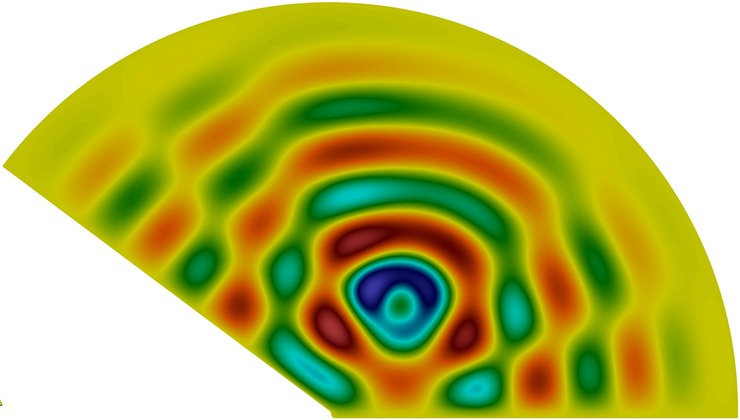}};
    \end{tikzpicture}%
}}

\caption{Real part of the effective model solutions and the errors w.r.t.\ the exact solution for $\Theta=2.5$ and $\peps=0.1$}
\label{tab: ME angle 2.5}

\makebox[\textwidth][c]{
\begin{tabular}{|c|cl<{\!}|cl<{\!}|}
    \cline{2-5}
    \multicolumn{1}{c|}{} &
    \multicolumn{2}{c|}{Computed solution} & 
    \multicolumn{2}{c|}{Difference with $u_\peps$} \\
    \hline
    $u_\peps$ &
    \imageSE & \colorbarSE &
    \multicolumn{2}{c|}{$\underset{\textstyle\diagup}\ $} \\[-3ex]
    & \multicolumn{2}{c|}{(exact solution)} & & \\
    \hline
    $\.u_{0,0}\.$ &
    \addgraphe{CL1_p0_l0} & \colorbarSE &
    \addgrapheE{errSE-CL1_p0} & \colorbar[min=-4, max=3.3, 1e-3] \\
    \hline
    $\sme3$ &
    \addgrapheE{ME3} & \colorbarSE &
    \addgrapheE{errSE-ME3} & \colorbar[min=-12, max=6.7, 1e-4] \\
    \hline
    $\sme5$ &
    \addgrapheE{ME5} & \colorbarSE &
    \addgrapheE{errSE-ME5} & \colorbar[min=-7, max=2.5, 1e-4] \\
    \hline
\end{tabular}}
\end{table}%
	\tikzset{external/export=false}%
	\tikzsetfigurename{main-figure}
	\tikzset{external/export=true}%
	\tikzsetfigurename{tab-ME_pi-}%
	\begin{table}[h]
\def\largeurimages{44mm}
\def\boite#1{\raisebox{-.5\height}{\hskip1.5mm \vbox{\vskip2.5mm\relax \hbox{#1} \vskip3mm}\!}}%
\def\addgraphe#1{\boite{\includegraphics[width=\largeurimages]{images/anglepi_#1.jpg}}}%
\def\addgrapheE#1{\addgraphe{#1_eps0,01}}%
\def\colorbarSE{\colorbar[min=-2.1, max=1.4, 1e-2]}
\def\colorbarheight{20mm}
\renewcommand{\arraystretch}{1.8}%

\def\imageSE{\boite{%
    \dimen0=\largeurimages\relax
    \dimen0=1.375\dimen0\relax
    \edef\tempo{\the\dimen0}%
    \begin{tikzpicture}[x=\largeurimages, y=\largeurimages]
        \clip (-0.5,{-0.5*0.1/4}) -- (0.5,{-0.5*0.1/4}) -- (0.5,0) arc (0:180:0.5) -- cycle;
        \draw (0,0.335) node{\includegraphics[width=\tempo]{images/anglepi_SE_eps0,1.jpg}};
    \end{tikzpicture}%
}}

\caption{Real part of the effective model solutions and the errors w.r.t.\ the exact solution for $\Theta=\pi$ and $\peps=0.1$}
\label{tab: ME angle pi}

\makebox[\textwidth][c]{
\begin{tabular}{|c|cl<{\!}|cl<{\!}|}
    \cline{2-5}
    \multicolumn{1}{c|}{} &
    \multicolumn{2}{c|}{Computed solution} & 
    \multicolumn{2}{c|}{Difference with $u_\peps$} \\
    \hline
    $u_\peps$ &
    \imageSE & \colorbarSE&
    \multicolumn{2}{c|}{$\underset{\textstyle\diagup}\ $} \\[-2.5ex]
    & \multicolumn{2}{c|}{(exact solution)} & & \\
    \hline
    $\.u_{0,0}\.$ &
    \addgraphe{CL1_p0_l0} & \colorbarSE &
    \addgrapheE{errSE-CL1_p0} & \colorbar[min=-3.9, max=3.1, 1e-3] \\
    \hline
    $\sme3$ &
    \addgrapheE{ME3} & \colorbarSE &
    \addgrapheE{errSE-ME3} & \colorbar[min=-3.4, max=4.7, 1e-4] \\
    \hline
    $\sme5$ &
    \addgrapheE{ME5} & \colorbarSE &
    \addgrapheE{errSE-ME5} & \colorbar[min=-5, max=10, 1e-5] \\
    \hline
\end{tabular}}
\end{table}%
	\tikzset{external/export=false}%
	\tikzsetfigurename{main-figure}
	\tikzset{external/export=true}%
	\tikzsetfigurename{tab-ME_4-}%
	\begin{table}[h]
\def\largeurimages{44mm}
\def\boite#1{\raisebox{-.5\height}{\hskip2mm \vbox{\vskip2.5mm\relax \hbox{#1} \vskip2mm}\!}}%
\def\addgraphe#1{\boite{\includegraphics[width=\largeurimages]{images/angle4_#1.jpg}}}%
\def\addgrapheE#1{\addgraphe{#1_eps0,01}}%
\def\colorbarSE{\colorbar[min=-2, max=1.4, 1e-2]}
\def\colorbarheight{20mm}
\renewcommand{\arraystretch}{1.8}%

\def\imageSE{\boite{%
    \dimen0=\largeurimages\relax
    \dimen0=1.375\dimen0\relax
    \edef\tempo{\the\dimen0}%
    \begin{tikzpicture}[x=\largeurimages, y=\largeurimages]
        \clip (0.5,-0.1) -- (0.5,0) arc (0:{4*180/3.14}:0.5) -- cycle;
        \draw (0,0.085) node{\includegraphics[width=\tempo]{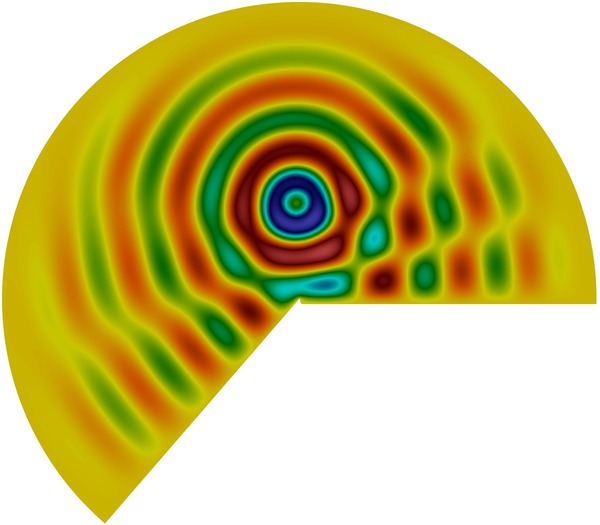}};
    \end{tikzpicture}%
}}

\caption{Real part of the effective model solutions and the errors w.r.t.\ the exact solution for $\Theta=4$ and $\peps=0.1$}
\label{tab: ME angle 4}

\makebox[\textwidth][c]{
\begin{tabular}{|c|cl<{\!}|cl<{\!}|}
    \cline{2-5}
    \multicolumn{1}{c|}{} &
    \multicolumn{2}{c|}{Computed solution} & 
    \multicolumn{2}{c|}{Difference with $u_\peps$} \\
    \hline
    $u_\peps$ &
    \imageSE & \colorbarSE&
    \multicolumn{2}{c|}{$\diagup$} \\[-3.5ex]
    & \multicolumn{2}{c|}{(exact solution)} & & \\
    \hline
    $\.u_{0,0}\.$ &
    \addgraphe{CL1_p0_l0} & \colorbarSE &
    \addgrapheE{errSE-CL1_p0} & \colorbar[min=-3.5, max=2.7, 1e-3] \\
    \hline
    $\sme3$ &
    \addgrapheE{ME3} & \colorbarSE &
    \addgrapheE{errSE-ME3} & \colorbar[min=-2.9, max=3.4, 1e-4] \\
    \hline
    $\sme5$ &
    \addgrapheE{ME5} & \colorbarSE &
    \addgrapheE{errSE-ME5} & \colorbar[min=-1.5, max=2.5, 1e-4] \\
    \hline
\end{tabular}}
\end{table}%
	\tikzset{external/export=false}%
	\tikzsetfigurename{main-figure}

	\tikzset{external/export=true}%
	\tikzsetfigurename{fig-erreurs_ME_2,5-}%
	\pgfplotsset{
    every axis/.append style={
        thick,
    },
    cycle list name=Dark2_avec_regressions,
    /pgfplots/table/x=epsilon,
    /pgfplots/table/y=nbCC,
    legend cell align=left,
}

\begin{figure}
\begin{tikzpicture}
\begin{axisME}{err_2_tot.dat}{$\OME$}
\courbeME[y=errH1_CL_1, fct=$u_{0,0}$, reg=0:5]
\courbeME[y=errH1_ME_3, fct=$\sme3$, reg=0:5]
\courbeME[y=errH1_ME_5, fct=$\sme5$, reg=0:5]
\end{axisME}
\end{tikzpicture}
\vspace{-1ex}

\begin{tikzpicture}
\begin{axisME}{err_2_loin.dat}{$\OME\cap\{r>1\}$}
\courbeME[y=errH1_CL_1, fct=$u_{0,0}$, reg=0:5]
\courbeME[y=errH1_ME_3, fct=$\sme3$, reg=0:5]
\courbeME[y=errH1_ME_5, fct=$\sme5$, reg=0:5]
\end{axisME}
\end{tikzpicture}
\caption{
    Numerical errors of the effective models for $\Theta=2.5$
}
\label{fig: graphe err ME 2.5}
\end{figure}%
	\tikzset{external/export=false}%
	\tikzsetfigurename{main-figure}
	\tikzset{external/export=true}%
	\tikzsetfigurename{fig-erreurs_ME_pi-}%
	\pgfplotsset{
    every axis/.append style={
        thick,
    },
    cycle list name=Dark2_avec_regressions,
    /pgfplots/table/x=epsilon,
    /pgfplots/table/y=nbCC,
    legend cell align=left,
}

\begin{figure}
\begin{tikzpicture}
\begin{axisME}{err_pi_tot.dat}{$\OME$}
\courbeME[y=errH1_CL_1, fct=$u_{0,0}$, reg=0:5]
\courbeME[y=errH1_ME_3, fct=$\sme3$, reg=0:5]
\courbeME[y=errH1_ME_5, fct=$\sme5$, reg=1:5]
\end{axisME}
\end{tikzpicture}
\vspace{-1ex}

\begin{tikzpicture}
\begin{axisME}{err_pi_loin.dat}{$\OME\cap\{r>1\}$}
\courbeME[y=errH1_CL_1, fct=$u_{0,0}$, reg=0:5]
\courbeME[y=errH1_ME_3, fct=$\sme3$, reg=0:5]
\courbeME[y=errH1_ME_5, fct=$\sme5$, reg=1:5]
\end{axisME}
\end{tikzpicture}
\caption{
    Numerical errors of the effective models for $\Theta=\pi$
}
\label{fig: graphe err ME pi}
\end{figure}%
	\tikzset{external/export=false}%
	\tikzsetfigurename{main-figure}
	\tikzset{external/export=true}%
	\tikzsetfigurename{fig-erreurs_ME_4-}%
	\pgfplotsset{
    every axis/.append style={
        thick,
    },
    cycle list name=Dark2_avec_regressions,
    /pgfplots/table/x=epsilon,
    /pgfplots/table/y=nbCC,
    legend cell align=left,
}

\begin{figure}
\begin{tikzpicture}
\begin{axisME}{err_4_tot.dat}{$\OME$}
\courbeME[y=errH1_CL_1, fct=$u_{0,0}$, reg=0:5]
\courbeME[y=errH1_ME_3, fct=$\sme3$, reg=0:5]
\courbeME[y=errH1_ME_5, fct=$\sme5$, reg=1:5]
\end{axisME}
\end{tikzpicture}
\vspace{-1ex}

\begin{tikzpicture}
\begin{axisME}{err_4_loin.dat}{$\OME\cap\{r>1\}$}
\courbeME[y=errH1_CL_1, fct=$u_{0,0}$, reg=0:5]
\courbeME[y=errH1_ME_3, fct=$\sme3$, reg=0:5]
\courbeME[y=errH1_ME_5, fct=$\sme5$, reg=1:5]
\end{axisME}
\end{tikzpicture}
\caption{
    Numerical errors of the effective models for $\Theta=4$
}
\label{fig: graphe err ME 4}
\end{figure}%
	\tikzset{external/export=false}%
	\tikzsetfigurename{main-figure}
\begin{subappendices}
	
\section{Appendix: proofs of the technical lemmas of \texorpdfstring{\cref{sec: lemmes de regularite elliptique ME}}{Section 9.5.2}}
\label{sec: preuves lemmes de regularite elliptique ME}

\begin{textesansboite}
This section contains the proofs of \cref{regularite elliptique pour ME,diffeo ME,regularite elliptique pour ME 2}.
\end{textesansboite}

\begin{demo}[Proof of \cref{regularite elliptique pour ME}]
\phantomsection\label{preuve: regularite elliptique pour ME}%
\def\recI{R}
\def\recJ{R'}
\def\J{J}
\def\ukhi{\tilde u}%
\def\fkhi{\tilde f}%
\def\gkhi{\tilde g}%
\def\dpl{\mathscr P}
\def\Lip{\operatorname{Lip}}%
The proof has two steps. First we prove that $u$ is $H^2 $ in $\recJ:=\J\times (0, \frac{\eta }{2} )$ and we bound $\|u\|_{H^2 (\recJ)}$ using finite differences. Then we rewrite the boundary condition as $\nu _{1} u = \varepsilon \,(g-\partial _{n} u+\dots )$ to deduce $\|u\|_{H^{1/2} (\J_{0} )} =O(\varepsilon )$. The first step is the longest part.\medskip

\underline{Step 1:} We focus on the second point (with Ventcel's condition), as the first one is simpler. We will explain in the end what should be modified for Robin's condition. Let $\dpl:=\mathbb{R}\times \mathbb{R}_+^*$. Let $\chi \in \mathcal{C}^2 (\overline{\dpl})$ be s.t.\ $\chi =1$ in $\recJ$, $\chi =0$ outside $\recI$, and $\chi $ and its derivatives are bounded. We also assume that $\partial _{y} \chi _{|\partial \dpl} =0$ (such a function $\chi $ is not difficult to construct). Let $\ukhi :\dpl\rightarrow \mathbb{C}$ be the extension of $\chi u$ by 0 outside $\recI$. We clearly have $\ukhi \in H^{1} (\dpl)$ and $\ukhi$ satisfies
\[\left\{\begin{array}{r@{}l@{}ll}
	\operatorname{div}(A \nabla \ukhi) &{}= \chi  f+ 2A \nabla  \chi \cdot \nabla  u+\operatorname{div}(A \nabla  \chi )\,u &{} =:\fkhi &\text{ in }\dpl \\[.5ex]
	\partial _{n} \ukhi + \frac{\nu _{1} }{\varepsilon } \ukhi - \varepsilon \,\nu _{2} \partial _{x} ^2 \ukhi &{}=\chi  g-\cancel{\partial _{y} \chi }\,u +2\varepsilon  \nu _{2} \,\partial _{x} \chi  \,\partial _{x} u +\varepsilon  \nu _{2} \partial _{x} ^2 \chi  \,u &{} = \chi g+\varepsilon \gkhi &\text{ on }\partial \dpl=\mathbb{R}\times \{0\}
\end{array}\right.\]
with $\fkhi \in L^2 (\dpl)$, $\chi g \in H^{1} (\partial \dpl)$ and $\gkhi:=2\nu _{2} \partial _{x} \chi  \,\partial _{x} u +\nu _{2} \partial _{x} ^2 \chi  \,u \in L^2 (\partial \dpl)$. Since $A$ is the identity matrix on $I_{0}$, the associated variational formulation is:
\begin{equation}{\label{eq: regularite elliptique pour ME preuve FV}}
	\int _{\dpl} A \nabla  \ukhi\cdot \nabla  \bar{v}+ \int _{\partial \dpl} \left(\frac{\nu _{1} }{\varepsilon } \,\ukhi \bar{v}+\varepsilon  \,\nu _{2} \,\partial _{x} \ukhi \,\partial _{x} \bar{v}\right) = -\int _{\dpl} \fkhi \bar{v}+\int _{\partial \dpl} (\chi g+\varepsilon \gkhi) \bar{v}.
\end{equation}
Let us denote $D_{h} \varphi :(x,y)\mapsto \frac{\varphi (x+h,y)-\varphi (x,y)}{h}$ for~any $h\in \mathbb{R}^*$ and $\varphi :\dpl\rightarrow \mathbb{C}$ (and similarly for $\varphi :\partial \dpl\rightarrow \mathbb{C}$). Taking $v:=D_{-h} D_{h} \ukhi$ in \cref{eq: regularite elliptique pour ME preuve FV} and doing discrete integrations by parts, we obtain:
\begin{multline*}
	\int _{\dpl} D_{h} (A\nabla  \ukhi)\cdot D_{h} \nabla  \bar\ukhi + \int _{\partial \dpl} \left( \frac{\nu _{1} }{\varepsilon } |D_{h} \ukhi|^2  +\varepsilon  \nu _{2} |D_{h} \partial _{x} \ukhi|^2  \right) \\
	= \int _{\dpl} \fkhi \cdot D_{-h} D_{h} \bar\ukhi -\int _{\partial \dpl} D_{h} (\chi g) \cdot  D_{h} \bar\ukhi -\int _{\partial \dpl} \varepsilon \gkhi \cdot D_{-h} D_{h} \bar\ukhi .
\end{multline*}
Now we have $D_{h} (A \varphi )=D_{h} A \cdot \varphi +A \cdot D_{h} \varphi $ with $\varphi :=\nabla \ukhi$. Passing the term $D_{h} A \cdot \varphi $ to the right-hand side and using $\|D_{h} A\|_{L^{\infty }} \leqslant \Lip(A)$, we get:
\begin{multline*}
	\int _{\dpl} A D_{h} \nabla  \ukhi \cdot  D_{h} \nabla  \bar\ukhi + \int _{\partial \dpl} \left( \frac{\nu _{1} }{\varepsilon } |D_{h} \ukhi |^2  +\varepsilon  \nu _{2} |D_{h} \partial _{x} \ukhi|^2  \right) \\
	\leqslant  \|\fkhi\|_{L^2 (\dpl)} \|D_{-h} D_{h} \ukhi\|_{L^2 (\dpl)} +\|D_{h} (\chi g)\|_{L^2 (\partial \dpl)} \|D_{h} \ukhi\|_{L^2 (\partial \dpl)} \\
	+\varepsilon  \.\|\gkhi\|_{L^2 (\partial \dpl)} \|D_{-h} D_{h} \ukhi\|_{L^2 (\partial \dpl)} + \Lip(A) \|\nabla \ukhi\|_{L^2 (\dpl)} \|D_{h} \nabla  \ukhi\|_{L^2 (\dpl)} .
\end{multline*}
Then we use the ellipticity assumption in the left-hand side, and we bound some $\|D_{\pm h} \cdot \|_{L^2 }$ by $\|\partial _{x} \cdot \|_{L^2 }$ in the right-hand side. This gives:
\begin{multline*}
	S:=\|D_{h} \nabla  \ukhi\|_{L^2 (\dpl)} ^2  + \frac{1}{\varepsilon } \|D_{h} \ukhi \|_{L^2 (\partial \dpl)} ^2  +\varepsilon  \.\|D_{h} \partial _{x} \ukhi\|_{L^2 (\partial \dpl)} ^2  \\
	\begin{array}{l@{}l@{}}
		\leqslant  C_{A,\nu _{1} ,\nu _{2} } \Big( & \|\fkhi\|_{L^2 (\dpl)} \|D_{h} \nabla \ukhi\|_{L^2 (\dpl)} +\|\partial _{x} (\chi g)\|_{L^2 (\partial \dpl)} \|D_{h} \ukhi\|_{L^2 (\partial \dpl)} \\[2\jot]
		& {} + \varepsilon  \.\|\gkhi\|_{L^2 (\partial \dpl)} \|D_{h} \partial _{x} \ukhi\|_{L^2 (\partial \dpl)}
		+ \|\nabla \ukhi\|_{L^2 (\dpl)} \|D_{h} \nabla \ukhi\|_{L^2 (\dpl)} \Big)\!\! \\[2\jot]
		\underset{\mathclap{\rm C.S.}}\leqslant \displaystyle C_{A,\nu _{1} ,\nu _{2} } \big(& \|\fkhi\|_{L^2 (\dpl)} ^2  +\varepsilon  \.\|\partial _{x} (\chi g)\|_{L^2 (\partial \dpl)} + \varepsilon  \.\|\gkhi\|_{L^2 (\partial \dpl)} ^2 + \|\nabla \ukhi\|_{L^2 (\dpl)} ^2  \big)^{1/2} \cdot S^{1/2} .
	\end{array}
\end{multline*}
Now, we simplify the factor $S^{1/2}$ on both sides, we take the limit $h\rightarrow 0$, and we use the estimates $\|\fkhi\|_{L^2 } \leqslant C_{\chi ,A} (\|f\|_{L^2 } +\|u\|_{H^{1} (\recI)})$, $\|\partial _{x} (\chi g)\|_{L^2 } \leqslant C_{\chi } \|g\|_{H^{1}}$, $\|\gkhi\|_{L^2 } \leqslant C_{\chi ,\nu _{2}} (\snh{u}(I_{0} )+\|u\|_{H^{1} (\recI)})$ and $\|\nabla  \ukhi\|_{L^2 (\dpl)} \leqslant C_{\chi } \|u\|_{H^{1} (\recI)}$. This gives:
\begin{align*}
	\|\partial _{x} \nabla  \ukhi\|_{L^2 (\dpl)} ^2  &\leqslant  \limsup\limits _{h\rightarrow 0} \|D_{h} \nabla  \ukhi\|_{L^2 (\dpl)} ^2  \\
	&\leqslant  C_{A,\nu _{1} ,\nu _{2} } \big(\|\fkhi\|_{L^2 (\dpl)} ^2  +\varepsilon \.\|\partial _{x} (\chi g)\|_{L^2 (\partial \dpl)} ^2 +\varepsilon \.\|\gkhi\|_{L^2 (\partial \dpl)} ^2  + \|\nabla \ukhi\|_{L^2 (\dpl)} ^2  \big) \\[\jot]
	&\leqslant  C_{I,J,\eta ,A,\nu _{1} ,\nu _{2} ,\varepsilon _{0} } \big(\|f\|_{L^2 (\recI)} ^2 + \varepsilon \.\|g\|_{H^{1} (I_{0} )} ^2   +\|u\|_{H^{1} (\recI)} ^2  + \varepsilon \.\snh{u}(I_{0} )^2  \big).
\end{align*}
Since $\ukhi_{|\recJ} =u_{|\recJ}$ and $\partial _{y} ^2 u = f-\partial _{x} ^2 u$ in $\recJ$, we deduce in particular $u_{|\recJ}\in H^2 $ and
\begin{equation}{\label{eq: regularite elliptique pour ME preuve H2 recJ}}
	\|u\|_{H^2 (\recJ)} \leqslant  C_{I,\J,\eta ,A,\nu _{1} ,\nu _{2} ,\varepsilon _{0} } \big(\|f\|_{L^2 (\recI)} + \sqrt{\varepsilon } \,\|g\|_{H^{1} (I_{0} )}  +\|u\|_{H^{1} (\recI)} + \sqrt{\varepsilon } \,\snh{u}(I_{0} ) \big).
\end{equation}
Similarly, in the Robin case (i.e.\ when $\partial _{n} u + \frac{\nu _{1} }{\varepsilon } u=g$ on $I_{0}$), one can prove $\|u\|_{H^2 (\recJ)} \leqslant  C_{I,\J,\eta ,A,\nu _{1} } \*\big(\|f\|_{L^2 (\recI)} + \sqrt{\varepsilon } \.\|g\|_{H^{1} (I_{0} )} +\|u\|_{H^{1} (\recI)} \big)$. It suffices to replace $\nu _{2}$ and $\gkhi$ by 0 and to ignore the terms that involve $\|D_{h} \partial _{x} \ukhi\|_{L^2 (\partial \dpl)}$.\medskip\smallskip

\underline{Step 2:} For the Robin case, the estimate on $\|u\|_{H^{1/2} (\J_{0} )}$ is a direct consequence of $u_{|\J_{0}} = \frac{\varepsilon }{\nu _{1}} (g-\partial _{n} u)$, the inequality $\|\partial _{n} u\|_{H^{1/2} (\J_{0} )} \leqslant  C_{\J,\eta } \|u\|_{H^2 (\recJ)}$ and step 1. The Ventcel case is less obvious because we have $u_{|\J_{0}} = \frac{\varepsilon }{\nu _{1}} (g-\partial _{n} u +\nu _{2} \partial _{x} ^2 u)$ with no control over $\|\partial _{x} ^2 u\|_{H^{1/2} (\J_{0} )}$ for now. To conclude in this case, we will do one more step of elliptic regularity on $\partial \dpl$.\\
\def\ukhi{w}%
Let $K$ be an open interval of $\mathbb{R}$ s.t.\ $\overline{J}\subset K$ and $\overline{K}\subset I$. We denote $K_{0} :=K\times \{0\}$. Let $\chi  \in  \mathcal{C}^2 (\partial \dpl)$ be a new truncation function s.t.\ $\chi _{|\J_{0}} =1$ and $\chi _{|\mathbb{R} \setminus K_{0} }=0$. Let $\ukhi:=\chi  u_{|K_{0}}$ extended by 0 outside $K_{0}$. The boundary condition of $u$ on $I_{0}$ implies:
\begin{equation}{\label{eq: regularite elliptique pour ME preuve pb sur ddpl}}
	-\varepsilon ^2 \nu _{2} \partial _{x} ^2 \ukhi + \nu _{1} \ukhi = \varepsilon  \chi  g-\varepsilon  \chi  \partial _{n} u -\varepsilon ^2 \nu _{2} (2 \partial _{x} \chi  \,\partial _{x} u +\partial _{x} ^2 \chi  \,u ) =: h \quad \text{ on }\partial \dpl.
\end{equation}
\Cref{eq: regularite elliptique pour ME preuve H2 recJ} of step 1 applied to $K$ instead of $J$ implies
\begin{align*}
	\|h\|_{H^{1/2} (\partial \dpl)} &\leqslant  C_{\chi ,\nu _{2}} \,\varepsilon \, \big(\|g\|_{H^{1/2} (I_{0} )} + \|\partial _{n} u\|_{H^{1/2} (K_{0} )} + \varepsilon \.\|u\|_{H^{3/2} (K_{0} )} \big)\\
	&\leqslant  C_{\J,K,\eta ,\nu _{2} ,\varepsilon _{0} } \,\varepsilon \, \big(\|g\|_{H^{1} (I_{0} )} + \|u\|_{H^2 (K\times (0,\eta /2))} \big)\\
	&\leqslant  C_{I,\J,\eta ,A,\nu _{1} ,\nu _{2} ,\varepsilon _{0} } \,\varepsilon \, \big(\|f\|_{L^2 (\recI)} + \|g\|_{H^{1} (I_{0} )}  +\|u\|_{H^{1} (\recI)} + \sqrt{\varepsilon } \,\snh{u}(I_{0} ) \big).
\end{align*}
Applying the one-dimensional Fourier transform on $\partial \dpl$ to \cref{eq: regularite elliptique pour ME preuve pb sur ddpl} gives $(\varepsilon ^2 \nu _{2} \xi ^2 +\nu _{1} )\hat\ukhi= \hat{h}$, so $\nu _{1} |\hat\ukhi|\leqslant | \hat{h}|$. Therefore\vspace{-2pt}
\[\smash{\|u\|_{H^{1/2} (\J_{0} )} \leqslant  \|\ukhi\|_{H^{1/2} (\partial \dpl)} = \big\| (1+\xi ^2 )\hat\ukhi \big\|_{L^2 (\mathbb{R})} \leqslant  \frac{1}{\nu _{1}} \.\big\| (1+\xi ^2 )\hat h \big\|_{L^2 (\mathbb{R})} = \frac{1}{\nu _{1}} \.\|h\|_{H^{1/2} (\partial \dpl)} .}\]
This concludes the proof.
\end{demo}

\begin{demo}[Proof of \cref{diffeo ME}]
\phantomsection\label{preuve: diffeo ME}%
\def\pol{\phi}%
\def\dompol{E}%
\def\I{(-\RME,\infty)}%
\def\J{I_\Theta}%
\def\dompsi{\I\times\J}%
\def\Opsi{\Omega'}%
\def\Gpsi{\Gamma'}%
\def\Spsi{\Sigma'}%
\def\psit{\psi_2}%
Let $(r,\theta )$ denote the polar coordinates with $\theta  \in (-\pi , \pi )$ and let $\rho :=r-\RME$. We define $\J:=(-\pi ,\pi )\cap (\Theta -2\pi ,\Theta )$, $\dompol:=\{(x,y)\in \mathbb{R}^2  \setminus (\mathbb{R}_-\times \{0\}) \mid \theta (x,y)\in \J\}$ and
\[\pol :\left( \begin{array}{ccc}
	\dompol & \rightarrow  & \dompsi\\[2pt]
	(x,y) & \mapsto  & (\rho ,\theta )
\end{array} \right) .\]
We will build $\varphi $ as $\psi \circ \pol$ with $\psi :\mathbb{R}^2  \rightarrow  \mathbb{R}^2 $ ($V_{1}$ and $V_{2}$ will be defined later). Let $\Opsi:=(\mathbb{R}_+^*)^2 $, $\Gpsi:=\mathbb{R}_+^*\times \{0\}$ and $\Spsi:=\{0\}\times \mathbb{R}_+^*$. It is clear that $\pol$ satisfies:
\begin{itemize}
	\item $\pol$ is a $\mathcal{C}^{\infty }$-diffeomorphism,
	\item $\pol(\PME)=0$, $\pol(\OME \cap \dompol)=\Opsi \cap \pol(\dompol)$, $\pol(\GME \cap \dompol)=\Gpsi \cap \pol(\dompol)$, $\pol(\SME \cap \dompol)=\Spsi \cap \pol(\dompol)$ (thanks to the definition of $\J$)
	\item $\forall (x,y) \in \Gpsi \cap \dompol, \; \diff\pol(x,y) = \mathrm{Diag} (1,1/x)$,
	\item $\forall (x,y) \in \Spsi \cap \dompol, \; \diff\pol(x,y) = \mathrm{Diag} (1,1/\RME)\circ R_{-\theta }$.
\end{itemize}
Hence, it suffices that $\psi $ satisfies:
\begin{enumerate}[label=(\roman*)]
	\item $\psi (0) = 0$, $\psi (\Opsi) = \Opsi$, $\psi (\Gpsi) = \Gpsi$, $\psi (\Spsi) = \Spsi$,
	\item $\psi $ is differentiable and $\diff \psi $ is Lipschitz,
	\item $\forall (\rho ,\theta ) \in \Gpsi, \; \diff \psi (\rho ,\theta ) = \mathrm{Diag} (1,\RME+\rho )$,
	\item $\forall (\rho ,\theta ) \in \Spsi, \; \diff \psi (\rho ,\theta ) = \mathrm{Diag} (1,\RME)$,
	\item there exist two neighborhoods $V',V_{2}$ of 0 in $\mathbb{R}^2 $ s.t.\ $\psi $ is a $\mathcal{C}^{1}$-diffeomorphism from $V'$ to $V_{2}$, and $V'\subset \pol(\dompol)$ (then it suffices to take $V_{1} :=\pol^{-1} (V')$).
\end{enumerate}
Let $\chi \in \mathcal{C}^{\infty } (\mathbb{R},[0,1])$ be equal to 0 in a vicinity of $-\infty $, and to 1 in a vicinity of $+\infty $. For~any $(\rho ,\theta )\in \mathbb{R}^2 $, we define:\vspace{-0.5ex}
\[\psi (\rho ,\theta ) := (\rho ,\psit(\rho ,\theta )) \qquad  \text{with}  \qquad  \psit(\rho ,\theta ) := \left\{\begin{array}{@{\;}ll@{}}
	\RME \theta  & \text{if } \rho \leqslant 0\\
	(\RME+\rho \,\chi (\ln \rho -\ln \theta ))\,\theta  & \text{if } \rho >0 \text{ and } \theta >0\\
	(\RME+\rho ) \theta  & \text{if } \rho >0 \text{ and } \theta =0\\
	-\psit(\rho ,-\theta ) & \text{if } \rho >0 \text{ and } \theta <0.
\end{array}\right.\]
It is clear that (i) is satisfied. Moreover, $\psi $ is differentiable in $\mathbb{R}^2  \setminus \{0\}$ because:
\begin{itemize}
	\item $\psit$ is $\mathcal{C}^{\infty }$ in $\mathbb{R}_+^*\times \mathbb{R}^* = \Opsi \cup  (\mathbb{R}_+^*\times \mathbb{R}_-^*)$,
	\item $\psit(\rho ,\theta )$ is equal to $\RME \theta $ in an open subset of $\Opsi$ that contains $(\mathbb{R}_-\times \mathbb{R}) \setminus \{0\}$,
	\item $\psit(\rho ,\theta )$ is equal to $(\RME +\rho )\theta $ in an open subset of  $\Opsi$ that contains $\Gpsi$,
	\item and these three open sets cover $\mathbb{R}^2  \setminus \{0\}$.
\end{itemize}
Let us show that $\diff \psi $ is Lipschitz on $\mathbb{R}^2  \setminus \{0\}$.
\begin{sousdemo}
	Let us show that $\diff\psit$ is Lipschitz in $\Opsi$, because the rest is easy.
	Let $(\rho ,\theta )\in \Opsi$. We have\vspace{-2pt}
	\[\!\Bigg\{ \begin{array}{@{\;}r@{\,=\,}l}
		\partial _{\rho } \psit(\rho ,\theta ) & \theta  \cdot  (\chi +\chi ')(\ln \rho -\ln \theta ) \\[3pt]
		\partial _{\theta } \psit(\rho ,\theta ) & \RME + \rho  \cdot  (\chi -\chi ')(\ln \rho -\ln \theta )
	\end{array} 
	\quad \  \mathrm{so} \quad \ 
	\left\{\begin{array}{@{\,}r@{\,=\,}l}
		\partial _{\rho } ^2  \psit(\rho ,\theta ) & \frac{\theta }{\rho } \cdot  (\chi '+\chi '')(\ln \rho -\ln \theta ) \\[3pt]
		\partial _{\rho } \partial _{\theta } \psit(\rho ,\theta ) & (\chi -\chi '')(\ln \rho -\ln \theta ) \\[4pt]
		\partial _{\theta } ^2 \psit(\rho ,\theta ) & \frac{\rho }{\theta } \cdot  (\chi ''-\chi ')(\ln \rho -\ln \theta ).
	\end{array} \right.\]
	Let $M>0$ be s.t.\ $\operatorname{supp}(\chi ')\subset [-M,M]$. We deduce $|\partial _{\rho } ^2  \psit(\rho ,\theta )| \leqslant  e^{M} \|\chi '+\chi ''\|_{L^{\infty }}$, $|\partial _{\rho } \partial _{\theta } \psit(\rho ,\theta )| \leqslant \|\chi -\chi ''\|_{L^{\infty }}$ and $|\partial _{\theta } ^2  \psit(\rho ,\theta )| \leqslant  e^{M} \|\chi ''-\chi '\|_{L^{\infty }}$. Thus $\diff ^2 \varphi $ is bounded on $\Opsi$, so $\diff \varphi $ is Lipschitz.
\end{sousdemo}
In~particular $\lim _{\bfx\rightarrow 0} \diff \varphi $ exists. This yields that $\psi $ is differentiable on $\mathbb{R}^2 $ and proves (ii). Moreover (iii) and (iv) are easy to verify. Finally, (v) is a direct consequence of the inverse function theorem, as $\diff \varphi (0)=\mathrm{Diag} (1,\RME)$ is invertible.
\end{demo}

\begin{demo}[Proof of \cref{regularite elliptique pour ME 2}]
\phantomsection\label{preuve: regularite elliptique pour ME 2}%
We only prove the second point, as the first one is simpler. First, we assume that $\RME=1$, then we will extend the result to any $\RME$ by a scaling method.\medskip

\underline{Step 1:} Let us assume that $\RME=1$, $0<\varepsilon \leqslant M$ and 
\begingroup\let\RME=1%
\[\left\{\begin{array}{r@{}ll}
	\Delta u &{}= f &\text{ in }\OME\\
	\mu _{0} \partial _{n} u + \frac{\mu _{1} }{\varepsilon } u - \varepsilon \,(\mu _{1} \partial _{x} ^2 +\omega ^2 \rho _{1} )u &{}=g &\text{ on }\GME\\
	\partial _{r} u &{}= h &\text{ on }\SME \text{ and } \text{at }  \PME
\end{array}\right.\]
for some $f\in L^2 (\OME)$. Let $\varphi :V_{1} \rightarrow V_{2}$ be the diffeomorphism given by \cref{diffeo ME}. Let $a>0$ be s.t.\ $[\RME,\RME+a)\cap \{0\} \subset V_{1}$. We define $I:=(\RME+ \frac{a}{5} ,\infty )$, $J:=(\RME+ \frac{2a}{5} ,\infty )$ and we take $\eta >0$ s.t.\ $I\times (0,\eta )\subset \OME$. Let $\tilde{g} := g-\varepsilon \,\omega ^2 \rho _{1} u_{|\GME}$. We apply \cref{regularite elliptique pour ME} on $u$ with these values of $I,J,\eta , \tilde{g}$, with $A$ the identity matrix, and with $\nu _{1} :=\mu _{1}$ and $\nu _{2} := \mu _{1} /3$. This gives
\begin{align*}
	\smash{\|u\|_{H^{1/2} ((\RME+ \frac{2a}{5} ,\infty )\times \{0\})} } &\leqslant  C_{M} \,\varepsilon \,\big(\|f\|_{L^2 (\OME)} + \| \tilde{g} \|_{H^{1} (\GME)} +\|u\|_{H^{1} (\OME)} +\sqrt{\varepsilon } \,\snhg{u} \big) \\
	&\leqslant  C_{M} \,\varepsilon \,\big( \|f\|_{L^2 (\OME)} + \|\red{g}\|_{H^{1} (\GME)} +\|u\|_{H^{1} (\OME)} + \sqrt{\varepsilon } \, \snhg{u} \big).
\end{align*}
It remains to treat the vicinity of $\PME$. Let $v:\OME\rightarrow \mathbb{C}$ be defined in polar coordinates by $v(r,\theta ):= r\cdot h(\RME,\theta )$. The function $u_{1} :=u+v$ satisfies
\[\left\{\begin{array}{@{\,}r@{}l@{\;}l}
	\Delta u_{1} &{}= f_{1} := f+\Delta v&\text{ in }\OME\\
	\mu _{0} \partial _{n} u_{1} + \frac{\mu _{1} }{\varepsilon } u_{1} - \varepsilon \.\mu _{1} \partial _{x} ^2 u_{1} &{}=g_{1} := \tilde{g} + \mu _{0} \partial _{\theta } h(\PME) + \frac{\mu _{1} }{\varepsilon } \.r \.h(\PME) &\text{ on }\GME\\
	\partial _{r} u_{1} &{}= 0 &\text{ on }\SME \text{ and } \text{at }  \PME
\end{array}\right.\]
Now let $u_{2} = u_{1} \circ \varphi ^{-1}$, $g_{2} = g_{1} \circ \varphi ^{-1}$, $f_{2} =(\det(\diff \varphi )^{-1} f_{1} )\circ \varphi ^{-1}$ and $A = [\det(\diff \varphi )^{-1} \cdot \diff \varphi \cdot \diff \varphi ^{\mathsf{T}} ]\circ \varphi ^{-1}$, where $\varphi ^{-1} :V_{2} \rightarrow V_{1}$ is the inverse function of $\varphi $ and $\det(\diff \varphi )^{-1} :x\mapsto  1/\det(\diff \varphi (x))$. The properties of $\varphi $ (see \cref{diffeo ME}) imply:
\begin{equation}{\label{regularite elliptique pour ME 2 preuve u2}}
\left\{\begin{array}{@{\,}r@{}l@{\;}l}
	\operatorname{div}(A \nabla  u_{2} ) &{}= f_{2} &\text{ in }(\mathbb{R}_+^*)^2  \cap V_{2} \\[1pt]
	\mu _{0} \partial _{n} u_{2} + \frac{\mu _{1} }{\varepsilon } u_{2} - \varepsilon \,\mu _{1} \partial _{x} ^2 u_{2} &{}=g_{2}  &\text{ on }(\mathbb{R}_+^*\times \{0\}) \cap  V_{2} \\[1pt]
	\partial _{x} u_{2} &{}= 0 &\text{ on }(\{0\}\times \mathbb{R}_+^*) \cap  V_{2} \text{ and } \text{at }  (0,0)
\end{array}\right.
\end{equation}
In~particular, we used the value of $\diff \varphi $ on $\partial \OME$, as it gives $\partial _{n} u_{2|(\mathbb{R}_+^*)^2 } = \partial _{n} u_{1|\partial \OME} \circ \varphi ^{-1}$ and $\partial _{x} ^2 u_{2|\mathbb{R}_+^*\times \{0\}} = \partial _{x} ^2 u_{1|\GME} \circ \varphi ^{-1}$. Let $b,c>0$ be s.t.\ $(b,0)=\varphi (\RME+ \frac{3a}{5} ,0)$ and $(c,0)=\varphi (\RME+ \frac{4a}{5} ,0)$. Let $\eta '>0$ be s.t.\ $[0,c]\times [0,\eta ']\subset V_{2}$. We extend $u_{2}$, $f_{2}$ and $g_{2}$ to $(-a',c)\times (0,\eta ')$ by setting:
\[\forall (x,y)\in (-c,0)\times (0,\eta '), \quad  u_{2} (x,y):=u_{2} (-x,y),\ f_{2} (x,y):=f_{2} (-x,y) \ \text{ and } \ g_{2} (x):=g_{2} (-x).\]
Thanks to \cref{regularite elliptique pour ME 2 preuve u2}, these extensions satisfy
\[\left\{\begin{array}{@{\,}r@{}l@{\;}l}
	\operatorname{div}(A \nabla  u_{2} ) &{}= f_{2} &\text{ in }R := (-c,c)\times (0,\eta ') \\[1pt]
	\mu _{0} \partial _{n} u_{2} + \frac{\mu _{1} }{\varepsilon } u_{2} - \varepsilon \,\mu _{1} \partial _{x} ^2 u_{2} &{}=g_{2}  &\text{ on }I_{0} := (-c,c)\times \{0\}
\end{array}\right.\]
Moreover $A$ is Lipschitz because $D\varphi $ is Lipschitz by \cref{diffeo ME}. A second application of \cref{regularite elliptique pour ME} yields $\|u_{2} \|_{H^{1/2} ((0,b)\times \{0\})} \leqslant  C_{M} \,\varepsilon \,\big(\|f_{2} \|_{L^2 (R)} + \|g_{2} \|_{H^{1} (I_{0} )} +\|u_{2} \|_{H^{1} (R)} +\sqrt{\varepsilon } \, \|u_{2} \|_{H^{1} (I_{0} )} \big)$. Thus:
\begin{align*}
	\|u\|_{H^{1/2} ((\RME,\RME+ \frac{3a}{5} )\times \{0\})} &\leqslant  \|u_{1} \|_{H^{1/2} ((\RME,\RME+ \frac{3a}{5} )\times \{0\})} + C |h(\PME)|\\
	&\leqslant  C \.\|u_{2} \|_{H^{1/2} ((0,b)\times \{0\})} + C |h(\PME)|\\
	&\leqslant  C_{M} \,\varepsilon \,\big(\|f_{1} \|_{L^2 (\OME \cap V_{1} )} + \|g_{1} \|_{H^{1} (\GME \cap V_{1} )} +\|u_{1} \|_{H^{1} (\OME \cap V_{1} )} +\sqrt{\varepsilon } \.\snh{u_{1} }(\GME \cap V_{1} ) \big) \\
	& \quad  + C \.|h(\PME)| \\
	&\leqslant  C_{M} \,\varepsilon \,\big(\|f\|_{L^2 (\OME)} + \|g\|_{H^{1} (\GME)} +\|h\|_{H^2 (\SME)} +\|u\|_{H^{1} (\OME)} +\sqrt{\varepsilon } \.\snhg{u} \big) \\
	& \quad + C_{M} \.|h(\PME)| 
\end{align*}
where we used in the last inequality:
\begin{itemize}
	\item $\|f_{1} \|_{L^2 (\OME \cap V_{1} )} \leqslant  \|f\|_{L^2 (\OME \cap V_{1} )} + C\.\|h\|_{H^2 (\SME)}$,
	\item $\|g_{1} \|_{H^{1} (\GME \cap V_{1} )} \leqslant \| \tilde{g} \|_{H^{1} (\GME \cap V_{1} )} + C\.|\partial _{\theta } h(\PME)|+ \frac{C}{\varepsilon } |h(\PME)|$ with $|\partial _{\theta } h(\PME)|\leqslant C \.\|h\|_{H^2 (\SME)}$,
	\item $\|u_{1} \|_{H^{1} (\OME \cap V_{1} )} \leqslant  \|u\|_{H^{1} (\OME \cap V_{1} )} + C\.\|h\|_{H^2 (\SME)}$,
	\item $\snh{u_{1} }(\GME \cap V_{1} ) \leqslant  \snh{u}(\GME \cap V_{1} ) + C\.|h(\PME)|$ with $|h(\PME)|\leqslant C \.\|h\|_{H^2 (\SME)}$. 
\end{itemize}
Finally we use $\|u\|_{H^{1/2} (\GME)} \leqslant  \|u\|_{H^{1/2} ((\RME,\RME+ \frac{3a}{5} )\times \{0\})} + \|u\|_{H^{1/2} ((\RME+ \frac{2a}{5} ,\infty )\times \{0\})}$.\medskip\smallskip
\endgroup

\underline{Step 2:} Let $\OME_{1} ,\RME_{1} ,\SME_{1} ,\PME_{1}$ be the values taken by $\OME,\RME,\SME,\PME$ when $\RME=1$. We assume $\RME \in (0,R]$ and $\varepsilon  \in (0,M\RME]$. We define the rescaled function $u_{1} :x\in \OME_{1} \mapsto u(\RME x)$ and similarly $g_{1} :=\RME g(\RME\varlibre)$ and $h_{1} :=\RME \.h(\RME\varlibre)$. Let $\varepsilon _{1} :=\varepsilon /\RME$. We have
\[\left\{\begin{array}{r@{}ll}
	\mu _{0} \Delta u_{1} &{}= f_{1} :=-\RME^2  \omega ^2 \rho _{0} u_{1} &\text{ in }\OME_{1} \\
	\mu _{0} \partial _{n} u_{1} + \frac{\mu _{1} }{\varepsilon _{1}} u_{1} - \varepsilon _{1} \,(\mu _{1} \partial _{x} ^2 +\omega ^2 \rho _{1} )u_{1} &{}=g_{1} &\text{ on }\GME_{1} \\
	\partial _{r} u_{1} &{}= h_{1} &\text{ on }\SME_{1} \text{ and } \text{at }  \PME_{1}
\end{array}\right.\]
Note that the boundary condition involves $\varepsilon _{1}$ instead of $\varepsilon $. Step 1 shows
\begin{multline*}
	\quad  \|u_{1} \|_{H^{1/2} (\GME_{1} )} \leqslant  C_{M} \,\varepsilon _{1} \Big(\|f_{1} \|_{L^2 (\OME_{1} )} + \|g_{1} \|_{H^{1} (\GME_{1} )} +\|h_{1} \|_{H^2 (\SME_{1} )} \\
	+\|u_{1} \|_{H^{1} (\OME_{1} )} +\sqrt{\varepsilon _{1}} \, \snh{u_{1} }(\GME_{1} ) \Big) + C_{M} \.|h_{1} (\PME_{1} )|. \quad 
\end{multline*}
Moreover, for any $\psi :\OME\rightarrow \mathbb{C}$ regular enough, $\psi _{1} :=\psi (\RME\varlibre)$ and $j\in \mathbb{N}$, we have $\|\psi _{1} \|_{L^2 (\OME_{1} )} = \RME^{-1} \|\psi \|_{L^2 (\OME)}$, $\|\nabla  \psi _{1} \|_{L^2 (\OME_{1} )} = \|\nabla  \psi \|_{L^2 (\OME)}$ and $\snh{\psi _{1} }(\GME_{1} ) = \sqrt{\RME} \.\snhg{\psi }$. Thus, by~definition of $\HME{j}$ in \cref{def: sobolev pondere pour ME}, we finally have:
\begin{multline*}
	\quad  \frac{\|u\|_{\HME{1/2}(\GME)}}{\sqrt{\RME}} \leqslant  C_{R,M} \frac{\varepsilon }{\RME} \bigg( \sqrt{\RME} \,\|g\|_{\HME1(\GME)} + \sqrt{\RME} \,\|h\|_{\HME2(\SME)} \\
	+ \frac{\|u\|_{H^{1} (\OME)}}{\RME}  +\sqrt{\varepsilon } \,\snhg{u} \bigg) + C_{M} \RME \.|h(\PME)|. \quad 
\end{multline*}
\end{demo}

\section{Appendix: behavior of the stability condition and the error estimate w.r.t.\ the number of modes}
\label{sec: annexe comportement en N des ME}

\def\black#1{\textcolor{black}{#1}}%
\def\darkred#1{\textcolor{red!75!gray}{#1}}%
\newenvironment{listeE}{\begin{enumerate}[leftmargin=2em,itemsep=2pt,labelsep=0.5em,topsep=2pt]}{\end{enumerate}}
\newenvironment{listeI}{\begin{itemize}[leftmargin=2.5em,itemsep=2pt,labelsep=0.5em,topsep=1pt,label={$-$}]}{\end{itemize}}
\newenvironment{nvdemo}{%
	\def\red##1{\textcolor{red!78!gray}{##1}}%
	\def\darkred##1{\textcolor{red!72!gray}{##1}}%
	\def\verydarkred##1{\textcolor{red!68!gray}{##1}}%
	\def\black##1{\textcolor{grisfonce}{##1}}%
	\begin{itemize}[label={},leftmargin=1em,topsep=-4pt]%
	\item[]\begin{demo}[Proof sketch]}{\end{demo}\end{itemize}}

\begin{textesansboite}
In this appendix, we explain how one can prove that the constant $\CstabME$ involved in the stability condition of \cref{def: condition de stabilite} can be taken as
\begin{equation}{\label{eq: comportement CstabME en N}}
	\CstabME = \frac{C_{\Rsource} }{N^{3}} .
\end{equation}
In addition, we prove \cref{norme opP}, which is a step of the demonstration of \cref{cvg des ME}, and we explain how the error estimate of this theorem can be improved into:
\begin{equation}{\label{eq: cvg des ME amelioree}}
	\red{\forall \eta >0,} \qquad  \|\sme-u_{\varepsilon } \|_{H^{1} (\OME)} \leqslant  C_{\fsource,R\red{,\eta }} \left[ \bigg( \frac{\varepsilon }{\RME} \bigg)^{m} \domi{\ln \varepsilon }^{\ell _{m}} N^{\red{2.5+\eta }} + \left( \frac{\RME}{R} \right)^{N\pst} \right].
\end{equation}
This requires modifying numerous intermediate results to make them more precise and thus increases the technicality of \cref{chap: decompo modale CL,chap: ME}. We do~not go into detail but only indicate how the different results should be modified, with sometimes hints in how to adapt the proofs. Changes are shown in red.

\begin{itemize}[leftmargin=1.45em,itemsep=4pt,labelsep=0.33em]
	\item Points 3 and 4 of \cref{norme-gamma des op R} become
	\begin{listeE}
		\item[3.] $\forall \varphi \in \AOprecis{d,n,\delta }, \ \ \normg{\partial _{y} \varphi _{|\Gamma } } \leqslant  C_{\delta ,n} \,\red{\domi{d}^{1/2} }\.\normg{\varphi }$,
    	\item[4.] $\forall \varphi \in \AGprecis{d,n}, \ \ \normg{\partial _{x} \varphi } \leqslant  C_{\delta ,n} \,\red{\domi{d}^{1/2} } \.\normg{\varphi }$ provided that $d\geqslant \delta $,
	\end{listeE}
	with $\domi{d}:=1+|d|$.	
	This version is more precise than the previous one. It is obtained by simply not bounding $\domi{d}^{1/2}$ by the geometric behavior $\eAMP^{d}$.

	\item The norms $\|.\|_{\eAMP}$ introduced in \cref{def: norme sur Amacro} are replaced, for~any $\eAMP>0$ and $n\in \mathbb{N}$, by 
	\[\|\Phi \|_{\eAMP\red{,n}} := \sum _{d\in \mathbb{R}} \red{\domi{d}^{-n/2} } \.\eAMP^{-d} \.\normg{\varphi _{d} } \ \in \ [0,+\infty ] \qquad  \text{for~any }  \Phi =(\varphi _{d} )_{d} \in \Amacro(\Omega ).\]
	In particular, we have $\|.\|_{\eAMP} =\|.\|_{\eAMP,0}$ for~any $\eAMP$.

	\item \Cref{norme des op R sur series formelles} becomes
	\begin{listeE}
		\item $\forall \Phi \in \AMprecis{n,\delta }(\Gamma ), \ \ \|\RGammaOmega(\Phi )\|_{\eAMP\red{,n}} \leqslant  C_{n} \|\Phi \|_{\eAMP\red{,n}}$,
		\item $\forall \Phi \in \AMprecis{n,\delta }(\Omega ), \ \ \|\ROmegaOmega(\Phi )\|_{\eAMP\red{,n}} \leqslant  C_{\delta ,n} \,\eAMP^{-2} \, \|\Phi \|_{\eAMP\red{,n}}$,
		\item $\forall \Phi \in \AMprecis{n,\delta }(\Omega ), \ \ \|\partial _{y} \Phi _{|\Gamma } \|_{\red{\eAMP,n+1}} \leqslant  C_{\delta ,n} \,\|\Phi \|_{\eAMP\red{,n}}$,
		\item $\forall \Phi \in \AMprecis{n,\delta }(\Gamma ), \ \ \|\partial _{x} \Phi \|_{\red{\eAMP,n+1}} \leqslant  C_{\delta ,n} \,\|\Phi \|_{\eAMP\red{,n}}$.
	\end{listeE}
	In particular, the increase in $\eAMP$ in points 3 and 4 is replaced by an increase in $n$, which is sharper.

	\item \Cref{norme de <Rlap>} becomes $\big\|\langle -k_{0} ^2 \ROmegaOmega\rangle (\Phi )\big\|_{\eAMP\red{,n}} \leqslant  C_{\delta ,n} \, \|\Phi \|_{\eAMP\red{,n}}$.

	\item \Cref{norme SF de mcl cas rationnel,norme SF de mcl cas general} become $\|\mcla_{d,n} \|_{\red{\eAMP,n}} \leqslant  C_{\delta ,n} \,\eAMP^{-d} \,\gamma (d)$.
	\begin{nvdemo}
		Start with $\|\mcla_{d,n} \|_{\eAMP\red{,0}} \leqslant  C_{\delta } \,\eAMP^{-d} \,\gamma (d)$ and use induction like in the original proof, but increase the norm index $n$ at each step instead of increasing $\eAMP$, using the previous modified results.
	\end{nvdemo}

	\item \Cref{majoration des fonctions de A par leur norme} becomes
	\[|\varphi (r,\theta )| + \red{ \frac{\black r}{\domi{q}} } \.|\nabla  \varphi (r,\theta )| \;\leqslant \; C_{\red{\delta ,}n} \,r^{d} \, (1+|{\ln r}|)^{2n+1} \|\varphi \|.\]

	\item \Cref{eq: vitesse de cvg des SF de norme finie} in \cref{cvg des SF de norme finie} becomes
	\providecommand{\restPhi}{}%
	\renewcommand{\restPhi}[1][q]{\rho_{#1}}%
	\[|\restPhi[q'](r,\theta )| + \red{ \frac{\black r}{\domi{q}} } \.|\nabla  \restPhi[q'](r,\theta )| \leqslant  C_{\delta ,n,R,\eAMP\red{,k}} \,\red{\domi{q}^{k/2} }\, \frac{(\red\eAMP r)^{q} }{\gamma (q)} \, (1+|{\ln r}|)^{2n+1} \|\Phi \|_{\eAMP\red{,k}} \]
	where $\restPhi[q']:=(\lim  \Phi ) -T_{<q'} (\Phi )$ as a reminder.

	\item \Cref{eq: vitesse de cvg des mcl} in \cref{convergence des mcl} becomes
	\providecommand{\restmcl}{}%
	\renewcommand{\restmcl}{\rho_{d,n,q}}%
	\[|\restmcl (r,\theta )| + \red{ \frac{\black r}{\domi{q}} } \.|\nabla  \restmcl (r,\theta )| \;\leqslant \; C_{\delta ,n,R,\eAMP} \,\red{\domi{q}^{n/2} }\, \frac{\gamma (d) \.\eAMP^{-d} }{\gamma (q) \.\red\eAMP^{-q} } \,r^{q} \, (1+|{\ln r}|)^{2n+1} \]
	where $\restmcl := \mcll_{d,n} - T_{<q} (\mcla_{d,n} )$.

	\item \Cref{estimation des mcl} becomes
	\[|\mcll_{d,n} (r,\theta )| + \red{ \frac{\black r}{\domi{d}} } |\nabla  \mcll_{d,n} (r,\theta )| \;\leqslant \; C_{\delta ,n,R} \,\red{\domi{d}^{n} }\.r^{d-n} \, (1+|{\ln r}|)^{2n+1} .\]
	\begin{nvdemo}
		Like in the original proof, apply \cref{convergence des mcl} with $q:=d-n$, and use the definition of $\gamma $ to write $\domi{d-n}^{n/2} \, \frac{\gamma (d)}{\gamma (d-n)} \leqslant C_{\delta ,n} \domi{d}^{n}$.
	\end{nvdemo}

	\item \Cref{estimation mme} becomes: if $\varepsilon  \,\domi{\ln \RME}^{3} \red{N} \leqslant M\RME$, then
	\setalignbox{$\|\partial _{\theta } ^{j} \partial _{r} (\mme_{n} - \mme0_{n} )\|_{L^{\infty } (\SME)}$}
	\begin{listeE}
		\item $\alignbox{\|\mme_{n} \|_{L^{\infty } (\SME)}} \leqslant  C_{R,M} \,(\RME)^{n\pst}$,
		\item $\alignbox{\|\partial _{\theta } ^{j} \partial _{r} \mme_{n} \|_{L^{\infty } (\SME)}} \leqslant  C_{R,M} \,\red{n^{j+1} }\.(\RME)^{n\pst-1}$,
		\item $\alignbox{\|\mme_{n} - \mme0_{n} \|_{L^{\infty } (\SME)}} \leqslant  C_{R,M} \,\varepsilon \.\red{n} \,(\RME)^{n\pst-1} \,\domi{\ln \RME}^{3}$,
		\item $\alignbox{\|\partial _{\theta } ^{j} \partial _{r} (\mme_{n} - \mme0_{n} )\|_{L^{\infty } (\SME)}} \leqslant  C_{R,M} \,\varepsilon \,\red{n^{j+2} }\.(\RME)^{n\pst-2} \,\domi{\ln \RME}^{3}$.
	\end{listeE}
	where the constants are independent of $N$.
	
	\item \Cref{opS - opS0} becomes: there exist two constants $\CdtnMEprecis1,\CdtnMEprecis2>0$ (independent of $N$) s.t., if $0\leqslant \varepsilon \leqslant \CdtnMEprecis1 \,\RME \,\domi{\ln \RME}^{-3} \red{N^{-2} }$, then we have:
	\begin{listeE}
		\item $\gram$ is invertible,
		\item $\displaystyle \forall j\in [\![0,2]\!],\;  \|\opS-\opSref\|_{\mathcal{L}(L^2 (\SME),\HME{\red{j}}(\SME))} \leqslant  \CdtnMEprecis2 \, \frac{\varepsilon  \.\red{N^{3+j} }}{\RME^{2}} \,\domi{\ln \RME}^{3} $,
		\item $\displaystyle \|\opP\|_{\mathcal{L}(L^2 (\SME),\mathbb{C})} \leqslant  \CdtnMEprecis2 \, \frac{\varepsilon  \.\red{N^{2.5} }}{\RME^{2.5}} \,\domi{\ln \RME}^{3} $.
	\end{listeE}
	\begin{nvdemo}
		We use the notations $a_{n}$, $\mmen_{n}$, $\gramn$, $V_{\varepsilon }$ and $V'_{\varepsilon }$ from the proof of \cref{opS - opS0}. We also introduce the following norms for any $(z_{n} )\in \mathbb{C}^{N}$ and $A\in \mathcal{M}_{N} (\mathbb{C})$:
		\[\|(z_{n} )\|_{\red2} :=\darkred{\bigg(\sum _{n=1} ^{N} |z_{n} |^2 \bigg)^{1/2} } \qquad  \text{ and } \qquad  \|A\|_{\red2} := \sup _{\mathbf z\in \mathbb{C}^{N} \setminus \{0\}} \frac{\|A{\bf z}\|_{\red2} }{\|{\bf z}\|_{\red2} }.\]
		Then step 1 of the proof becomes (under the assumption $\varepsilon \,\red{N}\.\domi{\ln \RME}^{3} \leqslant M\RME$):
		\setalignbox{$\|\partial _{\theta ^{j}} \partial _{r} (\mmen_{n} - \mmen0_{n} )\|_{L^{\infty } (\SME)}$}
		\begin{listeI}
			\item $\alignbox{\|\mmen_{n} \|_{L^{\infty } (\SME)}} \leqslant  C_{R,M} \,(\RME)^{-1/2}$,
			\item $\alignbox{\|\partial _{\theta } ^{j} \partial _{r} \mmen_{n} \|_{L^{\infty } (\SME)}} \leqslant  C_{R,M} \,\red{N^{j+1} }\.(\RME)^{-3/2}$,
			\item $\alignbox{\|\mmen_{n} - \mmen0_{n} \|_{L^{\infty } (\SME)}} \leqslant  C_{R,M} \,\varepsilon \,(\RME)^{-3/2} \,\red{N}\.\domi{\ln \RME}^{3}$,
			\item $\alignbox{\|\partial _{\theta } ^{j} \partial _{r} (\mmen_{n} - \mmen0_{n} ) \|_{L^{\infty } (\SME)}} \leqslant  C_{R,M} \,\varepsilon \,(\RME)^{-5/2} \,\red{N^{j+2} }\.\domi{\ln \RME}^{3}$.
		\end{listeI}
		Thus we have $|\langle \mmen_{n} ,\mmen_{n'} \rangle -\langle \mmen0_{n} ,\mmen0_{n'} \rangle | \leqslant  C_{R,M} \,\varepsilon \, \red{N} \.\RME^{-1} \,\domi{\ln \RME}^{3}$, which implies \[\|\gramn-\gramn0\|_{\red{2}} \leqslant  \bigg(\sum _{1\leqslant n,n'\leqslant N} \langle \mmen_{n} ,\mmen_{n'} \rangle -\langle \mmen0_{n} ,\mmen0_{n'} \rangle  \bigg)^{1/2} \leqslant  C_{R,M} \,\varepsilon \,\red{N^2 } \.\RME^{-1} \,\domi{\ln \RME}^{4} .\]
		Hence there exists $\CdtnMEprecis1 >0$ s.t., if $\varepsilon  \leqslant  \CdtnMEprecis1 \,\RME \domi{\ln \RME}^{-3} \red{N^{-2} }$ then we have $\|\gramn-\gramn0\|_{\infty } \leqslant  \frac{1}{2}$, which implies $\|\gramn^{-1} \|_{\infty } \leqslant  2$ and $\|\gramn^{-1} -\gramn0^{-1} \|_{\infty } \leqslant  C_{R,N} \, \varepsilon  \,\red{N^2 } \.\RME^{-1} \,\domi{\ln \RME}^{3}$.\\
		Next, step 2 is modified into:
		\setalignbox{$\|\partial _{\theta } ^{j} (V'_{\varepsilon } -V'_{0} )(\bfx)\|_{\red2}$}
		\begin{listeI}
			\item $\|V_{\red0} (v)\|_{\red2} \darkred{ {}= \left( \sum\limits _{n=1} ^{N} |\langle v,\mmen_{n} \rangle |^2  \right)^{\!1/2} = \left( \sum\limits _{n=1} ^{N} \left| \left\langle v, \sqrt{ \frac{2}{\RME\Theta } } \sin \!\left(n\pst \theta \right) \right\rangle  \right|^2  \right)^{\!1/2} } \leqslant  \|v\|_{L^2 (\SME)}$
			(the orthonormality of the modes $\mmen0_{n}$ avoids a factor $N^{1/2}$ due to the number of terms of the sum)
			\item $\alignbox{\|V_{\varepsilon } (v)-V_{0} (v)\|_{\red2} } \leqslant  C_{R} \,\varepsilon \,\red{N^{1.5} }\.\RME^{-1} \,\domi{\ln \RME}^{3} \,\|v\|_{L^2 (\SME)}$
			\item $\alignbox{\|\partial _{\theta } ^{j} V'_{\red0} (\bfx)\|_{\red2} } \leqslant  C_{R} \.\red{N^{j+1} }\.\RME^{-3/2} $ similarly to $\|V_{0} (v)\|_{2}$
			\item $\alignbox{\|\partial _{\theta } ^{j} (V'_{\varepsilon } -V'_{0} )(\bfx)\|_{\red2} } \leqslant  C_{R} \,\varepsilon \,\red{N^{j+2.5} }\.\RME^{-5/2} \,\domi{\ln \RME}^{3}$
		\end{listeI}
		Hence, using $\varepsilon  \leqslant  \CdtnMEprecis1 \,\RME \domi{\ln \RME}^{-3} N^{-2}$, we obtain:
		\begin{listeI}
			\item $\|V_{\varepsilon } (v)\|_{\red2} \darkred{ {}\leqslant \|V_{0} (v)\|_{2} + \|V_{\varepsilon } (v)-V_{0} (v)\|_{2} }\leqslant  C_{R} \|v\|_{L^2 (\SME)}$
			\item $\|\partial _{\theta } ^{j} V'_{\varepsilon } (\bfx)\|_{\red2} \darkred{ {} \leqslant \|\partial _{\theta } ^{j} V'_{0} (\bfx)\|_{2} + \|\partial _{\theta } ^{j} (V'_{\varepsilon } -V'_{0} )(\bfx)\|_{2} } \leqslant  C_{R} \,\red{N^{j+1} }\.\RME^{-3/2}$
		\end{listeI}
		Thus, we have
		\def\alignbox#1{\begingroup
			\setbox0=\hbox{$\big|\partial _{\theta } ^{j} (\opS-\opSref)(v)(\bfx)\big|$}
			\dimen0=\linewidth
			\advance\dimen0 by -\wd0
			\advance\dimen0 by -4.5em
			\hbox to \dimen0{${}#1$}\endgroup}
		\begin{align*}
			\big|\partial _{\theta } ^{j} (\opS-\opSref)(v)(\bfx)\big| &\alignbox{\leqslant  \|V_{\varepsilon } (v)-V_{0} (v)\|_{\red2} \,\|\gramn^{-1} \|_{\red2} \,\|\partial _{\theta } ^{j} V'_{\varepsilon } (\bfx)\|_{\red2} \hfill}\\
			&\alignbox{\hfill {} + \|V_{0} (v)\|_{\red2} \,\|\gramn^{-1} -\gramn0^{-1} \|_{\red2} \,\|\partial _{\theta } ^{j} V'_{\varepsilon } (\bfx)\|_{\red2} }\\
			& \alignbox{\hfill {} + \|V_{0} (v)\|_{\red2} \,\|\gramn0^{-1} \|_{\red2} \,\|\partial _{\theta } ^{j} (V'_{\varepsilon } -V'_{0} )(\bfx)\|_{\red2} }\\
			&\leqslant  C_{R} \,\varepsilon \,\red{N^{3+j} }\.\RME^{-5/2} \,\domi{\ln \RME}^{3} \,\|v\|_{L^2 (\SME)}
		\end{align*}
		\def\bfz{{\bf z}}\def\bfw{{\bf w}}%
		Note that using the norm 2 instead of $\infty $ avoids an additional factor $N$ here as we have, for any $\bfw,\bfz \in \mathbb{C}^{N}$, $|\bfw^{\mathsf{T}} \bfz|\leqslant \|\bfw\|_{2} \|\bfz\|_{2}$ but only $|\bfw^{\mathsf{T}} \bfz|\leqslant N \|\bfw\|_{\infty } \|\bfz\|_{\infty }$. Similarly, we have:
		\[|\opP(v)| \leqslant  \|V_{\varepsilon } (v)\|_{\red2} \,\|\gramn^{-1} \|_{\red2} \,\|(V'_{\varepsilon } -V'_{0} )(\bfx)\|_{\red2} \leqslant  C_{R} \,\varepsilon \,\red{N^{2.5} }\.\RME^{-5/2} \,\domi{\ln \RME}^{3} \,\|v\|_{L^2 (\SME).}\]
	\end{nvdemo}

	\item The stability condition of \cref{def: condition de stabilite} becomes ``$\varepsilon \leqslant \red\CstabMEprecis\, \RME \,\domi{\ln \RME}^{-3} \.\red{N^{-3} }$'' if $m=3$ and ``$\varepsilon \leqslant \red\CstabMEprecis\, \RME^{4/3} \,\domi{\ln \RME}^{-2} \.\red{N^{-3} }$'' if $m=5$, for some constant $\CstabMEprecis>0$.
	
	\item \Cref{caractere bien pose des ME augmentes} stays the same and its proof needs to be slightly modified in order to take $\CstabME = C_{\Rsource} /N^{3}$ as in \cref{eq: comportement CstabME en N}. 
	\begin{nvdemo}
		\def\cME#1{c^{(#1)}_\peps}%
		The case $m=3$ is easy to adapt. It suffices mainly to replace each $\CdtnME1$ by $\CdtnMEprecis1 /N^2 $, $\CdtnME2$ by $\CdtnMEprecis2 /N^{3}$ and $\CstabME$ by $\CstabMEprecis /N^{3}$.\\
		As for $m=5$, using the above generalization of \cref{opS - opS0}, \cref{eq: caractere bien pose des ME augmentes preuve m=5 forme c} becomes
		\[|\cME3(u,u)| \leqslant  C_{\Rsource,N} \,\domi{\ln \RME}^{3} \red{N^{3} } \left( \frac{\varepsilon }{\RME} + \frac{\varepsilon ^{3/2} }{\RME^2 } \right) \|u\|_{\VME5,\varepsilon } ^2  \]
		where, in the estimate of $\|\opP\|_{\mathcal{L}(L^2 (\SME),\mathbb{C})}$, we simply bounded $N^{2.5}$ by $N^{3}$. Then, \cref{eq: caractere bien pose des ME augmentes preuve m=5 etude condition stab} becomes
		\begin{align*}
			\varepsilon \leqslant  \red{ \frac{\black{M}}{N^{3}} } \, \frac{\RME^{4/3} }{\domi{\ln \RME}^2 } \ &\Rightarrow \ \varepsilon ^{3/2} \leqslant  \left( \red{ \frac{\black{M}}{N^{3}} } \right)^{3/2} \frac{\RME^2 }{\domi{\ln \RME}^{3}} \; \text{ and } \; \varepsilon \leqslant  \red{ \frac{\black{M}}{N^{3}} } \,C_{\Rsource} \frac{\RME}{\domi{\ln \RME}^{3}}  \\
			&\Rightarrow \ \domi{\ln \RME}^{3} \red{N^{3} } \left( \frac{\varepsilon }{\RME}  + \frac{\varepsilon ^{3/2} }{\RME^2 } \right) \leqslant  C_{\Rsource} \left( M+ \red{ \frac{\black{M^{3/2} }}{N^{1.5}} }\right) \leqslant  C_{\Rsource} \,(M+ M^{3/2} ).
		\end{align*}
		Hence the stability condition permits controlling $\cME3(u,u)$ in order to prove that the effective models are uniformly coercive, just like in the proof of \cref{caractere bien pose des ME augmentes}.
	\end{nvdemo}
\end{itemize}

Moreover, the stability condition of \cref{def: condition de stabilite} is one of the assumptions of \cref{cvg des ME}. For this theorem to still hold with $\CstabME = C_{\Rsource} /N^{3}$, we need two additional small modifications.
\begin{itemize}
	\item One can take the stability constant $C^{\mathrm{dir}} _{\Rsource,N}$ of \cref{ME avec Dirichlet sur Gamma} in the form \red{$C^{\mathrm{dir}} _{\Rsource,N} = C_{\Rsource} / N^{3}$}. This requires modifying the proof of \cref{ME avec Dirichlet sur Gamma} in the same way as \cref{caractere bien pose des ME augmentes}.
	\item \Cref{rq: ME avec Dirichlet sur Gamma} becomes: by reducing the value of \red{the constant $C_{\Rsource} $ in \cref{eq: comportement CstabME en N}}, one can assume that the stability condition of \cref{def: condition de stabilite} implies $\varepsilon \leqslant  C^{\mathrm{dir}} _{\Rsource,N} \,\RME \,\domi{\ln \RME}^{-3}$ (which is the assumption of \cref{ME avec Dirichlet sur Gamma}).
\end{itemize}
\end{textesansboite}

\begin{textesansboite}
Now that we have justified \cref{eq: comportement CstabME en N}, let us prove \cref{norme opP} and \cref{eq: cvg des ME amelioree}. We continue the previous list of modifications.

\begin{itemize}[leftmargin=1.45em,itemsep=4pt,labelsep=0.33em]
	\item \Cref{rq: condition de stabilite et CdtnME} becomes: the stability condition implies $\varepsilon  \leqslant \red{\CdtnMEprecis1} \,\RME \,\domi{\ln \RME}^{-3} \red{N^{-2} }$ and $\red{\CdtnMEprecis2} \.\varepsilon  \.\*\RME^{-1} \domi{\ln \RME}^{3} \red{N^{3} } \leqslant  C_{\Rsource}$.

	\item According to this new version of \cref{rq: condition de stabilite et CdtnME}, the stability condition implies the assumption ``$\varepsilon  \leqslant  \CdtnMEprecis1 \.\RME \.\domi{\ln \RME}^{-3} \.N^{-2}$'' of the above generalization of \cref{opS - opS0}. Thus, we have $\|\opP\|_{\mathcal{L}(L^2 (\SME),\mathbb{C})} \leqslant  \CdtnMEprecis2 \, \varepsilon  \.N^{2.5} \.\RME^{-2.5} \,\domi{\ln \RME}^{3} $, which proves \cref{norme opP}. The rest of this list explains how to prove \cref{eq: cvg des ME amelioree}.
	\phantomsection\label{page: preuve norme opP}

	\item \Cref{norme opS} becomes: under the stability condition of \cref{def: condition de stabilite}, we have for~any $j\in \{0,1,2\}$ and $\red{s\in [0, \frac{1}{2} )}$:
	\[\|\opS\|_{\mathcal{L}(\red{\HME{s}}(\SME),\HME{j}(\SME))} \leqslant  C_{\Rsource\red{,s}} \,N^{j+1\red{-s}} /\RME .\]
	\begin{nvdemo}
		Like in the original proof, $\|\opS-\opSref\|_{\mathcal{L}(\red{\HME{s}}(\SME),\HME{j}(\SME))} \leqslant  \|\opS-\opSref\|_{\mathcal{L}(L^2 (\SME),\HME{j}(\SME))}$ is controlled using \cref{opS - opS0}. Moreover, for any $v\in \red{\HME{s}}(\SME)$ and $\bfx \in \SME$, we have:
		\[\opSref(v)(\bfx) = \sum _{n=1} ^{N} \langle v,\varphi _{n} \rangle _{L^2 (\SME)} \, \partial _{r} \varphi _{n} (\bfx) \qquad  \text{with}  \quad  \varphi _{n|\SME} :\theta  \mapsto  \sqrt{ \frac{2}{\RME\Theta } } \sin (n\txtpst \theta ) \.\beta _{n}\]
		and $|\beta _{n} |=1$ and $\|\partial _{\theta } ^{j} \partial _{r} \varphi _{n} \|_{L^2 (\SME)} \leqslant  C_{\Rsource} N^{j+1} /\RME$. Now, for any $v\in \HME{s}(\SME)$, let us define
		\verydarkred{\[ \tilde{v} :\theta \in \mathbb{R} \mapsto  \left\{\begin{array}{ll}
			v(\RME, \theta  \,\mathrm{mod}\, 2\Theta ) & \text{if } (\theta  \,\mathrm{mod}\, 2\Theta ) \in [0,\Theta ] \\[1pt]
			-v(\RME, (-\theta ) \,\mathrm{mod}\, 2\Theta ) & \text{otherwise}.
		\end{array}\right.\]
		$\tilde{v}$ is $2\Theta $-periodic and, since $v\in \HME{s}(\SME)$ and $s< \frac{1}{2}$, we have $v\in H^{s} (\mathbb{R})$. Moreover, $\tilde{v}$ is odd, so we have for~any $n$: $\langle v,\varphi _{n} \rangle _{L^2 (\SME)} = \frac{1}{2} \RME^{1/2} \.c_{n} ( \tilde{v} )$ where $c_{n}$ denotes the $n$-th Fourier coefficient associated to a sinus (using Fourier series with sin and cos rather than with complex exponentials). Using $\|n^{s} \.c_{n} ( \tilde{v} )\|_{\ell ^2 (\mathbb{N}^*)} \leqslant  C \| \tilde{v} \|_{H^{s} (\mathbb{R})} \leqslant  C_{s} \.\RME^{-1/2} \|v\|_{\HME{s}(\mathbb{R})}$, we deduce}
		\begin{multline*}
			\big\|\partial _{\theta } ^{j} \opSref(v) \big\|_{L^2 (\SME)} ^2  \.=\. \sum _{n=1} ^{N} |\langle v,\varphi _{n} \rangle |^2  \, \|\partial _{\theta } ^{j} \partial _{r} \varphi _{n} \|_{L^2 (\SME)} ^2  \\
			\leqslant \. \left( C_{\Rsource} \frac{N^{j+1} }{\RME} \right)^2  \cdot  \red{N^{-2s} } \sum _{n=1} ^{N} \red{n^{2s} } |\langle v,\varphi _{n} \rangle |^2  
			\.\leqslant \. C_{\Rsource,s} \left( \frac{N^{j+1\red{-s}} }{\RME} \right)^2  \|v\|_{\red{\HME{s}}(\SME)} ^2 .
		\end{multline*}
	\end{nvdemo}

	\item In order to show \cref{eq: cvg des ME amelioree}, the proof of \cref{cvg des ME} page~\pageref{preuve: cvg des ME} is modified as follows.
	\begin{nvdemo}
		Similarly to \cref{eq: norme L2 SME vs H1 OME}, one can show $\|u\|_{\red{\HME{1/2}}(\SME)} \leqslant  C \.\RME^{1/2} \|u\|_{H^{1} (\OME)}$ for any $u\in H^{1} (\OME)$. Hence we have \red{for~any $s\in [0, \frac{1}{2} )$}:
		\begin{align*}
			\big\|\opS(\emeg_{|\SME} )\big\|_{\HME2(\SME)} &\leqslant  \|\opS\|_{\mathcal{L}(\red{\HME{s}}(\SME),\HME2(\SME))} \,\|\emeg\|_{\red{\HME{s}}(\SME)} \\
			&\leqslant  C_{\Rsource} \frac{N^{3\red{-s}} }{\RME} \cdot  \sqrt{\RME} \.\|\emeg\|_{H^{1} (\OME)} \\
			&\leqslant  C_{\fsource} \.\errmdl[3/2][m-1] N^{3\red{-s}} \\
			&= C_{\fsource} \.\errmdl[3/2][m-1] N^{\red{2.5+\eta }} 
		\end{align*}
		taking \red{$\eta = \frac{1}{2} -s$}. The rest of the proof stays the same, simply replacing each $N^{3}$ with \red{$N^{2.5+\eta }$}.
	\end{nvdemo}
\end{itemize}
\end{textesansboite}

\end{subappendices}

\chapter{Conclusion and perspectives}
\label{chap: conclusion}
\begin{textesansboite}
Let us summarise the main contributions of this thesis and list some future research directions.
\end{textesansboite}

\begin{textesansboite}
In \crefrange{chap: SF}{chap: construction du DA}, we provided and rigorously justified an asymptotic expansion at any order for a Helmholtz problem in a Dirichlet corner partially covered by a thin layer. This case is much more complex than previous works on the Poisson equation at any order or on the Helmholtz equation at order 2, because of the intricacy of the corner singularities. We handled this complexity by introducing a new system of algebraic tools to obtain the matching conditions, including:
\begin{itemize}
    \item an explicit basis of the corner singularities,
    \item operators that allow to describe neatly the corner singularities appearing in the fields' behaviors in terms of the fundamental singularities $\phi_{d} =r^{d} \sin (d \theta )$,
    \item and formal series, which permit rigorous complex operations on the fields' behaviors, like applying infinite series of operators, without having to deal with truncations of asymptotic behaviors.
\end{itemize}
This allows to replace cumbersome case-specific calculations by abstract generalizable ones. In addition, we developed a symbolic computation module that allows an exact, very fast and memory-thrifty computation of the matching coefficients and the fields' asymptotic behaviors w.r.t.\ $\bfr$ (the generalized distance to the corner).\\

Moreover, we showed in \cref{chap: complements sur le DA} that the asymptotic expansion simplifies in many cases, namely at any order when $\pst$ is irrational and at the first orders when $\pst$ is rational. The simplifications include simpler matching conditions and the absence of powers of $\ln\varepsilon $. However the general expression of the asymptotic expansion remains necessary in some important cases, like the angle $\Theta =\pi $.\\

In \cref{chap: validation numerique du DA}, we proposed an enriched Galerkin method to compute the corner profiles, based on complements at infinity. We carried out the numerical analysis of the method and provided numerical results. In addition, we illustrated numerically the validity of our asymptotic expansion up to  order~5 for various angles.\\

Finally, in \cref{chap: ME}, we proposed high order effective models of  order~3 and 5 w.r.t.\ $\varepsilon $. Our approach involves a novel DtN condition on an artificial circular boundary near the corner, that takes into account the singular behavior of the far fields at the corner. This relies on a modal decomposition of the far fields, that we proved in \cref{chap: decompo modale CL} using specially designed weighted norms on spaces of formal series. We proved the well-posedness and the uniform stability of the effective models, and established precise error estimates w.r.t.\ $\varepsilon $, the DtN radius $\RME$ and the number $N$ of modes of the DtN, that show the convergence of the models. In addition, we provided numerical results that illustrate the theoretical results. We noticed that the choice of $N$ and $\RME$ can be complicated in practice and needs to be investigated. Moreover, by~taking $\omega =0$, our method also provides high order effective models for the Poisson equation, which did~not exist in the literature either.
\end{textesansboite}

\begin{textesansboite}
There remain a few questions on the studied setting that we did~not have time to address. First, it would be interesting to determine precisely the range of angles $\Theta $ for which the naive effective models -- which replace the layer by the Robin, resp.\ Ventcel, effective boundary condition on $\Gamma $ entirely with no special treatment of the corner -- converge at  order~3, resp.\ 5. In these cases, it~is better to use the naive model than ours because the former is simpler to implement. We expect that this concerns the small angles because these are the angles for which the first far fields are the most regular at the corner. \cite{CalCosDauVia06,Via03} proved an optimal order of convergence of $\varepsilon ^{\min (3,\pst)}$ on the entire domain for the naive effective model of  order~3. We conjecture that the convergence far from the corner is rather in $\varepsilon ^{\min (3,2\pst)}$ because the first far field that has a non-trivial corner condition is $u_{2\pst,0}$ (whereas the first non-zero corner field is $S_{\pst,0}$). This would mean that the naive model converges at  order~3 far from the corner when $\Theta  \leqslant  \frac{2\pi }{3}$. Similarly, we conjecture that the naive effective model with the Ventcel condition converges at order $\varepsilon ^{\min (5,2\pst)}$ far from the corner, which is an $\varepsilon ^{5}$ when $\Theta  \leqslant  \frac{2\pi }{5}$. As~in \cite{Via03} for infinite layers, those rates of convergence could be proven by building an asymptotic expansion of the naive models when $\varepsilon \rightarrow 0$, and noticing that its first terms coincide with the far fields $u_{p,\ell }$.\\

Secondly, we think that a DtN condition similar to our effective models might be derived to compute the corner profiles. This DtN condition would help truncate the domain $\Omega _{1}$ at infinity and would replace the complements at infinity that we introduced in \cref{chap: validation numerique du DA}. This could reduce the computational costs because the matrix to be inverted in the finite element method would be more sparse. Indeed, a DtN condition leads to a dense part of the matrix on a 1D subdomain, whereas complements at infinity induce to a dense part of the matrix on a 2D subdomain. However, like in \cref{chap: ME}, the error analysis would probably be more complex than for the complements at infinity. In addition, we do~not know what error estimate we would obtain on the computation of the coefficients $\sigma _{q} (\ccDA_{d,n} )$.
\end{textesansboite}

\begin{textesansboite}
Next, we aim to extend our work to other settings. We already explained in \cref{eq: generalisation a une double couche} how to generalize the asymptotic expansion of $u_{\varepsilon }$ to a Dirichlet corner \emph{entirely} covered by a layer, and the rest of the thesis can similarly be adapted to this case without difficulty.\\

In a future work, we intend to treat also the case of a Neumann boundary condition. This case seems to work overall similarly to Dirichlet. For instance, in \cref{chap: raccord algebrique} one essentially needs to replace the spaces $\mathcal{A}_{d} (\Omega )$ by
\[\widetilde{\mathcal{A}}_{d} (\Omega ) := \operatorname{Span}\!\big(\big\{ \red\Re [(\alpha z)^{q} \, \overline{\alpha z}^{k} \log (\alpha z)^{j} ] \;\big|\; (q,k,j)\in \mathbb{R}\times \mathbb{N}^2 , q+k=d \big\}\big)\]
and $\ROmegaPi$, $\RGammaPi$ and $\RcouchePi$ by the appropriate operators. The rest of the reasoning to obtain the matching conditions is unchanged thanks to the abstract general method we proposed. A noticeable difference with the Dirichlet case is that the functions $\phi_{d} =r^{d} \sin (d\theta )$ are replaced by $r^{d} \cos (d\theta )$ ($d\in \TZ*$), plus the functions $(r,\theta ) \mapsto  1$ and $(r,\theta ) \mapsto  \ln r$ (which both belong to $\widetilde{\mathcal{A}}_{0} (\Omega )$). In addition, the construction of the corner fields is less obvious because the Poisson problem they satisfy requires a compatibility condition with the source term. This influences the matching conditions. We have~not yet investigated how to adapt the results of the other chapters to the Neumann case, but we suspect that the expansion involves powers of $\ln \varepsilon $ for any angle $\Theta $ because of the presence of the singularity $\ln r$.\\

On the longer term, it would be of interest to study 3D cases, because they are closer to real applications. The possible situations, of increasing difficulty, include an edge, a conical point and a polyhedral point. One can consider the Helmholtz equation or, more complicated, the Maxwell equations. The case of Helmholtz in an edge seems to be simply a tensorization of our 2D work along the direction of the edge, but the other ones are less obvious to treat. In particular they require finding new appropriate expressions of the corner singularity basis functions in order to adapt our method.\\

Finally, we would like to study finite periodic layers. This is much more challenging because an infinite periodic layer alone already requires to use a matched asymptotic expansion or a multiscale expansion, see \cref{sec: etat de l'art}. An approach used in \cite{DelSchSem16,SemDelSch18} combines a multiscale expansion for the layer and a matched asymptotic expansion near the corner. Another possibility that we are considering is to use matched expansions on both the layer and the corner. This would involve three main zones (note that the far and layer zones cannot be easily merged into one like we did), three matching zones where two types of fields must match, and one zone where all three types of fields must match, see \cref{fig: dvlpmt raccorde couche periodique finie}. We suspect that the matching between the far and corner zones would be similar to the work of this thesis and also involve the corner singularities of $\mathcal{A}(\Omega )$. As~for the layer zone and its matching with the far zone, it could be treated like infinite periodic layers. In particular, one would need to compute numerically some layer profiles, in addition to the corner layers. The third matching -- between the corner zone and the layer -- is not clear for now.
\end{textesansboite}

	\tikzset{external/export=true}%
	\tikzsetfigurename{fig-couche_periodique_finie-}%
	\begin{figure}[h]
\def\epsFig{0.6} 
\def\widthFig{4.2} 
\def\heightFig{3.2} 
\def\depthFig{1.8} 
\def\ThetaFig{225} 
\def\TFig{1.05} 
\def\lci{1.35} 
\def\lce{2.1} 
\def\lri{0.6} 
\def\lre{1.3} 
\pgfmathsetmacro{\tri}{asin(\lri/\lce)}
\pgfmathsetmacro{\tre}{asin(\lre/\lce)}
\def\plt#1#2{{#1*cos(#2)},{#1*sin(#2)}}
\newcommand{\cellulePeriodicite}[1]{
	\begin{scope}[shift={({#1*\epsFig*\TFig+0.1},0)}, xscale={\epsFig*\TFig}, yscale={-\epsFig}]
		\fill [color=white!54!gray] plot [smooth cycle] coordinates {(0.1,0.2) (0.3,0.2) (0.5,0.1) (0.7,0.3) (0.8,0.5) (0.7,0.9) (0.4,0.8) (0.3,0.9) (0.2,0.7) (0.1,0.5)} ;
		\fill [color=white!36!gray] plot [smooth cycle] coordinates {(0.3,0.25) (0.5,0.2) (0.6,0.3) (0.7,0.5) (0.7,0.8) (0.4,0.7) (0.3,0.7) (0.2,0.5)} ;
		\fill [color=white!16!gray] plot [smooth cycle] coordinates {(0.3,0.4) (0.5,0.3) (0.65,0.65) (0.4,0.6) (0.3,0.6) (0.3,0.5)} ;
	\end{scope}}

\begin{tikzpicture}[scale=0.9]
	\coordinate (finLambda) at ({-\depthFig*cot(\ThetaFig)}, -\depthFig);
	\coordinate (sommetCouche) at ({-\epsFig*cot(\ThetaFig)}, -\epsFig);
	\fill [color=white!85!gray] (0,0) -- (\widthFig,0) -- (\widthFig,-\epsFig) -- (sommetCouche) -- cycle; 
	\begin{scope}
		\clip (0,0) rectangle (\widthFig,-\epsFig) ; 
		\fill [color=white!75!gray] (-\widthFig,0) rectangle (\widthFig,-\epsFig) ;
		\foreach \i in {0,1,...,6} {\cellulePeriodicite{\i}}
	\end{scope}
	\fill [color=cyan, opacity=0.3] (\plt\lce\tre) arc (\tre:\ThetaFig:\lce) -- (finLambda) -- (-\widthFig,-\depthFig) -- (-\widthFig,\heightFig) -- (\widthFig,\heightFig) -- (\widthFig,\lre)  -- cycle; 
	\fill [color=green, opacity=0.3] (\lce,0) rectangle (\widthFig,-\epsFig); 
	\fill [color=green, opacity=0.3] (\lce,0) arc (0:\tri:\lce) -- (\widthFig,\lri) -- (\widthFig,0) -- cycle; 
	\fill [color=red, opacity=0.3] (\lci,0) arc (0:\ThetaFig:\lci) -- (sommetCouche) -- (\lci,-\epsFig) -- cycle; 
	\foreach \i in {5,15,...,95}{
		\pgfmathsetmacro{\ri}{(\lci*(100-\i+5)+\lce*(\i-5))/100}
		\pgfmathsetmacro{\re}{(\lci*(100-\i-5)+\lce*(\i+5))/100}
		\pgfmathsetmacro{\ti}{(\tri*(100-\i+5)+\tre*(\i-5))/100}
		\pgfmathsetmacro{\te}{(\tri*(100-\i-5)+\tre*(\i+5))/100}
		\fill [color=cyan!\i !red, opacity=0.3] (\plt\ri\tre) -- (\plt\re\tre) arc (\tre:\ThetaFig:\re) -- (\plt\ri\ThetaFig) arc (\ThetaFig:\tre:\ri) -- cycle; 
		\fill [color=green!\i !red, opacity=0.3] (\ri,0) -- (\re,0) arc (0:\tri:\re) -- (\plt\ri\tri) arc (\tri:0:\ri) -- cycle; 
		\fill [color=green!\i !red, opacity=0.3] (\ri,0) rectangle (\re,-\epsFig); 
		\fill [color=cyan, opacity=0.3*\i/100] (\plt\lce\ti) arc (\ti:\te:\lce) -- (\widthFig,{\lce*sin(\te)}) -- (\widthFig,{\lce*sin(\ti)}) -- cycle;
		\fill [color=green, opacity=0.3*(1-\i/100)] (\plt\lce\ti) arc (\ti:\te:\lce) -- (\widthFig,{\lce*sin(\te)}) -- (\widthFig,{\lce*sin(\ti)}) -- cycle;
		\foreach \j in {5,15,...,95}{ 
			\pgfmathsetmacro{\ti}{(\tri*(100-\j+5)+\tre*(\j-5))/100}
			\pgfmathsetmacro{\te}{(\tri*(100-\j-5)+\tre*(\j+5))/100}
			\pgfmathsetmacro{\op}{0.1 * \i/100 * (1-\j/100)}
			\fill [color=cyan!\j !green!\i !red, opacity=0.3*\j/100] (\plt\re\ti) arc (\ti:\te:\re) -- (\plt\ri\te) arc (\te:\ti:\ri) -- cycle;
			\fill [color=green, opacity=0.3 * \i/100 * (1-\j/100)] (\plt\re\ti) arc (\ti:\te:\re) -- (\plt\ri\te) arc (\te:\ti:\ri) -- cycle;
			\fill [color=red, opacity=0.3 * (1-\i/100) * (1-\j/100)] (\plt\re\ti) arc (\ti:\te:\re) -- (\plt\ri\te) arc (\te:\ti:\ri) -- cycle;
		}
	}
	\draw [thick] (finLambda) -- (sommetCouche) -- (\widthFig,-\epsFig);
	\draw (sommetCouche) -- (0,0) -- (\widthFig,0);
	\draw [thick,dashed] (\lci,-\epsFig) -- (\lci,0) arc (0:\ThetaFig:\lci);
	\draw [thick,dashed] (\lce,-\epsFig) -- (\lce,0) arc (0:\tri:\lce);
	\draw [thick,dashed] (\plt\lce\tre) arc (\tre:\ThetaFig:\lce);
	\draw [thick,dashed] (\plt\lce\tre) -- (\widthFig,\lre);
	\draw [thick,dashed] (\plt\lce\tri) -- (\widthFig,\lri);
	\draw [>=stealth,->] (-\widthFig-0.4,0) -- (\widthFig+0.4,0) node [right] {$x$};
	\draw [>=stealth,->] (0,-\depthFig) -- (0,{\heightFig+0.3}) node [right] {$y$};
	\draw [color=black!50!blue] (-\widthFig*0.7,\heightFig*0.63) node {far zone};
	\draw [color=black!60!green,<-,>=latex] (\widthFig*0.77,-\epsFig*0.6) to[bend right=20] ++(-0.05,-0.6) node[below]{layer zone};
	\draw [color=black!50!purple] (0,\lci*0.4) node {\parbox{20mm}{\centering corner\\ zone}};
	\draw [<->,>=latex] (\widthFig*0.9,0) -- node[midway,right]{$\varepsilon$} (\widthFig*0.9,-\epsFig);
	\draw [<-,>=latex] ({(\lci+\lce)/2*cos(67)},{(\lci+\lce)/2*sin(67)}) to[bend left=20] ++(0.5,0.55) node[right]{matching zones};
	\draw [<-,>=latex] ({(\lce+\widthFig)/2},{(\lri+\lre)/2}) to[bend right=20] ++(0,0.95);
\end{tikzpicture}

\caption{Zones of a hypothetical matched asymptotic expansion for a finite periodic layer}
\label{fig: dvlpmt raccorde couche periodique finie}
\end{figure}
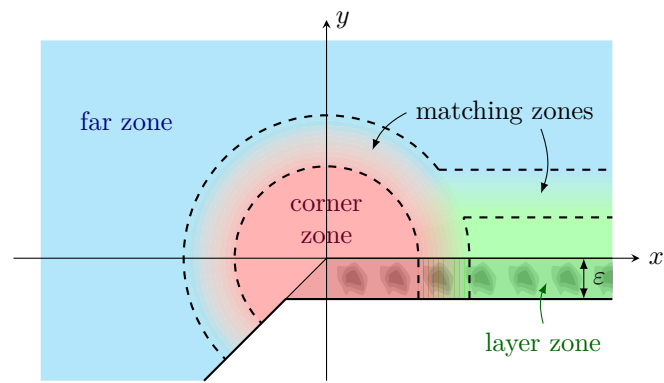%
	\tikzset{external/export=false}%
	\tikzsetfigurename{main-figure}


\printnomenclature 
\label{nomenclature}
\listoffigures
\listoftables
\printbibliography[heading=bibintoc, title=References]

\newpage\pagestyle{empty}\ 
\newpage
\tikzifexternalizing{}{\includepdf{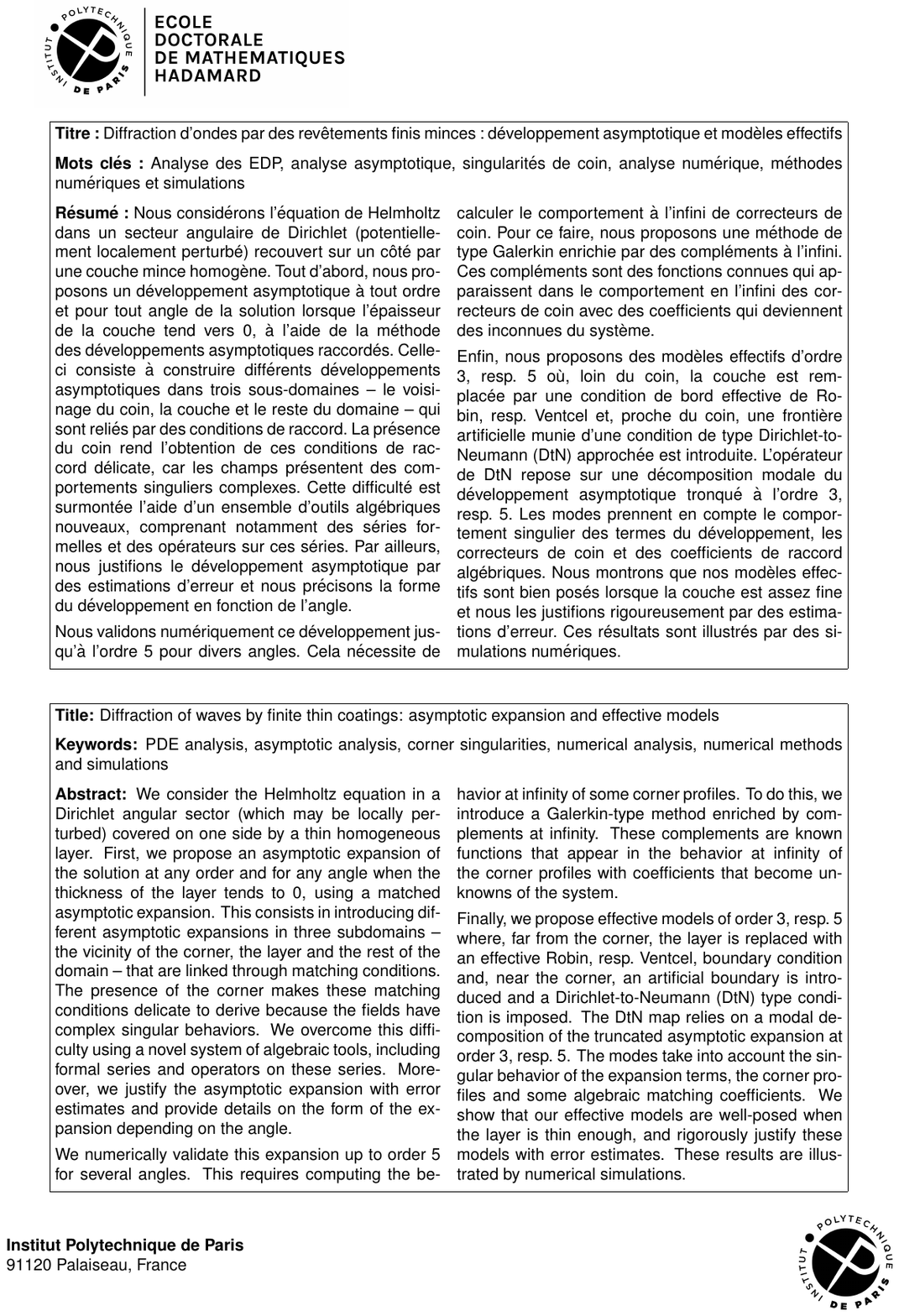}}

\end{document}